\documentclass[11pt,a4paper,twoside,openany]{book}

\usepackage[T1]{fontenc}
\usepackage{lmodern}
\usepackage{iftex}
\ifXeTeX
\usepackage{fontspec}
\newfontfamily\BookJapaneseFont{ipaexm.ttf}
\newcommand{\JapaneseName}[1]{{\BookJapaneseFont #1}}
\newfontfamily\BookChineseFont{FandolSong-Regular.otf}
\newcommand{\ChineseName}[1]{{\BookChineseFont #1}}
\else
\PackageError{book-preface}{XeLaTeX is required for the Japanese names}%
{Select XeLaTeX as the compiler, as specified in the first line.}
\fi

\usepackage{indentfirst}
\usepackage[a4paper,top=27mm,bottom=28mm,left=25mm,right=25mm,headheight=14pt,headsep=8mm]{geometry}
\usepackage{microtype}
\usepackage{tikz-cd}
\usepackage{amsmath,amssymb,amsfonts,amsthm,mathtools,amscd}
\usepackage{graphicx,float}
\usetikzlibrary{arrows.meta}
\usepackage[all]{xy}
\usepackage[table]{xcolor}
\IfFileExists{bbm.sty}{\usepackage{bbm}}{}
\providecommand{\mathbbm}[1]{\mathbb{#1}}
\usepackage{ragged2e}
\definecolor{tablehead}{RGB}{235,228,244}
\definecolor{tablerule}{RGB}{217,209,226}
\definecolor{earlytablehead}{RGB}{224,233,241}
\definecolor{earlytablerule}{RGB}{217,217,217}

\usepackage{booktabs,array,caption,longtable,hhline,refcount}
\usepackage{enumitem,needspace,titlesec,fancyhdr}
\usepackage{xurl}
\usepackage{hyperref}
\definecolor{inkblue}{RGB}{30,57,83}
\hypersetup{colorlinks=true,linkcolor=green,citecolor=blue,urlcolor=inkblue,
	pdftitle={Homotopy groups of special orthogonal groups in degrees 23 to 29}}

\allowdisplaybreaks[4]
\setlist[description]{font=\normalfont\bfseries,leftmargin=0pt,labelsep=.45em,
	topsep=2pt,itemsep=2pt,parsep=1pt,style=nextline}
\renewcommand{\arraystretch}{1.2}
\titleformat{\section}{\Large\bfseries}{\thesection.}{.55em}{}
\titleformat{\subsection}{\large\bfseries}{\thesubsection.}{.55em}{}
\titlespacing{\section}{0pt}{22pt plus 4pt minus 3pt}{10pt}
\titlespacing{\subsection}{0pt}{16pt plus 3pt minus 2pt}{8pt}
\fancypagestyle{plain}{\fancyhf{}\fancyfoot[C]{\small\thepage}}

\numberwithin{equation}{chapter}

\newtheoremstyle{paper}{9pt}{7pt}{\itshape}{}{\bfseries}{.}{.6em}{}
\theoremstyle{paper}
\newtheorem{thm}{Theorem}[chapter]

\newtheorem{prop}[thm]{Proposition}

\newtheorem{lem}[thm]{Lemma}

\theoremstyle{definition}
\newtheorem{defn}[thm]{Definition}
\newtheorem{rem}[thm]{Remark}
\let\originalrem\rem
\let\endoriginalrem\endrem
\renewenvironment{rem}
{\pushQED{\qed}\originalrem}
{\popQED\endoriginalrem}
\newtheorem{exam}[thm]{Example}

\titleformat{\chapter}[display]{\normalfont\bfseries}{\large\chaptertitlename\ \thechapter}{8pt}{\huge}
\titlespacing{\chapter}{0pt}{12pt}{24pt}

\newtheorem{conj0}[thm]{Conjecture}
\setlist[itemize]{leftmargin=1.8em,topsep=4pt,itemsep=4pt,parsep=2pt}
\setlist[enumerate]{label=\textup{(\arabic*)},leftmargin=2em,topsep=4pt,itemsep=4pt}
\usepackage{tocloft}
\makeatletter
\newenvironment{mainleft}{\@fleqntrue\@mathmargin=0pt\relax}{}
\makeatother

\makeatletter
\newcommand{\ResolveBookPhoto}[2]{%
	\def#1{#2}%
	\let\BookPhotoFound\@empty
	\@for\BookPhotoExt:=jpg,JPG,jpeg,JPEG,png,PNG,pdf,PDF\do{%
		\ifx\BookPhotoFound\@empty
		\IfFileExists{#2.\BookPhotoExt}{%
			\edef#1{#2.\BookPhotoExt}%
			\def\BookPhotoFound{yes}%
		}{}%
		\fi}}
\makeatother
\ResolveBookPhoto{\PrefacePhotoFile}{p1karrowbig}
\ResolveBookPhoto{\GroupPhotoFile}{p1k}
\ResolveBookPhoto{\DetailPhotoFile}{jtx}
\ResolveBookPhoto{\RecentPhotoFile}{p13}
\ResolveBookPhoto{\MoriPhotoFile}{mowz}
\ResolveBookPhoto{\MoriAnnotatedPhotoFile}{moyz}

\newcommand{\BookPhoto}[2][]{%
	\IfFileExists{#2}{\includegraphics[#1]{#2}}{%
		\PackageWarning{book-photos}{Photograph #2 not found}%
		\fbox{\parbox[c][18mm][c]{.75\linewidth}{\centering\small
				Photograph file missing: \texttt{\expandafter\detokenize\expandafter{#2}}}}}}

\newcommand{\FitMath}[1]{\begingroup
	\setbox0=\hbox{$\displaystyle #1$}%
	\ifdim\wd0>\linewidth
	\resizebox{\linewidth}{!}{\usebox0}%
	\else\usebox0\fi\endgroup}
\newcommand{\MainWhere}[2]{\begingroup
	\setbox0=\hbox{#1 #2}%
	\ifdim\wd0>\linewidth
	#1\par\nobreak{\raggedleft #2\par}%
	\else #1 #2\par\fi\endgroup}

\definecolor{internalpurple}{RGB}{128,0,160}

\newcommand{\GroupLayout}[2]{\begingroup
	\setbox0=\hbox{$\displaystyle #1$}%
	\ifdim\wd0>\linewidth \FitMath{#2}%
	\else\usebox0\fi\endgroup}

\newcommand{\HighlightMath}[1]{\begingroup
	\setbox0=\hbox{$\displaystyle #1$}%
	\ifdim\wd0>\dimexpr\linewidth-1em\relax
	\resizebox{\dimexpr\linewidth-1em\relax}{!}{\usebox0}%
	\else\usebox0\fi\endgroup}

\newcommand{\hc}{\circ}
\newcommand{\n}{\{ }
\newcommand{\nn}{\} }
\newcommand{\ca}{\eta}
\newcommand{\q}{\Delta}

\newcommand{\cg}{\sigma}
\newcommand{\lt}{\varepsilon}
\newcommand{\ch}{ \bar{\nu}  }

\newcommand{\ck}{\zeta}
\newcommand{\cn}{\kappa}
\newcommand{\co}{\rho}

\newcommand{\cq}{\bar{\mu}}
\newcommand{\cs}{\bar{\zeta}}
\newcommand{\ccs}{\bar{\sigma}}
\newcommand{\ct}{\bar{\kappa}}
\newcommand{\kj}{\phi}

\newcommand{\w}{\omega}

\newcommand{\wq}{\infty}
\newcommand{\af}{\alpha}
\newcommand{\sk}{\mathrm{sk}}

\newcommand{\ind}{\mathrm{Ind}}

\newcommand{\m}{\;\mathrm{mod}\,}
\newcommand{\pa}{\partial}

\newcommand{\x}{\langle}
\newcommand{\xx}{\rangle}

\newcommand{\e}{\iota}

\newcommand{\dyd}{\supseteq}
\newcommand{\xyd}{\subseteq}

\newcommand{\ty}{ \equiv}

\newcommand{\lo}{\Omega}

\newcommand{\cz}{\,\exists\;}

\newcommand{\jia}{\oplus}
\newcommand{\z}{\mathbb{Z}}

\newcommand{\s}{\Sigma}

\newcommand{\StemEntry}[2]{\par\addvspace{5pt}\noindent
	\begingroup\setbox0=\hbox{#1}%
	\ifdim\wd0>\linewidth
	\setbox0=\hbox{#2}%
	\ifdim\wd0>\linewidth\resizebox{\linewidth}{!}{\usebox0}%
	\else\usebox0\fi
	\else\usebox0\fi\endgroup\par}
\theoremstyle{definition}
\newtheorem{exercise}[thm]{Exercise}

\DeclareRobustCommand{\uprightsemicolon}{%
	\ifmmode\mathpunct{\text{\normalfont;}}\else\textnormal{;}\fi}

\begin{document}
	\frontmatter
	\begin{titlepage}
		\centering
		
		\vspace*{24mm}
		{\large\scshape Unstable Homotopy Theory\par}
		
		\vspace{18mm}
		{\fontsize{28}{34}\selectfont\bfseries
			The Art of Toda Brackets\par}
		
		\vspace{18pt}
		
		{\fontsize{20}{27}\selectfont\itshape
			A Journey through the Homotopy Groups\\[4pt]
			of $\mathit{SO}(n)$ in Degrees 23-29\par}
		
		\vfill
		{\large
			\begin{tabular}{@{}l@{\hspace{12mm}}l@{}}
				Juxin Yang   & Yoshihiro Hirato \\[6pt]
				Jinho Lee    & Toshiyuki Miyauchi \\[6pt]
				Juno Mukai   & Mariko Ohara
			\end{tabular}\par}
		
		\vspace{14mm}
		{\normalsize October 2026\par}
	\end{titlepage}
	
	\clearpage
	\thispagestyle{empty}
	\null
	\newpage
	
	\pagenumbering{roman}
	\chapter*{Preface}
	\addcontentsline{toc}{chapter}{Preface}
	\markboth{Preface}{Preface}
	\section*{Part I\quad Recollections of Hirosi Toda}
	\addcontentsline{toc}{section}{Part I\quad Recollections of Hirosi Toda}
	
	Professor Hirosi Toda (\JapaneseName{戸田宏}; Jan.\,1928~--~Feb.\,2026) played a formative role in
	my mathematical life. I would like to begin this book with a few
	recollections of his ideas, his teaching, and the conversations through
	which I came to know his mathematics.
	
	% Two portraits accompanying the recollections in Part I.
	% Keep to1.png and to2.png in the same folder as the LaTeX source.
	\begin{figure}[H]
		\centering
		\captionsetup{font=small,justification=centering,
			singlelinecheck=false,skip=5pt}
		\begin{minipage}[t]{0.44\textwidth}
			\centering
			\includegraphics[width=38mm,height=50mm]{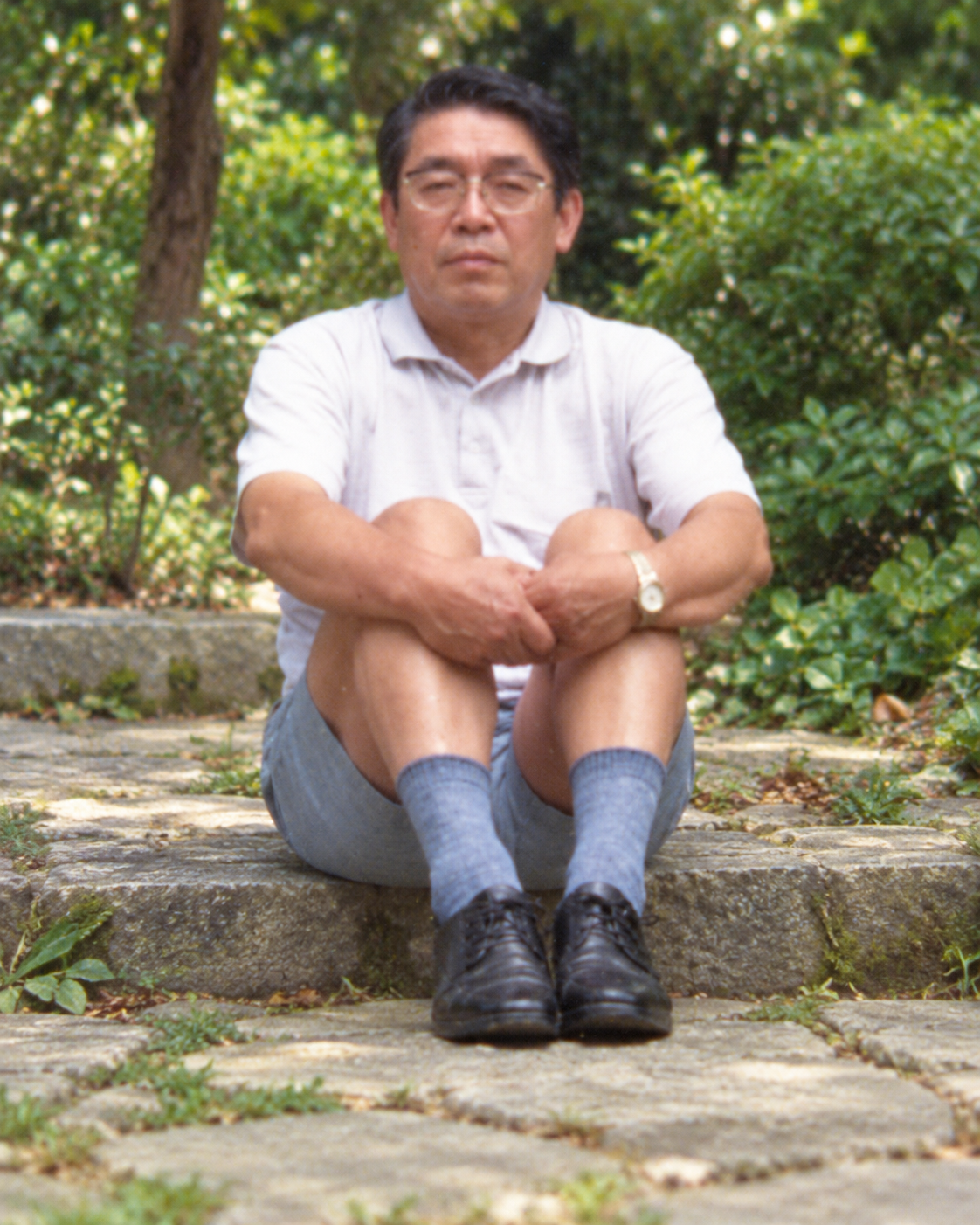}
			\caption*{H.~Toda, photographed in 1988.}
		\end{minipage}\hspace{0.04\textwidth}%
		\begin{minipage}[t]{0.44\textwidth}
			\centering
			\includegraphics[height=50mm]{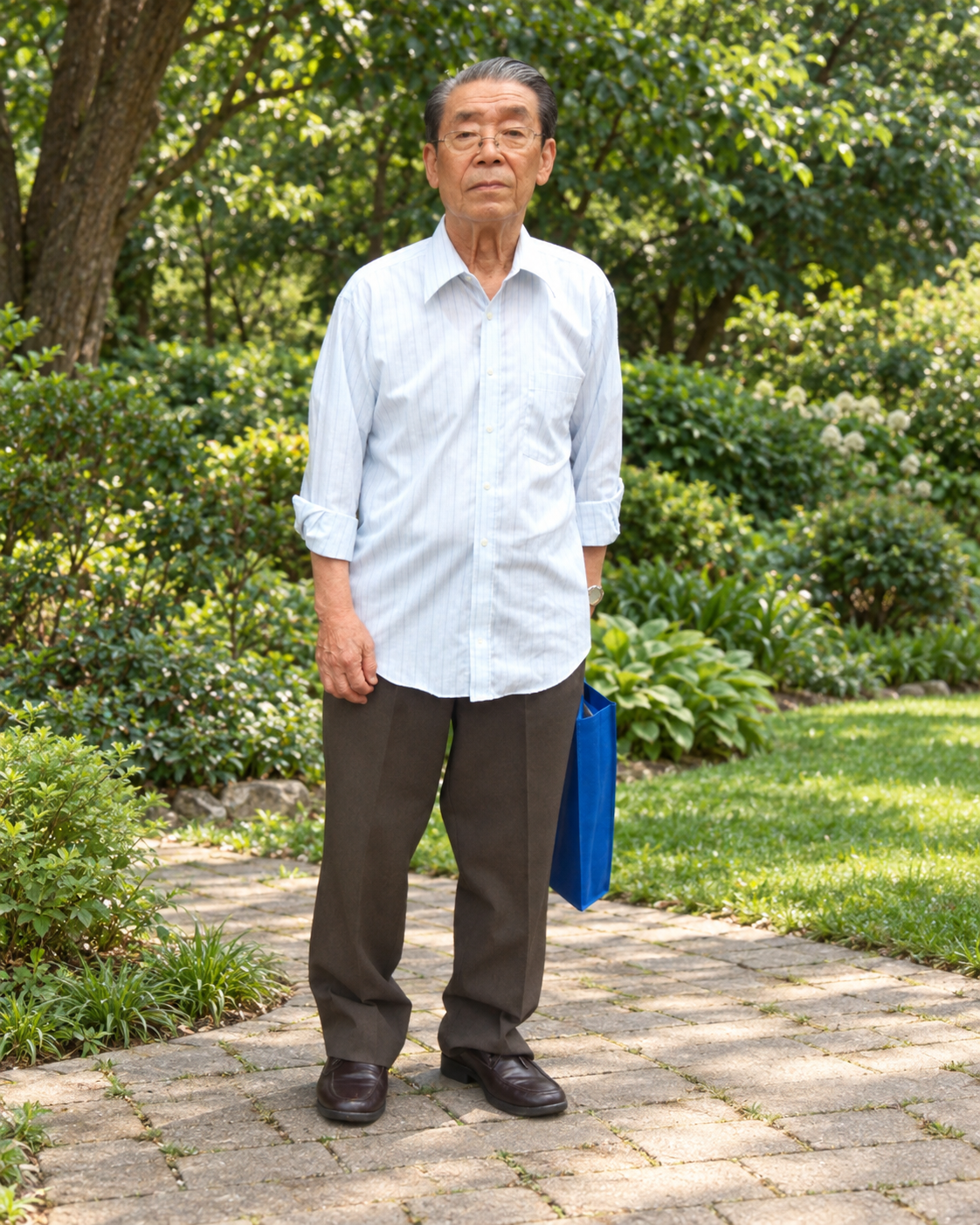}
			\caption*{H.~Toda, photographed in the 2010s.}
		\end{minipage}
	\end{figure}

	In the second year of my master's studies at Kyushu University, I read
	his ``Red Book'', \emph{Composition Methods in Homotopy Groups of Spheres},
	under the guidance of Professor Tatsuji Kudo (\JapaneseName{工藤達二}).
	In 1964, I became a research assistant in the Faculty of Engineering at
	Kyushu University. I occupied one half of a large office divided down
	the middle by a partition, and covered the partition wall with
	handwritten excerpts from the Red Book. That was the year of the Tokyo
	Olympics, and the autumn meeting of the Mathematical Society of Japan
	was held at Kyushu University. The topology programme included a
	lecture by Mamoru Mimura (\JapaneseName{三村護}), entitled
	``The $(n+i)$-th Homotopy Groups of $n$-Spheres for $i=21$ and $22.$ ''
	Professor Toda and Professor Masahiro Sugawara (\JapaneseName{菅原正博})
	visited my office on that occasion. As I recall, Mimura was then in
	France, and Professor Toda delivered the lecture on his behalf,
	discussing the homotopy groups of spheres in the $21$- and $22$-stems.
	
	A year or two later, Professor Toda returned to Kyushu University to
	give an intensive lecture course, based principally on his article
	``An Introduction to Homotopy Theory''
	(\emph{S\=ugaku} \textbf{15}, no.~3 (1964), pp.~141--155).
	
	After moving to Osaka University, I had the privilege of attending a
	weekly seminar at Kyoto University. During the period when I was
	studying the homotopy groups of mod~$2$ Moore spaces, I once
	demonstrated to Franklin Paul Peterson of MIT, who was visiting at the
	time, how a Toda bracket could be used to determine an extension of
	$\mathbb Z/2$ by $\mathbb Z/2.$ 
	
	It was also during this period that Goro Nishida
	(\JapaneseName{西田吾郎}) introduced me to the Kahn--Priddy theorem.
	
	The Kansai Topology Seminar met at various universities, and we would
	play Go before the sessions began. I played Black, with a handicap of
	five stones or more.
	
	After moving to Shinshu University, I worked with Hideyuki Kachi
	(\JapaneseName{可知偉行}) to determine the homotopy groups
	$\pi_{17}(\mathit{SO}(n))$ and $\pi_{18}(\mathit{SO}(n))$ of the
	rotation groups $\mathit{SO}(n).$  We presented our results at a meeting
	of the Mathematical Society of Japan and received comments from
	Professor Toda. He drew up tables of the calculations and seemed to
	have mastered the results more thoroughly than we, the authors, had
	ourselves.
	
	Professor Toda was a seemingly inexhaustible source of ideas, and at
	times I struggled to keep pace with his thinking. At homotopy theory
	symposia at Nagasaki University and in Kinosaki, he explained to me a
	method for determining the homotopy groups of
	$M^3=\Sigma\mathbb RP^2,$   the suspension of the real projective plane.
	The method struck me as groundbreaking. Through Dai Tamaki
	(\JapaneseName{玉木大}) of the Faculty of Science at Shinshu
	University, however, I learned that Frederick Ronald Cohen had also
	recognized the idea, and that his student at the time, Jie Wu (\ChineseName{吴杰}), had
	already incorporated it into his doctoral thesis. Then I read Professor
	Toda's paper ``Order of the Identity Class of a Suspension Space''
	(\emph{Annals of Mathematics} \textbf{78}, no.~2 (1963)) over and
	over again.
	
	His article ``The Stable Homotopy Groups of Spheres''
	(\emph{S\=ugaku} \textbf{10}, no.~2 (1959), pp.~87--97) is somewhat
	demanding, but I regard it as a particularly characteristic
	expression of his mathematical insight, revealing both the depth of
	his expertise and the direction of his interests. 
	It discusses the \(2\)-primary components of the stable homotopy groups of spheres through the \(14\)th stem. About twelve years
	ago, an algebraic topology conference was held at Nankai University
	in Tianjin, China. With Professor Toda in the audience, I gave a
	lecture on the determination of the $2$-primary components of
	$\pi_{n+32}(S^n),$   the homotopy groups of spheres in the $32$-stem.
	
	Shichiro Oka (\JapaneseName{岡七郎}) passed away around 1984. On
	Professor Toda's recommendation, I was asked at short notice to serve
	the remainder of Oka's term as a councillor. Oka had written many
	papers with Professor Toda on the $p$-primary components of homotopy
	groups, and had generously guided me in the study of unstable Adams
	maps and related topics.
	
	In calculations with Toda brackets, it is the determination of the
	indeterminacy that requires particular patience. Consider the element
	$\varepsilon$ in the $2$-primary part of the eighth  stem homotopy group,
	$\pi_8^S(S^{0}:2)=\mathrm{span}\{\bar\nu,\varepsilon\}\cong(\mathbb Z/2)^2,$  
	where $\bar\nu+\varepsilon=\eta\sigma.$ 
	One may express $\varepsilon$ by a Toda bracket, for example as
	$\varepsilon\in\x\eta,2\nu,\nu\xx.$ 
	Its indeterminacy is
	$\eta\circ\pi_7^S+\pi_5^S\circ\nu=\mathrm{span}\{\eta\sigma\},$  
	where $\mathrm{span}\{\eta\sigma\}$ denotes the subgroup generated by
	$\eta\sigma.$  Thus
	$\x\eta,2\nu,\nu\xx=\varepsilon+\mathrm{span}\{\eta\sigma\}.$ 
	Using a mod~$2$ Moore space, one has
	$\bar\eta\widetilde\eta=\pm2\nu$ and can write
	$\varepsilon=\x\eta\bar\eta,\widetilde\eta,\nu\xx.$ 
	
	There is no general method for keeping such indeterminacy calculations
	to a minimum. Nevertheless, I believe it is important to determine
	whether the elements involved lie in the image of the
	$J$-homomorphism.
	
	Professor Toda's ``Unstable $3$-Primary Homotopy Groups of Spheres'',
	published at Himeji Dokkyo University, is a paper of nearly eighty
	pages. I have understood very little of it so far, and would like to
	study it further.
	
	Professor Toda was one of the world's most influential homotopy
	theorists and a central figure in the development of the subject.
	Through my association with him, I also became acquainted, in some
	measure, with the work of John Frank Adams, Mark Edward Mahowald,
	Ioan Mackenzie James, and others, and had opportunities for brief
	conversations with Mark Edward Mahowald, Ioan Mackenzie James,
	Frederick Ronald Cohen, and Brayton Gray.
	
	I remain deeply grateful for the opportunity to have known Professor
	Toda and for his guidance and kindness over so many years.
	
	\par\addvspace{1.25\baselineskip}
	\Needspace{4\baselineskip}
	\begin{flushright}
		\textsc{Juno Mukai} (\JapaneseName{向井純夫})\\[3pt]
		Professor Emeritus, Shinshu University, Japan\\[3pt]
		\textit{Written in October 2026}
	\end{flushright}
	
	\clearpage
	\section*{Part II\quad About This Book}
	\addcontentsline{toc}{section}{Part II\quad About This Book}

	This book has several fundamental aims, one of which is to serve as a reference book for computing unstable homotopy groups using  methods of the Toda school.
	The Toda bracket serves the art of constructing homotopy liftings and extensions of maps.  It plays a
	fundamental role in dealing with composition relations of homotopy classes.  Appendix~\ref{tlqdzj} collects results
	on the $2$-primary components of unstable homotopy groups\footnote{Here, by ``unstable'' we  mean ``not necessarily stable''.} of spheres
	in stems $0$--$32$ and the known strictly unstable groups in stem $33$, with
	explicit generators. Section~\ref{sec:sphere-history} describes the sources and the scope of these tables.
	Appendix~\ref{app:oguchi} explains the naming conventions and records composition relations.

	The book mainly presents calculations of the $2$-primary components of the homotopy groups
	$\pi_i(\mathit{SO}(n))$ for all $23\leq i\leq 28$ and all $n,$  
	as well as     some $\pi_{29}(\mathit{SO}(n)),$  
	with explicit generators, while retaining two possible group structures for $\pi_{28}(\mathit{SO}(11)).$  In degree $29,$   the cases $7\leq n\leq14$ remain under investigation. The odd-primary components of the homotopy groups
	$\pi_i(\mathit{SO}(n))$ in the above range can be computed with little effort; these computations essentially amount to looking up the relevant entries in existing tables. See  Section \ref{sec:odd-method} and Section \ref{sec:odd23}. 
	
	As a prelude to the computation of $\pi_i(\mathit{SO}(n)),$   we present an introduction to the theory of Toda brackets, supplemented by carefully chosen exercises. Through this exposition, we hope to familiarize readers with their fundamental properties and the essential techniques involved in working with them.

	The study of groups $\pi_i(\mathit{SO}(n))$ grew out of geometric questions about vector
	bundles and spheres. An oriented real vector bundle over a sphere is
	described by a clutching map into a special orthogonal group. Adding a
	trivial line to the bundle gives the inclusion into the next group.
	Bott's 1959 calculation revealed a remarkable simplification after
	stabilization: the homotopy groups of the orthogonal group repeat with
	period eight. His proof used the geometry of geodesics on symmetric
	spaces \cite{Bo}.
	
	For a fixed value of $n,$   however, the answer need not be stable. Kervaire's
	1960 paper illustrates the geometric significance of this distinction:
	it relates nonstable homotopy groups of special orthogonal groups to
	the question of which spheres are parallelizable \cite{Ke}. To calculate
	these groups, one must also understand the homotopy groups of spheres
	and the homomorphisms in the exact sequence of the fibration
	$\mathit{SO}(n)\to\mathit{SO}(n+1)\to S^n.$ 
	Thus progress on spheres and progress on classical groups have remained
	closely connected.
	
	For reference, write
	\begin{equation}\label{wlsh1023}
		\mathbf{Stem}_{k}=\{\pi_{n+k}(S^{n})\otimes\mathbb Z_{(2)}\}_{2\leq n\leq k+2}
	\end{equation}
	
	\noindent   for the indicated family of $2$-localized homotopy groups of spheres.
	Toda's 1962 monograph established the classical framework
	($\mathbf{Stem}_{k}$ with $0\leq k\leq 19$).
	His joint work with Mimura in 1963 ($\mathbf{Stem}_{20}$)
	and Mimura's computations in 1965
	($\mathbf{Stem}_{21}$ and $\mathbf{Stem}_{22}$)
	extended these results.
	Mimura, Mori and Oda advanced the calculations in 1975
	($\mathbf{Stem}_{23}$ and $\mathbf{Stem}_{24}$),
	followed by Oda in 1979
	($\mathbf{Stem}_{25}$ to $\mathbf{Stem}_{31}$).
	Miyauchi and Mukai computed the thirty-second stem in 2017
	($\mathbf{Stem}_{32}$).  Section~\ref{sec:sphere-history} gives the precise ranges and sources.
	
	This sequence of work supplied more than tables of abstract groups.
	Toda's composition methods gave explicit generators and relations,
	including relations expressed by higher compositions. These make it
	possible to use the  exact sequence in concrete calculations.
	Mimura's work on classical and exceptional Lie groups brought the same methods to
	fibrations involving classical and exceptional Lie groups \cite{Mi2}.
	Kachi then calculated the $2$-primary components in degrees $15$ and $16,$   with generators \cite{Ka}. Kachi and Mukai treated
	degrees $17$--$18$ and $19$--$20$ \cite{KM1,KM2}, and Hirato, Kachi
	and Mukai continued with degrees $21$--$22$ \cite{HKM}.
	The present book continues this line of calculation.
	
	The central difficulty is often to determine a relation, rather than
	to find the order of a group.  A short
	exact sequence with kernel and quotient both isomorphic to
	$\mathbb Z/2,$   for example, can have middle group $\mathbb Z/4$ or
	$\mathbb Z/2\oplus\mathbb Z/2.$  Distinguishing these possibilities
	requires additional information. Toda brackets, the unstable
	$J$-homomorphism and composition relations provide that information
	in the calculations below. We record the generators and their images
	under inclusion so that the results can also be used in further work.
	
	Chapter~\ref{ch:toda} provides the theoretical basis for this approach,
	beginning with extensions and coextensions. Its treatment of
	composition, suspension, and additivity makes explicit the hypotheses
	under which the relevant bracket formulas hold. The discussion is
	integral, so that the general constructions are presented before their
	use in the localized calculations. For applications, a central task
	is to follow the indeterminacy through the homomorphisms involved:
	an ambiguity in the original bracket may disappear after applying a
	homomorphism, yielding a definite relation in the target group.
	
	Matrix Toda brackets are used extensively in this book. As explained
	in Section~\ref{sec:mtb}, they are three-fold Toda brackets formed
	from maps involving wedges, with the component maps arranged in
	matrices. This formulation allows us to use relations in which a sum
	of composites vanishes, even when the individual composites need not
	vanish. In higher-dimensional homotopy groups, such relations are
	particularly important for constructing generators, computing their
	images under connecting homomorphisms, and determining extensions
	in  exact sequences. Three examples illustrate how the
	method enters the calculations.
	
	\noindent$\bullet$ First, Equation~(\ref{m1:rho13}) expresses $\rho_{13}$ by
	$$  \rho_{13}\in\left\{\bigl(\sigma_{13}, \Sigma\theta\bigr),\,\begin{pmatrix}8\sigma_{20}\\\eta^2_{25}\end{pmatrix},2\iota_{27}\right\}_4\xyd\pi_{28}(S^{13}).$$\noindent 
	To check that a matrix Toda bracket is well-defined, we proceed
	just as for an ordinary Toda bracket, using the usual rules of
	matrix multiplication to compose its entries. This description of $\rho_{13}$ makes it possible to apply
	$P$ and then a connecting homomorphism; it is used in
	Lemma~\ref{m1:mtb2}(2) to determine the connecting images of
	$P(\rho_{13}).$ 
	
	\noindent$\bullet$  Second, the proof of Lemma~\ref{m1:mtb2}(1) uses
	\[\FitMath{
		\bar{\kappa}'\in \left\{\nu_6,\bigl(\eta_9, 2(\sigma_9\eta_{16})\bigr),\,\begin{pmatrix}\bar{\varepsilon}_{10}\\\bar{\nu}_{17}\end{pmatrix}\right\}_2\xyd\pi_{26}(S^{6}).
	}\] \noindent 
	Here the two components in the middle entry allow the relations
	involving $\bar\varepsilon$ and $\bar\nu$ to be used together.
	The resulting bracket provides the connecting images of
	$\bar\kappa'$ in $\pi_{25}(\mathit{SU}(3))$ and
	$\pi_{25}(\mathit{SO}(6)).$ 
	
	\noindent$\bullet$  Third, in the proof of Proposition~\ref{m1:10}, the lift
	$[\bar\varepsilon_9]\in\pi_{24}(\mathit{SO}(10))$ is chosen as
	\[\FitMath{
		[\bar{\varepsilon}_9]\in\left\{\bigl(\Delta(\nu_{10}), {[\iota_7]_{10}}\bigr),\,\begin{pmatrix}\nu_{12}\\\varepsilon_7 + \sigma'\eta_{14}\end{pmatrix},\bar\nu_{15}\right\}.
	}\]\noindent
	Projection identifies its image with $\bar\varepsilon_9,$   while
	composition properties of the matrix Toda bracket give the relation
	$2[\bar\varepsilon_9]=[\iota_7]_{10}\nu_7\kappa_{10}.$ 
	Thus the bracket supplies both a representative of a generator
	and information about its order. These examples explain why matrix
	Toda brackets are an important part of the method, rather than
	merely a convenient way of writing formulas.
	
	Four-fold Toda brackets also play an essential role in the
	calculation of $\pi_{28}(G_2)$ in Proposition~\ref{m3:p28g2}.
	The relevant expressions are
	$$\{\nu_5^{2},2\iota_{11},\eta_{11}\sigma_{12},\sigma_{19}\}_{1}\;\;\text{and}\;\;\{\varepsilon_3,2\iota_{11},\eta_{11}\sigma_{12},\sigma_{19}\}_{2},$$\noindent
	which describe $\phi_5$ and a choice of $\phi',$   with the
	indeterminacies specified in Lemma~\ref{m3:phi5}.
	We also use $\{j\varepsilon_3,2\iota_{11},\eta_{11}\sigma_{12},
	\sigma_{19}\}_{2},$   where $j:S^3\hookrightarrow\mathit{SU}(3)$
	is the inclusion; its construction is given in
	Lemma~\ref{m3:tsd}(2). Together with the indeterminacy calculation
	in Lemma~\ref{m3:ghc}, these constructions yield
	$\Delta(\phi_6)\equiv j_*\phi'\m2[\bar\rho''']$
	in Proposition~\ref{m3:drgpp}(2). This connecting image is needed
	to determine the cokernel in the  exact sequence of homotopy groups induced by the fibration
	$$\mathit{SU}(3)\longrightarrow G_2\longrightarrow S^6$$\noindent
	and hence to obtain $\pi_{28}(G_2).$    The meaning of the four-fold notation in these calculations
	deserves particular attention. We use it as shorthand for a
	three-fold Toda bracket with suitable extensions and coextensions
	chosen and fixed. For example, Lemma~\ref{m3:tsd}(1) specifies
	$$\{\varepsilon_3,2\iota_{11},\eta_{11}\sigma_{12},\sigma_{19}\}_{2}
	:=\allowbreak
	\{\varepsilon_3,\Sigma^2\overline{2\iota_9},
	\Sigma^2\widetilde{\sigma_{17}}\}_{2},$$ \noindent
	where $\overline{2\iota_9}$ is an extension of $2\iota_9$ over
	$C_{\eta_9\sigma_{10}}$ and $\widetilde{\sigma_{17}}$ is a
	coextension of $\sigma_{17}$ into that cofiber. The expressions
	with $\nu_5^2$ and $j\varepsilon_3$ are interpreted in the same
	way, with the suspension indices stated above.
	The rigorous definition of an unstable four-fold Toda bracket is
	considerably more involved: it requires compatible choices of
	null-homotopies and a careful treatment of the resulting
	extensions and coextensions. Section~\ref{sec:higher-toda}
	explains the distinction and points to the formal definition in
	\cite[Section~6]{O2015};
	Proposition~\ref{prop:admissible-null-quadruple} gives the
	admissibility criterion used in Lemma~\ref{m3:tsd}.
	In particular, our shorthand does not assert independence of
	these choices. The proofs use the specified three-fold brackets
	and explicitly control their indeterminacies.
	
	The introductory chapters explain the notation and recall the
	constructions used in the proofs. The chapters organized by degree
	present the calculations, while the appendices collect results on
	homotopy groups of spheres. Chapter~\ref{ch:geometry} discusses the
	interpretation in terms of vector bundles. Cases still under
	investigation are explicitly indicated in the text.
	
	Readers interested primarily in Toda brackets may read
	Chapter~\ref{ch:toda} together with selected proofs in
	Chapter~\ref{ch:fundamental} and the chapters devoted to individual
	degrees. This offers a way to study the methods through their
	applications: the general properties explain the steps of a proof,
	while the calculation shows why a particular bracket was chosen and
	how its indeterminacy was handled.
	
	We assume familiarity with based homotopy groups and the homotopy
	exact sequence of a fibration. A reader may begin with a single degree
	and consult the introductory chapters and appendices as needed.
	Unless explicitly stated otherwise, \textbf{spaces and maps are
		localized at $2.$ }
	
	The photographs below recall   mathematicians of the Toda school whose work shaped the history of calculations of the unstable homotopy groups of spheres. Recall that when introducing $\mathbf{Stem}_{k}$ in Equation \eqref{wlsh1023}, we briefly reviewed the groups calculated by these mathematicians. Section~\ref{sec:sphere-history} will give the precise ranges and sources.

	\medskip
	\noindent\begin{minipage}{\linewidth}
		{\centering
			\BookPhoto[width=1.2\linewidth,height=.36\textheight,keepaspectratio]{\PrefacePhotoFile}
			\par}
		\smallskip
		\noindent\hspace*{.05\linewidth}\begin{minipage}{.90\linewidth}
			% Normal LaTeX paragraph justification; final line stays left-aligned.
			\footnotesize\normalfont
			\setlength{\leftskip}{0pt}%
			\setlength{\rightskip}{0pt}%
			\setlength{\parfillskip}{0pt plus 1fil}%
			\setlength{\parindent}{0pt}%
			\spaceskip=0pt\relax
			\xspaceskip=0pt\relax
			\pretolerance=100\relax
			\tolerance=1000\relax
			\emergencystretch=1em\relax
			\noindent
			Taken in August 1988, this group photograph records the International
			Conference on Homotopy Theory, held at the Kinosaki Conference Hall
			in Kinosaki, Hyogo, Japan, on 19--24 August 1988 to celebrate
			H.~Toda's sixtieth birthday. H.~Toda, M.~Mimura and N.~Oda are
			identified in the photograph. Their work on compositions and the
			homotopy groups of spheres \cite{T,MT1, Mi1,MMO,Oda1979} underlies the calculations in this book.\par
		\end{minipage}
	\end{minipage}
	\par\medskip
	\par\medskip
	\noindent The photograph below introduces M.~Mori, whose joint work with
	M.~Mimura and N.~Oda \cite{MMO} determined the $23$- and $24$-stem homotopy
	groups of spheres. Their results form part of the foundation for the
	calculations in this book. The photograph was taken on Shikanoshima,
	Fukuoka, Japan, in 1973, a few days after the close of the International
	Conference on Manifolds and Related Topics in Topology. Organized by
	M.~F.~Atiyah and Kunihiko Kodaira, the conference was held from
	10 to 17 April 1973.
	\par\nopagebreak\smallskip
	\noindent\begin{minipage}{\linewidth}
		\centering
		\BookPhoto[width=.554\linewidth,height=.264\textheight,keepaspectratio]{\MoriPhotoFile}
		\par\smallskip
		\begin{minipage}{.84\linewidth}
			\normalsize\justifying\noindent
			The dashed arrow identifies M.~Mori, seated at the lower right.
			The sole figure of European or American appearance is M.~F.~Atiyah.\par
		\end{minipage}
	\end{minipage}
	\par\medskip
	
	\noindent The following photograph, taken in 2025, shows Juno Mukai and Toshiyuki Miyauchi.
	Their joint work determines  the $32$-stem \cite{MiM1}.

	\par\smallskip
	\noindent\makebox[\linewidth][c]{\BookPhoto[width=.468\linewidth,height=.221\textheight,keepaspectratio]{\RecentPhotoFile}}%
	\par\nopagebreak\smallskip
	{\footnotesize\justifying\noindent
		\centerline{ Juno Mukai and his former student Toshiyuki Miyauchi (2025).}\par}

	\subsection*{Authorship and Contributions}
	
	This book brings together calculations by Yoshihiro Hirato, Jinho Lee, Toshiyuki Miyauchi, Juno Mukai and Mariko Ohara in degrees $23$--$27,$   presented in \cite{r2324} (2021) and \cite{r2567} (2022), and by Juxin Yang, Juno Mukai and Toshiyuki Miyauchi in degrees $28$--$29,$   presented in \cite{r289} (2023). In the latter work, Yang carried out the calculations under Mukai's guidance, and Miyauchi verified them in full.
	In 2023, under the guidance of Juno Mukai, Juxin Yang independently
	recomputed the groups $\pi_i(\mathit{SO}(n))$ for $i=25, 26,27$
	and $n\geq2,$   verifying the corresponding results of
	\cite{r2567} without relying on its calculations.

	The present volume brings together the material from these preprints, compiled and organized by Juxin Yang in 2026. It covers the calculations of $\pi_i(\mathit{SO}(n))$ for $23\leq i\leq28,$   while retaining two possible group structures for $\pi_{28}(\mathit{SO}(11)),$   and selected results in degree $29,$   where the cases $7\leq n\leq14$ remain under investigation.  Chapters~\ref{ch:history}--\ref{xzzj} and~\ref{ch:geometry}, together with Appendix~\ref{app:oguchi}, were also written by Yang, while Appendix~\ref{tlqdzj} was written jointly by Lee and Yang. These chapters and appendices provide the background and theoretical framework for the calculations. The tables below give a more detailed account.
	
% Division of responsibilities: local table formatting only.
\begingroup
\small
\setlength{\tabcolsep}{4pt}
\renewcommand{\arraystretch}{1.2}
\setlength{\LTleft}{0pt}
\setlength{\LTright}{\fill}
\setlength{\LTpre}{8pt}
\setlength{\LTpost}{8pt}
\arrayrulecolor{earlytablerule}

% The three tasks apply only to Chapters 6--12.
\begin{longtable}{@{}
		>{\raggedright\arraybackslash}p{35mm}
		>{\raggedright\arraybackslash}p{49mm}
		>{\raggedright\arraybackslash}p{33mm}
		>{\raggedright\arraybackslash}
		p{\dimexpr\linewidth-117mm-6\tabcolsep\relax}@{}}
	
	\caption*{\textbf{Responsibilities for
			Chapters~\ref{ch:degree23}--\ref{ch:degree29}}}\\
	\toprule
	\rowcolor{earlytablehead}
	\textbf{Chapters\newline and source} &
	\textbf{Calculations\newline and proofs} &
	\textbf{\LaTeX{} preparation \newline  and preliminary verification} &
	\textbf{Comprehensive verification\newline of results}\\
	\midrule
	\endfirsthead
	
	\multicolumn{4}{@{}l@{}}{\small\itshape
		Responsibilities for
		Chapters~\ref{ch:degree23}--\ref{ch:degree29}
		(continued)
	}\\[4pt]
	\toprule
	\rowcolor{earlytablehead}
	\textbf{Chapters\newline and source} &
	\textbf{Calculations\newline and proofs} &
	\textbf{\LaTeX{} preparation \newline  and preliminary verification} &
	\textbf{Comprehensive verification\newline of results}\\
	\midrule
	\endhead
	
	\midrule
	\multicolumn{4}{r@{}}{\footnotesize\itshape
		Continued on the next page
	}\\
	\endfoot
	
	\bottomrule
	\multicolumn{4}{@{}p{\linewidth}@{}}{\footnotesize
		\LaTeX{} preparation refers to the cited preprints.
		\textsuperscript{*}Ohara assisted with the group extension in
		$\pi_{23}(\mathit{G}_2).$ 
		Yang's verification in degrees $25$-$27$ was carried out
		by independent calculation in 2023 under Mukai's guidance.
	}\\
	\endlastfoot
	
	Chapters~\ref{ch:degree23}--\ref{ch:degree24}\par
	\cite{r2324} &
	Miyauchi and Mukai\textsuperscript{*} &
	Miyauchi &
	Hirato and Lee\\[6pt]
	
	\rowcolor{blue!3}
	Chapters~\ref{ch:degree25}--\ref{ch:degree26}\par
	\cite{r2567} &
	Hirato, Miyauchi,\newline Mukai and Ohara &
	Lee &
	Miyauchi and Yang\\[6pt]
	
	Chapter~\ref{ch:degree27}\par
	\cite{r2567} &
	Hirato &
	Lee &
	Yang\\[6pt]
	
	\rowcolor{blue!3}
	Chapters~\ref{ch:degree28}--\ref{ch:degree29}\par
	\cite{r289} &
	Yang   &
	Yang &
	Mukai and   Miyauchi\\[6pt]
	
\end{longtable}

% Overall responsibility: no subdivision into three tasks.
\begin{longtable}{@{}
		>{\raggedright\arraybackslash}p{50mm}
		>{\raggedright\arraybackslash}
		p{\dimexpr\linewidth-50mm-2\tabcolsep\relax}@{}}
	
	\caption*{\textbf{Overall authorship and preparation}}\\
	\toprule
	\rowcolor{earlytablehead}
	\textbf{Chapters or scope} &
	\textbf{Responsibility}\\
	\midrule
	\endfirsthead
	
	\multicolumn{2}{@{}l@{}}{\small\itshape
		Overall authorship and preparation (continued)
	}\\[4pt]
	\toprule
	\rowcolor{earlytablehead}
	\textbf{Chapters or scope} &
	\textbf{Responsibility}\\
	\midrule
	\endhead
	
	\midrule
	\multicolumn{2}{r@{}}{\footnotesize\itshape
		Continued on the next page
	}\\
	\endfoot
	
	\bottomrule
	\endlastfoot
	
	Chapters~\ref{ch:history}--~\ref{grdl1102}&
	Written by Yang.\\[6pt]
	\rowcolor{blue!3}
	Chapter~\ref{ch:fundamental} &
	Jointly prepared by all authors.\\[6pt]
	
	Chapter~\ref{ch:geometry},\;
	Appendix~\ref{app:oguchi}  &
	Summarized by Yang.\\[6pt]
	
	\rowcolor{blue!3}
	Appendix~\ref{tlqdzj}&
	Summarized jointly by Lee and Yang.\\[6pt]
	
	Whole volume (2026) &
	Compiled and organized by Yang.\\[6pt]
	
\end{longtable}
\endgroup

	The resulting manuscript was reviewed by the other authors. In particular,
	Professor Juno Mukai, now in his eighties, has brought undiminished
	enthusiasm and meticulous care to every stage of this work, from the
	writing of \cite{r2324}, which corresponds to
	Chapters~\ref{ch:degree23} and~\ref{ch:degree24} of this book,
	to the writing and organization of the present volume.
	Throughout this process, he has carefully reviewed successive drafts
	and generously shared his knowledge and experience with the other
	authors, whose work has benefited greatly from his patient guidance. Note that this book is self-contained in the sense that its mathematical
	proofs do not use \cite{r2324,r2567,r289}.

	\subsection*{A Note on the Use of AI}
	
	The mathematical work presented in this book was \emph{essentially complete
		by 2023} \cite{r2324,r2567,r289}, \emph{before the recent advances in AI-assisted mathematical research. }\footnote{By 2026, when the material was assembled into the present volume,
		the role of AI in mathematics had changed considerably.
		OpenAI had announced solutions to numerous longstanding open problems,
		including a solution to the Navier--Stokes existence and smoothness
		problem, one of the seven Millennium Prize Problems, accompanied
		by a formalization in Lean.}
	
	In the preparation of this volume, ChatGPT~6 was used to improve the
	\LaTeX{} typesetting and the layout of the resulting PDF,
	correct typographical errors, polish the English, and check the consistency of notation.
	Its assistance was particularly valuable in checking the subscripts
	of generators of homotopy groups, such as
	$\nu_{4+n},$   $\cg_{8+n},$   $\lt_{3+n},\mu_{3+n}$ and $\cn_{7+n},$  
	and their compatibility in compositions. ChatGPT~6 made a substantial contribution to checking the subscripts
	of homotopy group generators. In this narrowly defined task,
	and only in this respect, its checks of the present manuscript
	proved more effective than those previously carried out by
	the authors.
	\par\addvspace{1.25\baselineskip}
	
	\begin{flushright}  \textsc{Juxin Yang}\\[2pt]
		\textsc{Yoshihiro Hirato}\\[2pt]
		\textsc{Jinho Lee}\\[2pt]
		\textsc{Toshiyuki Miyauchi}\\[2pt]
		\textsc{Juno Mukai}\\[2pt]
		\textsc{Mariko Ohara}\\[6pt]
		\textit{Written in October 2026}
	\end{flushright}

	\clearpage
	\chapter*{Authors}
	\addcontentsline{toc}{chapter}{Authors}
	\markboth{Authors}{Authors}
	\begingroup
	\small
	\setlength{\parindent}{0pt}
	\setlength{\parskip}{0pt}
	\newcommand{\AuthorEntry}[3]{%
		\begin{minipage}{\linewidth}
			{\normalsize\bfseries #1}\par\smallskip
			#2\par
			\href{mailto:#3}{\nolinkurl{#3}}\par
		\end{minipage}\par\addvspace{1.15\baselineskip}}    \AuthorEntry{Juxin Yang}{%
		School of Mathematical Sciences\\
		Dalian University of Technology\\
		No.~2 Linggong Road, Ganjingzi District\\
		Dalian, Liaoning 116024, China;\\
		Department of Mathematics\\
		Stockholm University\\
		SE-106 91 Stockholm, Sweden}{juxinyanghomotopy@gmail.com}
	\AuthorEntry{Yoshihiro Hirato}{%
		Department of General Education\\
		Nagano National College of Technology\\
		716 Tokuma\\
		Nagano 381-8550, Japan}{hirato@nagano-nct.ac.jp}
	\AuthorEntry{Jinho Lee}{%
		Department of Energy Science,\\ Sungkyunkwan University (SKKU),
		Suwon, Republic of Korea;\\
		Quantum AI, \\Block S, Sejong, Republic of Korea%
	}{champang@skku.edu}

	\AuthorEntry{Toshiyuki Miyauchi}{%
		Department of Applied Mathematics, Faculty of Science\\
		Fukuoka University\\
		8-19-1 Nanakuma, Jonan-ku\\
		Fukuoka 814-0180, Japan}{tmiyauchi@fukuoka-u.ac.jp}
	
	\AuthorEntry{Juno Mukai}{%
		Shinshu University\\
		3-1-1 Asahi\\
		Matsumoto, Nagano 390-8621, Japan}{jmukai@shinshu-u.ac.jp}
	
	\AuthorEntry{Mariko Ohara}{%
		National Institute of Technology, Oshima College\\
		1091-1 Komatsu\\
		Suo-Oshima, Oshima, Yamaguchi 742-2193, Japan}{primarydecomposition@gmail.com}

	\endgroup
	\clearpage
	\begingroup
	\hypersetup{linkcolor=black,linktoc=all}
	\renewcommand{\cftchappagefont}{\bfseries\color{blue}\hypersetup{linkcolor=blue}}
	\renewcommand{\cftsecpagefont}{\color{blue}\hypersetup{linkcolor=blue}}
	\renewcommand{\cftsubsecpagefont}{\color{blue}\hypersetup{linkcolor=blue}}
	\tableofcontents
	\endgroup
	\mainmatter
	% Start every chapter on a new page; openany avoids inserted blank pages.
	\titleclass{\chapter}{top}
	\chapter{Historical Background and the Scope of the Calculations}\label{ch:history}
	Throughout the book, all spaces and maps are pointed. Moreover, spaces and maps are assumed to be \textbf{localized at 2} unless otherwise stated.  However, since localization does not in general preserve manifold structures, any discussion involving the manifold structure of a space or a Lie group structure is to be understood as referring to its original, unlocalized form; for clarity, we often emphasize in such cases that we are “working integrally.”

	How much of a vector bundle is determined by its stable class, and what information is lost when
	trivial directions are added? For bundles over a sphere, these questions lead directly to the homotopy
	groups of special orthogonal groups. The dimension of the sphere fixes the homotopy degree, while the fiber dimension of
	the bundle selects the special orthogonal group in which the calculation takes place.
	
	We begin with this geometric interpretation working integrally. It explains why the unstable groups deserve study even
	though Bott periodicity gives a concise description of their eventual stable values. We then place the
	calculations in their historical setting and explain why the generators, their orders, and their
	images under inclusion are retained throughout the book. In this book,  spaces and maps are all pointed.
	
	\section*{Vector Bundles}
	Work integrally. An oriented real vector bundle with fiber $\mathbb R^n$ over $S^{k+1}$ is trivial over each hemisphere. Its global
	information is therefore concentrated along the equator: there the two trivializations differ by a map
	into the special orthogonal group. Passing from this description to homotopy classes gives the following
	classification.     For $k\geq1,$   the clutching construction gives the natural bijection  \[
	\begin{gathered}
		\pi_k(\mathit{SO}(n))
		\\[-4pt]
		\rotatebox[origin=c]{270}{\scalebox{1.8}[1]{$\cong$}}
		\\[-3pt]
		\bigl\{
		\text{oriented real vector bundles over }S^{k+1}\text{ with fiber }\mathbb R^n
		\bigr\}
		\big/
		\text{orientation-preserving isomorphism}.
	\end{gathered}
	\]
	\noindent
	Here $\mathit{SO}(n)$ denotes the ordinary, unlocalized
	special orthogonal group. A bundle is obtained by gluing two trivial bundles over
	the hemispheres using a map $S^k\to\mathit{SO}(n)$\uprightsemicolon{}
	two oriented bundles are isomorphic exactly when their clutching maps
	are homotopic \cite{St}. Under our standing 2-localization convention, our
	$\pi_k(\mathit{SO}(n))$ is obtained by tensoring the integral homotopy group with $\mathbb Z_{(2)}.$  Free cyclic summands become $\mathbb Z_{(2)}$; for finite groups, this localization retains precisely the $2$-primary summands.

	\section{From Bundles to Homotopy Groups}
	
	The two statements above describe different stages of the same problem. Clutching keeps the fiber dimension
	fixed\uprightsemicolon{} stabilization allows it to increase. To understand what happens between these stages, we study
	the inclusions $\mathit{SO}(n)\to\mathit{SO}(n+1)$ and follow their action on individual homotopy
	classes. Chapter~\ref{ch:geometry} develops this interpretation through vector bundles, sections, and
	Thom spaces.
	
	The general theory of fiber bundles provides the geometric setting \cite{St}. Early calculations for
	special orthogonal groups include those of Sugawara \cite{Su1953,Su1954} and Kervaire \cite{Ke}. Bott periodicity
	determines the stable groups \cite{Bo}\uprightsemicolon{} the unstable groups contain additional classes and
	extension phenomena that require separate calculations. The distinction between these two problems is
	essential: a class can disappear under stabilization while remaining nonzero in the homotopy group of
	a particular special orthogonal group.
	The calculations collected here belong to the tradition of composition methods developed by Toda
	\cite{T}. They use explicit generators in the homotopy groups of spheres, relations among their
	composites, and higher compositions to determine not only the orders of groups but also the orders of
	chosen lifts. Work on spheres by Mimura, Toda, Mori, Oda, and others supplies the input for these
	arguments \cite{MT1,Mi1,MMO,Oda1979}.
	\section{The Sequence of Calculations}
	The unstable calculations build on Kachi's work \cite{Ka}, the degree $17$--$18$ and degree
	$19$--$20$ calculations of Kachi and Mukai \cite{KM1,KM2}, and the degree $21$--$22$ calculations of
	Hirato, Kachi and Mukai \cite{HKM}. Tables of classical groups and homogeneous spaces provide useful
	comparisons \cite{Lun}. The small-$n$ cases are closely tied to computations for $\mathit{SU}(3),$  
	$\mathit{SU}(4),$   and $\mathit{Sp}(2)$ \cite{MT2,Mi2,F,Oda1978}.
	The present volume continues this sequence through degrees $23$--$29.$  Chapters~\ref{ch:degree23},
	\ref{ch:degree24}, and \ref{ch:degree25} contain the degree $23$--$25$ calculations and their
	odd-primary components. Chapters~\ref{ch:degree26} and \ref{ch:degree27} treat degrees $26$ and $27.$ 
	Chapters~\ref{ch:degree28} and \ref{ch:degree29} contain the final calculations. The preliminary
	formulas are collected in Chapter~\ref{ch:fundamental}.
	Early contributions of Sugawara (1953--1954) \cite{Su1953,Su1954} precede this systematic sequence.
	The table distinguishes the stable groups, the families $\pi_{n+r}(\mathit{SO}(n)),$  
	and calculations for specified degrees or values of $n.$  The stated ranges and any unresolved extensions
	are part of the scope of each result. The exceptional case $n=2$ is elementary:
	$\mathit{SO}(2)\cong S^1,$   so $\pi_1(\mathit{SO}(2))=\mathbb Z_{(2)}$ and $\pi_i(\mathit{SO}(2))=0$ for $i\geq2.$ 
	\begingroup
	\small
	\setlength{\LTpre}{6pt}
	\setlength{\LTpost}{6pt}
	\setlength{\tabcolsep}{5pt}
	\renewcommand{\arraystretch}{1.25}
	\arrayrulecolor{earlytablerule}
	\begin{longtable}{@{}
			>{\raggedright\arraybackslash}p{13mm}
			>{\raggedright\arraybackslash}p{34mm}
			>{\raggedright\arraybackslash}p{24mm}
			>{\raggedright\arraybackslash}p{\dimexpr\linewidth-71mm-6\tabcolsep\relax}
			@{}}
		\caption{Calculations of Homotopy Groups of Special Orthogonal Groups}
		\label{tab:special orthogonal-history}\\
		\toprule
		\rowcolor{earlytablehead}
		\textbf{Year} &
		\textbf{Authors} &
		\textbf{Degree $i$} &
		\textbf{Results for $\mathit{SO}(n)$}\\
		\midrule
		\endfirsthead
		
		\multicolumn{4}{@{}l@{}}{
			\small\itshape Table~\thetable\ continued
		}\\[4pt]
		\toprule
		\rowcolor{earlytablehead}
		\textbf{Year} &
		\textbf{Authors} &
		\textbf{Degree $i$} &
		\textbf{Results for $\mathit{SO}(n)$}\\
		\midrule
		\endhead
		
		\midrule
		\multicolumn{4}{r@{}}{
			\footnotesize\itshape Continued on the next page
		}\\
		\endfoot
		
		\bottomrule
		\endlastfoot
		
		1942 & G.~W.~Whitehead & $1\leq i\leq5$ &
		The groups $\pi_i(\mathit{SO}(n))$ for all $n\geq2$
		in these degrees, with explicit generators
		\cite{Wh42}.\\[5pt]
		
		\rowcolor{blue!3}
		1953 & M.~Sugawara & $i=6,7,8$ &
		Calculations of $\pi_i(\mathit{SO}(n))$ in these degrees
		\cite{Su1953}\uprightsemicolon{}
		the results in degree $7$ were corrected in
		\cite{Su1954}.\\[5pt]
		
		1954 & M.~Sugawara & $i=9,10,11$ &
		The groups $\pi_9(\mathit{SO}(n))$ for all $n\geq3$
		\uprightsemicolon{}
		partial results in degrees $10$ and $11$, with explicitly
		stated ambiguities\uprightsemicolon{}
		corrections to the earlier results in degree $7$
		\cite{Su1954}.\\[5pt]
		
		\rowcolor{blue!3}
		1959 & R.~Bott & All $i\geq1$ &
		The groups $\pi_i(\mathit{SO})$ and their period-eight
		pattern\uprightsemicolon{}
		these give $\pi_i(\mathit{SO}(n))$ for $n\geq i+2$
		\cite{Bo}.\\[5pt]
		
		1960 & M.~A.~Kervaire & $i=n+r,$ $-1\leq r\leq4$ &
		Tables of $\pi_{n+r}(\mathit{SO}(n))$ for $-1\leq r\leq4$,
		with an explicitly stated ambiguity\uprightsemicolon{}
		relations for the connecting homomorphism
		\cite{Ke}.\\[5pt]
		
		\rowcolor{blue!3}
		1968 & H.~Kachi & $i=15,16$ &
		The groups $\pi_i(\mathit{SO}(n))$ for all $n\geq3$
		in these degrees, with explicit generators
		\cite{Ka}.\\[5pt]
		
		1999 & H.~Kachi and J.~Mukai & $i=17,18$ &
		The groups $\pi_i(\mathit{SO}(n))$ for all $n\geq3$
		in these degrees, with explicit generators
		\cite{KM1}.\\[5pt]
		
		\rowcolor{blue!3}
		2000 & H.~Kachi and J.~Mukai & $i=19,20$ &
		The groups $\pi_i(\mathit{SO}(n))$ for all $n\geq3$
		in these degrees, with explicit generators
		\cite{KM2}.\\[5pt]
		
		2006 & Y.~Hirato, H.~Kachi and J.~Mukai & $i=21,22$ &
		The groups $\pi_i(\mathit{SO}(n))$ for all $n\geq3$
		in these degrees, with explicit generators
		\cite{HKM}.\\[5pt]
	\end{longtable}
	\endgroup

	Table~\ref{tab:related-lie-group-history} records six related studies in chronological order.
	Their connection with the special orthogonal groups is made precise by the coverings and
	fibrations in Section~\ref{lqzh}. For $i\geq2,$   the coverings
	$\mathit{Sp}(2)\cong\mathit{Spin}(5)\to\mathit{SO}(5)$ and
	$\mathit{SU}(4)\cong\mathit{Spin}(6)\to\mathit{SO}(6)$ induce isomorphisms on $\pi_i.$ 
	The fibration $\mathit{SU}(3)\to G_2\to S^6,$   together with the $2$-local splitting
	$\mathit{Spin}(7)\simeq G_2\times S^7,$   also connects the calculations for
	$\mathit{SU}(3)$ with those for $\mathit{SO}(7).$

	\begingroup
	\small
	\setlength{\LTpre}{6pt}
	\setlength{\LTpost}{6pt}
	\setlength{\tabcolsep}{5pt}
	\renewcommand{\arraystretch}{1.2}
	\arrayrulecolor{earlytablerule}
	\begin{longtable}{@{}>{\raggedright\arraybackslash}p{30mm}>{\raggedright\arraybackslash}p{62mm}>{\raggedright\arraybackslash}p{\dimexpr\linewidth-92mm-4\tabcolsep\relax}@{}}
		\caption{Calculations for Other Lie Groups and Their Applications to $\mathit{SO}(n)$}
		\label{tab:related-lie-group-history}\\
		\toprule
		\rowcolor{earlytablehead}
		\textbf{Year and reference}
		& \textbf{Results}
		& \textbf{Consequences for $\mathit{SO}(n)$}\\
		\midrule
		\endfirsthead
		
		\multicolumn{3}{@{}l@{}}{\small\itshape Table~\thetable\ continued}\\[4pt]
		\toprule
		\rowcolor{earlytablehead}
		\textbf{Year and reference}
		& \textbf{Results}
		& \textbf{Consequences for $\mathit{SO}(n)$}\\
		\midrule
		\endhead
		
		\midrule
		\multicolumn{3}{r@{}}{\footnotesize\itshape Continued on the next page}\\
		\endfoot
		\bottomrule
		\endlastfoot
		
		\textbf{1964}\par M.~Mimura and H.~Toda\par \cite{MT2}
		&
		Homotopy groups of $\mathit{SU}(3),$   $\mathit{SU}(4),$  
		and $\mathit{Sp}(2)$ for $i\leq23,$   with generators.
		&
		$\pi_i(\mathit{SO}(5))$ and $\pi_i(\mathit{SO}(6))$
		for $2\leq i\leq23.$ 
		\\[7pt]
		
		\rowcolor{blue!3}
		\textbf{1964}\par M.~Mimura and H.~Toda\par \cite{MT3}
		&
		Groups $\pi_i(\mathit{Sp}(n))$ for $i\leq23,$  
		selected cases with $i=24,$   and the families
		$i=4n+2,4n+3,4n+4.$ 
		&
		Relations in $\pi_i(\mathit{SO}(5))$ for $i=21,22.$ 
		\\[7pt]
		
		\textbf{1964}\par K.~\^{O}guchi\par \cite{Og}
		&
		Generators and composition relations in
		$\pi_i(\mathit{SU}(n))$ and $\pi_i(\mathit{Sp}(n))$
		for $i\leq22.$ 
		&
		Explicit generators and composition relations in
		$\pi_i(\mathit{SO}(5))$ and $\pi_i(\mathit{SO}(6))$
		for $14\leq i\leq22.$ 
		\\[7pt]
		
		\rowcolor{blue!3}
		\textbf{1967}\par M.~Mimura\par \cite{Mi2}
		&
		Homotopy groups of $G_2,$   $F_4,$   $\mathit{Spin}(7),$  
		and $\mathit{Spin}(9)$ for $i\leq23.$ 
		&
		$\pi_i(\mathit{SO}(7))$ and $\pi_i(\mathit{SO}(9))$
		for $2\leq i\leq22$;    
		$\pi_{23}(\mathit{SO}(7))$
		and $\pi_{23}(\mathit{SO}(9)),$  
		with unresolved group extensions \\[7pt]
		
		\textbf{1978}\par Y.~Furukawa\par \cite{F}
		&
		Homotopy groups of $\mathit{SU}(3),$   $\mathit{SU}(4),$  
		and $\mathit{Sp}(2)$ in degrees $24$--$27.$ 
		&
		$\pi_i(\mathit{SO}(5))$ and $\pi_i(\mathit{SO}(6))$
		for $24\leq i\leq27.$ 
		\\[7pt]
		
		\rowcolor{blue!3}
		\textbf{1978}\par N.~Oda\par \cite{Oda1978}
		&
		Homotopy groups of $\mathit{SU}(3),$   $\mathit{SU}(4),$  
		and $\mathit{Sp}(2)$ in degrees $24$--$36,$  
		with generators and composition relations.
		&
		$\pi_i(\mathit{SO}(5))$ and $\pi_i(\mathit{SO}(6))$
		for $24\leq i\leq36,$   with generators and composition relations.
		\\[7pt]
	\end{longtable}
	\endgroup
	\section{Why Generators and Relations Matter}
	
	The algebraic structure of homotopy groups extends beyond their
	decomposition into cyclic summands. Composition of maps relates
	elements in different dimensions: for a based space $X,$  
	classes $\alpha\in\pi_q(X)$ and $\beta\in\pi_r(S^q)$ determine
	$\alpha\circ\beta\in\pi_r(X).$ 
	These composition operations are associative.
	The homotopy groups of spheres provide operations on the
	homotopy groups of other spaces, and composition relations among
	homotopy classes of maps between spheres yield corresponding relations in those groups.
	In the stable homotopy groups of spheres, composition gives a
	graded-commutative ring structure. Unstable composition retains
	more information, although it need not be additive in both
	variables.
	
	Explicit generators allow us to describe this additional structure. Composition relations describe how generators interact and
	allow one calculation to be used in another.
	When successive composites vanish, Toda brackets provide
	secondary compositions, with a specified indeterminacy.
	These operations are central to Toda's composition methods
	\cite{T} and to the calculations in this book.
	
	For the special orthogonal groups, these methods are applied
	to the fibration
	$\mathit{SO}(n)\xrightarrow{i_{n+1}}
	\mathit{SO}(n+1)\xrightarrow{p_{n+1}}S^n.$ 
	Its homotopy group exact sequence relates the homotopy groups of $\mathit{SO}(n)$ and $\mathit{SO}(n+1)$ to those of $S^n.$ 
	Determining the kernels and cokernels gives a short exact
	sequence, whose extension must then be resolved.
	This requires relations among explicit lifts and the  elements  
	in the image of the inclusion.
	Composition relations, Toda brackets, and the unstable
	$J$-homomorphism supply the information needed to determine
	these lifts and their orders.
	The resulting calculations therefore describe both the groups
	and the behavior of their generators under the maps and
	operations used throughout the book.

	\section{Groups $\pi_k(\mathit{SO}(n))$ in the Stable or the  Metastable Range}
	For $k\geq1$ and $n\geq k+2,$   the inclusion into the stable special orthogonal group induces an isomorphism on
	$\pi_k.$  Indeed, in the fibration above, both $\pi_{k+1}(S^n)$ and $\pi_k(S^n)$ vanish when $k+1<n.$ 
	Iterating these isomorphisms reaches the stable group. Bott periodicity \cite{Bo} then gives the
	following localized values.
	\begin{center}
		\begin{tabular}{@{}cccccccc@{}}
			\toprule
			$k$ & $23$ & $24$ & $25$ & $26$ & $27$ & $28$ & $29$\\
			\midrule
			$\pi_k(\mathit{SO})$ & $\mathbb Z_{(2)}$ & $\mathbb Z/2$ & $\mathbb Z/2$ & $0$ & $\mathbb Z_{(2)}$ & $0$ & $0$\\
			$\pi_k(\mathit{SO}(n))$ is  stable when $n\geq$ & $25$ & $26$ & $27$ & $28$ & $29$ & $30$ & $31$\\
			\bottomrule
		\end{tabular}
	\end{center}
	The unstable tables must eventually agree with these values. They also describe the information
	that the stable groups cannot retain: unstable torsion, the changing orders of lifts, and the precise
	stage at which a class vanishes.
	
	This   gives a destination for each calculation. The work in the later chapters is to
	describe the groups for intermediate values of $n$ and the maps joining them. Before taking up those calculations, we
	develop the Toda-bracket constructions used to determine relations among the generators.
	
	The table in \cite[pp.~258--261]{Lun} gives the groups
	$\pi_{n-2+k}(\mathit{SO}(n))$ for $n>9$ and $1\leq k\leq14,$  
	subject to the additional restrictions on $n$ in individual entries.
	They use $K$-theory, $\operatorname{Im}(J)$-theory, and other methods,
	but do not provide explicit generating sets of the kind given by Toda.
	The groups in the higher stable range follow from Bott periodicity.
	For groups covered by these results, our remaining task is to give explicit generating sets.
	Even this is a nontrivial problem.
	
	More precisely, for $i=23,24,\ldots,29,$   we seek to determine
	$\pi_i(\mathit{SO}(n))$ together with explicit generators.
	Applying the restrictions in the cited tables with $k=i-n+2,$  
	and adjoining the stable values, gives the ranges below.
	The complementary ranges require calculations or sources beyond these tables;
	this does not assert that all such groups are unknown elsewhere.
	
	\begingroup
	\renewcommand{\arraystretch}{1.15}\setlength{\tabcolsep}{6pt}
	\setlength{\LTleft}{0pt}\setlength{\LTright}{0pt}
	\setlength{\LTpre}{5pt}\setlength{\LTpost}{6pt}
	\setlength{\LTcapwidth}{\linewidth}\arrayrulecolor{earlytablerule}
	\begin{longtable}{>{\centering\arraybackslash}p{23mm}>{\centering\arraybackslash}p{85mm}>{\centering\arraybackslash}p{\dimexpr\linewidth-108mm-6\tabcolsep\relax}}
		\hline
		\rowcolor{earlytablehead}Degree $i$ & Range beyond the cited coverage & Known Range \\\hline
		\endfirsthead
		\multicolumn{3}{c}{\small\itshape Table~\thetable\ continued}\\
		\hline
		\rowcolor{earlytablehead}Degree $i$ & Range beyond the cited coverage & Known Range\\\hline
		\endhead
		\hline
		\endfoot
		\endlastfoot
		$23$ & $n\leq12$ & $n\geq13$\\
		\rowcolor{blue!3}
		$24$ & $n\leq13$ & $n\geq14$\\
		$25$ & $n\leq13$ & $n\geq14$\\
		\rowcolor{blue!3}
		$26$ & $n\leq14$ & $n\geq15$\\
		$27$ & $n\leq14$ & $n\geq15$\\
		\rowcolor{blue!3}
		$28$ & $n\leq15$ & $n\geq16$\\
		$29$ & $n\leq16$ & $n\geq17$\\*
		\hline
		\caption{Ranges for $\pi_i(\mathit{SO}(n))$ covered by the table in \cite[pp.~258--261]{Lun}  and Bott periodicity, and their complementary ranges. In the known range,  explicit generators still require further work.}\label{tab:known-range-so}
	\end{longtable}
	\endgroup
	
	The equal lower bounds for degrees $24$ and $25$ come from the restrictions in the individual entries.
	At $n=13,$   these degrees correspond to $k=13$ and $k=14,$   respectively, in the column $n\equiv5\pmod8$;
	both entries require $n>13.$  At $n=14,$   the entries for both degrees apply, so both ranges begin at $n=14.$ 
	For degrees $26$ and $27,$   at $n=14$ the former corresponds to $k=14$ in the column $n\equiv6\pmod8,$   whose entry requires $n>14,$  
	whereas the latter would require $k=15,$   outside the tables.
	At $n=15,$   they correspond to $k=13$ and $k=14$ in the column $n\equiv7\pmod8,$   and both entries apply.
	Thus both of these ranges begin at $n=15.$

	\chapter{Toda Brackets}\label{ch:toda}
	
	The Toda bracket is a tool for constructing homotopy liftings and extensions of maps. It plays a
	fundamental role in studying composition relations among homotopy classes.  A Toda bracket begins when two successive composites
	vanish. Choices of their null-homotopies can then be joined to produce
	a new class one dimension higher. The choices matter: in general the
	answer is a coset.

	Both Toda brackets and Massey products
	have their origins in the mathematical
	developments of the early 1950s.
	The published record begins, among
	the works considered here, with
	Massey's announcement in 1951
	\cite{Massey1951}, followed by
	Toda's study of triple compositions
	in 1952 \cite{Toda1952}.
	While this establishes a clear
	chronological order among the
	publications, it leaves open the
	historically distinct questions
	of when the respective constructions
	were first conceived and whether
	either development influenced the other.
	
	Toda's 1952 paper is remarkable not only for
	introducing the triple Toda bracket, but also
	for demonstrating its effectiveness in the
	computation of homotopy groups of spheres.
	In Chapter~5(ii), under the assumptions
	$\alpha\beta=0$ and $\beta\gamma=0$, Toda
	constructed the bracket $\{\alpha,\beta,\gamma\}$
	by joining two chosen nullhomotopies and
	carefully described the resulting indeterminacy;
	see \cite[p.~70, (5.5)--(5.6)]{Toda1952}.
	He subsequently established a connection between
	this construction and maps into a cell attachment
	\cite[p.~72, Lemma~(5.9)]{Toda1952}.
	
	An illuminating application is provided by
	the relations
	\[
	\mathbf{x}\in\{\eta_r,2\iota_{r+1},\eta_{r+1}\},
	\qquad
	2\mathbf{x}=\eta_r\eta_{r+1}\eta_{r+2}\ne0
	\quad(r\ge3),
	\]
	which yield an element
	$\mathbf{x}\in\pi_{r+3}(S^r)$ of order $4$
	\cite[p.~72, Lemma~(5.10)]{Toda1952}.
	Thus, long before the appearance of his
	celebrated monograph \cite{T}, Toda had
	already demonstrated how triple brackets
	could serve as a powerful means of constructing
	homotopy classes and determining their orders.

	In this chapter, we \textbf{work integrally}. Recall that throughout this book,  spaces and maps are pointed. For trivial homotopy classes $\mathbf{x}$ and $\mathbf{y}$, following Toda, we write $$\mathbf{x}=\mathbf{y}=0$$  even though their domains or  codomains may not be  the same.
	
	\section{Extensions and Coextensions}\label{sec:toda-cones}
	
	A null-homotopy carries more information than the assertion that a composite is zero. The mapping cone
	allows us to use that information as part of a map: an extension continues a map across an attached
	cone, while a coextension uses the two cones of a suspension. Keeping these choices visible will
	explain the indeterminacy of the Toda bracket.
	
	For the sequence of based spaces and maps
	\[\begin{aligned}
		&L\xrightarrow{f}M\xrightarrow{g}N, \qquad gf\simeq *,
	\end{aligned}\]
	let $C_f=M\cup_f CL$ be the mapping cone of $f.$ 
	A map $\bar g:C_f\to N$ is called an
	\emph{extension of $g$ with respect to $f$} if the diagram
	\[
	\begin{tikzcd}[column sep=large,row sep=large]
		M \arrow[r,"g"] \arrow[d,hook,"i_f"'] & N \\
		C_f \arrow[ur,"\bar g"'] &
	\end{tikzcd}
	\]
	commutes strictly. We also denote this map by
	$\operatorname{ext}_f(g).$ 
	
	The coextension is the complementary construction. One half of the suspension uses $f,$   and the other
	uses a chosen null-homotopy of $gf.$  The signs in the construction determine the suspension convention
	used in the bracket below.
	
	A \emph{coextension of $f$ with respect to $g$} is a map
	\[\begin{aligned}
		&\widetilde f :\Sigma L=C^+L\cup    C^-L\longrightarrow C_g
	\end{aligned}\]
	constructed as follows. Use the reduced cone functor
	$C(-)=(-)\wedge I,$   where $I=[0,1]$ has basepoint $1.$ 
	Define $\widetilde f ^+$ by the composite
	\[\begin{aligned}
		&C^+L\xrightarrow{\mathrm{id}}CL \xrightarrow{ C(f) }CM \xrightarrow{j_{CM}}N\vee CM \xrightarrow{q}C_g,
	\end{aligned}\]
	where $q$ is the quotient map defining $C_g=N\cup_g CM.$  Note that $C(f)=f\wedge \mathrm{id}_{I}.$ 
	Choose an extension $\overline{gf}:CL\to N$ of $gf$
	with respect to $\mathrm{id}_L,$   and define $\widetilde f ^-$ by the composition
	\[\begin{aligned}
		&C^-L\xrightarrow{\mathrm{id}}CL \xrightarrow{\overline{gf}}N \xrightarrow{j_N}C_g.
	\end{aligned}\]
	The restrictions $\widetilde f ^+|_L$ and $\widetilde f ^-|_L$ agree,
	so these maps give
	\[\begin{aligned}
		&\widetilde f =-(\widetilde f ^+\cup   \widetilde f ^-)=(\widetilde f ^+\cup   \widetilde f ^-)\hc[-1]: C^+L\cup    C^-L\longrightarrow C_g.
	\end{aligned}\]
	
	\noindent   Here, $[-1]$ denotes the map of degree $-1$ which is given by $$\s L=L\wedge S^{1}\to \s L=L\wedge S^{1},\;\ell\wedge e^{\theta\sqrt{-1}}\mapsto \ell\wedge e^{-\theta\sqrt{-1}},\;(\theta\in [0,2\pi)).$$\noindent
	We also denote this map by $\operatorname{coext}_g(f).$  With this convention, $$\begin{aligned}
		&(C_{g}\stackrel{proj}\twoheadrightarrow \Sigma M)\circ \operatorname{coext}_g(f)=-\Sigma f \in [ \Sigma L, \Sigma M ].
	\end{aligned}$$
	
	In general, neither the extension nor the coextension is unique
	up to homotopy. Suppose that
	$\alpha\in[L,M],$   $\beta\in[M,N],$   and $\beta\alpha=0.$ 
	For representative maps $a\in\alpha$ and $b\in\beta,$   the notation
	$\operatorname{ext}_{\alpha}(\beta)$ denotes the homotopy class
	of a chosen extension $\operatorname{ext}_a(b).$ 
	Likewise, $\operatorname{coext}_{\beta}(\alpha)$ denotes the
	homotopy class of a chosen coextension
	$\operatorname{coext}_b(a).$ 
	Thus both the representatives and the choices of extensions
	and coextensions are allowed to vary. Denote by $\operatorname{Ext}_\alpha(\beta)$   the set of all
	extensions of $\beta$  with respect to $\alpha$\uprightsemicolon{}  denote  by  $\operatorname{Coext}_\beta(\alpha)$
	the set of all coextensions of $\alpha$ with respect to $\beta.$ 
	
	We can now combine the two constructions. For three composable classes, the vanishing of the first
	adjacent composite permits an extension, and the vanishing of the second permits a coextension.
	Composing these maps produces the bracket value\uprightsemicolon{} varying the choices produces the bracket as a set.
	
	It is worth noting that, in calculations involving Toda brackets, maps are often written with arrows pointing from right to left, as in
	\[
	Y \xleftarrow{f} X.
	\] This was the practice Toda himself followed when demonstrating
	calculations involving Toda brackets at the blackboard.  Prominent mathematicians such as Mimura, $\mathrm{\hat{O}}$guchi, Oda, and $\mathrm{\overline{O}}$shima also follow this convention.
	
	For a map or homotopy class $f\colon Y\leftarrow X,$   we use the identification $C_{\s f}=\s C_{f}.$ \footnote{This identification requires care, as it can affect signs. We follow
		Toda's convention, in which these sign issues are accounted for.
		Thus the properties of Toda brackets stated in his book \cite{T}
		hold exactly, not merely up to sign.
		See \cite[Equation~(1.16) and Proposition 1.7, p.~13]{T}.}
	\begin{defn}\label{def:toda}
		Let $n\geq0$ and consider homotopy classes
		\[\begin{aligned}
			&Z\xleftarrow{\alpha}\Sigma^nY \xleftarrow{\Sigma^n\beta}\Sigma^nX \xleftarrow{\Sigma^n\gamma}\Sigma^nW,
		\end{aligned}\]
		where $\alpha\Sigma^n\beta=0$ and $\beta\gamma=0.$ 
		The Toda bracket indexed by $n,$   denoted by
		$\{\alpha,\Sigma^n\beta,\Sigma^n\gamma\}_n,$  
		is the set of all compositions
		\[\begin{aligned}
			&(-1)^n \operatorname{ext}_{\Sigma^n\beta}(\alpha) \circ\Sigma^n\operatorname{coext}_{\beta}(\gamma).
		\end{aligned}\]
		More explicitly, it is the set of homotopy classes
		\begin{equation*}\notag
			\begin{aligned}
				&(-1)^n\bigl[ \operatorname{ext}_{\Sigma^n b}(a) \circ\Sigma^n\operatorname{coext}_b(c) \bigr] \in[\Sigma^{n+1}W,Z],
			\end{aligned}
		\end{equation*}
		where $a\in\alpha,$   $b\in\beta,$   and $c\in\gamma,$  
		and all choices of extensions and coextensions are admitted.
		The bracket $\{\alpha,\beta,\gamma\}_0$ is abbreviated to
		$\{\alpha,\beta,\gamma\}.$ 
	\end{defn}
	\begin{rem} For a fixed $n,$   the above Toda bracket is a function of the triad $(\af, \beta,\gamma)$; more specifically speaking, for a given  integer $n,$   the  Toda bracket $\n \af,\s^{n}\beta,\s^{n}\gamma\nn_{n}$   depends only on $(\af, \beta,\gamma)$ but not $(\af, \s^{n}\beta,\s^{n}\gamma).$   It is necessary to point out that  even if\; $$\s^{n}\beta'=\s^{n}\beta, \s^{n}\gamma'=\s^{n}\gamma \; \text{and}\; \beta'\gamma'=0,$$ \noindent we have
		$$\n \af,\s^{n}\beta',\s^{n}\gamma'\nn_{n}\neq\n \af,\s^{n}\beta,\s^{n}\gamma\nn_{n}\;\text{in general}.$$
		
		However, in accordance with Toda's convention, the suspension symbols $\s^n$ are often left implicit in brackets of the form $\n-,\s^n-,\s^n-\nn_n$ with $n\geq1.$  For example, in Toda's book \cite{T}, $\n\ca_6,2\e_7,\lt_7\nn_1$ means $\n\ca_6,\s(2\e_6),\s\lt_6\nn_1$; the third entry is understood to arise from $\lt_6,$   although $\lt_6+2\ch_6\neq \lt_6$ and $\s(\lt_6+2\ch_6)=\lt_7.$
	\end{rem}

	It is noteworthy that the  Toda bracket indexed by $n$ with $n\geq1$ makes better use of the desuspension properties of homotopy classes and plays a critical role in the theory of computing unstable homotopy groups. By presenting the following three arguments, we substantiate  in the unstable range the essential significance of the condition $n\geq1$ on the indexing number.\begin{rem}\label{xbh}
		
		\begin{itemize} 
			\item [\rm(1)]   It proves to be advantageous for defining new elements in $\pi_{*}(S^{i}).$  For example, recall from Toda's work \cite[Lemma 6.1, p. 51]{T} that  $\lt_{3}\in\n\ca_{3},\s\nu', \nu_{7} \nn_{1}.$  The presence of the subscript 1 enables us to utilize the $H$-formula, $$H(\lt_{3})\in H\n\ca_{3},\s\nu', \nu_{7} \nn_{1}=-P^{-1}(\ca_{2}\hc\nu')\hc\nu_{8}=\nu_{5}^{2},$$ thereby  confirming $\lt_{3}$ to be a lifting of $\nu_{5}^{2}.$  This, in fact, constitutes the fundamental method for constructing new elements of homotopy groups of spheres in \cite{T}. The vast majority of new elements in \cite{T} are constructed using this method.

			\item [\rm(2)]   It is useful to consider the connecting homomorphism. Suppose we aim to determine $\pi_{k}(X),$   with $\pi_{*}(Y)$ and $\pi_{*}(F)$ being already known, where $F\rightarrow X \rightarrow Y$ is a homotopy fibration.
			We need to transform the long exact sequence
			\[\begin{tikzcd}[column sep=0.4cm]\
				&&\pi_{k+1}(Y) \arrow[r,"\pa "] & \pi_{k}(F) \arrow[r] & \pi_{k}(X)\arrow[r] & \pi_{k}(Y) \arrow[r,r,"\pa "] & \pi_{k-1}(F) 
			\end{tikzcd}\]
			
			into a short exact sequence\[\begin{tikzcd}[column sep=0.4cm]
				&0\arrow[r,] & \mathrm{coker}(\pa\colon \pi_{k+1}(Y)\rightarrow \pi_{k}(F) ) \arrow[r] & \pi_{k}(X)\arrow[r] & \mathrm{ker}(\pa \colon \pi_{k}(Y)\rightarrow \pi_{k-1}(F) ) \arrow[r] & 0.
			\end{tikzcd}\]
			
			For $\mathbf{x}\in\n\af, \s^{n}\beta,\s^{n}\gamma \nn_{n}\xyd\pi_{k+1}(Y),$   the relation $$\pa (\mathbf{x})\in \pm \n \pa (\af), \s^{n-1}\beta,\s^{n-1}\gamma \nn_{n-1}$$   holds for $n\geq1$ but  not for $n=0.$  Similarly, when we consider the image of $$\mathbf{x}\in\n\af, \s^{n}\beta,\s^{n}\gamma \nn_{n}\xyd\pi_{k+2}(S^{2m+1}) $$\noindent under the homomorphism $P:\pi_{k+2}(S^{2m+1})\rightarrow\pi_{k}(S^{m}) $    in the 2-local   EHP sequence, the   relation $$P(\mathbf{x})\in \pm \n P(\af), \s^{n-2}\beta,\s^{n-2}\gamma \nn_{n-2}$$   holds for $n\geq2$ but  not for $n\in\n0, 1\nn.$ 
			
			\item [\rm(3)]   When studying unstable homotopy groups, we can trace the calculation back to the dimension before the birth of certain troublesome elements, thereby reducing interference. For instance, notice the proof of determining $P(\lt_{13})$ from \cite[p. 78]{T}. It tells us that, after localization at 2, $$P(\lt_{13})\in \n \nu_{6},\ca_{9}, 2\e_{10}\nn_{6}\hc(-\lt_{11})=\nu_{6}\hc \s^{6}\n\ca_{3},2\e_{4}, \lt_{4}
			\nn\xyd \nu_{6}\hc\s^{6}\pi_{13} (S^{3}).$$ Computing $\nu_{6}\hc\s^{6}\pi_{13} (S^{3})$ is easier than 
			computing $\nu_{6}\hc\pi_{19} (S^{9})$ due to the fact that some unstable elements in $\pi_{13} (S^{3})$ are killed by the iterated suspension functor $\s^{6}.$  An alternative perspective is that formula (\ref{ind1022}) reveals it is typically more practical to calculate the indeterminacy of the Toda bracket indexed by $n\geq1,$   as some troublesome elements usually appear in $[\s^{n+1} W,\s^{n} Y]$ but not in $[\s W, Y].$

		\end{itemize}
	\end{rem}

	\begin{lem}[Indeterminacy {\cite{T}}]
		\label{lem:toda-indeterminacy}
		Under the hypotheses of Definition~\ref{def:toda}, the bracket is a
		double coset of the subgroups
		\[\begin{aligned}
			&[\Sigma^{n+1}X,Z]\circ\Sigma^{n+1}\gamma, \qquad \alpha\circ\Sigma^n[\Sigma W,Y].
		\end{aligned}\]
		in $[\Sigma^{n+1}W,Z].$  If this group is abelian, the bracket is
		a coset with indeterminacy
		\begin{equation}\label{ind1022}
			\begin{aligned}
				&\operatorname{Ind}\{\alpha, \Sigma^n\beta, \Sigma^n\gamma\}_n =\alpha\circ\Sigma^n[\Sigma W,Y] +[\Sigma^{n+1}X,Z]\circ\Sigma^{n+1}\gamma.
			\end{aligned}
		\end{equation}
		In particular, the abelian conclusion holds if $n\geq1,$   if
		$W=\Sigma W',$   or if $Z=\Omega Z'.$ 
		\qed
	\end{lem}
	
	The two summands have a geometric meaning. Changing the
	null-homotopy of $\beta\gamma$ changes the coextension and contributes
	the first summand of \eqref{ind1022}. Changing the
	null-homotopy of $\alpha\circ \Sigma^n\beta$ changes the extension and
	contributes the second. Thus finding one representative determines
	the whole bracket only after both contributions have been evaluated.

	Recall that $\operatorname{Ext}_\alpha(\beta)$  is  the set of all
	extensions of $\beta$  with respect to $\alpha,$   and that   $\operatorname{Coext}_\beta(\alpha)$
	is  the set of all coextensions of $\alpha$ with respect to $\beta.$    The extensions and coextensions themselves can also be expressed in terms of Toda brackets. The following statements are well-known to experts.
	
	\begin{prop} \label{dsjd}Let $Z\xleftarrow{\alpha}{Y}\xleftarrow{\beta}{X}$ be given, where $\alpha\beta=0,$   and let  
		\begin{equation*}
			\begin{tikzcd}
				\cdots
				& \s C_{\alpha} \arrow[l]
				& \s Z \arrow[l, "{-\s j}"']
				& \s Y \arrow[l, "{-\s\af}"']
				& C_{\alpha} \arrow[l, "q"']
				& Z \arrow[l, "j"', "\dyd"]
				& Y \arrow[l, "\af"']
			\end{tikzcd}
		\end{equation*}
		\noindent be the canonical cofiber sequence induced by $\af.$ Then  \[
		\mathrm{Coext}_{\alpha}(\beta)
		\subseteq
		\{j,\alpha,\beta\},\]   
		\vspace{-1em}$$  q\hc\mathrm{Coext}_{\alpha}(\beta)
		=\n-\s\beta\nn,
		\,\mathrm{(a \;singleton\; set),}        $$
		\vspace{-1em} \[
		\s\mathrm{Coext}_{\alpha}(\beta)
		\xyd\pm\mathrm{Coext}_{\s\alpha}(\s\beta),\quad \s\mathrm{Ext}_{\beta}(\alpha)\xyd\pm\mathrm{Ext}_{\s\beta}(\s\alpha).
		\] 
		\noindent   Moreover,   we have $$\mathrm{id}_{\s Y}\in \n q,j,\af\nn,\quad\mathrm{id}_{\s Z}\in  \n-\s\af,q,j\nn\;\text{   and   } \;\;\mathrm{id}_{\s C_{\af}}\in\n \s j,\s\af,q\nn.   $$
		If $\beta$ is a suspended homotopy class, denote by $\s^{-1}\beta\colon \s ^{-1}Y\leftarrow \s ^{-1}X$  one of its desuspensions. Then $$\begin{aligned}
			&\pm\mathrm{Ext}_{\beta}(\alpha)\xyd\{ \alpha, \beta, ( X \stackrel{ proj.}\twoheadleftarrow C_{\Sigma^{-1}\beta}) \}.
		\end{aligned}$$
		
	\end{prop}
	
	\begin{proof}
		The map $\mathrm{id}\colon C_{\af}\leftarrow C_{\af}$ is an extension of $j$ with respect to $\af$; namely, we can take $\mathrm{ext}_{\af}(j)=\mathrm{id}.$ So the first relation holds. Similarly, we may take  $$\mathrm{coext}_{\beta}(  X \stackrel{ proj.}\twoheadleftarrow C_{\Sigma^{-1}\beta})=(\mathrm{id}\colon  C_{\beta}\leftarrow C_{\beta})$$ \noindent up to the identification $C_{\beta}=\s C_{\s^{-1}\beta}.$  
		So the last relation holds. Note that 
		the signs  come from the switching maps $T|_{C(\s X)}$ given by $$\s C_{f}=\s(Y\cup_{f} CX)= (\s Y)\cup_{\s f} (\s CX)\xrightarrow[\simeq]{T} (\s Y)\cup_{\s f}  C(\s X),$$
		$$\s (CX)=(X\wedge I)\wedge S^{1}\xrightarrow[\cong]{T|_{C(\s X)}} C(\s X)=(X\wedge S^{1})\wedge I,$$ \noindent where $f\colon X\to Y$ is a map.  The remaining relations are straightforward and are left to the reader.   \end{proof}

	The following result was essentially observed by Spanier in 1962
	\cite[Corollary~(3.4)$^{\prime}$]{Spanier1962}. However, Spanier used the
	``co-Toda bracket'', defined in terms of loop spaces and fiber
	sequences, which is the ``co-lifts $\circ$ lifts'' analogue  
	of the classical Toda bracket. Note that the two kinds of brackets are the same up to sign and the natural isomorphisms $[-,\lo-]\cong[\s-,-].$ 
	We provide a proof without using Spanier's ``co-Toda bracket''.

	\begin{prop}\label{fby1}
		Let $f:Y\to X$ be a map and let
		\[
		X\xleftarrow{\,f}Y\xleftarrow{\,j}F
		\xleftarrow{\,\partial}\Omega X
		\]
		be the canonical fiber sequence induced by $f.$  Then
		\[
		\operatorname{ev}\in\{f,j,\partial\}
		\subseteq[\Sigma\Omega X,X],
		\]
		where $\operatorname{ev}:\Sigma\Omega X\to X$ is the evaluation map.
	\end{prop}
	
	\begin{proof} Throughout this proof, whether a morphism denotes an actual map or a
		homotopy class of a map will be clear from the context.

		Note that all maps and homotopies mentioned here are based. Write $I=[0,1],$   with basepoint $1.$ 
		Observe that $F=\{(y,\lambda)\in Y\times X^{\mathrm{I}}:
		\lambda(0)=*_X,\ \lambda(1)=f(y)\}$ with
		\[
		j:F\to Y,\quad j(y,\lambda)=y\quad  \text{and}\quad 
		\partial:\Omega X\to F,
		\partial(\omega)=(*_Y,\omega).
		\]
		The basepoint of $F$ consists of $*_Y$ and the constant
		path at $*_X.$   Recall that the  cone and mapping
		cone are
		\[
		CF=F\wedge I=\frac{F\times I}
		{F\times\{1\}\cup\{*_F\}\times I},
		\quad
		C_j=\frac{Y\vee CF}{[u,0]\sim j(u),\;(u\in F)}.
		\]
		Here,   the symbol $[u,s]$ denotes the image of $(u,s)\in F\times I$
		under the quotient map $F\times I\to CF,$   and also under the
		composite $F\times I\to CF\to C_j.$  In particular, in $C_j,$   we observe that 
		$   [u,0]=j(u),\;[u,1]=*,$   and $[*_F,s]=*.$
		Similarly,
		\[
		\Sigma\Omega X=
		\frac{\Omega X\times I}
		{\Omega X\times\{0,1\}\cup\{*_{\Omega X}\}\times I},
		\]
		and $[\omega,t]$ denotes the image of $(\omega,t)$ under the
		quotient map $\Omega X\times I\to\Sigma\Omega X.$ 
		Note that   the evaluation map $\operatorname{ev}:\Sigma\Omega X\to X$  is given by $\operatorname{ev}([\omega,t])=\omega(t).$ Since $\lambda(1)=f(y)$ and $\lambda(0)=*_X,$   the map
		\[
		H:F\times I\to X,\;\;
		((y,\lambda),s)\mapsto\lambda(1-s)
		\]
		is a homotopy from $fj$ to $*.$  Thus $H$ determines
		an extension $\operatorname{ext}_j(f):C_j\to X$ given by
		\[
		\operatorname{ext}_j(f)(y)=f(y), \;(y\in Y)\;
		\text{and}\;\operatorname{ext}_j(f)([(y,\lambda),s])=\lambda(1-s),\;\;([(y,\lambda),s]\in CF).
		\]
		The composite $j\partial:\Omega X\to Y$ is strictly constant.
		Using the constant null-homotopy of $j\partial,$   we take a coextension  $ \operatorname{coext}_j(\partial):\Sigma\Omega X\to C_j$ given by
		\[  \qquad
		\operatorname{coext}_j(\partial)([\omega,t])
		=
		\begin{cases}
			[\partial(\omega),1-2t],&0\le t\le\frac12,\\
			*,&\frac12\le t\le1.
		\end{cases}
		\]
		\noindent   The composite of the chosen extension and coextension is given by
		\[
		\operatorname{ext}_j(f)\circ\operatorname{coext}_j(\partial):
		\Sigma\Omega X\to X,\quad
		[\omega,t]\mapsto\omega(\min\{2t,1\}).
		\]  \noindent   We now verify that
		$\operatorname{ext}_j(f)\circ\operatorname{coext}_j(\partial)\simeq\operatorname{ev}.$  Consider the function
		\[
		(\Omega X\times I)\times I\to X,\quad
		((\omega,t),r)\mapsto
		\omega\!\left((1-r)\min\{2t,1\}+rt\right).
		\]
		The displayed function is continuous because evaluation is
		continuous and $(1-r)\min\{2t,1\}+rt$ is a continuous
		function of $(t,r).$   Namely, it is a map.
		Then, 
		we have a homotopy 
		$$G\colon (\Sigma\Omega X)\times I\to X,\quad ([\omega,t],r)\mapsto
		\omega\!\left((1-r)\min\{2t,1\}+rt\right)$$ with $G(-,0)=\operatorname{ext}_j(f)\circ\operatorname{coext}_j(\partial)$ and $G(-,1)=\operatorname{ev}.$ 
		Therefore, $    \operatorname{ext}_j(f)\circ\operatorname{coext}_j(\partial)
		=\operatorname{ev}$ as  homotopy classes.
		Hence the result holds.
	\end{proof}

	\section{Common  Terminology for Dealing with Toda Brackets}

	We introduce the following conventions, which are used in works of the Toda school such as
	\cite{T} and the subsequent papers \cite{MT1,Mi1,MMO,Oda1979,MiM1}.
	\begin{rem}[Conventions]
		\label{rem:toda-conventions}
		The following conventions apply whenever the ambient sets of homotopy classes $[X,Y]$ including Toda brackets as their subsets
		have the abelian group structures. Let $G$ be an abelian group.
		\begin{enumerate}
			\item The notation
			$$x\in\{\alpha, \Sigma^n\beta,  \Sigma^n\gamma\}_n\bmod A$$
			means that the bracket contains $x$ and is the coset $x+A.$ 
			Here $A$ is the subgroup which the Toda bracket is subject to. If
			$A$ is generated by elements $a_\lambda\;(\lambda\in\Lambda),$   we also write
			``$\bmod \;a_\lambda\;\ (\lambda\in\Lambda)$'' for ``$\bmod \,A$''. Note that, for a coset
			$a\in G/H$ with $H\le G,$   $H$ is conventionally called 
			the indeterminacy of $a$; we write $\operatorname{Ind}(a)=H.$ 
			If $f$ and $\Sigma g$ are composable with the bracket, then we see that
			\[\begin{aligned}
				&\operatorname{Ind}\bigl( f\circ\{\alpha, \Sigma^n\beta, \Sigma^n\gamma\}_n\circ\Sigma g \bigr) =f\circ\operatorname{Ind}\{\alpha, \Sigma^n\beta, \Sigma^n\gamma\}_n \circ\Sigma g.
			\end{aligned}\]

			\item Let $K_i\leq G,$   and let $a_i\in G/K_i.$ 
			The notation in a chain
			\[
			\begin{aligned}
				&g\in a_{j_1}\subseteq a_{j_2}\subseteq\cdots\subseteq a_{j_m} \supseteq a_{k_1}\supseteq\cdots\supseteq a_{k_r}\ni g' \qquad\bmod A,
			\end{aligned}
			\]
			indicates that $A=\operatorname{Ind}(a_{j_m}),$  
			the indeterminacy of the largest coset in the chain. Consequently
			$g\equiv g'\bmod A$\uprightsemicolon{} the indeterminacies of the smaller cosets
			need not equal $A.$  Similarly, the notation in a chain
			\[
			\begin{aligned}
				&g\in a_{j_1}\subseteq a_{j_2}\subseteq\cdots\subseteq a_{j_m}  \qquad\bmod A,
			\end{aligned}
			\] \vspace{-1.15em}
			or      \[
			\begin{aligned}
				& a_{j_m} \supseteq a_{k_1}\supseteq\cdots\supseteq a_{k_r}\ni g' \qquad\bmod A,
			\end{aligned}
			\]
			
			indicates that $A=\operatorname{Ind}(a_{j_m}),$   \footnote{Clearly, if ``$\m A$ '' is written in the form ``$\bmod \;a_\lambda,\;(\lambda\in\Lambda)$'', the equation $A=\operatorname{Ind}(a_{j_m})$ should be read as  ``$\operatorname{Ind}(a_{j_m})$  is generated by elements $a_\lambda,\; \lambda\in\Lambda$''.}
			the indeterminacy of the largest coset in the chain.    
			
			\item For $g,g'\in G$ and a subgroup $B\leq G,$  
			$g\equiv g'\bmod B$ means $g-g'\in B.$  For $U\subseteq G,$  
			$x\equiv U\bmod B$ means that $x\equiv y\bmod B$ for some
			$y\in U;$  and $0\equiv U$ means $0\in U.$ 
			
			\item A Toda bracket $\{\alpha,\s^{n}\beta,\s^{n}\gamma\}_{n}$ containing only one element is identified with
			that element. If its unique element is $a,$   we write
			$$\{\alpha,\s^{n}\beta,\s^{n}\gamma\}_{n}=a.$$ 
			
			\item In the $p$-local setting, if the domain of the Toda bracket
			$\{\alpha,\Sigma^n\beta,\Sigma^n\gamma\}_n$ is a double suspension,
			we regard the bracket as a coset in the corresponding quotient
			$\z_{(p)}$-module. The indeterminacy subgroup $A$ in~(1),
			determined by Lemma~\ref{lem:toda-indeterminacy},
			becomes a $\z_{(p)}$-submodule.
		\end{enumerate}
	\end{rem}

	\section{Composition and Suspension}\label{sec:toda-naturality}
	
	The rules below compare sets of possible values. Inclusions must
	remain inclusions unless the corresponding indeterminacies have
	also been compared. Each assertion includes the vanishing conditions needed to define its brackets.
	
	\begin{prop}[Composition Rules {\cite{T}}]
		\label{prop:toda-composition}
		Let $n\geq0,$   with the sources and targets chosen so that the
		displayed compositions exist.
		\begin{enumerate}
			\item Assume that $\alpha\circ \Sigma^n\beta=\beta\gamma=0.$  If one of $\alpha,$   $\beta,$   or $\gamma$ is zero, then
			$0\in\{\alpha,\Sigma^n\beta,\Sigma^n\gamma\}_n.$ 
			
			\item If $\alpha\circ \Sigma^n\beta=\beta\gamma=0,$   then
			\[\begin{aligned}
				&\{\alpha,\Sigma^n\beta,\Sigma^n\gamma\}_n\circ \Sigma^{n+1}\delta \subseteq\{\alpha,\Sigma^n\beta,\Sigma^n(\gamma\delta)\}_n.
			\end{aligned}\]
			
			\item If $\alpha\circ \Sigma^n\beta=\beta\gamma\delta=0,$   then
			\[\begin{aligned}
				&\{\alpha,\Sigma^n\beta,\Sigma^n(\gamma\delta)\}_n \subseteq\{\alpha,\Sigma^n(\beta\gamma),\Sigma^n\delta\}_n.
			\end{aligned}\]
			
			\item If $\alpha\circ \Sigma^n\beta\circ  \Sigma^n\gamma=\gamma\delta=0,$   then
			\[\begin{aligned}
				&\{\alpha\circ \Sigma^n\beta,\Sigma^n\gamma,\Sigma^n\delta\}_n \subseteq\{\alpha,\Sigma^n(\beta\gamma),\Sigma^n\delta\}_n.
			\end{aligned}\]
			
			\item If $\beta \circ\Sigma^n\gamma=\gamma\delta=0,$   then
			\[\begin{aligned}
				&\alpha\circ\{\beta,\Sigma^n\gamma,\Sigma^n\delta\}_n \subseteq\{\alpha\beta,\Sigma^n\gamma,\Sigma^n\delta\}_n.
			\end{aligned}\]
		\end{enumerate}
		\qed
	\end{prop}
	\begin{rem}
		For the above proposition, a useful mnemonic is that \emph{moving a factor of a composite closer to the
			middle entry enlarges the bracket.} In $\n1,\s^{n}2,\s^{n}3\nn_n,$  
		the middle entry is $2.$ 
	\end{rem}

	\begin{prop}[Adjunction and Suspension
		{\cite{T}}]
		\label{prop:toda-suspension}
		Let $\Omega_0:[\Sigma V,Z]\to[V,\Omega Z]$ be the suspension--loop
		adjunction. If $\alpha\in[\Sigma^{n+1}Y,Z]$ and
		$\alpha \Sigma^{n+1}\beta=\beta\gamma=0,$   then
		\begin{equation*}
			\begin{aligned}
				&-\Omega_0\{\alpha,\Sigma^{n+1}\beta,\Sigma^{n+1}\gamma\}_{n+1} =\{\Omega_0\alpha,\Sigma^n\beta,\Sigma^n\gamma\}_n.
			\end{aligned}
		\end{equation*}
		If $\alpha\in[\Sigma^nY,Z]$ and
		$\alpha\circ \Sigma^n\beta=\beta\gamma=0,$   then
		\begin{equation*}
			\begin{aligned}
				&-\Sigma\{\alpha,\Sigma^n\beta,\Sigma^n\gamma\}_n \subseteq\{\Sigma\alpha,\Sigma^{n+1}\beta,\Sigma^{n+1}\gamma\}_{n+1}.
			\end{aligned}
		\end{equation*}
		Both statements hold for $n\geq0.$ 
		\qed
	\end{prop}
	
	The minus sign records the interchange of the two suspension
	coordinates. Thus suspending a chosen bracket value yields its
	negative in the bracket with the next index. It does not, by itself,
	show that every value of the latter bracket desuspends.
	
	The index also keeps track of how much desuspended data are retained.
	For $n\geq m\geq0,$   whenever the bracket on the left is defined,
	\begin{equation*}
		\begin{aligned}
			&\{\alpha,\Sigma^n\beta,\Sigma^n\gamma\}_n \subseteq \{\alpha,\Sigma^m(\Sigma^{n-m}\beta),\Sigma^m(\Sigma^{n-m}\gamma)\}_m.
		\end{aligned}
	\end{equation*}
	Indeed, the original null-homotopy of $\beta\gamma$ supplies a
	suspended null-homotopy on the right, but the right-hand bracket
	permits other choices as well \cite[(1.15), p.~10]{T}.
	
	\section{The Symmetry and the Jacobi Identity}\label{sec:toda-relations}
	
	When three successive composites vanish, a bracket can be moved
	across the remaining factor with a precise sign. 
	\begin{prop}[The Symmetry
		\cite{T}]
		\label{prop:toda-four}
		Let
		$$V\xrightarrow{\delta}W\xrightarrow{\gamma}X
		\xrightarrow{\beta}Y$$ be given, and let $\alpha\in[\Sigma^nY,Z]$ satisfy
		$\alpha\circ \Sigma^n\beta=\beta\gamma=\gamma\delta=0.$  Then
		\begin{equation*}
			\begin{aligned}
				&\{\alpha,\Sigma^n\beta,\Sigma^n\gamma\}_n\circ \Sigma^{n+1}\delta =(-1)^{n+1}\alpha\circ \Sigma^n\{\beta,\gamma,\delta\}.
			\end{aligned}
		\end{equation*}
		\qed
	\end{prop}
	
	Here both sides are the full sets of composites obtained from the
	indicated brackets. For $n=0,$   this reads
	$\{\alpha,\beta,\gamma\}\circ \Sigma\delta
	=-\alpha\circ\{\beta,\gamma,\delta\}.$ 
	It is useful when a composite of a bracket is easier to calculate
	than the bracket itself\uprightsemicolon{} the equality still concerns all allowed
	values, rather than two independently chosen representatives.

	\begin{exercise}
		
		Let the sequence
		\[
		X_0
		\xleftarrow{\alpha_1}\Sigma^n X_1
		\xleftarrow{\Sigma^n\alpha_2}\Sigma^n X_2
		\xleftarrow{\Sigma^n\alpha_3}\Sigma^n X_3
		\xleftarrow{\Sigma^n\alpha_4}\Sigma^n X_4
		\]

		\noindent   of spaces and homotopy classes of maps
		be given, where $\af_{1}\s^{n}\af_{2}=\af_{2}\af_{3}=\af_{3}\af_{4}=0$ and  \begin{equation}\notag
			0\in A:=\n \af_{1},\s^{n}\af_{2},\s^{n}\af_{3}\nn_{n},\qquad    0\in B:=\n \af_{2},\af_{3},\af_{4}\nn.
		\end{equation}
		\noindent Suppose that \begin{equation}\notag
			\mathrm{Ind}(A)=0\;\;\mathrm{or}\;\;\mathrm{Ind}(B)=0.
		\end{equation}
		Show that  there are two well-defined Toda brackets $$\{\overline{\alpha_1},\Sigma^n\widetilde{\alpha_3},
		\Sigma^n(-\Sigma\alpha_4)\}_n\;\text{and}\; \{\alpha_1,\Sigma^n\overline{\alpha_2},
		\Sigma^n\widetilde{\alpha_4}\}_n$$
		\noindent  for \vspace{-0.4em}$$\text{some}\; \overline{\alpha_1}  \in\operatorname{Ext}_{\Sigma^n\alpha_2}(\alpha_1)
		\subseteq[C_{\Sigma^n\alpha_2},X_0],\quad  \text{some  }    \widetilde{\alpha_3}\in   \operatorname{Coext}_{\alpha_2}(\alpha_3)
		\subseteq[\Sigma X_3,C_{\alpha_2}],$$
		\vspace{-0.4em}
		$$  \text{some}\;   \overline{\alpha_2}\in \operatorname{Ext}_{\alpha_3}(\alpha_2)
		\subseteq[C_{\alpha_3},X_1], \quad    \text{and some  }\widetilde{\alpha_4}\in \operatorname{Coext}_{\alpha_3}(\alpha_4)
		\subseteq[\Sigma X_4,C_{\alpha_3}] .$$
		
		\noindent   Hint: use the above proposition and Proposition \ref{dsjd}.
	\end{exercise}

	The next relation concerns \emph{compatible} choices in three
	brackets. Its existential formulation prevents a common mistake:
	one cannot choose all three inner bracket values independently.
	
	\needspace{3\baselineskip}
	\begin{prop}[The Jacobi Identity
		{\cite{T}}]
		\label{prop:toda-five}  Consider the sequence of  spaces and homotopy classes
		\[
		\begin{aligned}
			&X\xleftarrow{\alpha_1}Y
			\xleftarrow{\alpha_2}Z
			\xleftarrow{\alpha_3}K
			\xleftarrow{\alpha_4}L
			\xleftarrow{\alpha_5}U.
		\end{aligned}
		\]
		Suppose that
		\[
		\begin{aligned}
			&\alpha_1\alpha_2
			=\alpha_2\alpha_3
			=\alpha_3\alpha_4
			=\alpha_4\alpha_5=0,
		\end{aligned}
		\]
		and that
		\[
		\begin{aligned}
			&0\in\alpha_1\circ\{\alpha_2,\alpha_3,\alpha_4\},
			\qquad
			0\in\{\alpha_2,\alpha_3,\alpha_4\}\circ\Sigma\alpha_5.
		\end{aligned}
		\]
		Then there exist choices
		\[
		\begin{aligned}
			&\xi_1\in\{\alpha_1,\alpha_2,\alpha_3\},\quad
			\xi_2\in\{\alpha_2,\alpha_3,\alpha_4\},\quad
			\xi_3\in\{\alpha_3,\alpha_4,\alpha_5\},
		\end{aligned}
		\]
		for which the following brackets are well-defined and satisfy
		\begin{equation}\label{eq:generalized-jacobi-n0}
			\begin{aligned}
				&0\in
				\{\xi_1,\Sigma\alpha_4,\Sigma\alpha_5\}
				+\{\alpha_1,\xi_2,\Sigma\alpha_5\}
				+\{\alpha_1,\alpha_2,\xi_3\}.
			\end{aligned}
		\end{equation}
		\noindent
		This relation may be remembered symbolically as
		\[
		\begin{gathered}
			\{\{\alpha_1,\alpha_2,\alpha_3\},
			\Sigma\alpha_4,\Sigma\alpha_5\}
			\\[-5pt]
			+\\[-5pt]
			\{\alpha_1,\{\alpha_2,\alpha_3,\alpha_4\},
			\Sigma\alpha_5\}
			\\[-5pt]
			+\\[-5pt]
			\{\alpha_1,\alpha_2,
			\{\alpha_3,\alpha_4,\alpha_5\}\}\\[-3pt]
			\rotatebox[origin=c]{90}{$\equiv$}\\[-3pt]
			0.
		\end{gathered}
		\]
		
		\qed
	\end{prop}

	In addition to the definitions of elements via Toda brackets,
	the Jacobi Identity  provides a valuable tool for detecting elements
	in a given Toda bracket. We provide an example as an exercise.

	\begin{exercise} Work 2-locally. Follow     Appendix~\ref{app:oguchi} and Appendix~\ref{tlqdzj}. Show $$\n\ca_{16}, 2\e_{17},16\rho_{17} \nn\ni\cq_{16}$$ in $\pi_{33}(S^{16}).$ Hint: use the Jacobi Identity  and    $16\co_{17}=\n 8\cg_{17},2\e_{24},8\cg_{24}\nn$ which is implied by   $16\co_{17}=\s^{12}\co^{\mathrm{IV}}$; the element $\cq_{3}$ is defined to be an element of $$\n \mu_{3},2\e_{12},8\cg_{12}\nn_{1}$$ \noindent and so $\cq_{16}\in \n \mu_{16},2\e_{25},8\cg_{25}\nn.$ Additionally, in view of the definition of $\mu_{n}$, we know that $\mu_{16}\in \n\ca_{16},2\e_{17},8\cg_{17}\nn.$
	\end{exercise}

	In 2023, Professor Mukai raised the question of whether the above result admits a generalization with
	any subscript \(n\geq0\), suitable for computing \(\pi_{33+i}(S^i)\). Juxin Yang subsequently
	obtained the desired generalization by a short argument.
	
	\begin{prop}[The Generalized Jacobi Identity]
		\label{prop:generalized-jacobi}
		Let $n\geq0$ and consider the sequence of  spaces and
		homotopy classes
		\[\begin{aligned}
			&X\xleftarrow{\alpha_1}\Sigma^nY \xleftarrow{\Sigma^n\alpha_2}\Sigma^nZ \xleftarrow{\Sigma^n\alpha_3}\Sigma^nK \xleftarrow{\Sigma^n\alpha_4}\Sigma^nL \xleftarrow{\Sigma^n\alpha_5}\Sigma^nU.
		\end{aligned}\]
		Suppose that
		\[\begin{aligned}
			&\alpha_1\circ\Sigma^n\alpha_2 =\alpha_2\alpha_3 =\alpha_3\alpha_4 =\alpha_4\alpha_5=0,
		\end{aligned}\]
		and that
		\[\begin{aligned}
			&0\in\alpha_1\circ\Sigma^n\{\alpha_2,\alpha_3,\alpha_4\}, \qquad 0\in\{\alpha_2,\alpha_3,\alpha_4\}\circ\Sigma\alpha_5.
		\end{aligned}\]
		Then there exist choices
		\[
		\begin{aligned}
			&\xi_1\in \{\alpha_1,\Sigma^n\alpha_2,\Sigma^n\alpha_3\}_n,\; \xi_2\in\{\alpha_2,\alpha_3,\alpha_4\},\; \xi_3\in\{\alpha_3,\alpha_4,\alpha_5\},
		\end{aligned}
		\]
		for which the following brackets are well-defined and satisfy
		\begin{equation}\label{eq:generalized-jacobi}
			\begin{aligned}
				&0\in{} (-1)^n \{\xi_1,\Sigma^{n+1}\alpha_4,\Sigma^{n+1}\alpha_5\}_n+\{\alpha_1,\Sigma^n\xi_2,\Sigma^{n+1}\alpha_5\}_n+\{\alpha_1,\Sigma^n\alpha_2,\Sigma^n\xi_3\}_n.
			\end{aligned}
		\end{equation}
		\noindent
		Following Toda, we may remember this relation symbolically as% Preamble: \usepackage{amsmath,graphicx}
		\[
		\begin{gathered}
			(-1)^n
			\{\{\alpha_1,\Sigma^n\alpha_2,\Sigma^n\alpha_3\}_n,
			\Sigma^{n+1}\alpha_4,\Sigma^{n+1}\alpha_5\}_n
			\\[-5pt]
			+\\[-5pt]
			\{\alpha_1,\Sigma^n\{\alpha_2,\alpha_3,\alpha_4\},
			\Sigma^{n+1}\alpha_5\}_n
			\\[-5pt]
			+\\[-5pt]
			\{\alpha_1,\Sigma^n\alpha_2,
			\Sigma^n\{\alpha_3,\alpha_4,\alpha_5\}\}_n\\[-3pt]
			\rotatebox[origin=c]{90}{$\equiv$}\\[-3pt]
			\;0.
		\end{gathered}
		\]\end{prop}
	
	\begin{proof}
		The result follows from the suspension--loop adjunction formula in  Proposition~\ref{prop:toda-suspension}
		and Toda's classical Jacobi
		identity.
		
		Write \vspace{-0.8em}
		\[
		\begin{aligned}
			&\operatorname{ad}^{n}: [\Sigma^nM,N]\longrightarrow[M,\Omega^nN],\qquad \operatorname{ad}^{-n}: [M,\Omega^nN]\longrightarrow[\Sigma^nM,N],
		\end{aligned}
		\]
		for the mutually inverse suspension--loop adjunction bijections.
		Applying adjunction to $\alpha_1$ gives the sequence
		\[\begin{aligned}
			&\Omega^nX \xleftarrow{\operatorname{ad}^{n}(\alpha_1)}Y \xleftarrow{\alpha_2}Z \xleftarrow{\alpha_3}K \xleftarrow{\alpha_4}L \xleftarrow{\alpha_5}U.
		\end{aligned}\]
		Toda's classical Jacobi identity supplies compatible choices
		in the inner brackets such that
		\[
		\begin{aligned}
			0&\in{} \bigl\{ \{\operatorname{ad}^{n}(\alpha_1),\alpha_2,\alpha_3\}, \Sigma\alpha_4,\Sigma\alpha_5 \bigr\} \\
			& {}+\bigl\{ \operatorname{ad}^{n}(\alpha_1), \{\alpha_2,\alpha_3,\alpha_4\}, \Sigma\alpha_5 \bigr\} +\bigl\{ \operatorname{ad}^{n}(\alpha_1),\alpha_2, \{\alpha_3,\alpha_4,\alpha_5\} \bigr\}.
		\end{aligned}
		\]
		Apply $\operatorname{ad}^{-n}$ to this relation.
		Repeated use of the adjunction formula in
		Proposition~\ref{prop:toda-suspension}  introduces the factor $(-1)^n$
		in the first summand and gives
		\eqref{eq:generalized-jacobi}.
	\end{proof}
	
	\section{The Additivity}\label{sec:toda-additivity}
	
	In an unstable calculation, distributivity requires attention to
	the source. Postcomposition induces a homomorphism on homotopy
	classes from a suspension\uprightsemicolon{} precomposition has that property when
	the precomposing map is itself a suspension. The hypotheses below
	retain exactly the suspension conditions needed for each slot.
	
	\begin{prop}[The Additivity {\cite{T}}]
		\label{prop:toda-additivity}
		Let $\alpha,\alpha'\in[\Sigma^nY,Z],$  
		$\beta,\beta'\in[X,Y],$   and $\gamma,\gamma'\in[W,X].$ 
		In each assertion, assume that the two brackets on the left are
		defined and that the displayed sums exist. In particular, an
		unsuspended sum is used in its suspension group structure\uprightsemicolon{} this
		condition is automatic for the positive-dimensional sphere sources
		in the later calculations.
		\begin{enumerate}
			\item If $n\geq1$ or $W=\Sigma W',$   then
			\[
			\begin{aligned}
				&\{\alpha,\Sigma^n\beta,\Sigma^n\gamma\}_n +\{\alpha,\Sigma^n\beta,\Sigma^n\gamma'\}_n\supseteq \{\alpha,\Sigma^n\beta,\Sigma^n(\gamma+\gamma')\}_n.
			\end{aligned}
			\]
			
			\item If $\gamma=\Sigma\gamma'$ with
			$X=\Sigma X'$ and $W=\Sigma W',$   then
			\[
			\begin{aligned}
				&\{\alpha,\Sigma^n\beta,\Sigma^n\gamma\}_n +\{\alpha,\Sigma^n\beta',\Sigma^n\gamma\}_n= \{\alpha,\Sigma^n(\beta+\beta'),\Sigma^n\gamma\}_n.
			\end{aligned}
			\]
			
			\item If $n\geq1,$   or if $\beta=\Sigma\beta'$ and
			$\gamma=\Sigma\gamma'$ with $Y=\Sigma Y',$   $X=\Sigma X'$ and $W=\Sigma W',$   then
			\[
			\begin{aligned}
				&\{\alpha,\Sigma^n\beta,\Sigma^n\gamma\}_n +\{\alpha',\Sigma^n\beta,\Sigma^n\gamma\}_n\supseteq \{\alpha+\alpha',\Sigma^n\beta,\Sigma^n\gamma\}_n.
			\end{aligned}
			\]
		\end{enumerate}
		\qed
	\end{prop}
	
	The middle slot gives an equality, while the outer slots generally
	give inclusions. For instance, the indeterminacy associated with
	$\alpha+\alpha'$ is contained in the sum of the two
	indeterminacies associated with $\alpha$ and $\alpha'.$ 
	This explains why an outer-slot sum can contain additional values.
	The formulas are used only with defined sums. They must not be
	read as introducing an addition on an arbitrary unsuspended
	homotopy set.
	
	\section{Hopf Invariants and Stable Brackets}\label{sec:toda-hopf-stable}
	
	The preceding rules describe how brackets behave under composition. We now use the Hopf invariant and
	stabilization to compare bracket calculations in different dimensions. These comparisons are useful
	precisely because information may be lost under suspension\uprightsemicolon{} the kernels and indeterminacies must
	remain visible.
	
	\subsection{The Hopf Invariant of a Toda  Bracket}
	
	The EHP sequence, reviewed in Section~\ref{sec:ehp-review}, distinguishes
	a class killed by suspension from a class that is already zero.
	This distinction is essential in the next result: only
	$\Sigma(\alpha\beta)=0$ is required.
	
	Work 2-locally. For $m\geq1,$   use the second James--Hopf homomorphism
	$H:\pi_q(S^{m+1})\to\pi_q(S^{2m+1})$ and the EHP homomorphism
	$P:\pi_{q+2}(S^{2m+1})\to\pi_q(S^m),$   with Toda's compatible
	sign convention. Toda denotes this homomorphism between homotopy groups of spheres by $\Delta$\uprightsemicolon{}
	we write $P$ to distinguish it from the connecting homomorphism for the fibration of special orthogonal groups
	in this book. For $u\in\ker \Sigma\subseteq\pi_k(S^m),$   let
	\begin{equation*}
		\begin{aligned}
			&P^{-1}(u)=\{v\in\pi_{k+2}(S^{2m+1}):P(v)=u\}.
		\end{aligned}
	\end{equation*}
	This notation denotes a whole inverse-image coset, not an inverse
	homomorphism or a preferred lift. Exactness gives its indeterminacy
	$H\pi_{k+2}(S^{m+1}).$ 
	
	\begin{prop}[The $H$-Formula
		{\cite{T}}]
		\label{prop:toda-hopf}
		Work 2-locally.     Let $\alpha\in[K,S^m],$   $\beta\in\pi_k(K)$ and
		$\gamma\in\pi_l(S^k),$   with
		$\Sigma(\alpha\beta)=\beta\gamma=0.$  Then, in the  EHP sequence,
		\begin{equation*}
			\begin{aligned}
				&H\{\Sigma\alpha,\Sigma\beta,\Sigma\gamma\}_1 =-P^{-1}(\alpha\beta)\circ \Sigma^2\gamma \subseteq\pi_{l+2}(S^{2m+1}).
			\end{aligned}
		\end{equation*}
		\qed
	\end{prop}
	
	The dimension check is useful: $\alpha\beta\in\pi_k(S^m),$  
	its $P$-lifts lie in $\pi_{k+2}(S^{2m+1}),$   and
	$\Sigma^2\gamma:S^{l+2}\to S^{k+2}.$  Thus the right side has exactly
	the degree of the Hopf image on the left. In practice, one first
	finds the $P$-lifts of $\alpha\beta$ and then composes them with
	$\Sigma^2\gamma.$  The remaining ambiguity is the displayed image of $H$\uprightsemicolon{}
	choosing one $P$-lift alone does not remove it.
	\begin{exercise} Work 2-locally. Follow     Appendix~\ref{app:oguchi} and Appendix~\ref{tlqdzj}. Assume that  $m\in \z_{+}$ is an odd number. Show the relation for the  Whitehead product: $$\pm[\e_{m+1},\e_{m+1}]\ty \n \s\e_{m}, \s[\e_{m},\e_{m}], \s(2\e_{2m-1})\nn_{1}\m \s\pi_{2m}(S^{m}).$$
		\noindent Hint: use the well-known fact for the Hopf invariant  $H([\e_{m+1},\e_{m+1}])=\pm2\e_{2m+1}$ and use the above $H$-formula;  consider the $\mathrm{EHP}$ sequence $\pi_{2m}(S^{m})\stackrel{\s}\longrightarrow\pi_{2m+1}(S^{m+1})\stackrel{H}\longrightarrow\pi_{2m+1}(S^{2m+1}) .$ 
	\end{exercise}
	
	\begin{exercise} \label{1028lx}Work 2-locally. Note that $\e_{m}$ is the identity on $S^{m}.$ Assume that $P(\s^{2}\mathbf{a})=0$ in $\pi_{k}(S^{m})$ where $\mathbf{a}\in\pi_{k}(S^{2m-1})$ with $k,\,m\in\z_{+}$.  
		\begin{itemize}
			\item[\rm (1)]  
			Show that\[
			H\n\s\e_{m},\s[\e_{m},\e_{m}],\s\mathbf{a}\nn_1
			\begin{cases}
				=\pm\n\s^{2}\mathbf{a}\nn, & \text{if } m\ty0\m2,\\[4pt]
				\ni\pm\s^{2}\mathbf{a}, & \text{if } m\ty1\m2.
			\end{cases}
			\]  \noindent Here, $ \n\s^{2}\mathbf{a}\nn$ denotes a singleton set.\footnote{Warning: in Toda's notation, $\n x\nn$ denotes
				the subgroup generated by $x$. In this book, we do not follow this convention.}
			
			\item[\rm (2)]  Show that for any $\mathrm{Lif}(\s^{2}\mathbf{a})\in \pi_{k+2}(S^{m+1})$ with $H(\mathrm{Lif}(\s^{2}\mathbf{a}))=\s^{2}\mathbf{a}$, there is a relation $$\pm\,\mathrm{Lif}(\s^{2}\mathbf{a})\ty \n\s\e_{m},\s[\e_{m},\e_{m}],\s\mathbf{a}\nn_1\m\s\pi_{k+1}(S^{m}).$$ 
			
		\end{itemize}
		\noindent Hint: use the well-known $P$-formulas given in \eqref{1028ppp} and employ the  $H$-formula  given in the above proposition.
	\end{exercise}

	\begin{rem}Work 2-locally.  
		Exercise~\ref{1028lx} and the relation
		$2\ca_{2}=\pm[\e_{2},\e_{2}]$
		offer a useful clue as to how a generator
		of $\pi_{6}(S^{3})$ might be constructed.
		They naturally suggest looking for an element
		$\nu'$ satisfying
		\[
		\nu'\in
		\n \ca_{3},2\e_{4},\ca_{4}\nn_1
		\xyd
		\n \s\e_{2},\s(\pm[\e_{2},\e_{2}]),\s\ca_3\nn_1.
		\]
		The element $\nu'$ so chosen lifts
		$\ca_{5}\in H\pi_{6}(S^{3})$
		and generates
		$\pi_{6}(S^{3})\cong\z/4$. See \cite{T}.  In fact, many new elements in the homotopy groups of spheres can be constructed in this way. \end{rem}
	
	\subsection{Brackets of the Form $\n1,2,1\nn$}

	\begin{prop}\label{121X}Work 2-locally. Let $k\geq3$. \begin{itemize}
			\item[\rm (1)]
			{\cite[Corollary 3.7]{T}}  Assume that $\mathbf x\in\pi_{n-1}(S^{k-1})$ $(n>1)$  satisfies $2\mathbf x=0$. Then, $$\n 2\e_k, \s\mathbf x, 2\e_n\nn_1\ni (\s\mathbf x)\ca_{n}.$$
			\item[\rm (2)]
			{\cite[Lemma 5.5]{T}} Assume that $\mathbf y\in\pi_{n-3}(S^{k-2})$ $(n>\max\n5,k+1\nn)$ satisfies $\mathbf y\hc \ca_{n-3}=0$. Then, $$\n \ca_k, \s^3\mathbf y, \ca_n \nn_3\ni(\s^2\mathbf y)\nu_{n-1}.$$\qed
		\end{itemize}
		
	\end{prop}

	Recall that in this chapter we work integrally unless otherwise stated.
	Write
	$\pi_d^S(S^0)=\operatorname*{colim}_r\pi_{r+d}(S^r)$
	for the stable homotopy group of spheres. Stable composition is
	bilinear and graded commutative:
	$\alpha\beta=(-1)^{kh}\beta\alpha$ for
	$\alpha\in\pi_k^S(S^0)$ and $\beta\in\pi_h^S(S^0).$ 
	
	Let $\alpha\in\pi_k^S(S^0),$   $\beta\in\pi_h^S(S^0)$ and
	$\gamma\in\pi_l^S(S^0)$ satisfy $\alpha\beta=\beta\gamma=0.$ 
	For sufficiently large $r,$   choose 
	$\alpha_r:S^{r+k}\to S^r,$  
	$\beta_{r+k}:S^{r+k+h}\to S^{r+k}$ and
	$\gamma_{r+k+h}:S^{r+k+h+l}\to S^{r+k+h}$
	such that $$\s^{\wq}\af_{r}=\af, \s^{\wq}\beta_{r+k}=\beta,\s^{\wq}\gamma_{r+k+h}=\gamma$$ \noindent and that
	$\alpha_r\beta_{r+k}= 0$ and
	$\beta_{r+k}\gamma_{r+k+h}=0.$ 
	Define the stable Toda bracket by
	\begin{equation*}
		\begin{aligned}
			&\langle\alpha,\beta,\gamma\rangle
			=\Sigma^\infty\bigl((-1)^{r-1}
			\{\alpha_r,\beta_{r+k},\gamma_{r+k+h}\}\bigr).
		\end{aligned}
	\end{equation*}
	The sign $(-1)^{r-1}$ ensures, by the suspension formula,
	that this definition is independent of a further increase
	of $r$ \cite[p.~32]{T}.
	The resulting bracket is a coset in $\pi_{k+h+l+1}^S(S^0)$
	with indeterminacy
	\[\begin{aligned}
		&\alpha\circ\pi_{h+l+1}^S(S^0)
		+\pi_{k+h+1}^S(S^0)\circ\gamma.
	\end{aligned}\]
	
	\noindent An important property of $\x \alpha,\beta,\gamma\xx$ defined above is graded symmetry, which tells us that $$\x \alpha,\beta,\gamma\xx=\pm\x \gamma,\beta,\af\xx. $$
	\noindent Another important property of $\x \alpha,\beta,\gamma\xx$ defined above is the Jacobi identity, which gives $$\x \alpha,\beta,\gamma\xx\pm
	\x \beta,\gamma, \alpha\xx\pm
	\x \gamma, \alpha,\beta\xx \ty0$$
	\noindent
	if we additionally have $\af\gamma=0.$  (See \cite{T}.)
	\needspace{3\baselineskip}
	\begin{prop}[Repeated Outer Class {\cite{T}}]
		\label{prop:toda-stable-repeat}
		Let $\alpha\in\pi_k^S(S^0)$ and $\beta\in\pi_h^S(S^0),$  
		with $\alpha\beta=0.$ 
		\begin{enumerate}
			\item If $k$ is odd, then
			\[\begin{aligned}
				&\langle\alpha,\beta,\alpha\rangle
				\cap(-1)^h\langle\beta,\alpha,2\alpha\rangle
				\neq\varnothing.
			\end{aligned}\]
			\item If $k$ is even, then
			\[\begin{aligned}
				&\langle\alpha,\beta,\alpha\rangle
				\cap\bigl(\beta\circ\pi_{2k+1}^S(S^0)\bigr)
				\neq\varnothing.
			\end{aligned}\]
		\end{enumerate}
		\qed
	\end{prop}
	
	\begin{exercise} Work 2-locally. Follow     Appendix~\ref{app:oguchi} and Appendix~\ref{tlqdzj}. Show that $$\x2\e,8\cg,2\e\xx=0$$ in $\pi_{8}^{S}(S^{0}).$ Hint: use the above proposition and $\pi_{8}^{S}(S^{0})\cong(\z/2)^{2}.$ 
	\end{exercise}
	
	\section{Matrix Toda Brackets}\label{sec:mtb}

	Suppose $f_{i}: T_{i}\rightarrow W\; (i=1,2)$ are maps or homotopy classes. By abuse of notation, the following composition is still denoted by $f_{1}\vee f_{2},$  
	\[
	\begin{tikzcd}[column sep=3em]
		W &
		W\vee W
		\arrow[l, "\nabla" {above}, "\text{folding}" {below}] &
		T_1\vee T_2
		\arrow[l, "f_1\vee f_2"]
	\end{tikzcd}
	\]
	\noindent   Note that $$f_{1}\vee f_{2}=f_{1}\hc (T_{1}\xleftarrow{proj}T_{1}\vee T_{2})+f_{2}\hc (T_{2}\xleftarrow{proj}T_{1}\vee T_{2})$$ if $T_{1}$ and $T_{2}$ are suspensions. 
	We use the identification $\s^{n}(V_{1}\vee V_{2})=\s^{n}V_{1}\vee \s^{n}V_{2},\;\,  (V_{1}, V_{2}\colon$ spaces). For  $k$-fold wedges, we make comparable assumptions.

	Given $\ell,r\in\z_{+},$   let $W, X_{k}, Y_{s}$ and $U$ be spaces ($k,s\in\z,$    $1\leq k\leq \ell,\;1\leq s\leq r$), where $U=\s^{2}U''$ is a double suspension.  For each $k$ and $s,$   let the  homotopy classes
	$a_{k},$   $b_{k,s}$ and $c_{s}$ be given with respect to the following sequences
	$$\xymatrix@C=1cm{
		W\ar@{<-}[r]^{a_{k}}&\s^{n}X_{k} \ar@{<-}[r]^{\;\,\s^{n}b_{k,s}} & \s^{n}Y_{s}\ar@{<-}[r]^{\s^{n} c_{s}}& \s^{n}U.} $$
	\noindent Moreover, suppose the homotopy classes $c_{s}$ are \text{suspensions}\footnote{Namely, $c_{s}=\s c_{s}'$ for some homotopy class $c_{s}'.$  } for all $1\leq s\leq r.$  Let $\bigvee_{p=1}^{\ell} X_{p}\stackrel{\mathrm{i}_{k}}\longleftarrow X_{k}$ and $\bigvee_{p=1}^{r} Y_{p}  \stackrel{\;\mathrm{j}_{s}}\longleftarrow Y_{s}$ be the inclusions.
	We introduce formal matrices to denote homotopy classes,\vspace{-0.45\baselineskip}
	$$\mathbbm{a}_{1\times \ell}=( 
	a_{1}, a_{2}, \cdots, a_{\ell} ):=a_{1}\vee a_{2}\vee\cdots\vee a_{\ell}\in [\s^{n}(\bigvee_{k=1}^{\ell} X_{k}),\; W ],$$\vspace{-0.61\baselineskip}
	$$\mathbbm{b}_{\bullet,s}=(b_{1,s},\; b_{2,s},\; \cdots, \;b_{\ell,s})^{T}:=\mathrm{i}_{1}\hc b_{1,s}+\mathrm{i}_{2}\hc b_{2,s}+\cdots+\mathrm{i}_{\ell}\hc b_{\ell,s}\in [Y_{s}, \bigvee_{k=1}^{\ell} X_{k}],$$
	\vspace{-0.31\baselineskip}
	$$\mathbbm{b}_{\ell\times r}=(b_{k,s}):=\mathbbm{b}_{\bullet,1}\vee \mathbbm{b}_{\bullet,2}\vee\cdots \vee \mathbbm{b}_{\bullet,r}\in [ 
	\bigvee_{s=1}^{r} Y_{s}   ,\bigvee_{k=1}^{\ell} X_{k}],$$\vspace{-0.31\baselineskip}
	$$\s^{n}\mathbbm{b}_{\ell\times r}:=(\s^{n}b_{k,s})=\s^{n}(\mathbbm{b}_{\bullet,1}\vee \mathbbm{b}_{\bullet,2}\vee\cdots \vee \mathbbm{b}_{\bullet,r})\in [\s^{n}( 
	\bigvee_{s=1}^{r} Y_{s}  )  ,\;\s^{n}(\bigvee_{k=1}^{\ell} X_{k})]$$\vspace{-0.31\baselineskip}
	$$\text{and}\qquad\mathbbm{c}_{r\times 1}=( c_{1}, c_{2}, \cdots, c_{r})^{T}:=\mathrm{j}_{1}\hc c_{1}+\mathrm{j}_{2}\hc c_{2}+\cdots+\mathrm{j}_{r}\hc c_{r}\in [U, \bigvee_{s=1}^{r} Y_{s}].\;\qquad\qquad\;$$
	\noindent A homotopy class is regarded as a $1 \times 1$ formal matrix.   We say that two formal matrices are equal if the corresponding entries  are equal.
	
	Through direct verification, it is evident that the composition of $\mathbbm{a}_{1\times \ell}$ and $\s^{n}\mathbbm{b}_{\ell\times r}$  as homotopy classes  corresponds to \emph{their multiplication as formal matrices}; 
	the composition of   $\mathbbm{b}_{\ell\times r}$ and $\mathbbm{c}_{r\times 1}$ as homotopy classes  corresponds to  \emph{their multiplication as formal matrices} because $c_{s}$ are suspensions for all $1\leq s\leq r.$  (Notice that $f$ being a suspension is a sufficient condition for $f^{*}(-)=(-)\hc f $ being  a homomorphism.) The formal matrices of trivial homotopy classes are denoted by $O.$  For any formal matrix $\mathbbm{x}_{\ell\times r},$   the notation $\mathbbm{x}_{\bullet,s}$ always denotes its $s$-$th$ column  vector. 
	
	\begin{rem}If $U$ is not assumed to be a double suspension, the condition that all $c_s$ are suspensions \textbf{DOES NOT}, in general, ensure that composition agrees with formal matrix multiplication.  In \cite{YMM}, the space $U$ should be assumed to be a double suspension, or to satisfy the weaker condition that the functor $[U,-]$ takes values in abelian groups. The omission of this hypothesis was an oversight on Juxin Yang's part. See the following example.
	\end{rem}

	\begin{exam} If we do not assume that the functor $[U,-]$ takes values in abelian groups,   the following example shows that this can fail even when the spaces
		before suspension are simply connected.
		
		Take $n=0,$   $\ell=r=2,$   and $W=*.$  Let
		$U'=Y_1'=Y_2'=S^2\times S^2,$  
		$U=Y_1=Y_2=\s(S^2\times S^2),$   and $X_1=X_2=S^3.$ 
		Let $$b_{1,1}: X_1=S^3\leftarrow U=\s(S^2\times S^2) \text{   and    } b_{2,1}:X_2=S^3\leftarrow U=\s(S^2\times S^2)$$ \noindent be the suspensions
		of the projections of $S^2\times S^2$ onto the first and second
		factors, respectively.  Take the zero classes
		$*\xleftarrow{a_1=a_2=0}S^3,$  
		the   elements  
		$$S^3\xleftarrow{b_{1,2}=-b_{1,1}}\s(S^2\times S^2) \text{   and    }S^3\xleftarrow{b_{2,2}=-b_{2,1}}\s(S^2\times S^2),$$
		and the suspended identity 
		$\s(S^2\times S^2)
		\xleftarrow{c_1=c_2=\s\mathrm{id}_{S^2\times S^2}}
		\s(S^2\times S^2).$ 
		Thus the formal matrices are
		\begin{align}
			\mathbbm{a}_{1\times2}&=(0,0),\qquad
			\mathbbm{b}_{2\times2}=
			\begin{pmatrix}
				b_{1,1}&-b_{1,1}\\
				b_{2,1}&-b_{2,1}
			\end{pmatrix},\qquad
			\mathbbm{c}_{2\times1}=
			\begin{pmatrix}\mathrm{id}_U\\\mathrm{id}_U\end{pmatrix}.
			\notag
		\end{align}
		Both $c_1$ and $c_2$ are suspensions of identity maps of the simply
		connected space $U'.$  Also,  we have $X_1=X_2=\s S^2,$   so every space
		before suspension in this construction is simply connected.
		
		As formal matrices, these satisfy
		$\mathbbm{a}_{1\times2}\hc\mathbbm{b}_{2\times2}=O$ and
		\begin{align}
			\mathbbm{b}_{2\times2}\hc\mathbbm{c}_{2\times1}
			&=\begin{pmatrix}b_{1,1}+b_{1,2}\\b_{2,1}+b_{2,2}\end{pmatrix}
			\notag=\begin{pmatrix}b_{1,1}-b_{1,1}\\b_{2,1}-b_{2,1}\end{pmatrix}
			=O.
			\notag
		\end{align}
		These cancellations use only the group axioms, not commutativity.
		
		We next calculate the actual composition using the realizations
		specified in this section. Recall that 
		$\mathrm{i}_k:X_k\to X_1\vee X_2$ and
		$\mathrm{j}_s:Y_s\to Y_1\vee Y_2$ are the inclusions. The two columns represent
		\begin{align}
			\mathbbm{b}_{\bullet,1}
			&=\mathrm{i}_1\hc b_{1,1}+\mathrm{i}_2\hc b_{2,1},
			\notag\\
			\mathbbm{b}_{\bullet,2}
			&=-\mathrm{i}_1\hc b_{1,1}-\mathrm{i}_2\hc b_{2,1},
			\notag
		\end{align}
		and $\mathbbm{c}_{2\times1}=\mathrm{j}_1+\mathrm{j}_2.$ 
		Thus, as homotopy classes, we obtain
		\begin{align}
			\mathbbm{b}_{2\times2}\hc\mathbbm{c}_{2\times1}
			&=(\mathbbm{b}_{\bullet,1}\vee\mathbbm{b}_{\bullet,2})
			\hc(\mathrm{j}_1+\mathrm{j}_2)
			\notag\\
			&=\mathbbm{b}_{\bullet,1}+\mathbbm{b}_{\bullet,2}
			\notag\\
			&=\mathrm{i}_1\hc b_{1,1}+\mathrm{i}_2\hc b_{2,1}
			-\mathrm{i}_1\hc b_{1,1}-\mathrm{i}_2\hc b_{2,1}\neq0.
			\notag
		\end{align}
		This expression is the group commutator
		of $\mathrm{i}_1\hc b_{1,1}$ and $\mathrm{i}_2\hc b_{2,1}.$ 
		In particular, termwise negation of the first column does not
		automatically give the inverse of its realization: the inverse
		of a sum reverses the order of its two terms. Note that it can be shown that $$\mathrm{i}_1\hc b_{1,1}+\mathrm{i}_2\hc b_{2,1}
		-\mathrm{i}_1\hc b_{1,1}-\mathrm{i}_2\hc b_{2,1}\neq0$$ \noindent by using  Samelson products. \qed
	\end{exam}

	Then, we are poised to present the ensuing definition.
	
	\begin{defn}(\cite{YMM})\label{lwndy}
		Using the above symbols, further  suppose that $\mathbbm{a}_{1\times \ell}\hc\s^{n}\mathbbm{b}_{\ell\times r}=O  $
		and  that $\mathbbm{b}_{\ell\times r}\hc\mathbbm{c}_{r\times 1}=O.$  The matrix Toda bracket  indexed by $n,$   $\n \mathbbm{a}_{1\times \ell},\; \s^{n}\mathbbm{b}_{\ell\times r}, \;\s^{n}\mathbbm{c}_{r\times 1}\nn_{n},$   is defined to be the \textit{3}-fold Toda bracket with respect to the following sequence\vspace{-0.2\baselineskip}$$\xymatrix@C=1.3cm{
			W\ar@{<-}[r]^{\mathbbm{a}_{1\times \ell}\qquad}&\s^{n}(\bigvee_{k=1}^{\ell} X_{k}) \ar@{<-}[r]^{\s^{n}\mathbbm{b}_{\ell\times r}} & \s^{n}( \bigvee_{s=1}^{r} Y_{s})\ar@{<-}[r]^{\qquad\s^{n}\mathbbm{c}_{r\times 1}}& \s^{n}U.} $$
		
	\end{defn}
	The notations we have introduced have the advantage that, in the above definition, \emph{the symbol ``\, $\hc$\, '' can be interpreted both as the composition of the homotopy classes and as the multiplication of the formal matrices.} The matrix Toda bracket $\n \mathbbm{a}_{1\times \ell},\; \mathbbm{b}_{\ell\times r}, \mathbbm{c}_{r\times 1}\nn_{0}$ is denoted by $\n \mathbbm{a}_{1\times \ell},\; \mathbbm{b}_{\ell\times r}, \mathbbm{c}_{r\times 1}\nn_{}$ for short.

	The following proposition presents an  efficient method for computing the indeterminacy of the matrix Toda bracket. 
	
	\begin{prop}[{\cite{YMM}}]\label{ind}
		Under the condition of Definition  \ref{lwndy},
		additionally  suppose $n\geq1$ and the spaces $X_{k}$ are suspensions for all $1\leq k\leq\ell.$  Then, $$\n \mathbbm{a}_{1\times \ell},\; \s^{n}\mathbbm{b}_{\ell\times r}, \;\s^{n}\mathbbm{c}_{r\times 1}\nn_{n}\in [\s^{n+1}U, W]/A,$$\noindent where\[A=\sum_{k=1}^{\ell} (a_{k}\hc \s^{n}[\s U, X_{k}]) +       \sum_{s=1}^{r}    ([\s^{n+1}Y_{s}, W]\hc\s^{n+1}c_{s}).\] \qed
	\end{prop}

	The properties in the following exercises are very useful.
	\begin{exercise}Under the condition of Definition  \ref{lwndy},
		additionally  suppose $n\geq1$ and the spaces $X_{k}$ are suspensions for all $1\leq k\leq\ell.$ 
		
		\begin{itemize}
			\item [\rm (1)] Take $\ell=2$ and $r=1.$

			\begin{itemize}
				\item [\rm (i)]     If $a_{1}\hc\s^{n}b_{1,1}=0,$   then, \begingroup
				\renewcommand{\arraystretch}{1.05}
				\setlength{\jot}{0pt}
				\[
				\begin{gathered}
					\left\{
					(a_1,a_2),\,
					\Sigma^n
					\begin{pmatrix}
						b_{1,1}\\
						b_{2,1}
					\end{pmatrix},\,
					\Sigma^n c_1
					\right\}_n
					\\[2pt]
					\vcenter{\hbox{%
							\rotatebox{90}{\scalebox{2.8}[1]{$\displaystyle=$}}%
					}}
					\\[2pt]
					\left\{
					a_1,\Sigma^n b_{1,1},\Sigma^n c_1
					\right\}_n
					+
					\left\{
					a_2,\Sigma^n b_{2,1},\Sigma^n c_1
					\right\}_n.
				\end{gathered}
				\]
				\endgroup
				
				\item [\rm (ii)]  If $c_{1}\hc d_{1}=0$ for some map $d_{1}\colon U\leftarrow \s^{2}U',$   then
				\begingroup
				\renewcommand{\arraystretch}{1.05}
				\setlength{\jot}{0pt}
				\[
				\begin{gathered}
					(-1)^{n+1}\left\{
					(a_1,a_2),\,
					\Sigma^n
					\begin{pmatrix}
						b_{1,1}\\
						b_{2,1}
					\end{pmatrix},\,
					\Sigma^n c_1
					\right\}_n
					\circ\Sigma^{n+1}d_1
					\\[2pt]
					\vcenter{\hbox{%
							\rotatebox{90}{\scalebox{2.8}[1]{$\displaystyle=$}}%
					}}
					\\[2pt]
					(a_1,a_2)\circ\Sigma^n
					\left\{
					\begin{pmatrix}
						b_{1,1}\\
						b_{2,1}
					\end{pmatrix},\,
					c_1,d_1
					\right\}
					\\[2pt]
					\vcenter{\hbox{%
							\rotatebox{-90}{\scalebox{2.8}[1]{$\displaystyle\subseteq$}}%
					}}
					\\[2pt]
					a_1\circ\Sigma^n\{b_{1,1},c_1,d_1\}
					+a_2\circ\Sigma^n\{b_{2,1},c_1,d_1\}.
				\end{gathered}
				\]
				\endgroup
			\end{itemize}
			\item [\rm (2)] Take $\ell=1$ and $r=2.$ \begin{itemize}
				\item [\rm (i)] If $b_{1,1}\hc c_{1}=0,$   then\vspace{-0.8em}\begingroup
				\renewcommand{\arraystretch}{1.05}
				\setlength{\jot}{0pt}
				\[
				\begin{gathered}
					\left\{
					a_1,\,
					\Sigma^n(b_{1,1},b_{1,2}),\,
					\Sigma^n
					\begin{pmatrix}
						c_1\\
						c_2
					\end{pmatrix}
					\right\}_n
					\\[2pt]
					\vcenter{\hbox{%
							\rotatebox{90}{\scalebox{2.8}[1]{$\displaystyle=$}}%
					}}
					\\[2pt]
					\left\{
					a_1,\Sigma^n b_{1,1},\Sigma^n c_1
					\right\}_n
					+
					\left\{
					a_1,\Sigma^n b_{1,2},\Sigma^n c_2
					\right\}_n.
				\end{gathered}
				\]
				\endgroup
				\item [\rm (ii)] If $   \begin{pmatrix}
					c_1\\
					c_2
				\end{pmatrix}\hc d=0$ for some map $d\colon U\leftarrow \s^{2}U'',$   then \vspace{-1em}\begingroup
				\renewcommand{\arraystretch}{1.05}
				\setlength{\jot}{0pt}
				\[
				\begin{gathered}
					(-1)^{n+1}\left\{
					a_1,\,
					\Sigma^n(b_{1,1},b_{1,2}),\,
					\Sigma^n
					\begin{pmatrix}
						c_1\\
						c_2
					\end{pmatrix}
					\right\}_n
					\circ\Sigma^{n+1}d
					\\[2pt]
					\vcenter{\hbox{%
							\rotatebox{90}{\scalebox{2.8}[1]{$\displaystyle=$}}%
					}}
					\\[2pt]
					a_1\circ\Sigma^n
					\left\{
					(b_{1,1},b_{1,2}),\,
					\begin{pmatrix}
						c_1\\
						c_2
					\end{pmatrix},\,
					d
					\right\}.
				\end{gathered}
				\]
				\endgroup
			\end{itemize}
		\end{itemize}
		Hint: use
		the definition of matrix Toda brackets and Propositions \ref{prop:toda-composition} and \ref{prop:toda-four}.
	\end{exercise}

	For Hopf-invariant calculations, work 2-locally.
	Take $n=1$ and $W=S^{m+1},$   and suppose that
	\[
	a_k=\Sigma f_k,\qquad f_k:X_k\longrightarrow S^m.
	\]
	Write
	\[
	\mathbbm{f}_{1\times\ell}=(f_1,\ldots,f_\ell),\qquad
	\mathbbm{b}_{\bullet,s}
	=(b_{1,s},\ldots,b_{\ell,s})^{\mathrm T},
	\]
	where $\mathbbm{b}_{\bullet,s}$ denotes the $s$-th column of
	$\mathbbm{b}_{\ell\times r}.$  Suppose that
	\[
	\Sigma\bigl(
	\mathbbm{f}_{1\times\ell}\circ
	\mathbbm{b}_{\ell\times r}
	\bigr)=O,\qquad
	\mathbbm{b}_{\ell\times r}\circ
	\mathbbm{c}_{r\times1}=O.
	\]
	Then \cite[Theorem~2]{YMM} gives the generalized $H$-formula
	\begin{equation}
		H\left\{
		\Sigma\mathbbm{f}_{1\times\ell},
		\Sigma\mathbbm{b}_{\ell\times r},
		\Sigma\mathbbm{c}_{r\times1}
		\right\}_1
		=
		-\sum_{s=1}^{r}
		P^{-1}\!\left(
		\mathbbm{f}_{1\times\ell}\circ
		\mathbbm{b}_{\bullet,s}
		\right)
		\circ\Sigma^2c_s.
	\end{equation}
	\noindent   Here
	\[
	\mathbbm{f}_{1\times\ell}\circ
	\mathbbm{b}_{\bullet,s}
	=
	\sum_{k=1}^{\ell}f_kb_{k,s}\in[Y_s,S^m],
	\]
	and $P^{-1}$ denotes the inverse-image set under the generalized
	EHP homomorphism
	\[
	P:[\Sigma^2Y_s,S^{2m+1}]
	\longrightarrow[Y_s,S^m].
	\]
	The right-hand side is the sum of the indicated subsets of
	$[\Sigma^2U,S^{2m+1}],$   with the inverse images chosen
	independently for each $s.$  The formula holds for arbitrary
	$\ell,r\geq1.$

	\section{The $n$-Fold Toda Brackets with $n\geq4$}\label{sec:higher-toda}
	
	To date, the most complete references on the theory of unstable
	$n$-fold Toda brackets for $n\geq4$ are \cite{O2015} and
	\cite{O2022}, both written by H.~{\=O}shima and K.~{\=O}shima.
	The paper \cite{O2015} studies the case $n=4,$   while
	\cite{O2022} studies all cases $n\geq4.$   (By ``unstable'',  we mean ``not necessarily stable''.)
	
	The theory of tertiary compositions,
	or $4$-fold Toda brackets, has its
	origins in Toda's pioneering work
	on unstable homotopy theory.
	In 1955, Toda first made use of
	the existence of such compositions
	\cite{Toda1955}, and he subsequently
	employed methods essentially based
	on them in constructing the elements
	$\mu_3\in\pi_{12}(S^3)$ and
	$\kappa_7\in\pi_{21}(S^7)$
	\cite{T}.
	Motivated by these developments,
	\^{O}guchi \cite{Oguchi1963}
	and Mimura \cite{Mimura1964}
	undertook the systematic study
	of tertiary compositions, developing
	their respective theories in 1963
	and 1964.
	However, both papers \cite{Oguchi1963} and \cite{Mimura1964} contain substantial
	mathematical errors, which were
	subsequently identified and corrected
	by H.~\={O}shima and K.~\={O}shima
	in \cite{O2015}.

	The theory of higher compositions was
	further developed by J.~Cohen in 1968,
	through his definition of $n$-fold
	Toda brackets for every $n\geq3$
	\cite{Cohen1968}.
	His construction provides a general
	setting for higher Toda brackets,
	although it may yield a larger set
	than the classical construction,
	even in the case $n=3$; see
	\cite[Appendix~B]{O2015}.
	While such generality is valuable,
	the resulting indeterminacy may
	occasionally make the bracket less
	effective for computations in
	unstable homotopy theory, where
	more precise information about
	individual homotopy classes
	is often required.

	A careful reassessment of the
	earlier theory of unstable tertiary
	compositions was eventually undertaken
	by H.~\={O}shima and K.~\={O}shima
	in \cite{O2015}.
	Their work revisited the formulations
	proposed by \^{O}guchi and Mimura,
	identified the difficulties in the
	original arguments, and provided
	the corrections necessary to address
	the mathematical gaps and errors.

	Let the sequence
	\[
	X_0
	\xleftarrow{\alpha_1}\Sigma^n X_1
	\xleftarrow{\Sigma^n\alpha_2}\Sigma^n X_2
	\xleftarrow{\Sigma^n\alpha_3}\Sigma^n X_3
	\xleftarrow{\Sigma^n\alpha_4}\Sigma^n X_4
	\]

	\noindent   of spaces and homotopy classes of maps
	be given, where $\af_{1}\s^{n}\af_{2}=\af_{2}\af_{3}=\af_{3}\af_{4}=0$ and  \begin{equation}\label{rtj}
		0\in A:=\n \af_{1},\s^{n}\af_{2},\s^{n}\af_{3}\nn_{n},\qquad    0\in B:=\n \af_{2},\af_{3},\af_{4}\nn.
	\end{equation}
	\noindent Suppose that \begin{equation}\label{y0yq}
		\mathrm{Ind}(A)=0\;\;\mathrm{or}\;\;\mathrm{Ind}(B)=0.
	\end{equation}
	
	\noindent Choose the extensions and coextensions below from a single admissible representative, and fix these compatible choices.
	In particular, $\widetilde{\alpha_3}$ and $\overline{\alpha_2}$ are induced by the same nullhomotopy of $\alpha_2\alpha_3.$ 
	With these choices, we define the Toda brackets $$L=\{\alpha_1,\Sigma^n\alpha_2,\Sigma^n\alpha_3,
	\Sigma^n\alpha_4\}_n^L\;\;\mathrm{and}\;\;R=\{\alpha_1,\Sigma^n\alpha_2,\Sigma^n\alpha_3,
	\Sigma^n\alpha_4\}_n^R $$ \noindent as follows:\vspace{-0.5em}
	\[
	\begin{aligned}
		L
		&:=\{\overline{\alpha_1},\Sigma^n\widetilde{\alpha_3},
		\Sigma^n(-\Sigma\alpha_4)\}_n
		\subseteq[\Sigma^{n+2}X_4,X_0], \\[6pt]
		R
		&:=\{\alpha_1,\Sigma^n\overline{\alpha_2},
		\Sigma^n\widetilde{\alpha_4}\}_n
		\subseteq[\Sigma^{n+2}X_4,X_0].
	\end{aligned}
	\]
	\noindent   We note \vspace{-0.8em}$$ \overline{\alpha_1}  \in\operatorname{Ext}_{\Sigma^n\alpha_2}(\alpha_1)
	\subseteq[C_{\Sigma^n\alpha_2},X_0],\quad       \widetilde{\alpha_3}\in   \operatorname{Coext}_{\alpha_2}(\alpha_3)
	\subseteq[\Sigma X_3,C_{\alpha_2}],$$
	\vspace{-0.8em}
	$$      \overline{\alpha_2}\in \operatorname{Ext}_{\alpha_3}(\alpha_2)
	\subseteq[C_{\alpha_3},X_1], \quad  \widetilde{\alpha_4}\in \operatorname{Coext}_{\alpha_3}(\alpha_4)
	\subseteq[\Sigma X_4,C_{\alpha_3}]. $$
	
	\noindent  In general, the construction \textbf{does not} justify $L=R,$   or even the weaker assertion $L=\pm R.$ 
	\noindent Instead, we have $$L\cap R\neq \emptyset.$$

	\noindent   See \cite[Theorem 4.9 and its proof]{O2015}.
	However, $L$ is not independent of the choices of
		$\overline{\alpha_1}$ and $\widetilde{\alpha_3}$  and thus does not depend only
	on $(\af_1,\af_2,\af_3,\af_4,n);$  a similar statement applies to $R.$  See \cite[Example 5.10]{O2015}.

	Mimura \cite{Mimura1964} incorrectly claimed that $L=R$ and did not
	prove that either of the Toda brackets
	\[
	\{\overline{\alpha_1},\Sigma^n\widetilde{\alpha_3},
	\Sigma^n(-\Sigma\alpha_4)\}_n
	\quad\text{or}\quad
	\{\alpha_1,\Sigma^n\overline{\alpha_2},
	\Sigma^n\widetilde{\alpha_4}\}_n
	\]
	is independent of the choices of extensions or coextensions. (He considered the case $n=0$; note that analogues of his statements apply to the case
	$n\geq1$ stated above.) 
	
	Unfortunately, many classical papers follow Mimura's definition,
	treating $L$ and $R$ as equal and independent of the choices of
	extensions and coextensions. For a specific calculation, their notation may be interpreted
	$$\{\alpha_1,\Sigma^n\alpha_2,\Sigma^n\alpha_3,\Sigma^n\alpha_4\}_n$$\noindent
	as shorthand for either
	$\{\overline{\alpha_1},\Sigma^n\widetilde{\alpha_3},
	\Sigma^n(-\Sigma\alpha_4)\}_n$
	or
	$\{\alpha_1,\Sigma^n\overline{\alpha_2},
	\Sigma^n\widetilde{\alpha_4}\}_n,$  
	with compatible extensions and coextensions chosen and fixed. Under this interpretation, the results on the homotopy groups $\pi_*(S^n)$ established by means of $4$-fold Toda brackets continue to hold.

	In concrete calculations, using $L$ or $R$ is  usually more convenient
	than using the formal definition of the $4$-fold Toda bracket given
	in \cite[Section~6, p.~56]{O2015}. We use $L$ or $R$ simply as shorthand for the corresponding
	$3$-fold Toda bracket, with compatible extensions and coextensions fixed and the relevant indeterminacy checked in each application.
	See Lemma~\ref{m3:tsd}, which makes the required choices explicitly in the calculation of $\pi_{28}(G_{2}).$ 
	For $n$-fold Toda brackets with $n\geq5,$   we may use ideas analogous
	to those used in the case $n=4,$   whereas the situation becomes
	more complicated. For further knowledge of these higher brackets,   see \cite{O2015,O2022}.

	\begin{rem} To ensure that the Toda brackets $L$ and $R$ are well-defined for some
		choices of extensions and coextensions, the hypothesis in
		Equation~(\ref{y0yq}) can be weakened to the requirement that
		$(\af_{1},\af_{2},\af_{3},\af_{4})_{n,0}$ be admissible.
		See Proposition~\ref{prop:admissible-null-quadruple}.
	\end{rem}

	\chapter{Notation and Methods of Calculation}\label{ch:notation}\label{m2:sec:preliminaries}
	The calculations combine three kinds of information: the boundary maps
	in the fibrations of special orthogonal groups, relations among lifts of elements of homotopy groups of spheres, and the images of these lifts under the unstable
	$J$-homomorphism.  We first fix the notation and explain how these
	pieces fit together.  
	\section{Notation, Compositions, and Lifts}\label{m2:subsec:notation}
	
	The notation is designed to keep track of where a class lives and how it was obtained. Subscripts on
	elements of homotopy groups of spheres specify their targets, while brackets around a class indicate a chosen lift. The
	conventions below make these distinctions precise\uprightsemicolon{} the exact sequence in
	Chapter~\ref{chb} explains how they enter a calculation.
	
	All spaces, maps and homotopies are based. We refer to a homotopy class of a map as a map  when no confusion arises. Unless otherwise stated, we \textbf{work
		$2$-locally}\uprightsemicolon{}  the subscripts  ``\((2)\)'' in the symbols for spaces \(X_{(2)}\), maps \(f_{(2)}\) and cells \(e^n_{(2)}\) are   omitted. The notation $\pi_k^n=\pi_k(S^n)$ stands for homotopy groups of spheres, and
	$U_k^3=\pi_k(\mathit{SU}(3)).$  The stable group $\pi_d^S=\pi_{d}^{S}(S^{0})$ is also understood  $2$-locally.
	
	For $r\geq1,$   $\mathbb Z/2^r$ denotes the cyclic group $\mathbb Z/2^r\mathbb Z,$   and
	$\frac{\mathbb Z}{2^r}\{x\}$ specifies both a cyclic summand and its generator.
	In $(\frac{\mathbb Z}{2^r})^{\,s}\{x_1,\ldots,x_s\},$   the listed elements generate
	the $s$ respective summands. We write
	$\mathrm{span}_{\mathbb Z_{(p)}}\{x_1,\ldots,x_s\}$ for the
	$\mathbb Z_{(p)}$-submodule generated by the listed elements.

	Under the standing $2$-local convention,
	$\mathrm{span}\{x_1,\ldots,x_s\}$ 
	denotes the $\mathbb Z_{(2)}$-submodule generated by $x_1,\ldots,x_s$\uprightsemicolon{} explicitly,
	\[
	\mathrm{span}\{x_1,\ldots,x_s\}
	=\left\{\sum_{i=1}^s r_i x_i : r_i\in\mathbb Z_{(2)}\right\}.
	\]

	\noindent
	For convenience, in proposition, lemma, or theorem statements we write
	$$\pi_{*}(X)=\mathrm{span}\{a_1,\ldots,a_r\}
	\cong G_1\oplus\cdots\oplus G_r$$
	to mean a direct sum decomposition into cyclic
	$\mathbb Z_{(2)}$-modules, with \textbf{$a_i$ generating the summand
		corresponding to $G_i.$ } Consecutive equal cyclic summands may be abbreviated by a power\uprightsemicolon{} the generators are still listed in the corresponding order.

	Odd-primary components are exceptions to the standing convention.   For $i\geq2,$   we write
	\[
	\pi_i(\mathit{SO}(n):\mathrm{odd})
	:=\{x\in\pi_i(\mathit{SO}(n)):
	mx=0\text{ for some positive odd integer }m\}.
	\] Here, clearly, $\mathit{SO}(n)$ should be understood as an integral space. This is the subgroup generated by all elements of odd order, or equivalently,
	the direct sum of the odd-primary components of the integral homotopy group.
	In a calculation at an odd prime $p,$   all spaces and maps are localized at $p$\uprightsemicolon{} a prime-specific space is written $X_{(p)}$ when more than one prime occurs in the same section.

	A Toda bracket
	with zero indeterminacy is identified with its unique element: if that
	element is $a,$   we write $\{\alpha,\beta,\gamma\}=a.$ 
	We write $\mathrm{ord}(x)$\;  (or simply    $\mathrm{ord}\,x$)  for the order of $x$  in a group. And we write $|A|$ for the cardinality of a set $A.$ 
	Recall that,  for elements $a,b,c,d,e,\cdots$ of the same homotopy group, the notation
	$a\equiv b\bmod c,d,e,\cdots$ means
	\[
	a-b\in\mathrm{span}\{c,d,e,\cdots\},
	\]
	where the right-hand side is the $\mathbb Z_{(2)}$-submodule generated by
	$c,d,e,\cdots.$  In an explicitly $p$-local section, use $\mathbb Z_{(p)}$ instead.
	
	The identity class of $S^n$ is $\iota_n.$    Let $\af$ and $\beta$ be homotopy classes so that $\af\hc\beta$ is defined and the domain and codomain of  $\beta$ are  suspensions; let $t\in\z.$
	
	Note that $(t\alpha)\circ\beta$ and
	$t(\alpha\beta)$ agree when $\beta$ is a suspension, and also for
	the Lie-group targets used here, but they must not be identified for
	arbitrary unstable maps between spheres.  We write $$t(\alpha\beta):=t(\alpha\beta)=\alpha\circ\beta\hc t\,\mathrm{id}.$$

	For $a<b,$   let $i_{a,b}:\mathit{SO}(a)\to\mathit{SO}(b)$ be
	the standard inclusion, and write $i_{n+1}=i_{n,n+1}.$ 
	The projection $p_{n+1}:\mathit{SO}(n+1)\to S^n$ sends a matrix
	to its last column (before 2-localization).  In the exact sequence induced by the fibration
	$\mathit{SO}(n)\stackrel{{i_{n+1}}}{\longrightarrow}\mathit{SO}(n+1)
	\stackrel{{p_{n+1}}}{\longrightarrow}S^n,$   the connecting homomorphism is denoted by 
	$$\Delta\colon \pi_{k}^n\to \pi_{k-1}(\mathit{SO}(n)).$$
	If $\mathbf x\in\pi_{k}^n$ satisfies $\Delta\mathbf x=0,$   then $[\mathbbm x]\in \pi_{k}(\mathit{SO}(n+1))$
	denotes a chosen lift, so that $(p_{n+1})_*[\mathbf x]=\mathbf x.$ 
	The image of $[\mathbf x]\in \pi_{k}(\mathit{SO}(n+1))$  in  $\pi_{k}(\mathit{SO}(b))$ ($b> n+1$) is written as
	$$[\mathbf x]_b:=(i_{n+1,b})_*[\mathbf x]\in\pi_{k}(\mathit{SO}(b)).$$
	
	\noindent When the context is clear, for convenience, we sometimes denote $p_{n+1},$   $i_{a,b}$ by $p, i,$   respectively. 
	
	The element  lying in $ \pi_{14}(\mathit{SO}(7)) $ revised by Kachi and Mukai in \cite{KM1} is denoted by
	$$[\mathrm{c}_{\mathrm{KM}}(\bar{\nu}_6+\varepsilon_6)],$$
	where $\mathrm{c}_{\mathrm{KM}}\in\z$ is a fixed odd integer and
	$(p_7)_*[\mathrm{c}_{\mathrm{KM}}(\bar{\nu}_6+\varepsilon_6)]
	=\mathrm{c}_{\mathrm{KM}}(\bar{\nu}_6+\varepsilon_6).$
	This notation refers throughout to the choice satisfying
	\eqref{m1:Jbn6e6} and \eqref{m1:let6sg}.
	The lift induced by $\mathcal{T}_*$ is denoted separately by
	$[\bar{\nu}_6+\varepsilon_6]$ in Section~\ref{lqzh}.
	
	\section{The Calculation of Unstable Homotopy Groups of Spheres}\label{sec:sphere-history}
	For $\pi_{i+k}(S^k),$   the \emph{stem} is $i$ and the sphere dimension is $k.$  Fixing $i$ and varying
	$k$ follows the suspension maps through the unstable groups and into the stable group.
	Table~\ref{tab:sphere-history} displays both indices, so a result about a stem cannot be mistaken for
	a result about one sphere. It records the main sequence of systematic computations. Complementary
	results on generator orders, individual groups and unresolved extensions
	are discussed below the table. The historical account in
	\cite{YangWuHHA} provides an overview\uprightsemicolon{} the original sources are cited in each row.

	For the reader's convenience, we summarize the known unstable homotopy
	groups of spheres in the ranges covered by Appendix~\ref{tlqdzj}, together with explicit generating sets,
	based on the references listed below.
	Here, ``unstable'' means ``not necessarily stable''.

	Recall that, by default, all statements are understood $2$-locally unless otherwise stated.
	\begingroup\small\setlength{\tabcolsep}{5pt}\renewcommand{\arraystretch}{1.35}\arrayrulecolor{earlytablerule}
	\begin{longtable}{|>{\raggedright\arraybackslash}p{18mm}|>{\raggedright\arraybackslash}p{44mm}|>{\raggedright\arraybackslash}p{24mm}|>{\raggedright\arraybackslash}p{54mm}|}
		\caption{Systematic Computations of the Unstable Homotopy Groups of Spheres $\pi_{i+k}^{k}=\pi_{i+k}(S^{k})$. The scope column describes the results established in the cited publications.}\label{tab:sphere-history}\\
		\hline
		\rowcolor{earlytablehead} \textbf{Year} & \textbf{Authors} & \textbf{Stem $i$} & \textbf{Scope}\\ \hline
		\endfirsthead
		\multicolumn{4}{c}{\small\itshape Systematic Computations of the Unstable Homotopy Groups of Spheres at 2 Continued}\\
		\hline
		\rowcolor{earlytablehead} \textbf{Year} & \textbf{Authors} & \textbf{Stem $i$} & \textbf{Scope}\\ \hline
		\endhead
		\hline
		\endfoot
		\endlastfoot
		\rowcolor{blue!3} 1962 & H.~Toda & $0\leq i\leq19$ &all  $k\geq2$ and all displayed $i,$     \cite{T}.\\ \hline
		1963 & M.~Mimura and H.~Toda & $i=20$ & all $k\geq2$ and the displayed $i$ \cite{MT1}.\\ \hline
		\rowcolor{blue!3} 1965 & M.~Mimura & $i=21,22$ & all  $k\geq2$ and all displayed $i,$    \cite{Mi1}.\\ \hline
		1975 & M.~Mimura, M.~Mori and N.~Oda & $i=23,24$ & all  $k\geq2$  and all displayed $i,$    \cite{MMO}.\\ \hline
		\rowcolor{blue!3} 1979 & N.~Oda & $25\leq i\leq31$ & all  $k\geq2$ and all displayed $i$ (except for $\pi_{40}^{9}$ and $\pi_{41}^{10},$   whose extension problems were left open in that paper and were resolved in \cite{IMM}), \cite{Oda1979}.\\ \hline
		\rowcolor{blue!3} 2017 & T.~Miyauchi and J.~Mukai & $i=32$ &all  $k\geq2$  and the displayed $i,$   \cite{MiM1}.\\ \hline
	\end{longtable}
	\endgroup
	The $0$-stem is elementary: $\pi_k(S^k)=\mathbb Z_{(2)}.$  For fixed positive $i,$   Freudenthal
	suspension gives the stable value once $k\geq i+2$\uprightsemicolon{} the lower-dimensional entries are the unstable
	part of the calculation. \^{O}guchi's 1964 paper \cite{Og} supplies a particularly convenient catalogue of
	generators and their relations in the classical range. Its reference material is collected in
	Appendix~\ref{app:oguchi}, together with the companion composition table. We repeatedly use these
	tables without a separate citation at each occurrence.
	
	The following complementary results refine the calculations summarized in the table
	and extend some of them to individual groups in the next stems.
	In 1986, Oda clarified how the named generators realize the abstract group decompositions.
	For the $18$-stem, Oda determined the orders of Toda's generators in
	$\pi_{28}^{10},$   $\pi_{29}^{11}$ and $\pi_{30}^{12}$ and gave explicit cyclic direct summands
	\cite[Theorem, pp.~231--232]{Oda}. In the first of these groups, the sign in
	$\xi''\pm\lambda''$ is chosen so that this element has order $2.$ 
	These decompositions are recorded in Appendix~\ref{app:stem-18}.
	
	In 2015, Inoue, Miyauchi and Mukai resolved, in the $31$-stem, the two cases left open in Oda's tables,
	$\pi_{40}^{9}$ and $\pi_{41}^{10}$ \cite[Theorem~1.1]{IMM}.
	They proved that $\bar\nu_9\nu_{17}\bar\kappa_{20}=0$ and
	$\delta_{10}\sigma_{34}=2\kappa_{10}^{*},$   with $\kappa_{10}^{*}$ of order $4.$ 
	Appendix~\ref{app:stem-31} incorporates these determinations.
	
	Oda's 1979 work also treats the $32$- and $33$-stems in sphere dimensions
	$2\leq k\leq8.$  For the $33$-stem, Proposition~9.20 gives the groups
	$\pi_{33+k}^{k}$ for $2\leq k\leq8,$   with unresolved extensions in dimensions $6,$   $7$ and $8$
	\cite[pp.~144--145]{Oda1979}. In 2017, Miyauchi and Mukai determined $\pi_{42}^{9},$  
	together with its six generators, in their study of the $32$-stem \cite{MiM1}.
	In 2025, Yang and Wu proved $2\kappa'_6=\nu_6\sigma_9\nu_{16}\bar\kappa_{19}$ and thereby
	resolved the three remaining extension problems in dimensions $6,$   $7$ and $8$
	\cite[Theorem~1.1 and Corollary~1.2]{YangWuHHA}.
	Combining these results, Appendix~\ref{app:stem-33} lists the groups $\pi_{33+k}^{k}$
	for $2\leq k\leq9,$   with generators matched to their cyclic summands.
	
	\medskip
	Table~\ref{tab:sphere-history} shows how these calculations developed:
	Toda's composition methods provided a framework for explicit calculations, while the work of
	Mimura, Oda and their collaborators extended the range of known homotopy groups of spheres. These computations
	supply the generators and relations used throughout the later chapters.
	
	Figure~\ref{fig:toda-mimura-oda} accompanies this history with a
	view of the people behind it. The full photograph preserves the
	setting of the gathering\uprightsemicolon{} the detail below identifies H.~Toda,
	M.~Mimura and N.~Oda within the group. Their names recur in the
	preceding account because their methods and calculations remain
	part of the working language of this book.
	
	\par\medskip
	\noindent\begin{minipage}{\linewidth}
		\centering
		\BookPhoto[width=1\linewidth,height=0.41\textheight,keepaspectratio]{\GroupPhotoFile}
		\par\smallskip
		{\footnotesize\justifying (a) The complete group photograph, taken in August 1988 at the
			International Conference on Homotopy Theory in Kinosaki, Hyogo, Japan
			(19--24 August), celebrating H.~Toda's sixtieth birthday.\par}
	\end{minipage}
	\par\medskip
	\noindent\begin{minipage}{\linewidth}
		\centering
		\BookPhoto[width=1.2\linewidth,height=.345\textheight,keepaspectratio]{\DetailPhotoFile}
		\par\smallskip
		{\footnotesize\justifying (b) A closer view with the three mathematicians identified.\par}
		\captionsetup{width=.90\linewidth}
		\captionof{figure}{H.~Toda, M.~Mimura and N.~Oda in a group photograph taken in
			August 1988.
			The full view records the gathering\uprightsemicolon{} the annotated detail identifies
			three mathematicians whose work is central to this book.}
		\label{fig:toda-mimura-oda}
	\end{minipage}
	\par\medskip

	\par\medskip
	\noindent The photograph below complements this account with a view of
	M.~Mori, one of the authors of \cite{MMO}. It is the group photograph
	shown in the Preface, reproduced here with handwritten identifications.
	The setting is Shikanoshima, Fukuoka, Japan, in 1973. The photograph was
	taken a few days after the close of the International Conference on
	Manifolds and Related Topics in Topology, which was organized by
	M.~F.~Atiyah and Kunihiko Kodaira and held from 10 to 17 April 1973.
	\par\nopagebreak\smallskip
	\noindent\begin{minipage}{\linewidth}
		\centering
		\BookPhoto[width=.616\linewidth,height=.294\textheight,keepaspectratio]{\MoriAnnotatedPhotoFile}
		\captionsetup{width=.88\linewidth,font=normalsize,skip=6pt,
			justification=justified,singlelinecheck=false}
		
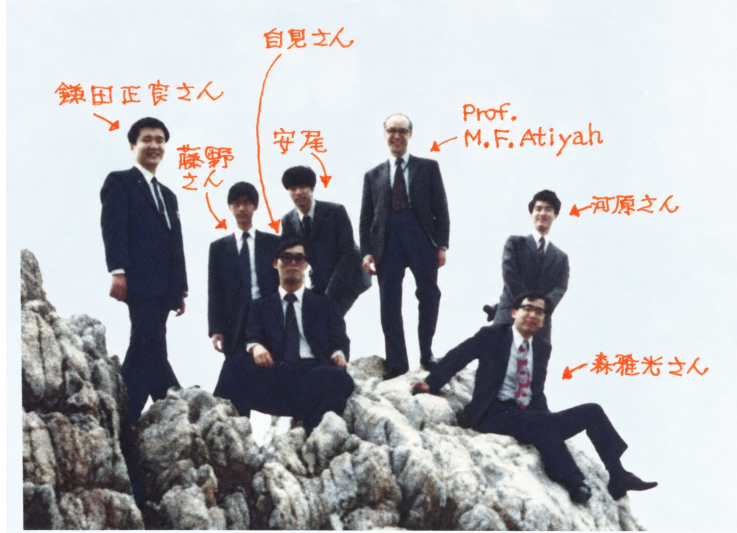
\captionof{figure}{Handwritten names and arrows identify the participants
			in the gathering. M.~Mori is seated at the lower right.}
		\label{fig:mori-annotated}
	\end{minipage}

	\chapter{The EHP Sequence and an Introduction to New Elements in Toda's Range}\label{xzzj}
	\par\medskip
	
	\section{A Review of the  EHP Sequence}\label{sec:ehp-review}
	We now turn to the EHP sequence, a principal tool in the calculations
	that follow. The term “new elements” in the title refers to those elements of \(\pi_{*}(S^{n})\) whose \(H\)-images are nontrivial.
	\par
	Suspension relates homotopy groups of spheres in successive dimensions, but unstable information can disappear
	along the way. The EHP sequence organizes this loss of information through the maps $H$ and $P.$  We
	recall it here because the calculations of the homotopy groups of special orthogonal groups repeatedly compare a relation with its
	suspended image or with its Hopf invariant.
	
	The suspension $\Sigma,$   the second James--Hopf homomorphism $H,$   and the Whitehead homomorphism $P$
	form the EHP sequence.  The James--Hopf map $H_{2}$ fits into the homotopy fibration
	\[\begin{aligned}
		&S^n\longrightarrow\Omega S^{n+1} \xrightarrow{H_2}\Omega S^{2n+1},\;n\geq1.
	\end{aligned}\]
	It gives  the exact sequence  
	\begin{equation*}       \begin{aligned}
			&\cdots\longrightarrow\pi_{q+1}(S^{2n+1}) \xrightarrow{P}\pi_{q-1}(S^n) \xrightarrow{\Sigma}\pi_q(S^{n+1}) \xrightarrow{H}\pi_q(S^{2n+1}) \xrightarrow{P}\pi_{q-2}(S^n)\longrightarrow\cdots .
		\end{aligned}
	\end{equation*}
	Then there is a short exact sequence
	$$0\to \s\pi_{q-1}(S^n)\stackrel{\xyd}\longrightarrow \pi_{q}(S^{n+1})\stackrel{H}\longrightarrow H\pi_{q}(S^{n+1})\to 0,$$
	\noindent where $$ \s\pi_{q-1}(S^n)\cong \mathrm{coker}(\pi_{q+1}(S^{2n+1}) \xrightarrow{P}\pi_{q-1}(S^n)),\quad H\pi_{q}(S^{n+1})=\mathrm{ker}(\pi_q(S^{2n+1}) \xrightarrow{P}\pi_{q-2}(S^n)).$$ The short exact sequence presented above forms the backbone of the computation of the homotopy groups of spheres, serving as the unifying thread that runs through the successive arguments and brings the various calculations together.

	For fixed \(k=q-(n+1)\), the homotopy groups \(\pi_{i+k}(S^i)\) may be computed successively as \(i\) increases, thereby providing a systematic procedure for determining these groups for all \(i\). As is well known, the computations need only be carried out for \(2\leq i\leq k+2\).  The transition from the computations in the $k$-stem to those in the \((k+1)\)-stem is facilitated by the isomorphism \(\pi_{\bullet\geq3}(S^3)\cong\pi_{\bullet\geq3}(S^2)\), induced by the Hopf map \(\ca_2\colon S^3\to S^2\).

	The following formulas play an important role in calculating the relevant kernels and cokernels:
	\begin{equation}\label{1028ppp}
		\begin{gathered}
			P(\iota_{2n+1})=\pm[\iota_n,\iota_n],
			\qquad
			P(\Sigma^{n+1}\alpha)=\pm[\iota_n,\alpha],\\
			P(\beta\hc\s^{2}\gamma)=P(\beta)\hc\gamma.
		\end{gathered}
	\end{equation}
	\noindent
	Here, $\alpha\in\pi_a(S^n)$,
	$\beta\in\pi_{b+2}(S^{2m+1})$ and
	$\gamma\in\pi_{c}(S^{b})$,
	where $a,b,c,m, n\in\z_+$.  Moreover, we have $$P\n\af,\s^{n}\beta,\s^n\gamma\nn_n\xyd \pm\n P\af,\s^{n-2}\beta,\s^{n-2}\gamma\nn_{n-2}$$
	with respect to $$S^{m}\stackrel{\af}\longleftarrow \s^{n}X\stackrel{\s^n \beta}\longleftarrow \s^{n}Y\stackrel{\s^n \gamma}\longleftarrow \s^{n}Z,\quad(n\geq2).$$
	
	Recall that
	that the $2$-localized group $\pi_{i}(S^{n})$ is also denoted by $\pi_{i}^{n}.$    The Hopf invariants and Whitehead products of the other elements are calculated by the following formulas:
	\begin{align*}
		H(\alpha\circ \Sigma\beta)&=H(\alpha)\circ \Sigma\beta
		&&\text{for }\alpha\in\pi_m^k,\ \beta\in\pi_{n-1}^{m-1}.\\
		H(\Sigma\alpha\circ\beta)&=\Sigma^k\alpha\circ \Sigma^m\alpha\circ H(\beta)
		&&\text{for }\alpha\in\pi_{m-1}^{k-1},\ \beta\in\pi_n^m.
	\end{align*}
	Additionally, $H([\e_{m},\e_{m}])=0$ if $m\in\z_{+}$ is odd, and
	$H([\e_{m},\e_{m}])=\pm2\e_{2m-1}$ if  $m\in\z_{+}$ is even.

	For $\alpha\in\pi_m^k$ and $\beta\in\pi_n^k,$   which may not be suspensions,
	\begin{equation*}
		[\alpha,\beta]=\begin{cases}
			[\iota_k,\iota_k]\circ \Sigma^{k-1}\beta\circ \Sigma^{n-1}\alpha,&\text{if }k\text{ is odd},\\
			(-1)^{m+n}[\iota_k,\iota_k]\circ \Sigma^{k-1}\beta\circ \Sigma^{n-1}\alpha,&\text{if }k\text{ is even}.
		\end{cases}
	\end{equation*}
	(See \cite[Corollary (7.4)]{BB}.)
	We recall the following results of P.\,J. Hilton.
	\textit{For an element $\alpha\in\pi_n^r,$   we have}
	
	\begin{enumerate}
		\item \textit{If $r$ is odd, then}
		\[
		\begin{aligned}
			&k\iota_r\circ\alpha = k\alpha\quad  \text{when } k\equiv 0,\ 1 \pmod{4}, \\
			&k\iota_r\circ\alpha = k\alpha+P\bigl(\Sigma^2H(\alpha)\bigr) \quad \text{when } k\equiv 2,\ 3 \pmod{4}.
		\end{aligned}
		\]
		
		\item \textit{If $r$ is even, then}
		\[\begin{aligned}
			&k\iota_r\circ\alpha = k\alpha \pm \frac{k(k-1)}{2}\, P\bigl(\Sigma^2H(\alpha)\bigr).
		\end{aligned}\]
	\end{enumerate}
	
	The morphism $P,$   denoted
	$\Delta$ in some  references such as Toda's book \cite{T}, is distinct from the connecting homomorphism $\Delta$ for the fibration of special orthogonal groups
	used in this book.

	\section{Generators in the Toda Range and a Working Reference}\label{sec:oguchi-working}
	Recall that the $2$-localized group $\pi_{i}(S^{n})$ is denoted by $\pi_{i}^{n}.$
	The generators in the Toda range are summarized below. The Hopf invariants and suspension properties are given in \cite{T,Og}.
	\begin{enumerate}[label=\textup{(\arabic*)},leftmargin=2.4em]
		\item \label{eq:oguchi-1}\(\begin{aligned}[t]
			&\iota_n\in\pi_n^n\quad(n\ge1).
		\end{aligned}\)
		\item \(\begin{aligned}[t]
			&\eta_n\in\pi_{n+1}^n\quad(n\ge2),\quad H(\eta_2)=\iota_3.
		\end{aligned}\)
		\item \(\begin{aligned}[t]
			&\nu'\in\pi_6^3,\quad H(\nu')=\eta_5,\quad2\nu'=\eta_3^3.\\
			&\nu_n\in\pi_{n+3}^n\quad(n\ge4),\quad H(\nu_4)=\iota_7,\quad2\nu_5=\Sigma^2\nu'.
		\end{aligned}\)
		\item \(\begin{aligned}[t]
			&\sigma'''\in\pi_{12}^5,\quad\sigma''\in\pi_{13}^6,\quad\sigma'\in\pi_{14}^7,\quad
			\sigma_n\in\pi_{n+7}^n\quad(n\ge8),\\
			&H(\sigma''')=4\nu_9,\quad H(\sigma'')=\eta_{11}^2,\quad H(\sigma')=\eta_{13},\\
			&H(\sigma_8)=\iota_{15},\quad2\sigma''=\Sigma\sigma''',\quad2\sigma'=\Sigma\sigma'',\quad2\sigma_9=\Sigma^2\sigma'.
		\end{aligned}\)
		\item \(\begin{aligned}[t]
			&\varepsilon_n\in\pi_{n+8}^n\quad(n\ge3),\quad H(\varepsilon_3)=\nu_5^2,\quad
			\bar\nu_n\in\pi_{n+8}^n\quad(n\ge6),\quad H(\bar\nu_6)=\nu_{11}.
		\end{aligned}\)
		\item \(\begin{aligned}[t]
			&\mu_n\in\pi_{n+9}^n\quad(n\ge3),\quad H(\mu_3)=\sigma'''.
		\end{aligned}\)
		\item \(\begin{aligned}[t]
			&\varepsilon'\in\pi_{13}^3,\quad H(\varepsilon')=\varepsilon_5,\quad
			2\varepsilon'=\eta_3^2 \varepsilon_5,\quad \Sigma^2\varepsilon'=\pm2\nu_5 \cg_8.
		\end{aligned}\)
		\item \(\begin{aligned}[t]
			&\mu'\in\pi_{14}^3,\quad H(\mu')=\mu_5,\quad2\mu'=\eta_3^2 \mu_5.\\
			&\zeta_n\in\pi_{n+11}^n\quad(n\ge5),\quad H(\zeta_5)=8\sigma_9,\quad2\zeta_5=\pm \Sigma^2\mu'.
		\end{aligned}\)
		\item \(\begin{aligned}[t]
			&\theta'\in\pi_{23}^{11},\quad\theta\in\pi_{24}^{12},\quad H(\theta')=\eta_{21}^2,\quad
			H(\theta)=\eta_{23},\quad \Sigma^2\theta'=0,\quad \Sigma^2\theta=0.
		\end{aligned}\)
		\item \(\begin{aligned}[t]
			&\kappa_n\in\pi_{n+14}^n\quad(n\ge7),\quad H(\kappa_7)=\bar\nu_{13}, \\
			&2\kappa_7=\bar\nu_7 \nu_{15}^2, (\text{using Mimura--Toda's revised definition}).\mbox{\footnotemark}
		\end{aligned}\)\footnotetext{The result on $H(\cn_{7})$ is improved by \cite{MiM1}.}
		\item \(\begin{aligned}[t]
			&\rho^{\mathrm{IV}}\in\pi_{20}^5,\quad\rho'''\in\pi_{21}^6,\quad\rho''\in\pi_{22}^7,\quad
			\rho'\in\pi_{24}^9,\quad\rho_n\in\pi_{n+15}^n\quad(n\ge13),\\
			&H(\rho^{\mathrm{IV}})=4\zeta_9,\quad H(\rho''')=\eta_{11} \mu_{12},\\
			&H(\rho'')\equiv\mu_{13}\m\nu_{13}^3, \eta_{13} \varepsilon_{14}\uprightsemicolon{}\quad H(\rho')=8\sigma_{17},\\
			&H(\rho_{13})=4\nu_{25},\quad2\rho'''= \Sigma\rho^{\mathrm{IV}},\\
			&2\rho''= \Sigma\rho''',\quad
			2\rho'= \Sigma^2\rho'',\quad2\rho_{13}=\Sigma^4\rho'.\\
			&\bar\varepsilon_n\in\pi_{n+15}^n\quad(n\ge3),\quad
			H(\bar\varepsilon_3)\equiv\nu_5 \sigma_8 \nu_{15}\mathbin{\mathrm{mod}}\nu_5 \eta_8 \mu_9.
		\end{aligned}\)
		\item \(\begin{aligned}[t]
			&\zeta'\in\pi_{22}^6,\quad H(\zeta')\equiv\zeta_{11}\mathbin{\mathrm{mod}}2\zeta_{11},\quad \Sigma^3\zeta'=0.\\
			&\omega_n\in\pi_{n+16}^n\quad(n\ge14),\quad H(\omega_{14})=\nu_{27}.\\
			&\eta^{*\prime}\in\pi_{31}^{15},\quad H(\eta^{*\prime})=\eta_{29}^2,\quad \Sigma^2\eta^{*\prime}\equiv0\mathbin{\mathrm{mod}} \Sigma^3\pi_{30}^{14}.\\
			&\eta_n^*\in\pi_{n+16}^n\quad(n\ge16),\quad H(\eta_{16}^*)=\eta_{31},\quad
			\eta_{18}^*\equiv\omega_{18}\mathbin{\mathrm{mod}}\sigma_{18} \mu_{25}.
		\end{aligned}\)
		\item \(\begin{aligned}[t]
			&\bar\mu_n\in\pi_{n+17}^n\quad(n\ge3),\quad H(\bar\mu_3)=\rho^{\mathrm{IV}}.\\
			&\bar\varepsilon'\in\pi_{20}^3,\quad H(\bar\varepsilon')=\bar\varepsilon_5,\quad
			2\bar\varepsilon'=\eta_3^2 \bar\varepsilon_5,\quad \Sigma^2\bar\varepsilon'=2\nu_5 \kappa_8.\\
			&P(\Sigma\theta)\in\pi_{23}^6,\quad H(P(\Sigma\theta))=\theta'.\\
			&\varepsilon_n^*\in\pi_{n+17}^n\quad(n\ge12),\quad H(\varepsilon_{12}^*)=\nu_{23}^2,\quad
			\varepsilon_{13}^*=\eta_{13} \omega_{14}.
		\end{aligned}\)
		\item \(\begin{aligned}[t]
			&\xi'',\lambda''\in\pi_{28}^{10},\quad\xi',\lambda'\in\pi_{29}^{11},\quad\lambda\in\pi_{31}^{13},\\
			&\xi_n\in\pi_{n+18}^n\quad(n\ge12),\quad\nu_n^*\in\pi_{n+18}^n\quad(n\ge16),\\
			&H(\xi'')=\nu_{19}^3+\eta_{19} \varepsilon_{20},\quad
			H(\lambda'')=\eta_{19} \varepsilon_{20}, \\
			&H(\xi')=\varepsilon_{21}+\bar\nu_{21},\quad H(\lambda')=\varepsilon_{21}, \\
			&H(\xi_{12})\equiv\sigma_{23}\mathbin{\mathrm{mod}}2\sigma_{23},\quad H(\lambda)=\nu_{25}^2,\quad
			H(\nu_{16}^*)\equiv\nu_{31}\mathbin{\mathrm{mod}}2\nu_{31},\\
			&2\xi'=\Sigma\xi'',\quad2\xi_{13}=\Sigma^2\xi',\quad2\lambda'=\Sigma\lambda'',\quad
			2\lambda=\Sigma^2\lambda',\quad2\nu_{17}^*=\Sigma^4\lambda,\\
			&4(\xi_{18}+\nu_{18}^*)=2(\xi_{19}+\nu_{19}^*)=\xi_{20}+\nu_{20}^*=0.\mbox{\footnotemark}
		\end{aligned}\)\footnotetext{The improved formulas for $H(\xi''),$   $H(\lambda''),$   $H(\xi')$ and $H(\lambda')$ are given in \cite{Oda}.}
		\item \(\begin{aligned}[t]
			&\bar\mu'\in\pi_{22}^3,\quad H(\bar\mu')\equiv\bar\mu_5\mathbin{\mathrm{mod}}\eta_5 \mu_6 \sigma_{15},\quad
			2\bar\mu'=\eta_3^2 \bar\mu_5.\\
			&\bar\zeta_n\in\pi_{n+19}^n\quad(n\ge5),\quad H(\bar\zeta_5)=8\rho',\quad2\bar\zeta_5=\Sigma^2\bar\mu'.\\
			&\bar\sigma_n\in\pi_{n+19}^n\quad(n\ge6),\quad H(\bar\sigma_6)\equiv\sigma_{11}^2\mathbin{\mathrm{mod}}2\sigma_{11}^2,\quad
			16\bar\sigma_6=\nu_6 \mu_9 \sigma_{18}.\\
			&\omega'\in\pi_{31}^{12},\quad H(\omega')=\varepsilon_{23}\ \text{or}\ \bar\nu_{23},\quad
			\Sigma^2\omega'=2\omega_{14} \nu_{30}.
		\end{aligned}\)
	\end{enumerate}
	Throughout this book, $\cn_7\in\pi_{21}(S^{7})$ follows the revised definition given in \cite{MT1}. Some of these elements can be chosen from Toda brackets \cite{T}:
	\begin{enumerate}[label=\textup{(\arabic*)},leftmargin=2.4em]
		\item \label{eq:oguchi-2}\(\begin{aligned}[t]
			&\ca_{3}^{2}\in\n 2\e_{3},\ca_{3},2\e_{4}\nn. 
		\end{aligned}\)
		\item \(\begin{aligned}[t]
			&\nu'\in\{\eta_3,2\iota_4,\eta_4\}_1.
		\end{aligned}\)
		\item \(\begin{aligned}[t]
			&\nu_{5}^{2}\in\{\eta_5,\nu_6,\eta_9\}.
		\end{aligned}\)
		\item \(\begin{aligned}[t]
			&\sigma'''\in\{\nu_5,8\iota_8,\nu_8\}_3\cap\n\nu_{5},2\nu_{8},4\e_{11}\nn_{3}.\mbox{\footnotemark}
		\end{aligned}\)\footnotetext{This bracket description follows from \cite[p.~84]{Oda1979}.}
		\item \(\begin{aligned}[t]
			&\varepsilon_3\in\{\eta_3,\Sigma\nu',\nu_7\}_1=\{\nu',\nu_6,\eta_9\},\\
			&\varepsilon_n\in\{\eta_n,2\iota_{n+1},\nu_{n+1}^2\}_{n-4}\quad(n\ge4),\\
			&\varepsilon_n\in\{\eta_n,\nu_{n+1}^2,2\iota_{n+7}\}_t,\quad
			\{\nu_n^2,2\iota_{n+6},\eta_{n+6}\}_t,\\
			&\{2\iota_n,\nu_n^2,\eta_{n+6}\}_t\quad(n\ge5),\quad(0\le t\le n-5).\\
			&\bar\nu_n\in\{\nu_n,\eta_{n+3},\nu_{n+4}\}_{ n-2}\quad(n\ge6).
		\end{aligned}\)
		\item \(\begin{aligned}[t]
			&\mu_n\in\{\eta_n,2\iota_{n+1},\Sigma^{n-4}\sigma'''\}_{n-4}+\mathrm{span}\{\nu_n^3\}\quad(n\ge4).
		\end{aligned}\)
		\item \(\begin{aligned}[t]
			&\varepsilon'\in\{\nu',2\nu_6,\nu_9\}_3.
		\end{aligned}\)
		\item \(\begin{aligned}[t]
			&\mu'\in\{\eta_3,2\iota_4,\mu_4\}_1,\quad
			\zeta_n\in\{\nu_n,8\iota_{n+3},\Sigma^{n-4}\sigma'\}_{n-4}\quad(n\ge6).
		\end{aligned}\)
		\item \(\begin{aligned}[t]
			&\theta'\in\{\sigma_{11},2\nu_{18},\eta_{21}\}_1,\quad
			\theta\in\{\sigma_{12},\nu_{19},\eta_{22}\}_1.
		\end{aligned}\)
		\item \(\begin{aligned}[t]&-\cn_7+ 2t\cg'\cg_{14}\in\n \nu_7,\;\s\mathrm{ext}_{2\e_{10}}(\ca_9),\;\s\mathrm{coext}_{2\e_{10}}(\ch_{10})\nn_1\text{\;for some\;}t\in\z, \text{some extension}\\&\text{and some coextension}.\\
			&\sigma_{11}^{2}\in\n \nu_{11},\cg_{14},\;\ell\nu_{21}\nn_{2} \quad(\ell:\text{odd}).\mbox{\footnotemark}
		\end{aligned}\)\footnotetext{This bracket description follows from  \cite[(2.8)]{MMO} and the fact that $\pi_{25}^{11}\stackrel{\s}\longrightarrow\pi_{26}^{12}$ is injective.} 
		\item \(\begin{aligned}[t]
			&\bar\varepsilon_n\in\{\varepsilon_n,2\iota_{n+8},\nu_{n+8}^2\}_{n+3}\quad(n\ge3),\quad
			\rho^{\mathrm{IV}}\in\{\sigma''',2\iota_{12},8\sigma_{12}\}_1,\\
			&\rho'''\in\{\sigma'',4\iota_{13},4\sigma_{13}\}_1,\quad
			\rho''\in\{\sigma',8\iota_{14},2\sigma_{14}\}_1,\quad
			\rho'\in\{\sigma_9,16\iota_{16},\sigma_{16}\}_1,\\
			&\rho_{13}\in\left\{\bigl(\sigma_{13}, \Sigma\theta\bigr),\,\begin{pmatrix}8\sigma_{20}\\\eta^2_{25}\end{pmatrix},2\iota_{27}\right\}_4\cap\left\{\bigl(\sigma_{13}, \Sigma\theta\hc\ca_{25}\bigr),\,\begin{pmatrix}8\sigma_{20}\\\eta_{26}\end{pmatrix},2\iota_{27}\right\}_4.\\&\pm\co_{17}\in\n \cg_{17},2\cg_{24},8\e_{31}\nn.\mbox{\footnotemark}
		\end{aligned}\)\footnotetext{The matrix bracket describing $\rho_{13}$ is established in \eqref{m1:rho13} and  \eqref{m001:rho13}  of this book.}
		\item \(\begin{aligned}[t]
			&\zeta'\in\{\sigma'',\varepsilon_{13},2\iota_{21}\}_1,\quad
			\eta^{*\prime}\in\{\sigma_{15},4\sigma_{22},\eta_{29}\}_1,\quad
			\eta_{16}^*\in\{\sigma_{16},2\sigma_{23},\eta_{30}\}_1.
		\end{aligned}\)
		\item \(\begin{aligned}[t]
			&\bar\mu_3\in\{\mu_3,2\iota_{12},8\sigma_{12}\}_1.
		\end{aligned}\)
		\item \(\begin{aligned}[t]
			&\nu_{16}^*\in\{\sigma_{16},2\sigma_{23},\nu_{30}\}_1,\quad
			\xi_{12}\in\{\sigma_{12},\nu_{19},\sigma_{22}\}_1.
		\end{aligned}\)
		\item \(\begin{aligned}[t]
			&\bar\mu'\in\{\mu',4\iota_{14},4\sigma_{14}\}_1,\quad
			\bar\zeta_5\in\{\zeta_5,8\iota_{16},2\sigma_{16}\}_1,\quad
			\bar\sigma_6\in\{\bar\nu_6+\varepsilon_6,\sigma_{14},\nu_{21}\}_1.
		\end{aligned}\)
	\end{enumerate}
	We list some Whitehead products of these generators below.
	\begin{align*}
		&[\iota_1,\iota_1]=[\iota_3,\iota_3]=[\iota_7,\iota_7]=0,\quad
		[\iota_5,\iota_5]=\nu_5 \eta_8,\\
		&[\iota_9,\iota_9]=\sigma_9 \eta_{16}+\varepsilon_9+\bar\nu_9,\quad
		[\iota_{11},\iota_{11}]=\sigma_{11} \nu_{18},\\
		&[\iota_{13},\iota_{13}]=\Sigma\theta,\quad[\iota_{15},\iota_{15}]=2\sigma_{15}^2,\\
		&[\iota_{17},\iota_{17}]\equiv\eta_{17}^*+\omega_{17}\mathbin{\mathrm{mod}}\sigma_{17} \mu_{24},\quad
		[\iota_{19},\iota_{19}]=\nu_{19}^*+\xi_{19}.\\
		&[\iota_2,\iota_2]=\pm2\eta_2,\quad[\iota_4,\iota_4]=\pm(\Sigma\nu'-2\nu_4),\\
		&[\iota_8,\iota_8]=\pm(\Sigma\sigma'-2\sigma_8),\quad[\iota_6,\eta_6]=0,\quad
		[\iota_6,\nu_6]=\pm2\bar\nu_6,\\
		&[\iota_6,\varepsilon_6]=[\iota_6,\bar\nu_6]=0,\quad[\iota_6,\mu_6]=0,\quad
		[\iota_6,\zeta_6]=\pm2\zeta'.\\
		&[\iota_{10},\eta_{10}]=2\sigma_{10} \nu_{17},\quad
		[\iota_{10},\nu_{10}^2]=\sigma_{10} \varepsilon_{17}=\sigma_{10} \bar\nu_{17},\\
		&[\iota_{10},\varepsilon_{10}]=[\iota_{10},\bar\nu_{10}]=0,\quad
		[\iota_{10},\mu_{10}]=2\sigma_{10} \zeta_{17}.\\
		&[\iota_{12},\eta_{12}]=\Sigma\theta',\quad[\iota_{12},\sigma_{12}]=\pm(\Sigma\xi'-2\xi_{12}).\\
		&[\iota_{14},\eta_{14}]=4\sigma_{14}^2,\quad[\iota_{14},\nu_{14}]=\pm2\omega_{14}.\\
		&[\iota_{16},\eta_{16}]\equiv \Sigma\eta^{*\prime}\mathbin{\mathrm{mod}} \Sigma^2\pi_{30}^{14},\quad
		[\iota_{16},\nu_{16}]=\pm(\Sigma^3\lambda-2\nu_{16}^*).\\
		&[\iota_{18},\eta_{18}]=2(\nu_{18}^*+\xi_{18}).
	\end{align*}

	\textbf{Standing Convention for the Reference Tables.}
	The \^{O}guchi material in Appendix~\ref{app:oguchi} and its composition table are convenient working
	tables. Throughout this book we frequently use their generator names and relations.

	\chapter{Boundaries,  Lifts, $J$-Images,  Diagrams for Lie Groups and So On}\label{chb}

	We often freely use the descriptions of the groups $\pi_{*}(S^{j})$
	given in Appendix~\ref{tlqdzj}, together with the composition
	relations among generators in Toda's range recorded in
	\^{O}guchi's table in Appendix~\ref{app:oguchi-relations}.
	
	Suppose that $\pi_k(\mathit{SO}(n))$ is known and we want to determine $\pi_k(\mathit{SO}(n+1)).$  The special orthogonal
	fibration separates this problem into the   elements   inherited from $\mathit{SO}(n)$ and the elements of homotopy groups of spheres that
	admit lifts. The exact sequence tells us which classes occur\uprightsemicolon{} relations among their representatives
	then determine how they fit together.
	
	For fixed degree $k,$   increasing $n$ gives the exact sequence
	\begin{equation}
		\begin{aligned}
			&\pi_{k+1}^n\stackrel{{\Delta}}{\longrightarrow}\pi_{k}(\mathit{SO}(n)) \stackrel{{(i_{n+1})_*}}{\longrightarrow}\pi_{k}(\mathit{SO}(n+1)) \stackrel{{(p_{n+1})_*}}{\longrightarrow}\pi_k^n \stackrel{{\Delta}}{\longrightarrow}\pi_{k-1}(\mathit{SO}(n)).\label{m2:eq:special orthogonal-sequence}
		\end{aligned}
	\end{equation}
	We refer to this sequence as sequence $(k)_{n}.$ 
	Its two boundary maps have different domains.  The first identifies
	which classes in $\pi_k(\mathit{SO}(n))$ vanish under inclusion into $\pi_k(\mathit{SO}(n+1))$\uprightsemicolon{} the second identifies
	which elements of homotopy groups of spheres lift.  Consequently, \eqref{m2:eq:special orthogonal-sequence}
	reduces the calculation to an extension
	\begin{equation}
		\begin{aligned}
			0&\longrightarrow \operatorname{coker}(\Delta\colon \pi_{k+1}^n\to \pi_{k}(\mathit{SO}(n))) \longrightarrow \pi_{k}(\mathit{SO}(n +1))\label{m2:eq:special orthogonal-extension} \\
			&\longrightarrow \operatorname{ker}(\Delta\colon \pi_k^n\to \pi_{k-1}(\mathit{SO}(n))) \longrightarrow0.
		\end{aligned}
	\end{equation}

	\section{$\Delta$-Images}\label{m2:subsec:boundaries}
	For convenience, we first summarize the images under $\Delta$ that will be discussed in this section. The results involving elements in the 7th stem are given in \cite[Lemma 1.2]{KM1} and \cite[(4.1)]{HKM}, and the results for $\q(\e_{4}),\;\q(\e_{15})$ are given in \cite{Ka}. All remaining entries in the table will be explained below.\par
	\begingroup
	\renewcommand{\arraystretch}{1.7}
	\setlength{\extrarowheight}{3pt}
	\setlength{\tabcolsep}{9pt}
	\setlength{\LTleft}{0pt}
	\setlength{\LTright}{0pt}
	\setlength{\LTcapwidth}{\textwidth}
	
	\begin{longtable}{@{}p{\textwidth}@{}}
		\caption{Values of the Connecting Homomorphism
			$\Delta\colon\pi_{\bullet+1}(S^{j})\to\pi_{\bullet}(SO(j))$}
		\label{m2:tab:boundary-values}\\
		\endfirsthead
		
		\caption[]{Values of the Connecting Homomorphism
			$\Delta$ (Continued)}\\
		\endhead
		
		\endfoot
		\endlastfoot
		
		\resizebox{\textwidth}{!}{%
			\begin{tabular}{|c|c|c|c|c|c|c|}
				\hline
				$\mathbf x$
				& $\iota_{4m+1}$
				& $\iota_{8k+3}$
				& $\e_4$
				& $\iota_5$
				& $\iota_7$
				& $\iota_{15}$
				\\
				\hline
				$\Delta(\mathbf x)$
				& $[\eta_{4m-1}]_{4m+1}$
				& $[\nu_{8k-1}]_{8k+3}$
				& $2[\e_3]-[\ca_2]_4$
				& $[\iota_3]_5\eta_3$
				& $0$
				& $[\bar{\nu}_6+\varepsilon_6]_{15}=[\mathrm{c}_{KM}(\bar{\nu}_6+\varepsilon_6)]_{15}$
				\\
				\hline
				$\mathbf x$
				& $\eta_{4m}$
				& $\eta_{4n+3}$
				& $\eta_{8k+2}$
				& $\eta_{4n+2}$
				& $\eta_{4n+2}^{2}$
				& $\eta_{8k+1}^{2}$
				\\
				\hline
				$\Delta(\mathbf x)$
				& $[\eta_{4m-2}^{\,2}]_{4m}$
				& $0$
				& $2[\nu_{8k-1}]_{8k+2}$
				& \textit{can be divisible by 2}
				& $0$
				& $4[\nu_{8k-1}]_{8k+1}$
				\\
				\hline
			\end{tabular}%
		}\\[14pt]
		
		\resizebox{\textwidth}{!}{%
			\begin{tabular}{|c|c|c|c|c|c|c|}
				\hline
				$\mathbf x$
				& $\eta_{10}$
				& $\nu_{8k-1}$
				& $\nu_4$
				& $\nu_5$
				& $\nu_6$
				& $\iota_9$
				\\
				\hline
				$\Delta(\mathbf x)$
				& $2[\iota_7]_{10}\nu_7$
				& $0$
				& $[\iota_3]\nu'\bmod [\eta_2]_4\nu'$
				& $0$
				& $2[\nu_5]$
				& $[\nu_5]_9+[\iota_7]_9\eta_7$
				\\
				\hline
				$\mathbf x$
				& \multicolumn{2}{c|}{$\s\sigma'$}
				& \multicolumn{2}{c|}{$\sigma_8$}
				& \multicolumn{2}{c|}{$\sigma_9$}
				\\
				\hline
				$\Delta(\mathbf x)$
				& \multicolumn{2}{c|}{$
					\begin{aligned}
						&2[\iota_7]\sigma'
						+4[\mathrm{c}_{\mathrm{KM}}(\bar{\nu}_6+\varepsilon_6)]_8\\
						&\quad+[ \nu_5]_8\nu_8^2
						-[\eta_5\varepsilon_6]_8
					\end{aligned}$}
				& \multicolumn{2}{c|}{$
					\begin{gathered}
						[\iota_7]\sigma'
						+c[\mathrm{c}_{\mathrm{KM}}(\bar{\nu}_6+\varepsilon_6)]_8\\
						\bmod\;[ \nu_5]_8\nu_8^2,\;
						[\eta_5\varepsilon_6]_8,\\
						\text{for some odd }c
					\end{gathered}$}
				& \multicolumn{2}{c|}{$
					\begin{aligned}
						&
						[\iota_7]_9(\bar{\nu}_7+
						\varepsilon_7)
					\end{aligned}$}
				\\[5pt]
				\hline
			\end{tabular}%
		}\\*[9pt]
		
		In this table, $m,n,k$ are integers satisfying
		$m\geq3,$ $n\geq2$ and $k\geq1.$ The number $\mathrm{c}_{\mathrm{KM}}\in \z$ is some odd integer due to Kachi and Mukai's revision for the lift of $\ch_6+\lt_6$ \cite{KM1}.
		\\
	\end{longtable}
	\endgroup
	\par
	We begin with the boundary formula that allows a suspended right factor to be carried through $\Delta$
	\cite[(2.2)]{MT2}. It is the basic reduction used when the element of a homotopy group of a sphere is given as a composite.
	\begin{equation}\notag
		\begin{aligned}
			&\Delta(\alpha\circ \Sigma\beta)=\Delta\alpha\circ\beta.
		\end{aligned}
	\end{equation}
	
	Suppose that we are given elements $\alpha\in\pi_k(S^n)$ and $\beta\in\pi_m(S^k)$ such
	that $\Delta\alpha=0.$  Since $\Delta(\alpha\beta)=
	\Delta({p_{n+1}}_*([\alpha]\beta))=0,$  
	we obtain
	\begin{equation}\notag
		\begin{aligned}
			&\Delta(\alpha\beta)=0\ (\Delta\alpha=0).
		\end{aligned}
	\end{equation}
	\noindent       
	By \cite[Theorem 4]{Wh42}, we have the relation
	\begin{equation}\notag
		\begin{aligned}
			&{p_n}_*\Delta\iota_n=(1+(-1)^n)\iota_{n-1}.
		\end{aligned}
	\end{equation}
	Therefore, for an
	element $\alpha\in\pi_k(S^n),$   we obtain  a formula
	\begingroup
	\renewcommand{\theHequation}{m1:DelS2}
	\begin{equation}
		\begin{aligned}
			&{p_{n+1}}_*\Delta(\Sigma\alpha)=0\ (n: \mbox{even})\uprightsemicolon{} \ {p_{n+1}}_*\Delta(\Sigma\alpha)=2\iota_n\circ\alpha\ (n: \mbox{odd}).\label{m1:DelS2}
		\end{aligned}
	\end{equation}
	\endgroup
	\noindent           For $n\geq1,$   the following formula applies whenever the Toda bracket on the left is defined \cite{MimuraOdaPeriodic}:
	\begin{equation}\notag
		\begin{aligned}
			&\Delta\{\alpha, \Sigma^n\beta, \Sigma^n\gamma\}_n\xyd\pm\{\Delta(\alpha),\Sigma^{n-1}\beta,\Sigma^{n-1}\gamma\}_{n-1}.
		\end{aligned}
	\end{equation}
	
	Recall that,  for an element $\alpha$ of a group, we denote by $\mathrm{ord}(\alpha)$ or       $\mathrm{ord}\,\alpha$       the order of $\alpha.$ 
	By \cite[(3.2)]{GM} and the fact that
	$2\pi_{4n+4}(\mathit{SO}(4n+4))=0$ \cite{Ke}, we deduce
	\begin{equation}\notag
		\begin{aligned}
			&\Delta(\eta_{4n+3})=0\uprightsemicolon{}\ \mathrm{ord}[\eta_{4n+3}]=2 \ (n\ge 2).
		\end{aligned}
	\end{equation}
	\noindent     
	In light of \cite[(4.2)]{GM} and the fact that $2\pi_{4n+4}(\mathit{SO}(4n+3))=0$ \cite{Ke}, we find
	\begin{equation}\notag
		\begin{aligned}
			&\Delta(\eta^2_{4n+2})=0\uprightsemicolon{}\ \mathrm{ord}[\eta^2_{4n+2}]=2.
		\end{aligned}
	\end{equation}
	
	\noindent       As a consequence of  \cite[(4.1), (4.4)]{GM}, we obtain
	\begin{equation}\notag
		\begin{aligned}
			&\Delta\iota_{4n+1} = [\eta_{4n-1}]_{4n+1}\ \mbox{for}\ n\geq 2.
		\end{aligned}
	\end{equation}
	
	\noindent       According to \cite[p. 161, Table]{Ke}, we know that $\pi_{4n+2}(\mathit{SO}(4n+2))\cong\mathbb Z/4$ for $n\ge 2.$ 
	Therefore, we conclude that
	\begin{equation}\notag
		\begin{aligned}
			&\Delta(\eta_{4n+2})\in 2\pi_{4n+2}(\mathit{SO}(4n+2))\ \mbox{for}\ n\ge 2.
		\end{aligned}
	\end{equation}
	\noindent     
	By \cite[p. 410, \mbox{the proof of} Lemma 3.5(1)]{GM}, we see that
	\begin{equation}\notag
		\begin{aligned}
			&\Delta(\nu_{8n-1})=0\uprightsemicolon{}\ \mathrm{ord}[\nu_{8n-1}]=8.
		\end{aligned}
	\end{equation}
	
	\noindent       In view of \cite[(4.1)-(4.5)]{GM}, we have
	\begin{equation}\notag
		\begin{aligned}
			&\Delta(\eta^2_{8n+1})=4[\nu_{8n-1}]_{8n+1},\ \Delta(\eta_{8n+2})=2[\nu_{8n-1}]_{8n+2},
		\end{aligned}
	\end{equation}
	\begin{equation}\notag
		\begin{aligned}
			&\Delta(\iota_{8n+3})=[\nu_{8n-1}]_{8n+3}.
		\end{aligned}
	\end{equation}
	
	\noindent       Recall from \cite[(1.8)]{HKM} that
	\begin{equation}\notag
		\begin{aligned}
			&\Delta\eta_{4n+4}=[\eta^2_{4n+2}]_{4n+4}\ \mbox{for}\ n\geq 2.
		\end{aligned}
	\end{equation}

	\section{The Constructions of  Lifts by Using $\e_{n}$}

	\begin{prop}\label{mdsyq}Let $n\geq2.$  With respect to the fiber sequence $$S^{n}\xleftarrow{p_{n+1}}\mathit{SO}(n+1)\xleftarrow{i_{n+1}}\mathit{SO}(n)\xleftarrow{\partial}\lo S^{n},$$ we have
		$$\{p_{n+1},i_{n+1},\Delta\iota_n\}\ni\pm\iota_n.$$
	\end{prop}
	\begin{proof} Let $\mathrm{i}_{E}\colon \lo S^{n}\leftarrow S^{n-1}$  be the classical inclusion.   Employing Proposition \ref{fby1}, we infer
		$$\{p_{n+1},i_{n+1},\partial\}\ni ( S^{n}\xleftarrow{\mathrm{ev}}\s\lo S^{n}).$$
		\noindent Then, \vspace{-1.6em}\begin{align}
			\{p_{n+1},i_{n+1},\Delta\iota_n\}
			&= \{p_{n+1},i_{n+1},\partial\hc(\mp\mathrm{i}_{E})\} \notag\\
			&\dyd \{p_{n+1},i_{n+1},\partial\}\hc(\pm\s\,\mathrm{i}_{E}) \notag\\
			&\ni \bigl(S^{n}\xleftarrow{\mathrm{ev}}\s\lo S^{n}\bigr)
			\hc\bigl(\s\lo S^{n}\xleftarrow{\pm\s\,\mathrm{i}_{E}}S^{n}\bigr) \notag\\
			&= \pm\e_{n}. \notag
		\end{align}
		\noindent   Hence the result holds.
	\end{proof}
	
	Using the above proposition, we are able to show the following lemma.
	\begin{lem}\label{m1:Dn4sg'}
		Suppose that $\alpha\in\pi_{n-2+k}^{n-2}$ satisfies
		$\Delta(\Sigma\alpha)=0.$  Then $$
		[\Sigma\alpha]\in{}
		\pm\{i_n,\Delta\iota_{n-1},\alpha\}\bmod\;
		i_{n*}\pi_{n-1+k}(\mathit{SO}(n-1))
		+\pi_{n-1}(\mathit{SO}(n))\circ\Sigma\alpha.
		$$\noindent In particular,$$
		\Delta\iota_{2n}\in{}
		\pm\{i_{2n},\Delta\iota_{2n-1},2\iota_{2n-2}\}\bmod\;
		i_{2n*}\pi_{2n-1}(\mathit{SO}(2n-1))
		+\pi_{2n-1}(\mathit{SO}(2n))\circ2\iota_{2n-1}.$$ 
	\end{lem}
	
	\begin{proof} Observe that
		\noindent
		$$\begin{aligned}[t]
			p_n\circ\bigl(-\{i_n,\Delta\iota_{n-1},\alpha\}\bigr)
			&=\{p_n,i_n,\Delta\iota_{n-1}\}\circ\Sigma\alpha\\
			&\ni\pm\Sigma\alpha
			\bmod\;
			p_n\circ\pi_{n-1}(\mathit{SO}(n))\circ\Sigma\alpha.
		\end{aligned}$$
		\noindent       Together with $p_{n*}[\Sigma\alpha]=\Sigma\alpha,$  
		the exact sequence  induced by the fibration
		$$\mathit{SO}(n-1)\to\mathit{SO}(n)\to S^{n-1}$$ gives $
		[\Sigma\alpha]\in{}
		\pm\{i_n,\Delta\iota_{n-1},\alpha\}
		\bmod\;
		i_{n*}\pi_{n-1+k}(\mathit{SO}(n-1))
		+\pi_{n-1}(\mathit{SO}(n))\circ\Sigma\alpha.
		$   The second assertion follows by taking
		$\alpha=2\iota_{2n-2}$ and using
		$2\iota_{2n-1}=p_{*}(\Delta\iota_{2n}).$ 
	\end{proof}
	
	\section{$J$-Images}\label{m2:subsec:j-toda}
	
	The $J$-homomorphism converts elements of homotopy groups of special orthogonal groups into elements of homotopy groups of spheres. Its compatibility with
	composition and with the boundary map lets us compare the two kinds of calculation. The signs relating $J$ to suspension and inclusion are retained explicitly below.
	Signs in related formulas are chosen consistently with the same maps and lifts.
	
	There is a canonical homomorphism $J\colon\pi_k(\mathit{SO}(n))\to \pi_{k+n}(S^n)$ defined by  Whitehead \cite{Wh}. We work integrally and regard an $n\times n$ matrix in $\mathit{SO}(n)$ as a linear map $\mathbb R^{n}\to\mathbb R^{n}.$ 
	The homomorphism $J$ is induced by the composition
	\[
	\begin{aligned}
		\mathit{SO}(n)\xyd \mathrm{Map}_{*}(\mathbb R^n\cup\n\wq\nn,\mathbb R^n\cup\n\wq\nn)_{1}=\mathrm{Map}_{*}(S^n,S^n)_{1}\xrightarrow{\simeq}\mathrm{Map}_{*}(S^n,S^n)_{0}
		=(\lo^n S^n)_0\xyd\lo^n S^n.
	\end{aligned}
	\]
	\noindent Here $\mathrm{Map}_{*}(S^n,S^n)_i$ denotes the component of based maps of degree $i.$  The homotopy equivalence is given by translation from the degree-$1$ component to the degree-$0$ component of the loop space.

	Now, work 2-locally. For $\alpha\in\pi_{m}(\mathit{SO}(n))$ and $\beta\in\pi_{k}(S^{m}),$   we have a formula
	\[\begin{aligned}
		&J(\alpha\beta )=J(\alpha )\Sigma^n\beta.
	\end{aligned}\]

	\noindent       The following diagram commutes up to sign: \[
	\begin{tikzcd}[column sep=large, row sep=large]
		\pi_k(\mathit{SO}(n))
		\arrow[r, "i_*"]
		\arrow[d, "J"']
		&
		\pi_k(\mathit{SO}(n+1))
		\arrow[r, "p_*"]
		\arrow[d, "J"']
		&
		\pi^n_k
		\arrow[r, "\Delta"]
		\arrow[d, "\Sigma^{n+1}"']
		&
		\pi_{k-1}(\mathit{SO}(n))
		\arrow[d, "J"]
		\\
		\pi^n_{k+n}
		\arrow[r, "\Sigma"']
		&
		\pi^{n+1}_{k+n+1}
		\arrow[r, "H"']
		&
		\pi^{2n+1}_{k+n+1}
		\arrow[r, "P"']
		&
		\pi^n_{k+n-1},
	\end{tikzcd}
	\]
	
	\noindent   where the lower sequence is the EHP sequence. This diagram is also a direct
	consequence of \cite[(11.4)]{T}. In particular, 
	\[
	\begin{aligned}
		&\Sigma J(\alpha)=\pm J(i_{n+1}\circ\alpha)
		\quad\text{for }\alpha\in\pi_k(\mathit{SO}(n)),\\
		&J\Delta(\gamma)=\pm[\gamma,\iota_n]
		\quad\text{for }\gamma\in\pi_k(S^n),\\
		&HJ(\beta)=\Sigma^n(p_n\circ\beta)
		\quad\text{for }\beta\in\pi_k(\mathit{SO}(n)).
	\end{aligned}
	\]
	\noindent 
	Note that $[\iota_n,[\iota_n,\iota_n]]=0.$  Hence, for $\gamma\in\pi_{k}^{n},$  
	by the formula \cite[Corollary (7.4)]{BB}, we have the relation
	\begin{equation}\notag
		\begin{aligned}
			&J\Delta(\gamma)=\pm[\iota_{n},\gamma] =\pm[\iota_n,\iota_n]\circ \Sigma^{n-1}\gamma \in\pi_{k+n-1}^{n}.\label{m1:JD}
		\end{aligned}
	\end{equation}

	Recall the relations from \cite[Lemma 1.1(ii)]{KM1}:
	\begingroup
	\renewcommand{\theHequation}{\ifnum\value{equation}=4 m1:Ji7\else\theHchapter.\arabic{equation}\fi}
	\begin{gather}
		\label{m1:Ji7}
		J[\iota_7]=\sigma_8, \\
		J[ \nu_5]=\bar\nu_6+\varepsilon_6,\notag
		\\
		\label{m1:Jet5ep}
		J[\eta_5\varepsilon_6]=-\sigma''\sigma_{13}+\bar{\nu}_6\nu^2_{14},\\
		\label{m1:Jeta6}
		J[\eta_6]=\sigma', \\
		\label{m1:Jbn6e6}
		J[\mathrm{c}_{\mathrm{KM}}(\bar\nu_6+\varepsilon_6)]=\sigma'\sigma_{14},\\
		\label{m1:let6sg}
		[\eta_6]\sigma'=4[\mathrm{c}_{\mathrm{KM}}(\bar\nu_6+\varepsilon_6)]+[ \nu_5]_7\nu^2_8+[\eta_5\varepsilon_6]_7.
	\end{gather}
	\endgroup
	The $\iota_{15}$ entry in Table~\ref{m2:tab:boundary-values} uses the same
	Kachi--Mukai choice. By \cite[Tables 2 and 3]{Ka},
	$\pi_{14}(\mathit{SO}(15))=\mathrm{span}\{\Delta\iota_{15}\}\cong\mathbb Z/2.$
	On the other hand, \eqref{m1:Jbn6e6} and Appendix~\ref{tlqdzj} give
	$J([\mathrm{c}_{\mathrm{KM}}(\bar{\nu}_6+\varepsilon_6)]_{15})
	=\pm2\sigma_{15}^2\neq0.$ Hence
	$\Delta\iota_{15}=[\mathrm{c}_{\mathrm{KM}}(\bar{\nu}_6+\varepsilon_6)]_{15}.$
	
	\section{Useful Diagrams of Lie Groups}\label{lqzh}
	
	In this section, we work integrally until Diagram~(\ref{tu1029}).
	
	There are well-known homeomorphisms  $$\mathit{SO}(4)\cong \mathit{SO}(3)\times S^{3},\quad \mathit{SO}(8)\cong \mathit{SO}(7)\times S^{7},$$
	$$\mathit{Spin}(4)\cong \mathit{Spin}(3)\times S^{3}, \quad\mathit{Spin}(8)\cong \mathit{Spin}(7)\times S^{7}$$
	\noindent of topological spaces.\footnote{The homeomorphism for $\mathit{SO}(4)$ or $\mathit{SO}(8)$ is not a homeomorphism of Lie groups. Neither of these two is  an $H$-space equivalence. See \cite[Chapter 19]{Hand}. Of course, this does not affect the calculation of homotopy groups.} Note that $\mathit{Spin}(3)=S^{3}.$ Moreover, $\mathit{SO}(3)\cong \mathbb RP^3,$   whose universal covering is $S^{3}.$ This tells us $$\pi_{*>1}(\mathit{SO}(3))\cong \pi_{*>1}(S^{3}).$$
	\noindent
	More precisely, we have
	$$\pi_{*>1}(\mathit{SO}(3))=[\ca_2]\hc\pi_{*>1}(S^3),$$ \noindent and
	$$\pi_*(\mathit{SO}(4))
	=i_{3,4}\hc\pi_*(\mathit{SO}(3))
	\jia[\e_3]\hc\pi_*(S^3)$$ induced by $\mathit{SO}(4)\cong \mathit{SO}(3)\times S^{3}$
	\cite{Ka}. An analogous description for $\mathit{SO}(8)$ is also given in \cite{Ka}.

	For small values of $n,$   we use fibrations involving classical and exceptional Lie groups.  The following diagram summarizes the fibrations listed
	in \cite{Mi2}.  In the following diagram, except  for  $\mathit{Spin}(7)\stackrel{{S^{15}}}{\longrightarrow}\mathit{Spin}(9),$  
	an arrow $$\begin{aligned}
		&G\stackrel{{Y}}{\longrightarrow}X
	\end{aligned}$$ means a principal $G$-bundle
	$$\begin{aligned}
		&G\to X\to Y,
	\end{aligned}$$ with the orbit space $X/G\cong Y.$  Here $V_{7,2}=V_{2}(\mathbb R^7)$ is the Stiefel manifold of orthonormal
	$2$-frames in $\mathbb{R}^7,$   and $\mathbb{O}P^2=F_4/\mathit{Spin}(9)$
	is the Cayley projective plane. The symbol  $   \mathit{Spin}(7)\stackrel{{S^{15}}}{\longrightarrow}\mathit{Spin}(9)$ means we have a fibration $$\mathit{Spin}(7)\to\mathit{Spin}(9)\to S^{15}.$$
	The covering
	$\mathit{Spin}(n)\to\mathit{SO}(n)$ induces an isomorphism on
	$\pi_k$ for $k\geq2,$   so these fibrations can be used directly in
	the degrees considered here.
	
	\par\smallskip\noindent
	\resizebox{\linewidth}{!}{
		\begin{tikzpicture}[
			>=Stealth,
			every node/.style={font=\large},
			every path/.style={semithick}
			]
			\node (S3) at (0,0)
			{$S^3=\mathit{Sp}(1)=\mathit{SU}(2)=\mathit{Spin}(3)$};
			\node (SU3) at (4,1.7)
			{$\mathit{SU}(3)$};
			\node (Sp2) at (4,-1.7)
			{$\mathit{Sp}(2)=\mathit{Spin}(5)$};
			\node (SU4) at (7,0)
			{$\mathit{SU}(4)=\mathit{Spin}(6)$};
			\node (G2) at (10,1.7)
			{$G_2$};
			\node (Spin7) at (11,0)
			{$\mathit{Spin}(7)$};
			\node (Spin9) at (14,0)
			{$\mathit{Spin}(9)$};
			\node (F4) at (14,-1.7)
			{$F_4$};
			\draw[->]
			(S3.north east)
			-- node[midway,above left] {$S^5$}
			(SU3.south west);
			\draw[->]
			(SU3.south east)
			-- node[midway,above right] {$S^7$}
			(SU4.north west);
			\draw[->]
			(S3.south east)
			-- node[midway,left=3pt] {$S^7$}
			([xshift=-6mm]Sp2.north);
			\draw[->]
			([xshift=6mm]Sp2.north)
			-- node[midway,right=3pt] {$S^5$}
			(SU4.south west);
			\draw[->]
			(SU3.east)
			-- node[above] {$S^6$}
			(G2.west);
			\draw[->]
			(G2.south east)
			-- node[right] {$S^7$}
			(Spin7.north);
			\draw[->]
			(SU4.east)
			-- node[below] {$S^6$}
			(Spin7.west);
			\draw[->]
			(Spin7.east)
			-- node[above] {$S^{15}$}
			(Spin9.west);
			\draw[->]
			(Spin9.south)
			-- node[right] {$\mathbb{O}P^2$}
			(F4.north);
			\draw[->]
			(S3.north) -- (0,2.8)
			-- node[above] {$V_{7,2}$} (9.4,2.8)
			-- (G2.north west);
		\end{tikzpicture}
	}
	\par\smallskip

	It should be pointed out that the fibration $\mathit{Spin}(7)\to\mathit{Spin}(9)\to S^{15}$
	is not resulting from the classical inclusion of Lie groups $\mathit{Spin}(7)\xyd\mathit{Spin}(9).$  (If so, we would have $\mathit{Spin}(9)/\mathit{Spin}(7)=V_{2}(\mathbb R^{9})\simeq S^{15}.$  It is impossible because $H_{7}(V_{2}(\mathbb R^{9});\z/2)=\z/2.$)

	Let $\mathcal{T}:G_2\to\mathit{SO}(7)$ be the composite
	$G_2\hookrightarrow\mathit{Spin}(7)\to\mathit{SO}(7),$  
	where the second map is the universal covering.  Under the canonical
	inclusions into $\mathit{Spin}(7),$   the intersection
	$G_2\cap\mathit{Spin}(6)$ is $\mathit{SU}(3)\uprightsemicolon{}$ equivalently,
	$G_2\times_{\mathit{Spin}(7)}\mathit{SU}(4)\cong\mathit{SU}(3)$
	\cite{A}.  Passing to the special orthogonal groups gives the following
	morphism of fibrations:
	\begin{equation}\label{m2:Y}
		\begin{tikzcd}[column sep=large,row sep=large]
			\mathit{SU}(3) \arrow[r,"i_G"] \arrow[d,"r_3"'] &
			G_2 \arrow[r,"p_G"] \arrow[d,"\mathcal{T}"] & S^6 \arrow[d,equal] \\
			\mathit{SO}(6) \arrow[r,"i_7"'] &
			\mathit{SO}(7) \arrow[r,"p_7"'] & S^6.
		\end{tikzcd}
	\end{equation}
	Similarly, the classical isomorphism 
	$\mathit{Spin}(5)\times_{S^{3}}\mathit{SU}(3)\cong\mathit{Spin}(6)$ \cite{A} together with the covering maps gives a morphism of fibrations:
	\begin{equation}\label{tu1029}
		\begin{tikzcd}
			\mathit{Spin}(3) \arrow[d] \arrow[r]
			& \mathit{SU}(3) \arrow[d, "\mathcal{D}"] \arrow[r]
			& S^{5} \arrow[d, equal] \\ \mathit{SO}(5) \arrow[r]
			& \mathit{SO}(6) \arrow[r]
			& S^{5}.
		\end{tikzcd}
	\end{equation}
	
	\noindent   We denote the connecting homomorphism of the top row of Diagram~\eqref{m2:Y} by
	$\Delta':\pi_{k+1}^6\to U_k^3\uprightsemicolon{}$ thus
	$(r_3)_*\Delta'=\Delta.$
	For the lift $\langle\bar{\nu}_6+\varepsilon_6\rangle\in\pi_{14}(G_2)$
	of \cite[Theorem 6.1]{Mi2}, Kachi \cite{Ka} defines
	$$[\bar{\nu}_6+\varepsilon_6]
	:=\mathcal{T}_*\langle\bar{\nu}_6+\varepsilon_6\rangle\in\pi_{14}(  \mathit{SO}(7) )$$
	and
	$[\bar{\nu}_6+\varepsilon_6]_{8+n}
	:=(i_{7,8+n})_*[\bar{\nu}_6+\varepsilon_6]$ for $n\geq0.$
	Diagram~\eqref{m2:Y} gives
	$(p_7)_*[\bar{\nu}_6+\varepsilon_6]=\bar{\nu}_6+\varepsilon_6.$
	In general, this  is different from Kachi and Mukai's revised element $[\bar{\nu}_6+\varepsilon_6]$ in \cite{KM1}. 
	For the revised element, we only have $$p_{*}[\bar{\nu}_6+\varepsilon_6]\ty\bar{\nu}_6+\varepsilon_6\m2(\bar{\nu}_6+\varepsilon_6).$$
	Recall that, to clarify, we use
	$[\mathrm{c}_{\mathrm{KM}}(\bar{\nu}_6+\varepsilon_6)]$
	to denote the  revised element $[\bar{\nu}_6+\varepsilon_6]$ where $\mathrm{c}_{\mathrm{KM}}\in\z$ is some odd number.
	
	The map
	$i_U:S^3=\mathit{SU}(2)\to\mathit{SU}(3)$ used later is the
	standard inclusion in the fibration over $S^5.$

	To transfer the calculations for $G_2$ to $\mathit{SO}(7),$ we use the following decomposition on
	homotopy groups. It separates the contribution induced by $\mathcal{T}$ from the contribution obtained by
	composition with $[\eta_6].$
	In fact, in view of $\mathit{Spin}(7)/G_2=S^{7}$ and
	$\pi_{7}(BG_{2})\cong\pi_{6}(G_{2})=0$ \,\cite{Mi2}, we have a (2-local) splitting $$\mathit{Spin}(7)\simeq G_2\times S^{7}.$$  The space
	$\mathit{Spin}(7)$ is the universal covering of $\mathit{SO}(7) ,$ leading to the following.
	\begin{prop}\label{so7fj}
		There are decompositions:
		\[ \begin{aligned}
			&\pi_{k}(\mathit{SO}(7)) =\mathcal{T}_{*}\pi_k (G_2) \oplus [\eta_6]\circ\pi^7_k,
		\end{aligned} \] for all $k\ge 2.$
		\qed
	\end{prop}

	As is well known, $\mathit{Sp}(2)\cong\mathit{Spin}(5).$ 
	The  covering map
	$\mathcal{C}:\mathit{Sp}(2)\to\mathit{SO}(5)$
	induces the Hopf map
	$\mathcal{H}:S^7=\mathit{Sp}(2)/\mathit{Sp}(1)
	\to\mathit{SO}(5)/\mathit{SO}(4)=S^4$ as can be checked on the $7$-th homotopy groups.
	We need the following result.
	
	\begin{prop}[{\cite[Lemma~1.2]{HKM}}]\label{sp2so5gx}
		Let $n\geq2.$  The fibrations
		$S^3\to\mathit{Sp}(2)\to S^7$ and
		$\mathit{SO}(4)\to\mathit{SO}(5)\to S^4$
		give rise to the following commutative diagram
		of exact sequences:
		\[
		\begin{tikzcd}[column sep=2.2em, row sep=3em]
			\pi_n(S^3)
			\arrow[r, "i_\ast"]
			\arrow[d, "{[\iota_3]_\ast}"']
			&
			\pi_n(\mathit{Sp}(2))
			\arrow[r, "p_\ast"]
			\arrow[d, "\mathcal{C}_\ast"', "\cong"]
			&
			\pi_n(S^7)
			\arrow[r, "\Delta_{\mathbb H}"]
			\arrow[d, "\mathcal{H}_\ast"']
			&
			\pi_{n-1}(S^3)
			\arrow[d, "{[\iota_3]_\ast}"]
			\\
			\pi_n(\mathit{SO}(4))
			\arrow[r, "i_\ast"]
			&
			\pi_n(\mathit{SO}(5))
			\arrow[r, "p_\ast"]
			&
			\pi_n(S^4)
			\arrow[r, "\Delta"]
			&
			\pi_{n-1}(\mathit{SO}(4)).
		\end{tikzcd}
		\]
		Here $\Delta_{\mathbb H}$ and $\Delta$ are the boundary homomorphisms
		of the first and second fibrations, respectively. Note that $\mathcal {H}\ty \nu_{4}\m2\nu_{4},\s \nu'.$
		\qed
	\end{prop}

	\section{Bott's and James' Theorems}\label{m3:sowq}
	
	As $n$ increases, the special orthogonal groups approach the stable groups discussed in
	Chapter~\ref{ch:history}. The stable $J$-homomorphism supplies a corresponding comparison with elements of the stable homotopy groups of spheres. We collect this information here for use at the upper end of the degree calculations.
	
	Working integrally, for every $n\in\mathbb{Z}_{+},$ we have a principal $\mathit{SO}(n)$-bundle
	$\mathit{SO}(n) \hookrightarrow\mathit{O}(n)\rightarrow S^{0}$;
	$\mathit{O}:=\mathop{\mathrm{colim}}\limits_{ n } \mathit{O}(n),$
	$\mathit{SO}:=\mathop{\mathrm{colim}}\limits_{ n } \mathit{SO}(n),\allowbreak \,\Omega^{\infty}S^{\infty}:=\mathop{\mathrm{colim}}\limits_{ n }\Omega^{n}S^{n}$; 
	there are also natural isomorphisms $\pi_{m}^{S}(S^{0})\cong\pi_{m}(\Omega^{\infty}S^{\infty}).$ Work 2-locally again. 
	The following proposition is given in \cite[Chapter 19]{Hand}.
	\needspace{3\baselineskip}\begin{prop}
		\leavevmode\par
		\begin{itemize}
			\item[\textup{(1)}] $B\mathit{O}\simeq(\Omega^{7}\mathit{SO})_0,$ where $B(-)$ denotes the classifying space
			functor and the subscript $0$ denotes the component of the constant loop.
			\item[\textup{(2)}] If $n\geq k+2$ and $k\geq1,$ then
			\[\pi _{k}(\mathit{SO}(n))\cong\pi _{k}(\mathit{SO})\cong\begin{cases}
				\mathbb Z_{(2)}, &\text{$k \equiv3,\allowbreak 7\;\bmod\;8$} \\
				\mathbb Z/2 , &\text{$k \equiv0,\allowbreak 1\;\bmod\;8$} \\
				0, &\text{else}\end{cases}.\]
		\end{itemize}
		\qed
	\end{prop}
	We refer to the natural isomorphism $[\Sigma^{k+2}X,\allowbreak Y]\cong [\Sigma^{2}X,\allowbreak \Omega^{k}Y],$ for
	$k\geq0,$ as the adjunction isomorphism. The following proposition follows from \cite[Theorem
	14.4]{Ja76}.
	\needspace{3\baselineskip}\begin{prop}\label{jmstm}
		[James' Theorem] For each $n\geq1,$ there is a homotopy-commutative diagram whose
		vertical maps induce the $J$-homomorphisms under adjunction:
		\[
		\begin{tikzcd}
			\mathit{SO}(n) \arrow[hookrightarrow]{r} \arrow[hookrightarrow]{d} & \mathit{SO}
			\arrow[hookrightarrow]{d}
			\\
			(\Omega^{n}S^{n})_{1} \arrow[hookrightarrow]{r} & (\Omega^{\infty}S^{\infty})_{1}.
		\end{tikzcd}.
		\]
		
		\noindent   Here, $(-)_{k}$ indicates the path-connected component of  degree $k$  maps\uprightsemicolon{} note that  $ (\Omega^{n}S^{n})_{1}\simeq     (\Omega^{n}S^{n})_{0}$ for all $n\geq1.$
		These diagrams induce isomorphisms if $r\leq 2n-2$:
		\[ \begin{aligned}
			&\pi_{r}(\mathit{SO},\allowbreak \mathit{SO}(n)) \;\cong\; \pi_{r}(\Omega^{\infty}S^{\infty},\allowbreak \Omega^{n}S^{n}) .
		\end{aligned} \]
		\qed
	\end{prop}

		\chapter{Classical Methods for Computing Homotopy Groups of Multi-Cell Complexes}\label{grdl1102}

	We work integrally throughout this section; no localization at
	a prime $p$ is assumed.
	
	If a space $X$ does not fit into a suitable fibration, such as the 2-local $EHP$-fibration $S^{n}\to\lo S^{n+1}\to\lo S^{2n+1}$ and the familiar fibrations of Lie groups  presented in Section~\ref{lqzh}, a different approach is needed to study its unstable homotopy groups.
	
	\section{The Homotopy Fiber of  a Pinch Map}
	
	Gray's relative James construction $J(X,A)$, denoted by
	$(X,A)_{\infty}$ in Gray's notation, provides a model for the
	homotopy fiber of the pinch map
	$X\cup CA\longrightarrow (X\cup CA)/X=\Sigma A$.
	
	When applying the relative James construction, we denote the
	mapping cylinder of a map $X\longrightarrow Y$ by $M_Y$.
	
	\begin{prop}[{\cite{Gray,hp2}}]\label{grdl}
		Let $X$ be a path-connected CW complex and let $A$ be a
		path-connected subcomplex of $X$.
		The inclusion $i\colon A\hookrightarrow X$ determines an
		inclusion of pairs $i\colon (A,A)\hookrightarrow (X,A)$.
		Then the following statements hold.
		
		\begin{enumerate}
			\item\label{grdl1a}
			$J(A,A)=J(A)$, where $J(A)$ is the ordinary James
			construction. There is a fiber sequence
			\[
			\begin{tikzcd}
				J(A)
				\arrow[r,hook,"J(i)"]
				& J(X,A)
				\arrow[r]
				& X\cup_i CA
				\arrow[r,"p"]
				& \Sigma A,
			\end{tikzcd}
			\]
			where $p$ is the pinch map and $J(i)$ is the inclusion
			induced by $i\colon (A,A)\hookrightarrow (X,A)$.
			
			Moreover, if
			\[
			K\xrightarrow{g}L\xrightarrow{j_L}C_g
			\xrightarrow{q}\Sigma K
			\]
			is a cofiber sequence of path-connected CW complexes
			$K$ and $L$, then there is a fiber sequence
			\[
			\begin{tikzcd}[column sep=0.8cm]
				J(K)
				\arrow[r,hook,"J(i_g)"]
				& J(M_L,K)
				\arrow[r]
				& C_g
				\arrow[r,"q"]
				& \Sigma K,
			\end{tikzcd}
			\]
			where $M_L$ is the mapping cylinder of $g$ and
			$J(i_g)$ is the inclusion induced by
			$i_g\colon (K,K)\hookrightarrow (M_L,K)$.
			
			\item\label{grdl2a}
			Suppose that $X/A$ is path-connected and that
			$\Omega\Sigma A$ and $X/A$ are CW complexes of finite
			type. For example, these conditions hold when $X$ is
			a finite simply connected CW complex and $A$ is a
			skeleton of $X$. Then, for any field $\mathbbm{k}$,
			\[
			\widetilde{H}_*(J(X,A);\mathbbm{k})
			\cong
			\widetilde{H}_*(X;\mathbbm{k})
			\otimes H_*(\Omega\Sigma A;\mathbbm{k}).
			\]
			
			\item\label{grdl3a}
			There is a filtration with $J_0(X,A)=*$ and
			$J_1(X,A)=X$ such that
			\[
			\begin{aligned}
				&J_0(X,A)\subseteq J_1(X,A)\subseteq
				J_2(X,A)\subseteq\cdots,\quad J(X,A)=\bigcup_{m\geq0}J_m(X,A).
			\end{aligned}
			\]
			
			\item
			If $X=\Sigma X'$ and $A=\Sigma A'$, then
			\[
			J_2(X,A)
			=X\cup_{[\mathrm{id}_X,i]}C(X\wedge A'),
			\]
			where $[\mathrm{id}_X,i]$ is the Whitehead product.
			
			\item\label{grdl5a}
			We have
			\[
			\begin{aligned}
				&J_m(X,A)/J_{m-1}(X,A)
				=X\wedge A^{\wedge(m-1)},\\
				&\Sigma J(X,A)
				\simeq
				\Sigma\bigvee_{m\geq0}X\wedge A^{\wedge m}.
			\end{aligned}
			\]
			
			\item
			There are relative James--Hopf invariants
			\[
			H_m\colon J(X,A)\longrightarrow
			J\bigl(X\wedge A^{\wedge(m-1)}\bigr),
			\]
			which are natural with respect to pairs.
		\end{enumerate}
		\qed
	\end{prop}
	
	The following lemma is useful for studying the skeletons of
	homotopy fibers. It follows directly from
	Proposition~\ref{grdl}\textup{\ref{grdl5a}}.
	
	\begin{lem}\label{j2}
		Let $m\geq2$ be an integer and let
		$X=Y\cup_g e^{m+1}$, where $g\colon S^m\longrightarrow Y$
		and $Y$ is a CW complex satisfying the following conditions:
		\begin{enumerate}[label=\textup{(\roman*)}]
			\item
			$Y$ is simply connected and not contractible;
			
			\item
			$Y$ is of finite type, and its cell structure is minimal,
			that is, consistent with its homology;
			
			\item
			$\dim(Y)\leq m$.
		\end{enumerate}
		Let $r\geq2$ be the dimension of the bottom cells of $Y$.
		For the homotopy fiber $J(M_Y,S^m)$ of the pinch map
		$X\longrightarrow S^{m+1}$, we have
		\[
		J_i(M_Y,S^m)
		=\mathrm{sk}_j\bigl(J(M_Y,S^m)\bigr),
		\qquad j=mi+r-1.
		\]
		\qed
	\end{lem}
	
	\begin{exercise}\label{nbdlx}Work 2-locally.  \begin{itemize}
			\item [\rm (1)]Let $P^{6}(2)=S^5\cup_2 e^{6}.$ Show that the homotopy fiber $F_{q}$ of the pinch map $q\colon P^{6}(2)\to S^{6}$ has a CW decomposition $$F_q=S^{5}\cup e^{10}\cup e^{15}\cup \cdots \cup e^{5i}\cup \cdots,$$
			and that $\sk_{14}(F_q)=S^{5}\cup_{[\e_5,2\e_5]} e^{10}=S^{5}\cup_{0} e^{10}=S^{5}\vee S^{10}$ up to homotopy.
			\item [\rm (2)]  Follow     Appendix~\ref{app:oguchi} and Appendix~\ref{tlqdzj}. Let $Z=S^{4}\cup_{\lt_{4}}e^{13}$ be the homotopy cofiber of $\lt_{4}\in\pi_{12}(S^{4}).$ Show that the homotopy fiber $F_{q'}$ of the pinch map $q'\colon Z\to S^{13}$ has a CW decomposition $$F_{q'}=S^{4}\cup e^{16}\cup e^{28}\cup \cdots \cup e^{4+12i}\cup \cdots,$$
			and that $\sk_{26}(F_{q'})=S^{4}\cup_{\s\nu'\hc\lt_7} e^{16}\not\simeq S^{4}\vee S^{16}.$
		\end{itemize}
		Hint: $F_{q}=J(M_{S^5}, S^{5})$,  $F_{q'}=J(M_{S^4}, S^{12})$, $[\e_{4},\lt_{4}]=\pm[\e_{4},\e_{4}]\hc\lt_{7}=\pm(2\nu_4-\s\nu')\hc\lt_{7}.$
	\end{exercise}
	
	\section{On Computing Homotopy Groups of Multi-Cell Complexes}

	We explain the philosophy behind the classical methods for calculations of $\pi_{\bullet}(Y\cup_g CX)$ when $Y\cup_g CX$ does not fit into a nice fibration such as the well-known fibrations of Lie groups given in  Section \ref{lqzh}.
	
	\begin{rem}\label{jssq}
		We aim to determine
		$\pi_{\bullet}(C_g)=\pi_{\bullet}(Y\cup_g CX)$,
		where $g\colon X\longrightarrow Y$ is a map between
		simply connected spaces. There is a cofiber sequence
		\[
		X\xrightarrow{g}Y
		\xrightarrow{j_Y}C_g
		\xrightarrow{q}\Sigma X.
		\]
		The homotopy fiber of $q$ is Gray's relative James
		construction, namely $F_q=J(M_Y,X)$.
		The fiber sequence
		$F_q\xrightarrow{i}C_g\xrightarrow{q}\Sigma X$
		gives an exact sequence
		\[
		\pi_{\bullet+1}(\Sigma X)
		\xrightarrow{\partial}\pi_{\bullet}(F_q)
		\xrightarrow{i_*}\pi_{\bullet}(C_g)
		\xrightarrow{q_*}\pi_{\bullet}(\Sigma X)
		\xrightarrow{\partial}\pi_{\bullet-1}(F_q).
		\]
		After determining the images of $\partial$, we obtain
		a short exact sequence
		\[
		0\longrightarrow\operatorname{im}(i_*)
		\hookrightarrow\pi_{\bullet}(C_g)
		\xrightarrow{q_*}\ker(\partial)
		\longrightarrow0,
		\]
		where
		$\operatorname{im}(i_*)
		\cong
		\operatorname{coker}\bigl(
		\partial\colon\pi_{\bullet+1}(\Sigma X)
		\longrightarrow\pi_{\bullet}(F_q)
		\bigr).$ Thus the problem reduces to determining the images of
		$\partial$ and resolving the group extension associated
		with the generators of $\ker(\partial)$.
		
		\begin{enumerate}
			\item
			To determine the images of $\partial$, we use
			\[
			\begin{aligned}
				\partial(\Sigma\mathbf{x})
				=i_Y\circ g\circ\mathbf{x},\quad
				\partial(\mathbf{x}\circ\Sigma\mathbf{y})
				=\partial(\mathbf{x})\circ\mathbf{y},
			\end{aligned}
			\]
			together with
			\[
			\begin{aligned}
				\partial
				\{\mathbf{x}_1,\Sigma^{k+1}\mathbf{x}_2,
				\Sigma^{k+1}\mathbf{x}_3\}_{k+1}\subseteq
				\pm
				\{\partial\mathbf{x}_1,\Sigma^k\mathbf{x}_2,
				\Sigma^k\mathbf{x}_3\}_k,
			\end{aligned}
			\]
			whenever the compositions are defined and the Toda
			bracket on the left is well-defined.
			Here
			\[
			i_Y\colon Y\hookrightarrow M_Y
			\hookrightarrow J(M_Y,X)=F_q
			\]
			is the canonical inclusion, and $j_Y=i\circ i_Y$.
			
			In applying these formulas, we need to express elements as
			linear combinations of the standard generators given in Appendix~\ref{tlqdzj} in view of Appendix~\ref{app:oguchi}.

			\item
			To resolve the group extension associated with the
			generators of $\ker(\partial)$, we use the following
			fact. Suppose that
			$\Sigma\mathbf{x}\in\ker(\partial)$ is a generator so that
			  $g\circ\mathbf{x}=0$ and $m\mathbf{x}=0$ for some integer $m\geq2$.
			Such an element often occurs in concrete calculations.
			There is a choice
			\[
			\operatorname{coext}_g(\mathbf{x})
			\in\operatorname{Coext}_g(\mathbf{x})
			\subseteq\{j_Y,g,\mathbf{x}\}
			\]
			satisfying
			\[
			q_*\bigl(\operatorname{coext}_g(\mathbf{x})\bigr)
			=-\Sigma\mathbf{x}.
			\]
			Consequently,
			\[
			-m\,\operatorname{coext}_g(\mathbf{x})
			\in j_Y\circ\{g,\mathbf{x},[m]\},
			\qquad [m]=m\,\mathrm{id}.
			\]
			It remains to analyze the Toda bracket
			$\{g,\mathbf{x},[m]\}$ using the fundamental properties
			of Toda brackets. If existing results (such as the Toda brackets given in \ref{sec:oguchi-working}) provide no information about $\{g,\mathbf{x},[m]\}$, this step often requires considerable ingenuity.
		\end{enumerate}
	\end{rem}

	Working 2-locally, recall from  Exercise \ref{nbdlx} that $P^{6}(2)=S^5\cup_2 e^{6}.$ The homotopy fiber $F_{q}$ of the pinch map $q\colon P^{6}(2)\to S^{6}$ has a CW decomposition $$F_q=S^{5}\cup e^{10}\cup e^{15}\cup \cdots \cup e^{5i}\cup \cdots,$$
	and  $\sk_{14}(F_q)=S^{5}\vee S^{10}$.
	
	\begin{exercise} Work 2-locally. Show that $$\pi_{9}(P^{6}(2))\cong\z/2\jia\z/2.$$ In fact, $\pi_{9}(P^{6}(2))\cong\z/2\jia\z/2$ before taking the 2-localization since the integral space $P^{6}(2)$ itself is 2-local.

		Hint: use the above remark and notice that $\n  2\e_5,4\nu_{5},2\e_{8}\nn=0$ implied by Proposition \ref{121X}. Note that the bracket $\n 2\e_5,4\nu_{5},2\e_{8}\nn$ lies outside the stable range in stem 4. With $\pi_{4}^S(S^0)=0$, stable homotopy groups provide no assistance in determining its value.
	\end{exercise}
	
	After reading the preceding remark, readers should find this exercise readily manageable with the help of $\n  2\e_5,4\nu_{5},2\e_{8}\nn=0$ and the appendixes, and we hope that it will deepen their interest in unstable homotopy groups. 
	
	\chapter{Preliminary Calculations}\label{ch:fundamental}
	Before turning to the principal computations considered in this volume, we must carry out several preliminary calculations. These calculations form the subject of the present chapter and are arranged so as to provide a systematic foundation for the mathematical arguments that follow. Although their immediate purpose is preparatory, they also help to clarify the relationships among the various ingredients that enter into the subsequent computations, thereby contributing to a more coherent presentation of the overall method.
	
	As in the previous chapter, we draw freely on the descriptions of the groups $\pi_{*}(S^{j})$
	presented in Appendix~\ref{tlqdzj} and on the relations governing
	compositions of generators in Toda's range collected in
	\^{O}guchi's table in Appendix~\ref{app:oguchi-relations}.

	\section{Composition Relations in Homotopy Groups of Spheres}
	
	The formulas for maps are useful only when the composites of maps between spheres on their right-hand sides are known.
	We begin with vanishing relations and then evaluate the brackets built from them. In each case, the
	null-homotopic composites establish that the bracket is defined before its value is calculated.

	We use the standard relation
	\begin{equation}
		\begin{aligned}
			&\eta^2_9\varepsilon_{11}=4(\nu_9\sigma_{12})=0.\label{m1:et^2e}
		\end{aligned}
	\end{equation}
	\noindent
	We also have
	\[\begin{aligned}
		&\nu_6\zeta_9=\zeta_6\nu_{17}=2(\sigma''\sigma_{13}).
	\end{aligned}\]
	Therefore, by Table~\ref{tab:oguchi-compositions} and $P(\iota_{17})=\pm(2\sigma_8-\Sigma\sigma')$ \cite[(5.16)]{T}, we obtain
	$$\nu_9\zeta_{12}=2(\Sigma^3 \sigma'')\sigma_{16} = 8\sigma^2_{9}.$$\noindent
	Hence, since the order of $\sigma^2_{14}$ is $8,$  
	we deduce
	\begin{equation}
		\begin{aligned}
			&\nu_9\zeta_{12}=\zeta_{9}\nu_{20}=8\sigma^2_{9} \quad\text{and}\quad \nu_{14}\zeta_{17}=0.\label{m1:n9z}
		\end{aligned}
	\end{equation}
	
	We consider the Toda bracket $\langle\eta^2,\varepsilon,\nu\rangle$ which is defined by the vanishing relations
	\eqref{m1:et^2e} and $\varepsilon\nu=0.$ 
	We note that $\langle\eta^2,\varepsilon,\nu\rangle
	\subseteq\pi^S_{14}=\frac{\mathbb Z}{2}\{\sigma^2\}\oplus\frac{\mathbb Z}{2}\{\kappa\}.$ 
	By the relations $\nu\eta=0$ and $\pi^S_{13}=0,$   we find
	$\langle\eta^2,\varepsilon,\nu\rangle\circ\eta=\eta^2\langle\varepsilon,\nu,\eta\rangle
	\subseteq\eta^2\circ\pi^S_{13}=0.$ 
	Hence, by $\eta\kappa=\bar\varepsilon\neq 0,$   we conclude that
	$\kappa \not\in \langle\eta^2,\varepsilon,\nu\rangle.$ 
	Moreover, by $\nu\sigma=0$ we see that
	$\langle\eta^2,\varepsilon,\nu\rangle\circ\sigma=\eta^2\circ\langle\varepsilon,\nu,\sigma\rangle.$ 
	By \cite[Lemma 2.2(2)]{MiM} and \cite[(3.9)-i)]{T}, it follows that $\langle\varepsilon,\nu,\sigma\rangle=0.$   Hence, by $\sigma^3 \neq 0$ \cite[Theorem A]{Mi1}, we obtain
	$\sigma^2 \not\in \langle\eta^2,\varepsilon,\nu\rangle$ and
	\begin{equation}
		\begin{aligned}
			&\langle\eta^2,\varepsilon,\nu\rangle=0.\label{m1:<et2,e,nu>}
		\end{aligned}
	\end{equation}
	
	We consider the Toda bracket $\{2\iota_{9}, \eta^2_{9}, \varepsilon_{11}\}$ which is well-defined
	by
	\eqref{m1:et^2e}.
	Let us recall from \cite[Lemma 9.1]{T} the relation $\zeta\in\langle 2\iota,\eta^2,\varepsilon\rangle.$ 
	Then, by the facts that $\pi^9_{20} = \frac{\mathbb Z}{8}\{\zeta_9\}\oplus\frac{\mathbb Z}{2}\{\bar\nu_9\nu_{17}\}$ and
	$\pi^S_{11}=\frac{\mathbb Z}{8}\{\zeta\}$ and the relation $\bar\nu\nu=0,$  
	there exists $a\in\{0,1\}$ satisfying
	$\zeta_9 + a\bar\nu_9\nu_{17}\in\{2\iota_{9}, \eta^2_{9}, \varepsilon_{11}\}.$ 
	Composing this relation on the right with $\nu_{20}$ and using
	$\varepsilon_{11}\nu_{19}=\nu_{11}\bar\nu_{14}=0,$   we have
	\[\begin{aligned}
		&\zeta_9\nu_{20} + a\bar\nu_9\nu^2_{17} \in\{2\iota_{9}, \eta^2_{9}, \varepsilon_{11}\}\circ \nu_{20} = -2\iota_{9} \circ \{\eta^2_{9}, \varepsilon_{11}, \nu_{19}\}.
	\end{aligned}\]
	Since $\zeta_9\nu_{20}=8\sigma^2_9$ \eqref{m1:n9z}, $\bar\nu_9\nu^2_{17}=2\kappa_9$ and
	$$\{\eta^2_{9}, \varepsilon_{11}, \nu_{19}\}
	\subseteq \operatorname{ker}(\Sigma^\infty: \pi^9_{23}\to\pi^S_{14})=\mathrm{span}\{2\sigma^2_9,2\kappa_9\}$$
	by \eqref{m1:<et2,e,nu>} and the standard group structures, we find
	\[\begin{aligned}
		&8\sigma^2_9 + 2a\kappa_9\in 2\iota_9\circ\mathrm{span}\{2\sigma^2_9,2\kappa_9\}= \mathrm{span}\{4\sigma^2_9\}.
	\end{aligned}\]
	This implies $a=0$ and
	\begin{equation}
		\begin{aligned}
			&\zeta_9\in\{2\iota_{9}, \eta^2_{9}, \varepsilon_{11}\}.\label{m1:zeta9}
		\end{aligned}
	\end{equation}
	Similarly, the Toda bracket $\{\eta^2_{9}, \varepsilon_{11}, 2\iota_{19}\}$ is well-defined.
	By the relation $\zeta\in\langle\eta^2,\varepsilon, 2\iota\rangle$ from \cite[Lemma 9.1]{T},
	we obtain
	\begin{equation*}
		\begin{aligned}
			&\zeta_9 + b\bar\nu_9\nu_{17} \in\{\eta^2_{9}, \varepsilon_{11}, 2\iota_{19}\} \text{\ for \ } b\in\{0,1\}.
		\end{aligned}
	\end{equation*}
	
	\subsection{Relations Involving $\nu,$    $\mu$ and $\sigma$}
	
	The next relations concern composites that recur in the small-$n$ calculations. Keeping their unstable
	subscripts is essential: a relation after stabilization does not by itself determine the corresponding
	composite in a lower-dimensional sphere.
	
	By $\sigma'\nu_{14}=x\nu_7\sigma_{10}$ for  some  odd $x,$   we conclude that
	\begin{equation}
		\begin{aligned}
			&\nu_7\sigma^2_{10}=0.\label{m1:n7sigsig}
		\end{aligned}
	\end{equation}
	By \eqref{m1:n9z},
	we see that
	\begin{equation}
		\begin{aligned}
			&\nu^2_6\zeta_{12}=0.\label{m1:n6nz}
		\end{aligned}
	\end{equation}
	
	The following congruences hold:
	\[\begin{aligned}
		&\mu_3\varepsilon_{12}\equiv\varepsilon_3\mu_{11}\ \bmod\ 2\bar{\varepsilon}' \ \mbox{and}\ \mu_3\varepsilon_{12}\equiv\eta_3\mu_4\sigma_{13}\ \bmod\ 2\bar{\varepsilon}'.
	\end{aligned}\]
	Therefore, by
	$\pm(2\nu_4-\Sigma\nu')=P(\iota_9)$ \cite[(5.8)]{T},
	$\Sigma\bar{\varepsilon}'=(\Sigma\nu')\kappa_7$ and
	the fact that the order of $\kappa_8$ is $4,$   we have\
	$2\Sigma^2\bar{\varepsilon}'=4(\nu_5\kappa_8)=0$ and
	\begin{equation}
		\begin{aligned}
			&\mu_5\varepsilon_{14}=\varepsilon_5\mu_{13}=\eta_5\mu_6\sigma_{15}.\label{m1:Ogu0}
		\end{aligned}
	\end{equation}
	
	\section{The $J$-Images for Small Values of $n$ and Connecting Values}
	
	For small values of $n,$   a chosen lift is often recognized by its $J$-image. We first record the reference values
	and then discuss their normalization. This separates a statement about a particular representative
	from one that is determined only up to a unit or an indeterminacy.
	
	\subsection{The Boundary of $\nu_4\sigma^\prime$}
	
	The next lemma records the ambiguity in the boundary formula $\Delta(\nu_4\sigma')=2[\iota_3]\varepsilon'$ \cite[Table 3]{Ka}.
	
	\begin{lem}
		$\Delta(\nu_4\sigma')\equiv 2[\iota_3]\varepsilon'\ \bmod\ 2[\eta_2]_4\varepsilon'.$ 
	\end{lem}
	\begin{proof}
		We use the relation\ $2\sigma'=\Sigma\sigma''$ in Table~\ref{tab:oguchi-compositions}. By
		the $\nu_4$ entry of Table~\ref{m2:tab:boundary-values}, we find
		$2\Delta(\nu_4\sigma')=\Delta(\nu_4(\Sigma\sigma''))=\Delta(\nu_4)\sigma''
		\equiv[\iota_3]\nu'\sigma''\ \bmod\ [\eta_2]_4\nu'\sigma''.$ 
		Since $\nu'\sigma''=0,$  
		it follows that $2\Delta(\nu_4\sigma')=0.$ 
		As a consequence of \cite[Table 2]{Ka}, we know
		\[\begin{aligned}
			&\pi_{13}(\mathit{SO}(4))={}\mathrm{span}\{[\iota_3]\varepsilon', [\eta_2]_4\varepsilon', [\iota_3]\eta_3\mu_4, [\eta_2]_4\eta_3\mu_4\}\cong(\mathbb Z/4)^2\oplus(\mathbb Z/2)^2.
		\end{aligned}\]
		Therefore, we can set
		$\Delta(\nu_4\sigma')=2a([\iota_3]\varepsilon')+2b([\eta_2]_4\varepsilon')+c [\iota_3]\eta_3\mu_4+d[\eta_2]_4\eta_3\mu_4$\ for\ $a, b, c, d\in\{0, 1\}.$  Let us take the $J$-image of this relation.
		By  the relation $[\iota_4,\nu_4]=\pm 2\nu^2_4$ \cite[p. 46]{T}, we conclude that
		$J\Delta(\nu_4\sigma')=[\iota_4,\nu_4\sigma']=4\nu_4^2\sigma_{10}.$  Then we see that
		$J([\iota_3]\varepsilon')=\nu_4(\Sigma^4\varepsilon')
		=\pm 2\nu^2_4\sigma_{10}.$  By  the relation $\nu'\nu_6\in\pi^3_9=0,$  
		we obtain $J([\eta_2]_4\varepsilon')=2(\Sigma\nu')\nu_7\sigma_{10}=0.$ 
		By similar arguments, we obtain $J([\iota_3]\eta_3\mu_4)=\nu_4\eta_7\mu_8$ and
		$J([\eta_2]_4\eta_3\mu_4)=(\Sigma\nu')\eta_7\mu_8.$  By inspecting the group structure of $\pi^4_{17},$   we obtain $a=1$ and $c=d=0.$ 
		This completes the proof.
	\end{proof}

	\subsection{The Lift of $\eta_{14}^{2}$ and Its $J$-Image}
	
	For $[\eta_{14}^{2}],$   the Hopf invariant gives the first comparison with a standard element in a homotopy group of a sphere. The
	remaining terms lie in the suspension image, so the calculation proceeds by examining that ambiguity.
	
	Since $HJ[\eta^2_{14}]=\eta^2_{29}= H(\eta^{*\prime})$ \cite[Lemma 12.14]{T},
	we deduce
	\[\begin{aligned}
		&J[\eta^2_{14}]\equiv \eta^{*\prime} \bmod \Sigma\pi^{14}_{30} = \mathrm{span}\{\omega_{15}, \sigma_{15}\mu_{22}\}.
	\end{aligned}\]
	As a consequence of               \cite[Lemma 2.10]{Og}, we have
	\[\begin{aligned}
		&[\iota_{16},\eta_{16}]\equiv \Sigma\eta^{*\prime} \bmod \sigma_{16}\mu_{23}.
	\end{aligned}\]
	
	Hence, since $[\eta^2_{14}]_{17} = 0$  and $\mathrm{ord}\omega=2,$   we obtain
	\begin{equation*}
		\begin{aligned}
			&J[\eta^2_{14}]\equiv \eta^{*\prime} \bmod \sigma_{15}\mu_{22}.
		\end{aligned}
	\end{equation*}
	We use the following: $H(\sigma')=\eta_{13},$  \ $2\sigma'=\Sigma\sigma'',$  \
	$2\sigma''=\Sigma\sigma'''$\ and\ $H(\sigma'')=\eta^2_{11}$ \cite[Lemma 5.14]{T}\uprightsemicolon{}
	$\pi^7_{14}=\frac{\mathbb Z}{8}\{\sigma'\},$   $\pi^6_{13}=\frac{\mathbb Z}{4}\{\sigma''\}$\ and\
	$\pi^5_{12}=\frac{\mathbb Z}{2}\{\sigma'''\}.$ 
	By \cite[Table 2]{Ka}, it follows that $\pi_{14}(\mathit{SO}(6))=\mathrm{span}\{[\eta_5\varepsilon_6],
	[ \nu_5]\nu^2_8\}\cong\mathbb Z/16\oplus\mathbb Z/2$ and
	\begin{gather}
		\begin{aligned}
			&{i_7}_*\pi_{14}(\mathit{SO}(6))={}\mathrm{span}\{[\eta_5\varepsilon_6]_7, [ \nu_5]_7\nu^2_8\} \cong\mathbb Z/8\oplus\mathbb Z/2.\label{m1:i7R614}
		\end{aligned}
	\end{gather}
	The two choices have the same projection after multiplying the  lift
	by $\mathrm{c}_{\mathrm{KM}}.$ Exactness therefore gives
	\begin{align*}
		[\mathrm{c}_{\mathrm{KM}}(\bar{\nu}_6+\varepsilon_6)]
		-\mathrm{c}_{\mathrm{KM}}[\bar{\nu}_6+\varepsilon_6]
		&\in\mathrm{span}\{[\eta_5\varepsilon_6]_7,[\nu_5]_7\nu_8^2\}.
	\end{align*}
	This comparison is preserved by the inclusion homomorphisms.
	In particular, the $J$-image in \eqref{m1:Jbn6e6} belongs to the
	Kachi--Mukai choice.

	\section{Zeta Lifts and Their Boundary Images}
	
	The $\zeta$-families show how a bracket description of an element of a homotopy group of a sphere leads to a relation for its
	lift. We first verify the null-homotopic composites that define the bracket and then compare its image under the
	boundary map. The  following lemma is useful.
	\begin{lem}\label{m1:zeta-lift-bracket}
		Suppose that $n\geq2.$ 
		\begin{itemize}
			\item[\textup{(1)}]
			The Toda bracket
			$\{\Delta\iota_{4n+2},\eta^2_{4n+1},\varepsilon_{4n+3}\}$
			is well-defined, and
			\begin{align}
				(p_{4n+2})_*
				\{\Delta\iota_{4n+2},\eta^2_{4n+1},\varepsilon_{4n+3}\}
				&=\{2\iota_{4n+1},\eta^2_{4n+1},\varepsilon_{4n+3}\}
				\notag\\
				&=\zeta_{4n+1}+\mathrm{span}\{2\zeta_{4n+1}\}.
				\notag
			\end{align}
			
			\item[\textup{(2)}]
			A lift $[\zeta_{4n+1}]\in\pi_{4n+12}(\mathit{SO}(4n+2))$
			can be chosen so that $(p_{4n+2})_*[\zeta_{4n+1}]=\zeta_{4n+1}$
			and
			$[\zeta_{4n+1}]\in
			\{\Delta\iota_{4n+2},\eta^2_{4n+1},\varepsilon_{4n+3}\}.$ 
		\end{itemize}
	\end{lem}
	
	\begin{proof}

		\begin{itemize}
			\item[\rm (1)]   Clearly,
			we have
			$\Delta\iota_{4n+2}\circ\eta^2_{4n+1}
			=\Delta(\eta^2_{4n+2})=0$ as well as $\eta^2_{4n+1}\circ\varepsilon_{4n+3}=0.$ 
			Thus $\{\Delta\iota_{4n+2},\eta^2_{4n+1},\varepsilon_{4n+3}\}$
			is well-defined.  Then we have
			\begin{align}
				(p_{4n+2})_*
				\{\Delta\iota_{4n+2},\eta^2_{4n+1},\varepsilon_{4n+3}\}
				&\subseteq
				\{p_{4n+2}\circ\Delta\iota_{4n+2},
				\eta^2_{4n+1},\varepsilon_{4n+3}\}
				\notag\\
				&=\{2\iota_{4n+1},\eta^2_{4n+1},\varepsilon_{4n+3}\}.
				\notag
			\end{align}
			By \eqref{m1:zeta9}, we infer 
			$\zeta_{4n+1}\in
			\{2\iota_{4n+1},\eta^2_{4n+1},\varepsilon_{4n+3}\}.$ 
			The standard group structures and composition relations show that
			the indeterminacy of this latter bracket is
			\begin{align}
				2\iota_{4n+1}\circ\pi^{4n+1}_{4n+12}
				+\pi^{4n+1}_{4n+4}\circ\varepsilon_{4n+4}
				&=\mathrm{span}\{2\zeta_{4n+1}\}.
				\notag
			\end{align}
			Hence this bracket is the coset
			$\zeta_{4n+1}+\mathrm{span}\{2\zeta_{4n+1}\}.$  Observe that 
			$$\ind\{\Delta\iota_{4n+2},\eta^2_{4n+1},\varepsilon_{4n+3}\}=\Delta\iota_{4n+2}\circ\pi^{4n+1}_{4n+12}
			+\pi_{4n+4}(\mathit{SO}(4n+2))\circ\varepsilon_{4n+4}.$$ \noindent
			The image of its first summand under
			$(p_{4n+2})_*$ is
			$2\iota_{4n+1}\circ\pi^{4n+1}_{4n+12},$  
			and the image of its second summand is contained in
			$\pi^{4n+1}_{4n+4}\circ\varepsilon_{4n+4}.$ 
			The preceding calculation therefore shows that the image of
			this indeterminacy is contained in
			$\mathrm{span}\{2\zeta_{4n+1}\}.$ 
			
			Conversely, the first summand contains
			$\Delta\iota_{4n+2}\circ\zeta_{4n+1},$   whose image is
			\begin{align}
				(p_{4n+2})_*
				(\Delta\iota_{4n+2}\circ\zeta_{4n+1})
				&=(p_{4n+2}\circ\Delta\iota_{4n+2})\circ\zeta_{4n+1}
				\notag\\
				&=2\iota_{4n+1}\circ\zeta_{4n+1}
				=2\zeta_{4n+1}.
				\notag
			\end{align}
			Thus the image of the indeterminacy is exactly
			$\mathrm{span}\{2\zeta_{4n+1}\},$   which is also the
			indeterminacy of
			$\{2\iota_{4n+1},\eta^2_{4n+1},\varepsilon_{4n+3}\}.$  Since one is contained in the other,
			they are equal. This proves the first assertion.
			\item[\rm (2)] By~\textup{(1)}, the fixed element $\zeta_{4n+1}$
			belongs to the image of
			$\{\Delta\iota_{4n+2},\eta^2_{4n+1},\varepsilon_{4n+3}\}$
			under $(p_{4n+2})_*.$  Choose an element of this bracket mapping
			to $\zeta_{4n+1}$ and denote it by $[\zeta_{4n+1}].$ 
			Then $(p_{4n+2})_*[\zeta_{4n+1}]=\zeta_{4n+1}$ and
			$[\zeta_{4n+1}]\in
			\{\Delta\iota_{4n+2},\eta^2_{4n+1},\varepsilon_{4n+3}\},$  
			as required.
		\end{itemize}

	\end{proof}

	\subsection{Connecting Images of Zeta Families}

	Recall that
	\begin{equation}
		\begin{aligned}
			&4\zeta_5=\eta^2_5\mu_7.\label{m1:4z5}
		\end{aligned}
	\end{equation}
	The preceding bracket construction yields a boundary formula and an order statement for the
	$\zeta$-lifts. Both are recorded in the next lemma because the degree calculations need the relation
	between the representatives, as well as their orders.
	
	\begin{lem} \label{m1:zeta1} Use the notation of the above lemma. We have
		\[\begin{aligned}
			&\Delta\zeta_{4n+2} = \pm 2[\zeta_{4n+1}] \quad\text{and}\quad \mathrm{ord}[\zeta_{4n+1}]=8 \text{\ for\ } n\ge 1.
		\end{aligned}\]
		In particular, $\Delta\zeta_6=2[\zeta_{5}]$ holds.
	\end{lem}
	\begin{proof}
		The assertion for $n=1$ follows from \cite[Proposition 4.1]{Ka}.
		
		Observe that \vspace{-1.2em}
		\begin{align*}
			- 2[\zeta_{4n+1}] &\in\{\Delta\iota_{4n+2}, \eta^2_{4n+1}, \varepsilon_{4n+3}\}\circ (-2\iota_{4n+12}) \\
			&=\Delta\iota_{4n+2}\circ\{\eta^2_{4n+1},\varepsilon_{4n+3}, 2\iota_{4n+11}\} \\
			&\ni\Delta\iota_{4n+2}\circ(\zeta_{4n+1}+ b\bar{\nu}_{4n+1}\nu_{4n+9})\quad(b\in\z) \\
			&= \Delta\zeta_{4n+2}+ \Delta(b\bar{\nu}_{4n+2}\nu_{4n+10}) = \Delta\zeta_{4n+2}.
		\end{align*}
		Since the indeterminacy of this relation is
		$$\begin{aligned}
			&\Delta\iota_{4n+2}\circ\pi_{4n+12}^{4n+1}\circ 2\iota_{4n+12}= \mathrm{span}\{\Delta\iota_{4n+2}\circ 2\zeta_{4n+1}\}=\mathrm{span}\{ 2\Delta\zeta_{4n+2}\},
		\end{aligned}$$
		we see that
		$2[\zeta_{4n+1}]\equiv -\Delta\zeta_{4n+2} \bmod 2\Delta\zeta_{4n+2}.$ 
		By the relations \eqref{m1:4z5},
		we deduce
		$4\Delta\zeta_{4n+2} = \Delta(\eta^2_{4n+2}\mu_{4n+4}) = 0.$ 
		Hence we obtain the first relation.
		We also have
		$
		8[\zeta_{4n+1}] = 4\Delta\zeta_{4n+2} = 0.
		$
		Since $\zeta_{4n+1}$ has order 8,
		$[\zeta_{4n+1}]$ has order at least 8.
		This leads to the second assertion, completing the proof.
	\end{proof}

	\section{The Orders of  Some Lifts}
	
	To determine the order of a lift defined by a Toda bracket, we compose the bracket with an appropriate
	degree map. The bracket identities then reduce the multiple to another bracket, whose vanishing must
	be established. The following relations prepare this step.
	
	Recall from \cite[Lemma 6.5]{T} the relation
	\begin{equation}
		\begin{aligned}
			&\mu_n\in\{\eta_n,2\iota_{n+1},\Sigma^{n-4}\sigma'''\}_{n-4}+\mathrm{span}\{\nu^3_n\} \text{\ \ for\ } n\ge 4.\label{m1:mu4}
		\end{aligned}
	\end{equation}
	According to \cite[Corollary 3.7]{T}, the Toda brackets
	$\{2\iota_n,8\sigma_n,2\iota_{n+7}\}_1$
	and $\{8\iota_n,2\sigma_n,8\iota_{n+7}\}_1$ are trivial
	for $n \ge 10.$  Since $\Sigma\pi^{n-1}_{n+7}=\pi^n_{n+8}$ for $n \geq 10,$  
	we find
	\begin{equation}
		\begin{aligned}
			&\{2\iota_n,8\sigma_n,2\iota_{n+7}\}=\{8\iota_n,2\sigma_n,8\iota_{n+7}\}=0 \text{\quad for } n\ge 10.\label{m1:2i8s2i}
		\end{aligned}
	\end{equation}
	
	\subsection{A Criterion for a Lift Annihilated by Two}
	
	The next criterion states this argument for lifts obtained from an $\eta$-bracket. Its 
	hypotheses specify the null-homotopic composites needed in the construction and in the calculation of the
	multiple.
	
	\begin{lem} \label{m1:ord}
		Let $\alpha\in\pi^{4n+2}_k$ with $n\ge 2$ be an element satisfying the relations
		$$2\alpha=0, \;\s^{3}(\nu_{4n-1}\af)=0,\; \Sigma^2(\alpha\eta_{k})=0\text{ and } \eta^2_{4n+2}(\Sigma^2\alpha)=0.$$\noindent
		Then, for each element $$\beta\in\{\eta_{4n+2},2\iota_{4n+3},\Sigma\alpha\}_1\subseteq\pi^{4n+2}_{k+2},$$\noindent
		there exists a lift $[\eta_{4n+2}(\Sigma\beta)]\in \pi_{k+3}(\mathit{SO}(4n+3))$ of $\eta_{4n+2}(\Sigma\beta)$
		which is annihilated by $2.$  If $\eta_{4n+2}(\Sigma\beta)\neq0,$   this lift has order exactly $2.$ 
	\end{lem}
	\begin{proof} We know that
		\[\begin{aligned}
			\Delta(\eta_{4n+2}(\Sigma\beta))&\in \pi_{4n+2}(\mathit{SO}(4n+2))\circ 2\iota_{4n+2}\circ\beta \\
			&\subseteq -\pi_{4n+2}(\mathit{SO}(4n+2))\circ\{2\iota_{4n+2}, \eta_{4n+2},2\iota_{4n+3}\}_1\circ(\Sigma^2\alpha) \\
			&=-\pi_{4n+2}(\mathit{SO}(4n+2))\circ\eta^2_{4n+2}(\Sigma^2\alpha)=0.
		\end{aligned}\]
		Therefore, $\eta_{4n+2}(\Sigma\beta)$ has a lift.

		Recalling \cite[(4.2)]{GM} and the fact that $2\pi_{4n+4}(\mathit{SO}(4n+3))=0$ \cite{Ke}, we infer
		\begin{equation}\notag
			\begin{aligned}
				&\Delta(\eta^2_{4n+2})=0\uprightsemicolon{}\ \mathrm{ord}[\eta^2_{4n+2}]=2.
			\end{aligned}
		\end{equation}
		In light of $\s^{3}(\nu_{4n-1}\af)=0,$   we calculate $$\ind (p_{4n+3}\hc\{[\eta^2_{4n+2}],2\iota_{4n+4},\Sigma^2\alpha\}_2
		)=\ind(\{\eta^2_{4n+2},2\iota_{4n+4},\Sigma^2\alpha\}_2).$$
		So,  $$p_{4n+3}\hc\{[\eta^2_{4n+2}],2\iota_{4n+4},\Sigma^2\alpha\}_2=\{\eta^2_{4n+2},2\iota_{4n+4},\Sigma^2\alpha\}_2.$$
		Then we may take
		\[\begin{aligned}
			&[\eta_{4n+2}(\Sigma\beta)]\in\{[\eta^2_{4n+2}],2\iota_{4n+4},\Sigma^2\alpha\}_2.
		\end{aligned}\]
		\noindent   Hence  we deduce that
		\[\begin{aligned}
			2[\eta_{4n+2}(\Sigma\beta)] &\in\{[\eta^2_{4n+2}],2\iota_{4n+4},\Sigma^2\alpha\}_2\circ 2\iota_{k+3} \\
			&=-[\eta^2_{4n+2}]\circ\{2\iota_{4n+4},\Sigma^2\alpha, 2\iota_{k+2}\}_2 \\
			&\ni[\eta^2_{4n+2}](\Sigma^2\alpha)\eta_{k+2}=0 \mod\ [\eta^2_{4n+2}]\circ \Sigma^2\pi^{4n+2}_{k+1}\circ 2\iota_{k+3}=0.
		\end{aligned}\]
		If $\eta_{4n+2}(\Sigma\beta)\neq0,$   its lift is nonzero and therefore has order exactly $2.$  This completes the proof.
	\end{proof}
	
	Note that $\mu_n\in\{\eta_n,2\iota_{n+1},8\sigma_{n+1}\}_1$ for $n \ge 9.$  In light of Lemma~\ref{m1:ord}, we obtain
	\begin{exam}\label{m1:shetep}
		$\mathrm{ord}[\eta_{4n+2}\mu_{4n+3}]=2$ for $n \ge 2.$ 
	\end{exam}

	To facilitate the arguments that follow, we  establish the following lemma, which will serve as a useful preliminary result.
	
	\begin{lem}\label{myjbu1026}Suppose $n\geq3.$ Let $M^n=S^{n-1}\cup_{2\iota_{n-1}}e^n$ be the mod~$2$
		Moore space, with inclusion $i'_n\colon S^{n-1}\to M^n$ 
		and projection $p'_n\colon M^n\to S^n.$     Let
		$$\overline{\;\ca_{n}\;}\in[M^{n+2},S^n]\cong\mathbb Z/4,\quad\widetilde{\;\ca_{n}\;} \in\pi_{n+2}(M^{n+1})\cong\mathbb Z/4$$
		be an extension and a coextension of $\eta_n,$   respectively. Then the following statements hold.
		\begin{itemize}
			\item[\textup{(1)}] $\overline{\;\ca_{n}\;}\in\{\eta_n,2\iota_{n+1},p'_{n+1}\},$  
			$\widetilde{\;\ca_{n}\;} \in\{i'_{n+1},2\iota_n,\eta_n\},$  
			$2\overline{\;\ca_{n}\;}=\eta_n^2p'_{n+2},$   and
			$2\widetilde{\;\ca_{n}\;} =i'_{n+1}\eta_n^2.$ 
			Moreover,
			\begingroup
			\renewcommand{\theHequation}{eeta}
			\begin{equation}\label{eeta}
				\overline{\;\ca_{n}\;}\hc\widetilde{\;\ca_{n+1}\;} 
				=\pm\Sigma^{n-3}\nu',
			\end{equation}
			\endgroup
			and
			\begin{equation}\label{ietap}
				2\,\mathrm{id}=i'_n\eta_{n-1}p'_n.
			\end{equation}
			\noindent In addition, $\Delta(\eta_{4n+2}\hc\overline{\;\ca_{4n+3}\;})$ is $0.$  Here, $\Delta\colon[M^{4n+5},S^{4n+2}]
			\to[M^{4n+4},\mathit{SO}(4n+2)]$ is the connecting homomorphism. Thus $\eta_{4n+2}\hc\overline{\;\ca_{4n+3}\;}$ admits a lift
			$$[\eta_{4n+2}\hc\overline{\;\ca_{4n+3}\;}]\in [M^{4n+5},\mathit{SO}(4n+3)].$$
			
			\item[\textup{(2)}]
			The elements $\varepsilon_k$ admit the bracket description
			\begin{equation}\label{vep}
				\varepsilon_k\in
				\{\eta_k\hc       \overline{\;\ca_{k+1}\;},\widetilde{\;\eta_{k+2}\;},\nu_{k+4}\}
				\qquad(k\geq5).
			\end{equation}
			Additionally, we have \vspace{-0.8em}
			\begin{align}
				\eta_{4n+2}\varepsilon_{4n+3}
				&=\varepsilon_{4n+2}\eta_{4n+10}\notag\\
				&=\eta_{4n+2}\hc\overline{\;\ca_{4n+3}\;}\circ
				\{\widetilde{\;\eta_{4n+4}\;},\nu_{4n+6},\eta_{4n+9}\}.
				\notag
			\end{align}
			
			\item[\textup{(3)}]A lift
			$[\eta_{4n+2}\varepsilon_{4n+3}]
			\in\pi_{4n+11}(\mathit{SO}(4n+3))$ may be chosen so that
			\begin{equation}\label{[eta-ep]}
				[\eta_{4n+2}\varepsilon_{4n+3}]
				\in[\eta_{4n+2}\hc\overline{\;\ca_{4n+3}\;}]\circ
				\{\widetilde{\;\eta_{4n+4}\;},\nu_{4n+6},\eta_{4n+9}\}.
			\end{equation}
			This choice satisfies
			\begin{equation}\label{2[eta-ep]}
				2[\eta_{4n+2}\varepsilon_{4n+3}]=0.
			\end{equation}
		\end{itemize}
	\end{lem}
	
	\begin{proof}
		\begin{itemize}
			\item[\textup{(1)}]
			The extension and coextension constructions give the bracket
			descriptions and the two doubling identities.
			Equation~\eqref{eeta} is \cite[(5.3)]{T}.
			By \cite{Toda1963}, we know \eqref{ietap} holds for $n\geq4.$ 
			Since $[M^3,M^3]$ has four elements and
			$\Sigma\colon[M^3,M^3]\to[M^4,M^4]$ is surjective,
			this homomorphism is an isomorphism.
			For $n=3,$   the two sides of \eqref{ietap} have the same
			suspension by the case $n=4,$   and hence are equal.
			This establishes the Moore space identities in the stated range.    The vanishing $\pi_{4n+3}(\mathit{U}(2n+1))=0$
			from \cite[Theorem~4.4]{To1} gives
			$\omega_{2n+1}(\mathbb C)\eta_{4n+2}=0.$ 
			It follows that
			\begin{align}
				\Delta(\eta_{4n+2}\hc\overline{\;\ca_{4n+3}\;})
				&=\Delta(\eta_{4n+2})\circ\overline{\;\ca_{4n+2}\;}\notag\\
				&=r_{2n+1}\omega_{2n+1}(\mathbb C)
				\circ\eta_{4n+2}^2p'_{4n+4}=0.\notag
			\end{align}
			The  exactness gives the required lift.
			
			\item[\textup{(2)}]
			The bracket description \eqref{vep} is given in
			\cite[Lemma~4.1]{GM}.  From \eqref{vep}, we infer
			\begin{align}
				\eta_{4n+2}\varepsilon_{4n+3}
				&=\varepsilon_{4n+2}\eta_{4n+10}\notag\\
				&\in\{\eta_{4n+2}   \hc\overline{\;\ca_{4n+3}\;},
				\widetilde{\;\eta_{4n+4}\;},\nu_{4n+6}\}
				\circ\eta_{4n+10}\notag\\
				&=\eta_{4n+2}   \hc\overline{\;\ca_{4n+3}\;}\circ
				\{\widetilde{\;\eta_{4n+4}\;},\nu_{4n+6},\eta_{4n+9}\}
				\notag\\
				&\quad\bmod\
				\eta_{4n+2}\    \hc\overline{\;\ca_{4n+3}\;}\circ
				\pi_{4n+10}(M^{4n+5})\circ\eta_{4n+10}.
				\notag
			\end{align}
		It is esay to infer that 
			$\pi_{4n+10}(M^{4n+5})
			=\mathrm{span}\{i'_{4n+5}\nu_{4n+4}^2\}.$ 
			Then we obtain
			\begin{align}
				\pi_{4n+10}(M^{4n+5})\circ\eta_{4n+10}
				&=\mathrm{span}\{i'_{4n+5}\nu_{4n+4}^2
				\eta_{4n+10}\}=0.\notag
			\end{align}
			The indeterminacy in the preceding calculation therefore
			vanishes, giving the equality in \textup{(2)}.
			
			\item[\textup{(3)}]
			Fix the lift $[\eta_{4n+2}  \hc\overline{\;\ca_{4n+3}\;}]$
			from \textup{(1)}.
			The equality in \textup{(2)} shows that every element of
			$[\eta_{4n+2}   \hc\overline{\;\ca_{4n+3}\;}]\circ
			\{\widetilde{\;\eta_{4n+4}\;},\nu_{4n+6},\eta_{4n+9}\}$
			projects to $\eta_{4n+2}\varepsilon_{4n+3}.$ 
			Choose one such element and denote it by
			$[\eta_{4n+2}\varepsilon_{4n+3}],$   as in \eqref{[eta-ep]}.
			To calculate its double, apply the Toda bracket composition
			rule once more:
			\begin{align}
				2[\eta_{4n+2}\varepsilon_{4n+3}]
				&\in[\ca_{4n+2}\hc\overline{\;\ca_{4n+3}\;} ]\circ
				\{\widetilde{\;\eta_{4n+4}\;},\nu_{4n+6},\eta_{4n+9}\}
				\circ2\iota_{4n+11}\notag\\
				&=[\ca_{4n+2}\hc\overline{\;\ca_{4n+3}\;} ]
				\hc\widetilde{\;\eta_{4n+4}\;}\circ
				\{\nu_{4n+6},\eta_{4n+9},2\iota_{4n+10}\}
				\notag\\
				&=0.\notag
			\end{align}
			\noindent   Thus the chosen lift is annihilated by $2,$   proving
			\eqref{2[eta-ep]}.
		\end{itemize}
	\end{proof}

	\begin{rem} The preceding lemma illustrates the essential role of the mod $2$ Moore spaces
		$S^{n-1}\cup_2 e^n$ in homotopy group calculations.
		Following Fred Cohen's notation, we write the mod $2$ Moore space
		$P^n(2)=M^n=S^{n-1}\cup_2 e^n$.
		Assume henceforth that $n\geq4$, so that $P^n(2)$ is a double suspension.
		It has been known since at least the 1960s that
		$[P^n(2),P^n(2)]\cong\z/4$,
		generated by $\mathrm{id}_{P^n(2)}$.
		Nevertheless, $\pi_i(P^n(2))$ contains a direct summand isomorphic to
		$\z/8$ for infinitely many values of $i$. See \cite{Wu03} and \cite{CW2013}.
		This tells us of an interesting phenomenon:
		\[
		4\pi_*(P^n(2))\neq0
		\qquad\text{in general}.
		\]
		\noindent In fact, for $\af\in \pi_{*}(\s^2 X)$ and $k\in\z$, in $\pi_{*}(\s ^2X)$, $$k(\mathrm{id})\hc\af\neq \af\hc k(\mathrm{id})\qquad\text{in general}.$$
		M.~G.~Barratt conjectured that
		\[
		8\pi_i(P^n(2))=0
		\qquad\text{for every }i\in\z_+.
		\]
		This remains open and is among the most difficult conjectures in unstable
		homotopy theory.
		Fred Cohen devoted much of his career to efforts to resolve it.  The odd-prime analogue of the above conjecture is more tractable;
		see the classic work of Cohen, Moore, and Neisendorfer \cite{CMN1979}.
	\end{rem}

	Related conjectures due to M.G. Barratt are as follows. 
	
	\begin{conj0}\label{conj:barratt} Work integrally. Let $p$ be a prime and let $X$ and $Y$ be CW complexes.
		\begin{itemize} 
			\item[\textup{(1)}]
			\textup{(The  Barratt conjecture).}
			If $f\colon\Sigma^2X\to Y$ is a map of order $p^r$
			in $[\Sigma^2X,Y]$ for some $r\in\z_{+}$, then, for all $k$, $p^{r+1}$ annihilates
			the image of
			$f_*\colon\pi_k(\Sigma^2X)\to\pi_k(Y)$. \item[\textup{(1)$'$}]
			\textup{(The strong  Barratt conjecture).}
			If $f\colon\Sigma^2X\to Y$ is a map of order $p^r$
			in $[\Sigma^2X,Y]$ for some $r\in\z_{+}$, then, 
			\[
			p^{r+1}\Omega^2f=0
			\quad\text{in\;\;}[\Omega^2\Sigma^2X,\Omega^2Y].
			\]   \item[\textup{(2)}]
			\textup{(The special Barratt conjecture).}
			If
			$p^r\mathrm{id}_{\Sigma^2X}=0$ in $[\Sigma^2X,\Sigma^2X]$ for some $r\in\z_{+}$, then
			\[
			p^{r+1}\pi_k(\Sigma^2X)=0
			\quad\text{for all }k.
			\]
			
			\item[\textup{(2)$'$}]
			\textup{(The strong special Barratt conjecture).}
			If
			$p^r\mathrm{id}_{\Sigma^2X}=0$ in $[\Sigma^2X,\Sigma^2X]$ for some $r\in\z_{+}$, then
			\[
			p^{r+1}\mathrm{id}_{\Omega^2\Sigma^2X}=0
			\quad\text{in }[\Omega^2\Sigma^2X,\Omega^2\Sigma^2X].
			\]
		\end{itemize}
	\end{conj0}
	
	\subsection{Epsilon Relations and Sigma Compositions}

	Recall that
	\begin{equation}
		\begin{aligned}
			&\varepsilon^2_3=\varepsilon_3\bar{\nu}_{11}=\eta_3\bar{\varepsilon}_4,\ \nu_5\sigma_{8}\nu^2_{15}=\eta_5\bar\varepsilon_6,\ \sigma'\eta_{14}\bar{\nu}_ {15} = \sigma'\nu^3_{14} = \nu_7\sigma_{10}\nu^2_{17} = \eta_7\bar{\varepsilon}_8,\ \eta_9\bar{\varepsilon}_{10}=0.\label{m1:epep}
		\end{aligned}
	\end{equation}
	
	Using $\pi^7_{27}=\mathrm{span}\{\bar{\kappa}_7\}\cong\mathbb Z/8$
	and the argument in the proof of \cite[(12.16)]{T}, we obtain
	\begin{equation}
		\begin{aligned}
			&H(\bar{\kappa}_7)\equiv\kappa_{13}\ \bmod\ 8\sigma^2_{13}.\label{m1:Hbk7}
		\end{aligned}
	\end{equation}
	We next record several relations that arise
	in the course of our calculations.
	\begin{lem} \label{m1:zesig}
		\begin{enumerate}
			\item $\bar{\nu}^2_6=\bar{\nu}_6\varepsilon_{14}=0.$ 
			\item $(\bar{\nu}_6+\varepsilon_6)\bar{\nu}_{14} =
			(\bar{\nu}_6+\varepsilon_6)\varepsilon_{14} = \eta_6\bar{\varepsilon}_7.$ 
			
			\item $\nu_5\sigma_8\bar{\nu}_{15} = \nu_5\sigma_8\varepsilon_{15}
			= \nu_5\bar{\varepsilon}_8.$ 
			\item $\zeta_5\sigma_{16}\in\{\nu_5\sigma_8, 8\iota_{15}, 2\sigma_{15}\}_1
			\ \bmod \nu_5\bar{\varepsilon}_8,\ 2\zeta_5\sigma_{16}.$ 
			
		\end{enumerate}
	\end{lem}
	\begin{proof}
		\begin{enumerate}[label=\textup{(\arabic*)},leftmargin=2em]
			\item These are standard vanishing relations. See   \^{O}guchi's table in Appendix~\ref{app:oguchi-relations}.

			\item By the relations (1) and \eqref{m1:epep}, we have (2).
			
			\item Assertion (3) is a standard composition identity. See     \^{O}guchi's table in Appendix~\ref{app:oguchi-relations}.
			
			\item We have $\pi^5_{23}=\mathrm{span}\{\zeta_5\sigma_{16},
			\nu_5\bar{\varepsilon}_8, \eta_5\bar{\mu}_6\}\cong\mathbb Z/8\oplus(\mathbb Z/2)^2$
			and, by \cite[Theorem 12.8]{T}, $\Sigma^2: \pi^5_{23}\to\pi^7_{25}$ is a split epimorphism and $\operatorname{ker} \Sigma^2=\mathrm{span}\{\nu_5\bar{\varepsilon}_8\}.$ 
			Then we see that
			\begin{align*}
				\Sigma^2\{\nu_5\sigma_8, 8\iota_{15}, 2\sigma_{15}\}_1 &\subseteq \{\nu_7\sigma_{10},8\iota_{17}, 2\sigma_{17}\}_1 \\
				&=\{x\sigma'\nu_{14},8\iota_{17}, 2\sigma_{17}\}_1 \\
				&\supseteq x\sigma'\circ\{\nu_{14},8\iota_{17}, 2\sigma_{17}\}_1 \\
				&\ni x\sigma'\zeta_{14} \\
				&=xy\zeta_7\sigma_{18} \mathrel{\mathrm{mod}}\;\begin{aligned}[t] &x\sigma'\nu_{14}\circ \Sigma\pi^{16}_{24} + \pi^7_{18}\circ 2\sigma_{18}=\mathrm{span}\{2\zeta_7\sigma_{18}\},
				\end{aligned}
			\end{align*}
			where $x$ and $y$ are odd.
			This leads to (4) and completes the proof.
		\end{enumerate}
	\end{proof}
	
	\subsection{Toda Brackets Involving $\mu$ and $\sigma$}
	
	Several later lifts involve the composites of $\mu$ or $\zeta$ with $\sigma.$  The following bracket
	descriptions provide representatives suitable for calculating their multiples. The congruences retain
	the indeterminacy that must be carried into those calculations.
	
	\begin{lem}\label{m1:m3sg}
		\begin{enumerate}
			\item
			$\mu_3\sigma_{12}$ is in
			$$\begin{aligned}
				&\{\varepsilon_3,\sigma_{11},16\iota_{18}\}_1 =\{\varepsilon_3,2\sigma_{11},8\iota_{18}\}_1 =\{\varepsilon_3,2\iota_{11},8\sigma_{11}\}_1 =\{\varepsilon_3,8\iota_{11},2\sigma_{11}\}_1,
			\end{aligned}$$
			which is a coset of $\mathrm{span}\n\eta_3\bar{\varepsilon}_4\nn\cong\mathbb Z/2.$
			
			\item
			Let $A=\mathrm{span}\{\eta_6\bar\varepsilon_7\}.$ 
			There exist $t\in\{0,1\}$ and $u\in\{1,3,5,7\}$ such that, with
			$\mathbf{x}=\zeta'+t\eta_6\bar\varepsilon_7,$   the following equalities of subsets of $\pi_{22}^{6}$ hold:
			$\{\bar\nu_6,8\iota_{14},2\sigma_{14}\}_1=\{\mathbf{x}\},$
			$$\{\bar\nu_6,8\sigma_{14},2\iota_{21}\}_1=\mathbf{x}+\mathrm{span}\{2\zeta'\},\;\{\bar\nu_6,\sigma_{14},16\iota_{21}\}_1
			=\{\bar\nu_6,2\sigma_{14},8\iota_{21}\}_1=\{u\mathbf{x}\}.$$\noindent
			Here $\{\mathbf{x}\}$ and $\{u\mathbf{x}\}$ denote singleton sets.

			\item
			With $A$ and $u$ as in (2),
			$\{\bar\nu_6+\varepsilon_6,\sigma_{14},16\iota_{21}\}_1
			=\{\bar\nu_6+\varepsilon_6,2\sigma_{14},8\iota_{21}\}_1
			=u\zeta'+\mu_6\sigma_{15}+A,$
			$\{\bar\nu_6+\varepsilon_6,8\iota_{14},2\sigma_{14}\}_1
			=\zeta'+\mu_6\sigma_{15}+A.$ 
			Moreover,
			$$\Sigma\zeta'+\mu_7\sigma_{16}\in\{\bar{\nu}_7+\varepsilon_7,\sigma_{15},16\iota_{22}\}_1
			=\{\bar{\nu}_7+\varepsilon_7,2\sigma_{15},8\iota_{22}\}_1=\{\bar{\nu}_7+\varepsilon_7,2\iota_{15},
			8\sigma_{15}\}_1\ \bmod \eta_7\bar{\varepsilon}_8.$$
		\end{enumerate}
	\end{lem}
	\begin{proof}
		\begin{enumerate}[label=\textup{(\arabic*)},leftmargin=2em]
			\item   As a consequence of          \cite[I-Proposition 3.2(1), its proof]{Oda1979}, we conclude that
			\[\begin{aligned}
				&\mu_3\sigma_{12}\in\{\varepsilon_3,2\iota_{11},8\sigma_{11}\}_1 \bmod \eta_3\bar{\varepsilon}_4.
			\end{aligned}\]
			We see that
			\[\begin{aligned}
				&\{\varepsilon_3,\sigma_{11},16\iota_{18}\}_1 \subseteq \{\varepsilon_3,2\sigma_{11},8\iota_{18}\}_1 \supseteq\{\varepsilon_3,2\iota_{11},8\sigma_{11}\}_1 \bmod\ \varepsilon_3\circ \Sigma\pi^{10}_{18} + 8\pi^3_{19}.
			\end{aligned}\]
			We know that $\Sigma\pi^{10}_{18}=\mathrm{span}\{\bar\nu_{11},\varepsilon_{11}\}$ and $8\pi^3_{19}=0.$  Therefore, by  \eqref{m1:epep}, we obtain (1), except for its last
			assertion. Since $\pi^3_{12}\cong(\mathbb Z/2)^2,$  
			we obtain
			$$\{\varepsilon_3,2\iota_{11},8\sigma_{11}\}_1
			\subseteq\{\varepsilon_3,8\iota_{11},2\sigma_{11}\}_1
			\bmod\ \varepsilon_3\circ \Sigma\pi^{10}_{18} + \pi^3_{12}\circ 2\sigma_{12}
			=\mathrm{span}\{\eta_3\bar{\varepsilon}_4\}.$$ This leads to (1).
			
			\item By \cite[Lemma~5.4(iv), pp.~156--157]{Mi2}, the ordinary bracket
			$\{\bar\nu_6,8\iota_{14},2\sigma_{14}\}$ is congruent to $\zeta'$ modulo $A.$ 
			Its nonempty subbracket with subscript $1$ therefore has its values in $\zeta'+A.$ 
			By the Stem~16 table in Appendix~B,
			$$\pi_{22}^{6}=\mathrm{span}\{\zeta',\mu_6\sigma_{15},\eta_6\bar\varepsilon_7\}
			\cong\mathbb Z/8\oplus(\mathbb Z/2)^2.$$
			Also, $\bar\nu_6\circ\Sigma\pi_{21}^{13}=0$ and $2\pi_{15}^{6}=0.$ 
			Thus the exact indeterminacy of $\{\bar\nu_6,8\iota_{14},2\sigma_{14}\}_1$ is zero,
			so this bracket is $\{\mathbf{x}\}$ for some $\mathbf{x}=\zeta'+t\eta_6\bar\varepsilon_7$ with $t\in\{0,1\}.$ 
			The bracket $\{\bar\nu_6,8\sigma_{14},2\iota_{21}\}_1$ contains $\mathbf{x}$ and has
			indeterminacy $\bar\nu_6\circ\Sigma\pi_{21}^{13}+2\pi_{22}^{6}
			=\mathrm{span}\{2\zeta'\},$   giving the second equality.
			
			The composition properties give
			$$\{\bar\nu_6,\sigma_{14},16\iota_{21}\}_1
			\subseteq\{\bar\nu_6,2\sigma_{14},8\iota_{21}\}_1
			\subseteq\{\bar\nu_6,8\sigma_{14},2\iota_{21}\}_1.$$
			The first two brackets have zero indeterminacy, since $16\pi_{22}^{6}=8\pi_{22}^{6}=0.$ 
			They are consequently the same singleton, whose value is $\mathbf{x}+2k\zeta'$ for some integer $k.$ 
			Since $2\mathbf{x}=2\zeta',$   this value is $u\mathbf{x}$ for an odd residue $u$ modulo $8.$ 
			This proves (2) without identifying distinct singleton brackets from containment in a larger coset.
			
			\item We have
			$(\bar\nu_6+\varepsilon_6)\circ\Sigma\pi_{21}^{13}=A,$  
			$8\pi_{22}^{6}=16\pi_{22}^{6}=0$ and $2\pi_{15}^{6}=0.$ 
			Hence each of the three brackets with first entry $\bar\nu_6+\varepsilon_6$ in (3)
			has exact indeterminacy $A.$ 
			Suspending (1) three times, and then using the inclusion into brackets with subscript $1,$  
			shows that the corresponding brackets with first entry $\varepsilon_6$ contain $\mu_6\sigma_{15}.$ 
			Their indeterminacy is also $A.$ 
			Additivity in the first entry, (2), and equality of these indeterminacies give the first two assertions in (3);
			the term $t\eta_6\bar\varepsilon_7$ is absorbed by $A.$ 
			
			For the three brackets in $\pi_{23}^{7},$   the common indeterminacy is
			$\mathrm{span}\{\eta_7\bar\varepsilon_8\}.$ 
			Suspending the first equality in (3) gives the stated value for the first of these brackets,
			since $2\Sigma\zeta'=0$ and $u$ is odd.
			The composition properties and equality of indeterminacies identify the first two of them,
			and \cite[(2.5)]{HKM} gives the third. This completes the proof.
		\end{enumerate}
	\end{proof}
	
	\section{Comparison Methods for Small Values of $n$}
	
	The special orthogonal sequence is not the only available comparison for small values of $n.$  Special descriptions of the
	small-$n$ groups can make a projection injective or separate the contribution of an exceptional group.
	We record these two useful forms of simplification.
	
	\subsection{The Inclusion from $\mathit{SO}(4)$ to $\mathit{SO}(5)$}
	
	The first lemma identifies conditions under which the image from the case $n=4$ vanishes. Exactness then
	turns the projection from $\mathit{SO}(5)$ into an injection, so a class can be detected by its image in the
	sphere.
	
	\begin{lem}\label{m1:236}
		\begin{enumerate}
			
			\item If $i_{5 \ast}\pi_{n}(\mathit{SO}(4))=0,$   where $n \ge 4,$   then $p_{5 \ast} :
			\pi_{n}(\mathit{SO}(5)) \to \pi_n^4$ is a monomorphism.
			\item If $\pi^3_n=\nu' \circ \pi^6_n$ for $n \ge 6,$   then $i_{5
				\ast}\pi_{n}(\mathit{SO}(4))=0.$ 
		\end{enumerate}
	\end{lem}
	\begin{proof}
		\begin{enumerate}[label=\textup{(\arabic*)},leftmargin=2em]
			\item Assertion (1) follows directly from $(n)_4.$ 
			
			\item Recall that
			$\pi_{n}(\mathit{SO}(4))=[\eta_2]_4 \circ \pi^3_n
			\oplus [\iota_3] \circ \pi^3_n$ (\cite[p.~14]{Ka}, \cite[p. 116]{St}).
			
			We note that $\pi_n(\mathit{SO}(3))\cong\pi_n(S^3),$   $\pi_n(\mathit{SO}(4))\cong\pi_n(\mathit{SO}(3))
			\oplus\pi_n(S^3).$ 
			It follows that $\pi_6(\mathit{SO}(5)) \cong \pi_6(\mathit{Sp}(2))=0.$ 
			Then $i_{5 \ast}([\eta_2]_4\nu')=[\eta_2]_5\nu'=0$ and $i_{5 \ast}([\iota_3]\nu')=[\iota_3]_5\nu'=0.$ 
		\end{enumerate}
	\end{proof}
	
	We
	know $\pi^3_6=\frac{\mathbb Z}{4}\{\nu'\},$  \ $\pi^3_7=\frac{\mathbb Z}{2}\{\nu'\eta_6\},$  
	$\pi^3_8=\frac{\mathbb Z}{2}\{\nu'\eta^2_6\},$   $\pi^3_9=\pi^3_{10}=0,$  \ $\pi^3_{15}
	=\nu'\circ\pi^6_{15},$   $\pi^3_{16}=\nu'\circ\pi^6_{16}$ and $\pi^3_{17}
	=\mathrm{span}\{\varepsilon_3\nu^2_{11}\},$   where $\varepsilon_3\nu^2_{11}
	=\nu'\bar{\nu}_6\nu_{14}$ and $\nu'\zeta_6=0.$  Therefore, by Lemma~\ref{m1:236}(2), we conclude that
	\begin{equation}
		\begin{aligned}
			&{i_5}_\ast \pi_{n}(\mathit{SO}(4))=0\ \mbox{for}\ 6\le n\le 10\ \mbox{and}\ 15\le n\le 17.\label{m1:i5R4n}
		\end{aligned}
	\end{equation}
	\subsection{The Splitting for $n=7$}

	\begin{prop}
		The following identity holds.
		$\mathcal{T}_{*}(\left<\eta_6^2\right>)=[ \nu_5]_7+[\eta_6]\eta_7.$
	\end{prop}
	\begin{proof}
		Observe that $\pi_8 (G_2)=\frac{\mathbb Z}{2}\{\left<\eta_6^2\right>\}$ \cite[Theorem 6.1]{Mi2} and
		$$\pi_{8}(\mathit{SO}(7))=\frac{\mathbb Z}{2}\{[ \nu_5]_7\}\oplus\frac{\mathbb Z}{2}\{[\eta_6]\eta_7\}$$ from \cite[Table 2]{Ka}.
		Using the diagram \eqref{m2:Y} and the relation
		${p_7}_*\mathcal{T}_{*}\left<\eta_6^2\right>=\eta_6^2={p_7}_*([\eta_6]\eta_7),$   we obtain
		$\mathcal{T}_{*}(\left<\eta_6^2\right>)\equiv [\eta_6]\eta_7 \bmod [ \nu_5]_7.$
		Hence, by Proposition~\ref{so7fj}, the asserted identity follows.
	\end{proof}
	
	\section{Relations in the Homotopy Groups of Spheres Used in Degrees 26 and 27}\label{m2:sec:sphere-relations}
	
	This section collects the composition relations and Toda-bracket
	calculations involving elements of homotopy groups of spheres that are used in the subsequent
	Lie-group calculations.  
	\subsection{Composition Relations}
	
	We begin with ordinary composites. These vanishing results will supply the null-homotopic compositions and
	remove terms in the bracket indeterminacies that follow.
	
	\begin{lem}\label{m2:lem:nu6-rho}
		The composition $\nu_6\rho'$ is zero.
		\begin{equation}
			\begin{aligned}
				&\nu_6\rho'=0.\label{m2:Osm}
			\end{aligned}
		\end{equation}
	\end{lem}
	\begin{proof}
		This is a standard vanishing relation in the Toda range. See    \^{O}guchi's table in Appendix~\ref{app:oguchi-relations}.\end{proof}
	For later reference, we include the following
	elementary observation concerning composition.
	\begin{lem}\label{m2:lem:nuprime-compositions}
		The compositions $\nu'\bar\sigma_6$ and $\nu'\bar\zeta_6$ vanish.
		\[\begin{aligned}
			&\nu'\bar\sigma_6=\nu'\bar\zeta_6=0.
		\end{aligned}\]
	\end{lem}
	\begin{proof}
		These are the composition relations in the homotopy groups of spheres used in
		\cite[pp.~313--314]{Mi1}.
	\end{proof}
	
	\begin{lem}\label{m2:lem:sigma-kappa-vanishing}
		The composition $\sigma''\kappa_{13}$ vanishes.
	\end{lem}
	\begin{proof}
		We use the homotopy groups of spheres
		$\pi^6_{27}=\frac{\mathbb Z}{2}\{\eta_6\bar\kappa_7\}\cong\mathbb Z/2$ and $\pi^7_{28}=\frac{\mathbb Z}{2}\{\eta_7\bar\kappa_8\}\oplus\frac{\mathbb Z}{2}\{\sigma'\kappa_{14}\}\cong(\mathbb Z/2)^2.$
		Since $\sigma''\kappa_{13}\in\pi^6_{27},$ we can write
		$\sigma''\kappa_{13}=x\eta_6\bar\kappa_7, \; x\in\{0,1\}.$
		Suspending gives
		$0=2\sigma'\kappa_{14}=\Sigma(\sigma''\kappa_{13}) =x(\eta_7\bar\kappa_8).$
		Thus $x=0,$   so $\sigma''\kappa_{13}=0.$ 
	\end{proof}
	We next record a useful relation
	concerning the composition under consideration.
	\begin{lem}\label{m2:lem:eta5-composition}
		\begin{gather}
			\FitMath{\eta_5\circ\pi^6_{22}={}\mathrm{span}\{\eta_5\mu_6\sigma_{15}\}.\label{m2:et5pi}}
		\end{gather}
	\end{lem}
	\begin{proof}
		The relevant group is
		$\pi^6_{22}=\mathrm{span}\{\zeta', \mu_6\sigma_{15}, \eta_6\bar{\varepsilon}_7\}\cong\mathbb Z/8\oplus(\mathbb Z/2)^2.$
		Using the relation $\eta_5\zeta'=0$ and
		$\eta^2_5\bar{\varepsilon}_7=\eta_5\nu_6\sigma_9\nu^2_{16}=0,$   we obtain the asserted relation.
	\end{proof}
	
	\subsection{Toda Bracket Evaluations}
	
	With the necessary composition relations in place, we turn to bracket values. The first calculation
	determines both a representative and zero indeterminacy\uprightsemicolon{} the distinction matters when its value is
	later used as an equality.
	
	\begin{lem}\label{m2:lem:rho-bracket}
		The Toda bracket $\{2\sigma_{11},8\iota_{18},2\sigma_{18}\}$
		has zero indeterminacy and contains $2\Sigma^2\rho'.$   That is,
		\begin{equation}
			\begin{aligned}
				&\{2\sigma_{11},8\iota_{18},2\sigma_{18}\} =2\Sigma^2\rho'.\label{m2:2e2rh'}
			\end{aligned}
		\end{equation}
	\end{lem}
	\begin{proof}
		The definition of $\rho''$ in \cite[p.~103]{T}, together with
		$2\rho'=\Sigma^2\rho''$ from \cite[Lemma 2.9]{Og}, gives
		$2\Sigma^2\rho'=\Sigma^4\rho''\in\{2\sigma_{11},8\iota_{18},2\sigma_{18}\}.$ 
		The indeterminacy is
		$2\sigma_{11}\circ\pi_{26}^{18}+\pi_{19}^{11}\circ2\sigma_{19}=0,$  
		as in the cited calculation.  Thus the unique value of the bracket is
		$2\Sigma^2\rho'.$ 
	\end{proof}
	
	First we need
	\begin{lem}\label{tx1023}
		\begin{enumerate}
			\item $\{\sigma'\eta_{14}+\varepsilon_7,\bar{\nu}_{15},2\iota_{23}\}_1
			=\nu_7\kappa_{10}+\mathrm{span}\{\sigma'\eta_{14}\mu_{15}+\eta_7\mu_8\sigma_{17}\}.$ 
			\item $\{2\iota_9,\eta_9,\bar{\varepsilon}_{10}\}=\nu_9\kappa_{12}.$ 
		\end{enumerate}
	\end{lem}
	\begin{proof}
		\begin{enumerate}[label=\textup{(\arabic*)},leftmargin=2em]
			\item We consider the Toda bracket
			$\{\sigma'\eta_{14}+\varepsilon_7,\bar{\nu}_{15},2\iota_{23}\}_1,$   which is well-defined by $$(\sigma'\eta_{14}+\varepsilon_7)\bar\nu_{15}=0$$ \noindent (see \eqref{m1:epep}).
			It has the indeterminacy
			$$(\sigma'\eta_{14}+\varepsilon_7)\circ \Sigma\pi^{14}_{23}+2\pi^{7}_{24}
			=(\sigma'\eta_{14}+\varepsilon_7)\circ \mathrm{span}\{\nu^3_{15},\mu_{15},\eta_{15}\varepsilon_{16}\}$$\noindent 
			by the standard group structures.
			We obtain $(\sigma'\eta_{14}+\varepsilon_7)\nu^3_{15}=0.$ 
			We deduce $$(\sigma'\eta_{14}+\varepsilon_7)\mu_{15}
			=\sigma'\eta_{14}\mu_{15}+\eta_7\mu_8\sigma_{17}$$ \noindent
			by \eqref{m1:Ogu0}.
			Observe that $(\sigma'\eta_{14}+\varepsilon_7)\eta_{15}\varepsilon_{16}=0$
			by \eqref{m1:et^2e} and
			$\varepsilon_7\eta_{15}\varepsilon_{16}=\eta_7\varepsilon^2_8=\eta_7\nu_8\sigma_{11}\nu^2_{18}=0$
			\eqref{m1:epep}.
			Hence,  the indeterminacy  of $\{\sigma'\eta_{14}+\varepsilon_7,\bar{\nu}_{15},2\iota_{23}\}_1$
			is $\mathrm{span}\{\sigma'\eta_{14}\mu_{15}+\eta_7\mu_8\sigma_{17}\}.$ 
			We note that
			$\Sigma^3\{\sigma'\eta_{14}+\varepsilon_7,\bar{\nu}_{15},2\iota_{23}\}_1
			\subseteq -\{\varepsilon_{10},\bar\nu_{18},2\iota_{26}\}.$ 
			By \cite[III-Proposition 6.1(2)]{Oda1979} and the fact that
			$\varepsilon_{10}\circ\pi^{18}_{27} + 2\pi^{10}_{27} = \mathrm{span}\{\eta_{10}\mu_{11}\sigma_{20}\},$  
			we have $\{\varepsilon_{10},\bar{\nu}_{18},2\iota_{26}\}
			=\nu_{10}\kappa_{13}+\mathrm{span}\{\eta_{10}\mu_{11}\sigma_{20}\}.$ 
			Since $$\operatorname{ker}(\Sigma^3: \pi^7_{24}\to\pi^{10}_{27}) = \mathrm{span}\{\sigma'\eta_{14}\mu_{15}\},$$\noindent  we can set
			\[\begin{aligned}
				&\nu_7\kappa_{10} + a\sigma'\eta_{14}\mu_{15} + b\eta_7\mu_8\sigma_{17} \in\{\sigma'\eta_{14}+\varepsilon_7,\bar{\nu}_{15},2\iota_{23}\}_1,
			\end{aligned}\]
			for $a,$   $b\in\{0, 1\}.$ 
			Compose this relation with $\eta_{24}$ on the right. Then, in $\pi^{7}_{25}=\frac{\mathbb Z}{8}\{\zeta_7\sigma_{18}\}\oplus\frac{\mathbb Z}{2}\{\eta_7\bar\mu_8\},$  
			\[\begin{aligned}
				4(a+b)\zeta_7\sigma_{18} &\in\{\sigma'\eta_{14}+\varepsilon_7,\bar{\nu}_{15},2\iota_{23}\}_1\circ\eta_{24} \\
				&=(\sigma'\eta_{14} +\varepsilon_7)\circ\{\bar{\nu}_{15},2\iota_{23},\eta_{23}\} \\
				&\subseteq (\sigma'\eta_{14}+\varepsilon_7)\circ \pi^{15}_{25}=0.
			\end{aligned}\]
			Here $\pi^{15}_{25}=\mathbb Z/2\{\eta_{15}\mu_{16}\}$ by Appendix~\ref{tlqdzj}. The sphere relations in Appendix~\ref{app:oguchi} and \eqref{m1:Ogu0} give $(\sigma'\eta_{14}+\varepsilon_7)\eta_{15}\mu_{16}=4\sigma'\zeta_{14}+4\zeta_7\sigma_{18}=0,$   since $\sigma'\zeta_{14}=u\zeta_7\sigma_{18}$ for an odd integer $u$ and $\zeta_7\sigma_{18}$ has order $8.$ 
			Hence, we obtain $a+b\equiv 0 \bmod 2.$ 
			This leads to (1).
			
			\item   As a consequence of           \cite[III-Proposition 6.1(1)]{Oda1979}, we obtain
			\[\begin{aligned}
				&\{2\iota_{10},\eta_{10},\bar{\varepsilon}_{11}\}\equiv\nu_{10}\kappa_{13}\ \bmod \eta_{10}\mu_{11}\sigma_{20}.
			\end{aligned}\]
			
			Applying \cite[Proposition 3.4(i)]{T} with $\alpha=2\iota_1,$   $\beta=\eta_6$ and
			$\gamma=\bar{\varepsilon}_3,$   we see that
			$\{2\iota_{10},\eta_{10},\bar{\varepsilon}_{11}\}$ and
			$\{\bar{\varepsilon}_{10},\eta_{25},2\iota_{26}\}$ have a common element.
			This leads to
			\begin{equation}
				\begin{aligned}
					&\{\bar{\varepsilon}_{10},\eta_{25},2\iota_{26}\}\equiv\nu_{10}\kappa_{13} \bmod\ \eta_{10}\mu_{11}\sigma_{20}.\label{m1:tbep10}
				\end{aligned}
			\end{equation}
			
			We have the decomposition
			\[\begin{aligned}
				&\pi^9_{26}={}\mathrm{span}\{\sigma_9\eta_{16}\mu_{17},\nu_9\kappa_{12}, \bar{\mu}_9, \eta_9\mu_{10}\sigma_{19}\}\cong(\mathbb Z/2)^4.
			\end{aligned}\]
			Therefore, we see that
			\[\begin{aligned}
				&\{2\iota_9,\eta_9,\bar{\varepsilon}_{10}\}\equiv\nu_9\kappa_{12}\ \bmod\ \sigma_9\eta_{16}\mu_{17}, \eta_9\mu_{10}\sigma_{19}.
			\end{aligned}\]
			Composing this equation on the left with $\eta_8,$   we obtain
			\[\begin{aligned}
				\eta_8\circ\{2\iota_9,\eta_9,\bar{\varepsilon}_{10}\} &= -\{\eta_8,2\iota_9,\eta_9\}\circ\bar{\varepsilon}_{11} \\
				&\ni \pm 2\nu_8\bar{\varepsilon}_{11} =0\ \bmod\ \eta_8 \circ 2\iota_9\circ\pi^{9}_{26}+\eta_8\circ\pi^9_{11}\circ\bar{\varepsilon}_{11}=0.
			\end{aligned}\]
			
			Using \eqref{m1:Ogu0} and \eqref{m1:4z5},
			we obtain
			\[\begin{aligned}
				&\eta_8\sigma_9\eta_{16}\mu_{17}=((\Sigma\sigma')\eta_{15}+\bar{\nu}_8+\varepsilon_8) \eta_{16}\mu_{17}=4(\Sigma\sigma')\zeta_{15}+4\zeta_8\sigma_{19}=0.
			\end{aligned}\]
			Notice that
			\[\begin{aligned}
				&\eta_8(\eta_9\mu_{10}\sigma_{19})=4\zeta_8\sigma_{19}\ne 0.
			\end{aligned}\]
			This implies
			\[\begin{aligned}
				&\{2\iota_9,\eta_9,\bar{\varepsilon}_{10}\}\equiv\nu_9\kappa_{12}\ \bmod\ \sigma_9\eta_{16}\mu_{17}.
			\end{aligned}\]
			
			Next, we record
			\[\begin{aligned}
				&\pi^9_{27}={}\mathrm{span}\{\sigma_9\zeta_{16}, \eta_9\bar{\mu}_{10}\}\cong\mathbb Z/8\oplus\mathbb Z/2.
			\end{aligned}\]
			
			Composing both sides of the preceding equation on the right with $\eta_{26}$ gives
			\[\begin{aligned}
				&\{2\iota_9,\eta_9,\bar{\varepsilon}_{10}\}\circ\eta_{26} = -2\iota_9\circ\{\eta_9,\bar{\varepsilon}_{10},\eta_{25}\},
			\end{aligned}\]
			$\nu_9\kappa_{12}\eta_{26}=\nu_9\eta_{12}\kappa_{13}=0$ and
			$\sigma_9\eta_{16}\mu_{17}\eta_{26}=4\sigma_9\zeta_{16}.$ 
			As a consequence of           \eqref{m1:tbep10}, we find
			\[\begin{aligned}
				&2\iota_9\circ \{\eta_9,\bar{\varepsilon}_{10},\eta_{25}\} =2\{\eta_9,\bar{\varepsilon}_{10},\eta_{25}\} =\{\eta_9,\bar{\varepsilon}_{10},\eta_{25}\}\circ2\iota_{27} =-\eta_9\circ\{\bar{\varepsilon}_{10}, 
				{}\eta_{25},2\iota_{26}\} \equiv 0 
			\end{aligned}\] $\bmod\ \eta^2_9\mu_{11}\sigma_{20}=4\zeta_9\sigma_{20} =8(\sigma_9\zeta_{16})=0.$
			This leads to (2) and completes the proof.
		\end{enumerate}
	\end{proof}

	We now turn our attention to another pair
	of Toda brackets, whose precise determination
	is given below.
	\begin{lem}\label{m2:lem:nu-kappa-bracket}
		The following Toda brackets have zero indeterminacy.
		\begin{equation}
			\begin{aligned}
				&\{2\iota_9,\eta_9,\bar\varepsilon_{10}\} =\nu_9\kappa_{12}.\label{m2:nu9k}
			\end{aligned}
		\end{equation}
		\[\begin{aligned}
			&\{2\iota_{10},\eta_{10},\bar\varepsilon_{11}\} =\nu_{10}\kappa_{13}.
		\end{aligned}\]
	\end{lem}
	\begin{proof}
		The first evaluation is Lemma~\ref{tx1023}(2).
		Suspension sends it to $\nu_{10}\kappa_{13}$ in the second bracket.
		The indeterminacy of that bracket is
		$2\iota_{10}\circ\pi_{27}^{10}
		+\pi_{12}^{10}\circ\bar\varepsilon_{12}=0$
		by the standard composition relations.
	\end{proof}
	For convenient reference, we collect below
	a number of relations that arise naturally
	in the present setting.
	\begin{lem}\label{m2:ods}
		The following relations hold.
		\begin{itemize}
			\item[\normalfont(1)] $\{2\iota_9,\eta_9,\Sigma\rho'\} = \sigma_9\eta_{16}\mu_{17}.$
			\item[\normalfont(2)] $\{\nu_7,\Sigma\rho',16\iota_{25}\} \subseteq \bar\zeta_7+\mathrm{span}\{2\bar\zeta_7\}.$
			This bracket has zero indeterminacy, so its unique value is an odd multiple of $\bar\zeta_7$.
			\item[\normalfont(3)]
			\[
			\begin{aligned}
				&\{\eta_5,\Sigma\rho^{\mathrm{IV}},2\iota_{21}\} =\bar\mu_5+\mathrm{span}\{\eta_5\mu_6\sigma_{15},2\nu_5\kappa_8\}, \\
				&\{\eta_9,\Sigma\rho',16\iota_{25}\} =\bar\mu_9+\mathrm{span}\{\eta_9\mu_{10}\sigma_{19}\}.
			\end{aligned}
			\]
			
			\item[\normalfont(4)]   $\n\nu_{7},\ca_{10}\cn_{11},2\e_{25}\nn\ni 0$
		\end{itemize}
	\end{lem}
	\begin{proof}
		\begin{enumerate}[label=\textup{(\arabic*)},leftmargin=2em]
			\item This is a direct consequence of \cite[Theorem 1]{Oda2}.
			
			\item The representative congruence is given by \cite[Theorem 1]{Oda2}. We compute the indeterminacy explicitly.
			By Appendix~\ref{app:stem-16}, $\pi^{10}_{26}=\mathrm{span}\{P(\sigma_{21}),\sigma_{10}\mu_{17}\}$.
			Since $P(\sigma_{21})=[\iota_{10},\iota_{10}]\sigma_{19}$ and $S^7$ is an H-space, naturality of the Whitehead product gives $\nu_7P(\sigma_{21})=[\nu_7,\nu_7]\sigma_{19}=0$.
			By Appendix~\ref{app:oguchi-relations}, $\sigma'\nu_{14}=x\nu_7\sigma_{10}$ for an odd $x$.
			Also, $\pi^{14}_{26}=0$ by Appendix~\ref{app:stem-12}, so $x\nu_7\sigma_{10}\mu_{17}=\sigma'\nu_{14}\mu_{17}=0$.
			As $x$ is odd, $\nu_7\sigma_{10}\mu_{17}=0$, and hence $\nu_7\circ\pi^{10}_{26}=0$.
			Appendix~\ref{app:stem-19} gives $\pi^7_{26}=\mathrm{span}\{\bar\zeta_7,\bar\sigma_7\}\cong\mathbb Z/8\oplus\mathbb Z/2$, so $16\pi^7_{26}=0$.
			Therefore $\ind\{\nu_7,\Sigma\rho',16\iota_{25}\}=\nu_7\circ\pi^{10}_{26}+16\pi^7_{26}=0$.
			Thus the unique value lies in the stated coset; the subgroup $\mathrm{span}\{2\bar\zeta_7\}$ records the representative congruence, not the indeterminacy.
			
			\item Recall that
			$\pi^5_{22} = \mathrm{span}\{\nu_5\kappa_8, \bar{\mu}_5, \eta_5\mu_6\sigma_{15}\} \cong \mathbb Z/4\oplus(\mathbb Z/2)^2.$
			
			The Toda bracket $\{\eta_5,\Sigma\rho^{\mathrm{IV}},2\iota_{21}\}$ is defined because $\eta_4\rho^{\mathrm{IV}} = 0.$ 
			By \eqref{m2:et5pi}, its indeterminacy is $\eta_5\circ\pi^6_{22}+2\pi^5_{22}
			=\mathrm{span}\{\eta_5\mu_6\sigma_{15},2\nu_5\kappa_8\}.$ Since $\langle\eta,2\iota,16\rho\rangle\ni\bar{\mu}$ \cite[I-Proposition 3.3(1)]{Oda1979},
			$\Sigma^\infty\rho^{\mathrm{IV}}=16\rho$  and
			$\pi^S_{16}=\mathrm{span}\{\mu\sigma, \omega\},$   we see that
			$$\Sigma^\infty\{\eta_5,\Sigma\rho^{\mathrm{IV}},2\iota_{21}\}\subseteq \langle\eta,16\rho,2\iota\rangle\supseteq\langle\eta,2\rho,16\iota\rangle \ni\bar{\mu}\m\eta\circ\pi^S_{16}=\mathrm{span}\{\eta\mu\sigma, \eta\omega\}.$$  \noindent  This implies the equation
			$\{\eta_5,\Sigma\rho^{\mathrm{IV}},2\iota_{21}\} = \bar{\mu}_5 +\mathrm{span}\{\eta_5\mu_6\sigma_{15},2\nu_5\kappa_8\}.$
			From this equation and $\Sigma^5\rho^{\mathrm{IV}} = 8\Sigma\rho',$ we note that$$\n\eta_9,\Sigma^5\rho^{\mathrm{IV}},2\iota_{25}\}=\{\eta_9,8\Sigma\rho',2\iota_{25}\}\supseteq\{\eta_9,\Sigma\rho',16\iota_{25}\}\bmod\bigl(\eta_9\circ\pi^{10}_{26}+2\pi^9_{26}\bigr),$$
			$$\eta_9\circ\pi^{10}_{26}+2\pi^9_{26}=\mathrm{span}\{\eta_9\mu_{10}\sigma_{19}\nn.$$
			\noindent   This implies that
			$\{\eta_9,\Sigma\rho',16\iota_{25}\} = \bar{\mu}_9 +\mathrm{span}\{\eta_9\mu_{10}\sigma_{19}\}.$ 
			This proves (3).
			\item   Clearly,
			$\x\nu,\ca\cn,2\e\xx=\x\nu,\cn\ca,2\e\xx\dyd\cn\hc\x\nu,\ca,2\e\xx=0.$
			So $\x\nu,\ca\cn,2\e\xx=\frac{\z}{4}     \n2\cs\nn.$ 
			Since the homomorphism $\pi_{19+7}^{7}\xrightarrow{\s^{\wq}}\pi_{19}^{S}$ is injective, we have
			$\n\nu_{7},\ca_{10}\cn_{11},2\e_{25}\nn\ni t(2\cs_{7})$ for some $t\in\z.$ 
			Observe that
			$$\ind\n\nu_{7},\ca_{10}\cn_{11},2\e_{25}\nn\ni2\cs_{7}.$$ \noindent Thus, $\n\nu_{7},\ca_{10}\cn_{11},2\e_{25}\nn\ni0.$ 
		\end{enumerate}
	\end{proof}
	We examine
	a particular Toda bracket, whose properties
	are established in the following lemma.
	\begin{lem}\label{m2:lem:mu10-sigma-bracket}
		The Toda bracket $\{\varepsilon_{10},8\iota_{18},2\sigma_{18}\}_1$
		has zero indeterminacy and satisfies
		\[\begin{aligned}
			&\{\varepsilon_{10},8\iota_{18},2\sigma_{18}\}_1 =\mu_{10}\sigma_{19}.
		\end{aligned}\]
	\end{lem}
	\begin{proof}
		By \cite[I-Proposition 3.2(1), p.~57]{Oda1979} and its proof, we find
		$\mu_3\sigma_{12}\in\{\varepsilon_3,2\iota_{11},8\sigma_{11}\} \ \bmod\
		\varepsilon^2_3=\eta_3\bar{\varepsilon}_4.$ Since $\varepsilon^2_{10}=
		\eta_{10}\bar{\varepsilon}_{11}=\nu_{10}\sigma_{13}\nu^2_{20}=0,$   we conclude that \[\begin{aligned}
			&\{\varepsilon_{10},2\iota_{18},8\sigma_{18}\}_1 \subseteq\{\varepsilon_{10},8\iota_{18},2\sigma_{18}\}_1 \ \bmod \ \varepsilon_{10}\circ \Sigma\pi^{17}_{25}+\pi^{10}_{19}\circ 2\sigma_{19}=0.
		\end{aligned}\] That is,
		$\{\varepsilon_{10},8\iota_{18},2\sigma_{18}\}_1=\mu_{10}\sigma_{19}.$
	\end{proof}
	
	\section{Stable Groups and Results on the Homotopy Groups of Spheres for Degrees 28 and 29}
	
	The calculations in degrees $28$ and $29$ require further results on the homotopy groups of spheres. We retain the same notation, but make explicit the
	role of coefficients on composites before listing the additional relations.

	\subsection{Relations among Generators of $\pi_{*}(S^{n})$}
	We recall the following relations.
	\needspace{3\baselineskip}\begin{lem}
		\label{m3:2ltp}
		\leavevmode\par
		\begin{itemize}
			\item[\textup{(1)}] $2\varepsilon'=\eta_{3}^{2}\varepsilon_5,$
			$\varepsilon_{5}\rho_{13}=0.$
			\item[\textup{(2)}] $\nu'\nu_{6}=0.$
			
			\item[\textup{(3)}] $\ca_{6}\cs_{7}=\ca_{6}\ccs_{7}=0.$
		\end{itemize}
	\end{lem}
	\begin{proof}
		\begin{itemize}
			\item[\textup{(1)}] $2\varepsilon'=\eta_{3}^{2}\varepsilon_5$ is the standard relation.
			The second relation is given in \cite[Proposition 3.1(4), p. 55]{Oda1979}.
			\item[\textup{(2)}] $\nu'\nu_{6}=0$ follows from $\nu'\nu_{6}\in\pi_{9}^{3}=0.$ 
			
			\item[\textup{(3)}]  By \cite{T}, we find $\cs_{7}\in\n \ck_{7},8\e_{18},2\cg_{18}\nn$. In light of \cite[Proposition 3.4(3), p.~59]{Oda1979},  we obtain $\ccs_{7}=\n \nu_{7},\ca_{10},\cg_{11}^{2}\nn$. Note that $\ca_{6}\ck_{7}=\nu_{9}\cg_{12}^{2}=0$. The remaining details of the argument are straightforward and may be readily supplied by the reader.
		\end{itemize}
	\end{proof}
	\needspace{3\baselineskip}Another useful observation emerges from
	the preceding analysis and is recorded below.\begin{prop}
		\label{m3:np}
		The following identity holds.
		$\nu'\bar{\nu}_{6}\kappa_{14}=\varepsilon_{3}\nu_{11}\kappa_{14}=\nu'\varepsilon_{6}\kappa_{14}.$
	\end{prop}
	\begin{proof}
		We have $\nu'\bar{\nu}_{6}=\varepsilon_{3}\nu_{11},$ and hence
		$\nu'\bar{\nu}_{6}\kappa_{14}=\varepsilon_{3}\nu_{11}\kappa_{14}.$ Lemma 6.1 of \cite{T} gives
		$H(\varepsilon_{3})=\nu_{5}^{2}.$ Using $\nu_{7}^{2}\kappa_{13}=4\bar{\kappa}_{7}$ and \cite[Lemma
		15.4]{MT1}, we obtain
		$$H(\varepsilon_{3}\nu_{11}\kappa_{14})=\nu_{5}^{2}\nu_{11}\kappa_{14}=\nu_{5}\nu_{8}^{2}\kappa_{14}=\nu_{5}\circ4\bar{\kappa}_{8}=4\nu_{5}\bar{\kappa}_{8}.$$ \noindent
		Moreover,
		\cite[Proposition 6.1(1), p. 66]{Oda1979} gives
		$H(\nu'\varepsilon_{6}\kappa_{14})=4\nu_{5}\bar{\kappa}_{8},$ and therefore
		\[\begin{aligned}
			&\nu'\varepsilon_{6}\kappa_{14}-\varepsilon_{3}\nu_{11}\kappa_{14} \in \operatorname{ker}(H) =\, \Sigma\pi_{27}^{2} =\frac{\mathbb Z}{2}\{2\phi'\}\oplus \frac{\mathbb Z}{2}\{\eta_{3}\bar{\mu}_{4}\sigma_{21}\}.
		\end{aligned}\]
		Thus there exist integers $x$ and $y$ such that
		\[\begin{aligned}
			&\nu'\varepsilon_{6}\kappa_{14}-\varepsilon_{3}\nu_{11}\kappa_{14}=2x\phi'+y\eta_{3}\bar{\mu}_{4}\sigma_{21}.
		\end{aligned}\]
		Applying $\,\Sigma^{\infty}$ gives $0=y\eta\bar{\mu}\sigma,$
		since $\pi_{25}^{S}(S^{0})=\frac{\mathbb Z}{2}\{\mu_{3 }\}\oplus \frac{\mathbb Z}{2}\{\eta\bar{\mu}\sigma\}.$
		Hence $y \equiv0\;\bmod\;2$ and
		$\nu'\varepsilon_{6}\kappa_{14}-\varepsilon_{3}\nu_{11}\kappa_{14}=2x\phi'.$ By \cite[(7.6), p.~69]{Oda1979},
		\[
		\Sigma^2\phi'\equiv2\phi''\m\nu_5\varepsilon_8\kappa_{16},\;\nu_5^2\bar{\sigma}_{11},\;\mu_{3,5},\;\eta_5\bar{\mu}_6\sigma_{23}.
		\]
		Each term in the modulus is annihilated by $2.$  Applying $\Sigma^2$ to the preceding equality therefore gives $0=4x\phi''.$  Since $\mathrm{ord}(\phi'')=8$ \cite{Oda1979}, we have $4x\equiv0\pmod8,$   so $x$ is even. It follows that
		$\nu'\varepsilon_{6}\kappa_{14}-\varepsilon_{3}\nu_{11}\kappa_{14}=2x\phi'=0\uprightsemicolon{}$ for
		$\mathrm{ord}(\phi')=4.$ See \cite{Oda1979}.
	\end{proof}
	
	\section{Comparison Formulas for  Some  Odd-Primary Components}\label{sec:odd-method}
	
	The odd-primary calculations use comparison formulas with symplectic groups and spheres. They are
	treated separately because the localization convention changes here. The formulas below provide the
	common input for the odd-torsion sections of the next three chapters.
	
	Exception to the standing convention. Throughout this section $p$ is an odd prime and all spaces and maps are localized at $p.$  The convention at $2$ resumes in the next chapter.

	Let $n\geq1$ and let $p$ be an odd prime. Work $p$-locally. The following formulas determine the odd-primary components
	in Chapters~\ref{ch:degree23}--\ref{ch:degree25}.
	\begin{gather}\label{m1:R-Sp:odd}
		\pi_k(\mathit{SO}(2n+1)) \cong \pi_k(\mathit{Sp}(n)) \ \cite[(3)]{Ha}\uprightsemicolon{}\\
		\label{m1:R-Sp:even}
		\pi_k(\mathit{SO}(2n+2)) \cong \pi_k(\mathit{SO}(2n+1))\oplus\pi_k(S^{2n+1}) \  \cite{BS}.
	\end{gather}
	Also $\pi_k(\mathit{SO}(1))=\pi_k(\mathit{SO}(2))=0$ for $k=23,24,25,$   since $\mathit{SO}(1)$ is a point and $\mathit{SO}(2)=S^1.$ 
	\chapter{The Homotopy Groups $\pi_{23}(\mathit{SO}(n))$}\label{ch:degree23}

	Degree $23$ is our first complete application of the method. The calculation begins with small-$n$
	groups and their generators of homotopy groups of spheres, passes through the comparison with $G_2,$   and then follows the
	surviving classes through successive inclusions. Bott periodicity supplies the endpoint:
	$\pi_{23}(\mathit{SO})\cong\mathbb Z_{(2)},$   reached for $n\geq25.$ 
	
	As usual, in our calculations, we draw freely on the descriptions of
	$\pi_{*}(S^{j})$ in Appendix~\ref{tlqdzj} and on the composition
	relations among generators in Toda's range tabulated in
	\^{O}guchi's table in Appendix~\ref{app:oguchi-relations}.

	\section{Main Theorem and Generators}
	
	The table is arranged by the value of $n.$  A generator with a subscript indicating $n$ denotes the image of the previously
	chosen class under inclusion\uprightsemicolon{} a newly bracketed element of a homotopy group of a sphere denotes a lift. These conventions allow
	successive rows to be compared, rather than read as unrelated abstract groups.
	
	\begin{thm} \label{m1:thm1}
		\textup{All spaces in this theorem are localized at $2.$  The displayed groups are obtained from the integral homotopy groups by tensoring with $\mathbb Z_{(2)}.$  Free cyclic summands become $\mathbb Z_{(2)}$; for finite groups, this localization retains precisely the $2$-primary summands. The homotopy groups are given as follows.   }\par
		\begin{enumerate}[leftmargin=0pt,label={},itemindent=0pt]
			\item \textup{The Cases $3\leq n\leq6.$ }\par
			\smallskip\noindent\(
			\GroupLayout{\pi_{23}(\mathit{SO}(3)) = \mathrm{span}\{[\eta_2]\nu'\bar{\mu}_6,\; [\eta_2]\nu'\eta_6\mu_7\sigma_{16}\}\cong (\mathbb Z/2)^{2}.}{\begin{aligned}\pi_{23}(\mathit{SO}(3)) = \mathrm{span}\{&[\eta_2]\nu'\bar{\mu}_6,\\
					&[\eta_2]\nu'\eta_6\mu_7\sigma_{16}\}\cong (\mathbb Z/2)^{2}.\end{aligned}}
			\)\par
			\smallskip\par\smallskip\noindent\(
			\GroupLayout{\pi_{23}(\mathit{SO}(4)) = \mathrm{span}\{[\eta_2]_4\nu'\bar{\mu}_6,\; [\eta_2]_4\nu'\eta_6\mu_7\sigma_{16},\; [\iota_3]\nu'\bar{\mu}_6,\; [\iota_3]\nu'\eta_6\mu_7\sigma_{16}\}\cong (\mathbb Z/2)^{4}.}{\begin{aligned}\pi_{23}(\mathit{SO}(4)) = \mathrm{span}\{&[\eta_2]_4\nu'\bar{\mu}_6,\; [\eta_2]_4\nu'\eta_6\mu_7\sigma_{16},\\
					&[\iota_3]\nu'\bar{\mu}_6,\; [\iota_3]\nu'\eta_6\mu_7\sigma_{16}\}\cong (\mathbb Z/2)^{4}.\end{aligned}}
			\)\par
			\smallskip\par\smallskip\noindent\(
			\GroupLayout{\pi_{23}(\mathit{SO}(5)) = \mathrm{span}\{[\nu_4\sigma'\mu_{14}],\; [\nu_4(\Sigma\zeta')],\; [\nu_4\eta_7\bar{\varepsilon}_8]\}\cong (\mathbb Z/2)^{3}.}{\begin{aligned}\pi_{23}(\mathit{SO}(5)) = \mathrm{span}\{&[\nu_4\sigma'\mu_{14}],\; [\nu_4(\Sigma\zeta')],\\
					&[\nu_4\eta_7\bar{\varepsilon}_8]\}\cong (\mathbb Z/2)^{3}.\end{aligned}}
			\)\par
			\smallskip\noindent\(
			\GroupLayout{\pi_{23}(\mathit{SO}(6)) \;= {} \mathrm{span}\{[\zeta_5]\sigma_{16},\; [ \nu_5]\bar{\varepsilon}_8,\; [\nu_4\sigma'\mu_{14}]_6,\; [\nu_4(\Sigma\zeta')]_6,\; [\nu_4\eta_7\bar{\varepsilon}_8]_6\}\cong \mathbb Z/8\oplus (\mathbb Z/2)^{4}.}{\begin{aligned}\pi_{23}(\mathit{SO}(6)) \;= {} \mathrm{span}\{&[\zeta_5]\sigma_{16},\; [ \nu_5]\bar{\varepsilon}_8,\; [\nu_4\sigma'\mu_{14}]_6,\\
					&[\nu_4(\Sigma\zeta')]_6,\; [\nu_4\eta_7\bar{\varepsilon}_8]_6\}\cong \mathbb Z/8\oplus (\mathbb Z/2)^{4}.\end{aligned}}
			\)\par
			\item \textup{The Cases $7\leq n\leq9.$ }\par
			\smallskip\noindent\(
			\GroupLayout{\pi_{23}(\mathit{SO}(7)) \;= {} \mathrm{span}\{[\zeta_5]_7\sigma_{16},\; [\nu_4\sigma'\mu_{14}]_7,\; [\nu_4(\Sigma\zeta')]_7,\; [\nu_4\eta_7\bar{\varepsilon}_8]_7,\; [\eta_6]\mu_7\sigma_{16},\; [P(\Sigma\theta)+\nu_6\kappa_9]\}\cong (\mathbb Z/2)^{5}\oplus \mathbb Z/4,}{\begin{aligned}\pi_{23}(\mathit{SO}(7)) \;= {} \mathrm{span}\{&[\zeta_5]_7\sigma_{16},\; [\nu_4\sigma'\mu_{14}]_7,\; [\nu_4(\Sigma\zeta')]_7,\; [\nu_4\eta_7\bar{\varepsilon}_8]_7,\\
					&[\eta_6]\mu_7\sigma_{16},\; [P(\Sigma\theta)+\nu_6\kappa_9]\}\cong (\mathbb Z/2)^{5}\oplus \mathbb Z/4,\end{aligned}}
			\)\par
			\nobreak
			{\raggedleft
				\noindent where $2[P(\Sigma\theta)+\nu_6\kappa_9]\equiv[ \nu_5]_7\bar{\varepsilon}_8
				\bmod
				[\nu_4(\Sigma\zeta')]_7, [\nu_4\eta_7\bar{\varepsilon}_8]_7.$\par}
			
			\par\smallskip\noindent\(
			\GroupLayout{\pi_{23}(\mathit{SO}(8)) \;= {} \mathrm{span}\{[\zeta_5]_8\sigma_{16},\; [\nu_4\sigma'\mu_{14}]_8,\; [\nu_4(\Sigma\zeta')]_8,\; [\nu_4\eta_7\bar{\varepsilon}_8]_8,\; [\eta_6]_8\mu_7\sigma_{16},\; [P(\Sigma\theta)+\nu_6\kappa_9]_8,\; [\iota_7]\sigma'\mu_{14},\; [\iota_7](\Sigma\zeta'),\; [\iota_7]\eta_7\bar{\varepsilon}_8,\; [\iota_7]\mu_7\sigma_{16}\}\cong (\mathbb Z/2)^{5}\oplus \mathbb Z/4\oplus (\mathbb Z/2)^{4}.}{\begin{aligned}\pi_{23}(\mathit{SO}(8)) \;= {} \mathrm{span}\{&[\zeta_5]_8\sigma_{16},\; [\nu_4\sigma'\mu_{14}]_8,\; [\nu_4(\Sigma\zeta')]_8,\; [\nu_4\eta_7\bar{\varepsilon}_8]_8,\; [\eta_6]_8\mu_7\sigma_{16},\; [P(\Sigma\theta)+\nu_6\kappa_9]_8,\\
					&[\iota_7]\sigma'\mu_{14},\; [\iota_7](\Sigma\zeta'),\; [\iota_7]\eta_7\bar{\varepsilon}_8,\; [\iota_7]\mu_7\sigma_{16}\}\cong (\mathbb Z/2)^{5}\oplus \mathbb Z/4\oplus (\mathbb Z/2)^{4}.\end{aligned}}
			\)\par
			\smallskip\par\smallskip\noindent\(
			\GroupLayout{\pi_{23}(\mathit{SO}(9)) = \mathrm{span}\{[\zeta_5]_9\sigma_{16},\; [\iota_7]_9\mu_7\sigma_{16},\; [P(\Sigma\theta)+\nu_6\kappa_9]_9\}\cong (\mathbb Z/2)^{2}\oplus \mathbb Z/4,}{\begin{aligned}\pi_{23}(\mathit{SO}(9)) = \mathrm{span}\{&[\zeta_5]_9\sigma_{16},\; [\iota_7]_9\mu_7\sigma_{16},\\
					&[P(\Sigma\theta)+\nu_6\kappa_9]_9\}\cong (\mathbb Z/2)^{2}\oplus \mathbb Z/4,\end{aligned}}
			\)\par
			\nobreak
			{\raggedleft
				\noindent where
				$2[P(\Sigma\theta)+\nu_6\kappa_9]_9=[ \nu_5]_9\bar{\varepsilon}_8.$ \par}
			
			\item \textup{The Cases $10\leq n\leq14.$ }\par
			\smallskip\noindent\(
			\GroupLayout{\pi_{23}(\mathit{SO}(10)) = \mathrm{span}\{[\sigma^2_9],\; [\bar{\nu}_9\nu^2_{17}],\; [P(\Sigma\theta)+\nu_6\kappa_9]_{10}\}\cong \mathbb Z/32\oplus (\mathbb Z/2)^{2}.}{\begin{aligned}\pi_{23}(\mathit{SO}(10)) = \mathrm{span}\{&[\sigma^2_9],\; [\bar{\nu}_9\nu^2_{17}],\\
					&[P(\Sigma\theta)+\nu_6\kappa_9]_{10}\}\cong \mathbb Z/32\oplus (\mathbb Z/2)^{2}.\end{aligned}}
			\)\par
			\smallskip\par\smallskip\noindent\(
			\GroupLayout{\pi_{23}(\mathit{SO}(11)) = \mathrm{span}\{[\sigma^2_9]_{11},\; [P(\Sigma\theta)+\nu_6\kappa_9]_{11}\}\cong (\mathbb Z/2)^{2}.}{\begin{aligned}\pi_{23}(\mathit{SO}(11)) = \mathrm{span}\{&[\sigma^2_9]_{11},\\
					&[P(\Sigma\theta)+\nu_6\kappa_9]_{11}\}\cong (\mathbb Z/2)^{2}.\end{aligned}}
			\)\par
			\smallskip\par\smallskip\noindent\(
			\GroupLayout{\pi_{23}(\mathit{SO}(12)) = \mathrm{span}\{[\theta'],\; [\sigma^2_9]_{12},\; [P(\Sigma\theta)+\nu_6\kappa_9]_{12}\}\cong (\mathbb Z/2)^{3},}{\begin{aligned}\pi_{23}(\mathit{SO}(12)) = \mathrm{span}\{&[\theta'],\; [\sigma^2_9]_{12},\\
					&[P(\Sigma\theta)+\nu_6\kappa_9]_{12}\}\cong (\mathbb Z/2)^{3},\end{aligned}}
			\)\par
			\nobreak
			{\raggedleft
				\noindent where $\Delta(\Sigma\theta')
				=[P(\Sigma\theta)+\nu_6\kappa_9]_{12}.$\par}
			
			\par\smallskip\noindent\(
			\GroupLayout{\pi_{23}(\mathit{SO}(13)) = \mathrm{span}\{[32[\iota_{12}, \iota_{12}]],\; [\sigma^2_9]_{13}\}\cong \mathbb Z_{(2)}\oplus \mathbb Z/2.}{\begin{aligned}\pi_{23}(\mathit{SO}(13)) = \mathrm{span}\{&[32[\iota_{12}, \iota_{12}]],\\
					&[\sigma^2_9]_{13}\}\cong \mathbb Z_{(2)}\oplus \mathbb Z/2.\end{aligned}}
			\)\par
			\smallskip\par\smallskip\noindent\(
			\GroupLayout{\pi_{23}(\mathit{SO}(14)) = \mathrm{span}\{[32[\iota_{12}, \iota_{12}]]_{14},\; [\sigma^2_9]_{14}\}\cong \mathbb Z_{(2)}\oplus \mathbb Z/2.}{\begin{aligned}\pi_{23}(\mathit{SO}(14)) = \mathrm{span}\{&[32[\iota_{12}, \iota_{12}]]_{14},\\
					&[\sigma^2_9]_{14}\}\cong \mathbb Z_{(2)}\oplus \mathbb Z/2.\end{aligned}}
			\)\par
			\item \textup{The Cases $n\geq15.$ }\par
			\smallskip\noindent\(
			\GroupLayout{\pi_{23}(\mathit{SO}(15)) = \mathrm{span}\{[32[\iota_{12}, \iota_{12}]]_{15},\; [\sigma^2_9]_{15},\; [\eta_{14}\varepsilon_{15}],\; [\eta_{14}\sigma_{15}]\eta_{22}\}\cong \mathbb Z_{(2)}\oplus (\mathbb Z/2)^{3}.}{\begin{aligned}\pi_{23}(\mathit{SO}(15)) = \mathrm{span}\{&[32[\iota_{12}, \iota_{12}]]_{15},\; [\sigma^2_9]_{15},\; [\eta_{14}\varepsilon_{15}],\\
					&[\eta_{14}\sigma_{15}]\eta_{22}\}\cong \mathbb Z_{(2)}\oplus (\mathbb Z/2)^{3}.\end{aligned}}
			\)\par
			\smallskip\par\smallskip\noindent\(
			\GroupLayout{\pi_{23}(\mathit{SO}(16)) \;= {} \mathrm{span}\{[32[\iota_{12},\iota_{12}]]_{16},\; [\sigma^2_9]_{16},\; [\eta_{14}\varepsilon_{15}]_{16},\; [\eta_{14}\sigma_{15}]_{16}\eta_{22},\; [\eta_{15}]\sigma_{16},\; [\varepsilon_{15}]\}\cong \mathbb Z_{(2)}\oplus (\mathbb Z/2)^{5}.}{\begin{aligned}\pi_{23}(\mathit{SO}(16)) \;= {} \mathrm{span}\{&[32[\iota_{12},\iota_{12}]]_{16},\; [\sigma^2_9]_{16},\; [\eta_{14}\varepsilon_{15}]_{16},\\
					&[\eta_{14}\sigma_{15}]_{16}\eta_{22},\; [\eta_{15}]\sigma_{16},\; [\varepsilon_{15}]\}\cong \mathbb Z_{(2)}\oplus (\mathbb Z/2)^{5}.\end{aligned}}
			\)\par
			\smallskip\par\smallskip\noindent\(
			\GroupLayout{\pi_{23}(\mathit{SO}(17)) = \mathrm{span}\{[32[\iota_{12}, \iota_{12}]]_{17},\; [\sigma^2_9]_{17},\; [\eta_{15}]_{17}\sigma_{16},\; [\varepsilon_{15}]_{17}\}\cong \mathbb Z_{(2)}\oplus (\mathbb Z/2)^{3}.}{\begin{aligned}\pi_{23}(\mathit{SO}(17)) = \mathrm{span}\{&[32[\iota_{12}, \iota_{12}]]_{17},\; [\sigma^2_9]_{17},\\
					&[\eta_{15}]_{17}\sigma_{16},\; [\varepsilon_{15}]_{17}\}\cong \mathbb Z_{(2)}\oplus (\mathbb Z/2)^{3}.\end{aligned}}
			\)\par
			\smallskip\par\smallskip\noindent\(
			\GroupLayout{\pi_{23}(\mathit{SO}(18)) = \mathrm{span}\{[32[\iota_{12}, \iota_{12}]]_{18},\; [\sigma^2_9]_{18},\; [\varepsilon_{15}]_{18}\}\cong \mathbb Z_{(2)}\oplus (\mathbb Z/2)^{2},\;\text{where}\;\Delta(\nu^2_{18}) = [\varepsilon_{15}]_{18}.}{\begin{aligned}\pi_{23}(\mathit{SO}(18)) = \mathrm{span}\{&[32[\iota_{12}, \iota_{12}]]_{18},\; [\sigma^2_9]_{18},\\
					&[\varepsilon_{15}]_{18}\}\cong \mathbb Z_{(2)}\oplus (\mathbb Z/2)^{2},\;\text{where}\;\Delta(\nu^2_{18}) = [\varepsilon_{15}]_{18}.\end{aligned}}
			\)\par
			\smallskip\par\smallskip\noindent\(
			\GroupLayout{\pi_{23}(\mathit{SO}(n)) = \mathrm{span}\{[32[\iota_{12}, \iota_{12}]]_n,\; [\sigma^2_9]_n\}\cong \mathbb Z_{(2)}\oplus \mathbb Z/2,}{\begin{aligned}\pi_{23}(\mathit{SO}(n)) = \mathrm{span}\{&[32[\iota_{12}, \iota_{12}]]_n,\\
					&[\sigma^2_9]_n\}\cong \mathbb Z_{(2)}\oplus \mathbb Z/2,\end{aligned}}
			\)\par
			\nobreak
			{\raggedleft
				\noindent for $n = 19, 20, 21,$   where $\Delta\nu_{21}=[\sigma^2_9]_{21}.$\par}
			
			\par\smallskip\noindent\(
			\GroupLayout{\pi_{23}(\mathit{SO}(24)) = \mathrm{span}\{\Delta(\iota_{24}),\; [32[\iota_{12},\iota_{12}]]_{24}\}\cong (\mathbb Z_{(2)})^{2}.}{\begin{aligned}\pi_{23}(\mathit{SO}(24)) = \mathrm{span}\{&\Delta(\iota_{24}),\\
					&[32[\iota_{12},\iota_{12}]]_{24}\}\cong (\mathbb Z_{(2)})^{2}.\end{aligned}}
			\)\par
			\smallskip\par\smallskip\noindent\(
			\FitMath{\pi_{23}(\mathit{SO}(22+n)) = \mathrm{span}\{[32[\iota_{12}, \iota_{12}]]_{22+n}\}\cong \mathbb Z_{(2)},\; (n\geq0, n\neq2).}
			\)\par
		\end{enumerate}
	\end{thm}
	
	\section{The Cases $3\leq n\leq6$}
	\subsection{The Cases $3\leq n\leq5$}
	
	We start with the results on the homotopy groups of spheres for the case $n=3$ and the resulting two contributions for $n=4.$  The passage to
	the case $n=5$ then illustrates the use of the small-$n$ comparison lemma: determining the image of the
	inclusion simplifies the projection to the sphere.
	
	We have $\pi^3_{23} = \frac{\mathbb Z}{2}\{\nu'\bar{\mu}_6\}
	\oplus\frac{\mathbb Z}{2}\{\nu'\eta_6\mu_7\sigma_{16}\}.$  Thus we obtain
	\begin{gather*}
		\begin{aligned}
			&\pi_{23}(\mathit{SO}(3))={}\frac{\mathbb Z}{2}\{[\eta_2]\nu'\bar{\mu}_6\} \oplus\frac{\mathbb Z}{2}\{[\eta_2]\nu'\eta_6\mu_7\sigma_{16}\}
		\end{aligned}
		\intertext{and}
		\begin{aligned}
			&\pi_{23}(\mathit{SO}(4))={}\frac{\mathbb Z}{2}\{[\eta_2]_4\nu'\bar{\mu}_6\} \oplus\frac{\mathbb Z}{2}\{[\eta_2]_4\nu'\eta_6\mu_7\sigma_{16}\} \oplus\frac{\mathbb Z}{2}\{[\iota_3]\nu'\bar{\mu}_6\} \oplus\frac{\mathbb Z}{2}\{[\iota_3]\nu'\eta_6\mu_7\sigma_{16}\}.
		\end{aligned}
	\end{gather*}
	Applying Lemma~\ref{m1:236} for $\pi_{23}^3,$   we obtain
	\begin{equation}
		\begin{aligned}
			&i_{5 \ast}(\pi_{23}(\mathit{SO}(4)))=0.\label{m1:i5R423}
		\end{aligned}
	\end{equation}
	
	We note from \cite[Theorem 5.1]{MT2} that
	$\pi_{23}(\mathit{Sp}(2)) = \frac{\mathbb Z}{2}\{[\sigma'\mu_{14}]\}
	\oplus\frac{\mathbb Z}{2}\{[\Sigma\zeta']\}\oplus\frac{\mathbb Z}{2}\{[\eta_7\bar{\varepsilon}_8]\}.$ 
	Therefore, by \cite[Lemma 1.2]{HKM},
	we obtain
	\[\begin{aligned}
		&\pi_{23}(\mathit{SO}(5))={}\frac{\mathbb Z}{2}\{[\nu_4\sigma'\mu_{14}]\} \oplus\frac{\mathbb Z}{2}\{[\nu_4(\Sigma\zeta')]\} \oplus\frac{\mathbb Z}{2}\{[\nu_4\eta_7\bar{\varepsilon}_8]\}.
	\end{aligned}\]
	We also give a direct proof that $\mathrm{ord}[\nu_4\sigma'\mu_{14}]=2.$ 
	Since $\mu_{14}\in\{\eta_{14},2\iota_{15},8\sigma_{15}\}_1
	\bmod \eta_{14}\varepsilon_{15}, \nu^3_{14}$ \eqref{m1:mu4} and \eqref{m1:i5R423},
	we obtain
	\begin{equation}
		\begin{aligned}
			&[\nu_4\sigma'\mu_{14}] \in \{[\nu_4\sigma'\eta_{14}], 2\iota_{15},8\sigma_{15}\}_1 \bmod [\nu_4(\Sigma\zeta')], [\nu_4\eta_7\bar{\varepsilon}_8].\label{m1:mR524}
		\end{aligned}
	\end{equation}
	Since $2[\nu_4 \sigma'\mu_{14}] \in \{[\nu_4\sigma'\eta_{14}], 2\iota_{15},
	8\sigma_{15}\}_1 \circ 2\iota_{23}=[\nu_4\sigma'\eta_{14}]\circ \{2\iota_{15},
	8\sigma_{15}, 2\iota_{22}\}_1 = 0$ \eqref{m1:2i8s2i},
	it follows that $\mathrm{ord} [\nu_4\sigma'\mu_{14}]=2.$ 
	
	\subsection{The Case $n=6$}
	
	For $n=6,$   the elements of homotopy groups of spheres that lift are determined by the kernel of the next boundary. The
	bracket representative of $\bar\zeta_5$ is recalled because its multiple controls a boundary relation
	needed in this step.
	
	We see that $\pi^5_{24}=\frac{\mathbb Z}{8}\{\bar{\zeta}_5\}\oplus\frac{\mathbb Z}{2}\{\nu_5\mu_8\sigma_{17}\},
	\pi^5_{23}=
	\frac{\mathbb Z}{8}\{\zeta_5\sigma_{16}\}\oplus\frac{\mathbb Z}{2}\{\nu_5\bar{\varepsilon}_8\}\oplus\frac{\mathbb Z}{2}\{\eta_5\bar{\mu}_6\}$
	and
	$$\operatorname{ker}(\Delta\colon \pi^5_{23}\to \pi_{22}(\mathit{SO}(5)))=\mathrm{span}\{\zeta_5\sigma_{16}, \nu_5\bar{\varepsilon}_8\}$$ from \cite[Theorem 0.2(i), (3.2)]{HKM}.
	
	Note from \cite[p.~137]{T} that
	\begin{equation}
		\begin{aligned}
			&\bar{\zeta}_5\in\{\zeta_5, 8\iota_{16}, 2\sigma_{16}\}_1\subseteq\pi^5_{24}.\label{m1:bz5}
		\end{aligned}
	\end{equation}
	Since $\pi_{16}(\mathit{SO}(5))\cong(\mathbb Z/2)^2$ \cite[Proposition 4.1]{Ka}, we deduce
	\[\begin{aligned}
		&\Delta\bar{\zeta}_5\in\{0, 8\iota_{15}, 2\sigma_{15}\}=\pi_{16}(\mathit{SO}(5))\circ 2\sigma_{16}=0.
	\end{aligned}\]
	Recall that $8[\ck_{5}]=0,\;2\bar\lt_{8}=0$ and that $[ \nu_5]$ is well-defined (see \cite{Ka}). Therefore, using the exact sequence $(23)_5,$   we obtain
	\begin{gather}
		\begin{aligned}
			&\pi_{23}(\mathit{SO}(6))={}\mathrm{span}\{[\zeta_5]\sigma_{16}, [ \nu_5]\bar{\varepsilon}_8, [\nu_4\sigma'\mu_{14}]_6, [\nu_4(\Sigma\zeta')]_6, [\nu_4\eta_7\bar{\varepsilon}_8]_6\}\cong\mathbb Z/8\oplus(\mathbb Z/2)^4.\label{m1:R623}
		\end{aligned}
	\end{gather}
	
	\section{The Cases $7\leq n\leq9$}
	
	\subsection{The Case $n=7$ and the Comparison Fibrations}
	
	The case $n=7$ introduces a comparison with the exceptional group $G_2.$  The diagram below places its
	fibration beside the special orthogonal fibration over the same sphere, so that their connecting images and
	chosen lifts can be compared.
	
	We need the following lemma.
	
	\begin{lem}
		Let $M^n=S^{n-1}\cup_{2\iota_{n-1}}e^n$ be the Moore
		space of type $(\mathbb Z/2,n-1),$   with inclusion
		$i'_n:S^{n-1}\to M^n$ and projection $p'_n:M^n\to S^n.$ 
		The following statements hold.
		\begin{itemize}
			\item[\textup{(1)}]
			We have $\pi_{8}(\mathit{SO}(7))
			=\mathrm{span}\{[ \nu_5]_7,[\eta_6]\eta_7\}
			\cong(\mathbb Z/2)^2$ and
			\begin{align*}
				\pi_{16}(\mathit{SO}(7))
				&=\mathrm{span}\{[\zeta_5]_7,
				[\nu_4\sigma'\eta_{14}]_7\eta_{15},
				[  \nu^2_4]_7\nu^2_{10},[ \nu_5]_7\bar{\nu}_8,\\
				&\phantom{=\mathrm{span}\{}[ \nu_5]_7\varepsilon_8,
				[\eta_6]\eta_7\varepsilon_8,[\eta_6]\mu_7\}
				\cong(\mathbb Z/2)^7;\\[3pt]
				\pi_{18}(\mathit{SO}(7))
				&=\mathrm{span}\{[2[\iota_6,\iota_6]]\sigma_{11},
				[\nu_4\zeta_7]_7,[ \nu_5]_7\sigma_8\nu_{15}\}\\
				&\cong\mathbb Z/16\oplus\mathbb Z/8\oplus\mathbb Z/2.
			\end{align*}
			Moreover,
			\begin{equation}
				\begin{aligned}
					[\zeta_5]\sigma_{16}
					&\in\{[ \nu_5]\sigma_8,8\iota_{15},2\sigma_{15}\}_1\\
					&\quad\bmod\ 2[\zeta_5]\sigma_{16},
					[ \nu_5]\bar{\varepsilon}_8,
					[\nu_4\sigma'\mu_{14}]_6, [\nu_4(\Sigma\zeta')]_6,
					[\nu_4\eta_7\bar{\varepsilon}_8]_6.
				\end{aligned}
				\label{m1:zesig1}
			\end{equation}
			The boundary homomorphisms in $(23)_6$ and $(24)_6$ satisfy
			\begin{align*}
				&\operatorname{ker}(\Delta\colon \pi^6_{23}
				\longrightarrow\pi_{22}(\mathit{SO}(6)))=\mathrm{span}\{P(\Sigma\theta)+\nu_6\kappa_9,
				\eta_6\mu_7\sigma_{16}\}\cong(\mathbb Z/2)^2;\\[3pt]
				&\Delta(\pi^6_{24})=\mathrm{span}\{2[\zeta_5]\sigma_{16}\};\\[3pt]
				&\Delta(\pi^6_{24}\circ p'_{24})=0.
			\end{align*}
			
			\item[\textup{(2)}]
			For the lift $[P(\Sigma\theta)+\nu_6\kappa_9]
			\in\pi_{23}(\mathit{SO}(7)),$   we have
			\begin{align*}
				2[P(\Sigma\theta)+\nu_6\kappa_9]
				&\equiv0\bmod\ [ \nu_5]_7\bar{\varepsilon}_8,
				[\nu_4\sigma'\mu_{14}]_7, [\nu_4(\Sigma\zeta')]_7,
				[\nu_4\eta_7\bar{\varepsilon}_8]_7.
			\end{align*}
			Consequently, this lift has order either $2$ or $4.$ 
			
			\item[\textup{(3)}]
			Choose representatives
			\begin{align*}
				\overline{P(\Sigma\theta)}
				&\in\{P(\Sigma\theta),2\iota_{23},p'_{23}\}_2,\;\;
				\overline{\nu_6\kappa_9}
				\in\{\nu_6\kappa_9,2\iota_{23},p'_{23}\}_2.
			\end{align*}
			Then \vspace{-1.4em}
			\begin{align*}
				&\{P(\Sigma\theta)+\nu_6\kappa_9,
				2\iota_{23},p'_{23}\}_2\subseteq\overline{P(\Sigma\theta)}
				+\overline{\nu_6\kappa_9}+\pi^6_{24}\circ p'_{24};\\[3pt]
				&(p_6)_*\Delta\bigl(\overline{P(\Sigma\theta)}
				+\overline{\nu_6\kappa_9}\bigr)
				=\nu_5\bar{\varepsilon}_8p'_{23}.
			\end{align*}
			
			\item[\textup{(4)}]
			We have $\mathrm{ord}[P(\Sigma\theta)+\nu_6\kappa_9]=4.$ 
		\end{itemize}
	\end{lem}
	
	\begin{proof}
		\begin{itemize}
			\item[\textup{(1)}]
			The three groups of $\mathit{SO}(7)$ follow from
			\cite[Table 2, Proposition 4.1]{Ka},
			\cite[Theorem 2]{KM1} and the relation
			$\bar{\nu}_7\nu_{15}=\eta_7\sigma_8\nu_{15}.$ 
			Lemma~\ref{m1:zesig}(4), together with the description of
			$\pi_{23}(\mathit{SO}(5)),$   gives \eqref{m1:zesig1}.  Recall that
			$\pi^7_{23}=\mathrm{span}\{\sigma'\mu_{14},\Sigma\zeta',
			\mu_7\sigma_{16},\eta_7\bar{\varepsilon}_8\}
			\cong(\mathbb Z/2)^4.$ 
			The asserted kernel in $(23)_6$ is given in
			\cite[p.~12-3]{HKM}; see also \cite[p.~12]{HKM}.
			Moreover,
			$$\pi^6_{24}=\mathrm{span}\{\zeta_6\sigma_{17},
			P(\Sigma\theta)\eta_{23},\eta_6\bar{\mu}_7\}
			\cong\mathbb Z/8\oplus(\mathbb Z/2)^2.$$
			\noindent As a consequence of           Lemma~\ref{m1:zeta1},  we get   $\Delta(\zeta_6\sigma_{17})=2[\zeta_5]\sigma_{16}.$     Note that $\Delta(\eta_6\bar{\mu}_7)=0.$  Since
			\begin{align*}
				P(\Sigma\theta)\eta_{23}
				&=(P(\Sigma\theta)+\nu_6\kappa_9)\eta_{23},
			\end{align*}
			the kernel formula gives
			$\Delta(P(\Sigma\theta)\eta_{23})=0.$ 
			The values on these three generators prove the asserted
			formula for $\Delta(\pi^6_{24}).$ 
			Moreover,  $2\iota_{23}\circ p'_{23}=0$ gives
			\begin{align*}
				\Delta(\pi^6_{24}\circ p'_{24})
				&\subseteq\mathrm{span}\{2[\zeta_5]\sigma_{16}
				\circ p'_{23}\}=0.
			\end{align*}
			Thus the boundary also vanishes on the indeterminacy
			that occurs for extensions to $M^{24}.$ 
			
			\item[\textup{(2)}]
			Consider the following diagram with exact rows, which
			commutes up to sign:
			\[
			\begin{tikzcd}[column sep=1.5em,row sep=large]
				\pi^{6}_{24}
				\arrow[r,"\Delta"] \arrow[d,"\Sigma^{7}"']
				& \pi_{23}(\mathit{SO}(6))
				\arrow[r,"{(i_7)_*}"] \arrow[d,"J"']
				& \pi_{23}(\mathit{SO}(7))
				\arrow[r,"{(p_7)_*}"] \arrow[d,"J"']
				& \pi^{6}_{23}
				\arrow[r,"\Delta"] \arrow[d,"\Sigma^{7}"']
				& \pi_{22}(\mathit{SO}(6)) \arrow[d,"J"]
				\\
				\pi^{13}_{31} \arrow[r,"P"']
				& \pi^{6}_{29} \arrow[r,"\Sigma"']
				& \pi^{7}_{30} \arrow[r,"H"']
				& \pi^{13}_{30} \arrow[r,"P"']
				& \pi^{6}_{28}.
			\end{tikzcd}
			\]
			Using this diagram, the kernel formula in \textup{(1)}
			and \eqref{m1:R623}, we obtain
			\begin{align*}
				2[P(\Sigma\theta)+\nu_6\kappa_9]
				&\equiv0\bmod\ [ \nu_5]_7\bar{\varepsilon}_8,
				[\nu_4\sigma'\mu_{14}]_7, [\nu_4(\Sigma\zeta')]_7,
				[\nu_4\eta_7\bar{\varepsilon}_8]_7.
			\end{align*}
			Each generator of the modulus is the image of an element
			of order $2$ in \eqref{m1:R623}. Thus every element of
			this modulus has order dividing $2,$   and hence
			$4[P(\Sigma\theta)+\nu_6\kappa_9]=0.$ 
			On the other hand, $    (p_7)_*[P(\Sigma\theta)+\nu_6\kappa_9]
			=P(\Sigma\theta)+\nu_6\kappa_9$
			has order $2$ by \textup{(1)}. The lift is therefore
			nonzero, and its order is either $2$ or $4.$ 
			
			\item[\textup{(3)}]
			Since $\Sigma^2[M^{22},S^{21}]
			=\mathrm{span}\{\eta_{23}p'_{24}\},$   the indeterminacy
			of $\{P(\Sigma\theta),2\iota_{23},p'_{23}\}_2$ is
			\begin{align*}
				&P(\Sigma\theta)\circ\Sigma^2[M^{22},S^{21}]
				+\pi^6_{24}\circ p'_{24}=\pi^6_{24}\circ p'_{24}.
			\end{align*}
			The same calculation with $\nu_6\kappa_9$ in place of
			$P(\Sigma\theta)$ gives the same indeterminacy for
			$\{\nu_6\kappa_9,2\iota_{23},p'_{23}\}_2.$ 
			The addition property of Toda brackets therefore gives
			\begin{align*}
				\{P(\Sigma\theta)+\nu_6\kappa_9,
				2\iota_{23},p'_{23}\}_2&\subseteq\{P(\Sigma\theta),2\iota_{23},p'_{23}\}_2
				+\{\nu_6\kappa_9,2\iota_{23},p'_{23}\}_2\\
				&=\overline{P(\Sigma\theta)}
				+\overline{\nu_6\kappa_9}+\pi^6_{24}\circ p'_{24}.
			\end{align*}
			Finally, applying $(p_6)_*(r_3)_*$ to
			\cite[Theorem 4.8]{MMOh} gives
			\begin{align*}
				(p_6)_*\Delta\bigl(\overline{P(\Sigma\theta)}
				+\overline{\nu_6\kappa_9}\bigr)
				&=\nu_5\bar{\varepsilon}_8p'_{23}.
			\end{align*}
			The vanishing of $\Delta(\pi^6_{24}\circ p'_{24})$
			in \textup{(1)} shows that this value is independent of the choices
			of the two bracket representatives.
			
			\item[\textup{(4)}]
			By \textup{(2)}, it remains to exclude order $2.$ 
			Suppose, for a contradiction, that
			$\mathrm{ord}[P(\Sigma\theta)+\nu_6\kappa_9]=2.$ 
			There would then be an extension
			\begin{align*}
				\overline{[P(\Sigma\theta)+\nu_6\kappa_9]}
				&\in\{[P(\Sigma\theta)+\nu_6\kappa_9],
				2\iota_{23},p'_{23}\}_2\\
				&\subseteq[M^{24},\mathit{SO}(7)].
			\end{align*}
			By naturality, its projection would belong to
			$\{P(\Sigma\theta)+\nu_6\kappa_9,2\iota_{23},p'_{23}\}_2.$ 
			Since $\Delta\circ(p_7)_*=0$ and
			$\Delta(\pi^6_{24}\circ p'_{24})=0$ by \textup{(1)},
			the relations in \textup{(3)} would then give
			\begin{align*}
				0&=(p_6)_*\Delta\bigl(\overline{P(\Sigma\theta)}
				+\overline{\nu_6\kappa_9}\bigr)
				=\nu_5\bar{\varepsilon}_8p'_{23}.
			\end{align*}
			The cofiber sequence
			$M^{23}\xrightarrow{p'_{23}}S^{23}
			\xrightarrow{2\iota_{23}}S^{23}$ would therefore imply
			\begin{align*}
				\nu_5\bar{\varepsilon}_8&\in2\pi^5_{23}
				=\mathrm{span}\{2(\zeta_5\sigma_{16})\}.
			\end{align*}
			This  contradicts the decomposition
			$\pi^5_{23}=\frac{\mathbb Z}{8}\{\zeta_5\sigma_{16}\}
			\oplus\frac{\mathbb Z}{2}\{\nu_5\bar{\varepsilon}_8\}
			\oplus\frac{\mathbb Z}{2}\{\eta_5\bar{\mu}_6\}.$  Thus the lift cannot have
			order $2,$   and \textup{(2)} yields
			$\mathrm{ord}[P(\Sigma\theta)+\nu_6\kappa_9]=4.$ 
		\end{itemize}
	\end{proof}
	
	\subsection{The Extension for $n=7$ and Generator Relations}
	
	Recall from \cite[Theorem 6.1]{Mi2} and \cite[Theorem 1.3]{MMOh} that
	\[\begin{aligned}
		&\pi_{23}(G_2)={}\mathrm{span}\{\langle P(\Sigma\theta)+\nu_6\kappa_9\rangle, \langle\eta_6\mu_7\rangle\sigma_{16}\}\cong\mathbb Z/4\oplus\mathbb Z/2,
	\end{aligned}\]
	where $2\langle P(\Sigma\theta)+\nu_6\kappa_9\rangle={i_G}_*[\nu_5\bar{\varepsilon}_8].$ 
	Since
	\[\begin{aligned}
		&\mathcal{T}_{*}\left<P(\Sigma\theta) + \nu_6\kappa_9\right> \equiv [P(\Sigma\theta) + \nu_6\kappa_9] \bmod {i_7}_* \pi_{23}(\mathit{SO}(6)),
	\end{aligned}\]
	we see that $\mathrm{ord}\left<P(\Sigma\theta) + \nu_6\kappa_9\right>=4.$ 
	Moreover, we find
	\[\begin{aligned}
		&\mathcal{T}_{*}{i_G}_*[\nu_5\bar{\varepsilon}_8] \equiv {i_7}_*{r_3}_* [\nu_5\bar{\varepsilon}_8] \equiv [ \nu_5]_7\bar{\varepsilon}_8 \bmod [\nu_4\sigma' \mu_{14}]_7, [\nu_4(\Sigma\zeta')]_7, [\nu_4\eta_7\bar{\varepsilon}_8]_7.
	\end{aligned}\]
	
	Proposition \ref{so7fj} tells us that
	
	\[ \begin{aligned}
		&\pi_{k}(\mathit{SO}(7)) \;\cong\; \mathcal{T}_{*}\pi_k (G_2) \oplus [\eta_6]\circ\pi^7_k,
	\end{aligned} \] for all $k\ge 2.$

	The comparison above supplies the input for the extension in the case $n=7.$  The next proposition gives its
	decomposition\uprightsemicolon{} the following argument determines the multiple of the new lift and its relation to the
	inherited classes.
	
	\begin{prop} \label{m1:R723}
		The structure of $\pi_{23}(\mathit{SO}(7))$ is as follows.
		\[\GroupLayout{\pi_{23}(\mathit{SO}(7)) \;= {} \mathrm{span}\{[\zeta_5]_7\sigma_{16},\; [\nu_4\sigma'\mu_{14}]_7,\; [\nu_4(\Sigma\zeta')]_7,\; [\nu_4\eta_7\bar{\varepsilon}_8]_7,\; [\eta_6]\mu_7\sigma_{16},\; [P(\Sigma\theta)+\nu_6\kappa_9]\}\cong (\mathbb Z/2)^{5}\oplus \mathbb Z/4,}{\begin{aligned}\pi_{23}(\mathit{SO}(7)) \;= {} \mathrm{span}\{&[\zeta_5]_7\sigma_{16},\; [\nu_4\sigma'\mu_{14}]_7,\; [\nu_4(\Sigma\zeta')]_7,\; [\nu_4\eta_7\bar{\varepsilon}_8]_7,\\
				&[\eta_6]\mu_7\sigma_{16},\; [P(\Sigma\theta)+\nu_6\kappa_9]\}\cong (\mathbb Z/2)^{5}\oplus \mathbb Z/4,\end{aligned}}\]
		where $2[P(\Sigma\theta)+\nu_6\kappa_9]
		\equiv[ \nu_5]_7\bar{\varepsilon}_8\ \bmod\
		[\nu_4\sigma'\mu_{14}]_7, [\nu_4(\Sigma\zeta')]_7, [\nu_4\eta_7\bar{\varepsilon}_8]_7.$ 
		\qed
	\end{prop}
	
	The relation in Proposition~\ref{m1:R723} will be refined after we determine the group structure of $\pi_{24}(\mathit{SO}(7)).$     By \eqref{m1:Ogu0} and \eqref{m1:4z5}, we deduce
	\begin{equation}
		\begin{aligned}
			&\eta_5\varepsilon_6\mu_{14} = 4\zeta_5\sigma_{16}.\label{m1:4z5sig}
		\end{aligned}
	\end{equation}
	
	\noindent   As a consequence of           \cite[Lemma 2.3(ii)]{KM1}, we  know
	\begin{equation}
		\begin{aligned}
			&[\mathrm{c}_{\mathrm{KM}}(\bar{\nu}_6+\varepsilon_6)]\eta_{14} \equiv [\eta_6](\bar{\nu}_7+\varepsilon_7) + [ \nu_5]_7\sigma_8 \bmod [\nu_4\sigma'\eta_{14}]_7.\label{m1:tr2}
		\end{aligned}
	\end{equation}
	
	We pause briefly to collect the relations
	that arise in the present setting.
	\begin{lem} \label{lem6d4}
		\begin{enumerate}
			\item $[\nu_4\sigma'\eta_{14}]\bar{\nu}_{15}
			= [\nu_4\eta_7\bar{\varepsilon}_8] = [  \nu^2_4]\sigma_{10}\nu^2_{17}$ and
			$[\nu_4\sigma'\eta_{14}]\varepsilon_{15}
			= [\nu_4(\Sigma\zeta')].$ 
			\item $[\eta_5\varepsilon_6]\nu^3_{14}
			= [\nu_4\eta_7\bar{\varepsilon}_8]_6,$  
			$[\eta_5\varepsilon_6]\eta_{14}\varepsilon_{15}
			= [\nu_4(\Sigma\zeta')]_6,$   $[\eta_6]\eta_7\bar{\varepsilon}_8
			= [\nu_4\eta_7\bar{\varepsilon}_8]_7$  and \\
			$[\eta_6]\Sigma\zeta' = [\nu_4(\Sigma\zeta')]_7.$ 
			
			\item $[\eta_5\varepsilon_6]\mu_{14}+[\nu_4\sigma'\mu_{14}]_6
			\equiv 4[\zeta_5]\sigma_{16} \
			\bmod [\nu_4(\Sigma\zeta')]_6, [\nu_4\eta_7\bar{\varepsilon}_8]_6$ and \\
			$[\eta_6]\sigma'\mu_{14}=[\eta_5\varepsilon_6]_7\mu_{14}
			\equiv[\nu_4\sigma'\mu_{14}]_7\ \bmod \
			[\nu_4\eta_7\bar{\varepsilon}_8]_7, [\nu_4(\Sigma\zeta')]_7.$ 
		\end{enumerate}
	\end{lem}
	\begin{proof}
		\begin{enumerate}[label=\textup{(\arabic*)},leftmargin=2em]
			\item We recall
			${i_5}_*\pi_{23}(\mathit{SO}(4))=0$ \eqref{m1:i5R423}. Therefore, by Lemma~\ref{m1:236}(1) and
			\eqref{m1:epep}, we obtain
			$[\nu_4\sigma'\eta_{14}]\bar{\nu}_{15} = [\nu_4\sigma'\eta_{14}\bar{\nu}_{15}]
			= [\nu_4\eta_7\bar{\varepsilon}_8]$ and
			$[  \nu^2_4]\sigma_{10}\nu^2_{17} = [\nu^2_4\sigma_{10}\nu^2_{17}]
			= [\nu_4\eta_7\bar{\varepsilon}_8].$ 
			Similarly, the standard composition relations give the last relation in (1).
			
			\item By \cite[Lemma 1.1(v)]{KM1}, we record
			\begin{equation}
				\begin{aligned}
					&[\eta_6]\sigma'=4[\mathrm{c}_{\mathrm{KM}}(\bar{\nu}_6+\varepsilon_6)]+[ \nu_5]_7\nu^2_8+[\eta_5\varepsilon_6]_7.\label{m1:et6s'}
				\end{aligned}
			\end{equation}
			In view of \cite[Lemma 1.2(iii)]{KM1},
			we know
			\begin{equation}
				\begin{aligned}
					&[\eta_5\varepsilon_6]\eta_{14} =4[ \nu_5]\sigma_8 + [\nu_4\sigma'\eta_{14}]_6,\ [\eta_5\varepsilon_6]_7\eta_{14} =[\nu_4\sigma'\eta_{14}]_7.\label{m1:km1}
				\end{aligned}
			\end{equation}
			
			As  consequences of
			\eqref{m1:km1},
			we get
			$$
			[\eta_5\varepsilon_6]\nu^3_{14}  = [\eta_5\varepsilon_6]\eta_{14}\bar{\nu}_{15}
			= 4[ \nu_5]\sigma_8\bar{\nu}_{15} + [\nu_4\sigma'\eta_{14}]_6\bar{\nu}_{15}
			=[\nu_4\sigma'\eta_{14}]_6\bar{\nu}_{15}$$\noindent
			and
			$
			[\eta_5\varepsilon_6]\eta_{14}\varepsilon_{15} = [\nu_4\sigma'\eta_{14}]_6\varepsilon_{15}.
			$
			Therefore, in light of (1), we obtain
			$$[\eta_5\varepsilon_6]\nu^3_{14} = [\nu_4\sigma'\eta_{14}]_6\bar{\nu}_{15}
			= [\nu_4\eta_7\bar{\varepsilon}_8]_6,$$ 
			$$[\eta_5\varepsilon_6]\eta_{14}\varepsilon_{15}
			= [\nu_4\sigma'\eta_{14}]_6\varepsilon_{15} = [\nu_4(\Sigma\zeta')]_6.$$ \noindent
			According to \eqref{m1:epep} and the relations
			$[\eta_6]\nu_7 = b[  \nu^2_4]_7$ for some odd $b$ (\cite[Lemma 1.1(iv)]{KM1}),
			$2[\eta_6]\eta_7\bar{\varepsilon}_8=0$ and (1), we obtain
			$$[\eta_6]\eta_7\bar{\varepsilon}_8 = [\eta_6]\nu_7\sigma_{10}\nu^2_{17}
			= [  \nu^2_4]_7\sigma_{10}\nu^2_{17} = [\nu_4\eta_7\bar{\varepsilon}_8]_7.$$ \noindent
			In light of \eqref{m1:et6s'}, we see that
			$[\eta_6]\sigma'\eta_{14}
			= [\eta_5\varepsilon_6]_7\eta_{14}.$ 
			Then, from  $\Sigma\zeta'=\sigma'\eta_{14}\varepsilon_{15}$
			and (1), we infer
			$[\eta_6]\Sigma\zeta'
			= [\eta_6]\sigma'\eta_{14}\varepsilon_{15}
			= [\eta_5\varepsilon_6]_7\eta_{14}\varepsilon_{15}
			= [\nu_4(\Sigma\zeta')]_7.$  This leads to (2).
			
			\item Next, we show $[\eta_5\varepsilon_6]_7\mu_{14}\equiv[\nu_4\sigma'\mu_{14}]_7\ \bmod\ [\nu_4\eta_7\bar{\varepsilon}_8]_7, [\nu_4(\Sigma\zeta')]_7.$ 
			By (1), \eqref{m1:mR524}, \eqref{m1:km1} and the fact that $\pi_{16}(\mathit{SO}(7))\cong(\mathbb Z/2)^7,$   we have
			\[\begin{aligned}
				[\eta_5\varepsilon_6]_7\mu_{14}&\in\{[\nu_4\sigma'\eta_{14}]_7,2\iota_{15},8\sigma_{15}\}_1\ni[\nu_4\sigma'\mu_{14}]_7
			\end{aligned}\]
			$\mathrm{mod}\; [\nu_4\sigma'\eta_{14}]_7\circ \Sigma\pi^{14}_{22} +\pi_{16}(\mathit{SO}(7))\circ 8\sigma_{16} =\mathrm{span}\{[\nu_4(\Sigma\zeta')]_7, [\nu_4\eta_7\bar{\varepsilon}_8]_7\}.$  According to the relation
			${p_6}_*([\eta_5\varepsilon_6]\mu_{14})={p_6}_\ast(4[\zeta_5]\sigma_{16})$
			from \eqref{m1:4z5sig}, we deduce
			\[\begin{aligned}
				&[\eta_5\varepsilon_6]\mu_{14}\equiv 4[\zeta_5]\sigma_{16} \bmod {i_6}_*\pi_{23}(\mathit{SO}(5)) = \mathrm{span}\{[\nu_4\sigma'\mu_{14}]_6, [\nu_4(\Sigma\zeta')]_6, [\nu_4\eta_7\bar{\varepsilon}_8]_6\}.
			\end{aligned}\]
			This implies the first relation in (3).     As a consequence of \eqref{m1:et6s'} and \eqref{m1:km1}, we find
			\[\begin{aligned}
				&[\eta_6]\sigma'\mu_{14}=[\eta_5\varepsilon_6]_7\mu_{14} \equiv [\nu_4\sigma'\mu_{14}]_7\ \bmod\ [\nu_4\eta_7\bar{\varepsilon}_8]_7, [\nu_4(\Sigma\zeta')]_7.
			\end{aligned}\]     This leads to the second assertion in (3) and completes the proof.
		\end{enumerate}
	\end{proof}
	
	We record the following remark.
	\[\begin{aligned}
		&[\eta_6]\circ\pi^7_{23}={}\mathrm{span}\{[\eta_6]\mu_7\sigma_{16}, [\nu_4\sigma'\mu_{14}]_7, [\nu_4(\Sigma\zeta')]_7, [\nu_4\eta_7\bar{\varepsilon}_8]_7\},
	\end{aligned}\]
	where $[\eta_6](\Sigma\zeta')=[\nu_4(\Sigma\zeta')]_7,$   $[\eta_6]\eta_7\bar{\varepsilon}_8=[\nu_4\eta_7\bar{\varepsilon}_8]_7$ and
	$$\begin{aligned}
		&[\eta_6]\sigma'\mu_{14}\equiv[\nu_4\sigma'\mu_{14}]_7\ \bmod\ [\nu_4(\Sigma\zeta')]_7, [\nu_4\eta_7\bar{\varepsilon}_8]_7.
	\end{aligned}$$
	
	By \cite[I-Proposition 3.1(2)]{Oda1979}, we obtain
	\[\begin{aligned}
		&\varepsilon_5\kappa_{13}=\eta^2_5\bar{\kappa}_7\ \mbox{and}\ \eta_5\varepsilon_6\kappa_{14}=4\nu_5\bar{\kappa}_8.
	\end{aligned}\]
	Therefore, we conclude that
	\begin{equation}
		\begin{aligned}
			&\bar{\nu}_6\bar{\varepsilon}_{14}=\nu^3_6\kappa_{15}=\varepsilon_6\bar{\varepsilon}_{14}.\label{m1:bn6+epbep}
		\end{aligned}
	\end{equation}
	In light of the relation $J[ \nu_5]=\bar{\nu}_6+\varepsilon_6$ \cite[Lemma 1.1(ii)]{KM1},
	we see that
	\[\begin{aligned}
		&J([ \nu_5]\bar{\varepsilon}_8)=J[ \nu_5]\circ \bar{\varepsilon}_{14} =(\bar{\nu}_6+\varepsilon_6)\bar{\varepsilon}_{14}.
	\end{aligned}\]
	Therefore, from the relation \eqref{m1:bn6+epbep}, we obtain
	\begin{equation}
		\begin{aligned}
			&J([ \nu_5]\bar{\varepsilon}_8)=0.\label{m1:Jbar}
		\end{aligned}
	\end{equation}
	
	We examine the element $J[P(\Sigma\theta)+\nu_6\kappa_9]\in\pi^7_{30}.$ 
	Since $HJ[\nu_4\sigma'\mu_{14}]=\Sigma^5(\nu_4\sigma'\mu_{14})=0,$   we obtain
	$J[\nu_4\sigma'\mu_{14}]\equiv 0\ \bmod\ \Sigma\pi^4_{27}=\mathrm{span}\{\nu_5\bar{\kappa}_8\}.$  We see that $[\nu_4\sigma'\mu_{14}]_\infty=0\in \pi_{23}(\mathit{SO})\cong\mathbb Z_{(2)}$ \cite{Bo} and
	$\mathrm{ord}(\nu\bar{\kappa})=8$ \cite[Theorem 1.1(a)]{MMO}. This implies
	\[\begin{aligned}
		&J[\nu_4\sigma'\mu_{14}]=0.
	\end{aligned}\]
	
	In view of the similar argument, we deduce
	\[\begin{aligned}
		&J[\nu_4(\Sigma\zeta')]=0\ \mbox{and}\ J[\nu_4\eta_7\bar{\varepsilon}_8]=0.
	\end{aligned}\]
	Hence, according to \eqref{m1:Jbar} and the fact that
	\[\begin{aligned}
		&2[P(\Sigma\theta)+\nu_6\kappa_9] \equiv[ \nu_5]_7\bar{\varepsilon}_8\ \bmod\ [\nu_4\sigma'\mu_{14}]_7, [\nu_4(\Sigma\zeta')]_7, [\nu_4\eta_7\bar{\varepsilon}_8]_7,
	\end{aligned}\]
	we obtain
	\[\begin{aligned}
		&2J[P(\Sigma\theta)+\nu_6\kappa_9]=0.
	\end{aligned}\]
	
	The suspension identities below are given in \cite[(3.2), (4.9)]{MMO}; we also use the displayed group structure.
	
	$$\begin{aligned}
		&2\bar{\rho}'=\Sigma\bar{\rho}'',\ 2\bar{\rho}_{11}=\Sigma^4\bar{\rho}',\ \Sigma\pi^6_{29} = \mathrm{span}\{2\bar{\rho}',\nu_7\bar{\kappa}_{10}, \phi_7\}\cong\mathbb Z/4\oplus\mathbb Z/8\oplus\mathbb Z/2,
	\end{aligned}$$
	\vspace{-1.1em}
	$$  \begin{aligned}
		&\pi^S_{23}={}\mathrm{span}\{\bar{\rho}, \nu\bar{\kappa}, \phi\}\cong\mathbb Z/16\oplus\mathbb Z/8\oplus\mathbb Z/2.
	\end{aligned}$$
	
	\noindent
	As a consequence of \eqref{m1:Hbk7}, we have
	\[\begin{aligned}
		&HJ[P(\Sigma\theta)+\nu_6\kappa_9] = \Sigma^7(P(\Sigma\theta)+\nu_6\kappa_9) = \nu_{13}\kappa_{16} = \kappa_{13}\nu_{27} = H(\bar{\kappa}_7\nu_{27}).
	\end{aligned}
	\]\noindent Therefore,
	we obtain
	\[\begin{aligned}
		&J[P(\Sigma\theta)+\nu_6\kappa_9]\equiv\bar{\kappa}_7\nu_{27} \bmod \Sigma\pi^6_{29} = \mathrm{span}\{2\bar{\rho}', \nu_7\bar{\kappa}_{10}, \phi_7\}.
	\end{aligned}\]
	Note that $2J[P(\Sigma\theta)+\nu_6\kappa_9]=0,$  
	$[P(\Sigma\theta)+\nu_6\kappa_9]_\infty=0,$   $\Sigma^\infty(\bar{\kappa}_7\nu_{27}
	- \nu_7\bar{\kappa}_{10})=0$ and $\Sigma^\infty: \mathrm{span}\{2\bar{\rho}',
	\nu_7\bar{\kappa}_{10}, \phi_7\}\to\pi^S_{23}$ is a monomorphism \cite[Theorem 1.1(a)]{MMO}. Thus,
	we obtain
	\begin{equation}
		\begin{aligned}
			&J[P(\Sigma\theta)+\nu_6\kappa_9]=\bar{\kappa}_7\nu_{27}-\nu_7\bar{\kappa}_{10}.\label{m1:JP}
		\end{aligned}
	\end{equation}
	
	\subsection{The Case $n=8$}
	
	Once the case $n=7$ is known, the splitting for $n=8$ adds the contribution from $S^7.$  We keep the
	generators of the two contributions separate, since their images must be followed in the next
	inclusion.
	
	Since $\pi_{23}(\mathit{SO}(8))\cong \pi_{23}(\mathit{SO}(7))\oplus\pi^7_{23}$ and
	$\pi^7_{23}=\mathrm{span}\{\sigma'\mu_{14}, \Sigma\zeta', \eta_7\bar{\varepsilon}_8,
	\mu_7\sigma_{16}\}\cong(\mathbb Z/2)^4,$   we have the following.
	\begin{prop} \label{m1:R8}
		\[
		\FitMath{\begin{aligned}
				\pi_{23}(\mathit{SO}(8)) ={}\mathrm{span}\bigl\{& [\zeta_5]_8\sigma_{16}, [\nu_4\sigma'\mu_{14}]_8, [\nu_4(\Sigma\zeta')]_8, [\nu_4\eta_7\bar{\varepsilon}_8]_8, [\eta_6]_8\mu_7\sigma_{16}, \\
				&[P(\Sigma\theta)+\nu_6\kappa_9]_8, [\iota_7]\sigma'\mu_{14}, [\iota_7](\Sigma\zeta'), [\iota_7]\eta_7\bar{\varepsilon}_8, [\iota_7]\mu_7\sigma_{16} \bigr\} \cong(\mathbb Z/2)^5\oplus\mathbb Z/4\oplus(\mathbb Z/2)^4.
		\end{aligned}}
		\] In addition, $[\iota_7]_9\sigma'=\pm[\mathrm{c}_{\mathrm{KM}}(\bar{\nu}_6+\varepsilon_6)]_9.$
		\qed
	\end{prop}
	
	We consider $J: \pi_{14}(\mathit{SO}(9))\to\pi^9_{23}.$ 
	It follows that $\pi_{14}(\mathit{SO}(9))=\mathrm{span}\{[\mathrm{c}_{\mathrm{KM}}(\bar{\nu}_6+\varepsilon_6)]_9, [ \nu_5]_9\nu^2_8\}
	$ which is isomorphic to $\mathbb Z/8\oplus\mathbb Z/2$ \cite[Table 2]{Ka} and $\pi^9_{23}=\mathrm{span}\{\sigma^2_9,\kappa_9\}
	\cong\mathbb Z/16\oplus\mathbb Z/4.$  Therefore, $J$ is a monomorphism.
	By \eqref{m1:Ji7} and \eqref{m1:Jbn6e6}, we see that
	\[\begin{aligned}
		&J([\iota_7]_9\sigma')=\pm2\sigma^2_9,\qquad J([\mathrm{c}_{\mathrm{KM}}(\bar{\nu}_6+\varepsilon_6)]_9)=\pm2\sigma^2_9.
	\end{aligned}\]\noindent    
	Hence, we obtain
	\[\begin{aligned}
		&[\iota_7]_9\sigma'=\pm[\mathrm{c}_{\mathrm{KM}}(\bar{\nu}_6+\varepsilon_6)]_9.
	\end{aligned}\]\noindent    
	Choosing the sign so that the image in $\pi_{14}(\mathit{SO}(9))$ is zero gives
	\[\begin{aligned}
		&[\iota_7]_8\sigma'\pm[\mathrm{c}_{\mathrm{KM}}(\bar{\nu}_6+\varepsilon_6)]_8\in\Delta\pi^8_{15}.
	\end{aligned}\]\noindent    
	We know that $\pi^8_{15}=\mathrm{span}\{\sigma_8, \Sigma\sigma'\}\cong\mathbb
	Z_{(2)}\oplus\mathbb Z/8.$  Therefore, according to \cite[Lemma 1.2(ii)]{KM1}, we
	deduce, with the inner sign chosen as above,
	\[\begin{aligned}
		&\Delta\sigma_8\equiv\pm\bigl([\iota_7]\sigma'\pm[\mathrm{c}_{\mathrm{KM}}(\bar{\nu}_6+\varepsilon_6)]_8\bigr)\ \bmod\ \Delta(\Sigma\sigma').
	\end{aligned}\]\noindent    
	We find $\Delta(\Sigma\sigma')=2[\iota_7]\sigma'-[\eta_6]_8\sigma'.$    As a consequence of \eqref{m1:et6s'}, we obtain
	\begin{equation}
		\begin{aligned}
			&\Delta\sigma_8\equiv x_1[\iota_7]\sigma'+x_2[\mathrm{c}_{\mathrm{KM}}(\bar{\nu}_6+\varepsilon_6)]_8\ \bmod\ [ \nu_5]_8\nu^2_8+[\eta_5\varepsilon_6]_8.\label{m1:Dsg82}
		\end{aligned}
	\end{equation}
	Here $x_1,x_2\in\z$ are odd integers.

	\subsection{Boundary Relations for the Passage to the Case $n=9$}
	
	To pass to the case $n=9,$   we need relations among the composites already represented for smaller values of $n.$  The
	next lemma identifies these relations with their remaining congruence terms.
	
	\begin{lem} \label{m1:nusiep}
		\begin{enumerate}
			\item $[ \nu_5]\sigma_8\bar{\nu}_{15}
			\equiv [ \nu_5]\sigma_8\varepsilon_{15}\equiv [ \nu_5]\bar{\varepsilon}_8
			\bmod [\nu_4\sigma'\mu_{14}]_6, [\nu_4(\Sigma\zeta')]_6,
			[\nu_4\eta_7\bar{\varepsilon}_8]_6.$ 
			\item $[\mathrm{c}_{\mathrm{KM}}(\bar{\nu}_6+\varepsilon_6)]\nu^3_{14}
			\equiv [\mathrm{c}_{\mathrm{KM}}(\bar{\nu}_6+\varepsilon_6)]\eta_{14}\varepsilon_{15}
			\equiv
			[ \nu_5]_7\bar{\varepsilon}_8\
			\bmod \
			\ [\nu_4\sigma'\mu_{14}]_7, [\nu_4(\Sigma\zeta')]_7,
			[\nu_4\eta_7\bar{\varepsilon}_8]_7.$ 
		\end{enumerate}
	\end{lem}
	\begin{proof}
		\begin{enumerate}[label=\textup{(\arabic*)},leftmargin=2em]
			\item By Lemma~\ref{m1:zesig}(3), $$[ \nu_5]\sigma_8\bar{\nu}_{15}
			- [ \nu_5]\bar{\varepsilon}_8\in {i_6}_*\pi_{23}(\mathit{SO}(5))$$ \noindent and
			$[ \nu_5]\sigma_8\varepsilon_{15}- [ \nu_5]\bar{\varepsilon}_8\in {i_6}_*\pi_{23}(\mathit{SO}(5)).$ 
			This leads to (1).
			
			\item According to \cite[Lemma 3.6(i)-(ii)]{HKM}, we obtain
			\begin{gather*}
				[\mathrm{c}_{\mathrm{KM}}(\bar{\nu}_6+\varepsilon_6)]\sigma_{14}\in\mathrm{span}\{[\nu_4\sigma'\sigma_{14}]_7, 2[\eta_6]\kappa_7\}
				\shortintertext{and}
				[\nu_4\sigma'\sigma_{14}]_6\eta_{21}=[\nu_4\sigma'\bar{\nu}_{14}]_6+
				[\nu_4\sigma'\varepsilon_{14}]_6.
			\end{gather*}
			Therefore, we see that
			\begin{align*}
				&[\mathrm{c}_{\mathrm{KM}}(\bar{\nu}_6+\varepsilon_6)]\eta^2_{14}\sigma_{16} \in\mathrm{span}\{[\nu_4\eta_7\bar{\varepsilon}_8]_7+[\nu_4(\Sigma\zeta')]_7\}, \\
				&[\mathrm{c}_{\mathrm{KM}}(\bar{\nu}_6+\varepsilon_6)]\nu^3_{14}\equiv[\mathrm{c}_{\mathrm{KM}}(\bar{\nu}_6+\varepsilon_6)]\eta_{14}\varepsilon_{15}\ \bmod\ [\nu_4\eta_7\bar{\varepsilon}_8]_7+[\nu_4(\Sigma\zeta')]_7,
			\end{align*} \noindent
			and
			$[\mathrm{c}_{\mathrm{KM}}(\bar{\nu}_6+\varepsilon_6)]\eta_{14}\varepsilon_{15}
			\equiv [\eta_6]\eta_7\bar{\varepsilon}_8 + [ \nu_5]_7\sigma_8\varepsilon_{15}
			\bmod [\nu_4\sigma'\eta_{14}]_7\varepsilon_{15} = [\nu_4(\Sigma\zeta')]_7.$ 
			Hence, in light of (1) and Lemma~\ref{lem6d4}(2), we have (2). This completes the
			proof.
		\end{enumerate}
	\end{proof}
	
	We continue by stating the following result, taking up its proof next before returning to the remainder of the discussion.
	\begin{lem} \label{m1:bnuepmu}
		$[\mathrm{c}_{\mathrm{KM}}(\bar{\nu}_6+\varepsilon_6)]\mu_{14}\equiv
		[\eta_6]\mu_7\sigma_{16} + [\zeta_5]_7\sigma_{16}
		\bmod \ [ \nu_5]_7\bar{\varepsilon}_8, [\nu_4\sigma'\mu_{14}]_7,
		[\nu_4(\Sigma\zeta')]_7, [\nu_4\eta_7\bar{\varepsilon}_8]_7.$ 
	\end{lem}
	\begin{proof}
		Note that $(\bar{\nu}_6+\varepsilon_6)\mu_{14}=\eta_6\mu_7\sigma_{16}.$  Therefore, we find
		\[\begin{aligned}
			&[\mathrm{c}_{\mathrm{KM}}(\bar{\nu}_6+\varepsilon_6)]\mu_{14}\equiv[\eta_6]\mu_7\sigma_{16}\ \bmod\ {i_7}_*\pi_{23}(\mathit{SO}(6)).
		\end{aligned}\]
		Hence, by \eqref{m1:R623}, we obtain
		\[\begin{aligned}
			&[\mathrm{c}_{\mathrm{KM}}(\bar{\nu}_6+\varepsilon_6)]\mu_{14}\equiv[\eta_6]\mu_7\sigma_{16}\ \bmod\ [\zeta_5]_7\sigma_{16}, [ \nu_5]_7\bar{\varepsilon}_8, [\nu_4\sigma'\mu_{14}]_7, [\nu_4\eta_7\bar{\varepsilon}_8]_7, [\nu_4(\Sigma\zeta')]_7.
		\end{aligned}\]\noindent To control the bracket comparison, we use the subgroup generated by the four elements in the statement.
		The calculation of $\pi_{16}(\mathit{SO}(7))$ preceding \eqref{m1:zesig1} gives
		$\pi_{16}(\mathit{SO}(7))\cong(\mathbb Z/2)^7.$ Hence
		\[
		\pi_{16}(\mathit{SO}(7))\circ2\sigma_{16}=0,
		\qquad 2[\zeta_5]_7\sigma_{16}=0.
		\]
		For each of the three brackets below, the full indeterminacy in $\pi_{23}(\mathit{SO}(7))$ is
		\[
		\mathrm{Ind}\{\alpha,8\iota_{15},2\sigma_{15}\}
		=\alpha\circ\pi^{15}_{23}
		+\pi_{16}(\mathit{SO}(7))\circ2\sigma_{16},
		\]
		where $\alpha$ denotes its first entry. Since
		$\pi^{15}_{23}=\mathrm{span}\{\bar\nu_{15},\varepsilon_{15}\},$
		Lemmas~\ref{m1:zesig}(2), \ref{lem6d4}(1),(2) and \ref{m1:nusiep}(1) give
		\begin{align*}
			[\eta_6](\bar\nu_7+\varepsilon_7)\circ\pi^{15}_{23}
			&=\mathrm{span}\{[\nu_4\eta_7\bar\varepsilon_8]_7\},\\
			[ \nu_5]_7\sigma_8\circ\pi^{15}_{23}
			&\subseteq\mathrm{span}\{[ \nu_5]_7\bar\varepsilon_8,
			[\nu_4\sigma'\mu_{14}]_7,\\
			&\hspace{4em}[\nu_4(\Sigma\zeta')]_7,
			[\nu_4\eta_7\bar\varepsilon_8]_7\},\\
			[\nu_4\sigma'\eta_{14}]_7\circ\pi^{15}_{23}
			&=\mathrm{span}\{[\nu_4\eta_7\bar\varepsilon_8]_7,
			[\nu_4(\Sigma\zeta')]_7\}.
		\end{align*}
		Thus all three indeterminacies lie in the subgroup generated by the four stated elements.
		This also holds when the first entry of the third bracket is multiplied by $x\in\{0,1\}.$
		
		By \eqref{m1:mu4} and Proposition~\ref{prop:toda-composition}(3),
		$\mu_{14}$ can be replaced by a value of
		$\{\eta_{14},8\iota_{15},2\sigma_{15}\}$ modulo
		$\mathrm{span}\{\nu_{14}^3,\eta_{14}\varepsilon_{15}\}.$
		Lemma~\ref{m1:nusiep}(2) shows that both correction terms, after left composition with
		$[\mathrm{c}_{\mathrm{KM}}(\bar\nu_6+\varepsilon_6)],$ lie in the same subgroup of $\pi_{23}(\mathit{SO}(7)).$
		We now compare the images of the elements and brackets in the quotient by this subgroup.
		Using \eqref{m1:tr2}, \eqref{m1:zesig1}, \eqref{m1:mR524},
		Lemma~\ref{m1:m3sg}(3) and Proposition~\ref{prop:toda-composition}(3),(5), we obtain
		\begin{align*}
			[\mathrm{c}_{\mathrm{KM}}(\bar{\nu}_6+\varepsilon_6)]\mu_{14}
			&\in\{[\mathrm{c}_{\mathrm{KM}}(\bar{\nu}_6+\varepsilon_6)]\eta_{14},8\iota_{15},2\sigma_{15}\}\\
			&=\{[\eta_6](\bar\nu_7+\varepsilon_7)+[ \nu_5]_7\sigma_8
			+x[\nu_4\sigma'\eta_{14}]_7,8\iota_{15},2\sigma_{15}\}\\
			&\subseteq\{[\eta_6](\bar\nu_7+\varepsilon_7),8\iota_{15},2\sigma_{15}\}
			+\{[ \nu_5]_7\sigma_8,8\iota_{15},2\sigma_{15}\}\\
			&\qquad+\{x[\nu_4\sigma'\eta_{14}]_7,8\iota_{15},2\sigma_{15}\}\\
			&\ni[\eta_6](\Sigma\zeta'+\mu_7\sigma_{16})
			+[\zeta_5]_7\sigma_{16}+x[\nu_4\sigma'\mu_{14}]_7,
			\qquad x\in\{0,1\}.
		\end{align*}
		The term $2[\zeta_5]_7\sigma_{16}$ arising from \eqref{m1:zesig1} is zero, as shown above.
		Each bracket in the sum has zero indeterminacy in this quotient, so the sum has only one value.
		Consequently, the two displayed elements are congruent modulo the four stated elements.
		Finally, Lemma~\ref{lem6d4}(2) gives
		$[\eta_6]\Sigma\zeta'=[\nu_4(\Sigma\zeta')]_7,$ which belongs to the modulus,
		as does $x[\nu_4\sigma'\mu_{14}]_7.$ This proves the asserted congruence.
	\end{proof}
	
	From \cite[p. 77]{T}, we see that
	\begin{equation}
		\begin{aligned}
			&\nu^4_5=P(\bar{\nu}_{11}).\label{m1:nu45}
		\end{aligned}
	\end{equation}
	
	\begin{lem} \label{m1:89}
		$[\nu_4\sigma'\mu_{14}]\eta_{23}= [\nu_4\sigma'\eta_{14}]\mu_{15},$  
		$[\nu_4(\Sigma\zeta')]\eta_{23}= [\nu_4\eta_7\bar{\varepsilon}_8]\eta_{23}=0$\
		and\ $[ \nu_5]\eta_8\bar{\varepsilon}_9=0.$  \
		
	\end{lem}
	\begin{proof}
		
		Since $\nu_4\sigma'\mu_{14}\eta_{23} = \nu_4\sigma'\eta_{14}\mu_{15},$  
		we obtain $[\nu_4\sigma'\mu_{14}]\eta_{23} - [\nu_4\sigma'\eta_{14}]\mu_{15}
		\in {i_5}_*\pi_{24}(\mathit{SO}(4)) = 0$ \eqref{m1:i5R424}. This leads to the first relation.
		
		As $\zeta'\eta_{22} = 0,$  
		$[\nu_4(\Sigma\zeta')]\eta_{23}\in {i_5}_*\pi_{24}(\mathit{SO}(4)) = 0.$  Since
		$\eta_7\bar{\varepsilon}_8 = \nu_7\sigma_{10}\nu^2_{17}$ from \eqref{m1:epep}, we deduce
		$\nu_4\eta_7\bar{\varepsilon}_8\eta_{23} = 0.$  This implies the relation
		$[\nu_4\eta_7\bar{\varepsilon}_8]\eta_{23} = 0.$  Also, $\nu_5\eta_8\bar{\varepsilon}_9 = 0.$ 
		Since $\eta_8\bar{\varepsilon}_9=(\Sigma\sigma')\nu^3_{15}$ and $2[ \nu_5]\sigma_8
		=[ \nu_5](\Sigma\sigma'),$   we find $[ \nu_5]\eta_8\bar{\varepsilon}_9
		=[ \nu_5](\Sigma\sigma')\nu^3_{15}=0.$ 
		This proves the third relation and completes the proof.
	\end{proof}
	
	\subsection{The Case $n=9$ and Its Inherited Generators}
	
	We now apply the boundary relations to the passage from the case $n=8$ to the case $n=9.$  The resulting
	generators are written using their earlier representatives, so the effect of inclusion remains
	visible.
	
	\begin{prop} \label{m1:R923}
		The structure of $\pi_{23}(\mathit{SO}(9))$ is as follows.
		$$\GroupLayout{\pi_{23}(\mathit{SO}(9)) =\mathrm{span}\{[P(\Sigma\theta)+\nu_6\kappa_9]_9,\; [\iota_7]_9\mu_7\sigma_{16},\; [\zeta_5]_9\sigma_{16}\}\cong \mathbb Z/4\oplus (\mathbb Z/2)^{2}.}{\begin{aligned}\pi_{23}(\mathit{SO}(9)) =\mathrm{span}\{&[P(\Sigma\theta)+\nu_6\kappa_9]_9,\; [\iota_7]_9\mu_7\sigma_{16},\\
				&[\zeta_5]_9\sigma_{16}\}\cong \mathbb Z/4\oplus (\mathbb Z/2)^{2}.\end{aligned}}$$
		
	\end{prop}
	\begin{proof}
		In the exact sequence $(23)_8$:
		\[\begin{aligned}
			&\pi^8_{24}\xrightarrow{\Delta}\pi_{23}(\mathit{SO}(8))\xrightarrow{i_*}\pi_{23}(\mathit{SO}(9))\xrightarrow{p_*}\pi^8_{23} \xrightarrow{\Delta}\pi_{22}(\mathit{SO}(8)),
		\end{aligned}\]
		$\Delta\colon  \pi^8_{23}\to \pi_{22}(\mathit{SO}(8))$ is a monomorphism and
		\[\begin{aligned}
			&\pi^8_{24}={}\mathrm{span}\{\sigma_8\nu^3_{15}, \sigma_8\eta_{15}\varepsilon_{16}, \sigma_8\mu_{15}, (\Sigma\sigma')\mu_{15}, \Sigma^2\zeta', \mu_8\sigma_{17}, \eta_8\bar{\varepsilon}_9\}\cong(\mathbb Z/2)^7.
		\end{aligned}\]
		According to Lemma~\ref{lem6d4}(3) and the relation $\Delta\iota_8 = 2[\iota_7]
		- [\eta_6]_8,$   we obtain the following:
		\begin{align*}
			&\Delta((\Sigma\sigma')\mu_{15}) = [\eta_6]_8\sigma'\mu_{14} \equiv 0\ \bmod\ [\nu_4\sigma'\mu_{14}]_8, [\nu_4(\Sigma\zeta')]_8, [\nu_4\eta_7\bar{\varepsilon}_8]_8\uprightsemicolon{} \\
			&\Delta(\Sigma^2\zeta') = [\eta_6]_8\Sigma\zeta' = [\nu_4(\Sigma\zeta')]_8\uprightsemicolon{}\quad \Delta(\mu_8\sigma_{17}) = [\eta_6]_8\mu_7\sigma_{16}\uprightsemicolon{} \\
			&\Delta(\eta_8\bar{\varepsilon}_9) = [\eta_6]_8\eta_7\bar{\varepsilon}_8 = [\nu_4\eta_7\bar{\varepsilon}_8]_8.
		\end{align*}
		In the following applications of \eqref{m1:Dsg82}, composition with
		$\eta_{14}\varepsilon_{15},$   $\nu_{14}^3,$   or $\mu_{14}$ gives classes of order at most $2.$ 
		Thus the odd coefficients $x_1,x_2$ act as $1,$   and signs do not affect these composites.
		In light of \eqref{m1:Dsg82}, we conclude that
		\[\begin{aligned}
			&\Delta(\sigma_8\eta_{15}\varepsilon_{16}) \equiv [\iota_7]\Sigma\zeta' + [\mathrm{c}_{\mathrm{KM}}(\bar{\nu}_6+\varepsilon_6)]_8\eta_{14}\varepsilon_{15} \bmod [\eta_5\varepsilon_6]_8\eta_{14}\varepsilon_{15}.
		\end{aligned}\]
		Therefore, by Lemmas \ref{lem6d4}(1) and \ref{m1:nusiep}(2),
		\[\begin{aligned}
			&\Delta(\sigma_8\eta_{15}\varepsilon_{16}) \equiv [\iota_7]\Sigma\zeta' + [\nu_4\eta_7\bar{\varepsilon}_8]_8 + [ \nu_5]_8\bar{\varepsilon}_8 \bmod [\nu_4\sigma'\mu_{14}]_8, [\nu_4(\Sigma\zeta')]_8, [\nu_4\eta_7\bar{\varepsilon}_8]_8.
		\end{aligned}\]
		Since $\sigma'\nu^3_{14} = \eta_7\bar{\varepsilon}_8$ \eqref{m1:epep} and
		$\nu^5_8 = 0$ \eqref{m1:nu45},
		we have
		\[\begin{aligned}
			&\Delta(\sigma_8\nu^3_{15})\equiv [\iota_7]\eta_7\bar{\varepsilon}_8 + [\mathrm{c}_{\mathrm{KM}}(\bar{\nu}_6+\varepsilon_6)]_8\nu^3_{14} \bmod [\eta_5\varepsilon_6]_8\nu^3_{14}.
		\end{aligned}\]
		In view of Lemma~\ref{m1:nusiep}, we know
		\[\begin{aligned}
			&[\mathrm{c}_{\mathrm{KM}}(\bar{\nu}_6+\varepsilon_6)]_8\nu^3_{14} = [\mathrm{c}_{\mathrm{KM}}(\bar{\nu}_6+\varepsilon_6)]_8\eta_{14}\bar{\nu}_{15} \equiv [\mathrm{c}_{\mathrm{KM}}(\bar{\nu}_6+\varepsilon_6)]_8\eta_{14}\varepsilon_{15} \bmod [\nu_4\sigma'\mu_{14}]_8, [\nu_4(\Sigma\zeta')]_8, [\nu_4\eta_7\bar{\varepsilon}_8]_8.
		\end{aligned}\]
		According to Lemma~\ref{lem6d4}(2), we know $[\eta_5\varepsilon_6]_8\nu^3_{14}
		= [\nu_4\eta_7\bar{\varepsilon}_8]_8.$  Hence we obtain
		\[\begin{aligned}
			&\Delta(\sigma_8\nu^3_{15})\equiv [\iota_7]\eta_7\bar{\varepsilon}_8 + [\nu_4\eta_7\bar{\varepsilon}_8]_8 + [ \nu_5]_8\bar{\varepsilon}_8 \bmod [\nu_4\sigma'\mu_{14}]_8, [\nu_4(\Sigma\zeta')]_8, [\nu_4\eta_7\bar{\varepsilon}_8]_8.
		\end{aligned}\]
		
		\noindent   Next, according to \eqref{m1:Dsg82},  we obtain
		\[\begin{aligned}
			&\Delta(\sigma_8\mu_{15})\equiv [\iota_7]\sigma'\mu_{14} + [\mathrm{c}_{\mathrm{KM}}(\bar{\nu}_6+\varepsilon_6)]_8\mu_{14} \bmod [\eta_5\varepsilon_6]_8\mu_{14}.
		\end{aligned}\]
		As a consequence of Lemma~\ref{m1:bnuepmu}, we get
		\[\begin{aligned}
			&[\mathrm{c}_{\mathrm{KM}}(\bar{\nu}_6+\varepsilon_6)]_8\mu_{14} \equiv [\eta_6]_8(\Sigma\zeta'+\mu_7\sigma_{16}) + [\zeta_5]_8\sigma_{16} \ \bmod [ \nu_5]_8\bar{\varepsilon}_8, [\nu_4\sigma'\mu_{14}]_8, [\nu_4(\Sigma\zeta')]_8, [\nu_4\eta_7\bar{\varepsilon}_8]_8.
		\end{aligned}\]
		As a consequence of           Lemma~\ref{lem6d4}(3),
		\[\begin{aligned}
			&[\eta_5\varepsilon_6]_8\mu_{14}\equiv [\nu_4\sigma'\mu_{14}]_8 \ \bmod \ [\nu_4(\Sigma\zeta')]_8, [\nu_4\eta_7\bar{\varepsilon}_8]_8.
		\end{aligned}\]
		Hence we obtain
		\[\begin{aligned}
			&\Delta(\sigma_8\mu_{15})\equiv [\iota_7]\sigma'\mu_{14} + [\zeta_5]_8\sigma_{16} + [\eta_6]_8\mu_7\sigma_{16} \bmod [ \nu_5]_8\bar{\varepsilon}_8, [\nu_4\sigma'\mu_{14}]_8, [\nu_4(\Sigma\zeta')]_8, [\nu_4\eta_7\bar{\varepsilon}_8]_8.
		\end{aligned}\]
		Thus we obtain the group $\pi_{23}(\mathit{SO}(9)),$   completing the proof.
	\end{proof}
	
	In the exact sequence $(23)_9$:
	\[\begin{aligned}
		&\pi^9_{24}\xrightarrow{\Delta}\pi_{23}(\mathit{SO}(9))\xrightarrow{i_*}\pi_{23}(\mathit{SO}(10)) \xrightarrow{p_*}\pi^9_{23}\xrightarrow{\Delta}\pi_{22}(\mathit{SO}(9)),
	\end{aligned}\]
	we record
	\[\begin{aligned}
		&\pi^9_{24}={}\frac{\mathbb Z}{16}\{\rho'\}\oplus\frac{\mathbb Z}{2}\{\sigma_9\bar{\nu}_{16}\} \oplus\frac{\mathbb Z}{2}\{\sigma_9\varepsilon_{16}\}\oplus\frac{\mathbb Z}{2}\{\bar{\varepsilon}_9\}.
	\end{aligned}\]
	
	Observe from \cite[(4.1)]{HKM} that
	\begin{equation}
		\begin{aligned}
			&\Delta\sigma_9 = [\iota_7]_9(\bar{\nu}_7+\varepsilon_7).\label{m1:Desig9}
		\end{aligned}
	\end{equation}
	Then, as a consequence of
	Lemma~\ref{m1:zesig} (2), we deduce
	\begin{equation}\notag
		\Delta(\sigma_9\bar{\nu}_{16})=\Delta(\sigma_9)\bar{\nu}_{15} = [\iota_7]_9(\bar{\nu}_7+\varepsilon_7)\bar{\nu}_{15} = [\iota_7]_9\eta_7\bar{\varepsilon}_8 =[ \nu_5]_9\bar{\varepsilon}_8,
	\end{equation}
	\begin{equation}\label{djx1024}
		\Delta(\sigma_9\varepsilon_{16})=\Delta(\sigma_9)\varepsilon_{15} = [\iota_7]_9(\bar{\nu}_7+\varepsilon_7)\varepsilon_{15} = [\iota_7]_9\eta_7\bar{\varepsilon}_8 =[ \nu_5]_9\bar{\varepsilon}_8.
	\end{equation}
	\noindent   That is,
	\begin{equation}
		\begin{aligned}
			&\Delta(\sigma_9\bar{\nu}_{16})=\Delta(\sigma_9\varepsilon_{16})=[ \nu_5]_9\bar{\varepsilon}_8.\label{m1:Desig9e}
		\end{aligned}
	\end{equation}
	
	\noindent Recall  from Proposition \ref{m1:R723} that $$2[P(\Sigma\theta)+\nu_6\kappa_9]
	\equiv[ \nu_5]_7\bar{\varepsilon}_8\ \bmod\
	[\nu_4\sigma'\mu_{14}]_7, [\nu_4(\Sigma\zeta')]_7, [\nu_4\eta_7\bar{\varepsilon}_8]_7.$$
	Combining this and Lemma \ref{lem6d4}, we deduce that  \begin{equation}\label{1024zmf}
		\Delta(\sigma_9\bar{\nu}_{16})=\Delta(\sigma_9\varepsilon_{16})=[ \nu_5]_9\bar{\varepsilon}_8=2[P(\Sigma\theta)+\nu_6\kappa_9]_9.
	\end{equation}

	Summarizing the above results, we obtain the following.
	
	\begin{lem} \label{m1:rel2}
		\begin{enumerate}
			\item
			\(\begin{aligned}[t]
				[\eta_5\varepsilon_6]_9\mu_{14}
				&=[\eta_5\varepsilon_6]_9\nu^3_{14}
				=[\eta_5\varepsilon_6]_9\eta_{14}\varepsilon_{15}=[\eta_6]_9\mu_7\sigma_{16}\\
				&
				=[\nu_4\sigma'\mu_{14}]_9
				=[\nu_4\eta_7\bar{\varepsilon}_8]_9
				=[\nu_4(\Sigma\zeta')]_9=0.
			\end{aligned}\)
			
			\item
			\(\begin{aligned}[t]
				[\mathrm{c}_{\mathrm{KM}}(\bar{\nu}_6+\varepsilon_6)]_9\nu^3_{14}
				&=[\mathrm{c}_{\mathrm{KM}}(\bar{\nu}_6+\varepsilon_6)]_9
				\eta_{14}\varepsilon_{15}\\
				&=[\nu_5]_9\bar{\varepsilon}_8
				=[\nu_5]_9\sigma_8\bar{\nu}_{15}=[\nu_5]_9\sigma_8\varepsilon_{15}
				=[\iota_7]_9\eta_7\bar{\varepsilon}_8=[\iota_7]_9\Sigma\zeta'.
			\end{aligned}\)
			
			\item
			\(\begin{aligned}[t]
				[\iota_7]_9\sigma'\mu_{14}
				&=[\mathrm{c}_{\mathrm{KM}}(\bar{\nu}_6+\varepsilon_6)]_9
				\mu_{14}\equiv[\zeta_5]_9\sigma_{16}
				\bmod [\nu_5]_9\bar{\varepsilon}_8.
			\end{aligned}\)
			
			\item
			\(\begin{aligned}[t]
				[\mathrm{c}_{\mathrm{KM}}(\bar{\nu}_6+\varepsilon_6)]_{10}\nu^3_{14}
				&=[\mathrm{c}_{\mathrm{KM}}(\bar{\nu}_6+\varepsilon_6)]_{10}
				\eta_{14}\varepsilon_{15}=[\nu_5]_{10}\bar{\varepsilon}_8
				=[\nu_5]_{10}\sigma_8\bar{\nu}_{15}\\
				&=[\nu_5]_{10}\sigma_8\varepsilon_{15}
				=[\iota_7]_{10}\eta_7\bar{\varepsilon}_8=[\iota_7]_{10}\Sigma\zeta'=0.
			\end{aligned}\)
			
			\item
			\(\begin{aligned}[t]
				[\iota_7]_{10}\sigma'\mu_{14}
				&=[\mathrm{c}_{\mathrm{KM}}(\bar{\nu}_6+\varepsilon_6)]_{10}
				\mu_{14}=[\zeta_5]_{10}\sigma_{16}.
			\end{aligned}\)
		\end{enumerate}
		\qed
	\end{lem}
	\section{The Cases $10\leq n\leq14$}
	\subsection{The Case $n=10$}
	For $n=10,$   the first task is to identify the surviving class and constrain the boundary of $\rho'.$ 
	These two pieces of information enter the extension calculation that follows.
	
	\begin{lem} \label{m1:Derh}
		\begin{enumerate}
			\item $[\iota_7]_{10}\mu_7\sigma_{16}\neq 0.$ 
			\item $\Delta\rho'\equiv 0 \bmod [\zeta_5]_9\sigma_{16},
			[ \nu_5]_9\bar{\varepsilon}_8.$ 
		\end{enumerate}
	\end{lem}
	\begin{proof}
		\begin{enumerate}[label=\textup{(\arabic*)},leftmargin=2em]
			\item We examine the image of a map $J\Delta\colon \pi_{24}^9\to\pi_{32}^9,$   where
			$\pi_{24}^9=\mathrm{span}\{\rho',\sigma_9\bar{\nu}_{16},\sigma_9\varepsilon_{16},\bar{\varepsilon}_9\}.$ 
			By the formula for $J\Delta$ in Section~\ref{m2:subsec:j-toda} and the facts that
			$2\rho_{13} = \Sigma^4\rho'$ \cite[Lemma 10.9]{T} and
			$2[\iota_9, \iota_9]=0,$   we find
			$J\Delta(\rho')=\pm[\iota_9, \iota_9]\circ\Sigma^8\rho'=\pm[\iota_9, \iota_9]\circ2\rho_{17}=0.$ 
			In light of the relations \eqref{m1:Desig9e} and \eqref{m1:Jbar}, we conclude that
			$J\Delta(\sigma_9\bar{\nu}_{16})=J\Delta(\sigma_9\varepsilon_{16})=J([ \nu_5]_9\bar{\varepsilon}_8)
			=(-1)^3\Sigma^3(J[ \nu_5]\bar{\varepsilon}_{14})=0.$ 
			Using the standard relations, we see that
			\begin{equation}
				\begin{aligned}
					&\Delta\bar{\varepsilon}_9 = [ \nu_5]_9\bar{\varepsilon}_8 + [\iota_7]_9\eta_7\bar{\varepsilon}_8 = 2[ \nu_5]_9\bar{\varepsilon}_8 = 0.\label{m1:Dbep9}
				\end{aligned}
			\end{equation}
			Hence we see that $J\Delta\pi^9_{24} = 0.$ 
			By \eqref{m1:Ji7}
			and \cite[Theorem 1.1(a)]{MMO}, we obtain
			\begin{equation}
				\begin{aligned}
					&J([\iota_7]_9\mu_7\sigma_{16}) = \sigma_9\mu_{16}\sigma_{25} = \sigma^2_9\mu_{23}\neq 0.\label{m1:s_9smu}
				\end{aligned}
			\end{equation}
			This implies $[\iota_7]_9\mu_7\sigma_{16}\not\in\Delta\pi^9_{24},$   and so
			(1) follows.
			
			\item By the relations $\rho'\in\{\sigma_9, 16\iota_{16}, \sigma_{16}\}_1$ \cite[p. 103]{T} and \eqref{m1:Desig9}, we obtain
			\begin{align*}
				\Delta\rho'&\in \Delta\{\sigma_9, 16\iota_{16}, \sigma_{16}\}_1 \subseteq\{\Delta\sigma_9, 16\iota_{15}, \sigma_{15}\} \\
				&=\{[\iota_7]_9(\bar{\nu}_7+\varepsilon_7), 16\iota_{15}, \sigma_{15}\} \\
				&\subseteq\{[\iota_7]_9,16(\bar{\nu}_7+\varepsilon_7),\sigma_{15}\} =\{[\iota_7]_9,0,\sigma_{15}\} \\
				&=[\iota_7]_9\circ\pi^7_{23}+ \pi_{16}(\mathit{SO}(9))\circ\sigma_{16}.
			\end{align*}
			Since $\pi^7_{23}=\mathrm{span}\{\sigma'\mu_{14}, \Sigma\zeta', \mu_7\sigma_{16},
			\eta_7\bar{\varepsilon}_8\},$   the standard composition relations imply
			\[\begin{aligned}
				&[\iota_7]_9\circ\pi^7_{23}={}\mathrm{span}\{[\zeta_5]_9\sigma_{16}, [\iota_7]_9\mu_7\sigma_{16}, [ \nu_5]_9\bar{\varepsilon}_8\}.
			\end{aligned}\]
			In light of the fact that
			$\pi_{16}(\mathit{SO}(9)) = \mathrm{span}\{[\zeta_5]_9,
			[ \nu_5]\bar{\nu}_8, [ \nu_5]\varepsilon_8, [\iota_7]_9\nu^3_7,
			[\iota_7]_9\eta_7\varepsilon_8, [\iota_7]_9\mu_7\}\cong(\mathbb Z/2)^6$
			\cite[Proposition 4.1]{Ka} and
			the standard composition relations,
			we have $$\pi_{16}(\mathit{SO}(9))\circ\sigma_{16} = \mathrm{span}\{[\zeta_5]_9\sigma_{16},
			[\iota_7]_9\mu_7\sigma_{16}\},$$\noindent  and hence we find
			\[\begin{aligned}
				&\Delta\rho'\equiv 0 \bmod [\zeta_5]_9\sigma_{16}, [ \nu_5]_9\bar{\varepsilon}_8, [\iota_7]_{9}\mu_7\sigma_{16}.
			\end{aligned}\]
			Applying the $J$-homomorphism to this relation,
			by the relations $J\Delta(\rho')=0,$   $J([ \nu_5]_9\bar{\varepsilon}_8)=0,$  
			$J([\iota_7]_9\mu_7\sigma_{16}) \neq 0$ \eqref{m1:s_9smu} and
			$J([\zeta_5]_9\sigma_{16})=J[\zeta_5]_9 \circ\sigma_{25}
			=\mu_9\sigma^2_{18}=0$ from \cite[Lemma 2.1(v)]{HKM} and \cite[(2.9)]{MMO},
			we have (2). This completes the proof.
		\end{enumerate}
	\end{proof}
	
	Before proceeding further, we record the
	relations that will be useful in what follows.
	\begin{lem}\label{m1:i710ms}
		\begin{enumerate}
			\item $[\zeta_5]_{10}\equiv [\iota_7]_{10}\mu_7 + 8\Delta\sigma_{10}
			\bmod [\iota_7]_{10}\nu^3_7.$ 
			\item
			$[\iota_7]_{10}\mu_7\sigma_{16}=8\Delta(\sigma^2_{10})\neq 0$\ and\
			$[\zeta_5]_{10}\sigma_{16}=0.$
		\end{enumerate}
	\end{lem}
	\begin{proof}
		\begin{enumerate}[label=\textup{(\arabic*)},leftmargin=2em]
			\item This is just \cite[(4.7), p.~28]{Ka}.
			
			\item As $[\iota_7]_{10}\eta_7 = [ \nu_5]_{10}$ and $[\iota_7]_{10}\Sigma\zeta'
			= 0,$   the boundary formula \eqref{m1:Desig9} together with Lemma~\ref{m1:m3sg}(3) gives
			\[\begin{aligned}
				&[\iota_7]_{10}\mu_7\sigma_{16} = [\iota_7]_{10}(\Sigma\zeta'+\mu_7\sigma_{16}) \in[\iota_7]_{10}\circ\{\bar{\nu}_7+\varepsilon_7, 2\iota_{15}, 8\sigma_{15}\} = \{[\iota_7]_{10}, \bar{\nu}_7+\varepsilon_7, 2\iota_{15}\} \circ 8\sigma_{16}.
			\end{aligned}\]
			
			As a consequence of           \eqref{m1:Desig9}, we see that $$\{[\iota_7]_{10}, \bar{\nu}_7+\varepsilon_7,
			2\iota_{15}\} \subseteq \{i_{10 \ast}, \Delta\sigma_9, 2\iota_{15}\}
			\supseteq \{i_{10 \ast}, \Delta(\iota_9), 2\sigma_8\} \ni
			\Delta(\sigma_{10}) $$ \noindent
			$\bmod \; i_{10 \ast}\pi_{16}(\mathit{SO}(9)) +
			2\pi_{16}(\mathit{SO}(10)).$ By \cite[Proposition 4.1]{Ka}, we find
			\begin{equation}\label{1026abc}
				\pi_{16}(\mathit{SO}(10))=\mathrm{span}\{\Delta(\sigma_{10}),[\iota_7]_{10}\nu_7^3,[\iota_7]_{10}\mu_7\}
				\cong\mathbb Z/16\oplus(\mathbb Z/2)^2.
			\end{equation}
			\noindent   Since $p_{10*}\Delta\sigma_{10}=2\sigma_9,$   exactness shows that $i_{10*}\pi_{16}(\mathit{SO}(9))$ is annihilated by $8.$  The same is true of $2\pi_{16}(\mathit{SO}(10)).$  Hence
			\[
			\bigl(i_{10*}\pi_{16}(\mathit{SO}(9))+2\pi_{16}(\mathit{SO}(10))\bigr)\circ8\sigma_{16}=0.
			\]
			Thus the preceding bracket calculation gives
			\[
			\{[\iota_7]_{10},\bar{\nu}_7+\varepsilon_7,2\iota_{15}\}\circ8\sigma_{16}
			=8\Delta(\sigma_{10}^2).
			\]
			Hence, by Lemma~\ref{m1:Derh}(1), we obtain
			$[\iota_7]_{10}\mu_7\sigma_{16}=8\Delta(\sigma^2_{10}) \neq 0.$ 
			This leads to the first half of (2).    In light of (1),
			we have the equation
			\[\begin{aligned}
				&[\zeta_5]_{10}\sigma_{16} + [\iota_7]_{10}\mu_7\sigma_{16} = 8\Delta(\sigma^2_{10}).
			\end{aligned}\]
			Substituting the first equality of (2), proved above, into this equation gives $[\zeta_5]_{10}\sigma_{16}=0,$ the second half of (2).
			This completes the proof.
		\end{enumerate}
	\end{proof}
	
	We need to examine the  lifts.
	\begin{lem}
		There exist lifts $[\sigma^2_9]$ and
		$[\bar{\nu}_9\nu^2_{17}]$ in
		$\pi_{23}(\mathit{SO}(10))$ with the following properties. The lift $[\sigma^2_9]$ can be chosen so that
		\refstepcounter{equation}\label{m1:[sig^2_9]}%
		$[\sigma^2_9]\in
		\{[\iota_7]_{10},\bar{\nu}_7+\varepsilon_7,\sigma_{15}\}$
		\textup{(\theequation);}
		for this choice, we may choose
		$\psi_{10}\in
		\{\sigma_{10},\bar{\nu}_{17}+\varepsilon_{17},\sigma_{25}\}_4$
		so that
		\refstepcounter{equation}\label{m1:psi10}%
		$J[\sigma^2_9]=\psi_{10}$
		\textup{(\theequation)}.
		
		The lift $[\bar{\nu}_9\nu^2_{17}]$ can be chosen so that
		\refstepcounter{equation}\label{m1:[bn9n^2]}%
		$[\bar{\nu}_9\nu^2_{17}]
		\in\{[[\iota_9,\eta_9]],\nu_{18},\eta_{21}\}_1$
		\textup{(\theequation)}.
	\end{lem}
	
	\begin{proof}
		Consider the exact sequence $(23)_9$:
		\[
		\pi^9_{24}\xrightarrow{\Delta}
		\pi_{23}(\mathit{SO}(9))\xrightarrow{i_*}
		\pi_{23}(\mathit{SO}(10))\xrightarrow{p_*}
		\pi^9_{23}\xrightarrow{\Delta}
		\pi_{22}(\mathit{SO}(9)),
		\]
		where
		$\pi^9_{23}=\mathrm{span}\{\sigma^2_9,\kappa_9\}
		\cong\mathbb Z/16\oplus\mathbb Z/4$.
		By \eqref{m1:Desig9},
		$\Delta(\sigma^2_9)
		=[\iota_7]_9(\bar{\nu}_7+\varepsilon_7)\sigma_{15}=0$.
		According to
		\cite[Lemma 5.4 and the proof of Lemma 5.5]{HKM},
		$\Delta\kappa_9
		=[ \nu_5]_9\kappa_8+[\iota_7]_9\bar{\varepsilon}_7$
		has order $2$.
		Hence
		$\Delta(\bar{\nu}_9\nu^2_{17})
		=2\Delta\kappa_9=0$.
		Exactness therefore gives lifts of
		$\sigma^2_9$ and $\bar{\nu}_9\nu^2_{17}$
		to $\pi_{23}(\mathit{SO}(10))$.
		
		To obtain the required choice of $[\sigma^2_9]$,
		we apply Proposition~\ref{mdsyq}:
		\begin{align*}
			{p_{10}}_*
			\{[\iota_7]_{10},\bar{\nu}_7+\varepsilon_7,\sigma_{15}\}
			&=-\{p_{10},[\iota_7]_{10},
			\bar{\nu}_7+\varepsilon_7\}\circ\sigma_{16}\\
			&\supseteq
			-\{p_{10},i_{10},\Delta\iota_9\}\circ\sigma^2_9
			\ni\pm\sigma^2_9.
		\end{align*}
		The description of $\pi_{16}(\mathit{SO}(10))$
		in \eqref{1026abc} shows that the indeterminacy of
		$\{p_{10},[\iota_7]_{10},
		\bar{\nu}_7+\varepsilon_7\}\circ\sigma_{16}$
		is
		$p_{10}\circ\pi_{16}(\mathit{SO}(10))\circ\sigma_{16}
		=\mathrm{span}\{2\sigma^2_9\}$.
		So, the displayed image contains both
		$\sigma^2_9$ and $-\sigma^2_9$.
		Thus a lift satisfying \eqref{m1:[sig^2_9]} can be chosen.
		
		By \cite[Lemma 5.1]{MiM} and \eqref{m1:Ji7},
		\begin{align*}
			J[\sigma^2_9]
			&\in J\{[\iota_7]_{10},
			\bar{\nu}_7+\varepsilon_7,\sigma_{15}\}\\
			&\subseteq
			\{J[\iota_7]_{10},
			\bar{\nu}_{17}+\varepsilon_{17},\sigma_{25}\}_{10}\\
			&=\{\sigma_{10},
			\bar{\nu}_{17}+\varepsilon_{17},\sigma_{25}\}_{10}.
		\end{align*}
		The last bracket is contained in
		$\{\sigma_{10},
		\bar{\nu}_{17}+\varepsilon_{17},\sigma_{25}\}_4$.
		In \cite[(4.27)]{MMO} and its proof,
		$\psi_{10}$ may be chosen as any element of this latter bracket.
		We may therefore choose $\psi_{10}$ to satisfy
		\eqref{m1:psi10}.
		
		For the second lift, the description of
		$\pi_{21}(\mathit{SO}(10))$ in \cite{HKM},
		together with the corresponding $J$-images, gives
		$[[\iota_9,\eta_9]]\circ\nu_{18}=0$.
		Since $\nu_{17}\eta_{20}=0$, the bracket
		$\{[[\iota_9,\eta_9]],\nu_{18},\eta_{21}\}_1$
		is defined.
		
		The relations
		$\bar{\nu}_9\nu^2_{17}=[\iota_9,\nu^2_9]$
		from \cite[p.~102]{T} and
		$\nu^2_{17}=\{\eta_{17},\nu_{18},\eta_{21}\}$
		from \cite[Lemma 5.12]{T} imply that
		$\{[\iota_9,\eta_9],\nu_{18},\eta_{21}\}$
		contains
		$[\iota_9,\iota_9]\circ\nu^2_{17}
		=\bar{\nu}_9\nu^2_{17}$.
		Its indeterminacy is
		$[\iota_9,\eta_9]\circ\pi^{18}_{23}
		+\pi^9_{22}\circ\eta_{22}=0$.
		Consequently,
		$\{[\iota_9,\eta_9],\nu_{18},\eta_{21}\}
		=\{\bar{\nu}_9\nu^2_{17}\}$.
		
		By naturality,
		$p_{10}\circ
		\{[[\iota_9,\eta_9]],\nu_{18},\eta_{21}\}_1
		\subseteq
		\{[\iota_9,\eta_9],\nu_{18},\eta_{21}\}
		=\{\bar{\nu}_9\nu^2_{17}\}$.
		The bracket on the left is nonempty, so any of its elements
		can be chosen as a lift $[\bar{\nu}_9\nu^2_{17}]$
		satisfying \eqref{m1:[bn9n^2]}.
	\end{proof}
	\begin{lem}
		$\{\zeta_{14}, \nu_{25}, \eta_{28}\}_1 \equiv 4\omega_{14} \bmod
		\sigma_{14}\mu_{21}.$ 
	\end{lem}
	\begin{proof}
		Since $\zeta_{13}\nu_{24}=8\sigma^2_{13}=P(\eta^2_{27})$ by
		\cite[(10.10)]{T}, we obtain $$H\{\zeta_{14}, \nu_{25}, \eta_{28}\}_1
		= -P^{-1}(8\sigma^2_{13}) \circ
		\eta_{29}=\{\eta^3_{27}\}=\{4\nu_{27}\}.$$ Together with $H(\omega_{14})=\nu_{27}$ and \cite[Proposition 2.6]{T}, this implies that $$\{\zeta_{14}, \nu_{25}, \eta_{28}\}_1 \equiv 4\omega_{14} \bmod
		\Sigma\pi^{13}_{29}=\mathrm{span}\{\sigma_{14}\mu_{21}\}.$$\noindent This completes the proof.
	\end{proof}
	
	We next present another explicit calculation,
	together with the relations relevant to its description.
	\begin{prop} \label{m1:R10}
		The structure of $\pi_{23}(\mathit{SO}(10))$ is as follows.
		\[ \GroupLayout{\pi_{23}(\mathit{SO}(10)) = \mathrm{span}\{[\sigma^2_9],\; [\bar{\nu}_9\nu^2_{17}],\; [P(\Sigma\theta)+\nu_6\kappa_9]_{10}\}\cong \mathbb Z/32\oplus (\mathbb Z/2)^{2},}{\begin{aligned}\pi_{23}(\mathit{SO}(10)) = \mathrm{span}\{&[\sigma^2_9],\; [\bar{\nu}_9\nu^2_{17}],\\
				&[P(\Sigma\theta)+\nu_6\kappa_9]_{10}\}\cong \mathbb Z/32\oplus (\mathbb Z/2)^{2},\end{aligned}} \]
		with relations $[ \nu_5]_{10}\bar{\varepsilon}_8=0,$  
		$\Delta\rho'\equiv [\zeta_5]_9\sigma_{16} \bmod
		[ \nu_5]_9\bar{\varepsilon}_8$
		and $16[\sigma^2_9] = 8\Delta(\sigma^2_{10})
		= [\iota_7]_{10}\mu_7\sigma_{16}.$ 
	\end{prop}
	\begin{proof}
		By $\pi_{24}^9=\mathrm{span}\{\rho',\sigma_9\bar{\nu}_{16},\sigma_9\varepsilon_{16},\bar{\varepsilon}_9\},$  
		\eqref{m1:Desig9e} and  $\Delta\bar{\varepsilon}_9 =0$ \eqref{m1:Dbep9},
		we deduce $$\Delta\pi^9_{24}=\mathrm{span}\{\Delta\rho',[ \nu_5]_9\bar{\varepsilon}_8\}.$$\noindent Hence,
		by Lemmas \ref{m1:Derh}(2) and \ref{m1:i710ms}(2), we obtain
		$[ \nu_5]_{10}\bar{\varepsilon}_8=0$
		and
		$\Delta\rho'\equiv [\zeta_5]_9\sigma_{16} \bmod
		[ \nu_5]_9\bar{\varepsilon}_8.$ 
		
		We determine the
		orders of $[\sigma^2_9]$ and $[\bar{\nu}_9\nu^2_{17}].$  By using the relation
		\eqref{m1:[sig^2_9]},
		$$
		\{\bar{\nu}_7+\varepsilon_7, \sigma_{15}, 16\iota_{22}\}
		\ni\mu_7\sigma_{16} + \Sigma\zeta'
		$$
		from Lemma~\ref{m1:m3sg}(3) and  the relation
		$[\iota_7]_{10}(\Sigma\zeta')
		= [\iota_7]_{10}\eta_7\bar{\varepsilon}_8= 0$
		together with $2\pi^7_{23}=0,$   we infer that
		\begin{align*}
			16[\sigma^2_9]&\in\{[\iota_7]_{10},\bar{\nu}_7+\varepsilon_7, \sigma_{15}\}\circ 16\iota_{23} \\
			&= [\iota_7]_{10}\circ\{\bar{\nu}_7+\varepsilon_7, \sigma_{15}, 16\iota_{22}\} \\
			&\ni[\iota_7]_{10}\circ(\mu_7\sigma_{16}+\Sigma\zeta') \\
			&= [\iota_7]_{10}\mu_7\sigma_{16} \ \bmod\ [\iota_7]_{10}\circ 16\pi^7_{23}=0.
		\end{align*}
		Thus we obtain $16[\sigma^2_9] = 8\Delta(\sigma^2_{10})
		= [\iota_7]_{10}\mu_7\sigma_{16}\neq 0$ in light of Lemma~\ref{m1:i710ms}.
		Hence,  the order of $[\sigma^2_9]$ is $32.$ 
		By \eqref{m1:[bn9n^2]}, we conclude that
		\[\begin{aligned}
			&2[\bar{\nu}_9\nu^2_{17}]\in\{[[\iota_9,\eta_9]], \nu_{18}, \eta_{21}\} \circ 2\iota_{23} = - [[\iota_9,\eta_9]]\circ\{\nu_{18}, \eta_{21}, 2\iota_{22}\} \subseteq[[\iota_9,\eta_9]]\circ\pi^{18}_{23} = 0.
		\end{aligned}\]
		So the order of $[\bar{\nu}_9\nu^2_{17}]$ is $2.$  This determines the group
		$\pi_{23}(\mathit{SO}(10)).$ 
		This completes the proof.
	\end{proof}
	
	In view of  \cite[p. 22]{MMO}, we see that
	\begin{equation}
		\begin{aligned}
			&P(\rho_{19})=\sigma^2_9\mu_{23}.\label{m1:sigsigmu}
		\end{aligned}
	\end{equation}
	
	\noindent   By
	$\eta_5\varepsilon_6\kappa_{14}=4\nu_5\bar\kappa_8$ \cite[I-Proposition 3.1(2)]{Oda1979}
	and $4\bar\kappa_7=\nu^2_7\kappa_{13}$ \cite[Lemma 15.4]{MT1}, we infer that
	\begin{equation}
		\begin{aligned}
			&\varepsilon_5\bar{\varepsilon}_{13}=\varepsilon_5\eta_{13}\kappa_{14} =\eta_5\varepsilon_6\kappa_{14}=4\nu_5\bar\kappa_8=\nu^3_5\kappa_{14}.\label{m1:ep5bep}
		\end{aligned}
	\end{equation}
	
	\subsection{The Case $n=11$}
	
	The next boundary formulas express the relations imposed when passing to the case $n=11.$  They also show
	how the chosen $\sigma^2$-lift contributes to the resulting cyclic decomposition.
	
	\begin{prop} \label{m1:R11}
		The structure of $\pi_{23}(\mathit{SO}(11))$ is as follows.
		$$\GroupLayout{\pi_{23}(\mathit{SO}(11)) =\mathrm{span}\{[P(\Sigma\theta)+\nu_6\kappa_9]_{11},\; [\sigma^2_9]_{11}\}\cong (\mathbb Z/2)^{2}.}{\begin{aligned}\pi_{23}(\mathit{SO}(11)) =\mathrm{span}\{&[P(\Sigma\theta)+\nu_6\kappa_9]_{11},\\
				&[\sigma^2_9]_{11}\}\cong (\mathbb Z/2)^{2}.\end{aligned}}$$
		Here, $\Delta(\sigma^2_{10}) = 2x[\sigma^2_9]$ for  some  odd $x,$   $\Delta\kappa_{10}
		\equiv [\bar{\nu}_9\nu^2_{17}] + [P(\Sigma\theta)+\nu_6\kappa_9]_{10}
		\bmod\ 16[\sigma^2_9].$ 
	\end{prop}
	\begin{proof}
		In the exact sequence $(23)_{10}$:
		\[\begin{aligned}
			&\pi^{10}_{24}\xrightarrow{\Delta}\pi_{23}(\mathit{SO}(10))\xrightarrow{i_*}\pi_{23}(\mathit{SO}(11)) \xrightarrow{p_*}\pi^{10}_{23}\xrightarrow{\Delta}\pi_{22}(\mathit{SO}(10)),
		\end{aligned}\]
		$\Delta\colon  \pi^{10}_{23}\to \pi_{22}(\mathit{SO}(10))$ is a monomorphism and $\pi^{10}_{24}
		= \frac{\mathbb Z}{16}\{\sigma^2_{10}\}\oplus\frac{\mathbb Z}{2}\{\kappa_{10}\}.$ 
		
		Since ${p_{10}}_*\Delta(\sigma^2_{10})=2\sigma^2_9$
		by  \eqref{m1:DelS2}, we deduce $$\Delta(\sigma^2_{10}) - 2[\sigma^2_9]\in{i_{10}}_*\pi_{23}(\mathit{SO}(9))
		= \mathrm{span}\{[P(\Sigma\theta)+\nu_6\kappa_9]_{10}, [\iota_7]_{10}\mu_7\sigma_{16}\}.$$
		We know that $J\Delta(\sigma^2_{10}) = [\iota_{10}, \sigma^2_{10}] = 2\psi_{10}$
		\cite[(4.27)]{MMO}.
		We also know \eqref{m1:JP} and \cite[(4.16)]{MMO}:
		\begin{equation}
			\begin{aligned}
				&J[P(\Sigma\theta)+\nu_6\kappa_9]_{10} =\bar{\kappa}_{10}\nu_{30} - \nu_{10}\bar{\kappa}_{13} =P(\kappa_{21})=J\Delta\kappa_{10}.\label{m1:JPeth10}
			\end{aligned}
		\end{equation}
		Hence we take $\Delta(\sigma^2_{10})\equiv 2[\sigma^2_9] \bmod
		[\iota_7]_{10}\mu_7\sigma_{16} = 16[\sigma^2_9]$ by Proposition~\ref{m1:R10}.
		This leads to the first relation.
		
		Since ${p_{10}}_*\Delta\kappa_{10}=2\kappa_9=\bar{\nu}_9\nu^2_{17}$
		by \eqref{m1:DelS2}, we have
		\begin{equation}
			\begin{aligned}
				&\Delta\kappa_{10}-[\bar{\nu}_9\nu^2_{17}]\in{i_{10}}_*\pi_{23}(\mathit{SO}(9)) =\mathrm{span}\{[P(\Sigma\theta)+\nu_6\kappa_9]_{10}, [\iota_7]_{10}\mu_7\sigma_{16}\}.\label{m1:Dk10}
			\end{aligned}
		\end{equation}
		In light of the relations \eqref{m1:[bn9n^2]},
		$J[[\iota_9,\eta_9]] = \pm \sigma_{10}\zeta_{17}$ \cite[p.~343]{KM1} and
		$\zeta_{17}\nu_{28}=\nu_{17}\zeta_{20}=0$ from \eqref{m1:n9z}, we conclude that
		\begin{align*}
			J[\bar{\nu}_9\nu^2_{17}]&\in J\{[[\iota_9,\eta_9]], \nu_{18}, \eta_{21}\} \\
			&\subseteq\{J[[\iota_9,\eta_9]], \nu_{28}, \eta_{31}\} \\
			&=\{\pm\sigma_{10}\zeta_{17}, \nu_{28}, \eta_{31}\} \\
			&\supseteq\pm\sigma_{10}\circ\{\zeta_{17}, \nu_{28}, \eta_{31}\}.
		\end{align*}
		We note that $\sigma_{10} \circ \Sigma^3\{\zeta_{14}, \nu_{25}, \eta_{28}\}_1
		\subseteq \mathrm{span}\{\sigma^2_{10}\mu_{24}\}=0$ by (\ref{m1:sigsigmu}).
		The decomposition
		$\pi^{14}_{30}
		= \frac{\mathbb Z}{8}\{\omega_{14}\}\oplus\frac{\mathbb Z}{2}\{\sigma_{14}\mu_{21}\},$  
		together with the relation
		$2\omega_{14}=P(\nu_{29})$ \cite[(12.20)]{T}
		and \eqref{m1:sigsigmu}, gives
		$\sigma_{10}\circ\{\zeta_{17}, \nu_{28}, \eta_{31}\}=0.$ 
		Consequently, using
		$\pi^{28}_{33}=0$ and
		$\pi^{10}_{32}=\mathrm{span}\{\sigma_{10}\rho_{17},\varepsilon_{10}\kappa_{18},\nu_{10}\bar\sigma_{13}\},$  
		we obtain
		\[\begin{aligned}
			&J[\bar{\nu}_9\nu^2_{17}]\equiv 0 \bmod\sigma_{10}\zeta_{17}\circ\pi^{28}_{33}+\pi^{10}_{32}\circ\eta_{32} =\mathrm{span}\{\sigma_{10}\rho_{17}\eta_{32}, \varepsilon_{10}\kappa_{18}\eta_{32}, \nu_{10}\bar\sigma_{13}\eta_{32}\}.
		\end{aligned}\]
		\noindent Since $\sigma_{10}\rho_{17}\eta_{32}=\sigma^2_{10}\mu_{24} = 0,$  
		$\varepsilon_{10}\kappa_{18}\eta_{32}=4\nu_{10}\bar\kappa_{13}$ by \eqref{m1:ep5bep} and
		$\nu_{10}\bar\sigma_{13}\eta_{32}=\nu_{10}\eta_{13}\bar\sigma_{14}=0$ by \cite[(2.7)]{MMO},
		we find $J[\bar{\nu}_9\nu^2_{17}]\equiv 0\ \bmod\ 4\nu_{10}\bar\kappa_{13}.$ 
		Hence, applying $J$ to \eqref{m1:Dk10}, we conclude that
		\[\begin{aligned}
			&J\Delta\kappa_{10}\equiv 0 \ \bmod\ J[P(\Sigma\theta)+\nu_6\kappa_9]_{10}, J([\iota_7]_{10}\mu_7\sigma_{16}), 4\nu_{10}\bar\kappa_{13}.
		\end{aligned}\]
		By \eqref{m1:s_9smu}, \eqref{m1:sigsigmu}, the compatibility of $J$ with suspension in Section~\ref{m2:subsec:j-toda}, and $\Sigma P=0$, we have
		\[
		\begin{aligned}
			J([\iota_7]_{10}\mu_7\sigma_{16})
			&=\pm\Sigma J([\iota_7]_9\mu_7\sigma_{16})\\
			&=\pm\Sigma(\sigma_9^2\mu_{23})=\pm\Sigma P(\rho_{19})=0.
		\end{aligned}
		\]
		Furthermore, using \eqref{m1:JPeth10}, this identity, and the relation
		$\bar{\kappa}_{10}\nu_{30} - \nu_{10}\bar{\kappa}_{13}\ne 4\nu_{10}\bar\kappa_{13}
		\in\pi^{10}_{33}$, we obtain
		\[\begin{aligned}
			&\Delta\kappa_{10} \equiv [\bar{\nu}_9\nu^2_{17}] +[P(\Sigma\theta)+\nu_6\kappa_9]_{10} \ \bmod\ [\iota_7]_{10}\mu_7\sigma_{16}.
		\end{aligned}\]
		
		\noindent   Thus, by Proposition~\ref{m1:R10}, we obtain the second relation and the group $\pi_{23}(\mathit{SO}(11)).$ 
		This completes the proof.\end{proof}
	
	The preceding considerations lead to the
	relations that deserve to be recorded explicitly.
	\begin{lem} \label{m1:fin}
		$[8\sigma_8]_{11}\sigma_{15}\eta_{22}=0$
		and\ \ $[\bar{\nu}_9\nu^2_{17}]_{11} = [P(\Sigma\theta)+\nu_6\kappa_9]_{11}
		=[[\iota_{10},\nu_{10}]]\eta_{22}.$ 
	\end{lem}
	\begin{proof}
		By using
		$\pi_{16}(\mathit{SO}(11)) =\mathrm{span}\{[\iota_7]_{11}\mu_7,[\iota_7]_{11}\nu^3_7\}$
		\cite[Proposition 4.1]{Ka},
		and Lemma~\ref{m1:i710ms}(2),
		we have $$[8\sigma_8]_{11}\sigma_{15}\eta_{22}
		= [8\sigma_8]_{11}\eta_{15}\sigma_{16}\in \pi_{16}(\mathit{SO}(11))\circ\sigma_{16}
		= 0.$$\noindent Hence, we deduce 
		\begin{equation}\label{1030bu}
			[8\sigma_8]_{11}\sigma_{15}\eta_{22}=0.
		\end{equation}
		As a consequence of           Proposition~\ref{m1:R11}, we find
		$[\bar{\nu}_9\nu^2_{17}]_{11}=[P(\Sigma\theta)+\nu_6\kappa_9]_{11}\neq0.$ 
		To obtain  $[\bar{\nu}_9\nu^2_{17}]_{11} = [[\iota_{10},\nu_{10}]]\eta_{22},$  
		we show that $[\eta_9\varepsilon_{10}]\nu_{18} = 0.$  Since
		$\eta_9\varepsilon_{10}\nu_{18} = 0$  and $\pi_{21}(\mathit{SO}(9))
		= \mathrm{span}\{[\iota_7]_9\kappa_7\}$ \cite[Theorem 0.1]{HKM}, we conclude that
		$$[\eta_9\varepsilon_{10}]\nu_{18}\in i_*\pi_{21}(\mathit{SO}(9))
		= \frac{\mathbb Z}{2}\{[\iota_7]_{10}\kappa_7\}.$$ \noindent  Thus we can set
		$[\eta_9\varepsilon_{10}]\nu_{18} = x[\iota_7]_{10}\kappa_7$ for $x = 0$ or
		$1.$  Since $J[\eta_9\varepsilon_{10}]\equiv \lambda'' \bmod
		\sigma_{10}\zeta_{17}, \eta_{10}\bar{\mu}_{11}$ \cite[Lemma 2.2]{HKM}, 
		$\zeta_{17}\nu_{28} = \nu_{17}\zeta_{20}=0$ by \eqref{m1:n9z}  and $\eta_{27}\nu_{28} = 0,$   we see that
		$\lambda''\nu_{28} = x\sigma_{10}\kappa_{17}.$  Therefore, since
		$\lambda'\nu_{29} = \sigma_{11}\kappa_{18}$ \cite[I-Proposition 3.1(1)]{Oda1979} and
		$\Sigma\lambda''= 2\lambda'$ \cite[Lemma 12.19]{T},
		we obtain $0 = x\sigma_{11}\kappa_{18}.$  Since
		$\sigma_{11}\kappa_{18}\neq 0,$   we obtain $x = 0$ and the relation
		\[\begin{aligned}
			&[\eta_9\varepsilon_{10}]\nu_{18} = 0.
		\end{aligned}\]
		In light of the relation
		$\Delta[\iota_{10},\iota_{10}]
		=[[\iota_9,\iota_9]\eta_{17}] + 2[\eta_9\varepsilon_{10}]$
		\cite[p.~343]{KM1}, we deduce
		\[\begin{aligned}
			&[[\iota_9,\eta_9]]_{11} = 2[\eta_9\varepsilon_{10}]_{11}.
		\end{aligned}\]
		Therefore, by the relation \eqref{m1:[bn9n^2]} and
		the fact that $2\pi_{23}(\mathit{SO}(11)) = 0$ from Proposition~\ref{m1:R11},
		we see that\vspace{-1.4em}
		\begin{align*}
			0\neq [\bar{\nu}_9\nu^2_{17}]_{11} &\in\{2[\eta_9\varepsilon_{10}]_{11}, \nu_{18},\eta_{21}\} \\
			&\subseteq\{[\eta_9\varepsilon_{10}]_{11}, 2\nu_{18}, \eta_{21}\} \\
			&\dyd \{[\eta_9\varepsilon_{10}]_{11}, \nu_{18},2\eta_{21}\}  \ni 0
		\end{align*}
		with the indeterminacy
		$[\eta_9\varepsilon_{10}]_{11}\circ\pi^{18}_{23}
		+\pi_{22}(\mathit{SO}(11))\circ\eta_{22}=\pi_{22}(\mathit{SO}(11))\circ\eta_{22}.$ 
		We note that $\pi_{22}(\mathit{SO}(11)) = \mathrm{span}\{[[\iota_{10},\nu_{10}]],
		[8\sigma_8]_{11}\sigma_{15}\}$ \cite[Theorem 0.2]{HKM}.
		Then, by $[8\sigma_8]_{11}\sigma_{15}\eta_{22}=0$ which is given in \eqref{1030bu},  
		we obtain $\pi_{22}(\mathit{SO}(11))\circ\eta_{22}
		= \mathrm{span}\{[[\iota_{10},\nu_{10}]]\eta_{22}\}.$  This leads to the relation, completing
		the proof.
	\end{proof}
	
	\subsection{The Cases $n=12,13$}
	
	The next step follows the surviving $\sigma^2$-class into the cases $n=12,13.$  Its order is
	constrained by the nonvanishing relation in a homotopy group of a sphere below, after which the exact sequences determine the
	additional generators.
	
	By the relations \eqref{m1:psi10} and
	$\psi_n\neq 0$ for $10\leq n\leq 21$ \cite[Theorem 1.1(a)]{MMO},
	we have
	\[\begin{aligned}
		&\mathrm{ord}[\sigma^2_9]_n=2\ (11\leq n\leq 21).
	\end{aligned}\]
	\noindent
	By the relation ${p_{12}}_\ast\Delta\theta = \theta'$ from \cite[Lemma 6.2(ii)]{HKM},
	we may take $[\theta']=\Delta\theta\in \pi_{23}(\mathit{SO}(12)).$ 
	We show the following.
	\begin{prop}\label{m1:R12}
		The structure of $\pi_{23}(\mathit{SO}(12))$ is as follows.
		\[ \GroupLayout{\pi_{23}(\mathit{SO}(12)) = \mathrm{span}\{[\theta'],\; [\bar{\nu}_9\nu^2_{17}]_{12},\; [\sigma^2_9]_{12}\}\cong (\mathbb Z/2)^{3}.}{\begin{aligned}\pi_{23}(\mathit{SO}(12)) = \mathrm{span}\{&[\theta'],\; [\bar{\nu}_9\nu^2_{17}]_{12},\\
				&[\sigma^2_9]_{12}\}\cong (\mathbb Z/2)^{3}.\end{aligned}} \]
		Here $[\theta']=\Delta\theta.$ 
	\end{prop}
	\begin{proof}
		In the exact sequence $(23)_{11}$:
		\[\begin{aligned}
			&\pi^{11}_{24}\xrightarrow{\Delta}\pi_{23}(\mathit{SO}(11))\xrightarrow{i_*}\pi_{23}(\mathit{SO}(12)) \xrightarrow{p_*}\pi^{11}_{23}\xrightarrow{\Delta} \pi_{22}(\mathit{SO}(11)),
		\end{aligned}\]
		we see that $\pi^{11}_{23}=\frac{\mathbb Z}{2}\{\theta'\}$ and
		$\pi^{11}_{24}= \frac{\mathbb Z}{2}\{\theta'\eta_{23}\}\oplus\frac{\mathbb Z}{2}\{\sigma_{11}\nu^2_{18}\}.$ 
		We record the relation
		$\eta_{11}\theta=\sigma_{11}\nu^2_{18}+\theta'\eta_{23}.$  Therefore, we can take
		$\pi^{11}_{24}=\mathrm{span}\{\theta'\eta_{23}, \eta_{11}\theta\}.$ 
		We conclude that $2[\theta']=\Delta(2\theta)= 0.$ 
		We see that  $\Delta(\theta'\eta_{23}) =\Delta\theta'\circ\eta_{22}= 0$ and
		$\Delta(\eta_{11}\theta)=0.$ 
		
		So Proposition~\ref{m1:R11} implies the conclusion.
	\end{proof}
	
	Recall  Mahowald's result: 
	\begin{equation}\label{md1020}
		\mathrm{ord}\Delta([\iota_{12},\iota_{12}])=11!/8\;\;\text{before localization at 2}.
	\end{equation}
	
	\noindent   See \cite[Theorem 1.3]{HKM}. Observe that
	$\mathrm{ord}\Delta([\iota_{12},\iota_{12}])=32.$ 
	The discussion proceeds with the following statement, after which we give its proof and then continue with the rest of the section.
	\begin{prop} \label{m1:R13}
		The structure of $\pi_{23}(\mathit{SO}(13))$ is as follows.
		\[ \GroupLayout{\pi_{23}(\mathit{SO}(13)) = \mathrm{span}\{[32[\iota_{12}, \iota_{12}]],\; [\sigma^2_9]_{13}\}\cong \mathbb Z_{(2)}\oplus \mathbb Z/2.}{\begin{aligned}\pi_{23}(\mathit{SO}(13)) = \mathrm{span}\{&[32[\iota_{12}, \iota_{12}]],\\
				&[\sigma^2_9]_{13}\}\cong \mathbb Z_{(2)}\oplus \mathbb Z/2.\end{aligned}} \]
		Here
		$\Delta(\Sigma\theta') = [\bar{\nu}_9\nu^2_{17}]_{12}
		= [P(\Sigma\theta)+\nu_6\kappa_9]_{12}
		= [[\iota_{10},\nu_{10}]]_{12}\eta_{22}.$ 
	\end{prop}
	\begin{proof}
		In the exact sequence $(23)_{12}$:
		\[\begin{aligned}
			&\pi^{12}_{24}\xrightarrow{\Delta}\pi_{23}(\mathit{SO}(12))\xrightarrow{i_*}\pi_{23}(\mathit{SO}(13)) \xrightarrow{p_*}\pi^{12}_{23}\xrightarrow{\Delta}\pi_{22}(\mathit{SO}(12)),
		\end{aligned}\]
		$$\begin{aligned}
			&\operatorname{ker}(\Delta\colon  \pi^{12}_{23}\to \pi_{22}(\mathit{SO}(12))) = \mathrm{span}\{32[\iota_{12}, \iota_{12}]\}.
		\end{aligned}$$
		Notice that\ $\mathrm{ord} P(\zeta_{25})=8$ \cite[p. 320]{Mi1} implies $\mathrm{ord}\Delta(\zeta_{12})=8.$ 
		It follows that $\pi^{12}_{24} = \frac{\mathbb Z}{2}\{\theta\}\oplus\frac{\mathbb Z}{2}\{\Sigma\theta'\}.$  $\Delta\theta = [\theta']$ by Proposition~\ref{m1:R12} and $J[\sigma^2_9]_{13} = \psi_{13}$ by
		\eqref{m1:psi10}.
		By the relations $\Sigma\theta' = [\iota_{12},\iota_{12}]\circ\eta_{23},$  
		\cite[Lemma 6.3(ii)]{HKM} and Lemma~\ref{m1:fin}, we deduce
		\[\begin{aligned}
			&\Delta(\Sigma\theta') = \Delta([\iota_{12}, \iota_{12}])\circ\eta_{22} \equiv [[\iota_{10},\nu_{10}]]_{12}\eta_{22} \bmod [8\sigma_8]_{12}\sigma_{15}\eta_{22} = 0.
		\end{aligned}\]
		In light of Lemma~\ref{m1:fin}, we obtain
		\[\begin{aligned}
			&\Delta(\Sigma\theta') = [[\iota_{10},\nu_{10}]]_{12}\eta_{22} = [P(\Sigma\theta)+\nu_6\kappa_9]_{12} = [\bar{\nu}_9\nu^2_{17}]_{12}.
		\end{aligned}\]
		This determines the group $\pi_{23}(\mathit{SO}(13)),$   completing the proof.
	\end{proof}

	We will use James' theorem (i.e., Proposition \ref{jmstm}) to examine the $J$-images, by which we will see that  $[32[\iota_{12}, \iota_{12}]]_{m}$ can not be divisible by 2 for each $m>13$.

	\begin{lem} 
		
		\begin{itemize}
			\item[\rm (1)]  $\pi_{23}(\lo^{\wq}S^{\wq},\lo^{13}S^{13})\cong\z/16$ and so $\pi_{23}(\mathit{SO},\mathit{SO}(13))\cong\z/16.$
				\item[\rm (2)]  The boundary homomorphism $\pa\colon\pi_{23}(\mathit{SO},\mathit{SO}(13))\to \pi_{22}(\mathit{SO}(13))$ is an isomorphism. The  homomorphism 
				\[
				\begin{array}{rccc}
					i_{*}=(i_{13,\wq})_{*}\colon \hspace{-0.85em}& \pi_{23}(\mathit{SO}(13))
					& \longrightarrow
					& \pi_{23}(\mathit{SO}) \\[-3pt]
					& \Vert && \Vert \\[-3pt]
					& \z_{(2)}\jia\z/2
					&& \z_{(2)}
				\end{array}
				\]
				is surjective.
					\item[\rm (3)] For the generator $[32[\iota_{12}, \iota_{12}]]\in\pi_{23}(\mathit{SO}(13))$, we have 
					$$ J([32[\iota_{12}, \iota_{12}]]_{\wq})\ty\bar\co\m 2\bar\co.$$And hence $[32[\iota_{12}, \iota_{12}]]_{m}$ can not be divisible by 2 for each $m>13$.
				In addition, we have	$$J[32[\iota_{12}, \iota_{12}]]\ty\bar\co_{13}\m 2\bar\co_{13},4\nu_{13}\ct_{16},\psi_{13},\phi_{13}.$$ 
		\end{itemize}
		
	\end{lem}

	\begin{proof}\begin{itemize}
			\item [\rm (1)] The isomorphism $\pi_{23}(\lo^{\wq}S^{\wq},\lo^{13}S^{13})\cong\z/16$ is easy to obtain by using the exact sequence.  Then, $\pi_{23}(\mathit{SO},\mathit{SO}(13))\cong\z/16$ follows from  Proposition \ref{jmstm}.
				\item [\rm (2)] Note that we have the exact sequence 
				\[
				\begin{array}{ccccc}
					\pi_{23}(\mathit{SO},\mathit{SO}(13))
					& \stackrel{\pa}\longrightarrow
					& \pi_{22}(\mathit{SO}(13))
					& \longrightarrow
					& \pi_{22}(\mathit{SO}) \\[-4pt]
					\Vert && \Vert && \Vert \\[-4pt]
					\z/16 && \z/16 && \,\;0.
				\end{array}
				\] Thus, this $\pa$ is an isomorphism. In view of the exactness, the displayed homomorphism $(i_{13,\wq})_{*}$ is surjective.
					\item [\rm (3)] Mahowald and Tangora's result \cite{imj} tells us that $$J(\pi_{23}(\mathit{SO}))=\mathrm{span}\n\bar\co\nn;$$ also see \cite[Remark 1.3]{MMO}.  Recall that $$\pi_{23}(\mathit{SO}(13)) = \mathrm{span}\{[32[\iota_{12}, \iota_{12}]],\; [\sigma^2_9]_{13}\}\cong \mathbb Z_{(2)}\oplus \mathbb Z/2.$$
				According to (2), we have $$ J([32[\iota_{12}, \iota_{12}]]_{\wq})\ty\bar\co\m 2\bar\co.$$
				In light of the relation $$\psi\ty 8\bar\co,4\nu\ct,\kj$$ from \cite[(4.39)]{MMO} and the group decompositions in stem 23, we infer the result.
		\end{itemize}
		
	\end{proof}

	Using the fact that $\pi_{23}(\mathit{SO})\cong \mathbb Z_{(2)}$ \cite{Bo} and the results in \cite[Table 1]{HoM},
	we have Table~\ref{m1:R23}.
	\par
	\begin{table}[H]
		\centering
		\small\setlength{\tabcolsep}{3pt}\begin{tabular}{cccccc}\toprule
			$\pi_{23}(\mathit{SO}(13))$ & $\mathbb Z_{(2)} \oplus \mathbb Z/2$
			& $\pi_{23}(\mathit{SO}(14))$ & $\mathbb Z_{(2)} \oplus \mathbb Z/2$
			& $\pi_{23}(\mathit{SO}(15))$ & $\mathbb Z_{(2)} \oplus (\mathbb Z/2)^3$  \\ \midrule
			$\pi_{23}(\mathit{SO}(16))$ & $\mathbb Z_{(2)} \oplus (\mathbb Z/2)^5$
			& $\pi_{23}(\mathit{SO}(17))$ & $\mathbb Z_{(2)} \oplus (\mathbb Z/2)^3$
			& $\pi_{23}(\mathit{SO}(18))$ & $\mathbb Z_{(2)} \oplus (\mathbb Z/2)^2$\\ \bottomrule
		\end{tabular}
		\caption{$\pi_{23}(\mathit{SO}(n))$ for $13\leq n\leq 18$}\label{m1:R23}
	\end{table}
	\par
	
	In addition, we have Table~\ref{m1:R23K} from \cite{Ke}.
	\par
	\begin{table}[H]
		\centering
		\small\setlength{\tabcolsep}{3pt}\begin{tabular}{cccccc}\toprule
			$\pi_{23}(\mathit{SO}(19))$ & $\mathbb Z_{(2)} \oplus \mathbb Z/2$
			& $\pi_{23}(\mathit{SO}(20))$ & $\mathbb Z_{(2)} \oplus \mathbb Z/2$
			& $\pi_{23}(\mathit{SO}(21))$ & $\mathbb Z_{(2)} \oplus \mathbb Z/2$  \\ \midrule
			$\pi_{23}(\mathit{SO}(22))$ & $\mathbb Z_{(2)}$
			& $\pi_{23}(\mathit{SO}(23))$ & $\mathbb Z_{(2)}$
			& $\pi_{23}(\mathit{SO}(24))$ & $\mathbb Z_{(2)}\oplus\mathbb Z_{(2)}$\\ \bottomrule
		\end{tabular}
		\caption{$\pi_{23}(\mathit{SO}(n))$ for $19\le n\le 24$}\label{m1:R23K}
	\end{table}
	\par
	
	By Tables~\ref{m1:R23} and \ref{m1:R23K}, we have the group structure of $\pi_{23}(\mathit{SO}(n))$ for $n \ge 13.$ 
	We now determine the generators and relations in
	$\pi_{23}(\mathit{SO}(n))$ for $n \ge 13.$ 
	
	\subsection{The Case $n=14$}
	
	For $n=14,$   injectivity of the outgoing connecting homomorphism removes the possibility of new lifts.
	The calculation is consequently determined by the image from $\mathit{SO}(13)$ and the incoming boundary
	relations.
	
	In the exact sequence $(23)_{13}$:
	\[\begin{aligned}
		&\pi^{13}_{24}\xrightarrow{\Delta}\pi_{23}(\mathit{SO}(13))\xrightarrow{i_*}\pi_{23}(\mathit{SO}(14)) \xrightarrow{p_*}\pi^{13}_{23}\xrightarrow{\Delta}\pi_{22}(\mathit{SO}(13)),
	\end{aligned}\]
	$\Delta\colon  \pi^{13}_{23}\to \pi_{22}(\mathit{SO}(13))$ is a monomorphism \cite[Lemma 6.7]{HKM},
	$\pi^{13}_{24}=\mathrm{span}\{\zeta_{13}\}$
	and
	$\Delta\zeta_{13} = 0$ by Lemma~\ref{m1:zeta1}. Thus we obtain the following.
	\begin{prop} \label{m1:R14}
		The structure of $\pi_{23}(\mathit{SO}(14))$ is as follows.
		\[ \GroupLayout{\pi_{23}(\mathit{SO}(14)) = \mathrm{span}\{[32[\iota_{12}, \iota_{12}]]_{14},\; [\sigma^2_9]_{14}\}\cong \mathbb Z_{(2)}\oplus \mathbb Z/2.}{\begin{aligned}\pi_{23}(\mathit{SO}(14)) = \mathrm{span}\{&[32[\iota_{12}, \iota_{12}]]_{14},\\
				&[\sigma^2_9]_{14}\}\cong \mathbb Z_{(2)}\oplus \mathbb Z/2.\end{aligned}} \]
		\qed
	\end{prop}
	\section{The Cases $n\geq15$}
	\subsection{The Cases $n=15,16$}
	
	We now use the order-two lift criterion from Chapter~\ref{ch:fundamental}. It provides a controlled
	choice of lift for the $\eta\varepsilon$-class before the exact sequence is used to assemble the
	groups.

	Recall  that $\q(\ca_{14})$ is divisible by 2. So, $\q(\ca_{14}\lt_{15})=\q(\ca_{14}\cg_{15}\ca_{22})=0.$  In view of \cite{Mukaiwhiteheadji}, we know $J\q(\mu_{14})=[\e_{14},\mu_{14}]\neq0.$
	Hence, $$p_{15*}(\pi_{23}(\mathit{SO}(15)) )=\mathrm{ker}(  \pi_{23}^{14}\stackrel{\q}\longrightarrow \pi_{22}(\mathit{SO}(14)) )=\mathrm{span}\n\ca_{14}\lt_{15}, \ca_{14}\cg_{15}\ca_{22}\nn\cong(\z/2)^{2}.$$
	
	\noindent As a consequence of           Lemma~\ref{myjbu1026}(3), we have $\mathrm{ord}[\eta_{14}\varepsilon_{15}]=2.$ 
	
	Note that $\{\eta_{15}, 2\iota_{16}, \nu^2_{16}\}_1\ni \lt_{15}$ from \cite[(6.1)]{T}. Since  $\mathrm{ord} [\eta_{15}]=2$ and $$ \ind ( p_{16}    \hc \{[\eta_{15}],2\iota_{16},\nu^2_{16}\}_1            )=\ind(        \{\eta_{15},2\iota_{16},\nu^2_{16}\}_1                    ),$$
	\noindent   we can take
	\[\begin{aligned}
		&[\varepsilon_{15}]\in\{[\eta_{15}],2\iota_{16},\nu^2_{16}\}_1.
	\end{aligned}\]
	In light of $\pi_{20}(\mathit{SO}(16)) = 0$ \cite[p. 161, Table]{Ke}, we obtain
	$$\pi_{20}(\mathit{SO}(16))\circ\nu_{20}=0,\quad\pi_{17}(\mathit{SO}(16))\circ\nu^2_{17}=0.$$\noindent Hence, we deduce
	\[\begin{aligned}
		&\{[\eta_{15}],2\iota_{16},\nu^2_{16}\}_1 = \{[\eta_{15}],2\nu_{16},\nu_{19}\}_1.
	\end{aligned}\]
	The preceding analysis prepares the way
	for a further determination of the groups under consideration.
	\begin{prop} \label{m1:R1523}
		The structures of $\pi_{23}(\mathit{SO}(m))$ for $m=15,16$ are as follows.
		\begin{enumerate}
			\item $\pi_{23}(\mathit{SO}(15)) = \mathrm{span}\{[32[\iota_{12}, \iota_{12}]]_{15},\; [\sigma^2_9]_{15},\; [\eta_{14}\varepsilon_{15}],\; [\eta_{14}\sigma_{15}]\eta_{22}\}\cong \mathbb Z_{(2)}\oplus (\mathbb Z/2)^{3}.$
			\item \[ \GroupLayout{\pi_{23}(\mathit{SO}(16)) \;= {} \mathrm{span}\{[32[\iota_{12}, \iota_{12}]]_{16},\; [\sigma^2_9]_{16},\; [\eta_{14}\varepsilon_{15}]_{16},\; [\eta_{14}\sigma_{15}]_{16}\eta_{22},\; [\eta_{15}]\sigma_{16},\; [\varepsilon_{15}]\}\cong \mathbb Z_{(2)}\oplus (\mathbb Z/2)^{5}.}{\begin{aligned}\pi_{23}(\mathit{SO}(16)) \;= {} \mathrm{span}\{&[32[\iota_{12}, \iota_{12}]]_{16},\; [\sigma^2_9]_{16},\; [\eta_{14}\varepsilon_{15}]_{16},\\
					&[\eta_{14}\sigma_{15}]_{16}\eta_{22},\; [\eta_{15}]\sigma_{16},\; [\varepsilon_{15}]\}\cong \mathbb Z_{(2)}\oplus (\mathbb Z/2)^{5}.\end{aligned}} \]
		\end{enumerate}
	\end{prop}
	\begin{proof}
		\begin{enumerate}[label=\textup{(\arabic*)},leftmargin=2em]
			\item In the exact sequence $(23)_{14}$:
			\[\begin{aligned}
				&\pi^{14}_{24}\xrightarrow{\Delta}\pi_{23}(\mathit{SO}(14))\xrightarrow{i_*}\pi_{23}(\mathit{SO}(15)) \xrightarrow{p_*}\pi^{14}_{23}\xrightarrow{\Delta}\pi_{22}(\mathit{SO}(14)),
			\end{aligned}\]
			we note
			$\pi^{14}_{24}=\frac{\mathbb Z}{2}\{\eta_{14}\mu_{15}\},$  
			$\operatorname{ker}\ (\Delta\colon\pi^{14}_{23}\to \pi_{22}(\mathit{SO}(14)))
			= \mathrm{span}\{\eta_{14}\varepsilon_{15}, \eta^2_{14}\sigma_{16}\}$ and
			$\Delta(\eta_{14}\mu_{15})=0.$  From the fact that
			$\Delta(\bar\nu_{14})=\Delta(\varepsilon_{14})$~\cite[p.~17]{HKM}, we know that a lift
			$[\eta_{14}\sigma_{15}]$ exists. Recall $2[\eta_{14}\varepsilon_{15}] = 0.$  This leads to (1).
			
			\item Since $J[\eta_{14}\sigma_{15}] \equiv x\sigma^{*\prime}\ \bmod\ \omega_{15}\nu^2_{31}$ for some odd $x$ \cite[the proof of Lemma 6.8]{HKM}, we conclude that
			$J([\eta_{14}\sigma_{15}]\eta_{22}) = \sigma^{*\prime}\eta_{37}.$  Since
			$H(\bar{\varepsilon}^{*\prime}) = \eta_{29}\varepsilon_{30}$ \cite[(3.4)]{MMO}, we see that
			$J[\eta_{14}\varepsilon_{15}]\equiv \bar{\varepsilon}^{*\prime} \bmod
			\Sigma\pi^{14}_{37}.$  This leads to (2), completing the proof.
		\end{enumerate}
	\end{proof}
	
	\subsection{The Cases $n=17,18$}
	
	We now turn to another explicit calculation,
	which we present in the following proposition.
	\begin{prop} \label{m1:R17}
		The structure of $\pi_{23}(\mathit{SO}(17))$ is as follows.
		\[ \GroupLayout{\pi_{23}(\mathit{SO}(17)) = \mathrm{span}\{[32[\iota_{12}, \iota_{12}]]_{17},\; [\sigma^2_9]_{17},\; [\eta_{15}]_{17}\sigma_{16},\; [\varepsilon_{15}]_{17}\}\cong \mathbb Z_{(2)}\oplus (\mathbb Z/2)^{3}.}{\begin{aligned}\pi_{23}(\mathit{SO}(17)) = \mathrm{span}\{&[32[\iota_{12}, \iota_{12}]]_{17},\; [\sigma^2_9]_{17},\\
				&[\eta_{15}]_{17}\sigma_{16},\; [\varepsilon_{15}]_{17}\}\cong \mathbb Z_{(2)}\oplus (\mathbb Z/2)^{3}.\end{aligned}} \]
		Here $\Delta(\eta_{16}\sigma_{17})
		\equiv [\eta_{14}\sigma_{15}]_{16}\eta_{22} \bmod [\sigma^2_9]_{16},$ and
		$\Delta\varepsilon_{16} \equiv [\eta_{14}\varepsilon_{15}]_{16} \bmod
		[\eta_{14}\sigma_{15}]_{16}\eta_{22}, [\sigma^2_9]_{16}.$ 
	\end{prop}
	\begin{proof}
		In the exact sequence $(23)_{16}$:
		\[\begin{aligned}
			&\pi^{16}_{24}\xrightarrow{\Delta}\pi_{23}(\mathit{SO}(16))\xrightarrow{i_*}\pi_{23}(\mathit{SO}(17)) \xrightarrow{p_*}\pi^{16}_{23}\xrightarrow{\Delta}\pi_{22}(\mathit{SO}(16)),
		\end{aligned}\]
		the second $\Delta$ is a monomorphism \cite[the proof of Lemma 6.4]{T} and
		$\pi^{16}_{24} = \frac{\mathbb Z}{2}\{\eta_{16}\sigma_{17}\}
		\oplus\frac{\mathbb Z}{2}\{\varepsilon_{16}\}.$ 
		By  \cite[(4.19)]{Ka} and the standard composition relations, we obtain
		$\Delta(\eta_{16}\sigma_{17})
		= [\eta^2_{14}]_{16}\sigma_{16}$ and
		\[\begin{aligned}
			&\Delta\varepsilon_{16}\in\Delta\{\eta_{16}, 2\iota_{17}, \nu^2_{17}\}_2 \xyd\{[\eta^2_{14}]_{16}, 2\iota_{16}, \nu^2_{16}\}_1 \ni[\eta_{14}\varepsilon_{15}]_{16}.
		\end{aligned}\]
		Hence, we deduce $\Delta(\varepsilon_{16}) \equiv
		[\eta_{14}\varepsilon_{15}]_{16} \bmod
		[\eta_{14}\sigma_{15}]_{16}\eta_{22}, [\sigma^2_9]_{16}.$ 
		This leads to the assertion, completing the proof.
	\end{proof}
	
	We now establish the following result.
	\begin{prop} \label{m1:R18}
		The structure of $\pi_{23}(\mathit{SO}(18))$ is as follows.
		\[ \GroupLayout{\pi_{23}(\mathit{SO}(18)) = \mathrm{span}\{[32[\iota_{12}, \iota_{12}]]_{18},\; [\sigma^2_9]_{18},\; [\varepsilon_{15}]_{18}\}\cong \mathbb Z_{(2)}\oplus (\mathbb Z/2)^{2}.}{\begin{aligned}\pi_{23}(\mathit{SO}(18)) = \mathrm{span}\{&[32[\iota_{12}, \iota_{12}]]_{18},\; [\sigma^2_9]_{18},\\
				&[\varepsilon_{15}]_{18}\}\cong \mathbb Z_{(2)}\oplus (\mathbb Z/2)^{2}.\end{aligned}} \]
		Here
		$[\varepsilon_{15}]_{18}=[\bar{\nu}_{15}]_{18}.$ 
	\end{prop}
	\begin{proof}
		
		Since $P(\nu^2_{35}) =\omega_{17}\nu^2_{33} \neq 0$ \cite[p.~323]{Mi1},
		$\Delta\colon \pi^{17}_{23}\to \pi_{22}(\mathit{SO}(17))$ is a monomorphism.
		Therefore, we find ${i_{18}}_*\pi_{23}(\mathit{SO}(17))=\pi_{23}(\mathit{SO}(18)).$ 
		By the relations
		$\Delta(\iota_{17}) = [\eta_{15}]_{17}$ and
		$[\eta_{15}]\sigma_{16}=[\bar{\nu}_{15}]+[\varepsilon_{15}]$
		we obtain
		$
		\Delta(\sigma_{17})=[\eta_{15}]_{17}\sigma_{16}=[\bar{\nu}_{15}]_{17}+[\varepsilon_{15}]_{17}.
		$
		This completes the proof.
	\end{proof}
	
	\subsection{The Cases $19\leq n\leq22$}
	
	We now establish the following result.
	\begin{prop}\label{m1:R19}
		The structure of $\pi_{23}(\mathit{SO}(19))$ is as follows.
		\[ \GroupLayout{\pi_{23}(\mathit{SO}(19)) = \mathrm{span}\{[32[\iota_{12}, \iota_{12}]]_{19},\; [\sigma^2_9]_{19}\}\cong \mathbb Z_{(2)}\oplus \mathbb Z/2.}{\begin{aligned}\pi_{23}(\mathit{SO}(19)) = \mathrm{span}\{&[32[\iota_{12}, \iota_{12}]]_{19},\\
				&[\sigma^2_9]_{19}\}\cong \mathbb Z_{(2)}\oplus \mathbb Z/2.\end{aligned}} \]
		Here $\Delta(\nu^2_{18}) = [\varepsilon_{15}]_{18}.$ 
	\end{prop}
	\begin{proof}
		It suffices to show the relation.
		By the relation $[\eta_{15}]_{17} = \Delta(\iota_{17}),$   and the
		fact that $\pi_{17}(\mathit{SO}(17)) = \pi_{16}(\mathit{SO}(17))\circ\eta_{16}$ \cite[Theorem
		1(iv)]{KM1}, we conclude that
		\begin{align*}
			[\varepsilon_{15}]_{18} &= {i_{16,18}}_*\{[\eta_{15}],2\iota_{16},\nu^2_{16}\} \\
			&\subseteq {i_{18}}_*\{\Delta(\iota_{17}),2\iota_{16},\nu^2_{16}\} \\
			&= -\{i_{18},\Delta(\iota_{17}),2\iota_{16}\}\circ\nu^2_{17} \\
			& \ni -\Delta(\iota_{18})\circ\nu^2_{17} \\
			&= \Delta(\nu^2_{18}) \bmod {i_{18}}_*\pi_{17}(\mathit{SO}(17))\circ\nu^2_{17} = 0.
		\end{align*}
		This implies  $\Delta(\nu^2_{18}) = [\varepsilon_{15}]_{18}.$ 
	\end{proof}
	
	We consider the exact sequence $(23)_{21}$:
	\[\begin{aligned}
		&\pi^{21}_{24}\xrightarrow{\Delta}\pi_{23}(\mathit{SO}(21)) \xrightarrow{{i}_*}\pi_{23}(\mathit{SO}(22))\xrightarrow{{p}_*}\pi^{21}_{23}.
	\end{aligned}\]
	As a consequence of          \cite[Table]{Ke}, $\pi_{23}(\mathit{SO}(22))\cong\mathbb Z_{(2)}.$  By Proposition~\ref{m1:R19} and the relation
	$J(\Delta(\eta^2_{21}))=[\iota_{21},\eta^2_{21}]\ne 0$ \cite{Mi1},
	we conclude that ${i_{22}}_*$ is a split epimorphism. Hence, we obtain the following.
	\begin{prop} \label{m1:R22}
		The structure of $\pi_{23}(\mathit{SO}(22))$ is as follows.
		$\Delta(\nu_{21}) = [\sigma^2_9]_{21}$ and
		\[ \FitMath{\pi_{23}(\mathit{SO}(22)) = \mathrm{span}\{[32[\iota_{12}, \iota_{12}]]_{22}\}\cong \mathbb Z_{(2)}.} \]
		\qed
	\end{prop}
	
	\subsection{The Final Unstable Groups and the Stable Groups}
	
	It follows that $\pi_{23}(\mathit{SO}(n))={}\mathbb Z_{(2)}\{[32[\iota_{12}, \iota_{12}]]_n\} (n = 23, n\geq 25),$  and $$\pi_{23}(\mathit{SO}(24))={}\mathbb Z_{(2)}\{[32[\iota_{12}, \iota_{12}]]_{24}\}\oplus\mathbb Z_{(2)}\{\Delta\iota_{24}\}.$$

	\section{Odd-Primary Components}\label{sec:odd23}
	\textbf{Exception to the standing convention.} In this section,
	$\pi_i(\mathit{SO}(n):\mathrm{odd})$ denotes the subgroup of all odd-order
	elements of the integral homotopy group, as defined in Chapter~\ref{ch:notation}.
	In the proof, the symbol  $\pi_{*}(X)$   
	is short for $\pi_{*}(X)\otimes\z_{\mathcal{P}}$ where $\mathcal{P}$
	is the set of odd primes and $\z_{\mathcal{P}}$ denotes the $\mathcal{P}$-localization of $\z.$   A subscript $(p)$ denotes localization at the indicated odd prime
	$p,$   and $\operatorname{Tor}\pi_i(X_{(p)})\xyd \pi_i(X_{(p)})$ is the subgroup of elements whose orders are finite.
	The default convention at $2$ resumes in the next chapter.
	\label{m1:odd}
	For an odd prime $p,$   use the comparison formulas in Section~\ref{sec:odd-method}. The nonzero odd-primary part is as follows.
	\begin{prop}
		For $n\geq3,$   we have
		\[
		\pi_{23}(\mathit{SO}(n):\mathrm{odd})\cong
		\begin{cases}
			\mathbb Z/3, & n=6,14,\\
			0, & \text{otherwise}.
		\end{cases}
		\]
	\end{prop}
	\begin{proof}
		By using the relation \eqref{m1:R-Sp:odd} and the fact that
		$\operatorname{Tor}\pi_{23}(\bigl(\mathit{Sp}(m)\bigr)_{(p)})=0$ from \cite{MT3} and \cite{Bo},
		we obtain $\operatorname{Tor}\pi_{23}(\bigl(\mathit{SO}(2m+1)\bigr)_{(p)})=0.$ 
		The nonzero odd-primary components of
		$\pi_{23}(S^{2m+1})$ are $\operatorname{Tor}\pi_{23}(\bigl(S^5\bigr)_{(3)})\cong \mathbb Z/3$ and $\operatorname{Tor}\pi_{23}(\bigl(S^{13}\bigr)_{(3)}) \cong \mathbb Z/3.$ 
		Thus by using the relation \eqref{m1:R-Sp:even}, we obtain the asserted decomposition.
	\end{proof}
	\chapter{The Homotopy Groups $\pi_{24}(\mathit{SO}(n))$}\label{ch:degree24}
	
	The degree-$24$ calculation uses the degree-$23$ groups as targets of the boundary maps. This
	dependence explains why the previous chapter retained its generators and relations: the abstract group
	structures alone would not identify the new kernels.

	As in the previous chapter, we usually freely appeal to the descriptions of the groups $\pi_{*}(S^{j})$
	set out in Appendix~\ref{tlqdzj} and to the composition relations
	among generators in Toda's range assembled in \^{O}guchi's table
	in Appendix~\ref{app:oguchi-relations}.
	
	\section{Main Theorem and Generators}
	
	\begin{thm} \label{m1:thm2}
		\textup{All spaces in this theorem are localized at $2.$  The displayed groups are obtained from the integral homotopy groups by tensoring with $\mathbb Z_{(2)}.$  Free cyclic summands become $\mathbb Z_{(2)}$; for finite groups, this localization retains precisely the $2$-primary summands. The homotopy groups are given as follows.   }\par
		\begin{enumerate}[leftmargin=0pt,label={},itemindent=0pt]
			\item \textup{The Cases $3\leq n\leq6.$ }\par
			\smallskip\noindent\(
			\FitMath{\pi_{24}(\mathit{SO}(3)) = \mathrm{span}\{[\eta_2]\nu'\eta_6\bar{\mu}_7\}\cong \mathbb Z/2.}
			\)\par
			\smallskip\par\smallskip\noindent\(
			\GroupLayout{\pi_{24}(\mathit{SO}(4)) = \mathrm{span}\{[\eta_2]_4\nu'\eta_6\bar{\mu}_7,\; [\iota_3]\nu'\eta_6\bar{\mu}_7\}\cong (\mathbb Z/2)^{2}.}{\begin{aligned}\pi_{24}(\mathit{SO}(4)) = \mathrm{span}\{&[\eta_2]_4\nu'\eta_6\bar{\mu}_7,\\
					&[\iota_3]\nu'\eta_6\bar{\mu}_7\}\cong (\mathbb Z/2)^{2}.\end{aligned}}
			\)\par
			\smallskip\par\smallskip\noindent\(
			\GroupLayout{\pi_{24}(\mathit{SO}(5)) = \mathrm{span}\{[  \nu^2_4]\kappa_{10},\; [\nu_4\sigma'\eta_{14}]\mu_{15}\}\cong (\mathbb Z/2)^{2},\; \text{where}\;[\nu_4\sigma'\eta_{14}]\mu_{15}=[\nu_4\sigma'\mu_{14}]\eta_{23}.}{\begin{aligned}\pi_{24}(\mathit{SO}(5)) = \mathrm{span}\{&[  \nu^2_4]\kappa_{10},\\
					&[\nu_4\sigma'\eta_{14}]\mu_{15}\}\cong (\mathbb Z/2)^{2},\; \text{where}\;[\nu_4\sigma'\eta_{14}]\mu_{15}=[\nu_4\sigma'\mu_{14}]\eta_{23}.\end{aligned}}
			\)\par
			\smallskip\noindent\(
			\GroupLayout{\pi_{24}(\mathit{SO}(6)) = \mathrm{span}\{[\bar{\zeta}_5],\; [ \nu_5]\mu_8\sigma_{17},\; [  \nu^2_4]_6\kappa_{10},\; [\nu_4\sigma'\eta_{14}]_6\mu_{15}\}\cong \mathbb Z/8\oplus (\mathbb Z/2)^{3}.}{\begin{aligned}\pi_{24}(\mathit{SO}(6)) = \mathrm{span}\{&[\bar{\zeta}_5],\; [ \nu_5]\mu_8\sigma_{17},\; [  \nu^2_4]_6\kappa_{10},\\
					&[\nu_4\sigma'\eta_{14}]_6\mu_{15}\}\cong \mathbb Z/8\oplus (\mathbb Z/2)^{3}.\end{aligned}}
			\)\par
			\item \textup{The Cases $7\leq n\leq9.$ }\par
			\smallskip\noindent\(
			\GroupLayout{\pi_{24}(\mathit{SO}(7)) \;= {} \mathrm{span}\{[\bar{\zeta}_5]_7,\; [ \nu_5]_7\mu_8\sigma_{17},\; [\eta_6]\bar{\mu}_7,\; [\eta_6]\eta_7\mu_8\sigma_{17},\; [P(\Sigma\theta)+\nu_6\kappa_9]\eta_{23},\; [  \nu^2_4]_7\kappa_{10},\; [\nu_4\sigma'\eta_{14}]_7\mu_{15}\}\cong (\mathbb Z/2)^{7}.}{\begin{aligned}\pi_{24}(\mathit{SO}(7)) \;= {} \mathrm{span}\{&[\bar{\zeta}_5]_7,\; [ \nu_5]_7\mu_8\sigma_{17},\; [\eta_6]\bar{\mu}_7,\; [\eta_6]\eta_7\mu_8\sigma_{17},\\
					&[P(\Sigma\theta)+\nu_6\kappa_9]\eta_{23},\; [  \nu^2_4]_7\kappa_{10},\; [\nu_4\sigma'\eta_{14}]_7\mu_{15}\}\cong (\mathbb Z/2)^{7}.\end{aligned}}
			\)\par
			\smallskip\par\smallskip\noindent\(
			\GroupLayout{\pi_{24}(\mathit{SO}(8)) \;= {} \mathrm{span}\{[\bar{\zeta}_5]_8,\; [\eta_6]_8\bar{\mu}_7,\; [\eta_6]_8\eta_7\mu_8\sigma_{17},\; [P(\Sigma\theta)+\nu_6\kappa_9]_8\eta_{23},\; [ \nu_5]_8\mu_8\sigma_{17},\; [  \nu^2_4]_8\kappa_{10},\; [\nu_4\sigma'\eta_{14}]_8\mu_{15},\; [\iota_7]\sigma'\eta_{14}\mu_{15},\; [\iota_7]\nu_7\kappa_{10},\; [\iota_7]\bar{\mu}_7,\; [\iota_7]\eta_7\mu_8\sigma_{17}\}\cong (\mathbb Z/2)^{11},}{\begin{aligned}\pi_{24}(\mathit{SO}(8)) \;= {} \mathrm{span}\{&[\bar{\zeta}_5]_8,\; [\eta_6]_8\bar{\mu}_7,\; [\eta_6]_8\eta_7\mu_8\sigma_{17},\; [P(\Sigma\theta)+\nu_6\kappa_9]_8\eta_{23},\; [ \nu_5]_8\mu_8\sigma_{17},\; [  \nu^2_4]_8\kappa_{10},\\
					&[\nu_4\sigma'\eta_{14}]_8\mu_{15},\; [\iota_7]\sigma'\eta_{14}\mu_{15},\; [\iota_7]\nu_7\kappa_{10},\; [\iota_7]\bar{\mu}_7,\; [\iota_7]\eta_7\mu_8\sigma_{17}\}\cong (\mathbb Z/2)^{11},\end{aligned}}
			\)\par
			\nobreak
			{\raggedleft
				\noindent where $\Delta(\nu_8\kappa_{11})=[  \nu^2_4]_8\kappa_{10}.$\par}
			
			\par\smallskip\noindent\(
			\GroupLayout{\pi_{24}(\mathit{SO}(9)) \;= {} \mathrm{span}\{[\bar{\zeta}_5]_9,\; [P(\Sigma\theta)+\nu_6\kappa_9]_9\eta_{23},\; [ \nu_5]_9\mu_8\sigma_{17},\; [\iota_7]_9\nu_7\kappa_{10},\; [\iota_7]_9\bar{\mu}_7,\; [\iota_7]_9\eta_7\mu_8\sigma_{17}\}\cong (\mathbb Z/2)^{6},}{\begin{aligned}\pi_{24}(\mathit{SO}(9)) \;= {} \mathrm{span}\{&[\bar{\zeta}_5]_9,\; [P(\Sigma\theta)+\nu_6\kappa_9]_9\eta_{23},\; [ \nu_5]_9\mu_8\sigma_{17},\\
					&[\iota_7]_9\nu_7\kappa_{10},\; [\iota_7]_9\bar{\mu}_7,\; [\iota_7]_9\eta_7\mu_8\sigma_{17}\}\cong (\mathbb Z/2)^{6},\end{aligned}}
			\)\par
			\nobreak
			{\raggedleft
				\noindent where $\Delta(\mu_9\sigma_{18})=[ \nu_5]_9\mu_8\sigma_{17}
				+[\iota_7]_9\eta_7\mu_8\sigma_{17}.$ \par}
			\noindent \textup{The Cases $10\leq n\leq14.$ }\par
			\smallskip\noindent\(
			\GroupLayout{\pi_{24}(\mathit{SO}(10)) \;= {} \mathrm{span}\{\Delta(\Sigma\rho'),\; [\sigma^2_9]\eta_{23},\; [\bar{\varepsilon}_9],\; [P(\Sigma\theta)+\nu_6\kappa_9]_{10}\eta_{23},\; [\iota_7]_{10}\bar{\mu}_7\}\cong \mathbb Z/16\oplus \mathbb Z/2\oplus \mathbb Z/4\oplus (\mathbb Z/2)^{2},}{\begin{aligned}\pi_{24}(\mathit{SO}(10)) \;= {} \mathrm{span}\{&\Delta(\Sigma\rho'),\; [\sigma^2_9]\eta_{23},\; [\bar{\varepsilon}_9],\; [P(\Sigma\theta)+\nu_6\kappa_9]_{10}\eta_{23},\\
					&[\iota_7]_{10}\bar{\mu}_7\}\cong \mathbb Z/16\oplus \mathbb Z/2\oplus \mathbb Z/4\oplus (\mathbb Z/2)^{2},\end{aligned}}
			\)\par
			\nobreak
			{\raggedleft
				\noindent where $8\Delta(\Sigma\rho')=[\bar{\zeta}_5]_{10}+[\iota_7]_{10}\bar{\mu}_7$  and
				$2[\bar{\varepsilon}_9]=[\iota_7]_{10}\nu_7\kappa_{10}.$\par}
			
			\par\smallskip\noindent\(
			\GroupLayout{\pi_{24}(\mathit{SO}(11)) = \mathrm{span}\{[\sigma^2_9]_{11}\eta_{23},\; [\bar{\varepsilon}_9]_{11},\; [P(\Sigma\theta)+\nu_6\kappa_9]_{11}\eta_{23},\; [\iota_7]_{11}\bar{\mu}_7\}\cong \mathbb Z/2\oplus \mathbb Z/4\oplus (\mathbb Z/2)^{2}.}{\begin{aligned}\pi_{24}(\mathit{SO}(11)) = \mathrm{span}\{&[\sigma^2_9]_{11}\eta_{23},\; [\bar{\varepsilon}_9]_{11},\; [P(\Sigma\theta)+\nu_6\kappa_9]_{11}\eta_{23},\\
					&[\iota_7]_{11}\bar{\mu}_7\}\cong \mathbb Z/2\oplus \mathbb Z/4\oplus (\mathbb Z/2)^{2}.\end{aligned}}
			\)\par
			\smallskip\par\smallskip\noindent\(
			\GroupLayout{\pi_{24}(\mathit{SO}(12)) \;= {} \mathrm{span}\{[\theta']\eta_{23},\; [\eta_{11}]\theta,\; [\sigma^2_9]_{12}\eta_{23},\; [\bar{\varepsilon}_9]_{12},\; [P(\Sigma\theta)+\nu_6\kappa_9]_{12}\eta_{23},\; [\iota_7]_{12}\bar{\mu}_7\}\cong (\mathbb Z/2)^{6},}{\begin{aligned}\pi_{24}(\mathit{SO}(12)) \;= {} \mathrm{span}\{&[\theta']\eta_{23},\; [\eta_{11}]\theta,\; [\sigma^2_9]_{12}\eta_{23},\; [\bar{\varepsilon}_9]_{12},\\
					&[P(\Sigma\theta)+\nu_6\kappa_9]_{12}\eta_{23},\; [\iota_7]_{12}\bar{\mu}_7\}\cong (\mathbb Z/2)^{6},\end{aligned}}
			\)\par
			\nobreak
			{\raggedleft
				\noindent where $\Delta((\Sigma\theta')\eta_{24})=[P(\Sigma\theta)+\nu_6\kappa_9]_{12}\eta_{23}.$\par}
			
			\par\smallskip\noindent\(
			\GroupLayout{\pi_{24}(\mathit{SO}(13)) = \mathrm{span}\{[\eta_{11}]_{13}\theta,\; [\sigma^2_9]_{13}\eta_{23},\; [\bar{\varepsilon}_9]_{13},\; [\iota_7]_{13}\bar{\mu}_7\}\cong (\mathbb Z/2)^{4},}{\begin{aligned}\pi_{24}(\mathit{SO}(13)) = \mathrm{span}\{&[\eta_{11}]_{13}\theta,\; [\sigma^2_9]_{13}\eta_{23},\\
					&[\bar{\varepsilon}_9]_{13},\; [\iota_7]_{13}\bar{\mu}_7\}\cong (\mathbb Z/2)^{4},\end{aligned}}
			\)\par
			\nobreak
			{\raggedleft
				\noindent where $\Delta(\Sigma\theta)=[\eta_{11}]_{13}\theta.$\par}
			
			\par\smallskip\noindent\(
			\GroupLayout{\pi_{24}(\mathit{SO}(14)) = \mathrm{span}\{[\zeta_{13}],\; [\sigma^2_9]_{14}\eta_{23},\; [\bar{\varepsilon}_9]_{14},\; [\iota_7]_{14}\bar{\mu}_7\}\cong \mathbb Z/8\oplus (\mathbb Z/2)^{3}.}{\begin{aligned}\pi_{24}(\mathit{SO}(14)) = \mathrm{span}\{&[\zeta_{13}],\; [\sigma^2_9]_{14}\eta_{23},\; [\bar{\varepsilon}_9]_{14},\\
					&[\iota_7]_{14}\bar{\mu}_7\}\cong \mathbb Z/8\oplus (\mathbb Z/2)^{3}.\end{aligned}}
			\)\par
			\item \textup{The Cases $n\geq15.$ }\par
			\smallskip\noindent\(
			\GroupLayout{\pi_{24}(\mathit{SO}(15)) = \mathrm{span}\{[\eta_{14}\mu_{15}],\; [\zeta_{13}]_{15},\; [\sigma^2_9]_{15}\eta_{23},\; [\bar{\varepsilon}_9]_{15},\; [\iota_7]_{15}\bar{\mu}_7\}\cong (\mathbb Z/2)^{5}.}{\begin{aligned}\pi_{24}(\mathit{SO}(15)) = \mathrm{span}\{&[\eta_{14}\mu_{15}],\; [\zeta_{13}]_{15},\; [\sigma^2_9]_{15}\eta_{23},\\
					&[\bar{\varepsilon}_9]_{15},\; [\iota_7]_{15}\bar{\mu}_7\}\cong (\mathbb Z/2)^{5}.\end{aligned}}
			\)\par
			\smallskip\par\smallskip\noindent\(
			\GroupLayout{\pi_{24}(\mathit{SO}(16)) \;= {} \mathrm{span}\{[\mu_{15}],\; [\varepsilon_{15}]\eta_{23},\; [\nu_{15}]\nu^2_{18},\; [\eta_{14}\mu_{15}]_{16},\; [\zeta_{13}]_{16},\; [\sigma^2_9]_{16}\eta_{23},\; [\bar{\varepsilon}_9]_{16},\; [\iota_7]_{16}\bar{\mu}_7\}\cong (\mathbb Z/2)^{8},}{\begin{aligned}\pi_{24}(\mathit{SO}(16)) \;= {} \mathrm{span}\{&[\mu_{15}],\; [\varepsilon_{15}]\eta_{23},\; [\nu_{15}]\nu^2_{18},\; [\eta_{14}\mu_{15}]_{16},\; [\zeta_{13}]_{16},\\
					&[\sigma^2_9]_{16}\eta_{23},\; [\bar{\varepsilon}_9]_{16},\; [\iota_7]_{16}\bar{\mu}_7\}\cong (\mathbb Z/2)^{8},\end{aligned}}
			\)\par
			\nobreak
			{\raggedleft
				\noindent where $\Delta(\eta_{16}\varepsilon_{17})=[\zeta_{13}]_{16},$  
				$\Delta(\nu^3_{16})\equiv[\bar{\varepsilon}_9]_{16} \bmod [\iota_7]_{16}\bar{\mu}_7$ and $\Delta(\mu_{16})\equiv[\eta_{14}\mu_{15}]_{16} \bmod
				\Delta(\eta_{16}\varepsilon_{17}), \Delta(\nu^3_{16}).$\par}
			
			\par\smallskip\noindent\(
			\GroupLayout{\pi_{24}(\mathit{SO}(17)) \;= {} \mathrm{span}\{[\mu_{15}]_{17},\; [\eta_{15}]_{17}\varepsilon_{16},\; [\eta_{15}]_{17}\bar{\nu}_{16},\; [\sigma^2_9]_{17}\eta_{23},\; [\iota_7]_{17}\bar{\mu}_7\}\cong (\mathbb Z/2)^{5}.}{\begin{aligned}\pi_{24}(\mathit{SO}(17)) \;= {} \mathrm{span}\{&[\mu_{15}]_{17},\; [\eta_{15}]_{17}\varepsilon_{16},\; [\eta_{15}]_{17}\bar{\nu}_{16},\\
					&[\sigma^2_9]_{17}\eta_{23},\; [\iota_7]_{17}\bar{\mu}_7\}\cong (\mathbb Z/2)^{5}.\end{aligned}}
			\)\par
			\smallskip\par\smallskip\noindent\(
			\GroupLayout{\pi_{24}(\mathit{SO}(18)) = \mathrm{span}\{\Delta(\sigma_{18}),\; [\nu_{15}]_{18}\nu^2_{18},\; [\iota_7]_{18}\bar{\mu}_7\}\cong \mathbb Z/16\oplus (\mathbb Z/2)^{2},}{\begin{aligned}\pi_{24}(\mathit{SO}(18)) = \mathrm{span}\{&\Delta(\sigma_{18}),\; [\nu_{15}]_{18}\nu^2_{18},\\
					&[\iota_7]_{18}\bar{\mu}_7\}\cong \mathbb Z/16\oplus (\mathbb Z/2)^{2},\end{aligned}}
			\)\par
			{\raggedleft where $8\Delta(\sigma_{18})=[\mu_{15}]_{18}.$ \par}

			\par\smallskip\noindent\(
			\GroupLayout{\pi_{24}(\mathit{SO}(19)) = \mathrm{span}\{[\nu_{15}]_{19}\nu^2_{18},\; [\iota_7]_{19}\bar{\mu}_7\}\cong (\mathbb Z/2)^{2}, \text{where}\; \Delta(\nu^2_{19})=[\nu_{15}]_{19}\nu^2_{18}.}{\begin{aligned}\pi_{24}(\mathit{SO}(19)) = \mathrm{span}\{&[\nu_{15}]_{19}\nu^2_{18},\\
					&[\iota_7]_{19}\bar{\mu}_7\}\cong (\mathbb Z/2)^{2}, \text{where}\; \Delta(\nu^2_{19})=[\nu_{15}]_{19}\nu^2_{18}.\end{aligned}}
			\)\par
			\smallskip\par\smallskip\noindent\(
			\GroupLayout{\pi_{24}(\mathit{SO}(22)) = \mathrm{span}\{\Delta(\nu_{22}),\; [\iota_7]_{22}\bar{\mu}_7\}\cong \mathbb Z/4\oplus \mathbb Z/2.}{\begin{aligned}\pi_{24}(\mathit{SO}(22)) = \mathrm{span}\{&\Delta(\nu_{22}),\\
					&[\iota_7]_{22}\bar{\mu}_7\}\cong \mathbb Z/4\oplus \mathbb Z/2.\end{aligned}}
			\)\par
			\smallskip\par\smallskip\noindent\(
			\GroupLayout{\pi_{24}(\mathit{SO}(23)) = \mathrm{span}\{[\eta^2_{22}],\; [\iota_7]_{23}\bar{\mu}_7\}\cong (\mathbb Z/2)^{2}.}{\begin{aligned}\pi_{24}(\mathit{SO}(23)) = \mathrm{span}\{&[\eta^2_{22}],\\
					&[\iota_7]_{23}\bar{\mu}_7\}\cong (\mathbb Z/2)^{2}.\end{aligned}}
			\)\par
			\smallskip\par\smallskip\noindent\(
			\GroupLayout{\pi_{24}(\mathit{SO}(24)) = \mathrm{span}\{[\eta_{23}],\; [\eta^2_{22}]_{24},\; [\iota_7]_{24}\bar{\mu}_7\}\cong (\mathbb Z/2)^{3}, \text{where}\;\Delta(\eta_{24})=[\eta^2_{22}]_{24}.}{\begin{aligned}\pi_{24}(\mathit{SO}(24)) = \mathrm{span}\{&[\eta_{23}],\; [\eta^2_{22}]_{24},\\
					&[\iota_7]_{24}\bar{\mu}_7\}\cong (\mathbb Z/2)^{3}, \text{where}\;\Delta(\eta_{24})=[\eta^2_{22}]_{24}.\end{aligned}}
			\)\par
			\smallskip\par\smallskip\noindent\(
			\GroupLayout{\pi_{24}(\mathit{SO}(25)) = \mathrm{span}\{[\eta_{23}]_{25},\; [\iota_7]_{25}\bar{\mu}_7\}\cong (\mathbb Z/2)^{2},\text{where}\;\Delta(\iota_{25})=[\eta_{23}]_{25}.}{\begin{aligned}\pi_{24}(\mathit{SO}(25)) = \mathrm{span}\{&[\eta_{23}]_{25},\\
					&[\iota_7]_{25}\bar{\mu}_7\}\cong (\mathbb Z/2)^{2},\text{where}\;\Delta(\iota_{25})=[\eta_{23}]_{25}.\end{aligned}}
			\)\par
			\smallskip\par\smallskip\noindent\(
			\FitMath{\pi_{24}(\mathit{SO}(20+n)) = \mathrm{span}\{[\iota_7]_{20+n}\bar{\mu}_7\}\cong \mathbb Z/2,\; (n=0, 1\; \text{or}\; n\geq 6).}
			\)\par
		\end{enumerate}
	\end{thm}
	
	\section{The Cases $3\leq n\leq6$}
	\subsection{The Cases $3\leq n\leq5$}
	
	As in degree $23,$   the initial calculation in the homotopy groups of spheres gives the results for $n=3,4.$  We then
	determine what the inclusion kills for $n=5.$  This makes the kernel of the connecting homomorphism the
	decisive part of that calculation.
	
	We see that $\pi^3_{24}
	= \frac{\mathbb Z}{2}\{\nu'\eta_6\bar{\mu}_7\}.$  This yields the groups
	$\pi_{24}(\mathit{SO}(3))$ and $\pi_{24}(\mathit{SO}(4))$:
	\[\begin{aligned}
		&\pi_{24}(\mathit{SO}(3))={}\mathrm{span}\{[\eta_2]\nu'\eta_6\bar{\mu}_7\}\cong\mathbb Z/2\uprightsemicolon{}\ \pi_{24}(\mathit{SO}(4))=\mathrm{span}\{[\eta_2]_4\nu'\eta_6\bar{\mu}_7,[\iota_3]\nu'\eta_6\bar{\mu}_7\}\cong(\mathbb Z/2)^2.
	\end{aligned}\]
	\noindent
	Applying Lemma~\ref{m1:236} for $\pi_{24}^3,$   we deduce
	\begin{equation}
		\begin{aligned}
			&i_{5 \ast}(\pi_{24}(\mathit{SO}(4)))=0.\label{m1:i5R424}
		\end{aligned}
	\end{equation}

	We proceed to the next case, where the relevant
	calculations can be assembled into a complete description.
	\begin{prop} The group $\pi_{24}(\mathit{SO}(5))$ is given as follows:
		\begin{equation}\label{m1:524}
			\pi_{24}(\mathit{SO}(5))={}\frac{\mathbb Z}{2}\{[  \nu^2_4]\kappa_{10}\} \oplus\frac{\mathbb Z}{2}\{[\nu_4\sigma'\eta_{14}]\mu_{15}\}.
		\end{equation}
	\end{prop}

	\begin{proof}
		In the exact sequence $(24)_4$:
		\[\begin{aligned}
			&\pi^4_{25}\xrightarrow{\Delta}\pi_{24}(\mathit{SO}(4))\xrightarrow{i_*}\pi_{24}(\mathit{SO}(5)) \xrightarrow{p_*}\pi^4_{24}\xrightarrow{\Delta}\pi_{23}(\mathit{SO}(4)),
		\end{aligned}\]
		we point out that $\operatorname{ker}(\Delta\colon \pi^4_{24}\to \pi_{23}(\mathit{SO}(4)))= \frac{\mathbb Z}{2}\{\nu^2_4\kappa_{10},
		\nu_4\sigma'\eta_{14}\mu_{15}\}\cong(\mathbb Z/2)^2$ and
		\[\begin{aligned}
			&\pi^4_{25}={}\frac{\mathbb Z}{2}\{(\Sigma\nu')\eta_7\bar{\mu}_8\} \oplus\frac{\mathbb Z}{2}\{\nu_4\eta_7\bar{\mu}_8\}\oplus\frac{\mathbb Z}{8}\{\nu_4\zeta_7\sigma_{18}\}.
		\end{aligned}\]
		Using $\nu'\zeta_6=0,$   we find
		\[\begin{aligned}
			&\Delta((\Sigma\nu')\eta_7\bar{\mu}_8) = [\eta_2]_4\nu'\eta_6\bar{\mu}_7, \Delta(\nu_4\eta_7\bar{\mu}_8)\equiv [\iota_3]\nu'\eta_6\bar{\mu}_7 \bmod [\eta_2]_4\nu'\eta_6\bar{\mu}_7
		\end{aligned}\]
		and $\Delta(\nu_4\zeta_7\sigma_{18})\equiv[\iota_3]\nu'\zeta_6\sigma_{17}=0\ \bmod\
		[\eta_2]_4\nu'\zeta_6\sigma_{17}=0.$  This completes the proof.
	\end{proof}
	
	\subsection{The Case $n=6$ and the Lift of $\bar\zeta_5$}
	
	The new issue for $n=6$ is the order of a lift of $\bar\zeta_5.$  Its Toda-bracket representative
	allows the eighth multiple to be calculated, while its projection supplies the comparison with the
	order of the element of a homotopy group of a sphere.
	
	We recall \eqref{m1:bz5}:
	$\bar{\zeta}_5\in\{\zeta_5, 8\iota_{16}, 2\sigma_{16}\}_1.$ 
	Since the order of $[\zeta_5]\in \pi_{16}(\mathit{SO}(6))$ is $8$ \cite[Proposition 4.1]{Ka}, we define a lift of $\bar{\zeta}_5$ as
	\begin{equation}
		\begin{aligned}
			&[\bar{\zeta}_5]\in\{[\zeta_5], 8\iota_{16}, 2\sigma_{16}\}_1\subseteq \pi_{24}(\mathit{SO}(6)).\label{m1:[bzeta_5]}
		\end{aligned}
	\end{equation}
	
	We show the following.
	\begin{prop} \label{m1:R61}
		The structure of $\pi_{24}(\mathit{SO}(6))$ is as follows.
		\[ \GroupLayout{\pi_{24}(\mathit{SO}(6)) = \mathrm{span}\{[\bar{\zeta}_5],\; [ \nu_5]\mu_8\sigma_{17},\; [\nu_4\sigma'\eta_{14}]_6\mu_{15},\; [  \nu^2_4]_6\kappa_{10}\}\cong \mathbb Z/8\oplus (\mathbb Z/2)^{3}.}{\begin{aligned}\pi_{24}(\mathit{SO}(6)) = \mathrm{span}\{&[\bar{\zeta}_5],\; [ \nu_5]\mu_8\sigma_{17},\; [\nu_4\sigma'\eta_{14}]_6\mu_{15},\\
				&[  \nu^2_4]_6\kappa_{10}\}\cong \mathbb Z/8\oplus (\mathbb Z/2)^{3}.\end{aligned}} \]
	\end{prop}
	\begin{proof}
		We consider the exact sequence $(24)_5$:
		\[\begin{aligned}
			&\pi^5_{25}\xrightarrow{\Delta}\pi_{24}(\mathit{SO}(5))\xrightarrow{i_*}\pi_{24}(\mathit{SO}(6)) \xrightarrow{p_*}\pi^5_{24}\xrightarrow{\Delta}\pi_{23}(\mathit{SO}(5)).
		\end{aligned}\]
		Since $\pi^{5}_{24}=\frac{\mathbb Z}{8}\{\bar{\zeta}_5\}\oplus\frac{\mathbb Z}{2}\{\nu_5\mu_8\sigma_{17}\},$   we conclude that $\Delta(\pi^5_{24})=0.$ 
		Recall that $$\pi^5_{25} = \mathrm{span}\{\nu^2_5\kappa_{11}, \nu_5\bar{\mu}_8,
		\nu_5\eta_8\mu_9\sigma_{18}\}$$\noindent
		and $\Delta\nu_5=0.$ 
		In light of \eqref{m1:[bzeta_5]} and \eqref{m1:2i8s2i}, we see that
		\[\begin{aligned}
			&8[\bar{\zeta}_5]\in\{[\zeta_5], 8\iota_{16}, 2\sigma_{16}\}_1\circ 8\iota_{24} = -[\zeta_5]\circ \Sigma\{8\iota_{15}, 2\sigma_{15}, 8\iota_{22}\}=0.
		\end{aligned}\]
		This implies $\mathrm{ord}[\bar{\zeta}_5]=8$ and completes the proof.
	\end{proof}
	
	\section{Boundary Relations for the Cases $7\leq n\leq9$}

	\subsection{Connecting Images and Relations Inherited from Degree 25}
	
	We use the composition relations
	\begin{equation}
		\begin{aligned}
			&\varepsilon_3\zeta_{11}=0,\ \bar{\nu}_6\zeta_{14}=\pm 8\bar{\sigma}_6.\label{m1:Ogep}
		\end{aligned}
	\end{equation}
	Therefore, we have
	\[\begin{aligned}
		&\pm8\bar{\sigma}_6=\bar{\nu}_6\zeta_{14} =(\bar{\nu}_6+\varepsilon_6)\zeta_{14}.
	\end{aligned}\]
	
	Hence, by
	\eqref{m1:iGU25}, we obtain
	\begin{equation}
		\begin{aligned}
			&\langle\bar{\nu}_6+\varepsilon_6\rangle\zeta_{14}=\pm 8\langle\bar{\sigma}_6\rangle.\label{m1:bnu+epzt}
		\end{aligned}
	\end{equation}
	\noindent   We need the following group from \cite[Table 2]{Ka}:
	\[\begin{aligned}
		&\pi_{7}(\mathit{SO}(5))={}\mathbb Z_{(2)}\{[4 \nu_4]\},\ \pi_{7}(\mathit{SO}(6))=\mathbb Z_{(2)}\{[\eta^2_5]\}.
	\end{aligned}\]
	\noindent We record the relations $J[4 \nu_4]=\sigma'''$ and $J[\eta^2_5]=\pm
	\sigma''$ \cite[Lemma 2.1(i)]{HKM}.     Hence, we obtain $J(2[\eta^2_5]-[4 \nu_4]_6)=0$ and
	$2[\eta^2_5]\equiv[4 \nu_4]_6 \bmod \operatorname{ker} J=\mathrm{span}\{4[\eta^2_5]\}.$  That is, 
	\begin{equation}
		\begin{aligned}
			& 2[\eta^2_5]\ty [4 \nu_4]_6\m 4[\eta^2_5].\label{m1:2[et25}
		\end{aligned}
	\end{equation}
	\noindent   Since
	\ $[  \nu^2_4]_6\sigma^2_{10}\in i_*\pi_{24}(\mathit{SO}(5))$ and $\pi_{24}(\mathit{SO}(5))\cong(\mathbb Z/2)^2,$   we
	deduce
	\begin{equation}
		\begin{aligned}
			&2[  \nu^2_4]_6\sigma^2_{10} = 0.\label{m1:2n24sg2}
		\end{aligned}
	\end{equation}
	
	\section{The Cases $7\leq n\leq9$}
	\subsection{The Case $n=7$}
	The various relations established thus far
	may now be brought together to yield
	the following explicit determination.
	\begin{prop} \label{m1:724}
		The structure of $\pi_{24}(\mathit{SO}(7))$ is as follows.
		\[ \GroupLayout{\pi_{24}(\mathit{SO}(7)) = \mathrm{span}\{[\bar{\zeta}_5]_7,\; [ \nu_5]_7\mu_8\sigma_{17},\; [\nu_4\sigma'\eta_{14}]_7\mu_{15},\; [  \nu^2_4]_7\kappa_{10},\; [P(\Sigma\theta)+\nu_6\kappa_9]\eta_{23},\; [\eta_6]\bar{\mu}_7,\; [\eta_6]\eta_7\mu_8\sigma_{17}\} \;\cong\; (\mathbb Z/2)^7.}{\begin{aligned}\pi_{24}(\mathit{SO}(7)) = \mathrm{span}\{&[\bar{\zeta}_5]_7,\; [ \nu_5]_7\mu_8\sigma_{17},\; [\nu_4\sigma'\eta_{14}]_7\mu_{15},\; [  \nu^2_4]_7\kappa_{10},\\
				&[P(\Sigma\theta)+\nu_6\kappa_9]\eta_{23},\; [\eta_6]\bar{\mu}_7,\; [\eta_6]\eta_7\mu_8\sigma_{17}\} \;\cong\; (\mathbb Z/2)^7.\end{aligned}} \]
		Here $\Delta\bar{\zeta}_6 = 2[\bar{\zeta}_5].$ 
	\end{prop}
	\begin{proof}
		Examine the exact sequence $(24)_6$:
		\[\begin{aligned}
			&\pi^6_{25}\xrightarrow{\Delta}\pi_{24}(\mathit{SO}(6))\xrightarrow{i_*}\pi_{24}(\mathit{SO}(7)) \xrightarrow{p_*}\pi^6_{24}\xrightarrow{\Delta}\pi_{23}(\mathit{SO}(6)).
		\end{aligned}\]
		We know that $\operatorname{ker}(\Delta\colon \pi^6_{24}\to \pi_{23}(\mathit{SO}(6)))= \mathrm{span}\{(P(\Sigma\theta)+\nu_6\kappa_9)\eta_{23},
		\eta_6\bar{\mu}_7, \eta^2_6\mu_8\sigma_{17}\}$ and
		$\pi^6_{25} = \frac{\mathbb Z}{8}\{\bar{\zeta}_6\}\oplus\frac{\mathbb Z}{32}\{\bar{\sigma}_6\}.$ 
		By using  Lemma~\ref{m1:zeta1}, we find
		\[\begin{aligned}
			&\Delta\bar{\zeta}_6\in\Delta\{\zeta_6, 8\iota_{17}, 2\sigma_{17}\}_2 \subseteq\{2[\zeta_5], 8\iota_{16}, 2\sigma_{16}\}_1.
		\end{aligned}\]
		Since $\mathit{SO}(6)$ is an H-group, we conclude that
		\[\begin{aligned}
			&\{2[\zeta_5], 8\iota_{16}, 2\sigma_{16}\}_1 \supseteq 2\,\mathrm{id}\circ\{[\zeta_5], 8\iota_{16}, 2\sigma_{16}\}_1 =\{[\zeta_5], 8\iota_{16}, 2\sigma_{16}\}_1 \circ 2\iota_{24}.
		\end{aligned}\]
		Moreover, by using  the relation \eqref{m1:[bzeta_5]}, we see that
		\[\begin{aligned}
			&\Delta\bar{\zeta}_6\in \{2[\zeta_5], 8\iota_{16}, 2\sigma_{16}\}_1 \ni 2[\bar{\zeta}_5] \bmod 2[\zeta_5]\circ\pi^{16}_{24} + \pi_{17}(\mathit{SO}(6))\circ 2\sigma_{17}.
		\end{aligned}\]
		As $\pi^{16}_{24}\cong(\mathbb Z/2)^2$ and $\pi_{17}(\mathit{SO}(6))
		= \frac{\mathbb Z}{8}\{[  \nu^2_4]_6\sigma_{10}\}\oplus(\mathbb Z/2)^3$ \cite{KM1},
		the indeterminacy is equal to $\mathrm{span}\{2[  \nu^2_4]_6\sigma^2_{10}\}=0$
		\eqref{m1:2n24sg2}.
		Hence the indeterminacy is trivial.
		This implies the relation
		\[\begin{aligned}
			&\Delta\bar{\zeta}_6 = 2[\bar{\zeta}_5].
		\end{aligned}\]
		\noindent 
		By Lemma~\ref{m3:fh} and the relation $(r_3)_*\Delta'=\Delta$,
		we obtain $\Delta(\bar{\sigma}_6)=0$ in $\pi_{24}(\mathit{SO}(6)).$
		Thus, by Proposition~\ref{m1:R61}, we obtain
		the group structure of $\pi_{24}(\mathit{SO}(7)).$ 
	\end{proof}

	We need the following lemma.
	\begin{lem}
		$2[P(\Sigma\theta)+\nu_6\kappa_9] \equiv[ \nu_5]_7\bar{\varepsilon}_8\ \bmod\ [\nu_4\eta_7\bar{\varepsilon}_8]_7, [\nu_4(\Sigma\zeta')]_7.$
	\end{lem}
	\begin{proof} Recall the relation
		$$2[P(\Sigma\theta)+\nu_6\kappa_9]
		\equiv[ \nu_5]_7\bar{\varepsilon}_8\ \bmod\
		[\nu_4\sigma'\mu_{14}]_7, [\nu_4(\Sigma\zeta')]_7, [\nu_4\eta_7\bar{\varepsilon}_8]_7.$$

		In light of \eqref{m1:epep}, we deduce $$[ \nu_5]_7\bar{\varepsilon}_8\eta_{23}=
		[ \nu_5]_7\nu_8\sigma_{11}\nu^2_{18}=[ \nu_5]_7(\Sigma\sigma')\nu^3_{15}
		=(2[ \nu_5]_7\sigma_8)\nu^3_{15}=0.$$    Since $[\nu_4\eta_7\bar{\varepsilon}_8]\eta_{23}\in{i_5}_*\pi_{24}(\mathit{SO}(4))=0$
		\eqref{m1:i5R424}, we find $[\nu_4\eta_7\bar{\varepsilon}_8]_7\eta_{23}=0.$  Since
		$\eta^2_{15}\varepsilon_{17}=4\nu_{15}\sigma_{18}=0,$   we have
		$[\nu_4(\Sigma\zeta')]\eta_{23}\in{i_5}_*\pi_{24}(\mathit{SO}(4))=0$ and
		$[\nu_4(\Sigma\zeta')]_7\eta_{23}=0.$ 
		
		Finally, we infer that $$[\nu_4\sigma'\mu_{14}]\eta_{23}\equiv[\nu_4\sigma'\eta_{14}]\mu_{15}\ \bmod\ {i_5}_*\pi_{24}(\mathit{SO}(4))=0.$$ 
		This and the group structure of $\pi_{24}(\mathit{SO}(7))$ imply $[\nu_4\sigma'\mu_{14}]\eta_{23}=[\nu_4\sigma'\eta_{14}]\mu_{15}$ and
		$[\nu_4\sigma'\mu_{14}]_7\eta_{23}=[\nu_4\sigma'\eta_{14}]_7\mu_{15}.$ 
		Hence, we obtain $2[P(\Sigma\theta)+\nu_6\kappa_9] \equiv[ \nu_5]_7\bar{\varepsilon}_8\ \bmod\ [\nu_4\eta_7\bar{\varepsilon}_8]_7, [\nu_4(\Sigma\zeta')]_7.
		$\end{proof}

	\subsection{The Case $n=8$}
	Since $\pi_{24}(\mathit{SO}(8))\cong \pi_{24}(\mathit{SO}(7))\oplus\pi^7_{24}$ and $\pi^7_{24}
	= \mathrm{span}\{\sigma'\eta_{14}\mu_{15}, \nu_7\kappa_{10}, \bar{\mu}_7,
	\eta_7\mu_8\sigma_{17}\}\cong(\mathbb Z/2)^4,$   we have the following.
	
	\begin{prop} \label{m1:8}The structure of $\pi_{24}(\mathit{SO}(8))$ is as follows.
		\[\FitMath{\begin{aligned}
				\pi_{24}(\mathit{SO}(8)) \;= {}& \mathrm{span}\{[\bar{\zeta}_5]_8,\; [ \nu_5]_8\mu_8\sigma_{17},\; [\nu_4\sigma'\eta_{14}]_8\mu_{15},\; [\nu_4^2]_8\kappa_{10},\; [\eta_6]_8\bar{\mu}_7,\; [\eta_6]_8\eta_7\mu_8\sigma_{17}\} \\
				&\mathord{\oplus}\,\mathrm{span}\{[P(\Sigma\theta)+\nu_6\kappa_9]_8\eta_{23},\; [\iota_7]\sigma'\eta_{14}\mu_{15},\; [\iota_7]\nu_7\kappa_{10},\; [\iota_7]\bar{\mu}_7,\; [\iota_7]\eta_7\mu_8\sigma_{17}\}\cong (\mathbb Z/2)^{11}.
		\end{aligned}}\]
		\qed
	\end{prop}
	
	\subsection{The Case $n=9$}
	We now establish the following result.
	\begin{prop} \label{m1:924}
		The structure of $\pi_{24}(\mathit{SO}(9))$ is as follows.
		\[\GroupLayout{\pi_{24}(\mathit{SO}(9)) \;= {} \mathrm{span}\{[\bar{\zeta}_5]_9,\; [ \nu_5]_9\mu_8\sigma_{17},\; [P(\Sigma\theta) +\nu_6\kappa_9]_9\eta_{23},\; [\iota_7]_9\nu_7\kappa_{10},\; [\iota_7]_9\bar{\mu}_7,\; [\iota_7]_9\eta_7\mu_8\sigma_{17}\}\cong (\mathbb Z/2)^{6},}{\begin{aligned}\pi_{24}(\mathit{SO}(9)) \;= {} \mathrm{span}\{&[\bar{\zeta}_5]_9,\; [ \nu_5]_9\mu_8\sigma_{17},\; [P(\Sigma\theta) +\nu_6\kappa_9]_9\eta_{23},\\
				&[\iota_7]_9\nu_7\kappa_{10},\; [\iota_7]_9\bar{\mu}_7,\; [\iota_7]_9\eta_7\mu_8\sigma_{17}\}\cong (\mathbb Z/2)^{6},\end{aligned}}\]
		with relations
		$[\nu_4\sigma'\eta_{14}]_9\mu_{15}
		=[  \nu^2_4]_9\kappa_{10}
		=0$ and
		$$[\mathrm{c}_{\mathrm{KM}}(\bar{\nu}_6+\varepsilon_6)]_9\eta_{14}\mu_{15}
		=[\iota_7]_9\sigma'\eta_{14}\mu_{15}=[ \nu_5]_9\mu_8\sigma_{17}.$$
	\end{prop}
	\begin{proof}
		In the exact sequence $(24)_8$:
		\[\begin{aligned}
			&\pi^8_{25}\xrightarrow{\Delta}\pi_{24}(\mathit{SO}(8))\xrightarrow{i_*}\pi_{24}(\mathit{SO}(9))\xrightarrow{p_*}\pi^8_{24}\stackrel{\Delta}\hookrightarrow \pi_{23}(\mathit{SO}(8)),
		\end{aligned}\]
		note that $\Delta\colon  \pi^8_{24}\to \pi_{23}(\mathit{SO}(8))$ is a monomorphism and that
		\[\begin{aligned}
			&\pi^8_{25} \!={}\! \mathrm{span}\{\sigma_8\eta_{15}\mu_{16}, (\Sigma\sigma')\eta_{15}\mu_{16}, \nu_8\kappa_{11}, \bar{\mu}_8, \eta_8\mu_9\sigma_{18}\}\!\cong\!(\mathbb Z/2)^5.
		\end{aligned}\]\noindent    By \eqref{m1:Dsg82}, with its odd integers $x_1,x_2\in\z,$   we obtain the following. Since $\eta_{14}\mu_{15}$ has order at most $2,$   the coefficients reduce to $1$ after composition:
		\[\begin{aligned}
			\Delta(\sigma_8\eta_{15}\mu_{16}) &\equiv x_1[\iota_7]\sigma'\eta_{14}\mu_{15} + x_2[\mathrm{c}_{\mathrm{KM}}(\bar{\nu}_6+\varepsilon_6)]_8\eta_{14}\mu_{15} \ \bmod\ ([ \nu_5]_8\nu^2_8 + [\eta_5\varepsilon_6]_8)\circ\eta_{14}\mu_{15} \\
			&\equiv [\iota_7]\sigma'\eta_{14}\mu_{15} + [\mathrm{c}_{\mathrm{KM}}(\bar{\nu}_6+\varepsilon_6)]_8\eta_{14}\mu_{15}\ \bmod\ [\eta_5\varepsilon_6]_8\eta_{14}\mu_{15}.
		\end{aligned}\]
		In view of Lemma~\ref{m1:bnuepmu} and the relation
		$\zeta'\eta_{22}=0,$   we deduce
		\[\begin{aligned}
			[\mathrm{c}_{\mathrm{KM}}(\bar{\nu}_6+\varepsilon_6)]_8\eta_{14}\mu_{15} &\equiv [\eta_6]_8(\Sigma\zeta' + \mu_7\sigma_{16})\eta_{23} + [\zeta_5]_8\sigma_{16}\eta_{23} \\
			&=[\eta_6]_8\eta_7\mu_8\sigma_{17} + [ \nu_5]_8\mu_8\sigma_{17} \\
			&\bmod [ \nu_5]_8\bar\varepsilon_8\eta_{23}, [\nu_4\sigma'\mu_{14}]_8\eta_{23}, [\nu_4(\Sigma\zeta')]_8\eta_{23}, [\nu_4\eta_7\bar{\varepsilon}_8]_8\eta_{23}.
		\end{aligned}\]
		Therefore, by Lemma~\ref{m1:89}, we find
		\[\begin{aligned}
			&[\mathrm{c}_{\mathrm{KM}}(\bar{\nu}_6+\varepsilon_6)]_8\eta_{14}\mu_{15} \equiv [\eta_6]_8\eta_7\mu_8\sigma_{17} + [ \nu_5]_8\mu_8\sigma_{17} \bmod [  \nu^2_4]_8\kappa_{10}, [\nu_4\sigma'\eta_{14}]_8\mu_{15}.
		\end{aligned}\]
		As a consequence of          \cite[Lemmas 1.1(v), 1.2(iii)]{KM1} and Lemma~\ref{m1:89}, we know 
		\[\begin{aligned}
			&[\eta_6]_8\sigma'\eta_{14}\mu_{15} = [\eta_5\varepsilon_6]_8\eta_{14}\mu_{15} = [\nu_4\sigma'\eta_{14}]_8\mu_{15}.
		\end{aligned}\]
		Hence, we conclude that
		\begin{equation}
			\begin{aligned}
				&\Delta(\sigma_8\eta_{15}\mu_{16}) \equiv [\iota_7]\sigma'\eta_{14}\mu_{15} + [\eta_6]_8\eta_7\mu_8\sigma_{17} + [ \nu_5]_8\mu_8\sigma_{17}\ \bmod [  \nu^2_4]_8\kappa_{10}, [\nu_4\sigma'\eta_{14}]_8\mu_{15}.\label{m1:Dsg8etm}
			\end{aligned}
		\end{equation}
		We note the relation $\Delta\iota_8 = 2[\iota_7] - [\eta_6]_8.$ 
		Since $[\eta_6]\nu_7 = x[  \nu^2_4]_7$ for some odd $x$ \cite[Lemma 1.1(iv)]{KM1},
		we obtain
		$\Delta(\nu_8\kappa_{11}) = [\eta_6]_8\nu_7\kappa_{10}=[  \nu^2_4]_8\kappa_{10}.$ 
		\noindent   Note that
		\begin{align*}
			&\Delta((\Sigma\sigma')\eta_{15}\mu_{16}) = [\eta_6]_8\sigma'\eta_{14}\mu_{15}= [\nu_4\sigma'\eta_{14}]_8\mu_{15}\uprightsemicolon{} \\
			&\Delta(\nu_8\kappa_{11}) = [  \nu^2_4]_8\kappa_{10}\uprightsemicolon{}\ \Delta\bar{\mu}_8 = [\eta_6]_8\bar{\mu}_7\uprightsemicolon{}\ \Delta(\eta_8\mu_9\sigma_{18}) = [\eta_6]_8\eta_7\mu_8\sigma_{17}.
		\end{align*}
		Thus we obtain the group $\pi_{24}(\mathit{SO}(9))$ and
		the equation
		\[\begin{aligned}
			&[\nu_4\sigma'\eta_{14}]_9\mu_{15}=[  \nu^2_4]_9\kappa_{10} =[\eta_6]_9\bar{\mu}_7 = [\eta_6]_9\eta_7\mu_8\sigma_{17}=0.
		\end{aligned}\]
		Furthermore, we obtain
		$[\mathrm{c}_{\mathrm{KM}}(\bar{\nu}_6+\varepsilon_6)]_9\eta_{14}\mu_{15}=[ \nu_5]_9\mu_8\sigma_{17}$
		by the relation
		\[\begin{aligned}
			&[\mathrm{c}_{\mathrm{KM}}(\bar{\nu}_6+\varepsilon_6)]_9\eta_{14}\mu_{15} \equiv [\eta_6]_9\eta_7\mu_8\sigma_{17} + [ \nu_5]_9\mu_8\sigma_{17} \bmod [  \nu^2_4]_9\kappa_{10}, [\nu_4\sigma'\eta_{14}]_9\mu_{15}.
		\end{aligned}\]\noindent    As a consequence of           \eqref{m1:Dsg8etm}, we  obtain
		$[\iota_7]_9\sigma'\eta_{14}\mu_{15}=[ \nu_5]_9\mu_8\sigma_{17}.$ 
		This completes the proof.
	\end{proof}
	
	\section{The Cases $10\leq n\leq14$}
	\subsection{The Cases $n=10,11$}
	\hspace{2ex}
	First we note that $\sigma_9(\bar{\nu}_{16} + \varepsilon_{16})
	= \sigma^2_9\eta_{23}.$  In the exact sequence $(24)_{9}$:
	\[\begin{aligned}
		&\pi^9_{25}\xrightarrow{\Delta}\pi_{24}(\mathit{SO}(9))\xrightarrow{i_*}\pi_{24}(\mathit{SO}(10)) \xrightarrow{p_*}\pi^9_{24}\xrightarrow{\Delta}\pi_{23}(\mathit{SO}(9)),
	\end{aligned}\]
	we see that $\operatorname{ker}\ \Delta = \mathrm{span}\{2\rho', \sigma^2_9\eta_{23}, \bar{\varepsilon}_9\}$
	from the proof of Proposition~\ref{m1:R10}.
	So there exist elements
	$\Delta(\Sigma\rho'),$   $[\sigma^2_9]\eta_{23},$  
	$[\bar{\varepsilon}_9] \in \pi_{24}(\mathit{SO}(10)).$ 
	
	\begin{lem} There is a relation\vspace{-0.8em}
		\begin{equation}
			\begin{aligned} 
				&\bar{\varepsilon}_5=\{2\nu_5,\nu_8, \bar\nu_{11}\}_1.\label{m1:tbep5}
			\end{aligned}
		\end{equation}
	\end{lem}
	\begin{proof}
		The Toda bracket $\{2\nu_5,\nu_8,\bar{\nu}_{11}\}_1$ is
		defined by $\nu_7\bar{\nu}_{10} = 0.$ 
		Since $\pi^7_{19}=0,$   we see that
		$\{2\nu_5,\nu_8,\bar{\nu}_{11}\}_1=\{2\nu_5,\nu_8,\bar{\nu}_{11}\}.$ 
		In light of \cite[Theorem 2(2)]{MaO} and \cite[(3.9)-i)]{T}, we obtain
		$\langle 2\nu, \nu,\bar{\nu} \rangle
		= \langle \bar{\nu},\nu,2\nu\rangle=\bar{\varepsilon}.$ 
		Therefore, since $\Sigma^\infty:\pi^5_{20}\to \pi^S_{15}$ is a monomorphism
		\cite[Theorems 10.5 and 10.10]{T}, we obtain the assertion.
	\end{proof}
	We next consider the following statement, beginning with its formulation and proceeding to its proof before taking up the next topic.
	\begin{lem}\label{m1:so10-image}
		The following relations hold: $[P(\Sigma\theta)+\nu_6\kappa_9]_{10}\eta_{23}=[\bar{\nu}_9\nu^2_{17}]\eta_{23};$ $   [\iota_7]_{10}\eta_7\mu_8\sigma_{17}=[ \nu_5]_{10}\mu_8\sigma_{17}=0.$ Moreover,
		\begin{equation*}
			{i_{10}}_*\pi_{24}(\mathit{SO}(9))=\mathrm{span}\{[\bar{\zeta}_5]_{10},[P(\Sigma\theta)+\nu_6\kappa_9]_{10}\eta_{23},[\iota_7]_{10}\nu_7\kappa_{10},[\iota_7]_{10}\bar{\mu}_7\}\cong(\mathbb Z/2)^4.
		\end{equation*}
	\end{lem}
	
	\begin{proof}
		By Proposition~\ref{m1:R11},
		$\Delta\kappa_{10}=[\bar{\nu}_9\nu^2_{17}]+[P(\Sigma\theta)+\nu_6\kappa_9]_{10}\bmod16[\sigma^2_9].$ 
		Additionally, \cite[p.~96]{KM2} gives
		\begin{equation}
			\Delta(\eta_{10})=2[\iota_7]_{10}\nu_7.\label{m1:Det10}
		\end{equation}
		Therefore,
		$\Delta\kappa_{10}\circ\eta_{23}=\Delta(\eta_{10}\kappa_{11})=2[\iota_7]_{10}\nu_7\kappa_{10}=0,$  
		which proves the first relation.    The relevant group has the decomposition
		$\pi^9_{25}=\mathrm{span}\{\sigma_9\nu^3_{16},\sigma_9\mu_{16},\sigma_9\eta_{16}\varepsilon_{17},\mu_9\sigma_{18}\}\cong(\mathbb Z/2)^4.$ 
		We calculate the boundaries of these generators in $(24)_9.$ 
		Using $\Delta\sigma_9=[\iota_7]_9(\bar{\nu}_7+\varepsilon_7)$
		from \eqref{m1:Desig9}, together with $2\kappa_7=\bar{\nu}_7\nu^2_{15}$ and $\bar{\nu}_6\mu_{14}=0,$ we obtain
		\begin{align}
			\Delta(\sigma_9\nu^3_{16})&=[\iota_7]_9(\bar{\nu}_7+\varepsilon_7)\nu^3_{15}=[\iota_7]_9\bar{\nu}_7\nu^3_{15}=2[\iota_7]_9\nu_7\kappa_{10}=0;\notag\\
			\Delta(\sigma_9\mu_{16})&=[\iota_7]_9(\bar{\nu}_7+\varepsilon_7)\mu_{15}=[\iota_7]_9\eta_7\mu_8\sigma_{17}.\notag
		\end{align}
		Thus exactness gives $[\iota_7]_{10}\eta_7\mu_8\sigma_{17}=0.$ 
		By Lemma~\ref{m1:zesig}(2) and \eqref{m1:epep},
		\begin{equation*}
			\Delta(\sigma_9\eta_{16}\varepsilon_{17})=\Delta(\sigma_9\varepsilon_{16}\eta_{24})=[\iota_7]_9(\bar{\nu}_7+\varepsilon_7)\varepsilon_{15}\eta_{23}=[\iota_7]_9\sigma'\nu^3_{14}\eta_{23}=0.
		\end{equation*}
		Finally, the standard relations give
		$\Delta(\mu_9\sigma_{18})=[\iota_7]_9\eta_7\mu_8\sigma_{17}+[ \nu_5]_9\mu_8\sigma_{17}.$ 
		Exactness now yields $[ \nu_5]_{10}\mu_8\sigma_{17}=0$ as well.
		Applying these four boundary values to the description of
		$\pi_{24}(\mathit{SO}(9))$ gives the asserted image.
	\end{proof}
	We now turn to a further case, bringing together
	the relevant information in the following lemma.
	\begin{lem}\label{m1:so10-bracket-values}
		We have $\pi_{17}(\mathit{SO}(10))\circ2\sigma_{17}=0.$ 
		The following Toda brackets have zero indeterminacy:
		\begin{align}
			\{[\zeta_5]_{10},8\iota_{16},2\sigma_{16}\}_1&=[\bar{\zeta}_5]_{10};\notag\\
			\{[\iota_7]_{10}\mu_7,8\iota_{16},2\sigma_{16}\}_1&=[\iota_7]_{10}\bar{\mu}_7;\notag\\
			\{[\iota_7]_{10}\nu^3_7,8\iota_{16},2\sigma_{16}\}_1&=0.\notag
		\end{align}
	\end{lem}
	
	\begin{proof}
		By \cite[Theorem 1.1(iii)]{KM1},
		$\pi_{17}(\mathit{SO}(10))=\mathrm{span}\{[\iota_7]_{10}\nu_7\sigma_{10},[\iota_7]_{10}\eta_7\mu_8\}.$ 
		Composition with $2\sigma_{17}$ annihilates both generators.
		Relation \eqref{m1:[bzeta_5]} gives
		\begin{equation*}
			\{[\zeta_5]_{10},8\iota_{16},2\sigma_{16}\}_1\ni[\bar{\zeta}_5]_{10}
		\end{equation*}
		with indeterminacy
		$[\zeta_5]_{10}\circ\Sigma\pi^{15}_{23}+\pi_{17}(\mathit{SO}(10))\circ2\sigma_{17}.$ 
		By Lemma~\ref{m1:i710ms}(2), $[\zeta_5]_{10}\sigma_{16}=0.$ 
		Together with $\eta_9\sigma_{10}=\bar{\nu}_9+\varepsilon_9,$  
		this gives $[\zeta_5]_{10}\eta_{16}\sigma_{17}=0$ and
		$[\zeta_5]_{10}\bar{\nu}_{16}=[\zeta_5]_{10}\varepsilon_{16}.$ 
		On composing the congruence in Lemma~\ref{m1:i710ms}(1) with $\varepsilon_{16},$ the terms involving $8\Delta\sigma_{10}$ and $[\iota_7]_{10}\nu_7^3$ vanish, since $8\varepsilon_{16}=0$ and $\nu_7^3\varepsilon_{16}=0.$ The latter relation follows by suspension from the composition relations in Table~\ref{tab:oguchi-compositions}.
		Together with \eqref{m1:Ogu0} and Lemma~\ref{m1:so10-image}, this gives
		\begin{equation*}
			[\zeta_5]_{10}\varepsilon_{16}=[\iota_7]_{10}\mu_7\varepsilon_{16}=[\iota_7]_{10}\eta_7\mu_8\sigma_{17}=0.
		\end{equation*}
		Hence $[\zeta_5]_{10}\circ\pi^{16}_{24}=0,$   and the first
		bracket has the stated value with zero indeterminacy.
		
		Next, \cite[p.~136]{T} gives
		$\{\mu_7,8\iota_{16},2\sigma_{16}\}_1\supseteq\{\mu_7,2\iota_{16},8\sigma_{16}\}\ni\bar{\mu}_7.$ 
		Thus the second bracket contains $[\iota_7]_{10}\bar{\mu}_7.$ 
		Recall that $\pi^{15}_{23}=\mathrm{span}\{\bar{\nu}_{15},\varepsilon_{15}\},$   that $\mu_7\bar{\nu}_{16}=0,$   and that
		$\mu_7\varepsilon_{16}=\eta_7\mu_8\sigma_{17}$ by \eqref{m1:Ogu0}.
		The indeterminacy is therefore
		\begin{equation*}
			[\iota_7]_{10}\mu_7\circ\Sigma\pi^{15}_{23}+\pi_{17}(\mathit{SO}(10))\circ2\sigma_{17}=\mathrm{span}\{[\iota_7]_{10}\eta_7\mu_8\sigma_{17}\}=0.
		\end{equation*}
		
		Finally, $\zeta_{13}=\{\nu_{13},16\iota_{16},\sigma_{16}\}$
		from \cite[(9.1)]{T}, together with \eqref{m1:n6nz}, gives
		\begin{equation*}
			\{[\iota_7]_{10}\nu^3_7,8\iota_{16},2\sigma_{16}\}_1\supseteq[\iota_7]_{10}\nu^2_7\circ\{\nu_{13},8\iota_{16},2\sigma_{16}\}_1\ni[\iota_7]_{10}\nu^2_7\zeta_{13}=0.
		\end{equation*}
		Its indeterminacy is
		$[\iota_7]_{10}\nu^3_7\circ\Sigma\pi^{15}_{23}+\pi_{17}(\mathit{SO}(10))\circ2\sigma_{17}=0,$  
		which proves the third assertion.
	\end{proof}
	Continuing our examination of the groups under consideration,
	we arrive at the following explicit description.
	\begin{lem}\label{m1:so10-epsilon-lift}
		We have
		\begin{equation}
			\Delta(\nu^2_{11})=[\iota_7]_{11}\nu^3_7\neq0\label{m1:Dn211}
		\end{equation}
		and
		\begin{equation}
			\Delta(\nu^2_{10})=[\iota_7]_{10}(\bar{\nu}_7+\eta_7\sigma_8)=[\iota_7]_{10}(\varepsilon_7+\sigma'\eta_{14}).\label{m1:Dnu^2}
		\end{equation}
		A lift $[\bar{\varepsilon}_9]$ may be chosen so that
		\begin{equation*}
			[\bar{\varepsilon}_9]\in\left\{\bigl(\Delta(\nu_{10}),[\iota_7]_{10}\bigr),\begin{pmatrix}\nu_{12}\\\varepsilon_7+\sigma'\eta_{14}\end{pmatrix},\bar{\nu}_{15}\right\}.
		\end{equation*}
		In addition,\vspace{-1em}
		\begin{align}
			\{\nu_{12},\bar{\nu}_{15},2\iota_{23}\}&=0;\notag\\
			\{\varepsilon_7+\sigma'\eta_{14},\bar{\nu}_{15},2\iota_{23}\}&=\nu_7\kappa_{10}+\mathrm{span}\{\sigma'\eta_{14}\mu_{15}+\eta_7\mu_8\sigma_{17}\}.\notag
		\end{align}
		For the chosen lift, $2[\bar{\varepsilon}_9]\equiv[\iota_7]_{10}\nu_7\kappa_{10}\bmod[\iota_7]_{10}\sigma'\eta_{14}\mu_{15},[\iota_7]_{10}\eta_7\mu_8\sigma_{17}.$
	\end{lem}
	
	\begin{proof}
		As a consequence of          \cite[p.~161, Table]{Ke} and \cite{Bo}, we know 
		$\pi_{17}(\mathit{SO}(13))=\frac{\mathbb Z}{2}\{[\iota_7]_{13}\eta_7\mu_8\}.$ 
		The sequence $(17)_{12}$ gives
		$\pi_{17}(\mathit{SO}(12))=\mathrm{span}\{[\iota_7]_{12}\eta_7\mu_8\}+\Delta(\pi^{12}_{18}).$ 
		Since $${p_{12}}_*\pi_{17}(\mathit{SO}(12))\subseteq\mathrm{span}\{2\nu^2_{11}\}=0,$$\noindent 
		the sequence $(17)_{11}$ shows that
		$\Delta\colon \pi^{11}_{17}\to\pi_{16}(\mathit{SO}(11))$ is injective,
		which gives \eqref{m1:Dn211}.
		
		In light of \cite[p.~22]{Ka}, we know
		$\Delta(\nu^2_{10})=[\iota_7]_{10}\bar{\nu}_7+x[ \nu_5]_{10}\sigma_8$
		for some $x\in\{0,1\}.$  Since $[ \nu_5]_{10}=[\iota_7]_{10}\eta_7,$  
		this becomes
		$\Delta(\nu^2_{10})=[\iota_7]_{10}(\bar{\nu}_7+x\eta_7\sigma_8).$ 
		Composing with $\eta_{15}$ and using
		$\bar{\nu}_7\eta_{15}=\nu^3_7,$   we obtain
		$0=[\iota_7]_{10}(\nu^3_7+x\eta_7\sigma_8\eta_{15}).$ 
		By \eqref{m1:Dn211}, $[\iota_7]_{10}\nu^3_7\neq0,$  
		so $x=1.$  This proves \eqref{m1:Dnu^2}.
		
		The standard relations give $\nu_{12}\bar{\nu}_{15}=0.$ 
		Moreover, $(\varepsilon_7+\sigma'\eta_{14})\bar{\nu}_{15}=0,$  
		because $\varepsilon_3\bar{\nu}_{11}=\eta_3\bar{\varepsilon}_4$
		and $\sigma'\eta_{14}\bar{\nu}_{15}=\eta_7\bar{\varepsilon}_8$
		by \eqref{m1:epep}. Together with \eqref{m1:Dnu^2},
		these relations show that the stated matrix Toda bracket is defined.
		Let $A$ be a representative. By
		\cite[Lemma 4.1(2)]{MiM1},
		\begin{equation*}
			p_{10\ast}(A)\in\left\{\bigl(2\nu_9,0\bigr),\begin{pmatrix}\nu_{12}\\\varepsilon_7+\sigma'\eta_{14}\end{pmatrix},\bar{\nu}_{15}\right\}=\{2\nu_9,\nu_{12},\bar{\nu}_{15}\}.
		\end{equation*}
		By \eqref{m1:tbep5},
		$\bar{\varepsilon}_9\in\{2\nu_9,\nu_{12},\bar{\nu}_{15}\}\bmod2\nu_9\circ\pi^{12}_{24}+\pi^9_{16}\circ\bar{\nu}_{16}=0.$ 
		Thus $A$ is a lift of $\bar{\varepsilon}_9,$   and we may
		choose $[\bar{\varepsilon}_9]=A.$ 
		
		For the two scalar brackets, we have
		\begin{equation*}
			\{\nu_{12},\bar{\nu}_{15},2\iota_{23}\}\supseteq\Sigma^3\{\nu_9,\bar{\nu}_{12},2\iota_{20}\}=0\bmod\nu_{12}\circ\pi^{15}_{24}+2\pi^{12}_{24}=0,
		\end{equation*}
		since $\pi^9_{21}=0.$ 
		Because $\Sigma\pi^{14}_{23}=\pi^{15}_{24},$   the brackets
		$\{\varepsilon_7+\sigma'\eta_{14},\bar{\nu}_{15},2\iota_{23}\}_1$
		and $\{\varepsilon_7+\sigma'\eta_{14},\bar{\nu}_{15},2\iota_{23}\}$
		have the same indeterminacy. 
		Therefore we have the asserted element of the second bracket.
		The matrix Toda bracket relation gives
		\begin{equation*}
			2[\bar{\varepsilon}_9]\in-\Delta(\nu_{10})\circ\{\nu_{12},\bar{\nu}_{15},2\iota_{23}\}+[\iota_7]_{10}\circ\{\varepsilon_7+\sigma'\eta_{14},\bar{\nu}_{15},2\iota_{23}\}.
		\end{equation*}
		The two scalar bracket values just established yield the
		stated congruence for $2[\bar{\varepsilon}_9].$ 
	\end{proof}
	We proceed with our analysis by examining a further
	case, whose outcome is recorded in the following proposition.
	\begin{prop}\label{m1:10}
		The structure of $\pi_{24}(\mathit{SO}(10))$ is as follows.
		\begin{equation*}
			\pi_{24}(\mathit{SO}(10))=\mathrm{span}\{\Delta(\Sigma\rho'),[\bar{\varepsilon}_9],[\sigma^2_9]\eta_{23},[\iota_7]_{10}\bar{\mu}_7,[P(\Sigma\theta)+\nu_6\kappa_9]_{10}\eta_{23}\}\cong\mathbb Z/16\oplus\mathbb Z/4\oplus(\mathbb Z/2)^3.
		\end{equation*}
		The generators satisfy
		$[P(\Sigma\theta)+\nu_6\kappa_9]_{10}\eta_{23}=[\bar{\nu}_9\nu^2_{17}]\eta_{23},$  
		$[\iota_7]_{10}\eta_7\mu_8\sigma_{17}=[ \nu_5]_{10}\mu_8\sigma_{17}=0,$  
		$8\Delta(\Sigma\rho')=[\bar{\zeta}_5]_{10}+[\iota_7]_{10}\bar{\mu}_7$
		and $2[\bar{\varepsilon}_9]=[\iota_7]_{10}\nu_7\kappa_{10}.$ 
	\end{prop}
	
	\begin{proof}
		The first two relations follow from Lemma~\ref{m1:so10-image}.
		We prove the two remaining relations and then identify the group.
		Recall that
		$\rho^{\mathrm{IV}}\in\{\sigma''',2\iota_{12},8\sigma_{12}\}_1\subseteq\{\sigma''',8\iota_{12},2\sigma_{12}\}_1,$  
		and $8\rho'=\Sigma^4\rho^{\mathrm{IV}}.$ 
		Hence $8\rho'\in\{8\sigma_9,8\iota_{16},2\sigma_{16}\}_1.$ 
		Suspending and applying the boundary homomorphism gives
		$8\Sigma\rho'\in\{8\sigma_{10},8\iota_{17},2\sigma_{17}\}_2$
		and $8\Delta(\Sigma\rho')\in\{8\Delta(\sigma_{10}),8\iota_{16},2\sigma_{16}\}_1.$ 
		By Lemma~\ref{m1:i710ms}(1),
		$8\Delta\sigma_{10}=[\zeta_5]_{10}+[\iota_7]_{10}\mu_7+y[\iota_7]_{10}\nu^3_7$
		for some $y\in\{0,1\}.$  Additivity and
		Lemma~\ref{m1:so10-bracket-values} therefore give
		\begin{align}
			8\Delta(\Sigma\rho')
			&\in\{[\zeta_5]_{10},8\iota_{16},2\sigma_{16}\}_1
			+\{[\iota_7]_{10}\mu_7,8\iota_{16},2\sigma_{16}\}_1+y\{[\iota_7]_{10}\nu^3_7,8\iota_{16},2\sigma_{16}\}_1\notag\\
			&=[\bar{\zeta}_5]_{10}+[\iota_7]_{10}\bar{\mu}_7.\notag
		\end{align}
		The right-hand side represents a singleton, so
		$8\Delta(\Sigma\rho')=[\bar{\zeta}_5]_{10}+[\iota_7]_{10}\bar{\mu}_7.$ 
		
		By \cite[(3.9)]{MMO},
		$H(\tilde{\varepsilon}_{10})=\bar{\varepsilon}_{19}$ and
		$2\tilde{\varepsilon}_{10}=\sigma_{10}\nu_{17}\kappa_{20}.$ 
		The first identity yields
		\begin{equation*}
			J[\bar{\varepsilon}_9]\equiv\tilde{\varepsilon}_{10}\bmod\Sigma\pi^9_{33}.
		\end{equation*}
		\noindent 
		For the lift chosen in Lemma~\ref{m1:so10-epsilon-lift},
		that lemma gives
		\begin{equation*}
			2[\bar{\varepsilon}_9]\equiv[\iota_7]_{10}\nu_7\kappa_{10}\bmod[\iota_7]_{10}\sigma'\eta_{14}\mu_{15},[\iota_7]_{10}\eta_7\mu_8\sigma_{17}.
		\end{equation*}
		\noindent   In light of Proposition~\ref{m1:924}
		and Lemma~\ref{m1:so10-image},
		\begin{equation*}
			[\iota_7]_{10}\sigma'\eta_{14}\mu_{15}=[ \nu_5]_{10}\sigma_8\mu_{15}=[ \nu_5]_{10}\mu_8\sigma_{17}=[\iota_7]_{10}\eta_7\mu_8\sigma_{17}=0.
		\end{equation*}
		\noindent   Thus $2[\bar{\varepsilon}_9]=[\iota_7]_{10}\nu_7\kappa_{10}.$ 
		The image in Lemma~\ref{m1:so10-image} shows that
		$[\bar{\varepsilon}_9]$ has order $4$ and
		$\Delta(\Sigma\rho')$ has order $16.$ 
		Combining these relations with that image and the exact
		sequence $(24)_9$ gives the stated decomposition.
	\end{proof}

	The calculations developed in the preceding discussion
	also shed light on the present case, leading to
	the following determination.
	
	\begin{prop} \label{m1:11}
		The structure of $\pi_{24}(\mathit{SO}(11))$ is as follows.
		\[ \GroupLayout{\pi_{24}(\mathit{SO}(11)) = \mathrm{span}\{[\bar{\varepsilon}_9]_{11},\; [\sigma^2_9]_{11}\eta_{23},\; [\iota_7]_{11}\bar{\mu}_7,\; [P(\Sigma\theta)+\nu_6\kappa_9]_{11}\eta_{23}\}\cong \mathbb Z/4\oplus (\mathbb Z/2)^{3}.}{\begin{aligned}\pi_{24}(\mathit{SO}(11)) = \mathrm{span}\{&[\bar{\varepsilon}_9]_{11},\; [\sigma^2_9]_{11}\eta_{23},\; [\iota_7]_{11}\bar{\mu}_7,\\
				&[P(\Sigma\theta)+\nu_6\kappa_9]_{11}\eta_{23}\}\cong \mathbb Z/4\oplus (\mathbb Z/2)^{3}.\end{aligned}} \]
		Here $[\bar{\zeta}_5]_{11}=[\iota_7]_{11}\bar{\mu}_7.$ 
	\end{prop}
	\begin{proof}
		In the exact sequence $(24)_{10}$:
		\[\begin{aligned}
			&\pi^{10}_{25}\xrightarrow{\Delta}\pi_{24}(\mathit{SO}(10))\xrightarrow{i_*}\pi_{24}(\mathit{SO}(11)) \xrightarrow{p_*}\pi^{10}_{24}\xrightarrow{\Delta}\pi_{23}(\mathit{SO}(10)),
		\end{aligned}\]
		we note that $\Delta\colon  \pi^{10}_{24}\to \pi_{23}(\mathit{SO}(10))$ is a monomorphism and $\pi^{10}_{25}
		= \frac{\mathbb Z}{16}\{\Sigma\rho'\}\oplus\frac{\mathbb Z}{2}\{\sigma_{10}\bar{\nu}_{17}\}
		\oplus\frac{\mathbb Z}{2}\{\bar{\varepsilon}_{10}\}.$ 
		Since  $\sigma_{10}\bar{\nu}_{17}=[\iota_{10},\nu^2_{10}]$
		\cite[(10.20)]{T} and $\Delta[\iota_{10},\nu_{10}]=0$ \cite[(4.4)]{HKM},
		we see that
		\[\begin{aligned}
			&\Delta(\sigma_{10}\bar{\nu}_{17})=\Delta([\iota_{10},\nu_{10}]\circ\nu_{22}) =\Delta[\iota_{10},\nu_{10}]\circ \nu_{21}=0.
		\end{aligned}\]
		
		Using \eqref{m1:Det10}
		and $\Delta(\eta_{10}\kappa_{11})=\Delta(\eta_{10})\kappa_{10}=2[\iota_7]_{10}\nu_7\kappa_{10}=0$ in the proof of Lemma~\ref{m1:so10-image}, we have
		$\Delta\bar{\varepsilon}_{10} = 0.$ 
		This completes the proof.
	\end{proof}
	
	Continuing our examination of the groups under consideration,
	we arrive at the following explicit description.
	\begin{lem}\label{m1:[bep9]11eta}
		\begin{enumerate}
			\item
			$[\bar{\varepsilon}_9]\eta_{24} \equiv 0\ \bmod\
			[\iota_7]_{10}\zeta_7\sigma_{18},
			[\iota_7]_{10}\eta_7\bar\mu_8.$ 
			\item
			$[\bar\varepsilon_9]_{11}\eta_{24}\equiv 0\ \bmod\ [\iota_7]_{11}\eta_7\bar\mu_8.$ 
		\end{enumerate}
	\end{lem}
	\begin{proof}
		\begin{enumerate}[label=\textup{(\arabic*)},leftmargin=2em]
			\item In light of the definition of $[\bar{\varepsilon}_9]$ and \cite[Lemma 4.1(2)]{MiM1},
			we deduce
			\begin{align}
				[\bar{\varepsilon}_9]\eta_{24}
				&\in
				\left\{
				\bigl(\Delta(\nu_{10}),[\iota_7]_{10}\bigr),\,
				\begin{pmatrix}
					\nu_{12}\\
					\varepsilon_7+\sigma'\eta_{14}
				\end{pmatrix},
				\bar\nu_{15}
				\right\}\circ\eta_{24}
				\notag\\
				&\subseteq
				\left\{
				\bigl(\Delta(\nu_{10}),[\iota_7]_{10}\bigr),\,
				\begin{pmatrix}
					\nu_{12}\\
					\varepsilon_7+\sigma'\eta_{14}
				\end{pmatrix},
				\nu^3_{15}
				\right\}
				\notag\\
				&\subseteq
				\left\{
				\bigl(\Delta(\nu_{10}),[\iota_7]_{10}\bigr),\,
				\begin{pmatrix}
					\nu^3_{12}\\
					0
				\end{pmatrix},
				\nu_{21}
				\right\}
				\notag\\
				&\xyd
				\{\Delta(\nu_{10}),\nu^3_{12},\nu_{21}\}
				+[\iota_7]_{10}\circ\pi^7_{25}.
				\notag
			\end{align}\noindent    Observe that
			\[\begin{aligned}
				\{\Delta(\nu_{10}),\nu^3_{12},\nu_{21}\} &=\{\Delta(\nu_{10}),\eta_{12}\bar{\nu}_{13},\nu_{21}\} \\
				&\supseteq\{\Delta(\nu_{10}), \eta_{12},0\} \\
				&\ni 0 \bmod\ \Delta(\nu_{10})\circ\pi^{12}_{25} + \pi_{22}(\mathit{SO}(10))\circ\nu_{22}.
			\end{aligned}\]
			Since $\pi^{12}_{25}=\mathrm{span}\{\theta\eta_{24}, (\Sigma\theta')\eta_{24}\}\cong(\mathbb Z/2)^2,$  
			$\Sigma\pi^{12}_{25}=\pi^{13}_{26}=\mathrm{span}\{(\Sigma\theta)\eta_{25}\}\cong\mathbb Z/2$ and $(\Sigma\theta)\eta_{25}=P(\eta_{27})$ \cite[(7.30)]{T}, we find
			$\Delta(\nu_{10})\circ\pi^{12}_{25}=\Delta\bigl(\nu_{10}\circ\Sigma\pi^{12}_{25}\bigr)
			=\mathrm{span}\{\Delta(\nu_{10}\circ P(\eta_{27}))\}=0.$  Therefore, we see that
			\[\begin{aligned}
				&\{\Delta(\nu_{10}),\nu^3_{12},\nu_{21}\}=\pi_{22}(\mathit{SO}(10))\circ\nu_{22}.
			\end{aligned}\]
			Since
			$\pi_{22}(\mathit{SO}(10))=\mathrm{span}\{[8\sigma_8]_{10}\sigma_{15},[\iota_7]_{10}\rho'',[\iota_7]_{10}\bar\varepsilon_7\}$
			\cite[Theorem 0.2(iii)]{HKM}, it follows  that
			\[\begin{aligned}
				&[\bar{\varepsilon}_9]\eta_{24}\in \pi_{22}(\mathit{SO}(10))\circ\nu_{22} + [\iota_7]_{10}\circ \pi^7_{25} = [\iota_7]_{10}\circ \pi^7_{25}.
			\end{aligned}\]
			Since $\pi^7_{25} = \mathrm{span}\{\zeta_7\sigma_{18},
			\eta_7\bar{\mu}_8\}\cong\mathbb Z/8\oplus\mathbb Z/2,$   we obtain
			\[\begin{aligned}
				&[\iota_7]_{10}\circ\pi^7_{25}={}\mathrm{span}\{[\iota_7]_{10}\zeta_7\sigma_{18}, [\iota_7]_{10}\eta_7\bar{\mu}_8\}.
			\end{aligned}\]
			This leads to (1).
			
			\item (1) implies (2) if we show $[\iota_7]_{11}\zeta_7\sigma_{18} = 0.$ 
			Note that
			$\rho''\nu_{22}=\sigma'\zeta_{14}.$
			We deduce that $\pi^{11}_{26} = \frac{\mathbb Z}{16}\{\Sigma^2\rho'\} \oplus
			\frac{\mathbb Z}{2}\{\bar\varepsilon_{11}\}.$  Using the standard vanishing relation,   we have
			\[\begin{aligned}
				&\Delta(\Sigma^2\rho')= [\iota_7]_{11}\nu_7(\Sigma\rho') =0.
			\end{aligned}\]
			The relation $\Delta(\eta_{11})=0$  implies that    $\Delta(\bar\varepsilon_{11})=\Delta(\eta_{11})\kappa_{11} = 0.$
			So it suffices to show $[\iota_7]_{12}\rho''\nu_{22}=0.$  By \cite[Lemma 6.3]{HKM} and the fact that $\mathrm{ord} P(\nu_{25})=4,$   we conclude that
			\[\begin{aligned}
				&[\iota_7]_{12}\rho''\nu_{22} =4[[\e_{10},\nu_{10}]]_{12}\nu_{22}=4\Delta(P(\nu_{25}))=0.
			\end{aligned}\]
			This leads to (2) and completes the proof.
		\end{enumerate}
	\end{proof}
	
	\subsection{The Cases $n=12,13$}
	We show the following.
	
	\begin{prop} \label{m1:12}
		The structure of $\pi_{24}(\mathit{SO}(12))$ is as follows.
		\[\GroupLayout{\pi_{24}(\mathit{SO}(12)) \;= {} \mathrm{span}\{[\theta']\eta_{23},\; [\eta_{11}]\theta,\; [\bar{\varepsilon}_9]_{12},\; [\sigma^2_9]_{12}\eta_{23},\; [\iota_7]_{12}\bar{\mu}_7,\; [P(\Sigma\theta)+\nu_6\kappa_9]_{12}\eta_{23}\}\cong (\mathbb Z/2)^{6}.}{\begin{aligned}\pi_{24}(\mathit{SO}(12)) \;= {} \mathrm{span}\{&[\theta']\eta_{23},\; [\eta_{11}]\theta,\; [\bar{\varepsilon}_9]_{12},\; [\sigma^2_9]_{12}\eta_{23},\\
				&[\iota_7]_{12}\bar{\mu}_7,\; [P(\Sigma\theta)+\nu_6\kappa_9]_{12}\eta_{23}\}\cong (\mathbb Z/2)^{6}.\end{aligned}}\]
	\end{prop}
	\begin{proof}
		In the exact sequence $(24)_{11}$:
		\[\begin{aligned}
			&\pi^{11}_{25}\xrightarrow{\Delta}\pi_{24}(\mathit{SO}(11))\xrightarrow{i_*}\pi_{24}(\mathit{SO}(12)) \xrightarrow{p_*}\pi^{11}_{24}\xrightarrow{\Delta}\pi_{23}(\mathit{SO}(11)),
		\end{aligned}\]
		we know that $\pi^{11}_{25} = \mathrm{span}\{\sigma^2_{11}, \kappa_{11}\}\cong\mathbb Z/16
		\oplus\mathbb Z/2$
		and
		\[\begin{aligned}
			&\operatorname{ker}\ (\Delta\colon  \pi^{11}_{24}\to \pi_{23}(\mathit{SO}(11))) = \pi^{11}_{24} = \mathrm{span}\{\theta'\eta_{23}, \eta_{11}\theta\}\cong(\mathbb Z/2)^2
		\end{aligned}\]
		from the proof of Proposition~\ref{m1:R12}.
		Using \eqref{m1:n7sigsig} and  Proposition~\ref{m1:10},
		we see that
		$\Delta(\sigma^2_{11}) = [\iota_7]_{11}\nu_7\sigma^2_{10} = 0$ and
		$
		\Delta\kappa_{11}=[\iota_7]_{11}\nu_7\kappa_{10}
		=2[\bar{\varepsilon}_9]_{11}.
		$
		So Proposition~\ref{m1:11} implies the conclusion.
	\end{proof}
	
	We continue our discussion by presenting the following statement, whose proof will be given before we move on to other matters.
	\begin{prop} \label{m1:13}
		The structure of $\pi_{24}(\mathit{SO}(13))$ is as follows.
		\[ \GroupLayout{\pi_{24}(\mathit{SO}(13)) = \mathrm{span}\{[\eta_{11}]_{13}\theta,\; [\bar{\varepsilon}_9]_{13},\; [\sigma^2_9]_{13}\eta_{23},\; [\iota_7]_{13}\bar{\mu}_7\}\cong (\mathbb Z/2)^{4}.}{\begin{aligned}\pi_{24}(\mathit{SO}(13)) = \mathrm{span}\{&[\eta_{11}]_{13}\theta,\; [\bar{\varepsilon}_9]_{13},\; [\sigma^2_9]_{13}\eta_{23},\\
				&[\iota_7]_{13}\bar{\mu}_7\}\cong (\mathbb Z/2)^{4}.\end{aligned}} \]
		Here
		$\Delta(\Sigma\theta)=[\eta_{11}]_{13}\theta.$ 
	\end{prop}
	\begin{proof}
		In the exact sequence $(24)_{12}$:
		\[\begin{aligned}
			&\pi^{12}_{25}\xrightarrow{\Delta}\pi_{24}(\mathit{SO}(12))\xrightarrow{i_*}\pi_{24}(\mathit{SO}(13)) \xrightarrow{p_*}\pi^{12}_{24}\xrightarrow{\Delta}\pi_{23}(\mathit{SO}(12)),
		\end{aligned}\]
		$\Delta\colon  \pi^{12}_{24}\to \pi_{23}(\mathit{SO}(12))$ is a monomorphism
		by the proof of Proposition~\ref{m1:R13} and
		\[\begin{aligned}
			&\pi^{12}_{25}={}\frac{\mathbb Z}{2}\{\theta\eta_{24}\}\oplus\frac{\mathbb Z}{2}\{(\Sigma\theta')\eta_{24}\}.
		\end{aligned}\]
		As $\Delta\theta = [\theta']$ and $\Delta(\Sigma\theta')
		= [P(\Sigma\theta)+\nu_6\kappa_9]_{12}$ by Proposition~\ref{m1:R13}, we obtain
		\[\begin{aligned}
			&\Delta(\theta\eta_{24}) = [\theta']\eta_{23} \quad\text{and}\quad \Delta((\Sigma\theta')\eta_{24}) = [P(\Sigma\theta)+\nu_6\kappa_9]_{12}\eta_{23}.
		\end{aligned}\]
		
		Since $\Delta\iota_{13} = [\eta_{11}]_{13}$
		\cite[Table 3]{Ka}, we deduce
		$\Delta(\Sigma\theta) = [\eta_{11}]_{13}\theta.$ 
		This completes the proof.
	\end{proof}

	\subsection{The Case $n=14$}
	Before moving further into the discussion, we state the following result and give its proof in the course of our presentation.
	\begin{lem} \label{m1:14}
		$\pi_{24}(\mathit{SO}(14)) = \mathrm{span}\{[\zeta_{13}],\; [\bar{\varepsilon}_9]_{14},\; [\sigma^2_9]_{14}\eta_{23},\; [\iota_7]_{14}\bar{\mu}_7\}\cong \mathbb Z/8\oplus (\mathbb Z/2)^{3},$   where $J[\zeta_{13}]
		\equiv \pm\zeta^* \bmod \Sigma\pi^{13}_{37}.$ 
	\end{lem}
	\begin{proof}
		In the exact sequence $(24)_{13}$:
		\[\begin{aligned}
			&\pi^{13}_{25}\xrightarrow{\Delta}\pi_{24}(\mathit{SO}(13))\xrightarrow{i_*}\pi_{24}(\mathit{SO}(14)) \xrightarrow{p_*}\pi^{13}_{24}\xrightarrow{\Delta}\pi_{23}(\mathit{SO}(13)),
		\end{aligned}\]
		we see that
		$\Delta\colon  \pi^{13}_{24}\to \pi_{23}(\mathit{SO}(13))$ is trivial by Lemma~\ref{m1:zeta1}
		and $\pi^{13}_{24} = \frac{\mathbb Z}{8}\{\zeta_{13}\}.$ 
		We  know $\pi^{13}_{25} = \frac{\mathbb Z}{2}\{\Sigma\theta\}.$ 
		By Lemma~\ref{m1:zeta1}, we find
		$\Delta\zeta_{14} = \pm 2[\zeta_{13}]$ and
		$\mathrm{ord}[\zeta_{13}]=8.$ 
		This determines the group $\pi_{24}(\mathit{SO}(14)).$  By \cite[(3.10)]{MMO}, $H(\zeta^*)
		= \pm\zeta_{27}.$  Thus we conclude that $J[\zeta_{13}]
		\equiv\pm\zeta^* \bmod \Sigma\pi^{13}_{37}.$  This completes the proof.
	\end{proof}
	
	Using the fact that $\pi_{24}(\mathit{SO})\cong \mathbb Z/2$ \cite{Bo} and the results in \cite[Table 1]{HoM},
	we have Table~\ref{m1:R24}.
	\par
	\begin{table}[H]
		\centering
		\small\setlength{\tabcolsep}{3pt}\begin{tabular}{cccccc}\toprule
			$\pi_{24}(\mathit{SO}(14))$ & $\mathbb Z/2 \oplus \mathbb Z/8 \oplus (\mathbb Z/2)^2$
			& $\pi_{24}(\mathit{SO}(15))$ & $\mathbb Z/2 \oplus (\mathbb Z/2)^4$
			& $\pi_{24}(\mathit{SO}(16))$ & $\mathbb Z/2 \oplus (\mathbb Z/2)^7$  \\ \midrule
			$\pi_{24}(\mathit{SO}(17))$ & $\mathbb Z/2 \oplus (\mathbb Z/2)^4$
			& $\pi_{24}(\mathit{SO}(18))$ & $\mathbb Z/2 \oplus \mathbb Z/16 \oplus \mathbb Z/2$
			& $\pi_{24}(\mathit{SO}(19))$ & $\mathbb Z/2 \oplus \mathbb Z/2$\\ \bottomrule
		\end{tabular}
		\caption{$\pi_{24}(\mathit{SO}(n))$ for $14\leq n\leq 19$}\label{m1:R24}
	\end{table}
	\par
	By \cite{Ke}, we  have Table~\ref{m1:R24K}.
	\par
	\begin{table}[H]
		\centering
		\small\setlength{\tabcolsep}{3pt}\begin{tabular}{cccccc}\toprule
			$\pi_{24}(\mathit{SO}(20))$ & $\mathbb Z/2$
			& $\pi_{24}(\mathit{SO}(21))$ & $\mathbb Z/2$
			& $\pi_{24}(\mathit{SO}(22))$ & $\mathbb Z/4 \oplus \mathbb Z/2$  \\ \midrule
			$\pi_{24}(\mathit{SO}(23))$ & $(\mathbb Z/2)^2$
			& $\pi_{24}(\mathit{SO}(24))$ & $(\mathbb Z/2)^3$
			& $\pi_{24}(\mathit{SO}(25))$ & $(\mathbb Z/2)^2$\\ \bottomrule
		\end{tabular}
		\caption{$\pi_{24}(\mathit{SO}(n))$ for $20\leq n\leq 25$}\label{m1:R24K}
	\end{table}
	\par
	
	In light of Tables~\ref{m1:R24} and \ref{m1:R24K}, we have the group structure of $\pi_{24}(\mathit{SO}(n))$ for $n \ge 14.$ 
	We now determine the generators and relations in
	$\pi_{24}(\mathit{SO}(n))$ for $n \ge 14.$ 
	
	\section{The Cases $n\geq15$}
	\subsection{The Cases $15\leq n\leq17$}
	
	We now turn to the following statement, which we present here before moving on to the rest of the discussion.
	
	\begin{lem} \label{m1:15} The group $   \pi_{24}(\mathit{SO}(15))$ is given by
		$$\GroupLayout{\pi_{24}(\mathit{SO}(15))=\mathrm{span}\{[\eta_{14}\mu_{15}],\; [\zeta_{13}]_{15},\; [\bar{\varepsilon}_9]_{15},\; [\sigma^2_9]_{15}\eta_{23},\; [\iota_7]_{15}\bar{\mu}_7\}\cong (\mathbb Z/2)^{5},}{\begin{aligned}\pi_{24}(\mathit{SO}(15))=\mathrm{span}\{&[\eta_{14}\mu_{15}],\; [\zeta_{13}]_{15},\; [\bar{\varepsilon}_9]_{15},\\
				&[\sigma^2_9]_{15}\eta_{23},\; [\iota_7]_{15}\bar{\mu}_7\}\cong (\mathbb Z/2)^{5},\end{aligned}}$$ where $J[\eta_{14}\mu_{15}]
		\equiv \mu^{*\prime} \bmod \Sigma\pi^{14}_{38}.$ 
	\end{lem}
	\begin{proof}
		In the exact sequence $(24)_{14}$:
		\[\begin{aligned}
			&\pi^{14}_{25}\xrightarrow{\Delta}\pi_{24}(\mathit{SO}(14))\xrightarrow{i_*}\pi_{24}(\mathit{SO}(15)) \xrightarrow{p_*}\pi^{14}_{24}\xrightarrow{\Delta}\pi_{23}(\mathit{SO}(14)),
		\end{aligned}\]
		$\Delta\colon  \pi^{14}_{24}\to \pi_{23}(\mathit{SO}(14))$ is zero, as shown in the proof of Proposition~\ref{m1:R1523} (1).
		By Lemma~\ref{m1:zeta1}, we have
		$\Delta\zeta_{14}=\pm 2[\zeta_{13}].$ 
		
		By Example~\ref{m1:shetep}, $\mathrm{ord}([\eta_{14}\mu_{15}])=2.$ 
		Since $HJ[\eta_{14}\mu_{15}] = \eta_{29}\mu_{30} = H(\mu^{*\prime})$ \cite[(3.11)]{MMO},
		we have the relation. This completes the proof.
	\end{proof}
	We now turn to a further case, bringing together
	the relevant information in the following lemma.
	
	\begin{lem} \label{m1:16} The group $   \pi_{24}(\mathit{SO}(16))$ is given by
		\[
		\begin{aligned}
			\pi_{24}(\mathit{SO}(16)) &={}\mathrm{span}\bigl\{ [\zeta_{13}]_{16}, [\sigma^2_9]_{16}\eta_{23}, [\iota_7]_{16}\bar{\mu}_7, [\bar{\varepsilon}_9]_{16}, [\eta_{14}\mu_{15}]_{16}, [\mu_{15}], [\eta_{15}]\varepsilon_{16}, [\eta_{15}]\bar{\nu}_{16} \bigr\} \\
			&\cong(\mathbb Z/2)^8.
		\end{aligned}
		\] Moreover, $J[\mu_{15}]\equiv \mu^*_{16}
		\bmod \Sigma\pi^{15}_{39}.$ 
	\end{lem}
	\begin{proof}
		In the exact sequence $(24)_{15}$:
		\[\begin{aligned}
			&\pi^{15}_{25}\xrightarrow{\Delta}\pi_{24}(\mathit{SO}(15))\xrightarrow{i_*}\pi_{24}(\mathit{SO}(16)) \xrightarrow{p_*}\pi^{15}_{24}\xrightarrow{\Delta}\pi_{23}(\mathit{SO}(15)),
		\end{aligned}\]
		we observe that $\Delta\colon  \pi^{15}_{24}\to \pi_{23}(\mathit{SO}(15))$ is trivial. Since
		$[\eta_{15}]$ has order $2$ \cite[Proposition 4.1]{Ka},
		the Toda bracket
		$\{[\eta_{15}], 2\iota_{16}, 8\sigma_{16}\}_1$ is well-defined.
		
		We define
		\[\begin{aligned}
			&[\mu_{15}]\in\{[\eta_{15}], 2\iota_{16}, 8\sigma_{16}\}_1,
		\end{aligned}\]
		and we observe that
		\begin{align*}
			&2[\mu_{15}]\in\{[\eta_{15}], 2\iota_{16}, 8\sigma_{16}\}_1\circ 2\iota_{24} =[\eta_{15}]\circ\{ 2\iota_{16}, 8\sigma_{16}, 2\iota_{23}\}_1.
		\end{align*}
		Moreover, 
		$\{ 2\iota_{16}, 8\sigma_{16}, 2\iota_{23}\}_1 \ni 8\sigma_{16}\eta_{23}=0$ \cite[Corollary 3.7]{T}
		and the indeterminacy is $2\iota_{16}\circ \pi^{16}_{24}+2\pi^{16}_{24}=0.$ 
		Hence, we obtain $2[\mu_{15}] = 0.$     By the relation
		$
		HJ[\mu_{15}]=\mu_{31}=H(\mu^\ast)$
		\cite[(3.11)]{MMO}, we have the relation.
		This completes the proof.
	\end{proof}
	
	It follows that $P(\nu^3_{33})\equiv\tilde{\varepsilon}_{16} \bmod\ \Sigma^7\pi^9_{33}$
	\cite[p. 39]{MMO}.
	Note that 
	$\Sigma^{n-9}\pi^9_{33}=\mathrm{span}\{\delta_n,\bar{\sigma}'_n, \bar{\mu}_n\sigma_{n+17}, \sigma_n\varepsilon^*_{n+7}\}$ for $n\ge 13.$ 
	According to           \cite[Lemma 12.15 ii)]{T} and \cite[I-Propositions 3.4(4), 3.5(9)]{Oda1979},
	we deduce $$\sigma_{11}\varepsilon^*_{18}=\sigma_{11}\omega_{18}\eta_{34}
	=\eta_{11}\sigma_{12}\omega_{19}=\eta_{11}\phi_{12}\equiv\delta_{11} \ \bmod\
	\bar{\mu}_{11}\sigma_{28}.$$ Hence, we obtain
	\begin{equation}
		\begin{aligned}
			&P(\nu^3_{33})\equiv\tilde{\varepsilon}_{16} \bmod\ \delta_{16},\bar{\sigma}'_{16}, \bar{\mu}_{16}\sigma_{33}.\label{m1:Pn333}
		\end{aligned}
	\end{equation}
	At this point, we pause to state and prove the following result before proceeding with the remainder of this section.
	
	\begin{prop} \label{m1:17}
		The structure of $\pi_{24}(\mathit{SO}(17))$ is as follows.
		\[ \GroupLayout{\pi_{24}(\mathit{SO}(17)) = \mathrm{span}\{[\mu_{15}]_{17},\; [\eta_{15}]_{17}\varepsilon_{16},\; [\eta_{15}]_{17}\bar{\nu}_{16},\; [\sigma^2_9]_{17}\eta_{23},\; [\iota_7]_{17}\bar{\mu}_7\}\cong (\mathbb Z/2)^{5}.}{\begin{aligned}\pi_{24}(\mathit{SO}(17)) = \mathrm{span}\{&[\mu_{15}]_{17},\; [\eta_{15}]_{17}\varepsilon_{16},\; [\eta_{15}]_{17}\bar{\nu}_{16},\\
				&[\sigma^2_9]_{17}\eta_{23},\; [\iota_7]_{17}\bar{\mu}_7\}\cong (\mathbb Z/2)^{5}.\end{aligned}} \]
		Here
		$\Delta\mu_{16}\equiv [\eta_{14}\mu_{15}]_{16} \bmod
		\Delta(\nu^3_{16}), \Delta(\eta_{16}\varepsilon_{17}),$  \
		$\Delta(\eta_{16}\varepsilon_{17}) = [\zeta_{13}]_{16}$\
		and\ $\Delta(\nu^3_{16})\equiv [\bar{\varepsilon}_9]_{16} \bmod
		[\iota_7]_{16}\bar{\mu}_7.$ 
	\end{prop}
	\begin{proof}
		In the exact sequence $(24)_{16}$:
		\[\begin{aligned}
			&\pi^{16}_{25}\xrightarrow{\Delta}\pi_{24}(\mathit{SO}(16))\xrightarrow{i_*}\pi_{24}(\mathit{SO}(17)) \xrightarrow{p_*}\pi^{16}_{24}\xrightarrow{\Delta}\pi_{23}(\mathit{SO}(16)),
		\end{aligned}\]
		$\Delta\colon \pi^{16}_{24}\to \pi_{23}(\mathit{SO}(16))$ is a monomorphism
		from the proof of Proposition~\ref{m1:R17}.
		By inspecting the calculation of $\pi^{16}_{40}$~\cite[(5.21)]{MMO} and
		\cite[(3.10), p.~39]{MMO}, we see that $P : \pi^{33}_{42} \to \pi^{16}_{40}$
		is a monomorphism. So is $\Delta : \pi^{16}_{25} \to \pi_{24}(\mathit{SO}(16)).$ 
		We know that
		\begin{align*}
			\Delta\mu_{16}&\in\Delta\{\eta_{16}, 2\iota_{17}, 8\sigma_{17}\}_1 \subseteq\{[\eta^2_{14}]_{16}, 2\iota_{16}, 8\sigma_{16}\} \\
			&\ni[\eta_{14}\mu_{15}]_{16}\ \bmod\ [\eta^2_{14}]_{16} \pi^{16}_{24} + \pi_{17}(\mathit{SO}(16))\circ 8\sigma_{17}.
		\end{align*}
		Since $\pi_{17}(\mathit{SO}(16))\cong(\mathbb Z/2)^3$ \cite[Theorem 1(iv)]{KM1}, we conclude that
		$\pi_{17}(\mathit{SO}(16))\circ 8\sigma_{17} = 0$ and $[\eta^2_{14}]_{16}\pi^{16}_{24}
		= \mathrm{span}\{[\eta^2_{14}]_{16}\bar{\nu}_{16}, [\eta^2_{14}]_{16}\varepsilon_{16}\}.$ 
		Thus we see that
		\[\begin{aligned}
			&\Delta\mu_{16}\equiv [\eta_{14}\mu_{15}]_{16} \bmod \Delta(\eta_{16}\bar{\nu}_{17}), \Delta(\eta_{16}\varepsilon_{17}).
		\end{aligned}\]
		
		We have $\Delta(\eta_{16}\varepsilon_{17})
		= [\eta^2_{14}]_{16}\varepsilon_{16}.$ 
		By Lemmas~\ref{m1:Dn4sg'} and~\ref{m1:zeta-lift-bracket}(2), together with \eqref{m1:et^2e}, we obtain the following calculation modulo $i_{15}\circ\pi_{16}(\mathit{SO}(14))\circ\varepsilon_{16}$.
		The other ambiguity in Lemma~\ref{m1:Dn4sg'} vanishes after composition with $\varepsilon_{16}$ because $\eta_{14}^2\varepsilon_{16}=0$.
		\[\begin{aligned}
			[\eta^2_{14}]\varepsilon_{16} &\in (-\{i_{15},\Delta(\iota_{14}),\eta^2_{13}\}) \circ \varepsilon_{16} \\
			&= i_{15}\circ \{\Delta(\iota_{14}), \eta^2_{13},\varepsilon_{15}\} \\
			&\ni [\zeta_{13}]_{15} \bmod i_{15}\circ \pi_{16}(\mathit{SO}(14)) \circ \varepsilon_{16}.
		\end{aligned}\]
		By the fact that $\pi_{16}(\mathit{SO}(14))=\mathrm{span}\{\Delta(\nu_{14}), [\iota_7]_{14}\mu_7\}$
		\cite[Proposition  4.1 and (4.11)]{Ka} and
		the relations $\mu_7\varepsilon_{16} = \eta_7\mu_8\sigma_{17}$ and
		$[\iota_7]_{10}\eta_7\mu_8\sigma_{17}=0$ (Proposition~\ref{m1:10}), we find
		${i_{15}}_*\pi_{16}(\mathit{SO}(14)) \circ \varepsilon_{16}=0.$  Applying $(i_{16})_*$ to the resulting equality in $\pi_{24}(\mathit{SO}(15)),$   we obtain
		\[\begin{aligned}
			&\Delta(\eta_{16}\varepsilon_{17}) = [\eta^2_{14}]_{16}\varepsilon_{16} = [\zeta_{13}]_{16}.
		\end{aligned}\]
		
		According to          the relation $\eta_{33}\bar{\nu}_{34}=\nu^3_{33}$  and
		\eqref{m1:Pn333}, we conclude that
		\[\begin{aligned}
			&J\Delta(\eta_{16}\bar{\nu}_{17})=P(\eta_{33}\bar{\nu}_{34}) =P(\nu^3_{33})\equiv \tilde{\varepsilon}_{16} \bmod \delta_{16},\bar{\sigma}'_{16}, \bar{\mu}_{16}\sigma_{33}.
		\end{aligned}\]
		We note from \cite[(3.9)]{MMO} that
		$H(\tilde{\varepsilon}_{10}) = \bar{\varepsilon}_{19}$
		and
		$2\tilde{\varepsilon}_{10} = \sigma_{10}\nu_{17}\kappa_{20}.$ 
		By using the EHP sequence and the relation $H(\tilde{\varepsilon}_{10})
		= \bar{\varepsilon}_{19} =HJ[\bar{\varepsilon}_9],$   we see that
		$
		J[\bar{\varepsilon}_9]\equiv \tilde{\varepsilon}_{10} \bmod \Sigma\pi^9_{33}.
		$
		Thus we obtain
		\[\begin{aligned}
			&J[\bar{\varepsilon}_9]_{12}\equiv \tilde{\varepsilon}_{12} \bmod \delta_{12}, \sigma_{12}\bar{\mu}_{19}, \bar{\sigma}'_{12}.
		\end{aligned}\]
		\noindent   Note that 
		$[\eta^2_{14}]\bar{\nu}_{16}\in{i_{15}}_*\{\Delta(\iota_{14}),\eta^2_{13},
		\bar{\nu}_{15}\}.$ 
		We obtain
		\[\begin{aligned}
			{p_{14}}_*\{\Delta(\iota_{14}),\eta^2_{13},\bar{\nu}_{15}\} &\subseteq\{2\iota_{13},\eta^2_{13},\bar{\nu}_{15}\} \\
			&\supseteq\{2\iota_{13},\eta_{13}, \nu^3_{14}\} \\
			&\supseteq\{2\iota_{13},\eta_{13},\nu_{14}\}\circ\nu^2_{18} \\
			&\ni 0 \m2\pi^{13}_{24}=\mathrm{span}\{2\zeta_{13}\}.
		\end{aligned}\]
		This implies
		\[\begin{aligned}
			&\{\Delta(\iota_{14}),\eta^2_{13},\bar{\nu}_{15}\}\subseteq\mathrm{span}\{2[\zeta_{13}]\}+{i_{14}}_*\pi_{24}(\mathit{SO}(13)).
		\end{aligned}\]
		Here ${i_{14}}_*\pi_{24}(\mathit{SO}(13))=\mathrm{span}\{[\bar{\varepsilon}_9]_{14},[\sigma^2_9]_{14}\eta_{23},[\iota_7]_{14}\bar{\mu}_7\}.$ 
		By Lemma~\ref{m1:zeta1}, $\Delta\zeta_{14}=\pm2[\zeta_{13}],$   so $(i_{15})_*(2[\zeta_{13}])=0.$ 
		After applying $(i_{15})_*$ and then $(i_{16})_*,$   the bracket therefore lies in $\mathrm{span}\{[\bar{\varepsilon}_9]_{16},[\sigma^2_9]_{16}\eta_{23},[\iota_7]_{16}\bar{\mu}_7\}.$ 
		Comparing the preceding $J$-images shows that the coefficient of $[\bar{\varepsilon}_9]_{16}$ is odd. Thus
		\[\begin{aligned}
			&\Delta(\eta_{16}\bar{\nu}_{17})\equiv [\bar{\varepsilon}_9]_{16} \bmod [\sigma^2_9]_{16}\eta_{23}, [\iota_7]_{16}\bar{\mu}_7.
		\end{aligned}\]
		\noindent   In light of Lemma~\ref{m1:[bep9]11eta}(2), we know that $[\bar{\varepsilon}_9]_{11}\eta_{24}
		\equiv 0\ \bmod\ [\iota_7]_{11}\eta_7\bar{\mu}_8.$ 
		Since $J([\sigma^2_9]_{16}\eta^2_{23})=\psi_{16}\eta^2_{39}=4\xi_{16}\sigma_{34}\ne 0$ \cite[Lemma 3.5]{MiM1}, we obtain
		\[\begin{aligned}
			&\Delta(\eta_{16}\bar{\nu}_{17})\equiv [\bar{\varepsilon}_9]_{16} \bmod [\iota_7]_{16}\bar{\mu}_7.
		\end{aligned}\]
		This completes the proof.
	\end{proof}
	\subsection{The Cases $18\leq n\leq21$}
	
	We turn our attention next to the following statement, which we shall prove here before proceeding to the next part.
	\begin{prop} \label{m1:18}
		The structure of $\pi_{24}(\mathit{SO}(18))$ is as follows.
		\[ \GroupLayout{\pi_{24}(\mathit{SO}(18)) = \mathrm{span}\{\Delta(\sigma_{18}),\; [\sigma^2_9]_{18}\eta_{23},\; [\iota_7]_{18}\bar{\mu}_7\}\cong \mathbb Z/16\oplus (\mathbb Z/2)^{2}.}{\begin{aligned}\pi_{24}(\mathit{SO}(18)) = \mathrm{span}\{&\Delta(\sigma_{18}),\; [\sigma^2_9]_{18}\eta_{23},\\
				&[\iota_7]_{18}\bar{\mu}_7\}\cong \mathbb Z/16\oplus (\mathbb Z/2)^{2}.\end{aligned}} \]
		Here
		$8\Delta(\sigma_{18}) = [\mu_{15}]_{18}.$ 
	\end{prop}
	\begin{proof}
		In the exact sequence $(24)_{17},$  
		it follows that $\operatorname{ker}\ (\Delta\colon  \pi^{17}_{24}\to \pi_{23}(\mathit{SO}(17))) = \mathrm{span}\{2\sigma_{17}\}$
		by the proof of Proposition~\ref{m1:R18}.
		Note that $\Delta(\sigma_{18})$ is a lift of $2\sigma_{17}.$ 
		Since $P(\sigma_{37})$ has order $16$ by \cite[(5.36)]{MMO},
		$\mathrm{ord}(\Delta(\sigma_{18}))=16.$ 
		
		Then we obtain
		\[\begin{aligned}
			&\Delta\varepsilon_{17} = [\eta_{15}]_{17}\varepsilon_{16} \quad\text{and}\quad \Delta\bar\nu_{17} = [\eta_{15}]_{17}\bar\nu_{16}.
		\end{aligned}\]
		Therefore we obtain
		$
		{i_{18}}_\ast \pi_{24}(\mathit{SO}(17))=\mathrm{span}\{[\mu_{15}]_{18}, [\sigma^2_9]_{18}\eta_{23},
		[\iota_7]_{18}\bar{\mu}_7\}\cong(\mathbb Z/2)^3.
		$
		
		In light of the definition of $[\mu_{15}],$  
		we deduce
		\begin{align*}
			[\mu_{15}]_{18} &\in i_{18}\circ i_{17}\circ \{[\eta_{15}], 2\iota_{16}, 8\sigma_{16}\}_1 \\
			&\subseteq i_{18}\circ\{[\eta_{15}]_{17}, 2\iota_{16}, 8\sigma_{16}\}_1 \\
			&= i_{18}\circ\{\Delta(\iota_{17}),2\iota_{16},8\sigma_{16}\}_1 \\
			&= -\{i_{18}, \Delta(\iota_{17}),2\iota_{16}\} \circ 8\sigma_{17} \ni 8\Delta(\sigma_{18}).
		\end{align*}
		The indeterminacy is $i_{18}\circ \pi_{17}(\mathit{SO}(17)) \circ 8\sigma_{17}=0$
		by \cite{Bo}. That is, $8\Delta(\sigma_{18}) = [\mu_{15}]_{18}.$ 
		This completes the proof.
	\end{proof}

	We note that the stable $J$-image in $\pi^S_{24}(S^0)$ is isomorphic to
	$\mathbb Z/2$ generated by $\sigma\bar{\mu} = \eta\bar{\rho}.$  Therefore, from the fact that
	$\pi^S_k(S^0) = 0 \;\;(k = 4, 5),$  
	by using the exact sequence $(24)_n$  $(n = 20, 21)$
	and \cite{Ke}, we obtain the following.
	
	\begin{lem} \label{m1:Ker1}
		$\pi_{24}(\mathit{SO}(n)) = \mathrm{span}\{[\iota_7]_n\bar{\mu}_{7}\}\cong \mathbb Z/2$ for $n = 20$ and $21.$ 
		\qed
	\end{lem}

	We find
	\[\begin{aligned}
		&\Delta(\nu^2_{19}) = [\nu_{15}]_{19}\nu^2_{18}.
	\end{aligned}\]
	According to          
	\cite[Theorem 1.1(b), (5.39)]{MMO}, we have
	$J\Delta(\nu^2_{19})=P(\nu^2_{39})\neq 0.$ 
	
	That is, $\Delta(\nu^2_{19})=[\nu_{15}]_{19}\nu^2_{18} \neq 0\in \pi_{24}(\mathit{SO}(19)).$ 
	Therefore, using the exact sequence $(24)_{19},$   we have the following.
	
	\begin{prop} \label{m1:19}
		The structure of $\pi_{24}(\mathit{SO}(19))$ is as follows.
		\[ \GroupLayout{\pi_{24}(\mathit{SO}(19)) = \mathrm{span}\{[\iota_7]_{19}\bar{\mu}_{7},\; [\nu_{15}]_{19}\nu^2_{18}\}\cong (\mathbb Z/2)^{2}.}{\begin{aligned}\pi_{24}(\mathit{SO}(19)) = \mathrm{span}\{&[\iota_7]_{19}\bar{\mu}_{7},\\
				&[\nu_{15}]_{19}\nu^2_{18}\}\cong (\mathbb Z/2)^{2}.\end{aligned}} \]
		Here $\Delta(\nu^2_{19})=[\nu_{15}]_{19}\nu^2_{18}.$ 
		\qed
	\end{prop}
	
	\subsection{The Case $n=22$ and the Stable Groups}
	The following result follows directly from the preceding calculations and overlaps with the results of \cite{Ke}.
	
	\begin{lem} \label{m1:Ker2}
		\begin{enumerate}
			\item $\pi_{24}(\mathit{SO}(22)) = \mathrm{span}\{\Delta(\nu_{22}),\; [\iota_7]_{22}\bar{\mu}_7\}\cong \mathbb Z/4\oplus \mathbb Z/2.$ 
			\item $\pi_{24}(\mathit{SO}(23)) = \mathrm{span}\{[\eta^2_{22}],\; [\iota_7]_{23}\bar{\mu}_7\}\cong (\mathbb Z/2)^{2}.$ 
			\item $\pi_{24}(\mathit{SO}(24)) = \mathrm{span}\{[\eta_{23}],\; [\eta^2_{22}]_{24},\; [\iota_7]_{24}\bar{\mu}_7\}\cong (\mathbb Z/2)^{3},$   where $\Delta\eta_{24}=[\eta^2_{22}]_{24}.$ 
			\item $\pi_{24}(\mathit{SO}(25))= \mathrm{span}\{[\eta_{23}]_{25},\; [\iota_7]_{25}\bar{\mu}_7\}\cong (\mathbb Z/2)^{2}.$ 
			\item $\pi_{24}(\mathit{SO}(n)) = \mathrm{span}\{[\iota_7]_n\bar{\mu}_7\}\cong \mathbb Z/2$ for $n\geq 26.$ 
		\end{enumerate}
		\qed
	\end{lem}
	
	\section{Spin Group Comparisons for the Cases $n=7,9$}
	Recall that
	$\pi_{*}(\mathit{Spin}(7)) \cong \pi_{*}(G_2) \oplus \pi^7_{*}.$ 
	
	By Proposition~\ref{m1:724}, we obtain
	\[\FitMath{\FitMath{\begin{aligned}
				\pi_{24}(\mathit{Spin}(7))={}\mathrm{span}\bigl\{& \langle P(\Sigma\theta)+\nu_6\kappa_9\rangle\eta_{23},\, \langle\eta^2_6\rangle\mu_8\sigma_{17},\, \langle\eta_6\bar\mu_7\rangle, \\
				&\sigma'\eta_{14}\mu_{15},\, \nu_7\kappa_{10},\,\bar\mu_7,\, \eta_7\mu_8\sigma_{17}\bigr\} \cong(\mathbb Z/2)^7.
	\end{aligned}}}\]
	\noindent   and
	$\pi_{24}(\mathit{Spin}(9))\cong(\mathbb Z/2)^6.$
	\section{Odd-Primary Components}\label{sec:odd24}
	\textbf{Exception to the standing convention.} In this section,
	$\pi_i(\mathit{SO}(n):\mathrm{odd})$ denotes the subgroup of all odd-order
	elements of the integral homotopy group, as defined in Chapter~\ref{ch:notation}.
	In the proof, the symbol $\pi_{*}(X)$  
	is short for $\pi_{*}(X)\otimes\z_{\mathcal{P}}$ where $\mathcal{P}$
	is the set of odd primes and $\z_{\mathcal{P}}$ denotes the $\mathcal{P}$-localization of $\z.$  A subscript $(p)$ denotes localization at the indicated odd prime
	$p,$   and $\operatorname{Tor}\pi_i(X_{(p)})\xyd \pi_i(X_{(p)})$ is the subgroup of elements whose orders are finite.
	The default convention at $2$ resumes in the next chapter.
	For an odd prime $p,$   use the comparison formulas in Section~\ref{sec:odd-method}. The nonzero odd-primary part is as follows.
	\begin{prop}
		For $n\geq3,$   we have
		\[
		\pi_{24}(\mathit{SO}(n):\mathrm{odd})\cong
		\begin{cases}
			\mathbb Z/33, & n=6,\\
			\mathbb Z/3, & n=7,8,12,22,\\
			\mathbb Z/15, & n=10,18,\\
			\mathbb Z/63, & n=14,\\
			0, & \text{otherwise}.
		\end{cases}
		\]
	\end{prop}
	\begin{proof}
		By using the relation \eqref{m1:R-Sp:odd}, we find
		$\operatorname{Tor}\pi_{24}(\bigl(\mathit{SO}(3)\bigr)_{(p)})=\operatorname{Tor}\pi_{24}(\bigl(\mathit{SO}(5)\bigr)_{(p)})=0$ from $\operatorname{Tor}\pi_{24}(\bigl(\mathit{Sp}(1)\bigr)_{(p)}) \cong \operatorname{Tor}\pi_{24}(\bigl(S^3\bigr)_{(p)}) =0$
		and $\operatorname{Tor}\pi_{24}(\bigl(\mathit{Sp}(2)\bigr)_{(p)}) =0$ \cite{F}.
		For the following integral input only, the subscript $\mathrm{int}$ denotes the spaces before
		localization. We observe that $\pi_{24}(\mathit{Sp}(2)_{\mathrm{int}}) \cong (\mathbb Z/2)^2$ and
		$\pi_{23}(\mathit{Sp}(2)_{\mathrm{int}}) \cong (\mathbb Z/2)^3$
		by \cite{F} and \cite{MT1}.
		By the following exact sequence
		\[\begin{aligned}
			&\pi_{24}(\mathit{Sp}(2)) \xrightarrow{i_\ast} \pi_{24}(\mathit{Sp}(3)) \xrightarrow{p_\ast} \pi_{24}(S^{11} ) \xrightarrow{\partial} \pi_{23}(\mathit{Sp}(2)),
		\end{aligned}\]
		we observe that $\operatorname{Tor}\pi_{24}(\bigl(\mathit{Sp}(3)\bigr)_{(p)}) \cong \operatorname{Tor}\pi_{24}(\bigl(S^{11}\bigr)_{(p)}).$ 
		Therefore, since $\pi_{24}(S^{11}_{\mathrm{int}})\cong \mathbb Z/6\oplus \mathbb Z/2,$  
		we see that $\operatorname{Tor}\pi_{24}(\bigl(\mathit{Sp}(3)\bigr)_{(3)}) \cong \mathbb Z/3$ and
		$\operatorname{Tor}\pi_{24}(\bigl(\mathit{Sp}(3)\bigr)_{(p)}) = 0$ for
		$p \ge 5.$ 
		Thus we obtain
		$$\operatorname{Tor}\pi_{24}(\bigl(\mathit{SO}(7)\bigr)_{(3)}) \cong \operatorname{Tor}\pi_{24}(\bigl(\mathit{Sp}(3)\bigr)_{(3)}) \cong \mathbb Z/3$$ and
		$\operatorname{Tor}\pi_{24}(\bigl(\mathit{SO}(7)\bigr)_{(p)})= 0$ for $p \ge 5.$ 
		Similarly, we have $\operatorname{Tor}\pi_{24}(\bigl(\mathit{SO}(2m+1)\bigr)_{(p)})=0$ for $m\ge 4$ by the relations
		$\operatorname{Tor}\pi_{24}(\bigl(\mathit{Sp}(4)\bigr)_{(p)}) =\operatorname{Tor}\pi_{24}(\bigl(\mathit{Sp}(5)\bigr)_{(p)}) =0$ from \cite{MT3} and $\operatorname{Tor}\pi_{24}(\bigl(\mathit{Sp}(m)\bigr)_{(p)}) =0 $ for $m \ge 6$
		from \cite{Bo}. Thus we obtain the asserted decomposition for $n=2m+1.$ 
		For $m\geq2,$   the odd-primary part of $\pi_{24}(S^{2m+1})$
		is isomorphic to $\mathbb Z/33$ for $m=2,$  
		$\mathbb Z/15$ for $m=4,8,$  
		$\mathbb Z/3$ for $m=5,10,$  
		and $\mathbb Z/63$ for $m=6,$   and vanishes otherwise.
		Combining these groups with \eqref{m1:R-Sp:even} and the results
		for $n=2m+1$ determines the corresponding groups for $n=2m+2.$ 
		This proves the assertion.
	\end{proof}
	\chapter{The Homotopy Groups $\pi_{25}(\mathit{SO}(n))$}\label{ch:degree25}
	
	The investigation of the homotopy groups in degree $25$ continues the general computational strategy employed in the preceding degrees, while requiring an additional ingredient at one of its initial stages. As before, the boundary homomorphisms provide the essential connection with the calculations already completed in the preceding degree, allowing the previously obtained results to serve as a foundation for the present determination. 
	
	The case $n=7$, however, is supplemented by an auxiliary computation associated with the exceptional Lie group $G_2$, which supplies the information required at this stage of the argument. Once this case has been established, the calculation proceeds through the successive inclusions of  $\mathit{SO}(n)$, tracing the relevant homotopy-theoretic information across the intermediate values of $n$. These successive steps form a coherent chain of calculations, leading from the initial cases toward the range in which stabilization occurs. Ultimately, the computation reaches the stable range, where the groups no longer vary with $n$. In particular, the homotopy group in degree $25$ is isomorphic to $\mathbb Z/2$ for every $n\geq27$, completing the determination in the stable range.
	
	As usual, we draw freely on the descriptions of the groups $\pi_{*}(S^{j})$
	presented in Appendix~\ref{tlqdzj}, as well as the relations
	governing compositions of generators in Toda's range collected
	in \^{O}guchi's table in Appendix~\ref{app:oguchi-relations}.
	
	\section{Main Theorem and Generators}
	
	The table retains the representatives chosen for small values of $n$ as their images are followed into groups with larger values of $n.$  This
	is useful when a later relation expresses a multiple of a new lift in terms of an inherited class.
	
	\begin{thm} \label{m1:thm3}
		\textup{All spaces in this theorem are localized at $2.$  The displayed groups are obtained from the integral homotopy groups by tensoring with $\mathbb Z_{(2)}.$  Free cyclic summands become $\mathbb Z_{(2)}$; for finite groups, this localization retains precisely the $2$-primary summands. The homotopy groups are given as follows.   }\par
		\begin{enumerate}[leftmargin=0pt,label={},itemindent=0pt]
			\item \textup{The Cases $3\leq n\leq6.$ }\par
			\smallskip\noindent\(
			\FitMath{\pi_{25}(\mathit{SO}(3)) = \mathrm{span}\{[\eta_2]\varepsilon_3\kappa_{11}\}\cong \mathbb Z/2.}
			\)\par
			\smallskip\par\smallskip\noindent\(
			\GroupLayout{\pi_{25}(\mathit{SO}(4)) = \mathrm{span}\{[\eta_2]_4\varepsilon_3\kappa_{11},\; [\iota_3]\varepsilon_3\kappa_{11}\}\cong (\mathbb Z/2)^{2}.}{\begin{aligned}\pi_{25}(\mathit{SO}(4)) = \mathrm{span}\{&[\eta_2]_4\varepsilon_3\kappa_{11},\\
					&[\iota_3]\varepsilon_3\kappa_{11}\}\cong (\mathbb Z/2)^{2}.\end{aligned}}
			\)\par
			\smallskip\par\smallskip\noindent\(
			\GroupLayout{\pi_{25}(\mathit{SO}(5)) = \mathrm{span}\{[\nu_4\zeta_7]\sigma_{18},\; [\iota_3]_5\varepsilon_3\kappa_{11}\}\cong \mathbb Z/8\oplus \mathbb Z/2,}{\begin{aligned}\pi_{25}(\mathit{SO}(5)) = \mathrm{span}\{&[\nu_4\zeta_7]\sigma_{18},\\
					&[\iota_3]_5\varepsilon_3\kappa_{11}\}\cong \mathbb Z/8\oplus \mathbb Z/2,\end{aligned}}
			\)\par
			\nobreak
			{\raggedleft
				\noindent where $[\nu_4\sigma'\zeta_{14}]\equiv [\nu_4\zeta_7]\sigma_{18} \bmod [\iota_3]_5\varepsilon_3\kappa_{11},\; 2[\nu_4\zeta_7]\sigma_{18}.$\par}
			
			\par\smallskip\noindent\(
			\GroupLayout{\pi_{25}(\mathit{SO}(6)) = \mathrm{span}\{[\nu_4\zeta_7]_6\sigma_{18},\; [ \nu_5]\nu_8\kappa_{11},\; [ \nu_5]\bar{\mu}_8,\; [ \nu_5]\eta_8\mu_9\sigma_{18}\}\cong \mathbb Z/8\oplus (\mathbb Z/2)^{3},}{\begin{aligned}\pi_{25}(\mathit{SO}(6)) = \mathrm{span}\{&[\nu_4\zeta_7]_6\sigma_{18},\; [ \nu_5]\nu_8\kappa_{11},\; [ \nu_5]\bar{\mu}_8,\\
					&[ \nu_5]\eta_8\mu_9\sigma_{18}\}\cong \mathbb Z/8\oplus (\mathbb Z/2)^{3},\end{aligned}}
			\)\par
			\nobreak
			{\raggedleft
				\noindent where $\Delta(\bar\kappa'
				)=[ \nu_5]\nu_8\kappa_{11}$
				and $\Delta(P(\rho_{13}))=[ \nu_5]\eta_8\mu_9\sigma_{18}+4[\nu_4\zeta_7]_6\sigma_{18}.$\par}
			
			\item \textup{The Cases $7\leq n\leq9.$ }\par
			\smallskip\noindent\(
			\GroupLayout{\pi_{25}(\mathit{SO}(7)) = \mathrm{span}\{[\bar{\sigma}_6],\; [\eta_6]\eta_7\bar{\mu}_8,\; [\nu_4\zeta_7]_7\sigma_{18},\; [ \nu_5]_7\bar{\mu}_8\}\cong \mathbb Z/32\oplus \mathbb Z/2\oplus \mathbb Z/8\oplus \mathbb Z/2,}{\begin{aligned}\pi_{25}(\mathit{SO}(7)) = \mathrm{span}\{&[\bar{\sigma}_6],\; [\eta_6]\eta_7\bar{\mu}_8,\; [\nu_4\zeta_7]_7\sigma_{18},\\
					&[ \nu_5]_7\bar{\mu}_8\}\cong \mathbb Z/32\oplus \mathbb Z/2\oplus \mathbb Z/8\oplus \mathbb Z/2,\end{aligned}}
			\)\par
			\nobreak
			{\raggedleft
				\noindent where $[\eta_6]\zeta_7\sigma_{18}=[\nu_4\zeta_7]_7\sigma_{18}$
				and  $4[\bar{\sigma}_6]=x[[2\e_{6},\e_{6}]]\sigma^2_{11}$ for some odd $x.$ 
				Moreover, we have $16[\bar{\sigma}_6]=4[[2\e_{6},\e_{6}]]\sigma^2_{11}\equiv
				[P(\Sigma\theta)+\nu_6\kappa_9]\eta^2_{23} \bmod 4[\nu_4\sigma'\zeta_{14}]_7.$ \par}
			
			\par\smallskip\noindent\(
			\GroupLayout{\pi_{25}(\mathit{SO}(8)) \;= {} \mathrm{span}\{[\bar{\sigma}_6]_8,\; [\eta_6]_8\eta_7\bar{\mu}_8,\; [\nu_4\zeta_7]_8\sigma_{18},\; [ \nu_5]_8\bar{\mu}_8,\; [\iota_7]\zeta_7\sigma_{18},\; [\iota_7]\eta_7\bar{\mu}_8\}\cong \mathbb Z/32\oplus \mathbb Z/2\oplus \mathbb Z/8\oplus \mathbb Z/2\oplus \mathbb Z/8\oplus \mathbb Z/2.}{\begin{aligned}\pi_{25}(\mathit{SO}(8)) \;= {} \mathrm{span}\{&[\bar{\sigma}_6]_8,\; [\eta_6]_8\eta_7\bar{\mu}_8,\; [\nu_4\zeta_7]_8\sigma_{18},\; [ \nu_5]_8\bar{\mu}_8,\\
					&[\iota_7]\zeta_7\sigma_{18},\; [\iota_7]\eta_7\bar{\mu}_8\}\cong \mathbb Z/32\oplus \mathbb Z/2\oplus \mathbb Z/8\oplus \mathbb Z/2\oplus \mathbb Z/8\oplus \mathbb Z/2.\end{aligned}}
			\)\par
			\smallskip\par\smallskip\noindent\(
			\GroupLayout{\pi_{25}(\mathit{SO}(9)) = \mathrm{span}\{[\bar{\sigma}_6]_9,\; [ \nu_5]_9\bar{\mu}_8,\; [\iota_7]_9\eta_7\bar{\mu}_8\}\cong \mathbb Z/32\oplus (\mathbb Z/2)^{2},}{\begin{aligned}\pi_{25}(\mathit{SO}(9)) = \mathrm{span}\{&[\bar{\sigma}_6]_9,\; [ \nu_5]_9\bar{\mu}_8,\\
					&[\iota_7]_9\eta_7\bar{\mu}_8\}\cong \mathbb Z/32\oplus (\mathbb Z/2)^{2},\end{aligned}}
			\)\par
			\nobreak
			{\raggedleft
				\noindent where
				$16[\bar\sigma_6]_9 = 2[\iota_7]_9\zeta_7\sigma_{18}=[\nu_4\zeta_7]_9\sigma_{18}
				=[\eta_6]_9\zeta_7\sigma_{18}=[P(\Sigma\theta)+\nu_6\kappa_9]_9\eta^2_{23}.$\par}
			
			\item \textup{The Cases $10\leq n\leq14.$ }\par
			\smallskip\noindent\(
			\GroupLayout{\pi_{25}(\mathit{SO}(10)) = \mathrm{span}\{[\bar{\sigma}_6]_{10},\; [\iota_7]_{10}\eta_7\bar{\mu}_8,\; [\sigma^2_9]\eta^2_{23},\; [\sigma_9\eta_{16}\varepsilon_{17}]\}\cong \mathbb Z/32\oplus (\mathbb Z/2)^{3},}{\begin{aligned}\pi_{25}(\mathit{SO}(10)) = \mathrm{span}\{&[\bar{\sigma}_6]_{10},\; [\iota_7]_{10}\eta_7\bar{\mu}_8,\\
					&[\sigma^2_9]\eta^2_{23},\; [\sigma_9\eta_{16}\varepsilon_{17}]\}\cong \mathbb Z/32\oplus (\mathbb Z/2)^{3},\end{aligned}}
			\)\par
			\nobreak
			{\raggedleft
				\noindent where $[\bar{\varepsilon}_9]\eta_{24}\equiv 0 \bmod  [\iota_7]_{10}\sigma'\zeta_{14}, [\iota_7]_{10}\eta_7\bar{\mu}_8.$\par}
			
			\par\smallskip\noindent\(
			\GroupLayout{\pi_{25}(\mathit{SO}(11)) = \mathrm{span}\{[\bar{\sigma}_6]_{11},\; [[\e_{10},\nu_{10}]]\nu_{22},\; [\iota_7]_{11}\eta_7\bar\mu_8,\; [\bar{\varepsilon}_{10}]\}\cong (\mathbb Z/4)^{2}\oplus (\mathbb Z/2)^{2},}{\begin{aligned}\pi_{25}(\mathit{SO}(11)) = \mathrm{span}\{&[\bar{\sigma}_6]_{11},\; [[\e_{10},\nu_{10}]]\nu_{22},\; [\iota_7]_{11}\eta_7\bar\mu_8,\\
					&[\bar{\varepsilon}_{10}]\}\cong (\mathbb Z/4)^{2}\oplus (\mathbb Z/2)^{2}.\end{aligned}}
			\)\par
			\nobreak
			
			\par\smallskip\noindent\(
			\GroupLayout{\pi_{25}(\mathit{SO}(12)) \;= {} \mathrm{span}\{[\sigma^2_{11}],\; [[\e_{10},\nu_{10}]]_{12}\nu_{22},\; [\bar{\sigma}_6]_{12},\; [\iota_7]_{12}\eta_7\bar\mu_8,\; [\bar{\varepsilon}_{10}]_{12}\}\cong \mathbb Z/16\oplus (\mathbb Z/4)^{2}\oplus (\mathbb Z/2)^{2},}{\begin{aligned}\pi_{25}(\mathit{SO}(12)) \;= {} \mathrm{span}\{&[\sigma^2_{11}],\; [[\e_{10},\nu_{10}]]_{12}\nu_{22},\; [\bar{\sigma}_6]_{12},\\
					&[\iota_7]_{12}\eta_7\bar\mu_8,\; [\bar{\varepsilon}_{10}]_{12}\}\cong \mathbb Z/16\oplus (\mathbb Z/4)^{2}\oplus (\mathbb Z/2)^{2},\end{aligned}}
			\)\par
			\nobreak
			{\raggedleft
				\noindent where $\Delta(P(\nu_{25}))=[[\e_{10},\nu_{10}]]_{12}\nu_{22}, \Delta(\kappa_{12})=[\bar{\varepsilon}_{10}]_{12}$ and
				$\Delta (\sigma^2_{12})\equiv2[\sigma^2_{11}] +  [\bar\sigma_{6}]_{12} \bmod 2[\bar\sigma_{6}]_{12}, \Delta(P(\nu_{25})).$\par}
			
			\par\smallskip\noindent\(
			\GroupLayout{\pi_{25}(\mathit{SO}(13+n)) = \mathrm{span}\{[\iota_7]_{13+n}\eta_7\bar{\mu}_8,\; [\sigma^2_{11}]_{13+n}\}\cong \mathbb Z/2\oplus \mathbb Z/8,\; (n=0\;\text{or}\; 1).}{\begin{aligned}\pi_{25}(\mathit{SO}(13+n)) = \mathrm{span}\{&[\iota_7]_{13+n}\eta_7\bar{\mu}_8,\\
					&[\sigma^2_{11}]_{13+n}\}\cong \mathbb Z/2\oplus \mathbb Z/8,\; (n=0\;\text{or}\; 1).\end{aligned}}
			\)\par
			
			\item \textup{The Cases $n\geq15.$ }\par
			\smallskip\noindent\(
			\GroupLayout{\pi_{25}(\mathit{SO}(15)) = \mathrm{span}\{[\iota_7]_{15}\eta_7\bar{\mu}_8,\; [\sigma^2_{11}]_{15},\; [\eta^2_{14}]\mu_{16}\}\cong \mathbb Z/2\oplus \mathbb Z/8\oplus \mathbb Z/2.}{\begin{aligned}\pi_{25}(\mathit{SO}(15)) = \mathrm{span}\{&[\iota_7]_{15}\eta_7\bar{\mu}_8,\; [\sigma^2_{11}]_{15},\\
					&[\eta^2_{14}]\mu_{16}\}\cong \mathbb Z/2\oplus \mathbb Z/8\oplus \mathbb Z/2.\end{aligned}}
			\)\par
			\smallskip\par\smallskip\noindent\(
			\GroupLayout{\pi_{25}(\mathit{SO}(16)) = \mathrm{span}\{[\iota_7]_{16}\eta_7\bar{\mu}_8,\; [\sigma^2_{11}]_{16},\; [\eta^2_{14}]_{16}\mu_{16},\; [\eta_{15}]\mu_{16}\}\cong \mathbb Z/2\oplus \mathbb Z/8\oplus (\mathbb Z/2)^{2},}{\begin{aligned}\pi_{25}(\mathit{SO}(16)) = \mathrm{span}\{&[\iota_7]_{16}\eta_7\bar{\mu}_8,\; [\sigma^2_{11}]_{16},\; [\eta^2_{14}]_{16}\mu_{16},\\
					&[\eta_{15}]\mu_{16}\}\cong \mathbb Z/2\oplus \mathbb Z/8\oplus (\mathbb Z/2)^{2},\end{aligned}}
			\)\par
			\nobreak
			{\raggedleft
				\noindent where $\Delta(\eta_{16}\mu_{17})=
				[\eta^2_{14}]_{16}\mu_{16}$ and $[\nu_{15}]\sigma_{18} \equiv[\sigma^2_{11}]_{16}\mod2[\sigma^2_{11}]_{16}.$\par}
			
			\par\smallskip\noindent\(
			\GroupLayout{\pi_{25}(\mathit{SO}(17)) = \mathrm{span}\{[\iota_7]_{17}\eta_7\bar{\mu}_8,\; [\sigma^2_{11}]_{17},\; [\eta_{15}]_{17}\mu_{16}\}\cong \mathbb Z/2\oplus \mathbb Z/8\oplus \mathbb Z/2,}{\begin{aligned}\pi_{25}(\mathit{SO}(17)) = \mathrm{span}\{&[\iota_7]_{17}\eta_7\bar{\mu}_8,\; [\sigma^2_{11}]_{17},\\
					&[\eta_{15}]_{17}\mu_{16}\}\cong \mathbb Z/2\oplus \mathbb Z/8\oplus \mathbb Z/2,\end{aligned}}
			\)\par
			\nobreak
			{\raggedleft
				\noindent where
				$\Delta(\mu_{17})=[\eta_{15}]_{17}\mu_{16}$ and $\Delta(\eta^2_{17}\sigma_{19})=4[\nu_{15}]_{17}\sigma_{18}=4[\sigma^2_{11}]_{17}.$\par}
			
			\par\smallskip\noindent\(
			\GroupLayout{\pi_{25}(\mathit{SO}(18)) = \mathrm{span}\{[\iota_7]_{18}\eta_7\bar{\mu}_8,\; [\sigma^2_{11}]_{18}\}\cong \mathbb Z/2\oplus \mathbb Z/4,}{\begin{aligned}\pi_{25}(\mathit{SO}(18)) = \mathrm{span}\{&[\iota_7]_{18}\eta_7\bar{\mu}_8,\\
					&[\sigma^2_{11}]_{18}\}\cong \mathbb Z/2\oplus \mathbb Z/4,\end{aligned}}
			\)\par
			{\raggedleft where $\Delta(\eta_{18}\sigma_{19})=2[\nu_{15}]_{18}\sigma_{18}=2[\sigma^2_{11}]_{18}.$ \par}

			\par\smallskip\noindent\(
			\GroupLayout{\pi_{25}(\mathit{SO}(19)) = \mathrm{span}\{[\iota_7]_{19}\eta_7\bar{\mu}_8,\; [\sigma^2_{11}]_{19}\}\cong (\mathbb Z/2)^{2}, \text{where}\;\Delta(\sigma_{19})=[\sigma^2_{11}]_{19} = [\nu_{15}]_{19}\sigma_{18}.}{\begin{aligned}\pi_{25}(\mathit{SO}(19)) = \mathrm{span}\{&[\iota_7]_{19}\eta_7\bar{\mu}_8,\\
					&[\sigma^2_{11}]_{19}\}\cong (\mathbb Z/2)^{2}, \text{where}\;\Delta(\sigma_{19})=[\sigma^2_{11}]_{19} = [\nu_{15}]_{19}\sigma_{18}.\end{aligned}}
			\)\par
			\smallskip\par\smallskip\noindent\(
			\FitMath{\pi_{25}(\mathit{SO}(20))=\mathrm{span}\{[\iota_7]_{20}\eta_7\bar{\mu}_8\}\cong\mathbb Z/2.}
			\)\par
			\smallskip\par\smallskip\noindent\(
			\GroupLayout{\pi_{25}(\mathit{SO}(23)) = \mathrm{span}\{[\iota_7]_{23}\eta_7\bar{\mu}_8,\; [\eta^2_{22}]\eta_{24}\}\cong (\mathbb Z/2)^{2}.}{\begin{aligned}\pi_{25}(\mathit{SO}(23)) = \mathrm{span}\{&[\iota_7]_{23}\eta_7\bar{\mu}_8,\\
					&[\eta^2_{22}]\eta_{24}\}\cong (\mathbb Z/2)^{2}.\end{aligned}}
			\)\par
			\smallskip\par\smallskip\noindent\(
			\GroupLayout{\pi_{25}(\mathit{SO}(24)) = \mathrm{span}\{[\iota_7]_{24}\eta_7\bar{\mu}_8,\; [\eta_{23}]_{24}\eta_{24},\; [\eta^2_{22}]_{24}\eta_{24}\}\cong (\mathbb Z/2)^{3},\;\;\text{where\;   }\Delta(\eta^2_{24})=[\eta^2_{22}]_{24}\eta_{24}.}{\begin{aligned}\pi_{25}(\mathit{SO}(24)) = \mathrm{span}\{&[\iota_7]_{24}\eta_7\bar{\mu}_8,\; [\eta_{23}]_{24}\eta_{24},\\
					&[\eta^2_{22}]_{24}\eta_{24}\}\cong (\mathbb Z/2)^{3},\;\;\text{where\;   }\Delta(\eta^2_{24})=[\eta^2_{22}]_{24}\eta_{24}.\end{aligned}}
			\)\par

			\par\smallskip\noindent\(
			\GroupLayout{\pi_{25}(\mathit{SO}(25)) = \mathrm{span}\{[\iota_7]_{25}\eta_7\bar{\mu}_8,\; [\eta_{23}]_{25}\eta_{24}\}\cong (\mathbb Z/2)^{2}, \text{where}\;\Delta(\eta_{25})=[\eta_{23}]_{25}\eta_{24}.}{\begin{aligned}\pi_{25}(\mathit{SO}(25)) = \mathrm{span}\{&[\iota_7]_{25}\eta_7\bar{\mu}_8,\\
					&[\eta_{23}]_{25}\eta_{24}\}\cong (\mathbb Z/2)^{2}, \text{where}\;\Delta(\eta_{25})=[\eta_{23}]_{25}\eta_{24}.\end{aligned}}
			\)\par
			\smallskip\par\smallskip\noindent\(
			\GroupLayout{\pi_{25}(\mathit{SO}(26)) = \mathrm{span}\{\Delta(\iota_{26}),\; [\iota_7]_{26}\eta_7\bar{\mu}_8\}\cong \mathbb Z_{(2)}\oplus \mathbb Z/2.}{\begin{aligned}\pi_{25}(\mathit{SO}(26)) = \mathrm{span}\{&\Delta(\iota_{26}),\\
					&[\iota_7]_{26}\eta_7\bar{\mu}_8\}\cong \mathbb Z_{(2)}\oplus \mathbb Z/2.\end{aligned}}
			\)\par
			\smallskip\par\smallskip\noindent\(
			\FitMath{\pi_{25}(\mathit{SO}(21+n)) = \mathrm{span}\{[\iota_7]_{21+n}\eta_7\bar{\mu}_8\}\cong \mathbb Z/2, \;(n=0, 1 \;\text{or} \;n\geq 6).}
			\)\par
		\end{enumerate}
	\end{thm}
	\section{The Cases $3\leq n\leq6$}
	Next, we shall determine $\pi_{25}(\mathit{SO}(n))$ for $n\leq 6.$ 
	
	\subsection{The Cases $3\leq n\leq5$}
	By
	\cite[Table 2]{Ka}, we know
	\begin{gather}
		\begin{aligned}
			&\pi_{11}(\mathit{SO}(5))={}\frac{\mathbb Z}{2}\{[\iota_3]_5\varepsilon_3\},\ \pi_{11}(\mathit{SO}(6))=\frac{\mathbb Z}{4}\{[ \nu_5]\nu_8\},\ \mbox{where}\ 2[ \nu_5]\nu_8=[\iota_3]_6\varepsilon_3,\label{m1:3epis-1}
		\end{aligned}
		\\[3pt]
		\begin{aligned}
			&\pi_{11}(\mathit{SO}(7))={}\frac{\mathbb Z}{2}\{[ \nu_5]_7\nu_8\}\oplus\mathbb Z_{(2)}\{[[2\e_{6},\e_{6}]]\}.\label{m1:3epis-2}
		\end{aligned}
	\end{gather}
	
	Since $\pi^3_{25}=\frac{\mathbb Z}{2}\{\varepsilon_3\kappa_{11}\},$   we have the following:
	$$   \pi_{25}(\mathit{SO}(3))={}\frac{\mathbb Z}{2}\{[\eta_2]\varepsilon_3\kappa_{11}\},\;\;     \pi_{25}(\mathit{SO}(4))={}\frac{\mathbb Z}{2}\{[\eta_2]_4\varepsilon_3\kappa_{11}\}\oplus\frac{\mathbb Z}{2}\{[\iota_3]\varepsilon_3\kappa_{11}\}.    $$
	\noindent   In the exact sequence $(25)_4$:
	\[\begin{aligned}
		&\pi^4_{26}\xrightarrow{\Delta}\pi_{25}(\mathit{SO}(4))\xrightarrow{i_*}\pi_{25}(\mathit{SO}(5))\xrightarrow{p_*}\pi^4_{25} \xrightarrow{\Delta}\pi_{24}(\mathit{SO}(4)),
	\end{aligned}\]
	recall that $$\begin{aligned}
		&\operatorname{ker}(\Delta\colon \pi^4_{25}\to \pi_{24}(\mathit{SO}(4)))=\frac{\mathbb Z}{8}\{\nu_4\zeta_7\sigma_{18}\}
	\end{aligned}$$
	from the proof of \eqref{m1:524} and
	$\pi^4_{26}=\mathrm{span}\{\nu_4\bar{\zeta}_7,\varepsilon_4\kappa_{12},
	\nu_4\bar{\sigma}_7\}\cong\mathbb Z/8\oplus(\mathbb Z/2)^2.$ 
	Noting that
	$2\iota_3\circ \varepsilon_3=2\varepsilon_3=0$ and that
	$\nu'\bar{\zeta}_6=\nu'\bar{\sigma}_6=0$ \cite[p.~313--4]{Mi1}, we see that
	$\Delta(\varepsilon_4\kappa_{12})
	=[\eta_2]_4\varepsilon_3\kappa_{11},$  
	$\Delta(\nu_4\bar{\zeta}_7)=0$ and
	$\Delta(\nu_4\bar{\sigma}_7)=0.$  Hence, we obtain
	\begin{gather}
		\begin{aligned}
			&\pi_{25}(\mathit{SO}(5))={}\frac{\mathbb Z}{8}\{[\nu_4\zeta_7]\sigma_{18}\}\oplus\frac{\mathbb Z}{2}\{[\iota_3]_5\varepsilon_3\kappa_{11}\}.\label{m1:255}
		\end{aligned}
	\end{gather}
	
	\subsection{The Case $n=6$}
	In the exact sequence $(25)_5$:
	\[\begin{aligned}
		&\pi^5_{26}\xrightarrow{\Delta}\pi_{25}(\mathit{SO}(5))\xrightarrow{i_*}\pi_{25}(\mathit{SO}(6))\xrightarrow{p_*}\pi^5_{25} \xrightarrow{\Delta}\pi_{24}(\mathit{SO}(5)),
	\end{aligned}\]
	it follows that $\operatorname{ker}\ (\Delta\colon \pi^5_{25}\to \pi_{24}(\mathit{SO}(5)))
	=\mathrm{span}\{\nu_5\nu_8\kappa_{11},\nu_5\bar{\mu}_8,\nu_5\eta_8\mu_9\sigma_{18}\}\cong(\mathbb Z/2)^3$
	from the proof of Proposition~\ref{m1:R61}
	and $\pi^5_{26}=\mathrm{span}\{\alpha,\nu_5\eta_8\bar{\mu}_9\}\cong(\mathbb Z/2)^2,$   where
	$\Sigma\alpha=\eta_6\bar{\kappa}_7$ \cite[p.~305--307]{Mi1}.
	We observe that $\Delta(\nu_5\eta_8\bar{\mu}_9)=0.$  In light of \eqref{m1:3epis-1},
	we find
	$[\iota_3]_6\varepsilon_3\kappa_{11}=2[ \nu_5]\nu_8\kappa_{11}=0.$  This implies the following.
	\begin{lem}\label{g1016:lem:delta-alpha}
		For the   element   $\alpha\in\pi_{26}^5$ fixed above, with $\Sigma\alpha=\eta_6\bar\kappa_7,$   we have
		$\Delta\alpha=[\iota_3]_5\varepsilon_3\kappa_{11}.$ \qed
	\end{lem}
	We therefore have
	\begin{equation}\label{m1:i4,6R425}
		{i_5}_*\pi_{25}(\mathit{SO}(4))={}\frac{\mathbb Z}{2}\{[\iota_3]_5\varepsilon_3\kappa_{11}\},\ {i_{4,6}}\pi_{25}(\mathit{SO}(4))=0,
	\end{equation}
	\begin{equation}\label{m1:R256}
		\pi_{25}(\mathit{SO}(6))={}\mathrm{span}\{[\nu_4\zeta_7]_6\sigma_{18},[ \nu_5]\nu_8\kappa_{11},[ \nu_5]\bar{\mu}_8,[ \nu_5]\eta_8\mu_9\sigma_{18}\}\cong\mathbb Z/8\oplus(\mathbb Z/2)^3.
	\end{equation}

	\section{The Auxiliary Degree 25 Calculation for $G_2$}\label{sec:degree25-g2}
	
	The passage to the case $n=7$ requires degree-$25$ information about $G_2.$  We obtain it through the
	fibration with fiber $\mathit{SU}(3)$ and base $S^6,$   retaining the connecting images and multiples of
	lifts needed in the subsequent comparison with special orthogonal groups.
	
	\subsection{Homotopy Groups of Spheres and Unitary Groups}
	We need the group structure of $\pi_{25}(G_2).$ 
	We list the following groups from \cite[Theorem 4.1]{MT2} and \cite[Theorem 1.1]{F}:
	\[\begin{gathered}
		\begin{aligned}
			&\pi_{11}(\mathit{SU}(3))={}\frac{\mathbb Z}{4}\{[\nu^2_5]\},
		\end{aligned}\\[3pt]
		\begin{aligned}
			&\pi_{23}(\mathit{SU}(3))={}\frac{\mathbb Z}{4}\{[2\iota_5]\zeta_5\sigma_{16}\} \oplus\frac{\mathbb Z}{2}\{[\nu_5\bar{\varepsilon}_8]\},
		\end{aligned}\\[3pt]
		\begin{aligned}
			&\pi_{24}(\mathit{SU}(3))={}\frac{\mathbb Z}{4}\{[2\iota_5]\bar{\zeta}_5\},
		\end{aligned}\\[3pt]
		\begin{aligned}
			&\pi_{25}(\mathit{SU}(3))={}\frac{\mathbb Z}{2}\{[\nu^2_5]\kappa_{11}\}\oplus\frac{\mathbb Z}{2}\{[\nu_5\eta_8\mu_9]\sigma_{18}\}.
		\end{aligned}
	\end{gathered}\]
	
	By \eqref{m2:Y} and
	\eqref{m1:3epis-2}, we obtain
	\begin{equation}
		\begin{aligned}
			&\mathcal{T}_{*}\langle [2\e_{6},\e_{6}]\rangle\equiv[[2\e_{6},\e_{6}]]\ \bmod\ [ \nu_5]_7\nu_8.\label{m1:hl2Pi13}
		\end{aligned}
	\end{equation}
	\subsection{The Canonical Inclusion and Connecting Images}
	
	The unitary generators must be compared with the chosen lifts to special orthogonal groups. Agreement after
	projection leaves a possible fiber term, so the next lemma records the correction terms explicitly.
	
	At this point, it is convenient to gather
	several related observations into a  lemma.
	\begin{lem}\label{m1:r3n5}
		\begin{enumerate}
			\item
			${r_3}_\ast[\nu^2_5]=\pm [ \nu_5]\nu_8.$ 
			\item
			${r_3}_\ast[\nu_5\eta_8\mu_9]=[ \nu_5]\eta_8\mu_9+4[\nu_4\zeta_7]_6 \ \bmod\ 2[ \nu_5]\sigma_8\nu_{15};$ moreover,
			$${r_3}_\ast([\nu_5\eta_8\mu_9]\sigma_{18})=[ \nu_5]\eta_8\mu_9\sigma_{18}
			+4[\nu_4\zeta_7]_6\sigma_{18}.$$
		\end{enumerate}
	\end{lem}
	\begin{proof}
		\begin{enumerate}[label=\textup{(\arabic*)},leftmargin=2em]
			\item By the equality ${p_6}_\ast({r_3}_\ast[\nu^2_5])=\nu^2_5={p_6}_\ast([ \nu_5]\nu_8)$ and \eqref{m1:3epis-1},
			we obtain the relation (1).
			
			\item Since
			$
			{p_6}_\ast ({r_3}_\ast[\nu_5\eta_8\mu_9]) = \nu_5\eta_8\mu_9 = {p_6}_\ast ([ \nu_5]\eta_8\mu_9),
			$
			we deduce
			\[\begin{aligned}
				&{r_3}_\ast([\nu_5\eta_8\mu_9]) - [ \nu_5]\eta_8\mu_9 \in {i_6}_\ast \pi_{18}(\mathit{SO}(5)) = \frac{\mathbb Z}{8}\{[\nu_4\zeta_7]_6\} \oplus \frac{\mathbb Z}{2}\{2[ \nu_5]\sigma_8\nu_{15}\}
			\end{aligned}\]
			\cite[Theorem 2(i), Lemma 3.4 (ii)]{KM1}.
			Since $[\nu_5\eta_8\mu_9]$ has order 2 \cite[Theorem 4.1]{MT2}, we find
			\[\begin{aligned}
				&{r_3}_\ast[\nu_5\eta_8\mu_9]=[ \nu_5]\eta_8\mu_9+4a[\nu_4\zeta_7]_6+2b[ \nu_5]\sigma_8\nu_{15}\ \mbox{for}\ a, b\in\{0, 1\}.
			\end{aligned}\]
			
			We note from \cite[p. 160]{Mi2} that
			\begin{equation}
				\begin{aligned}
					&8\langle [2\e_{6},\e_{6}]\rangle\sigma_{11} = {i_G}_{\ast}[\nu_5\eta_8\mu_9].\label{m1:8<Pi13>s}
				\end{aligned}
			\end{equation}
			Therefore,        according to        \eqref{m2:Y}, \eqref{m1:hl2Pi13} and the fact that $2[ \nu_5]_7\sigma_8\nu_{15}=0,$   we conclude that
			\[\begin{aligned}
				&8[[2\e_{6},\e_{6}]]\sigma_{11}= 8\mathcal{T}_{*}\langle [2\e_{6},\e_{6}]\rangle\sigma_{11}=\mathcal{T}_{*}{i_G}_{\ast}[\nu_5\eta_8\mu_9] ={i_7}_*{r_3}_*[\nu_5\eta_8\mu_9] =[ \nu_5]_7\eta_8\mu_9+4a[\nu_4\zeta_7]_7.
			\end{aligned}\]
			From \cite[Lemma 3.7(ii)]{KM1}, we obtain $a=1.$  This implies
			the first relation in (2).
			Hence, by the relation $\nu_{15}\sigma_{18}=0,$   we have the second relation in (2). This completes the proof.
		\end{enumerate}
	\end{proof}
	
	By \eqref{m1:8<Pi13>s}, we obtain
	\[\begin{aligned}
		&{i_G}_*([\nu_5\eta_8\mu_9]\sigma_{18})= 8\langle[2\e_{6},\e_{6}]\rangle\sigma^2_{11}.
	\end{aligned}\]
	\noindent  
	Observe that $\Sigma\theta=P(\iota_{27})$ \cite[(7.30)]{T} and
	$8\sigma^2_{13}=P(\eta^2_{27})$ \cite[(10.10)]{T}. Therefore, we obtain
	\[\begin{aligned}
		&P(8\sigma^2_{13}) =P(P(\eta^2_{27}))=P(P(\iota_{27}))\eta^2_{23} =P(\Sigma\theta)\eta^2_{23}=(P(\Sigma\theta)+\nu_6\kappa_9)\eta^2_{23}.
	\end{aligned}\]
	\noindent
	Hence, we deduce
	\begin{gather}\label{m1:42P13}
		4\langle[2\e_{6},\e_{6}]\rangle\sigma^2_{11}\equiv\langle P(\Sigma\theta)+\nu_6\kappa_9\rangle \eta^2_{23}\ \bmod\
		{i_G}_*\pi_{25}(\mathit{SU}(3))
		\shortintertext{and}
		\label{m1:diamond}
		{i_G}_*([\nu_5\eta_8\mu_9]\sigma_{18})=8\langle[2\e_{6},\e_{6}]\rangle\sigma^2_{11}=0.
	\end{gather}
	\noindent   
	The group is
	\[\begin{aligned}
		&\pi^6_{26}={}\frac{\mathbb Z}{32}\{P(\rho_{13})\}\oplus\frac{\mathbb Z}{4}\{\bar{\kappa}'\}.
	\end{aligned}\]
	\noindent       
	In light of \eqref{m2:Y}, we have
	\[\begin{aligned}
		&\Delta={r_3}_*\Delta': \pi^6_{26}\to \pi_{25}(\mathit{SO}(6)).
	\end{aligned}\]
	\noindent   
	In \cite[p. 12]{Ba2}, Barratt defined  $\rho_{13}\in\pi^{13}_{28}$ by a matrix
	Toda bracket:
	\[\FitMath{
		\left\{\bigl(\sigma_{13}, \Sigma\theta\bigr),\,\begin{pmatrix}8\sigma_{20}\\\eta^2_{25}\end{pmatrix},2\iota_{27}\right\}.
	}\]
	\noindent       We now follow Toda's definition of $\rho_{13}.$  The stable relation is $\x \cg,8\cg,2\e \xx=\rho+\mathrm{span}\n 2\rho\nn.$ 
	This follows from $\pm\rho\in\x \cg,2\cg,8\e \xx\subseteq\x \cg,8\cg,2\e \xx$ and the indeterminacy $\sigma\pi_8^S+2\pi_{15}^S=\mathrm{span}\{2\rho\}.$ 
	Since $\Sigma\theta$ stabilizes to zero, the stable image of the matrix bracket in \eqref{m1:rho13} lies in this coset. Its indeterminacy contains $2\pi^{13}_{28},$   whose stable image is $2\pi_{15}^S.$ 
	By Appendix~\ref{tlqdzj}, stabilization sends the generators $\rho_{13}$ and $\bar\varepsilon_{13}$ of $\pi^{13}_{28}\cong\mathbb Z/32\oplus\mathbb Z/2$ to $\rho$ and $\bar\varepsilon,$   respectively, and is an isomorphism onto $\pi_{15}^S.$ 
	The stable image of the matrix bracket is therefore the whole coset $\rho+\mathrm{span}\{2\rho\}.$  Injectivity of stabilization now gives
	\begin{equation}\label{m1:rho13}
		\rho_{13}\in\left\{\bigl(\sigma_{13}, \Sigma\theta\bigr),\,\begin{pmatrix}8\sigma_{20}\\\eta^2_{25}\end{pmatrix},2\iota_{27}\right\}_4.
	\end{equation}
	\noindent 	Similarly, we infer \begin{equation}\label{m001:rho13}
		\rho_{13}\in\left\{\bigl(\sigma_{13}, \Sigma\theta\hc\ca_{25}\bigr),\,\begin{pmatrix}8\sigma_{20}\\\eta_{26}\end{pmatrix},2\iota_{27}\right\}_4.
	\end{equation}

	\noindent       
	We know that $\pi^3_{25}=\frac{\mathbb Z}{2}\{\varepsilon_3\kappa_{11}\}$
	\ and\ ${i_U}_*\varepsilon_3=2[\nu^2_5]$ \cite[(4.1)]{MT2}. Therefore, ${i_U}_*:
	\pi^3_{25}\to\pi_{25}(\mathit{SU}(3))$ is trivial and
	\begin{equation}
		\begin{aligned}
			&{p_U}_*: \pi_{25}(\mathit{SU}(3))\to\pi^5_{25}\ \mbox{is a monomorphism}.\label{m1:pUmono}
		\end{aligned}
	\end{equation}
	
	The discussion proceeds with the following statement, after which we give its proof and then continue with the rest of the section.
	\begin{lem}\label{m1:mtb2}
		\begin{enumerate}
			\item $\Delta'(\bar\kappa')=[\nu^2_5]\kappa_{11}$\ and \
			$\Delta(\bar\kappa')=[ \nu_5]\nu_8\kappa_{11}.$ 
			\item $\Delta'(P(\rho_{13}))=[\nu_5\eta_8\mu_9]\sigma_{18}$\ and \
			$\Delta(P(\rho_{13}))=[ \nu_5]\eta_8\mu_9\sigma_{18}+4[\nu_4\zeta_7]_6\sigma_{18}.$ 
		\end{enumerate}
	\end{lem}
	\begin{proof}
		\begin{enumerate}[label=\textup{(\arabic*)},leftmargin=2em]
			\item By the proof of \cite[Lemma 15.3]{MT1}, we conclude that
			\[\begin{aligned}
				&\bar{\kappa}'\in\{\nu_6,\eta_9,\bar{\varepsilon}_{10}\}.
			\end{aligned}\]
			Since $\eta_7\bar{\varepsilon}_8=\sigma'(\eta_{14}\bar{\nu}_{15})$
			\eqref{m1:epep}, we can define $\bar{\kappa}'$ as
			\[\FitMath{
				\bar{\kappa}'\in \left\{\nu_6,\bigl(\eta_9, 2(\sigma_9\eta_{16})\bigr),\,\begin{pmatrix}\bar{\varepsilon}_{10}\\\bar{\nu}_{17}\end{pmatrix}\right\}_2.
			}\]
			Since $\Delta'(\nu_6)=[2\iota_5]\nu_5$ \cite[Proposition 6.3]{Mi2}, we see that
			\[\FitMath{
				\Delta'(\bar{\kappa}')\in \Delta'\left\{\nu_6,\bigl(\eta_9, 2(\sigma_9\eta_{16})\bigr),\,\begin{pmatrix}\bar{\varepsilon}_{10}\\\bar{\nu}_{17}\end{pmatrix}\right\}_2
				\subseteq\left\{[2\iota_5]\nu_5,\bigl(\eta_8, (\Sigma\sigma')\eta_{15}\bigr),\,\begin{pmatrix}\bar{\varepsilon}_9\\\bar{\nu}_{16}\end{pmatrix}\right\}_1.
			}\]
			Therefore, we obtain
			\[\FitMath{\begin{split}
					{p_U}_*\Delta'(\bar{\kappa}')\in\left\{2\nu_5,\bigl(\eta_8, (\Sigma\sigma')\eta_{15}\bigr),\,\begin{pmatrix}\bar{\varepsilon}_9\\\bar{\nu}_{16}\end{pmatrix}\right\}_1
					&\supseteq\nu_5\circ\left\{2\iota_8,\bigl(\eta_8, (\Sigma\sigma')\eta_{15}\bigr),\,\begin{pmatrix}\bar{\varepsilon}_9\\\bar{\nu}_{16}\end{pmatrix}\right\}_1\\
					&\bmod\ 2\nu_5\circ \Sigma\pi^7_{24}+\pi^5_{10}\circ\bar{\varepsilon}_{10}
					+\pi^5_{17}\circ\bar{\nu}_{17}.
			\end{split}}\]
			From the fact that $\pi^7_{24}\cong(\mathbb Z/2)^4,$   we note that
			$2\nu_5\circ \Sigma\pi^7_{24}=0.$  Observe that $$\pi^5_{10}\circ\bar{\varepsilon}_{10}=
			\mathrm{span}\{\nu_5\eta^2_8\bar{\varepsilon}_{10}\}=0.$$  Clearly, $\pi^5_{17}=\mathrm{span}\{\nu^4_5,
			\nu_5\mu_8, \nu_5\eta_8\varepsilon_9\}.$  Therefore, we conclude that
			$\pi^5_{17}\circ\bar{\nu}_{17}=0.$  This implies
			\[\FitMath{
				{p_U}_*\Delta'(\bar{\kappa}')=\nu_5\circ\left\{2\iota_8,\bigl(\eta_8, (\Sigma\sigma')\eta_{15}\bigr),\,\begin{pmatrix}\bar{\varepsilon}_9\\\bar{\nu}_{16}\end{pmatrix}\right\}_1.
			}\]
			
			The group decomposes as
			\[\begin{aligned}
				&\pi^8_{25}={}\mathrm{span}\{\sigma_8\eta_{15}\mu_{16}\}\oplus \Sigma\pi^7_{24}\cong(\mathbb Z/2)^5.
			\end{aligned}\]
			
			We have
			\[\FitMath{\begin{split}
					H\left\{2\iota_8,\bigl(\eta_8, (\Sigma\sigma')\eta_{15}\bigr),\,\begin{pmatrix}\bar{\varepsilon}_9\\\bar{\nu}_{16}\end{pmatrix}\right\}_1
					\subseteq&
					\left\{0,\bigl(\eta_8, (\Sigma\sigma')\eta_{15}\bigr),\,\begin{pmatrix}\bar{\varepsilon}_9\\\bar{\nu}_{16}\end{pmatrix}\right\}_1\ni 0\\
					&
					\qquad\bmod\ 0\circ\bar{\varepsilon}_{10}+\pi^{15}_{17}\circ\bar{\nu}_{17}=0.
			\end{split}}\]
			This implies
			\[\FitMath{
				\left\{2\iota_8,\bigl(\eta_8, (\Sigma\sigma')\eta_{15}\bigr),\,\begin{pmatrix}\bar{\varepsilon}_9\\\bar{\nu}_{16}\end{pmatrix}\right\}_1\subseteq \Sigma\pi^7_{24}.
			}\]
			
			It follows that
			\[\FitMath{\begin{split}
					\eta_7\circ\left\{2\iota_8,\bigl(\eta_8, (\Sigma\sigma')\eta_{15}\bigr),\,\begin{pmatrix}\bar{\varepsilon}_9\\\bar{\nu}_{16}\end{pmatrix}\right\}_1&=
					-\{\eta_7,2\iota_8,\eta_8\}_1\circ\bar{\varepsilon}_{10}
					-\{\eta_7,2\iota_8,(\Sigma\sigma')\eta_{15}\}_1\circ\bar{\nu}_{17}\\
					&\subseteq \pi^7_{10}\circ\bar{\varepsilon}_{10}-\pi^7_{17}\circ\bar{\nu}_{17}
					=0.
			\end{split}}\]
			Then we  have $\eta_7(\eta_8\mu_9\sigma_{18})=4\zeta_7\sigma_{18}\ne 0$ and $\eta_7\bar{\mu}_8\ne 0.$  Therefore, according to        Lemma~\ref{tx1023}(2), we obtain
			\[\FitMath{
				\left\{2\iota_8,\bigl(\eta_8, (\Sigma\sigma')\eta_{15}\bigr),\,\begin{pmatrix}\bar{\varepsilon}_9\\\bar{\nu}_{16}\end{pmatrix}\right\}_1\ni\nu_8\kappa_{11}\ \bmod\
				(\Sigma\sigma')\eta_{15}\mu_{16}.
			}\]
			Since $2(\nu_5\sigma_8)=\nu_5(\Sigma\sigma'),$   we find
			\[\begin{aligned}
				&{p_U}_*\Delta'(\bar{\kappa}')\equiv\nu^2_5\kappa_{11}\ \bmod\ \nu_5(\Sigma\sigma')\eta_{15}\mu_{16}=0.
			\end{aligned}\]
			This and \eqref{m1:pUmono} lead to the first half of (1). The second half of Assertion (1) follows directly from the first and
			Lemma~\ref{m1:r3n5}(1).
			
			\item  By \eqref{m1:rho13}, we conclude that
			\[\FitMath{
				P(\rho_{13})\in P\left\{\bigl(\sigma_{13}, \Sigma\theta\bigr),\,\begin{pmatrix}8\sigma_{20}\\\eta^2_{25}\end{pmatrix},2\iota_{27}\right\}_4
				\subseteq
				\left\{\bigl(P(\sigma_{13}), P(\Sigma\theta)\bigr),\,\begin{pmatrix}8\sigma_{18}\\\eta^2_{23}\end{pmatrix},2\iota_{25}\right\}_2
			}\]
			and
			\[\FitMath{
				\Delta'(P(\rho_{13}))\in\Delta'\left\{\bigl(P(\sigma_{13}), P(\Sigma\theta)\bigr),\,\begin{pmatrix}8\sigma_{18}\\\eta^2_{23}\end{pmatrix},2\iota_{25}\right\}_2
				\subseteq
				\left\{\bigl(\Delta'P(\sigma_{13}), \Delta'P(\Sigma\theta)\bigr),\,\begin{pmatrix}8\sigma_{17}\\\eta^2_{22}\end{pmatrix},2\iota_{24}\right\}_1.
			}\]
			According to        \cite[Corollary 5.3, Proposition 6.3]{Mi2}, we get 
			\[\begin{aligned}
				&\Delta'P(\sigma_{13})=[\nu_5\eta^2_8]\sigma_{10}\ \mbox{and}\ \Delta'P(\Sigma\theta)=[2\iota_5]\nu_5\kappa_8.
			\end{aligned}\]
			Since $\pi^5_{25}\cong(\mathbb Z/2)^3,$   we obtain
			\[\FitMath{\begin{split}
					{p_U}_*(\Delta'(P(\rho_{13}))&\in\left\{\bigl(\nu_5(\bar{\nu}_8+\varepsilon_8)\eta_{16}, 2\nu_5\kappa_8\bigr),\,\begin{pmatrix}8\sigma_{17}\\\eta^2_{22}\end{pmatrix},2\iota_{24}\right\}_1\\
					&
					=\{\nu_5(\bar{\nu}_8+\varepsilon_8)\eta_{16},8\sigma_{17},2\iota_{24}\}_1
					+\{2\nu_5\kappa_8,\eta^2_{22},2\iota_{24}\}_1\\
					&
					\ni\nu_5(\bar{\nu}_8+\varepsilon_8)\mu_{16}+\nu_5\kappa_8\circ 4\nu_{22}
					=\nu_5\eta_8\mu_9\sigma_{18}\\
					&\qquad
					\bmod\ \nu_5(\bar{\nu}_8+\varepsilon_8)\eta_{16}\circ \Sigma\pi^{16}_{24}+2\nu_5\kappa_8\circ \Sigma\pi^{21}_{24}
					+ 2\pi^5_{25}=0.
			\end{split}}\]
			This and \eqref{m1:pUmono} lead to the first half of (2).
			
			The second half of (2) is a direct consequence of the first and  Lemma~\ref{m1:r3n5}(2). This completes the proof.
		\end{enumerate}
	\end{proof}
	
	\subsection{Choices of Lifts and the Group Decomposition}
	
	As a consequence of Lemma~\ref{m1:mtb2}, we deduce
	\begin{equation}
		\begin{aligned}
			&{i_G}_*\pi_{25}(\mathit{SU}(3))=0.\label{m1:iGU25}
		\end{aligned}
	\end{equation}
	
	In light of \cite[p. 151]{T} and \cite[I-Proposition 6.4(4)]{Oda1979}, we have
	\begin{equation}
		\begin{aligned}
			&2\bar{\sigma}_6 = yP(\sigma^2_{13}) \text{\ for  some  odd\ } y.\label{m1:2bs6}
		\end{aligned}
	\end{equation}
	
	\needspace{3\baselineskip}\begin{lem}\label{m3:fh} \leavevmode\par
		\begin{itemize}
			\item[\textup{(1)}] In $\pi_{*}(G_{2}),$  for some odd integer $t_{1} \in\z$, the lift of $\bar{\sigma}_{6}$ may
			be chosen as $$\begin{aligned}
				&\langle\bar{\sigma}_{6}\rangle\in\{ \langle \bar{\nu} _{6}+\varepsilon_{6}\rangle,\allowbreak \sigma_{14},\allowbreak t_{1} \nu_{21}\}_{1}.
			\end{aligned}$$
			\item[\textup{(2)}]
			For the homomorphisms $\Delta\colon \pi_{*+1}(S^{6})\rightarrow \pi_{*}(\mathit{SU}(3)),$ we note that
			$\Delta(\bar{\sigma}_{6})=0.$
		\end{itemize}
	\end{lem}
	\begin{proof}
		\begin{itemize}
			\item[\textup{(1)}]
			By \cite{Mi2}, we see that
			$\langle \bar{\nu} _{6}+\varepsilon_{6}\rangle$ is well-defined, and  that $\langle \bar{\nu}
			_{6}+\varepsilon_{6}\rangle\circ \sigma_{14}\in\pi_{21}(G_{2})=0.$
			Thus, $\{ \langle \bar{\nu} _{6}+\varepsilon_{6}\rangle,\allowbreak  \sigma_{14},\allowbreak \nu_{21}\}_{1}$ is
			well-defined.
			By \cite[Proposition 3.4(3), p. 59]{Oda1979}, applying $\Sigma^{\infty}$ gives, for some odd $t_1,$  
			$t_1\bar{\sigma}_{6}\in\{ \bar{\nu} _{6}+\varepsilon_{6},\allowbreak \sigma_{14},\allowbreak \nu_{21} \}_{1}\;\bmod\;
			8\bar{\sigma}_{6}.$
			Hence (1) follows.
			
			\item[\textup{(2)}] This follows immediately from (1).
		\end{itemize}
	\end{proof}
	
	The next step in our presentation is to state the following result and then give its proof in the paragraphs below.
	\begin{lem}\label{m1:G2}
		
		$\pi_{25}(G_2)=
		\mathrm{span}\{\langle\bar{\sigma}_6\rangle,\; \langle\eta^2_6\rangle\bar{\mu}_8\}\cong \mathbb Z/32\oplus \mathbb Z/2$
		with relations
		${i_G}_*[\nu_5\eta_8\mu_9]\sigma_{18}=8\langle2[\iota_6,\iota_6]\rangle\sigma^2_{11}=0$ and
		$16\langle\bar{\sigma}_6\rangle
		= \langle P(\Sigma\theta)+\nu_6\kappa_9\rangle\eta^2_{23}.$ 
	\end{lem}
	\begin{proof}
		The first relation is \eqref{m1:diamond}.
		
		We consider  the exact sequence:
		\[\begin{aligned}
			&\pi^6_{26}\xrightarrow{\Delta'}\pi_{25}(\mathit{SU}(3))\xrightarrow{{i_G}_*}\pi_{25}(G_2) \xrightarrow{{p_G}_*}\pi^6_{25}\xrightarrow{\Delta'}\pi_{24}(\mathit{SU}(3)).
		\end{aligned}\]
		
		According to        \eqref{m1:iGU25},  we know ${i_G}_*=0.$ 
		Therefore, we conclude that
		\[\begin{aligned}
			&\pi_{25}(G_2)\cong\operatorname{ker}(\Delta':\pi^6_{25}\to\pi_{24}(\mathit{SU}(3))).
		\end{aligned}\] \noindent  
		Moreover, by \eqref{m1:2bs6},
		we see that
		\begin{equation}
			\begin{aligned}
				&4\langle\bar{\sigma}_6\rangle = y\langle [2\e_{6},\e_{6}]\rangle\sigma^2_{11}\ \text{\ for  some  odd\ } y.\label{m1:heart}
			\end{aligned}
		\end{equation}  \noindent   Since $8\langle2[\iota_6,\iota_6]\rangle\sigma_{11}={i_G}_*([\nu_5\eta_8\mu_9])$ \cite[p.160]{Mi2}, we obtain
		\[\begin{aligned}
			&0={i_G}_*([\nu_5\eta_8\mu_9]\sigma_{18})= 8\langle2[\iota_6,\iota_6]\rangle\sigma^2_{11}=32\langle\bar{\sigma}_6\rangle.
		\end{aligned}\]
		\noindent  This implies $\mathrm{ord}\langle\bar{\sigma}_6\rangle=32.$ 
		By \cite[Proposition 6.2]{Mi2}, we obtain $\Delta'(\bar{\zeta}_6)=[2\iota_5]\bar{\zeta}_5.$  This and the
		group structure of $\pi_{24}(\mathit{SU}(3))$ yield
		$\pi_{25}(G_2).$    By Lemma~\ref{m1:mtb2},
		$\Delta':\pi^6_{26}\to \pi_{25}(\mathit{SU}(3))$ is an epimorphism, and hence we obtain  the group structure of $\pi_{25}(G_2).$   According to        \eqref{m1:42P13} and \eqref{m1:iGU25},
		we deduce
		\begin{equation}
			\begin{aligned}
				&16\langle\bar{\sigma}_6\rangle = 4\langle2[\iota_6,\iota_6]\rangle\sigma^2_{11}=\langle P(\Sigma\theta)+\nu_6\kappa_9\rangle\eta^2_{23}.\label{m1:16<bs6>}
			\end{aligned}
		\end{equation}
		This completes the proof.
		
	\end{proof}
	\section{The Cases $7\leq n\leq9$}
	\subsection{The Case $n=7$ and Its Composition Relations}
	
	We first compare two expressions for a composite involving $\nu,$   $\sigma,$   and $\mu.$  The
	corresponding relation on lifts retains a fiber ambiguity, which must be followed through the
	inclusions.
	
	\begin{lem}\label{m1:n5musg}
		\begin{enumerate}
			\item
			$\nu_5\sigma_8\mu_{15} = \nu_5\mu_8\sigma_{17}.$ 
			\item $[ \nu_5]\sigma_8\mu_{15}\equiv[ \nu_5]\mu_8\sigma_{17}\ \bmod \
			[  \nu^2_4]_6\kappa_{10}, [\nu_4\sigma'\eta_{14}]_6\mu_{15}$\ and \\
			$[ \nu_5]_9\sigma_8\mu_{15}=[ \nu_5]_9\mu_8\sigma_{17}.$ 
		\end{enumerate}
	\end{lem}
	\begin{proof}
		\begin{enumerate}[label=\textup{(\arabic*)},leftmargin=2em]
			\item Recall that $\nu_5\sigma_8\eta_{15} = \nu_5\varepsilon_8$ from the stem $11$ entry of Table~\ref{tab:oguchi-compositions} and its footnote citing \cite[p.~152]{T}.
			In view of  Lemma~\ref{m1:m3sg}(1), we have \vspace{-1.1em}
			\begin{align*}
				\nu_5\sigma_8\mu_{15} &\in \nu_5\sigma_8\circ\{\eta_{15},8\iota_{16},2\sigma_{16}\}_1 \\
				&\subseteq\{\nu_5\varepsilon_8,8\iota_{16},2\sigma_{16}\}_1 \\
				&\supseteq\nu_5\circ\{\varepsilon_8,8\iota_{16}, 2\sigma_{16}\}_1 \\
				&\ni\nu_5\mu_8\sigma_{17} \bmod \nu_5\varepsilon_8\circ \Sigma\pi^{15}_{23} + \pi^5_{17}\circ 2\sigma_{17} + \mathrm{span}\{\nu_5\eta_8\bar{\varepsilon}_9\}=0.
			\end{align*}
			This leads to (1).
			
			\item By (1), \eqref{m1:524} and Proposition~\ref{m1:924}, we have (2). This completes the proof.
		\end{enumerate}
	\end{proof}
	
	Recall the relation
	\begin{equation}
		\begin{aligned}
			&2\varepsilon'=\eta^2_3\varepsilon_5.\label{m1:2ep'}
		\end{aligned}
	\end{equation}
	
	Using the standard relations, we obtain $\nu_4\sigma'\zeta_{14}=x\nu_4\zeta_7\sigma_{18}$
	for some odd $x.$ 
	By Lemma~5.2, the boundary composition formula in Section~4.1, and \eqref{m1:2ep'},
	we have
	$$\Delta(\nu_4\sigma'\zeta_{14})=\Delta(\nu_4\sigma')\zeta_{13}=2([\iota_3]+a[\eta_2]_4)\varepsilon'\zeta_{13}=([\iota_3]+a[\eta_2]_4)\eta^2_3\varepsilon_5\zeta_{13}=0,$$
	where the last equality uses the suspension of \eqref{m1:Ogep}. Therefore, we define $[\nu_4\sigma'\zeta_{14}].$ 
	
	According to        \eqref{m1:i4,6R425},
	we find
	\begin{equation}
		\begin{aligned}
			&[\nu_4\sigma'\zeta_{14}]\equiv x[\nu_4\zeta_7]\sigma_{18}\ \bmod\ [\iota_3]_5\varepsilon_3\kappa_{11}\ \mbox{for  some  odd}\ x.\label{m1:n4sg'zt}
		\end{aligned}
	\end{equation}  \noindent  Then we conclude that
	\begin{gather}
		\begin{aligned}
			&{i_7}_*\pi_{25}(\mathit{SO}(6))={}\mathrm{span}\{[\nu_4\zeta_7]_7\sigma_{18}, [ \nu_5]_7\bar{\mu}_8\}\cong\mathbb Z/8\oplus\mathbb Z/2.\label{m1:i7R625}
		\end{aligned}
	\end{gather}        Using the lift chosen in Lemma~\ref{m3:fh}, we define $[\bar\sigma_6]$ in $\pi_{25}(\mathit{SO}(7))$ as follows:
	\begin{equation}
		\begin{aligned}
			&[\bar\sigma_6]=\mathcal{T}_{*}\langle\bar{\sigma}_6\rangle \in \{[\bar{\nu}_6+\varepsilon_6],\sigma_{14},t_1\nu_{21}\}_1.\label{m1:[bsig6]}
		\end{aligned}
	\end{equation}

	The calculations developed in the preceding discussion
	also shed light on the present case, leading to
	the following determination.
	\begin{prop}\label{m1:257}
		
		The structure of $\pi_{25}(\mathit{SO}(7))$ is as follows.
		\[ \GroupLayout{\pi_{25}(\mathit{SO}(7)) = \mathrm{span}\{[\nu_4\zeta_7]_7\sigma_{18},\; [ \nu_5]_7\bar{\mu}_8,\; [\eta_6]\eta_7\bar{\mu}_8,\; [\bar{\sigma}_6]\} \cong\mathbb Z/8 \oplus(\mathbb Z/2)^2 \oplus\mathbb Z/32.}{\begin{aligned}\pi_{25}(\mathit{SO}(7)) = \mathrm{span}\{&[\nu_4\zeta_7]_7\sigma_{18},\; [ \nu_5]_7\bar{\mu}_8,\\
				&[\eta_6]\eta_7\bar{\mu}_8,\; [\bar{\sigma}_6]\} \cong\mathbb Z/8 \oplus(\mathbb Z/2)^2 \oplus\mathbb Z/32.\end{aligned}} \]
		Here $\Delta(P(\rho_{13}))=[ \nu_5]\eta_8\mu_9\sigma_{18}+4[\nu_4\zeta_7]_6\sigma_{18}.$ 
		
	\end{prop}
	\begin{proof}
		In the exact sequence $(25)_6$:
		\[\begin{aligned}
			&\pi^6_{26}\xrightarrow{\Delta}\pi_{25}(\mathit{SO}(6))\xrightarrow{i_*}\pi_{25}(\mathit{SO}(7))\xrightarrow{p_*}\pi^6_{25} \xrightarrow{\Delta}\pi_{24}(\mathit{SO}(6)),
		\end{aligned}\]
		we obtain
		$\operatorname{ker}(\Delta\colon \pi^6_{25}\to \pi_{24}(\mathit{SO}(6)))
		=\mathrm{span}\{\bar{\sigma}_6,\eta^2_6\bar{\mu}_8\}\cong\mathbb Z/32\oplus\mathbb Z/2$
		by Lemma~\ref{m3:fh}, the relation $(r_3)_*\Delta'=\Delta$,
		and the identity $\Delta\bar{\zeta}_6=2[\bar{\zeta}_5]$ from Proposition~\ref{m1:724}.
		Moreover, $\pi^6_{26}=\mathrm{span}\{P(\rho_{13}),\bar{\kappa}'\}\cong\mathbb Z/32\oplus\mathbb Z/4.$   Recall from \eqref{m1:i7R625} that
		\begin{gather}
			\begin{aligned}
				&{i_7}_*\pi_{25}(\mathit{SO}(6))={}\mathrm{span}\{[\nu_4\zeta_7]_7\sigma_{18}, [ \nu_5]_7\bar{\mu}_8\}\cong\mathbb Z/8\oplus\mathbb Z/2.
			\end{aligned}
		\end{gather} 
		In light of \eqref{m1:16<bs6>}, we see that
		\begin{gather}\label{m1:bsig11}
			4[2[\iota_6,\iota_6]]\sigma^2_{11}= [P(\Sigma\theta)+\nu_6\kappa_9]\eta^2_{23}=16[\bar{\sigma}_6]\ \bmod\ 4[\nu_4\zeta_7]_7\sigma_{18}
			\shortintertext{and}
			\label{m1:bsig12}
			8[2[\iota_{6},\iota_{6}]]\sigma^2_{11}=32[\bar{\sigma}_6]=0.
		\end{gather}  
		By Lemma~\ref{m1:mtb2}(2),
		we deduce
		\[\begin{aligned}
			&\Delta(P(\rho_{13})) = [ \nu_5]\eta_8\mu_9\sigma_{18} + 4[\nu_4\zeta_7]_6\sigma_{18}.
		\end{aligned}\]
		\noindent   By \eqref{m1:bsig12} and \cite[Lemma 3.7(ii)]{KM1}, we find
		\[\begin{aligned}
			&0=8[2[\iota_6,\iota_6]]\sigma^2_{11} = [ \nu_5]_7\eta_8\mu_9\sigma_{18} + 4[\nu_4\zeta_7]_7\sigma_{18}.
		\end{aligned}\]
		\noindent   This leads to the relation and completes the proof.
	\end{proof}
	
	According to       \eqref{m1:bnu+epzt}, we conclude that
	\begin{equation}
		\begin{aligned}
			&\pm8[\bar{\sigma}_6]=[\bar{\nu}_6+\varepsilon_6]\zeta_{14}.\label{m1:Ogu1}
		\end{aligned}
	\end{equation}
	
	\subsection{The Case $n=8$}
	Since $\pi_{25}(\mathit{SO}(8))\cong \pi_{25}(\mathit{SO}(7))\oplus \pi^7_{25}$ and
	$\pi^7_{25}=\mathrm{span}\{\zeta_7\sigma_{18},\eta_7\bar{\mu}_8\}\cong\mathbb Z/8\oplus\mathbb Z/2,$   we obtain the following.
	\begin{prop}\label{g1016:prop:25-8}
		The structure of $\pi_{25}(\mathit{SO}(8))$ is as follows.
		\[ \GroupLayout{\pi_{25}(\mathit{SO}(8)) = \mathrm{span}\{[\nu_4\zeta_7]_8\sigma_{18},\; [ \nu_5]_8\bar{\mu}_8,\; [\eta_6]_8\eta_7\bar{\mu}_8,\; [\bar{\sigma}_6]_8,\; [\iota_7]\zeta_7\sigma_{18},\; [\iota_7]\eta_7\bar{\mu}_8\} \;\cong\; \mathbb Z/8 \oplus(\mathbb Z/2)^2 \oplus\mathbb Z/32 \oplus\mathbb Z/8 \oplus\mathbb Z/2.}{\begin{aligned}\pi_{25}(\mathit{SO}(8)) = \mathrm{span}\{&[\nu_4\zeta_7]_8\sigma_{18},\; [ \nu_5]_8\bar{\mu}_8,\; [\eta_6]_8\eta_7\bar{\mu}_8,\; [\bar{\sigma}_6]_8,\\
				&[\iota_7]\zeta_7\sigma_{18},\; [\iota_7]\eta_7\bar{\mu}_8\} \;\cong\; \mathbb Z/8 \oplus(\mathbb Z/2)^2 \oplus\mathbb Z/32 \oplus\mathbb Z/8 \oplus\mathbb Z/2.\end{aligned}} \]
		\qed
	\end{prop}
	
	\subsection{The Case $n=9$}
	We show the following.
	\begin{prop}\label{m1:R259}
		The structure of $\pi_{25}(\mathit{SO}(9))$ is as follows.
		\begin{enumerate}
			\item
			$\pi_{25}(\mathit{SO}(9)) = \mathrm{span}\{[ \nu_5]_9\bar{\mu}_8, [\iota_7]_9\eta_7\bar{\mu}_8, [\bar{\sigma}_6]_9\} \cong(\mathbb Z/2)^2 \oplus\mathbb Z/32.$  Moreover, 
			$$\pm 8[\bar\sigma_6]_9 = [\iota_7]_9\zeta_7\sigma_{18},$$$$
			[P(\Sigma\theta)+\nu_6\kappa_9]_9\eta^2_{23}=
			16[\bar\sigma_6]_9 = 2[\iota_7]_9\zeta_7\sigma_{18}
			=[\eta_6]_9\zeta_7\sigma_{18}=[\nu_4\zeta_7]_9\sigma_{18}.$$
			\item ${i_{10}}_*\pi_{25}(\mathit{SO}(9))=\mathrm{span}\{[\bar{\sigma}_6]_{10},[\iota_7]_{10}\eta_7\bar{\mu}_8\}\cong\mathbb Z/32\oplus\mathbb Z/2.$  
		\end{enumerate}
	\end{prop}
	\begin{proof}
		\begin{enumerate}[label=\textup{(\arabic*)},leftmargin=2em]
			\item In the exact sequence $(25)_8$:
			\[\begin{aligned}
				\pi^8_{26}&\xrightarrow{\Delta}\pi_{25}(\mathit{SO}(8))
				\xrightarrow{i_*}\pi_{25}(\mathit{SO}(9))\xrightarrow{p_*}\pi^8_{25}
				\xrightarrow{\Delta}\pi_{24}(\mathit{SO}(8)),
			\end{aligned}\]
			$\Delta\colon \pi^8_{25}\to\pi_{24}(\mathit{SO}(8))$ is a monomorphism
			by the proof of Proposition~\ref{m1:924}. Hence $i_*$ is surjective.
			Note that 
			$\pi^8_{26}=\mathrm{span}\{\sigma_8\zeta_{15},\zeta_8\sigma_{19},
			\eta_8\bar{\mu}_9\}\cong(\mathbb Z/8)^2\oplus\mathbb Z/2.$
			We compute the $\q$-images of all three generators.
			
			By \cite[Lemma 3.4(iii)]{KM1},
			$[\nu_4\zeta_7]_7\equiv[\eta_6]\zeta_7
			\bmod 8[[2\iota_6,\iota_6]]\sigma_{11}.$
			Composing with $\sigma_{18}$ and using \eqref{m1:bsig12} gives
			$[\nu_4\zeta_7]_7\sigma_{18}=[\eta_6]\zeta_7\sigma_{18}.$
			Since $\Delta\iota_8=2[\iota_7]-[\eta_6]_8,$ we obtain
			$$  \Delta(\zeta_8\sigma_{19})
			=2[\iota_7]\zeta_7\sigma_{18}
			-[\nu_4\zeta_7]_8\sigma_{18},\quad
			\Delta(\eta_8\bar{\mu}_9)
			=[\eta_6]_8\eta_7\bar{\mu}_8.$$
			To compute the remaining boundary, we retain the distinction between
			$[\mathrm{c}_{\mathrm{KM}}(\bar{\nu}_6+\varepsilon_6)]$
			and $[\bar{\nu}_6+\varepsilon_6].$
			In view of the comparison following \eqref{m1:i7R614} and
			$\nu_8^2\zeta_{14}=0,$  we get 
			\begin{align*}
				[\mathrm{c}_{\mathrm{KM}}(\bar{\nu}_6+\varepsilon_6)]_8\zeta_{14}
				&\equiv\mathrm{c}_{\mathrm{KM}}
				[\bar{\nu}_6+\varepsilon_6]_8\zeta_{14}\equiv\pm8\mathrm{c}_{\mathrm{KM}}[\bar{\sigma}_6]_8
				\bmod [\eta_5\varepsilon_6]_8\zeta_{14},
			\end{align*}
			where the second congruence uses \eqref{m1:Ogu1}.
			The suspension of $\varepsilon_3\zeta_{11}=0$ in \eqref{m1:Ogep}
			implies $\eta_5\varepsilon_6\zeta_{14}=0.$ Thus
			$[\eta_5\varepsilon_6]\zeta_{14}\in(i_6)_*\pi_{25}(\mathit{SO}(5)).$
			By \eqref{m1:255}, \eqref{m1:R256}, and
			$[\iota_3]_6\varepsilon_3\kappa_{11}=0,$ this image is
			$\mathrm{span}\{[\nu_4\zeta_7]_6\sigma_{18}\}.$
			Consequently, \eqref{m1:Dsg82} gives
			\begin{align*}
				\Delta(\sigma_8\zeta_{15})
				&=x_1[\iota_7]\sigma'\zeta_{14}
				\pm8x_2\mathrm{c}_{\mathrm{KM}}[\bar{\sigma}_6]_8+a[\nu_4\zeta_7]_8\sigma_{18}
			\end{align*}
			for some $a\in\z,$ with the odd integers $x_1,x_2$ of
			\eqref{m1:Dsg82}. 
			
			Observe that 
			$\sigma'\zeta_{14}=x\zeta_7\sigma_{18}$ for some odd integer $x.$
			Start with the direct sum in Proposition~\ref{g1016:prop:25-8}.
			The  boundary formulas give
			$[\eta_6]_9\eta_7\bar{\mu}_8=0$ and
			$[\nu_4\zeta_7]_9\sigma_{18}=2[\iota_7]_9\zeta_7\sigma_{18}.$
			The quotient $\pi_{25}(\mathit{SO}(8))/\mathrm{span}\{\Delta(\eta_8\bar{\mu}_9),
			\Delta(\zeta_8\sigma_{19})\}$
			is the direct sum of the cyclic submodules generated by 
			$$[\nu_5]_8\bar{\mu}_8+W,\;
			[\iota_7]\eta_7\bar{\mu}_8+W,\;
			[\bar{\sigma}_6]_8+W, \text{ \;and  \;}
			[\iota_7]\zeta_7\sigma_{18}+W,$$
			whose orders in this quotient are $2,2,32,8,$ respectively; here, $W=\mathrm{span}\{\Delta(\eta_8\bar{\mu}_9),
			\Delta(\zeta_8\sigma_{19})\}.$

			Since $i_*\Delta(\sigma_8\zeta_{15})=0$, we obtain
			\[
			(x_1x+2a)[\iota_7]_9\zeta_7\sigma_{18}
			\pm8x_2\mathrm{c}_{\mathrm{KM}}[\bar{\sigma}_6]_9=0.
			\]Observe that $x_1x+2a$ and $x_2\mathrm{c}_{\mathrm{KM}}$ are odd. So,     $$[\iota_7]_9\zeta_7\sigma_{18}\ty
			8[\bar{\sigma}_6]_9\m2[\iota_7]_9\zeta_7\sigma_{18}.$$
			It follows that
			\begin{align*}
				\pi_{25}(\mathit{SO}(9))
				&=\mathrm{span}\{[\nu_5]_9\bar{\mu}_8,
				[\iota_7]_9\eta_7\bar{\mu}_8,[\bar{\sigma}_6]_9\}\\*
				&\cong(\mathbb Z/2)^2\oplus\mathbb Z/32.
			\end{align*}
			In addition,
			$[\iota_7]_9\zeta_7\sigma_{18}=\pm8[\bar{\sigma}_6]_9.$ Doubling gives
			\begin{align*}
				16[\bar{\sigma}_6]_9=2[\iota_7]_9\zeta_7\sigma_{18}=[\eta_6]_9\zeta_7\sigma_{18}
				=[\nu_4\zeta_7]_9\sigma_{18}.
			\end{align*}
			Finally, \eqref{m1:bsig11} gives
			$[P(\Sigma\theta)+\nu_6\kappa_9]_9\eta_{23}^2
			=16[\bar{\sigma}_6]_9,$ since 
			$4[\nu_4\zeta_7]_9\sigma_{18}$ is  zero. This proves (1).
			
			\item In the exact sequence $(25)_{9}$:
			\[\begin{aligned}
				&\pi^9_{26}\xrightarrow{\Delta}\pi_{25}(\mathit{SO}(9))\xrightarrow{i_*}\pi_{25}(\mathit{SO}(10))\xrightarrow{p_*}\pi^9_{25} \xrightarrow{\Delta}\pi_{24}(\mathit{SO}(9)),
			\end{aligned}\]
			we know $\pi^9_{26}=\mathrm{span}\{\sigma_9\eta_{16}\mu_{17},\nu_9\kappa_{12},\bar{\mu}_9,\eta_9\mu_{10}\sigma_{19}\}\cong(\mathbb Z/2)^4.$ 
			\noindent  
			According to        \eqref{m1:Desig9}, the relation $\bar{\nu}_7\mu_{15}=0,$  
			\eqref{m1:4z5sig} and (1), we infer that
			\[\begin{aligned}
				&\Delta(\sigma_9\eta_{16}\mu_{17}) = [\iota_7]_9(\bar{\nu}_7+\varepsilon_7)\eta_{15}\mu_{16} = [\iota_7]_9(\bar{\nu}_7\mu_{15}\eta_{24}+\eta_{7}\varepsilon_{8}\mu_{16}) = 4[\iota_7]_9\zeta_7\sigma_{18}=32[\bar{\sigma}_6]_9=0.
			\end{aligned}\]
			Using $\eta^2_7\mu_9=4\zeta_7$ \eqref{m1:4z5}, Proposition~\ref{m1:257},
			(1), \eqref{m1:Ogu0} and \eqref{m1:bsig12}, we find
			\begin{align*}
				&\Delta(\eta_9\mu_{10}\sigma_{19}) =[ \nu_5]_9\eta_8\mu_9\sigma_{18}+[\iota_7]_9\eta^2_7\mu_9\sigma_{18} =4[\nu_4\zeta_7]_9\sigma_{18}+32[\bar{\sigma}_6]_9 =0, \\
				&\Delta(\nu_9\kappa_{12})=[ \nu_5]_9\nu_8\kappa_{11}=0 \quad\text{and}\quad \Delta(\bar{\mu}_9)=[\iota_7]_9\eta_7\bar{\mu}_8+[ \nu_5]_9\bar{\mu}_8.
			\end{align*}
			This leads to the group ${i_{10}}_*\pi_{25}(\mathit{SO}(9))$ and the relation
			$[\iota_7]_{10}\eta_7\bar{\mu}_8=[ \nu_5]_{10}\bar{\mu}_8.$ 
			So we infer (2) and complete the proof.
		\end{enumerate}
	\end{proof}
	
	\section{Spin Group Comparisons and Additional Boundary Relations}
	\subsection{The Degree 25 Spin Groups}
	We also need degree-$25$ information for the spin groups. The next proposition records generators for
	$\mathit{Spin}(7)$ in a form suitable for the fibration with total space $\mathit{Spin}(9).$
	
	\begin{prop} \label{m2:spin7} The structure of $\pi_{25}(\mathit{Spin}(7))$ is as follows.
		\[\begin{aligned}
			&\pi_{25}(\mathit{Spin}(7)) = \mathrm{span}\bigl\{\langle\bar\sigma_6\rangle, \, \langle\eta^2_6\rangle\bar\mu_8, \, \zeta_7\sigma_{18}, \,\eta_7\bar\mu_8\bigr\} \cong\mathbb Z/32 \oplus\mathbb Z/2 \oplus\mathbb Z/8 \oplus\mathbb Z/2.
		\end{aligned}\]
		\qed
	\end{prop}
	
	We consider the elements $\langle \bar{\nu}_6 + \varepsilon_6 \rangle
	\eta_{14}\mu_{15}$ and $\langle \eta^2_6 \rangle \mu_8 \sigma_{17}$ in
	$\pi_{24}(G_2).$
	Since $${p_G}_\ast(\langle \bar{\nu}_6 + \varepsilon_6 \rangle \eta_{14}\mu_{15}) =
	\eta_6\varepsilon_7\mu_{15}=\eta^2_6\mu_8\sigma_{17},$$ \noindent   we have $\langle \bar{\nu}_6 + \varepsilon_6
	\rangle \eta_{14}\mu_{15} \equiv \langle \eta^2_6 \rangle \mu_8 \sigma_{17} \bmod
	{i_G}_{\ast}\pi_{24}(\mathit{SU}(3))=0.$ This implies that $\langle \bar{\nu}_6 + \varepsilon_6
	\rangle \eta_{14}\mu_{15}$ is equal to $\langle \eta^2_6 \rangle \mu_8 \sigma_{17}.$ We now pass from
	$\mathit{Spin}(7)$ to $\mathit{Spin}(9).$ The exact sequence of the fibration over $S^{15}$ determines
	which degree-$25$ classes survive, while the preceding composition relation helps identify the
	relevant boundary.
	
	\begin{prop}\label{m2:spin9}
		The structure of $\pi_{25}(\mathit{Spin}(9))$ is as follows.
		\[ \pi_{25}(\mathit{Spin}(9)) \cong \mathrm{span}\{\langle\bar{\sigma}_6\rangle,\; \langle\eta^2_6\rangle\bar{\mu}_8,\; \eta_7\bar{\mu}_8\} \cong \z/32 \oplus(\mathbb Z/2)^2.\]
	\end{prop}
	\begin{proof}
		Consider  the exact sequence
		\[\begin{aligned}
			&\pi^{15}_{k+1}\stackrel{{\Delta}}{\longrightarrow}\pi_k(\mathit{Spin}(7))\stackrel{{i_*}}{\longrightarrow}\pi_k(\mathit{Spin}(9))\stackrel{{p_*}}{\longrightarrow}\pi^{15}_k\stackrel{{\Delta}}{\longrightarrow}\pi_{k-1}(\mathit{Spin}(7)).
		\end{aligned}\]
		Recall from \cite[Proposition 9.2]{Mi2} the relation
		$\Delta\iota_{15}=x\langle\bar{\nu}_6+\varepsilon_6\rangle+y\sigma'+zi_*[\nu^2_5]\nu_{11},$ where $x$
		and $y$ are odd integers. By \cite[p.175]{Mi2},
		$\Delta : \pi^{15}_{24} \to \pi_{23}(\mathit{Spin}(7))$
		is
		injective. Consequently,
		$i_*:\pi_{24}(\mathit{Spin}(7))\to\pi_{24}(\mathit{Spin}(9))$
		is surjective.
		By the equation $\langle \bar{\nu}_6 + \varepsilon_6 \rangle \eta_{14}\mu_{15}= \langle \eta^2_6
		\rangle \mu_8 \sigma_{17},$ we see that $\Delta(\eta_{15}\mu_{16})= \langle \bar{\nu}_6 + \varepsilon_6
		\rangle \eta_{14}\mu_{15} + \sigma'\eta_{14}\mu_{15} = \langle \eta^2_6 \rangle \mu_8
		\sigma_{17}+\sigma'\eta_{14}\mu_{15}.$ Hence, we get
		$\pi_{24}(\mathit{Spin}(9)) \cong (\mathbb Z/2)^6.$
		
		Next, we consider $\Delta\zeta_{15}=x\langle\bar{\nu}_6+\varepsilon_6\rangle\zeta_{14}+y\sigma'\zeta_{14}.$
		The standard composition relations give
		\[\begin{aligned}
			&\pm8\bar{\sigma}_6=\bar{\nu}_6\zeta_{14} =(\bar{\nu}_6+\varepsilon_6)\zeta_{14}.
		\end{aligned}\]
		Therefore,
		$\langle\bar{\nu}_6+\varepsilon_6\rangle\zeta_{14}\equiv \pm 8\langle\bar{\sigma}_6\rangle\ \bmod\ {i_G}_*\pi_{25}(\mathit{SU}(3)).$
		Hence we obtain  $\langle\bar{\nu}_6+\varepsilon_6\rangle\zeta_{14}= \pm 8\langle\bar{\sigma}_6\rangle$
		by  ${i_G}_*\pi_{25}(\mathit{SU}(3))=0$ from \eqref{m1:iGU25}.
		This induces $\Delta\zeta_{15}=y\sigma'\zeta_{14}\pm8x\langle\bar{\sigma}_6\rangle
		,$ and thus, we get the group $\pi_{25}(\mathit{Spin}(9)).$
	\end{proof}
	
	\subsection{The $J$-Images for Small Values of $n$}
	We need the following group from \cite[Table 2]{Ka}:
	$$\pi_{7}(\mathit{SO}(5))=\mathbb Z_{(2)}\{[4 \nu_4]\},\ \pi_{7}(\mathit{SO}(6))=\mathbb Z_{(2)}\{[\eta^2_5]\}.$$
	Recall the relations $J[4 \nu_4]=\sigma'''$ and $J[\eta^2_5]=\pm
	\sigma''$ \cite[Lemma 2.1(i)]{HKM}. Hence, we obtain $J(2[\eta^2_5]-[4 \nu_4]_6)=0$ and
	$2[\eta^2_5]\equiv[4 \nu_4]_6 \bmod \operatorname{ker} J=\mathrm{span}\{4[\eta^2_5]\}.$ Consequently,
	$[4 \nu_4]_6=2u[\eta^2_5]$ for some $u\in\mathbb Z_{(2)}^{\times}.$ 
	
	Observe that $$\pi^{10}_{25}=\mathrm{span}\{\Sigma\rho', \sigma_{10}\bar\nu_{17},
	\bar\varepsilon_{10}\} \cong \mathbb Z/16 \oplus ((\mathbb Z/2)^2).$$ 
	We deduce $2[\bar\varepsilon_{10}] \in \{[\iota_7]_{11}\nu_7, 2\iota_{10}, \kappa_{10}\} \circ
	2\iota_{25} \subseteq \{[\iota_7]_{11}\nu_7, 2\iota_{10}, 0\}= [\iota_7]_{11}\nu_7 \circ
	\pi^{10}_{25}.$ Suspending \eqref{m2:Osm} gives $\nu_7\Sigma\rho'=\Sigma(\nu_6\rho')=0.$
	The stem $10$ relation in Table~\ref{tab:oguchi-compositions} gives $\sigma'\nu_{14}=x\nu_7\sigma_{10}$ for some odd integer $x.$ Thus $x\nu_7\sigma_{10}\bar\nu_{17}=\sigma'\nu_{14}\bar\nu_{17}=0,$ so $\nu_7\sigma_{10}\bar\nu_{17}=0$ because multiplication by $x$ is invertible in the 2-local group.
	The composition relations in the same table also give $\nu_7\bar\varepsilon_{10}=\nu_7\eta_{10}\kappa_{11}=0.$ Hence $[\iota_7]_{11}\nu_7\circ\pi^{10}_{25}=0,$ and therefore $2[\bar\varepsilon_{10}]=0.$
	\section{The Cases $10\leq n\leq14$}
	\subsection{The Case $n=10$ and the Required Lift Relations}
	
	For the case $n=10,$   a relation involving $P(\bar\kappa_{13})$ is needed before the extension can be
	resolved. We establish this relation in a homotopy group of a sphere first and then use it to compare the multiples of the
	lifts to special orthogonal groups.
	
	We show \begin{lem} $2P(\bar{\kappa}_{13}) = \bar{\nu}_6\nu_{14}\kappa_{17}.\label{m1:2Pbk_13}$
	\end{lem}

	\begin{proof}
		By \cite[Propositions 2.7, 6.2]{T} and \cite[Lemma 15.4]{MT1}, we see that
		\[\begin{aligned}
			&H(2P(\bar{\kappa}_{13})) = 4\iota_{11}\circ \bar{\kappa}_{11}= 4\bar{\kappa}_{11} \text{\ and\ } H(\bar{\nu}_6\nu_{14}\kappa_{17}) =\nu^2_{11}\kappa_{17}= 4\bar{\kappa}_{11}.
		\end{aligned}\]
		Thus we infer
		\[\begin{aligned}
			&2P(\bar{\kappa}_{13}) \equiv \bar{\nu}_6\nu_{14}\kappa_{17}\,\,\bmod\,\Sigma\pi^5_{30},
		\end{aligned}\]
		where
		\[\begin{aligned}
			&\Sigma\pi^5_{30}={}\mathrm{span}\{\Sigma\phi'',\nu^2_6\bar{\sigma}_{12}, \mu_{3,6}, \eta_6\bar{\mu}_7\sigma_{24}\} \cong \mathbb Z/8 \oplus (\mathbb Z/2)^3.
		\end{aligned}\]
		In light of the relations in \cite[I-Proposition 3.1(2)]{Oda1979}, we obtain
		\[\begin{aligned}
			&\bar{\nu}_7\kappa_{15}\nu_{29} = \sigma'\eta_{14}\kappa_{15}\nu_{29} + \varepsilon_7\kappa_{15} \nu_{29} = \sigma'\eta_{14}\nu_{15}\kappa_{18} + \varepsilon_7\nu_{15}\kappa_{18}= 0.
		\end{aligned}\]
		The facts that
		$
		\Sigma(2P(\bar{\kappa}_{13})) = 0
		$
		and that
		$
		\Sigma \colon \Sigma\pi^5_{30} \to \Sigma\pi^6_{31}
		$
		is an isomorphism
		\cite[I-(7.7),(7.12)]{Oda1979}
		yield the desired assertion.
		
	\end{proof}
	
	Proposition~\ref{m1:R259}(2) gives ${i_{10}}_*\pi_{25}(\mathit{SO}(9))$ and the relation
	$[\iota_7]_{10}\eta_7\bar{\mu}_8=[ \nu_5]_{10}\bar{\mu}_8.$

	Notice that the element $[\nu_{15}]\in \pi_{18}(\mathit{SO}(16))$ has order $8$ \cite{Ke,KM1}. Recall the relation from Proposition~\ref{m1:R259} (1):
	\begin{gather}\label{m1:r1}
		16[\bar\sigma_6]_9 = 2[\iota_7]_9\zeta_7\sigma_{18} =
		[\eta_6]_9\zeta_7\sigma_{18}= [\nu_4\zeta_7]_9\sigma_{18}
		\shortintertext{and}
		\label{m1:8bs}
		[\iota_7]_{9}\zeta_7\sigma_{18} = \pm 8[\bar{\sigma_6}]_{9}.
	\end{gather}
	\noindent
	By \cite[p.~166]{T}, we deduce
	\begin{equation}
		\begin{aligned}
			&\sigma_{12}\zeta_{19}=8P(\sigma_{25})=16\xi_{12}.\label{m1:s10z}
		\end{aligned}
	\end{equation}
	
	By the proof of \cite[I-Proposition 4.1(4)]{Oda1979}, we find
	$\phi'''\in\{\sigma'\sigma_{14},\sigma_{21},\nu_{28}\}_1.$ 
	The indeterminacy of this Toda bracket is
	$\sigma'\sigma_{14}\circ \Sigma\pi^{20}_{31}
	+\pi^7_{29}\circ\nu_{29}.$ 
	By \eqref{m1:s10z}, we conclude that
	$\sigma'\sigma_{14}\circ \Sigma\pi^{20}_{31}=\mathrm{span}\{\sigma'\sigma_{14}\zeta_{21}\}=0.$ 
	We use the group structure
	$$\pi^7_{29}=\mathrm{span}\{\sigma'\rho_{14}, \bar{\nu}_7\kappa_{15}, \varepsilon_7\kappa_{15}, \nu_7\bar{\sigma}_{10}\}\cong\mathbb Z/8\oplus(\mathbb Z/2)^3.$$    \noindent  
	Note that
	$\varepsilon_7\kappa_{15}\nu_{29}=\varepsilon_7\nu_{15}\kappa_{18}=0.$ 
	By the proof of Lemma~\ref{m1:2Pbk_13}, we obtain $\bar{\nu}_7\kappa_{15}\nu_{29}=\bar{\nu}_7\nu_{15}\kappa_{18}=0.$ 
	By \cite[I-Proposition 3.5(1)]{Oda1979}, we deduce $\nu_7\bar{\sigma}_{10}\nu_{29}
	=\nu^2_7\bar{\sigma}_{13}.$  By \cite[II-Proposition 2.1(4)]{Oda1979}, we have
	$\sigma'\rho_{14}\equiv x\rho''\sigma_{22}\ \bmod\ \bar{\nu}_7\kappa_{15}+\varepsilon_7\kappa_{15}$ for  some  odd $x.$  Then $\sigma'\rho_{14}\nu_{29}=0.$ 
	This implies that the indeterminacy of $\{\sigma'\sigma_{14},\sigma_{21},\nu_{28}\}_1$ is $\nu^2_7\bar{\sigma}_{13}.$  Therefore, we conclude that
	$\phi'''\in\{\sigma'\sigma_{14},\sigma_{21},\nu_{28}\}_1 \ \bmod\ \nu^2_7\bar{\sigma}_{13}.$ 
	
	The comparison following \eqref{m1:i7R614}, the normalizations
	\eqref{m1:Jet5ep} and \eqref{m1:Jbn6e6}, and
	$2\kappa_7=\bar\nu_7\nu_{15}^2$ give
	$$J([\bar{\nu} _{6}+\varepsilon_{6}])\equiv \sigma'\sigma_{14}
	\bmod 2\kappa_7,\;2\sigma'\sigma_{14}.$$
	To control this correction in the defining bracket of
	$[\bar\sigma_6],$ observe that
	$\bar\nu_7\nu_{15}\kappa_{18}=0$ by the proof of
	Lemma~\ref{m1:2Pbk_13}, while $\nu_{15}\sigma_{18}=0.$
	Appendix~\ref{tlqdzj} gives
	$\pi_{32}^{18}=\mathrm{span}\{\sigma_{18}^2,\kappa_{18}\}.$
	Thus
	\[
	0\in\bar\nu_7\nu_{15}\circ\{\nu_{18},\sigma_{21},\nu_{28}\}_1
	\subseteq\{2\kappa_7,\sigma_{21},\nu_{28}\}_1.
	\]
	The first indeterminacy term is zero, since
	$2\kappa_7\zeta_{21}=\bar\nu_7\nu_{15}^2\zeta_{21}=0$
	by \eqref{m1:n6nz}; the second is
	$\pi_{29}^7\circ\nu_{29}=\mathrm{span}\{\nu_7^2\bar\sigma_{13}\},$
	as computed above. Hence
	$\{2\kappa_7,\sigma_{21},\nu_{28}\}_1
	=\mathrm{span}\{\nu_7^2\bar\sigma_{13}\}.$
	By Lemma~\ref{m3:fh} and the definition $[\bar\sigma_6]=\mathcal{T}_*\langle\bar\sigma_6\rangle$,
	applying $J$ and using additivity in the first entry of the Toda bracket gives
	\begin{equation*}
		J[\bar\sigma_6]\equiv x\phi'''
		\ \bmod\ \nu_7^2\bar\sigma_{13}\qquad(x:\text{ odd}).
	\end{equation*}
	The odd factor $t_1$ is absorbed into $x$ and does not change the indeterminacy subgroup.
	\noindent
	In light of the definition of $\nu^*_{16}$ \cite[p. 153]{T} and that of $\phi''',$   we obtain
	\begin{align*}
		\sigma_9\nu^*_{16}&\in\sigma_9\circ\{\sigma_{16},2\sigma_{23},\nu_{30}\}_1
		\subseteq\{\sigma^2_9,2\sigma_{23},\nu_{30}\}_1
		\shortintertext{and}
		\Sigma^2\phi'''&\in\{2\sigma^2_9,\sigma_{23},\nu_{30}\}_1
		\subseteq\{\sigma^2_9,2\sigma_{23},\nu_{30}\}_1.
	\end{align*}
	Therefore, we deduce
	\[\begin{aligned}
		&\sigma_9\nu^*_{16}\equiv \Sigma^2\phi'''\ \bmod\ \sigma^2_9\circ \Sigma\pi^{22}_{33} + \pi^9_{31}\circ\nu_{31}.
	\end{aligned}\]
	We know that $\sigma^2_9\circ \Sigma\pi^{22}_{33}=\sigma^2_9\circ\mathrm{span}\{\zeta_{23}\}=0,$  
	$\pi^9_{31}\circ\nu_{31}=\mathrm{span}\{\nu^2_9\bar{\sigma}_{15}\}$ and
	$\nu^2_9\bar\sigma_{15}=2(\sigma_9\nu^{\ast}_{16})$ \cite[(3.6)]{MiM1}.
	Thus \cite[I-Proposition 5.1(2)]{Oda1979} gives the first relation below.
	The corrected congruence for $J[\bar\sigma_6]$ then gives the second:
	$x$ is odd, and $\sigma_9\nu^*_{16}$ has order $4$ by
	Appendix~\ref{tlqdzj}.
	\begin{gather}
		\Sigma^2\phi'''\equiv\sigma_9\nu^*_{16}\ \bmod\ 2\sigma_9\nu^*_{16} \notag
		\shortintertext{and}
		\label{m1:[bsg69}
		J([\bar\sigma_6]_9)= \pm \sigma_9\nu^*_{16}.
	\end{gather}
	
	By \cite[Proposition]{Ba} and by the relation $[\iota_{11}, \sigma^2_{11}]=\sigma_{11}\nu_{18}\sigma^2_{21}=0$ \cite[Lemma 12.14]{T}, we find
	$\psi_{11}\eta_{34} = \eta_{11}\psi_{12} + [\iota_{11},\sigma^2_{11}] =
	\eta_{11}\psi_{12}.$ 
	Therefore, by the EHP sequence, we conclude that
	\[\begin{aligned}
		&\psi_{10}\eta_{33} - \eta_{10}\psi_{11} \in P\pi^{21}_{36} = \mathrm{span}\{P(\bar\varepsilon_{21}), P(\rho_{21})\}.
	\end{aligned}\]
	Then we deduce that
	\[\begin{aligned}
		&P(\bar\varepsilon_{21}) =P(\eta_{21})\kappa_{20}= (2\sigma_{10}\nu_{17})\kappa_{20} =0.
	\end{aligned}\] \noindent   Since the order of
	$P(\rho_{21})$
	is 32 \cite[Theorem 1.1]{MMO}, we get
	\[\begin{aligned}
		&\psi_{10}\eta_{33} + \eta_{10}\psi_{11} = 16xP(\rho_{21}) \,\,\,\,\,(x \in \{0,1\})
	\end{aligned}\]
	and
	$
	\psi_{10}\eta^2_{33} = \eta_{10}\psi_{11}\eta_{34} = \eta^2_{10}\psi_{12}.
	$
	By the relation
	$2(\sigma_9\nu^{\ast}_{16})=\nu^2_9\bar\sigma_{15}=\eta^2_9\psi_{11}$ \cite[Lemma 3.5]{MiM1},
	we have
	\begin{equation}
		\begin{aligned}
			&\psi_{10}\eta^2_{33} = \eta^2_{10}\psi_{12} = 2\sigma_{10}\nu^\ast_{17} = \nu^2_{10}\bar\sigma_{16}.\label{m1:ps10et2}
		\end{aligned}
	\end{equation}  \noindent   As a consequence of  \eqref{m1:heart}, we have 
	\begin{equation}
		\begin{aligned}
			&4[\bar{\sigma}_6] =x[[2\e_{6},\e_{6}]]\sigma^2_{11}\ \mbox{for  some  odd}\ x.\label{m1:4bs}
		\end{aligned}
	\end{equation}
	
	We now take up the following result, presenting its statement first and then turning to the proof that follows.
	\begin{prop}\label{m1:25-10}
		The structure of $\pi_{25}(\mathit{SO}(10))$ is as follows.
		\[ \GroupLayout{\pi_{25}(\mathit{SO}(10)) = \mathrm{span}\{[\bar{\sigma}_6]_{10},\; [\iota_7]_{10}\eta_7\bar{\mu}_8,\; [\sigma^2_9]\eta^2_{23},\; [\sigma_9\eta_{16}\varepsilon_{17}]\}\cong \mathbb Z/32\oplus (\mathbb Z/2)^{3}.}{\begin{aligned}\pi_{25}(\mathit{SO}(10)) = \mathrm{span}\{&[\bar{\sigma}_6]_{10},\; [\iota_7]_{10}\eta_7\bar{\mu}_8,\; [\sigma^2_9]\eta^2_{23},\\
				&[\sigma_9\eta_{16}\varepsilon_{17}]\}\cong \mathbb Z/32\oplus (\mathbb Z/2)^{3},\end{aligned}} \] where $[\cg_{9}\lt_{16}\ca_{24}]+8[\ccs_{6}]_{10}\ty\n [P(\Sigma\theta)+\nu_{6}\cn_{9}]_{10},\;2\e_{23},\ca_{23}\nn\m[\cg_{9}^{2}]\ca_{23}^{2}.$
	\end{prop}
	\begin{proof}
		Consider  the exact sequence
		\[\begin{aligned}
			&\pi^9_{26} \xrightarrow{ \Delta} \pi_{25}(\mathit{SO}(9)) \xrightarrow{i_\ast } \pi_{25}(\mathit{SO}(10)) \xrightarrow{ p_\ast} \pi^9_{25} \xrightarrow{\Delta}\pi_{24}(\mathit{SO}(9)),
		\end{aligned}\]
		where $ \pi^9_{26} = (\frac{\mathbb Z}{2})^4\{\sigma_9\eta_{16}\mu_{17},
		\nu_9\kappa_{12}, \bar{\mu}_9, \eta_9\mu_{10}\sigma_{19}\} .$  By the proof of Proposition~\ref{m1:R259}(2),
		we know
		\[\begin{aligned}
			&\Delta(\sigma_9\eta_{16}\mu_{17}) = \Delta(\nu_9\kappa_{12}) = \Delta(\eta_9\mu_{10}\sigma_{19}) = 0.
		\end{aligned}\]
		Since $ \Delta \iota_9 = [ \nu_5]_9
		+ [\iota_7]_9\eta_7 ,$   we find
		\begin{align*}
			\Delta \bar{\mu}_9 &=
			[ \nu_5]_9\bar{\mu}_8 + [\iota_7]_9\eta_7\bar{\mu}_8
			\shortintertext{and}
			\Delta(\mu_9\sigma_{18}) &=
			[ \nu_5]_9\mu_8\sigma_{17}+[\iota_7]_9\eta_7\mu_8\sigma_{17}.
		\end{align*}
		As a consequence of \eqref{m1:Desig9}, 
		we get 
		$\Delta(\sigma_9\mu_{16}) =\allowbreak  [\iota_7]_9\eta_7\mu_8\sigma_{17}=
		[\iota_7]_9(\bar\nu_7 + \varepsilon_7)\mu_{15}.$ 
		By the calculation of $\pi_{24}(\mathit{SO}(10)),$   we see that $$\operatorname{ker}
		(\Delta : \pi^9_{25} \to \pi_{24}(\mathit{SO}(9)))=\mathrm{span}\{\sigma_9\nu^3_{16},
		\sigma_9\eta_{16}\varepsilon_{17}\}.$$
		\noindent
		In view of  Proposition~\ref{m1:R259}(2), we have
		a short exact sequence
		\[\begin{aligned}
			&0\to \operatorname{coker} \Delta \xrightarrow{i_\ast } \pi_{25}(\mathit{SO}(10)) \xrightarrow{ p_\ast} (\frac{\mathbb Z}{2})^2\{\sigma^2_9\eta^2_{23}, \sigma_9\eta_{16}\varepsilon_{17}\} \to 0.
		\end{aligned}\]
		Consider  lifts of $\sigma_9\eta_{16}\varepsilon_{17}.$  By Equation \eqref{1024zmf}, we know $\Delta(\cg_{9}\lt_{16})=2[P(\Sigma\theta)+\nu_{6}\cn_{9}]_{9}.$  Then we deduce that 
		\begin{align}
			\n i_{9,10},\Delta\e_{9},\cg_{8}\lt_{15}\ca_{23}\nn&\xyd
			\n i_{9,10},\;\Delta(\cg_{9}\lt_{16}),\ca_{23}\nn
			\notag\\
			&=
			\n i_{9,10},\;2[P(\Sigma\theta)+\nu_{6}\cn_{9}]_{9},\ca_{23}\nn
			\notag\\
			&\dyd
			\n [P(\Sigma\theta)+\nu_{6}\cn_{9}]_{10},\;2\e_{23},\ca_{23}\nn.
			\notag
		\end{align} \noindent$\ind\n i_{9,10},\;\Delta(\cg_{9}\lt_{16}),\ca_{23}\nn=i_{9,10_{*}}\pi_{25}(\mathit{SO}(9))+\pi_{24}(\mathit{SO}(10))\hc\ca_{24}.$  Take an element $$\mathbf{x}\in\n [P(\Sigma\theta)+\nu_{6}\cn_{9}]_{10},\;2\e_{23},\ca_{23}\nn .$$
		\noindent   
		It follows that 
		$$p_{*}(\mathbf{x})\ty\cg_{9}\lt_{16}\ca_{24}\m p_{*}( i_{9,10_{*}}\pi_{25}(\mathit{SO}(9))+\pi_{24}(\mathit{SO}(10))\hc\ca_{24})=\mathrm{span}\n \cg_{9}^{2}\ca_{23}^{2}\nn.$$\noindent
		So, $\cz\mathbf{y}=\mathbf{x}+t[\cg_{9}^{2}]\ca_{23}^{2}$ for some $t\in\z$ such that $p_{*}(\mathbf{y})=\cg_{9}\lt_{16}\ca_{24}.$
		Thus, 
		\begin{align}
			2\mathbf{y}
			&=2\mathbf{x}
			\notag\\
			&\in
			\n [P(\Sigma\theta)+\nu_{6}\cn_{9}]_{10},\;2\e_{23},\ca_{23}\nn
			\hc2\e_{25}
			\notag\\
			&=
			-[P(\Sigma\theta)+\nu_{6}\cn_{9}]_{10}
			\hc\n 2\e_{23},\ca_{23},2\e_{24}\nn
			\notag\\
			&=
			[P(\Sigma\theta)+\nu_{6}\cn_{9}]_{10}\ca_{23}^{2}
			\notag\\
			&=
			16[\ccs_{6}]_{10},\;
			(\mathrm{Prop.}\,\ref{m1:R259}).
			\notag
		\end{align}\noindent
		Define $[\cg_{9}\lt_{16}\ca_{24}]=\mathbf{y}-8[\ccs_{6}]_{10}=\mathbf{x}+t[\cg_{9}^{2}]\ca_{23}^{2}-8[\ccs_{6}]_{10}.$  Then, $$p_{*}([\cg_{9}\lt_{16}\ca_{24}])=\cg_{9}\lt_{16}\ca_{24},\qquad2[\cg_{9}\lt_{16}\ca_{24}]=0.$$  \noindent In addition,  $[\cg_{9}\lt_{16}\ca_{24}]+8[\ccs_{6}]_{10}\ty\n [P(\Sigma\theta)+\nu_{6}\cn_{9}]_{10},\;2\e_{23},\ca_{23}\nn\m[\cg_{9}^{2}]\ca_{23}^{2}.$ Hence the results hold.\end{proof}
	
	\subsection{The Case $n=11$}
	
	The passage to the case $n=11$ depends on the incoming boundary of $P(\sigma_{21}).$  The following formula
	expresses it in the generators already chosen for $n=10,$   making the quotient relation explicit.
	
	\begin{lem}\label{m1:DPs}
		$
		\Delta(P(\sigma_{21})) = [\sigma^2_9]\eta^2_{23} + 2x[\bar\sigma_6]_{10}
		$ for some odd $x.$ 
	\end{lem}
	
	\begin{proof}
		By \cite[p.~342-3]{KM1}, \eqref{m1:8bs} and \eqref{m1:4bs}, we have
		\begin{align*}
			\Delta(P(\sigma_{21})) &= [P(\eta_{19})]\sigma_{18} + 2[\eta_9\varepsilon_{10}]\sigma_{18} \\
			&= [P(\eta_{19})]\sigma_{18} + [\iota_7]_{10}\zeta_7\sigma_{18} + 2[[2\e_{6},\e_{6}]]_{10}\sigma^2_{11} \\
			&= [P(\eta_{19})]\sigma_{18} +16y[\bar\sigma_6]_{10},
		\end{align*}
		for $ y \in \{0,1\}. $  Since $P(\eta_{19})\sigma_{18}=(\eta_9\sigma_{10}+\sigma_9\eta_{16})\eta_{17}\sigma_{18}=\sigma^2_9\eta^2_{23},$  
		we see that
		\[\begin{aligned}
			&[P(\eta_{19})]\sigma_{18}\equiv[\sigma^2_9]\eta^2_{23}\ \bmod\ {i_{10}}_*\pi_{25}(\mathit{SO}(9)).
		\end{aligned}\]
		By Proposition~\ref{m1:R259}(2), we obtain
		\[\begin{aligned}
			&[P(\eta_{19})]\sigma_{18}=[\sigma^2_9]\eta^2_{23}+2a[\bar{\sigma}_6]_{10} +b[\iota_7]_{10}\eta_7\bar{\mu}_8\ \mbox{for}\ a\in\mathbb{Z}\ \mbox{and}\ b\in\{0, 1\}.
		\end{aligned}\] \noindent  We know that $J[P(\eta_{19})]=\pm\sigma_{10}\zeta_{17}$ \cite[p. 343]{KM1}.
		In light of \eqref{m1:psi10}, \eqref{m1:ps10et2} and \eqref{m1:[bsg69}, we obtain
		$0=
		2\sigma_{10}\nu_{17}^{\ast}\pm 2a\sigma_{10}\nu_{17}^{\ast}+b\sigma_{10}\eta_{17}\bar{\mu}_{18}.$  This implies that $a$ is odd and $b=0.$ 
		Taking $x=a+8y,$   which is odd, gives the stated formula. This completes the proof.
	\end{proof}
	
	Since $\kappa_{11}$ is a suspension and $2\kappa_{11}=0$ by Appendix~\ref{app:stem-14}, we have $([2\e_{6},\e_{6}])\circ\kappa_{11}=P(\iota_{13})\circ(2\kappa_{11})=0$. Hence we deduce $\langle
	[2\e_{6},\e_{6}]\rangle\circ\kappa_{11}\in{i_G}_*\pi_{25}(\mathit{SU}(3))=0$
	\eqref{m1:iGU25}. Therefore, by \eqref{m1:hl2Pi13} and Lemma~\ref{m1:mtb2}(1), we obtain
	\begin{equation}
		\begin{aligned}
			&\langle [2\e_{6},\e_{6}]\rangle\circ\kappa_{11}=0,\ [[2\e_{6},\e_{6}]]\circ\kappa_{11}=0.\label{m1:2P13kp11}
		\end{aligned}
	\end{equation}
	\noindent
	Recall that $\Delta(\iota_{11}) = [\iota_7]_{11}\nu_7.$ 
	By  \eqref{m1:Det10}, we find
	\begin{align*}
		{p_{11}}_*\{[\iota_7]_{11}\nu_7, 2\iota_{10}, \kappa_{10}\} &= - \{p_{11},[\iota_7]_{11}\nu_7, 2\iota_{10}\}\circ\kappa_{11} \\
		&\supseteq - \{p_{11}, i_{11},\Delta(\eta_{10})\}\circ\kappa_{11} \\
		&\ni\eta_{10}\kappa_{11} \\
		&=\bar{\varepsilon}_{10} \ \bmod \ {p_{11}}_*\pi_{11}(\mathit{SO}(11))\circ\kappa_{11}.
	\end{align*}    \noindent   Since $\pi_{11}(\mathit{SO}(11))=\mathrm{span}\{[[2\e_{6},\e_{6}]]_{11}\}\cong\mathbb Z_{(2)},$   we conclude that
	$\pi_{11}(\mathit{SO}(11))\circ\kappa_{11}=0$ \eqref{m1:2P13kp11}.
	Therefore, we can take
	\begin{equation}
		\begin{aligned}
			&[\bar\varepsilon_{10}]\in\{\Delta(\iota_{11}), 2\iota_{10}, \kappa_{10}\}.\label{m1:bep10}
		\end{aligned}
	\end{equation}
	\noindent       Note that $2[\bar\varepsilon_{10}]\in -\{2\iota_{\mathit{SO}(11)},\Delta(\iota_{11}),
	2\iota_{10}\}\circ\kappa_{11}\subseteq \pi_{11}(\mathit{SO}(11))\circ\kappa_{11}=0.$ 
	That is,
	\begin{equation}
		\begin{aligned}
			&\mathrm{ord}[\bar\varepsilon_{10}]=2.\label{m1:bep102}
		\end{aligned}
	\end{equation}
	\noindent   
	It follows that $\sigma_{10}\varepsilon_{17}=P(\nu^2_{21}).$ 
	By the relation $\Delta([\e_{10},\nu_{10}])=0$ \cite[(4.4)]{HW53}, we can set
	\begin{equation}\label{hecanset}
		[\sigma_{10}\varepsilon_{17}]=[[\e_{10},\nu_{10}]]\nu_{22}.
	\end{equation}
	We show the following.
	\begin{prop}\label{m1:25-11}
		The structure of $\pi_{25}(\mathit{SO}(11))$ is as follows.
		\[ \GroupLayout{\pi_{25}(\mathit{SO}(11)) = \mathrm{span}\{[\bar\sigma_6]_{11},\; [[\e_{10},\nu_{10}]]\nu_{22},\; [\iota_7]_{11}\eta_7\bar\mu_8,\; [\bar\varepsilon_{10}]\}\cong (\mathbb Z/4)^{2}\oplus (\mathbb Z/2)^{2}.}{\begin{aligned}\pi_{25}(\mathit{SO}(11)) = \mathrm{span}\{&[\bar\sigma_6]_{11},\; [[\e_{10},\nu_{10}]]\nu_{22},\; [\iota_7]_{11}\eta_7\bar\mu_8,\\
				&[\bar\varepsilon_{10}]\}\cong (\mathbb Z/4)^{2}\oplus (\mathbb Z/2)^{2}.\end{aligned}} \]
		Moreover,  $2[[\e_{10},\nu_{10}]]\nu_{22} \equiv [\sigma_9\eta_{16}\varepsilon_{17}]_{11}\m 2[\bar\sigma_6]_{11},\;[\iota_7]_{11}\eta_7\bar\mu_8.$ 
	\end{prop}
	\begin{proof}
		We recall\ $\mathrm{ord}[\bar\varepsilon_{10}]=2$ \eqref{m1:bep102}.
		Consider  the exact sequence
		\[\begin{aligned}
			&\pi^{10}_{26} \xrightarrow{\Delta} \pi_{25}(\mathit{SO}(10)) \xrightarrow{i_\ast} \pi_{25}(\mathit{SO}(11)) \xrightarrow{p_\ast} \pi^{10}_{25} \xrightarrow{\Delta} \pi_{24}(\mathit{SO}(10)),
		\end{aligned}\]
		where
		$
		\pi^{10}_{26} = \frac{\mathbb Z}{16}\{P(\sigma_{21})\} \oplus
		\frac{\mathbb Z}{2}\{\sigma_{10}\mu_{17}\}.$ 
		By \cite[p.~343]{KM1}, (\ref{m1:r1}) and \eqref{m1:4bs}, we see that
		\[\begin{aligned}
			&\Delta(\mu_{10}\sigma_{19})=2[\iota_7]_{10}\zeta_7\sigma_{18}+4[[2\e_{6},\e_{6}]]_{10}\sigma^2_{11} = 16[\bar{\sigma}_6]_{10}+16[\bar{\sigma}_6]_{10}=0.
		\end{aligned}\]
		\noindent Recall from Lemma \ref{m1:DPs} that $$
		\Delta(P(\sigma_{21})) = [\sigma^2_9]\eta^2_{23} + 2x[\bar\sigma_6]_{10}
		$$ for some odd $x.$  
		Hence, $$i_{10,11_{*}}\pi_{25}(\mathit{SO}(10))=\mathrm{span\n
			[\ccs_{6}]_{11},
			[\e_{7}]_{11}\ca_{7}\cq_{8}, [\cg_{9}\ca_{16}\lt_{17}]_{11} 
			\nn\cong\z/4\jia(\z/2})^{2}.$$ \noindent
		In light of  the calculation of $\pi_{24}(\mathit{SO}(11)),$   we know that $$\operatorname{ker} (\Delta \colon
		\pi^{10}_{25} \to \pi_{24}(\mathit{SO}(10)))=\mathrm{span}\{\sigma_{10}\bar\nu_{17},
		\bar\varepsilon_{10}\}.$$
		\noindent        Recall  from Proposition \ref{m1:25-10}  that $$[\cg_{9}\lt_{16}\ca_{24}]+8[\ccs_{6}]_{10}\ty\n [P(\Sigma\theta)+\nu_{6}\cn_{9}]_{10},\;2\e_{23},\ca_{23}\nn\m[\cg_{9}^{2}]\ca_{23}^{2}.$$ \noindent  Equivalently, for some $t\in\mathbb Z,$   we have
		$$[\sigma_9\varepsilon_{16}\eta_{24}]_{11}
		+2t[\bar\sigma_6]_{11}
		\in\{[P(\Sigma\theta)+\nu_6\kappa_9]_{11},
		2\iota_{23},\eta_{23}\}.$$\noindent Lemma  \ref{m1:fin} tells us that  $[P(\Sigma\theta)+\nu_6\kappa_9]_{11}
		=[[\e_{10},\nu_{10}]]\eta_{22}.$ 
		Then,   by  
		$\langle\eta,2\iota,\eta\rangle\ni2\nu\bmod4\nu,$  
		we have\vspace{-1em}
		\begin{align*}
			\{[P(\Sigma\theta)+\nu_6\kappa_9]_{11},
			2\iota_{23},\eta_{23}\}&=\{[[\e_{10},\nu_{10}]]\eta_{22},
			2\iota_{23},\eta_{23}\}\\
			&\supseteq [[\e_{10},\nu_{10}]]\circ
			\{\eta_{22},2\iota_{23},\eta_{23}\}\\
			&\ni 2[[\e_{10},\nu_{10}]]\nu_{22}\bmod I.
		\end{align*}
		Here,
		$I={i_{10,11}}_*\pi_{24}(\mathit{SO}(10))\circ\eta_{24}
		\xyd\mathrm{span}\{[\iota_7]_{11}\eta_7\bar\mu_8,
		2[\bar\sigma_6]_{11}\}$.
		To verify this inclusion, use the generators of $\pi_{24}(\mathit{SO}(11))$ in Proposition~\ref{m1:11}. Their composites with $\eta_{24}$ satisfy
		\[
		\begin{aligned}
			[\bar\varepsilon_9]_{11}\eta_{24}
			&\in\mathrm{span}\{[\iota_7]_{11}\eta_7\bar\mu_8\},\\
			[\sigma_9^2]_{11}\eta_{23}^2
			&=-2x[\bar\sigma_6]_{11},\\
			[\iota_7]_{11}\bar\mu_7\eta_{24}
			&=[\iota_7]_{11}\eta_7\bar\mu_8,\\
			[P(\Sigma\theta)+\nu_6\kappa_9]_{11}\eta_{23}^2
			&=16[\bar\sigma_6]_{11}=0.
		\end{aligned}
		\]
		The first line follows from Lemma~\ref{m1:[bep9]11eta}(2).
		The second follows by applying $(i_{11})_*$ to Lemma~\ref{m1:DPs} and using exactness, with the same odd $x$.
		The third uses the Stem~17 relation in Appendix~\ref{app:oguchi-relations}.
		The fourth follows from Proposition~\ref{m1:R259}(1) and the order $4$ of $[\bar\sigma_6]_{11}$ established above.
		Since ${i_{10,11}}_*\pi_{24}(\mathit{SO}(10))$ is a subgroup of $\pi_{24}(\mathit{SO}(11))$, these four calculations give the asserted inclusion for $I$.
		Consequently,
		$$2[[\e_{10},\nu_{10}]]\nu_{22}
		\equiv[\sigma_9\varepsilon_{16}\eta_{24}]_{11}
		\bmod\mathrm{span}\{[\iota_7]_{11}\eta_7\bar\mu_8,
		2[\bar\sigma_6]_{11}\}.$$
		\noindent   We may choose the lift
		$[\eta_{10}\kappa_{11}]$ so that
		$[\eta_{10}\kappa_{11}]
		\in\{\Delta(\iota_{11}),2\iota_{10},\kappa_{10}\}.$ 
		It then follows that
		\begin{align*}
			2[\eta_{10}\kappa_{11}]
			&\in\{\Delta(\iota_{11}),2\iota_{10},0\}\\
			&=\{[\iota_7]_{11}\nu_7,2\iota_{10},0\}\\
			&=[\iota_7]_{11}\nu_7\circ
			\mathrm{span}\{\Sigma\rho',\eta_{10}\kappa_{11}\}\\
			&=0.
		\end{align*}\vspace{-1.5em} Hence the result holds.
	\end{proof}
	\subsection{The Case $n=12$}
	
	A Toda bracket now provides a representative whose projection can be compared with the desired element of a homotopy group of a sphere. The calculation first establishes that projection and then determines the relation satisfied by
	the resulting lift.
	
	Since $\Delta(\sigma_{11})=[\iota_7]_{11}\nu_7\sigma_{10}=x[\iota_7]_{11}\sigma'\nu_{14}$ for some odd $x,$ we consider the Toda bracket $\{[\iota_7]_{12}\sigma',\nu_{14},\sigma_{17}\}.$
	By Proposition~\ref{mdsyq}, $\{p_{12},i_{12},\Delta(\iota_{11})\}\ni\pm\iota_{11}.$ Composing with $\sigma_{11}$ and then with $\sigma_{18},$ and using Proposition~\ref{prop:toda-composition}(2),(3), gives
	\[\begin{aligned}
		{p_{12}}_*\{x[\iota_7]_{12}\sigma',\nu_{14},\sigma_{17}\} &= -\{p_{12},x[\iota_7]_{12}\sigma',\nu_{14}\}\circ\sigma_{18} \\
		&\supseteq -\{p_{12},i_{12},\Delta(\sigma_{11})\}\circ\sigma_{18} \\
		&\ni\pm\sigma^2_{11}.
	\end{aligned}\]
	Multiplication by the odd integer $x$ is an equivalence after localization at 2. Naturality of the bracket under this equivalence, including its indeterminacy, shows that the projection of $\{[\iota_7]_{12}\sigma',\nu_{14},\sigma_{17}\}$ contains an odd multiple of $\sigma_{11}^2.$ Choose an odd integer $t_0$ that inverts this multiple modulo the order of $\sigma_{11}^2.$ By right-composition naturality, we can therefore choose a lift with ${p_{12}}_*[\sigma_{11}^2]=\sigma_{11}^2$ such that
	\begin{equation}
		\begin{aligned}
			&[\sigma^2_{11}]\in\{[\iota_7]_{12}\sigma',\nu_{14},t_{0}\sigma_{17}\}.\label{m1:sgm211}
		\end{aligned}
	\end{equation}
	With the relevant preliminary results established,
	we may now present the outcome of the corresponding
	calculation.
	\begin{prop}\label{m1:25-12}
		The structure of $\pi_{25}(\mathit{SO}(12))$ is as follows.
		\[ \GroupLayout{\pi_{25}(\mathit{SO}(12)) = \mathrm{span}\{[\sigma^2_{11}],\; [[\e_{10},\nu_{10}]]_{12}\nu_{22},\; [\bar\sigma_6]_{12},\; [\iota_7]_{12}\eta_7\bar\mu_8,\; [\bar\varepsilon_{10}]_{12}\}\cong \mathbb Z/16\oplus (\mathbb Z/4)^{2}\oplus (\mathbb Z/2)^{2}.}{\begin{aligned}\pi_{25}(\mathit{SO}(12)) = \mathrm{span}\{&[\sigma^2_{11}],\; [[\e_{10},\nu_{10}]]_{12}\nu_{22},\; [\bar\sigma_6]_{12},\\
				&[\iota_7]_{12}\eta_7\bar\mu_8,\; [\bar\varepsilon_{10}]_{12}\}\cong \mathbb Z/16\oplus (\mathbb Z/4)^{2}\oplus (\mathbb Z/2)^{2}.\end{aligned}} \]
		Here $\Delta(P(\nu_{25}))=[[\e_{10},\nu_{10}]]_{12}\nu_{22}$ and
		$\Delta(\kappa_{12})=[\bar\varepsilon_{10}]_{12}.$ 
	\end{prop}
	
	\begin{proof}
		By \cite[Lemma 6.3(ii)]{HKM}, we deduce $\Delta([\iota_{12},\nu_{12}])
		=[[\e_{10},\nu_{10}]]_{12}\nu_{22}.$    Consider  the exact sequence
		\[\begin{aligned}
			&\pi^{11}_{26} \xrightarrow{ \Delta} \pi_{25}(\mathit{SO}(11)) \xrightarrow{i_\ast } \pi_{25}(\mathit{SO}(12)) \xrightarrow{p_\ast } \pi^{11}_{25} \xrightarrow{\Delta } \pi_{24}(\mathit{SO}(11)),
		\end{aligned}\]
		where $ \pi^{11}_{26} = \frac{\mathbb Z}{16}\{\Sigma^2\rho'\} \oplus
		\frac{\mathbb Z}{2}\{\bar\varepsilon_{11}\}. $ Using the standard vanishing relation,
		\cite[Table 3]{Ka},
		we see that
		\begin{align*}
			&\Delta(\Sigma^2\rho')= [\iota_7]_{11}\nu_7(\Sigma\rho') =0 \qquad\text{and}\qquad \Delta(\bar\varepsilon_{11})=\Delta(\eta_{11})\kappa_{11} = 0.
		\end{align*}
		\noindent       Recall that
		$ \Delta(\kappa_{11})=2[\bar\varepsilon_9]_{11} $ from the calculation of $\pi_{24}(\mathit{SO}(12)).$ 
		By \eqref{m1:sgm211},
		we have
		\[\begin{aligned}
			&16[\sigma^2_{11}] \in - [\iota_7]_{12}\sigma'\circ\{\nu_{14}, t_{0}\sigma_{17},16\iota_{24}\}\ni -t_{0} [\iota_7]_{12}\sigma'\zeta_{14} = 8t_{0}[\bar\sigma_6]_{12} = 0 \ \bmod\ 0.
		\end{aligned}\]
		\noindent       
		In view of \eqref{m1:bep10} and \eqref{m1:2P13kp11}, we infer that
		\[\begin{aligned}
			[\bar{\varepsilon}_{10}]_{12} &\in {i_{12}}_*\{\Delta(\iota_{11}),2\iota_{10},\kappa_{10}\} \\
			&= - \{i_{12},\Delta(\iota_{11}), 2\iota_{10}\}\circ\kappa_{11} \ni\Delta(\iota_{12})\circ\kappa_{11} \\
			&= \Delta(\kappa_{12})\ \bmod\ {i_{12}}_*\pi_{11}(\mathit{SO}(11))\circ\kappa_{11}=0.
		\end{aligned}\]
		This implies $[\bar{\varepsilon}_{10}]_{12}=\Delta(\kappa_{12})$\ and
		completes the proof.
	\end{proof}
	
	Employing  \cite[Theorem 1(a), I-Proposition 3.5(3)]{Oda1979}, we obtain
	\begin{equation}
		\begin{aligned}
			&\mathrm{ord}(\xi_{12}\sigma_{30})=8,\ 2(\xi_{12}\sigma_{30})=\sigma_{12}\xi_{19}\ \ \mbox{and}\ \ 2\sigma_{12}\nu^{\ast}_{19}=2\sigma_{12}\xi_{19}.\label{m1:2xs12}
		\end{aligned}
	\end{equation}
	\noindent From \cite[II-Proposition 2.1(6)]{Oda1979}, we know
	\begin{equation}
		\begin{aligned}
			&\xi'\sigma_{29}\equiv\sigma_{11}\nu_{18}^{\ast}\ \bmod\ 2\sigma_{11}\nu_{18}^{\ast}.\label{m1:x's29}
		\end{aligned}
	\end{equation}  \noindent  According to  \cite[Lemma 3.5(1)]{MiM}, $H(P(\lambda))=x\lambda'+2y\xi',$   where
	$x$ is odd. In light of \cite[Lemma 12.19]{T}, we have  $\Sigma^2\lambda'=2\lambda$ and
	$\Sigma^2\xi'=2\xi_{13}.$ 
	
	We show

	\begin{lem}$\lambda'\sigma_{29}=0.$ \end{lem}
	\begin{proof}
		By \cite[III-Proposition 2.2(2)]{Oda1979}, it follows that $\lambda\sigma_{31}=0.$ 
		Therefore, we obtain
		\[\begin{aligned}
			&0=x\lambda'\sigma_{29}+2y\xi'\sigma_{29}.
		\end{aligned}\]
		By suspending this relation twice, we conclude that
		\[\begin{aligned}
			&0=2x\lambda\sigma_{31}+4y\xi_{13}\sigma_{31}=4y\xi_{13}\sigma_{31}.
		\end{aligned}\]
		Since $\mathrm{ord}(\xi_{13}\sigma_{31})=8$ \cite[Theorem 1(a)]{Oda1979},
		$y$ is even. We see that $\xi'\sigma_{29}\equiv\sigma_{11}\nu_{18}^{\ast}\ \bmod\
		2\sigma_{11}\nu_{18}^{\ast}$ \eqref{m1:x's29}. Therefore, by
		\eqref{m1:ps10et2}, we see that $4\xi'\sigma_{29}=0.$  This implies the desired relation.
	\end{proof}
	
	Since $HJ[\sigma^2_{11}]=\sigma^2_{23}=H(\xi_{12}\sigma_{30}),$   we know
	$J[\sigma^2_{11}] \equiv \xi_{12}\sigma_{30} \bmod \Sigma\pi^{11}_{36}.$  In light of  \eqref{m1:sgm211}, we obtain\vspace{-0.8em}
	\begin{align}
		\pm J[\sigma^2_{11}]
		&\in \{2\sigma^2_{12},\nu_{26},t_{0}\sigma_{29}\}_1 \notag\\
		&\supseteq 2\sigma_{12}\circ
		\{\sigma_{19},\nu_{26},t_{0}\sigma_{29}\}_1 \notag\\
		&\ni 2t_{0}\sigma_{12}\xi_{19}
		\quad\bmod\ 
		2\sigma^2_{12}\circ\Sigma\pi^{25}_{36}
		+\pi^{12}_{30}\circ\sigma_{30}, \notag
	\end{align}
	\noindent where $t_{0}\in\z$ is some odd number.
	Since $\sigma_{19}\zeta_{26}=0$ \eqref{m1:s10z}, we deduce $2\sigma^2_{12}\circ \Sigma\pi^{25}_{36}=0.$  From the fact that
	$$\begin{aligned}
		&\pi^{12}_{30}=\mathrm{span}\{\xi_{12}, \Sigma\lambda', \Sigma\xi', \eta_{12}\bar{\mu}_{13}\},
	\end{aligned}$$  we obtain
	$\pi^{12}_{30}\circ\sigma_{30}=\mathrm{span}\{\xi_{12}\sigma_{30}, (\Sigma\xi')\sigma_{30}, \eta_{12}\bar{\mu}_{13}\sigma_{30}\}.$  This implies
	\begin{equation}
		\begin{aligned}
			&J[\sigma^2_{11}] \equiv \xi_{12}\sigma_{30} \bmod\ \sigma_{12}\nu_{19}^{\ast}, \eta_{12}\bar{\mu}_{13}\sigma_{30}, 2\xi_{12}\sigma_{30}.\label{m1:Jsg^211}
		\end{aligned}
	\end{equation}
	We shall also need the following congruence
	in the subsequent calculations.
	\begin{lem}\label{m1:1112}
		$\Delta (\sigma^2_{12})\equiv2[\sigma^2_{11}] + x [\bar\sigma_{6}]_{12}\ \bmod \ \Delta(P(\nu_{25})),$   where $x$ is odd.
	\end{lem}
	\begin{proof}
		Since $p_{12}(\Delta(\sigma^2_{12})-2[\sigma^2_{11}])=0,$  
		we find $\Delta(\sigma^2_{12}) - 2[\sigma^2_{11}] \in
		i_{12 \ast}\pi_{25}(\mathit{SO}(11)).$  Next, by the relation $P(\iota_{39})=\nu_{19}^{\ast}+\xi_{19}$
		\cite[Corollary 12.25]{T}, we obtain $$J\Delta(\sigma^2_{12})=  \pm  [\iota_{12},
		\sigma^2_{12}]=  \pm  \sigma_{12}[\iota_{19},
		\iota_{19}]=  \pm  \sigma_{12}(\xi_{19}+\nu^{\ast}_{19}).$$  Propositions \ref{m1:25-11} and \ref{m1:25-12} imply $$\Delta(\sigma^2_{12}) =
		2[\sigma^2_{11}] +x[\bar\sigma_6]_{12}+y[[\e_{10},\nu_{10}]]_{12}\nu_{22} +z\Delta(\kappa_{12})
		+w[\iota_7]_{12}\eta_{7}\bar\mu_{8},$$\noindent where $x, y \in \{0, 1, 2, 3\}$
		and $z, w \in \{0, 1\}.$ 
		From the fact that $J\Delta(P(\nu_{25}))=0$ and
		that    $J\Delta(\kappa_{12})=P(\kappa_{25}),$   we deduce  $$  \pm  \sigma_{12}(\xi_{19}+\nu^{\ast}_{19})=J\Delta(\sigma^2_{12}) =
		2J([\sigma^2_{11}]) +x J([\bar\sigma_6]_{12})+zP(\kappa_{25})
		+w\eta_{12}\bar\mu_{13}\sigma_{30}.$$
		\noindent
		In view of  \cite[I-Proposition 4.3(1)]{Oda1979}, we conclude that $P(\kappa_{25})=\Sigma K_1.$ 
		By \eqref{m1:[bsg69}, \eqref{m1:2xs12} and \eqref{m1:Jsg^211}, we see that
		\[\begin{aligned}
			&\sigma_{12}(\xi_{19}+\nu^{\ast}_{19})=a\sigma_{12}\xi_{19}\pm x\sigma_{12}\nu_{19}^{\ast}+z\Sigma K_1 +w\eta_{12}\bar\mu_{13}\sigma_{30}\ \mbox{for  some  odd}\ a.
		\end{aligned}\]
		Employing the  structure of $\pi^{12}_{37},$   we obtain $z=w=0$ and find that $x$ is odd.
		This completes the proof.   \end{proof}
	
	\subsection{The Cases $n=13,14$}
	
	The next two cases are described using the images of the existing generators. The relation following
	the decomposition specifies how the $\bar\sigma$-class is absorbed into a multiple of the chosen
	$\sigma^2$-lift.
	
	Let us now turn to the following statement and present its proof before proceeding with the next part of our discussion.
	\begin{prop}\label{g1016:prop:25-13}
		The structure of $\pi_{25}(\mathit{SO}(13))$ is as follows.
		\[ \GroupLayout{\pi_{25}(\mathit{SO}(13)) = \mathrm{span}\{[\iota_7]_{13}\eta_7\bar{\mu}_8,\; [\sigma^2_{11}]_{13}\}\cong \mathbb Z/2\oplus \mathbb Z/8.}{\begin{aligned}\pi_{25}(\mathit{SO}(13)) = \mathrm{span}\{&[\iota_7]_{13}\eta_7\bar{\mu}_8,\\
				&[\sigma^2_{11}]_{13}\}\cong \mathbb Z/2\oplus \mathbb Z/8.\end{aligned}} \]
		Here\ $2[\sigma^2_{11}]_{13}=-x[\bar\sigma_6]_{13}$\ and\ $x$ is the same odd integer as in Lemma~\ref{m1:1112}.
	\end{prop}
	
	\begin{proof}
		
		Consider  the exact sequence
		\[\begin{aligned}
			&\pi^{12}_{26} \xrightarrow{ \Delta} \pi_{25}(\mathit{SO}(12)) \xrightarrow{i_\ast } \pi_{25}(\mathit{SO}(13)) \xrightarrow{p_\ast } \pi^{12}_{25} \xrightarrow{ \Delta} \pi_{24}(\mathit{SO}(12)),
		\end{aligned}\]
		where
		$\pi^{12}_{26}=\frac{\mathbb Z}{16}\{\sigma^2_{12}\}\oplus\frac{\mathbb Z}{4}\{P(\nu_{25})\}\oplus\frac{\mathbb Z}{2}\{\kappa_{12}\}.$ 
		In light of  the calculation of $\pi_{24}(\mathit{SO}(13)),$   we know that  $$\Delta : \pi^{12}_{25} \to \pi_{24}(\mathit{SO}(12))$$ is a
		monomorphism. As a consequence of  Proposition~\ref{m1:25-12}, we find 
		$\Delta(P(\nu_{25}))=[[\e_{10},\nu_{10}]]_{12}\nu_{22}.$ 
		From Proposition~\ref{m1:25-12}, we infer 
		$\Delta(\kappa_{12})=[\bar\varepsilon_{10}]_{12}.$ 
		Then,   by Lemma~\ref{m1:1112}, we obtain the result.   \end{proof}
	The preceding discussion prepares the way
	for the following determination.
	\begin{prop}
		The structure of $\pi_{25}(\mathit{SO}(14))$ is as follows. $$\GroupLayout{\pi_{25}(\mathit{SO}(14)) =\mathrm{span}\{[\iota_7]_{14}\eta_7\bar{\mu}_8,\; [\sigma_{11}^2]_{14}\}\cong \mathbb Z/2\oplus \mathbb Z/8.}{\begin{aligned}\pi_{25}(\mathit{SO}(14)) =\mathrm{span}\{&[\iota_7]_{14}\eta_7\bar{\mu}_8,\\
				&[\sigma_{11}^2]_{14}\}\cong \mathbb Z/2\oplus \mathbb Z/8.\end{aligned}}$$
	\end{prop}
	\begin{proof}
		In the exact sequence $(25)_{13}$
		\[\begin{aligned}
			&\pi^{13}_{26} \xrightarrow{ \Delta} \pi_{25}(\mathit{SO}(13)) \xrightarrow{{i_{14}}_\ast } \pi_{25}(\mathit{SO}(14)) \xrightarrow{{p_{14}}_\ast } \pi^{13}_{25}\xrightarrow{ \Delta}\pi_{24}(\mathit{SO}(13)),
		\end{aligned}\]
		recall that ${p_{14}}_\ast$ is trivial since the $\Delta$ on the right is a monomorphism.
		Consider the equation
		$p_{14 \ast}\Delta(P(\iota_{29}))=P(\eta_{27})=(\Sigma\theta)\eta_{25} \neq
		0$~(\cite[Lemma 1.4]{HKM})\uprightsemicolon{}  $\Delta$ on the left
		is trivial and hence,
		${i_{14}}_{\ast} : \pi_{25}(\mathit{SO}(13)) \to \pi_{25}(\mathit{SO}(14))$ is a monomorphism.
		This completes the proof.
		
	\end{proof}
	\section{The Cases $n\geq15$}
	
	\subsection{The Cases $n=15,16$}
	
	For $n=15,16,$   the inherited generators are joined by a composite involving the chosen
	lift of $\eta_{14}^{2}.$  We determine the decomposition while retaining these representatives for the
	next inclusions.
	
	\begin{prop}
		The structure of $\pi_{25}(\mathit{SO}(15))$ is as follows.        $$\GroupLayout{\pi_{25}(\mathit{SO}(15))=\mathrm{span}\{[\iota_7]_{15}\eta_7\bar{\mu}_8,\; [\eta^2_{14}]\mu_{16},\; [\sigma^2_{11}]_{15}\}\cong (\mathbb Z/2)^{2}\oplus \mathbb Z/8.}{\begin{aligned}\pi_{25}(\mathit{SO}(15))=\mathrm{span}\{&[\iota_7]_{15}\eta_7\bar{\mu}_8,\; [\eta^2_{14}]\mu_{16},\\
				&[\sigma^2_{11}]_{15}\}\cong (\mathbb Z/2)^{2}\oplus \mathbb Z/8.\end{aligned}}$$
	\end{prop}
	\begin{proof}
		Consider  the exact sequence
		\[\begin{aligned}
			&0= \pi^{14}_{26} \xrightarrow{ \Delta} \pi_{25}(\mathit{SO}(14)) \xrightarrow{i_\ast } \pi_{25}(\mathit{SO}(15)) \xrightarrow{p_\ast } \pi^{14}_{25} \xrightarrow{ \Delta} \pi_{24}(\mathit{SO}(14)).
		\end{aligned}\] \noindent   In light of  Lemma~\ref{m1:zeta1}, we have
		$\Delta(\zeta_{14})=\pm 2[\zeta_{13}].$ 
		By the calculation of $\pi_{24}(\mathit{SO}(15)),$   we know that $\operatorname{ker}
		(\pi^{14}_{25} \to \pi_{24}(\mathit{SO}(14))) =
		\mathrm{span}\{4\zeta_{14}\}=\mathrm{span}\{\eta^2_{14}\mu_{16}\} \cong \mathbb Z/2.$  Therefore, there exists
		$[\eta^2_{14}]\mu_{16} \in \pi_{25}(\mathit{SO}(15)).$  \end{proof}
	
	Using the fact that $\pi_{25}(\mathit{SO})\cong \mathbb Z/2$ \cite{Bo} and the results in \cite[Table 1]{HoM},
	we have Table~\ref{m1:R25}.
	\par
	\begin{table}[H]
		\centering
		\small\setlength{\tabcolsep}{3pt}\begin{tabular}{cccccc}\toprule
			$\pi_{25}(\mathit{SO}(15))$ & $\mathbb Z/2 \oplus \mathbb Z/2 \oplus \mathbb Z/8$
			& $\pi_{25}(\mathit{SO}(16))$ & $\mathbb Z/2 \oplus (\mathbb Z/2)^2 + \mathbb Z/8$
			& $\pi_{25}(\mathit{SO}(17))$ & $\mathbb Z/2 \oplus \mathbb Z/8 \oplus \mathbb Z/2$  \\ \midrule
			$\pi_{25}(\mathit{SO}(18))$ & $\mathbb Z/2 \oplus \mathbb Z/4$
			& $\pi_{25}(\mathit{SO}(19))$ & $\mathbb Z/2 \oplus \mathbb Z/2$
			& $\pi_{25}(\mathit{SO}(20))$ & $\mathbb Z/2$\\ \bottomrule
		\end{tabular}
		\caption{$\pi_{25}(\mathit{SO}(n))$ for $15\leq n\leq 20$}\label{m1:R25}
	\end{table}
	\par
	By \cite{Ke}, we          have Table~\ref{m1:R25K}.
	\par
	\begin{table}[H]
		\centering
		\small\setlength{\tabcolsep}{3pt}\begin{tabular}{cccccc}\toprule
			$\pi_{25}(\mathit{SO}(21))$ & $\mathbb Z/2$
			& $\pi_{25}(\mathit{SO}(22))$ & $\mathbb Z/2$
			& $\pi_{25}(\mathit{SO}(23))$ & $(\mathbb Z/2)^2$  \\ \midrule
			$\pi_{25}(\mathit{SO}(24))$ & $(\mathbb Z/2)^3$
			& $\pi_{25}(\mathit{SO}(25))$ & $(\mathbb Z/2)^2$
			& $\pi_{25}(\mathit{SO}(26))$ & $\mathbb Z_{(2)} \oplus \mathbb Z/2$\\ \bottomrule
		\end{tabular}
		\caption{$\pi_{25}(\mathit{SO}(n))$ for $21\leq n\leq 26$}\label{m1:R25K}
	\end{table}
	\par
	
	By Tables~\ref{m1:R25} and  \ref{m1:R25K}, we have the group structure of $\pi_{25}(\mathit{SO}(n))$ for $n \ge 15.$ 
	We now determine the generators and relations in
	$\pi_{25}(\mathit{SO}(n))$ for $n \ge 15.$ 
	
	\begin{prop}\label{m1:16-25}
		The structure of $\pi_{25}(\mathit{SO}(16))$ is as follows.
		\[ \GroupLayout{\pi_{25}(\mathit{SO}(16)) = \mathrm{span}\{[\nu_{15}]\sigma_{18},\; [\eta_{15}]\mu_{16},\; [\eta^2_{14}]_{16}\mu_{16},\; [\iota_7]_{16}\eta_7\bar{\mu}_8\} \cong \mathbb Z/8 \oplus (\mathbb Z/2)^3.}{\begin{aligned}\pi_{25}(\mathit{SO}(16)) = \mathrm{span}\{&[\nu_{15}]\sigma_{18},\; [\eta_{15}]\mu_{16},\; [\eta^2_{14}]_{16}\mu_{16},\\
				&[\iota_7]_{16}\eta_7\bar{\mu}_8\} \cong \mathbb Z/8 \oplus (\mathbb Z/2)^3.\end{aligned}} \]
		Here,  $[\sigma^2_{11}]_{16}=a[\nu_{15}]\sigma_{18}$ for some odd $a.$ 
	\end{prop}
	\begin{proof}
		Consider  the exact sequence
		\[\begin{aligned}
			&\pi^{15}_{26} \xrightarrow{ \Delta} \pi_{25}(\mathit{SO}(15)) \xrightarrow{i_\ast } \pi_{25}(\mathit{SO}(16)) \xrightarrow{p_\ast } \pi^{15}_{25} \xrightarrow{ \Delta} \pi_{24}(\mathit{SO}(15)),
		\end{aligned}\]
		where $\pi^{15}_{26}=\frac{\mathbb Z}{8}\{\zeta_{15}\}.$ 
		Since $\zeta_{15} \in\{\nu_{15}, 8\iota_{18}, 2\sigma_{18}\}_1\m0$
		\cite[p.~59]{T} and $8[\nu_{15}]=0,$   we may take $$[\zeta_{15}] \in \{[\nu_{15}], 8\iota_{18},
		2\sigma_{18}\}_1.$$     \noindent   Then, we conclude that  $\Delta(\zeta_{15})=\Delta({p_{16}}_*[\zeta_{15}])=0.$      Therefore, $i_{\ast}$ is a monomorphism.  We have a lift $[\eta_{15}]\mu_{16}$ of $\eta_{15}\mu_{16} \in
		\pi^{15}_{25}.$  We see that
		\[\begin{aligned}
			&[\nu_{15}]\sigma_{18}\in {i_{16}}_*\{\Delta(\iota_{15}),\nu_{14},\sigma_{17}\}.
		\end{aligned}\]
		Recall that $\Delta(\iota_{15})=[\iota_7]_{15}\sigma'.$  Therefore, by \eqref{m1:sgm211},  for some odd $t_{0}$, we obtain
		\[\begin{aligned}
			&[\nu_{15}]\sigma_{18}\in {i_{12,16}}_*\{[\iota_7]_{12}\sigma',\nu_{14}, t_{0}\sigma_{17}\}\ni[\sigma^2_{11}]_{16} \ \bmod\ {i_{16}}_*\pi_{18}(\mathit{SO}(15))\circ\sigma_{18}.
		\end{aligned}\]
		Note  that
		${i_{16}}_*\pi_{18}(\mathit{SO}(15))\circ\sigma_{18}=\mathrm{span}\{[\nu^2_{12}]_{16}\sigma_{18}\}.$ 
		By \cite[Lemma 3.5(1)]{GM}, we deduce $[\nu^2_{12}]_{16}\sigma_{18}= 2[\nu_{15}]\sigma_{18}.$  This
		completes the proof.    \end{proof}
	
	We find $H(J[\nu_{15}])=\nu_{31}=H(c\nu_{16}^{\ast}),$   where $c\in\mathbb Z$ is  odd. Therefore, we conclude that
	$$\begin{aligned}
		&J[\nu_{15}]\equiv c\nu_{16}^{\ast}\ \bmod\ \Sigma\pi^{15}_{33},
	\end{aligned}$$
	where \ $\pi^{15}_{33}=\mathrm{span}\{\xi_{15}, \Sigma^2\lambda, \eta_{15}\bar{\mu}_{16}\}
	\cong(\mathbb Z/8)^2\oplus\mathbb Z/2.$     Since $[\nu_{15}]_\infty\in \pi_{18}(\mathit{SO})=0$ \cite{Bo}, $\nu^{\ast}=-\xi$\
	and $\Sigma^\infty\lambda=2\nu^{\ast}$ \cite[Theorem 12.22]{T}, we obtain
	\[\begin{aligned}
		&J([\nu_{15}])=c\nu^{\ast}_{16}+x\xi_{16}+y\Sigma^3\lambda,
	\end{aligned}\]
	where $c-x+2y \equiv 0 \bmod 8,$   $c, x$ are odd  and $c, x, y \in[0,7]\cap \mathbb{Z}.$  Note that this equation is an improvement of
	\cite[Lemma 2.2 (ii), the last]{HKM}.
	
	\subsection{The Cases $17\leq n\leq20$}
	
	The larger-$n$ calculation uses lifts of $\nu_{15}$ and $\eta_{15}.$  Their composites give convenient
	representatives to follow through the cases $17\leq n\leq20.$ 
	
	We show the following.
	\begin{prop}
		The structure of $\pi_{25}(\mathit{SO}(17))$ is as follows.
		\[ \GroupLayout{\pi_{25}(\mathit{SO}(17)) = \mathrm{span}\{[\iota_7]_{17}\eta_7\bar{\mu}_8,\; [\nu_{15}]_{17}\sigma_{18},\; [\eta_{15}]_{17}\mu_{16}\}\cong \mathbb Z/2\oplus \mathbb Z/8\oplus \mathbb Z/2.}{\begin{aligned}\pi_{25}(\mathit{SO}(17)) = \mathrm{span}\{&[\iota_7]_{17}\eta_7\bar{\mu}_8,\; [\nu_{15}]_{17}\sigma_{18},\\
				&[\eta_{15}]_{17}\mu_{16}\}\cong \mathbb Z/2\oplus \mathbb Z/8\oplus \mathbb Z/2.\end{aligned}} \]
		Here, $[\eta_{15}]_{17}\mu_{16}=\Delta(\mu_{17}).$ 
	\end{prop}
	\begin{proof}
		
		We have  the exact sequence
		\[\begin{aligned}
			&\pi^{17}_{26} \xrightarrow{ \Delta} \pi_{25}(\mathit{SO}(17)) \xrightarrow{i_\ast } \pi_{25}(\mathit{SO}(18)) \xrightarrow{p_\ast } \pi^{17}_{25} \xrightarrow{ \Delta} \pi_{24}(\mathit{SO}(17)).
		\end{aligned}\]
		We see that
		$\pi^{17}_{26}= \frac{\mathbb Z}{2}\{\mu_{17}\} \oplus \frac{\mathbb Z}{2}\{\eta_{17}\varepsilon_{18}\} \oplus
		\frac{\mathbb Z}{2}\{\eta^2_{17}\sigma_{19}\}$ and $\Delta : \pi^{17}_{25} \to
		\pi_{24}(\mathit{SO}(17))$ is a monomorphism by the calculation of
		$\pi_{24}(\mathit{SO}(17))$ and $\pi_{24}(\mathit{SO}(18)).$ 
		By the $\eta^2_{8k+1}$ entry of Table~\ref{m2:tab:boundary-values} with $k=2$ and the boundary composition formula, we have
		$\Delta(\eta^2_{17}\sigma_{19})=4[\nu_{15}]_{17}\sigma_{18}.$ 
		
		To determine the generators and their orders, we also use the exact sequence
		\[
		\pi^{16}_{26}\xrightarrow{\Delta}
		\pi_{25}(\mathit{SO}(16))\xrightarrow{i_*}
		\pi_{25}(\mathit{SO}(17))\xrightarrow{p_*}
		\pi^{16}_{25}\xrightarrow{\Delta}
		\pi_{24}(\mathit{SO}(16)).
		\]
		The last homomorphism is injective by the proof of
		Proposition~\ref{m1:17}. Hence $i_*$ is surjective.
		The stem-$10$ calculation in Appendix~\ref{tlqdzj} gives
		$\pi^{16}_{26}=\frac{\mathbb Z}{2}\{\eta_{16}\mu_{17}\}$.
		By Table~\ref{m2:tab:boundary-values} and the boundary composition formula,
		\[
		\Delta(\eta_{16}\mu_{17})
		=\Delta(\eta_{16})\mu_{16}
		=[\eta_{14}^2]_{16}\mu_{16}.
		\]
		By Proposition~\ref{m1:16-25}, this element generates
		a direct summand of order $2$ in
		$\pi_{25}(\mathit{SO}(16))$. Exactness therefore gives
		\[
		\pi_{25}(\mathit{SO}(17))
		\cong\pi_{25}(\mathit{SO}(16))
		\big/\mathrm{span}\{[\eta_{14}^2]_{16}\mu_{16}\}.
		\]
		The images of the other three generators are
		$[\iota_7]_{17}\eta_7\bar\mu_8$,
		$[\nu_{15}]_{17}\sigma_{18}$ and
		$[\eta_{15}]_{17}\mu_{16}$, respectively.
		They generate independent cyclic summands of orders $2$, $8$ and $2$,
		which proves the stated decomposition.
		Finally, the $\iota_{4m+1}$ entry of
		Table~\ref{m2:tab:boundary-values}, with $m=4$, gives
		$\Delta(\iota_{17})=[\eta_{15}]_{17}$. Thus
		$\Delta(\mu_{17})=\Delta(\iota_{17})\mu_{16}
		=[\eta_{15}]_{17}\mu_{16}$.	\end{proof}
	We now bring together the relevant calculations
	to obtain the following description.
	\begin{prop}
		The structure of $\pi_{25}(\mathit{SO}(18))$ is as follows.
		\[ \GroupLayout{\pi_{25}(\mathit{SO}(18)) = \mathrm{span}\{[\iota_7]_{18}\eta_7\bar{\mu}_8,\; [\nu_{15}]_{18}\sigma_{18}\}\cong \mathbb Z/2\oplus \mathbb Z/4.}{\begin{aligned}\pi_{25}(\mathit{SO}(18)) = \mathrm{span}\{&[\iota_7]_{18}\eta_7\bar{\mu}_8,\\
				&[\nu_{15}]_{18}\sigma_{18}\}\cong \mathbb Z/2\oplus \mathbb Z/4.\end{aligned}} \]
		Moreover,  $\Delta(\eta^2_{17}\sigma_{19})=\Delta(\nu^3_{17})=4[\nu_{15}]_{17}\sigma_{18}.$ 
	\end{prop}
	\begin{proof} Note that $\Delta(\eta^2_{17}\sigma_{19})=4[\nu_{15}]_{17}\sigma_{18}.$  In light of Proposition~\ref{m1:16-25}, we obtain
		$$J([\sigma^2_{11}]_{18})\ty J([\nu_{15}]_{18}\sigma_{18})\ty \xi_{18}\sigma_{36}\m2\xi_{18}\sigma_{36}.$$\noindent
		This element has order $4$ by \cite[Theorem 1(a)]{Oda1979}.
		
		Now, we consider  the exact sequence
		\[\begin{aligned}
			&\pi^{17}_{26} \xrightarrow{ \Delta} \pi_{25}(\mathit{SO}(17)) \xrightarrow{i_\ast } \pi_{25}(\mathit{SO}(18)) \xrightarrow{p_\ast} \pi^{17}_{25} \xrightarrow{ \Delta} \pi_{24}(\mathit{SO}(17)),
		\end{aligned}\]
		where $\pi^{17}_{26}=\mathrm{span}\{\nu^3_{17}, \mu_{17}, \eta_{17}\varepsilon_{18}\}\cong(\z/2)^3.$ 
		By the calculation of $\pi_{24}(\mathit{SO}(17)),$   the right-hand homomorphism $\Delta$ is a monomorphism.
		Consequently, $i_*$ is surjective. The preceding proposition and
		Table~\ref{m1:R25} give
		\[
		|\operatorname{Im}\Delta|
		=\frac{|\pi_{25}(\mathit{SO}(17))|}
		{|\pi_{25}(\mathit{SO}(18))|}
		=\frac{32}{8}=4.
		\]
		The elements
		$\Delta(\mu_{17})=[\eta_{15}]_{17}\mu_{16}$ and
		$\Delta(\eta_{17}^2\sigma_{19})=4[\nu_{15}]_{17}\sigma_{18}$
		are independent elements of order $2$ by the preceding proposition. So, 
		\[
		\operatorname{Im}\Delta
		=\mathrm{span}\{[\eta_{15}]_{17}\mu_{16},\,
		4[\nu_{15}]_{17}\sigma_{18}\}
		\cong(\mathbb Z/2)^2.
		\]
		Taking the quotient of the preceding decomposition by this subgroup
		shows that $[\iota_7]_{18}\eta_7\bar\mu_8$ and
		$[\nu_{15}]_{18}\sigma_{18}$ generate independent cyclic summands
		of orders $2$ and $4$, respectively. This proves the asserted
		decomposition of $\pi_{25}(\mathit{SO}(18))$.
		
		It remains to examine $\Delta(\nu_{17}^3)$.
		We first prove that $J$ is injective on $\operatorname{Im}\Delta$.
		Consider the EHP sequence
		\[
		\pi^{35}_{44}\xrightarrow{P}\pi^{17}_{42}
		\xrightarrow{\Sigma}\pi^{18}_{43}.
		\]
		Recall that	$|\pi^{17}_{42}|=64$ and $|\pi^{18}_{43}|=16$.
		Moreover, the generators $\xi_{18}\sigma_{36}$,
		$\mu_{3,18}$ and $\eta_{18}\bar\mu_{19}\sigma_{36}$ of
		$\pi^{18}_{43}$ are the suspensions of
		$\xi_{17}\sigma_{35}$, $\mu_{3,17}$ and
		$\eta_{17}\bar\mu_{18}\sigma_{35}$, respectively.
		Thus $\Sigma$ is surjective and
		$|\operatorname{Im}P|=|\ker\Sigma|=64/16=4$.
		The homomorphism 
		$\Sigma^{18}\colon\pi^{17}_{26}\to\pi^{35}_{44}$
		is an isomorphism.
		The formula $J\Delta=\pm P\Sigma^{18}$  consequently implies
		\[
		J(\operatorname{Im}\Delta)=\operatorname{Im}P\text{\;\;and\;\;}
		|J(\operatorname{Im}\Delta)|=4=|\operatorname{Im}\Delta|.
		\]
		Therefore $\ker J\cap\operatorname{Im}\Delta=0$.
		
		By our Section~\ref{sec:oguchi-working} from \cite{T},
		$[\iota_{17},\iota_{17}]
		\equiv\eta_{17}^*+\omega_{17}\bmod\sigma_{17}\mu_{24}$. 
		The formula for $J\Delta$ now yields
		\[
		J\Delta(\eta_{17}\varepsilon_{18})
		=\pm[\iota_{17},\iota_{17}]\eta_{33}\varepsilon_{34}=0.
		\]
		As $\Delta(\eta_{17}\varepsilon_{18})\in\operatorname{Im}\Delta$,
		the injectivity (namely, $\ker J\cap\operatorname{Im}\Delta=0$) implies
		$\Delta(\eta_{17}\varepsilon_{18})=0$.
		Recall that
		$\eta_{17}^2\sigma_{19}
		=\eta_{17}\varepsilon_{18}+\nu_{17}^3$.
		Applying $\Delta$, we infer $\Delta(\nu_{17}^3)
		=\Delta(\eta_{17}^2\sigma_{19})
		=4[\nu_{15}]_{17}\sigma_{18}.$	\end{proof}
	
	We now introduce the following result into the discussion, stating it first and then presenting its proof in what follows.
	\begin{prop}
		The structure of $\pi_{25}(\mathit{SO}(n))$ for $19\leq n\leq20$ is as follows.
		\begin{enumerate}
			\item $\pi_{25}(\mathit{SO}(19)) = \mathrm{span}\{[\iota_7]_{19}\eta_7\bar{\mu}_8,\; [\sigma^2_{11}]_{19}\}\cong (\mathbb Z/2)^{2}.$
			Here, $[\nu_{15}]_{19}\sigma_{18} =[\sigma^2_{11}]_{19}.$ 
			\item   $\pi_{25}(\mathit{SO}(20)) = \mathrm{span}\{[\iota_7]_{20}\eta_7\bar{\mu}_8\}\cong \mathbb Z/2.$
		\end{enumerate}
	\end{prop}
	\begin{proof}
		\begin{enumerate}[label=\textup{(\arabic*)},leftmargin=2em]
			\item In the exact sequence $(25)_{18}$:
			\[\begin{aligned}
				&\pi^{18}_{26} \xrightarrow{ \Delta} \pi_{25}(\mathit{SO}(18)) \xrightarrow{i_\ast } \pi_{25}(\mathit{SO}(19)) \xrightarrow{p_\ast } \pi^{18}_{25},
			\end{aligned}\]
			we see that the right-hand $\Delta$ is a monomorphism by the proof of Proposition~\ref{m1:18}. By \cite[Lemma
			6.2(1)]{GM}, it follows that $\Delta(\varepsilon_{18})=0.$ 
			We know $\Delta(\eta_{18}\sigma_{19})=2[\nu_{15}]_{18}\sigma_{18}.$  This leads to (1).
			
			\item In the exact sequence $(25)_{19}$:
			\[\begin{aligned}
				&\pi^{19}_{26} \xrightarrow{ \Delta} \pi_{25}(\mathit{SO}(19)) \xrightarrow{i_\ast } \pi_{25}(\mathit{SO}(20)) \xrightarrow{p_\ast } \pi^{19}_{25},
			\end{aligned}\]
			we note  $\Delta(\nu^2_{19})=[\nu_{15}]_{19}\nu^2_{18} \neq 0,$   because
			$P(\nu^2_{39})\ne 0$ \cite[(5.39)]{MMO}.
			By  Proposition~\ref{m1:16-25},
			we have $\Delta(\sigma_{19})=[\nu_{15}]_{19}\sigma_{18}=[\sigma^2_{11}]_{19}.$ 
			This leads to (2) and completes the proof.
		\end{enumerate}
	\end{proof}
	
	\subsection{The Case $n=21$ and the Stable Groups}
	
	The remaining cases connect the unstable calculations with the stable value. The next proposition
	records explicit representatives and can be compared with Table~\ref{m1:R25K}.
	
	\begin{prop}\label{m1:27}
		The structures of $\pi_{25}(\mathit{SO}(21+n)),$   $\pi_{25}(\mathit{SO}(23)),$  
		$\pi_{25}(\mathit{SO}(24)),$   $\pi_{25}(\mathit{SO}(25))$ and $\pi_{25}(\mathit{SO}(26))$ are as follows.
		\begin{enumerate}
			\item
			$\pi_{25}(\mathit{SO}(21+n)) = \mathrm{span}\{[\iota_7]_{21+n}\eta_7\bar{\mu}_8\}\cong \mathbb Z/2$ for $n=0, 1$\ or $n\ge 6.$ 
			\item
			$\pi_{25}(\mathit{SO}(23)) = \mathrm{span}\{[\eta^2_{22}]\eta_{24},\; [\iota_7]_{23}\eta_7\bar{\mu}_8\}\cong (\mathbb Z/2)^{2}.$
			\item
			$\pi_{25}(\mathit{SO}(24)) = \mathrm{span}\{[\eta^2_{22}]_{24}\eta_{24},\; [\eta_{23}]\eta_{24},\; [\iota_7]_{24}\eta_7\bar{\mu}_8\}\cong (\mathbb Z/2)^{3}.$
			Here $\Delta(\eta^2_{24})=[\eta^2_{22}]_{24}\eta_{24}.$ 
			\item
			$\pi_{25}(\mathit{SO}(25)) = \mathrm{span}\{[\eta_{23}]_{25}\eta_{24},\; [\iota_7]_{25}\eta_7\bar{\mu}_8\}\cong (\mathbb Z/2)^{2}.$
			Here $\Delta(\eta_{25})=[\eta_{23}]_{25}\eta_{24}.$ 
			\item
			$\pi_{25}(\mathit{SO}(26)) = \mathrm{span}\{\Delta(\iota_{26}),\; [\iota_7]_{26}\eta_7\bar{\mu}_8\}\cong \mathbb Z_{(2)}\oplus \mathbb Z/2.$
		\end{enumerate}
		\qed
	\end{prop}
	
	\section{Odd-Primary Components}\label{sec:odd25}
	\textbf{Exception to the standing convention.} In this section,
	$\pi_i(\mathit{SO}(n):\mathrm{odd})$ denotes the subgroup of all odd-order
	elements of the integral homotopy group, as defined in Chapter~\ref{ch:notation}.
	In the proof, the symbol $\pi_{*}(X)$  
	is short for $\pi_{*}(X)\otimes\z_{\mathcal{P}}$ where $\mathcal{P}$
	is the set of odd primes and $\z_{\mathcal{P}}$ denotes the $\mathcal{P}$-localization of $\z.$  A subscript $(p)$ denotes localization at the indicated odd prime
	$p,$   and  $\operatorname{Tor}\pi_i(X_{(p)})\xyd \pi_i(X_{(p)})$ is the subgroup of elements whose orders are finite.
	The default convention at $2$ resumes in the next chapter.
	For an odd prime $p,$   use the comparison formulas in Section~\ref{sec:odd-method}. The nonzero odd-primary part is as follows.
	\begin{prop}
		For $n\geq3,$   we have
		\[
		\pi_{25}(\mathit{SO}(n):\mathrm{odd})\cong
		\begin{cases}
			\mathbb Z/105, & n=3,5,\\
			(\mathbb Z/105)^2, & n=4,\\
			(\mathbb Z/3)^2\oplus\mathbb Z/5\oplus\mathbb Z/7, & n=6,\\
			\mathbb Z/21, & n=7,\\
			(\mathbb Z/3)^2\oplus\mathbb Z/7, & n=8,\\
			\mathbb Z/3, & n=9,10,16,\\
			0, & \text{otherwise}.
		\end{cases}
		\]
	\end{prop}
	\begin{proof}
		We list the odd-primary components:
		$\pi_{25}(S^3) \cong \mathbb Z/105$\uprightsemicolon{}
		$\pi_{25}(\mathit{Sp}(2)) \cong \mathbb Z/105$ ~\cite[Theorem 1.1]{F}\uprightsemicolon{}
		$\pi_{25}(S^5) \cong \mathbb Z/3$\uprightsemicolon{}
		$\pi_{25}(S^7) \cong \mathbb Z/3$\uprightsemicolon{}
		$\pi_{25}(S^{15}) \cong \mathbb Z/3.$ 
		
		We consider the exact sequence induced by the fibration:
		$\mathit{Sp}(3) \to \mathit{Sp}(3)/\mathit{Sp}(2)=S^{11}$:
		\[\begin{aligned}
			&\pi_{26}(S^{11}) \xrightarrow{\Delta} \pi_{25}(\mathit{Sp}(2)) \xrightarrow{i_{\ast}} \pi_{25}(\mathit{Sp}(3)) \xrightarrow{p_{\ast}}\pi_{25}(S^{11})=0.
		\end{aligned}\]
		In light of \cite[(3.14)]{Os}, we obtain $\operatorname{ker} \Delta \cong \mathbb Z/3.$  This
		implies $\pi_{25}(\mathit{Sp}(3)) \cong \mathbb Z/21.$  We know $\pi_{25}(\mathit{Sp}(4)) \cong \mathbb Z/3$~\cite[Theorem 5]{Os}.
		
		Now we have the following:
		$\pi_{25}(\mathit{SO}(3))
		\cong \pi_{25}(S^3) \cong \mathbb Z/3\oplus \mathbb Z/5 \oplus \mathbb Z/7,$  \\
		$\pi_{25}(\mathit{SO}(4)) \cong \pi_{25}(\mathit{SO}(3)) \oplus \pi_{25}(S^3)
		\cong (\mathbb Z/3)^2\oplus (\mathbb Z/5)^2 \oplus (\mathbb Z/7)^2,$  \\
		$\pi_{25}(\mathit{SO}(5)) \cong \pi_{25}(\mathit{Sp}(2)) \cong \mathbb Z/3 \oplus \mathbb Z/5 \oplus \mathbb Z/7,$  \\
		$\pi_{25}(\mathit{SO}(6)) \cong \pi_{25}(\mathit{SO}(5)) \oplus \pi_{25}(S^5) \cong (\mathbb Z/3)^2
		\oplus \mathbb Z/5 \oplus \mathbb Z/7,$  \\
		$\pi_{25}(\mathit{SO}(7)) \cong \pi_{25}(\mathit{Sp}(3)) \cong \mathbb Z/3 \oplus \mathbb Z/7,$  \\
		$\pi_{25}(\mathit{SO}(8)) \cong \pi_{25}(\mathit{SO}(7))\oplus \pi_{25}(S^7) \cong (\mathbb{Z}/3)^2 \oplus \mathbb Z/7,$  \\
		$\pi_{25}(\mathit{SO}(9)) \cong \pi_{25}(\mathit{Sp}(4))
		\cong\mathbb Z/3,$  \\
		$\pi_{25}(\mathit{SO}(10)) \cong \pi_{25}(\mathit{SO}(9))\cong \mathbb Z/3,$  \\
		$\pi_{25}(\mathit{SO}(11))\cong\pi_{25}(\mathit{Sp}(5))=0$ \cite{MimuraMetastable},\\
		$\pi_{25}(\mathit{SO}(12))\cong\pi_{25}(\mathit{SO}(11))\oplus\pi_{25}(S^{11})=0$ \cite{T},\\
		$\pi_{25}(\mathit{SO}(13))\cong\pi_{25}(\mathit{Sp}(6))=0$ \cite{Bo},\\
		$\pi_{25}(\mathit{SO}(14))\cong\pi_{25}(\mathit{SO}(13))\oplus\pi_{25}(S^{13})=0$ \cite{Bo}, \cite{T},\\
		$\pi_{25}(\mathit{SO}(15)) \cong \pi_{25}(\mathit{Sp}(7))=0$ \cite{Bo},\\
		$\pi_{25}(\mathit{SO}(16)) \cong \pi_{25}(\mathit{SO}(15))\oplus\pi_{25}(S^{15})\cong\mathbb Z/3$
		\cite{Bo}, \cite{T}.
		\noindent The same argument gives the groups
		$\pi_{25}(\mathit{SO}(n):\mathrm{odd})=0\ \mbox{for}\ n \ge 17.$ 
		This completes the proof.
	\end{proof}
	\chapter{The Homotopy Groups $\pi_{26}(\mathit{SO}(n))$}\label{ch:degree26}
	\label{m2:chapter}
	
	The computation of $\pi_{26}(\mathit{SO}(n))$ presents an interesting contrast between the unstable and stable situations. On the one hand, nontrivial homotopy classes may appear for finite values of $n$, and their precise orders must be determined as part of the calculation. On the other hand, the stable homotopy group in this degree is trivial, so every such class is necessarily annihilated at some later stage by the homomorphisms induced by the inclusions $\mathit{SO}(n)\hookrightarrow\mathit{SO}(n+1)$. 
	
	The principal task is thus to understand how these unstable classes behave as $n$ increases, rather than merely to identify the groups in which they first occur. The resulting calculation records the orders of the relevant elements and explains their eventual disappearance, thereby connecting the unstable computations with the vanishing of the stable group. 
	
	The investigation is initiated by an auxiliary computation for the exceptional Lie group $G_2$, after which the groups corresponding to $3\leq n\leq13$ are examined in detail. The remaining cases are then treated successively for larger values of $n$, completing the passage from the unstable range to the stable one.
	
	We make free use of the descriptions of the groups $\pi_{*}(S^{j})$
	in Appendix~\ref{tlqdzj} and of the composition relations among
	generators in Toda's range listed in \^{O}guchi's table in
	Appendix~\ref{app:oguchi-relations}.
	
	\section{Main Theorem and Generators}
	
	The theorem records both the cyclic factors and the chosen representatives. The case-by-case proofs
	explain the relations between them\uprightsemicolon{} in particular, a larger order for a lift must be justified by a
	nonzero multiple in the image from the preceding value of $n.$ 
	The lift $[\nu_9\kappa_{12}]$ used in the theorem is fixed below,
	after Lemma~\ref{m2:Dnuka}, with
	$2[\nu_9\kappa_{12}]=\Delta(\nu_{10}\kappa_{13})
	=[\nu_6^2\kappa_{12}]_{10}$.
	The elements $[\nu_9\kappa_{12}]_n$ are the images of this fixed lift.
	
	\begingroup
	\setlength{\abovedisplayskip}{6pt plus 1pt minus 1pt}
	\setlength{\belowdisplayskip}{6pt plus 1pt minus 1pt}
	\setlength{\parindent}{1.5em}
	\begin{thm}\label{m2:thm1}
		\textup{All spaces in this theorem are localized at $2.$  The displayed groups are obtained from the integral homotopy groups by tensoring with $\mathbb Z_{(2)}.$  Free cyclic summands become $\mathbb Z_{(2)}$; for finite groups, this localization retains precisely the $2$-primary summands. The homotopy groups are given as follows.   }\par
		\begin{itemize}[leftmargin=0pt,itemindent=0pt]
			\item[\normalfont(1)] \textup{The Cases $3\leq m\leq13.$ }\par
			\smallskip\noindent\(
			\FitMath{\pi_{26}(\mathit{SO}(3)) = \mathrm{span}\{[\eta_2]\bar{\alpha}\}\cong \mathbb Z/4, \bar\alpha \in \{ \varepsilon_3, 2\iota_{11}, \kappa_{11}\}_1,\; 2\bar \alpha=\eta_3 \varepsilon_4 \kappa_{12}.}
			\)
			
			\par\smallskip\noindent\(
			\GroupLayout{\pi_{26}(\mathit{SO}(4)) = \mathrm{span}\{[\eta_2]_4\bar{\alpha},\; [\iota_3]\bar{\alpha}\}\cong (\mathbb Z/4)^{2}.}{\begin{aligned}\pi_{26}(\mathit{SO}(4)) = \mathrm{span}\{&[\eta_2]_4\bar{\alpha},\\
					&[\iota_3]\bar{\alpha}\}\cong (\mathbb Z/4)^{2}.\end{aligned}}
			\)\par
			\smallskip\noindent\(
			\GroupLayout{\pi_{26}(\mathit{SO}(5)) = \mathrm{span}\{[\nu_4\bar{\zeta}_7],\; [\nu_4\bar{\sigma}_7]\}\cong \mathbb Z/8\oplus \mathbb Z/2.}{\begin{aligned}\pi_{26}(\mathit{SO}(5)) = \mathrm{span}\{&[\nu_4\bar{\zeta}_7],\\
					&[\nu_4\bar{\sigma}_7]\}\cong \mathbb Z/8\oplus \mathbb Z/2.\end{aligned}}
			\)\par
			\par\smallskip\noindent\(
			\GroupLayout{\pi_{26}(\mathit{SO}(6)) = \mathrm{span}\{[\nu_4\bar{\zeta}_7]_6,\; [\nu_4\bar{\sigma}_7]_6,\; [ \nu_5]\eta_8\bar{\mu}_9\}\cong \mathbb Z/8\oplus (\mathbb Z/2)^{2}.}{\begin{aligned}\pi_{26}(\mathit{SO}(6)) = \mathrm{span}\{&[\nu_4\bar{\zeta}_7]_6,\; [\nu_4\bar{\sigma}_7]_6,\\
					&[ \nu_5]\eta_8\bar{\mu}_9\}\cong \mathbb Z/8\oplus (\mathbb Z/2)^{2}.\end{aligned}}
			\)\par
			\smallskip\par\smallskip\noindent\(
			\GroupLayout{\pi_{26}(\mathit{SO}(7)) = \mathrm{span}\{[2P(\rho_{13})],\; [\nu_4\bar{\zeta}_7]_7,\; [\nu_4\bar{\sigma}_7]_7,\; [\nu^2_6\kappa_{12}]\}\cong \mathbb Z/32\oplus \mathbb Z/8\oplus (\mathbb Z/2)^{2}.}{\begin{aligned}\pi_{26}(\mathit{SO}(7)) = \mathrm{span}\{&[2P(\rho_{13})],\; [\nu_4\bar{\zeta}_7]_7,\; [\nu_4\bar{\sigma}_7]_7,\\
					&[\nu^2_6\kappa_{12}]\}\cong \mathbb Z/32\oplus \mathbb Z/8\oplus (\mathbb Z/2)^{2}.\end{aligned}}
			\)\par
			\smallskip
			
			\par\smallskip\noindent\(
			\GroupLayout{\pi_{26}(\mathit{SO}(8)) = \mathrm{span}\{[2P(\rho_{13})]_8,\; [\nu_4\bar{\zeta}_7]_8,\; [\iota_7]\bar{\zeta}_7,\; [\nu_4\bar{\sigma}_7]_8,\; [\nu^2_6\kappa_{12}]_8,\; [\iota_7]\bar{\sigma}_7\}\cong \mathbb Z/32\oplus (\mathbb Z/8)^{2}\oplus (\mathbb Z/2)^{3}.}{\begin{aligned}\pi_{26}(\mathit{SO}(8)) = \mathrm{span}\{&[2P(\rho_{13})]_8,\; [\nu_4\bar{\zeta}_7]_8,\; [\iota_7]\bar{\zeta}_7,\; [\nu_4\bar{\sigma}_7]_8,\\
					&[\nu^2_6\kappa_{12}]_8,\; [\iota_7]\bar{\sigma}_7\}\cong \mathbb Z/32\oplus (\mathbb Z/8)^{2}\oplus (\mathbb Z/2)^{3}.\end{aligned}}
			\)\par
			\smallskip
			
			\par\smallskip\noindent\(
			\GroupLayout{\pi_{26}(\mathit{SO}(9)) = \mathrm{span}\{[2P(\rho_{13})]_9,\; [\iota_7]_9\bar{\zeta}_7,\; [\nu^2_6\kappa_{12}]_9,\; [\iota_7]_9\bar{\sigma}_7\}\cong \mathbb Z/32\oplus \mathbb Z/8\oplus (\mathbb Z/2)^{2}.}{\begin{aligned}\pi_{26}(\mathit{SO}(9)) = \mathrm{span}\{&[2P(\rho_{13})]_9,\; [\iota_7]_9\bar{\zeta}_7,\; [\nu^2_6\kappa_{12}]_9,\\
					&[\iota_7]_9\bar{\sigma}_7\}\cong \mathbb Z/32\oplus \mathbb Z/8\oplus (\mathbb Z/2)^{2}.\end{aligned}}
			\)\par
			\smallskip
			
			\par\smallskip\noindent\(
			\GroupLayout{\pi_{26}(\mathit{SO}(10)) = \mathrm{span}\{[\sigma_9\eta_{16}\mu_{17}],\; [\eta_9\mu_{10}\sigma_{19}],\; [\nu_9\kappa_{12}],\; [\iota_7]_{10}\bar\sigma_7\}\cong \mathbb Z/64\oplus \mathbb Z/8\oplus \mathbb Z/4\oplus \mathbb Z/2,}{\begin{aligned}\pi_{26}(\mathit{SO}(10)) = \mathrm{span}\{&[\sigma_9\eta_{16}\mu_{17}],\; [\eta_9\mu_{10}\sigma_{19}],\; [\nu_9\kappa_{12}],\\
					&[\iota_7]_{10}\bar\sigma_7\}\cong \mathbb Z/64\oplus \mathbb Z/8\oplus \mathbb Z/4\oplus \mathbb Z/2.\end{aligned}}
			\)
			\par\smallskip\noindent\(
			\GroupLayout{\pi_{26}(\mathit{SO}(11)) = \mathrm{span}\{[\sigma_9\eta_{16}\mu_{17}]_{11},\; [\mu_{10}\sigma_{19}],\; [\nu_9\kappa_{12}]_{11},\; [\iota_7]_{11}\bar\sigma_7\}\cong \mathbb Z/64\oplus \mathbb Z/8\oplus (\mathbb Z/2)^{2},}{\begin{aligned}\pi_{26}(\mathit{SO}(11)) = \mathrm{span}\{&[\sigma_9\eta_{16}\mu_{17}]_{11},\; [\mu_{10}\sigma_{19}],\; [\nu_9\kappa_{12}]_{11},\\
					&[\iota_7]_{11}\bar\sigma_7\}\cong \mathbb Z/64\oplus \mathbb Z/8\oplus (\mathbb Z/2)^{2},\end{aligned}}
			\)\par
			\smallskip
			{\raggedleft
				where $2[\mu_{10}\sigma_{19}]\equiv x[\eta_9\mu_{10}\sigma_{19}]_{11}\m16[\sigma_9\eta_{16}\mu_{17}]_{11}$ for some odd $x.$ 
				\par}
			
			\par\smallskip\noindent\(
			\GroupLayout{\pi_{26}(\mathit{SO}(12)) \;= {} \mathrm{span}\{[\Sigma^2\rho'],\; \ell_{0}[\sigma_9\eta_{16}\mu_{17}]_{12}-2[\Sigma^2\rho'],\; 2\ell_{0}[\sigma_9\eta_{16}\mu_{17}]_{12} +[\mu_{10}\sigma_{19}]_{12}-4[\Sigma^2\rho'],\; [\nu_9\kappa_{12}]_{12},\; [\iota_7]_{12}\bar\sigma_7,\; [\eta_{11}]\kappa_{12}\}\cong \mathbb Z/128\oplus \mathbb Z/16\oplus \mathbb Z/4\oplus (\mathbb Z/2)^{3}.}{\begin{aligned}\pi_{26}(\mathit{SO}(12)) \;= {} \mathrm{span}\{&[\Sigma^2\rho'],\; \ell_{0}[\sigma_9\eta_{16}\mu_{17}]_{12}-2[\Sigma^2\rho'],\; 2\ell_{0}[\sigma_9\eta_{16}\mu_{17}]_{12} +[\mu_{10}\sigma_{19}]_{12}-4[\Sigma^2\rho'],\\
					&[\nu_9\kappa_{12}]_{12},\; [\iota_7]_{12}\bar\sigma_7,\; [\eta_{11}]\kappa_{12}\}\cong \mathbb Z/128\oplus \mathbb Z/16\oplus \mathbb Z/4\oplus (\mathbb Z/2)^{3},\end{aligned}}
			\)\par
			\noindent
			{\raggedleft
				
				where $\ell_{0}$ is some odd integer and $16[\Sigma^2\rho']=8\ell_{0}[\sigma_9\eta_{16}\mu_{17}]_{12}+4[\mu_{10}\sigma_{19}]_{12}.$ 
				\par}   
			\par\smallskip
			
			\par\smallskip\noindent\(
			\GroupLayout{\pi_{26}(\mathit{SO}(13)) = \mathrm{span}\{[\Sigma^2\rho']_{13},\; [8\sigma^2_{12}],\; [\eta_{11}]_{13}\kappa_{12},\; [\iota_7]_{13}\bar\sigma_7\}\cong \mathbb Z/128\oplus \mathbb Z/8\oplus (\mathbb Z/2)^{2},}{\begin{aligned}\pi_{26}(\mathit{SO}(13)) = \mathrm{span}\{&[\Sigma^2\rho']_{13},\; [8\sigma^2_{12}],\; [\eta_{11}]_{13}\kappa_{12},\\
					&[\iota_7]_{13}\bar\sigma_7\}\cong \mathbb Z/128\oplus \mathbb Z/8\oplus (\mathbb Z/2)^{2},\end{aligned}}
			\)\par
			\smallskip
			{\raggedleft
				where
				$2[8\sigma_{12}^2]\equiv2\ell_{0}[\sigma_9\eta_{16}\mu_{17}]_{13}+[\mu_{10}\sigma_{19}]_{13}-4[\Sigma^2\rho']_{13}\m[\eta_{11}]_{13}\kappa_{12},[\iota_7]_{13}\bar\sigma_7.$
				\par}
			
			\item[\normalfont(2)] \textup{The Cases $14\leq m\leq17.$ }\par
			\smallskip\noindent\(
			\GroupLayout{\pi_{26}(\mathit{SO}(14)) = \mathrm{span}\{[(\Sigma\theta)\eta_{25}],\; [8\sigma^2_{12}]_{14}\}\cong \mathbb Z/256\oplus \mathbb Z/8.}{\begin{aligned}\pi_{26}(\mathit{SO}(14)) = \mathrm{span}\{&[(\Sigma\theta)\eta_{25}],\\
					&[8\sigma^2_{12}]_{14}\}\cong \mathbb Z/256\oplus \mathbb Z/8.\end{aligned}}
			\)\par
			\smallskip\noindent\(
			\FitMath{\pi_{26}(\mathit{SO}(15)) = \mathrm{span}\{[8\sigma^2_{12}]_{15}\}\cong \mathbb Z/8.}
			\)\par
			\smallskip\noindent\(
			\GroupLayout{\pi_{26}(\mathit{SO}(16)) = \mathrm{span}\{[\zeta_{15}],\; [8\sigma^2_{12}]_{16}\}\cong (\mathbb Z/8)^{2}.}{\begin{aligned}\pi_{26}(\mathit{SO}(16)) = \mathrm{span}\{&[\zeta_{15}],\\
					&[8\sigma^2_{12}]_{16}\}\cong (\mathbb Z/8)^{2}.\end{aligned}}
			\)\par
			\smallskip\noindent\(
			\FitMath{\pi_{26}(\mathit{SO}(17)) = \mathrm{span}\{[\zeta_{15}]_{17}\}\cong \mathbb Z/8.}
			\)\par
			\item[\normalfont(3)] \textup{The Cases $m\geq18.$ }\par
			\smallskip\noindent
			\MainWhere{$\pi_{26}(\mathit{SO}(18))
				=\mathrm{span}\{[\eta_{17}\varepsilon_{18}]\}\cong \mathbb Z/8,$  }{where $2[\eta_{17}\varepsilon_{18}]\ty[\zeta_{15}]_{18}\m 2[\zeta_{15}]_{18}.$ }
			\par\smallskip\noindent
			\MainWhere{$\pi_{26}(\mathit{SO}(19))
				=\mathrm{span}\{[\varepsilon_{18}]\}\cong \mathbb Z/8,$  }{where $2[\varepsilon_{18}]\ty[\eta_{17}\varepsilon_{18}]_{19}\m2[\eta_{17}\varepsilon_{18}]_{19}.$ }
			\par\smallskip\noindent
			$\pi_{26}(\mathit{SO}(20))
			=\mathrm{span}\{\Delta\sigma_{20},\; [\varepsilon_{18}]_{20}-2\Delta\sigma_{20}\}\cong \mathbb Z/16\oplus \mathbb Z/4.$ 
			\par\smallskip\noindent
			\MainWhere{$\pi_{26}(\mathit{SO}(21))
				=\mathrm{span}\{[\nu^2_{20}]\}\cong \mathbb Z/8,$  }{where $J([\nu^2_{20}])\ty M'_2\m 2 M'_2.$ }
			\par\smallskip\noindent
			$\pi_{26}(\mathit{SO}(22))
			=\mathrm{span}\{[\nu^2_{20}]_{22}\}\cong \mathbb Z/8.$ 
			\par\smallskip\noindent
			$\pi_{26}(\mathit{SO}(23))
			=\mathrm{span}\{[\nu^2_{20}]_{23}\}\cong \mathbb Z/8.$ 
			\par\smallskip\noindent
			\MainWhere{$\pi_{26}(\mathit{SO}(24))
				=\mathrm{span}\{[\nu^2_{20}]_{24},\; [\nu_{23}]\}\cong (\mathbb Z/8)^{2},$  }{where $\Delta\nu_{24}\equiv
				2[\nu_{23}]-[\nu^2_{20}]_{24}\bmod 2[\nu^2_{20}]_{24}.$ }
			\par\smallskip\noindent
			\MainWhere{$\pi_{26}(\mathit{SO}(25))
				=\mathrm{span}\{[\nu_{23}]_{25}\}\cong \mathbb Z/8,$  }{where $\Delta(\eta^2_{25})=4[\nu_{23}]_{25}.$ }
			\par\smallskip\noindent
			\MainWhere{$\pi_{26}(\mathit{SO}(26))
				=\mathrm{span}\{[\nu_{23}]_{26}\}\cong \mathbb Z/4,$  }{where $\Delta\eta_{26}=2[\nu_{23}]_{26}.$ }
			\par\smallskip\noindent
			\MainWhere{$\pi_{26}(\mathit{SO}(27))
				=\mathrm{span}\{\Delta(\iota_{27})\}\cong \mathbb Z/2,$  }{where $\Delta(\iota_{27})=[\nu_{23}]_{27}.$ }
			\par\smallskip\noindent
			$\pi_{26}(\mathit{SO}(m))=0$ for $m\geq28.$ 
			
		\end{itemize}
	\end{thm}
	\endgroup
	\section{Auxiliary Fibrations and the $G_2$ Calculations}\label{m2:sec:small-n}
	
	\subsection{Unitary and Symplectic Input}
	We need the group structures of $\pi_k(G_2)$ for $k=25,26.$
	We list the following groups from \cite[Theorem 4.1]{MT2} and \cite[Theorem 1.1]{F}:
	\[\begin{gathered}
		\begin{aligned}
			&\pi_{11}(\mathit{SU}(3))={}\frac{\mathbb Z}{4}\{[\nu^2_5]\},
		\end{aligned}\\[3pt]
		\begin{aligned}
			&\pi_{26}(\mathit{SU}(3))={}\frac{\mathbb Z}{2}\{[\nu_5\eta_8\bar\mu_9]\}\oplus \frac{\mathbb Z}{2}\{i_{\ast}\bar\alpha \},
		\end{aligned}\\[3pt]
		\begin{aligned}
			&\pi_{27}(\mathit{SU}(3))={}\frac{\mathbb Z}{4}\{[\rho^{\mathrm{IV}}]\sigma_{20}\}\oplus \frac{\mathbb Z}{2}\{[\nu_5\bar\sigma_8]\}
		\end{aligned}
	\end{gathered}\]
	and
	\[\begin{gathered}
		\begin{aligned}
			&\pi_{26}(\mathit{Sp}(2))={}\frac{\mathbb Z}{8}\{[\bar{\zeta}_7]\}\oplus\frac{\mathbb Z}{2}\{[\bar{\sigma}_7]\},
		\end{aligned}\\[3pt]
		\begin{aligned}
			&\pi_{27}(\mathit{Sp}(2))={}\frac{\mathbb Z}{2}\{i_*\delta_3\}\oplus \frac{\mathbb Z}{2}\{i_*\bar{\mu}_3\sigma_{20}\} \oplus\frac{\mathbb Z}{2}\{[4\bar{\kappa}_7]\}.
		\end{aligned}
	\end{gathered}\]
	
	We use the maps $r_3,$   $i_G,$   $p_G,$   and $\mathcal{T}$ in diagram~\eqref{m2:Y}, with the choices of lifts fixed in Section~\ref{lqzh}.
	
	\subsection{Boundary Maps in the $G_2$ Fibration}
	To use the fibration $\mathit{SU}(3)\to G_2\to S^6,$ we first isolate a vanishing composition in
	$\mathit{SU}(3).$ This relation will remove one of the indeterminacy terms in the Toda-bracket
	calculation below.
	
	\begin{lem}\label{m2:[s5kp}
		$[\sigma''']\kappa_{12} =0.$
	\end{lem}
	\begin{proof}
		By \cite[p.162]{Mi2}, we conclude that
		$\Delta(\sigma'') \equiv [\sigma'''] \bmod 2[\sigma'''] .$
		This implies
		$[\sigma''']\kappa_{12} = \Delta(\sigma'' \kappa_{13}) .$
		By Lemma~\ref{m2:lem:sigma-kappa-vanishing}, the composition of maps between spheres
		$\sigma''\kappa_{13}$ vanishes, and hence so does its boundary.
		This proves the assertion.
	\end{proof}
	
	We next identify the connecting image of $\nu_6^2.$ Expressing this image as a multiple of the chosen
	generator of $\pi_{11}(\mathit{SU}(3))$ makes it possible to compare the subsequent Toda brackets in
	the same group.
	
	\begin{lem}\label{m2:[2]n25}
		$\Delta'(\nu^2_6)=[2\iota_5]\nu^2_5=2[\nu^2_5] \in \pi_{11}(\mathit{SU}(3)).$
	\end{lem}
	\begin{proof}
		By \cite[Proposition 6.2]{Mi2}, we see that $\Delta'(\nu^2_6)=[2\iota_5]\nu^2_5.$
		We obtain $2[\nu^2_5]={i_U}_*\varepsilon_3$ \cite[(4.1)]{MT1}.
		Since $[2\iota_5]\in\{i,\eta_3,2\iota_4\},$   the bracket
		composition rule gives
		\[
		\begin{aligned}
			&[2\iota_5]\eta_5 \in\{i_U,\eta_3,2\iota_4\}\circ\eta_5 =i\circ\{\eta_3,2\iota_4,\eta_4\} \equiv i_*\nu'\pmod{i_*\pi_5^3\circ\eta_5},
		\end{aligned}
		\]
		where $i_*\pi_5^3\circ\eta_5=\mathrm{span}\{2i_*\nu'\}=0.$ 
		Consequently,
		\[
		\begin{aligned}
			&[2\iota_5]\nu_5^2 \in[2\iota_5]\circ\{\eta_5,\nu_6,\eta_9\} \subseteq\{{i_U}_*\nu',\nu_6,\eta_9\} \supseteq i_U\circ\{\nu',\nu_6,\eta_9\} \ni {i_U}_*\varepsilon_3.
		\end{aligned}
		\]
		The comparison is modulo ${i_U}_*\pi_{10}^3\circ\eta_{10}=0.$ 
		This proves the stated identity.
	\end{proof}
	
	The preceding relations allow us to compute the boundary of $\eta_6\bar\kappa_7.$ The point is to
	identify the resulting class with the image of $\bar\alpha$ and to check that the bracket
	indeterminacy does not leave an additional ambiguity.
	
	\begin{lem}\label{m2:lem:eta6-kappa-boundary}
		$ \Delta'(\eta_6\bar\kappa_7) =i_{U\ast} \bar{\alpha}.$
	\end{lem}
	\begin{proof}
		By \cite[Lemma 6.1]{Mi1} and Lemma~\ref{m2:[2]n25}, we have
		\begin{align*}
			\notag \Delta'(\eta_6\bar\kappa_7)=\Delta'(\Sigma\alpha)=[2\iota_5]\alpha &\in  [2\iota_5] \circ \{\nu^2_5,2\iota_{11},\kappa_{11}\}_1\\
			\notag &\subseteq   \{[2\iota_5]\nu^2_5,2\iota_{11},\kappa_{11}\}_1\\
			\notag &=  \{2[\nu^2_5],2\iota_{11},\kappa_{11}\}_1.
		\end{align*}
		On the other hand, we see that $i_{U\ast}(\bar{\alpha}) \in
		\{i_{U\ast}(\varepsilon_3),2\iota_{11},\kappa_{11}\}_1 =\{2[\nu^2_5],2\iota_{11},\kappa_{11}\}_1$ by
		\cite[(3.1)]{MMO} and \cite[(4.1)]{MT2}. This implies that $\Delta'(\eta_6\bar\kappa_7) \equiv
		i_{U\ast}\bar{\alpha} \,\,\bmod\, \operatorname{Ind}\{2[\nu^2_5],2\iota_{11},\kappa_{11}\}_1.$
		Consider $$\operatorname{Ind}\{2[\nu^2_5],2\iota_{11},\kappa_{11}\}_1 = U^3_{12} \circ \kappa_{12} +
		i_{U\ast}(\varepsilon_3) \circ \Sigma\pi^{10}_{25},$$\noindent where
		$
		U^3_{12} = \frac{\mathbb Z}{4}\{[\sigma''']\}
		$ \cite[Theorem 4.1]{MT2} and
		$
		\Sigma\pi^{10}_{25} = \frac{\mathbb Z}{16}\{\Sigma^2\rho'\} \jia \frac{\mathbb Z}{2}\{\bar{\varepsilon}_{11}\}
		.$ 
		Since
		$
		i_{U\ast}(\varepsilon_3)=2[\nu^2_5],
		$
		we find
		$
		i_{U\ast}( \varepsilon_3\bar{\varepsilon}_{11})= 2[\nu^2_5]\bar{\varepsilon}_{11}=0.
		$
		By \eqref{m2:Osm}, we conclude that
		$$2[\nu^2_5]\Sigma^2\rho'=\Delta'(\nu^2_6(\Sigma^3\rho'))=\Delta'(\nu_6\circ \Sigma^3(\nu_6\rho'))=0.$$\noindent This implies $\operatorname{Ind}\{2[\nu^2_5],2\iota_{11},\kappa_{11}\}_1
		=0$ and completes the proof.
	\end{proof}
	
	\subsection{The Degree 26 Group}
	We now collect the information about $G_2$ needed in degrees $26$ and $27.$ Besides the group
	decompositions, we record the multiples of the degree-$26$ lift that enter the extension calculations
	for the special orthogonal groups.
	
	\begin{prop}\label{m2:G2}
		\(\pi_{26}(G_2)=\mathrm{span}\{\left<2P(\rho_{13}) \right>,\; \left<\nu^2_6 \kappa_{12} \right>\}\cong \mathbb Z/32\oplus \mathbb Z/2  .\)
		The generators satisfy, for an odd integer $x,$  
		\[
		\begin{aligned}
			&2\left<2P(\rho_{13})\right>=x\left<[2\e_{6},\e_{6}]\right>\Sigma^2\rho',\quad 16\left<2P(\rho_{13})\right>=8\left<[2\e_{6},\e_{6}]\right>\Sigma^2\rho' =i_{G*}[\nu_5\eta_8\bar{\mu}_9].
		\end{aligned}
		\]
	\end{prop}
	\begin{proof}For $\pi_{26}(G_2),$ we use the  exact sequence
		\[\begin{aligned}
			&\pi^6_{27} \xrightarrow{ \Delta'} U^3_{26} \xrightarrow{i_{G\ast} } \pi_{26}(G_2) \xrightarrow{ p_{G\ast}} \pi^6_{26} \xrightarrow{\Delta' } U^3_{25},
		\end{aligned}\]
		where
		\[
		\begin{aligned}
			\pi^6_{27}
			&=\frac{\mathbb Z}{2}\{\eta_6\bar\kappa_7\},
			\qquad&
			U^3_{26}
			&=\left(\frac{\mathbb Z}{2}\right)^2
			\{[\nu_5\eta_8\bar\mu_9],i_{U\ast}\bar\alpha\},
			\\[6pt]
			\pi^6_{26}
			&=\frac{\mathbb Z}{32}\{P(\rho_{13})\}
			\oplus\frac{\mathbb Z}{4}\{\bar\kappa'\},
			\qquad&
			U^3_{25}
			&=\left(\frac{\mathbb Z}{2}\right)^2
			\{[\nu^2_5]\kappa_{11},[\nu_5\eta_8\mu_9]\sigma_{18}\}.
		\end{aligned}
		\]
		In light of Lemma~\ref{m2:lem:eta6-kappa-boundary}, we have a short exact sequence
		\[\begin{aligned}
			&0 \to \frac{\mathbb Z}{2}\{[\nu_5\eta_8\bar{\mu}_9]\} \xrightarrow{i_{G\ast} } \pi_{26}(G_2) \xrightarrow{ p_{G\ast}} \frac{\mathbb Z}{16}\{2P(\rho_{13})\} \oplus \frac{\mathbb Z}{2}\{2\bar{\kappa}'\} \to 0.
		\end{aligned}\]
		By \cite[Lemma 10.9]{T} and \cite[Lemma 15.3]{MT1}, $2P(\rho_{13}) = P(\Sigma^4 \rho')=P(\iota_{13})
		\Sigma^2\rho' \,\uprightsemicolon{}\,\, 2\bar{\kappa}' =\nu^2_6 \kappa_{12}.$ Thus we can rewrite
		\[\begin{aligned}
			&0 \to \frac{\mathbb Z}{2}\{[\nu_5\eta_8\bar{\mu}_9]\} \xrightarrow{i_{G\ast} } \pi_{26}(G_2) \xrightarrow{ p_{G\ast}} \frac{\mathbb Z}{16}\{P(\iota_{13}) \Sigma^2\rho'\} \oplus \frac{\mathbb Z}{2}\{\nu^2_6 \kappa_{12}\} \to 0.
		\end{aligned}\]
		The bracket description of $\rho'$ in
		$\{\sigma_9,16\iota_{16},\sigma_{16}\}_1$
		\cite[p.~103]{T} also makes the multiple explicit.  Using
		$P(\nu_{21}^2)=\sigma_{10}\bar\nu_{17}
		=\sigma_{10}\varepsilon_{17}$ \cite[(10.20)]{T}, we obtain
		\[
		\begin{aligned}
			&8\Sigma^2\rho' \in\{\sigma_{11},16\iota_{18},8\sigma_{18}\}_1 \supseteq\{8\sigma_{11},2\iota_{18},8\sigma_{18}\}_1\m\sigma_{11}\circ \Sigma\pi_{25}^{17}+\pi_{19}^{11}\circ8\sigma_{19}=0.
		\end{aligned}
		\]
		
		\noindent

		Recall that $i_{G\ast}[\nu_5\eta_8\mu_9]=8 \langle [2\e_{6},\e_{6}] \rangle \sigma_{11}$ from ~\cite[Theorem
		6.1]{Mi2}.  Since
		$8\rho'=\Sigma^4(\rho^{\mathrm{IV}}) = \{8\sigma_9, 2\iota_{16}, 8\sigma_{16}\}_1$
		\cite[p.103]{T} and $\pi_{19}(G_{2})\cong\z/2$ \cite[Theorem 6.1]{Mi2}, we obtain
		\begin{align*}
			\left<[2\e_{6},\e_{6}] \right> \circ \Sigma^2(8\rho')
			\notag   &\in  \{\left<[2\e_{6},\e_{6}] \right>\circ 8\sigma_{11}, 2\iota_{18}, 8\sigma_{18}\} \\
			\notag   &=  \{i_{G\ast}[\nu_5\eta_8\mu_9], 2\iota_{18}, 8\sigma_{18}\} \\
			\notag   &\supseteq   i_{G\ast}\{[\nu_5\eta_8\mu_9], 2\iota_{18}, 8\sigma_{18}\} \\
			\notag   &\ni  i_{G\ast}[\nu_5\eta_8\bar\mu_9]\m i_{G\ast}[\nu_5\eta_8\mu_9]\circ \pi^{18}_{26}+\pi_{19}(G_{2})\circ 8\sigma_{19} =0.
		\end{align*}
		\noindent  It follows that $\left<[2\e_{6},\e_{6}] \right> \circ \Sigma^2(8\rho')= i_{G\ast}[\nu_5\eta_8\bar{\mu}_9].$  Thus,  $2\left<2P(\rho_{13})\right>=x\left<[2\e_{6},\e_{6}]\right>\Sigma^2\rho'$ for  some  odd
		$x$ .  Observe that    
		\[\begin{aligned}
			&(p_G)_*\bigl(\langle[2\e_{6},\e_{6}]\rangle \Sigma^2\rho'\bigr) =[2\e_{6},\e_{6}]\Sigma^2\rho'=4P(\rho_{13}).
		\end{aligned}\]
		Thus $2\langle2P(\rho_{13})\rangle$ and
		$\langle[2\e_{6},\e_{6}]\rangle \Sigma^2\rho'$ have the same projection.
		By the short exact sequence above, their difference lies in
		$\mathrm{span}\{i_{G\ast}[\nu_5\eta_8\bar\mu_9]\}.$ 
		In addition,  $8 \langle [2\e_{6},\e_{6}] \rangle \circ \Sigma^2 \rho' = i_{G \ast}[\nu_5 \eta_8\bar \mu_9].$
		Thus, we find
		$2 \langle 2P(\rho_{13}) \rangle  \equiv \langle [2\e_{6},\e_{6}] \rangle \circ \Sigma^2 \rho' \bmod 8\langle [2\e_{6},\e_{6}] \rangle \circ \Sigma^2 \rho'.$
		
		We choose $$\langle\nu_6^2\kappa_{12}\rangle:=\langle P(\Sigma\theta)+\nu_6\kappa_9\rangle\nu_{23},$$  using the lift from Chapter~\ref{ch:degree23}.
		Its projection is $\nu_6^2\kappa_{12},$   since $P(\Sigma\theta)\nu_{23}=0$ and $\nu_6\kappa_9\nu_{23}=\nu_6^2\kappa_{12}.$ 
		The relation $$2\langle P(\Sigma\theta)+\nu_6\kappa_9\rangle=i_{G*}[\nu_5\bar\varepsilon_8]$$ was established there.
		The composition relations in Appendix~\ref{app:oguchi} give $$\bar\varepsilon_8\nu_{23}=\eta_8\kappa_9\nu_{23}=\eta_8\nu_9\kappa_{12}=0.$$ 
		Thus $[\nu_5\bar\varepsilon_8]\nu_{23}$ projects to zero in $\pi_{26}^5.$ 
		In the displayed decomposition of $U^3_{26},$   the generator $[\nu_5\eta_8\bar\mu_9]$ has nonzero projection $\nu_5\eta_8\bar\mu_9,$   so $[\nu_5\bar\varepsilon_8]\nu_{23}\in\mathrm{span}\{i_{U*}\bar\alpha\}.$ 
		Lemma~\ref{m2:lem:eta6-kappa-boundary} and exactness give $i_{G*}i_{U*}\bar\alpha=0.$ 
		Consequently, $2\langle\nu_6^2\kappa_{12}\rangle=i_{G*}([\nu_5\bar\varepsilon_8]\nu_{23})=0.$ 
		Since its projection is $\nu_6^2\kappa_{12}=2\bar\kappa'\neq0,$   this lift has order exactly $2.$ 
		Its projection is independent of $2P(\rho_{13}).$  Together with the order-$32$ relation proved above and the short exact sequence, this gives the asserted direct sum and its generators.
		
		Since ${i_G}_*([\nu_5\eta_8\mu_9]\sigma_{18})=0,$ we can write $[\nu_5\eta_8\mu_9]\sigma_{18} =
		\Delta' (xP(\rho_{13}) + y \bar\kappa').$ 
		For the choices fixed in Lemma~\ref{m1:mtb2}(1),
		\[\FitMath{\Delta'(\bar\kappa')=[\nu_5^2]\kappa_{11}.}\]    \noindent       By Lemma~\ref{m1:r3n5}(1),
		we obtain \(\FitMath{r_{3*}([\nu_5^2])=\pm[ \nu_5]\nu_8\in\pi_{11}(\mathit{SO}(6)).}\)  \end{proof}
	
	\subsection{The Boundary of $P(\rho_{13})$}
	The remaining boundary calculation concerns $P(\rho_{13}).$ Identifying its image explicitly will let
	us compare the connecting maps for the two fibrations in diagram~\eqref{m2:Y}.
	
	\begin{lem}\label{m2:D'Prh13}
		$
		\Delta'(P(\rho_{13})) = [\nu_5\eta_8\mu_9]\sigma_{18}.
		$
	\end{lem}
	\begin{proof}
		By the relation ${i_G}_*[\nu_5\eta_8\mu_9]\sigma_{18}=0$ from Proposition~\ref{m2:G2},
		$
		[\nu_5\eta_8\mu_9]\sigma_{18}
		$
		is contained in the image of
		$
		\Delta':\pi^6_{26} \to \pi_{25}(\mathit{SU}(3))
		,$  
		where
		$\pi_{25}(\mathit{SU}(3)) = \frac{\mathbb Z}{2}\{[\nu^2_5]\kappa_{11}\}\oplus\frac{\mathbb Z}{2}\{[\nu_5\eta_8\mu_9]\sigma_{18}\}.$ Lemma~\ref{m1:mtb2}(2) gives
		\[\FitMath{\Delta'(P(\rho_{13}))=[\nu_5\eta_8\mu_9]\sigma_{18}.}\]
		In particular,
		\[\FitMath{\Delta'(P(\rho_{13}))\equiv[\nu_5\eta_8\mu_9]\sigma_{18}\pmod{\mathrm{span}\{[\nu_5^2]\kappa_{11}\}}.}\] Applying $J{r_3}_*$ to the above relation, we conclude that $J(\Delta P(\rho_{13})) \equiv
		J({r_3}_*([\nu_5\eta_8\mu_9]\sigma_{18})) \,\,\,\bmod \,\,\,J({r_3}_*([\nu^2_5]\kappa_{11})).$ Using
		the relation
		$H(J({r_3}_{\ast}([\nu^2_5]\kappa_{11})))=4\bar{\kappa}_{11} \neq 0,$
		we see that
		$
		\Delta'(P(\rho_{13})) = [\nu_5\eta_8\mu_9]\sigma_{18}.
		$
	\end{proof}
	
	Set $x=0.$  For the choices of Lemma~\ref{m1:mtb2}(1), Lemma~\ref{m1:r3n5}(1) gives
	\[\FitMath{\begin{aligned}
			\Delta'(\bar\kappa')&=[\nu_5^2]\kappa_{11},\\
			\Delta(\bar\kappa')&=r_{3*}\Delta'(\bar\kappa')=[ \nu_5]\nu_8\kappa_{11}.
	\end{aligned}}\]\noindent       Applying $r_{3 \ast}$ to this relation and using $r_{3 \ast}([\nu^2_5])=\pm [ \nu_5]\nu_8,$   we obtain
	$\Delta (\bar\kappa' +xP(\rho_{13}))=r_{3 \ast}[\nu^2_5]\kappa_{11} = [ \nu_5]\nu_8\kappa_{11}.$
	Since the right-hand side is in the image of $\Delta,$ we have
	\[\begin{aligned}
		&[ \nu_5]_7\nu_8\kappa_{11}=0 .
	\end{aligned}\]
	
	The splitting for $n=7$ is given in Proposition~\ref{so7fj}. The required degree 25
	spin-group calculations are collected in Chapter~\ref{ch:degree25}.
	\section{The Cases $3\leq n\leq13$}\label{m2:sec:degree26-low}
	\subsection{The Cases $3\leq n\leq6$}
	
	We begin the degree-$26$ calculation for $n=3,4.$  These groups provide the initial generators
	whose images will be followed through the inclusions into special orthogonal groups with larger values of $n.$ 
	
	\begin{prop}\label{m2:3426}
		The group structures of $\pi_{26}(\mathit{SO}(m))$ for $m=3,4$ are given as follows. $\pi_{26}(\mathit{SO}(3)) = \mathrm{span}\{[\eta_2]\bar{\alpha}\}\cong \mathbb Z/4 ,$      $   \pi_{26}(\mathit{SO}(4))  = \mathrm{span}\{[\eta_2]_4\bar{\alpha},\; [\iota_3]\bar{\alpha}\}\cong (\mathbb Z/4)^{2}.$ 
		Moreover,
		\[
		\begin{aligned}
			&2([\eta_2]\bar\alpha)=[\eta_2]\eta_3\varepsilon_4\kappa_{12},\qquad 2([\eta_2]_4\bar\alpha)=[\eta_2]_4\eta_3\varepsilon_4\kappa_{12}, \\
			&2([\iota_3]\bar\alpha)=[\iota_3]\eta_3\varepsilon_4\kappa_{12}.
		\end{aligned}
		\]
		Here $\bar\alpha\in\{\varepsilon_3,2\iota_{11},\kappa_{11}\}_1$ and $2\bar\alpha=\eta_3\varepsilon_4\kappa_{12}.$ 
		\qed
	\end{prop}
	
	At $n=5,$ the exact sequence over $S^4$ relates the calculation to the elements $\nu_4\bar\zeta_7$
	and $\nu_4\bar\sigma_7.$ The next proposition identifies their lifts and the resulting cyclic
	summands.
	
	\begin{prop}\label{m2:526}
		The structure of $\pi_{26}(\mathit{SO}(5))$ is as follows.
		\[ \GroupLayout{\pi_{26}(\mathit{SO}(5)) = \mathrm{span}\{[\nu_4\bar{\zeta}_7],\; [\nu_4\bar{\sigma}_7]\}\cong \mathbb Z/8\oplus \mathbb Z/2,}{\begin{aligned}\pi_{26}(\mathit{SO}(5)) = \mathrm{span}\{&[\nu_4\bar{\zeta}_7],\\
				&[\nu_4\bar{\sigma}_7]\}\cong \mathbb Z/8\oplus \mathbb Z/2,\end{aligned}} \]
		where $8[\nu_4\bar\zeta_7]=0$ and $2[\nu_4\bar\sigma_7]=0.$
		
	\end{prop}
	\begin{proof}
		Consider  the exact sequence
		\[\begin{aligned}
			&\pi^4_{27} \xrightarrow{ \Delta} \pi_{26}(\mathit{SO}(4)) \xrightarrow{i_\ast } \pi_{26}(\mathit{SO}(5)) \xrightarrow{p_\ast } \pi^4_{26} \xrightarrow{\Delta } \pi_{25}(\mathit{SO}(4)),
		\end{aligned}\]
		where
		$
		\pi^4_{27} = \frac{\mathbb Z}{8}\{\nu_4\bar{\kappa}_7\} \oplus \frac{\mathbb Z}{4}\{\Sigma\bar{\alpha}\}.
		$
		By Lemma~\ref{m2:lem:nuprime-compositions} and \cite[pp.~313--314]{Mi1}, we deduce
		$
		\Delta(\nu_4 \bar{\sigma}_7) = \Delta(\nu_4 \bar{\zeta}_7) = 0.
		$
		From the relation
		$
		\Delta(\iota_4) = 2[\iota_3] - [\eta_2]_4
		$
		\cite[Table 3]{Ka}, we see that $\Delta(\varepsilon_4 \kappa_{12}) = [\eta_2]_4 \varepsilon_3\kappa_{11}$ and $ \Delta(\Sigma\bar{\alpha}) = (2[\iota_3] - [\eta_2]_4)\bar{\alpha}.$ By \cite[(5.13)]{T},
		$(J\circ \Delta)(\nu_4\bar{\kappa}_7) = \pm2 \nu^2_4\bar{\kappa}_{10}.$
		By \cite[p.19]{MMO},
		$$J([\iota_3]\bar{\alpha}) = \nu_4(\Sigma^4\bar{\alpha}) = \pm2\nu^2_4 \bar{\kappa}_{10}.$$\noindent
		Since $J[\eta_2] = \nu'$ and $ P(\bar\kappa_9)=\pm2\nu_4\bar\kappa_7
		\pm \Sigma\bar\alpha $ \cite[(4.5)]{MMO}, we conclude that
		$$J([\eta_2]_4 \bar{\alpha}) = -(\Sigma\nu') \Sigma^4\bar{\alpha} = \pm2(\Sigma\nu')(\nu_7\bar{\kappa}_{10}) =0.$$ \noindent
		Then we infer that
		\begin{equation}
			\begin{aligned}
				&\Delta(\nu_4\bar{\kappa}_7) \equiv\pm  [\iota_3]\bar{\alpha}\,\, \bmod\, [\eta_2]_4\bar{\alpha}.\label{m2:Dn4bk}
			\end{aligned}
		\end{equation}
		In light of the preceding boundary calculations, $\Delta\colon  \pi^4_{27} \to \pi_{26}(\mathit{SO}(4))$ is surjective.
		This completes the proof.
	\end{proof}
	
	Using the  group structure of  $\pi_{26}(\mathit{SO}(4))$ and \eqref{m2:Dn4bk}, we obtain
	\begin{equation}
		\begin{aligned}
			&{i_5}_*\pi_{26}(\mathit{SO}(4))=0.\label{m2:iR426}
		\end{aligned}
	\end{equation}
	
	Having determined the image of $\pi_{26}(\mathit{SO}(4)),$ we turn to the fibration over $S^5.$ The
	calculation keeps track of the   elements   inherited from $n=5$ and the additional class represented by
	$[ \nu_5]\eta_8\bar\mu_9.$
	
	\begin{prop}\label{m2:626} The structure of $\pi_{26}(\mathit{SO}(6))$ is as follows.
		\[ \GroupLayout{\pi_{26}(\mathit{SO}(6)) = \mathrm{span}\{[\nu_4\bar{\zeta}_7]_6,\; [\nu_4\bar{\sigma}_7]_6,\; [ \nu_5]\eta_8\bar{\mu}_9\}\cong \mathbb Z/8\oplus (\mathbb Z/2)^{2},}{\begin{aligned}\pi_{26}(\mathit{SO}(6)) = \mathrm{span}\{&[\nu_4\bar{\zeta}_7]_6,\; [\nu_4\bar{\sigma}_7]_6,\\
				&[ \nu_5]\eta_8\bar{\mu}_9\}\cong \mathbb Z/8\oplus (\mathbb Z/2)^{2},\end{aligned}} \] where
		$2[ \nu_5]\eta_8\bar\mu_9=0.$
	\end{prop}
	\begin{proof}
		Consider  the exact sequence
		\[\begin{aligned}
			&\pi^5_{27} \xrightarrow{ \Delta} \pi_{26}(\mathit{SO}(5)) \xrightarrow{i_\ast } \pi_{26}(\mathit{SO}(6)) \xrightarrow{p_\ast } \pi^5_{26} \xrightarrow{\Delta } \pi_{25}(\mathit{SO}(5)),
		\end{aligned}\]
		where
		$
		\pi^5_{27} = (\frac{\mathbb Z}{2})^3\{\varepsilon_5\kappa_{13},\nu_5\bar{\zeta}_8,\nu_5\bar{\sigma}_8\}.
		$
		Since
		$
		\Delta \nu_5 = 0
		$
		\cite[Table 3]{Ka},
		we deduce
		$\Delta(\nu_5\bar{\zeta}_8) = \Delta(\nu_5\bar{\sigma}_8) = 0\,\uprightsemicolon{}\, \Delta(\nu_5\eta_8\bar{\mu}_9) = 0.$
		Since
		$
		\Delta \iota_5 = [\iota_3]_5 \eta_3
		$
		\cite[Table 3]{Ka}, we have
		$$\Delta (\varepsilon_5 \kappa_{13} ) = [\iota_3]_5\eta_3\varepsilon_4\kappa_{12} = [\iota_3]_5\circ 2\bar{\alpha}.$$ 
		Since $\Delta : \pi^4_{27} \to \pi_{26}(\mathit{SO}(4))$ is a surjective homomorphism,
		$[\iota_3]_5\circ 2\bar{\alpha} = i_\ast([\iota_3]\circ
		2\bar{\alpha}) =0.$ 
		For the  element  $\alpha\in\pi_{26}^5$ fixed in Lemma~\ref{g1016:lem:delta-alpha}, we note that
		\[\FitMath{\Sigma\alpha=\eta_6\bar\kappa_7,\quad\Delta\alpha=[\iota_3]_5\varepsilon_3\kappa_{11}.}\]\noindent       Thus we have a short exact sequence
		\[\begin{aligned}
			&0 \to \frac{\mathbb Z}{8}\{[\nu_4\bar{\zeta}_7]\} \oplus \frac{\mathbb Z}{2}\{[\nu_4\bar{\sigma}_7]\} \xrightarrow{i_\ast } \pi_{26}(\mathit{SO}(6)) \xrightarrow{p_\ast } \frac{\mathbb Z}{2}\{\nu_5\eta_8\bar{\mu}_9\} \to 0.
		\end{aligned}\]
		Since the lift $[ \nu_5]\eta_8\bar{\mu}_9 \in \pi_{26}(\mathit{SO}(6))$ of
		$\nu_5\eta_8\bar{\mu}_9$ has order 2, we see that the homotopy group $\pi_{26}(\mathit{SO}(6)) $ is equal to $ \frac{\mathbb Z}{8}\{[\nu_4\bar{\zeta}_7]_6\} \oplus (\frac{\mathbb Z}{2})^2\{[\nu_4\bar{\sigma}_7]_6, [ \nu_5]\eta_8\bar{\mu}_9\}.$ \end{proof}

	\subsection{The Cases $7\leq n\leq9$}
	
	For the passage to $n=7,$ we return to the degree-$26$ calculation for $G_2.$ We give the lift
	relations explicitly here, since their coefficients control the extension involving $2P(\rho_{13}).$
	
	Recall that the structure of $\pi_{26}(G_2)$ is as follows:
	\[ \GroupLayout{\pi_{26}(G_2) = \mathrm{span}\{\left<2P(\rho_{13}) \right>, \x\nu^2_6 \kappa_{12} \xx\}\cong \mathbb Z/32\oplus \mathbb Z/2.}{\begin{aligned}\pi_{26}(G_2) = \mathrm{span}\{&\left<2P(\rho_{13}) \right>,\\
			&\x\nu^2_6 \kappa_{12} \xx\}\cong \mathbb Z/32\oplus \mathbb Z/2.\end{aligned}} \] The generators satisfy, for an odd integer $x,$  
	\begin{equation}\label{g2bm1024}
		2\left<2P(\rho_{13})\right>=x\left<[2\e_{6},\e_{6}]\right>\Sigma^2\rho',\qquad 16\left<2P(\rho_{13})\right>=8\left<[2\e_{6},\e_{6}]\right>\Sigma^2\rho' =i_{G*}[\nu_5\eta_8\bar{\mu}_9].
	\end{equation}
	\noindent The preceding calculations can now be assembled in the exact sequence for
	$\mathit{SO}(6)\to\mathit{SO}(7)\to S^6.$ The following proposition expresses the result using lifts of elements of homotopy groups of spheres and images of the generators from $n=6.$
	
	\begin{prop}\label{m2:726} The group $\pi_{26}(\mathit{SO}(7))$  is as follows.
		\[\GroupLayout{\pi_{26}(\mathit{SO}(7)) = \mathrm{span}\{[2P(\rho_{13})],\; [\nu_4\bar\zeta_7]_7,\; [\nu_4\bar\sigma_7]_7,\; [\nu_6^2\kappa_{12}]\}\cong \mathbb Z/32\oplus \mathbb Z/8\oplus (\mathbb Z/2)^{2}.}{\begin{aligned}\pi_{26}(\mathit{SO}(7)) = \mathrm{span}\{&[2P(\rho_{13})],\; [\nu_4\bar\zeta_7]_7,\; [\nu_4\bar\sigma_7]_7,\\
				&[\nu_6^2\kappa_{12}]\}\cong \mathbb Z/32\oplus \mathbb Z/8\oplus (\mathbb Z/2)^{2}.\end{aligned}}\]
		The group extension relations are
		\[
		\begin{aligned}
			&16[2P(\rho_{13})]=[ \nu_5]_7\eta_8\bar\mu_9+4[\nu_4\bar\zeta_7]_7,\qquad 2[\nu_6^2\kappa_{12}]=0.
		\end{aligned}
		\]
	\end{prop}
	\begin{proof}
		Consider  the exact sequence
		\[\begin{aligned}
			&\pi^6_{27} \xrightarrow{ \Delta} \pi_{26}(\mathit{SO}(6)) \xrightarrow{i_\ast } \pi_{26}(\mathit{SO}(7)) \xrightarrow{p_\ast } \pi^6_{26} \xrightarrow{\Delta } \pi_{25}(\mathit{SO}(6)),
		\end{aligned}\]
		where
		$
		\pi^6_{27} = \frac{\mathbb Z}{2}\{\eta_6\bar{\kappa}_7\}.
		$
		Since
		$
		\Delta(\eta_6) = 0
		$
		\cite[Lemma 1.1(i)]{KM1}, we obtain
		$
		\Delta (\eta_6\bar{\kappa}_7) = 0.
		$
		By Lemma~\ref{m1:mtb2}, we deduce $\Delta P(\rho_{13}) =
		[ \nu_5]\eta_8\mu_9\sigma_{18}+4[\nu_4\zeta_7]_6\sigma_{18}\ \mbox{and} \ \Delta
		\bar{\kappa}'=[ \nu_5]\nu_8\kappa_{11}.$ By \cite[Lemma 15.3]{MT1}, we have an exact sequence
		\[\begin{aligned}
			&0 \to \frac{\mathbb Z}{8}\{[\nu_4\bar{\zeta}_7]\} \oplus (\frac{\mathbb Z}{2})^2\{[\nu_4\bar{\sigma}_7]_6, [ \nu_5]\eta_8\bar{\mu}_9\} \xrightarrow{i_\ast } \pi_{26}(\mathit{SO}(7)) \xrightarrow{p_\ast } \frac{\mathbb Z}{16}\{2P(\rho_{13})\} \oplus \frac{\mathbb Z}{2}\{\nu^2_6\kappa_{12}\} \to 0.
		\end{aligned}\]
		The two short exact sequences fit into the following commutative
		diagram.  The groups and generators in the outer terms are those
		specified above.
		\[
		\begin{tikzcd}[column sep=5mm,row sep=large]
			0\arrow[r]&
			\frac{\mathbb Z}{2}\{[\nu_5\eta_8\bar\mu_9]\}
			\arrow[r,"i_{G*}"]\arrow[d,"r_{3*}"']&
			\pi_{26}(G_2)\arrow[r,"p_{G*}"]\arrow[d,"\mathcal{T}_{*}"]&
			\frac{\mathbb Z}{16}\{2P(\rho_{13})\}\oplus\frac{\mathbb Z}{2}\{\nu_6^2\kappa_{12}\}
			\arrow[r]\arrow[d,equal]&0\\
			0\arrow[r]&\pi_{26}(\mathit{SO}(6))\arrow[r,"i_*"]&
			\pi_{26}(\mathit{SO}(7))\arrow[r,"p_*"]&
			\frac{\mathbb Z}{16}\{2P(\rho_{13})\}\oplus\frac{\mathbb Z}{2}\{\nu_6^2\kappa_{12}\}
			\arrow[r]&0
		\end{tikzcd}
		\]
		Using the homomorphism ${p_6}_\ast\colon \pi_{26}(\mathit{SO}(6)) \to \pi^5_{26}, $ we obtain
		$$r_{3\ast}[\nu_5\eta_8\bar\mu_9] \equiv [ \nu_5]\eta_8\bar\mu_9 \,\,\, \bmod
		\,\,\,4[\nu_4\bar\zeta_7]_6, [\nu_4\bar\sigma_7]_6.$$ Then we
		conclude that $16[ 2P(\rho_{13}) ] \equiv [ \nu_5]_7\eta_8\bar{\mu}_9 \,\,\, \bmod \,\,\,4[\nu_4\bar\zeta_7]_7,
		[\nu_4\bar\sigma_7]_7$ by Equation \eqref{g2bm1024}. From \cite[(8.1),(8.6)]{Mi1} and \cite[Lemmas 10.6, 12.4]{T}, we infer that $$(H\circ
		J)[\nu_4\bar{\zeta}_7] = \nu_9\bar{\zeta}_{12}=8\sigma_9\rho_{16}=8\rho'\sigma_{24} = H(\bar{\zeta}_5
		\sigma_{24}).$$ This implies that
		$$J[\nu_4\bar{\zeta}_7] \equiv \bar{\zeta}_5\sigma_{24} \,\, \bmod \,\Sigma\pi^4_{30}.$$ \noindent 
		In light of \cite[(4.6)]{MMO}, we have
		$$(H \circ J)[\nu_4\bar{\sigma}_7] = \nu_9\bar\sigma_{12}=\bar\sigma_9\nu_{28}=H(\phi_5\nu_{28}).$$ \noindent 
		It follows that
		$J[\nu_4\bar{\sigma}_7] \equiv \phi_5\nu_{28} \,\, \bmod \, \Sigma\pi^4_{30}.$
		By \cite[Lemma 15.4]{MT1}, we deduce
		$(H \circ J)[\nu^2_6\kappa_{12}] = \nu^2_{13}\kappa_{19} = 4\bar{\kappa}_{13}.$
		By \cite[Lemma 2.1 (ii)]{HKM}, \cite[I-Proposition 3.5 (2)]{Oda1979},
		\cite[(16.6)]{MT1} and \cite[Lemma 12.4]{T}, we find $$J([ \nu_5]_7\eta_8\bar\mu_9) = (\bar\nu_7+\varepsilon_7)\eta_{15}\bar\mu_{16}=\eta^2_7\bar\mu_9\sigma_{26}=4\bar\zeta_7\sigma_{26}.$$\noindent  Since
		$
		J(16[2P(\rho_{13})]) \in 16\pi^7_{33}=0,
		$ we conclude that
		\[\begin{aligned}
			&16[ 2P(\rho_{13}) ] = [ \nu_5]_7\eta_8\bar{\mu}_9+4[\nu_4\bar\zeta_7]_7,
		\end{aligned}\]
		\noindent by comparing these $J$-images. Noting that $\nu^2_6\kappa_{12} = \nu_6\kappa_9\nu_{23}$
		and
		$$\theta \nu_{24} \in \sigma_{12} \circ \{\nu_{19},\eta_{22},\nu_{23}\}_1=\sigma_{12}\bar\nu_{19}=0,$$ \noindent we define $[\nu^2_6\kappa_{12}]$
		as
		$[\nu^2_6\kappa_{12}]=[P(\Sigma\theta)+\nu_6\kappa_9]\nu_{23}.$ According to the relation
		$$2[P(\Sigma\theta)+\nu_6\kappa_9]\equiv[ \nu_5]_7\bar\varepsilon_8\ \bmod\ [\nu_4(\Sigma\zeta')]_7,
		[\nu_4\eta_7\bar{\varepsilon}_8]_7$$
		\noindent from  the case $n=7$ in Theorem~\ref{m1:thm1},
		we deduce that
		$$2[\nu^2_6\kappa_{12}]=(2[P(\Sigma\theta)+\nu_6\kappa_9])\nu_{23}\equiv[ \nu_5]_7\bar\varepsilon_8\nu_{23}=0\
		\bmod\ [\nu_4(\Sigma\zeta')]_7\nu_{23}, [\nu_4\eta_7\bar{\varepsilon}_8]_7\nu_{23}.$$\noindent Here,
		$[\nu_4(\Sigma\zeta')]_7\nu_{23}, [\nu_4\eta_7\bar{\varepsilon}_8]_7\nu_{23}\in
		{i_{4,7}}_*\pi_{26}(\mathit{SO}(4))=0$ \eqref{m2:iR426}.
		This completes the proof.
	\end{proof}
	
	Since $ \pi_{26}(\mathit{SO}(8)) \cong \pi_{26}(\mathit{SO}(7)) \oplus \pi^7_{26} $ and $ \pi^7_{26}
	= \frac{\mathbb Z}{8}\{\bar{\zeta}_7\} \oplus \frac{\mathbb Z}{2}\{\bar{\sigma}_7\}, $ we obtain the splitting used below.

	\begin{prop}\label{m2:826} The structure of $\pi_{26}(\mathit{SO}(8))$ is as follows.
		\[\GroupLayout{\pi_{26}(\mathit{SO}(8)) \;= {} \mathrm{span}\{[2P(\rho_{13})]_8,\; [\nu_4\bar\zeta_7]_8,\; [\iota_7]\bar\zeta_7,\; [\nu_4\bar\sigma_7]_8,\; [\nu_6^2\kappa_{12}]_8,\; [\iota_7]\bar\sigma_7\}\cong \mathbb Z/32\oplus (\mathbb Z/8)^{2}\oplus (\mathbb Z/2)^{3}.}{\begin{aligned}\pi_{26}(\mathit{SO}(8)) \;= {} \mathrm{span}\{&[2P(\rho_{13})]_8,\; [\nu_4\bar\zeta_7]_8,\; [\iota_7]\bar\zeta_7,\; [\nu_4\bar\sigma_7]_8,\\
				&[\nu_6^2\kappa_{12}]_8,\; [\iota_7]\bar\sigma_7\}\cong \mathbb Z/32\oplus (\mathbb Z/8)^{2}\oplus (\mathbb Z/2)^{3}.\end{aligned}}\]
		The homomorphism induced by the inclusion and the chosen lift of
		$\iota_7$ is an isomorphism:
		\[\begin{aligned}
			&(i_8)_*+([\iota_7])_*: \pi_{26}(\mathit{SO}(7))\oplus\pi_{26}^7\xrightarrow{\cong}\pi_{26}(\mathit{SO}(8)), \qquad ([\iota_7])_*(\alpha)=[\iota_7]\circ\alpha.
		\end{aligned}\]
		Under this isomorphism, the summands generated by
		$[\iota_7]\bar\zeta_7$ and $[\iota_7]\bar\sigma_7$ are the images of
		$\mathrm{span}\{\bar\zeta_7\}\cong \mathbb Z/8$ and
		$\mathrm{span}\{\bar\sigma_7\}\cong \mathbb Z/2,$   respectively.
		\qed
	\end{prop}
	
	Before moving to $n=9,$ we compare two descriptions of the   elements   involving $\bar\sigma_7$ and
	$\bar\zeta_7.$ These comparisons express the boundary formulas in the same notation as the generators
	already chosen.
	
	\begin{lem}\label{m2:lem:rank7-sigma-zeta}
		\begin{itemize}
			\item[\normalfont(1)]$[\eta_6]\bar{\sigma}_7 = [\nu_4\bar{\sigma}_7]_7.$
			\item[\normalfont(2)]
			$[\eta_6]\bar{\zeta}_7\ty [\nu_4\bar{\zeta}_7]_7\m8[[2\e_{6},\e_{6}]]\Sigma^2\rho'.$
		\end{itemize}
	\end{lem}
	
	\begin{proof}
		\begin{enumerate}[label=\textup{(\arabic*)},leftmargin=2em]
			\item By the definition of $\bar{\sigma}_6$ \cite[p.138]{T}, we obtain $$\nu_4\bar\sigma_7 \in
			\nu_4\circ\{\nu_7, \bar\nu_{10}+\varepsilon_{10}, \sigma_{18}\}_1 \subseteq\{\nu^2_4,
			\bar\nu_{10}+\varepsilon_{10}, \sigma_{18}\}_1.$$ \noindent By \cite[Lemma 1.1(ii)]{KM1} and by the relation
			$\nu_5\varepsilon _8=\nu_5\sigma_8\eta_{15},$   we deduce $J([  \nu^2_4]\eta_{10}) = J([\iota_3]_5\varepsilon_3) = \nu_5\varepsilon_8.$
			Then we obtain $[  \nu^2_4](\bar\nu_{10} + \varepsilon_{10}) = [  \nu^2_4]\eta_{10}\sigma_{11} =
			[\iota_3]_5\varepsilon_3\sigma_{11}=0.$
			Then we can take
			$[\nu_4\bar\sigma_7] \in \{[  \nu^2_4], \bar\nu_{10}+\varepsilon_{10}, \sigma_{18}\}_1.$
			Since
			$
			J[  \nu^2_4]_7 = J([\eta_6]\nu_7) = \sigma'\nu_{14} = a\nu_7\sigma_{10}
			$
			\cite[Lemma 1.1(iv)]{KM1}, we find
			$
			[  \nu^2_4]_7 = a[\eta_6]\nu_7
			$
			for an odd $a.$ 
			It follows that
			\begin{align*}
				\notag [\nu_4\bar\sigma_7]_7 &\in  \{[  \nu^2_4]_7, \bar\nu_{10}+\varepsilon_{10}, \sigma_{18}\}_1\\
				\notag &=  \{a[\eta_6]\nu_7,\bar\nu_{10}+\varepsilon_{10}, \sigma_{18}\}_1\\
				\notag &\ni  [\eta_6]\bar\sigma_7 \,\,\, \bmod \,\,\, [\eta_6]\nu_7\circ \Sigma\pi^9_{25}+\pi_{19}(\mathit{SO}(7)) \circ \sigma_{19},
			\end{align*}
			where
			$
			\nu_7\circ \Sigma\pi^9_{25} = \mathrm{span}\{\nu_7\sigma_{10}\nu^3_{17}, \nu_7\sigma_{10}\mu_{17}, \nu_7\sigma_{10}\eta_{17}\varepsilon_{18}, \nu_7\mu_{10}\sigma_{19}\}.
			$
			The standard vanishing relations give
			$
			\nu_7\sigma_{10}\nu^3_{17} = \nu_7\sigma_{10}\eta_{17}\varepsilon_{18} =0.
			$
			We have
			$
			\nu_7\mu_{10}\sigma_{19} = \nu_7\sigma_{10}\mu_{17} = 0.
			$
			By \cite[Theorem 2.1]{KM2}, we obtain
			$
			\pi_{19}(\mathit{SO}(7)) \circ \sigma_{19} = 0.
			$
			It follows that
			$[\nu_4\bar\sigma_7]_7 = [\eta_6]\bar\sigma_7.$
			
			\item By  \cite[Theorem 2]{KM1}, we deduce
			$[\nu_4\bar\zeta_7]_7 \in \{[\nu_4\zeta_7]_7, 8\iota_{18}, 2\sigma_{18}\}_1.$
			Using the relation \cite[Lemma 3.4(iii)]{KM1}
			$$[\nu_4\zeta_7]_7 \equiv [\eta_6]\zeta_7\m8[[2\e_{6},\e_{6}]]\sigma_{11},$$ \noindent      we find
			\begin{align*}
				\notag \{[\nu_4\zeta_7]_7, 8\iota_{18}, 2\sigma_{18}\}_1 &=  \{[\eta_6]\zeta_7+8a[[2\e_{6},\e_{6}]]\sigma_{11}, 8\iota_{18}, 2\sigma_{18}\}_1\\
				\notag &\subseteq  \{[\eta_6]\zeta_7 , 8\iota_{18}, 2\sigma_{18}\}_1 + \{ 8a[[2\e_{6},\e_{6}]]\sigma_{11}, 8\iota_{18}, 2\sigma_{18}\}_1\\
				\notag &\ni  [\eta_6]\bar\zeta_7 + 8a[[2\e_{6},\e_{6}]]\Sigma^2\rho'\;\;(a\in\z)
			\end{align*}
			with indeterminacy
			$
			[\eta_6]\zeta_7\circ \Sigma\pi^{17}_{25} + \pi_{19}(\mathit{SO}(7))\circ 2\sigma_{19} + 8a[[2\e_{6},\e_{6}]]\sigma_{11} \circ \Sigma\pi^{17}_{25}.
			$
			By \cite[Theorem 2.1]{KM2} and the standard sphere composition relations, the indeterminacy is 0. So, 
			$[\nu_4\bar\zeta_7] _7=[\eta_6]\bar\zeta_7 + 8a[[2\e_{6},\e_{6}]]\Sigma^2\rho'.$

		\end{enumerate}
	\end{proof}
	
	We now apply those comparisons to the fibration over $S^8.$ Computing the boundaries of $\bar\zeta_8$
	and $\bar\sigma_8$ identifies the relations imposed on the $n=8$ generators.
	
	\begin{prop}\label{g1016:prop:26-9} The structure of $\pi_{26}(\mathit{SO}(9))$ is as follows.
		\[ \GroupLayout{\pi_{26}(\mathit{SO}(9)) = \mathrm{span}\{[2P(\rho_{13})]_9,\; [\iota_7]_9\bar{\zeta}_7,\; [\nu^2_6\kappa_{12}]_9,\; [\iota_7]_9\bar{\sigma}_7\}\cong \mathbb Z/32\oplus \mathbb Z/8\oplus (\mathbb Z/2)^{2},}{\begin{aligned}\pi_{26}(\mathit{SO}(9)) = \mathrm{span}\{&[2P(\rho_{13})]_9,\; [\iota_7]_9\bar{\zeta}_7,\; [\nu^2_6\kappa_{12}]_9,\\
				&[\iota_7]_9\bar{\sigma}_7\}\cong \mathbb Z/32\oplus \mathbb Z/8\oplus (\mathbb Z/2)^{2}.\end{aligned}} \] Additionally, $2[\iota_7]_9\bar\zeta_7\ty[\nu_4\bar\zeta_7]_9\m8[[2\e_{6},\e_{6}]]\Sigma^2\rho'; \;[\nu_4\bar\sigma_7]_9=0.$
	\end{prop}
	
	\begin{proof}
		Consider  the exact sequence $(26)_8$:
		\[\begin{aligned}
			&\pi^8_{27} \xrightarrow{ \Delta} \pi_{26}(\mathit{SO}(8)) \xrightarrow{ i_\ast} \pi_{26}(\mathit{SO}(9)) \xrightarrow{ p_\ast} \pi^8_{26} \xrightarrow{\Delta} \pi_{25}(\mathit{SO}(8)),
		\end{aligned}\]
		where
		$
		\pi^8_{27} = \frac{\mathbb Z}{8}\{\bar{\zeta}_8\} \oplus \frac{\mathbb Z}{2}\{\bar{\sigma}_8\}.
		$
		Using the relation
		$\Delta \iota_8 = 2[\iota_7]-[\eta_6]_8
		$
		and Lemma~\ref{m2:lem:rank7-sigma-zeta}, we conclude that
		\begin{align*}
			\Delta\bar{\zeta}_8 &=  2[\iota_7]\bar{\zeta}_7-[\eta_6]_8\bar{\zeta}_7 \ty2[\iota_7]\bar{\zeta}_7 - [\nu_4\bar{\zeta}_7]_8\m8[[2\e_{6},\e_{6}]]\Sigma^2\rho',\\
			\Delta\bar{\sigma}_8    &=  [\eta_6]_8\bar{\sigma}_7 = [\nu_4\bar{\sigma}_7]_8.
		\end{align*}\noindent       By Propositions~\ref{g1016:prop:25-8} and~\ref{m1:R259}(1), with $i_*=(i_9)_*,$  
		\[\FitMath{\begin{aligned}
				\operatorname{Im}\Delta&=\ker\bigl(i_*:\pi_{25}(\mathit{SO}(8))\to\pi_{25}(\mathit{SO}(9))\bigr),\\
				|\operatorname{Im}\Delta|&=|\pi_{26}^8|.
		\end{aligned}}\]
		\noindent       Hence, $\ker\bigl(\Delta\colon \pi_{26}^8\to\pi_{25}(\mathit{SO}(8))\bigr)=0.$    This completes the
		proof.  \end{proof}
	
	\subsection{The Cases $n=10,11$}
	
	The $n=10$ calculation requires a relation between a composite of a lift of $\sigma_9^2$ and an
	inherited generator. The following identity permits these two expressions to be used interchangeably
	in the subsequent argument.
	
	\begin{lem}
		$ [\sigma^2_9]\nu_{23} = [\iota_7]_{10}\bar\sigma_7. $
	\end{lem}
	
	\begin{proof}
		By \cite[Proposition 1.4]{T}, \eqref{m1:[sig^2_9]} and
		\cite[I-Proposition 3.4 (3)]{Oda1979}, we see that $$[\sigma^2_9]\nu_{23} \in \{[\iota_7]_{10},\bar\nu_7+\varepsilon_7,\sigma_{15}\} \circ \nu_{23} = [\iota_7]_{10} \circ \{\bar\nu_7+\varepsilon_7,\sigma_{15},\nu_{22}\}=[\iota_7]_{10}\bar\sigma_7.$$ This completes the proof.
	\end{proof}
	
	We next examine the $J$-image of $[\eta_9\varepsilon_{10}].$ Identifying this image is the first step
	in comparing the Toda-bracket lift used below with the corresponding element of a homotopy group of a sphere.
	
	\begin{lem}\label{m2:J[et9ep}
		For an odd $x,$  
		$J[\eta_9\varepsilon_{10}] =x \lambda''.$
	\end{lem}
	
	\begin{proof}
		By \cite[Proposition 3 (3)]{Oda}, we have
		$J[\eta_9\varepsilon_{10}] \equiv \lambda'' \,\,\, \bmod \,\,\Sigma\pi^9_{27}.$
		Recall that
		$$\begin{aligned}
			&\Sigma\pi^9_{27} = \frac{\mathbb Z}{4}\{\sigma_{10}\zeta_{17}\} \oplus \frac{\mathbb Z}{2}\{\eta_{10}\bar{\mu}_{11}\}
		\end{aligned}$$
		and that
		$2\lambda'' \equiv \sigma_{10}\zeta_{17} \,\,\,\bmod\,\,2\sigma_{10}\zeta_{17}.$
		So,
		$J[\eta_9\varepsilon_{10}] \equiv \lambda '' \,\,\,\bmod\,\,2\lambda'', \eta_{10}\bar{\mu}_{11}.$
		Note that $ [\eta_9\varepsilon_{10}]_\infty \in \pi_{18}(\mathit{SO})=0, $ $
		\Sigma^\infty(\eta_{10}\bar{\mu}_{11})=\eta\bar{\mu} $ and $
		\Sigma^7\lambda''=0 $ by \cite[Lemma 12.18, 12.19]{T}, where $ \pi^S_{18} = \mathrm{span}\{\nu^\ast,\eta\bar{\mu}\} .$  This completes the proof.
	\end{proof}
	
	By the fourth suspension of \eqref{m1:Ogu0},
	$\eta_9\mu_{10}\sigma_{19} = \mu_9\varepsilon_{18} = \varepsilon_9\mu_{17}.$
	Since $$\begin{aligned}
		&\varepsilon_9\mu_{17} \in \varepsilon_9 \circ \{\eta_{17},2\iota_{18}, 8\sigma_{18}\}_1 \subseteq \{\eta_9\varepsilon_{10}, 8\iota_{18},2\sigma_{18}\}_1,
	\end{aligned}$$ we define
	\begin{equation}
		\begin{aligned}
			&[\eta_9\mu_{10}\sigma_{19}] \in \{[\eta_9\varepsilon_{10}], 8\iota_{18},2\sigma_{18}\}_1.\label{m2:etms}
		\end{aligned}
	\end{equation}
	With the lift in~\eqref{m2:etms} fixed, we can compare its $J$-image with $\tau''.$ The congruence
	below records precisely the ambiguity coming from suspension.
	
	\begin{lem}
		$ J[\eta_9\mu_{10}\sigma_{19}] \equiv \tau''
		\,\,\,\bmod\,\,\Sigma\pi^9_{35}. $
	\end{lem}
	
	\begin{proof}
		By Lemma~\ref{m2:J[et9ep}, for an odd integer $x,$  
		\[
		\begin{aligned}
			&J\{[\eta_9\varepsilon_{10}],8\iota_{18},2\sigma_{18}\}_1 \subseteq\{J[\eta_9\varepsilon_{10}],8\iota_{28},2\sigma_{28}\}_1 =\{x\lambda'',8\iota_{28},2\sigma_{28}\}_1.
		\end{aligned}
		\]
		It follows from \cite[Proposition 3(3)]{Oda} that
		\begin{align*}
			H\{x\lambda'',8\iota_{28},2\sigma_{28}\}_1
			\notag &\subseteq   \{H(x\lambda''),8\iota_{28},2\sigma_{28}\}_1 \\
			\notag &=  \{\eta_{19}\varepsilon_{20},8\iota_{28},2\sigma_{28}\}_1 \\
			\notag &\supseteq   \varepsilon_{19}\circ \{\eta_{27},8\iota_{28}, 2\sigma_{28}\}_1\\
			\notag &\supseteq   \varepsilon_{19}\circ \{\eta_{27},2\iota_{28},8\sigma_{28}\}_1\\
			\notag &\ni   \varepsilon_{19}\mu_{27}
		\end{align*}
		with indeterminacy $ \eta_{19}\varepsilon_{20}\circ \Sigma\pi^{27}_{35} +
		\pi^{19}_{29}\circ 2\sigma_{29} = 0. $ By \cite[I-Proposition 4.2]{Oda1979}
		and the standard composition relations, we deduce $ H(\tau'')=  \varepsilon_{19}\mu_{27}.
		$ This completes the proof.
	\end{proof}
	
	The $J$-image alone does not settle the extension in the special orthogonal group. We therefore compute a
	multiple of the lift itself and record an annihilating multiple\uprightsemicolon{} these relations will be used in the
	$n=10$ decomposition.
	
	\begin{lem}\label{m2:2[et9ms}
		The chosen lift satisfies $8[\eta_9\mu_{10}\sigma_{19}]=0$ and
		\[\begin{aligned}
			&2[\eta_9\mu_{10}\sigma_{19}] =[\iota_7]_{10}\bar\zeta_7+2[[2\e_{6},\e_{6}]]_{10}\Sigma^2\rho'.
		\end{aligned}\]
	\end{lem}
	\begin{proof}
		Since $\{8\iota_{17},2\sigma_{17},8\iota_{24}\}=0,$   we obtain
		\[
		\begin{aligned}
			&8\{[\eta_9\varepsilon_{10}],8\iota_{18},2\sigma_{18}\}_1 =-[\eta_9\varepsilon_{10}]\circ \Sigma\{8\iota_{17},2\sigma_{17},8\iota_{24}\} =0.
		\end{aligned}
		\]
		This proves the first assertion.
		
		Since $2[\eta_9\varepsilon_{10}]=[\iota_7]_{10}\zeta_7+2[[2\e_{6},\e_{6}]]_{10}\sigma_{11}$ \cite[p. 342]{KM1},
		we find
		\[\begin{aligned}
			&2[\eta_9\mu_{10}\sigma_{19}] \in \{2[\eta_9\varepsilon_{10}], 8\iota_{18}, 2\sigma_{18}\}_1 = \{[\iota_7]_{10}\zeta_7 + 2[[2\e_{6},\e_{6}]]_{10}\sigma_{11}, 8\iota_{18}, 2\sigma_{18} \},
		\end{aligned}\]
		which is contained in
		$\{[\iota_7]_{10} \zeta_7 , 8\iota_{18}, 2\sigma_{18} \}  +  \{2[[2\e_{6},\e_{6}]]_{10}\sigma_{11}, 8\iota_{18}, 2\sigma_{18} \}.$ Notice that $$\{[\iota_7]_{10} \zeta_7 , 8\iota_{18}, 2\sigma_{18} \} \ni[\iota_7]_{10}\bar{\zeta}_7 \
		\bmod\ [\iota_7]_{10} \zeta_7\circ\pi^{18}_{26} +\pi_{19}(\mathit{SO}(10))\circ 2\sigma_{19}.$$ Then we see that
		$\pi_{19}(\mathit{SO}(10))\circ 2\sigma_{19}=\mathrm{span}\{\Delta(\sigma_{10}\nu_{17})\circ
		2\sigma_{19}\}=0$ \cite[Theorem 2.1, p. 96]{KM2}. Using 
		$\zeta_5\bar{\nu}_{16}=\zeta_7\varepsilon_{18}=0,$   we obtain $[\iota_7]_{10} \zeta_7\circ\pi^{18}_{26}=0.$ That is,
		$\{[\iota_7]_{10}\zeta_7,8\iota_{18},2\sigma_{18}\}=[\iota_7]_{10}\bar\zeta_7.$ By \eqref{m2:2e2rh'}, we deduce
		\[\begin{aligned}
			&\{2[[2\e_{6},\e_{6}]]_{10}\sigma_{11}, 8\iota_{18}, 2\sigma_{18} \} \supseteq [[2\e_{6},\e_{6}]]_{10}\circ\{2\sigma_{11},8\iota_{18},2\sigma_{18}\} =2[[2\e_{6},\e_{6}]]_{10}\Sigma^2\rho'
		\end{aligned}\]
		$\bmod\ 2[[2\e_{6},\e_{6}]]_{10}\sigma_{11}\circ\pi^{18}_{26}+\pi_{19}(\mathit{SO}(10))\circ 2\sigma_{19}.$
		We have $2[[2\e_{6},\e_{6}]]_{10}\sigma_{11}\circ\pi^{18}_{26}=0$ and
		$\pi_{19}(\mathit{SO}(10))\circ 2\sigma_{19}=0.$ This implies
		$\{2[[2\e_{6},\e_{6}]]_{10}\sigma_{11},8\iota_{18},2\sigma_{18}\}=2[[2\e_{6},\e_{6}]]_{10}\Sigma^2\rho'.$
		This proves the second assertion.
	\end{proof}
	
	Two further composition identities are needed to control the   element   involving $\nu\kappa.$ The first is
	a vanishing relation, and the second identifies the corresponding connecting image in terms of a
	previously chosen lift.
	
	\begin{lem}\label{m2:Dnuka}
		\begin{itemize}
			\item[\normalfont(1)] $[\bar\nu_9\nu^2_{17}]\nu_{23} = [P(\eta_{19})]\bar\nu_{18}=0.$
			\item[\normalfont(2)] $\Delta(\nu_{10}\kappa_{13}) = [\nu^2_6\kappa_{12}]_{10}.$
		\end{itemize}
	\end{lem}
	
	\begin{proof}
		\begin{enumerate}[label=\textup{(\arabic*)},leftmargin=2em]
			\item By \eqref{m1:[bn9n^2]}, we conclude that
			\begin{align*}
				\notag [\bar\nu_9\nu^2_{17}] \nu_{23} &\in  \{[P(\eta_{19})], \nu_{18}, \eta_{21}\} \circ \nu_{23}\\
				\notag &=  [P(\eta_{19})] \circ \{\nu_{18}, \eta_{21}, \nu_{22}\} \\
				\notag &\ni  [P(\eta_{19})]\bar\nu_{18} .
			\end{align*}
			Using the relation $\Delta(P(\iota_{21}))=[P(\eta_{19})]+2[\eta_9\varepsilon_{10}]$
			\cite[p. 343]{KM1}, we see that $\Delta(P(\bar{\nu}_{21}))=[P(\eta_{19})]\bar{\nu}_{18}.$ Since
			$
			\Sigma:\pi^{10}_{27} \to \pi^{11}_{28}
			$
			is an isomorphism, we obtain $P(\bar\nu_{21})=0.$ This proves (1).
			
			\item We record a relation
			$
			\Delta(\kappa_{10}) \equiv [\bar\nu_9\nu^2_{17}]+[P(\Sigma\theta)+\nu_6\kappa_9]_{10}
			\bmod\ 16[\sigma_9^2]
			$
			from Proposition~\ref{m1:R11}.
			Since $\nu_{23}$ is a suspension of order $8$,
			right composition with $\nu_{23}$ is additive, and
			$(16[\sigma_9^2])\nu_{23}=16([\sigma_9^2]\nu_{23})
			=[\sigma_9^2](16\nu_{23})=0.$
			Thus the congruence becomes an equality after right composition with $\nu_{23}$.
			Using $\kappa_{10}\nu_{24}=\nu_{10}\kappa_{13}$, the boundary composition formula,
			and (1), we deduce $\Delta(\nu_{10}\kappa_{13})=[\bar\nu_9\nu^2_{17}]\nu_{23} +
			[P(\Sigma\theta)+\nu_6\kappa_9]_{10}\nu_{23}=[\nu^2_6\kappa_{12}]_{10}.$ This proves (2).
		\end{enumerate}
	\end{proof}
	
	We quote from \cite[Lemma 2.4]{KM1}:
	\begin{equation}
		\begin{aligned}
			&\Delta(\eta_{10})=2[\iota_7]_{10}\nu_7.\label{m2:Det10}
		\end{aligned}
	\end{equation}
	\noindent In light of Lemma~\ref{m2:lem:nu-kappa-bracket}, we define a matrix Toda bracket
	$$
	\left\{
	\bigl({{\Delta(\iota_{10})}}, {{2[\iota_7]_{10}}}\bigr),\,
	\begin{pmatrix}{{\eta_9}}\\{{\nu_7}}\end{pmatrix}\!,\,
	{\bar\varepsilon_{10}}\right\}{}.
	$$
	Take
	$\mathbf{y}\in\left\{
	\bigl({{\Delta(\iota_{10})}}, {{2[\iota_7]_{10}}}\bigr),\,
	\begin{pmatrix}{{\eta_9}}\\{{\nu_7}}\end{pmatrix}\!,\,
	{\bar\varepsilon_{10}}\right\}{}.$
	By \eqref{m2:nu9k}, we find
	\begin{align*}
		\notag {p_{10}}_*(\mathbf{y})
		&\in \left\{
		\bigl({{p_{10}(\Delta(\iota_{10}))}}, {{p_{10}([\iota_7]_{10})}}\bigr),\,
		\begin{pmatrix}{{\eta_9}}\\{{2\nu_7}}\end{pmatrix}\!,\,
		{\bar\varepsilon_{10}}\right\}{}\\
		\notag &\subseteq \left\{
		\bigl({{2\iota_9}}, {{0}}\bigr),\,
		\begin{pmatrix}{{\eta_9}}\\{{2\nu_7}}\end{pmatrix}\!,\,
		{\bar\varepsilon_{10}}\right\}{}\\
		\notag &\subseteq
		\{2\iota_9,\eta_9,\bar\varepsilon_{10}\}
		=\nu_9\kappa_{12}.
	\end{align*}
	Note that
	$2\mathbf{y} \in
	-\Delta(\iota_{10}) \circ
	\{\eta_9,\bar\varepsilon_{10},2\iota_{25}\}
	+2[\iota_7]_{10} \circ
	\{\nu_7,\bar\varepsilon_{10},2\iota_{25}\}.$
	The bracket $\{\nu_7,\bar\varepsilon_{10},2\iota_{25}\}$
	contains $0$: the relation
	$\bar\varepsilon_{10}=\eta_{10}\kappa_{11}$
	reduces the comparison to $\nu_7\circ\pi_{26}^{10}=0.$
	Its indeterminacy is
	$\nu_7\circ\pi_{26}^{10}+2\pi_{26}^7
	=\mathrm{span}\{2\bar\zeta_7\}.$ Therefore
	\[
	\begin{aligned}
		&\{\nu_7,\bar\varepsilon_{10},2\iota_{25}\}
		=\mathrm{span}\{2\bar\zeta_7\}.
	\end{aligned}
	\]
	\noindent
	In light of \cite[Lemma 15.1]{MT1},
	we have
	\[
	\begin{aligned}
		&\{\eta_9,\bar\varepsilon_{10},2\iota_{25}\}
		\subseteq\nu_9\kappa_{12}
		+\mathrm{span}\{\eta_9\mu_{10}\sigma_{19},P(\mu_{19})\}.
	\end{aligned}
	\]
	\noindent
	It follows that
	$$\Delta(\iota_{10}) \circ
	\{\eta_9,\bar\varepsilon_{10},2\iota_{25}\}
	\equiv \Delta(\nu_{10}\kappa_{13})
	\ \bmod\ \Delta(\eta_{10}\mu_{11}\sigma_{20}),
	\Delta(\iota_{10})P(\mu_{19}).$$ \noindent
	By \eqref{m2:Det10}, we obtain
	$\Delta(\eta_{10}\mu_{11}\sigma_{20})=0$
	and
	$\Delta(\sigma_{10}\eta_{17}\mu_{18})
	=\Delta(\eta_{10}\sigma_{11}\mu_{18})=0.$
	Hence we deduce
	$$
	2\mathbf{y}
	\equiv\Delta(\nu_{10}\kappa_{13})
	\ \bmod\ 4[\iota_7]_{10}\bar{\zeta}_7.
	$$
	\noindent The indeterminacy of the matrix Toda bracket containing
	$\mathbf{y}$ includes $2[\iota_7]_{10}\circ\pi_{26}^{7}$.
	Adding an integral multiple of $2[\iota_7]_{10}\bar{\zeta}_7$
	to $\mathbf{y}$ leaves it in the same bracket and does not change
	its image under $(p_{10})_*$.
	Its double changes by the corresponding multiple of
	$4[\iota_7]_{10}\bar{\zeta}_7$.
	The preceding congruence therefore allows us to choose
	\[
	\begin{aligned}
		[\nu_9\cn_{12}]
		&\in\mathbf{y}+\mathrm{span}\{2[\iota_7]_{10}\bar{\zeta}_7\}\\
		&\subseteq\left\{
		\bigl({{\Delta(\iota_{10})}}, {{2[\iota_7]_{10}}}\bigr),\,
		\begin{pmatrix}{{\eta_9}}\\{{\nu_7}}\end{pmatrix}\!,\,
		{\bar\varepsilon_{10}}\right\}
	\end{aligned}
	\]
	so that
	$$2[\nu_9\cn_{12}]=\Delta(\nu_{10}\cn_{13})
	=[\nu_6^2\cn_{12}]_{10}$$ in view of Lemma~\ref{m2:Dnuka}(2).
	In particular, $(p_{10})_*([\nu_9\cn_{12}])=\nu_9\cn_{12}$.
	We fix this choice throughout. For $n\geq11$,
	$[\nu_9\cn_{12}]_n=(i_{10,n})_*[\nu_9\cn_{12}]$.
	\noindent Let $R_1=J[\nu_9\cn_{12}].$
	By the relation $4\sigma_{10}\bar{\zeta}_{17}=0$
	\cite[I-(8.14)]{Oda1979}, we find
	$2R_1\equiv P(\nu_{21}\kappa_{24})
	\ \bmod 4\sigma_{10}\bar{\zeta}_{17}
	=P(\nu_{21}\kappa_{24}).$
	We have
	\[\FitMath{
		J[\nu_9\cn_{12}]
		\in J\left\{
		\bigl({{\Delta(\iota_{10})}}, {{2[\iota_7]_{10}}}\bigr),\,
		\begin{pmatrix}{{\eta_9}}\\{{\nu_7}}\end{pmatrix}\!,\,
		{\bar\varepsilon_{10}}\right\}{}
		\subseteq
		\left\{
		\bigl({{P(\iota_{21})}}, {{\sigma_{10}}}\bigr),\,
		\begin{pmatrix}{{\eta_{19}}}\\{{2\nu_{17}}}\end{pmatrix}\!,\,
		{\bar\varepsilon_{20}}\right\}_1
	}\]
	and
	\[\FitMath{
		H\left\{
		\bigl({{P(\iota_{21})}}, {{\sigma_{10}}}\bigr),\,
		\begin{pmatrix}{{\eta_{19}}}\\{{2\nu_{17}}}\end{pmatrix}\!,\,
		{\bar\varepsilon_{20}}\right\}_1
		\subseteq
		\left\{
		\bigl({{2\iota_{19}}}, {{0}}\bigr),\,
		\begin{pmatrix}{{\eta_{19}}}\\{{2\nu_{17}}}\end{pmatrix}\!,\,
		{\bar\varepsilon_{20}}\right\}_1
		=\{2\iota_{19},\eta_{19},\bar\varepsilon_{20}\}
		=\{\nu_{19}\kappa_{22}\}.
	}\]
	This defines $R_1\in\pi_{36}^{10}$  in
	\cite[I-Proposition 4.3(2)]{Oda1979}.

	We now combine the lift relations, boundary calculations, and Toda-bracket evaluations for $n=10.$
	The proposition records both the cyclic factors and the connecting images that relate their generators
	to the elements of homotopy groups of spheres.
	
	\begin{prop}\label{m2:10} The structure of $\pi_{26}(\mathit{SO}(10))$ is as follows.
		\[\GroupLayout{\pi_{26}(\mathit{SO}(10)) = \mathrm{span}\{[\sigma_9\eta_{16}\mu_{17}],\; [\eta_9\mu_{10}\sigma_{19}],\; [\nu_9\kappa_{12}],\; [\iota_7]_{10}\bar\sigma_7\}\cong \mathbb Z/64\oplus \mathbb Z/8\oplus \mathbb Z/4\oplus \mathbb Z/2.}{\begin{aligned}\pi_{26}(\mathit{SO}(10)) = \mathrm{span}\{&[\sigma_9\eta_{16}\mu_{17}],\; [\eta_9\mu_{10}\sigma_{19}],\; [\nu_9\kappa_{12}],\\
				&[\iota_7]_{10}\bar\sigma_7\}\cong \mathbb Z/64\oplus \mathbb Z/8\oplus \mathbb Z/4\oplus \mathbb Z/2.\end{aligned}}\]
		The three new lifts satisfy
		\[
		\begin{aligned}
			&2[\sigma_9\eta_{16}\mu_{17}]\equiv x[2P(\rho_{13})]_{10}\quad  (\mathrm{mod}\;[\iota_7]_{10}\bar\zeta_7+4[2P(\rho_{13})]_{10},\ [\nu_6^2\kappa_{12}]_{10},\ [\iota_7]_{10}\bar\sigma_7), \\
			&32[\sigma_9\eta_{16}\mu_{17}]=-4[\iota_7]_{10}\bar\zeta_7=16[2P(\rho_{13})]_{10},\qquad 2[\eta_9\mu_{10}\sigma_{19}]=[\iota_7]_{10}\bar\zeta_7+2[[2\e_{6},\e_{6}]]_{10}\Sigma^2\rho', \\
			&2[\nu_9\kappa_{12}]=\Delta(\nu_{10}\kappa_{13})=[\nu_6^2\kappa_{12}]_{10}.
		\end{aligned}
		\]
		Here $x$ is odd.  Moreover,
		\[\begin{aligned}
			&2[2P(\rho_{13})]_{10}\equiv[[2\e_{6},\e_{6}]]_{10}\Sigma^2\rho'\m(i_{6,10})_*\pi_{26}(\mathit{SO}(6));\quad \Delta\bar\mu_{10}=4[\eta_9\mu_{10}\sigma_{19}].
		\end{aligned}\]
	\end{prop}
	\begin{proof}
		Consider  the exact sequence
		\[\begin{aligned}
			&\pi^{9}_{27} \xrightarrow{ \Delta} \pi_{26}(\mathit{SO}(9)) \xrightarrow{i_\ast } \pi_{26}(\mathit{SO}(10)) \xrightarrow{p_\ast } \pi^{9}_{26} \xrightarrow{\Delta} \pi_{25}(\mathit{SO}(9)).
		\end{aligned}\]
		We observe that $\pi^9_{27} \cong \frac{\mathbb Z}{2}\{\bar\mu_9\eta_{26}\} \oplus
		\frac{\mathbb Z}{8}\{\sigma_9\zeta_{16}\}$ and $$\pi^9_{26} \cong \frac{\mathbb Z}{2}\{\sigma_9\eta_{16}\mu_{17}\} \oplus \frac{\mathbb Z}{2}\{\nu_9\kappa_{12}\} \oplus \frac{\mathbb Z}{2}\{\eta_9\mu_{10}\sigma_{19}\} \oplus \frac{\mathbb Z}{2}\{\bar\mu_9\},$$ where the right-hand $\Delta$ sends $\sigma_9\eta_{16}\mu_{17},$
		$\nu_9\kappa_{12}$ and $\eta_9\mu_{10}\sigma_{19}$ to zero in
		$\pi_{25}(\mathit{SO}(9)).$
		
		Note that
		$\varepsilon_7\zeta_{15}=0$ and $\bar\nu_7 \zeta_{15}=0.$ So, $$\Delta(\sigma_9 \zeta_{16}) = \Delta(\sigma_9) \zeta_{15} = [\iota_7]_9 (\bar\nu_7 + \varepsilon_7) \zeta_{15} = 0.$$ \noindent $\Delta(\bar\mu_9\eta_{26}) = \Delta(\iota_9)\bar\mu_8\eta_{25}=([ \nu_5]_9 + [\iota_7]_9 \eta_7)\bar\mu_8\eta_{25} = [ \nu_5]_9 \bar\mu_8 \eta_{25} + [\iota_7]_9 \eta_7 \bar\mu_8 \eta_{25}.$
		By Lemma~\ref{m2:2[et9ms}, we have
		$2[\eta_9\mu_{10}\sigma_{19}]=[\iota_7]_{10}\bar\zeta_7
		+2[[2\e_{6},\e_{6}]]_{10}\Sigma^2\rho'.$
		Since $2\rho_{13}=\Sigma^4 \rho',$   we have $2P(\rho_{13})=P(\iota_{13})\Sigma^2 \rho',$   and hence
		$$2[2P(\rho_{13})] \equiv [[2\e_{6},\e_{6}]]\Sigma^2\rho' \bmod i_{\ast}\pi_{26}(\mathit{SO}(6)).$$
		
		\noindent   From this, we conclude that $4[2P(\rho_{13})]_{10} \equiv 2[[2\e_{6},\e_{6}]]_{10}\Sigma^2\rho' \bmod 2[\nu_4\bar\zeta_7]_{10},$   where $2[\nu_4\bar\zeta_7]_{10}=4[\iota_7]_{10}\bar\zeta_7.$  Next consider
		$p_{10}\circ \left\{
		\bigl({{\Delta\iota_{10}}}, {{2[\iota_7]_{10}}}\bigr),\,\begin{pmatrix}{{\eta_9}}\\{{\nu_7}}\end{pmatrix}\!,\,{\Sigma\rho'}\right\}{}.$ 
		By Lemma~\ref{m2:ods}(1), this is a subset of
		\[\FitMath{
			\left\{
			\bigl({{2\iota_9}}, {{0}}\bigr),\,\begin{pmatrix}{{\eta_9}}\\{{2\nu_7}}\end{pmatrix}\!,\,{\Sigma\rho'}\right\}{}
			=\{2\iota_9,\eta_9,\Sigma\rho'\}=\sigma_9\eta_{16}\mu_{17}.
		}\]
		Thus, we can take a lift\vspace{-0.7em}
		\[\FitMath{
			[\sigma_9 \eta_{16} \mu_{17}] \in
			\left\{
			\bigl({{\Delta\iota_{10}}}, {{[\iota_7]_{10}}}\bigr),\,\begin{pmatrix}{{\eta_9}}\\{{2\nu_7}}\end{pmatrix}\!,\,{\Sigma\rho'}\right\}.
		}\]
		Then,
		\[\FitMath{
			16 [\sigma_9 \eta_{16} \mu_{17}] \in
			\left\{
			\bigl({{\Delta\iota_{10}}}, {{[\iota_7]_{10}}}\bigr),\,\begin{pmatrix}{{\eta_9}}\\{{2\nu_7}}\end{pmatrix}\!,\,{\Sigma\rho'}\right\}
			\circ 16 \iota_{26},
		}\]
		where the right-hand side is included in
		$-\Delta(\iota_{10}) \{\eta_9, \Sigma\rho', 16 \iota_{25}\} + [\iota_7]_{10}\{2\nu_7, \Sigma\rho', 16 \iota_{25}\}.$
		\noindent   
		It follows that \[\begin{aligned}
			&32[\sigma_9 \eta_{16} \mu_{17}] \in [\iota_7]_{10} \{2\nu_7, \Sigma\rho', 16\iota_{25} \} \circ 2\iota_{26} \supseteq [\iota_7]_{10} 2\iota_7 \{\nu_7, \Sigma\rho', 16\iota_{25}\} \circ 2\iota_{26} \ni 4[\iota_7]_{10}\bar\zeta_7.
		\end{aligned}\] Consequently, we get $32[\sigma_9 \eta_{16}\mu_{17}]=-4[\iota_7]_{10}\bar\zeta_7 = 16[2P(\rho_{13})]_{10}.$
		Thus, we deduce $$2[\sigma_9 \eta_{16} \mu_{17}] \equiv x[2P(\rho_{13})]_{10} \bmod \{[\iota_7]_{10}\bar\zeta_7+ 4[2P(\rho_{13})]_{10}, [\nu^2_6\kappa_{12}]_{10}, [\iota_7]_{10} \bar\sigma_7 \}.$$
		\indent From    \cite[pp.~343--344]{KM1}, 
		we have
		\[
		\Delta\mu_{10}
		= 2\bigl(
		[\iota_7]_{10}\zeta_7
		+ 2[P(2\iota_{13})]_{10}\sigma_{11}
		\bigr)
		= 4[\eta_9\varepsilon_{10}].
		\]
		\noindent   Then,  using \eqref{m2:etms},  we infer that
		\[
		\begin{aligned}
			\pm\Delta\bar{\mu}_{10}
			&\in \q\{\mu_{10},\,8\iota_{19},\,2\sigma_{19}\}_{1} \\ 
			&\xyd \{\Delta\mu_{10},\,8\iota_{18},\,2\sigma_{18}\} \\
			&= \{4[\eta_9\varepsilon_{10}],\,8\iota_{18},\,2\sigma_{18}\} \\
			&\supseteq 4\{[\eta_9\varepsilon_{10}],\,8\iota_{18},\,2\sigma_{18}\} \\
			&\ni 4[\eta_9\mu_{10}\cg_{19}]\m0.
		\end{aligned}
		\]\noindent     Hence $\Delta\bar{\mu}_{10}=4[\eta_9\mu_{10}\cg_{19}].$  \end{proof}
	
	In light of Lemma~\ref{m2:lem:mu10-sigma-bracket}, we may  choose the lift
	$[\mu_{10}\sigma_{19}]$ in the following Toda bracket:  
	$[\mu_{10}\sigma_{19}]\in\{[\varepsilon_{10}],8\iota_{18},2\sigma_{18}\}_1.$
	In addition, $\mathrm{ord}[\mu_{10}\sigma_{19}]=8.$ 
	Note that the   element  
	$[\varepsilon_{10}]$ has order $8$ by \cite[Theorem 2(iii)]{KM1}.
	
	At $n=11,$ the chosen lift of $\mu_{10}\sigma_{19}$ must be compared with the image of the
	$n=10$ generators. The following proposition gives the resulting decomposition together with the
	relation needed for the next inclusion.
	
	\begin{prop}\label{g1016:prop:26-11} The structure of $\pi_{26}(\mathit{SO}(11))$ is as follows.
		\[ \GroupLayout{\pi_{26}(\mathit{SO}(11)) = \mathrm{span}\{[\sigma_9 \eta_{16} \mu_{17}]_{11},\; [\mu_{10}\sigma_{19}],\; [\nu_9\kappa_{12}]_{11},\; [\iota_7]_{11} \bar\sigma_7\}\cong \mathbb Z/64\oplus \mathbb Z/8\oplus (\mathbb Z/2)^{2},}{\begin{aligned}\pi_{26}(\mathit{SO}(11)) = \mathrm{span}\{&[\sigma_9 \eta_{16} \mu_{17}]_{11},\; [\mu_{10}\sigma_{19}],\; [\nu_9\kappa_{12}]_{11},\\
				&[\iota_7]_{11} \bar\sigma_7\}\cong \mathbb Z/64\oplus \mathbb Z/8\oplus (\mathbb Z/2)^{2},\end{aligned}} \] where
		$2[\mu_{10}\sigma_{19}]\equiv x[\eta_9\mu_{10}\sigma_{19}]_{11} \ \bmod\
		16[\sigma_9\eta_{16}\mu_{17}]_{11}$ for some odd $x.$ 
	\end{prop}
	\begin{proof}
		Consider  the exact sequence
		\[\begin{aligned}
			&\pi^{10}_{27} \xrightarrow{ \Delta} \pi_{26}(\mathit{SO}(10)) \xrightarrow{i_\ast } \pi_{26}(\mathit{SO}(11)) \xrightarrow{p_\ast } \pi^{10}_{26}.
		\end{aligned}\]
		The fixed lift satisfies
		$2[\nu_9\kappa_{12}]_{11}
		=(i_{11})_*\Delta(\nu_{10}\kappa_{13})=0$.
		By Lemma~\ref{m2:Dnuka}(2), Proposition~\ref{m2:10},
		$\Delta(\eta_{10}\mu_{11}\sigma_{20})=0$,
		and the description of $\pi_{27}^{10}$ in Appendix~\ref{tlqdzj},
		\[
		\operatorname{Im}\bigl(\Delta\colon \pi_{27}^{10}
		\to\pi_{26}(\mathit{SO}(10))\bigr)
		=\mathrm{span}\{2[\nu_9\kappa_{12}],
		4[\eta_9\mu_{10}\sigma_{19}]\}.
		\]
		The ordered decomposition in Proposition~\ref{m2:10} shows that
		$[\nu_9\kappa_{12}]$ does not belong to this image.
		Hence $[\nu_9\kappa_{12}]_{11}\ne0$, and
		$\mathrm{ord}\,[\nu_9\kappa_{12}]_{11}=2$.
		The image of $\pi_{26}(\mathit{SO}(11))$ in $\pi^{10}_{26}$ is
		$\frac{\mathbb Z}{2}\{\mu_{10}\sigma_{19}\}.$
		
		By \cite[Theorem 2.(iii)]{KM1}, we find $$2[\varepsilon_{10}]\in
		{i_{11}}_*\pi_{18}(\mathit{SO}(10))=\mathrm{span}\{[P(\eta_{19})]_{11},
		[\eta_9\varepsilon_{10}]_{11}\}.$$
		Since $\Delta P(\iota_{21})=[P(\eta_{19})]+2[\eta_9\varepsilon_{10}]$
		\cite[p. 343]{KM1}, we conclude that $[P(\eta_{19})]_{11}=-2[\eta_9\varepsilon_{10}]_{11}.$
		This implies $2[\varepsilon_{10}]=x[\eta_9\varepsilon_{10}]_{11}$ for an integer $x.$ Applying $J$ and using
		$J[\eta_9\varepsilon_{10}]=y\lambda''$ for some odd $y$ (Lemma~\ref{m2:J[et9ep})
		and $J[\varepsilon_{10}]=\lambda'+2\xi'$
		\cite[Lemma 2.2(ii)]{HKM}, we have $2(\lambda'+2\xi')=xy(\Sigma\lambda'').$
		Recall that $2\lambda'=\Sigma\lambda''$ \cite[Lemma 12.19]{T} and $4\lambda'=4\xi'.$  This gives the relation $6\lambda'=2(xy)\lambda'.$
		This implies that $x$ is odd.
		
		Notice that $$2[\mu_{10}\sigma_{19}]\in\{2[\varepsilon_{10}],8\iota_{18},2\sigma_{18}\}_1
		=\{x[\eta_9\varepsilon_{10}]_{11},8\iota_{18},2\sigma_{18}\}_1\ni
		x[\eta_9\mu_{10}\sigma_{19}]_{11}$$\noindent $\mathrm{mod}\;2[\varepsilon_{10}]\circ \Sigma\pi^{17}_{25} +
		\pi_{19}(\mathit{SO}(11))\circ 2\sigma_{19}.$ Then  we have  $2[\varepsilon_{10}]\circ
		\Sigma\pi^{17}_{25}=0.$ By \cite[Theorem 2.1]{KM2}, we conclude that $\pi_{19}(\mathit{SO}(11))\circ
		2\sigma_{19}=\mathrm{span}\{2[32P(\iota_{21})]\sigma_{19}\}.$ This gives the intermediate congruence.
		The refinement to the stated modulus is verified at the beginning of the following proof.
	\end{proof}
	
	\begin{lem}\label{m2:eq:zeta-in-rank11}
		For some odd integer $y_3,$
		\[
		[\iota_7]_{11}\bar\zeta_7
		=8y_3[\sigma_9\eta_{16}\mu_{17}]_{11}
		+4[\mu_{10}\sigma_{19}].
		\]
	\end{lem}
	
	\begin{proof}
		We retain the lift $[\mu_{10}\sigma_{19}]$ chosen before
		Proposition~\ref{g1016:prop:26-11}.
		The bracket computation in the proof of that proposition gives
		\[
		2[\mu_{10}\sigma_{19}]
		=x[\eta_9\mu_{10}\sigma_{19}]_{11}
		+2y[32P(\iota_{21})]\sigma_{19}
		\]
		for integers $x,y,$ with $x$ odd.
		We first justify the stated refinement of its modulus using the image of
		$\pi_{26}(\mathit{SO}(10))$ in $\pi_{26}(\mathit{SO}(11)).$
		Clearly,  $16P(\sigma_{21})=0.$
		So,
		\[
		(p_{11})_*\bigl([32P(\iota_{21})]\sigma_{19}\bigr)
		=32P(\iota_{21})\sigma_{19}
		=32P(\sigma_{21})=0.
		\]
		Exactness and Proposition~\ref{m2:10} therefore give
		\[
		\begin{aligned}
			{}[32P(\iota_{21})]\sigma_{19}
			&\in(i_{10,11})_*\pi_{26}(\mathit{SO}(10))\\
			&=\mathrm{span}\{[\sigma_9\eta_{16}\mu_{17}]_{11},\;
			[\eta_9\mu_{10}\sigma_{19}]_{11},\;[\nu_9\kappa_{12}]_{11},\;
			[\iota_7]_{11}\bar\sigma_7\}.
		\end{aligned}
		\]
		The relation $\Delta\bar\mu_{10}=4[\eta_9\mu_{10}\sigma_{19}]$
		in Proposition~\ref{m2:10} implies
		$4[\eta_9\mu_{10}\sigma_{19}]_{11}=0.$
		The fixed choice gives $2[\nu_9\kappa_{12}]_{11}=0$.
		Also, $2[\iota_7]_{11}\bar\sigma_7
		=[\iota_7]_{11}(2\bar\sigma_7)=0$.
		It follows that, for some integers $x_1,y_1$ with $x_1$ odd,
		\[
		2[\mu_{10}\sigma_{19}]
		=x_1[\eta_9\mu_{10}\sigma_{19}]_{11}
		+2y_1[\sigma_9\eta_{16}\mu_{17}]_{11}.
		\]
		The coefficient $x_1$ is odd because it differs from $x$ by an even integer.
		Multiplication by $4,$ together with
		$8[\mu_{10}\sigma_{19}]=0,$ now gives
		\[
		8y_1[\sigma_9\eta_{16}\mu_{17}]_{11}=0.
		\]
		Since $[\sigma_9\eta_{16}\mu_{17}]_{11}$ has order $64,$
		we have $y_1=8y_2$ for some integer $y_2.$
		Thus
		\[
		\begin{aligned}
			2[\mu_{10}\sigma_{19}]
			&=x_1[\eta_9\mu_{10}\sigma_{19}]_{11}
			+16y_2[\sigma_9\eta_{16}\mu_{17}]_{11},\\
			4[\mu_{10}\sigma_{19}]
			&=2[\eta_9\mu_{10}\sigma_{19}]_{11}
			+32y_2[\sigma_9\eta_{16}\mu_{17}]_{11}.
		\end{aligned}
		\]
		The first equality proves the refinement
		\[
		2[\mu_{10}\sigma_{19}]
		\equiv x_1[\eta_9\mu_{10}\sigma_{19}]_{11}
		\pmod{16[\sigma_9\eta_{16}\mu_{17}]_{11}},
		\qquad x_1\ \text{odd}.
		\]
		Here the modulus is generated by the image under $(i_{10,11})_*$
		of $16[\sigma_9\eta_{16}\mu_{17}].$
		The second equality gives the intermediate comparison
		\[
		2[\eta_9\mu_{10}\sigma_{19}]_{11}
		\equiv4[\mu_{10}\sigma_{19}]
		\pmod{32[\sigma_9\eta_{16}\mu_{17}]_{11}}.
		\]
		
		The relations in Proposition~\ref{m2:10} also give
		\[
		\begin{aligned}
			{}[\iota_7]_{10}\bar\zeta_7
			&=2[\eta_9\mu_{10}\sigma_{19}]
			-4(2x_2+1)[P(2\rho_{13})]_{10},\\
			4[P(2\rho_{13})]_{10}
			&=8x_3[\sigma_9\eta_{16}\mu_{17}]
		\end{aligned}
		\]
		for integers $x_2,x_3,$ with $x_3$ odd.
		Applying $(i_{10,11})_*$ and using the exact expression for
		$4[\mu_{10}\sigma_{19}]$ obtained above, we get
		\[
		\begin{aligned}
			{}[\iota_7]_{11}\bar\zeta_7
			&=2[\eta_9\mu_{10}\sigma_{19}]_{11}
			-8(2x_2+1)x_3[\sigma_9\eta_{16}\mu_{17}]_{11}\\
			&=4[\mu_{10}\sigma_{19}]
			+8\bigl(-(2x_2+1)x_3-4y_2\bigr)
			[\sigma_9\eta_{16}\mu_{17}]_{11}.
		\end{aligned}
		\]
		The integer $y_3=-(2x_2+1)x_3-4y_2$ is odd.
		This proves the asserted equality after absorbing the error into
		its odd coefficient.
	\end{proof}
	
	\subsection{The Cases $n=12,13$}
	
	The next two cases require more detailed extension calculations.  We first
	identify the multiple of the new lift at $n=12,$   then choose generators
	that exhibit the cyclic factors.  At $n=13,$   we quotient by the boundary
	image before adding the new contribution from the homotopy groups of spheres.

	\begin{lem}\label{1023dxb}
		The lift $[\nu_{9}\cn_{12}]$ fixed after Lemma~\ref{m2:Dnuka} satisfies
		\[
		[\nu_{9}\cn_{12}]\in\left\{
		(\Delta\e_{10},2[\e_{7}]_{10}),\,
		\begin{pmatrix}\ca_{9}\\\nu_{7}\end{pmatrix},\,
		\ca_{10}\cn_{11}
		\right\}\xyd\pi_{26}(\mathit{SO}(10)).
		\]
		Moreover, $2[\nu_{9}\cn_{12}]=[\nu_{6}^{2}\cn_{12}]_{10}$
		and $\mathrm{ord}\,[\nu_{9}\cn_{12}]=4$.
	\end{lem}
	\begin{proof}
		Since $\bar\varepsilon_{10}=\ca_{10}\cn_{11}$,
		this is the matrix Toda bracket used to choose $[\nu_{9}\cn_{12}]$.
		The fixed lift satisfies
		$2[\nu_{9}\cn_{12}]=\Delta(\nu_{10}\cn_{13})$.
		By Lemma~\ref{m2:Dnuka}(2),
		$\Delta(\nu_{10}\cn_{13})=[\nu_{6}^{2}\cn_{12}]_{10}$.
		The order of $[\nu_{9}\cn_{12}]$ is $4$ by Proposition~\ref{m2:10}.
	\end{proof}

	To examine the group extension problem\footnote{Seriously speaking, it is a $\z_{(2)}$-module  extension problem.} for $\pi_{26}(\mathit{SO}(12)),$   let us do an exercise from homological algebra.
	\begin{lem}\label{xlx1023}
		Consider a short exact sequence of $\z_{(2)}$-modules
		\[
		0\longrightarrow\z/64\oplus\z/8\xrightarrow{j}B\xrightarrow{q}\z/16\longrightarrow0.
		\]
		Suppose that there is an element $b\in B$ such that
		$$16b=j(8,0)+j(0,4),\;\;q(b)=1.$$
		\noindent Then $$B=\mathrm{span}\{b,j(1,0)-2b,j(2,1)-4b\}\cong\z/128\oplus\z/16\oplus\z/4,$$
		where the three listed generators correspond to the cyclic summands in the indicated order.
	\end{lem}
	
	\begin{proof}   We use the description of $\z_{(2)}$-module extensions in terms of
		pushouts and the functor $\mathrm{Ext}_{1}(-,-).$  This is standard material in
		homological algebra textbooks; see, for example,
		\cite[Lemma~7.28 and Theorem~7.30]{Rotman}.
		The sequence
		\[
		\cdots\longrightarrow0\longrightarrow\z_{(2)}\xrightarrow{\times16}\z_{(2)}\xrightarrow{\mathrm{proj.}}\z/16\longrightarrow0
		\]
		is a projective resolution of $\z/16.$ 
		Define  the homomorphisms 
		$$f\colon\z_{(2)}\to\z/64\oplus\z/8,\quad 1\mapsto (8,4),  $$
		$$g\colon \z_{(2)}\to B,\quad 1\mapsto b.$$
		\noindent The hypotheses give the following commutative diagram with exact rows:
		\[
		\begin{tikzcd}[column sep=large,row sep=large]
			0\arrow[r]
			&\z_{(2)}\arrow[r,"\times16"]\arrow[d,"f"']
			&\z_{(2)}\arrow[r,"\mathrm{proj.}"]\arrow[d,"g"']
			&\z/16\arrow[r]\arrow[d,equal]
			&0\\
			0\arrow[r]
			&\z/64\oplus\z/8\arrow[r,"j"']
			&B\arrow[r,"q"']
			&\z/16\arrow[r]
			&0.
		\end{tikzcd}
		\]
		\noindent   Let $$\operatorname{cls}(x,y,t)=(x,y,t)+\mathrm{span}\{(8,4,-16)\nn\in  \bigl(\z/64\oplus\z/8\oplus\z_{(2)}\bigr)/\mathrm{span}\{(8,4,-16)\}.$$ \noindent   The left square is a pushout. Consequently, we have
		\[B\cong
		\bigl(\z/64\oplus\z/8\oplus\z_{(2)}\bigr)/\mathrm{span}\{(8,4,-16) \nn\cong\z/128\oplus\z/16\oplus\z/4,
		\]
		\noindent generated by $\mathrm{cls}(0,0,1), \mathrm{cls}(1,0,-2)$ and $\mathrm{cls}(2,1,-4),$   respectively.  The defining relation gives $$128\operatorname{cls}(0,0,1)=0,\;
		16\operatorname{cls}(1,0,-2)=0 \text{ and } 4\operatorname{cls}(2,1,-4)=0.$$
		Hence the result holds. \end{proof}
	
	Yoshihiro Hirato (the second author) first computed the groups $\pi_{26}(\mathit{SO}(n))$ for all $n\geq 3$ in 2022, and Juxin Yang independently recomputed them in 2023 under the guidance of Juno Mukai.  The proof presented here in the case $n=12$ is due to Yang and improves upon Hirato's result.  Hirato's result for $\pi_{26}(\mathit{SO}(12))$ requires two odd parameters $\af,\beta$, while the displayed result involves only a single odd parameter $\ell_0$.
	
	\begin{prop}\label{m2:lem:rank12-extension}
		For an odd integer $\ell_{0},$  
		\[\GroupLayout{\pi_{26}(\mathit{SO}(12)) = {} \mathrm{span}\{[\Sigma^2\rho'],\; \ell_{0}[\sigma_9\eta_{16}\mu_{17}]_{12}-2[\Sigma^2\rho'],\; 2\ell_{0}[\sigma_9\eta_{16}\mu_{17}]_{12} +[\mu_{10}\sigma_{19}]_{12}-4[\Sigma^2\rho'],\; [\nu_9\kappa_{12}]_{12},\; [\iota_7]_{12}\bar\sigma_7,\; [\eta_{11}]\kappa_{12}\}\cong \mathbb Z/128\oplus \mathbb Z/16\oplus \mathbb Z/4\oplus (\mathbb Z/2)^{3}.}{\begin{aligned}\pi_{26}(\mathit{SO}(12)) \;= {} \mathrm{span}\{&[\Sigma^2\rho'],\; \ell_{0}[\sigma_9\eta_{16}\mu_{17}]_{12}-2[\Sigma^2\rho'],\; 2\ell_{0}[\sigma_9\eta_{16}\mu_{17}]_{12} +[\mu_{10}\sigma_{19}]_{12}-4[\Sigma^2\rho'],\\
				&[\nu_9\kappa_{12}]_{12},\; [\iota_7]_{12}\bar\sigma_7,\; [\eta_{11}]\kappa_{12}\}\cong \mathbb Z/128\oplus \mathbb Z/16\oplus \mathbb Z/4\oplus (\mathbb Z/2)^{3}.\end{aligned}}\]
		The group extension relations are
		\[
		\begin{aligned}
			&16[\Sigma^2\rho']=8\ell_{0}[\sigma_9\eta_{16}\mu_{17}]_{12}+4[\mu_{10}\sigma_{19}]_{12},\qquad 2[\eta_{11}]\kappa_{12}=0.
		\end{aligned}
		\]
		The lifts may be chosen compatibly so that
		\begingroup
		\renewcommand{\theHequation}{m2:eq:rank12-kappa-boundary}
		\begin{equation}
			\begin{aligned}
				&\Delta\bar\varepsilon_{12} =[\eta_{10}^2]_{12}\kappa_{12} \ty[\nu_9\kappa_{12}]_{12}\m 64[\Sigma^2\rho'].\label{m2:eq:rank12-kappa-boundary}
			\end{aligned}
		\end{equation}
		\endgroup
	\end{prop}
	
	\begin{proof}
		We separate the exact-sequence calculation, the extension relation,
		and the change of generators.
		
		\smallskip\noindent\emph{The Exact Sequence.}
		
		Consider
		\[\begin{aligned}
			&\pi_{27}^{11}\xrightarrow{\Delta}\pi_{26}(\mathit{SO}(11)) \xrightarrow{(i_{12})_*}\pi_{26}(\mathit{SO}(12)) \xrightarrow{(p_{12})_*}\pi_{26}^{11} \xrightarrow{\Delta}\pi_{25}(\mathit{SO}(11)).
		\end{aligned}\]
		The homotopy groups of spheres are
		\[\begin{aligned}
			&\pi_{27}^{11}={}\frac{\mathbb Z}{2}\{\mu_{11}\sigma_{20}\},\qquad \pi_{26}^{11}=\frac{\mathbb Z}{16}\{\Sigma^2\rho'\} \oplus\frac{\mathbb Z}{2}\{\eta_{11}\kappa_{12}\}.
		\end{aligned}\]
		Since $\Delta\iota_{11}=[\iota_7]_{11}\nu_7,$   the first boundary is
		$\Delta(\mu_{11}\sigma_{20})
		=[\iota_7]_{11}\nu_7\mu_{10}\sigma_{19}=0.$
		Notice  that $[\iota_7]\nu_7\mu_{10}\in \pi_{19}(\mathit{SO}(8)),$   and
		$\pi_{19}(\mathit{SO}(8))=\frac{\mathbb Z}{2}\{[ \nu_5]_8\bar\nu_8\nu_{16}\}$
		by \cite[Theorem~2.1]{KM2}\uprightsemicolon{} its composite with $\sigma_{19}$ is
		zero because $\nu_{16}\sigma_{19}=0.$ 
		Thus $(i_{12})_*$ is injective in degree $26.$ 
		For the image of the fixed lift, we obtain
		\[
		2[\nu_9\kappa_{12}]_{12}
		=(i_{12})_*\bigl(2[\nu_9\kappa_{12}]_{11}\bigr)=0,
		\qquad [\nu_9\kappa_{12}]_{12}\ne0.
		\]
		Thus $\mathrm{ord}\,[\nu_9\kappa_{12}]_{12}=2$.
		
		The other boundary is also zero.  The relation $\nu_6\rho'=0$
		gives
		$\Delta(\Sigma^2\rho')=[\iota_7]_{11}\nu_7\Sigma\rho'=0.$
		In addition, $\Delta\eta_{11}=0.$   So,  $[\eta_{11}]\kappa_{12}$
		lifts $\eta_{11}\kappa_{12}.$   This lift has order exactly $2$:
		it is killed by $2\kappa_{12}=0,$   and its projection is nonzero.
		Consequently there is a short exact sequence
		\begin{equation}
			\begin{aligned}
				&0\longrightarrow \pi_{26}(\mathit{SO}(11)) \xrightarrow{(i_{12})_*}\pi_{26}(\mathit{SO}(12)) \xrightarrow{(p_{12})_*} \frac{\mathbb Z}{16}\{\Sigma^2\rho'\}\oplus\frac{\mathbb Z}{2}\{\eta_{11}\kappa_{12}\} \longrightarrow0.\label{m2:eq:rank12-short-exact}
			\end{aligned}
		\end{equation}
		The injected group has the four cyclic factors of orders
		$64,8,2,2$ in the preceding $n=11$ decomposition.
		
		\smallskip\noindent\emph{The Extension Problem.}
		
		Since $(i_{12})_*\Delta\iota_{11}=0$ and $\nu_6\rho'=0,$   the
		bracket $\{[\iota_7]_{12},\nu_7,\Sigma\rho'\}_1$ is defined.
		Its projection consists of odd multiples of $\Sigma^2\rho'.$ 
		Note that
		$\{p_{12},i_{12},\Delta\iota_{11}\}\ni\e_{11}$ and  its composite with
		$\Sigma^2\rho'$ has indeterminacy $p_{12*}\pi_{11}(\mathit{SO}(12))\hc\s^{2}\co'=\mathrm{span}\{2\Sigma^2\rho'\}.$ 
		Choose an odd $t$ and an element
		$[\Sigma^2\rho']\in\{[\iota_7]_{12},\nu_7,t\Sigma\rho'\}_1$
		after normalizing the projection to be $\Sigma^2\rho'.$ 
		The bracket composition rule gives
		$16[\Sigma^2\rho']
		\in-[\iota_7]_{12}\circ\{\nu_7,t\Sigma\rho',16\iota_{25}\}_1.$
		The latter Toda bracket in the homotopy groups of spheres gives an odd multiple of $\bar\zeta_7.$ 
		This follows by stabilizing, using the stable bracket relations of
		\cite[p.~58, Proposition~3.3(2)--(3)]{Oda1979}, and using the injectivity
		of stabilization on $\pi_{26}^7\uprightsemicolon{}$ the vanishing $\nu_6\rho'=0$
		ensures that the indicated unstable bracket is defined.
		Together with Lemma~\ref{m2:eq:zeta-in-rank11}, this yields
		\begin{equation}
			\begin{aligned}
				&16[\Sigma^2\rho'] =8\ell_{0}[\sigma_9\eta_{16}\mu_{17}]_{12} +4[\mu_{10}\sigma_{19}]_{12},\qquad \ell_{0}\text{:  odd}.\label{m2:eq:rank12-sixteen}
			\end{aligned}
		\end{equation}
		\noindent In light of the method of Lemma \ref{xlx1023}, we see that $$\pi_{26}(\mathit{SO}(12))\cong \z/128\jia\z/16\jia\z/4\jia(\z/2)^{3}$$ generated by the corresponding generators displayed in the proposition.

		\smallskip\noindent\emph{The Boundary.}
		
		Our composition formula and the boundary table give
		$$\Delta\bar\varepsilon_{12}
		=\Delta(\eta_{12}\kappa_{13})
		=[\eta_{10}^2]_{12}\kappa_{12}\in (i_{11,12})_*\bigl(\{i_{10,11},\Delta\iota_{10},\eta_9^2\}\circ\kappa_{12}\bigr).$$ 
		\noindent
		We introduce the following notation:
		\[
		\begin{aligned}
			\mathbf{a}_{7}
			&= [\s^{2}\co'],\\
			\mathbf{a}_{4}
			&= \ell_{0}[\cg_{9}\ca_{16}\mu_{17}]_{12}
			-2[\s^{2}\co'],\\
			\mathbf{a}_{2}
			&= 2\ell_{0}[\cg_{9}\ca_{16}\mu_{17}]_{12}
			+[\mu_{10}\cg_{19}]_{12}
			-4[\s^{2}\co'],\\
			\mathbf{a}_{1}
			&= [\nu_{9}\cn_{12}]_{12},\\
			\mathbf{b}_{1}
			&= [\e_{7}]_{12}\ccs_{7},\\
			\mathbf{c}_{1}
			&= [\ca_{11}]\cn_{12}.
		\end{aligned}
		\]
		Here, the subscripts indicate the exponents of $2$
		in the orders of the corresponding elements.
		For instance, $\mathrm{ord}(\mathbf{a}_{7})=2^{7}.$ 
		Observe that
		$\Delta(\eta_{12}\kappa_{13})
		=\Delta(\eta_{12})\circ\kappa_{12}
		=[\eta_{10}^{2}]_{12}\kappa_{12}.$ 
		In view of
		$\pi_{12}(\mathit{SO}(11))
		=\mathrm{span}\{[\eta_{10}^{2}]\}\cong\z/2$
		from \cite{Ka}, we have
		$[\eta_{10}^{2}]
		\in\{i_{10,11},\Delta\iota_{10},\eta_{9}^{2}\}.$ 
		Hence,
		\begin{align}
			[\eta_{10}^{2}]_{12}\kappa_{12}
			&\in i_{11,12}\circ
			\{i_{10,11},\Delta\iota_{10},\eta_{9}^{2}\}
			\circ\kappa_{12}
			\notag\\
			&\subseteq
			\{i_{10,12},\Delta\iota_{10},\eta_{9}^{2}\}
			\circ\kappa_{12}.
			\notag
		\end{align}
		\noindent
		Recall that
		$\pi_{12}(\mathit{SO}(10))
		=\mathrm{span}\{\Delta\nu_{10}\}$
		from \cite{Ka}. In light of Lemma \ref{m2:eq:zeta-in-rank11}, we deduce that
		$$2[\iota_{7}]_{12}\bar{\zeta}_{7}
		=16z[\sigma_{9}\eta_{16}\mu_{17}]_{12}
		=16t(\mathbf{a}_{4}+2\mathbf{a}_{7})
		=32t\mathbf{a}_{7},$$
		where $t$ and $z$ are some odd integers.  So, 
		\begin{align}
			[\nu_{9}\kappa_{12}]_{12}
			&\in i_{10,12}\circ
			\left\{
			\bigl(\Delta\iota_{10},\,2[\iota_{7}]_{10}\bigr),\,
			\begin{pmatrix}
				\eta_{9}\\
				\nu_{7}
			\end{pmatrix},\,
			\eta_{10}\kappa_{11}
			\right\}
			\notag\\
			&\subseteq i_{10,12}\circ
			\left\{
			\bigl(\Delta\iota_{10},\,2[\iota_{7}]_{10}\bigr),\,
			\begin{pmatrix}
				\eta_{9}^{2}\\
				0
			\end{pmatrix},\,
			\kappa_{11}
			\right\}
			\notag\\
			&\xyd i_{10,12}\circ
			\{\Delta\iota_{10},\eta_{9}^{2},\kappa_{11}\}
			+i_{10,12}\circ
			\{2[\iota_{7}]_{10},0,\kappa_{11}\}
			\notag\\
			&\supseteq [\eta_{10}^{2}]_{12}\kappa_{12}+0
			\notag\\
			&=[\eta_{10}^{2}]_{12}\kappa_{12};
			\notag
		\end{align}
		here,
		$$i_{10,12}\circ
		\Bigl(
		(\Delta\iota_{10})\circ\pi_{26}^{9}
		+\pi_{12}(\mathit{SO}(10))\circ\kappa_{12}
		+2[\iota_{7}]_{10}\circ\pi_{26}^{7}
		\Bigr)=\mathrm{span}\{2[\iota_{7}]_{12}\bar{\zeta}_{7}\}
		=\mathrm{span}\{2^{5}\mathbf{a}_{7}\}.$$ 
		
		\noindent   Consequently,	$$[\eta_{10}^{2}]_{12}\kappa_{12}
		\equiv[\nu_{9}\kappa_{12}]_{12}
		=\mathbf{a}_{1}\m2^{5}\mathbf{a}_{7}.$$
		Since $2[\eta_{10}^{2}]_{12}\kappa_{12}=2\mathbf{a}_{1}=0$
		and $\mathrm{ord}(\mathbf{a}_{7})=2^{7}$,
		the coefficient of $2^{5}\mathbf{a}_{7}$ in
		$[\eta_{10}^{2}]_{12}\kappa_{12}-\mathbf{a}_{1}$ is even.
		Therefore,$$	[\eta_{10}^{2}]_{12}\kappa_{12}
		\equiv[\nu_{9}\kappa_{12}]_{12}
		=\mathbf{a}_{1}\m2^{6}\mathbf{a}_{7}.$$
		\noindent       Then the result holds.
	\end{proof}

	Proposition~\ref{m1:25-12} and Lemma~\ref{m1:1112} give
	\[\FitMath{\begin{aligned}
			\pi_{26}^{12}&=\frac{\mathbb Z}{16}\{\sigma_{12}^2\}\oplus\frac{\mathbb Z}{4}\{P(\nu_{25})\}\oplus\frac{\mathbb Z}{2}\{\kappa_{12}\},\\
			\Delta(\sigma_{12}^2)&\equiv2[\sigma_{11}^2]+u[\bar\sigma_6]_{12}\pmod{\mathrm{span}\{\Delta(P(\nu_{25}))\}},\quad u\text{ odd},\\
			\Delta(P(\nu_{25}))&=[[\e_{10},\nu_{10}]]_{12}\nu_{22},\qquad
			\Delta(\kappa_{12})=[\bar\varepsilon_{10}]_{12}.
	\end{aligned}}\]
	The direct sum in Proposition~\ref{m1:25-12} therefore yields
	\begin{equation}\label{m2:eq:rank13-required-kernel}
		\FitMath{\ker\bigl(\Delta\colon \pi_{26}^{12}\to\pi_{25}(\mathit{SO}(12))\bigr)=\frac{\mathbb Z}{2}\{8\sigma_{12}^2\}.}
	\end{equation}\noindent     
	\begin{lem}\label{m2:lem:rank13-inherited-image}
		Let $\ell_{0}$ be the odd integer in
		Proposition~\ref{m2:lem:rank12-extension}, and write
		$$([\mu_{10}\sigma_{19}]_{13})'
		=2\ell_{0}[\sigma_9\eta_{16}\mu_{17}]_{13}
		+[\mu_{10}\sigma_{19}]_{13}
		-4[\Sigma^2\rho']_{13}.$$
		Then
		\vspace{-1em}   \begin{align}
			(i_{13})_*\pi_{26}(\mathit{SO}(12))={}
			&\frac{\mathbb Z}{128}\{[\Sigma^2\rho']_{13}\}
			\oplus
			\frac{\mathbb Z}{4}\{([\mu_{10}\sigma_{19}]_{13})'\}
			\notag\\
			&\mathord{\oplus}\,
			\frac{\mathbb Z}{2}\{[\iota_7]_{13}\bar\sigma_7\}
			\oplus
			\frac{\mathbb Z}{2}\{[\eta_{11}]_{13}\kappa_{12}\}.
			\label{m2:eq:rank13-inherited}
		\end{align}
	\end{lem}
	\begin{proof}
		We use the six generators
		$\mathbf{a}_{7},\mathbf{a}_{4},\mathbf{a}_{2},
		\mathbf{a}_{1},\mathbf{b}_{1},\mathbf{c}_{1}$
		defined in the proof of
		Proposition~\ref{m2:lem:rank12-extension}.
		Their orders are $128,16,4,2,2,2,$   respectively.
		
		Consider the exact sequence
		$\pi_{27}^{12}\xrightarrow{\Delta}
		\pi_{26}(\mathit{SO}(12))
		\xrightarrow{(i_{13})_*}\pi_{26}(\mathit{SO}(13)),$  
		where
		$$\pi_{27}^{12}
		=\frac{\mathbb Z}{16}\{\Sigma^3\rho'\}
		\oplus\frac{\mathbb Z}{2}\{\eta_{12}\kappa_{13}\}.$$
		By \eqref{m2:eq:rank12-kappa-boundary},
		$\Delta(\eta_{12}\kappa_{13})
		=\mathbf{a}_{1}+64t\mathbf{a}_{7}$
		for some $t\in\{0,1\}.$ 
		Hence
		$\operatorname{ord}(\Delta(\eta_{12}\kappa_{13}))=2.$ 
		
		Note that
		$(p_{12})_*\Delta(\Sigma^3\rho')=2\Sigma^2\rho'.$ 
		Moreover,
		$J\Delta(\Sigma^3\rho')=\pm2P(\rho_{25})$
		has order $16,$   by \cite[p.~75, (8.20)]{Oda1979}
		and the direct sum decomposition of $\pi_{38}^{12}.$ 
		Together with $16\Delta(\Sigma^3\rho')=0,$   this gives
		$\operatorname{ord}(\Delta(\Sigma^3\rho'))=16.$ 
		Since $\operatorname{ord}(\mathbf{a}_{7})=128,$  
		the equality $16\Delta(\Sigma^3\rho')=0$ implies
		$$\Delta(\Sigma^3\rho')
		=8x_{7}\mathbf{a}_{7}+x_{4}\mathbf{a}_{4}
		+x_{2}\mathbf{a}_{2}+x_{1}\mathbf{a}_{1}
		+y_{1}\mathbf{b}_{1}+z_{1}\mathbf{c}_{1}$$
		for some integers $x_{7},x_{4},x_{2},x_{1},y_{1},z_{1}\in\z.$ 
		By definition,
		$(p_{12})_*\mathbf{a}_{7}=\Sigma^2\rho',$  
		$(p_{12})_*\mathbf{a}_{4}=-2\Sigma^2\rho',$  
		$(p_{12})_*\mathbf{a}_{2}=-4\Sigma^2\rho',$  
		$(p_{12})_*\mathbf{a}_{1}
		=(p_{12})_*\mathbf{b}_{1}=0,$  
		and
		$(p_{12})_*\mathbf{c}_{1}=\eta_{11}\kappa_{12}.$ 
		Therefore
		$$2\Sigma^2\rho'
		=(8x_{7}-2x_{4}-4x_{2})\Sigma^2\rho'
		+z_{1}\eta_{11}\kappa_{12}.$$
		The decomposition
		$\pi_{26}^{11}
		=\frac{\mathbb Z}{16}\{\Sigma^2\rho'\}
		\oplus\frac{\mathbb Z}{2}\{\eta_{11}\kappa_{12}\}$
		gives
		$8x_{7}-2x_{4}-4x_{2}\equiv2\pmod{16}$
		and $z_{1}\equiv0\pmod2.$ 
		Thus $x_{4}\equiv1\pmod2,$   $z_{1}\mathbf{c}_{1}=0,$   and
		$\Delta(\Sigma^3\rho')
		=8x_{7}\mathbf{a}_{7}+x_{4}\mathbf{a}_{4}
		+x_{2}\mathbf{a}_{2}+x_{1}\mathbf{a}_{1}
		+y_{1}\mathbf{b}_{1}.$ 
		
		\noindent Consider
		$\mathbf{a}_{1}
		=\Delta(\eta_{12}\kappa_{13})-64t\mathbf{a}_{7}$
		and
		$$x_{4}\mathbf{a}_{4}
		=\Delta(\Sigma^3\rho')
		-x_{1}\Delta(\eta_{12}\kappa_{13})
		-8(x_{7}-8tx_{1})\mathbf{a}_{7}
		-x_{2}\mathbf{a}_{2}-y_{1}\mathbf{b}_{1}.$$
		Since $x_{4}\in\mathbb Z$ is odd, we see that
		$\pi_{26}(\mathit{SO}(12))
		=\mathrm{span}\{
		\mathbf{a}_{7},\mathbf{a}_{2},
		\mathbf{b}_{1},\mathbf{c}_{1}\}
		+\Delta\pi_{27}^{12}.$ 
		
		For $m,n\in\mathbb Z,$   suppose that
		$m\Delta(\Sigma^3\rho')
		+n\Delta(\eta_{12}\kappa_{13})
		\in\mathrm{span}\{
		\mathbf{a}_{7},\mathbf{a}_{2},
		\mathbf{b}_{1},\mathbf{c}_{1}\}.$ 
		Comparing the $\mathbf{a}_{4}$- and
		$\mathbf{a}_{1}$-coefficients gives
		$mx_{4}\equiv0\pmod{16}$
		and $mx_{1}+n\equiv0\pmod2.$ 
		Since $x_{4}$ is odd, $16\mid m$ and $2\mid n.$ 
		Consequently,
		$m\Delta(\Sigma^3\rho')
		+n\Delta(\eta_{12}\kappa_{13})=0.$ 
		Hence
		$\mathrm{span}\{
		\mathbf{a}_{7},\mathbf{a}_{2},
		\mathbf{b}_{1},\mathbf{c}_{1}\}
		\cap\Delta\pi_{27}^{12}=0,$   and therefore
		$\pi_{26}(\mathit{SO}(12))
		=\mathrm{span}\{
		\mathbf{a}_{7},\mathbf{a}_{2},
		\mathbf{b}_{1},\mathbf{c}_{1}\}
		\oplus\Delta\pi_{27}^{12}.$ 
		
		By exactness, the restriction
		$(i_{13})_*\colon
		\mathrm{span}\{
		\mathbf{a}_{7},\mathbf{a}_{2},
		\mathbf{b}_{1},\mathbf{c}_{1}\}
		\longrightarrow
		(i_{13})_*\pi_{26}(\mathit{SO}(12))$
		is an isomorphism.
		Since
		$(i_{13})_*\mathbf{a}_{7}=[\Sigma^2\rho']_{13},$  
		$(i_{13})_*\mathbf{a}_{2}=([\mu_{10}\sigma_{19}]_{13})',$  
		$(i_{13})_*\mathbf{b}_{1}=[\iota_7]_{13}\bar\sigma_7,$  
		and
		$(i_{13})_*\mathbf{c}_{1}=[\eta_{11}]_{13}\kappa_{12},$  
		we obtain \eqref{m2:eq:rank13-inherited}.
	\end{proof}

	\begin{lem}\label{m2:lem:rank13-toda-order}
		For the order-$8$ lift $[\nu_{12}^2]\in\pi_{18}(\mathit{SO}(13))$ of \cite[Theorem~2]{KM1}, the Toda bracket $\{[\nu_{12}^2],8\iota_{18},2\sigma_{18}\}_1$ is well-defined.
		Every element of this bracket has order $8,$   and the $J$-image of its fourth multiple is $4\Sigma\tau^{\mathrm{IV}}\ne0.$ 
		Moreover,
		\begin{equation}
			\{4\lambda,8\iota_{31},2\sigma_{31}\}=4\Sigma\tau^{\mathrm{IV}}\in\pi_{39}^{13}.
			\label{m2:eq:rank13-four-lambda}
		\end{equation}
	\end{lem}
	
	\begin{proof}
		The order of $[\nu_{12}^2]$ and the relation $8\iota_{18}\circ2\sigma_{18}=0$ show that the bracket is defined.
		Its elements are killed by $8,$   by the bracket composition rule and the stable relation $\langle8\iota,2\sigma,8\iota\rangle=0,$   which follows from \cite[(3.10)]{T}.
		For the $J$-image, we use $J[\nu_{12}^2]=\lambda+2\xi_{13}$ from \cite[Lemma~2.2(ii), p.~7]{HKM}.
		
		To establish \eqref{m2:eq:rank13-four-lambda}, start with $H\lambda''=\eta_{19}\varepsilon_{20}.$ The Hopf calculation in the proof of the lemma following \eqref{m2:etms} gives $H\{\lambda'',8\iota_{28},2\sigma_{28}\}_1=\{\varepsilon_{19}\mu_{27}\}=\{\eta_{19}\sigma_{20}\mu_{27}\}$, with zero indeterminacy in the Hopf image.
		By \cite[I-Proposition~4.2(1), p.~63]{Oda1979}, $H\tau''=\eta_{19}\sigma_{20}\mu_{27}.$ Exactness of the EHP sequence therefore gives $$\{\lambda'',8\iota_{28},2\sigma_{28}\}_1\subseteq\tau''+\Sigma\pi_{35}^{9}.$$
		After three suspensions, the bracket maps into $$\{\Sigma^3\lambda'',8\iota_{31},2\sigma_{31}\}_4\subseteq\{4\lambda,8\iota_{31},2\sigma_{31}\},$$ where $\Sigma^3\lambda''=4\lambda$ and $\Sigma^3\tau''=4\Sigma\tau^{\mathrm{IV}}$; see \cite[pp.~158--159]{T} and \cite[p.~76, (8.22)]{Oda1979}.
		The resulting error lies in $\Sigma^4\pi_{35}^{9}.$ 
		By \cite[pp.~74--75, (8.12), (8.14), (8.18)]{Oda1979}, this subgroup is contained in the order-$2$ summands of $\pi_{39}^{13}.$ 
		
		The bracket in \eqref{m2:eq:rank13-four-lambda} has zero indeterminacy and contains four times an element of $\{\lambda,8\iota_{31},2\sigma_{31}\}.$ 
		The group description $\pi_{39}^{13}=\frac{\mathbb Z}{8}\{\Sigma\tau^{\mathrm{IV}}\}\oplus(\mathbb Z/2)^4$ shows that a fourfold multiple has no nonzero component in those order-$2$ summands.
		The error must therefore vanish, proving \eqref{m2:eq:rank13-four-lambda}.
		
		Together with the formula for $J[\nu_{12}^2],$   this evaluation shows that the fourth multiple of every element of $\{[\nu_{12}^2],8\iota_{18},2\sigma_{18}\}_1$ has $J$-image $4\Sigma\tau^{\mathrm{IV}}\ne0.$ 
		Since these elements are killed by $8,$   each has order exactly $8.$ 
	\end{proof}
	
	\begin{lem}\label{m2:lem:rank13-initial-lift}
		Assume \eqref{m2:eq:rank13-required-kernel}.
		An element of $\{[\nu_{12}^2],8\iota_{18},2\sigma_{18}\}_1$ can be chosen as a lift $[8\sigma_{12}^2]$ in $\pi_{26}(\mathit{SO}(13)).$ 
		This lift has order $8$ and satisfies $J(4[8\sigma_{12}^2])=4\Sigma\tau^{\mathrm{IV}}\ne0.$ 
		For some integers $u,v,$   with $v$ odd, its double satisfies
		\begin{equation}
			2[8\sigma_{12}^2]\equiv32u[\Sigma^2\rho']_{13}+v([\mu_{10}\sigma_{19}]_{13})'\m [\iota_7]_{13}\bar\sigma_7,\ [\eta_{11}]_{13}\kappa_{12}.
			\label{m2:eq:rank13-double-before-choice}
		\end{equation}
	\end{lem}
	
	\begin{proof}
		In light of Lemma~\ref{m2:lem:rank13-toda-order}, every element of the indicated bracket has order $8$ and has a fourth multiple with nonzero $J$-image $4\Sigma\tau^{\mathrm{IV}}.$ 
		No such element belongs to the subgroup in \eqref{m2:eq:rank13-inherited}.
		For an order-$8$ element of that subgroup, the coefficient on the order-$128$ factor must be divisible by $16.$ 
		Its fourth multiple is consequently a multiple of $64[\Sigma^2\rho']_{13},$   whose $J$-image is zero; all the other factors are killed by $4.$ 
		
		By exactness and \eqref{m2:eq:rank13-required-kernel}, the projection of a bracket element is therefore $8\sigma_{12}^2.$ 
		Choose one as the initial lift $[8\sigma_{12}^2].$ 
		Its double has order $4$ and lies in the inherited subgroup.
		The coefficient of its order-$128$ generator must thus be divisible by $32,$   which gives \eqref{m2:eq:rank13-double-before-choice}.
		
		It remains to prove that $v$ is odd.
		The relevant $J$-images follow from \cite[p.~63, Proposition~4.2(1), and pp.~75--76, (8.18)--(8.22)]{Oda1979}.
		On the cyclic $\frac{\mathbb Z}{8}\{\Sigma\tau^{\mathrm{IV}}\}$-summand, the coefficient of $J([\Sigma^2\rho']_{13})$ is even, whereas that of $J(([\mu_{10}\sigma_{19}]_{13})')$ is twice an odd integer.
		The other two generators of \eqref{m2:eq:rank13-inherited} have images killed by $2.$ 
		If $v$ were even, doubling the $J$-image of \eqref{m2:eq:rank13-double-before-choice} would therefore give zero.
		This contradicts $J(4[8\sigma_{12}^2])=4\Sigma\tau^{\mathrm{IV}}\ne0.$ 
		Hence $v$ is odd.
	\end{proof}
	
	\begin{prop}\label{m2:lem:rank13-conditional}
		The structure of $\pi_{26}(\mathit{SO}(13))$ is as follows.
		\[
		\pi_{26}(\mathit{SO}(13))=\mathrm{span}\{[\Sigma^2\rho']_{13},\;[8\sigma_{12}^2],\;[\eta_{11}]_{13}\kappa_{12},\;[\iota_7]_{13}\bar\sigma_7\}\cong\mathbb Z/128\oplus\mathbb Z/8\oplus(\mathbb Z/2)^2.
		\]
		For the odd integer $\ell_{0}$ of Proposition~\ref{m2:lem:rank12-extension}, the lift $[8\sigma_{12}^2]$ can be chosen so that
		\begin{equation}
			2[8\sigma_{12}^2]\equiv2\ell_{0}[\sigma_9\eta_{16}\mu_{17}]_{13}+[\mu_{10}\sigma_{19}]_{13}-4[\Sigma^2\rho']_{13}\m[\eta_{11}]_{13}\kappa_{12},\ [\iota_7]_{13}\bar\sigma_7.
			\label{m2:eq:rank13-normalized-double}
		\end{equation}
	\end{prop}
	
	\begin{proof}
		Start with the lift provided by Lemma~\ref{m2:lem:rank13-initial-lift}, and use the integers $u,v$ in \eqref{m2:eq:rank13-double-before-choice}.
		Replace this lift by $[8\sigma_{12}^2]-16u[\Sigma^2\rho']_{13}.$ 
		The subtracted element belongs to the inherited subgroup and is killed by $8.$ 
		Thus the projection remains $8\sigma_{12}^2,$   and the new lift is still killed by $8.$ 
		Its double has coefficient $v$ on the order-$4$ factor and has no component on the order-$128$ factor, modulo the two order-$2$ factors.
		
		Since $v$ is odd, multiply the new lift by an odd integer inverse to $v$ modulo $4.$ 
		Its projection remains $8\sigma_{12}^2,$   because that element has order $2.$ 
		The resulting double has coefficient $1$ on $([\mu_{10}\sigma_{19}]_{13})'.$ 
		Expanding the definition of this generator gives \eqref{m2:eq:rank13-normalized-double}.
		The double has order exactly $4,$   so the normalized lift has order exactly $8.$ 
		This normalization is not asserted to preserve membership in the original Toda bracket.
		
		The normalized relation expresses the order-$4$ generator of the inherited subgroup as twice the chosen lift, modulo the two order-$2$ generators.
		Consequently, \eqref{m2:eq:rank13-inherited} and exactness show that the four elements in the proposition generate $\pi_{26}(\mathit{SO}(13)).$ 
		Their orders divide $128,8,2,2,$   respectively.
		The same inherited subgroup description and \eqref{m2:eq:rank13-required-kernel} give $$|\pi_{26}(\mathit{SO}(13))|=128\cdot4\cdot2\cdot2\cdot2=2^{12}.$$ Hence the resulting surjection from $\mathbb Z/128\oplus\mathbb Z/8\oplus(\mathbb Z/2)^2$ has source and target of the same finite order and is an isomorphism.
	\end{proof}
	\section{The Cases $n\geq14$}\label{m2:sec:degree26-high}
	\subsection{The Cases $14\leq n\leq17$}
	
	We begin the larger-$n$ calculation with $\mathit{SO}(14).$ Here the lift of $\Sigma\theta\eta_{25}$
	and its relation to the preceding generators determine the large cyclic factor in degree $26.$
	
	\begin{prop}\label{g1016:prop:26-14}
		The structure of $\pi_{26}(\mathit{SO}(14))$ is as follows.
		\[ \GroupLayout{\pi_{26}(\mathit{SO}(14)) = \mathrm{span}\{[ (\Sigma\theta )\eta_{25}],\; [8\sigma^2_{12}]_{14}\}\cong \mathbb Z/256\oplus \mathbb Z/8,}{\begin{aligned}\pi_{26}(\mathit{SO}(14)) = \mathrm{span}\{&[ (\Sigma\theta )\eta_{25}],\\
				&[8\sigma^2_{12}]_{14}\}\cong \mathbb Z/256\oplus \mathbb Z/8,\end{aligned}} \]
		where
		$2[(\Sigma\theta)\eta_{25}]\equiv[\Sigma^2\rho']_{14}\m[8\sigma^2_{12}]_{14}, 2[\Sigma^2\rho']_{14}.$
	\end{prop}
	\begin{proof}
		Consider the  exact sequence
		\[\begin{aligned}
			&\pi^{13}_{27} \stackrel{\Delta}{\longrightarrow} \pi_{26}(\mathit{SO}(13)) \stackrel{i_\ast}{\longrightarrow} \pi_{26}(\mathit{SO}(14)) \stackrel{p_\ast}{\longrightarrow} \pi^{13}_{26} \stackrel{\Delta}{\longrightarrow} \pi_{25}(\mathit{SO}(13)).
		\end{aligned}\]
		Here $$\pi_{26}(\mathit{SO}(13)) = \frac{\mathbb Z}{128}\{[\Sigma^2 \rho']_{13}\} \oplus
		\frac{\mathbb Z}{8}\{[8\sigma^2_{12}]\} \oplus \frac{\mathbb Z}{2}\{[\eta_{11}]_{13}\kappa_{12}\} \oplus
		\frac{\mathbb Z}{2}\{[\iota_7]_{13}\bar{\sigma_7}\},$$ $\pi^{13}_{27}=\frac{\mathbb Z}{16}\{\sigma^2_{13}\} \oplus
		\frac{\mathbb Z}{2}\{\kappa_{13}\}.$ Recall that $\Delta(\iota_{13})=[\eta_{11}]_{13}$ \cite[Table 3]{Ka}. Then, $$\Delta(\kappa_{13})=[\eta_{11}]_{13}\kappa_{12},\quad\Delta(\sigma^2_{13})=[\eta_{11}]_{13}\sigma^2_{12}.$$ Now we observe that
		$J(\Delta(\sigma^2_{13}))=\sigma_{13}\bar{\sigma_{20}}=J([\iota_7]_{13}\bar{\sigma_7}).$ We must compare these classes in $\pi_{26}(\mathit{SO}(13))$ itself. From $J([\eta_{11}]_{13})
		=\Sigma \theta,$ we find $J([\eta_{11}]_{13}\kappa_{12})=(\Sigma\theta)
		\kappa_{25},$   a class of order $2$ in $\pi_{39}^{13}$
		according to \cite{Oda1979}.
		
		Using Sugawara's relation $$\Sigma^{17}( p_{14} \circ \Delta(\sigma_{14}))=2\Sigma^{16}\sigma_{14}$$  (\cite[Lemma 3.10]{Ka})
		and the fact that $\pi^{13}_{20}$ is in the stable range, we obtain $\Delta(2\sigma^2_{13})=0.$
		Thus,
		$\Delta(\sigma^2_{13}) \equiv [\iota_7]_{13} \bar{\sigma_7} \bmod 64[\Sigma^2
		\rho']_{13}, 4[8\sigma^2_{12}]_{13}.$ From the
		calculation of $\pi_{25}(\mathit{SO}(14)),$  we know the morphism $\Delta\colon\pi^{13}_{26}=\frac{\mathbb Z}{2}\{ (\Sigma\theta )\eta_{25}\} \to
		\pi_{25}(\mathit{SO}(13))$ is $0.$ Using the relation $$\mathrm{ord} \Delta (P(\iota_{29}))=256$$ from 
		\cite[Theorem 1.3]{HKM} and $p_{14 \ast}(\Delta
		(P(\iota_{29})))=(\Sigma\theta)\eta_{25}$ from \cite[ Lemma 1.4]{HKM},
		we infer $\mathrm{ord} [(\Sigma\theta)\eta_{25}]=256,$ and $2[(\Sigma\theta)\eta_{25}]\equiv[\Sigma^2\rho']_{14}\m[8\sigma^2_{12}]_{14}, 2[\Sigma^2\rho']_{14}.$
	\end{proof}
	
	The group becomes simpler after the next inclusion. The following proposition records the surviving
	cyclic group at $n=15$ in terms of the lift already chosen at $n=13.$
	
	\begin{prop}
		The structure of $\pi_{26}(\mathit{SO}(15))$ is as follows.
		\[ \FitMath{\pi_{26}(\mathit{SO}(15)) = \mathrm{span}\{[8\sigma^2_{12}]_{15}\}\cong \mathbb Z/8,} \]
		where $[8\sigma_{12}^2]_{15}=(i_{15})_*[8\sigma_{12}^2]_{14}$ and $8[8\sigma_{12}^2]_{15}=0.$   \qed
	\end{prop}
	
	Using the relation $\pi_{26}(\mathit{SO})=0$ \cite{Bo} and the results in \cite[Table 1]{HoM},
	we have Table~\ref{m2:R26}.
	\par
	\begin{table}[H]
		\centering
		\caption{$\pi_{26}(\mathit{SO}(n))$ for $16\leq n\leq 21$}\label{m2:R26}
		\begin{tabular*}{\textwidth}{@{\extracolsep{\fill}}clclcl@{}}
			\toprule
			\textbf{Group} & \textbf{Structure} & \textbf{Group} & \textbf{Structure} & \textbf{Group} & \textbf{Structure} \\
			\midrule
			$\pi_{26}(\mathit{SO}(16))$ & $\mathbb Z/8 \oplus \mathbb Z/8$
			& $\pi_{26}(\mathit{SO}(17))$ & $\mathbb Z/8$
			& $\pi_{26}(\mathit{SO}(18))$ & $\mathbb Z/8$  \\ \midrule
			$\pi_{26}(\mathit{SO}(19))$ & $\mathbb Z/8$
			& $\pi_{26}(\mathit{SO}(20))$ & $\mathbb Z/16 \oplus \mathbb Z/4$
			& $\pi_{26}(\mathit{SO}(21))$ & $\mathbb Z/8$\\ \bottomrule
		\end{tabular*}
		
	\end{table}
	\par
	We have Table~\ref{m2:R26K} from \cite{Ke}.
	\par
	\begin{table}[H]
		\centering
		\caption{$\pi_{26}(\mathit{SO}(n))$ for $22\le n\le 27$}\label{m2:R26K}
		\begin{tabular*}{\textwidth}{@{\extracolsep{\fill}}clclcl@{}}
			\toprule
			\textbf{Group} & \textbf{Structure} & \textbf{Group} & \textbf{Structure} & \textbf{Group} & \textbf{Structure} \\
			\midrule
			$\pi_{26}(\mathit{SO}(22))$ & $\mathbb Z/8$
			& $\pi_{26}(\mathit{SO}(23))$ & $\mathbb Z/8$
			& $\pi_{26}(\mathit{SO}(24))$ & $\mathbb Z/8 \oplus \mathbb Z/8$  \\ \midrule
			$\pi_{26}(\mathit{SO}(25))$ & $\mathbb Z/8$
			& $\pi_{26}(\mathit{SO}(26))$ & $\mathbb Z/4$
			& $\pi_{26}(\mathit{SO}(27))$ & $\mathbb Z/2$\\ \bottomrule
		\end{tabular*}
		
	\end{table}
	\par
	
	By Tables~\ref{m2:R26} and \ref{m2:R26K}, we have the group structure of $\pi_{26}(\mathit{SO}(n))$ for $n \ge 16.$
	We now determine generators and their relations in
	$\pi_{26}(\mathit{SO}(n))$ for $n \ge 16.$
	
	We obtain $\Delta(\nu_{8n-1})=0$ \cite[p. 410]{GM}.
	The lift $[\nu_{8n+7}]$ has order $8$ \cite[p.~161, Table]{Ke}. We therefore choose
	a lift $[\zeta_{8n+7}]$ of $\zeta_{8n+7}$ as follows:
	$[\zeta_{8n+7}]\in\{[\nu_{8n+7}],8\iota_{8n+10},2\sigma_{8n+10}\}_1.$
	
	The displayed bracket gives a convenient choice of $[\zeta_{15}].$ We use this lift to describe the
	new contribution from the homotopy groups of spheres at $n=16$ and to determine its order within the extension.
	
	\begin{prop}
		The structure of $\pi_{26}(\mathit{SO}(16))$ is as follows.
		\[ \GroupLayout{\pi_{26}(\mathit{SO}(16)) = \mathrm{span}\{[\zeta_{15}],\; [8\sigma^2_{12}]_{16}\}\cong (\mathbb Z/8)^{2},}{\begin{aligned}\pi_{26}(\mathit{SO}(16)) = \mathrm{span}\{&[\zeta_{15}],\\
				&[8\sigma^2_{12}]_{16}\}\cong (\mathbb Z/8)^{2},\end{aligned}} \]
		where $[\zeta_{15}]\in\{[\nu_{15}],8\iota_{18},2\sigma_{18}\}_1 $ and $ 8[\zeta_{15}]=0.$   \end{prop}
	\begin{proof}
		Consider the  exact sequence
		\[\begin{aligned}
			&\pi^{15}_{27} \stackrel{\Delta}{\longrightarrow} \pi_{26}(\mathit{SO}(15)) \stackrel{i_\ast}{\longrightarrow} \pi_{26}(\mathit{SO}(16)) \stackrel{p_\ast}{\longrightarrow} \pi^{15}_{26} \stackrel{\Delta}{\longrightarrow} \pi_{25}(\mathit{SO}(15)).
		\end{aligned}\]
		Here $\pi^{15}_{27} =0$ and the generator $\zeta_{15}$ of
		$\pi_{26}^{15}$ goes to zero in $\pi_{25}(\mathit{SO}(15)).$
		
		Let us examine $\{[\nu_{15}], 8\iota_{18}, 2\sigma_{18}\} \subseteq
		\pi_{26}(\mathit{SO}(16)),$ where $$[\nu_{15}] \in \pi_{18}(\mathit{SO}(16))=\frac{\mathbb Z}{8}\{[\nu_{15}]\}
		+\frac{\mathbb Z}{8}\{[\nu^2_{12}]_{16}\}.$$
		
		In light of \cite[p.54]{T}, we deduce $$p_{*}\{[\nu_{15}], 8\iota_{18}, 2\sigma_{18}\}
		\subseteq \{\nu_{15},8\iota_{18},2\sigma_{18}\}=\zeta_{15}.$$ \noindent Recall that $\nu_{15}\pi^{18}_{26}=\mathrm{span}\{\nu_{15}\bar\nu_{18},
		\nu_{15}\varepsilon_{18}\}=0$ and $0=\pi^{15}_{19} \circ 2\sigma_{19}.$
		Therefore, we can take $[\zeta_{15}] \in \{[\nu_{15}], 8\iota_{18},
		2\sigma_{18}\}.$ Then, we conclude that $$8[\zeta_{15}] \in [\nu_{15}]\{8\iota_{18}, 2 \sigma_{18},
		8\iota_{25}\} \subseteq [\nu_{15}]\{2\iota_{18},16\sigma_{18},4\iota_{25}\}=0.$$ \noindent Thus, we obtain $8[\zeta_{15}]=0.$
	\end{proof}
	To reach $n=17,$ we compute the relation imposed by the boundary of $\zeta_{16}.$ This expresses
	the contribution inherited from $[8\sigma_{12}^2]$ in terms of the surviving $\zeta$-lift.
	
	\begin{prop}
		The structure of $\pi_{26}(\mathit{SO}(17))$ is as follows.
		\[ \FitMath{\pi_{26}(\mathit{SO}(17)) = \mathrm{span}\{[\zeta_{15}]_{17}\}\cong \mathbb Z/8,} \]
		where $\Delta\zeta_{16}=2[\zeta_{15}]+a[8\sigma^2_{12}]_{16}$ for some odd integer $a.$  Moreover,  $2[\zeta_{15}]_{17}+a[8\sigma^2_{12}]_{17}=0.$ 
	\end{prop}
	
	\begin{proof}
		Consider the  exact sequence
		\[\begin{aligned}
			&\pi^{16}_{27} \stackrel{\Delta}{\longrightarrow} \pi_{26}(\mathit{SO}(16)) \stackrel{i_\ast}{\longrightarrow} \pi_{26}(\mathit{SO}(17)) \stackrel{p_\ast}{\longrightarrow} \pi^{16}_{26}.
		\end{aligned}\]
		In view of  \cite{Su1953,Su1954}, we see that $\Delta(\zeta_{16})=2[\zeta_{15}]
		+a[8\sigma^2_{12}]_{16},$ so that
		$\Delta(4\zeta_{16})=4a[8\sigma^2_{12}]_{16}$ holds.
		
		We now show that $a$ is odd.
		Consider the images of $J : \pi_{26}(\mathit{SO}(16)) \to \pi^{16}_{42}$ and $P:
		\pi^{33}_{44} \to \pi^{16}_{42}.$ By \cite[p.77]{Oda1979},
		we obtain $J\Delta(4\zeta_{16})=4P(\zeta_{33}) = 4\Sigma^4\tau^{\mathrm{IV}}$ and
		$$J(4[8\sigma_{12}^2]_{16})\equiv 4\Sigma^4\tau^{\mathrm{IV}} \bmod
		\nu^2_{16}\bar\kappa_{22}, \eta_{16}\mu_{3, 17}.$$
		Thus $\Delta\zeta_{16}=2[\zeta_{15}]+a[8\sigma_{12}^2]_{16}$ with $a$ odd. This comparison determines the parity of $a,$   since $4a\equiv4\pmod8$ is equivalent to $a\equiv1\pmod2.$   \end{proof}
	
	\subsection{The Cases $18\leq n\leq23$}
	
	At $n=18,$ we choose a generator lying over $\eta_{17}\varepsilon_{18}.$ The doubling relation
	below records how the preceding $\zeta$-class is represented in this new choice of generator.
	
	\begin{prop}
		The structure of $\pi_{26}(\mathit{SO}(18))$ is as follows.
		\[ \FitMath{\pi_{26}(\mathit{SO}(18)) = \mathrm{span}\{[\eta_{17}\varepsilon_{18}]\}\cong \mathbb Z/8,} \] where $2[\eta_{17}\varepsilon_{18}]\ty[\zeta_{15}]_{18}\m 2[\zeta_{15}]_{18}.$
	\end{prop}
	
	\begin{proof}
		Consider the  exact sequence
		\[\begin{aligned}
			&\pi^{17}_{27} \stackrel{\Delta}{\longrightarrow} \pi_{26}(\mathit{SO}(17)) \stackrel{i_\ast}{\longrightarrow} \pi_{26}(\mathit{SO}(18)) \stackrel{p_\ast}{\longrightarrow} \pi^{17}_{26} \stackrel{\Delta}{\longrightarrow} \pi_{25}(\mathit{SO}(17)).
		\end{aligned}\]
		We first compute $\Delta(\eta_{17}\mu_{18})=4[\zeta_{15}]_{17}.$ By using $J
		: \pi_{26}(\mathit{SO}(17)) \to \pi^{17}_{43}$ and $P : \pi^{35}_{45} \to
		\pi^{17}_{43},$ we have
		$P(\eta_{35}\mu_{36})=2\Sigma^5\tau^{\mathrm{IV}}=4\Sigma M^{(1)}_1$ from
		\cite[p.77]{Oda1979}.
		Additionally, by using $H : \pi^{16}_{42} \to \pi^{31}_{42},$ we
		find $H(M_1^{(1)}) = \zeta_{31}.$ Thus $J([\zeta_{15}]_{17}) =
		\Sigma M^{(1)}_1 \bmod \Sigma^2\pi^{15}_{41}.$ Therefore, we conclude that
		$\Delta(\eta_{17}\mu_{18})=4[\zeta_{15}]_{17}.$
		
		Next, we establish
		$2[\eta_{17}\varepsilon_{18}]\ty[\zeta_{15}]_{18}\m2[\zeta_{15}]_{18}.$
		Consider the images of $$J : \pi_{26}(\mathit{SO}(18)) \to \pi^{18}_{44}$$ and $H
		\colon \pi^{18}_{44} \to\pi^{35}_{44}.$ As a consequence of $H(M^{\prime\prime\prime}_2) \equiv
		\eta_{35}\varepsilon_{36} \bmod \eta^2_{35}\sigma_{37}$ \cite[p.~77]{Oda1979} and $[\e_{17},\nu^{3}_{17}]=[\e_{17}, \eta_{17}\varepsilon_{18}+\eta^2_{17}\sigma_{19} ]\neq0$  \cite{Mukaiwhiteheadji}, we deduce  
		$$H(M^{\prime\prime\prime}_2) =
		\eta_{35}\varepsilon_{36}.$$    \noindent   Then we see that
		$J([\eta_{17}\varepsilon_{18}]) \equiv M^{\prime\prime\prime}_2 \bmod
		\Sigma\pi_{43}^{17},$ so that $J(2[\eta_{17}\varepsilon_{18}]) \equiv
		2M^{\prime\prime\prime}_2 \bmod 4M_2^{\prime\prime\prime}$ from \cite[p.77]{Oda1979}. In addition,
		we know $$J([\zeta_{15}]_{18}) \equiv 2M_2^{\prime\prime\prime} \bmod \eta_{18}\mu_{3,
			19}, \nu^2_{18}\bar\kappa_{24}.$$ \noindent Therefore, we obtain
		$2[\eta_{17}\varepsilon_{18}]\ty[\zeta_{15}]_{18}\m 2[\zeta_{15}]_{18}.$
	\end{proof}

	We now turn our attention to the computation of \(\pi_{26}(\mathit{SO}(19))\), which will be carried out in the following proposition.
	
	\begin{prop}
		The structure of $\pi_{26}(\mathit{SO}(19))$ is as follows.
		\[ \FitMath{\pi_{26}(\mathit{SO}(19)) = \mathrm{span}\{[\varepsilon_{18}]\}\cong \mathbb Z/8,} \] where $2[\varepsilon_{18}]\ty[\eta_{17}\varepsilon_{18}]_{19}\m2[\eta_{17}\varepsilon_{18}]_{19}.$
	\end{prop}
	\begin{proof}
		Consider the  exact sequence
		\[\begin{aligned}
			&\pi^{18}_{27} \stackrel{\Delta}{\longrightarrow} \pi_{26}(\mathit{SO}(18)) \stackrel{i_\ast}{\longrightarrow} \pi_{26}(\mathit{SO}(19)) \stackrel{p_\ast}{\longrightarrow} \pi^{18}_{26} \stackrel{\Delta}{\longrightarrow} \pi_{25}(\mathit{SO}(18)).
		\end{aligned}\]
		Here, $\pi_{25}(\mathit{SO}(18)) \cong \frac{\mathbb Z}{4}\{[\nu_{15}]_{18}\sigma_{18}\} + \frac{\mathbb Z}{2}\{[\iota_7]_{18}\eta_7\bar\mu_8\};$ the boundary of
		$\varepsilon_{18}$ is zero, whereas the boundary of
		$\eta_{18}\sigma_{19}$ is $2[\nu_{15}]_{18}\sigma_{18}.$  Note that $[\nu_{15}]_{18}\sigma_{18} =
		\pm[\sigma^2_{11}]_{18}.$ 
		
		Consider the map $J : \pi_{26}(\mathit{SO}(18)) \to \pi^{18}_{44}.$ Since
		$J([\eta_{17}\varepsilon_{18}]) \equiv M^{\prime\prime\prime}_2 \bmod
		\mathrm{span}\{\nu^2_{18}\bar\kappa_{24}, \eta_{18}\mu_{3,19}\},$ the map $J$ is
		injective.
		
		In light of \cite[p.67, p.77]{Oda1979}, we deduce $P(\eta^2_{37}\sigma_{39}) =
		P(\eta_{37}\varepsilon_{38})=0$ and $P(\mu_{37})=4M_2^{\prime\prime\prime}.$
		Therefore, we find $\Delta(\mu_{18})=4[\eta_{17}\varepsilon_{18}]$ and $\Delta(\eta^2_{18}\sigma_{20})=\Delta(\eta_{18}\varepsilon_{19})=0.$
		Since $J([\varepsilon_{18}]) \equiv M^{\prime\prime}_2 \bmod \Sigma\pi^{18}_{44},$ we
		obtain $2[\varepsilon_{18}]\ty[\eta_{17}\varepsilon_{18}]_{19}\m2[\eta_{17}\varepsilon_{18}]_{19}$ by the
		determination of $\pi_{26}(\mathit{SO}(18)).$
	\end{proof}
	
	The $n=20$ calculation involves both a boundary class and an inherited $\varepsilon$-class. We
	describe a choice of generators adapted to their relation, so that the two cyclic factors can be read
	off explicitly.
	
	\begin{prop}
		The structure of $\pi_{26}(\mathit{SO}(20))$ is as follows.
		\[ \GroupLayout{\pi_{26}(\mathit{SO}(20)) = \mathrm{span}\{[2\sigma_{19}],\; [\varepsilon_{18}]_{20} - 2[2\sigma_{19}]\}\cong \mathbb Z/16\oplus \mathbb Z/4,}{\begin{aligned}\pi_{26}(\mathit{SO}(20)) = \mathrm{span}\{&[2\sigma_{19}],\\
				&[\varepsilon_{18}]_{20} - 2[2\sigma_{19}]\}\cong \mathbb Z/16\oplus \mathbb Z/4,\end{aligned}} \] where
		$[2\sigma_{19}]=\Delta\sigma_{20}$ and
		$8[2\sigma_{19}]=4[\varepsilon_{18}]_{20}.$
	\end{prop}
	\begin{proof}
		Consider the  exact sequence
		\[\begin{aligned}
			&\pi^{19}_{27} \stackrel{\Delta}{\longrightarrow} \pi_{26}(\mathit{SO}(19)) \stackrel{i_\ast}{\longrightarrow} \pi_{26}(\mathit{SO}(20)) \stackrel{p_\ast}{\longrightarrow} \pi^{19}_{26} \stackrel{\Delta}{\longrightarrow} \pi_{25}(\mathit{SO}(19)),
		\end{aligned}\]
		and consider the map $J : \pi_{26}(\mathit{SO}(19)) \to \pi^{19}_{45}.$ Since
		$J([\varepsilon_{18}]) \equiv M^{\prime\prime}_2,$ the map $J$ is a
		monomorphism. By \cite[p.67]{Oda1979},
		$P(\eta_{39}\sigma_{40})=P(\varepsilon_{39})=0.$ Thus $\Delta :
		\pi^{19}_{27} \to \pi_{26}(\mathit{SO}(19))$ is a zero map.
		
		Next, consider $(p_{20})_*:\pi_{26}(\mathit{SO}(20))\longrightarrow\pi_{26}^{19},$  
		$\Sigma^{20}:\pi_{26}^{19}=\frac{\mathbb Z}{16}\{\sigma_{19}\}\xrightarrow{\cong}\pi_{46}^{39}=\frac{\mathbb Z}{16}\{\sigma_{39}\}$
		and $H : \pi^{20}_{46}
		\to\pi_{46}^{39}=\frac{\mathbb Z}{16}\{\sigma_{39}\}.$ We note that $H(J\Delta\sigma_{20})=H(P\sigma_{41})=H(P\iota_{41})\sigma_{39}=2\sigma_{39}.$ Moreover, $(p_{20})_*(\Delta\sigma_{20})=2\sigma_{19},$   so $\Delta\sigma_{20}=[2\sigma_{19}]$
		is a lift of $2\sigma_{19}.$  To determine its order, use
		$\pm J\Delta\sigma_{20}=[\iota_{20},\sigma_{20}].$ 
		For even $n\geq10,$     $[\iota_n,\sigma_n]$ has order $16$; see \cite{Mukaiwhiteheadji}. Thus $J\Delta\sigma_{20}$ has order $16.$ 
		Since $16\sigma_{20}=0,$   it follows that $\Delta\sigma_{20}$ has order exactly $16.$ 
		
		Write $a=[\varepsilon_{18}]_{20}$ and $b=\Delta\sigma_{20}.$ 
		The vanishing of the incoming boundary proved above shows that $a$ has order $8$ and
		$\ker(p_{20})_*=\mathrm{span}\{a\}.$ 
		Moreover, $J\Delta\sigma_{19}=[\iota_{19},\sigma_{19}]\ne0$ by  \cite{Mukaiwhiteheadji},
		whereas $2\Delta\sigma_{19}=0$ because $2\sigma_{19}=(p_{20})_*b.$ 
		Hence $\operatorname{Im}(p_{20})_*=\mathrm{span}\{2\sigma_{19}\}\cong\mathbb Z/8,$  
		and the group is generated by $a,b$ and has order $64.$ 
		Now $8b$ is a nonzero element of order $2$ in $\mathrm{span}\{a\},$  
		so $8b=4a.$  Consequently, $a-2b$ has order $4$: its fourth multiple is zero,
		and its projection is $-4\sigma_{19},$   which has order $4.$ 
		The elements $b,a-2b$ generate the group, and the product of their orders is $64.$ 
		They therefore give the asserted direct sum decomposition.
	\end{proof}
	
	We treat the cases $n=21,22,23$ together, following the lift of $\nu_{20}^2$ through the
	successive inclusions. The relevant homotopy groups of spheres simplify these exact sequences and make the behavior
	of this class particularly transparent.
	
	\begin{prop}
		The structure of $\pi_{26}(\mathit{SO}(m))$ for $21\leq m\leq23$ is as follows.
		\[ \begin{aligned}
			&\pi_{26}(\mathit{SO}(21)) = \mathrm{span}\{[\nu^2_{20}]\} \cong \mathbb Z/8,\; \\
			&\pi_{26}(\mathit{SO}(22)) = \mathrm{span}\{[\nu^2_{20}]_{22}\} \cong \mathbb Z/8 , \\
			&\pi_{26}(\mathit{SO}(23)) = \mathrm{span}\{[\nu^2_{20}]_{23}\} \cong\mathbb Z/8,
		\end{aligned} \]
		where $2[\nu^2_{20}]=x[\varepsilon_{18}]_{21}$ for some odd $x.$ 
	\end{prop}
	\begin{proof}
		For $\pi_{26}(\mathit{SO}(21)),$ it follows directly from \cite[p.421 Lemma
		6.2]{GM}.
		Since $\Delta\sigma_{20}=[2\sigma_{19}],$   the $n=20$ decomposition gives
		$(i_{21})_*\pi_{26}(\mathit{SO}(20))=\frac{\mathbb Z}{4}\{[\varepsilon_{18}]_{21}\}.$  The resulting extension
		by $\frac{\mathbb Z}{2}\{\nu_{20}^2\}$ is cyclic of order $8.$  Hence the double of $[\nu_{20}^2]$ is an odd
		multiple of $[\varepsilon_{18}]_{21},$   as stated.
		
		Consider  the exact sequence
		\[\begin{aligned}
			&\pi^{21}_{27} \xrightarrow{ \Delta} \pi_{26}(\mathit{SO}(21)) \xrightarrow{i_\ast } \pi_{26}(\mathit{SO}(22)) \xrightarrow{p_\ast } \pi^{21}_{26}=0.
		\end{aligned}\]
		The displayed vanishing shows that $i_{\ast}$ is surjective. Thus, we have
		$\pi_{26}(\mathit{SO}(22))=\mathrm{span}\{[\nu^2_{20}]_{22}\} \cong \mathbb Z/8 $ by \cite{Ke}.     Since $\pi_{27}^{22}=0$ and $\pi_{26}^{22}=0,$ we obtain
		$\pi_{26}(\mathit{SO}(23))=\mathrm{span}\{[\nu^2_{20}]_{23}\} \cong
		\mathbb Z/8.$
	\end{proof}
	
	\subsection{The Cases $n\geq24$}
	
	Near the stable range, the calculation is governed by lifts of $\nu_{23}$ and the connecting image of
	$\nu_{24}.$ 
	\begin{prop}
		The structure of $\pi_{26}(\mathit{SO}(m))$ for $24\leq m\leq26$ is as follows.
		\[ \FitMath{\pi_{26}(\mathit{SO}(26)) = \mathrm{span}\{[\nu_{23}]_{26}\} \cong \mathbb Z/4 .} \]  The conjectural identification is $\xi=[\nu_{23}]_{26}.$ Here, $\xi:=r_{13}\omega_{13}(\mathbb{C}).$
		\[ \begin{aligned}
			&\pi_{26}(\mathit{SO}(25)) = \mathrm{span}\{[\nu_{23}]_{25}\} \cong \mathbb Z/8 , \\
			&\pi_{26}(\mathit{SO}(24)) = \mathrm{span}\{\Delta\nu_{24},\; [\nu_{23}]\}\cong (\mathbb Z/8)^{2},
		\end{aligned} \]
		where
		\[
		\begin{aligned}
			&8[\nu_{23}]=0,\qquad 2[\nu_{23}]-\Delta\nu_{24}=x[\nu_{20}^2]_{24},\qquad x\text{: odd}, \\
			&\Delta(\eta_{25}^2)=4[\nu_{23}]_{25},\qquad [\nu_{23}]_{26}=u\xi,\qquad u\text{ odd}, \\
			&\Delta\eta_{26}=2\xi.
		\end{aligned}
		\]
	\end{prop}
	\begin{proof}
		We see that $\xi=r_{13}\omega_{13}(\mathbb{C})=\tau'_{26}$ by \cite[p.411
		(4.2)]{GM}.
		At $n=24,$   the order-$8$ lift and the relation
		\[\begin{aligned}
			&2[\nu_{23}]-\Delta\nu_{24}=x[\nu_{20}^2]_{24},\qquad x\text{: odd},
		\end{aligned}\]
		follow from \cite[p.~410, (3.11), and the proof of Lemma~3.5]{GM}.
		Thus $\Delta\nu_{24}$ and $[\nu_{23}]$ form an alternative generating set of the two order-$8$ summands in Theorem~\ref{m2:thm1}.
		Consider  the exact sequence
		\[\begin{aligned}
			&\pi^{25}_{27} \xrightarrow{ \Delta} \pi_{26}(\mathit{SO}(25)) \xrightarrow{i_\ast } \pi_{26}(\mathit{SO}(26)) \xrightarrow{p_\ast } \pi^{25}_{26} \xrightarrow{ \Delta} \pi_{25}(\mathit{SO}(25)).
		\end{aligned}\]
		The right $\Delta$ is injective by \cite[Theorem~1]{Ke}, so $p_*=0$ and $i_*$ is surjective. The left
		boundary has nonzero image by \cite[Theorem~2]{Ke}\uprightsemicolon{} it is therefore the order-$2$ subgroup of
		$\pi_{26}(\mathit{SO}(25))=\frac{\mathbb Z}{8}\{[\nu_{23}]_{25}\},$   giving
		$\Delta(\eta_{25}^2)=4[\nu_{23}]_{25}.$
		It follows that $[\nu_{23}]_{26}$ generates $\pi_{26}(\mathit{SO}(26))$ and equals $u\xi$ for some odd
		$u.$  The nonzero order-$2$ boundary at the next value of $n$ is $\Delta\eta_{26}=2\xi.$  The equality $u=1$ is
		the conjectural identification stated above.
	\end{proof}
	
	\chapter{The Homotopy Groups $\pi_{27}(\mathit{SO}(n))$}\label{ch:degree27}\label{m2:sec:degree27}
	
	Degree $27$ has stable group $\mathbb Z_{(2)},$   reached for $n\geq29.$  The unstable calculation
	must therefore account for the free part as well as the torsion. We first study the cases with small values of $n,$   then
	the lifted classes for $15\leq n\leq17,$   and finally the larger-$n$ information.

	The descriptions of the groups $\pi_{*}(S^{j})$ in
	Appendix~\ref{tlqdzj} and the composition relations among
	generators in Toda's range recorded in \^{O}guchi's table
	in Appendix~\ref{app:oguchi-relations} are used freely
	in our arguments.
	
	\section{Main Theorem and Generators}

	\begingroup
	\setlength{\abovedisplayskip}{6pt plus 1pt minus 1pt}
	\setlength{\belowdisplayskip}{6pt plus 1pt minus 1pt}
	\setlength{\parindent}{1.5em}
	\begin{thm} \label{m2:thm2}
		\textup{All spaces in this theorem are localized at $2.$  The displayed groups are obtained from the integral homotopy groups by tensoring with $\mathbb Z_{(2)}.$  Free cyclic summands become $\mathbb Z_{(2)}$; for finite groups, this localization retains precisely the $2$-primary summands. The homotopy groups are given as follows.   }\par
		\begin{itemize}[leftmargin=0pt,itemindent=0pt]
			\item[\normalfont(1)] \textup{The Cases $3\leq n\leq6.$ }\par
			\smallskip\noindent\(
			\GroupLayout{\pi_{27}(\mathit{SO}(3)) = \mathrm{span}\{[\eta_2]\delta_3,\; [\eta_2]\cq _3\sigma_{20},\; [\eta_2]\varepsilon'\kappa_{13}\}\cong (\mathbb Z/2)^{3}.}{\begin{aligned}\pi_{27}(\mathit{SO}(3)) = \mathrm{span}\{&[\eta_2]\delta_3,\; [\eta_2]\cq _3\sigma_{20},\\
					&[\eta_2]\varepsilon'\kappa_{13}\}\cong (\mathbb Z/2)^{3}.\end{aligned}}
			\)\par
			\smallskip\par\smallskip\noindent\(
			\GroupLayout{\pi_{27}(\mathit{SO}(4)) \;= {} \mathrm{span}\{[\eta_2]_4\delta_3,\; [\eta_2]_4\cq _3\sigma_{20},\; [\eta_2]_4\varepsilon'\kappa_{13},\; [\iota_3]\delta_3,\; [\iota_3]\cq _3\sigma_{20},\; [\iota_3]\varepsilon'\kappa_{13}\}\cong (\mathbb Z/2)^{6}.}{\begin{aligned}\pi_{27}(\mathit{SO}(4)) \;= {} \mathrm{span}\{&[\eta_2]_4\delta_3,\; [\eta_2]_4\cq _3\sigma_{20},\; [\eta_2]_4\varepsilon'\kappa_{13},\\
					&[\iota_3]\delta_3,\; [\iota_3]\cq _3\sigma_{20},\; [\iota_3]\varepsilon'\kappa_{13}\}\cong (\mathbb Z/2)^{6}.\end{aligned}}
			\)\par
			\smallskip
			
			\par\par\smallskip\noindent\(
			\GroupLayout{\pi_{27}(\mathit{SO}(5)) = \mathrm{span}\{[\iota_3]_5\delta_3,\; [\iota_3]_5\cq _3\sigma_{20},\; [\nu_4^2]\nu_{10}\kappa_{13}\}\cong (\mathbb Z/2)^{3}.}{\begin{aligned}\pi_{27}(\mathit{SO}(5)) = \mathrm{span}\{&[\iota_3]_5\delta_3,\; [\iota_3]_5\cq _3\sigma_{20},\\
					&[\nu_4^2]\nu_{10}\kappa_{13}\}\cong (\mathbb Z/2)^{3}.\end{aligned}}
			\)\par
			\par\smallskip\noindent\(
			\GroupLayout{\pi_{27}(\mathit{SO}(6)) = \mathrm{span}\{[\eta_5^2]\ct _7,\; [\rho^{\mathrm{IV}}]\sigma_{20},\; [ \nu_5]\ccs_8\}\cong (\mathbb Z/4)^{2}\oplus \mathbb Z/2.}{\begin{aligned}\pi_{27}(\mathit{SO}(6)) = \mathrm{span}\{&[\eta_5^2]\ct _7,\; [\rho^{\mathrm{IV}}]\sigma_{20},\\
					&[ \nu_5]\ccs_8\}\cong (\mathbb Z/4)^{2}\oplus \mathbb Z/2.\end{aligned}}
			\)\par
			\item[\normalfont(2)] \textup{The Cases $7\leq n\leq14.$ }\par
			\smallskip\noindent\(
			\GroupLayout{\pi_{27}(\mathit{SO}(7)) = \mathrm{span}\{[\eta_6]\ct _7,\; [ \nu_5]_7\ccs_8\}\cong \mathbb Z/8\oplus \mathbb Z/2.}{\begin{aligned}\pi_{27}(\mathit{SO}(7)) = \mathrm{span}\{&[\eta_6]\ct _7,\\
					&[ \nu_5]_7\ccs_8\}\cong \mathbb Z/8\oplus \mathbb Z/2.\end{aligned}}
			\)\par
			\smallskip\par\smallskip\noindent\(
			\GroupLayout{\pi_{27}(\mathit{SO}(8)) = \mathrm{span}\{[\iota_7]\ct _7,\; [\eta_6]_8\ct _7,\; [ \nu_5]_8\ccs_8\}\cong (\mathbb Z/8)^{2}\oplus \mathbb Z/2.}{\begin{aligned}\pi_{27}(\mathit{SO}(8)) = \mathrm{span}\{&[\iota_7]\ct _7,\; [\eta_6]_8\ct _7,\\
					&[ \nu_5]_8\ccs_8\}\cong (\mathbb Z/8)^{2}\oplus \mathbb Z/2.\end{aligned}}
			\)\par
			\smallskip
			
			\par\par\smallskip\noindent\(
			\GroupLayout{\pi_{27}(\mathit{SO}(9)) = \mathrm{span}\{[\iota_7]_9\ct _7,\; [ \nu_5]_9\ccs_8\}\cong \mathbb Z/8\oplus \mathbb Z/2.}{\begin{aligned}\pi_{27}(\mathit{SO}(9)) = \mathrm{span}\{&[\iota_7]_9\ct _7,\\
					&[ \nu_5]_9\ccs_8\}\cong \mathbb Z/8\oplus \mathbb Z/2.\end{aligned}}
			\)\par
			\par\smallskip\noindent\(
			\GroupLayout{\pi_{27}(\mathit{SO}(10)) = \mathrm{span}\{[\sigma_9\zeta_{16}],\; [\iota_7]_{10}\ct _7\}\cong (\mathbb Z/8)^{2}.}{\begin{aligned}\pi_{27}(\mathit{SO}(10)) = \mathrm{span}\{&[\sigma_9\zeta_{16}],\\
					&[\iota_7]_{10}\ct _7\}\cong (\mathbb Z/8)^{2}.\end{aligned}}
			\)\par
			\par\smallskip\noindent\(
			\GroupLayout{\pi_{27}(\mathit{SO}(11)) = \mathrm{span}\{[\iota_7]_{11}\ct _7,\; [\eta_{10}\mu_{11}]\sigma_{20}\}\cong \mathbb Z/8\oplus \mathbb Z/2.}{\begin{aligned}\pi_{27}(\mathit{SO}(11)) = \mathrm{span}\{&[\iota_7]_{11}\ct _7,\\
					&[\eta_{10}\mu_{11}]\sigma_{20}\}\cong \mathbb Z/8\oplus \mathbb Z/2.\end{aligned}}
			\)\par
			\smallskip\par\smallskip\noindent\(
			\GroupLayout{\pi_{27}(\mathit{SO}(12)) = \mathrm{span}\{[\iota_7]_{12}\ct _7,\; [\eta_{10}\mu_{11}]_{12}\sigma_{20},\; [\mu_{11}]\sigma_{20}\}\cong \mathbb Z/4\oplus (\mathbb Z/2)^{2}.}{\begin{aligned}\pi_{27}(\mathit{SO}(12)) = \mathrm{span}\{&[\iota_7]_{12}\ct _7,\; [\eta_{10}\mu_{11}]_{12}\sigma_{20},\\
					&[\mu_{11}]\sigma_{20}\}\cong \mathbb Z/4\oplus (\mathbb Z/2)^{2}.\end{aligned}}
			\)\par
			\smallskip
			
			\par\par\smallskip\noindent\(
			\GroupLayout{\pi_{27}(\mathit{SO}(13)) = \mathrm{span}\{[\iota_7]_{13}\ct _7,\; [\mu_{11}]_{13}\sigma_{20}\}\cong \mathbb Z/4\oplus \mathbb Z/2.}{\begin{aligned}\pi_{27}(\mathit{SO}(13)) = \mathrm{span}\{&[\iota_7]_{13}\ct _7,\\
					&[\mu_{11}]_{13}\sigma_{20}\}\cong \mathbb Z/4\oplus \mathbb Z/2.\end{aligned}}
			\)\par
			\par\smallskip\noindent\(
			\GroupLayout{\pi_{27}(\mathit{SO}(14)) = \mathrm{span}\{[2\iota_{13}]\sigma_{13}^2,\; [\iota_7]_{14}\ct _7\}\cong \mathbb Z/8\oplus \mathbb Z/2.}{\begin{aligned}\pi_{27}(\mathit{SO}(14)) = \mathrm{span}\{&[2\iota_{13}]\sigma_{13}^2,\\
					&[\iota_7]_{14}\ct _7\}\cong \mathbb Z/8\oplus \mathbb Z/2.\end{aligned}}
			\)\par
			\item[\normalfont(3)] \textup{The Cases $n\geq15.$ }\par
			\smallskip\noindent\(
			\FitMath{\pi_{27}(\mathit{SO}(15+n)) = \mathrm{span}\{[256P(\iota_{29})]_{15+n}\}\cong \mathbb Z_{(2)},\;\;0\leq n\leq 3.}
			\)\par
			\par\smallskip\noindent\(
			\GroupLayout{\pi_{27}(\mathit{SO}(19)) = \mathrm{span}\{[256P(\iota_{29})]_{19},\; [\eta^2_{18}]\sigma_{20},\; [\varepsilon_{18}]\eta_{26}\}\cong \mathbb Z_{(2)}\oplus (\mathbb Z/2)^{2}.}{\begin{aligned}\pi_{27}(\mathit{SO}(19)) = \mathrm{span}\{&[256P(\iota_{29})]_{19},\; [\eta^2_{18}]\sigma_{20},\\
					&[\varepsilon_{18}]\eta_{26}\}\cong \mathbb Z_{(2)}\oplus (\mathbb Z/2)^{2}.\end{aligned}}
			\)\par
			\smallskip\par\smallskip\noindent\(
			\GroupLayout{\pi_{27}(\mathit{SO}(20)) \;= {} \mathrm{span}\{[256P(\iota_{29})]_{20},\; [\eta^2_{18}]_{20}\sigma_{20},\; [\varepsilon_{18}]_{20}\eta_{26},\; [\eta_{19}]\sigma_{20},\; [\varepsilon_{19}]\}\cong \mathbb Z_{(2)}\oplus (\mathbb Z/2)^{4}.}{\begin{aligned}\pi_{27}(\mathit{SO}(20)) \;= {} \mathrm{span}\{&[256P(\iota_{29})]_{20},\; [\eta^2_{18}]_{20}\sigma_{20},\\
					&[\varepsilon_{18}]_{20}\eta_{26},\; [\eta_{19}]\sigma_{20},\; [\varepsilon_{19}]\}\cong \mathbb Z_{(2)}\oplus (\mathbb Z/2)^{4}.\end{aligned}}
			\)\par
			\smallskip
			
			\par\par\smallskip\noindent\(
			\GroupLayout{\pi_{27}(\mathit{SO}(21)) = \mathrm{span}\{[256P(\iota_{29})]_{21},\; [\eta_{19}]_{21}\sigma_{20},\; [\varepsilon_{19}]_{21}\}\cong \mathbb Z_{(2)}\oplus (\mathbb Z/2)^{2}.}{\begin{aligned}\pi_{27}(\mathit{SO}(21)) = \mathrm{span}\{&[256P(\iota_{29})]_{21},\; [\eta_{19}]_{21}\sigma_{20},\\
					&[\varepsilon_{19}]_{21}\}\cong \mathbb Z_{(2)}\oplus (\mathbb Z/2)^{2}.\end{aligned}}
			\)\par
			\smallskip\noindent\(
			\GroupLayout{\pi_{27}(\mathit{SO}(22)) = \mathrm{span}\{[256P(\iota_{29})]_{22},\; [\nu^2_{21}]\}\cong \mathbb Z_{(2)}\oplus \mathbb Z/4.}{\begin{aligned}\pi_{27}(\mathit{SO}(22)) = \mathrm{span}\{&[256P(\iota_{29})]_{22},\\
					&[\nu^2_{21}]\}\cong \mathbb Z_{(2)}\oplus \mathbb Z/4.\end{aligned}}
			\)\par
			\smallskip\noindent\(
			\GroupLayout{\pi_{27}(\mathit{SO}(23)) = \mathrm{span}\{[256P(\iota_{29})]_{23},\; [\nu^2_{21}]_{23}\}\cong \mathbb Z_{(2)}\oplus \mathbb Z/2.}{\begin{aligned}\pi_{27}(\mathit{SO}(23)) = \mathrm{span}\{&[256P(\iota_{29})]_{23},\\
					&[\nu^2_{21}]_{23}\}\cong \mathbb Z_{(2)}\oplus \mathbb Z/2.\end{aligned}}
			\)\par
			\smallskip\noindent\(
			\GroupLayout{\pi_{27}(\mathit{SO}(24)) = \mathrm{span}\{[256P(\iota_{29})]_{24},\; [\eta_{23}]\nu_{24}\}\cong \mathbb Z_{(2)}\oplus \mathbb Z/2.}{\begin{aligned}\pi_{27}(\mathit{SO}(24)) = \mathrm{span}\{&[256P(\iota_{29})]_{24},\\
					&[\eta_{23}]\nu_{24}\}\cong \mathbb Z_{(2)}\oplus \mathbb Z/2.\end{aligned}}
			\)\par
			\smallskip\noindent\(
			\GroupLayout{\pi_{27}(\mathit{SO}(25)) = \mathrm{span}\{[256P(\iota_{29})]_{25},\; \Delta(\nu_{25})\}\cong \mathbb Z_{(2)}\oplus \mathbb Z/2.}{\begin{aligned}\pi_{27}(\mathit{SO}(25)) = \mathrm{span}\{&[256P(\iota_{29})]_{25},\\
					&\Delta(\nu_{25})\}\cong \mathbb Z_{(2)}\oplus \mathbb Z/2.\end{aligned}}
			\)\par
			\par\smallskip\noindent\(
			\FitMath{\pi_{27}(\mathit{SO}(26+n)) = \mathrm{span}\{[256P(\iota_{29})]_{26+n}\}\cong \mathbb Z_{(2)},\qquad n=0,1\;\; \text{or}\;\; n\geq3.}
			\)\par
			\smallskip\noindent\(
			\GroupLayout{\pi_{27}(\mathit{SO}(28)) = \mathrm{span}\{[256P(\iota_{29})]_{28},\; \Delta(\iota_{28})\}\cong  \mathbb Z_{(2)}\jia\z_{(2)}.}{\begin{aligned}\pi_{27}(\mathit{SO}(28)) = \mathrm{span}\{&[256P(\iota_{29})]_{28},\\
					&\Delta(\iota_{28})\}\cong \mathbb Z_{(2)}\jia\z_{(2)}.\end{aligned}}
			\)\par
		\end{itemize}
	\end{thm}
	\endgroup
	\section{The Cases $3\leq n\leq 14$}\label{m2:subsec:degree27-low}
	
	\subsection{The Cases $3\leq n\leq6$}
	
	The initial calculations for $n=3,4$ use the homotopy groups of $S^3.$  The subsequent exact
	sequences compare these classes with the generators for $n=5,6,$   providing the input for the
	next set of cases.
	
	\begin{prop}\label{m2:2734}
		The degree-$27$ groups for $n=3,4$ are given as follows.
		\[\begin{gathered}
			\GroupLayout{\pi_{27}(\mathit{SO}(3)) = \mathrm{span}\{[\eta_2]\delta_3,\; [\eta_2]\bar\mu_3\sigma_{20},\; [\eta_2]\varepsilon'\kappa_{13}\}\cong (\mathbb Z/2)^{3},}{\begin{aligned}\pi_{27}(\mathit{SO}(3)) = \mathrm{span}\{&[\eta_2]\delta_3,\; [\eta_2]\bar\mu_3\sigma_{20},\\
					&[\eta_2]\varepsilon'\kappa_{13}\}\cong (\mathbb Z/2)^{3},\end{aligned}}\\[3pt]
			\GroupLayout{\pi_{27}(\mathit{SO}(4)) \;= {} \mathrm{span}\{[\eta_2]_4\delta_3,\; [\eta_2]_4\bar\mu_3\sigma_{20},\; [\eta_2]_4\varepsilon'\kappa_{13},\; [\iota_3]\delta_3,\; [\iota_3]\bar\mu_3\sigma_{20},\; [\iota_3]\varepsilon'\kappa_{13}\}\cong (\mathbb Z/2)^{6}.}{\begin{aligned}\pi_{27}(\mathit{SO}(4)) \;= {} \mathrm{span}\{&[\eta_2]_4\delta_3,\; [\eta_2]_4\bar\mu_3\sigma_{20},\; [\eta_2]_4\varepsilon'\kappa_{13},\\
					&[\iota_3]\delta_3,\; [\iota_3]\bar\mu_3\sigma_{20},\; [\iota_3]\varepsilon'\kappa_{13}\}\cong (\mathbb Z/2)^{6}.\end{aligned}}
		\end{gathered}\]
		For each $g\in\{\delta_3,\bar\mu_3\sigma_{20},\varepsilon'\kappa_{13}\},$   we have
		\[\begin{aligned}
			&2([\eta_2]g)=0,\qquad 2([\eta_2]_4g)=0,\qquad 2([\iota_3]g)=0.
		\end{aligned}\]
		The $n=4$ decomposition is the splitting induced by $i_*+[\iota_3]_*.$\qed 
	\end{prop}

	With the degree-$27$ groups for $n=3,4$ in hand, we proceed to $n=5.$ The exact sequence
	over $S^4$ determines the surviving images and the additional class involving $\nu_{10}\kappa_{13}.$
	
	\begin{prop}\label{m2:275}
		The structure of $\pi_{27}(\mathit{SO}(5))$ is as follows.
		\[ \GroupLayout{\pi_{27}(\mathit{SO}(5)) = \mathrm{span}\{[\iota_3]_5\delta_3,\; [\iota_3]_5\bar\mu_3\sigma_{20},\; [\nu_4^2]\nu_{10}\kappa_{13}\}\cong (\mathbb Z/2)^{3}.}{\begin{aligned}\pi_{27}(\mathit{SO}(5)) = \mathrm{span}\{&[\iota_3]_5\delta_3,\; [\iota_3]_5\bar\mu_3\sigma_{20},\\
				&[\nu_4^2]\nu_{10}\kappa_{13}\}\cong (\mathbb Z/2)^{3}.\end{aligned}} \]
		
	\end{prop}
	
	\begin{proof}
		Consider  the exact sequence $(27)_4$:
		$$\begin{aligned}
			&\pi^4_{28}\stackrel{{\Delta}}{\longrightarrow}\pi_{27}(\mathit{SO}(4))\stackrel{{i_*}}{\longrightarrow}\pi_{27}(\mathit{SO}(5))\stackrel{{p_*}}{\longrightarrow}\pi^4_{27}\stackrel{{\Delta}}{\longrightarrow}\pi_{26}(\mathit{SO}(4)).
		\end{aligned}$$\noindent    
		The boundary calculation in the proof of Proposition~\ref{m2:526} gives
		\[\FitMath{\begin{aligned}
				\pi_{27}^4&=\frac{\mathbb Z}{8}\{\nu_4\bar\kappa_7\}\oplus\frac{\mathbb Z}{4}\{\Sigma\bar\alpha\},\\
				\operatorname{Im}\Delta&=\pi_{26}(\mathit{SO}(4))\cong(\mathbb Z/4)^2,\qquad
				|\ker\Delta|=2,\\
				\ker\bigl(\Delta\colon \pi_{27}^4\to\pi_{26}(\mathit{SO}(4))\bigr)
				&=\mathrm{span}\{4\nu_4\bar\kappa_7\}
				=\mathrm{span}\{\nu_4^3\kappa_{13}\}\cong\mathbb Z/2.
		\end{aligned}}\]
		Here $4\nu_4\bar\kappa_7=\nu_4^3\kappa_{13}$ by \cite[Lemma~15.4]{MT1}.
		Note that $$\pi^4_{28}=\mathrm{span}\{\delta_4, \cq _4\sigma_{21},
		\Sigma\varepsilon'\kappa_{14}, \nu_4\sigma'\kappa_{14},
		\nu_4\eta_7\ct _8\}\cong(\mathbb Z/2)^5.$$ \noindent In light of the relation
		$\Delta(\iota_4)=2[\iota_3]-[\eta_2]_4$,
		we see that $$\Delta(\delta_4)=[\eta_2]_4\delta_3,\
		\Delta(\cq _4\sigma_{21})=[\eta_2]_4\cq _3\sigma_{20},\ \text{ and} \
		\Delta(\Sigma\varepsilon'\kappa_{14})=[\eta_2]_4\varepsilon'\kappa_{13}.$$ By \cite[Lemma 1.4]{KM2}, we
		obtain
		$\Delta(\nu_4\sigma'\kappa_{14})=0.$
		According to \cite[Lemma 1.1]{KM1} and \cite[Proposition 3.1 (1)]{Oda1979}, we have
		$$\Delta(\nu_4\eta_7\ct _8) =([\iota_3]+a[\eta_2]_4)\nu'\eta_6\ct _7
		=([\iota_3]+a[\eta_2]_4)\varepsilon'\kappa_{13}\ (0\leq a \leq 3).$$ Since $\nu_{10}\kappa_{13}$ has order 2, we find
		$2[\nu_4^2]\nu_{10}\kappa_{13}=0.$
		This proves the proposition.
	\end{proof}
	We need the follwing lemma.
	\begin{lem}The relation $\sigma''\nu_{13}=\pm2\nu_6\sigma_9$ in Appendix~\ref{app:oguchi-relations} permits the choice of an odd integer $x_3\in\z$ such that
		\begin{equation}
			\begin{aligned}
				&\sigma''\nu_{13}=2x_3\nu_6\sigma_9.\label{m2:s''n}
			\end{aligned}
		\end{equation}
		\[\begin{aligned}
			&[\eta^2_5]\nu_7=2x_4[  \nu^2_4]_6\ \text{for an odd integer }x_4\in\z.
		\end{aligned}\]
	\end{lem}
	\begin{proof}
		Since $[\eta^2_5]\nu_7\in{i_6}_*\pi_{10}(\mathit{SO}(5))=\frac{\mathbb Z}{8}\{[  \nu^2_4]_6\}$ \cite[Table
		2]{Ka}, we see that $[\eta^2_5]\nu_7=y[  \nu^2_4]_6$ for an integer $y.$
		We obtain $J[\nu_4^2]=\nu_5\sigma_8$ \cite[Lemma 1.1(i)]{KM1} and $J[\eta_5^2]=\sigma''$ \cite[Lemma
		2.1(i)]{HKM}. The inclusion formula in Section~\ref{m2:subsec:j-toda} gives $J([\nu_4^2]_6)=\pm\nu_6\sigma_9.$ 
		Consequently, \eqref{m2:s''n} gives $y\equiv\pm2x_3\pmod8,$   since $\nu_6\sigma_9$ has order $8$ by Appendix~\ref{tlqdzj}.
		There is therefore an odd integer $x_4\in\z$ with $y\equiv2x_4\pmod8$ and $x_4\equiv\pm x_3\pmod4.$ 
		Since $[\nu_4^2]_6$ has order $8,$   this proves the asserted equality.
		This completes the proof.
	\end{proof}
	
	The degree-$27$ extension at $n=6$ is determined by two doubling relations. Recording them
	alongside the decomposition explains how the generators from $n=5$ appear inside the new cyclic
	summands.
	
	\begin{prop}\label{m2:276}
		The group $\pi_{27}(\mathit{SO}(6))$ is as follows.
		\[\GroupLayout{\pi_{27}(\mathit{SO}(6)) = \mathrm{span}\{[\eta_5^2]\bar\kappa_7,\; [\rho^{\mathrm{IV}}]\sigma_{20},\; [ \nu_5]\bar\sigma_8\}\cong (\mathbb Z/4)^{2}\oplus \mathbb Z/2.}{\begin{aligned}\pi_{27}(\mathit{SO}(6)) = \mathrm{span}\{&[\eta_5^2]\bar\kappa_7,\; [\rho^{\mathrm{IV}}]\sigma_{20},\\
				&[ \nu_5]\bar\sigma_8\}\cong (\mathbb Z/4)^{2}\oplus \mathbb Z/2.\end{aligned}}\]
		The new lifts satisfy
		\[
		\begin{aligned}
			&2([\eta_5^2]\bar\kappa_7)=[\nu_4^2]_6\nu_{10}\kappa_{13},\qquad 2([\rho^{\mathrm{IV}}]\sigma_{20})=[\iota_3]_6\bar\mu_3\sigma_{20}, \\
			&2([ \nu_5]\bar\sigma_8)=0.
		\end{aligned}
		\]
	\end{prop}
	
	\begin{proof}
		Consider  the exact sequence $(27)_5$:
		$$\pi^5_{28}\stackrel{{\Delta}}{\longrightarrow}\pi_{27}(\mathit{SO}(5))\stackrel{{i_*}}{\longrightarrow}\pi_{27}(\mathit{SO}(6))\stackrel{{p_*}}{\longrightarrow}\pi^5_{27}\stackrel{{\Delta}}{\longrightarrow}\pi_{26}(\mathit{SO}(5)).$$ 
		In light of the proof of Proposition~\ref{m2:626},
		\[\FitMath{\Delta(\nu_5\bar\sigma_8)=\Delta(\nu_5\bar\zeta_8)=\Delta(\varepsilon_5\kappa_{13})=0.}\]
		Thus
		\[\FitMath{\begin{aligned}
				\ker\bigl(\Delta\colon \pi_{27}^5\to\pi_{26}(\mathit{SO}(5))\bigr)
				&=\pi_{27}^5\\
				&=\mathrm{span}\{\nu_5\bar\sigma_8,\nu_5\bar\zeta_8,\varepsilon_5\kappa_{13}\}\\
				&=\mathrm{span}\{\nu_5\bar\sigma_8,\rho^{\mathrm{IV}}\sigma_{20},\eta_5^2\bar\kappa_7\}\cong(\mathbb Z/2)^3.
		\end{aligned}}\]
		The alternative generators follow from \cite[Proposition~3.1(2)]{Oda1979} and \cite[(2.4)]{MMO}.
		
		We find
		$\pi^5_{28}=\mathrm{span}\{\nu_5\ct _8, \overline{\rho}''', \phi_5\}\cong\mathbb Z/8\oplus(\mathbb Z/2)^2.$ 
		Using the relation $\Delta(\nu_5)=0$ \cite[Lemma 1.1(i)]{KM1},
		we have $\Delta(\nu_5\ct _8)=0.$
		By \cite[(3.2)]{MMO}, \cite[p.94]{KM2} and the fact that
		$\pi_{20}(\mathit{SO}(5))\cong(\mathbb Z/2)^3$ \cite[Theorem 3.2]{KM2}, we see that
		$$\Delta(\overline{\rho}''')\in \{\Delta(\rho^{\mathrm{IV}}), 2\iota_{19},
		8\sigma_{19}\}=\{0,2\iota_{19},8\sigma_{19}\}=0.$$ Here the bracket has indeterminacy
		$\pi_{20}(\mathit{SO}(5))\circ8\sigma_{20}=0.$ 
		Thus $\Delta(\overline{\rho}''')=0.$
		By \cite[II-Proposition 4.4 (3)]{Oda1979},
		we obtain $J\Delta(\phi_5)=P(\phi_{11})\equiv\nu_5\delta_8\bmod \nu_5\cq _8\sigma_{25}.$
		Using the relation $J[\iota_3]=\nu_4$ \cite[(2.2)]{HKM}, we deduce $J([\iota_3]_5\delta_3)=\nu_5\delta_8$ and
		$J([\iota_3]_5\cq _3\sigma_{20})=\nu_5\cq _8\sigma_{25}.$
		In view of  the relation $J[\nu_4^2]=\nu_5\sigma_8$ \cite[Lemma 1.1(i)]{KM1} and
		$\sigma'''\ct _{12}=\nu_5\sigma_8\nu_{15}\kappa_{18}$ \cite[III-Proposition 2.4 (2)]{Oda1979}, we find
		$J([\nu_4^2]\nu_{10}\kappa_{13})=\sigma'''\ct _{12}.$
		Thus
		$\Delta(\phi_5)=[\iota_3]_5\delta_3 \bmod [\iota_3]_5\cq _3\sigma_{20}.$
		
		Employing the relation
		$\nu^2_7\kappa_{13}=4\bar{\kappa}_7$ \cite[p. 48]{MT1},
		$\sigma''\eta_{13}=4\bar\nu_6,$   and
		$\pi_{12}(\mathit{SO}(6))=\frac{\mathbb Z}{4}\{\Delta\sigma''\},$   we obtain\vspace{-1em}\begin{align*}
			2[\eta_5^2]\bar\kappa_7
			&= [\eta_5^2](\Sigma\bar\kappa') \notag\\
			&\in [\eta_5^2]\circ\{\nu_7,\eta_{10},\bar\varepsilon_{11}\}_1 \notag\\
			&\subseteq \{[\eta_5^2]\nu_7,\eta_{10},\bar\varepsilon_{11}\}_1 \notag\\
			&= \{2x_4[\nu_4^2]_6,\eta_{10},\bar\varepsilon_{11}\}_1 \notag\\
			&\supseteq x_4[\nu_4^2]_6\circ \Sigma\{2\iota_9,\eta_9,\bar\varepsilon_{10}\} \notag\\
			&\ni x_4[\nu_4^2]_6\nu_{10}\kappa_{13}. \notag
		\end{align*}
		The comparison is modulo
		\(\begin{aligned}
			&[\eta_5^2]\nu_7\circ \Sigma\pi_{26}^9+\pi_{12}(\mathit{SO}(6))\circ\bar\varepsilon_{12}={}\mathrm{span}\{4[\eta_5^2]\bar\kappa_7\}.
		\end{aligned}\)
		The   element   $[\nu_4^2]_6\nu_{10}\kappa_{13}$ has order at most $2,$   so multiplication by the odd integer $x_4$ leaves it unchanged.
		The displayed congruence has the form $2[\eta_5^2]\bar\kappa_7\equiv[\nu_4^2]_6\nu_{10}\kappa_{13}$ modulo $4[\eta_5^2]\bar\kappa_7.$ 
		Doubling shows that an odd multiple of $4[\eta_5^2]\bar\kappa_7$ is zero; the group is $2$-primary, so $4[\eta_5^2]\bar\kappa_7=0.$ 
		Thus we conclude that
		$2[\eta_5^2]\ct _7=[\nu_4^2]_6\nu_{10}\kappa_{13}.$
		Using the relation $2[\rho^{\mathrm{IV}}]=[\iota_3]_6\bar{\mu}_3$ \cite[Theorem 3.2]{KM2},
		we see that
		$2[\rho^{\mathrm{IV}}]\sigma_{20}=[\iota_3]_6\cq _3\sigma_{20}.$
		This completes the proof.
	\end{proof}
	
	We next use the fibration over $S^6$ to pass to $n=7.$ The lift involving $\bar\kappa_7$ accounts
	for the principal cyclic summand, and its double is compared with an inherited class.
	
	\subsection{The Cases $7\leq n\leq9$}
	
	The auxiliary $G_2$ information and the contribution from the homotopy groups of spheres now enter the calculation for $n=7.$  We
	retain the chosen representatives while passing through the cases $n=8,9,$   so their subsequent
	images remain identifiable.
	
	\begin{prop}\label{m2:277}
		The group $\pi_{27}(\mathit{SO}(7))$ is as follows.
		\[\GroupLayout{\pi_{27}(\mathit{SO}(7)) = \mathrm{span}\{[\eta_6]\bar\kappa_7,\; [ \nu_5]_7\bar\sigma_8\}\cong \mathbb Z/8\oplus \mathbb Z/2.}{\begin{aligned}\pi_{27}(\mathit{SO}(7)) = \mathrm{span}\{&[\eta_6]\bar\kappa_7,\\
				&[ \nu_5]_7\bar\sigma_8\}\cong \mathbb Z/8\oplus \mathbb Z/2.\end{aligned}}\]
		The extension relation is
		$2([\eta_6]\bar\kappa_7)=[\eta_5^2]_7\bar\kappa_7.$
	\end{prop}
	
	\begin{proof}
		Consider  the exact sequence $(27)_6$:
		$$\begin{aligned}
			&\pi^6_{28}\stackrel{{\Delta}}{\longrightarrow}\pi_{27}(\mathit{SO}(6))\stackrel{{i_*}}{\longrightarrow}\pi_{27}(\mathit{SO}(7))\stackrel{{p_*}}{\longrightarrow}\pi^6_{27}\stackrel{{\Delta}}{\longrightarrow}\pi_{26}(\mathit{SO}(6)).
		\end{aligned}$$\noindent    
		By the proof of Proposition~\ref{m2:726},
		\[\FitMath{\pi_{27}^6=\frac{\mathbb Z}{2}\{\eta_6\bar\kappa_7\},\qquad\Delta(\eta_6\bar\kappa_7)=0.}\]
		Consequently,
		\[\FitMath{\ker\bigl(\Delta\colon \pi_{27}^6\to\pi_{26}(\mathit{SO}(6))\bigr)=\frac{\mathbb Z}{2}\{\eta_6\bar\kappa_7\}.}\]\noindent      
		We have $\pi^6_{28}=\mathrm{span}\{\rho'''\sigma_{21}, \nu_6\ccs_9, \bar{\nu}_6\kappa_{14}, \varepsilon_6\kappa_{14}\} \cong\mathbb Z/4\oplus(\mathbb Z/2)^3.$ 
		By \cite[p.106]{KM2}, we deduce
		$\Delta(\rho'''\sigma_{21})=[\rho^{\mathrm{IV}}]\sigma_{20}$.
		By \cite[Table 3]{Ka}, we find
		$\Delta(\nu_6\ccs_9)=2[ \nu_5]\ccs_8=0.$
		Using the relations $\Delta(\bar{\nu}_6)=\Delta(\varepsilon_6)=2[\nu_4^2]_6\nu_{10}$ and $2\kappa_{13}=0,$ we conclude that
		$\Delta(\bar{\nu}_6\kappa_{14})=\Delta(\varepsilon_6\kappa_{14})=0.$
		Employing the relation $2[\eta_6]=[\eta_5^2]_7,$ we see that
		$2[\eta_6]\ct _7=[\eta_5^2]_7\ct _7.$
		This proves the proposition.
	\end{proof}
	
	The homomorphism
	$(i_8)_*+([\iota_7])_*:\pi_{27}(\mathit{SO}(7))\oplus\pi_{27}^7\to \pi_{27}(\mathit{SO}(8))$
	is an isomorphism.  Here $([\iota_7])_*(\alpha)=[\iota_7]\circ\alpha.$ 
	Applying it to the generators gives the following decomposition.
	
	\begin{prop}\label{m2:278}
		The structure of $\pi_{27}(\mathit{SO}(8))$ is as follows.
		\[\GroupLayout{\pi_{27}(\mathit{SO}(8)) = \mathrm{span}\{[\iota_7]\bar\kappa_7,\; [\eta_6]_8\bar\kappa_7,\; [ \nu_5]_8\bar\sigma_8\}\cong (\mathbb Z/8)^{2}\oplus \mathbb Z/2.}{\begin{aligned}\pi_{27}(\mathit{SO}(8)) = \mathrm{span}\{&[\iota_7]\bar\kappa_7,\; [\eta_6]_8\bar\kappa_7,\\
				&[ \nu_5]_8\bar\sigma_8\}\cong (\mathbb Z/8)^{2}\oplus \mathbb Z/2.\end{aligned}}\]
		The decomposition is induced by the isomorphism
		$(i_8)_*+([\iota_7])_*:\pi_{27}(\mathit{SO}(7))\oplus\pi_{27}^7\to \pi_{27}(\mathit{SO}(8)).$ 	\qed
	\end{prop}
	
	By the relation $J[\iota_7]=\sigma_8$ \cite[Lemma 1.2(ii)]{KM1}, we obtain
	$J([\iota_7]_n\bar{\kappa}_7)=\pm\sigma_n\bar{\kappa}_{n+7}\ \text{for}\ n\geq 8.$
	The order of $[\iota_7]_n\bar\kappa_7$ divides $8,$   since
	$8\bar\kappa_7=0.$   Its $J$-image has order $8$ for
	$8\leq n\leq11$ by \cite[Theorem 2(a)]{Oda1979}.  Therefore
	$\mathrm{ord}([\iota_7]_n\bar{\kappa}_7)=8\ \text{for}\ 8\leq n\leq 11.$
	
	The preceding order calculation identifies a class that survives the next inclusion. We now determine
	the full group at $n=9$ by computing the boundary of $\bar\kappa_8.$
	
	\begin{prop}\label{m2:279}
		The structure of $\pi_{27}(\mathit{SO}(9))$ is as follows.
		\[\GroupLayout{\pi_{27}(\mathit{SO}(9)) = \mathrm{span}\{[\iota_7]_9\bar\kappa_7,\; [ \nu_5]_9\bar\sigma_8\}\cong \mathbb Z/8\oplus \mathbb Z/2.}{\begin{aligned}\pi_{27}(\mathit{SO}(9)) = \mathrm{span}\{&[\iota_7]_9\bar\kappa_7,\\
				&[ \nu_5]_9\bar\sigma_8\}\cong \mathbb Z/8\oplus \mathbb Z/2.\end{aligned}}\]
		The boundary and the induced relation are
		\[
		\begin{aligned}
			&\Delta\bar\kappa_8=2[\iota_7]\bar\kappa_7-[\eta_6]_8\bar\kappa_7,\qquad 2[\iota_7]_9\bar\kappa_7=[\eta_6]_9\bar\kappa_7.
		\end{aligned}
		\]
	\end{prop}
	
	\begin{proof}
		Consider  the exact sequence $(27)_8$:
		$$\pi^8_{28}\stackrel{{\Delta}}{\longrightarrow}\pi_{27}(\mathit{SO}(8))\stackrel{{i_*}}{\longrightarrow}\pi_{27}(\mathit{SO}(9))\stackrel{{p_*}}{\longrightarrow}\pi^8_{27}\stackrel{{\Delta}}{\longrightarrow}\pi_{26}(\mathit{SO}(8)).$$ 
		In light of Proposition~\ref{g1016:prop:26-9},
		\[\FitMath{\begin{aligned}
				\pi_{27}^8&=\frac{\mathbb Z}{8}\{\bar\zeta_8\}\oplus\frac{\mathbb Z}{2}\{\bar\sigma_8\},\\
				\Delta\bar\zeta_8&\ty2[\iota_7]\bar\zeta_7-[\nu_4\bar\zeta_7]_8\m8[[2\e_{6},\e_{6}]]_8\s^{2}\rho',\\
				\Delta\bar\sigma_8&=[\nu_4\bar\sigma_7]_8.
		\end{aligned}}\]
		The direct sum in Proposition~\ref{m2:826} gives
		\[\FitMath{\begin{aligned}
				\operatorname{Im}\Delta=\frac{\mathbb Z}{8}\{\Delta\bar\zeta_8\}\oplus\frac{\mathbb Z}{2}\{\Delta\bar\sigma_8\},\;\;\;
				\ker\bigl(\Delta\colon \pi_{27}^8\to\pi_{26}(\mathit{SO}(8))\bigr)=0.
		\end{aligned}}\]\noindent       We note that
		$\pi^8_{28}=\mathrm{span}\{\ct _8\}\cong\mathbb Z/8.$ 
		By the relation $\Delta(\iota_8)=2[\iota_7]-[\eta_6]_8$ \cite[Table 3]{Ka}, we find
		$\Delta(\ct _8)=\Delta(\iota_8)\ct _7=2[\iota_7]\ct _7-[\eta_6]_8\ct _7.$
		This proves the proposition.
	\end{proof}
	
	At $n=10,$ a lift of $\sigma_9\zeta_{16}$ joins the inherited $\bar\kappa$-class. The key extension
	question is whether this lift has order $8\uprightsemicolon{}$ the proof uses its $J$-image to settle that point.
	
	\subsection{The Cases $n=10,11$}
	
	For $n=10,11,$   the order of each chosen lift must be compared with its image under
	projection. The following decomposition records the representatives used for this comparison.
	
	\begin{prop}\label{m2:2710}
		The structure of $\pi_{27}(\mathit{SO}(10))$ is as follows.
		\[ \GroupLayout{\pi_{27}(\mathit{SO}(10)) = \mathrm{span}\{[\sigma_9\zeta_{16}],\; [\iota_7]_{10}\bar\kappa_7\}\cong (\mathbb Z/8)^{2}.}{\begin{aligned}\pi_{27}(\mathit{SO}(10)) = \mathrm{span}\{&[\sigma_9\zeta_{16}],\\
				&[\iota_7]_{10}\bar\kappa_7\}\cong (\mathbb Z/8)^{2}.\end{aligned}} \]
	\end{prop}
	
	\begin{proof}
		Consider  the exact sequence $(27)_9$:
		\[\begin{aligned}
			&\pi^9_{28}\stackrel{{\Delta}}{\longrightarrow}\pi_{27}(\mathit{SO}(9))\stackrel{{i_*}}{\longrightarrow}\pi_{27}(\mathit{SO}(10))\stackrel{{p_*}}{\longrightarrow}\pi^9_{27}\stackrel{{\Delta}}{\longrightarrow}\pi_{26}(\mathit{SO}(9)).
		\end{aligned}\]
		\noindent       The proof of Proposition~\ref{m2:10} gives
		\[\FitMath{\begin{aligned}
				\pi_{27}^9&=\frac{\mathbb Z}{8}\{\sigma_9\zeta_{16}\}\oplus\frac{\mathbb Z}{2}\{\bar\mu_9\eta_{26}\},\\
				\Delta(\sigma_9\zeta_{16})&=0,\qquad\Delta(\bar\mu_9\eta_{26})\ne0.
		\end{aligned}}\]
		Hence $\ker\bigl(\Delta\colon \pi_{27}^9\to\pi_{26}(\mathit{SO}(9))\bigr)=\frac{\mathbb Z}{8}\{\sigma_9\zeta_{16}\}.$     Note that
		$\pi^9_{28}=\mathrm{span}\{\overline{\zeta}_9, \ccs_9\}\cong\mathbb Z/8\oplus\mathbb Z/2.$ 
		Using the relation $\Delta(\iota_9)=[ \nu_5]_9+[\iota_7]_9\eta_7$
		\cite[Table 3]{Ka} and $\eta_6\bar{\sigma}_7=0$ \cite[(2.7)]{MMO}, we see that $$\Delta(\ccs_9)=\Delta(\iota_9)\ccs_8 =[ \nu_5]_9\ccs_8+[\iota_7]_9\eta_7\ccs_8=[ \nu_5]_9\ccs_8.$$ By the definition of \ $\overline{\zeta}_5$ \cite[p.137]{T} and \cite[p.96]{KM2}, we obtain $\Delta(\overline{\zeta}_9)\in \{\Delta(\zeta_9), 8\iota_{19}, 2\sigma_{19}\}=\{0, 2\iota_{19}, 8\sigma_{19}\}=\pi_{20}(\mathit{SO}(9))\circ 8\sigma_{20}.$ Since $\pi_{20}(\mathit{SO}(9))\cong\mathbb Z/2$ \cite[Theorem 2.1]{KM2},
		we deduce
		$\Delta(\overline{\zeta}_9)=0.$
		By \cite[Lemma 1.1(ii)]{KM1}, we find
		\begin{equation}
			\begin{aligned}
				&J([\iota_7]_{10}\ct _7)=\pm\sigma_{10}\ct _{17}.\label{m2:J2710}
			\end{aligned}
		\end{equation}
		By  using the relations $$8[\sigma_9\zeta_{16}]\in
		i_*\pi_{27}(\mathit{SO}(9))=\mathrm{span}\{[\iota_7]_{10}\ct _7\},\;\;J(8[\sigma_9\zeta_{16}])=0\in\pi^{10}_{37}\cong \mathbb Z/4\oplus(\mathbb Z/8)^2$$ \noindent and
		(\ref{m2:J2710}), whose right-hand side has order $8,$   we see that $J$ is injective on this cyclic subgroup. Hence we infer
		$8[\sigma_9\zeta_{16}]=0.$
		This proves the proposition.
	\end{proof}
	
	For $n=11,$ we follow the inherited $\bar\kappa$-class and examine a lift of
	$\eta_{10}\mu_{11}\sigma_{20}.$ The next proposition identifies the two contributions and their
	orders.
	
	\begin{prop}\label{m2:2711}
		The structure of $\pi_{27}(\mathit{SO}(11))$ is as follows.
		\[ \GroupLayout{\pi_{27}(\mathit{SO}(11)) = \mathrm{span}\{[\iota_7]_{11}\bar\kappa_7,\; [\eta_{10}\mu_{11}]\sigma_{20}\}\cong \mathbb Z/8\oplus \mathbb Z/2.}{\begin{aligned}\pi_{27}(\mathit{SO}(11)) = \mathrm{span}\{&[\iota_7]_{11}\bar\kappa_7,\\
				&[\eta_{10}\mu_{11}]\sigma_{20}\}\cong \mathbb Z/8\oplus \mathbb Z/2.\end{aligned}} \]
	\end{prop}
	
	\begin{proof}
		Consider  the exact sequence $(27)_{10}$:
		\[\begin{aligned}
			&\pi^{10}_{28}\stackrel{{\Delta}}{\longrightarrow}\pi_{27}(\mathit{SO}(10))\stackrel{{i_*}}{\longrightarrow}\pi_{27}(\mathit{SO}(11))\stackrel{{p_*}}{\longrightarrow}\pi^{10}_{27}\stackrel{{\Delta}}{\longrightarrow}\pi_{26}(\mathit{SO}(10)).
		\end{aligned}\]
		\noindent  By Lemma~\ref{m2:Dnuka}(2), Proposition~\ref{m2:10}, and \eqref{m2:Det10}, we have 
		\[\FitMath{\begin{aligned}
				\pi_{27}^{10}&=\mathrm{span}\{\nu_{10}\kappa_{13},\bar\mu_{10},\eta_{10}\mu_{11}\sigma_{20}\}\cong(\mathbb Z/2)^3,
		\end{aligned}}\]
		$$\Delta(\nu_{10}\kappa_{13})=[\nu_6^2\kappa_{12}]_{10}=2[\nu_9\kappa_{12}],\;
		\Delta\bar\mu_{10}=4[\eta_9\mu_{10}\sigma_{19}],\;
		\Delta(\eta_{10}\mu_{11}\sigma_{20})=0.$$
		These images satisfy $\mathrm{span}\{\Delta(\nu_{10}\kappa_{13}),4[\eta_9\mu_{10}\sigma_{19}]\}\cong(\mathbb Z/2)^2.$  In the ordered decomposition of Proposition~\ref{m2:10}, their projections onto the cyclic summands generated by $[\nu_9\kappa_{12}]$ and $[\eta_9\mu_{10}\sigma_{19}]$ are twice the former generator and four times the latter generator, respectively. Both images have order $2,$   and these projections show that they are independent. Therefore
		\[\FitMath{\ker\bigl(\Delta\colon \pi_{27}^{10}\to\pi_{26}(\mathit{SO}(10))\bigr)=\frac{\mathbb Z}{2}\{\eta_{10}\mu_{11}\sigma_{20}\}.}\]
		
		\noindent       Recall that
		$$\pi^{10}_{28}=\mathrm{span}\{\xi'', \xi''\pm\lambda'', \eta_{10}\cq _{11}\}\cong\mathbb Z/8\oplus(\mathbb Z/2)^2.$$ 
		Using the relation $\Delta(\eta_{10})=2[\iota_7]_{10}\nu_7$ \cite[p.336]{KM1}, we see that
		$\Delta(\eta_{10}\cq _{11})=2[\iota_7]_{10}\nu_7\cq _{10}=0.$
		Recall  the relation $4\xi''=2\sigma_{10}\zeta_{17}$ from \cite{T}. 
		Then by \cite[(4.7)]{Ka} and the fact that $4\sigma_9\zeta_{16}\neq 0,$ we deduce
		$$\Delta(4\xi'')=\Delta(\sigma_{10})2\zeta_{16}=2[2\sigma_9]\zeta_{16}\neq 0.$$ \noindent 
		Hence $\Delta(\xi'')$ has order 8.
		In light of \cite[Lemma 12.19, Corollary 12.25]{T}, we find $$J\Delta(\xi'')=P(4\xi_{21})=P(4\nu^*_{21})=0.$$ According to (\ref{m2:J2710}), $J$ is injective on $\mathrm{span}\{[\iota_7]_{10}\ct _7\}.$ 
		Thus the induced homomorphism $$\ker (J\colon \pi_{27}(\mathit{SO}(10))\to \pi_{37}^{10})\hookrightarrow  \frac{\pi_{27}(\mathit{SO}(10))}{\mathrm{span}\{[\iota_7]_{10}\ct _7\}}:=W\cong\z/8$$ is injective.
		The quotient module $W$ is   generated by the   class of $[\sigma_9\zeta_{16}].$ 
		Since $\Delta(\xi'')\in\ker J$ has order $8,$   its image is a generator of $W$. Consequently,
		$$\Delta(\xi'')\equiv x_5[\sigma_9\zeta_{16}]\bmod[\iota_7]_{10}\ct _7$$ for some odd integer $x_5\in\z.$ 
		Moreover, $\ker J=\mathrm{span}\{\Delta(\xi'')\},$ in view of  $|\ker J|\leq |W|=8$. 
		By \cite[Lemma 12.18, Lemma 12.19]{T}, we see that
		$$J\Delta(\lambda'')=P(\Sigma^{11}\lambda'')=P(8\nu^*_{21})=0.$$
		Therefore $\Delta(\lambda'')\in\mathrm{span}\{\Delta(\xi'')\},$   and  $\Delta(\eta_{10}\cq _{11})=0$ implies that $\operatorname{Im}\Delta=\mathrm{span}\{\Delta(\xi'')\}.$ Hence we obtain
		$$\operatorname{coker}(\Delta\colon \pi^{10}_{28}\to \pi_{27}(\mathit{SO}(10)))=\frac{\mathbb Z}{8}\{[\iota_7]_{10}\ct _7\}.$$ \noindent
		Since $[\eta_{10}\mu_{11}]\in \pi_{20}(\mathit{SO}(11))\cong(\mathbb Z/2)^3$ \cite[Theorem 3.2]{KM2}, we deduce
		$2[\eta_{10}\mu_{11}]\sigma_{20}=0.$
		This proves the proposition.
	\end{proof}

	We need the following lemma.

	\begin{lem}\label{m2:[zt9sg}
		$[\zeta_9]\sigma_{20}\in\mathrm{span}\{\Delta(\sigma_{10}\zeta_{17})\}$.
	\end{lem}
	\begin{proof}
		By Lemma~\ref{m1:zeta-lift-bracket}(2), with $n=2$, we have $[\zeta_9]\in\{\Delta\iota_{10},\eta^2_9,\varepsilon_{11}\}.$ Since
		$\pi^9_{27}=\frac{\mathbb Z}{8}\{\sigma_9\zeta_{16}\}\oplus\frac{\mathbb Z}{2}\{\eta_9\bar{\mu}_{10}\},$   we conclude that $$[\zeta_9]\sigma_{20}\in\{\Delta\iota_{10},\eta^2_9,\varepsilon_{11}\} \circ\sigma_{20}=-\Delta\iota_{10}\circ\{\eta^2_9,\varepsilon_{11},\sigma_{19}\} \subseteq \Delta\iota_{10}\circ\pi^9_{27}.$$ \noindent  Since $\Delta(\eta_{10})=2[\iota_7]_{10}\nu_7,$ we see that
		$\Delta(\eta_{10}\bar{\mu}_{11})=0.$
		This completes the proof.
	\end{proof}
	
	From \cite[III-Proposition 2.2(5)]{Oda1979}, we obtain
	$\bar{\kappa}_7\sigma_{27}\equiv 0 \ \bmod\ 2\sigma'\bar{\kappa}_{14}, \bar{\nu}_7\bar{\sigma}_{15}.$
	In light of \cite[II-Proposition 4.4(6), (8)]{Oda1979}, we deduce
	$P(\bar{\sigma}_{19})=\bar{\nu}_9\bar{\sigma}_{17}$
	and $P(\nu_{23}\kappa_{26})=4\sigma_{11}\bar{\kappa}_{18}.$
	Thus
	\begin{equation}
		\begin{aligned}
			&\bar{\kappa}_{11}\sigma_{31}\equiv 0 \bmod\ P(\nu_{23}\kappa_{26}).\label{m2:kappabarsig}
		\end{aligned}
	\end{equation}
	
	The passage to $n=12$ changes the order of the inherited $\bar\kappa$-class and introduces a class
	represented by $[\mu_{11}]\sigma_{20}.$ Computing the connecting images distinguishes these
	contributions from the   element   already present at $n=11.$
	
	\subsection{The Cases $12\leq n\leq14$}
	
	We next follow the inherited classes through the cases $12\leq n\leq14.$  
	\begin{prop}\label{m2:2712}
		The group $\pi_{27}(\mathit{SO}(12))$ is as follows.
		\[\GroupLayout{\pi_{27}(\mathit{SO}(12)) = \mathrm{span}\{[\iota_7]_{12}\bar\kappa_7,\; [\eta_{10}\mu_{11}]_{12}\sigma_{20},\; [\mu_{11}]\sigma_{20}\}\cong \mathbb Z/4\oplus (\mathbb Z/2)^{2}.}{\begin{aligned}\pi_{27}(\mathit{SO}(12)) = \mathrm{span}\{&[\iota_7]_{12}\bar\kappa_7,\; [\eta_{10}\mu_{11}]_{12}\sigma_{20},\\
				&[\mu_{11}]\sigma_{20}\}\cong \mathbb Z/4\oplus (\mathbb Z/2)^{2}.\end{aligned}}\]
		The new lift and the boundary satisfy
		\[
		\begin{aligned}
			&2([\mu_{11}]\sigma_{20})=0,\qquad \Delta(\nu_{11}\kappa_{14})=4[\iota_7]_{11}\bar\kappa_7.
		\end{aligned}
		\]
	\end{prop}
	
	\begin{proof}
		Consider  the exact sequence $(27)_{11}$:
		\[\begin{aligned}
			&\pi^{11}_{28}\stackrel{{\Delta}}{\longrightarrow}\pi_{27}(\mathit{SO}(11))\stackrel{{i_*}}{\longrightarrow}\pi_{27}(\mathit{SO}(12))\stackrel{{p_*}}{\longrightarrow}\pi^{11}_{27}\stackrel{{\Delta}}{\longrightarrow}\pi_{26}(\mathit{SO}(11)).
		\end{aligned}\]
		\noindent   
		By the proof of Proposition~\ref{m2:lem:rank12-extension}, we have
		$\pi_{27}^{11}=\frac{\mathbb Z}{2}\{\mu_{11}\sigma_{20}\},\quad\Delta(\mu_{11}\sigma_{20})=0.$
		Thus
		\[\FitMath{\ker\bigl(\Delta\colon \pi_{27}^{11}\to\pi_{26}(\mathit{SO}(11))\bigr)=\frac{\mathbb Z}{2}\{\mu_{11}\sigma_{20}\}.}\]

		\noindent       The relevant group is
		$\pi_{28}^{11}
		=\mathrm{span}\{\eta_{11}\sigma_{12}\mu_{19},\nu_{11}\kappa_{14},\bar\mu_{11}\}
		\cong(\mathbb Z/2)^3.$
		By \cite[Lemma 1.1 (i)]{KM1} and \cite[ Lemma 15.4, (16.6)]{MT1}, we conclude that
		\[
		\begin{aligned}
			&\Delta(\eta_{11}\sigma_{12}\mu_{19})=[\iota_7]_{11}\nu_7\eta_{10}\sigma_{11}\mu_{18}=0,\qquad \Delta(\nu_{11}\kappa_{14})=[\iota_7]_{11}\nu_7^2\kappa_{13}=4[\iota_7]_{11}\bar\kappa_7, \\
			&\Delta\bar\mu_{11}=[\iota_7]_{11}\nu_7\bar\mu_{10}=0.
		\end{aligned}
		\]
		Since $[\mu_{11}]\in \pi_{20}(\mathit{SO}(12))\cong(\mathbb Z/2)^5$   \cite[Theorem 3.2]{KM2}, we have
		$2[\mu_{11}]\sigma_{20}=0.$
		This proves the proposition.
	\end{proof}
	
	At $n=13,$ the boundary of $\mu_{12}\sigma_{21}$ removes one of the preceding order-$2$ classes. We
	record that boundary explicitly together with the remaining generators.
	
	\begin{prop}\label{m2:2713}
		The group $\pi_{27}(\mathit{SO}(13))$ is as follows.
		\[\GroupLayout{\pi_{27}(\mathit{SO}(13)) = \mathrm{span}\{[\iota_7]_{13}\bar\kappa_7,\; [\mu_{11}]_{13}\sigma_{20}\}\cong \mathbb Z/4\oplus \mathbb Z/2.}{\begin{aligned}\pi_{27}(\mathit{SO}(13)) = \mathrm{span}\{&[\iota_7]_{13}\bar\kappa_7,\\
				&[\mu_{11}]_{13}\sigma_{20}\}\cong \mathbb Z/4\oplus \mathbb Z/2.\end{aligned}}\]
		The boundary and the surviving lift satisfy
		\[
		\begin{aligned}
			&\Delta(\mu_{12}\sigma_{21})=[\eta_{10}\mu_{11}]_{12}\sigma_{20},\qquad 2([\mu_{11}]_{13}\sigma_{20})=0.
		\end{aligned}
		\]
	\end{prop}
	
	\begin{proof}
		Consider  the exact sequence $(27)_{12}$:
		\[\begin{aligned}
			&\pi^{12}_{28}\stackrel{{\Delta}}{\longrightarrow}\pi_{27}(\mathit{SO}(12))\stackrel{{i_*}}{\longrightarrow}\pi_{27}(\mathit{SO}(13))\stackrel{{p_*}}{\longrightarrow}\pi^{12}_{27}\stackrel{{\Delta}}{\longrightarrow}\pi_{26}(\mathit{SO}(12)).
		\end{aligned}\]
		\noindent By Lemma~\ref{m2:lem:rank13-inherited-image} and \eqref{m2:eq:rank12-kappa-boundary},
		\[\FitMath{\begin{aligned}
				\pi_{27}^{12}&=\frac{\mathbb Z}{16}\{\Sigma^3\rho'\}\oplus\frac{\mathbb Z}{2}\{\eta_{12}\kappa_{13}\},\\
				\Delta(\eta_{12}\kappa_{13})&\equiv[\nu_9\kappa_{12}]_{12}\pmod{64[\Sigma^2\rho']},\\
				\operatorname{Im}\Delta&=\frac{\mathbb Z}{16}\{\Delta(\Sigma^3\rho')\}\oplus\frac{\mathbb Z}{2}\{\Delta(\eta_{12}\kappa_{13})\}.
		\end{aligned}}\]
		Consequently,
		\[\FitMath{\ker\bigl(\Delta\colon \pi_{27}^{12}\to\pi_{26}(\mathit{SO}(12))\bigr)=0.}\]
		\noindent       We note that
		$\pi^{12}_{28}=\mathrm{span}\{\mu_{12}\sigma_{21}\}\cong\mathbb Z/2.$ 
		By \cite[Lemma 3.9 (2)]{KM2} and Lemma \ref{m2:[zt9sg}, we deduce $\Delta(\mu_{12}\sigma_{21})\equiv
		[\eta_{10}\mu_{11}]_{12}\sigma_{20} \bmod [\zeta_9]_{12}\sigma_{20}=0.$ This completes the proof.
	\end{proof}
	
	We now examine $n=14,$ expressing the proposed generators as connecting images of $\sigma_{14}^2$
	and $\kappa_{14}.$ The accompanying identities compare these boundary classes with the inherited
	lifts.
	
	\begin{prop}\label{m2:2714}
		The group $\pi_{27}(\mathit{SO}(14))$ is as follows.
		\[\GroupLayout{\pi_{27}(\mathit{SO}(14)) = \mathrm{span}\{\Delta(\sigma_{14}^2),\; \Delta\kappa_{14}\}\cong \mathbb Z/8\oplus \mathbb Z/2.}{\begin{aligned}\pi_{27}(\mathit{SO}(14)) = \mathrm{span}\{&\Delta(\sigma_{14}^2),\\
				&\Delta\kappa_{14}\}\cong \mathbb Z/8\oplus \mathbb Z/2.\end{aligned}}\]
		The new lift and the boundary classes satisfy
		\[
		\begin{aligned}
			&8\Delta(\sigma_{14}^2)=[\mu_{11}]_{14}\sigma_{20}=0,\qquad \Delta\kappa_{14}=[\iota_7]_{14}\bar\kappa_7, \\
			&\Delta\bar\varepsilon_{13}=2[\iota_7]_{13}\bar\kappa_7.
		\end{aligned}
		\]
	\end{prop}
	
	\begin{proof}
		Consider  the exact sequence $(27)_{13}$:
		\[\begin{aligned}
			&\pi^{13}_{28}\stackrel{{\Delta}}{\longrightarrow}\pi_{27}(\mathit{SO}(13))\stackrel{{i_*}}{\longrightarrow}\pi_{27}(\mathit{SO}(14))\stackrel{{p_*}}{\longrightarrow}\pi^{13}_{27}\stackrel{{\Delta}}{\longrightarrow}\pi_{26}(\mathit{SO}(13)).
		\end{aligned}\]
		
		\noindent   The proof of Proposition~\ref{g1016:prop:26-14} and the decomposition in Proposition~\ref{m2:lem:rank13-conditional} give
		\[\FitMath{\begin{aligned}
				\pi_{27}^{13}&=\frac{\mathbb Z}{16}\{\sigma_{13}^2\}\oplus\frac{\mathbb Z}{2}\{\kappa_{13}\},\\
				\operatorname{Im}\Delta&=\frac{\mathbb Z}{2}\{\Delta(\sigma_{13}^2)\}\oplus\frac{\mathbb Z}{2}\{\Delta\kappa_{13}\}.
		\end{aligned}}\]
		Hence
		\[\FitMath{\ker\bigl(\Delta\colon \pi_{27}^{13}\to\pi_{26}(\mathit{SO}(13))\bigr)=\frac{\mathbb Z}{8}\{2\sigma_{13}^2\}.}\]
		\noindent   We find
		$\pi^{13}_{28}=\mathrm{span}\{\rho_{13}, \overline{\varepsilon}_{13}\}\cong\mathbb Z/32\oplus\mathbb Z/2.$ 
		By \cite[Lemma 3.8(5)]{KM2}, we conclude that
		$J[\mu_{11}]\equiv \beta''' \bmod \mathrm{span}\{\Sigma\theta'\varepsilon_{24}, \Sigma\beta''\}$ and $$J([\mu_{11}]_{13}\sigma_{20})\equiv(\Sigma\beta''')\sigma_{33} \bmod \mathrm{span}\{(\Sigma^2\theta')\varepsilon_{25}\sigma_{33}, (\Sigma^2\beta'')\sigma_{33}\}.$$ In light of \cite[(7.30)]{T}, we see that $(\Sigma^2\theta')\varepsilon_{25}\sigma_{33}=0.$
		By \cite[I-(2.3), p. 95]{Oda1979}, we obtain $$\Sigma^2(\beta'')\sigma_{33}=(\Sigma^2C^{(1)}_1)\sigma_{33}=0,\;\;(\Sigma\beta''')\sigma_{33}=(\Sigma A^{(1)}_1)\sigma_{33}.$$
		Thus
		$J([\mu_{11}]_{13}\sigma_{20})=(\Sigma A^{(1)}_1)\sigma_{33}.$	Notice that  $\Delta(\bar{\varepsilon}_{13})=[\eta_{11}]_{13}\bar{\varepsilon}_{12}.$
		Using $\Delta\iota_{11}=[\iota_7]_{11}\nu_7$ and
		$2\bar\kappa_7\in\{\nu_7,\eta_{10},\bar\varepsilon_{11}\},$   we have
		\begin{align*}
			[\eta_{11}]\bar\varepsilon_{12}
			&\in (i_{12})_*\{\Delta\iota_{11},\eta_{10},\bar\varepsilon_{11}\} \notag\\
			&= (i_{12})_*\{[\iota_7]_{11}\nu_7,\eta_{10},\bar\varepsilon_{11}\} \notag\\
			&\ni 2[\iota_7]_{12}\bar\kappa_7. \notag
		\end{align*}
		The comparison is modulo $(i_{12})_*\pi_{12}(\mathit{SO}(11))\circ\bar\varepsilon_{12}.$ 
		We note that $\Delta(\eta_{12})=[\eta^2_{10}]_{12}$ \cite[(1.8)]{HKM},
		$\pi_{12}(\mathit{SO}(11))=\frac{\mathbb Z}{2}\{[\eta^2_{10}]\}$ \cite[Table 3]{Ka} and $\eta_9\bar{\varepsilon}_{10}=0.$  Therefore,
		$${i_{12}}_*\pi_{12}(\mathit{SO}(11))\circ\bar{\varepsilon}_{12}=\mathrm{span}\{\Delta(\eta_{12}\bar{\varepsilon}_{13})\}=0.$$ This implies
		$[\eta_{11}]\bar{\varepsilon}_{12}=2[\iota_7]_{12}\bar{\kappa}_7$ and
		$\Delta(\bar{\varepsilon}_{13})=2[\iota_7]_{13}\bar{\kappa}_7.$ 	Since $J\Delta(\rho_{13})=P(\rho_{27})\equiv(\Sigma A^{(1)}_1)\sigma_{33} \bmod 2\sigma_{13}\ct _{20}$
		(\cite[p.96]{Oda1979}), we see that
		$\Delta(\rho_{13})\equiv [\mu_{11}]_{13}\sigma_{20} \bmod 2[\iota_7]_{13}\ct _7.$
		This implies $[\mu_{11}]_{14}\sigma_{20}=0.$
		
		In view of \cite[p. 111]{KM2}, we obtain $[\mu_{11}]_{14}=8\Delta(\sigma_{14}).$
		Thus $[\mu_{11}]_{14}\sigma_{20}=8\Delta(\sigma^2_{14})=0.$
		Hence we obtain the group structure.
		Since ${p_{14}}_*\Delta\kappa_{14}=2\kappa_{13}=0,$ we deduce
		$$\Delta\kappa_{14}\in {i_{14}}_*\pi_{27}(\mathit{SO}(13))=\mathrm{span}\{[\iota_7]_{14}\bar{\kappa}_7\}.$$ \noindent 
		By \cite[II-Proposition~4.4(10), p.~92]{Oda1979}, $P(\kappa_{29})=\sigma_{14}\bar{\kappa}_{21}-\bar{\kappa}_{14}\sigma_{34}.$ Suspending \eqref{m2:kappabarsig} three times and using $\Sigma P=0$ gives $\bar{\kappa}_{14}\sigma_{34}=0.$ Hence, by the formula relating $J$ and $\Delta$ in Section~\ref{m2:subsec:j-toda}, $J\Delta\kappa_{14}=P(\kappa_{29})=\sigma_{14}\bar{\kappa}_{21}.$
		This leads to the second relation and completes the proof.
	\end{proof}
	
	\section{The Cases $15\leq n\leq 17$}\label{m2:subsec:degree27-middle}
	The cases $n=15,16$ introduce the free cyclic summand of the degree-$27$ calculation. We describe it
	using the same lifted multiple of $P(\iota_{29}),$ so that its images under the inclusions remain
	explicit.
	
	\begin{prop}\label{m2:2715}
		The homotopy groups for $n=15,16$ are given as follows.
		\[\FitMath{\pi_{27}(\mathit{SO}(15+n)) = \mathrm{span}\{[256P(\iota_{29})]_{15+n}\}\cong \mathbb Z_{(2)},\qquad n=0,1.}\]
		The lift and its image under inclusion satisfy
		\[
		\begin{aligned}
			&(p_{15})_*[256P(\iota_{29})]_{15}=256P(\iota_{29}),\qquad (i_{16})_*[256P(\iota_{29})]_{15}=[256P(\iota_{29})]_{16}.
		\end{aligned}
		\]
		\qed
	\end{prop}
	
	By \cite[I-(8.28)]{Oda1979}, we conclude that
	$P(\zeta_{33})=\pm 2M^{(1)}_1 + x\Sigma^4\tau^{\mathrm{IV}}\ \text{for an odd integer }x.$
	This implies $\mathrm{ord}(\Delta\zeta_{16})=8.$
	
	This boundary calculation also allows the degree-$27$ generator to be followed to $n=17.$ The
	proposition records its image before we turn to the larger-$n$ tables.
	
	\begin{prop}\label{m2:2717}
		We have \[ \FitMath{\pi_{27}(\mathit{SO}(17)) = \mathrm{span}\{[256P(\iota_{29})]_{17}\}\cong \mathbb Z_{(2)}.} \]
		Additionally,  $\Delta\zeta_{16}=2[\zeta_{15}]+a[8\sigma_{12}^2]_{16}$ for some odd integer $a.$ 
		The element  $\Delta\zeta_{16}$ has order $8.$ 
		\qed
	\end{prop}
	
	\section{The Cases $n\geq 18$}\label{m2:subsec:degree27-high}
	Using the relation $\pi_{27}(\mathit{SO})\cong \mathbb Z_{(2)}$ \cite{Bo} and the results in \cite[Table 1]{HoM},
	we have Table~\ref{m2:R27}.
	\par
	\begin{table}[H]
		\centering
		\caption{$\pi_{27}(\mathit{SO}(n))$ for $18\leq n\leq 23$}\label{m2:R27}
		\begin{tabular*}{\textwidth}{@{\extracolsep{\fill}}clclcl@{}}
			\toprule
			\textbf{Group} & \textbf{Structure} & \textbf{Group} & \textbf{Structure} & \textbf{Group} & \textbf{Structure} \\
			\midrule
			$\pi_{27}(\mathit{SO}(18))$ & $\mathbb Z_{(2)}$
			& $\pi_{27}(\mathit{SO}(19))$ & $\mathbb Z_{(2)} \oplus (\mathbb Z/2)^2$
			& $\pi_{27}(\mathit{SO}(20))$ & $\mathbb Z_{(2)} \oplus (\mathbb Z/2)^4$  \\ \midrule
			$\pi_{27}(\mathit{SO}(21))$ & $\mathbb Z_{(2)}\oplus (\mathbb Z/2)^2$
			& $\pi_{27}(\mathit{SO}(22))$ & $\mathbb Z_{(2)} \oplus \mathbb Z/4$
			& $\pi_{27}(\mathit{SO}(23))$ & $\mathbb Z_{(2)} \oplus \mathbb Z/2$\\ \bottomrule
		\end{tabular*}
		
	\end{table}
	\par
	
	In view of \cite{Ke}, we  have Table~\ref{m2:R27K}.
	\par
	\begin{table}[H]
		\centering
		\caption{$\pi_{27}(\mathit{SO}(n))$ for $24\leq n\leq 29$}\label{m2:R27K}
		\begin{tabular*}{\textwidth}{@{\extracolsep{\fill}}clclcl@{}}
			\toprule
			\textbf{Group} & \textbf{Structure} & \textbf{Group} & \textbf{Structure} & \textbf{Group} & \textbf{Structure} \\
			\midrule
			$\pi_{27}(\mathit{SO}(24))$ & $\mathbb Z_{(2)}\oplus\mathbb Z/2$
			& $\pi_{27}(\mathit{SO}(25))$ & $\mathbb Z_{(2)}\oplus\mathbb Z/2$
			& $\pi_{27}(\mathit{SO}(26))$ & $\mathbb Z_{(2)}$     \\ \midrule
			$\pi_{27}(\mathit{SO}(27))$ & $\mathbb Z_{(2)}$
			& $\pi_{27}(\mathit{SO}(28))$ & $\mathbb Z_{(2)} \oplus \mathbb Z_{(2)}$
			& $\pi_{27}(\mathit{SO}(29))$ & $\mathbb Z_{(2)}$ \\ \bottomrule
		\end{tabular*}
		
	\end{table}
	\par
	
	Tables~\ref{m2:R27} and~\ref{m2:R27K} give the abstract group structures
	for $n\geq18.$   The corresponding generators and their relations
	are listed in Theorem~\ref{m2:thm2}.
	
	\medskip
	
	\chapter{The Homotopy Groups $\pi_{28}(\mathit{SO}(n))$}\label{ch:degree28}
	\label{m3:chapter}
	
	In degree $28,$   the stable group vanishes. The argument approaches the unstable groups from both
	directions: small-$n$ calculations use $G_2$ and spin-group comparisons, while the final part works
	downwards through the cases $n=15,14,13,12.$

	The   element   for $n=14$ is also the source of the geometric examples in Chapter~\ref{ch:geometry}.
	Following its boundary representative and $J$-image explains how an explicit group generator becomes a
	bundle and a Thom attaching map.
	
	As usual, in our calculations, we often draw freely on the descriptions of
	$\pi_{*}(S^{j})$ in Appendix~\ref{tlqdzj} and on the composition
	relations among generators in Toda's range tabulated in
	\^{O}guchi's table in Appendix~\ref{app:oguchi-relations}.

	\section{Main Theorem and Generators}

	\begin{mainleft}
		\begin{thm}\label{m3:thm28}
			\textup{All spaces in this theorem are localized at $2.$  The displayed groups are obtained from the integral homotopy groups by tensoring with $\mathbb Z_{(2)}.$  Free cyclic summands become $\mathbb Z_{(2)}$; for finite groups, this localization retains precisely the $2$-primary summands. The homotopy groups are given as follows. }\par
			\textup{The Cases $3\leq n\leq6.$ }\par
			\smallskip\noindent\(
			\GroupLayout{\pi_{28}( \mathit{SO}(3)) = \mathrm{span}\{[\eta_{2}] \phi',\; [\eta_{2}] \nu' \varepsilon_{6}\kappa_{14},\; [\eta_{2}] \mu_{3,3},\; [\eta_{2}] \eta_{3}\bar{\mu}_{4}\sigma_{21}\}\cong \mathbb Z/4\oplus (\mathbb Z/2)^{3}.}{\begin{aligned}\pi_{28}( \mathit{SO}(3)) = \mathrm{span}\{&[\eta_{2}] \phi',\; [\eta_{2}] \nu' \varepsilon_{6}\kappa_{14},\; [\eta_{2}] \mu_{3,3},\\
					&[\eta_{2}] \eta_{3}\bar{\mu}_{4}\sigma_{21}\}\cong \mathbb Z/4\oplus (\mathbb Z/2)^{3}.\end{aligned}}
			\)\par
			\smallskip\noindent\(
			\GroupLayout{\pi_{28}( \mathit{SO}(4)) = \mathrm{span}\{[\iota_{3}] \phi',\; [\iota_{3}] \nu' \varepsilon_{6}\kappa_{14},\; [\iota_{3}] \mu_{3,3},\; [\iota_{3}] \eta_{3}\bar{\mu}_{4}\sigma_{21},\; [\eta_{2}]_{4} \phi',\; [\eta_{2}] _{4}\nu' \varepsilon_{6}\kappa_{14},\; [\eta_{2}]_{4} \mu_{3,3},\; [\eta_{2}] _{4} \eta_{3}\bar{\mu}_{4}\sigma_{21}\} \, \;\cong\; \mathbb Z/4 \oplus(\mathbb Z/2)^{3} \oplus\mathbb Z/4 \oplus(\mathbb Z/2)^{3}.}{\begin{aligned}\pi_{28}( \mathit{SO}(4)) = \mathrm{span}\{&[\iota_{3}] \phi',\; [\iota_{3}] \nu' \varepsilon_{6}\kappa_{14},\; [\iota_{3}] \mu_{3,3},\; [\iota_{3}] \eta_{3}\bar{\mu}_{4}\sigma_{21},\; [\eta_{2}]_{4} \phi',\; [\eta_{2}] _{4}\nu' \varepsilon_{6}\kappa_{14},\\
					&[\eta_{2}]_{4} \mu_{3,3},\; [\eta_{2}] _{4} \eta_{3}\bar{\mu}_{4}\sigma_{21}\} \, \;\cong\; \mathbb Z/4 \oplus(\mathbb Z/2)^{3} \oplus\mathbb Z/4 \oplus(\mathbb Z/2)^{3}.\end{aligned}}
			\)\par
			\smallskip\noindent\(
			\GroupLayout{\pi_{28}(\mathit{SO}(5)) = \mathrm{span}\{[\iota_{3}]_{5}\phi',\; [\iota_{3}]_{5}\mu_{3,3},\; [\iota_{3}]_{5}\eta_{3}\bar{\mu}_{4}\sigma_{21},\; [ \nu_{4}\sigma'\kappa_{14} ]\}\cong \mathbb Z/4\oplus (\mathbb Z/2)^{3}.}{\begin{aligned}\pi_{28}(\mathit{SO}(5)) = \mathrm{span}\{&[\iota_{3}]_{5}\phi',\; [\iota_{3}]_{5}\mu_{3,3},\; [\iota_{3}]_{5}\eta_{3}\bar{\mu}_{4}\sigma_{21},\\
					&[ \nu_{4}\sigma'\kappa_{14} ]\}\cong \mathbb Z/4\oplus (\mathbb Z/2)^{3}.\end{aligned}}
			\)\par
			\smallskip\noindent\(
			\GroupLayout{\pi_{28}(\mathit{SO}(6)) = \mathrm{span}\{[\iota_{3}]_{6}\phi',\; [\nu_{4}\sigma'\kappa_{14}]_{6},\; [\bar{\rho}'''],\; [ \nu_{5}]\bar{\kappa}_{8}\}\cong (\mathbb Z/2)^{2}\oplus \mathbb Z/4\oplus \mathbb Z/8.}{\begin{aligned}\pi_{28}(\mathit{SO}(6)) = \mathrm{span}\{&[\iota_{3}]_{6}\phi',\; [\nu_{4}\sigma'\kappa_{14}]_{6},\; [\bar{\rho}'''],\\
					&[ \nu_{5}]\bar{\kappa}_{8}\}\cong (\mathbb Z/2)^{2}\oplus \mathbb Z/4\oplus \mathbb Z/8.\end{aligned}}
			\)\par
			\medskip\noindent\textup{The Cases $7\leq n\leq11.$ }\par
			\smallskip\noindent\(
			\GroupLayout{\pi_{28}(\mathit{SO}(7)) \;= {} \mathrm{span}\{[\bar{\nu} _{6}+\varepsilon_{6}]\kappa_{14},\; [\bar{\sigma}_{6} ]\nu_{25},\; [ \eta_{6}]\eta_{7}\bar{\kappa}_{8},\; [ \nu_5]_{7}\bar{\kappa}_{8},\; [\nu_{4}\sigma'\kappa_{14}]_{7}\}\cong (\mathbb Z/2)^{5}.}{\begin{aligned}\pi_{28}(\mathit{SO}(7)) \;= {} \mathrm{span}\{&[\bar{\nu} _{6}+\varepsilon_{6}]\kappa_{14},\; [\bar{\sigma}_{6} ]\nu_{25},\; [ \eta_{6}]\eta_{7}\bar{\kappa}_{8},\\
					&[ \nu_5]_{7}\bar{\kappa}_{8},\; [\nu_{4}\sigma'\kappa_{14}]_{7}\}\cong (\mathbb Z/2)^{5}.\end{aligned}}
			\)\par
			\smallskip\noindent\(
			\GroupLayout{\pi_{28}(\mathit{SO}(8)) = \mathrm{span}\{[\bar{\nu} _{6}+\varepsilon_{6}]_{8}\kappa_{14},\; [\bar{\sigma}_{6} ] _{8}\nu_{25},\; [ \eta_{6}]_{8}\eta_{7}\bar{\kappa}_{8},\; [ \nu_5]_{8}\bar{\kappa}_{8},\; [\nu_{4}\sigma'\kappa_{14}]_{8},\; [\iota_{7}]\eta_{7}\bar{\kappa}_{8},\; [ \iota_{7}]\sigma'\kappa_{14}\} \; \;\cong\; (\mathbb Z/2)^{7}.}{\begin{aligned}\pi_{28}(\mathit{SO}(8)) = \mathrm{span}\{&[\bar{\nu} _{6}+\varepsilon_{6}]_{8}\kappa_{14},\; [\bar{\sigma}_{6} ] _{8}\nu_{25},\; [ \eta_{6}]_{8}\eta_{7}\bar{\kappa}_{8},\; [ \nu_5]_{8}\bar{\kappa}_{8},\\
					&[\nu_{4}\sigma'\kappa_{14}]_{8},\; [\iota_{7}]\eta_{7}\bar{\kappa}_{8},\; [ \iota_{7}]\sigma'\kappa_{14}\} \; \;\cong\; (\mathbb Z/2)^{7}.\end{aligned}}
			\)\par
			\smallskip\noindent\(
			\GroupLayout{\pi_{28}(\mathit{SO}(9)) = \mathrm{span}\{[\bar{\nu} _{6}+\varepsilon_{6}]_{9}\kappa_{14},\; [\bar{\sigma}_{6}] _{9} \nu_{25},\; [ \nu_5]_{9}\bar{\kappa}_{8},\; [\iota_{7}]_{9}\eta_{7}\bar{\kappa}_{8}\}\cong (\mathbb Z/2)^{4}.}{\begin{aligned}\pi_{28}(\mathit{SO}(9)) = \mathrm{span}\{&[\bar{\nu} _{6}+\varepsilon_{6}]_{9}\kappa_{14},\; [\bar{\sigma}_{6}] _{9} \nu_{25},\\
					&[ \nu_5]_{9}\bar{\kappa}_{8},\; [\iota_{7}]_{9}\eta_{7}\bar{\kappa}_{8}\}\cong (\mathbb Z/2)^{4}.\end{aligned}}
			\)\par
			\smallskip\noindent\(
			\GroupLayout{\pi_{28}(\mathit{SO}(10)) = \mathrm{span}\{[\bar{\zeta}_{9}],\; [\bar{\nu} _{6}+\varepsilon_{6}]_{10}\kappa_{14},\; [\bar{\sigma}_{6}]_{10}\nu_{25},\; [ \nu_5]_{10}\bar{\kappa}_{8}\}\cong \mathbb Z/8\oplus (\mathbb Z/2)^{3}.}{\begin{aligned}\pi_{28}(\mathit{SO}(10)) = \mathrm{span}\{&[\bar{\zeta}_{9}],\; [\bar{\nu} _{6}+\varepsilon_{6}]_{10}\kappa_{14},\; [\bar{\sigma}_{6}]_{10}\nu_{25},\\
					&[ \nu_5]_{10}\bar{\kappa}_{8}\}\cong \mathbb Z/8\oplus (\mathbb Z/2)^{3}.\end{aligned}}
			\)\par
			
			For $\pi_{28}(\mathit{SO}(11)),$   the remaining possibilities are:
			\[\GroupLayout{(i_{10,11})_*\pi_{28}(\mathit{SO}(10))={}\mathrm{span}\{[\bar\zeta_9]_{11},\; [\mathrm{c}_{\mathrm{KM}}(\bar\nu_6+\varepsilon_6)]_{11}\kappa_{14},\; [ \nu_5]_{11}\bar\kappa_8\}\cong (\mathbb Z/2)^{3}.}{\begin{aligned}(i_{10,11})_*\pi_{28}(\mathit{SO}(10))={}\mathrm{span}\{&[\bar\zeta_9]_{11},\; [\mathrm{c}_{\mathrm{KM}}(\bar\nu_6+\varepsilon_6)]_{11}\kappa_{14},\\
					&[ \nu_5]_{11}\bar\kappa_8\}\cong (\mathbb Z/2)^{3}.\end{aligned}}\]
			\[\GroupLayout{p_*\pi_{28}(\mathit{SO}(11))={}\mathrm{span}\{\eta_{10}\bar\mu_{11},\; 4t\xi''+(\xi''\pm\lambda'')\}\cong (\mathbb Z/2)^{2},}{\begin{aligned}p_*\pi_{28}(\mathit{SO}(11))={}\mathrm{span}\{&\eta_{10}\bar\mu_{11},\\
					&4t\xi''+(\xi''\pm\lambda'')\}\cong (\mathbb Z/2)^{2},\end{aligned}}\]
			{\raggedleft where $t\in\mathbb Z$ and $\mathrm{ord}(\xi''\pm\lambda'')=2.$ \par} \noindent The full group is given by
			
			\par\smallskip\noindent\(
			\GroupLayout{\pi_{28}(\mathit{SO}(11)) = \mathrm{span}\{[\bar\zeta_9]_{11},\; [\mathrm{c}_{\mathrm{KM}}(\bar\nu_6+\varepsilon_6)]_{11}\kappa_{14},\; [ \nu_5]_{11}\bar\kappa_8,\; [\eta_{10}\bar\mu_{11}],\; [4t\xi''+(\xi''\pm\lambda'')]\},}{\begin{aligned}\pi_{28}(\mathit{SO}(11)) = \mathrm{span}\{&[\bar\zeta_9]_{11},\; [\mathrm{c}_{\mathrm{KM}}(\bar\nu_6+\varepsilon_6)]_{11}\kappa_{14},\; [ \nu_5]_{11}\bar\kappa_8,\\
					&[\eta_{10}\bar\mu_{11}],\; [4t\xi''+(\xi''\pm\lambda'')]\},\end{aligned}}
			\)\par{\raggedleft which is isomorphic to  $(\mathbb Z/2)^5$ or $(\mathbb Z/2)^3\oplus\mathbb Z/4.$ \par}
			\noindent   \textup{The Cases $n\geq12.$ }\par
			\smallskip\noindent\(
			\begin{aligned}
				&\pi_{28}(\mathit{SO}(12)) \;\cong\; (\mathbb Z/2)^{6} \text{generated by the elements}\; [\eta_{11}]\sigma_{12}\mu_{19},\; [\bar\mu_{11}] \;\text{and elements in}\;  \Delta\pi_{29}^{12}.
			\end{aligned}
			\)
			\par\smallskip\noindent\(
			\GroupLayout{\pi_{28}(\mathit{SO}(13)) = \mathrm{span}\{[\eta_{11}]_{13}\sigma_{12}\mu_{19},\; [\bar{\mu}_{11}]_{13}\}\cong (\mathbb Z/2)^{2}.}{\begin{aligned}\pi_{28}(\mathit{SO}(13)) = \mathrm{span}\{&[\eta_{11}]_{13}\sigma_{12}\mu_{19},\\
					&[\bar{\mu}_{11}]_{13}\}\cong (\mathbb Z/2)^{2}.\end{aligned}}
			\)\par
			\smallskip\noindent\(
			\FitMath{\pi_{28}(\mathit{SO}(14)) = \mathrm{span}\{\Delta\rho_{14}\}\cong \mathbb Z/32.}
			\)\par
			\smallskip\noindent\(
			\begin{aligned}
				&\pi_{28}(\mathit{SO}(15)) = 0.
			\end{aligned}
			\)\par
			\smallskip\noindent\(
			\begin{aligned}
				&\pi_{28}(\mathit{SO}(16)) = \pi_{28}(\mathit{SO}(17))=0.
			\end{aligned}
			\)\par
			\smallskip\noindent\(
			\FitMath{\pi_{28}(\mathit{SO}(18)) = \mathrm{span}\{[\zeta_{17}]\}\cong \mathbb Z/8.}
			\)\par
			\smallskip\noindent\(
			\GroupLayout{\pi_{28}(\mathit{SO}(19)) = \mathrm{span}\{[\zeta_{17}]_{19},\; [\eta_{18}\mu_{19}]\}\cong (\mathbb Z/2)^{2}.}{\begin{aligned}\pi_{28}(\mathit{SO}(19)) = \mathrm{span}\{&[\zeta_{17}]_{19},\\
					&[\eta_{18}\mu_{19}]\}\cong (\mathbb Z/2)^{2}.\end{aligned}}
			\)\par
			\smallskip\noindent\(
			\GroupLayout{\pi_{28}(\mathit{SO}(20)) = \mathrm{span}\{[\zeta_{17}]_{20},\; [\eta_{18}\mu_{19}]_{20},\; [\eta_{19}] \bar{\nu} _{20},\; [\eta_{19}]\varepsilon_{20},\; [\mu_{19}]\}\cong (\mathbb Z/2)^{5}.}{\begin{aligned}\pi_{28}(\mathit{SO}(20)) = \mathrm{span}\{&[\zeta_{17}]_{20},\; [\eta_{18}\mu_{19}]_{20},\; [\eta_{19}] \bar{\nu} _{20},\\
					&[\eta_{19}]\varepsilon_{20},\; [\mu_{19}]\}\cong (\mathbb Z/2)^{5}.\end{aligned}}
			\)\par
			\smallskip\noindent\(
			\GroupLayout{\pi_{28}(\mathit{SO}(21)) = \mathrm{span}\{[\eta_{19}]_{21} \bar{\nu} _{20},\; [\eta_{19}]_{21}\varepsilon_{20},\; [\mu_{19}]_{21}\}\cong (\mathbb Z/2)^{3}.}{\begin{aligned}\pi_{28}(\mathit{SO}(21)) = \mathrm{span}\{&[\eta_{19}]_{21} \bar{\nu} _{20},\; [\eta_{19}]_{21}\varepsilon_{20},\\
					&[\mu_{19}]_{21}\}\cong (\mathbb Z/2)^{3}.\end{aligned}}
			\)\par
			\smallskip\noindent\(
			\FitMath{\pi_{28}(\mathit{SO}(22)) = \mathrm{span}\{\Delta(\sigma_{22})\}\cong \mathbb Z/16.}
			\)\par
			\smallskip\noindent\(
			\begin{aligned}
				&\pi_{28}(\mathit{SO}(23)) = 0.
			\end{aligned}
			\)\par
			\smallskip\noindent\(
			\begin{aligned}
				&\pi_{28}(\mathit{SO}(24)) = \pi_{28}(\mathit{SO}(25))=0.
			\end{aligned}
			\)\par
			\smallskip\noindent\(
			\FitMath{\pi_{28}(\mathit{SO}(26)) = \mathrm{span}\{[2\nu_{25}]\}\cong \mathbb Z/4.}
			\)\par
			\smallskip\noindent\(
			\FitMath{\pi_{28}(\mathit{SO}(27)) = \mathrm{span}\{[\eta_{26}^{2}]\}\cong \mathbb Z/2.}
			\)\par
			\smallskip\noindent\(
			\GroupLayout{\pi_{28}(\mathit{SO}(28)) = \mathrm{span}\{[\eta_{26}^{2}]_{28},\; [\eta_{27}]\}\cong (\mathbb Z/2)^{2}.}{\begin{aligned}\pi_{28}(\mathit{SO}(28)) = \mathrm{span}\{&[\eta_{26}^{2}]_{28},\\
					&[\eta_{27}]\}\cong (\mathbb Z/2)^{2}.\end{aligned}}
			\)\par
			\smallskip\noindent\(
			\FitMath{\pi_{28}(\mathit{SO}(29)) = \mathrm{span}\{[\eta_{27}]_{29}\}\cong \mathbb Z/2.}
			\)\par
			\smallskip\noindent\(
			\begin{aligned}
				&\pi_{28}(\mathit{SO}(30+n)) = 0\,(n\geq0).
			\end{aligned}
			\)\par
			
		\end{thm}
	\end{mainleft}
	
	\section{The Cases $3\leq n\leq6$}
	
	We begin with the boundary and lift relations needed for the cases $3\leq n\leq6.$  These relations
	prepare the small-$n$ decompositions and the comparison with the unitary group used in the following
	section.
	
	\needspace{3\baselineskip}\begin{prop}\label{m3:dn4}
		The following relations hold.
		\begin{itemize}
			\item[\textup{(1)}]
			$\Delta(\nu_4\sigma'\rho_{14})
			=\Delta(\nu_4^2\bar{\sigma}_{10})=0.$ 
			
			\item[\textup{(2)}]
			$\Delta(\nu_4\bar{\nu}_7\kappa_{15})
			=\Delta(\nu_4\varepsilon_7\kappa_{15})
			\equiv[\iota_3]\nu'\varepsilon_6\kappa_{14}
			\bmod [\eta_2]_4\nu'\varepsilon_6\kappa_{14}.$ 
			
			\item[\textup{(3)}]
			$\Delta(\Sigma\phi')=2[\iota_3]\phi'-[\eta_2]_4\phi',$  
			$\Delta(\Sigma\nu'\circ\varepsilon_7\kappa_{15})
			=[\eta_2]_4\nu'\varepsilon_6\kappa_{14}.$ 
			
			\item[\textup{(4)}]
			$\Delta(\mu_{3,4})=[\eta_2]_4\mu_{3,3},$  
			$\Delta(\eta_4\bar{\mu}_5\sigma_{22})
			=[\eta_2]_4\eta_3\bar{\mu}_4\sigma_{21}.$ 
			
			\item[\textup{(5)}]
			$\Delta(\delta_4)=[\eta_2]_4\delta_3,$  
			$\Delta(\bar{\mu}_4\sigma_{21})
			=[\eta_2]_4\bar{\mu}_3\sigma_{20},$  
			$\Delta(\Sigma\varepsilon'\circ\kappa_{14})
			=[\eta_2]_4\varepsilon'\kappa_{13},$  
			$\Delta(\nu_4\sigma'\kappa_{14})=0.$ 
			
			\item[\textup{(6)}]
			$\Delta(\nu_4\eta_7\bar{\kappa}_8)
			\equiv[\iota_3]\varepsilon'\kappa_{13}
			\bmod [\eta_2]_4\varepsilon'\kappa_{13}.$ 
			
			\item[\textup{(7)}]
			$\Delta(\delta_5)\equiv2[\iota_3]_5\phi'
			\bmod [\iota_3]_5\eta_3\bar{\mu}_4\sigma_{21}.$ 
			
			\item[\textup{(8)}]
			$\Delta(\bar{\mu}_5\sigma_{22})
			=[\iota_3]_5\eta_3\bar{\mu}_4\sigma_{21},$  
			$\Delta(\nu_5\eta_8\bar{\kappa}_9)
			=\Delta(\nu_5\sigma_8\kappa_{15})=0.$ 
		\end{itemize}
	\end{prop}
	
	\begin{proof}
		Write $\Delta\nu_{4}=[\iota_{3}]\nu'+t[\eta_{2}]_{4}\nu'$ for
		some $t\in\mathbb Z.$
		\begin{itemize}
			\item[\textup{(1)}] For the first assertion, we note that
			$\Delta(\nu_{4}\sigma'\rho_{14})=
			2[\iota_{3}]\varepsilon'\circ\rho_{13}=[\iota_{3}]\circ2\varepsilon'\rho_{13}.$ Then, by Lemma
			\ref{m3:2ltp}(1), we find
			$[\iota_{3}]\circ2\varepsilon'\rho_{13}=[\iota_{3}]\circ(\eta_{3}^{2}\varepsilon_{5})\rho_{13}=[\iota_{3}]\circ\eta_{3}^{2}(\varepsilon_{5}\rho_{13})$=0,
			and hence $\Delta(\nu_{4}\sigma'\rho_{14})=0.$ For the second assertion, by Lemma~\ref{m3:2ltp}, we conclude that
			\[\begin{aligned}
				&\Delta(\nu_{4}\nu_{7}\bar{\sigma}_{10}) =([\iota_{3}]\nu' +t[\eta_{2}]_{4}\nu')\nu_{6}\bar{\sigma}_{9} =([\iota_{3}] +t[\eta_{2}]_{4})(\nu'\nu_{6})\bar{\sigma}_{9} =0.
			\end{aligned}\]
			\item[\textup{(2)}] In light of Proposition~\ref{m3:np}, we know 
			\[
			\begin{aligned}
				\Delta(\nu_{4} \bar{\nu} _{7}\kappa_{15}) &={} ([\iota_{3}]\nu' \, +t[\eta_{2}]_{4}\nu') \bar{\nu} _{6}\kappa_{14} {} \\
				&={} [\iota_{3}]\nu' \bar{\nu} _{6}\kappa_{14} +t[\eta_{2}]_{4}\nu' \bar{\nu} _{6}\kappa_{14} {} \\
				&={} [\iota_{3}]\nu'\varepsilon_{6}\kappa_{14} +t[\eta_{2}]_{4}\nu'\varepsilon_{6}\kappa_{14},
			\end{aligned}
			\]
			\[\begin{aligned}
				&\Delta(\nu_{4}\varepsilon_{7}\kappa_{15}) =([\iota_{3}]\nu' +t[\eta_{2}]_{4}\nu')\varepsilon_{6}\kappa_{14} =[\iota_{3}]\nu'\varepsilon_{6}\kappa_{14} +t[\eta_{2}]_{4}\nu'\varepsilon_{6}\kappa_{14}.
			\end{aligned}\]
			\item[\rm (3)]  For the first assertion, we see that
			\[\begin{aligned}
				&\Delta(\, \Sigma\phi') = \Delta(\iota_{4}\circ\, \Sigma\phi') = (2[\iota_{3}]-[\eta_{2}]_{4})\phi' =2[\iota_{3}]\phi'-[\eta_{2}]_{4}\phi'.
			\end{aligned}\]
			For the second
			assertion, we note that
			\[\begin{aligned}
				&\Delta( \, \Sigma\nu'\circ\varepsilon_{7}\kappa_{15} ) = \Delta( \iota_{4}\circ\, \Sigma\nu'\circ\varepsilon_{7}\kappa_{15} ) = (2[\iota_{3}]-[\eta_{2}]_{4})\nu'\varepsilon_{6}\kappa_{14} =[\eta_{2}]_{4}\nu'\varepsilon_{6}\kappa_{14}\uprightsemicolon{}
			\end{aligned}\]
			\item[\textup{(4)}] For the first assertion, we observe that $\Delta(\mu_{3,\allowbreak 4})=
			\Delta(\iota_{4}\circ\mu_{3,\allowbreak 4})= (2[\iota_{3}]-[\eta_{2}]_{4})\mu_{3,\allowbreak 3}=
			[\eta_{2}]_{4}\mu_{3,\allowbreak 3}.$
			For the second assertion, we have
			\[\begin{aligned}
				&\Delta(\eta_{4}\bar{\mu}_{5}\sigma_{22}) = \Delta(\iota_{4}\circ\eta_{4}\bar{\mu}_{5}\sigma_{22}) = (2[\iota_{3}]-[\eta_{2}]_{4})\eta_{3}\bar{\mu}_{4}\sigma_{21} =[\eta_{2}]_{4}\eta_{3}\bar{\mu}_{4}\sigma_{21}.
			\end{aligned}\]
			\item[\textup{(5)}] The four formulas follow from the characteristic-element formula in Table~\ref{m2:tab:boundary-values} by the same argument as in (4).
			\item[\textup{(6)}] The  statement \cite[Proposition 3.1 (2), p. 79]{Oda1979} gives $\nu'\eta_{6}\bar{\kappa}_{7}=\varepsilon'\kappa_{13},$ so
			\begin{align*}
				\Delta(\nu_{4}\eta_{7}\bar{\kappa}_{8}) {}&={} ([\iota_{3}]\nu' +t[\eta_{2}]_{4}\nu')\eta_{6}\bar{\kappa}_{7} {} \\
				&={} ([\iota_{3}]+t[\eta_{2}]_{4})\nu'\eta_{6}\bar{\kappa}_{7} \\
				&=([\iota_{3}]+t[\eta_{2}]_{4})\varepsilon'\kappa_{13} \\
				&=[\iota_{3}]\varepsilon'\kappa_{13}+t[\eta_{2}]_{4}\varepsilon'\kappa_{13}.
			\end{align*}
			\item[\textup{(7)}] The  third line from the bottom of \cite[p.~68]{Oda1979} gives $\eta_{3}\delta_{4}
			\equiv2\phi'\;\bmod\;\eta_{3}\bar{\mu}_{4}\sigma_{21}.$ Then $\Delta(\delta_{5})=[\iota_{3}]_{5}\eta_{3}\delta_{4}
			\equiv2[\iota_{3}]_{5}\phi'\;\bmod\;[\iota_{3}]_{5}\eta_{3}\bar{\mu}_{4}\sigma_{21}.$
			\item[\textup{(8)}] Table~\ref{m2:tab:boundary-values} gives $\Delta\iota_5=[\iota_3]_5\eta_3$ and $\Delta\nu_5=0.$ The suspended-right-factor formula in Section~\ref{m2:subsec:boundaries} gives $\Delta(\bar\mu_5\sigma_{22})=(\Delta\iota_5)\bar\mu_4\sigma_{21}=[\iota_3]_5\eta_3\bar\mu_4\sigma_{21}.$ The rule $\Delta\alpha=0\Rightarrow\Delta(\alpha\beta)=0$ in the same section gives $\Delta(\nu_5\eta_8\bar\kappa_9)=\Delta(\nu_5\sigma_8\kappa_{15})=0.$
		\end{itemize}
	\end{proof}
	
	The following proposition follows immediately from $\mathit{SO}(3)=\mathbb{R}P^{3}$ and $$\pi_{j}(
	\mathit{SO}(4))=[\iota_{3}] \circ\pi_{j}^{3}\oplus [\eta_{2}]_{4} \circ\pi_{j}^{3} \cong
	\pi_{j}^{3}\oplus\pi_{j}^{3}\qquad(j\geq2).$$ \noindent See \cite[p. 14]{Ka} for details.
	\needspace{3\baselineskip}\begin{prop}
		The homotopy groups are given as follows.
		\[\GroupLayout{\pi_{28}( \mathit{SO}(3)) = \mathrm{span}\{[\eta_{2}] \phi',\; [\eta_{2}] \nu' \varepsilon_{6}\kappa_{14},\; [\eta_{2}] \mu_{3,3},\; [\eta_{2}] \eta_{3}\bar{\mu}_{4}\sigma_{21}\}\cong \mathbb Z/4\oplus (\mathbb Z/2)^{3}.}{\begin{aligned}\pi_{28}( \mathit{SO}(3)) = \mathrm{span}\{&[\eta_{2}] \phi',\; [\eta_{2}] \nu' \varepsilon_{6}\kappa_{14},\; [\eta_{2}] \mu_{3,3},\\
				&[\eta_{2}] \eta_{3}\bar{\mu}_{4}\sigma_{21}\}\cong \mathbb Z/4\oplus (\mathbb Z/2)^{3}.\end{aligned}}\]
		\begin{gather*}
			\GroupLayout{\pi_{28}( \mathit{SO}(4)) \;= {} \mathrm{span}\{[\iota_{3}] \phi',\; [\iota_{3}] \nu' \varepsilon_{6}\kappa_{14},\; [\iota_{3}] \mu_{3,3},\; [\iota_{3}] \eta_{3}\bar{\mu}_{4}\sigma_{21},\; [\eta_{2}]_{4} \phi',\; [\eta_{2}] _{4}\nu' \varepsilon_{6}\kappa_{14},\; [\eta_{2}]_{4} \mu_{3,3},\; [\eta_{2}] _{4} \eta_{3}\bar{\mu}_{4}\sigma_{21}\}\cong \mathbb Z/4\oplus (\mathbb Z/2)^{3}\oplus \mathbb Z/4\oplus (\mathbb Z/2)^{3}.}{\begin{aligned}\pi_{28}( \mathit{SO}(4)) \;= {} \mathrm{span}\{&[\iota_{3}] \phi',\; [\iota_{3}] \nu' \varepsilon_{6}\kappa_{14},\; [\iota_{3}] \mu_{3,3},\; [\iota_{3}] \eta_{3}\bar{\mu}_{4}\sigma_{21},\; [\eta_{2}]_{4} \phi',\\
					&[\eta_{2}] _{4}\nu' \varepsilon_{6}\kappa_{14},\; [\eta_{2}]_{4} \mu_{3,3},\; [\eta_{2}] _{4} \eta_{3}\bar{\mu}_{4}\sigma_{21}\}\cong \mathbb Z/4\oplus (\mathbb Z/2)^{3}\oplus \mathbb Z/4\oplus (\mathbb Z/2)^{3}.\end{aligned}}
		\end{gather*}
		\qed
	\end{prop}
	\needspace{3\baselineskip}\begin{lem}
		\label{m3:f284} For $\mathrm{(28)}_{4},$ $\operatorname{coker}(
		\Delta)=\mathrm{span}\{[\iota_{3}]\phi',\; [\iota_{3}]\mu_{3, 3},\; [\iota_{3}]\eta_{3}\bar{\mu}_{4}\sigma_{21}\}\cong \mathbb Z/4\oplus (\mathbb Z/2)^{2}.$
		Equivalently, $i_{5_{*}} (\pi_{28}(\mathit{SO}(4)))=\mathrm{span}\{[\iota_{3}]_{5}\phi',\; [\iota_{3}]_{5}\mu_{3, 3},\; [\iota_{3}]_{5}\eta_{3}\bar{\mu}_{4}\sigma_{21}\}\cong \mathbb Z/4\oplus (\mathbb Z/2)^{2}.$
	\end{lem}
	\begin{proof}
		This follows immediately from Proposition~\ref{m3:dn4}(1)--(4).
	\end{proof}
	\needspace{3\baselineskip}The preceding calculations provide the necessary groundwork
	for the following result.\begin{prop}
		\label{m3:r285}
		The structure of $\pi_{28}(\mathit{SO}(5))$ is as follows.
		\[\GroupLayout{\pi_{28}(\mathit{SO}(5)) = \mathrm{span}\{[\iota_{3}]_{5}\phi',\; [\iota_{3}]_{5}\mu_{3,3},\; [\iota_{3}]_{5}\eta_{3}\bar{\mu}_{4}\sigma_{21},\; [ \nu_{4}\sigma'\kappa_{14} ]\}\cong \mathbb Z/4\oplus (\mathbb Z/2)^{3}.}{\begin{aligned}\pi_{28}(\mathit{SO}(5)) = \mathrm{span}\{&[\iota_{3}]_{5}\phi',\; [\iota_{3}]_{5}\mu_{3,3},\; [\iota_{3}]_{5}\eta_{3}\bar{\mu}_{4}\sigma_{21},\\
				&[ \nu_{4}\sigma'\kappa_{14} ]\}\cong \mathbb Z/4\oplus (\mathbb Z/2)^{3}.\end{aligned}}\]
	\end{prop}
	\begin{proof}
		By \cite[Theorem~5.1(A), p.~85]{Oda1978},
		$\pi_{28}(\mathit{SO}(5))\cong\pi_{28}(\mathit{Sp}(2))\cong\mathbb Z/4\oplus(\mathbb Z/2)^{3}.$
		Lemma~\ref{m3:f284} and Proposition~\ref{m3:dn4}(5) then give the group structure of $\pi_{28}(\mathit{SO}(5)).$
	\end{proof}
	\needspace{3\baselineskip}\begin{lem}
		\label{m3:f285} For $\mathrm{(28)}_{5},$ $i_{6_{*}} (\pi_{28}(\mathit{SO}(5)))=\mathrm{span}\{[\iota_{3}]_{6}\phi',\; [\iota_{3}]_{6}\mu_{3, 3},\; [\nu_{4}\sigma'\kappa_{14}]_{6}\}\cong (\mathbb Z/2)^{3}.$
	\end{lem}
	\begin{proof}
		The structure of $\pi_{29}^{5}$ is known, and $\pi_{28}(\mathit{SO}(5))$ is given
		in Proposition~\ref{m3:r285}. The image follows immediately from Proposition~\ref{m3:dn4}(7)--(8).
	\end{proof}
	In the following lemma, the same notation is used for maps and their homotopy classes.
	\needspace{3\baselineskip}\begin{lem}
		\label{m3:qr} Let $i_{S^{3},\allowbreak \mathit{SU}(3)}: S^{3}\hookrightarrow \mathit{SU}(3)$ and
		$i_{\mathit{SU}(3),\allowbreak \mathit{SU}(4)}: \mathit{SU}(3)\hookrightarrow \mathit{SU}(4)$ be the inclusions,
		and let $q_{_{\mathit{SU}(4),\allowbreak \mathit{SO}(6)}}: \mathit{SU}(4)\rightarrow \mathit{SO}(6)$ be the
		covering map. Then $$\begin{aligned}
			&q_{_{\mathit{SU}(4),\allowbreak \mathit{SO}(6)}} \circ i_{\mathit{SU}(3),\allowbreak \mathit{SU}(4)} \circ i_{S^{3},\allowbreak \mathit{SU}(3)} =\pm[\iota_{3}]_{6}\in \pi_{3}( \mathit{SO}(6)).
		\end{aligned}$$
	\end{lem}
	\begin{proof}
		We know that $\pi_{3}(\mathit{SO}(6))=\mathbb Z_{(2)}\{ [\iota_{3}]_{6}\}$ \cite[p.~18]{Ka},
		$\pi_{3}(\mathit{SU}(4))=\mathbb Z_{(2)}\{ \,i_{S^{3},\allowbreak \mathit{SU}(4)} \},$  
		where $i_{S^{3},\allowbreak \mathit{SU}(4)}:=i_{\mathit{SU}(3),\allowbreak \mathit{SU}(4)}\circ i_{S^{3},\allowbreak \mathit{SU}(3)}: S^{3}\hookrightarrow \mathit{SU}(4)$ \cite{MT2}. Since
		$q_{_{\mathit{SU}(4),\allowbreak \mathit{SO}(6)}}$ induces isomorphisms on the $k$th homotopy groups for
		$k\geq3,$ the lemma follows.
	\end{proof}
	\needspace{3\baselineskip}\begin{prop}
		\label{m3:r286}
		The structure of $\pi_{28}(\mathit{SO}(6))$ is as follows.
		\[\GroupLayout{\pi_{28}(\mathit{SO}(6)) = \mathrm{span}\{[\iota_{3}]_{6}\phi',\; [\nu_{4}\sigma'\kappa_{14}]_{6},\; [ \bar{\rho}''' ],\; [ \nu_{5}]\bar{\kappa}_{8}\}\cong (\mathbb Z/2)^{2}\oplus \mathbb Z/4\oplus \mathbb Z/8.}{\begin{aligned}\pi_{28}(\mathit{SO}(6)) = \mathrm{span}\{&[\iota_{3}]_{6}\phi',\; [\nu_{4}\sigma'\kappa_{14}]_{6},\; [ \bar{\rho}''' ],\\
				&[ \nu_{5}]\bar{\kappa}_{8}\}\cong (\mathbb Z/2)^{2}\oplus \mathbb Z/4\oplus \mathbb Z/8.\end{aligned}}\]
		In addition, $ [ \bar{\rho}''' ]:=\mathcal{D}_*[ \bar{\rho}''' ]^{U}$, where $[ \bar{\rho}''' ]^{U}=[ \bar{\rho}''' ]\in\pi_*(\mathit{SU}(3)).$ Here, $$\mathcal{D}\colon \mathit{SU}(3)\xyd \mathit{SU}(4)\stackrel{\cong}\longrightarrow \mathit{Spin}(6)\twoheadrightarrow \mathit{SO}(6)$$ is the map in Diagram~\eqref{tu1029}.
	\end{prop}
	\begin{proof}
		First, we claim that, for $\mathrm{(28)}_{5},$ $\operatorname{ker}(
		\Delta)=\frac{\mathbb Z}{8}\{\nu_{5}\bar{\kappa}_{8}\}\oplus \frac{\mathbb Z}{2}\{\bar{\rho}'''\}.$
		We note that
		$$\begin{aligned}
			&\pi_{28}(\mathit{SO}(6))\cong\pi_{28}(\mathit{SU}(4))\cong(\mathbb Z/2)^{2}\oplus\mathbb Z/4\oplus\mathbb Z/8
		\end{aligned}$$
		and
		$i_{5,\allowbreak 6*} (\pi_{28}(\mathit{SO}(5)))=\frac{\mathbb Z}{2}\{[\iota_{3}]_{6}\phi'\}\oplus \frac{\mathbb Z}{2}\{[\iota_{3}]_{6}\mu_{3, 3}\}\oplus \frac{\mathbb Z}{2}\{[\nu_{4}\sigma'\kappa_{14}]_{6}\}$ by Lemma~\ref{m3:f285}. For $\mathrm{(28)}_{5},$ exactness gives
		$|\operatorname{ker}( \Delta)|=16.$ Observe that $$\operatorname{ker}( \Delta)
		\subseteq\pi_{28}^{5}\cong\mathbb Z/8\oplus(\mathbb Z/2)^{2}$$
		and $\frac{\mathbb Z}{8}\{\nu_{5}\bar{\kappa}_{8}\}\subseteq\operatorname{ker}( \Delta).$  So
		$\operatorname{ker}( \Delta)\cong\mathbb Z/8\oplus\mathbb Z/2.$
		In light of \cite[Theorem 6.1 A]{Oda1978}, we find
		$$\bar{\rho}''' \in\operatorname{Im}(p_{*}:\pi_{28}(\mathit{SU}(4))\rightarrow\pi_{28}(S^{5})\oplus\pi_{28}
		(S^{7})\;).$$ \noindent This means $$\bar{\rho}'''
		\in\operatorname{Im}\mathbf{(}p_{*}:\pi_{28}(\mathit{SU}(4))\rightarrow\pi_{28}(S^{5})\mathbf{)}.$$
		\noindent   Therefore, $\bar{\rho}''' \in\operatorname{ker}(\partial_{*}).$
		Consider the morphism of fiber sequences
		\[
		\xymatrix{
			\Omega S^{5} \ar[d]^{\mathrm{id}}\ar[r]^{\partial} & \mathrm{\mathit{Spin}(5)} \ar[d]^{q} \ar[r] &
			\mathrm{\mathit{SU}(4)}\ar[d] \ar[r]^{}& S^{5}\ar[d]^{\mathrm{id}} \\
			\Omega S^{5} \ar[r]^{ \Delta}&\mathit{SO}(5) \ar[r] &\mathit{SO}(6) \ar[r]&S^{5}.}
		\]
		Then,   $\Delta\bar{\rho}''' =q_{*}\partial_{*}\bar{\rho}''' =0$ and $\bar{\rho}''' \in\operatorname{ker}(
		\Delta).$ Since the kernel has order $16$ in $\pi_{28}^{5}\cong\mathbb Z/8\oplus(\mathbb Z/2)^{2},$   we obtain $\operatorname{ker}( \Delta)=\frac{\mathbb Z}{8}\{\nu_{5}\bar{\kappa}_{8}\}\oplus \frac{\mathbb Z}{2}\{\bar{\rho}'''\}.$
		Then the group structures $i_{5,\allowbreak 6*}(\pi_{28}(\mathit{SO}(5)))\cong (\mathbb Z/2)^{3},\allowbreak
		\;\pi_{28}(\mathit{SO}(6))\cong(\mathbb Z/2)^{2}\oplus\mathbb Z/4\oplus\mathbb Z/8$
		suggest that $\mathrm{ord}([ \nu_5]\bar{\kappa}_{8})=8$ and $\mathrm{ord}[\bar{\rho}''']=4.$
		
		By \cite[Theorem 4.1 B), p. 83]{Oda1978},
		we obtain $i_{S^{3},\allowbreak \mathit{SU}(3)}\circ\mu_{3,\allowbreak 3}=2 [\bar{\rho}''' ]\in\pi_{28}(\mathit{SU}(3)).$
		Lemma~\ref{m3:qr} then gives the following commutative diagram:
		\[
		\begin{tikzcd}[row sep=1.2cm, column sep=1.8cm]
			& & \mathit{SU}(3)
			\arrow[d, "i"']
			\\
			S^{28}
			\arrow[r, "\mu_{3,3}"']
			\arrow[urr, "{2[\bar{\rho}''']}"]
			& S^{3}
			\arrow[r, "i"']
			\arrow[ur, "i"']
			\arrow[dr, "{\pm[\iota_{3}]_{6}}"']
			& \mathit{SU}(4)
			\arrow[d, "q"', "(\text{covering})"]
			\\
			& & \mathit{SO}(6).
		\end{tikzcd}
		\]  \noindent   Thus, $[\iota_{3}]_{6}\mu_{3,\allowbreak 3}= \pm q\circ i \circ2[\bar{\rho}'''] \in \pi_{28}(\mathit{SO}(6))$ is divisible by
		2.
		This proves the proposition.
	\end{proof}
	
	\section{The Auxiliary Degree 28 Calculation for $G_2$}\label{sec:degree28-g2}
	
	To pass from the special orthogonal groups with small values of $n$ to the case $n=7,$   we need explicit generators for
	$\pi_{28}(G_2).$  We state the decomposition first, then construct the unitary-group input, calculate
	the connecting images, and choose the required lifts.
	
	\begin{thm} Proposition~\ref{m3:p28g2} gives the following decomposition.  $$\pi_{28}(G_{2})=\mathrm{span}\{\langle \bar{\nu}
		_{6}+\varepsilon_{6}\rangle\kappa_{14},\; \langle \bar{\sigma}_{6} \nu_{25}\rangle,\; \langle
		\eta_{6}^{2}\rangle\bar{\kappa}_{8}\}\cong (\mathbb Z/2)^{3}.$$
	\end{thm}
	
	Recall from \cite[Theorem 4.1]{MT2} that $\pi_{20} (\mathit{SU}(3))=\frac{\mathbb Z}{4}\{ [\rho^{\mathrm{IV}}]\}\oplus\frac{\mathbb Z}{2}\{i_{*}\bar{\varepsilon}'\}.$
	\subsection{Unitary Generators and Toda Brackets}
	
	The unitary calculation begins with a choice of generator expressed by a Toda bracket. Such a
	representative is useful because compositions and multiples can then be compared using the bracket
	identities from Chapter~\ref{ch:toda}.
	
	\needspace{3\baselineskip}\begin{lem}
		\label{m3:fp20} In $\pi_{20} (\mathit{SU}(3))=\mathrm{span}\{[\rho^{\mathrm{IV}}],\; i_{*}\bar{\varepsilon}'\}\cong \mathbb Z/4\oplus \mathbb Z/2,$ the generator $[\rho^{\mathrm{IV}}]$ may be
		chosen so that $[\rho^{\mathrm{IV}}]\in\{ [\sigma'''],\allowbreak 4\iota_{12},\allowbreak  4\sigma_{12}\}_{1}.$
	\end{lem}
	\begin{proof}
		By \cite[Theorem 4.1]{MT2}, $[\sigma''']\in \pi_{12} (\mathit{SU}(3))$ has
		order 4, so $\{ [\sigma'''],\allowbreak 4\iota_{12},\allowbreak  4\sigma_{12}\}_{1}$ is well-defined. Choose
		$[\rho^{\mathrm{IV}}]\in\{ [\sigma'''],\allowbreak 4\iota_{12},\allowbreak  4\sigma_{12}\}_{1}.$ By \cite[p. 103]{T},
		$\rho^{\mathrm{IV}}\in\{ \sigma''',\allowbreak 2\iota_{12},\allowbreak  8\sigma_{12}\}_{1}\subseteq\{ \sigma''',\allowbreak 4\iota_{12},\allowbreak
		4\sigma_{12}\}_{1}.$ It follows that
		\[\begin{aligned}
			&p_{*}([\rho^{\mathrm{IV}}]) \in\{ \sigma''', 4\iota_{12}, 4\sigma_{12}\}_{1} \ni \rho^{\mathrm{IV}}\; \bmod\; \sigma'''\varepsilon_{12} =0, \; \sigma''' \bar{\nu} _{12} =0, \; 4\varepsilon_{5}\sigma_{13} =0\uprightsemicolon{}
		\end{aligned}\]
		for $\sigma'''\varepsilon_{12}= \sigma''' \bar{\nu}
		_{12}=0.$ 
		This proves the assertion.
	\end{proof}
	\needspace{3\baselineskip}\begin{lem}
		\label{m3:2g0}
		\leavevmode\par
		\begin{itemize}
			\item[\textup{(1)}]
			$\GroupLayout{\{ 2\iota_{9}, \bar{\nu} _{9},\sigma_{17}\}_{1} =\{2\iota_{9},\varepsilon_{9},\sigma_{17}\}_{1}=\mathrm{span}\{\sigma_{9}\eta_{16}^{2}\sigma_{18},\; \mu_{9}\sigma_{18}\}\cong (\mathbb Z/2)^{2},\qquad 0 \in\{2\iota_{9},\eta_{9}\sigma_{10},\sigma_{17}\}_{1}.\;}{\begin{aligned}\{ 2\iota_{9}, \bar{\nu} _{9},\sigma_{17}\}_{1} =\{2\iota_{9},\varepsilon_{9},\sigma_{17}\}_{1}=\mathrm{span}\{&\sigma_{9}\eta_{16}^{2}\sigma_{18},
					\mu_{9}\sigma_{18}\}\cong (\mathbb Z/2)^{2},\quad 0 \in\{2\iota_{9},\eta_{9}\sigma_{10},\sigma_{17}\}_{1}.\;\end{aligned}}$
			\item[\textup{(2)}]
			$\begin{aligned}
				&0 \in\{ \varepsilon_{3},2\iota_{11},\eta_{11}\sigma_{12} \}_{1}=\{ \varepsilon_{3},2\iota_{11},\eta_{11}\sigma_{12} \}_{2}.
			\end{aligned}$
		\end{itemize}
	\end{lem}
	\begin{proof}
		\begin{enumerate}[label=\textup{(\arabic*)},leftmargin=2em]
			\item We note that $\eta_{9}\bar{\varepsilon}_{10}= \bar{\nu} _{6}\sigma_{14}=0.$ 
			The structure of $\pi_{16+7}^{7}$ shows that there exist $x_{1},\allowbreak x_{2},\allowbreak x_{3},\allowbreak x_{4}\in\mathbb Z_{(2)}$ such that
			$x_{1}\sigma'\mu_{14}+x_{2}\sigma'\eta_{14}\varepsilon_{15}+x_{3}\mu_{7}\sigma_{16}+x_{4}\eta_{7}\bar{\varepsilon}_{8}\in\{
			2\iota_{7},\allowbreak  \bar{\nu} _{7},\allowbreak \sigma_{15}\}_{1}.$
			Applying $\Sigma^{2},$ we obtain
			$x_{3}\mu_{9}\sigma_{18}\in\{ 2\iota_{9},\allowbreak  \bar{\nu} _{9},\allowbreak \sigma_{17}\}_{1}.$ Notice that 
			\[\begin{aligned}
				&\mathrm{Ind}\{ 2\iota_{9}, \bar{\nu} _{9}, \sigma_{17}\}_{1} =\pi_{18}^{9}\circ\sigma_{18} =\mathrm{span}\{\sigma_{9}\eta_{16}^{2}\sigma_{18}, \; \mu_{9}\sigma_{18}\} \ni x_{3}\mu_{9}\sigma_{18}.
			\end{aligned}\]
			Then we
			deduce $$\{ 2\iota_{9},\allowbreak  \bar{\nu} _{9},\allowbreak \sigma_{17}\}_{1}=\frac{\mathbb Z}{2}\{
			\sigma_{9}\eta_{16}^{2}\sigma_{18}\}\oplus\frac{\mathbb Z}{2}\{\mu_{9}\sigma_{18}\}.$$ The structure of
			$\pi_{16+3}^{3}$ gives $\{2\iota_{3},\allowbreak \varepsilon_{3},\allowbreak \sigma_{11}\}=t\eta_{3}\bar{\varepsilon}_{4}+
			\mathrm{span}\{ \mu_{3}\sigma_{12}\}$ for some $t\in\mathbb Z,$
			and $\eta_{9}\bar{\varepsilon}_{10}=0,$
			and hence $$\{ 2\iota_{9},\allowbreak \varepsilon_{9},\allowbreak \sigma_{17}\}_{1}=\frac{\mathbb Z}{2}\{
			\sigma_{9}\eta_{16}^{2}\sigma_{18}\}\oplus\frac{\mathbb Z}{2}\{\mu_{9}\sigma_{18}\}.$$
			The remaining assertion follows from $\ca_{9}\cg_{10}=\ch_9+\lt_{9}$ and $2\e_{9}\hc\pi_{25}^{9}=2\e_{9}\hc\s\pi_{24}^{8}=0$.
			
			\item Observe that $\varepsilon_{3}\nu_{11}^{3}=2\bar{\varepsilon}'=\varepsilon'\sigma_{13}.$ 
			We consider $\{ \varepsilon_{3},\allowbreak 2\iota_{11},\allowbreak \eta_{11} \}_{1}.$ Notice that
			$$H\{ \varepsilon_{3},\allowbreak 2\iota_{11},\allowbreak \eta_{11} \}_{1}\subseteq \{ \nu_{5}^{2} ,\allowbreak 2\iota_{11},\allowbreak \eta_{11}
			\}_{1}=\varepsilon_{5}.$$\noindent
			By \cite[p. 67]{T}, $\varepsilon'+t\mu_{3}\eta_{12}\in\{ \varepsilon_{3},\allowbreak 2\iota_{11},\allowbreak \eta_{11}
			\}_{1}$ for some $t\in\mathbb Z.$
			Then $$(\varepsilon'+t\mu_{3}\eta_{12})\sigma_{13}=\varepsilon_{3}\nu_{11}^{3}+t\mu_{3}\eta_{12}\sigma_{13}\in\{
			\varepsilon_{3},\allowbreak 2\iota_{11},\allowbreak \eta_{11}\sigma_{12} \}_{1}.$$ \noindent
			The indeterminacy of $\{ \varepsilon_{3},\allowbreak 2\iota_{11},\allowbreak \eta_{11}\sigma_{12} \}_{1}$ is
			$\varepsilon_{3}\circ \mathrm{span}\{ \varepsilon_{11}\eta_{19},\allowbreak \nu_{11}^{3},\allowbreak \mu_{11} \}+
			\mathrm{span}\{ \mu_{3}\eta_{12}\sigma_{13} \}.$ Thus
			$(\varepsilon'+t\mu_{3}\eta_{12})\sigma_{13}=2\bar{\varepsilon}'+t\mu_{3}\eta_{12}\sigma_{13}$
			belongs to the indeterminacy of $\{ \varepsilon_{3},\allowbreak 2\iota_{11},\allowbreak \eta_{11}\sigma_{12} \}_{1},$ and
			hence $0\in\{ \varepsilon_{3},\allowbreak 2\iota_{11},\allowbreak \eta_{11}\sigma_{12} \}_{1}.$ It is evident that
			$\{ \varepsilon_{3},2\iota_{11},\eta_{11}\sigma_{12} \}_{1}
			=\{ \varepsilon_{3},2\iota_{11},\eta_{11}\sigma_{12} \}_{2}$.
		\end{enumerate}
	\end{proof}
	
	We recall several results on higher Toda brackets.
	
	Work integrally. Let $n_1,n_2\geq 0,$   and consider homotopy classes
	\begin{equation}\label{eq:null-quadruple}
		\begin{aligned}
			&X_0 \xleftarrow{\alpha_1} \Sigma ^{n_1}X_1, \qquad X_1 \xleftarrow{\alpha_2} \Sigma ^{n_2}X_2, \qquad X_2 \xleftarrow{\alpha_3} X_3 \xleftarrow{\alpha_4} X_4, \\
			&\alpha_1\circ \Sigma^{n_1}\alpha_2=0, \qquad \alpha_2\circ \Sigma^{n_2}\alpha_3=0, \qquad \alpha_3\circ\alpha_4=0.
		\end{aligned}
	\end{equation}
	These data form a \emph{null quadruple}, denoted by
	$(\alpha_1,\alpha_2,\alpha_3,\alpha_4)_{n_1,n_2}.$ 
	
	Choose representatives $a_i\in\alpha_i$ and based
	null-homotopies
	\[\begin{aligned}
		&A_1:a_1\circ \Sigma^{n_1}a_2\simeq *, \qquad A_2:a_2\circ \Sigma^{n_2}a_3\simeq *, \qquad A_3:a_3\circ a_4\simeq *.
	\end{aligned}\]
	For a null homotopy $H:f\circ g\simeq *,$   write
	$[f,H,g]$ and $(f,H,g)$ for the associated extension
	and coextension, respectively. If
	$g:X\to Y$ and $f:Y\to Z,$   their domains and
	codomains are
	\[\begin{aligned}
		&[f,H,g]:Y\cup_g CX\longrightarrow Z, \qquad (f,H,g):\Sigma X\longrightarrow Z\cup_f CY.
	\end{aligned}\]
	Here $\widetilde \Sigma^{\,n}H$ denotes the suspended
	homotopy, given by
	\[\begin{aligned}
		&(\widetilde \Sigma^{\,n}H)(x\wedge s,t) =H(x,t)\wedge s.
	\end{aligned}\]
	
	The representative
	$(a_1,a_2,a_3,a_4\uprightsemicolon{}A_1,A_2,A_3)_{n_1,n_2}$
	is called \emph{admissible} if
	\begin{align*}
		&[a_1,A_1,\Sigma^{n_1}a_2]\circ \bigl(\Sigma^{n_1}a_2,\widetilde \Sigma^{\,n_1}A_2, \Sigma^{n_1+n_2}a_3\bigr) \simeq *, \\
		&[a_2,A_2,\Sigma^{n_2}a_3]\circ \bigl(\Sigma^{n_2}a_3,\widetilde \Sigma^{\,n_2}A_3, \Sigma^{n_2}a_4\bigr) \simeq *.
	\end{align*}
	The null quadruple is called \emph{admissible}
	if it has an admissible representative.

	The following is given in   \cite[Proposition~4.1, p.~26]{O2015}.
	\begin{prop}
		\label{prop:admissible-null-quadruple}
		Work integrally. For the null quadruple \eqref{eq:null-quadruple},
		define subgroups of $[\Sigma^{n_2+1}X_3,X_1]$ by
		\begin{align*}
			&K_1 =(\Sigma^{n_1})^{-1}\left( (\alpha_{1*})^{-1}\left( (\Sigma^{n_1+n_2+1}\alpha_3)^* [\Sigma^{n_1+n_2+1}X_2,X_0]\right)\right), \\
			&K_2 =\left((\Sigma^{n_2+1}\alpha_4)^*\right)^{-1} \left(\alpha_{2*} \left(\Sigma^{n_2}[\Sigma X_4,X_2]\right)\right).
		\end{align*}
		Here $\alpha_{i*}$ denotes postcomposition,
		$(\Sigma^r\alpha_i)^*$ denotes precomposition, and
		each superscript ``$-1$'' denotes an inverse image.
		Equivalently,
		\begin{align*}
			&K_1 =\left\{ \gamma\in[\Sigma ^{n_2+1}X_3,X_1]\ :\
			\alpha_1\circ \Sigma^{n_1}\gamma
			\in [\Sigma ^{n_1+n_2+1}X_2,X_0]
			\circ \Sigma^{n_1+n_2+1}\alpha_3 \right\}, \\
			&K_2 =\left\{ \gamma\in[\Sigma ^{n_2+1}X_3,X_1]\ :\
			\gamma\circ \Sigma^{n_2+1}\alpha_4
			\in \alpha_2\circ \Sigma^{n_2}[\Sigma X_4,X_2] \right\}.
		\end{align*}
		If
		\[\begin{aligned}
			&0\in\{\alpha_1,\Sigma^{n_1}\alpha_2,\Sigma^{n_1+n_2}\alpha_3\}_{n_1},\quad 0\in\{\alpha_2,\Sigma^{n_2}\alpha_3,\Sigma^{n_2}\alpha_4\}_{n_2}, \quad K_1+K_2=[\Sigma ^{n_2+1}X_3,X_1],
		\end{aligned}\]
		then $(\alpha_1,\alpha_2,\alpha_3,\alpha_4)_{n_1,n_2}$
		is admissible.
		\qed
	\end{prop}
	
	\needspace{3\baselineskip}\begin{lem}
		\label{m3:tsd}
		\leavevmode\par
		\begin{itemize}
			\item[\textup{(1)}] Use the notation above and work $2$-locally. Then the null quadruple $(
			\varepsilon_{3},\allowbreak 2\iota_{9},\allowbreak \eta_{9}\sigma_{10},\allowbreak \sigma_{17})_{2,0}$ is admissible.   Consequently, there is an extension $\overline{2\iota_{9}}\in \mathrm{Ext}_{\eta_{9}\sigma_{10}}(2\iota_{9})$ and a coextension $\widetilde{\sigma_{17}}\in\mathrm{Coext}_{\eta_{9}\sigma_{10}}(\sigma_{17})$   such that the Toda bracket
			$$\begin{aligned}
				&\{ \varepsilon_{3},\allowbreak \,\Sigma^{2}\overline{2\iota_{9}},\allowbreak \,\Sigma^{2}\widetilde{\sigma_{17}} \}_{2} \end{aligned}$$ is well-defined. We write
			$$\begin{aligned}
				&\{ \varepsilon_{3},\allowbreak 2\iota_{11},\allowbreak \eta_{11}\sigma_{12},\allowbreak \sigma_{19} \}_2=\{ \varepsilon_{3},\allowbreak \,\Sigma^{2}\overline{2\iota_{9}},\allowbreak \,\Sigma^{2}\widetilde{\sigma_{17}} \}_{2}.
			\end{aligned}$$
			\item[\textup{(2)}] Let  $j:S^{3}\hookrightarrow
			\mathit{SU}(3)$ be the inclusion and use the notation of assertion (1). Then the Toda bracket
			$$\begin{aligned}
				&\{ j\varepsilon_{3},\allowbreak \,\Sigma^{2}\overline{2\iota_{9}},\allowbreak \,\Sigma^{2}\widetilde{\sigma_{17}} \}_{2} \end{aligned}$$ is well-defined. We write
			$$\begin{aligned}
				&\{ j\varepsilon_{3},\allowbreak 2\iota_{11},\allowbreak \eta_{11}\sigma_{12},\allowbreak \sigma_{19} \}_2=\{j \varepsilon_{3},\allowbreak \,\Sigma^{2}\overline{2\iota_{9}},\allowbreak \,\Sigma^{2}\widetilde{\sigma_{17}} \}_{2}.
			\end{aligned}$$
		\end{itemize}
	\end{lem}
	\begin{proof}
		\begin{itemize}
			\item[\textup{(1)}]  Take
			$(X_{0}, \allowbreak \Sigma^{2} X_{1},\allowbreak \Sigma^{2}X_{2},\allowbreak \Sigma^{2}X_{3},\allowbreak \Sigma^{2}X_{4})=(S^{3},\allowbreak S^{11},\allowbreak S^{11},\allowbreak S^{19},\allowbreak S^{26}),\allowbreak $
			and take
			$$\begin{aligned}
				&(\alpha_{1},\allowbreak \Sigma^{2}\alpha_{2},\allowbreak \Sigma^{2}\alpha_{3},\allowbreak \Sigma^{2}\alpha_{4})=( \varepsilon_{3},\allowbreak 2\iota_{11},\allowbreak \eta_{11}\sigma_{12},\allowbreak \sigma_{19} ).
			\end{aligned}$$
			
			(For convenience, here  $\Sigma^{2}\alpha_{3}=\eta_{11}\sigma_{12}$ means
			$\alpha_{3}=\eta_{9}\sigma_{10}$\uprightsemicolon{} similar conventions apply to the remaining elements.)
			Now, use Proposition~\ref{prop:admissible-null-quadruple}. In this case, $$\begin{aligned}
				&K_{1}=(\Sigma^{2})^{-1} ((\varepsilon_{3})_{*})^{-1} (\pi_{12}^{3}\circ \eta_{12}\sigma_{13}),
			\end{aligned}$$
			which is equal to $\mathrm{span}\{\nu_{9}^3+\ca_{9}\lt_{10}, \mu_{9}+\ell \nu_{9}^{3}, [\iota_{9}, \eta_9]  \}$ for some $\ell\in \mathbb Z\uprightsemicolon{}$
			and $$\begin{aligned}
				&K_2= ((\sigma_{18})^{*})^{-1} ( 2\iota_9 \circ \pi_{25}^9),
			\end{aligned}$$
			which is equal to $\mathrm{span}\{   \nu_{9}^{3}, \eta_9\varepsilon_{10} \}.$
			Therefore, $K_{1}+K_{2}=\pi_{18}^{9}=[\Sigma X_{3}, X_{1}].$ 
			\item[\textup{(2)}] The assertion follows from (1).
		\end{itemize}
	\end{proof}
	\subsection{Unitary Connecting Images and Cofiber Calculations}
	
	We now calculate the unitary connecting images. The cofiber descriptions keep track of the extension
	and coextension choices used in these computations, so the resulting relations can be applied to the
	$G_2$ fibration.
	
	\needspace{3\baselineskip}\begin{lem}
		\label{m3:phi5} Use the symbols in Lemma \ref{m3:tsd}.
		\leavevmode\par
		\begin{itemize}
			
			\item[\textup{(1)}] $\{ \sigma''',\allowbreak \;\eta_{12}\sigma_{13},\allowbreak \;\sigma_{20}\}=0.$
			
			\item[\textup{(2)}]
			$\begin{aligned}
				&\phi_{5} \equiv\{ \nu_{5}^{2},2\iota_{11},\eta_{11}\sigma_{12},\sigma_{19} \}_{1}=\{ \nu_{5}^{2},\; \Sigma^{2}\overline{2\iota_{9}},\;\Sigma^{2}\widetilde{\sigma_{17} } \}_{1}\;\bmod\;4\nu_{5}\bar{\kappa}_{8}.
			\end{aligned}$
			\item[\textup{(3)}]
			$\begin{aligned}
				&\pm\phi' \equiv\{ \varepsilon_{3},2\iota_{11},\eta_{11}\sigma_{12},\sigma_{19} \}_{2}=\{ \varepsilon_{3},\,\Sigma^{2}\overline{2\iota_{9}},\,\Sigma^{2}\widetilde{\sigma_{17}} \}_{2}\;\bmod\; \mu_{3,3}, \nu'\varepsilon_{6}\kappa_{14},\;\eta_{3}\bar{\mu}_{4}\sigma_{21}.
			\end{aligned}$
			
			Thus, the generator $\phi'$ may be chosen so that $\phi'\in\{
			\varepsilon_{3},\allowbreak \,\Sigma^{2}\overline{2\iota_{9}},\allowbreak \,\Sigma^{2}\widetilde{\sigma_{17}} \}_{2}.$
			
		\end{itemize}
	\end{lem}
	\begin{proof}
		\begin{itemize}
			\item[\textup{(1)}] Since $\Sigma^{\infty}\colon\pi_{28+i}^{12+i}\rightarrow \pi_{16}^{S}$ ($i=0,1$) is injective
			and $\langle \sigma^{2},\allowbreak \eta,\allowbreak 2\iota \rangle=0$ by \cite[Proposition 3.3(6), p. 58]{Oda1979},
			we infer
			$\{ 2\iota_{12+i},\allowbreak \eta_{12+i},\allowbreak \sigma_{13+i}^{2} \}=0\subseteq \pi_{28+i}^{12+i}$ ($i=0,1$).
			Observe that  $$\begin{aligned}
				&\Sigma\{ \sigma''',\allowbreak \;\eta_{12}\sigma_{13},\allowbreak \;\sigma_{20} \}\subseteq-\{ 2\sigma'',\allowbreak \;\eta_{13}\sigma_{14},\allowbreak \;\sigma_{21} \}\subseteq\pi_{23+6}^{6},
			\end{aligned}$$
			$\Sigma:\pi_{23+5}^{5}\rightarrow\pi_{23+6}^{6}$ is injective, and $$\begin{aligned}
				&\{ 2\sigma'',\allowbreak \;\eta_{13}\sigma_{14},\allowbreak \;\sigma_{21} \}\supseteq\sigma''\circ\{ 2\iota_{13},\allowbreak \;\eta_{13}\sigma_{14},\allowbreak \;\sigma_{21} \}=0\;\bmod\;0.
			\end{aligned}$$
			Here we use $\zeta'\sigma_{22}=\kappa_{7}\sigma_{21}=0\;$ from \cite[Proposition 2.1(2)(11), pp.~81--82]{Oda1979} to compute $\mathrm{Ind}\{ 2\sigma'',\allowbreak \;\eta_{13}\sigma_{14},\allowbreak \;\sigma_{21} \}=0$. This proves the
			assertion.
			\item[\textup{(2)}] In light of \cite[the proof of (4.6), p. 19 and p. 20]{MMO}, we find $\{
			\nu_{5}^{2},\allowbreak 2\iota_{11},\allowbreak \eta_{11}\sigma_{12}\}_{1}=0.$
			Lemma~\ref{m3:2g0}  gives $\{2\iota_{9},\allowbreak \eta_{9}\sigma_{10},\allowbreak \sigma_{17} \}\ni0.$
			Arguing as in the proof of \cite[(4.6)]{MMO}, we see that
			$$\begin{aligned}
				&\{ \nu_{5}^{2},\allowbreak 2\iota_{11},\allowbreak \eta_{11}\sigma_{12},\allowbreak \sigma_{19} \}_{i}:=\{\nu_{5}^{2}, \Sigma^{2}\overline{2\iota_{9}},\allowbreak \,\Sigma^{2}\widetilde{\sigma_{17}}\}_{i},\;(i\leq2)
			\end{aligned}$$ are well-defined and $H\{
			\nu_{5}^{2},\allowbreak 2\iota_{11},\allowbreak \eta_{11}\sigma_{12},\allowbreak \sigma_{19} \}_{1}=\bar{\sigma}_{9}.$  Let $q\colon C_{\eta_{12}\sigma_{13}}\to S^{21}$ be the pinch map that collapses the bottom sphere $S^{12}.$    We have
			\begin{align*}
				2\iota_{5}\circ\{ \nu_{5}^{2},2\iota_{11},\eta_{11}\sigma_{12},\sigma_{19} \}_{1} {}&\subseteq{} \{2\nu_{5}^{2}, \Sigma^{2}\overline{2\iota_{9}},\allowbreak \,\Sigma^{2}\widetilde{\sigma_{17}}\}_{1} \\
				&={} \{0, \Sigma^{2}\overline{2\iota_{9}},\allowbreak \,\Sigma^{2}\widetilde{\sigma_{17}}\}_{1} \\
				&={} [C_{\eta_{12}\sigma_{13}}, S^5]\circ \Sigma^{3}\widetilde{\sigma_{17}} \\
				&={} \mathrm{span}\{\mu_{5}\sigma_{14}q, \eta_{5}^{2}\kappa_{7}q, \mathrm{ext}_{\eta_{12}\sigma_{13}}(\sigma''') \}\circ \Sigma^{3}\widetilde{\sigma_{17}} \\
				&={} \mathrm{span}\{ \eta_{5}^{2}\kappa_{7}\sigma_{21} ,\mu_{5}\sigma_{14}^{2},0\} ={} 0.
			\end{align*}
			
			The relation $\kappa_{7}\sigma_{21}=0$ is given in \cite[II-Proposition~2.1, pp.~81--82]{Oda1979}. The identities $H(\phi_{5})=\bar{\sigma}_{9}$ and $H ( \bar{\rho}''')=4\bar{\zeta}_{9}$ are given in \cite[(3.2)]{MMO}. We also have $\pi_{28}^{5}=\frac{\mathbb Z}{8}\{\nu_{5}\bar{\kappa}_{8}\}\oplus \frac{\mathbb Z}{2}\{\bar{\rho}'''\}\oplus \frac{\mathbb Z}{2}\{\phi_{5}\}.$ Therefore,
			$$\phi_{5} \equiv\{ \nu_{5}^{2},\allowbreak 2\iota_{11},\allowbreak \eta_{11}\sigma_{12},\allowbreak \sigma_{19}
			\}_{1}\;\bmod\;4\nu_{5}\bar{\kappa}_{8}.$$
			This proves the assertion.
			\item[\textup{(3)}] By Lemma~\ref{m3:tsd},  we conclude that $\{
			\varepsilon_{3},\allowbreak 2\iota_{11},\allowbreak \eta_{11}\sigma_{12},\allowbreak \sigma_{19} \}_2=\{
			\varepsilon_{3},\allowbreak \,\Sigma^{2}\overline{2\iota_{9}},\allowbreak \,\Sigma^{2}\widetilde{\sigma_{17}} \}_{2}.$
			By  assertion (2), we obtain
			\begin{align*}
				H\{ \varepsilon_{3},\,\Sigma^{2}\overline{2\iota_{9}},\,\Sigma^{2}\widetilde{\sigma_{17}} \}_{2}&\subseteq{} \{ \nu_{5}^{2},\,\Sigma^{2}\overline{2\iota_{9}},\,\Sigma^{2}\widetilde{\sigma_{17}} \}_{2} {} \\
				&\subseteq{} \{\nu_{5}^{2},2\iota_{11},\eta_{11}\sigma_{12},\sigma_{19} \}_{1} {} \\
				&\equiv{} \phi_{5}\;\bmod\;4\nu_{5}\bar{\kappa}_{8}.
			\end{align*}
			By \cite[Proposition 6.1(1), p. 66]{Oda1979}, we have $H(\mu_{3,\allowbreak 3})= \bar{\rho}'''$ and $H(\nu'\varepsilon_{6}\kappa_{14})=4\nu_{5}\bar{\kappa}_{8}$. Moreover, \cite[p.~68]{Oda1979} gives $H(\phi') \equiv\phi_{5}\;\bmod\;4\nu_{5}\bar{\kappa}_{8},\allowbreak \;\bar{\rho}'''$, proving the assertion.
		\end{itemize}
	\end{proof}
	
	The cofiber sequence
	\[
	\begin{tikzcd}[column sep=0.85cm]
		S^{20} \arrow[r, "\eta_{12}\sigma_{13}"]
		& S^{12} \arrow[r, "\bar{i}_{12}"]
		& C_{\eta_{12}\sigma_{13}} \arrow[r, "\bar{p}_{12}"]
		& S^{21} \arrow[r, "\eta_{13}\sigma_{14}"]
		& S^{13},
	\end{tikzcd}
	\]
	where $\bar{i}_{12}$ and $\bar{p}_{12}$ are the inclusion and pinch map, respectively, induces the
	exact sequence
	\[
	\begin{tikzcd}[column sep=1cm]
		\pi_{13}(\mathit{SU}(3))
		\arrow[r, "(\eta_{13}\sigma_{14})^{*}"]
		& \pi_{21}(\mathit{SU}(3))
		\arrow[r, "\bar{p}_{12}^{*}"]
		& {[C_{\eta_{12}\sigma_{13}}, \mathit{SU}(3)]}
		\arrow[r, "\bar{i}_{12}^{*}"]
		& \pi_{12}(\mathit{SU}(3))
		\arrow[r, "(\eta_{12}\sigma_{13})^{*}"]
		& \pi_{20}(\mathit{SU}(3)),
	\end{tikzcd}
	\]
	where $j:S^{3}\hookrightarrow \mathit{SU}(3)$ is the inclusion. We obtain the following result:
	\needspace{3\baselineskip}\begin{lem}
		\label{m3:tlj}
		\leavevmode\par
		\begin{itemize}
			\item[\textup{(1)}]
			$\FitMath{\mathrm{ Im}\bar{p}_{12}^{*} =\mathrm{span}\{j\mu'\sigma_{14}\bar{p}_{12}\}\cong \mathbb Z/2,\quad \operatorname{ker}(\eta_{12}\sigma_{13})^{*}=\mathrm{span}\{[\sigma''']\}\cong \mathbb Z/4.}$
			\item[\textup{(2)}]
			$\GroupLayout{[ C_{\eta_{12}\sigma_{13}} , \mathit{SU}(3)] =\mathrm{span}\{\overline{[\sigma'''] },\; j\mu'\sigma_{14}\bar{p}_{12}\}\cong\mathbb Z/8\;or\;\mathbb Z/4 \oplus\mathbb Z/2.}{\begin{aligned}[ C_{\eta_{12}\sigma_{13}} , \mathit{SU}(3)] =\mathrm{span}\{&\overline{[\sigma'''] },\\
					&j\mu'\sigma_{14}\bar{p}_{12}\}\cong\mathbb Z/8\;or\;\mathbb Z/4 \oplus\mathbb Z/2.\end{aligned}}$
		\end{itemize}
		Here $\,\overline{[\sigma'''] }\;\in\{\,[\sigma'''],\allowbreak \;\eta_{12}\sigma_{13}
		,\allowbreak \;\bar{p}_{11}  \}$ has order 4 or 8, and $\bar{p}_{11}:C_{\eta_{11}\sigma_{12}}\rightarrow
		S^{20}$ is the pinch map.
	\end{lem}
	\begin{proof}
		\begin{itemize}
			\item[\textup{(1)}] According to \cite[Theorem 4.1]{MT2}, $\pi_{13}(\mathit{SU}(3))=\frac{\mathbb Z}{2}\{ j\varepsilon' \}$ and $\pi_{21}(\mathit{SU}(3))=\frac{\mathbb Z}{2}\{ j\mu'\sigma_{14} \}$. We also have $\varepsilon'\sigma_{13}=2\bar{\varepsilon}'.$ Then $$(\eta_{13}\sigma_{14})^{*} (j\varepsilon' )= j\varepsilon' \circ\sigma_{13}\eta_{20}=j(\varepsilon'
			\sigma_{13})\eta_{20}=2j\bar{\varepsilon}'\eta_{20}=0.$$
			The same result, \cite[Theorem 4.1]{MT2}, also gives $\pi_{12}(\mathit{SU}(3))=\frac{\mathbb Z}{4}\{[\sigma'''] \}$ and $\pi_{20}(\mathit{SU}(3))=\frac{\mathbb Z}{4}\{ [\rho^{\mathrm{IV}}] \}\oplus \frac{\mathbb Z}{2}\{ j\bar{\varepsilon}'\}$.
			Thus $[\sigma''']\eta_{12}=x (j\varepsilon')$ for some $x\in\mathbb Z,$ and consequently
			\[\begin{aligned}
				&(\eta_{12}\sigma_{13})^{*}([\sigma''']) =[\sigma''']\eta_{12}\sigma_{13} =x (j\varepsilon')\sigma_{13} =x( j(\varepsilon'\sigma_{13})) =2x j\bar{\varepsilon}' =0.
			\end{aligned}\]
			This proves the
			assertion.
			\item[\textup{(2)}] By (1), it remains to show that there exists $\,\overline{[\sigma'''] }\;\in\{\,[\sigma'''],\allowbreak \;\eta_{12}\sigma_{13} ,\allowbreak \;\bar{p}_{11}  \}$ such that
			$$\bar{i}_{12}^{*}(\,\overline{[\sigma'''] }\,)= \pm[\sigma'''].$$ First,
			$(\eta_{12}\sigma_{13})^{*}([\sigma'''])=0$ and the cofiber sequence associated with
			$\eta_{11}\sigma_{12}$ show that the Toda bracket $\{\,[\sigma'''],\allowbreak \;\eta_{12}\sigma_{13}
			,\allowbreak \;\bar{p}_{11}  \}$ is well-defined. Hence
			\begin{align*}
				\bar{i}_{12}^{*}\{[\sigma'''],\;\eta_{12}\sigma_{13},\;\bar{p}_{11}\;\}&={} \{[\sigma'''],\;\eta_{12}\sigma_{13},\;\bar{p}_{11}\;\}\circ\bar{i}_{12} {} \\
				&={} -[\sigma''']\circ\{\eta_{12}\sigma_{13},\;\bar{p}_{11},\;i_{11}\} {} \\
				&\ni{} -[\sigma''']\circ\mathrm{id} \\
				&=-[\sigma'''],
			\end{align*}
			where $i_{11}:S^{11}\to C_{\eta_{11}\sigma_{12}}$ is the inclusion in this cofiber sequence and $\Sigma i_{11}=\bar{i}_{12}.$  Therefore,
			there exists $\,\overline{[\sigma'''] }\;\in\{\,[\sigma'''],\allowbreak \;\eta_{12}\sigma_{13}
			,\allowbreak \;\bar{p}_{11}  \}$ such that $\bar{i}_{12}^{*}(\,\overline{[\sigma'''] }\,)= -[\sigma'''].$
		\end{itemize}
	\end{proof}
	\needspace{3\baselineskip}\begin{lem}
		\label{m3:inp15}
		\leavevmode\par
		\begin{itemize}
			\item[\textup{(1)}] In $\pi_{15}( \mathit{SU}(3)),$
			$[\sigma''']\circ\nu_{12}=2[2\iota_{5}]\nu_{5}\sigma_{8}.$
			\item[\textup{(2)}] In $\pi_{28}^{3},$ $\mu'\sigma_{14}^{2}
			\equiv0\;\bmod\;\nu'\varepsilon_{6}\kappa_{14}.$
		\end{itemize}
	\end{lem}
	\begin{proof}
		\begin{itemize}
			\item[\textup{(1)}] We use $\pi_{15}( \mathit{SU}(3))=\frac{\mathbb Z}{4}\{ [2\iota_{5}]\nu_{5}\sigma_{8} \}$ from \cite[Theorem 4.1]{MT2}, together with $\sigma'''\circ\nu_{12}=4\nu_{5}\sigma_{8}.$ 
			Write $[\sigma''']\nu_{12}=x[2\iota_{5}]\nu_{5}\sigma_{8}$ for some $x\in\mathbb Z.$
			Applying $p_{*},$
			we obtain $4\nu_{5}\sigma_{8}=2x\nu_{5}\sigma_{8}.$ Since $\mathrm{ord}(\nu_{5}\sigma_{8})=8,$ we find $2x
			\equiv4\;\bmod\;8$ and hence $x \equiv2\;\bmod\;4.$
			\item[\textup{(2)}] Observe that $\Sigma^{2}\phi' \equiv2\phi''\;\bmod\;\nu_{5}\varepsilon_{8}\kappa_{16},\allowbreak \nu_{5}^{2}\bar{\sigma}_{11},\allowbreak \mu_{3,\allowbreak 5},\allowbreak \eta_{5}\bar{\mu}_{6}\sigma_{23}$ by \cite[(7.6), p. 69]{Oda1979},
			and $\Sigma^{2}\mu'=\pm2 \zeta_{5}.$ Hence $\mu'\sigma_{14}^{2}\in \mathrm{Ker }(\Sigma^{2})=\frac{\mathbb Z}{2}\{
			\nu'\varepsilon_{6}\kappa_{14} \},$ proving the assertion.
		\end{itemize}
	\end{proof}
	
	Recall that $\widetilde{\sigma_{17}}\in\pi_{25}(C_{\eta_{9}\sigma_{10}})$ is the coextension of
	$\sigma_{17}.$
	\needspace{3\baselineskip}\begin{lem}
		\label{m3:ghc} $[ C_{\eta_{12}\sigma_{13}} ,\allowbreak  \mathit{SU}(3) ] \circ
		\Sigma^{3}\widetilde{\sigma_{17}} \equiv0\;\bmod\;2[\bar{\rho}'''].$
	\end{lem}
	\begin{proof}
		Lemma~\ref{m3:tlj}(2) gives $[ C_{\eta_{12}\sigma_{13}} ,\allowbreak  \mathit{SU}(3) ]
		=\mathrm{span}\{\,\overline{[\sigma'''] }\,,\allowbreak \; j\mu'\sigma_{14}\bar{p}_{12} \},$ where $$\;
		\overline{[\sigma'''] }\;\in\{\,[\sigma'''],\allowbreak \;\eta_{12}\sigma_{13} ,\allowbreak \;\bar{p}_{11}  \}.$$
		It therefore suffices to show that $$\{[\sigma'''],\allowbreak \;\eta_{12}\sigma_{13} ,\allowbreak \;\bar{p}_{11}
		\}\circ \Sigma^{3}\widetilde{\sigma_{17}} \equiv0\;\bmod\;2[\bar{\rho}''']=j_{*}\mu_{3,\allowbreak 3},\text{and that }\;j\mu'\sigma_{14}\bar{p}_{12}\circ \Sigma^{3}\widetilde{\sigma_{17}}=0.$$ \noindent Consider  the exact sequence\vspace{-1em}
		\[
		\xymatrix@C=0.5cm{
			\pi_{28}^{3} \ar[r]^{j_{*}\qquad} & \pi_{28}(\mathit{SU}(3)) \ar[r]^{\qquad p_{*}} & \pi_{28}^{5}.
		}
		\]
		We know that	$\operatorname{Im} j_{*}=\frac{\mathbb Z}{2}\{j_{*}\mu_{3, 3}\}\oplus \frac{\mathbb Z}{2}\{j_{*}\phi'\}$ by the
		fifth row of \cite[Proposition 3.2]{Oda1978}, and $j_{*}\mu_{3,\allowbreak 3}=2[\bar{\rho}''']$ by \cite[Theorem
		B)]{Oda1978}.
		Since $\bar{p}_{11}\circ \Sigma^{2}\widetilde{\sigma_{17}}\in\pi_{27}^{20},$ we conclude that
		$\bar{p}_{11}\circ \Sigma^{2}\widetilde{\sigma_{17}}=t\sigma_{20}$ for some $t\in\mathbb Z.$
		Consequently,
		\begin{align*}
			&\{[\sigma'''],\;\eta_{12}\sigma_{13} ,\;\bar{p}_{11} \}\circ \Sigma^{3}\widetilde{\sigma_{17}} {}\subseteq{} - \{[\sigma'''],\;\eta_{12}\sigma_{13} , \bar{p}_{11} \Sigma^{2}\widetilde{\sigma_{17}}\;\}=- \{ [\sigma'''],\;\eta_{12}\sigma_{13},\;t\sigma_{20}\}.
		\end{align*}
		In light of Lemma~\ref{m3:phi5}(1), the bracket $\{\sigma''',\eta_{12}\sigma_{13},\sigma_{20}\}$
		contains only zero. For every integer $t,$   composition with $t\iota_{28}$ gives
		$0=\{\sigma''',\eta_{12}\sigma_{13},t\sigma_{20}\}.$ 	Consequently,
		\begin{align*}
			&p_{*}\{\,[\sigma'''],\;\eta_{12}\sigma_{13} ,\;\bar{p}_{11} \}\circ \Sigma^{3}\widetilde{\sigma_{17}} {}\subseteq{} - \{ p_{*}[\sigma'''],\;\eta_{12}\sigma_{13} , t\sigma_{20}\}= \{ \sigma''',\;\eta_{12}\sigma_{13} ,\;t\sigma_{20}\}=0.
		\end{align*} Hence,
		\[\begin{aligned}
			&\{\,[\sigma'''], \; \eta_{12}\sigma_{13} , \; \bar{p}_{11} \}\circ \Sigma^{3}\widetilde{\sigma_{17}} \subseteq \operatorname{Im}(j_{*}) =\frac{\mathbb Z}{2}\{2[\bar{\rho}''']\}\oplus \frac{\mathbb Z}{2}\{j_{*}\phi'\}.
		\end{aligned}\]
		Now, suppose $x(2[\bar{\rho}'''])
		+y (j_{*}\phi' ) \in\{\,[\sigma'''],\allowbreak \;\eta_{12}\sigma_{13} ,\allowbreak \;\bar{p}_{11}  \}\circ
		\Sigma^{3}\widetilde{\sigma_{17}},$ where $x,\allowbreak y\in\mathbb Z.$
		Then,
		\begin{align*}
			&x(2[\bar{\rho}''']) +y (j_{*}\phi' ) \in \{\,[\sigma'''],\;\eta_{12}\sigma_{13} ,\;\bar{p}_{11} \}\circ \Sigma^{3}\widetilde{\sigma_{17}}\subseteq- \{ \,[\sigma'''], \;\eta_{12}\sigma_{13} ,\;t\sigma_{20} \}.
		\end{align*}
		Thus $(x(2[\bar{\rho}''']) +y (j_{*}\phi' ))\circ\nu_{28}^{2} \in- \{
		\,[\sigma'''],\allowbreak \;\eta_{12}\sigma_{13} ,\allowbreak \;t\sigma_{20} \}\circ\nu_{28}^{2}.$ It follows that
		\begin{align*}
			y(j_{*}\phi'\nu_{28}^{2})&\in{} [\sigma''']\circ\{\eta_{12}\sigma_{13},\;t\sigma_{20},\;\nu_{27}\}\circ\nu_{31} {} \\
			&={} [\sigma''']\circ t\bar{\sigma}_{12}\circ\nu_{31} {} \\
			&={} t([\sigma''']\nu_{12})\bar{\sigma}_{15} {}={} 2t[2\iota_{5}]\nu_{5}\sigma_{8}\circ\bar{\sigma}_{15} {} \\
			&={} t[2\iota_{5}]\nu_{5}\sigma_{8}\circ2\bar{\sigma}_{15}=0,
		\end{align*} by Lemma~\ref{m3:inp15}(1) and $\{\eta_{12}\sigma_{13}
		,\allowbreak \;\sigma_{20}\;,\allowbreak \nu_{27}\}=\bar{\sigma}_{12}\allowbreak \uprightsemicolon{}$ see \cite[Proposition 3.4, p. 59]{Oda1979}. By \cite[Theorem 4.1]{Oda1978}, $j_{*}\phi'\nu_{28}^{2}\in\pi_{34}(\mathit{SU}(3))$ has order 2. Therefore, $y \equiv0\;\bmod\;2.$ Thus, $\{\,[\sigma'''],\allowbreak \;\eta_{12}\sigma_{13} ,\allowbreak \;\bar{p}_{11}  \}\circ \Sigma^{3}\widetilde{\sigma_{17}}\subseteq
		\mathrm{span}\{ \,2[\bar{\rho}'''] \; \}.$
		Hence,
		$$\{\,[\sigma'''],\allowbreak \;\eta_{12}\sigma_{13} ,\allowbreak \;\bar{p}_{11}  \}\circ
		\Sigma^{3}\widetilde{\sigma_{17}} \equiv0\;\bmod\;2[\bar{\rho}''']=j_{*}\mu_{3,\allowbreak 3}.$$\noindent
		Next, we show that $j\mu'\sigma_{14}\bar{p}_{12}\circ \Sigma^{3}\widetilde{\sigma_{17}}=0.$
		Clearly, $\bar{p}_{12}\circ \Sigma^{3}\widetilde{\sigma_{17}}=r\sigma_{21}$ for  $r=\pm1.$
		Lemma~\ref{m3:inp15}(2) gives $\mu'\sigma_{14}^{2}=\ell\nu'\varepsilon_{6}\kappa_{14}$ for some
		$\ell\in\mathbb Z.$
		By \cite[Proposition 3.3]{Oda1978}, $\nu'\varepsilon_{6}\kappa_{14}$ belongs to the image of $\Delta\uprightsemicolon{}$
		equivalently, $j_{*}(\nu'\varepsilon_{6}\kappa_{14})=0.$
		Therefore, \[
		\begin{aligned}
			&j\mu'\sigma_{14}\bar{p}_{12}\circ \Sigma^{3}\widetilde{\sigma_{17}} =j\mu'\sigma_{14}\circ r\sigma_{21} =rj\mu'\sigma_{14}^{2} =(r\ell) j\nu'\varepsilon_{6}\kappa_{14} =(r\ell) j_{*}(\nu'\varepsilon_{6}\kappa_{14}) = 0.
		\end{aligned}
		\]
		This proves the assertion.
	\end{proof}
	
	For the exact sequence\vspace{-1em}
	$$\xymatrix@C=0.42cm{
		\pi_{29}^{6} \ar[r]^{\Delta\qquad\,} & \pi_{28}(\mathit{SU}(3)) \ar[r]^{\quad i_{*}} &
		\pi_{28}(G_{2}) \ar[r]^{\quad p_{*}} & \pi_{28}^{6} \ar[r]^{\Delta\qquad\,} &
		\pi_{27}(\mathit{SU}(3))    },$$
	
	\noindent we will use the following group structures; the unitary groups are given in \cite{Mi2} and \cite{Oda1978}:
	\begin{gather*}
		\begin{aligned}
			&\pi_{29}^{6}={}\mathrm{span}\{ \nu_{6}\bar{\kappa}_{9},\bar{\rho}'' ,\phi_{6}\}\oplus \mathrm{span}\{ P(\lambda),P(\xi_{13}) \}\cong(\mathbb Z/8 \oplus\mathbb Z/4\oplus\mathbb Z/2) \oplus(\mathbb Z/8\oplus\mathbb Z/4),
		\end{aligned}
		\\[3pt]
		\begin{aligned}
			&\pi_{28}(\mathit{SU}(3))={}\frac{\mathbb Z}{4}\{[\bar{\rho}''']\}\oplus \frac{\mathbb Z}{4}\{[2\iota_{5}]\nu_{5}\bar{\kappa}_{8}\} \oplus \frac{\mathbb Z}{2}\{j_{*}\phi'\},
		\end{aligned}
	\end{gather*} where our notation $j_{*}$ corresponds to $i_{*}$ in \cite{Oda1978}.
	\[\begin{aligned}
		&\pi_{28}^{6}={}\frac{\mathbb Z}{4}\{\rho'''\sigma_{21}\}\oplus \frac{\mathbb Z}{2}\{\bar{\nu} _{6}\kappa_{14}\} \oplus \frac{\mathbb Z}{2}\{\varepsilon_{6}\kappa_{14}\} \oplus \frac{\mathbb Z}{2}\{\nu_{6}\bar{\sigma}_{9}\}.\allowbreak \end{aligned}\] We will use $\varepsilon_{6}\kappa_{14}=\eta^{2}_{6}\bar{\kappa}_{8}$ from \cite[I-Proposition~3.1(2), p.~54]{Oda1979} and
	$$\pi_{27}(\mathit{SU}(3))=\frac{\mathbb Z}{4}\{[\rho^{\mathrm{IV}}] \sigma_{20}\}\oplus \frac{\mathbb Z}{2}\{[\nu_{5}\bar{\sigma}_{8}]\}.$$

	\subsection{Connecting Images in the $G_2$ Fibration}
	
	The preceding unitary relations supply the target information for the $G_2$ boundary map. We use them
	to identify its images and hence the elements of homotopy groups of spheres for which lifts can be chosen.
	
	\needspace{3\baselineskip}\begin{prop}
		\label{m3:drgpp}
		\leavevmode\par
		\begin{itemize}
			\item[\textup{(1)}]
			$\begin{aligned}
				&\Delta(\bar{\rho}'') \equiv\pm[\bar{\rho}''']\;\bmod\; j_{*}\phi';\;\,\Delta( \nu_{6}\bar{\kappa}_{9}) =[2\iota_{5}]\nu_{5}\bar{\kappa}_{8}.\;
			\end{aligned}$
			\item[\textup{(2)}]
			$\begin{aligned}
				&\Delta( \phi_{6}) \equiv j_{*}\phi' \;\bmod\;2[\bar{\rho}'''].
			\end{aligned}$
			\item[\textup{(3)}]
			$\begin{aligned}
				&\Delta(\rho'''\sigma_{21}) =[\rho^{\mathrm{IV}}] \sigma_{20},
			\end{aligned}$
			
			$\Delta(\mathrm{span}\{ \bar{\nu}
			_{6}\kappa_{14},\allowbreak \varepsilon_{6}\kappa_{14},\allowbreak \nu_{6}\bar{\sigma}_{9} \}\;) \equiv
			0\;\bmod\;2[\rho^{\mathrm{IV}}] \sigma_{20}.$
		\end{itemize}
	\end{prop}
	
	\begin{proof}
		\begin{itemize}
			\item[\textup{(1)}] We find $\Delta(\iota_{6})=[2\iota_{5}],\allowbreak $ and \cite[p. 21]{MMO} gives
			$2\bar{\rho}''=\Sigma\bar{\rho}'''.$ Thus,
			$$\begin{aligned}
				&2\Delta( \bar{\rho}'')=\Delta( 2 \bar{\rho}'')=\Delta(\Sigma\bar{\rho}''' )=[2\iota_{5}]\bar{\rho}'''.
			\end{aligned}$$
			Examining the image under $p_{*}$ gives
			$\Delta(\bar{\rho}'') \equiv\pm[\bar{\rho}''']\;\bmod\; j_{*}\phi' .\;$ It is obvious that
			$\Delta( \nu_{6}\bar{\kappa}_{9})=[2\iota_{5}]\nu_{5}\bar{\kappa}_{8}\;$ .
			\item[\textup{(2)}]  Lemma~\ref{m3:phi5} gives
			\begin{align*}
				&\phi_{5} \equiv\{ \nu_{5}^{2},2\iota_{11},\eta_{11}\sigma_{12},\sigma_{19} \}_{1}=\{ \nu_{5}^{2},\; \Sigma^{2}\overline{2\iota_{9}},\;\Sigma^{2}\widetilde{\sigma_{17} } \}_{1}\;\bmod\;4\nu_{5}\bar{\kappa}_{8} \\
				&\phi' \in\{ \varepsilon_{3},2\iota_{11},\eta_{11}\sigma_{12},\sigma_{19} \}_{2}=\{ \varepsilon_{3},\,\Sigma^{2}\overline{2\iota_{9}},\,\Sigma^{2}\widetilde{\sigma_{17}} \}_{2}.
			\end{align*}
			Recall that $\{ j\varepsilon_{3}
			,\allowbreak 2\iota_{11},\allowbreak \eta_{11}\sigma_{12},\allowbreak \sigma_{19} \}_{2}$ is defined in Lemma \ref{m3:tsd}, where $j:S^{3}\hookrightarrow
			\mathit{SU}(3)$ is the inclusion. Hence
			\[\begin{aligned}
				&\Delta(\phi_{6}) =[2\iota_{5}]\phi_{5} \equiv \{ [2\iota_{5}] \nu_{5}^{2}, \; \Sigma^{2}\overline{2\iota_{9}}, \; \Sigma^{2}\widetilde{\sigma_{17} } \}_{1}\; \bmod\; 4\Delta(\nu_{6}\bar{\kappa}_{9}) =0.
			\end{aligned}\]
			Since $[2\iota_{5}]\nu_{5}^{2}=
			j_{*}\varepsilon_{3}=2[\nu_{5}^{2}]=\Delta(\nu_{6}^{2})$ \cite[p. 162]{Mi2}, we obtain
			\begin{align*}
				\Delta(\phi_{6})&\in{} \{ 2[\nu_{5}^{2}],\; \Sigma^{2}\overline{2\iota_{9}},\;\Sigma^{2}\widetilde{\sigma_{17} } \}_{1} {} \\
				&={} \{ j\varepsilon_{3},\; \Sigma^{2}\overline{2\iota_{9}}, \;\Sigma^{2}\widetilde{\sigma_{17} } \}_{1} {} \\
				&\supseteq{} j\circ\{ \varepsilon_{3},\; \Sigma^{2}\overline{2\iota_{9}},\;\Sigma^{2}\widetilde{\sigma_{17} } \}_{1} {}\ni{} j\phi'.
			\end{align*}
			In light of Lemma~\ref{m3:ghc},
			\begin{align*}
				\mathrm{Ind}\{ 2[\nu_{5}^{2}],\; \Sigma^{2}\overline{2\iota_{9}},\;\Sigma^{2}\widetilde{\sigma_{17} } \}_{1}&= 2[\nu_{5}^{2}]\circ \Sigma\pi_{27}^{10}+[ C_{\eta_{12}\sigma_{13}}, \;\mathit{SU}(3)]\circ \Sigma^{3}\widetilde{\sigma_{17} } {} \\
				&={} 0+[ C_{\eta_{12}\sigma_{13}},\;\mathit{SU}(3)]\circ \Sigma^{3}\widetilde{\sigma_{17} } {} \\
				&\subseteq{} \mathrm{ span} \{ 2[\bar{\rho}''']\}.
			\end{align*} This proves the assertion.
			\item[\textup{(3)}] Toda's bracket description in \cite[p.~103]{T} gives $\rho'''\in\{\sigma'',\allowbreak 4\iota_{13},\allowbreak 4\sigma_{13}\}_{1}$.
			Then, by Lemma~\ref{m3:fp20},
			\begin{align*}
				\Delta(\rho'''\sigma_{21})&\in{} \Delta(\{\sigma'',4\iota_{13},4\sigma_{13}\}_{1}\circ\sigma_{21}) {} \\
				&\subseteq{} \{ \Delta\sigma'',4\iota_{12}, 4\sigma_{12}\}_{1}\circ\sigma_{20} {} \\
				&={} \{\; [\sigma'''],4\iota_{12},4\sigma_{12}\;\}_{1}\circ\sigma_{20} {}\ni{} [\rho^{\mathrm{IV}}] \sigma_{20}.
			\end{align*}
			\[\begin{aligned}
				&\mathrm{ Ind}(\{\; [\sigma'''], 4\iota_{12}, 4\sigma_{12}\; \}_{1}\circ\sigma_{20}) =\pi_{13}(\mathit{SU}(3))\circ4\sigma_{13}^{2} +[\sigma''']\circ\pi_{20}^{12}\circ\sigma_{20} = 0,
			\end{aligned}\]
			where we use $\pi_{13}(\mathit{SU}(3))\cong\mathbb Z/2$ as computed in \cite{MT2}. Thus,
			$\Delta(\rho'''\sigma_{21})= [\rho^{\mathrm{IV}}] \sigma_{20}.$
			For the second assertion, $\Delta( \bar{\nu} _{6})=j_{*}\varepsilon'$ by \cite[Proposition
			6.3]{Mi2}\uprightsemicolon{} our notation $j_{*}$ corresponds to $i_{*}$ in that reference.
			Together with $\Delta(\iota_{6})=[2\iota_{5}],$   this gives
			\[\begin{aligned}
				&\Delta(\mathrm{span}\{ \bar{\nu} _{6}\kappa_{14}, \varepsilon_{6}\kappa_{14}, \nu_{6}\bar{\sigma}_{9} \}) = \mathrm{span}\{ j_{*}\varepsilon'\kappa_{13}, [2\iota_{5}]\varepsilon_{5}\kappa_{13}, [2\iota_{5}]\nu_{5}\bar{\sigma}_{8} \}.
			\end{aligned}\]
			For any $x\in\mathbb Z,$ examining the image under $p_{*}$ gives $j_{*}\varepsilon'\kappa_{13}\neq
			[\nu_{5}\bar{\sigma}_{8}]+x [\rho^{\mathrm{IV}}] \sigma_{20},$ since $p_{*}j_{*}=0.$ It follows that
			$j_{*}\varepsilon'\kappa_{13} \equiv0\;\bmod\;2[\rho^{\mathrm{IV}}] \sigma_{20}.$ Similarly,
			$$[2\iota_{5}]\varepsilon_{5}\kappa_{13},\allowbreak \,[2\iota_{5}]\nu_{5}\bar{\sigma}_{8} \equiv
			0\;\bmod\;2[\rho^{\mathrm{IV}}] \sigma_{20}.$$\noindent Hence $\Delta(\mathrm{span}\{ \bar{\nu}
			_{6}\kappa_{14},\allowbreak \varepsilon_{6}\kappa_{14},\allowbreak \nu_{6}\bar{\sigma}_{9} \}) \equiv
			0\;\bmod\;2[\rho^{\mathrm{IV}}] \sigma_{20}.$
		\end{itemize}
	\end{proof}
	\subsection{Chosen Lifts and the Degree 28 Decomposition}
	
	Having determined the relevant connecting images, we use the lift $\langle\bar\sigma_6\rangle$
	chosen in Lemma~\ref{m3:fh} to determine the degree-$28$ decomposition.
	
	\needspace{3\baselineskip}\begin{prop}
		\label{m3:p28g2} The structure of $\pi_{28}(G_{2})$ is as follows.
		\[ \GroupLayout{\pi_{28}(G_{2}) = \mathrm{span}\{\langle \bar{\nu} _{6}+\varepsilon_{6}\rangle\kappa_{14},\; \langle \bar{\sigma}_{6} \nu_{25}\rangle,\; \langle \eta_{6}^{2}\rangle\bar{\kappa}_{8}\}\cong (\mathbb Z/2)^{3}.}{\begin{aligned}\pi_{28}(G_{2}) = \mathrm{span}\{&\langle \bar{\nu} _{6}+\varepsilon_{6}\rangle\kappa_{14},\; \langle \bar{\sigma}_{6} \nu_{25}\rangle,\\
				&\langle \eta_{6}^{2}\rangle\bar{\kappa}_{8}\}\cong (\mathbb Z/2)^{3}.\end{aligned}} \]
	\end{prop}
	\begin{proof}
		By Proposition~\ref{m3:drgpp}(3), we have $$\operatorname{Im}(\Delta\colon  \pi_{28}^{6}\rightarrow
		\pi_{27}(\mathit{SU}(3)))=\frac{\mathbb Z}{4}\{ [\rho^{\mathrm{IV}}] \sigma_{20} \},\allowbreak $$ \noindent
		and $\pi_{28}^{6}\cong\mathbb Z/4\oplus(\mathbb Z/2)^{3}.$
		Consequently, $\operatorname{ker}(\Delta\colon \pi_{28}^{6}\rightarrow
		\pi_{27}(\mathit{SU}(3)))\cong(\mathbb Z/2)^{3}.$
		By Proposition~\ref{m3:drgpp}(1),(2), we conclude that $\operatorname{coker}(\Delta\colon  \pi_{29}^{6}\rightarrow
		\pi_{28}(\mathit{SU}(3))\;)=0,$ so, $$\begin{aligned}
			&\pi_{28}(G_{2})\cong \operatorname{ker}(\Delta\colon  \pi_{28}^{6}\rightarrow \pi_{27}(\mathit{SU}(3))\;)\cong(\mathbb Z/2)^{3}.
		\end{aligned}$$
		
		Next, we use the following facts.
		In view of  Mimura's definition of $\langle \;-\;\rangle$ \cite[p. 160]{Mi2},
		the elements $\langle \bar{\nu} _{6}+\varepsilon_{6}\rangle$ and $\langle \eta_{6}^{2}\rangle$ are
		well-defined; see \cite[Theorem 6.1]{Mi2}. Using the lift chosen in Lemma~\ref{m3:fh},
		we define $\langle\bar{\sigma}_{6}\nu_{25}\rangle:=\langle\bar{\sigma}_{6}\rangle\nu_{25}\in \pi_{28}(G_{2}),$
		which satisfies $p_{*}\langle\bar{\sigma}_{6}\nu_{25}\rangle=\bar{\sigma}_{6}\nu_{25}.$
		Consider  the exact sequence
		\[
		\xymatrix@C=0.6cm{
			\pi_{28}(G_{2}) \ar[r]^{\quad p_{*}} & \pi_{28}^{6} \ar[r]^{} &
			\pi_{27}(\mathit{SU}(3)) .}
		\]
		The group has the following descriptions:
		\begin{align}
			\pi_{28}^{6}
			&= \mathrm{span}\{
			\rho'''\sigma_{21},\;
			\bar{\nu}_{6}\kappa_{14},\;
			\varepsilon_{6}\kappa_{14},\;
			\nu_{6}\bar{\sigma}_{9}
			\} \notag \\
			&= \mathrm{span}\{
			\rho'''\sigma_{21},\;
			(\bar{\nu}_{6}+\varepsilon_{6})\kappa_{14},\;
			\eta_{6}^{2}\bar{\kappa}_{8},\;
			\bar{\sigma}_{6}\nu_{25}
			\} \notag \\
			&\cong \mathbb Z/4\oplus(\mathbb Z/2)^{3}
			\uprightsemicolon{} \notag
		\end{align}
		here we use $\varepsilon_{6}\kappa_{14}=\eta^{2}_{6}\bar{\kappa}_{8}\;$ from \cite[Proposition 3.1(2), p.~54]{Oda1979} and $\nu_{6}\bar{\sigma}_{9}=\bar{\sigma}_{6}\nu_{25}$ from \cite[Proposition 3.5(1), p. 60]{Oda1979}.
		The three elements $(\bar{\nu} _{6}+\varepsilon_{6})\kappa_{14} ,\allowbreak \; \bar{\sigma}_{6} \nu_{25}$ and $ \eta_{6}^{2}\bar{\kappa}_{8}$ therefore
		belong to
		\[\begin{aligned}
			&\operatorname{Im}(p_{*}: \pi_{28}(G_{2})\rightarrow \pi_{28}^{6} ) =\operatorname{ker}(\Delta\colon  \pi_{28}^{6}\rightarrow \pi_{27}(\mathit{SU}(3))) \cong(\mathbb Z/2)^{3}.
		\end{aligned}\]
		If $a \langle \bar{\nu} _{6}+\varepsilon_{6}\rangle\kappa_{14} +b \langle \bar{\sigma}_{6}
		\nu_{25}\rangle +c \langle \eta_{6}^{2}\rangle\bar{\kappa}_{8}=0$ , $a,\allowbreak b,\allowbreak c\in\mathbb Z,$ applying
		$p_{*}$ gives
		$$\begin{aligned}
			&a( \bar{\nu} _{6}+\varepsilon_{6})\kappa_{14} +b \bar{\sigma}_{6} \nu_{25}+c \eta_{6}^{2}\bar{\kappa}_{8}=0\in\pi_{28}^{6},\allowbreak
		\end{aligned}$$
		and thus, $a \equiv b \equiv c \equiv0\;\bmod\;2.$
		This shows that
		$\mathrm{span}\{\langle \bar{\nu} _{6}+\varepsilon_{6}\rangle\kappa_{14} \}+\mathrm{span}\{ \langle \bar{\sigma}_{6}
		\nu_{25}\rangle\} +\mathrm{span}\{ \langle \eta_{6}^{2}\rangle\bar{\kappa}_{8} \}$
		is a direct sum isomorphic to $(\mathbb Z/2)^{3}.$ This proves the assertion.
	\end{proof}
	
	\section{The Cases $7\leq n\leq11$}
	
	We now transfer the $G_2$ information to the special orthogonal groups and use spin-group comparisons where they simplify the calculation.
	
	\subsection{The Cases $n=7,8$}
	
	The comparison of fibrations carries a lift in $G_2$ to a lift in $\mathit{SO}(7).$  The next lemma
	states this implication explicitly before it is applied to the generators just constructed.
	
	\needspace{3\baselineskip}\begin{lem}
		\label{m3:gjh} If $\alpha\in\pi_{*}^{6}$ and $\langle \alpha\rangle\in\pi_{*}(G_{2})$ is
		well-defined, then $[ \alpha]\in\pi_{*}(\mathit{SO}(7))$ is also well-defined. Equivalently, for
		$\alpha\in\pi_{*}^{6},$ $\Delta_{G_{2}}(\alpha)=0$ implies $\Delta_{\mathit{SO}(7)}(\alpha)=0,$
		where
		$\Delta_{G_{2}}: \pi_{n+1}^{6}\rightarrow\pi_{n}(\mathit{SU}(3))$ and $\Delta_{\mathit{SO}(7)}:
		\pi_{n+1}^{6}\rightarrow\pi_{n}(\mathit{SO}(6))$ are the boundary homomorphisms.
		\qed
	\end{lem}
	
	The exact sequence  \[
	\pi_{29}^{6}
	\longrightarrow \pi_{28}(\mathit{SO}(6))
	\xrightarrow{i_*} \pi_{28}(\mathit{SO}(7))
	\xrightarrow{p_*} \pi_{28}^{6}
	\longrightarrow \pi_{27}(\mathit{SO}(6))
	\]
	gives rise to the short exact sequence
	\[
	0 \longrightarrow \operatorname{Im}(i_*)
	\longrightarrow \pi_{28}(\mathit{SO}(7))
	\xrightarrow{p_*} \operatorname{Im}(p_*)
	\longrightarrow 0.
	\]
	Recall that
	$$\begin{aligned}
		&\pi_{28}(\mathit{SO}(7))=\mathcal{T}_{*}\pi_{28}(G_{2})\jia[\eta_{6}]\circ\pi_{28}^{7}\cong\pi_{28}(G_{2})\oplus\pi_{28}^{7},
	\end{aligned}$$
	where $\mathcal{T}:\,G_{2}\rightarrow \mathit{SO}(7)$ is the composition
	\[
	G_2
	\xrightarrow[\text{inclusion}]{i'}
	\mathit{Spin}(7)
	\xrightarrow[\text{covering}]{q}
	\mathit{SO}(7).
	\]
	We assert that
	\needspace{3\baselineskip}\begin{lem}
		\label{m3:i37}
		$[\iota_{3}]_{7}\phi'=[\bar{\rho}''']_7=0\in\pi_{28}(\mathit{SO}(7)).$
	\end{lem}
	\begin{proof}
		Proposition~\ref{m3:r286} gives
		\[\begin{aligned}
			&\pi_{28}(\mathit{SO}(6))={}\frac{\mathbb Z}{2}\{[\iota_{3}]_{6}\phi'\}\oplus \frac{\mathbb Z}{2}\{[\nu_{4}\sigma'\kappa_{14}]_{6}\} \oplus \frac{\mathbb Z}{4}\{[\bar{\rho}''']\} \oplus \frac{\mathbb Z}{8}\{[ \nu_{5}]\bar{\kappa}_{8}\}.
		\end{aligned}\]
		Since $\pi_{28}(G_{2})$ is, up to isomorphism, a direct summand of $\pi_{28}(\mathit{SO}(7)),$ and
		$i_{\mathit{SU}(3)_{*}}( j_{*}\phi')=0\in\pi_{28}(G_{2})$
		by Proposition~\ref{m3:drgpp}(1)--(2), where $i_{\mathit{SU}(3)}: \mathit{SU}(3)\hookrightarrow G_{2},$
		we obtain $[\iota_{3}]_{7}\phi'=0.$  Since $[\bar{\rho}''']^{U}\in\pi_{28}(\mathit{SU}(3))$ maps to zero in $\pi_{28}(G_{2})$, it follows that $[\bar{\rho}''']_7=0$; see Proposition~\ref{m3:drgpp}(1),(2).
	\end{proof}
	\needspace{3\baselineskip}\begin{prop}
		The homotopy group $  \pi_{28}(\mathit{SO}(7))$ is given by
		\[\GroupLayout{\pi_{28}(\mathit{SO}(7)) = \mathrm{span}\{[\bar{\nu} _{6}+\varepsilon_{6}]\kappa_{14},\; [\bar{\sigma}_{6}] \nu_{25},\; [\eta_{6}]\eta_{7}\bar{\kappa}_{8},\; [ \nu_5]_{7}\bar{\kappa}_{8},\; [\nu_{4}\sigma'\kappa_{14}]_{7}\} {} \cong{} (\mathbb Z/2)^{5}.}{\begin{aligned}\pi_{28}(\mathit{SO}(7)) = \mathrm{span}\{&[\bar{\nu} _{6}+\varepsilon_{6}]\kappa_{14},\; [\bar{\sigma}_{6}] \nu_{25},\; [\eta_{6}]\eta_{7}\bar{\kappa}_{8},\\
				&[ \nu_5]_{7}\bar{\kappa}_{8},\; [\nu_{4}\sigma'\kappa_{14}]_{7}\} {} \cong{} (\mathbb Z/2)^{5}.\end{aligned}}\]
	\end{prop}
	\begin{proof}
		Proposition~\ref{m3:p28g2} gives $\pi_{28}(G_{2})=\frac{\mathbb Z}{2}\{\langle \bar{\nu}
		_{6}+\varepsilon_{6}\rangle\kappa_{14}\}\oplus \frac{\mathbb Z}{2}\{\langle \bar{\sigma}_{6}
		\nu_{25}\rangle\}\oplus \frac{\mathbb Z}{2}\{\langle \eta_{6}^{2}\rangle\bar{\kappa}_{8}\}.$ Notice that 
		$\pi_{28}^{7}=\frac{\mathbb Z}{2}\{\eta_{7}\bar{\kappa}_{8}\}\oplus \frac{\mathbb Z}{2}\{\sigma'\kappa_{14}\}.$
		Proposition~\ref{so7fj} gives $\pi_{28}(\mathit{SO}(7))=\mathcal{T}_{*}\pi_{28}(G_2)\oplus[\eta_6]\circ\pi_{28}^{7}.$ Combining this decomposition with the groups above gives $$\pi_{28}(\mathit{SO}(7))\cong(\mathbb Z/2)^5.$$ The required lifts exist by Lemma~\ref{m3:gjh}, and the group is generated by $[\bar{\nu} _{6}+\varepsilon_{6}]\kappa_{14},\allowbreak \; [\bar{\sigma}_{6} ]\nu_{25},\allowbreak \; \mathcal{T}_{*}(\langle
		\eta_{6}^{2}\rangle\bar{\kappa}_{8})$ and
		$[\eta_{6}]\eta_{7}\bar{\kappa}_{8},\allowbreak \;[\eta_{6}]\sigma'\kappa_{14}.$ Consider  the exact sequence
		\[
		\begin{tikzcd}[column sep=0.3cm]
			\pi_{29}^{6} \arrow[r]
			& \pi_{28}(\mathit{SO}(6)) \arrow[rr, "i_{*}"]
			&& \pi_{28}(\mathit{SO}(7)) \arrow[rr, "p_{*}"]
			&& \pi_{28}^{6} \arrow[r]
			& \pi_{27}(\mathit{SO}(6)).
		\end{tikzcd}
		\]
		It induces an exact sequence
		\[
		\begin{tikzcd}[column sep=0.3cm]
			0 \arrow[r]
			& \operatorname{Im}(i_{*}) \arrow[rr]
			&& \pi_{28}(\mathit{SO}(7)) \arrow[rr, "p_{*}"]
			&& \operatorname{Im}(p_{*}) \arrow[r]
			& 0.
		\end{tikzcd}
		\]
		We know
		\begin{gather*}
			\FitMath{\begin{aligned}
					\operatorname{Im}(p_{*})={}p_{*}\mathrm{span}\{& [\bar{\nu} _{6}+\varepsilon_{6}]\kappa_{14} ,\; [\bar{\sigma}_{6} \nu_{25}] ,\; \mathcal{T}_{*}(\langle \eta_{6}^{2}\rangle\bar{\kappa}_{8}),\; [\eta_{6}]\eta_{7}\bar{\kappa}_{8}, \\
					&[\eta_{6}]\sigma'\kappa_{14} \} \,=\frac{\mathbb Z}{2}\{( \bar{\nu} _{6}+\varepsilon_{6})\kappa_{14}\}\oplus \frac{\mathbb Z}{2}\{\bar{\sigma}_{6} \nu_{25}\} \oplus \frac{\mathbb Z}{2}\{\varepsilon_{6}\kappa_{14}\}.
			\end{aligned}}
		\end{gather*}
		Observe that $\varepsilon_{6}\kappa_{14}=\eta_{6}^{2}\bar{\kappa}_{8},\allowbreak \;\bar{\sigma}_{6}
		\nu_{25}=\nu_{6}\bar{\sigma}_{9},$
		as shown in the proof of Proposition~\ref{m3:p28g2}, and that $\quad \mathcal{T}_{*}(\langle
		\eta_{6}^{2}\rangle\bar{\kappa}_{8}) \equiv[\eta_{6}]\eta_{7}\bar{\kappa}_{8}\;\bmod\;
		\operatorname{Im}(i_{*}).$
		Hence
		$\operatorname{Im}(i_{*})\cong(\mathbb Z/2)^{2}.$
		Thus, by Lemma~\ref{m3:i37}, $$\operatorname{Im}(i_{*})=
		\frac{\mathbb Z}{2}\{[ \nu_5]_{7}\bar{\kappa}_{8}\}\oplus
		\frac{\mathbb Z}{2}\{[\nu_{4}\sigma'\kappa_{14}]_{7}\}.$$
		The image of $\mathrm{span}\{ [\bar{\nu} _{6}+\varepsilon_{6}]\kappa_{14} ,\allowbreak \; [\bar{\sigma}_{6}
		]\nu_{25},\allowbreak \;[\eta_{6}]\eta_{7}\bar{\kappa}_{8}\} \,\subseteq \pi_{28}(\mathit{SO}(7)) $ under
		$p_{*}$ is $$\mathrm{span} \{ ( \bar{\nu} _{6}+\varepsilon_{6})\kappa_{14} ,\allowbreak \; \bar{\sigma}_{6} \nu_{25},\allowbreak \;
		\varepsilon_{6}\kappa_{14} \}= \operatorname{Im}(p_{*}) \cong(\mathbb Z/2)^{3}.$$ This proves the
		assertion.
	\end{proof}
	The group $\pi_{28}(\mathit{SO}(7))$ was determined above, and the decomposition of $\pi_{28}^{7}$ was recalled above.
	Moreover, $[\iota_{7}]$ is well-defined \cite{Ka}.
	The following proposition now follows immediately.
	
	The isomorphism $\pi_{28}(\mathit{SO}(8))\cong \pi_{28}(\mathit{SO}(7))\oplus \pi_{28}^{7}$
	is induced by the fibration
	$\mathit{SO}(7)\rightarrow \mathit{SO}(8)\rightarrow S^{7}.$
	\needspace{3\baselineskip}\begin{prop}
		The structure of $\pi_{28}(\mathit{SO}(8))$ is as follows.
		\begin{gather*}
			\GroupLayout{\pi_{28}(\mathit{SO}(8)) \;= {} \mathrm{span}\{[\bar{\nu} _{6}+\varepsilon_{6}]_{8}\kappa_{14},\; [\bar{\sigma}_{6}] _{8} \nu_{25},\; [ \eta_{6}]_{8}\eta_{7}\bar{\kappa}_{8},\; [ \nu_5]_{8}\bar{\kappa}_{8},\; [\nu_{4}\sigma'\kappa_{14}]_{8},\; [\iota_{7}]\eta_{7}\bar{\kappa}_{8},\; [\iota_{7}]\sigma'\kappa_{14}\}\cong (\mathbb Z/2)^{7}.}{\begin{aligned}\pi_{28}(\mathit{SO}(8)) \;= {} \mathrm{span}\{&[\bar{\nu} _{6}+\varepsilon_{6}]_{8}\kappa_{14},\; [\bar{\sigma}_{6}] _{8} \nu_{25},\; [ \eta_{6}]_{8}\eta_{7}\bar{\kappa}_{8},\; [ \nu_5]_{8}\bar{\kappa}_{8},\\
					&[\nu_{4}\sigma'\kappa_{14}]_{8},\; [\iota_{7}]\eta_{7}\bar{\kappa}_{8},\; [\iota_{7}]\sigma'\kappa_{14}\}\cong (\mathbb Z/2)^{7}.\end{aligned}}
		\end{gather*}
		\qed
	\end{prop}

	\subsection{The Case $n=9$ and Its Spin Group Comparison}
	We next determine $\pi_{28}(\mathit{Spin}(9)).$  Recall  that the fibration
	$G_{2}\stackrel{j_{_{G_{2}}}\,}\longrightarrow\mathit{Spin}(7)\stackrel{q\,}\rightarrow S^{7}$ is trivial.
	
	Let $j_{S^{7}}\colon S^{7}\rightarrow \mathit{Spin}(7)\simeq G_{2}\times S^{7}$ be the (2-local)
	cross-section of $q,$ so that $q\circ j_{S^{7}}=\e_{7}$.
	
	For the fiber sequence $$   \Omega S^{15} \to  \mathit{Spin}(7)\to  \mathit{Spin}(9)
	\to  S^{15},   $$   in view of  \cite{Mi2}, we find
	\[\begin{aligned}
		&\Delta_{\mathit{Spin}}(\iota_{15}) \equiv xj_{G_{2}}\circ\langle \bar{\nu} _{6} +\varepsilon_{6}\rangle +yj_{S^{7}}\circ\sigma'\; \bmod\; j_{G_{2}}\circ i_{\mathit{SU}(3)}\circ[\nu_{5}^{2}]\nu_{11},
	\end{aligned}\]
	where $x,\allowbreak  y$ are odd integers and
	$i_{\mathit{SU}(3)}: \mathit{SU}(3)\hookrightarrow G_{2}$ is the standard inclusion
	\cite[Proposition 9.1]{Mi2}.
	Moreover, \cite[p. 174]{Mi2} gives $\Delta_{\mathit{Spin}}(\sigma_{15})=j_{S^{7}}\circ\sigma'\sigma_{14}$.
	\needspace{3\baselineskip}\begin{prop}
		\label{m3:p28spin9}The structure of $\pi_{28}(\mathit{Spin}(9))$ is as follows.
		\[ \begin{aligned}
			&\pi_{28}(\mathit{Spin}(9)) \;\cong\; (\mathbb Z/2)^{4}.
		\end{aligned} \]
	\end{prop}
	\begin{proof}
		Since $\pi_{28}^{15}=0,$   $\pi_{28}(\mathit{Spin}(9))\cong \operatorname{coker} (
		\Delta_{\mathit{Spin}}: \pi_{29}^{15}\rightarrow \pi_{28}(\mathit{Spin}(7))).$ We have $\pi_{29}^{15}=\frac{\mathbb Z}{4}\{ \sigma_{15}^{2}\}\oplus\frac{\mathbb Z}{2}\{\kappa_{15}\}$, and \cite[p. 55]{Oda1979} gives $\sigma'\sigma_{14}^{2}=0$.
		Then we deduce $$\Delta_{\mathit{Spin}}( \sigma_{15}^{2})=( \Delta_{\mathit{Spin}}(
		\sigma_{15}))\sigma_{21}=j_{S^{7}}\circ\sigma'\sigma_{14}\sigma_{21}=0,\vspace{-1em}$$
		\begin{align*}
			\Delta_{\mathit{Spin}}( \kappa_{15}) {}&={} ( \Delta_{\mathit{Spin}}( \iota_{15}))\kappa_{14} \\
			&\equiv xj_{G_{2}}\circ\langle \bar{\nu} _{6}+\varepsilon_{6}\rangle\kappa_{14} +yj_{S^{7}}\circ\sigma' \kappa_{14} \;\bmod\; j_{G_{2}}\circ i_{\mathit{SU}(3)}\circ[\nu_{5}^{2}]\nu_{11}\kappa_{14},
		\end{align*} where $x,\allowbreak  y$ are odd integers. Hence $\Delta_{\mathit{Spin}}( \kappa_{15}) \equiv
		yj_{S^{7}}\circ\sigma' \kappa_{14}\;\bmod\; j_{G_{2}}\circ \pi_{28}(G_{2}).$ Since the element  $\sigma'
		\kappa_{14}\in\pi_{28}^{7}$ has order 2 and the fibration  $ G_{2}\to\mathit{Spin}(7)\to S^{7}$ is trivial at the
		prime 2, with $$\pi_{28} (\mathit{Spin}(7))\cong\pi_{28} (\mathit{SO}(7))\cong(\mathbb Z/2)^{5},$$\noindent 
		it follows that $\Delta_{\mathit{Spin}}( \kappa_{15})$ has order 2.
		Consequently, $$\operatorname{coker} ( \Delta_{\mathit{Spin}}: \pi_{29}^{15}\rightarrow
		\pi_{28}(\mathit{Spin}(7)))\cong(\mathbb Z/2)^{4}.$$ Hence, $\pi_{28}(\mathit{Spin}(9))\cong
		\operatorname{coker} ( \Delta_{\mathit{Spin}}: \pi_{29}^{15}\rightarrow
		\pi_{28}(\mathit{Spin}(7)))\cong(\mathbb Z/2)^{4}.$
	\end{proof}
	
	Consider  the exact sequence
	\[
	\begin{tikzcd}[column sep=0.4cm]
		\pi_{29}^{8} \arrow[r, "\Delta"]
		& \pi_{28}(\mathit{SO}(8)) \arrow[r, "i_{*}"]
		& \pi_{28}(\mathit{SO}(9)) \arrow[r, "p_{*}"]
		& \pi_{28}^{8} \arrow[r, "\Delta"]
		& \pi_{27}(\mathit{SO}(8)).
	\end{tikzcd}
	\]
	We summarize the following facts:
	$\pi_{29}^{8}=\mathrm{span}\{
	\sigma_{8}^{3},\allowbreak \sigma_{8}\kappa_{15},\allowbreak \eta_{8}\bar{\kappa}_{9},\allowbreak (\Sigma\sigma')\kappa_{15}
	\}=\mathbb Z/4\oplus(\mathbb Z/2)^{3},$  
	$\pi_{28}^{8}=\frac{\z}{8}\{\bar{\kappa}_{8}\},$   and 
	$\pi_{27}(\mathit{SO}(8)) =\mathrm{span}\{ [\eta_{6}]_{8}\bar{\kappa}_{7},\allowbreak [ \nu_5]_{8}
	\bar{\sigma}_{8},\allowbreak \;[\iota_{7}]\bar{\kappa}_{7} \}=\mathbb Z/8\oplus\mathbb Z/2\oplus\mathbb Z/8.$ 
	
	The following lemma determines $\operatorname{Im}(i_{*})$ and $\operatorname{Im}(p_{*}).$
	\needspace{3\baselineskip}\begin{lem}
		\label{m3:zca8}
		\leavevmode\par
		\begin{itemize}
			\item[\textup{(1)}] $\Delta(\eta_{8}\bar{\kappa}_{9})= [
			\eta_{6}]_{8}\eta_{7}\bar{\kappa}_{8},$
			$\Delta(\sigma_{8}\kappa_{15}) \equiv[\iota_{7}]\sigma'\kappa_{14} \;\bmod\;
			i_{7,\allowbreak 8}\circ\pi_{28}(\mathit{SO}(7)).$
			\item[\textup{(2)}] $\Delta\bar{\kappa}_{8}=2[\iota_{7}]\bar{\kappa}_{7}-[\eta_{6}]_{8}\bar{\kappa}_{7}.$
		\end{itemize}
	\end{lem}
	\begin{proof}
		\begin{enumerate}[label=\textup{(\arabic*)},leftmargin=2em]
			\item The identity $\Delta(\iota_{8})=2[\iota_{7}]-[\eta_{6}]_{8}$ gives the assertion for
			$\Delta(\eta_{8}\bar{\kappa}_{9}),$ while
			$$\begin{aligned}
				\Delta(\sigma_8)&\equiv x_1[\iota_7]\sigma'
				+x_2[\mathrm{c}_{\mathrm{KM}}(\bar{\nu}_6+\varepsilon_6)]_8\m [\nu_5]_8\nu_8^2,[\eta_5\varepsilon_6]_8
			\end{aligned}$$ from \eqref{m1:Dsg82}, with $x_1,x_2$ odd,
			gives the assertion for $\Delta(\sigma_8\kappa_{15})$ because
			$\sigma'\kappa_{14}$ has order $2.$
			
			\item $\Delta(\iota_{8})=2[\iota_{7}]-[\eta_{6}]_{8}$ gives
			$\Delta\bar{\kappa}_{8}=2[\iota_{7}]\bar{\kappa}_{7}-[\eta_{6}]_{8}\bar{\kappa}_{7}.$
		\end{enumerate}
	\end{proof}
	\needspace{3\baselineskip}\begin{lem} \label{m3:zca5}
		\leavevmode\par
		\begin{itemize}
			\item[\textup{(1)}] $[\eta_{5}\varepsilon_{6}]_{8}\kappa_{14} \equiv0\;\bmod\;
			[\nu_{4}\sigma'\kappa_{14}]_{8};$
			$[ \nu_5]_{8}\nu_{8}^{2}\kappa_{14} \equiv0\;\bmod\; [\nu_{4}\sigma'\kappa_{14}]_{8}.$
			\item[\textup{(2)}] $[\eta_{5}\varepsilon_{6}]_{8}\sigma_{14}^{2}
			\equiv0\;\bmod\;[\nu_{4}\sigma'\kappa_{14}]_{8};\allowbreak $
			$[ \nu_5]_{8}\nu_{8}^{2}\sigma_{14}^{2}=0.$
			\item[\textup{(3)}] $\Delta(\sigma_{8}^{3}) \equiv0\;\bmod\;[\nu_{4}\sigma'\kappa_{14}]_{8};$
			$\Delta((\Sigma\sigma')\kappa_{15}) \equiv0\;\bmod\; [\nu_{4}\sigma'\kappa_{14}]_{8}.$
		\end{itemize}
	\end{lem}
	\begin{proof} \begin{itemize}
			\item[\textup{(1)}]
			For the first assertion, by \cite[(3.1), p. 10]{MMO} and
			\cite[(4.5), p. 19]{MMO}, we have
			$\eta_{5}\varepsilon_{6}\kappa_{14}=4\nu_{5}\bar{\kappa}_{8}.$
			Recall from Proposition~\ref{m3:r286} that \[\begin{aligned} &\pi_{28}(\mathit{SO}(6))={}\frac{\mathbb Z}{2}\{[\iota_{3}]_{6}\phi'\}\oplus \frac{\mathbb Z}{2}\{[\nu_{4}\sigma'\kappa_{14}]_{6}\} \oplus \frac{\mathbb Z}{4}\{[ \bar{\rho}''' ]\} \oplus \frac{\mathbb Z}{8}\{[ \nu_{5}]\bar{\kappa}_{8}\}. \end{aligned}\] Then  we infer $$[\eta_{5}\varepsilon_{6}]\kappa_{14}\ty4[ \nu_5]\bar{\kappa}_{8}\m
			[\iota_{3}]_{6}\phi',\,[\nu_{4}\sigma'\kappa_{14}]_{6},\,2[ \bar{\rho}''' ].$$
			
			\noindent Note that $2  \pi_{28}(\mathit{SO}(7))=0$.    Applying $i_{6,\allowbreak 7*}$ and using $[\iota_{3}]_{7}\phi'=[ \bar{\rho}''' ]_7=0$ from Lemma~\ref{m3:i37}, 
			we have     $$[\eta_{5}\varepsilon_{6}]_{7}\kappa_{14} \equiv4[ \nu_5]_7\bar{\kappa}_{8}=0\m[\nu_{4}\sigma'\kappa_{14}]_{7}.$$
			For the second assertion,  observe that
			$[ \nu_5]_{8}\nu_{8}^{2}\kappa_{14}=i_{6,\allowbreak 8}\circ[ \nu_5]\nu_{8}^{2}\kappa_{14},$ and
			$2[ \nu_5]\nu_{8}^{2}\kappa_{14}=0$ in $\pi_{28}(\mathit{SO}(6)).$
			Hence,
			$$[ \nu_5]\nu_{8}^{2}\kappa_{14}=x_{1} [\iota_{3}]_{6}\phi'+x_{2}[\nu_{4}\sigma'\kappa_{14}]_{6}
			+2x_{3}[ \bar{\rho}''' ]+4x_{4}[ \nu_5]\bar{\kappa}_{8}$$
			for some $x_{1},\allowbreak  x_{2},\allowbreak x_{3},\allowbreak x_{4}\in\mathbb Z.$ Applying $i_{6,\allowbreak 8*}$ and using  Lemma~\ref{m3:i37}, we obtain the relation 
			$[ \nu_5]_{8}\nu_{8}^{2}\kappa_{14}=x_{2}[\nu_{4}\sigma'\kappa_{14}]_{8}.$
			Therefore, $[ \nu_5]_{8}\nu_{8}^{2}\kappa_{14} \equiv0\;\bmod\; [\nu_{4}\sigma'\kappa_{14}]_{8}.$
			\item[\textup{(2)}] First, we claim that $i_{5,\allowbreak 8}
			\circ \pi_{28}(\mathit{SO}(5))=\mathrm{span}\{ [\nu_{4}\sigma'\kappa_{14}]_{8}\}.$
			From Proposition~\ref{m3:r285}, $$\pi_{28}(\mathit{SO}(5))=\mathrm{span}\{
			[\iota_{3}]_{5}\phi',\allowbreak \;[\iota_{3}]_{5}\mu_{3,\allowbreak 3},\allowbreak
			\;[\iota_{3}]_{5}\eta_{3}\bar{\mu}_{4}\sigma_{21},\allowbreak [\nu_{4}\sigma'\kappa_{14}]
			\}.$$ \noindent
			By Lemma~\ref{m3:i37}, $[\iota_{3}]_{7}\phi'=0.$ Moreover, \cite[Theorem 4.1 B)]{Oda1978} gives $j_{*}\mu_{3,\allowbreak 3}=2[\bar{\rho}''']\in\pi_{28}(\mathit{SU}(3))$. We also have $j_{*}\eta_{3}\bar{\mu}_{4}\sigma_{21} =0\in\pi_{28}(\mathit{SU}(3))$, since $\eta_{3}\bar{\mu}_{4}\sigma_{21}\in \operatorname{Im}(\Delta)$ by \cite[Proposition 3.2, the fifth row]{Oda1978}. Hence,
			$$\begin{aligned}
				&i_{5,\allowbreak 8} \circ \pi_{28}(\mathit{SO}(5))=\mathrm{span}\{ [\nu_{4}\sigma'\kappa_{14}]_{8}\}.
			\end{aligned}$$
			
			For the first assertion, $[\eta_{5}\varepsilon_{6}]\sigma_{14}\in \operatorname{Im}(i_{5,\allowbreak 6*}),$ and hence
			$$[\eta_{5}\varepsilon_{6}]_{8}\sigma_{14}^{2}\in i_{5,\allowbreak 8}
			\circ \pi_{28}(\mathit{SO}(5))=\mathrm{span}\{ [\nu_{4}\sigma'\kappa_{14}]_{8}\}.$$\noindent
			The second assertion follows from $\nu_{11}\sigma_{14}=0.$
			\item[\textup{(3)}] By \eqref{m1:Dsg82} and
			$\sigma'\sigma_{14}^2=0,$ we get 
			\begin{align*}
				\Delta(\sigma_8^3)
				&\equiv x_2[\mathrm{c}_{\mathrm{KM}}(\bar{\nu}_6+\varepsilon_6)]_8
				\sigma_{14}^2\m  [\nu_5]_8\nu_8^2\sigma_{14}^2,
				[\eta_5\varepsilon_6]_8\sigma_{14}^2.
			\end{align*}
			By definition,
			$[\bar{\nu}_6+\varepsilon_6]\sigma_{14}
			=\mathcal{T}_*(\langle\bar{\nu}_6+\varepsilon_6\rangle\sigma_{14})=0,$
			since $\pi_{21}(G_2)=0$ \cite[Theorem 6.1]{Mi2}.
			The comparison following \eqref{m1:i7R614} therefore yields
			\begin{align*}
				[\mathrm{c}_{\mathrm{KM}}(\bar{\nu}_6+\varepsilon_6)]_8\sigma_{14}^2
				&\equiv0\bmod [\eta_5\varepsilon_6]_8\sigma_{14}^2,
				[\nu_5]_8\nu_8^2\sigma_{14}^2.
			\end{align*}
			Part (2) now gives
			$\Delta(\sigma_8^3)\equiv0
			\bmod [\nu_4\sigma'\kappa_{14}]_8.$
			For the second assertion, in view of the formula
			$\Delta(\Sigma\sigma')=2[\iota_7]\sigma'-[\eta_6]_8\sigma'$
			and \eqref{m1:et6s'}, we  obtain
			\begin{align*}
				\Delta((\Sigma\sigma')\kappa_{15})
				&=2[\iota_7]\sigma'\kappa_{14}
				-4[\mathrm{c}_{\mathrm{KM}}(\bar{\nu}_6+\varepsilon_6)]_8\kappa_{14}-[\nu_5]_8\nu_8^2\kappa_{14}
				-[\eta_5\varepsilon_6]_8\kappa_{14}\\*
				&=[\nu_5]_8\nu_8^2\kappa_{14}
				+[\eta_5\varepsilon_6]_8\kappa_{14}.
			\end{align*}
			Here $2\pi_{28}(\mathit{SO}(8))=0$ by the preceding decomposition.
			Part (1) proves the second assertion.
		\end{itemize}
	\end{proof}
	
	\begin{prop}
		The structure of $\pi_{28}(\mathit{SO}(9))$ is as follows.
		\[\GroupLayout{\pi_{28}(\mathit{SO}(9)) = \mathrm{span}\{[\bar{\nu} _{6}+\varepsilon_{6}]_{9}\kappa_{14},\; [\bar{\sigma}_{6}] _{9} \nu_{25},\; [ \nu_5]_{9}\bar{\kappa}_{8},\; [\iota_{7}]_{9}\eta_{7}\bar{\kappa}_{8}\}\cong (\mathbb Z/2)^{4}.}{\begin{aligned}\pi_{28}(\mathit{SO}(9)) = \mathrm{span}\{&[\bar{\nu} _{6}+\varepsilon_{6}]_{9}\kappa_{14},\; [\bar{\sigma}_{6}] _{9} \nu_{25},\\
				&[ \nu_5]_{9}\bar{\kappa}_{8},\; [\iota_{7}]_{9}\eta_{7}\bar{\kappa}_{8}\}\cong (\mathbb Z/2)^{4}.\end{aligned}}\]
	\end{prop}
	\begin{proof}
		Consider  the exact sequence
		\[
		\begin{tikzcd}[column sep=0.4cm]
			\pi_{29}^{8}
			\arrow[r, "\Delta"]
			& \pi_{28}(\mathit{SO}(8))
			\arrow[r, "i_{*}"]
			& \pi_{28}(\mathit{SO}(9))
			\arrow[r, "p_{*}"]
			& \pi_{28}^{8}
			\arrow[r, "\Delta"]
			& \pi_{27}(\mathit{SO}(8)).
		\end{tikzcd}
		\]
		\noindent   In light of Lemma~\ref{m3:zca8}, we find that $\Delta\bar{\kappa}_{8}\in \pi_{27}(\mathit{SO}(8))$ has order $8.$  Thus,
		$$\operatorname{ker} ( \Delta : \pi_{28}^{8}\rightarrow \pi_{27}(\mathit{SO}(8)) )=0.$$

		\noindent   Recall that     $\pi_{28}(\mathit{SO}(8)) $ is equal to $$\mathrm{span}\{[\bar{\nu} _{6}+\varepsilon_{6}]_{8}\kappa_{14},\; [\bar{\sigma}_{6}] _{8} \nu_{25},\; [ \eta_{6}]_{8}\eta_{7}\bar{\kappa}_{8},\; [ \nu_5]_{8}\bar{\kappa}_{8},[\nu_{4}\sigma'\kappa_{14}]_{8},\; [\iota_{7}]\eta_{7}\bar{\kappa}_{8},\; [\iota_{7}]\sigma'\kappa_{14}\}\cong (\mathbb Z/2)^{7}.
		$$  \noindent   
		According to Lemma~\ref{m3:zca5}(3), we get  $$\Delta(\sigma_{8}^{3}) \equiv0\;\bmod\;[\nu_{4}\sigma'\kappa_{14}]_{8}.$$ Lemma~\ref{m3:zca8}(1) gives $$\Delta(\eta_{8}\bar{\kappa}_{9})=[\eta_{6}]_{8}\eta_{7}\bar{\kappa}_{8}.$$ Furthermore, Lemma~\ref{m3:zca5}(3) yields $$\Delta((\Sigma\sigma')\kappa_{15}) \equiv0\;\bmod\;[\nu_{4}\sigma'\kappa_{14}]_{8},$$and Lemma~\ref{m3:zca8}(1) implies $$\Delta(\sigma_{8}\kappa_{15}) \equiv [\iota_{7}]\sigma'\kappa_{14}\;\bmod\; i_{7,\allowbreak 8*}\pi_{28}(\mathit{SO}(7)).$$
		
		\noindent Write $\mathbf{a}^{C}=\mathbf{a}+ \operatorname{Im} ( \Delta : \pi_{29}^{8}\rightarrow \pi_{28}(\mathit{SO}(8)) )\in\operatorname{coker} ( \Delta) $, the class of $\mathbf{a}\in\pi_{28}(\mathit{SO}(8))$.       It follows that $$\mathrm{span} \{( [\bar{\nu} _{6}+\varepsilon_{6}]_{8}\kappa_{14})^{\mathrm{C}},\; ([\bar{\sigma}_{6}] _{8} \nu_{25} )^{\mathrm{C}},\; ([ \nu_5]_{8}\bar{\kappa}_{8})^{\mathrm{C}} , ([\iota_{7}]\eta_{7}\bar{\kappa}_{8})^{\mathrm{C}} \} \cong(\z/2)^{4}.$$ Additionally,
		\begin{align*}
			\operatorname{coker} ( \Delta : \pi_{29}^{8}\rightarrow \pi_{28}(\mathit{SO}(8)) ) {}&={} \mathrm{span}\{ ([\nu_{4}\sigma'\kappa_{14}]_{8} )^{\mathrm{C}}\} \\
			& {}\quad\jia\mathrm{span} \{( [\bar{\nu} _{6}+\varepsilon_{6}]_{8}\kappa_{14})^{\mathrm{C}},\; ([\bar{\sigma}_{6}] _{8} \nu_{25} )^{\mathrm{C}},\; ([ \nu_5]_{8}\bar{\kappa}_{8})^{\mathrm{C}} , ([\iota_{7}]\eta_{7}\bar{\kappa}_{8})^{\mathrm{C}} \}  \\
			&\cong{} (\mathbb Z/2)^{5}\,\text{or }\,(\mathbb Z/2)^{4}.
		\end{align*}
		\noindent
		Proposition~\ref{m3:p28spin9} tells us that $$\pi_{28}(\mathit{SO}(9))
		\cong\pi_{28}(\mathit{Spin}(9)) \cong(\mathbb Z/2)^{4}.$$\noindent Therefore, $([\nu_{4}\sigma'\kappa_{14}]_{8} )^{\mathrm{C}}=0.$
		Otherwise, $\operatorname{coker} ( \Delta : \pi_{29}^{8}\rightarrow \pi_{28}(\mathit{SO}(8)) )$ would be
		isomorphic to $(\mathbb Z/2)^{5}.$
		Thus, $\operatorname{coker} ( \Delta : \pi_{29}^{8}\rightarrow \pi_{28}(\mathit{SO}(8)) ) $ is equal to
		\begin{align*}
			\mathrm{span} \{( [\bar{\nu} _{6}+\varepsilon_{6}]_{8}\kappa_{14})^{\mathrm{C}},\; ([\bar{\sigma}_{6}] _{8} \nu_{25} )^{\mathrm{C}},\; ([ \nu_5]_{8}\bar{\kappa}_{8})^{\mathrm{C}} , ([\iota_{7}]\eta_{7}\bar{\kappa}_{8})^{\mathrm{C}} \} \cong{} (\mathbb Z/2)^{4}.
		\end{align*}
		Hence, $$\pi_{28}(\mathit{SO}(9))= \mathrm{span}\{[\bar{\nu} _{6}+\varepsilon_{6}]_{9}\kappa_{14},\; [\bar{\sigma}_{6}] _{9} \nu_{25},\; [ \nu_5]_{9}\bar{\kappa}_{8},\; [\iota_{7}]_{9}\eta_{7}\bar{\kappa}_{8}\}\cong (\mathbb Z/2)^{4}.$$
		Moreover, Lemma~\ref{m3:zca5}(1) and
		$[\nu_4\sigma'\kappa_{14}]_9=0$ show that
		$[\eta_5\varepsilon_6]_9\kappa_{14}
		=[\nu_5]_9\nu_8^2\kappa_{14}=0.$
		The comparison following \eqref{m1:i7R614} and the oddness of
		$\mathrm{c}_{\mathrm{KM}}$ therefore give
		\begin{align*}
			[\mathrm{c}_{\mathrm{KM}}(\bar{\nu}_6+\varepsilon_6)]_9\kappa_{14}
			&=[\bar{\nu}_6+\varepsilon_6]_9\kappa_{14}.
		\end{align*}
		The same equality holds after every further inclusion.
		
		This proves the assertion.
	\end{proof}
	\subsection{The Case $n=10$ and the Image from $\mathit{SO}(9)$}
	
	The calculation for $n=10$ begins by determining exactly which generators for $n=9$ survive inclusion.
	The following lemma also records the boundary relation needed to compare that image with the new
	elements of homotopy groups of spheres.
	
	\needspace{3\baselineskip}\begin{lem}
		\label{m3:k289} $i_{9,\allowbreak 10*} (\pi_{28}(\mathit{SO}(9)))=\mathrm{span}\{[ \bar{\nu}
		_{6}+\varepsilon_{6}]_{10}\kappa_{14},\; [\bar{\sigma}_{6}]_{10}\nu_{25},\; [ \nu_5]_{10}\bar{\kappa}_{8}\}\cong (\mathbb Z/2)^{3},$
		$\Delta(\bar{\sigma}_{9})=[ \nu_5]_{9}\circ\bar{\sigma}_{8}\neq0.$
	\end{lem}
	\begin{proof}
		The stem-$20$ table in Section~\ref{app:stem-20} gives $\pi_{29}^{9}=\frac{\mathbb Z}{8}\{\bar{\kappa}_{9}\}.$ 
		Using $\Delta\iota_{9}=[ \nu_5]_{9}+[\iota_{7}]_{9}\eta_{7}$ and the boundary composition formula, we obtain
		$\Delta\bar{\kappa}_{9}=[ \nu_5]_{9}\bar{\kappa}_{8}+[\iota_{7}]_{9}\eta_{7}\bar{\kappa}_{8}.$ 
		Exactness therefore gives $$i_{9,\allowbreak 10*}(\pi_{28}(\mathit{SO}(9)))\cong\pi_{28}(\mathit{SO}(9))/\mathrm{span}\{\Delta\bar{\kappa}_{9}\}.$$ 
		In the preceding decomposition of $\pi_{28}(\mathit{SO}(9))$ into four cyclic groups of order $2,$   this relation identifies the third and fourth generators.
		The images of the first three generators consequently remain independent and have order $2,$   so
		$i_{9,\allowbreak 10*} (\pi_{28}(\mathit{SO}(9)))=\frac{\mathbb Z}{2}\{[\bar{\nu} _{6}+\varepsilon_{6}]_{10}\kappa_{14}\}\oplus \frac{\mathbb Z}{2}\{[\bar{\sigma}_{6}]_{10}\nu_{25}\}\oplus \frac{\mathbb Z}{2}\{[ \nu_5]_{10}\bar{\kappa}_{8}\}.$
		
		For the boundary $\pi_{28}^{9}\to\pi_{27}(\mathit{SO}(9)),$   the relation $\eta_{7}\bar{\sigma}_{8}=0$ gives
		$$\Delta(\bar{\sigma}_{9})=([ \nu_5]_{9}+[\iota_{7}]_{9}\eta_{7})\bar{\sigma}_{8}=[ \nu_5]_{9}\bar{\sigma}_{8}.$$
		Since $J[ \nu_5]= \bar{\nu} _{6}+\varepsilon_{6}$ and $\varepsilon_{5}\bar{\sigma}_{13}=0$ by \cite[Proposition 2.2(7), p. 84]{Oda1979},
		we obtain
		$J([ \nu_5]_{9}\bar{\sigma}_{8})=( \bar{\nu} _{9}+\varepsilon_{9})\bar{\sigma}_{17}= \bar{\nu}
		_{9}\bar{\sigma}_{17}\neq0\in \pi_{36}^{9}.$
		This proves the assertion.
	\end{proof}
	\needspace{3\baselineskip}The preceding considerations lead to a further
	result concerning the choice and properties of
	certain homotopy classes.\begin{lem}\label{csdy1023}
		\leavevmode\par
		\begin{itemize}
			\item[\textup{(1)}] For every $n\geq2,$ $\bar{\zeta}_{4n+1}$ admits a lift to
			$\pi_{*}(\mathit{SO}(4n+2)),$ which may be chosen as $$[\bar{\zeta}_{4n+1}]\in  \{ [\zeta_{4n+1}],\allowbreak \;
			8\iota_{4n+12},\allowbreak \;2\sigma_{4n+12} \;\}_{1}.$$
			\item[\textup{(2)}] For the lift $[\bar{\zeta}_{4n+1}]$ above,
			$\mathrm{ord}([\bar{\zeta}_{4n+1}])=8,$ $\forall n\geq2.$
		\end{itemize}
	\end{lem}
	\begin{proof} Let $n\geq2.$
		\begin{itemize}
			\item[\textup{(1)}]
			By Lemma~\ref{m1:zeta1}, $ \{ [\zeta_{4n+1}],\allowbreak \; 8\iota_{4n+12},\allowbreak \;2\sigma_{4n+12} \;\}_{1}$ is
			well-defined.
			As a consequence of \cite[p. 137]{T}, we find
			$\bar{\zeta}_{5}\in \{ \zeta_{5},\allowbreak \; 8\iota_{16},\allowbreak \;2\sigma_{16} \;\}_{1}.$
			Notice that \begin{align*}
				&p_{*} \{ [\zeta_{4n+1}],\; 8\iota_{4n+12},\;2\sigma_{4n+12} \;\}_{1} \subseteq \{ \zeta_{4n+1},\; 8\iota_{4n+12},\;2\sigma_{4n+12} \;\}_{1}\ni \bar{\zeta}_{4n+1},\\
				&\mathrm{Ind} \{ \zeta_{4n+1},\; 8\iota_{4n+12},\;2\sigma_{4n+12} \;\}_{1} =\zeta_{4n+1}\circ\pi_{4n+20}^{4n+12}+\pi_{4n+13}^{4n+1}\circ2\sigma_{4n+13}=0.
			\end{align*}
			\noindent This proves the assertion.
			
			\item[\textup{(2)}] By \cite[(3.10)]{T}, we conclude that $\langle 8\iota,\allowbreak 2\sigma,\allowbreak 8\iota\rangle=0.$ Hence,
			\begin{align*}
				&\{ [\zeta_{4n+1}],\; 8\iota_{4n+12},\;2\sigma_{4n+12} \;\}_{1}\circ8\iota_{4n+20}={} [\zeta_{4n+1}]\circ\{ 8\iota_{4n+12},\;2\sigma_{4n+12} ,\;8\iota_{4n+19}\}_{1}= 0.
			\end{align*}
			Since $p_*[\bar\zeta_{4n+1}]=\bar\zeta_{4n+1}$ and $\bar\zeta_{4n+1}$ has order $8,$   the lift also has order $8.$ 
			
		\end{itemize}
	\end{proof}
	\needspace{3\baselineskip}\begin{prop}
		The structure of $\pi_{28}(\mathit{SO}(10))$ is as follows.
		\[\GroupLayout{\pi_{28}(\mathit{SO}(10)) = \mathrm{span}\{[\bar{\zeta}_{9}],\; [\bar{\nu} _{6}+\varepsilon_{6}]_{10}\kappa_{14},\; [\bar{\sigma}_{6}]_{10}\nu_{25},\; [ \nu_5]_{10}\bar{\kappa}_{8}\}\cong \mathbb Z/8\oplus (\mathbb Z/2)^{3}.}{\begin{aligned}\pi_{28}(\mathit{SO}(10)) = \mathrm{span}\{&[\bar{\zeta}_{9}],\; [\bar{\nu} _{6}+\varepsilon_{6}]_{10}\kappa_{14},\; [\bar{\sigma}_{6}]_{10}\nu_{25},\\
				&[ \nu_5]_{10}\bar{\kappa}_{8}\}\cong \mathbb Z/8\oplus (\mathbb Z/2)^{3}.\end{aligned}}\]
	\end{prop}
	\begin{proof}
		By Lemmas~\ref{m3:k289} and~\ref{csdy1023}, we see that
		$$i_{9,\allowbreak 10*} (\pi_{28}(\mathit{SO}(9)))=\frac{\mathbb Z}{2}\{[\bar{\nu} _{6}+\varepsilon_{6}]_{10}\kappa_{14}\}\oplus \frac{\mathbb Z}{2}\{[\bar{\sigma}_{6}]_{10}\nu_{25}\}\oplus \frac{\mathbb Z}{2}\{[ \nu_5]_{10}\bar{\kappa}_{8}\}$$ and that 
		$\operatorname{ker}( \Delta\colon \pi_{28}^{9}\rightarrow \pi_{27}(\mathit{SO}(9)))=\frac{\mathbb Z}{8}\{ \bar{\zeta}_{9}\};$ in addition, 
		$[\bar{\zeta}_{9}]$ has order 8. This proves the assertion.
	\end{proof}
	\subsection{The Case $n=11$ and the Remaining Connecting Images}
	
	For the case $n=11,$   we first need a controlled lift of $\eta_{10}\bar\mu_{11}.$  Its bracket description
	permits the order calculation\uprightsemicolon{} the remaining connecting images then determine the group.
	
	\needspace{3\baselineskip}\begin{lem}
		\label{m3:zkhca10cq11}
		$[\eta_{10}\bar{\mu}_{11}]\in \{ [\eta_{10}\mu_{11}],\allowbreak  2\iota_{20},\allowbreak  8\sigma_{20} \},$
		$\mathrm{ord}([\eta_{10}\bar{\mu}_{11}])=2.$
	\end{lem}
	\begin{proof}
		The lemma follows from the bracket description $\bar{\mu}_{3}\in \{ \mu_{3},\allowbreak 2\iota_{12},\allowbreak 8\sigma_{12}\}_{1}$ in \cite{T}, the order formula $\mathrm{ord}([\eta_{10}\mu_{11}])=2$ in \cite{KM2}, and $\langle 2\iota,\allowbreak 8\sigma,\allowbreak 2\iota\rangle=0.$
	\end{proof}
	\needspace{3\baselineskip}To prepare for the calculations that follow, we collect
	the relevant information in a form suited to
	the subsequent arguments.\begin{lem}
		\leavevmode\par
		\begin{itemize}
			\item[\textup{(1)}]     $   (i_{10,11})_*\pi_{28}(\mathit{SO}(10))=\mathrm{span}\{[\bar\zeta_9]_{11},\; [\mathrm{c}_{\mathrm{KM}}(\bar\nu_6+\varepsilon_6)]_{11}\kappa_{14},\; [ \nu_5]_{11}\bar\kappa_8\}\cong (\mathbb Z/2)^{3}.$
			\item[\textup{(2)}]     $   \pi_{28}(\mathit{SO}(11))=\mathrm{span}\{[\bar\zeta_9]_{11},\;[\mathrm{c}_{\mathrm{KM}}(\bar\nu_6+\varepsilon_6)]_{11}\kappa_{14},\;[ \nu_5]_{11}\bar\kappa_8,  [\eta_{10}\bar\mu_{11}], [4t\xi''+(\xi''\pm\lambda'')]\}$
			which is isomorphic to $(\mathbb Z/2)^{\,5}\;\text{or}\;(\mathbb Z/2)^{\,3}\oplus\mathbb Z/4.$
			Here, $t\in\mathbb Z$ and $\mathrm{ord}(\xi''\pm\lambda'')=2.$ 
		\end{itemize}
	\end{lem}
	\begin{proof}
		\begin{itemize}
			\item[\textup{(1)}]     In light of    $p_{*}( \Delta\bar{\sigma}_{10})=2\iota_{9}\circ\bar{\sigma}_{9}=0,$ we obtain
			\[
			\begin{aligned}
				&\Delta\bar{\sigma}_{10}{}\equiv{} 0\; \bmod\,i_{9,10_{*}}\pi_{28}(\mathit{SO}(9)) =\mathrm{span} \{ [\bar{\sigma}_{6}]_{10}\nu_{25}, \; [\bar{\nu} _{6}+\varepsilon_{6}]_{10}\kappa_{14} , \; [ \nu_5]_{10}\bar{\kappa}_{8} \; \}\cong (\mathbb Z/2)^{3}.
			\end{aligned}
			\]
			By \eqref{m1:Jbn6e6},
			$J([\mathrm{c}_{\mathrm{KM}}(\bar{\nu}_6+\varepsilon_6)])=\sigma'\sigma_{14},$
			and $J([\nu_5])=\bar{\nu}_6+\varepsilon_6.$
			The equality of the $\kappa_{14}$-composites proved for
			$\mathit{SO}(9)$ remains valid in $\mathit{SO}(10).$ Hence
			\begin{align*}
				J([\bar{\nu}_6+\varepsilon_6]_{10}\kappa_{14})
				&=J([\mathrm{c}_{\mathrm{KM}}(\bar{\nu}_6+\varepsilon_6)]_{10}\kappa_{14})\\
				&=\pm\Sigma^3(\sigma'\sigma_{14}\kappa_{21})=0.
			\end{align*}
			We also have
			$J( [ \nu_5]_{10}\bar{\kappa}_{8} )=( \bar{\nu} _{10}+\varepsilon_{10})\bar{\kappa}_{18}.$
			Moreover, \cite[p. 102]{Oda1979} gives $$J( \Delta\bar{\sigma}_{10})= P(\bar{\sigma}_{21})=\sigma_{10}\nu_{17}^{*}\nu_{35}\neq0.$$ Thus $J( \Delta\bar{\sigma}_{10})=\sigma_{10}\nu_{17}^{*}\nu_{35}$, and $\Delta\bar{\sigma}_{10}\in
			\frac{\mathbb Z}{2}\{[\bar{\sigma}_{6}]_{10}\nu_{25}\}\oplus \frac{\mathbb Z}{2}\{[ \bar{\nu}
			_{6}+\varepsilon_{6}]_{10}\kappa_{14}\}\oplus \frac{\mathbb Z}{2}\{[ \nu_5]_{10}\bar{\kappa}_{8}\}.$
			The elements $\sigma_{10}\nu_{17}^{*}\nu_{35}$ and $( \bar{\nu}
			_{10}+\varepsilon_{10})\bar{\kappa}_{18}$ are independent in $\pi_{28+10}^{10}.$
			We may therefore write $$\Delta\bar{\sigma}_{10}= x [\bar{\sigma}_{6}]_{10}\nu_{25}+y [ \bar{\nu}
			_{6}+\varepsilon_{6}]_{10}\kappa_{14} +z[ \nu_5]_{10}\bar{\kappa}_{8}$$\noindent for some
			$x,\allowbreak y,\allowbreak z\in\mathbb Z.$ If $x$ were even, the $J$-image of the two sides would differ. Thus, $x$ is
			odd.
			Thus,
			$\Delta\bar{\sigma}_{10}\equiv [\bar{\sigma}_{6}]_{10}\nu_{25} \;\bmod\,[ \bar{\nu}
			_{6}+\varepsilon_{6}]_{10}\kappa_{14} ,\allowbreak \;[ \nu_5]_{10}\bar{\kappa}_{8}.$ 
			Recall that
			$$[\bar{\zeta}_{9}]\in  \{ [\zeta_{9}],\allowbreak \; 8\iota_{20},\allowbreak \;2\sigma_{20} \}_{1}$$\noindent from Lemma
			\ref{csdy1023} and
			$\Delta(\zeta_{10})=\pm2[\zeta_{9}],$
			from Lemma~\ref{m1:zeta1}.
			Employing the definition of $[\bar{\zeta}_{9}]$ given in Lemma \ref{csdy1023} (1), we infer 
			$\Delta(\bar{\zeta}_{10})=\pm2[\bar{\zeta}_{9}].$  Therefore, $$\begin{aligned}
				&i_{10,\allowbreak 11_{*}}(\pi_{28}(\mathit{SO}(10)))=\frac{\mathbb Z}{2}\{[\bar{\zeta}_{9}]_{11}\}\oplus \frac{\mathbb Z}{2}\{[\mathrm{c}_{\mathrm{KM}}(\bar{\nu}_6+\varepsilon_6)]_{11}\kappa_{14}\}\oplus \frac{\mathbb Z}{2}\{[ \nu_5]_{11}\bar{\kappa}_{8}\}.
			\end{aligned}$$
			\item[\textup{(2)}] Note that $\Delta(4t\xi''+(\xi''\pm\lambda''))=0$ for some $t\in\mathbb Z$ by using the result of $\pi_{27}(\mathit{SO}(11)).$   Part (1) and
			Lemma~\ref{m3:zkhca10cq11} establish (2).
		\end{itemize}
	\end{proof}
	We use the decompositions
	$\pi_{29}^{11}=\frac{\mathbb Z}{8}\{\xi'\}\oplus\frac{\mathbb Z}{4}\{\xi'+\lambda'\}\oplus\frac{\mathbb Z}{2}\{\eta_{11}\bar{\mu}_{12}\}$ and
	\[\begin{aligned}
		&\pi_{30}^{12}={}\frac{\mathbb Z}{32}\{\xi_{12}\} \oplus\frac{\mathbb Z}{4}\{ \Sigma\xi' +4\xi_{12}\} \oplus\frac{\mathbb Z}{4}\{ \Sigma\xi' +\Sigma\lambda'\} \oplus\frac{\mathbb Z}{2}\{\eta_{12}\bar{\mu}_{13}\}
	\end{aligned}\]
	from \cite[Lemma 2.18, p. 78]{Og}.
	Employing the definition of $\xi'$ in \cite[Lemma 2.18, p. 78]{Og}, we show the following.
	\needspace{3\baselineskip}.\begin{prop}
		The following congruence holds:
		$p_{12*}(\Delta\xi_{12})\equiv-\xi'\pmod{4\xi'}.$ 
		Write $p_{12*}(\Delta\xi_{12})=u\xi',$   where $u$ is an odd integer satisfying $u\equiv-1\pmod4.$  We define $$[\xi']=u^{-1}\Delta\xi_{12},$$ where $u^{-1}\in\mathbb Z_{(2)}$ is the inverse of $u.$  Then $p_{12*}[\xi']=\xi'.$ 
	\end{prop}
	\begin{proof}
		By \cite[Lemma 2.18 (2), p. 78]{Og}, we find
		$-\xi'\in\{ 2\sigma_{11},\allowbreak \nu_{18},\allowbreak \sigma_{21}\}_{1}.$
		By \cite{T}, $$\begin{aligned}
			&\xi_{12}\in \{ \sigma_{12},\allowbreak \nu_{19},\allowbreak \sigma_{22}\}_{1}.
		\end{aligned}$$ Then $$p_{12*}\Delta\xi_{12}\in p_{12*}\Delta\{\sigma_{12},\allowbreak \nu_{19},\allowbreak \sigma_{22}\}_{1}\subseteq \{
		2\sigma_{11},\allowbreak \nu_{18},\allowbreak \sigma_{21}\}\ni-\xi'.$$
		The proof of \cite[Lemma 2.18(2), p.~78]{Og} gives $$p_{12*}\Delta\xi_{12}
		\equiv-\xi'\;\bmod\; 4\xi'.$$ Since $\xi'$ has order $8,$   this congruence gives $p_{12*}(\Delta\xi_{12})=u\xi'$ for an odd $u\equiv-1\pmod4.$  Multiplication by $u^{-1}$ therefore gives a lift projecting exactly to $\xi',$   as claimed.
	\end{proof}
	\needspace{3\baselineskip}\begin{lem} $p_{12*} \Delta\varepsilon_{12}^{*} \equiv0\;\bmod\;
		\eta_{11}\sigma_{12}\mu_{19},$ where $\varepsilon_{12}^{*}$ is given by \cite[Lemma 2.15(1), p.
		75]{Og}.
	\end{lem}
	\begin{proof}
		As in \cite[Lemma 2.15(1), p. 75]{Og},
		we redefine the element by taking $\varepsilon_{12}^{*}\in \{ \theta ,\allowbreak \;\nu_{24},\allowbreak \;\eta_{27} \}_{1},\allowbreak $
		and use $\Delta\theta=[\theta']\in \pi_{23}(\mathit{SO}(12))$ from Proposition~\ref{m1:R12}.
		So,\begin{align*}
			&p_{12*} \Delta\varepsilon_{12}^{*} \in p_{12*} \Delta \{ \theta ,\;\nu_{24},\;\eta_{27} \}_{1}\subseteq p_{12*} \{ [\theta'],\;\nu_{23},\;\eta_{26} \}\subseteq \{ \theta', \nu_{23},\;\eta_{26} \}, \\
			&\Sigma^{\infty} \{ \theta',\;\nu_{23},\;\eta_{26} \} \subseteq\langle 0,\;\nu,\;\eta \rangle=\frac{\mathbb Z}{2}\{ \eta\omega \}\oplus\frac{\mathbb Z}{2}\{ \eta\sigma\mu\}.
		\end{align*}
		Observe that	$\Sigma^{\infty}: \pi_{28}^{11}\rightarrow \pi_{17}^{S}$ is injective.
		Thus,
		$p_{12*} \Delta\varepsilon_{12}^{*} \equiv0\;\bmod\; \eta_{11}\sigma_{12}\mu_{19}.$
	\end{proof}
	\needspace{3\baselineskip}

	The preceding calculations lead to a further identity,
	which we state explicitly below.

	\begin{prop}
		The following identity holds.
		$\Delta(\pi_{29}^{11})=\mathrm{span}\{\Delta(\xi'+\lambda' )\}\cong \mathbb Z/2.$
	\end{prop}
	\begin{proof}
		We know that $\pi_{28}(\mathit{SO}(11))\cong (\mathbb Z/2)^{5} \;or\; (\mathbb Z/2)^{3}\oplus\mathbb Z/4$
		and that $\pi_{28}(\mathit{SO}(12))\cong(\mathbb Z/2)^{6},$   which is generated by $[\eta_{11}]\sigma_{12}\mu_{19} \;,\allowbreak
		\,[\bar{\mu}_{11}],\allowbreak \; \Delta \pi_{29}^{12}.$
		Notice that
		\[\begin{aligned}
			&\Delta \pi_{29}^{12}={}\frac{\mathbb Z}{2}\{\Delta(\varepsilon_{12}^{*})\}\oplus \frac{\mathbb Z}{2}\{[2\iota_{11}]\sigma_{11}\eta_{18}\mu_{19}\} \oplus \frac{\mathbb Z}{2}\{[2\iota_{11}]\nu_{11}\kappa_{14}\}\oplus \frac{\mathbb Z}{2}\{[2\iota_{11}]\bar{\mu}_{11}\}.
		\end{aligned}\]
		Thus,  we infer that \begin{gather*}
			\begin{aligned}
				p_{12*}( \pi_{28}(\mathit{SO}(12)) )&={}p_{12*}(\mathrm{span}\{ [\eta_{11}]\sigma_{12}\mu_{19} ,\,[\bar{\mu}_{11}]\;\} +\Delta \pi_{29}^{12}\;) \\
				&=\mathrm{span} \{ \eta_{11}\sigma_{12}\mu_{19} , \bar{\mu}_{11},\;p_{12*} \Delta(\varepsilon_{12}^{*})\;\}\cong {(\mathbb Z/2)^{2},}
			\end{aligned}
			\\[3pt]
			\begin{aligned}
				&|\;\operatorname{ker}(p_{12*})\;| =| \;\pi_{28}(\mathit{SO}(12))|\;/ \; |\;p_{12*}( \pi_{28}(\mathit{SO}(12)) )\;|=\frac{2^{6}}{ 2^{2}}=2^{4},
			\end{aligned}
			\\[3pt]
			\begin{aligned}
				&\frac{\pi_{28}(\mathit{SO}(11)) }{\operatorname{ker}( i_{12*})} \cong\operatorname{Im}( i_{12*})=\operatorname{ker}(p_{12*}).
			\end{aligned}
		\end{gather*}
		Then, we have\vspace{-1.1em}
		\begin{align*}
			| \Delta(\pi_{29}^{11})| &= | \operatorname{ker}( i_{12*}) |=| \;\pi_{28}(\mathit{SO}(11))\;|\;/ \; |\;\operatorname{Im}( i_{12*})\;| \\
			&=| \;\pi_{28}(\mathit{SO}(11))\;|\;/ \; |\operatorname{ker}(p_{12*})|=\frac{2^{5}}{ 2^{4}}=2,
		\end{align*}
		\noindent
		$\Delta(\eta_{11}\bar{\mu}_{12}) = \Delta\xi'=0.$  Therefore, $| \Delta(\pi_{29}^{11})|=|\,\Delta(\mathrm{ span}\{\,\xi'+\lambda' \; \})\,|=2,$   and
		$\Delta(\pi_{29}^{11})$ is equal to $\frac{\mathbb Z}{2}\{ \Delta(\xi'+\lambda' ) \}.$
	\end{proof}
	
	\section{The Cases $n\geq12$}
	The remaining calculations proceed downwards from the stable range. The starting calculation for $n=15$ is followed by the cases $n=14,13,12.$ 
	\subsection{The Stable Range and the Case $n=15$}
	
	The groups $\pi_{28}(\mathit{SO}(k ))$\;$(k\geq16)$ are covered, without explicit generators, by
	\cite[Tables, pp. 258--261]{Lun} together with Bott periodicity in the higher stable range.     As in the cases $23\leq i\leq 27,$   the results for $\pi_{i}(\mathit{SO}(k))$ with $i=28$ and $k\geq16$ in the
	theorem of this chapter can be obtained by  easy calculations
	using \cite[Tables, pp.~258--261]{Lun}, the propositions in
	Section~\ref{m3:sowq}, and \cite{Oda1979}.

	\needspace{3\baselineskip}\begin{prop}\label{r2815}$\pi_{28}(\mathit{SO}(15))=0.$
	\end{prop}
	\begin{proof}
		Consider  the diagram  with exact rows which commutes up to sign:
		\[
		\xymatrix@C=.7cm{
			0= \pi_{30}(\mathit{SO}) \ar[d] \ar[r] &
			\pi_{30}(\mathit{SO},\;\mathit{SO}(16))\ar[d]^{ J }_{(\text{by\;James'\, theorem})\,\approx}
			\ar[r]^{\,d_{_{30}}}_{\cong} & \pi_{29}(\mathit{SO}(16)) \ar[d]_{ J } \ar[r]^{} &
			\pi_{29}(\mathit{SO})=0 \ar[d] \\
			\pi_{30}^{S}(S^{0}) \ar[r]^{} & \pi_{30}(\Omega^{\infty}S^{\infty}, \Omega^{16}S^{16})
			\ar[r]^{\qquad \partial_{_{30}}\;}_{\quad\;\;onto} & \pi_{29}(\Omega^{16}S^{16}) \ar[r] &
			\pi_{29}^{S}(S^{0})=0.   }
		\]
		It follows that $J: \pi_{29}(\mathit{SO}(16))\rightarrow \pi_{29}(\Omega^{16}S^{16})$ is surjective.
		We next consider the diagram  with exact rows which commutes up to sign:
		\[
		\xymatrix@C=.7cm{
			\pi_{29}(\mathit{SO}(16)) \ar[d]_{ J }^{onto} \ar[r]^{p_{_{*}}} &
			\pi_{29}(S^{15}) \ar[d]_{\Sigma^{16}} ^{onto} \ar[r]^{\Delta\quad}_{onto\;} &
			\pi_{28}(\mathit{SO}(15))  \ar[d]_{} \ar[r] &\pi_{28}(\mathit{SO}(16))=0 \ar[d] \\
			\pi_{29}(\Omega^{16}S^{16})  \ar[r]^{H}_{onto} &\pi_{29}(\Omega^{16}S^{31})
			\ar[r]^{P}&\pi_{28+15}(S^{15}) \ar[r]^{\Sigma}_{\cong} &\pi_{28+16}(S^{16}).   }
		\]
		The homomorphism $\Sigma^{16}\circ p_{_{*}}=\pm(H\hc J)$ is surjective. We recall from Toda's book the following  homomorphism:
		$$  \s^{16}\colon  \frac{\z}{4}\n \cg_{15}^{2}\nn\oplus\frac{\z}{2}\n\cn_{15}\nn= \pi_{29}(S^{15}) \longrightarrow \pi_{29}(\Omega^{16}S^{31})= \frac{\z}{2}\n \cg_{31}^{2}\nn\oplus\frac{\z}{2}\n\cn_{31}\nn.$$\noindent
		Consequently,
		$p_{_{*}}\colon \pi_{29}(\mathit{SO}(16))\rightarrow\pi_{29}(S^{15})$ is surjective. In view of
		$\pi_{28}(\mathit{SO}(16))=0$ from \cite[Tables, pp. 258--261]{Lun}, we obtain 
		$$\pi_{28}(\mathit{SO}(15))=\operatorname{coker}(p_{_{*}}\colon
		\pi_{29}(\mathit{SO}(16))\rightarrow\pi_{29}(S^{15}))=0.$$
	\end{proof}
	\subsection{The Case $n=14$}
	
	The calculation now descends to the case $n=14.$  Since the group for the next value of $n$ has already been determined, the
	exact sequence can be used to identify the boundary generator and its order.
	
	\needspace{3\baselineskip}\begin{prop}\label{r2814}
		The structure of $\pi_{28}( \mathit{SO}(14))$ is as follows.
		\[ \FitMath{\pi_{28}( \mathit{SO}(14)) = \mathrm{span}\{\Delta\rho_{14}\}\cong \mathbb Z/32.} \]
	\end{prop}
	\begin{proof}
		We have shown that $\pi_{28}(\mathit{SO}(15))=0.$
		Thus, $\pi_{28}(\mathit{SO}(14))=\Delta \pi_{29}^{14}.$
		The following diagram commutes up to sign:
		\[
		\xymatrix@C=1.4cm{
			\pi_{29}^{14}\ar[d]_{\Sigma^{15}}^{\approx} \ar[r]^{\quad\Delta} & \pi_{28}(\mathit{SO}(14)) \ar[d]^{ J }\\
			\pi_{44}^{29}\ar[r]^{\quad P} & \pi_{42}^{14},   }
		\]
		\[\begin{aligned}
			&\pi_{42}^{14}={}\frac{\mathbb Z}{32}\{ P(\rho_{29}) \}\oplus\frac{\mathbb Z}{2}\{ \varepsilon_{14}\bar{\kappa}_{22}\}.
		\end{aligned}\]
		Then $\Delta\rho_{14}\in \pi_{28}(\mathit{SO}(14))$ has order 32. Moreover,
		$\Delta(\eta_{14}\kappa_{15})=\Delta(\eta_{14})\circ\kappa_{14}=(2\alpha)\circ\kappa_{14}=0$ for
		some $\alpha,$ since $\pi_{14}(\mathit{SO}(14))\cong\mathbb Z/4$ \cite{Ka}.
		Hence, $\pi_{28}(\mathit{SO}(14))=\frac{\mathbb Z}{32}\{ \Delta\rho_{14} \}.$
	\end{proof}
	\subsection{The Case $n=13$}
	
	With the result for the case $n=14$ available, we return to the case $n=13$ and compare its generators with their images
	under inclusion. The relations in the next statement record the information needed for this step.
	
	\needspace{3\baselineskip}\begin{prop}
		\label{m3:r2813}  The homotopy groups are given as follows.
		\[\GroupLayout{\pi_{28}( \mathit{SO}(13)) = \mathrm{span}\{[\eta_{11}]_{13}\sigma_{12}\mu_{19},\; [\bar{\mu}_{11}]_{13}\}\cong (\mathbb Z/2)^{2}.}{\begin{aligned}\pi_{28}( \mathit{SO}(13)) = \mathrm{span}\{&[\eta_{11}]_{13}\sigma_{12}\mu_{19},\\
				&[\bar{\mu}_{11}]_{13}\}\cong (\mathbb Z/2)^{2}.\end{aligned}}\]
	\end{prop}
	\begin{proof}
		Consider  the exact sequence
		\[
		\xymatrix@C=0.2cm{
			\pi_{29}^{13} \ar[rr]^{\Delta\qquad} && \pi_{28}(\mathit{SO}(13)) \ar[rr]^{i_{*}} &&\pi_{28}(\mathit{SO}(14))  \ar[rr]^{\qquad p_{*}} &&
			\pi_{28}^{13} }
		\]
		$\pi_{28}(\mathit{SO}(14))=\frac{\mathbb Z}{32}\{ \Delta\rho_{14} \}$ , $p_{*}( \Delta\rho_{14}) =2\rho_{13}.$
		Thus, we obtain  the exact sequence
		\[
		\xymatrix@C=0.27cm{
			0 \ar[r] &\Delta(\pi_{29}^{13}) \ar[rr]^{\subseteq\;\;} && \pi_{28}(\mathit{SO}(13)) \ar[rr]^{i_{*}\quad} &&
			\frac{\mathbb Z}{2}\{ 16\Delta(\rho_{14}) \} \ar[r] & 0 .}
		\]
		$\Delta(\pi_{29}^{13})$ is generated by
		$\Delta(\sigma_{13}\mu_{20})=[\eta_{11}]_{13}\sigma_{12}\mu_{19}.$
		In light of   \cite[Lemma 2.1(iv)]{HKM}, we obtain
		$\pm J([\eta_{11}]_{13})=\pm\Sigma J([\eta_{11}])=\Sigma A_{1} .$ Then $J([\eta_{11}]_{13}\sigma_{12}\mu_{19})=(\Sigma A_{1})\sigma_{25}\mu_{32}\in\pi_{28+13}^{13}$ is a
		generator of order 2.
		Therefore, $[\eta_{11}]_{13}\sigma_{12}\mu_{19}\neq0.\allowbreak $
		Thus, we obtain the exact sequence
		\[
		\xymatrix@C=0.25cm{
			0 \ar[r] &\frac{\mathbb Z}{2}\{ [\eta_{11}]_{13}\sigma_{12}\mu_{19} \}\ar[rr]^{\quad\subseteq} &&
			\pi_{28}(\mathit{SO}(13)) \ar[rr]^{i_{*}} && \frac{\mathbb Z}{2}\{ 16\Delta(\rho_{14}) \} \ar[r] & 0 .}
		\]
		Consider $\{ [\mu_{11}],\allowbreak 2\iota_{20},\allowbreak 8\sigma_{20} \}_{1}\subseteq \pi_{28}(\mathit{SO}(12)).$ Since $\mathrm{ord}
		[\mu_{11}]=2,$ the bracket is well-defined. By the definition of $\bar{\mu}_{3},$
		$p_{*}\{ [\mu_{11}],\allowbreak 2\iota_{20},\allowbreak 8\sigma_{20} \}_{1}\subseteq \{ \mu_{11},\allowbreak 2\iota_{20},\allowbreak 8\sigma_{20}
		\}_{1}\ni \bar{\mu}_{11}.$
		Let $B=\{[\mu_{11}],2\iota_{20},8\sigma_{20}\}_1.$  Since $2\pi_{21}(\mathit{SO}(12))=0,$   the second summand of $\mathrm{Ind}\,B$ vanishes. The corresponding summand of the sphere bracket also vanishes, since $8\pi_{21}^{11}=0.$  Consequently,
		\[
		\begin{aligned}
			p_*(\mathrm{Ind}\,B)&=p_*([\mu_{11}]\circ\pi_{28}^{20})
			=\mu_{11}\circ\pi_{28}^{20}=\mathrm{Ind}\{\mu_{11},2\iota_{20},8\sigma_{20}\}_1.
		\end{aligned}
		\]
		The inclusion above and equality of the indeterminacy subgroups imply
		\[
		p_*B=\{\mu_{11},2\iota_{20},8\sigma_{20}\}_1\ni\bar{\mu}_{11}.
		\]
		We may therefore choose $[\bar{\mu}_{11}]\in \{ [\mu_{11}],\allowbreak 2\iota_{20},\allowbreak 8\sigma_{20} \}_{1}$ such
		that $p_{*}([\bar{\mu}_{11}])=\bar{\mu}_{11}.$
		Since $2\pi_{21}(\mathit{SO}(12))=0,$ we find
		\[
		\begin{aligned}
			&2[\bar{\mu}_{11}] = 2\mathrm{id}\circ[\bar{\mu}_{11}] \in\{ 2[\mu_{11}], 2\iota_{20}, 8\sigma_{20} \}_{1} =\{ 0, 2\iota_{20}, 8\sigma_{20} \}_{1} =\pi_{21}(\mathit{SO}(12))\circ8\sigma_{21} =0,
		\end{aligned}
		\]Thus, $\mathrm{ord}[\bar{\mu}_{11}]=2.$
		Then,
		\[\begin{aligned}
			&[\bar{\mu}_{11}]_{14} \in i_{12,14} \circ \{ \,[\mu_{11}], 2\iota_{20}, 8\sigma_{20} \}_{1} \subseteq \{ \,[\mu_{11}]_{14}, \; 2\iota_{20}, 8\sigma_{20} \}_{1}.
		\end{aligned}\]
		Toda's formulas in \cite{T} give $\rho'\in \{ \sigma_{9},\allowbreak 16\iota_{16},\allowbreak \sigma_{16} \}_{1}$ and $16\rho_{14}=8\Sigma^{5}\rho'$. We also use $[\mu_{11}]_{14}=8[2\iota_{13}]\sigma_{13}$ from \cite[p. 23]{KM2}, together with $\pi_{21}(\mathit{SO}(14))\cong\mathbb Z/4.$ Thus, \begin{align*}
			\Delta(16\rho_{14})&={} \Delta(8\Sigma^{5}\rho') {}\in{} \Delta \{\,\sigma_{14}, 16\iota_{21},8\sigma_{21} \}_{2} {} \\
			&\subseteq{} \{\,\Delta\sigma_{14}, 16\iota_{20},8\sigma_{20} \}_{1} {} \\
			&={} \{\,[2\iota_{13}]\sigma_{13}, \;16\iota_{20},8\sigma_{20} \}_{1} {} \\
			&\supseteq{} \{\,8[2\iota_{13}]\sigma_{13}, \;2\iota_{20},8\sigma_{20} \}_{1} {} \\
			&={} \{ \,[\mu_{11}]_{14},\;2\iota_{20},8\sigma_{20} \}_{1} \ni [\bar{\mu}_{11}]_{14},
		\end{align*}
		and
		\[\begin{aligned}
			&\mathrm{Ind}\{\,[2\iota_{13}]\sigma_{13}, \; 16\iota_{20}, 8\sigma_{20} \}_{1} =\pi_{21}(\mathit{SO}(14))\circ8\sigma_{21} +[2\iota_{13}]\sigma_{13}\circ\pi_{28}^{20} =0 .
		\end{aligned}\]
		Therefore, $[\bar{\mu}_{11}]_{14}= 16\Delta(\rho_{14})$ has order 2.
		By this, we conclude that $\mathrm{ord} [\bar{\mu}_{11}]_{13}=2.$ Finally,
		$\pi_{28}(\mathit{SO}(13))=\frac{\mathbb Z}{2}\{[\eta_{11}]_{13}\sigma_{12}\mu_{19}\}\oplus
		\frac{\mathbb Z}{2}\{[\bar{\mu}_{11}]_{13}\}.$
	\end{proof}
	\subsection{The Case $n=12$}
	
	The case $n=12$ is the last step of this downward calculation. Here we combine the information from the case $n=13$ with the relevant homotopy groups of spheres to determine the stated abstract decomposition.
	
	\needspace{3\baselineskip}\begin{prop}
		The structure of $\pi_{28}( \mathit{SO}(12))$ is as follows.
		\[ \begin{aligned}
			&\pi_{28}( \mathit{SO}(12)) \;\cong\; (\mathbb Z/2)^{6}.
		\end{aligned} \] This group is generated by
		$[\eta_{11}]\sigma_{12}\mu_{19} ,\allowbreak  \;[\bar{\mu}_{11}],\allowbreak \; \Delta \pi_{29}^{12}.$
	\end{prop}
	\begin{proof}
		Consider  the exact sequence\[\xymatrix{
			\pi_{29}^{12} \ar[r]^{\Delta\qquad} & \pi_{28}(\mathit{SO}(12)) \ar[r]^{i_{12,13_{*}}} & \pi_{28}(\mathit{SO}(13)) \ar[r]^{} & 0
			.}
		\]
		Using $\Sigma C_{1}\circ\sigma_{24}\mu_{31}=\Sigma^{2}F_{1}^{\;(1)}$ from \cite[Proposition 3.2(1) (7), p. 89]{Oda1979}, we see that 
		\[\begin{aligned}
			&P(\bar{\mu}_{25})\equiv \Sigma C_{1}^{\;(2)}\;\bmod\;P(\eta_{25}\sigma_{26}\mu_{33}).
		\end{aligned}\]
		Then,
		by \cite[Proposition 4.5(8), p. 93]{Oda1979},
		we have, in $\pi_{28+12}^{12} , \;
		J\Delta\pi_{29}^{12}=P\Sigma^{13}\pi_{29}^{12}\cong(\mathbb Z/2)^{4}.$
		Since $\pi_{29}^{12}\cong(\mathbb Z/2)^{4},$ we obtain
		$\Delta\pi_{29}^{12}\cong(\mathbb Z/2)^{4}.$
		Finally, we obtain the exact sequence
		\[
		\begin{tikzcd}[
			column sep=1em,
			row sep=1em
			]
			0 \arrow[r]
			&
			\Delta\pi_{29}^{12}
			\arrow[rr, "\subseteq"]
			\arrow[d, equal]
			&&
			\pi_{28}(\mathit{SO}(12))
			\arrow[rr, "i_{12,13*}"]
			&&
			\pi_{28}(\mathit{SO}(13))
			\arrow[r]
			\arrow[d, equal]
			&
			0.
			\\
			&
			(\mathbb Z/2)^4
			&&&&
			(\mathbb Z/2)^2
		\end{tikzcd}
		\]
		\noindent
		By Proposition~\ref{m3:r2813}, the image of $\frac{\mathbb Z}{2}\{[\eta_{11}]\sigma_{12}\mu_{19}\}\oplus
		\frac{\mathbb Z}{2}\{[\bar{\mu}_{11}]\}$ under $i_{12,\allowbreak 13_{*}}$ is
		$\pi_{28}(\mathit{SO}(13)).$
		The proof of Proposition~\ref{m3:r2813} also shows that $[\bar{\mu}_{11}]$ is well-defined and has
		order 2.
		Hence, $\pi_{28}(\mathit{SO}(12))\cong(\mathbb Z/2)^{6}$ is generated by
		$[\eta_{11}]\sigma_{12}\mu_{19},\allowbreak \;[\bar{\mu}_{11}],\allowbreak \; \Delta \pi_{29}^{12}.$
	\end{proof}
	\chapter{Some Homotopy Groups $\pi_{29}(\mathit{SO}(n))$}\label{ch:degree29}
	
	Having considered the homotopy groups of $\mathit{SO}(n)$ in the preceding degrees, we now direct our attention toward degree $29$. The objective of the present chapter is to determine $\pi_{29}(\mathit{SO}(n))$ for certain values of $n$, thereby extending the range of computations treated in this volume. Although the discussion is restricted to selected cases, these calculations provide a natural continuation of the preceding investigations. The results will be presented within the same general mathematical framework, with the aim of maintaining a coherent and systematic account of the homotopy groups considered throughout the book.

	We often freely use both the descriptions of the groups $\pi_{*}(S^{j})$
	provided in Appendix~\ref{tlqdzj} and the composition relations
	among generators in Toda's range found in \^{O}guchi's table
	in Appendix~\ref{app:oguchi-relations}.

	\section{Main Theorem and Generators}
	
	The following theorem summarizes the results of this chapter.
	\begin{mainleft}
		\begin{thm}\label{m3:zdlr29}
			\textup{All spaces in this theorem are localized at $2.$  The displayed groups are obtained from the integral homotopy groups by tensoring with $\mathbb Z_{(2)}.$  Free cyclic summands become $\mathbb Z_{(2)}$; for finite groups, this localization retains precisely the $2$-primary summands. The homotopy groups are given as follows.   }\par
			\textup{The Cases $3\leq n\leq6.$ }\par
			\smallskip\noindent\(
			\GroupLayout{\pi_{29}(\mathit{SO}(3)) = \mathrm{span}\{[\eta_{2}]\bar{\mu}'\sigma_{22},\; [\eta_{2}]\nu'\phi_{6},\; [\eta_{2}]\bar{\alpha}\nu_{26},\; [\eta_{2}]\eta_{3}\mu_{3,4}\}\cong \mathbb Z/4\oplus (\mathbb Z/2)^{3}.}{\begin{aligned}\pi_{29}(\mathit{SO}(3)) = \mathrm{span}\{&[\eta_{2}]\bar{\mu}'\sigma_{22},\; [\eta_{2}]\nu'\phi_{6},\; [\eta_{2}]\bar{\alpha}\nu_{26},\\
					&[\eta_{2}]\eta_{3}\mu_{3,4}\}\cong \mathbb Z/4\oplus (\mathbb Z/2)^{3}.\end{aligned}}
			\)\par
			\smallskip\noindent\(
			\GroupLayout{\pi_{29}(\mathit{SO}(4)) \;= {} \mathrm{span}\{[\eta_{2}]_{4}\bar{\mu}'\sigma_{22},\; [\eta_{2}]_{4}\nu'\phi_{6},\; [\eta_{2}]_{4}\bar{\alpha}\nu_{26},\; [\eta_{2}]_{4}\eta_{3}\mu_{3,4},\; [\iota_{3}]\bar{\mu}'\sigma_{22},\; [\iota_{3}]\nu'\phi_{6},\; [\iota_{3}]\bar{\alpha}\nu_{26},\; [\iota_{3}]\eta_{3}\mu_{3,4}\}\cong \mathbb Z/4\oplus (\mathbb Z/2)^{3}\oplus \mathbb Z/4\oplus (\mathbb Z/2)^{3}.}{\begin{aligned}\pi_{29}(\mathit{SO}(4)) \;= {} \mathrm{span}\{&[\eta_{2}]_{4}\bar{\mu}'\sigma_{22},\; [\eta_{2}]_{4}\nu'\phi_{6},\; [\eta_{2}]_{4}\bar{\alpha}\nu_{26},\; [\eta_{2}]_{4}\eta_{3}\mu_{3,4},\; [\iota_{3}]\bar{\mu}'\sigma_{22},\\
					&[\iota_{3}]\nu'\phi_{6},\; [\iota_{3}]\bar{\alpha}\nu_{26},\; [\iota_{3}]\eta_{3}\mu_{3,4}\}\cong \mathbb Z/4\oplus (\mathbb Z/2)^{3}\oplus \mathbb Z/4\oplus (\mathbb Z/2)^{3}.\end{aligned}}
			\)\par
			\smallskip\noindent\(
			\GroupLayout{\pi_{29}(\mathit{SO}(5)) = \mathrm{span}\{[\nu_{4}\rho'']\sigma_{22},\; [\iota_{3}]_{5}\eta_{3}\mu_{3,4},\; [\nu_{4}^{2}]\bar{\sigma}_{10},\; [\nu_{4}\sigma'\eta_{14}]\kappa_{15}\}\cong \mathbb Z/16\oplus (\mathbb Z/2)^{3}.}{\begin{aligned}\pi_{29}(\mathit{SO}(5)) = \mathrm{span}\{&[\nu_{4}\rho'']\sigma_{22},\; [\iota_{3}]_{5}\eta_{3}\mu_{3,4},\; [\nu_{4}^{2}]\bar{\sigma}_{10},\\
					&[\nu_{4}\sigma'\eta_{14}]\kappa_{15}\}\cong \mathbb Z/16\oplus (\mathbb Z/2)^{3}.\end{aligned}}
			\)\par
			
			\smallskip\noindent\(
			\GroupLayout{\pi_{29}(\mathit{SO}(6)) \;= {} \mathrm{span}\{[\nu_{4}\rho'']_{6}\sigma_{22},\; [\nu_{4}^{2}]_{6}\bar{\sigma}_{10},\; [\nu_{4}\sigma'\eta_{14}]_{6}\kappa_{15},\; [ \nu_5]\sigma_{8}\kappa_{15},\; [ \nu_5]\eta_{8}\bar{\kappa}_{9}\}\cong \mathbb Z/16\oplus (\mathbb Z/2)^{4}.}{\begin{aligned}\pi_{29}(\mathit{SO}(6)) \;= {} \mathrm{span}\{&[\nu_{4}\rho'']_{6}\sigma_{22},\; [\nu_{4}^{2}]_{6}\bar{\sigma}_{10},\; [\nu_{4}\sigma'\eta_{14}]_{6}\kappa_{15},\\
					&[ \nu_5]\sigma_{8}\kappa_{15},\; [ \nu_5]\eta_{8}\bar{\kappa}_{9}\}\cong \mathbb Z/16\oplus (\mathbb Z/2)^{4}.\end{aligned}}
			\)\par
			
			The groups $\pi_{29}(\mathit{SO}(n))$ for $7\leq n\leq14$ remain under investigation.\par
			\textup{The Cases $n\geq15.$ }\par
			\smallskip\noindent\(
			\GroupLayout{\pi_{29}(\mathit{SO}(15)) = \mathrm{span}\{[\eta_{14}\sigma_{15}]\sigma_{22},\; [\eta_{14}\kappa_{15}]\}\cong (\mathbb Z/2)^{2}.}{\begin{aligned}\pi_{29}(\mathit{SO}(15)) = \mathrm{span}\{&[\eta_{14}\sigma_{15}]\sigma_{22},\\
					&[\eta_{14}\kappa_{15}]\}\cong (\mathbb Z/2)^{2}.\end{aligned}}
			\)\par
			\smallskip\noindent\(
			\GroupLayout{\pi_{29}(\mathit{SO}(16)) = ( \Delta\pi_{30}^{16}) \oplus \mathrm{span}\{[\sigma_{15}]\sigma_{22},\; [\kappa_{15}]\} \cong(\mathbb Z/2)^{2} \oplus(\mathbb Z/4 \oplus\mathbb Z/2).}{\begin{aligned}\pi_{29}(\mathit{SO}(16)) = ( \Delta\pi_{30}^{16}) \oplus \mathrm{span}\{&[\sigma_{15}]\sigma_{22},\\
					&[\kappa_{15}]\} \cong(\mathbb Z/2)^{2} \oplus(\mathbb Z/4 \oplus\mathbb Z/2).\end{aligned}}
			\)\par
			\smallskip\noindent\(
			\GroupLayout{\pi_{29}(\mathit{SO}(17)) = \mathrm{span}\{[\kappa_{15}]_{17},\; [\sigma_{15}]_{17}\sigma_{22}\}\cong \mathbb Z/2\oplus \mathbb Z/4.}{\begin{aligned}\pi_{29}(\mathit{SO}(17)) = \mathrm{span}\{&[\kappa_{15}]_{17},\\
					&[\sigma_{15}]_{17}\sigma_{22}\}\cong \mathbb Z/2\oplus \mathbb Z/4.\end{aligned}}
			\)\par
			\smallskip\noindent\(
			\GroupLayout{\pi_{29}(\mathit{SO}(18)) = \mathrm{span}\{[\kappa_{15}]_{18},\; [\sigma_{15}]_{18}\sigma_{22}\}\cong \mathbb Z/2\oplus \mathbb Z/4.}{\begin{aligned}\pi_{29}(\mathit{SO}(18)) = \mathrm{span}\{&[\kappa_{15}]_{18},\\
					&[\sigma_{15}]_{18}\sigma_{22}\}\cong \mathbb Z/2\oplus \mathbb Z/4.\end{aligned}}
			\)\par
			\smallskip\noindent\(
			\GroupLayout{\pi_{29}(\mathit{SO}(19)) = \mathrm{span}\{[\eta_{18}^{2}]\mu_{20},\; [\kappa_{15}]_{19},\; [\sigma_{15}]_{19}\sigma_{22}\}\cong (\mathbb Z/2)^{2}\oplus \mathbb Z/4.}{\begin{aligned}\pi_{29}(\mathit{SO}(19)) = \mathrm{span}\{&[\eta_{18}^{2}]\mu_{20},\; [\kappa_{15}]_{19},\\
					&[\sigma_{15}]_{19}\sigma_{22}\}\cong (\mathbb Z/2)^{2}\oplus \mathbb Z/4.\end{aligned}}
			\)\par
			\smallskip\noindent\(
			\GroupLayout{\pi_{29}(\mathit{SO}(20)) = \mathrm{span}\{[\eta_{19}]\mu_{20},\; [\eta_{18}^{2}]_{20}\mu_{20},\; [\kappa_{15}]_{20},\; [\sigma_{15}]_{20}\sigma_{22}\}\cong (\mathbb Z/2)^{3}\oplus \mathbb Z/4.}{\begin{aligned}\pi_{29}(\mathit{SO}(20)) = \mathrm{span}\{&[\eta_{19}]\mu_{20},\; [\eta_{18}^{2}]_{20}\mu_{20},\; [\kappa_{15}]_{20},\\
					&[\sigma_{15}]_{20}\sigma_{22}\}\cong (\mathbb Z/2)^{3}\oplus \mathbb Z/4.\end{aligned}}
			\)\par
			\smallskip\noindent\(
			\GroupLayout{\pi_{29}(\mathit{SO}(21)) = \mathrm{span}\{[\nu_{20}^{2}]\nu_{26},\; [\sigma_{15}]_{21}\sigma_{22},\; [\eta_{19}]_{21}\mu_{20}\}\cong (\mathbb Z/4)^{2}\oplus \mathbb Z/2.}{\begin{aligned}\pi_{29}(\mathit{SO}(21)) = \mathrm{span}\{&[\nu_{20}^{2}]\nu_{26},\; [\sigma_{15}]_{21}\sigma_{22},\\
					&[\eta_{19}]_{21}\mu_{20}\}\cong (\mathbb Z/4)^{2}\oplus \mathbb Z/2.\end{aligned}}
			\)\par
			\smallskip\noindent\(
			\GroupLayout{\pi_{29}(\mathit{SO}(22)) = \mathrm{span}\{[\nu_{20}^{2}]_{22}\nu_{26},\; [\sigma_{15}]_{22}\sigma_{22}\}\cong (\mathbb Z/4)^{2}.}{\begin{aligned}\pi_{29}(\mathit{SO}(22)) = \mathrm{span}\{&[\nu_{20}^{2}]_{22}\nu_{26},\\
					&[\sigma_{15}]_{22}\sigma_{22}\}\cong (\mathbb Z/4)^{2}.\end{aligned}}
			\)\par
			\smallskip\noindent\(
			\GroupLayout{\pi_{29}(\mathit{SO}(23)) = \mathrm{span}\{[\nu_{20}^{2}]_{23}\nu_{26},\; [\sigma_{15}]_{23}\sigma_{22}\}\cong (\mathbb Z/2)^{2}.}{\begin{aligned}\pi_{29}(\mathit{SO}(23)) = \mathrm{span}\{&[\nu_{20}^{2}]_{23}\nu_{26},\\
					&[\sigma_{15}]_{23}\sigma_{22}\}\cong (\mathbb Z/2)^{2}.\end{aligned}}
			\)\par
			\smallskip\noindent\(
			\GroupLayout{\pi_{29}(\mathit{SO}(24)) = \mathrm{span}\{[\nu_{20}^{2}]_{24}\nu_{26},\; [\nu_{23}]\nu_{26}\}\cong (\mathbb Z/2)^{2}.}{\begin{aligned}\pi_{29}(\mathit{SO}(24)) = \mathrm{span}\{&[\nu_{20}^{2}]_{24}\nu_{26},\\
					&[\nu_{23}]\nu_{26}\}\cong (\mathbb Z/2)^{2}.\end{aligned}}
			\)\par
			\smallskip\noindent\(
			\FitMath{\pi_{29}(\mathit{SO}(25)) = \mathrm{span}\{[\nu_{23}]_{25}\nu_{26}\}\cong \mathbb Z/2.}
			\)\par
			\smallskip\noindent\(
			\FitMath{\pi_{29}(\mathit{SO}(26)) = \mathrm{span}\{[\nu_{23}]_{26}\nu_{26}\}\cong \mathbb Z/2.}
			\)\par
			\smallskip\noindent\(
			\GroupLayout{\pi_{29}(\mathit{SO}(27)) = \mathrm{span}\{[\nu_{23}]_{27}\nu_{26},\; [\eta_{26}^{2}]\eta_{28}\}\cong (\mathbb Z/2)^{2}.}{\begin{aligned}\pi_{29}(\mathit{SO}(27)) = \mathrm{span}\{&[\nu_{23}]_{27}\nu_{26},\\
					&[\eta_{26}^{2}]\eta_{28}\}\cong (\mathbb Z/2)^{2}.\end{aligned}}
			\)\par
			\smallskip\noindent\(
			\GroupLayout{\pi_{29}(\mathit{SO}(28)) = \mathrm{span}\{[\eta_{27}]\eta_{28},\; [\eta_{26}^{2}]_{28}\eta_{28}\}\cong (\mathbb Z/2)^{2}.}{\begin{aligned}\pi_{29}(\mathit{SO}(28)) = \mathrm{span}\{&[\eta_{27}]\eta_{28},\\
					&[\eta_{26}^{2}]_{28}\eta_{28}\}\cong (\mathbb Z/2)^{2}.\end{aligned}}
			\)\par
			\smallskip\noindent\(
			\FitMath{\pi_{29}(\mathit{SO}(29)) = \mathrm{span}\{[\eta_{27}]_{29}\eta_{28}\}\cong \mathbb Z/2.}
			\)\par
			\smallskip\noindent\(
			\FitMath{\pi_{29}(\mathit{SO}(30)) = \mathrm{span}\{[2\iota_{29}]\}\cong \mathbb Z_{(2)}.}
			\)\par
			\smallskip\noindent\(
			\begin{aligned}
				&\pi_{29}(\mathit{SO}(31+n)) = 0\,(n\geq0).
			\end{aligned}
			\)\par
			
		\end{thm}
	\end{mainleft}
	
	\section{The Cases $3\leq n\leq6$}
	
	The degree-$29$ homotopy group of $S^3$ supplies the initial generators for $n=3.$  We retain their composites
	explicitly, since the later unitary and $G_2$ calculations use these same classes.
	
	The group is
	\[\begin{aligned}
		&\pi_{29}^{3}={}\frac{\mathbb Z}{4}\{\bar{\mu}'\sigma_{22}\}\oplus \frac{\mathbb Z}{2}\{\nu'\phi_{6}\} \oplus \frac{\mathbb Z}{2}\{\bar{\alpha}\nu_{26}\} \oplus \frac{\mathbb Z}{2}\{\eta_{3}\mu_{3, 4}\}.
	\end{aligned}\]
	The following proposition follows from $\mathit{SO}(3)=\mathbb{R}P^{3}$ and $$\pi_{j}(
	\mathit{SO}(4))=[\iota_{3}] \circ\pi_{j}^{3}\oplus [\eta_{2}]_{4} \circ\pi_{j}^{3} \cong
	\pi_{j}^{3}\oplus\pi_{j}^{3}$$\noindent for $j\geq2.$ 
	See \cite[p. 14]{Ka} for details.
	\needspace{3\baselineskip}\begin{prop}
		The structure of $\pi_{29}(\mathit{SO}(m))$ for $3\leq m\leq4$ is as follows.
		\begin{gather*}
			\GroupLayout{\pi_{29}(\mathit{SO}(3)) = \mathrm{span}\{[\eta_{2}]\bar{\mu}'\sigma_{22},\; [\eta_{2}]\nu'\phi_{6},\; [\eta_{2}]\bar{\alpha}\nu_{26},\; [\eta_{2}]\eta_{3}\mu_{3,4}\}\cong \mathbb Z/4\oplus (\mathbb Z/2)^{3}}{\begin{aligned}\pi_{29}(\mathit{SO}(3)) = \mathrm{span}\{&[\eta_{2}]\bar{\mu}'\sigma_{22},\; [\eta_{2}]\nu'\phi_{6},\; [\eta_{2}]\bar{\alpha}\nu_{26},\\
					&[\eta_{2}]\eta_{3}\mu_{3,4}\}\cong \mathbb Z/4\oplus (\mathbb Z/2)^{3}\end{aligned}},\\[3pt]
			\GroupLayout{\pi_{29}(\mathit{SO}(4)) \;= {} \mathrm{span}\{[\eta_{2}]_{4}\bar{\mu}'\sigma_{22},\; [\eta_{2}]_{4}\nu'\phi_{6},\; [\eta_{2}]_{4}\bar{\alpha}\nu_{26},\; [\eta_{2}]_{4}\eta_{3}\mu_{3,4},\; [\iota_{3}]\bar{\mu}'\sigma_{22},\; [\iota_{3}]\nu'\phi_{6},\; [\iota_{3}]\bar{\alpha}\nu_{26},\; [\iota_{3}]\eta_{3}\mu_{3,4}\}\cong \mathbb Z/4\oplus (\mathbb Z/2)^{3}\oplus \mathbb Z/4\oplus (\mathbb Z/2)^{3}.}{\begin{aligned}\pi_{29}(\mathit{SO}(4)) \;= {} \mathrm{span}\{&[\eta_{2}]_{4}\bar{\mu}'\sigma_{22},\; [\eta_{2}]_{4}\nu'\phi_{6},\; [\eta_{2}]_{4}\bar{\alpha}\nu_{26},\; [\eta_{2}]_{4}\eta_{3}\mu_{3,4},\; [\iota_{3}]\bar{\mu}'\sigma_{22},\\
					&[\iota_{3}]\nu'\phi_{6},\; [\iota_{3}]\bar{\alpha}\nu_{26},\; [\iota_{3}]\eta_{3}\mu_{3,4}\}\cong \mathbb Z/4\oplus (\mathbb Z/2)^{3}\oplus \mathbb Z/4\oplus (\mathbb Z/2)^{3}.\end{aligned}}
		\end{gather*}
		\qed
	\end{prop}
	
	Recall the principal bundle
	$S^{3}\stackrel
	{i}\hookrightarrow \mathit{Sp}(2)\twoheadrightarrow S^{7}.$
	By \cite[Theorem 5.1 A)]{Oda1978},
	\[\begin{aligned}
		&\pi_{29}(\mathit{Sp}(2))={}\frac{\mathbb Z}{16}\{[\rho'']\sigma_{22}\}\oplus \frac{\mathbb Z}{2}\{[\sigma'\eta_{14}]\kappa_{15}\} \oplus \frac{\mathbb Z}{2}\{[\nu_{7}]\bar{\sigma}_{10}\}\oplus \frac{\mathbb Z}{2}\{i_{*}\eta_{3}\mu_{3,4}\},
	\end{aligned}\]
	where $[x]$ denotes a lift of $x \in \pi_{*}(S^{7}).$ By Proposition \ref{sp2so5gx}, we obtain the following proposition.
	\needspace{3\baselineskip}\begin{prop}
		The structure of $\pi_{29}(\mathit{SO}(5))$ is as follows.
		\[\GroupLayout{\pi_{29}(\mathit{SO}(5)) = \mathrm{span}\{[\nu_{4}\rho'']\sigma_{22},\; [\iota_{3}]_{5}\eta_{3}\mu_{3,4},\; [\nu_{4}^{2}]\bar{\sigma}_{10},\; [\nu_{4}\sigma'\eta_{14}]\kappa_{15}\}\cong \mathbb Z/16\oplus (\mathbb Z/2)^{3}.}{\begin{aligned}\pi_{29}(\mathit{SO}(5)) = \mathrm{span}\{&[\nu_{4}\rho'']\sigma_{22},\; [\iota_{3}]_{5}\eta_{3}\mu_{3,4},\; [\nu_{4}^{2}]\bar{\sigma}_{10},\\
				&[\nu_{4}\sigma'\eta_{14}]\kappa_{15}\}\cong \mathbb Z/16\oplus (\mathbb Z/2)^{3}.\end{aligned}}\]
		\qed
	\end{prop}
	
	Proposition \ref{sp2so5gx} and  \cite[Theorem 5.1 A)]{Oda1978} also imply the following lemma.
	
	\needspace{3\baselineskip}\begin{lem}
		In $\pi_{29}(\mathit{SO}(5)),$
		$[\iota_{3}]_{5}\mu'\rho_{14}=8[\nu_{4}\sigma'\rho_{14}]= 8[\nu_{4}\rho'']\sigma_{22}.$
		\qed
	\end{lem}
	
	Recall from Section~\ref{lqzh} that the universal covering of $\mathit{SO}(6)$ is $\mathit{Spin}(6)=\mathit{SU}(4).$
	Oda gives $\pi_{29}(\mathit{SU}(4))$ in \cite{Oda1978}. Then we have the following.
	\needspace{3\baselineskip}\begin{prop}
		The structure of $\pi_{29}(\mathit{SO}(6))$ is as follows.
		\[\GroupLayout{\pi_{29}(\mathit{SO}(6)) = \mathrm{span}\{[\nu_{4}\rho'']_{6}\sigma_{22},\; \allowbreak \; [\nu_{4}^{2}]_{6}\bar{\sigma}_{10},\; \allowbreak \;[\nu_{4}\sigma'\eta_{14}]_{6}\kappa_{15},\; \allowbreak [ \nu_5]\sigma_{8}\kappa_{15},\; \allowbreak \;[ \nu_5]\eta_{8}\bar{\kappa}_{9}\} \cong\mathbb Z/16 \oplus(\mathbb Z/2)^{4}.}{\begin{aligned}\pi_{29}(\mathit{SO}(6)) = \mathrm{span}\{&[\nu_{4}\rho'']_{6}\sigma_{22},\; \allowbreak \; [\nu_{4}^{2}]_{6}\bar{\sigma}_{10},\; \allowbreak \;[\nu_{4}\sigma'\eta_{14}]_{6}\kappa_{15},\\
				&\allowbreak [ \nu_5]\sigma_{8}\kappa_{15},\; \allowbreak \;[ \nu_5]\eta_{8}\bar{\kappa}_{9}\} \cong\mathbb Z/16 \oplus(\mathbb Z/2)^{4}.\end{aligned}}\]
		\qed
	\end{prop}

	\section{The Cases $n\geq15$}
	
	The groups $\pi_{29}(\mathit{SO}(k ))$\;$(k\geq17)$ are covered, without explicit generators, by
	\cite[Tables, pp. 258--261]{Lun} together with Bott periodicity in the higher stable range. Each table entry must be used with its stated restrictions.     As in the cases $23\leq i\leq 27,$   the results for $\pi_{i}(\mathit{SO}(k))$ with $i=29$ and $k\geq17$ in the
	theorem of this chapter can be obtained by  easy calculations
	using \cite[Tables, pp.~258--261]{Lun}, the propositions in
	Section~\ref{m3:sowq}, and \cite{Oda1979}. The groups $\pi_{29}(\mathit{SO}(k))$ for $k=15,16$ can be readily
	determined from the result for $\pi_{29}(\mathit{SO}(17)),$   using
	the methods of Propositions~\ref{r2814} and~\ref{r2815}.
	We omit the proofs.

	\section{Auxiliary Results for $\pi_{29}(G_2)$}
	
	This section records partial information relevant to $\pi_{29}(\mathit{G}_2)$.
	It does not determine $\pi_{29}(\mathit{G}_2)$ or $\pi_{29}(\mathit{SO}(k))$ for $7\leq k\leq14$.
	These cases remain outside the calculations completed in this volume.

	Consider  the exact sequence
	\begin{equation*}
		\begin{tikzcd}[]
			\pi_{30}^{6}
			\arrow[r, "\Delta"]
			& \pi_{29}(\mathit{SU}(3))
			\arrow[r, "i_{*}"]
			& \pi_{29}(G_{2})
			\arrow[r, "p_{*}"]
			& \pi_{29}^{6}
			\arrow[r, "\Delta"]
			& \pi_{28}(\mathit{SU}(3)).
		\end{tikzcd}
	\end{equation*}
	The following group structures will be used; the unitary groups are given in \cite{Mi2} and \cite{Oda1978}:
	\[\FitMath{\begin{aligned}
			\pi_{30}^{6}= {}& \frac{\mathbb Z}{8}\{\bar{\zeta}^{\prime}_{6}\} \oplus\frac{\mathbb Z}{2}\{\delta_{6}\} \oplus\frac{\mathbb Z}{2}\{\bar{\mu}_{6}\sigma_{23}\} \oplus\frac{\mathbb Z}{2}\{\nu_{6}\sigma_{9}\kappa_{16}\} \\
			&\mathord{\oplus}\,\frac{\mathbb Z}{2}\{\bar{\sigma}^{\prime}_{6}\} \oplus\frac{\mathbb Z}{2}\{P(\lambda\eta_{31})\} \oplus\frac{\mathbb Z}{2}\{P(\xi_{13}\eta_{31})\},
	\end{aligned}}\]
	\begin{gather*}
		\begin{aligned}
			&\pi_{29}(\mathit{SU}(3))={}\frac{\mathbb Z}{2}\{j_{*}\bar{\mu}'\sigma_{22}\} \oplus\frac{\mathbb Z}{4}\{[\nu_{5}\sigma_{8}\kappa_{15}]\},
		\end{aligned}
		\\[3pt]
		\begin{aligned}
			&\pi_{29}^{6}={}\mathrm{span}\{\nu_{6}\bar{\kappa}_{9},\bar{\rho}'',\phi_{6}\} \oplus\mathrm{span}\{P(\lambda),P(\xi_{13})\} \cong(\mathbb Z/8\oplus\mathbb Z/4\oplus\mathbb Z/2) \oplus(\mathbb Z/8\oplus\mathbb Z/4),
		\end{aligned}
		\\[3pt]
		\begin{aligned}
			&\pi_{28}(\mathit{SU}(3))={}\frac{\mathbb Z}{4}\{[\bar{\rho}''']\} \oplus\frac{\mathbb Z}{4}\{[2\iota_{5}]\nu_{5}\bar{\kappa}_{8}\} \oplus\frac{\mathbb Z}{2}\{j_{*}\phi'\}.
		\end{aligned}
	\end{gather*}
	
	Here our symbol $j_{*}$ corresponds to $i_{*}$ in \cite{Oda1978}. The proof of
	Proposition~\ref{m3:p28g2} shows that
	$\operatorname{coker}(\Delta\colon  \pi_{29}^{6}\rightarrow \pi_{28}(\mathit{SU}(3)))=0,$
	which implies that
	$\Delta\colon  \pi_{29}^{6}\rightarrow \pi_{28}(\mathit{SU}(3))$ is surjective.  Consequently,
	$$|\operatorname{ker}(\Delta\colon  \pi_{29}^{6}\rightarrow
	\pi_{28}(\mathit{SU}(3)))|=\frac{8\times4\times2\times8\times4}{4\times4\times2}=64.$$ \noindent Recall that $i_{_{S^{3},\allowbreak G_{2}}}$ denotes the inclusion $S^{3}\hookrightarrow G_{2}.$
	We obtain the following lemma.
	\needspace{3\baselineskip}\begin{lem}
		\leavevmode\par
		\begin{itemize}
			\item[\textup{(1)}] $\Delta(\bar{\mu}_{6}\sigma_{23})=j_{*}\bar{\mu}'\sigma_{22},$
			$2(i_{*}[\nu_{5}\sigma_{8}\kappa_{15}])=0\in\pi_{29}(G_{2}).$
			\item[\textup{(2)}] $\mathrm{ord}( \Delta(\bar{\zeta}^{\prime}_{\,6}))\leq2,$
			\[
			\GroupLayout{\Delta \mathrm{span}\{\nu_{6}\sigma_{9}\kappa_{16},\; \delta_{6},\; \bar{\sigma}^{\prime}_{\,6},\; P(\lambda\eta_{31}),\; P(\xi_{13}\eta_{31})\} \subseteq \mathrm{span}\{j_{*}\bar{\mu}'\sigma_{22},\; 2[\nu_{5}\sigma_{8}\kappa_{15}]\}\cong (\mathbb Z/2)^{2}.}{\begin{aligned}\Delta \mathrm{span}\{&\nu_{6}\sigma_{9}\kappa_{16},\; \delta_{6},\; \bar{\sigma}^{\prime}_{\,6},\; P(\lambda\eta_{31}),\\
					&P(\xi_{13}\eta_{31})\} \subseteq \mathrm{span}\{j_{*}\bar{\mu}'\sigma_{22},\; 2[\nu_{5}\sigma_{8}\kappa_{15}]\}\cong (\mathbb Z/2)^{2}.\end{aligned}}
			\]
			\item[\textup{(3)}] $\operatorname{Im}(i_{*}: \pi_{29}(\mathit{SU}(3))\rightarrow\pi_{29}(G_{2})
			)=\mathrm{span}\{i_{*}[\nu_{5}\sigma_{8}\kappa_{15}]\}\cong \mathbb Z/2.$
		\end{itemize}
	\end{lem}
	\begin{proof}
		\begin{itemize}
			\item[\textup{(1)}] Proposition 6.3 of \cite{Mi2} gives
			$\Delta\bar{\mu}_{6}=j_{*}\bar{\mu}',$ proving the first assertion. For the second
			assertion, the proof for $i$=29 in \cite[p. 84]{Oda1978} gives
			$$2[\nu_{5}\sigma_{8}\kappa_{15}]=j_{*}\bar{\alpha}\nu_{26}\in\pi_{29}(\mathit{SU}(3)).$$ \noindent Recall from Section~\ref{m2:sec:small-n} that  $i_{_{S^{3},\allowbreak G_{2}}}\circ\bar{\alpha}=0\in\pi_{26}(G_{2}).$ 
			Then, in $\pi_{29}(G_{2}),$   $$2(i_{*}[\nu_{5}\sigma_{8}\kappa_{15}])=(i_{_{S^{3},\allowbreak G_{2}}}\circ\bar{\alpha})\circ\nu_{26}=0.$$
			
			\item[\textup{(2)}] By the proof of \cite[(3.8)]{MMO}, we find
			$\bar{\zeta}'_{\,6}\in\{ \zeta',\allowbreak  8\iota_{22},\allowbreak 2\sigma_{22}\}_{1}.$
			In view of  \cite[Proposition 6.3]{Mi2}, we have $\Delta\zeta'=j_{*}\mu'\sigma_{14}.$
			Then,
			\begin{align*}
				&\Delta\bar{\zeta}'_{\,6} \in\{\Delta\zeta',\;8\iota_{21},\;2\sigma_{21}\} =\{j_{*}\mu'\sigma_{14},\;8\iota_{21},\;2\sigma_{21}\} \supseteq j_{*}\{\mu'\sigma_{14},\;8\iota_{21},\;2\sigma_{21}\}, \\
				&\mathrm{Ind}\{ j_{*}\mu'\sigma_{14}, 8\iota_{21},2\sigma_{21}\} =j_{*}\mu'\sigma_{14}\circ\pi_{29}^{21}+\pi_{22}(\mathit{SU}(3))\circ2\sigma_{22}=0.
			\end{align*}
			\noindent Here we use the relation $2\pi_{22}(\mathit{SU}(3))=0$ from \cite{MT2}. Thus,
			\[\begin{aligned}
				&\Delta\bar{\zeta}'_{\,6} =j_{*}\{\mu'\sigma_{14}, 8\iota_{21}, 2\sigma_{21}\} \in\operatorname{Im}(j_{*}) =\frac{\mathbb Z}{2}\{ j_{*}\bar{\mu}'\sigma_{22}\} \oplus\frac{\mathbb Z}{2}\{ 2[\nu_{5}\sigma_{8}\kappa_{15}]\},
			\end{aligned}\]
			and hence
			$\mathrm{ord}( \Delta\bar{\zeta}'_{\,6})\leq2.$ The remaining part of (2) is immediate.
			\item[\textup{(3)}] If $[\nu_{5}\sigma_{8}\kappa_{15}]\in \operatorname{Im}(\Delta),$   then the structure of $\pi_{30}^{6}$ would imply that $\Delta\bar{\zeta}'_{\,6}$ has order 4, a contradiction. Parts (1) and (2) therefore give $\operatorname{coker}( \Delta\colon 
			\pi_{30}^{6}\rightarrow\pi_{29}(\mathit{SU}(3)))=\frac{\mathbb Z}{2}\{
			[\nu_{5}\sigma_{8}\kappa_{15}]+\operatorname{Im}( \Delta)\},$ proving the assertion.
		\end{itemize}
	\end{proof}
	\chapter{Geometric Applications}\label{ch:geometry}
	
	This chapter explains how an element of $\pi_k(\mathit{SO}(n))$
	defines a vector bundle over $S^{k+1},$   and how maps between homotopy
	groups answer questions about that bundle. The three maps we shall
	use are the inclusion $i_{n,n+1}\colon      \mathit{SO}(n)\to\mathit{SO}(n+1),$  
	the last-column projection $p_n\colon     \mathit{SO}(n)\to S^{n-1},$   and
	$J\colon     \pi_k(\mathit{SO}(n))\to\pi_{n+k}(S^n).$ 
	We shall explain their geometric meanings by constructing the relevant
	bundles, sections, and attaching maps.
	
	With the  exception of Proposition~\ref{jhyy1027}, we do not claim any originality for the results collected in this chapter. These results are standard in their respective areas and are well documented in the mathematical literature, with comprehensive treatments available in established textbooks. However, as this book is intended primarily for specialists in homotopy theory rather than geometry, we believe it is worthwhile to recall the basic geometric notions and results that are relevant to our investigation. Their inclusion is intended to provide readers with the necessary background without requiring frequent consultation of external references.
	
	The reader is assumed to know homotopy groups, fibrations, CW complexes,
	and mapping cones. We give the vector-bundle definitions needed here,
	including the meaning of a frame, a pullback, a trivial summand, and a
	Thom space. Our main application is a family indexed by
	$a\in\mathbb Z/32$: its rank-$14$ bundles $\xi_a\to S^{29}$ satisfy
	$\xi_a\oplus\varepsilon^1\cong\varepsilon^{15},$   but
	$\xi_a\cong\varepsilon^{14}$ only for $a=0.$  We shall also determine
	which $\xi_a$ admit a nowhere-zero section and describe their Thom
	spaces as mapping cones of elements of $\pi_{42}(S^{14}).$
	
	In this chapter, we \textbf{work integrally}. Homotopy groups are integral;
	localization at $2$ is indicated by the subscript $(2).$ 
	All vector bundles are real and oriented, and their isomorphisms cover
	the identity map of the stated base and preserve the fiber orientations.
	We write $\cong$ for a bundle isomorphism or a homeomorphism, with the
	objects specifying which is meant, and $\simeq$ for a homotopy
	equivalence. The bases used below are compact CW complexes, usually
	spheres. For bundles over $S^{k+1}$ we assume $k\geq1$ and $n\geq2.$ 
	The symbol $0\in\pi_k(\mathit{SO}(n))$
	is represented by the constant map with value $I_n,$   the identity
	matrix, rather than by the zero matrix.

	\section{Vector Bundles in Local Coordinates}
	
	\subsection{Fibers, Trivializations, and Transition Matrices}
	
	\begin{defn}[Vector Bundles and Local Trivializations]
		A rank-$n$ vector bundle $\xi$ over $B$ consists of a space $E(\xi),$  
		a continuous map $\pi\colon     E(\xi)\to B,$   and a real vector-space structure
		of dimension $n$ on each fiber $\xi_x=\pi^{-1}(x).$ 
		The required local condition is the following. Every $x\in B$ has
		an open neighborhood $U$ and a homeomorphism
		$\phi_U:\pi^{-1}(U)\to U\times\mathbb R^n$ such that
		$\operatorname{pr}_1\phi_U=\pi$ and
		$\phi_U|_{\xi_y}:\xi_y\to\{y\}\times\mathbb R^n$ is linear for
		every $y\in U.$  Such a map $\phi_U$ is a \emph{local trivialization}.
		The space $B$ is the base, $E(\xi)$ the total space, and
		$\operatorname{rank}\xi=n$ is the fiber dimension. 
	\end{defn}
	
	\begin{defn}[Transition Functions]
		On an overlap $U\cap V,$   the change of coordinates has the form
		$\phi_V\phi_U^{-1}(x,v)=(x,g_{VU}(x)v),$   where
		$g_{VU}\colon     U\cap V\to\mathit{GL}_n(\mathbb R)$ is continuous.
		Here $\mathit{GL}_n(\mathbb R)$ denotes the group of invertible real
		$n\times n$ matrices.
		The maps $g_{VU}$ are called the transition functions,
		and their values are the transition matrices.
	\end{defn}
	
	\begin{rem}
		The relations
		$g_{UU}=I_n,$   $g_{UV}=g_{VU}^{-1},$   and
		$g_{WU}=g_{WV}g_{VU}$ on triple overlaps express consistency of the
		fiber coordinates. Thus a bundle is described by local copies of
		$\mathbb R^n$ together with specified linear identifications on
		overlaps; the matrices $g_{VU}(x)$ describe these identifications; see
		\cite[Section 1.1, pp.~14--15]{HatcherVB}.
	\end{rem}
	
	\begin{defn}[Trivial Bundles]
		We use $\varepsilon_B^r$ to denote the trivial real vector
		bundle of rank $r$ over $B$: its projection is
		$B\times\mathbb R^r\to B$, $(x,v)\mapsto x$.
		When $B$ is fixed, we suppress the subscript.
	\end{defn}
	
	\begin{defn}[Orientations and Bundle Isomorphisms]
		An orientation of $\xi$ means that its trivializations can be chosen
		with $\det g_{VU}(x)>0$ on all overlaps. This condition identifies
		the usual orientation of $\mathbb R^n$ consistently in each fiber.
		If $\xi$ and $\eta$ are bundles over $B,$   a bundle isomorphism is a
		homeomorphism $F\colon     E(\xi)\to E(\eta)$ satisfying
		$\pi_\eta F=\pi_\xi$ whose restriction
		$F_x\colon     \xi_x\to\eta_x$ is a linear isomorphism for every $x.$ 
		In this chapter each $F_x$ also preserves orientation. A global
		trivialization is an isomorphism $\xi\cong\varepsilon_B^n.$
	\end{defn}
	
	The existence of a homeomorphism between total spaces alone would
	not supply the conditions $\pi_\eta F=\pi_\xi$ and linearity of $F_x.$ 
	For orientations and oriented isomorphisms, see
	\cite[Section 1.2, pp.~24--25]{HatcherVB}.

	\subsection{Metrics and Frames}
	
	\begin{defn}[Bundle Metrics]
		A metric on $\xi$ is a continuous choice of positive definite inner
		product $\langle-,-\rangle_x$ on $\xi_x.$
	\end{defn}
	
	\begin{prop}
		Every vector bundle over the bases considered here admits a metric.
	\end{prop}
	
	\begin{proof}
		Choose
		a finite trivializing cover $\{U_i\}$ and functions
		$\rho_i\colon     B\to[0,1]$ with support in $U_i$ and
		$\sum_i\rho_i(x)=1.$  If $v_i$ denotes the coordinates of
		$v\in\xi_x$ in the $i$th trivialization, define
		$\langle v,w\rangle_x=\sum_i\rho_i(x)(v_i\cdot w_i).$ 
		Each summand is extended by zero outside $U_i.$ 
		At least one weight is positive at every $x,$   so
		$\langle v,v\rangle_x>0$ for $v\neq0.$ 
		The functions $\rho_i$ are supplied by a partition of unity
		subordinate to the chosen trivializing cover.
	\end{proof}
	
	\begin{defn}[Frames]
		A frame in $\xi_x$ is an ordered basis $(v_1,\ldots,v_n)$ of that
		fiber. It is orthonormal when $\langle v_i,v_j\rangle_x=\delta_{ij},$  
		and oriented when its ordered basis gives the chosen orientation.
	\end{defn}
	
	\begin{lem}
		A local trivialization supplies continuous local basis sections;
		Gram--Schmidt changes them into continuous orthonormal sections.
	\end{lem}
	
	\begin{proof}
		Starting from a local basis $(w_1,\ldots,w_n),$   define
		$u_1=w_1/\|w_1\|$ and, recursively,
		$z_j=w_j-\sum_{i<j}\langle w_j,u_i\rangle u_i$ and
		$u_j=z_j/\|z_j\|.$  Linear independence gives $z_j\neq0,$   so every
		denominator is positive and each $u_j$ varies continuously.
		The positive normalization factors preserve the orientation of the
		ordered basis. This construction is also used in
		\cite[Section 1.1, Inner Products]{HatcherVB}.
	\end{proof}
	
	Two oriented orthonormal frames differ by a unique matrix in
	$$\mathit{SO}(n)=\{A\in M_n(\mathbb R)\colon     A^{\mathsf T}A=I_n,
	\ \det A=1\}.$$ Consequently we can take all transition matrices in
	$\mathit{SO}(n).$  This does not require the bundle to be trivial:
	the matrices can still vary with $x$ on the overlaps.
	It also does not add a new invariant to the oriented bundle.
	For example, Gram--Schmidt writes every
	$A\in\mathit{GL}_n^+(\mathbb R)$ as $A=QR,$   with
	$Q\in\mathit{SO}(n)$ and $R$ upper triangular with positive diagonal.
	The path $Q((1-t)R+tI_n),$   $0\leq t\leq1,$   stays invertible and
	ends at $Q.$  It gives a deformation retraction
	$\mathit{GL}_n^+(\mathbb R)\simeq\mathit{SO}(n).$ 
	This explains why the clutching classification can use orthogonal
	matrices; compare \cite[Section 1.2, Clutching Functions]{HatcherVB}.

	\section{From a Homotopy Class to a Bundle over a Sphere}
	
	\subsection{The Gluing Construction and Its Classification}
	
	\begin{defn}[Clutching Maps and Their Bundles]
		Write $S^{k+1}=D_+^{k+1}\cup D_-^{k+1},$   with intersection $S^k.$ 
		For a based map $\alpha:S^k\to\mathit{SO}(n),$   define $E_\alpha$
		as the quotient of
		$$(D_+^{k+1}\times\mathbb R^n)\amalg
		(D_-^{k+1}\times\mathbb R^n)$$ by
		$(x,v)_+\sim(x,\alpha(x)v)_-$ for $x\in S^k.$ 
		Projection to the first coordinate descends  to the map
		$\pi_\alpha:E_\alpha\to S^{k+1}.$ 
		Its fiber over any point is $\mathbb R^n$ with the vector-space
		operations induced by these coordinates.
		Denote this bundle by
		$\xi_\alpha.$  The map $\alpha$ is called its \textbf{clutching map};
		the word refers to the specified gluing along $S^k.$  The homotopy class $[\alpha]\in\pi_k(\mathit{SO}(n))$
		is called its clutching class.
	\end{defn}
	
	There is a local trivialization across the equator as well as over
	the interiors of the hemispheres. On a collar with coordinates
	$(x,t)\in S^k\times(-\delta,\delta),$   use the coordinate $v$ on
	the positive side and $\alpha(x)^{-1}w$ when $w$ is the negative-side
	coordinate. At $t=0$ these agree because $w=\alpha(x)v.$ 
	Thus $\pi_\alpha$ is a vector bundle. Its transition matrices have
	determinant $1,$   so it is oriented.
	
	\begin{prop}\label{geo:clutching}
		The assignment $[\alpha]\mapsto[\xi_\alpha]$ is a bijection from
		$\pi_k(\mathit{SO}(n))$ to the set of oriented rank-$n$ bundle
		isomorphism classes over $S^{k+1}.$  In particular,
		$\xi_\alpha\cong\varepsilon^n$ if and only if $[\alpha]=0.$ 
	\end{prop}
	\begin{proof}
		A vector bundle over a contractible CW complex is trivial
		\cite[Theorem 1.6 and Corollary 1.8]{HatcherVB}. For a bundle over
		$S^{k+1},$   choose oriented orthonormal trivializations over $D_+$
		and $D_-.$  Their change of coordinates on $S^k$ is a map
		$\alpha:S^k\to\mathit{SO}(n),$   and the original bundle is
		isomorphic to $\xi_\alpha.$  This proves surjectivity.
		
		An isomorphism between two glued bundles is represented on the
		hemispheres by maps $u_\pm:D_\pm^{k+1}\to
		\mathit{GL}_n^+(\mathbb R).$  Compatibility on the equator gives
		$\alpha'(x)=u_-(x)\alpha(x)u_+(x)^{-1}.$ 
		Since the disks are contractible and
		$\mathit{GL}_n^+(\mathbb R)$ is path connected, each $u_\pm$ is
		homotopic to the constant identity map. Thus $\alpha'$ and $\alpha$
		are homotopic in $\mathit{GL}_n^+(\mathbb R),$   and the deformation
		retraction above gives a homotopy in $\mathit{SO}(n).$ 
		Conversely, gluing with a homotopy
		$H:S^k\times I\to\mathit{SO}(n)$ gives a bundle over
		$S^{k+1}\times I$; its restrictions at the two ends are isomorphic
		by \cite[Proposition 1.7]{HatcherVB}.
		
		These arguments initially use unbased maps. Choose $x_0\in S^k.$ 
		Multiplying a map by the inverse of its value at $x_0$ normalizes it
		to have value $I_n$ there. If a homotopy $H$ has based endpoints,
		then $H'(x,t)=H(x_0,t)^{-1}H(x,t)$ is a based homotopy with the same
		endpoints. Hence the resulting set of classes is precisely
		$\pi_k(\mathit{SO}(n)).$  The constant identity clutching map glues
		to $S^{k+1}\times\mathbb R^n,$   proving the last assertion.\end{proof}
	
	\subsection{Pullback and Degree Maps}
	
	\begin{defn}[Pullback Bundles]
		For a map $f:X\to B$ and a bundle $\pi:E(\xi)\to B,$   define
		$$E(f^*\xi)=\{(x,v)\in X\times E(\xi):f(x)=\pi(v)\},$$
		with projection $(x,v)\mapsto x.$  Its fiber is
		$(f^*\xi)_x=\xi_{f(x)}.$
		This is the \emph{pullback bundle}.
	\end{defn}
	
	\begin{rem}
		If
		$\phi_U(v)=(\pi(v),v_U)$ is a trivialization of $\xi,$   then
		$(x,v)\mapsto(x,v_U)$ trivializes $f^*\xi$ over $f^{-1}(U).$ 
		The transition matrices become $g_{VU}(f(x)).$
	\end{rem}
	
	\begin{lem}
		$(f\circ h)^*\xi\cong h^*(f^*\xi).$ 
		Homotopic maps $f_0\simeq f_1$ give isomorphic bundles
		$f_0^*\xi\cong f_1^*\xi$.
	\end{lem}
	\begin{proof}See
		\cite[Proposition 1.5 and Theorem 1.6]{HatcherVB}.
	\end{proof}
	
	\begin{prop}
		If $h:S^{k+1}\to S^{k+1}$ has degree $d,$   then
		$h^*\xi_\alpha\cong\xi_{d\alpha}.$
	\end{prop}
	
	\begin{proof}
		Choose a
		degree-$d$ map $g:S^k\to S^k.$  Its suspension $\Sigma g$ has degree
		$d,$   so $h\simeq\Sigma g.$  The suspension respects the two
		hemispheres and restricts to $g$ on the equator. Therefore the
		clutching map of $(\Sigma g)^*\xi_\alpha$ is $\alpha\circ g.$ 
		Precomposition with $g$ multiplies the   element   in
		$\pi_k(\mathit{SO}(n))$ by $d,$   giving
		$[\alpha\circ g]=d[\alpha]$ and the asserted bundle isomorphism.
	\end{proof}
	
	\section{Trivial Summands and Nonzero Sections}
	
	\subsection{The Direct Sum Operation and Stable Triviality}
	
	\begin{defn}[Whitney Sums and Trivial Summands]
		For bundles $\xi$ and $\eta$ over the same base $B,$   their Whitney
		sum $\xi\oplus\eta$ has total space
		$$\{(v,w)\in E(\xi)\times E(\eta):\pi_\xi(v)=\pi_\eta(w)\},$$  
		projection $(v,w)\mapsto\pi_\xi(v),$   and fiber
		$(\xi\oplus\eta)_x=\xi_x\oplus\eta_x.$  A decomposition $\xi\cong\eta\oplus\varepsilon^r$
		is said to exhibit a trivial summand of rank $r$ in $\xi$.
	\end{defn}
	
	\begin{rem}
		If the summands have transition matrices $A(x)$ and $C(x),$   the
		sum has transition matrix $\operatorname{diag}(A(x),C(x))$;
		see \cite[Section 1.1, Direct Sums]{HatcherVB}.
		In particular, $\operatorname{rank}(\xi\oplus\eta)$ is equal to $
		\operatorname{rank}\xi+\operatorname{rank}\eta.$
	\end{rem}
	
	\begin{lem}
		For $N\geq n,$   define $i_{n,N}(A)=\operatorname{diag}(A,I_{N-n}).$    Then
		$\xi_{(i_{n,N})_*\alpha}\cong
		\xi_\alpha\oplus\varepsilon^{N-n}.$
	\end{lem}
	
	\begin{proof}
		In the gluing construction, the last $N-n$ coordinates glue by the identity, which is why they
		form a trivial summand.
	\end{proof}
	
	Fixed-rank addition has a different meaning:
	$\xi_{\alpha+\delta}$ has rank $n$ when
	$\alpha,\delta\in\pi_k(\mathit{SO}(n)),$   whereas
	$\xi_\alpha\oplus\xi_\delta$ has rank $2n.$

	\begin{defn}[Stable Triviality]
		A bundle $\xi$ of rank $n$ is \emph{stably trivial} if
		$\xi\oplus\varepsilon^r\cong\varepsilon^{n+r}$ for some integer
		$r\geq0.$  The number $r$ is allowed to exceed $1.$
	\end{defn}
	With this terminology in place, we may formulate the
	following characterization of stable triviality.
	\begin{prop}\label{geo:stable-trivial}
		The bundle $\xi_\alpha$ is stably trivial if and only if the image
		of $\alpha$ in
		$\pi_k(\mathit{SO})$
		is zero. In that case,
		$\xi_\alpha\oplus\varepsilon^{N-n}\cong\varepsilon^N$ for every
		$N\geq\max\{n,k+2\}.$ 
	\end{prop}
	\begin{proof}
		The block-matrix formula identifies vanishing at a finite rank $N$
		with $\xi_\alpha\oplus\varepsilon^{N-n}\cong\varepsilon^N.$ 
		An element of the direct limit is zero exactly when it becomes zero
		at some finite stage. For the bound, the fibration
		$\mathit{SO}(m)\to\mathit{SO}(m+1)\to S^m$ has
		$\pi_{k+1}(S^m)=\pi_k(S^m)=0$ when $m\geq k+2.$ 
		Exactness then makes the inclusion an isomorphism on $\pi_k.$ 
		Thus the stable image already determines the class in every rank
		$N\geq k+2.$  The next section gives the explicit bundle model of
		this fibration.
	\end{proof}
	
	\subsection{A Section and the Direct Summand That It Determines}
	
	\begin{defn}[Sections and Nowhere-Zero Sections]
		A section of $\pi:E(\xi)\to B$ is a continuous map
		$s:B\to E(\xi)$ with $\pi s=\operatorname{id}_B.$ 
		It is nowhere zero if $s(x)\neq0_x\in\xi_x$ for all $x.$    For a chosen metric, a unit section has length one in each fiber.
		Sections are independent if their values are linearly independent
		in every fiber. A global frame is a family of sections whose values
		form a basis in every fiber.
	\end{defn}
	With these definitions in place, we record the following
	basic result for later use.
	\begin{prop}
		A rank-$n$ oriented vector bundle $\xi$ admits a nowhere-zero section if and only if $\xi\cong\eta\oplus\varepsilon^1$ for some oriented rank-$(n-1)$ bundle $\eta$.
	\end{prop}

	\begin{proof}
		Given a nowhere-zero section and a chosen metric, replace it by the unit section
		$u(x)=s(x)/\|s(x)\|.$  Its span consists of the vectors $t u(x),$  
		and $(x,t)\mapsto t u(x)$ identifies this span with $\varepsilon^1.$ 
		Define the complementary fibers by
		$\eta_x=\{v\in\xi_x:\langle v,u(x)\rangle_x=0\}.$
		
		The fibers $\eta_x$ form a rank-$(n-1)$ vector bundle.
		One local justification is to choose, near a point $x_0,$   sections
		$w_1,\ldots,w_{n-1}$ such that
		$(w_1(x_0),\ldots,w_{n-1}(x_0),u(x_0))$ is a basis.
		The same ordered family remains a basis on a smaller neighborhood
		because its determinant is nonzero at $x_0$ and varies continuously.
		Orthogonal projection to $u(x)^\perp$ and Gram--Schmidt give a
		continuous local basis for $\eta_x.$  
		
		The splitting isomorphism is explicit:
		$\Phi_x:\eta_x\oplus\mathbb R\to\xi_x$ sends
		$(v,t)$ to $v+t u(x),$   and its inverse sends $w$ to
		$(w-\langle w,u(x)\rangle_xu(x),\langle w,u(x)\rangle_x).$ 
		These maps are continuous in $x.$  Orient $\eta_x$ so that an
		oriented basis of $\eta_x,$   followed by $u(x),$   is an oriented
		basis of $\xi_x.$  Then $\Phi$ preserves orientation.
		Conversely, an isomorphism $\eta\oplus\varepsilon^1\cong\xi$
		sends the section $x\mapsto(0_x,1)$ to a nowhere-zero section.
	\end{proof}
	
	The preceding proposition provides a useful starting point
	for the following observation.
	\begin{prop}
		For $1\leq r<n,$  
		\begin{equation}\label{geo:sections}
			\xi_\alpha\text{ admits $r$ independent sections}
			\;\Longleftrightarrow\;
			\alpha\in\operatorname{Im}\bigl(\pi_k(\mathit{SO}(n-r))
			\to\pi_k(\mathit{SO}(n))\bigr).
		\end{equation}
		Note that the  canonical inclusion on the right is $A\mapsto\operatorname{diag}(A,I_r).$
	\end{prop}
	
	\begin{proof}
		For $r$ everywhere independent sections, Gram--Schmidt supplies
		unit sections $u_1,\ldots,u_r$ satisfying
		$\langle u_i(x),u_j(x)\rangle_x=\delta_{ij}.$ 
		The projection onto their orthogonal complement is given by 
		$$w\mapsto w-\sum_{j=1}^r\langle w,u_j(x)\rangle_xu_j(x).$$
		The same argument gives $\xi\cong\eta\oplus\varepsilon^r$ with
		$\operatorname{rank}\eta=n-r.$
	\end{proof}
	
	\begin{rem}
		When $n-r=1,$   the complementary oriented line bundle is trivial:
		its positively oriented unit vector in each fiber defines a
		continuous global section, and $\mathit{SO}(1)=\{1\}.$ 
		This rank-one case is also explained in
		\cite[Section 1.2, p.~25]{HatcherVB}.
		For $r=n,$   a family of $n$ independent sections is a global frame
		and gives the trivialization
		$(x,t_1,\ldots,t_n)\mapsto\sum_jt_j u_j(x).$
	\end{rem}
	
	\begin{defn}[Reduction of the Structure Group]
		Sometimes~\eqref{geo:sections} is described as a reduction of the
		structure group to $\mathit{SO}(n-r).$  In this statement that phrase
		has a precise coordinate meaning: one can choose trivializations
		whose transition matrices have the form
		$\operatorname{diag}(A(x),I_r).$  The last $r$ basis vectors then
		agree on overlaps and define the $r$ global sections.
	\end{defn}
	
	\subsection{The Section Obstruction as an Element of a  Homotopy Group of a Sphere}
	
	\begin{defn}[The Section Obstruction]
		Define $p_n:\mathit{SO}(n)\to S^{n-1}$ by $p_n(A)=Ae_n,$   where
		$e_n$ is the last standard unit vector. 
		For the bundle with clutching class $\alpha$, the element
		$(p_n)_*\alpha\in\pi_k(S^{n-1})$ is its section obstruction.
	\end{defn}
	The significance of this obstruction is expressed
	in the following proposition.
	\begin{prop}
		For the bundle $\xi_\alpha$ over $S^{k+1}$,
		\begin{equation}\label{geo:section-test}
			\xi_\alpha\text{ admits a nowhere-zero section}
			\;\Longleftrightarrow\;
			(p_n)_*\alpha=0\in\pi_k(S^{n-1}).
		\end{equation}
	\end{prop}
	
	\begin{proof}
		For a clutching map $\alpha,$  
		choose the constant unit section $e_n$ on $D_+^{k+1}.$ 
		Its coordinates on the other side of the equator are
		$x\mapsto\alpha(x)e_n=p_n\alpha(x).$ 
		They extend to a map $D_-^{k+1}\to S^{n-1}$ exactly when
		$[p_n\alpha]=0\in\pi_k(S^{n-1}).$  Such an extension glues to a
		global unit section by the defining relation of $\xi_\alpha.$
		
		Conversely, write a global unit section in the two hemisphere
		coordinates as $$u_\pm:D_\pm^{k+1}\to S^{n-1}.$$ 
		On their boundary, $u_-(x)=\alpha(x)u_+(x).$ 
		Since $D_+$ is contractible and $S^{n-1}$ is path connected,
		$u_+$ is homotopic to the constant map $e_n.$  Multiplication by
		$\alpha(x)$ gives a boundary homotopy from $u_-|_{S^k}$ to
		$p_n\alpha.$  The first of these maps extends over $D_-,$   so
		$p_n\alpha$ is nullhomotopic.
	\end{proof}

	There are now two distinct tests, each with a specified target.
	The equality $(i_{n,n+1})_*\alpha=0$ in
	$\pi_k(\mathit{SO}(n+1))$ means
	$\xi_\alpha\oplus\varepsilon^1\cong\varepsilon^{n+1}.$ 
	The equality $(p_n)_*\alpha=0$ in $\pi_k(S^{n-1})$ means
	$\xi_\alpha\cong\eta\oplus\varepsilon^1$ for some $\eta.$ 
	The first statement adds a trivial line to the bundle; the second
	finds a trivial line inside it. Neither equality can be substituted
	for the other.

	\section{The Orthogonal Fibration and Pullbacks of the Tangent Bundle}
	
	\subsection{The Tangent Bundle of a Sphere without Manifold Notation}
	
	\begin{defn}[The Tangent Bundle of a Sphere]
		Let $S^n$ be the unit sphere in $\mathbb R^{n+1}.$
		Its tangent space at $x$ is
		$$x^\perp=
		\{v\in\mathbb R^{n+1}:\langle x,v\rangle=0\}.$$
		
		\noindent
		The tangent bundle is
		$$TS^n=\{(x,v)\in S^n\times\mathbb R^{n+1}:
		\langle x,v\rangle=0\},$$
		with projection $(x,v)\mapsto x.$
		This coincides with the standard tangent bundle of $S^n$. Let $p\in S^n.$ 
		Choose the orientation on $T_pS^n=\n(p,v)\in S^n\times\mathbb R^{n+1}:(p,v)\in TS^n\nn$ by declaring
		$(v_1,\ldots,v_n)$ positive when
		$\det(v_1,\ldots,v_n,p)=1$ for an orthonormal frame.
		This specifies the order in which the normal vector is appended.
	\end{defn}
	
	\begin{rem}
		For an explicit explanation, a differentiable curve $c(t)$ on $S^n$
		with $c(0)=x$ satisfies the relation 
		$0=(d/dt)_{t=0}\langle c(t),c(t)\rangle=2\langle x,c'(0)\rangle.$ 
		Conversely, for $0\neq v\in x^\perp,$   the curve
		$$c(t)=\cos(\|v\|t)x+\sin(\|v\|t)v/\|v\|$$ lies on $S^n$ and has
		$c'(0)=v.$  The zero vector is the velocity of the constant curve.
		Thus the velocity-vector description gives exactly $x^\perp.$
		
		Its local triviality can be checked directly. Fix $x_0$ and an
		orthonormal basis $w_1,\ldots,w_n$ of $x_0^\perp.$ 
		For $x$ near $x_0,$   the projected vectors
		$w_j(x)=w_j-\langle w_j,x\rangle x$ remain independent because
		their Gram matrix is the identity at $x_0$ and has nonzero
		determinant nearby. Applying Gram--Schmidt gives continuous
		orthonormal vectors $v_1(x),\ldots,v_n(x)$ in $x^\perp.$ 
		The map $(x,t_1,\ldots,t_n)\mapsto(x,\sum_jt_jv_j(x))$ is a local
		trivialization of $TS^n.$
	\end{rem}
	We now record a simple but useful observation.
	\begin{prop}
		With the orientation specified above,
		$TS^n\oplus\varepsilon^1\cong\varepsilon^{n+1}.$
	\end{prop}
	
	\begin{proof}
		The map
		$TS^n\oplus\varepsilon^1\to\varepsilon^{n+1}$ given by
		$(x,v,t)\mapsto(x,v+tx)$ is orientation-preserving.
		Its inverse is
		$(x,w)\mapsto(x,w-\langle w,x\rangle x,\langle w,x\rangle).$
	\end{proof}
	
	\subsection{Frames and the Last-Column Fibration}
	
	\begin{defn}[Oriented Frame Bundles and Principal Bundles]
		For an oriented metric bundle $\xi,$   its oriented orthonormal frame
		bundle $F^+(\xi)$ has points $(x;v_1,\ldots,v_n),$   where
		$(v_1,\ldots,v_n)$ is an oriented orthonormal basis of $\xi_x.$ 
		The projection sends this point to $x.$ 
		If $v$ denotes the matrix with these basis vectors as columns, the
		right action of $A\in\mathit{SO}(n)$ is $v\mapsto vA.$ 
		It is free and transitive on the frames over each fixed $x$:
		for two such frames there is a unique matrix $A$ relating them.
		A local orthonormal frame identifies $F^+(\xi)|_U$ with
		$U\times\mathit{SO}(n).$  A bundle with this local description and
		right action is called a principal $\mathit{SO}(n)$-bundle. \end{defn}

	We next consider a standard example that provides
	a concrete illustration of these notions.
	\begin{prop}
		For $TS^n,$   an element of $\mathit{SO}(n+1)$ can be written as
		$A=(v_1,\ldots,v_n,x),$   and its first $n$ columns are exactly a
		frame of $T_xS^n$ with the orientation just specified.
		Thus $F^+(TS^n)$ is identified with $\mathit{SO}(n+1),$   and its
		projection is $p_{n+1}(A)=Ae_{n+1}=x.$ 
		Its fiber over $e_{n+1}$ consists of
		$\operatorname{diag}(B,1),$   $B\in\mathit{SO}(n).$
	\end{prop}
	
	\begin{proof}
		For completeness, the local product description is also explicit.
		On a neighborhood $U$ with the frame constructed above, let
		$Q(x)=(v_1(x),\ldots,v_n(x),x)\in\mathit{SO}(n+1).$ 
		Observe     the map 
		$$   U\times\mathit{SO}(n)\to p_{n+1}^{-1}(U),\quad (x,B)\mapsto       Q(x)\operatorname{diag}(B,1).  $$
		Its inverse sends $A$ to $x=p_{n+1}(A)$ and the unique block $B$
		in $Q(x)^{\mathsf T}A=\operatorname{diag}(B,1).$ 
		This realizes the fibration
		$$\mathit{SO}(n)\xrightarrow{i_{n,n+1}}\mathit{SO}(n+1)
		\xrightarrow{p_{n+1}}S^n$$ as a concrete locally trivial bundle.
		Its fiber-bundle and lifting properties are the classical ones
		described in \cite{St}.
	\end{proof}
	
	\subsection{Why the Boundary Class Is a Clutching Class}
	
	\begin{defn}[The Boundary Homomorphism]
		Let $f:S^{k+1}\to S^n$ be based, with the basepoint of $S^n$ equal
		to $e_{n+1}.$  Over each hemisphere choose a lift
		$Q_\pm:D_\pm^{k+1}\to\mathit{SO}(n+1)$ satisfying
		$p_{n+1}Q_\pm=f|_{D_\pm}.$  Such lifts exist: the pullback of
		$TS^n$ over a disk is trivial, and a trivialization followed by
		Gram--Schmidt gives the first $n$ columns of $Q_\pm.$ 
		Their last columns are both $f(x).$ 
		On the equator there is consequently a unique
		$\alpha_f:S^k\to\mathit{SO}(n)$ such that
		$Q_+(x)=Q_-(x)\operatorname{diag}(\alpha_f(x),1).$
		
		The boundary homomorphism
		$\Delta\colon \pi_{k+1}(S^n)\to\pi_k(\mathit{SO}(n))$ can be
		defined by comparing lifts on the two hemispheres.
		We use the comparison $Q_+=Q_-\operatorname{diag}(\alpha_f,1),$  
		so that $\Delta[f]=[\alpha_f].$
	\end{defn}
	
	\begin{prop}
		For the boundary convention just specified,
		$f^*(TS^n)\cong\xi_{\Delta[f]}.$
	\end{prop}
	
	\begin{proof}
		The fiber coordinate $v\in\mathbb R^n$ on the positive side
		represents the vector $Q_+(x)(v,0)\in f(x)^\perp.$ 
		The displayed matrix identity gives
		$Q_+(x)(v,0)=Q_-(x)(\alpha_f(x)v,0).$ 
		Thus the negative-side coordinate is $\alpha_f(x)v,$   exactly the
		clutching convention used for $\xi_{\alpha_f}.$ 
		We have proved $f^*(TS^n)\cong\xi_{\alpha_f}.$
	\end{proof}
	
	This is the frame-bundle interpretation of the boundary construction;
	see \cite{St} and \cite[Section 1.2, Clutching Functions]{HatcherVB}.
	If the comparison is reversed, the boundary map changes to
	$-\Delta.$  The cyclic subgroups and orders used below are unchanged
	by that replacement.

	\begin{prop}
		$\xi_\alpha\oplus\varepsilon^1\cong\varepsilon^{n+1}$ if and
		only if $\alpha=\Delta[f]$ for some $f:S^{k+1}\to S^n.$ 
		Hence the bundles trivialized by one additional line are precisely
		the bundles isomorphic to these pullbacks $f^*(TS^n).$
	\end{prop}
	
	\begin{proof}
		The pullback trivialization is
		$f^*(TS^n)\oplus\varepsilon^1\to\varepsilon^{n+1},$  
		$(x,v,t)\mapsto(x,v+t f(x)).$ 
		Conversely, exactness gives
		$\ker((i_{n,n+1})_*)=\operatorname{Im}\Delta$ in
		$\pi_k(\mathit{SO}(n)).$
	\end{proof}
	
	\section{The Integral Family of Rank Fourteen Bundles}
	\label{geo:rank14-family}
	
	\subsection{Passing from a Localized Class to an Actual Clutching Map}
	
	\begin{defn}[The Integral $2$-Primary Subgroup]
		For an abelian group $A,$   define
		$$A[2^\infty]=\{a\in A:2^r a=0\text{ for some }r\geq0\}$$ and
		$A_{(2)}=A\otimes\mathbb Z_{(2)}.$ 
		An element of $A_{(2)}$ is represented by $a/s$ with $s$ odd.
	\end{defn}
	
	The groups calculated earlier are localized, while a bundle is
	constructed from an integral map $S^k\to\mathit{SO}(n).$ 
	The following algebraic fact makes the passage precise.
	\begin{rem}\label{geo:integral-torsion}
		Localization induces a natural isomorphism
		$A[2^\infty]\to\operatorname{Tor}(A_{(2)}),$
		where $\operatorname{Tor}$ denotes the subgroup of elements
		of finite order. This isomorphism preserves the order of
		each element. Consequently, every finite-order element of
		$A_{(2)}$ has a unique preimage in $A[2^\infty]$ under
		localization. The uniqueness is asserted within
		$A[2^\infty]$, not among all elements of $A$.
		
		Let $F$ be a homomorphism from $A$ to an abelian group.
		For $a\in A[2^\infty]$, the element $F(a)$ also has
		$2$-power order. Hence
		$\operatorname{ord}(F(a))
		=\operatorname{ord}((F(a))_{(2)})$;
		in particular,
		$F(a)=0$ if and only if $(F(a))_{(2)}=0$.
		The integral inclusion, projection, boundary, and
		$J$ homomorphisms used here commute with localization,
		so this observation applies to their restrictions to
		the integral $2$-primary subgroups.\end{rem}
	
	\subsection{The Thirty-Two Isomorphism Classes}
	
	Let $r\in\pi_{29}(S^{14})[2^\infty]$ be the integral lift of
	$\rho_{14}$. Appendix~\ref{tlqdzj}, Section~\ref{app:stem-15},
	tells us that $\operatorname{ord}(r)=32$.
	For the boundary homomorphism
	$\Delta\colon\pi_{29}(S^{14})\to\pi_{28}(\mathit{SO}(14))$,
	define 
	\begin{equation}\label{gs1028}
		\beta=\Delta r\in\pi_{28}(\mathit{SO}(14)).
	\end{equation}
	\noindent   Proposition~\ref{r2814} and its proof give
	$\pi_{28}(\mathit{SO}(14))[2^\infty]
	=\langle\beta\rangle\cong\mathbb Z/32$.

	\begin{defn}[The Bundles $\xi_a$]
		For $a\in\mathbb Z/32$, write $$\xi_a=\xi_{a\beta}$$ for the
		oriented real vector bundle of rank $14$ over $S^{29}$
		whose clutching class is $a\beta$.
		More explicitly, choose a based map
		$\alpha\colon S^{28}\to\mathit{SO}(14)$ with $[\alpha]=a\beta$.
		Writing $S^{29}=D_+^{29}\cup D_-^{29}$ with
		$D_+^{29}\cap D_-^{29}=S^{28}$, construct $E(\xi_a)$
		as the quotient of
		$$(D_+^{29}\times\mathbb R^{14})\amalg
		(D_-^{29}\times\mathbb R^{14})$$
		by the identifications
		$(x,v)_+\sim(x,\alpha(x)v)_-$
		for $x\in S^{28}$ and $v\in\mathbb R^{14}$.
		The bundle projection is
		$\pi_{\xi_a}\colon E(\xi_a)\to S^{29}$,
		$[(x,v)_\pm]\mapsto x$.
		Thus $\xi_a$ has total space $E(\xi_a)$, base space $S^{29}$,
		structure group $\mathit{SO}(14)$, and fibers
		$(\xi_a)_x=\pi_{\xi_a}^{-1}(x)\cong\mathbb R^{14}$.
	\end{defn}
	
	By Proposition~\ref{geo:clutching}, the oriented isomorphism
	class of $\xi_a$ depends only on $[\alpha]=a\beta$.
	Since $32\beta=0$, this isomorphism class is independent
	of the chosen integer representative of $a$ and of the
	chosen clutching map representing $a\beta$.

	Choose $f_a:S^{29}\to S^{14}$ representing $ar.$ 
	Since $\Delta(ar)=a\beta,$   the preceding section supplies the
	explicit model $$E(\xi_a)=\{(x,v)\in S^{29}\times\mathbb R^{15}:
	\langle v,f_a(x)\rangle=0\},$$ with projection $(x,v)\mapsto x.$ 
	Its fiber over $x$ is the hyperplane $f_a(x)^\perp.$ 
	This formula means $$\xi_a\cong f_a^*(TS^{14});$$ it gives a bundle
	model for each specified homotopy class without requiring a
	coordinate formula for a representative of $\rho_{14}.$

	\begin{prop}\label{geo:rank14-bundles}
		The bundles $\xi_a,$   $a\in\mathbb Z/32,$   are pairwise
		nonisomorphic as oriented bundles over $S^{29}.$ 
		Each satisfies $\xi_a\oplus\varepsilon^1\cong\varepsilon^{15}.$ 
		For a degree-$d$ map $h:S^{29}\to S^{29},$  
		$h^*\xi_a\cong\xi_{da}$; hence
		$h^*\xi_1\cong\varepsilon^{14}$ if and only if $32\mid d.$ 
	\end{prop}
	\begin{proof}
		By Proposition~\ref{geo:clutching}, $\xi_a\cong\xi_{a'}$ is
		equivalent to $a\beta=a'\beta.$  Since $\beta$ has order $32,$  
		this is equivalent to $a=a'$ in $\mathbb Z/32.$ 
		In the hyperplane model the isomorphism
		$\xi_a\oplus\varepsilon^1\to\varepsilon^{15}$ is
		$(x,v,t)\mapsto(x,v+t f_a(x)).$ 
		Its inverse is
		$(x,w)\mapsto(x,w-\langle w,f_a(x)\rangle f_a(x),
		\langle w,f_a(x)\rangle).$ 
		The degree formula for pullbacks gives
		$h^*\xi_a\cong\xi_{da},$   and $d\beta=0$ exactly when $32\mid d.$ 
	\end{proof}
	
	\begin{rem}
		These $32$ bundles exhaust the classes with clutching element in
		$\pi_{28}(\mathit{SO}(14))[2^\infty].$ 
		No assertion about an odd-primary subgroup is needed.
		For example, $\xi_1$ is nontrivial, its pullback by a degree-$16$
		map is the nontrivial bundle $\xi_{16},$   and its pullback by a
		degree-$32$ map is $\xi_0\cong\varepsilon^{14}.$ 
		Each statement follows from the corresponding equality or
		inequality among $\beta,$   $16\beta,$   and $32\beta.$
	\end{rem}
	
	\subsection{Exactly Which Members Admit a Nonzero Section}
	
	Recall  that    $\beta=\Delta r\in\pi_{28}(\mathit{SO}(14))$ from Equation \eqref{gs1028}.
	Let $c=(p_{14})_*\beta\in\pi_{28}(S^{13})$.
	We have
	\[
	p_*\colon \pi_{28}(\mathit{SO}(14))_{(2)}
	\to\pi_{28}(S^{13})_{(2)},
	\qquad
	\Delta\rho_{14}\mapsto 2\rho_{13}.
	\]
	Since $\rho_{13}$ has order $32$ by
	Appendix~\ref{tlqdzj}, Section~\ref{app:stem-15}, we obtain
	\[
	\operatorname{ord}\bigl(p_*(\Delta\rho_{14})\bigr)=16,
	\qquad
	\ker\bigl(p_*|_{\langle\Delta\rho_{14}\rangle}\bigr)
	=\langle16\Delta\rho_{14}\rangle.
	\]  \noindent   The sign of the boundary comparison does not change this kernel.
	Remark~\ref{geo:integral-torsion} therefore gives the integral
	statements $\operatorname{ord}(c)=16$ and
	$\ker((p_{14})_*|_{\langle\beta\rangle})
	=\langle16\beta\rangle.$
	
	\begin{prop}\label{geo:rank14-sections}
		For $a\in\mathbb Z/32,$   the bundle $\xi_a$ admits a nowhere-zero
		section if and only if $a=0$ or $a=16.$ 
	\end{prop}
	\begin{proof}
		Its section obstruction is the explicitly defined element
		$(p_{14})_*(a\beta)=ac\in\pi_{28}(S^{13}).$ 
		Since $\operatorname{ord}(c)=16,$   we have $ac=0$ if and only if
		$16\mid a.$  The two such residue classes modulo $32$ are $0$
		and $16.$  Applying~\eqref{geo:section-test}, we deduce the result.
	\end{proof}
	
	\begin{rem}
		In particular, $\xi_{16}$ has a trivial line summand although
		$\xi_{16}\not\cong\varepsilon^{14},$   since $16\beta\neq0.$ 
		For every $a\notin\{0,16\},$   the element $ac$ prevents a section
		while the explicit map in Proposition~\ref{geo:rank14-bundles}
		still trivializes $\xi_a\oplus\varepsilon^1.$ 
		Thus this single family distinguishes the conditions
		$(p_{14})_*(a\beta)=0$ and $(i_{14,15})_*(a\beta)=0.$
	\end{rem}
	
	\section{The Thom Space and Its Attaching Map}
	
	\subsection{Disk Bundles, Sphere Bundles, and the Quotient}
	
	\begin{defn}[Disk and Sphere Bundles]
		For a metric bundle $\xi,$   define
		$$D(\xi)=\{v\in E(\xi):\|v\|\leq1\},\;\;S(\xi)=\{v\in E(\xi):\|v\|=1\}.$$ 
		Their projections to $B$ have fibers $D^n$ and $S^{n-1},$  
		respectively. Orthonormal trivializations identify the restricted
		pairs with $(U\times D^n,U\times S^{n-1}).$
	\end{defn}
	The preceding definitions provide the setting
	for another standard construction.
	\begin{defn}[Thom Spaces]
		The \textbf{Thom space} is the based quotient
		$$T(\xi)=D(\xi)/S(\xi),$$   with all points of $S(\xi)$ identified
		to one basepoint. Thus boundaries from different fibers are
		identified with each other, as well as being collapsed within
		each fiber. 
	\end{defn}
	With this definition in place, it is useful to record
	a few elementary observations.
	\begin{prop}
		For the trivial bundle,
		$$T(\varepsilon_B^n)=(B\times D^n)/(B\times\partial D^n)
		\cong B_+\wedge S^n,$$ where $B_+=B\amalg\{*\}$ has a disjoint
		basepoint and $S^n=D^n/\partial D^n.$ 
		
		A bundle isomorphism induces a homeomorphism of Thom spaces.
	\end{prop}
	
	\begin{proof}
		This formula follows from the definition of the smash product:
		the subsets $B_+\times\{*\}$ and $\{*\}\times S^n$ are collapsed.

		For a bundle isomorphism $F$ that is not an isometry, the fiber map
		$v\mapsto\|v\|F(v)/\|F(v)\|$ for $v\neq0,$   with $0\mapsto0,$  
		maps disk to disk and sphere to sphere and realizes the induced
		quotient homeomorphism. The same radial construction compares
		Thom spaces formed from different metrics.
	\end{proof}
	
	\subsection{The Attaching Map on Both Boundary Faces}

	Let $m=k+1,$   and write the base sphere as $S^m=D^m/\partial D^m,$  
	with the collapsed boundary as its $0$-cell. A model for the bundle
	with clutching map $\alpha:S^{m-1}\to\mathit{SO}(n)$ is obtained
	from $(D^m\times\mathbb R^n)\amalg\mathbb R^n$ by identifying
	$(x,v)$ with $\alpha(x)v$ in the second copy whenever
	$x\in\partial D^m.$  The second copy is the fiber over the
	$0$-cell. In a collar of $\partial D^m,$   the coordinate
	$w=\alpha(x)v$ gives the trivialization at that basepoint.
	Here $D^m$ represents the positive-dimensional cell of $S^m,$  
	and its boundary is glued to the remaining fiber by $\alpha$; compare
	\cite[Lemma 4.8 and its proof]{HatcherVB}.

	\begin{defn}[The Attaching Map]
		After taking fiber disks and collapsing the sphere bundle, the
		second disk becomes $D^n/\partial D^n=S^n.$ 
		The first part contributes the disk $D^m\times D^n,$   of dimension
		$m+n.$  Its boundary is
		$$(\partial D^m\times D^n)\cup(D^m\times\partial D^n).$$
		Let $q\colon D^n\to D^n/\partial D^n$ be the quotient map.
		The attaching map $a_\alpha\colon\partial(D^m\times D^n)\to S^n$
		is specified by two formulas:
		$a_\alpha(x,v)=q(\alpha(x)v)$ on $\partial D^m\times D^n,$  
		and $a_\alpha(x,v)=*$ on $D^m\times\partial D^n.$
	\end{defn}
	
	\begin{rem}
		On the intersection, $\|v\|=1$ implies
		$\|\alpha(x)v\|=1,$   so the first formula also has value $*.$ 
		The two maps therefore agree and define a continuous map.
	\end{rem}
	We wish to draw the reader's attention to the following proposition, which we find particularly appealing for its elegant formulation and the mathematical insight it provides.
	
	\begin{prop}
		The Thom space of $\xi_\alpha$ has the following description:
		\begin{equation}\label{geo:thom}
			T(\xi_\alpha)\simeq S^n\cup_{J(\alpha)}e^{n+k+1},
			\qquad J(\alpha)\in\pi_{n+k}(S^n).
		\end{equation}
	\end{prop}
	
	\begin{proof}
		The boundary of $D^m\times D^n$ is a sphere of dimension
		$m+n-1=n+k.$  The action of $\mathit{SO}(n)$ on $D^n$ used in
		the formula for $a_\alpha$ is the disk description of the unstable
		$J$-homomorphism \cite{Wh}. Orient the characteristic disk so
		that $[a_\alpha]=J(\alpha)$ in the convention of
		Section~\ref{m2:subsec:j-toda}.
	\end{proof}
	
	This real-bundle construction uses exactly the two boundary faces
	specified above; the complex-bundle version is
	\cite[Lemma 4.8]{HatcherVB}.
	For a homotopy theorist, the resulting object is the mapping cone
	of the map $a_\alpha\colon S^{n+k}\to S^n.$ 
	The fiber dimension $n$ determines the bottom cell, and adding the
	base-cell dimension $k+1$ gives the top-cell dimension $n+k+1.$
	
	\subsection{The Exact Condition for a Wedge Decomposition}
	
	\begin{lem}\label{geo:thom-wedge}
		For $g\in\pi_{n+k}(S^n),$   where $n\geq2$ and $k\geq1,$  
		$S^n\cup_g e^{n+k+1}\simeq S^n\vee S^{n+k+1}$ if and only if
		$g=0.$  Consequently,
		$T(\xi_\alpha)\simeq S^n\vee S^{n+k+1}$ if and only if
		$J(\alpha)=0.$ 
	\end{lem}
	\begin{proof}
		If $g=0,$   its mapping cone is the stated wedge.
		Conversely, let $X=S^n\cup_g e^{n+k+1}$ and let
		$F:X\to S^n\vee S^{n+k+1}$ be a homotopy equivalence.
		Denote the bottom-cell inclusion by $j:S^n\to X$ and the
		projection of the wedge onto its first summand by
		$q:S^n\vee S^{n+k+1}\to S^n.$ 
		The inequality $n+k+1\geq n+2$ shows that $j$ and $q$ induce
		isomorphisms on $H_n(-;\mathbb Z).$ 
		Therefore $qFj:S^n\to S^n$ has degree $1$ or $-1$ and has a
		homotopy inverse $u:S^n\to S^n.$ 
		The characteristic map of the top cell extends $jg$ over
		$D^{n+k+1},$   so $jg\simeq*.$ 
		It follows that $g\simeq u qFjg\simeq*.$ 
		Apply this result to $g=J(\alpha)$ for the last assertion.
	\end{proof}
	
	For the constant identity clutching map, the boundary formula
	for $a_{I_n}$ extends over all of $D^m\times D^n$ by
	$(x,v)\mapsto q(v).$  Hence $[a_{I_n}]=0$ and
	$T(\varepsilon_{S^{k+1}}^n)\simeq S^n\vee S^{n+k+1}.$ 
	In general it is $J(\alpha),$   not the full class $\alpha,$   that
	this wedge criterion tests. Formula~\eqref{geo:thom} alone does
	not imply injectivity of $J.$

	\subsection{Adding a Line as a Quotient-Space Suspension}
	We next record another useful property of Thom spaces.
	\begin{prop}
		For a vector bundle $\xi$,
		$T(\xi\oplus\varepsilon^1)\cong\Sigma T(\xi)$.
	\end{prop}
	
	\begin{proof}
		Let $I=[-1,1].$ 
		For the direct sum, use the fiber disk
		$\{(v,t):\|v\|\leq1,\ |t|\leq1\},$   whose boundary is
		$\{\|v\|=1\}\cup\{|t|=1\}.$ 
		This disk is radially homeomorphic to the usual Euclidean disk.
		Explicitly, for $(v,t)\neq(0,0),$   multiplication by
		$\max\{\|v\|,|t|\}/\sqrt{\|v\|^2+t^2}$ sends it to that disk;
		the origin maps to the origin. This map and its inverse vary
		continuously with the base point and send boundary to boundary.
		
		The quotient defining its Thom space is therefore
		$(D(\xi)\times I)/((S(\xi)\times I)\cup
		(D(\xi)\times\partial I)).$ 
		Collapsing the indicated subsets identifies this space with
		$$(D(\xi)/S(\xi))\wedge(I/\partial I)
		=T(\xi)\wedge S^1=\Sigma T(\xi).$$
		This proves the claimed homeomorphism, with $\Sigma$ the reduced
		suspension. It is also the usual compatibility of Thom spaces
		with a trivial line; see
		\cite[Section 4.1, Thom Spaces]{HatcherVB} and \cite{Wh}.
	\end{proof}
	We may now bring together the preceding ideas
	to obtain the following conclusion.
	\begin{prop}
		If $\xi_\alpha\oplus\varepsilon^1\cong\varepsilon^{n+1}$, then
		$\Sigma T(\xi_\alpha)\simeq S^{n+1}\vee S^{n+k+2}$ and
		$\Sigma J(\alpha)=0$.
	\end{prop}
	
	\begin{proof}
		Suspending a mapping cone suspends its attaching map. Applying
		this to~\eqref{geo:thom} gives
		$\Sigma T(\xi_\alpha)\simeq
		S^{n+1}\cup_{\Sigma J(\alpha)}e^{n+k+2}.$ 
		If $\xi_\alpha\oplus\varepsilon^1\cong\varepsilon^{n+1},$  
		the preceding quotient calculation and the trivial-bundle case
		give $\Sigma T(\xi_\alpha)\simeq S^{n+1}\vee S^{n+k+2}.$ 
		Lemma~\ref{geo:thom-wedge}, applied after suspension, then gives
		$\Sigma J(\alpha)=0.$  This proof needs no choice of sign in a
		comparison between $J$ in two successive ranks.
	\end{proof}
	
	\section{An Attaching Class of Order Thirty-Two Killed by Suspension}
	\label{geo:order32-thom}
	
	\begin{defn}[The Attaching Class $\gamma$]
		Return to $\beta=\Delta r\in\pi_{28}(\mathit{SO}(14))$ and define
		$\gamma=J(\beta)\in\pi_{42}(S^{14}).$
	\end{defn}
	
	\begin{rem}
		The proof of Proposition~\ref{r2814} identifies its localized
		image with $P\rho_{29}$ up to sign and shows that this image has
		order $32.$  Since $32\beta=0,$   additivity of $J$ gives
		$32\gamma=0.$  Remark~\ref{geo:integral-torsion} therefore implies
		the integral equality $\operatorname{ord}(\gamma)=32.$ 
		Only this order, and not the choice between the two signs, is
		used in the following proposition.
	\end{rem}
	The preceding discussion culminates in the following
	result, which brings together several of the ideas
	developed above.
	\begin{prop}\label{jhyy1027}
		For $a\in\mathbb Z/32,$   the Thom space of $\xi_a$ satisfies
		$$T(\xi_a)\simeq S^{14}\cup_{a\gamma}e^{43},$$ where
		$\operatorname{ord}(a\gamma)=32/\gcd(32,a).$ 
		It is homotopy equivalent to $S^{14}\vee S^{43}$ if and only if
		$a=0.$  For every $a,$  
		$\Sigma T(\xi_a)\simeq S^{15}\vee S^{44}.$ 
		In particular, $\gamma\neq0$ has order $32$ but
		$\Sigma\gamma=0\in\pi_{43}(S^{15}).$ 
	\end{prop}
	\begin{proof}
		The clutching class of $\xi_a$ is $a\beta.$ 
		Formula~\eqref{geo:thom} and the equality
		$J(a\beta)=aJ(\beta)=a\gamma$ give the two-cell description.
		In a cyclic group of order $32,$   the order of $a\gamma$ is
		$32/\gcd(32,a).$  For $a=0$ we use $\gcd(32,0)=32,$   so that
		the zero element has order $1.$ 
		Lemma~\ref{geo:thom-wedge} gives the asserted equivalence between
		the wedge decomposition and $a=0.$
		
		Proposition~\ref{geo:rank14-bundles} gives
		$\xi_a\oplus\varepsilon^1\cong\varepsilon^{15}.$ 
		Thus $\Sigma T(\xi_a)\cong T(\xi_a\oplus\varepsilon^1)
		\cong T(\varepsilon_{S^{29}}^{15})
		\simeq S^{15}\vee S^{44}.$ 
		For $a=1,$   the suspended two-cell description has attaching map
		$\Sigma\gamma:S^{43}\to S^{15}.$ 
		The wedge criterion shows that this map is nullhomotopic, while
		the order calculation gives $\gamma\neq0$ before suspension.
	\end{proof}
	
	For $a\in\{1,2,4,8,16\}$, the Thom space $T(\xi_a)$
	has the same integral homology as $S^{14}\vee S^{43}$,
	but
	\[
	T(\xi_a)\not\simeq S^{14}\vee S^{43},
	\]
	since its attaching map $a\gamma$ is nonzero.  The bundle classification does not imply that the Thom spaces are
	pairwise inequivalent. To see the distinction explicitly, let
	$R:D^{43}\to D^{43}$ be a linear reflection and let
	$r=R|_{S^{42}},$   so that $\deg r=-1.$ 
	For a based map $g:S^{42}\to S^{14},$   the identity on $S^{14}$
	together with $R$ on the top disk gives a homeomorphism
	$S^{14}\cup_{g\circ r}D^{43}\cong S^{14}\cup_g D^{43}.$ 
	Since $[g\circ r]=-[g],$   we have
	$T(\xi_a)\simeq T(\xi_{-a}).$ 
	For example, $\xi_1\not\cong\xi_{31}$ as oriented bundles,
	while $T(\xi_1)\simeq T(\xi_{31}).$ 
	The order-$32$ attaching class and its vanishing suspension
	therefore give the stated wedge obstruction without a claim that
	all $32$ Thom spaces have different homotopy types.
	
	\appendix
	\chapter{Reference Tables}\label{app:oguchi}
	This appendix serves as a reference for $\pi_{*}(S^{n})$ (after localization at 2). 
	It consists of two parts: 
	a table of relations among generators in Toda's range
	(stems $0$--$19$)  mainly from \^{O}guchi's paper \cite{Og};
	and a table of relations in stems $20$ and above,
	mainly  collected from Oda's work \cite{Oda1979}.

	The homotopy groups in 0-32 stems and some homotopy groups in 33 stem are given in Appendix~\ref{tlqdzj}, with clear generating sets.
	\section{Reading the Tables}

	Recall that we work 2-locally.
	Roughly speaking, if $\mathbf{x}$ represents a Greek letter,
	$\mathbf{x}_n$ denotes one of the generators of
	$\pi_{n+r}(S^n)$ for some $r$\uprightsemicolon{}
	the subscript $n$ indicates the codomain of $\mathbf{x}_n.$ 
	Moreover, $\mathbf{x}_{n+k}:=\Sigma^k\mathbf{x}_n,$  
	$\mathbf{x}:=\Sigma^\infty\mathbf{x}_n,$  
	and $\mathbf{x}_n^\ell$ is an abbreviation for
	$$\mathbf{x}_n\circ\mathbf{x}_{n+r}\circ\cdots
	\circ\mathbf{x}_{n+(\ell-1)r},\,\; (\ell \;\text{factors}).$$
	The uses of $\overline{\mathbf{x}}_n$ and $\mathbf{x}_n^*$
	are similar to those above.

	In $\pi_{j+r}(S^j)$ (not a stable homotopy group),
	if a generator is written without a subscript,
	then this generator  does not survive in the stable homotopy group
	$\pi_r^S(S^0)$ or its $\Sigma^\infty$-image is divisible by $2.$  For instance, consider $\theta\in\pi_{24}(S^{12})$ and
	$\sigma'''\in\pi_{12}(S^5).$ 
	Their $\Sigma^\infty$-images satisfy
	$\Sigma^\infty\theta=0,$   $\Sigma^\infty\sigma'''=8\sigma$
	having order $2.$ 
	
	If a generator is written with a subscript,
	then this generator    survives in the stable homotopy group
	$\pi_r^S(S^0)$ and  its $\Sigma^\infty$-image is not divisible by $2,$   except for the following elements ($n\geq0$):
	$$  \text{the elements}\; \cg_{16+n}^{*}\;\text{ in Stem 22,  the elements}\; \bar{\nu}^\ast_{16+n},\;\bar{\lt}^*_{16+n}\;\text{ in Stem 23, } $$$$\text{the elements}\; \cs_{6+n}',\ccs_{6+n}'\text{ in Stem 24,}\quad\text{the elements}\; \bar{\bar{\lt}}_{18+n}\text{ in Stem 32}.$$
	\noindent Each displayed element does not strictly survive in the stable homotopy group
	$\pi_r^S(S^0),$   or its $\Sigma^\infty$-image is divisible by $2.$  \footnote{In fact, these  are  oversights in \cite{Mi1,MMO,MiM1} whose authors departed from Toda's naming conventions. These conventions were strictly followed in \cite{MT1,Oda1979}.}

	An advantage of this naming convention is
	that we can examine the commutativity of unstable compositions of elements in homotopy groups of spheres
	conveniently:
	\[
	\mathbf{x}_n\circ\mathbf{y}_i
	=\pm\mathbf{y}_n\circ\mathbf{x}_j
	\quad\text{for some }i,j
	\quad\text{if }n\geq a+b,
	\]
	where the family $\{\mathbf{x}_k\}_{k\geq r}$ has a member $\mathbf{x}_a$ and
	the family  $\{\mathbf{y}_\ell\}_{\ell\geq s }$ has a member $\mathbf{y}_b.$  In particular, 
	if
	the family $\{\mathbf{x}_k\}_{k\geq r}$ was born in $\pi_{*}(S^{a})$ 
	and the family  $\{\mathbf{y}_\ell\}_{\ell\geq s }$  was born in $\pi_{*}(S^{b})$ (namely, $r=a,\;s=b$),  then the above commutativity holds. It should be pointed out that the homotopy group containing the above $\mathbf{x}_n\circ\mathbf{y}_i$ with 
	$n= a+b$  \emph{may not be in the stable range.} For example, in $\pi_{14}(S^{5}),$   $$\lt_{5}\ca_{13}=\pm(\ca_{5}\lt_{6})=\ca_{5}\lt_{6},$$ whereas  $\pi_{14}(S^{5})$ is not 
	in the stable range.
	
	Let $\af$ and $\beta$ be homotopy classes so that $\af\hc\beta$ is defined and the domain and codomain of $\beta$ are  suspensions, and $t\in\z.$ Note that $(t\alpha)\circ\beta$ and
	$t(\alpha\beta)$ agree when $\beta$ is a suspension,  but they must not be identified for
	arbitrary unstable homotopy classes.  We write $$t\alpha\beta=:t(\alpha\beta)=\alpha\circ\beta\hc t\,\mathrm{id}.$$
	
	\noindent	Additionally, for trivial homotopy classes $\mathbf{x}$ and $\mathbf{y}$, following Toda, we write $$\mathbf{x}=\mathbf{y}=0$$  even though their domains or  codomains may not be  the same.

	Table~\ref{tab:oguchi-compositions} and Table~\ref{tab:oda-relations}  are arranged by the stem of the resulting composite, rather than by
	a single target dimension.
	
	In Table~\ref{tab:oguchi-compositions},
	the element $\kappa_7\in\pi_{21}(S^7)$ follows Mimura--Toda's revised definition in \cite{MT1}. (In fact, in 1963, with the publication of \cite{MT1}, the Toda school, including Toda himself, stopped using Toda's original $\cn_7$. It has not been used by the school since.)  The elements $\xi',\xi'', \lambda'', \lambda$ in the 18th stem follow Toda's definitions in \cite{T}.

	For clarity, in  Table~\ref{tab:oguchi-compositions},  we attach the subscript ``$\mathrm{Og}$'' to each of the four elements  $\xi',\xi'', \lambda'', \lambda$ 
	which follow \^{O}guchi's revised definitions in \cite[Lemma 2.18]{Og}.  \emph{Each of these four revised elements has the same order as the corresponding original element and
		the same $H$-image.}\footnote{The equality $H(\lambda''_{\mathrm{Og}})=H(\lambda'')$ does not follow directly from \cite{T} and \cite{Og}; its justification requires \cite[Proposition~3(3)]{Oda}.}

	These revisions due to \^{O}guchi make valuable contributions to the results of \cite{T}. Toda resolves the group extensions for $\pi_{28+i}(S^{10+i})$ $(0\leq i\leq2)$ and provides generating sets for the three groups, but does not specify generators for the individual direct summands. In \cite[Lemma 2.18]{Og}, \^{O}guchi then gives 
	\[
	\begin{aligned}
		\pi_{28}(S^{10})
		&=\mathrm{span}\n
		\xi''_{\mathrm{Og}},\,
		\xi''_{\mathrm{Og}}-\lambda''_{\mathrm{Og}},\,
		\ca_{10}\cq_{11}
		\nn
		\cong\mathbb Z/8\oplus(\mathbb Z/2)^2,\\[4pt]
		\pi_{29}(S^{11})
		&=\mathrm{span}\n
		\xi'_{\mathrm{Og}},\,
		\xi'_{\mathrm{Og}}-\lambda',\,
		\ca_{11}\cq_{12}
		\nn
		\cong\mathbb Z/8\oplus\mathbb Z/4\oplus\mathbb Z/2,\\[4pt]
		\pi_{30}(S^{12})
		&=\mathrm{span}\n
		\xi_{12},\,
		\Sigma\xi'_{\mathrm{Og}}-\Sigma\lambda',\,
		\Sigma\xi'_{\mathrm{Og}}+4\xi_{12},\,
		\ca_{12}\cq_{13}
		\nn\\
		&\cong\mathbb Z/32\oplus(\mathbb Z/4)^2\oplus\mathbb Z/2.
	\end{aligned}
	\]
	Here, the generators and their corresponding direct summands are listed in matching positions.

	Toda's original definitions of $\xi',\xi'',\lambda'',\lambda$
	are particularly well suited to constructing some new elements via
	periodicity in unstable homotopy groups. Accordingly, works such as
	\cite{Oda1979} adopt Toda's original definitions
	in preference to Oguchi's revised ones. Generating sets for the groups $\pi_{28+i}(S^{10+i})$ $(0\leq i\leq2)$, with generators specified for each direct summand using Toda's original definitions of $\xi',\xi'',\lambda'',\lambda$, were obtained by Oda in 1986 \cite[Theorem, pp.~231--232]{Oda}, more than two decades after Toda's seminal work of 1962 and \^{O}guchi's refinements published in 1964.

	It is worth pointing out the following.
	Since its first appearance, \^{O}guchi's table \cite{Og} has stood
	the test of time. The points requiring correction or clarification
	are limited to a very small number of readily apparent typographical
	errors, which we indicate in footnotes, and a notational ambiguity
	in the $18$-stem entries. Specifically, in \^{O}guchi's table \cite{Og},  the same symbols
	$\xi',\xi'', \lambda'', \lambda$  are used both for \^{O}guchi's elements
	and for Toda's original elements, without consistently distinguishing
	\^{O}guchi's $\xi',\xi'', \lambda'', \lambda$  from Toda's $\xi',\xi'', \lambda'', \lambda$ .This distinction is essential, as some relations have been established only for Toda's choices, while others have been established only for \^{O}guchi's choices. In the present version, we resolve this
	ambiguity by explicitly distinguishing the two sets of elements
	in accordance with the proofs of the $18$-stem relations in \cite{Og}.

	% Each stem is one outer cell; its parts form the smaller cells to the right.
	\definecolor{appstemshade}{RGB}{226,215,239}
	\newlength{\apprelationwidth}
	\newcommand{\AppStemRow}[3]{%
		\cellcolor{appstemshade}\(\pi^n_{n+#1}\) &
		\multicolumn{2}{@{}l@{}}{%
			\begin{tabular}[c]{>{\columncolor{tablehead!30}\centering\arraybackslash}m{13mm}|p{\apprelationwidth}}
				#2
			\end{tabular}%
		}#3\\\hline
	}
	\newcommand{\AppPartRow}[2]{%
		\textit{part~#1} &
		\raisebox{0pt}[\dimexpr\height+1pt\relax][\dimexpr\depth+1pt\relax]{%
			\(\begin{aligned}[c]#2\end{aligned}\)}\\
	}
	
	\section{A Table of Relations in Toda's Range--\^{O}guchi's  Table}\label{app:oguchi-relations}

	When a relation in \^{O}guchi's table has been improved, we give the improved version and cite its source in a footnote.

	Frequently used by the authors of \cite{Mi1,MMO,Oda1979,Ka,MiM1} and many other works, \^{O}guchi's table is a classic and convenient reference tool.

	\begingroup
	% 经典淡蓝色表头与浅灰色表格线
	\definecolor{tablehead}{RGB}{224,233,241}
	\definecolor{tablerule}{RGB}{217,217,217}
	\small
	\setlength{\tabcolsep}{5pt}
	\renewcommand{\arraystretch}{1.2}
	\setlength{\LTleft}{0pt}
	\setlength{\LTright}{0pt}
	\arrayrulecolor{tablerule}
	
	\renewcommand{\thetable}{\thesection.1}
	\begin{longtable}{>{\centering\arraybackslash}p{18mm}|p{132mm}}
		\caption{\^{O}guchi's table}\label{tab:oguchi-compositions}\\
		\hline
		\rowcolor{tablehead}\textbf{Stem} & \textbf{Composition relations}\\\hline
		\endfirsthead
		\multicolumn{2}{c}{\small\itshape \^{O}guchi's table Continued}\\
		\hline
		\rowcolor{tablehead}\textbf{Stem} & \textbf{Composition relations}\\\hline
		\endhead
		\hline
		\endfoot
		\endlastfoot
		
		$\pi_{n+3}(S^{n})$ & \(\begin{aligned}[t]&\eta_3^3=2\nu'.\end{aligned}\)\\
		& \(\begin{aligned}[t]&\Sigma^2\nu'=2\nu_5.\end{aligned}\)\\
		\hline
		
		$\pi_{n+4}(S^{n})$ & \(\begin{aligned}[t]&\eta_3 \nu_4=\nu' \eta_6.\end{aligned}\)\\
		\hline
		
		$\pi_{n+6}(S^{n})$ & \(\begin{aligned}[t]&\nu' \nu_6=0.\end{aligned}\)\\
		\hline
		
		$\pi_{n+7}(S^{n})$ & \(\begin{aligned}[t]&\Sigma\sigma'''=2\sigma'',\quad
			\Sigma\sigma''=2\sigma',\quad
			\Sigma^2\sigma'=2\sigma_9.\end{aligned}\)\\
		\hline
		
		$\pi_{n+8}(S^{n})$ & \(\begin{aligned}[t]&\eta_4 \sigma'''=\eta_5 \sigma''=\sigma''' \eta_{12}=0,\quad
			\eta_6 \sigma'=\sigma'' \eta_{13}=4\bar\nu_6.\end{aligned}\)\\
		& \(\begin{aligned}[t]&\eta_7 \sigma_8=\sigma' \eta_{14}+\varepsilon_7+\bar\nu_7,\quad
			\eta_9 \sigma_{10}=\varepsilon_9+\bar\nu_9.\end{aligned}\)\\
		& \(\begin{aligned}[t]&\eta_{10} \sigma_{11}=\sigma_{10} \eta_{17}=\varepsilon_{10}+\bar\nu_{10}.\end{aligned}\)\\
		\hline
		
		$\pi_{n+9}(S^{n})$ & \(\begin{aligned}[t]&\varepsilon_{3} \eta_{11}=\eta_{3} \varepsilon_{4}.\end{aligned}\)\\
		& \(\begin{aligned}[t]&\eta_5 \bar\nu_6=\nu_5^3,\quad
			\bar\nu_{6} \eta_{14}=\eta_{6} \bar\nu_{7}=\nu_{6}^3.\end{aligned}\)\\
		\hline
		
		$\pi_{n+10}(S^{n})$ & \(\begin{aligned}[t]&\eta_3^2 \varepsilon_5=2\varepsilon'.\end{aligned}\)\\
		& \(\begin{aligned}[t]&2(\nu_5 \sigma_8)=\nu_5 (\Sigma\sigma')=\pm \Sigma^2\varepsilon'.\end{aligned}\)\\
		& \(\begin{aligned}[t]&\nu' \sigma''=0,\quad
			\Sigma\nu'\circ \sigma'=2\Sigma\varepsilon',\quad
			\sigma''' \nu_{12}=4(\nu_5 \sigma_8).\end{aligned}\)\\
		& \(\begin{aligned}[t]&\sigma'' \nu_{13}=\pm2\nu_6 \sigma_9,\quad
			\sigma' \nu_{14}=x\nu_7 \sigma_{10}\quad\text{for some odd integer }x.\end{aligned}\)\\
		& \(\begin{aligned}[t]&\nu_9 \sigma_{12}=\pm2\sigma_9\nu_{16}.\end{aligned}\)\\
		& \(\begin{aligned}[t]&\mu_{3} \eta_{12}=\eta_{3} \mu_{4}.\end{aligned}\)\\
		\hline
		
		$\pi_{n+11}(S^{n})$ & \(\begin{aligned}[t]&\eta_3^2 \mu_5=2\mu',\quad
			\Sigma^2\mu'=\pm2\zeta_5.\end{aligned}\)\\
		& \(\begin{aligned}[t]&\nu' \bar\nu_6=\varepsilon_3 \nu_{11},\quad
			\nu_6 \bar\nu_9=\nu_6 \varepsilon_9=2\bar\nu_6 \nu_{14}.\end{aligned}\)\\
		& \(\begin{aligned}[t]&\varepsilon' \eta_{13}=\nu' \varepsilon_6.\end{aligned}\)\\
		& \(\begin{aligned}[t]&\nu_5 \sigma_8 \eta_{15}=\nu_5 \varepsilon_8.\mbox{\footnotemark}\end{aligned}\)\footnotetext{The equation $\nu_5 \sigma_8 \eta_{15}=\nu_5 \varepsilon_8$ is given in \cite[p. 152]{T}. Yet, \cite{Og} neglected it and gave a weaker relation.}\\
		\hline
		
		$\pi_{n+12}(S^{n})$ & \(\begin{aligned}[t]&\eta_4 \zeta_5\equiv (\Sigma\nu') \mu_7\mathbin{\mathrm{mod}} (\Sigma\nu') \eta_7 \varepsilon_8,\quad
			\eta_{5} \zeta_{6}=0.\end{aligned}\)\\
		& \(\begin{aligned}[t]&\zeta_6 \eta_{17}=8[\iota_6,\iota_6] \sigma_{11},\quad
			\zeta_{7} \eta_{18}=0.\end{aligned}\)\\
		& \(\begin{aligned}[t]&\mu' \eta_{14}=\nu' \mu_6.\end{aligned}\)\\
		& \(\begin{aligned}[t]&\mu_3 \nu_{12}=\nu' \eta_6 \varepsilon_7,\quad
			\nu_6 \mu_9=8[\iota_6,\iota_6] \sigma_{11}.\end{aligned}\)\\
		\hline
		
		$\pi_{n+13}(S^{n})$ & \(\begin{aligned}[t]&\varepsilon' \nu_{13}=0.\end{aligned}\)\\
		& \(\begin{aligned}[t]&\eta_{10} \theta'=0,\quad
			\eta_{11} \theta=\sigma_{11} \nu_{18}^2+\theta' \eta_{23}.\end{aligned}\)\\
		\hline
		
		$\pi_{n+14}(S^{n})$ & \(\begin{aligned}[t]&2\kappa_7=\bar\nu_7 \nu_{15}^2.\end{aligned}\)\\
		& \(\begin{aligned}[t]&\theta' \eta_{23}^2=0,\quad
			\theta \eta_{24}^2=2[\iota_{12},\nu_{12}]+8\sigma_{12}^2\;(\text{this has been improved by using}\;\s).\end{aligned}\)\\
		& \(\begin{aligned}[t]&\nu' \zeta_6=\mu' \nu_{14}=0,\quad
			\sigma''' \sigma_{12}=\nu_5 \zeta_8.\end{aligned}\)\\
		& \(\begin{aligned}[t]&\zeta_{5} \nu_{16}\ty\nu_{5} \zeta_{8}\m \nu_{5}\ch_{8}\nu_{16};\quad
			\nu_6 \zeta_9=\ck_{6}\nu_{17}=2\sigma'' \sigma_{13}.\mbox{\footnotemark}\end{aligned}\)\footnotetext{The equation $\nu_6 \zeta_9=\ck_{6}\nu_{17}=2\sigma'' \sigma_{13}$ is the improved result given in \cite[Lemma 2.3]{KM2}; N. Oda taught the authors of \cite{KM2} how to make this improvement.}\\
		\hline
		
		$\pi_{n+15}(S^{n})$ & \(\begin{aligned}[t]&\Sigma\rho^{\mathrm{IV}}=2\rho''',\quad
			\Sigma\rho'''=2\rho'',\quad
			\Sigma^2\rho''=2\rho',\quad
			\Sigma^4\rho'=2\rho_{13}.\end{aligned}\)\\
		& \(\begin{aligned}[t]&\eta_{6} \kappa_{7}=\bar\varepsilon_{6},\quad
			\kappa_7 \eta_{21}=\sigma' \bar\nu_{14}+\bar\varepsilon_7.\end{aligned}\)\\
		& \(\begin{aligned}[t]&\kappa_{9} \eta_{23}=\bar\varepsilon_{9},\quad
			\nu_4^5=0.\end{aligned}\)\\
		& \(\begin{aligned}[t]&\nu_8 \theta'\equiv (\Sigma\sigma') \varepsilon_{15} \mod{(\Sigma\sigma') \eta_{15}\sigma_{16}}.\end{aligned}\)\\
		& \(\begin{aligned}[t]&\nu_9 \theta=\sigma_9 \varepsilon_{16}\;\text{or}\;\sigma_9 \bar\nu_{16}\uprightsemicolon{}\quad
			\theta' \nu_{23}=\theta \nu_{24}=0.\end{aligned}\)\\
		& \(\begin{aligned}[t]&\sigma''' \varepsilon_{12}=\sigma''' \bar\nu_{12}=\sigma'' \varepsilon_{13}
			=\sigma'' \bar\nu_{13}=[\iota_6,\iota_6] \sigma_{11} \nu_{18}=0.\end{aligned}\)\\
		& \(\begin{aligned}[t]&\varepsilon_{3} \sigma_{11}=0,\quad
			\bar\nu_{6} \sigma_{14}=0.\end{aligned}\)\\
		& \(\begin{aligned}[t]&\sigma_{10} \varepsilon_{17}=\sigma_{10} \bar\nu_{17}=[\iota_{10},\nu_{10}^2].\end{aligned}\)\\
		\hline
		
		$\pi_{n+16}(S^{n})$ & \(\begin{aligned}[t]&\eta_4 \rho^{\mathrm{IV}}=\rho^{\mathrm{IV}} \eta_{20}=\eta_5 \rho'''=0,\quad
			\rho''' \eta_{21}=\eta_6 \rho''=4\zeta'.\end{aligned}\)\\
		& \(\begin{aligned}[t]&\rho'' \eta_{22}=\sigma' \mu_{14}.\end{aligned}\)\\
		& \(\begin{aligned}[t]&\rho' \eta_{24}\equiv\sigma_9 \mu_{16}+\mu_9 \sigma_{18}
			\mathbin{\mathrm{mod}} \sigma_9 \varepsilon_{16} \eta_{24},\sigma_9 \nu_{16}^3.\end{aligned}\)\\
		& \(\begin{aligned}[t]&\eta_{12} \rho_{13}=\sigma_{12} \mu_{19},\;(\s\rho')\ca_{25}=8[\e_{10},\cg_{10}].\end{aligned}\)\\
		& \(\begin{aligned}[t]&\eta_{13} \rho_{14}=\rho_{13} \eta_{28}=\mu_{13} \sigma_{22}
			=\sigma_{13} \mu_{20}.\end{aligned}\)\\
		& \(\begin{aligned}[t]&\sigma''' \mu_{12}=0,\quad
			\sigma'' \mu_{13}=4\zeta',\quad
			\mu_{10} \sigma_{19}=\sigma_{10} \mu_{17}+8[\iota_{10},\sigma_{10}].\end{aligned}\)\\
		& \(\begin{aligned}[t]&\sigma' \eta_{14} \varepsilon_{15}=\Sigma\zeta'.\end{aligned}\)\\
		& \(\begin{aligned}[t]&\nu_{5} \sigma_{8} \nu_{15}^2=\eta_{5} \bar\varepsilon_{6}.\end{aligned}\)\\
		& \(\begin{aligned}[t]&\sigma_{10} \nu_{17}^3=\sigma_{10} \eta_{17} \varepsilon_{18}=0.\end{aligned}\)\\
		& \(\begin{aligned}[t]&\varepsilon_{3}^2=\varepsilon_{3} \bar\nu_{11}=\eta_{3} \bar\varepsilon_{4}
			=\bar\varepsilon_{3} \eta_{18}.\end{aligned}\)\\
		& \(\begin{aligned}[t]&\bar\nu_{6}^2=\bar\nu_{6} \varepsilon_{14}=0.\end{aligned}\)\\
		& \(\begin{aligned}[t]&\eta_{18}^*\equiv\omega_{18}\mathbin{\mathrm{mod}}\sigma_{18} \mu_{25}.\end{aligned}\)\\
		\hline
		
		$\pi_{n+17}(S^{n})$ & \(\begin{aligned}[t]&\Sigma\bar\varepsilon'=(\Sigma\nu') \kappa_7,\quad
			\eta_3^2 \bar\varepsilon_5=2\bar\varepsilon'.\end{aligned}\)\\
		& \(\begin{aligned}[t]&\eta_{13} \omega_{14}=\varepsilon_{13}^*,\quad
			\eta_{14} \omega_{15}=\omega_{14} \eta_{30}=\varepsilon_{14}^*.\end{aligned}\)\\
		& \(\begin{aligned}[t]&\eta_{18}^* \eta_{34}\equiv\varepsilon_{18}^*\mathbin{\mathrm{mod}}\sigma_{18} \eta_{25} \mu_{26}.\end{aligned}\)\\
		& \(\begin{aligned}[t]&(\Sigma^2\eta^{*\prime}) \eta_{33}\equiv0\mathbin{\mathrm{mod}}\sigma_{17} \eta_{24} \mu_{25}.\end{aligned}\)\\
		& \(\begin{aligned}[t]&\eta_{14} \eta^{*\prime}\equiv0\mathbin{\mathrm{mod}}\sigma_{14} \eta_{21} \mu_{22}.\end{aligned}\)\\
		& \(\begin{aligned}[t]&\eta_{15} \eta_{16}^*\equiv\eta^{*\prime} \eta_{31}+\varepsilon_{15}^*
			\mathbin{\mathrm{mod}}\sigma_{15} \eta_{22} \mu_{23}.\end{aligned}\)\\
		& \(\begin{aligned}[t]&\zeta' \eta_{22}=\eta_5 \zeta'=0.\end{aligned}\)\\
		& \(\begin{aligned}[t]&\kappa_{7} \nu_{21}=\nu_{7} \kappa_{10}.\end{aligned}\)\\
		& \(\begin{aligned}[t]&\nu' \bar\nu_6 \nu_{14}^2=\varepsilon' \sigma_{13}
			=\varepsilon_3 \nu_{11}^3=2\bar\varepsilon'.\end{aligned}\)\\
		& \(\begin{aligned}[t]&\nu_5 \bar\nu_8 \nu_{16}^2=2(\nu_5 \kappa_8).\end{aligned}\)\\
		& \(\begin{aligned}[t]&\nu_5 (\Sigma\sigma') \sigma_{15}=0,\quad
			\nu_5 \sigma_8^2=2\nu_5 \kappa_8\quad\text{or}\quad0.\end{aligned}\)\\
		& \(\begin{aligned}[t]&\mu_3 \varepsilon_{12}\equiv\varepsilon_3 \mu_{11}\mathbin{\mathrm{mod}}2\bar\varepsilon';\quad
			\mu_3 \varepsilon_{12}\equiv\eta_3 \mu_4 \sigma_{13}\mathbin{\mathrm{mod}}2\bar\varepsilon'.\end{aligned}\)\\
		& \(\begin{aligned}[t]&\mu_3 \bar\nu_{12}\equiv0\mathbin{\mathrm{mod}}2\bar\varepsilon',\quad
			\mu_{6} \bar\nu_{15}=\bar\nu_{6} \mu_{14}=0.\end{aligned}\)\\
		& \(\begin{aligned}[t]&[\iota_6,\iota_6] \theta'=0.\end{aligned}\)\\
		\hline
		
		$\pi_{n+18}(S^{n})$ & \(\begin{aligned}[t]&\Sigma\xi''_{}=2\xi'_{},\quad
			\Sigma^2\xi'_{}=2\xi_{13},\quad
			\Sigma\lambda''_{}=2\lambda',\quad
			\Sigma^2\lambda'=2\lambda_{}.\end{aligned}\)\\
		& \(\begin{aligned}[t]&\Sigma^4\lambda_{}=2\nu_{17}^*,\quad4(\nu_{18}^{*}+\xi_{18})
			=2(\nu_{19}^{*}+\xi_{19})
			=\nu_{20}^{*}+\xi_{20}
			=0
			.\end{aligned}\)\hspace{-1em}\footnote{The relations $4(\nu_{18}^{*}+\xi_{18})
			=2(\nu_{19}^{*}+\xi_{19})
			=\nu_{20}^{*}+\xi_{20}
			=0$ are given in \cite[Corollary 12.25, p.~169]{T}. The corresponding entry in \^{O}guchi's table \cite{Og} contains a typo.}\\
		& \(\begin{aligned}[t]&2\xi''_{\mathrm{Og}}=2\lambda''_{\mathrm{Og}}=\pm\sigma_{10} \zeta_{17}.\end{aligned}\)\\
		& \(\begin{aligned}[t]&2\xi''_{}\ty2\lambda''_{}\ty\sigma_{10} \zeta_{17}\m2\sigma_{10} \zeta_{17};\quad
			4\xi'_{}=4\lambda'=\sigma_{11} \zeta_{18}.\end{aligned}\)\hspace{-1em}\footnote{The relations satisfied by Toda's original elements
			$\xi',\xi'',\lambda',\lambda''$ are recorded in \cite{Oda}.
		}\\
		& \(\begin{aligned}[t]&4\Sigma\xi'_{}=16\xi_{12}=\sigma_{12} \zeta_{19}.\end{aligned}\)\\
		& \(\begin{aligned}[t]&\eta_{3} \bar\mu_{4}=\bar\mu_{3} \eta_{20}.\end{aligned}\)\\
		& \(\begin{aligned}[t]&\bar\varepsilon' \eta_{20}=\bar\varepsilon_3 \nu_{18}=\nu' \bar\varepsilon_6.\end{aligned}\)\\
		& \(\begin{aligned}[t]&\eta^{*\prime} \eta_{31}^2=4\Sigma^2\lambda_{};\quad
			\eta_{16}^* \eta_{32}^2\equiv4\nu_{16}^*\mathbin{\mathrm{mod}}4\Sigma^3\lambda_{};\quad
			\eta_{17}^* \eta_{33}^2=4\nu_{17}^*.\end{aligned}\)\\
		& \(\begin{aligned}[t]&\varepsilon_{13}^* \eta_{30}=\eta_{13} \varepsilon_{14}^*
			\equiv4\xi_{13}\mathbin{\mathrm{mod}}4\lambda_{}.\end{aligned}\)\\
		& \(\begin{aligned}[t]&\nu' \rho'''=(\Sigma\nu') \rho''=0,\quad
			\nu_5 \Sigma\rho''=0.\end{aligned}\)\\
		& \(\begin{aligned}[t]&\nu_6 \rho'=0,\quad
			\nu_{10} \rho_{13}\equiv0\mathbin{\mathrm{mod}}2\sigma_{10} \zeta_{17}.\mbox{\footnotemark}\end{aligned}\)\footnotetext{The equations $\nu_5\s \rho''=0$ and $\nu_6 \rho'=0$ are the improved result given in \cite{Os}.}\\
		& \(\begin{aligned}[t]&\rho^{\mathrm{IV}} \nu_{20}\equiv4\zeta_5 \sigma_{16}\mathbin{\mathrm{mod}}\nu_5 \bar\varepsilon_8;\quad
			\rho''' \nu_{21}=\pm2\zeta_6 \sigma_{17}.\end{aligned}\)\\
		& \(\begin{aligned}[t]&\rho'' \nu_{22}=\sigma' \zeta_{14}=x\zeta_7 \sigma_{18}
			\quad\text{for some odd integer }x.\end{aligned}\)\\
		& \(\begin{aligned}[t]&\rho' \nu_{24}\equiv-\sigma_9 \zeta_{16}\mathbin{\mathrm{mod}}4\sigma_9 \zeta_{16},\quad
			\rho_{13} \nu_{28}=4\lambda.\end{aligned}\)\\
		& \(\begin{aligned}[t]&\nu_5 \sigma_8 \varepsilon_{15}=\nu_5 \sigma_8 \bar\nu_{15}
			=\nu_5 \bar\varepsilon_8.\end{aligned}\)\hspace{-1em}\footnote{The relation $\nu_5 \sigma_8 \varepsilon_{15}=\nu_5 \sigma_8 \bar\nu_{15}
			=\nu_5 \bar\varepsilon_8$ is the improved result given in \cite{Oda2}.}\\
		& \(\begin{aligned}[t]&\varepsilon' \varepsilon_{13}=\varepsilon' \bar\nu_{13}=\nu' \bar\varepsilon_6.\end{aligned}\)\\
		& \(\begin{aligned}[t]&4\sigma_{10} \zeta_{17}=2\sigma_{11} \zeta_{18}=\sigma_{13} \zeta_{20}=0.\end{aligned}\)\\
		& \(\begin{aligned}[t]&\sigma''' \zeta_{12}=4\zeta_5 \sigma_{16},\quad
			\sigma'' \zeta_{13}=\pm2\zeta_6 \sigma_{17}.\end{aligned}\)\\
		& \(\begin{aligned}[t]&\mu_3^2\equiv\eta_3 \bar\mu_4\mathbin{\mathrm{mod}}2\mu' \sigma_{14}.\end{aligned}\)\\
		\hline
		
		$\pi_{n+19}(S^{n})$ &\(
		\begin{aligned}[t]
			&\;\;\text{In the following part of this block, each formula containing }
			\xi',\xi'',\lambda''\text{ or }\lambda\text{ remains true}
		\end{aligned}
		\)\\
		&\(
		\begin{aligned}[t]
			&\text{when these elements are replaced by }
			\xi'_{\mathrm{Og}},\xi''_{\mathrm{Og}},\lambda''_{\mathrm{Og}}
			\text{ and }\lambda_{\mathrm{Og}}\text{, respectively.}
		\end{aligned}
		\)\\& \(\begin{aligned}[t]&\eta_3^2 \bar\mu_5=2\bar\mu',\quad
			\Sigma^2\bar\mu'=2\bar\zeta_5.\end{aligned}\)\\
		& \(\begin{aligned}[t]&\Sigma^2\omega'=2\omega_{14} \nu_{30},\quad
			\nu_6 \mu_9 \sigma_{18}=16\bar\sigma_6.\end{aligned}\)\\
		& \(\begin{aligned}[t]&\xi'' \eta_{28}=\lambda''_{} \eta_{28}=\eta_{10} \xi'_{}=\eta_{10} \lambda'=0.\end{aligned}\)\\
		& \(\begin{aligned}[t]&\eta_9 \lambda''_{}=\eta_9 \xi''_{}=0;\quad
			\eta_{11} \xi_{12}=\xi'_{} \eta_{29}.\end{aligned}\)\\
		& \(\begin{aligned}[t]&\nu_{16}^* \eta_{34}=0,\quad
			\eta_{17} \nu_{18}^*=\omega_{17} \nu_{33},\quad
			\lambda_{} \eta_{31}\equiv \Sigma\omega'\mathbin{\mathrm{mod}}\xi_{13} \eta_{31}.\end{aligned}\)\\
		& \(\begin{aligned}[t]&\eta_{12} \lambda_{}\equiv0\mathbin{\mathrm{mod}}(\Sigma\xi'_{}) \eta_{30},\;(\Sigma\lambda' )\eta_{30}.\end{aligned}\)\\
		& \(\begin{aligned}[t]&P(\Sigma\theta) \eta_{23}^2\equiv0\mathbin{\mathrm{mod}}16\bar\sigma_6.\end{aligned}\)\\
		& \(\begin{aligned}[t]&\nu' \zeta'=\nu' \eta_6 \bar\varepsilon_7=\varepsilon_3 \zeta_{11}=0.\end{aligned}\)\\
		& \(\begin{aligned}[t]&\zeta' \nu_{22}=\pm8\bar\sigma_6,\quad
			\nu_{11} \omega_{14}\equiv\lambda' \eta_{29}\mathbin{\mathrm{mod}}\xi'_{} \eta_{29}.\end{aligned}\)\\
		& \(\begin{aligned}[t]&\eta^{*\prime} \nu_{31}=\eta_{16}^* \nu_{32}=0.\end{aligned}\)\\
		& \(\begin{aligned}[t]&\nu_{12} \eta^{*\prime}\equiv0\mathbin{\mathrm{mod}} (\Sigma\xi'_{} )\eta_{30}, \Sigma\lambda' \eta_{30}.\end{aligned}\)\\
		& \(\begin{aligned}[t]&\nu_{13} \eta_{16}^*\equiv \Sigma\omega'\mathbin{\mathrm{mod}}\xi_{13} \eta_{31}.\end{aligned}\)\\
		& \(\begin{aligned}[t]&\sigma''' \Sigma\theta'=\sigma'' \s\theta=0,\quad
			\theta' \sigma_{23}=\xi'_{} \eta_{29}.\mbox{\footnotemark}\end{aligned}\)\footnotetext{The paper \cite{Og} has typos in the relation $\sigma''' \Sigma\theta'=\sigma'' \s\theta=0.$  We provide the revised version based on the proof in \cite{Og}.}\\
		& \(\begin{aligned}[t]&\theta \sigma_{24}\equiv\xi_{12} \eta_{30}
			\mathbin{\mathrm{mod}} (\Sigma\xi'_{}) \eta_{30},\;(\Sigma\lambda' )\eta_{30}.\end{aligned}\)\\
		& \(\begin{aligned}[t]&\bar\nu_6 \zeta_{14}=\pm8\bar\sigma_6.\end{aligned}\)\\
		& \(\begin{aligned}[t]&\zeta_5 \eta_{16} \sigma_{17}=\nu_5 \mu_8 \sigma_{17}.\end{aligned}\)\\
		& \(\begin{aligned}[t]&\varepsilon' \mu_{13}=\nu' \mu_6 \sigma_{15}.\end{aligned}\)\\
		& \(\begin{aligned}[t]&[\iota_6,\iota_6] \sigma_{11}^2\equiv2\bar\sigma_6\mathbin{\mathrm{mod}}4\bar\sigma_6.\end{aligned}\)\\
		\hline
	\end{longtable}
	\endgroup

	\section{A Table of Relations Beyond Toda's Range}\label{app:oda-relations}
	Table~\ref{tab:oda-relations} collects and arranges by stem the relations
	from thirty-one items of results: sixteen from \cite{Oda1979} and fifteen
	from \cite{MiM1}.\footnote{In stem 20, the relation $\ca_{6}\cs_{7}=\ca_{6}\ccs_{7}=0$ is not from \cite{Oda1979} or \cite{MiM1}. It is easy and given in our Lemma~\ref{m3:2ltp}(3).} The sixteen items from \cite{Oda1979} are listed below.
	Their references use Oda's printed page numbers and distinguish the three
	parts of the paper, in which proposition numbers begin afresh.
	\begin{itemize}
		\item Part~I: Proposition~3.1 (pp.~54--55), Proposition~3.5 (p.~60),
		Proposition~5.1 (p.~64), and formulas~(8.13)--(8.14) (p.~74).
		\item Part~II: Proposition~2.1 (pp.~81--82), Proposition~2.2 (p.~84),
		Proposition~2.3 (p.~86), Proposition~2.4 (p.~87), and
		formulas~(6.3)--(6.4) (pp.~98--99).
		\item Part~III: Proposition~2.1 (p.~107), Proposition~2.2 (p.~108),
		Proposition~2.3 (p.~109), Proposition~2.4 (p.~111),
		Proposition~2.5 (pp.~112--113), Proposition~2.6 (p.~114), and
		Propositions~3.1--3.2 (p.~116), taken together as the final group.
	\end{itemize}The following fifteen items are taken from \cite{MiM1}. Page numbers for
	this source refer to its 69-page PDF, numbered~1--69, rather than to the
	journal pagination~319--387. (In fact, the published paper \cite{MiM1} does not have printed page numbers.) Each numbered item below is one item.
	\begin{enumerate}[itemsep=1pt,parsep=0pt]
		\item Formula~(3.20) (p.~22).
		\item Lemma~3.10 (p.~23).
		\item Formula~(3.23) (p.~23).
		\item Lemma~3.11 (p.~24).
		\item Formulas~(3.36) and~(3.38) (p.~26).
		\item Formulas~(3.39) and~(3.40) (p.~27).
		\item Formula~(3.41) and Lemma~3.12 (p.~28).
		\item Formula~(3.44) (p.~30).
		\item Formula~(3.45) and Lemma~3.14 (p.~32).
		\item Lemma~5.3 (p.~46).
		\item Lemma~5.6 (p.~51).
		\item Lemma~5.7 (p.~53).
		\item Lemma~6.1 (p.~57).
		\item Lemma~6.2 (p.~59).
		\item Lemma~7.3 (p.~66).
	\end{enumerate}
	Only relations in stems $k\geq20$ are included: an unstable term in stem~$k$
	lies in $\pi_{n+k}(S^{n}),$   and a stable term lies in $\pi_k^S(S^0).$ 
	The entries span stems~20--32. Whitehead products are placed according to
	the dimension of the resulting class. When one source group contains
	relations in different stems, its entries appear in the corresponding
	stem blocks. Equations, congruences and alternatives are retained with
	their original qualifications.

	The generators and their dimensions follow Appendix~\ref{tlqdzj}, including
	its dictionary of Oda's notation. In particular, Oda's symbols
	$A_1,C_1,B_1,D_1,F_1$ and their decorated variants are retained whenever
	replacing them by a more familiar generator would turn a congruence into
	an unwarranted equality. The epsilon generator in stem~24 is written
	$\tilde\varepsilon_n,$   as in Subsection~\ref{app:stem-24}; it is distinct
	from $\bar{\varepsilon}_n$ in stem~15.
	Each coefficient \textbf{$x$ or $y$ is odd and is local to its displayed relation}.
	A list after $\bmod$ denotes the subgroup generated by the listed terms;
	a sum or difference within one listed term remains a single generator.

	In this table,
	the element $\kappa_7$ follows Mimura--Toda's revised definition in \cite{MT1}\uprightsemicolon{} the elements $\xi',\xi'', \lambda'', \lambda$ in the 18th stem follow Toda's original definition in \cite{T}, not \^{O}guchi's revised definition.

	In the row marked ``stably'', dimension subscripts are suppressed as in
	Appendix~\ref{tlqdzj}; the subscripts in $\mu_3$ and $\rho_3$ are retained
	as part of the generator names, with $\mu_3\in\pi_{25}^S(S^0)$ and
	$\rho_3\in\pi_{31}^S(S^0).$

	The document linked below provides precise references for each relation,
	including page numbers and the relevant propositions, lemmas, or formulas,
	allowing readers to locate the proofs easily:
	\href{https://drive.google.com/file/d/1pRnBoFZZ-mWxzT8RADWRJ9GCLGFoiEmh/view?usp=drive_link}{CLICK TO DOWNLOAD.}

	The elementary relations beyond Toda's range are collected in the following table.
	\begingroup
	% Parts remain together; a stem may continue on the following page.
	% Formula cells have explicit math delimiters and do not use nested tables.
	\small\setlength{\tabcolsep}{3.5pt}\renewcommand{\arraystretch}{1.2}
	\setlength{\extrarowheight}{2pt}\setlength{\jot}{1.2pt}
	\renewcommand{\thetable}{\thesection.1}
	% Add equal clearance above and below each formula's actual box.
	\newcommand{\OdaRelationCell}[1]{%
		\begingroup
		\setbox0=\hbox{#1}%
		\rule[-\dimexpr\dp0+3pt\relax]{0pt}{\dimexpr\ht0+\dp0+6pt\relax}%
		\usebox0%
		\endgroup
	}
	\setlength{\parskip}{0pt}
	\setlength{\LTleft}{0pt}\setlength{\LTright}{0pt}
	\setlength{\LTpre}{6pt}\setlength{\LTpost}{6pt}
	\setlength{\LTcapwidth}{\linewidth}
	\definecolor{odapurple}{RGB}{232,224,242}
	\definecolor{odapurplepart}{RGB}{232,224,242}
	\definecolor{odapurplebody}{RGB}{255,255,255}
	\definecolor{odagreen}{RGB}{218,232,224}
	\definecolor{odagreenpart}{RGB}{218,232,224}
	\definecolor{odagreenbody}{RGB}{255,255,255}
	\definecolor{odahead}{RGB}{224,233,241}
	\definecolor{odarule}{RGB}{207,207,215}
	\arrayrulecolor{odarule}
	% The page references repeat a stem label only at a genuine page break.
	% As with ordinary cross-references, compile until references have settled.
	\providecommand{\OdaPageStem}[2]{}
	\renewcommand{\OdaPageStem}[2]{%
		\ifnum#2=1\relax
		{\color{black}\bfseries #1}%
		\else
		\ifnum\getpagerefnumber{oda:part:#1:#2}=\getpagerefnumber{oda:part:#1:\number\numexpr#2-1\relax}\relax
		\else
		{\color{black}\bfseries #1}\par\textup{\scriptsize(cont.)}%
		\fi
		\fi
	}
	\begin{longtable}{>{\centering\arraybackslash}m{16mm}|>{\centering\arraybackslash}m{13mm}|m{\dimexpr\linewidth-29mm-6\tabcolsep-2\arrayrulewidth\relax}}
		\caption{Relations beyond Toda's range}\label{tab:oda-relations}\\
		\hline
		\rowcolor{odahead}\textbf{Stem} & \textbf{Part} & \cellcolor{white}\textbf{Composition relations}\\\hline
		\endfirsthead
		\multicolumn{3}{c}{\small\itshape Relations beyond Toda's range continued}\\
		\hline
		\rowcolor{odahead}\textbf{Stem} & \textbf{Part} & \cellcolor{white}\textbf{Composition relations}\\\hline
		\endhead
		\hline
		\endfoot
		\hline
		\endlastfoot
		
		% Stem 20, part 1: 4 rows.
		\cellcolor{odapurple}\leavevmode\color{black}\OdaPageStem{20}{1}
		& \cellcolor{odapurplepart}\textup{part~1}
		& \cellcolor{odapurplebody}\leavevmode\label{oda:part:20:1}%
		\OdaRelationCell{\(\begin{aligned}
				% Oda1979: I.3.1(3), p.55
				&\varepsilon_{12}^{*}\nu_{29}=\theta\bar\nu_{24}=A_1\bar\nu_{24},\quad \theta'\bar\nu_{23}=0.\\
				% Oda1979: II.2.1(7), p.82
				&C_1\varepsilon_{23}=C_1\eta_{23}\sigma_{24}=\xi'\eta_{29}^{2}=\Sigma F_1.\\
				% Oda1979: II.2.2(1)--(2), p.84
				&\bar\zeta_5\eta_{24}\equiv\zeta_5\mu_{16}\equiv\nu_5\bar\mu_8\m \nu_5\eta_8\mu_9\sigma_{18}.\\
				% Oda1979: II.2.3(1), p.86
				&\bar\mu_3\nu_{20}=\mu_3\zeta_{12}=\nu'\eta_6\mu_7\sigma_{16};\;\quad \ca_{6}\cs_{7}=\ca_{6}\ccs_{7}=0.     %   Lemma~\ref{m3:2ltp}(3) of this book.
			\end{aligned}\)}\\
		\hline
		% Stem 21, part 1: 4 rows.
		\cellcolor{odagreen}\leavevmode\color{black}\OdaPageStem{21}{1}
		& \cellcolor{odagreenpart}\textup{part~1}
		& \cellcolor{odagreenbody}\leavevmode\label{oda:part:21:1}%
		\OdaRelationCell{\(\begin{aligned}
				% Oda1979: I.3.1(1), p.54
				&\xi'\nu_{29}=0.\\
				% Oda1979: I.3.1(1), p.54
				&\nu_{10}\lambda=\sigma_{10}\kappa_{17},\quad \lambda'\nu_{29}=\sigma_{11}\kappa_{18}.\\
				% Oda1979: I.3.1(7), p.55
				&\sigma'''\sigma_{12}^{2}=\sigma''\sigma_{13}^{2}=\sigma'\sigma_{14}^{2}=0.\\
				% Oda1979: II.2.1(2), p.81
				&\kappa_7\sigma_{21}=0,\quad \nu_9\xi_{12}=\sigma_9^3,\quad \xi_{12}\nu_{30}=\sigma_{12}^3.
			\end{aligned}\)}\\
		\hhline{>{\arrayrulecolor{odagreen}}->{\arrayrulecolor{odarule}}--}
		% Stem 21, part 2: 4 rows.
		\cellcolor{odagreen}\leavevmode\color{black}\OdaPageStem{21}{2}
		& \cellcolor{odagreenpart}\textup{part~2}
		& \cellcolor{odagreenbody}\leavevmode\label{oda:part:21:2}%
		\OdaRelationCell{\(\begin{aligned}
				% Oda1979: II.2.3(3), p.86
				&\nu_{13}\nu_{16}^{*}=\pm\lambda\nu_{31}+\sigma_{13}^3.\\
				% Oda1979: II.2.4(3), p.87
				&\nu_5\sigma_8\zeta_{15}=0.\\
				% Oda1979: III.2.1(1), p.107
				&C_1\mu_{23}=\beta''\eta_{31}=C_1^{(2)}\eta_{31},\quad C_1\eta_{23}\varepsilon_{24}=0.\\
				% Oda1979: III.2.4(2), p.111
				&\bar\kappa_7\eta_{27}=\sigma'\kappa_{14}+\eta_7\bar\kappa_8.
			\end{aligned}\)}\\
		\hline
		% Stem 22, part 1: 4 rows.
		\cellcolor{odapurple}\leavevmode\color{black}\OdaPageStem{22}{1}
		& \cellcolor{odapurplepart}\textup{part~1}
		& \cellcolor{odapurplebody}\leavevmode\label{oda:part:22:1}%
		\OdaRelationCell{\(\begin{aligned}
				% Oda1979: I.3.1(1), p.54
				&\nu_{10}\lambda\eta_{31}=\nu_{10}^{2}\eta_{16}^{*}=0.\\
				% Oda1979: I.3.1(2), p.54
				&\sigma'\eta_{14}\kappa_{15}=\bar\nu_7\kappa_{15}+\varepsilon_7\kappa_{15}.\\
				% Oda1979: I.3.1(2),(6), pp.54--55
				&\alpha\eta_{26}=\varepsilon_5\kappa_{13}=\eta_5^2\bar\kappa_7.\\
				% Oda1979: I.3.1(6), p.55
				&\eta_6\bar\kappa_7\eta_{27}=\eta_6^2\bar\kappa_8.
			\end{aligned}\)}\\
		\hhline{>{\arrayrulecolor{odapurple}}->{\arrayrulecolor{odarule}}--}
		% Stem 22, part 2: 4 rows.
		\cellcolor{odapurple}\leavevmode\color{black}\OdaPageStem{22}{2}
		& \cellcolor{odapurplepart}\textup{part~2}
		& \cellcolor{odapurplebody}\leavevmode\label{oda:part:22:2}%
		\OdaRelationCell{\(\begin{aligned}
				% Oda1979: I.3.5(1), p.60
				&\bar\sigma_6\nu_{25}=\nu_6\bar\sigma_9.\\
				% Oda1979: II.2.1(3), p.81
				&\bar\zeta_5\nu_{24}=\nu_5\bar\zeta_8=\zeta_5^2.\\
				% Oda1979: II.2.1(4), pp.81--82
				&\sigma''\rho_{13}=\pm\rho'''\sigma_{21}.\\
				% Oda1979: II.2.1(4), p.81
				&\sigma'\rho_{14}\equiv x\rho''\sigma_{22}\m \bar\nu_7\kappa_{15}+\varepsilon_7\kappa_{15}\quad(x\text{: odd}).
			\end{aligned}\)}\\
		\hhline{>{\arrayrulecolor{odapurple}}->{\arrayrulecolor{odarule}}--}
		% Stem 22, part 3: 4 rows.
		\cellcolor{odapurple}\leavevmode\color{black}\OdaPageStem{22}{3}
		& \cellcolor{odapurplepart}\textup{part~3}
		& \cellcolor{odapurplebody}\leavevmode\label{oda:part:22:3}%
		\OdaRelationCell{\(\begin{aligned}
				% Oda1979: II.2.1(4), p.82
				&\sigma_9\rho_{16}\equiv x\rho'\sigma_{24}\m \sigma_9\bar\varepsilon_{16}\quad(x\text{: odd}).\\
				% Oda1979: II.2.1(4), p.82
				&(\Sigma\rho')\sigma_{25}=t\sigma_{10}\rho_{17}\quad(t\text{: odd}).\\
				% Oda1979: II.2.1(4), p.82
				&\sigma_{13}\rho_{20}\equiv2\rho_{13}\sigma_{28}\m 4\rho_{13}\sigma_{28}.\\
				% Oda1979: II.2.1(4), p.82
				&\rho_{17}\sigma_{32}=8\sigma_{17}^{*}.
			\end{aligned}\)}\\
		\hhline{>{\arrayrulecolor{odapurple}}->{\arrayrulecolor{odarule}}--}
		% Stem 22, part 4: 4 rows.
		\cellcolor{odapurple}\leavevmode\color{black}\OdaPageStem{22}{4}
		& \cellcolor{odapurplepart}\textup{part~4}
		& \cellcolor{odapurplebody}\leavevmode\label{oda:part:22:4}%
		\OdaRelationCell{\(\begin{aligned}
				% Oda1979: II.2.1(4), p.82
				&4\sigma_{14}\rho_{21}=2\sigma_{15}\rho_{22}=\sigma_{17}\rho_{24}=\rho_{21}\sigma_{36}=0.\\
				% Oda1979: II.2.1(5), p.82
				&\nu_5\bar\zeta_8=\rho^{\mathrm{IV}}\sigma_{20},\quad \nu_6\bar\zeta_9=2\rho'''\sigma_{21},\quad \nu_{13}\bar\zeta_{16}=0.\\
				% Oda1979: III.2.3(1), p.109
				&C_1\eta_{23}\mu_{24}=8\sigma_{11}\rho_{18}.\\
				% Oda1979: III.2.3(4), p.109
				&\rho_{14}\sigma_{29}=2x\sigma^{*\prime\prime},\quad \rho_{15}\sigma_{30}=4y\sigma^{*\prime}\quad(x,y\text{ odd}).
			\end{aligned}\)}\\
		\hhline{>{\arrayrulecolor{odapurple}}->{\arrayrulecolor{odarule}}--}
		% Stem 22, part 5: 1 rows.
		\cellcolor{odapurple}\leavevmode\color{black}\OdaPageStem{22}{5}
		& \cellcolor{odapurplepart}\textup{part~5}
		& \cellcolor{odapurplebody}\leavevmode\label{oda:part:22:5}%
		\OdaRelationCell{\(\begin{aligned}
				% Oda1979: III.2.5(7), p.113
				&A_2\eta_{40}^{2}=2P(\nu_{41})\pm4\sigma_{20}^{*}.
			\end{aligned}\)}\\
		\hline
		% Stem 23, part 1: 4 rows.
		\cellcolor{odagreen}\leavevmode\color{black}\OdaPageStem{23}{1}
		& \cellcolor{odagreenpart}\textup{part~1}
		& \cellcolor{odagreenbody}\leavevmode\label{oda:part:23:1}%
		\OdaRelationCell{\(\begin{aligned}
				% Oda1979: I.3.1(2), p.54
				&\eta_5\varepsilon_6\kappa_{14}=4\nu_5\bar\kappa_8.\\
				% Oda1979: I.3.1(4), p.55
				&\rho^{\mathrm{IV}}\varepsilon_{20}=\rho^{\mathrm{IV}}\bar\nu_{20}=\varepsilon_5\rho_{13}=\bar\nu_6\rho_{14}=0.\\
				% Oda1979: I.3.5(8), p.60
				&2\iota_5\circ\phi_5=0.\\
				% Oda1979: II.2.1(11), p.82
				&\zeta'\sigma_{22}=0.
			\end{aligned}\)}\\
		\hhline{>{\arrayrulecolor{odagreen}}->{\arrayrulecolor{odarule}}--}
		% Stem 23, part 2: 4 rows.
		\cellcolor{odagreen}\leavevmode\color{black}\OdaPageStem{23}{2}
		& \cellcolor{odagreenpart}\textup{part~2}
		& \cellcolor{odagreenbody}\leavevmode\label{oda:part:23:2}%
		\OdaRelationCell{\(\begin{aligned}
				% Oda1979: III.2.1(2), p.107
				&\varepsilon_{12}^{*}\nu_{29}^{2}=0,\quad \sigma_9\Sigma\eta^{*\prime}\equiv0\m \sigma_9^2\mu_{23}.\\
				% Oda1979: III.2.1(3), p.107
				&P(\nu_{43})=\omega_{21}\sigma_{37}-\sigma_{21}\omega_{28}.\\
				% Oda1979: III.2.1(4), p.107
				&\omega_{14}\sigma_{30}\equiv\psi_{14}\m \phi_{14},\;4\nu_{14}\bar\kappa_{17}.\\
				% Oda1979: III.2.2(1), p.108
				&\rho'\varepsilon_{24}=\sigma_9^2\mu_{23},\quad \rho'\bar\nu_{24}=0.
			\end{aligned}\)}\\
		\hhline{>{\arrayrulecolor{odagreen}}->{\arrayrulecolor{odarule}}--}
		% Stem 23, part 3: 2 rows.
		\cellcolor{odagreen}\leavevmode\color{black}\OdaPageStem{23}{3}
		& \cellcolor{odagreenpart}\textup{part~3}
		& \cellcolor{odagreenbody}\leavevmode\label{oda:part:23:3}%
		\OdaRelationCell{\(\begin{aligned}
				% Oda1979: III.2.5(5), p.113
				&\eta_{15}\sigma_{16}^{*}=\eta^{*\prime}\sigma_{31}+\sigma_{15}\omega_{22}+\omega_{15}\sigma_{31}.\\
				% Oda1979: III.2.5(5), p.113
				&\sigma_{16}^{*}\eta_{38}\equiv\eta_{16}^{*}\sigma_{32}+\sigma_{16}\omega_{23}\m  P(\sigma_{33}\eta_{40}).
			\end{aligned}\)}\\
		\hline
		% Stem 24, part 1: 4 rows.
		\cellcolor{odapurple}\leavevmode\color{black}\OdaPageStem{24}{1}
		& \cellcolor{odapurplepart}\textup{part~1}
		& \cellcolor{odapurplebody}\leavevmode\label{oda:part:24:1}%
		\OdaRelationCell{\(\begin{aligned}
				% Oda1979: I.3.1(2), p.54
				&\bar\alpha\eta_{26}=\varepsilon'\kappa_{13}=\nu'\eta_6\bar\kappa_7=\eta_3\nu_4\bar\kappa_7.\\
				% Oda1979: I.3.1(2), p.54
				&\zeta'\varepsilon_{22}=\zeta'\bar\nu_{22}=0.\\
				% Oda1979: I.3.1(5), p.55; I.3.5(2), p.60
				&\eta_3\mu_4\kappa_{13}=\bar\varepsilon_3\mu_{18}=\bar\varepsilon'\sigma_{20}=\varepsilon'\sigma_{13}^{2}=0.\\
				% Oda1979: I.3.1(6), p.55
				&\alpha\nu_{26}\equiv\nu_5\sigma_8\kappa_{15}\m \nu_5\eta_8\bar\kappa_9.
			\end{aligned}\)}\\
		\hhline{>{\arrayrulecolor{odapurple}}->{\arrayrulecolor{odarule}}--}
		% Stem 24, part 2: 4 rows.
		\cellcolor{odapurple}\leavevmode\color{black}\OdaPageStem{24}{2}
		& \cellcolor{odapurplepart}\textup{part~2}
		& \cellcolor{odapurplebody}\leavevmode\label{oda:part:24:2}%
		\OdaRelationCell{\(\begin{aligned}
				% Oda1979: I.3.5(9), p.60
				&\eta_4\phi_5\equiv\delta_4\m \bar\mu_4\sigma_{21},\;(\Sigma\varepsilon')\kappa_{14}.\\
				% Oda1979: I.3.5(9), p.60
				&\phi_5\eta_{28}\equiv\delta_5\m \bar\mu_5\sigma_{22},\;\nu_5\eta_8\bar\kappa_9.\\
				% Oda1979: I.3.5(9), p.60
				&\varepsilon_6\omega_{14}\equiv\eta_6\phi_7\m \nu_6\sigma_9\kappa_{16}.\\
				% Oda1979: II.2.1(10), p.82
				&\eta_9\psi_{10}\equiv\bar\sigma'_9+\delta_9\m \bar\mu_9\sigma_{26},\;\sigma_9^2\eta_{23}\mu_{24}.
			\end{aligned}\)}\\
		\hhline{>{\arrayrulecolor{odapurple}}->{\arrayrulecolor{odarule}}--}
		% Stem 24, part 3: 4 rows.
		\cellcolor{odapurple}\leavevmode\color{black}\OdaPageStem{24}{3}
		& \cellcolor{odapurplepart}\textup{part~3}
		& \cellcolor{odapurplebody}\leavevmode\label{oda:part:24:3}%
		\OdaRelationCell{\(\begin{aligned}
				% Oda1979: III.2.1(3), p.107
				&P(\nu_{39}^2)=\omega_{19}\bar\nu_{35}-\bar\nu_{19}\omega_{27}.\\
				% Oda1979: III.2.1(3), p.107
				&P(\nu_{33}^3)=\omega_{16}\varepsilon_{32}-\varepsilon_{16}\omega_{24}.\\
				% Oda1979: III.2.1(4), p.107
				&\omega_{14}\varepsilon_{30}\equiv\tilde\varepsilon_{14}\m \bar\mu_{14}\sigma_{31},\;\delta_{14},\;\bar\sigma'_{14}.\\
				% Oda1979: III.2.1(4), p.107
				&\omega_{14}\sigma_{30}\eta_{37}=\omega_{14}\varepsilon_{30}+\omega_{14}\bar\nu_{30}\equiv\bar\sigma'_{14}\m \delta_{14},\;\bar\mu_{14}\sigma_{31}.
			\end{aligned}\)}\\
		\hhline{>{\arrayrulecolor{odapurple}}->{\arrayrulecolor{odarule}}--}
		% Stem 24, part 4: 4 rows.
		\cellcolor{odapurple}\leavevmode\color{black}\OdaPageStem{24}{4}
		& \cellcolor{odapurplepart}\textup{part~4}
		& \cellcolor{odapurplebody}\leavevmode\label{oda:part:24:4}%
		\OdaRelationCell{\(\begin{aligned}
				% Oda1979: III.2.1(4), p.107
				&\varepsilon_{12}^{*}\sigma_{29}\equiv\bar\sigma'_{12}\m \delta_{12},\;\bar\mu_{12}\sigma_{29}.\\
				% Oda1979: III.2.2(1), p.108
				&\bar\rho_9\eta_{32}\equiv\sigma_9\bar\mu_{16}\m \sigma_9^2\mu_{23}\eta_{32},\;\sigma_9\varepsilon_{16}^{*},\;\phi_9\eta_{32}.\\
				% Oda1979: III.2.2(1), p.108
				&\bar\mu_4\sigma_{21}\equiv\mu_4\rho_{13}\m (\Sigma\varepsilon')\kappa_{14}.\\
				% Oda1979: III.2.2(1), p.108
				&\bar\mu_5\sigma_{22}=\mu_5\rho_{14}.
			\end{aligned}\)}\\
		\hhline{>{\arrayrulecolor{odapurple}}->{\arrayrulecolor{odarule}}--}
		% Stem 24, part 5: 4 rows.
		\cellcolor{odapurple}\leavevmode\color{black}\OdaPageStem{24}{5}
		& \cellcolor{odapurplepart}\textup{part~5}
		& \cellcolor{odapurplebody}\leavevmode\label{oda:part:24:5}%
		\OdaRelationCell{\(\begin{aligned}
				% Oda1979: III.2.2(4), p.108
				&P(16\rho_{21})=(\Sigma\rho')\mu_{25}=\bar\rho_{10}\eta_{33}+\eta_{10}\phi_{11}=\sigma_{10}\bar\mu_{17}+\bar\mu_{10}\sigma_{27}.\\
				% Oda1979: III.2.2(4), p.108
				&\bar\rho_9\eta_{32}\equiv\sigma_9\bar\mu_{16}\m \sigma_9^2\mu_{23}\eta_{32},\;\sigma_9\varepsilon_{16}^{*}+\eta_9\sigma_{10}\omega_{17}.\\
				% Oda1979: III.2.2(4), p.108
				&\bar\rho_{10}\eta_{33}=\sigma_{10}\bar\mu_{17}.\\
				% Oda1979: III.2.4(2), p.111
				&\nu_5(\Sigma\sigma')\kappa_{15}=0.
			\end{aligned}\)}\\
		\hhline{>{\arrayrulecolor{odapurple}}->{\arrayrulecolor{odarule}}--}
		% Stem 24, part 6: 2 rows.
		\cellcolor{odapurple}\leavevmode\color{black}\OdaPageStem{24}{6}
		& \cellcolor{odapurplepart}\textup{part~6}
		& \cellcolor{odapurplebody}\leavevmode\label{oda:part:24:6}%
		\OdaRelationCell{\(\begin{aligned}
				% Oda1979: III.2.6(3), p.114
				&\bar\rho'''\eta_{28}=\rho^{\mathrm{IV}}\mu_{20}=\sigma'''\bar\mu_{12}=0.\\
				% Oda1979: III.2.6(5), p.114
				&\bar\nu_7\omega_{15}\equiv0\m \nu_7\sigma_{10}\kappa_{17},\;\bar\zeta'_7,\quad \bar\nu_9\omega_{17}=0.
			\end{aligned}\)}\\
		\hline
		% Stem 25, part 1: 4 rows.
		\cellcolor{odagreen}\leavevmode\color{black}\OdaPageStem{25}{1}
		& \cellcolor{odagreenpart}\textup{part~1}
		& \cellcolor{odagreenbody}\leavevmode\label{oda:part:25:1}%
		\OdaRelationCell{\(\begin{aligned}
				% Oda1979: I.3.1(7), p.55
				&\nu_5^2\bar\zeta_{11}=2\nu_5(\Sigma\sigma')\rho_{15}=0.\\
				% Oda1979: I.3.1(7), p.55
				&\zeta_5\sigma_{16}^{2}=0.\\
				% Oda1979: I.3.1(8), p.55
				&\bar\nu_6\bar\mu_{14}=\eta_9^2\sigma_{11}\omega_{18}=0.\\
				% Oda1979: I.3.1(8), p.55; I.3.5(3), p.60
				&\sigma_9\Sigma^3\lambda=2(\sigma_9\nu_{16}^{*}).
			\end{aligned}\)}\\
		\hhline{>{\arrayrulecolor{odagreen}}->{\arrayrulecolor{odarule}}--}
		% Stem 25, part 2: 4 rows.
		\cellcolor{odagreen}\leavevmode\color{black}\OdaPageStem{25}{2}
		& \cellcolor{odagreenpart}\textup{part~2}
		& \cellcolor{odagreenbody}\leavevmode\label{oda:part:25:2}%
		\OdaRelationCell{\(\begin{aligned}
				% Oda1979: I.3.1(8), p.55
				&\sigma_9\eta_{16}^{2}\omega_{18}=4\sigma_9\xi_{16}.\\
				% Oda1979: I.3.1(8), p.55; III.2.2(3), p.108
				&\eta_9\sigma_{10}\bar\mu_{17}=\eta_9\bar\mu_{10}\sigma_{27}=0.\\
				% Oda1979: I.3.5(2), p.60
				&\bar\mu_3\bar\nu_{20}=0.\\
				% Oda1979: I.3.5(2), p.60
				&\varepsilon_3\bar\mu_{11}=\bar\mu_3\varepsilon_{20}=\mu_3^2\sigma_{21}=\eta_3\bar\mu_4\sigma_{21}.
			\end{aligned}\)}\\
		\hhline{>{\arrayrulecolor{odagreen}}->{\arrayrulecolor{odarule}}--}
		% Stem 25, part 3: 4 rows.
		\cellcolor{odagreen}\leavevmode\color{black}\OdaPageStem{25}{3}
		& \cellcolor{odagreenpart}\textup{part~3}
		& \cellcolor{odagreenbody}\leavevmode\label{oda:part:25:3}%
		\OdaRelationCell{\(\begin{aligned}
				% Oda1979: I.3.5(3), p.60
				&\sigma_{12}\xi_{19}=2\xi_{12}\sigma_{30}.\\
				% Oda1979: I.3.5(3), p.60
				&\xi_{13}\sigma_{31}\equiv\nu_{13}\sigma_{16}^{*}\m (\lambda+2\xi_{13})\sigma_{31}.\\
				% Oda1979: I.3.5(3), p.60
				&\nu_{16}^{*}\sigma_{34}=4\xi_{16}\sigma_{34}=2\sigma_{16}\xi_{23}.\\
				% Oda1979: I.3.5(3), p.60
				&\nu_{19}^{*}\sigma_{37}=\sigma_{19}\xi_{26}=2\xi_{19}\sigma_{37}=0.
			\end{aligned}\)}\\
		\hhline{>{\arrayrulecolor{odagreen}}->{\arrayrulecolor{odarule}}--}
		% Stem 25, part 4: 4 rows.
		\cellcolor{odagreen}\leavevmode\color{black}\OdaPageStem{25}{4}
		& \cellcolor{odagreenpart}\textup{part~4}
		& \cellcolor{odagreenbody}\leavevmode\label{oda:part:25:4}%
		\OdaRelationCell{\(\begin{aligned}
				% Oda1979: I.3.5(4), p.60
				&B_1\eta_{32}\varepsilon_{33}=\omega_{16}\eta_{32}\varepsilon_{33}=\eta_{16}^{*}\eta_{32}\varepsilon_{33}=0.\\
				% Oda1979: I.3.5(4), p.60
				&\eta_9\varepsilon_{10}\omega_{18}=0.\\
				% Oda1979: I.3.5(4), p.60
				&\varepsilon_9\varepsilon_{17}^{*}=\varepsilon_9^{*}\varepsilon_{26}=\bar\nu_7\varepsilon_{15}^{*}=0.\\
				% Oda1979: I.3.5(4), p.60
				&\eta_{14}\omega_{15}\bar\nu_{31}=\omega_{14}\eta_{30}\bar\nu_{31}=\omega_{14}\nu_{30}^{3}.
			\end{aligned}\)}\\
		\hhline{>{\arrayrulecolor{odagreen}}->{\arrayrulecolor{odarule}}--}
		% Stem 25, part 5: 4 rows.
		\cellcolor{odagreen}\leavevmode\color{black}\OdaPageStem{25}{5}
		& \cellcolor{odagreenpart}\textup{part~5}
		& \cellcolor{odagreenbody}\leavevmode\label{oda:part:25:5}%
		\OdaRelationCell{\(\begin{aligned}
				% Oda1979: I.3.5(4), p.60
				&\eta^{*\prime}\eta_{31}\varepsilon_{32}=\eta_{15}\eta_{16}^{*}\varepsilon_{32}=0.\\
				% Oda1979: I.3.5(4), p.60
				&\eta_{16}\omega_{17}\bar\nu_{33}=\eta_{16}\omega_{17}\eta_{33}\sigma_{34}=4\xi_{16}\sigma_{34}.\\
				% Oda1979: I.5.1(2), p.64
				&\sigma_9\nu_{16}^{*}\equiv\Sigma^2\phi'''\m \nu_9^2\bar\sigma_{15}.\\
				% Oda1979: I.5.1(3), p.64
				&\eta_9^2\phi_{11}=\sigma_9\eta_{16}^{*}\eta_{32}^{2}=4(\sigma_9\nu_{16}^{*})=0.
			\end{aligned}\)}\\
		\hhline{>{\arrayrulecolor{odagreen}}->{\arrayrulecolor{odarule}}--}
		% Stem 25, part 6: 4 rows.
		\cellcolor{odagreen}\leavevmode\color{black}\OdaPageStem{25}{6}
		& \cellcolor{odagreenpart}\textup{part~6}
		& \cellcolor{odagreenbody}\leavevmode\label{oda:part:25:6}%
		\OdaRelationCell{\(\begin{aligned}
				% Oda1979: I.5.1(3), p.64; III.2.2(3), p.108
				&\eta_9\delta_{10}=0.\\
				% Oda1979: II.2.1(6), p.82
				&\xi'\sigma_{29}\equiv\sigma_{11}\nu_{18}^{*}\m 2\sigma_{11}\nu_{18}^{*}.\\
				% Oda1979: II.2.1(7), p.82
				&\eta_5\bar\sigma'_6\equiv\nu_5^2\bar\sigma_{11}\m \eta_5\bar\mu_6\sigma_{23}.\\
				% Oda1979: II.2.1(7), p.82
				&\bar\sigma'_6\eta_{30}\equiv\nu_6^2\bar\sigma_{12}\m \eta_6\bar\mu_7\sigma_{24},\;P((\Sigma A_1)\varepsilon_{25}).
			\end{aligned}\)}\\
		\hhline{>{\arrayrulecolor{odagreen}}->{\arrayrulecolor{odarule}}--}
		% Stem 25, part 7: 4 rows.
		\cellcolor{odagreen}\leavevmode\color{black}\OdaPageStem{25}{7}
		& \cellcolor{odagreenpart}\textup{part~7}
		& \cellcolor{odagreenbody}\leavevmode\label{oda:part:25:7}%
		\OdaRelationCell{\(\begin{aligned}
				% Oda1979: II.2.3(2), p.86; II.(6.4), p.98
				&2\phi'''\equiv\nu_7^2\bar\sigma_{13}\m 2\Sigma^2\phi''.\\
				% Oda1979: II.2.3(2), p.86
				&2\Sigma^2\phi'''\equiv\nu_9^2\bar\sigma_{15}\m 4\sigma_9\xi_{16}.\\
				% Oda1979: II.2.3(3), p.86
				&\eta_5\phi_6\eta_{29}=4\phi''.\\
				% Oda1979: II.2.3(3), p.86; II.(6.4), p.98
				&\sigma'\Sigma\lambda\equiv0\m 4\Sigma^2\phi''.
			\end{aligned}\)}\\
		\hhline{>{\arrayrulecolor{odagreen}}->{\arrayrulecolor{odarule}}--}
		% Stem 25, part 8: 4 rows.
		\cellcolor{odagreen}\leavevmode\color{black}\OdaPageStem{25}{8}
		& \cellcolor{odagreenpart}\textup{part~8}
		& \cellcolor{odagreenbody}\leavevmode\label{oda:part:25:8}%
		\OdaRelationCell{\(\begin{aligned}
				% Oda1979: II.2.3(4), p.86
				&(\Sigma\nu')\sigma'\rho_{14}=0.\\
				% Oda1979: II.2.3(4), p.86
				&2\iota_5\circ\phi''\equiv2\phi''\m \nu_5\varepsilon_8\kappa_{16}.\\
				% Oda1979: II.2.3(4), p.86
				&2\phi''\equiv\pm\Sigma^2\phi'\m \nu_5\varepsilon_8\kappa_{16},\;\nu_5^2\bar\sigma_{11},\;\eta_5\bar\mu_6\sigma_{23}.\\
				% Oda1979: I.5.1(1), p.64; II.(6.4), pp.98--99
				&\sigma'\xi_{14}\equiv\Sigma^2\phi''\m 2\Sigma^2\phi'',\;2\phi''',\;\eta_7\bar\mu_8\sigma_{25}.
			\end{aligned}\)}\\
		\hhline{>{\arrayrulecolor{odagreen}}->{\arrayrulecolor{odarule}}--}
		% Stem 25, part 9: 3 rows.
		\cellcolor{odagreen}\leavevmode\color{black}\OdaPageStem{25}{9}
		& \cellcolor{odagreenpart}\textup{part~9}
		& \cellcolor{odagreenbody}\leavevmode\label{oda:part:25:9}%
		\OdaRelationCell{\(\begin{aligned}
				% Oda1979: III.2.2(2), p.108
				&\lambda\sigma_{31}=0.\\
				% Oda1979: III.2.2(2), p.108
				&\eta_{13}\bar\sigma'_{14}\equiv4\xi_{13}\sigma_{31}\m \eta_{13}\bar\mu_{14}\sigma_{31}.\\
				% Oda1979: III.2.2(6), p.108
				&\sigma_{16}^{*}\nu_{38}=\sigma_{16}\nu_{23}^{*}=\pm2\xi_{16}\sigma_{34}.
			\end{aligned}\)}\\
		\hline
		% Stem 26, part 1: 4 rows.
		\cellcolor{odapurple}\leavevmode\color{black}\OdaPageStem{26}{1}
		& \cellcolor{odapurplepart}\textup{part~1}
		& \cellcolor{odapurplebody}\leavevmode\label{oda:part:26:1}%
		\OdaRelationCell{\(\begin{aligned}
				% Oda1979: I.3.1(9), p.55
				&\varepsilon'\sigma_{13}\mu_{20}=0.\\
				% Oda1979: I.3.1(9), p.55
				&\zeta_5\rho_{16}=\pm\bar\zeta_5\sigma_{24}.\\
				% Oda1979: I.3.5(5), p.60
				&2\bar\sigma_6\sigma_{25}=\nu_6\phi_9,\quad \phi_9\nu_{32}=0.\\
				% Oda1979: I.3.5(6), p.60
				&\eta_9\sigma_{10}\nu_{17}^{*}=\bar\sigma_{10}\sigma_{29}=0.
			\end{aligned}\)}\\
		\hhline{>{\arrayrulecolor{odapurple}}->{\arrayrulecolor{odarule}}--}
		% Stem 26, part 2: 4 rows.
		\cellcolor{odapurple}\leavevmode\color{black}\OdaPageStem{26}{2}
		& \cellcolor{odapurplepart}\textup{part~2}
		& \cellcolor{odapurplebody}\leavevmode\label{oda:part:26:2}%
		\OdaRelationCell{\(\begin{aligned}
				% Oda1979: I.3.5(6), p.60
				&\eta_9\sigma_{10}\xi_{17}=\bar\sigma_9\sigma_{28}.\\
				% Oda1979: I.3.5(6), p.60
				&\theta'\sigma_{23}^{2}=0,\quad A_1\sigma_{24}^{2}=\theta\sigma_{24}^{2}=\sigma_{12}\bar\sigma_{19}.\\
				% Oda1979: I.3.5(7), p.60
				&\eta_9^2\sigma_{11}\bar\mu_{18}=0.\\
				% Oda1979: I.(8.13), p.74
				&2\tau''=\pm\sigma_{10}\bar\zeta_{17}.
			\end{aligned}\)}\\
		\hhline{>{\arrayrulecolor{odapurple}}->{\arrayrulecolor{odarule}}--}
		% Stem 26, part 3: 4 rows.
		\cellcolor{odapurple}\leavevmode\color{black}\OdaPageStem{26}{3}
		& \cellcolor{odapurplepart}\textup{part~3}
		& \cellcolor{odapurplebody}\leavevmode\label{oda:part:26:3}%
		\OdaRelationCell{\(\begin{aligned}
				% Oda1979: I.(8.14), p.74
				&\sigma'\bar\zeta_{14}=x\bar\zeta_7\sigma_{26}=y\zeta_7\rho_{18}\quad(x,y\text{: odd}).\\
				% Oda1979: I.(8.14), p.74
				&\bar\zeta_9\sigma_{28}=\pm2\sigma_9\bar\zeta_{16}.\\
				% Oda1979: I.(8.14), p.74
				&2\bar\zeta_{10}\sigma_{29}=4\sigma_{10}\bar\zeta_{17}=\bar\zeta_{11}\sigma_{30}=0.\\
				% Oda1979: I.(8.14), p.74
				&2\sigma_{11}\bar\zeta_{18}=\sigma_{13}\bar\zeta_{20}=0.
			\end{aligned}\)}\\
		\hhline{>{\arrayrulecolor{odapurple}}->{\arrayrulecolor{odarule}}--}
		% Stem 26, part 4: 4 rows.
		\cellcolor{odapurple}\leavevmode\color{black}\OdaPageStem{26}{4}
		& \cellcolor{odapurplepart}\textup{part~4}
		& \cellcolor{odapurplebody}\leavevmode\label{oda:part:26:4}%
		\OdaRelationCell{\(\begin{aligned}
				% Oda1979: II.2.1(8), p.82
				&\sigma'''\bar\zeta_{12}=\rho^{\mathrm{IV}}\zeta_{20}=\bar\rho'''\nu_{28}=4\zeta_5\rho_{16}.\\
				% Oda1979: II.2.1(8), p.82
				&\sigma'\bar\sigma_{14}=\phi_7\nu_{30}.\\
				% Oda1979: II.2.1(9), p.82
				&\eta_7\sigma_8\xi_{15}=\bar\sigma_7\sigma_{26},\quad \eta_7\Sigma\phi'''=\phi_7\nu_{30}.\\
				% Oda1979: II.2.1(9), p.82
				&\eta_8\sigma_9\nu_{16}^{*}=\phi_8\nu_{31}.
			\end{aligned}\)}\\
		\hhline{>{\arrayrulecolor{odapurple}}->{\arrayrulecolor{odarule}}--}
		% Stem 26, part 5: 3 rows.
		\cellcolor{odapurple}\leavevmode\color{black}\OdaPageStem{26}{5}
		& \cellcolor{odapurplepart}\textup{part~5}
		& \cellcolor{odapurplebody}\leavevmode\label{oda:part:26:5}%
		\OdaRelationCell{\(\begin{aligned}
				% Oda1979: II.(6.3), p.98
				&\varepsilon_4\xi_{12}=0,\quad \bar\nu_6\xi_{14}=\pm\bar\sigma_6\sigma_{25}.\\
				% Oda1979: III.2.2(3), p.108
				&\varepsilon_5\lambda=\zeta_{11}\rho_{22}=0,\quad \sigma_{10}\bar\sigma_{17}=\psi_{10}\nu_{33}.\\
				% Oda1979: III.2.2(3), p.108
				&\nu_7\psi_{10}\equiv\sigma'\bar\sigma_{14}\m \bar\sigma_7\sigma_{26},\;\bar\zeta_7\sigma_{26}.
			\end{aligned}\)}\\
		\hline
		% Stem 27, part 1: 4 rows.
		\cellcolor{odagreen}\leavevmode\color{black}\OdaPageStem{27}{1}
		& \cellcolor{odagreenpart}\textup{part~1}
		& \cellcolor{odagreenbody}\leavevmode\label{oda:part:27:1}%
		\OdaRelationCell{\(\begin{aligned}
				% Oda1979: I.3.1(10), p.55
				&\eta_3^2\bar\mu_5\varepsilon_{22}=\eta_3^2\bar\mu_5\bar\nu_{22}=0.\\
				% Oda1979: II.2.2(4), p.84
				&\bar\zeta_5\bar\nu_{24}=0.\\
				% Oda1979: II.2.2(4), p.84
				&\bar\zeta_5\varepsilon_{24}=\bar\zeta_5\eta_{24}\sigma_{25}=\zeta_5\mu_{16}\sigma_{25}=\nu_5\bar\mu_8\sigma_{25}.\\
				% Oda1979: II.2.2(5), p.84
				&\varepsilon_3\bar\sigma_{11}\equiv\nu'\bar\sigma'_6\m \nu'\bar\mu_6\sigma_{23}.
			\end{aligned}\)}\\
		\hhline{>{\arrayrulecolor{odagreen}}->{\arrayrulecolor{odarule}}--}
		% Stem 27, part 2: 4 rows.
		\cellcolor{odagreen}\leavevmode\color{black}\OdaPageStem{27}{2}
		& \cellcolor{odagreenpart}\textup{part~2}
		& \cellcolor{odagreenbody}\leavevmode\label{oda:part:27:2}%
		\OdaRelationCell{\(\begin{aligned}
				% Oda1979: II.2.2(5), p.84
				&\varepsilon_3\bar\zeta_{11}=0.\\
				% Oda1979: II.2.2(5), p.84
				&\varepsilon'\varepsilon_{13}^{*}=\varepsilon'\eta_{13}\omega_{14}=\nu'\varepsilon_6\omega_{14}\equiv\nu'\delta_6\m \nu'\bar\mu_6\sigma_{23}.\\
				% Oda1979: II.2.2(6), p.84
				&\varepsilon'\bar\mu_{13}=\nu'\bar\mu_6\sigma_{23}.\\
				% Oda1979: II.2.2(7), p.84
				&\varepsilon_5\bar\sigma_{13}=\nu_6\bar\mu_9\sigma_{26}=\nu_6\bar\sigma'_9=\nu_6\delta_9=0.
			\end{aligned}\)}\\
		\hhline{>{\arrayrulecolor{odagreen}}->{\arrayrulecolor{odarule}}--}
		% Stem 27, part 3: 4 rows.
		\cellcolor{odagreen}\leavevmode\color{black}\OdaPageStem{27}{3}
		& \cellcolor{odagreenpart}\textup{part~3}
		& \cellcolor{odagreenbody}\leavevmode\label{oda:part:27:3}%
		\OdaRelationCell{\(\begin{aligned}
				% Oda1979: II.2.2(7), p.84
				&\bar\sigma'_6\nu_{30}=\bar\nu_6\bar\sigma_{14}.\\
				% Oda1979: II.2.2(8), p.84
				&\nu_5\eta_8\bar\rho_9=\nu_5\bar\mu_8\sigma_{25}.\\
				% Oda1979: II.2.2(10), p.84
				&A_1\bar\varepsilon_{24}=2\sigma_{12}\bar\kappa_{19},\quad C_1\bar\varepsilon_{23}=4\sigma_{11}\bar\kappa_{18}=P(\nu_{23}\kappa_{26}).\\
				% Oda1979: III.2.1(1), p.107
				&A_1\rho_{24}=A_1^{(1)}\sigma_{32},\quad C_1\rho_{23}=C_1^{(1)}\sigma_{31}.
			\end{aligned}\)}\\
		\hhline{>{\arrayrulecolor{odagreen}}->{\arrayrulecolor{odarule}}--}
		% Stem 27, part 4: 4 rows.
		\cellcolor{odagreen}\leavevmode\color{black}\OdaPageStem{27}{4}
		& \cellcolor{odagreenpart}\textup{part~4}
		& \cellcolor{odagreenbody}\leavevmode\label{oda:part:27:4}%
		\OdaRelationCell{\(\begin{aligned}
				% Oda1979: III.2.2(5), p.108
				&\bar\kappa_7\sigma_{27}\equiv0\m 2\sigma'\bar\kappa_{14},\;\bar\nu_7\bar\sigma_{15}.\\
				% Oda1979: III.2.2(5), p.108
				&2\bar\kappa_9\sigma_{29}=0.\\
				% Oda1979: III.2.4(2), p.111
				&\nu_5\eta_8\bar\kappa_9\nu_{29}=0.\\
				% Oda1979: III.2.4(2), p.111
				&\alpha\nu_{26}^{2}=\sigma'''\bar\kappa_{12}=\nu_5\sigma_8\kappa_{15}\nu_{29}.
			\end{aligned}\)}\\
		\hline
		% Stem 28, part 1: 4 rows.
		\cellcolor{odapurple}\leavevmode\color{black}\OdaPageStem{28}{1}
		& \cellcolor{odapurplepart}\textup{part~1}
		& \cellcolor{odapurplebody}\leavevmode\label{oda:part:28:1}%
		\OdaRelationCell{\(\begin{aligned}
				% Oda1979: II.2.2(1)--(2), p.84
				&\zeta_{3,5}\eta_{32}\equiv\bar\zeta_5\mu_{24}\equiv\zeta_5\bar\mu_{16}\equiv\nu_5\mu_{3,8}\m \nu_5\eta_8\bar\mu_9\sigma_{26}.\\
				% Oda1979: II.2.3(1), p.86
				&\mu_{3,3}\nu_{28}=\bar\mu_3\zeta_{20}=\mu_3\bar\zeta_{12}=\nu'\eta_6\bar\mu_7\sigma_{24}.\\
				% Oda1979: II.2.3(1), p.86
				&\varepsilon'\eta_{13}\bar\mu_{14}=\nu'\varepsilon_6\bar\mu_{14}=\nu'\eta_6\bar\mu_7\sigma_{24}.\\
				% Oda1979: II.2.3(2), p.86
				&2\nu_5\Sigma\phi'''=\nu_5^3\bar\sigma_{14}.
			\end{aligned}\)}\\
		\hhline{>{\arrayrulecolor{odapurple}}->{\arrayrulecolor{odarule}}--}
		% Stem 28, part 2: 4 rows.
		\cellcolor{odapurple}\leavevmode\color{black}\OdaPageStem{28}{2}
		& \cellcolor{odapurplepart}\textup{part~2}
		& \cellcolor{odapurplebody}\leavevmode\label{oda:part:28:2}%
		\OdaRelationCell{\(\begin{aligned}
				% Oda1979: II.2.3(2), p.86
				&\nu_5\Sigma^3\phi''\equiv0\m \nu_5\eta_8\bar\mu_9\sigma_{26},\;\nu_5^3\bar\sigma_{14}.\\
				% Oda1979: II.2.3(3), p.86
				&\nu_6^3\bar\sigma_{15}=\nu_7\Sigma^3\phi'''=0.\\
				% Oda1979: II.2.3(4), p.86
				&(\Sigma\phi')\nu_{29}=(\Sigma\nu')\phi''',\quad 2\phi'''\nu_{32}=0.\\
				% Oda1979: II.2.3(4), p.86
				&\nu_5\Sigma\phi'''\equiv\phi''\nu_{30}\m \nu_5^3\bar\sigma_{14}.
			\end{aligned}\)}\\
		\hhline{>{\arrayrulecolor{odapurple}}->{\arrayrulecolor{odarule}}--}
		% Stem 28, part 3: 4 rows.
		\cellcolor{odapurple}\leavevmode\color{black}\OdaPageStem{28}{3}
		& \cellcolor{odapurplepart}\textup{part~3}
		& \cellcolor{odapurplebody}\leavevmode\label{oda:part:28:3}%
		\OdaRelationCell{\(\begin{aligned}
				% Oda1979: II.2.4(4), p.87
				&8P(\bar\rho_{13})=\nu_6\mu_{3,9}.\\
				% Oda1979: II.(6.3), p.98
				&\varepsilon'\xi_{13}=\varepsilon'\lambda=0.\\
				% Oda1979: III.2.2(6), p.108
				&\eta_5\zeta_{3,6}=0.\\
				% Oda1979: III.2.2(6), p.108
				&\sigma_{16}^{*}\nu_{38}^{2}=\nu_{13}\sigma_{16}^{*}\nu_{38}=0.
			\end{aligned}\)}\\
		\hhline{>{\arrayrulecolor{odapurple}}->{\arrayrulecolor{odarule}}--}
		% Stem 28, part 4: 3 rows.
		\cellcolor{odapurple}\leavevmode\color{black}\OdaPageStem{28}{4}
		& \cellcolor{odapurplepart}\textup{part~4}
		& \cellcolor{odapurplebody}\leavevmode\label{oda:part:28:4}%
		\OdaRelationCell{\(\begin{aligned}
				% Oda1979: III.2.4(3), p.111
				&\bar\kappa_7\varepsilon_{27}=\sigma'\eta_{14}\bar\kappa_{15}+\varepsilon_7\bar\kappa_{15}.\\
				% MiM1: group 1; (3.20), p.22
				&\kappa_7^2=\sigma'\eta_{14}\bar\kappa_{15}+\varepsilon_7\bar\kappa_{15}.\\
				% MiM1: group 1; (3.20), p.22
				&\kappa_9^2=\varepsilon_9\bar\kappa_{17}.
			\end{aligned}\)}\\
		\hline
		% Stem 29, part 1: 4 rows.
		\cellcolor{odagreen}\leavevmode\color{black}\OdaPageStem{29}{1}
		& \cellcolor{odagreenpart}\textup{part~1}
		& \cellcolor{odagreenbody}\leavevmode\label{oda:part:29:1}%
		\OdaRelationCell{\(\begin{aligned}
				% Oda1979: II.2.1(6), p.82
				&\sigma'\rho_{14}\sigma_{29}=0.\\
				% Oda1979: II.2.4(1), p.87
				&\varepsilon'\bar\sigma_{13}=\varepsilon'\bar\zeta_{13}=0.\\
				% Oda1979: II.2.4(2), p.87
				&\nu_5\phi_8\nu_{31}=\bar\rho'''\nu_{28}^2=0.\\
				% Oda1979: II.2.4(3), p.87
				&\nu_5\sigma_8\bar\zeta_{15}=\nu_5\bar\zeta_8\sigma_{27}=\nu_5\zeta_8\rho_{19}=0.
			\end{aligned}\)}\\
		\hhline{>{\arrayrulecolor{odagreen}}->{\arrayrulecolor{odarule}}--}
		% Stem 29, part 2: 4 rows.
		\cellcolor{odagreen}\leavevmode\color{black}\OdaPageStem{29}{2}
		& \cellcolor{odagreenpart}\textup{part~2}
		& \cellcolor{odagreenbody}\leavevmode\label{oda:part:29:2}%
		\OdaRelationCell{\(\begin{aligned}
				% Oda1979: II.2.4(3), p.87
				&\nu_5\eta_8\mu_9\omega_{18}=0.\\
				% Oda1979: III.2.4(1), p.111
				&\phi_5\nu_{28}^{2}\equiv\nu_5\sigma_8\bar\sigma_{15}\m \nu_5\bar\kappa_8\nu_{28}^{2}-\nu_5^2\bar\kappa_{11}\nu_{31}.\\
				% Oda1979: III.2.4(2), p.111
				&(\Sigma^2\bar\alpha)\nu_{28}^{2}=0.\\
				% Oda1979: III.2.4(2), p.111
				&\mu_3\bar\kappa_{12}\equiv\bar\alpha\nu_{26}^{2}\m \eta_3\varepsilon_4\bar\kappa_{12}.
			\end{aligned}\)}\\
		\hhline{>{\arrayrulecolor{odagreen}}->{\arrayrulecolor{odarule}}--}
		% Stem 29, part 3: 4 rows.
		\cellcolor{odagreen}\leavevmode\color{black}\OdaPageStem{29}{3}
		& \cellcolor{odagreenpart}\textup{part~3}
		& \cellcolor{odagreenbody}\leavevmode\label{oda:part:29:3}%
		\OdaRelationCell{\(\begin{aligned}
				% Oda1979: III.2.4(3), p.111
				&\bar\varepsilon_3\kappa_{18}=\bar\alpha\nu_{26}^{2}+\eta_3\varepsilon_4\bar\kappa_{12}.\\
				% Oda1979: III.2.4(3), p.111
				&\kappa_7^2\eta_{35}=\sigma'\varepsilon_{14}\kappa_{22}+\bar\varepsilon_7\kappa_{22}.\\
				% Oda1979: III.2.4(3), p.111
				&\eta_7\bar\kappa_8\varepsilon_{28}=\eta_7\varepsilon_8\bar\kappa_{16}=\bar\varepsilon_7\kappa_{22}.\\
				% Oda1979: III.2.4(3), p.111
				&\eta_6\bar\kappa_7\varepsilon_{27}=\eta_6\varepsilon_7\bar\kappa_{15}.
			\end{aligned}\)}\\
		\hhline{>{\arrayrulecolor{odagreen}}->{\arrayrulecolor{odarule}}--}
		% Stem 29, part 4: 4 rows.
		\cellcolor{odagreen}\leavevmode\color{black}\OdaPageStem{29}{4}
		& \cellcolor{odagreenpart}\textup{part~4}
		& \cellcolor{odagreenbody}\leavevmode\label{oda:part:29:4}%
		\OdaRelationCell{\(\begin{aligned}
				% Oda1979: III.2.4(4), p.111
				&\kappa_7\rho_{21}=0,\quad \sigma'\sigma_{14}\rho_{21}=0.\\
				% Oda1979: III.2.4(4), p.111
				&\rho''\kappa_{22}\ty0\m \nu_7\cg_{10}\ccs_{17},\;\sigma'\varepsilon_{14}\kappa_{22}.\\
				% Oda1979: III.2.4(4), pp.111--112
				&(\Sigma\sigma')\sigma^{*\prime}=-2\Sigma\delta''.\\
				% Oda1979: III.2.4(5), p.111
				&C_1\eta_{23}\omega_{24}=(\Sigma^2\rho')\sigma_{26}^{2}=8\sigma_{11}\sigma_{18}^{*}.
			\end{aligned}\)}\\
		\hhline{>{\arrayrulecolor{odagreen}}->{\arrayrulecolor{odarule}}--}
		% Stem 29, part 5: 4 rows.
		\cellcolor{odagreen}\leavevmode\color{black}\OdaPageStem{29}{5}
		& \cellcolor{odagreenpart}\textup{part~5}
		& \cellcolor{odagreenbody}\leavevmode\label{oda:part:29:5}%
		\OdaRelationCell{\(\begin{aligned}
				% Oda1979: III.2.4(5), p.111
				&4P(\nu_{25}^{*})=(\Sigma C_1)\eta_{24}\omega_{25}.\\
				% Oda1979: III.2.4(6), p.111
				&2(\sigma_{16}^{*}\sigma_{38})=0.\\
				% Oda1979: III.2.5(7), p.113
				&A_2\eta_{40}^{2}\sigma_{42}=0.\\
				% Oda1979: III.3.1, p.116
				&2\delta''=-\Sigma\delta',\quad \Sigma^2\delta''=2(\sigma_9\sigma_{16}^{*}).
			\end{aligned}\)}\\
		\hline
		% Stem 30, part 1: 4 rows.
		\cellcolor{odapurple}\leavevmode\color{black}\OdaPageStem{30}{1}
		& \cellcolor{odapurplepart}\textup{part~1}
		& \cellcolor{odapurplebody}\leavevmode\label{oda:part:30:1}%
		\OdaRelationCell{\(\begin{aligned}
				% Oda1979: III.2.2(5), p.108
				&2\nu_5\bar\kappa_8\sigma_{28}=0.\\
				% Oda1979: III.2.4(2), p.111
				&\eta_3\mu_4\bar\kappa_{13}\equiv0\m 2\varepsilon'\bar\kappa_{13}.\\
				% Oda1979: III.2.4(3), p.111
				&\eta_6\bar\kappa_7\eta_{27}\varepsilon_{28}=0.\\
				% Oda1979: III.2.5(1), p.112
				&\rho^{\mathrm{IV}}\rho_{20}=\bar\rho'''\sigma_{28},\quad \nu_5\zeta_{3,8}\equiv\rho^{\mathrm{IV}}\rho_{20}\m 2\pi_{35}(S^{5}).
			\end{aligned}\)}\\
		\hhline{>{\arrayrulecolor{odapurple}}->{\arrayrulecolor{odarule}}--}
		% Stem 30, part 2: 4 rows.
		\cellcolor{odapurple}\leavevmode\color{black}\OdaPageStem{30}{2}
		& \cellcolor{odapurplepart}\textup{part~2}
		& \cellcolor{odapurplebody}\leavevmode\label{oda:part:30:2}%
		\OdaRelationCell{\(\begin{aligned}
				% Oda1979: III.2.5(1), p.112
				&\nu_5\zeta_{3,8}=\zeta_5\bar\zeta_{16}=\bar\zeta_5\zeta_{24}=\zeta_{3,5}\nu_{32}=\sigma'''\bar\rho_{12}.\\
				% Oda1979: III.2.5(2), p.112
				&2\rho'\rho_{24}=\pm2\bar\rho_9\sigma_{32}.\\
				% Oda1979: III.2.5(2), p.112
				&2\bar\rho_9\sigma_{32}\equiv2x\sigma_9\bar\rho_{16}\m 4\sigma_9\nu_{16}\bar\kappa_{19}\quad(x\text{: odd}).\\
				% Oda1979: III.2.5(2),(6), pp.112--113
				&(\Sigma^2\rho')\rho_{26}\equiv\sigma_{11}\bar\rho_{18}\m 8\sigma_{11}\bar\rho_{18},\;\sigma_{11}\psi_{18}.
			\end{aligned}\)}\\
		\hhline{>{\arrayrulecolor{odapurple}}->{\arrayrulecolor{odarule}}--}
		% Stem 30, part 3: 4 rows.
		\cellcolor{odapurple}\leavevmode\color{black}\OdaPageStem{30}{3}
		& \cellcolor{odapurplepart}\textup{part~3}
		& \cellcolor{odapurplebody}\leavevmode\label{oda:part:30:3}%
		\OdaRelationCell{\(\begin{aligned}
				% Oda1979: III.2.5(2),(6), pp.112--113
				&\sigma_{13}\bar\rho_{20}=2x\rho_{13}^{2}\quad(x\text{: odd}).\\
				% Oda1979: III.2.5(2), p.113
				&\bar\rho_{10}\sigma_{33}\equiv \sigma_{10}\bar\rho_{17}\m 2\sigma_{10}\bar\rho_{17},\;2\sigma_{10}\nu_{17}\bar\kappa_{20},\;\sigma_{10}\omega_{17}\sigma_{33},\;\sigma_{10}(\Sigma B_1)\sigma_{33}.\\
				% Oda1979: III.2.5(2), p.113
				&16\rho_{16}^{2}=8\rho_{18}^{2}=4\rho_{22}^{2}=2\rho_{26}^{2}=0.\\
				% Oda1979: III.2.5(2), p.113
				&8\sigma_{16}\bar\rho_{23}=4\sigma_{18}\bar\rho_{25}=0.
			\end{aligned}\)}\\
		\hhline{>{\arrayrulecolor{odapurple}}->{\arrayrulecolor{odarule}}--}
		% Stem 30, part 4: 4 rows.
		\cellcolor{odapurple}\leavevmode\color{black}\OdaPageStem{30}{4}
		& \cellcolor{odapurplepart}\textup{part~4}
		& \cellcolor{odapurplebody}\leavevmode\label{oda:part:30:4}%
		\OdaRelationCell{\(\begin{aligned}
				% Oda1979: III.2.5(3), p.113
				&\rho'\eta_{24}\kappa_{25}\equiv0\m 4\sigma_9\nu_{16}\bar\kappa_{19}.\\
				% Oda1979: III.2.5(3), p.113
				&\nu_5\bar\kappa_8\sigma_{28}\equiv0\m 4\nu_5\sigma_8\bar\kappa_{15},\;\nu_5\bar\nu_8\bar\sigma_{16}.\\
				% Oda1979: III.2.5(4), p.113
				&\nu'\zeta_{3,6}=0.\\
				% Oda1979: III.2.5(5), p.113
				&\sigma_9\eta_{16}\sigma_{17}^{*}\equiv\sigma_9\psi_{16}\m \sigma_9\phi_{16},\;4\sigma_9\nu_{16}\bar\kappa_{19}.
			\end{aligned}\)}\\
		\hhline{>{\arrayrulecolor{odapurple}}->{\arrayrulecolor{odarule}}--}
		% Stem 30, part 5: 4 rows.
		\cellcolor{odapurple}\leavevmode\color{black}\OdaPageStem{30}{5}
		& \cellcolor{odapurplepart}\textup{part~5}
		& \cellcolor{odapurplebody}\leavevmode\label{oda:part:30:5}%
		\OdaRelationCell{\(\begin{aligned}
				% Oda1979: III.2.5(5), p.113
				&\eta_9\sigma_{10}\sigma_{17}^{*}=\phi_9\sigma_{32}.\\
				% Oda1979: III.2.5(6), p.113
				&\bar\nu_9\nu_{17}\bar\sigma_{20}=0,\quad 2\psi_{10}\sigma_{33}=\sigma_{10}\phi_{17}.\\
				% Oda1979: III.2.5(6), p.113
				&\sigma_{11}\phi_{18}=\sigma_{13}\psi_{20}=0.\\
				% Oda1979: III.2.5(6), p.113
				&\sigma_{11}\psi_{18}=\phi_{11}\sigma_{34}=\sigma_{11}\omega_{18}\sigma_{34}.
			\end{aligned}\)}\\
		\hhline{>{\arrayrulecolor{odapurple}}->{\arrayrulecolor{odarule}}--}
		% Stem 30, part 6: 4 rows.
		\cellcolor{odapurple}\leavevmode\color{black}\OdaPageStem{30}{6}
		& \cellcolor{odapurplepart}\textup{part~6}
		& \cellcolor{odapurplebody}\leavevmode\label{oda:part:30:6}%
		\OdaRelationCell{\(\begin{aligned}
				% Oda1979: III.2.5(8), p.113
				&\sigma_{16}^{*}\varepsilon_{38}=\sigma_{16}^{*}\bar\nu_{38}=0.\\
				% Oda1979: III.3.2(1), p.116
				&8\theta^{\mathrm{VII}}=\pm4\sigma_{12}\bar\rho_{19}.\\
				% Oda1979: III.3.2(1), p.116
				&2\theta^{\mathrm{VI}}\equiv\rho_{14}^{2}\m 2\rho_{14}^{2}.\\
				% Oda1979: III.3.2(1), p.116
				&2\theta^{\mathrm{V}}\equiv\Sigma\theta^{\mathrm{VI}}\m 2\Sigma\theta^{\mathrm{VI}},\;16\rho_{15}^{2}.
			\end{aligned}\)}\\
		\hhline{>{\arrayrulecolor{odapurple}}->{\arrayrulecolor{odarule}}--}
		% Stem 30, part 7: 4 rows.
		\cellcolor{odapurple}\leavevmode\color{black}\OdaPageStem{30}{7}
		& \cellcolor{odapurplepart}\textup{part~7}
		& \cellcolor{odapurplebody}\leavevmode\label{oda:part:30:7}\label{oda:part:30:8}%
		\OdaRelationCell{\(\begin{aligned}
				% Oda1979: III.3.2(1), p.116
				&32\theta^{\mathrm{IV}}\equiv2\sigma_{16}\bar\rho_{23}\m 4\sigma_{16}\bar\rho_{23}.\\
				% Oda1979: III.3.2(2), p.116
				&2\theta'_{23}=\Sigma\theta''.\\
				% MiM1: group 3; (3.23), p.23
				&\omega_{14}\sigma_{30}^{2}=\psi_{14}\sigma_{37}.\\
				% MiM1: group 4; Lemma 3.11, p.24
				&\bar\sigma_6\zeta_{25}=0.
			\end{aligned}\)}\\
		\hline
		% Stem 31, part 1: 4 rows.
		\cellcolor{odagreen}\leavevmode\color{black}\OdaPageStem{31}{1}
		& \cellcolor{odagreenpart}\textup{part~1}
		& \cellcolor{odagreenbody}\leavevmode\label{oda:part:31:1}%
		\OdaRelationCell{\(\begin{aligned}
				% Oda1979: III.2.1(1), p.107
				&B_1\rho_{32}=B_1^{(1)}\sigma_{40},\quad D_1\rho_{31}=D_1^{(1)}\sigma_{39}.\\
				% Oda1979: III.2.6(1), p.114
				&\bar\mu_3\sigma_{20}^{2}=\mu_3\sigma_{12}\rho_{19}=0.\\
				% Oda1979: III.2.6(1), p.114
				&\bar\rho_9\varepsilon_{32}=\rho'\sigma_{24}\mu_{31}=\sigma_9^2\bar\mu_{23}=\sigma_9\eta_{16}\bar\rho_{17}.\\
				% Oda1979: III.2.6(1), p.114
				&\varepsilon_3\psi_{11}=\delta_3\sigma_{27},\quad \varepsilon_3\phi_{11}=\varepsilon_3\bar\rho_{11}=0.
			\end{aligned}\)}\\
		\hhline{>{\arrayrulecolor{odagreen}}->{\arrayrulecolor{odarule}}--}
		% Stem 31, part 2: 4 rows.
		\cellcolor{odagreen}\leavevmode\color{black}\OdaPageStem{31}{2}
		& \cellcolor{odagreenpart}\textup{part~2}
		& \cellcolor{odagreenbody}\leavevmode\label{oda:part:31:2}%
		\OdaRelationCell{\(\begin{aligned}
				% Oda1979: III.2.6(2), p.114
				&\bar\nu_6\psi_{14}\equiv\bar\sigma'_6\sigma_{30}\m  P((\Sigma A_1)\sigma_{25}^{2}).\\
				% Oda1979: III.2.6(2), p.114
				&\bar\nu_6\phi_{14}=0,\quad \zeta'\rho_{22}\equiv\bar\zeta'_6\sigma_{30}\equiv\bar\nu_6\bar\rho_{14}\equiv0\m 8\pi_{37}(S^{6}).\\
				% Oda1979: III.2.6(3), p.114
				&2\nu_5(\Sigma\phi''')\nu_{33}=0.\\
				% Oda1979: III.2.6(3), p.114
				&\bar\rho'''\varepsilon_{28}=\bar\rho'''\bar\nu_{28}=0,\quad \mu_3\varepsilon_{12}\kappa_{20}=0.
			\end{aligned}\)}\\
		\hhline{>{\arrayrulecolor{odagreen}}->{\arrayrulecolor{odarule}}--}
		% Stem 31, part 3: 4 rows.
		\cellcolor{odagreen}\leavevmode\color{black}\OdaPageStem{31}{3}
		& \cellcolor{odagreenpart}\textup{part~3}
		& \cellcolor{odagreenbody}\leavevmode\label{oda:part:31:3}%
		\OdaRelationCell{\(\begin{aligned}
				% Oda1979: III.2.6(4), p.114
				&\sigma_9\varepsilon_{16}^{*}\sigma_{33}\equiv\sigma_9\bar\sigma'_{16}\m \sigma_9\delta_{16},\;\sigma_9^2\bar\mu_{23}.\\
				% Oda1979: III.2.6(5), p.114
				&\bar\nu_7\omega_{15}\sigma_{31}=0.\\
				% Oda1979: III.2.6(5), p.114
				&\varepsilon_6\omega_{14}\sigma_{30}=\delta_6\sigma_{30}=\eta_6\phi_7\sigma_{30}.\\
				% Oda1979: III.2.6(5), p.114
				&\bar\sigma'_7\sigma_{31}=\bar\nu_7\psi_{15}=\delta_{11}\sigma_{35}=\rho_{13}\omega_{28}=0.
			\end{aligned}\)}\\
		\hhline{>{\arrayrulecolor{odagreen}}->{\arrayrulecolor{odarule}}--}
		% Stem 31, part 4: 4 rows.
		\cellcolor{odagreen}\leavevmode\color{black}\OdaPageStem{31}{4}
		& \cellcolor{odagreenpart}\textup{part~4}
		& \cellcolor{odagreenbody}\leavevmode\label{oda:part:31:4}%
		\OdaRelationCell{\(\begin{aligned}
				% Oda1979: III.2.6(6), p.114
				&\eta_9\sigma_{10}\bar\rho_{17}\equiv0\m 8\Sigma^3\pi_{37}(S^{6}).\\
				% Oda1979: III.2.6(6), p.114
				&\eta_9\sigma_{10}\phi_{17}=\sigma_{10}^2\bar\mu_{24}=0.\\
				% Oda1979: III.2.6(6), p.114
				&\sigma_{10}\delta_{17}=\sigma_{10}\bar\mu_{17}\sigma_{34}=0,\quad \sigma_{10}\bar\sigma'_{17}=\delta_{10}\sigma_{34}.\\
				% Oda1979: III.2.6(7), p.114
				&\sigma_9\Sigma^2\zeta^{*}=\sigma_9\Sigma D_1^{(1)}=0.
			\end{aligned}\)}\\
		\hhline{>{\arrayrulecolor{odagreen}}->{\arrayrulecolor{odarule}}--}
		% Stem 31, part 5: 3 rows.
		\cellcolor{odagreen}\leavevmode\color{black}\OdaPageStem{31}{5}
		& \cellcolor{odagreenpart}\textup{part~5}
		& \cellcolor{odagreenbody}\leavevmode\label{oda:part:31:5}%
		\OdaRelationCell{\(\begin{aligned}
				% Oda1979: III.2.6(7), p.114
				&\sigma_9\tilde\varepsilon_{16}\equiv0\m \sigma_9\bar\sigma'_{16},\;\sigma_9\delta_{16},\;\sigma_9^2\bar\mu_{23}.\\
				% Oda1979: III.2.6(8), p.114
				&\eta_{16}^{*}\bar\varepsilon_{32}=\bar\varepsilon_{10}\omega_{25}=\omega_{17}\bar\varepsilon_{33}=(\Sigma B_1)\bar\varepsilon_{33}=0.\\
				% Oda1979: III.2.6(8), p.114
				&A_1\bar\sigma_{24}=C_1\bar\sigma_{23}=0.
			\end{aligned}\)}\\
		\hhline{>{\arrayrulecolor{odagreen}}->{\arrayrulecolor{odarule}}--}
		% Stem 31, part 6: 4 rows.
		\cellcolor{odagreen}\leavevmode\color{black}\OdaPageStem{31}{6}
		& \cellcolor{odagreenpart}\textup{part~6}
		& \cellcolor{odagreenbody}\leavevmode\label{oda:part:31:6}%
		\OdaRelationCell{\(\begin{aligned}
				% MiM1: group 5; (3.36), p.26
				&\eta_{16}^{*}\bar\varepsilon_{32}=\eta_{16}^{*}\eta_{32}\kappa_{33}=0.\\
				% MiM1: group 14; Lemma 6.2, p.59
				&\eta_{23}\theta'_{24}=\omega_{23}^{*}.\\
				% MiM1: group 14; Lemma 6.2, p.59
				&\theta_{24}\eta_{54}=\omega_{24}^{*}.\\
				% MiM1: group 15; Lemma 7.3(1), p.66
				&[\iota_{18},\kappa_{18}]\equiv\Sigma^3\kappa^{*\prime}\m [\iota_{18},\sigma_{18}^{2}].
			\end{aligned}\)}\\
		\hline
		% Stem 32, part 1: 4 rows.
		\cellcolor{odapurple}\leavevmode\color{black}\OdaPageStem{32}{1}
		& \cellcolor{odapurplepart}\textup{part~1}
		& \cellcolor{odapurplebody}\leavevmode\label{oda:part:32:1}%
		\OdaRelationCell{\(\begin{aligned}
				% MiM1: group 2; Lemma 3.10, p.23
				&\kappa_7\nu_{21}^{*}=0\ \text{or}\ 4\sigma'\xi_{14}\sigma_{32}.\\
				% MiM1: group 2; Lemma 3.10, p.23
				&\kappa_9\nu_{23}^{*}=0.\\
				% MiM1: group 5; (3.38), p.26
				&\bar\mu_3\rho_{20}=\mu_{3,3}\sigma_{28}.\\
				% MiM1: group 6; (3.39), p.27
				&\bar\mu_{13}\rho_{30}=\sigma_{13}\mu_{3,20}=\rho_{13}\bar\mu_{28}.
			\end{aligned}\)}\\
		\hhline{>{\arrayrulecolor{odapurple}}->{\arrayrulecolor{odarule}}--}
		% Stem 32, part 2: 4 rows.
		\cellcolor{odapurple}\leavevmode\color{black}\OdaPageStem{32}{2}
		& \cellcolor{odapurplepart}\textup{part~2}
		& \cellcolor{odapurplebody}\leavevmode\label{oda:part:32:2}%
		\OdaRelationCell{\(\begin{aligned}
				% MiM1: group 6; (3.40), p.27
				&\eta\rho_3=\bar\mu\rho=\mu_3\sigma\ne0\quad\text{(stably)}.\\
				% MiM1: group 7; (3.41), p.28
				&\xi_{12}\kappa_{30}=0.\\
				% MiM1: group 7; Lemma 3.12, p.28
				&\kappa_{10}^{*}\eta_{41}=0,\quad\eta_{11}\kappa_{12}^{*}=0.\\
				% MiM1: group 8; (3.44), p.30
				&\alpha_3^{\mathrm{V}}\eta_{44}=0\ \text{or}\ \sigma_{13}\mu_{3,20}.
			\end{aligned}\)}\\
		\hhline{>{\arrayrulecolor{odapurple}}->{\arrayrulecolor{odarule}}--}
		% Stem 32, part 3: 4 rows.
		\cellcolor{odapurple}\leavevmode\color{black}\OdaPageStem{32}{3}
		& \cellcolor{odapurplepart}\textup{part~3}
		& \cellcolor{odapurplebody}\leavevmode\label{oda:part:32:3}%
		\OdaRelationCell{\(\begin{aligned}
				% MiM1: group 9; (3.45), p.32
				&\eta_{13}\omega_{14}^{**}=\eta_{13}\omega_{14}^{*}.\\
				% MiM1: group 9; Lemma 3.14, p.32
				&\xi_{13}\sigma_{31}^{2}\equiv\eta_{13}\omega_{14}^{*}\m\sigma_{13}\mu_{3,20}.\\
				% MiM1: group 10; Lemma 5.3, p.46
				&2\Sigma\bar\zeta^{*}=0,\quad\pm2\bar\zeta^{*}=[\iota_{14},\bar\zeta_{14}].\\
				% MiM1: group 10; Lemma 5.3, p.46
				&[\iota_{14},\bar\sigma_{14}]=0.
			\end{aligned}\)}\\
		\hhline{>{\arrayrulecolor{odapurple}}->{\arrayrulecolor{odarule}}--}
		% Stem 32, part 4: 4 rows.
		\cellcolor{odapurple}\leavevmode\color{black}\OdaPageStem{32}{4}
		& \cellcolor{odapurplepart}\textup{part~4}
		& \cellcolor{odapurplebody}\leavevmode\label{oda:part:32:4}%
		\OdaRelationCell{\(\begin{aligned}
				% MiM1: group 11; Lemma 5.6, p.51
				&\omega_{14}^{2}=\eta_{14}\omega_{15}^{*}\equiv\xi_{14}\sigma_{32}^{2}\m\sigma_{14}\mu_{3,21}.\\
				% MiM1: group 12; Lemma 5.7, p.53; corrected subscript, see the note after the table
				&\omega_{14}\eta_{30}^{*}=\omega_{14}^{*}\eta_{45}.\\
				% MiM1: group 13; Lemma 6.1, p.57
				&\omega_{20}^{2}=(\eta_{20}^{*})^{2}=0.\\
				% MiM1: group 14; Lemma 6.2, p.59
				&\theta_{24}\eta_{54}^{2}=0,\;(\Sigma^3\kappa^{*\prime})\eta_{49}=0.
			\end{aligned}\)}\\
	\end{longtable}
	\endgroup

	\noindent\textit{Subscripts and source corrections.}
	Four misprints in the cited statements of \cite{Oda1979}
	are corrected here: $\xi_{16}\sigma_{34}$ for
	$\xi_{16}\sigma_{37}$ in Part~I, Proposition~3.5(3), p.~60
	(compare the proof on p.~61); $\xi'\eta_{29}^2$ for
	$\xi'\eta_{30}^2$ in Part~II, Proposition~2.1(7), p.~82
	(compare p.~83); $-2\Sigma\delta''$ for $-2\delta''$ in
	Part~III, Proposition~2.4(4), p.~111 (compare p.~112); and
	$\omega_{15}\sigma_{31}$ for $\omega_{15}\sigma_{23}$ in
	Part~III, Proposition~2.5(5), p.~113. In the last case,
	$\omega_{15}\in\pi_{31}(S^{15})$ forces the subscript $31.$ 
	In \cite[Lemma~5.7, p.~53]{MiM1}, the printed factor $\eta_{16}^{*}$ in
	$\omega_{14}\eta_{16}^{*}$ is read as $\eta_{30}^{*}$: the maps
	$\omega_{14}\colon S^{30}\to S^{14}$ and
	$\eta_{30}^{*}\colon S^{46}\to S^{30}$ compose, and the proof uses
	$\{2\iota_{30},\sigma_{30}^{2},\eta_{44}\}.$  Thus the entry is
	$\omega_{14}\eta_{30}^{*}=\omega_{14}^{*}\eta_{45}$ in $\pi_{46}(S^{14}).$

	\chapter[Homotopy Groups of Spheres]{Homotopy Groups of Spheres in Stems 0--32 and Known Unstable Groups in Stem 33}\label{tlqdzj}
	\begingroup\small
Recall that we work 2-locally unless otherwise stated.	So, $$\pi_{i}(S^{n})\cong\pi_i(S_{\mathrm{int}}^n)\otimes\mathbb Z_{(2)}$$ and
	$\pi_r^S(S^0)=\mathrm{colim}_{k}\pi_{k+r}(S^{k})
	=(\mathrm{colim}_{k}\pi_{k+r}(S^{k}_{\mathrm{int}}))\otimes\mathbb Z_{(2)}$,  the stable $2$-local group. Here, $S_{\mathrm{int}}^n$ means the integral  (unlocalized) $n$-sphere.
	
	An ordered decomposition
	$$\pi_{i}(S^{n})=\langle a_1,\ldots,a_r\rangle\cong A_1\oplus\cdots\oplus A_r,$$
	where each $A_j$ is a cyclic $\mathbb Z_{(2)}$-module,  indicates that 
	each $a_j$ generates 
	the  $A_j$-summand of $\pi_{i}(S^{n}).$

	Note that throughout  this book, the connecting homomorphism $\pi_{\bullet+2}(S^{2n+1})\to \pi_{\bullet}(S^{n})$ in the (2-local)  EHP sequence is denoted by $P,$   rather than Toda's $\Delta.$ 
	For stems $0$--$32,$   stable groups are listed at the end of each stem\uprightsemicolon{} their generators are obtained
	from the range $n\geq k+2$ by suppressing dimension subscripts.
	Family indices, as in $\mu_3,$   $\zeta_3,$   $\rho_3$ and $\eta_4,$   are retained.  
	
	The \textbf{2-localized} homotopy groups are given as follows.
	
	\section{The Toda Stems}
	\subsection{Stems 0--6}\label{app:stem-1}
	\subsubsection*{Stem 0}
	
	\StemEntry{\(\displaystyle \pi_{n}(S^{n})=\langle \e_n\rangle\cong \mathbb Z_{(2)}\quad(n\ge1).\)}{\(\displaystyle \begin{aligned}[t]\pi_{n}(S^{n})&=\langle \e_n\rangle\\&\cong \mathbb Z_{(2)}\quad(n\ge1).\end{aligned}\)}
	\StemEntry{\(\displaystyle \pi^S_{0}(S^0)=\langle \e\rangle\cong \mathbb Z_{(2)}.\)}{\(\displaystyle \pi^S_{0}(S^0)=\langle \e\rangle\cong \mathbb Z_{(2)}.\)}

	\subsubsection*{Stem 1}
	\StemEntry{\(\displaystyle \pi_{3}(S^{2})=\langle \eta_2\rangle\cong \mathbb Z_{(2)}.\)}{\(\displaystyle \begin{aligned}[t]\pi_{3}(S^{2})&=\langle \eta_2\rangle\\&\cong \mathbb Z_{(2)}.\end{aligned}\)}
	\StemEntry{\(\displaystyle \pi_{n+1}(S^{n})=\langle \eta_n\rangle\cong \mathbb Z/2\quad(n\ge3).\)}{\(\displaystyle \begin{aligned}[t]\pi_{n+1}(S^{n})&=\langle \eta_n\rangle\\&\cong \mathbb Z/2\quad(n\ge3).\end{aligned}\)}
	\StemEntry{\(\displaystyle \pi^S_{1}(S^0)=\langle \eta\rangle\cong \mathbb Z/2.\)}{\(\displaystyle \begin{aligned}[t]\pi^S_{1}(S^0)&=\langle \eta\rangle\\&\cong \mathbb Z/2.\end{aligned}\)}
	\subsubsection*{Stem 2}
	\StemEntry{\(\displaystyle \pi_{n+2}(S^{n})=\langle \eta^2_n\rangle\cong \mathbb Z/2\quad(n\ge2).\)}{\(\displaystyle \begin{aligned}[t]\pi_{n+2}(S^{n})&=\langle \eta^2_n\rangle\\&\cong \mathbb Z/2\quad(n\ge2).\end{aligned}\)}
	\StemEntry{\(\displaystyle \pi^S_{2}(S^0)=\langle \eta^2\rangle\cong \mathbb Z/2.\)}{\(\displaystyle \begin{aligned}[t]\pi^S_{2}(S^0)&=\langle \eta^2\rangle\\&\cong \mathbb Z/2.\end{aligned}\)}
	\subsubsection*{Stem 3}
	\StemEntry{\(\displaystyle \pi_{5}(S^{2})=\langle \eta^3_2\rangle\cong \mathbb Z/2.\)}{\(\displaystyle \begin{aligned}[t]\pi_{5}(S^{2})&=\langle \eta^3_2\rangle\\&\cong \mathbb Z/2.\end{aligned}\)}
	\StemEntry{\(\displaystyle \pi_{6}(S^{3})=\langle \nu'\rangle\cong \mathbb Z/4.\)}{\(\displaystyle \begin{aligned}[t]\pi_{6}(S^{3})&=\langle \nu'\rangle\\&\cong \mathbb Z/4.\end{aligned}\)}
	\StemEntry{\(\displaystyle \pi_{7}(S^{4})=\langle \nu_4, \Sigma \nu'\rangle\cong \mathbb Z_{(2)} \oplus \mathbb Z/4.\)}{\(\displaystyle \begin{aligned}[t]\pi_{7}(S^{4})&=\langle \nu_4, \Sigma \nu'\rangle\\&\cong \mathbb Z_{(2)} \oplus \mathbb Z/4.\end{aligned}\)}
	\StemEntry{\(\displaystyle \pi_{n+3}(S^{n})=\langle \nu_n\rangle\cong \mathbb Z/8\quad(n\ge5).\)}{\(\displaystyle \begin{aligned}[t]\pi_{n+3}(S^{n})&=\langle \nu_n\rangle\\&\cong \mathbb Z/8\quad(n\ge5).\end{aligned}\)}
	\StemEntry{\(\displaystyle \pi^S_{3}(S^0)=\langle \nu\rangle\cong \mathbb Z/8.\)}{\(\displaystyle \begin{aligned}[t]\pi^S_{3}(S^0)&=\langle \nu\rangle\\&\cong \mathbb Z/8.\end{aligned}\)}
	\subsubsection*{Stem 4}
	\StemEntry{\(\displaystyle \pi_{6}(S^{2})=\langle \eta_2 \nu'\rangle\cong \mathbb Z/4.\)}{\(\displaystyle \begin{aligned}[t]\pi_{6}(S^{2})&=\langle \eta_2 \nu'\rangle\\&\cong \mathbb Z/4.\end{aligned}\)}
	\StemEntry{\(\displaystyle \pi_{7}(S^{3})=\langle \nu' \eta_6\rangle\cong \mathbb Z/2.\)}{\(\displaystyle \begin{aligned}[t]\pi_{7}(S^{3})&=\langle \nu' \eta_6\rangle\\&\cong \mathbb Z/2.\end{aligned}\)}
	\StemEntry{\(\displaystyle \pi_{8}(S^{4})=\langle \nu_4 \eta_7, (\Sigma \nu') \eta_7\rangle\cong (\mathbb Z/2)^{2}.\)}{\(\displaystyle \begin{aligned}[t]\pi_{8}(S^{4})&=\langle \nu_4 \eta_7, (\Sigma \nu') \eta_7\rangle\\&\cong (\mathbb Z/2)^{2}.\end{aligned}\)}
	\StemEntry{\(\displaystyle \pi_{9}(S^{5})=\langle \nu_5 \eta_8\rangle\cong \mathbb Z/2.\)}{\(\displaystyle \begin{aligned}[t]\pi_{9}(S^{5})&=\langle \nu_5 \eta_8\rangle\\&\cong \mathbb Z/2.\end{aligned}\)}
	\StemEntry{\(\displaystyle \pi_{n+4}(S^{n})=0\quad(n\ge6).\)}{\(\displaystyle \begin{aligned}[t]\pi_{n+4}(S^{n})&=0\quad(n\ge6).\end{aligned}\)}
	\StemEntry{\(\displaystyle \pi^S_{4}(S^0)=0.\)}{\(\displaystyle \pi^S_{4}(S^0)=0.\)}
	
	\subsubsection*{Stem 5}
	\StemEntry{\(\displaystyle \pi_{7}(S^{2})=\langle \eta_2 \nu' \eta_6\rangle\cong \mathbb Z/2.\)}{\(\displaystyle \begin{aligned}[t]\pi_{7}(S^{2})&=\langle \eta_2 \nu' \eta_6\rangle\\&\cong \mathbb Z/2.\end{aligned}\)}
	\StemEntry{\(\displaystyle \pi_{8}(S^{3})=\langle \nu' \eta^2_6\rangle\cong \mathbb Z/2.\)}{\(\displaystyle \begin{aligned}[t]\pi_{8}(S^{3})&=\langle \nu' \eta^2_6\rangle\\&\cong \mathbb Z/2.\end{aligned}\)}
	\StemEntry{\(\displaystyle \pi_{9}(S^{4})=\langle \nu_4 \eta^2_7, (\Sigma \nu') \eta^2_7\rangle\cong (\mathbb Z/2)^{2}.\)}{\(\displaystyle \begin{aligned}[t]\pi_{9}(S^{4})&=\langle \nu_4 \eta^2_7, (\Sigma \nu') \eta^2_7\rangle\\&\cong (\mathbb Z/2)^{2}.\end{aligned}\)}
	\StemEntry{\(\displaystyle \pi_{10}(S^{5})=\langle \nu_5 \eta^2_8\rangle\cong \mathbb Z/2.\)}{\(\displaystyle \begin{aligned}[t]\pi_{10}(S^{5})&=\langle \nu_5 \eta^2_8\rangle\\&\cong \mathbb Z/2.\end{aligned}\)}
	\StemEntry{\(\displaystyle \pi_{11}(S^{6})=\langle P(\iota_{13})\rangle\cong \mathbb Z_{(2)}.\)}{\(\displaystyle \begin{aligned}[t]\pi_{11}(S^{6})&=\langle P(\iota_{13})\rangle\\&\cong \mathbb Z_{(2)}.\end{aligned}\)}
	\StemEntry{\(\displaystyle \pi_{n+5}(S^{n})=0\quad(n\ge7).\)}{\(\displaystyle \begin{aligned}[t]\pi_{n+5}(S^{n})&=0\quad(n\ge7).\end{aligned}\)}
	\StemEntry{\(\displaystyle \pi^S_{5}(S^0)=0.\)}{\(\displaystyle \pi^S_{5}(S^0)=0.\)}
	\subsubsection*{Stem 6}
	\StemEntry{\(\displaystyle \pi_{8}(S^{2})=\langle \eta_2 \nu' \eta^2_6\rangle\cong \mathbb Z/2.\)}{\(\displaystyle \begin{aligned}[t]\pi_{8}(S^{2})&=\langle \eta_2 \nu' \eta^2_6\rangle\\&\cong \mathbb Z/2.\end{aligned}\)}
	\StemEntry{\(\displaystyle \pi_{9}(S^{3})=0.\)}{\(\displaystyle \begin{aligned}[t]\pi_{9}(S^{3})&=0.\end{aligned}\)}
	\StemEntry{\(\displaystyle \pi_{10}(S^{4})=\langle \nu^2_4\rangle\cong \mathbb Z/8.\)}{\(\displaystyle \begin{aligned}[t]\pi_{10}(S^{4})&=\langle \nu^2_4\rangle\\&\cong \mathbb Z/8.\end{aligned}\)}
	\StemEntry{\(\displaystyle \pi_{n+6}(S^{n})=\langle \nu^2_n\rangle\cong \mathbb Z/2\quad(n\ge5).\)}{\(\displaystyle \begin{aligned}[t]\pi_{n+6}(S^{n})&=\langle \nu^2_n\rangle\\&\cong \mathbb Z/2\quad(n\ge5).\end{aligned}\)}
	\StemEntry{\(\displaystyle \pi^S_{6}(S^0)=\langle \nu^2\rangle\cong \mathbb Z/2.\)}{\(\displaystyle \begin{aligned}[t]\pi^S_{6}(S^0)&=\langle \nu^2\rangle\\&\cong \mathbb Z/2.\end{aligned}\)}
	
	\subsection{Stems 7--9}\label{app:stem-7}
	\subsubsection*{Stem 7}
	\StemEntry{\(\displaystyle \pi_{9}(S^{2})=0.\)}{\(\displaystyle \begin{aligned}[t]\pi_{9}(S^{2})&=0.\end{aligned}\)}
	\StemEntry{\(\displaystyle \pi_{10}(S^{3})=0.\)}{\(\displaystyle \begin{aligned}[t]\pi_{10}(S^{3})&=0.\end{aligned}\)}
	\StemEntry{\(\displaystyle \pi_{11}(S^{4})=0.\)}{\(\displaystyle \begin{aligned}[t]\pi_{11}(S^{4})&=0.\end{aligned}\)}
	\StemEntry{\(\displaystyle \pi_{12}(S^{5})=\langle \sigma'''\rangle\cong \mathbb Z/2.\)}{\(\displaystyle \begin{aligned}[t]\pi_{12}(S^{5})&=\langle \sigma'''\rangle\\&\cong \mathbb Z/2.\end{aligned}\)}
	\StemEntry{\(\displaystyle \pi_{13}(S^{6})=\langle \sigma''\rangle\cong \mathbb Z/4.\)}{\(\displaystyle \begin{aligned}[t]\pi_{13}(S^{6})&=\langle \sigma''\rangle\\&\cong \mathbb Z/4.\end{aligned}\)}
	\StemEntry{\(\displaystyle \pi_{14}(S^{7})=\langle \sigma'\rangle\cong \mathbb Z/8.\)}{\(\displaystyle \begin{aligned}[t]\pi_{14}(S^{7})&=\langle \sigma'\rangle\\&\cong \mathbb Z/8.\end{aligned}\)}
	\StemEntry{\(\displaystyle \pi_{15}(S^{8})=\langle \sigma_8, \Sigma \sigma'\rangle\cong \mathbb Z_{(2)} \oplus \mathbb Z/8.\)}{\(\displaystyle \begin{aligned}[t]\pi_{15}(S^{8})&=\langle \sigma_8, \Sigma \sigma'\rangle\\&\cong \mathbb Z_{(2)} \oplus \mathbb Z/8.\end{aligned}\)}
	\StemEntry{\(\displaystyle \pi_{n+7}(S^{n})=\langle \sigma_n\rangle\cong \mathbb Z/16\quad(n\ge9).\)}{\(\displaystyle \begin{aligned}[t]\pi_{n+7}(S^{n})&=\langle \sigma_n\rangle\\&\cong \mathbb Z/16\quad(n\ge9).\end{aligned}\)}
	\StemEntry{\(\displaystyle \pi^S_{7}(S^0)=\langle \sigma\rangle\cong \mathbb Z/16.\)}{\(\displaystyle \begin{aligned}[t]\pi^S_{7}(S^0)&=\langle \sigma\rangle\\&\cong \mathbb Z/16.\end{aligned}\)}
	
	\subsubsection*{Stem 8}
	\StemEntry{\(\displaystyle \pi_{10}(S^{2})=0.\)}{\(\displaystyle \begin{aligned}[t]\pi_{10}(S^{2})&=0.\end{aligned}\)}
	\StemEntry{\(\displaystyle \pi_{n+8}(S^{n})=\langle \varepsilon _n\rangle\cong \mathbb Z/2\quad(n = 3, 4, 5).\)}{\(\displaystyle \begin{aligned}[t]\pi_{n+8}(S^{n})&=\langle \varepsilon _n\rangle\\&\cong \mathbb Z/2\quad(n = 3, 4, 5).\end{aligned}\)}
	\StemEntry{\(\displaystyle \pi_{14}(S^{6})=\langle \bar\nu_6, \varepsilon _6\rangle\cong \mathbb Z/8 \oplus \mathbb Z/2.\)}{\(\displaystyle \begin{aligned}[t]\pi_{14}(S^{6})&=\langle \bar\nu_6, \varepsilon _6\rangle\\&\cong \mathbb Z/8 \oplus \mathbb Z/2.\end{aligned}\)}
	\StemEntry{\(\displaystyle \pi_{15}(S^{7})=\langle \sigma' \eta_{14}, \bar\nu_7, \varepsilon _7\rangle\cong (\mathbb Z/2)^{3}.\)}{\(\displaystyle \begin{aligned}[t]\pi_{15}(S^{7})&=\langle \sigma' \eta_{14}, \bar\nu_7, \varepsilon _7\rangle\\&\cong (\mathbb Z/2)^{3}.\end{aligned}\)}
	\StemEntry{\(\displaystyle \pi_{16}(S^{8})=\langle \sigma_8 \eta_{15}, (\Sigma \sigma') \eta_{15}, \bar\nu_8, \varepsilon _8\rangle\cong (\mathbb Z/2)^{4}.\)}{\(\displaystyle \begin{aligned}[t]\pi_{16}(S^{8})&=\langle \sigma_8 \eta_{15}, (\Sigma \sigma') \eta_{15}, \bar\nu_8, \varepsilon _8\rangle\\&\cong (\mathbb Z/2)^{4}.\end{aligned}\)}
	\StemEntry{\(\displaystyle \pi_{17}(S^{9})=\langle \sigma_9 \eta_{16}, \bar\nu_9, \varepsilon _9\rangle\cong (\mathbb Z/2)^{3}.\)}{\(\displaystyle \begin{aligned}[t]\pi_{17}(S^{9})&=\langle \sigma_9 \eta_{16}, \bar\nu_9, \varepsilon _9\rangle\\&\cong (\mathbb Z/2)^{3}.\end{aligned}\)}
	\StemEntry{\(\displaystyle \pi_{n+8}(S^{n})=\langle \bar\nu_n, \varepsilon _n\rangle\cong (\mathbb Z/2)^{2}\quad(n \ge 10).\)}{\(\displaystyle \begin{aligned}[t]\pi_{n+8}(S^{n})&=\langle \bar\nu_n, \varepsilon _n\rangle\\&\cong (\mathbb Z/2)^{2}\quad(n \ge 10).\end{aligned}\)}
	\StemEntry{\(\displaystyle \pi^S_{8}(S^0)=\langle \bar\nu, \varepsilon\rangle\cong \mathbb Z/2 \oplus \mathbb Z/2.\)}{\(\displaystyle \begin{aligned}[t]\pi^S_{8}(S^0)&=\langle \bar\nu, \varepsilon\rangle\\&\cong \mathbb Z/2 \oplus \mathbb Z/2.\end{aligned}\)}
	
	\subsubsection*{Stem 9}
	\StemEntry{\(\displaystyle \pi_{11}(S^{2})=\langle \eta_2 \varepsilon _3\rangle\cong \mathbb Z/2.\)}{\(\displaystyle \begin{aligned}[t]\pi_{11}(S^{2})&=\langle \eta_2 \varepsilon _3\rangle\\&\cong \mathbb Z/2.\end{aligned}\)}
	\StemEntry{\(\displaystyle \pi_{12}(S^{3})=\langle \mu_3, \eta_3 \varepsilon _4\rangle\cong (\mathbb Z/2)^{2}.\)}{\(\displaystyle \begin{aligned}[t]\pi_{12}(S^{3})&=\langle \mu_3, \eta_3 \varepsilon _4\rangle\\&\cong (\mathbb Z/2)^{2}.\end{aligned}\)}
	\StemEntry{\(\displaystyle \pi_{n+9}(S^{n})=\langle \nu^3_n, \mu_n, \eta_n \varepsilon _{n+1}\rangle\cong (\mathbb Z/2)^{3}\quad(n = 4, 5, 6).\)}{\(\displaystyle \begin{aligned}[t]\pi_{n+9}(S^{n})&=\langle \nu^3_n, \mu_n, \eta_n \varepsilon _{n+1}\rangle\\&\cong (\mathbb Z/2)^{3}\quad(n = 4, 5, 6).\end{aligned}\)}
	\StemEntry{\(\displaystyle \pi_{16}(S^{7})=\langle \sigma' \eta^2_{14}, \nu^3_7, \mu_7, \eta_7 \varepsilon _8\rangle\cong (\mathbb Z/2)^{4}.\)}{\(\displaystyle \begin{aligned}[t]\pi_{16}(S^{7})&=\langle \sigma' \eta^2_{14}, \nu^3_7, \mu_7, \eta_7 \varepsilon _8\rangle\\&\cong (\mathbb Z/2)^{4}.\end{aligned}\)}
	\StemEntry{\(\displaystyle \pi_{17}(S^{8})=\langle \sigma_8 \eta^2_{15}, (\Sigma \sigma') \eta^2_{15}, \nu^3_8, \mu_8, \eta_8 \varepsilon _9\rangle\cong (\mathbb Z/2)^{5}.\)}{\(\displaystyle \begin{aligned}[t]\pi_{17}(S^{8})&=\langle \sigma_8 \eta^2_{15}, (\Sigma \sigma') \eta^2_{15}, \nu^3_8, \mu_8, \eta_8 \varepsilon _9\rangle\\&\cong (\mathbb Z/2)^{5}.\end{aligned}\)}
	\StemEntry{\(\displaystyle \pi_{18}(S^{9})=\langle \sigma_9 \eta^2_{16}, \nu^3_9, \mu_9, \eta_9 \varepsilon _{10}\rangle\cong (\mathbb Z/2)^{4}.\)}{\(\displaystyle \begin{aligned}[t]\pi_{18}(S^{9})&=\langle \sigma_9 \eta^2_{16}, \nu^3_9, \mu_9, \eta_9 \varepsilon _{10}\rangle\\&\cong (\mathbb Z/2)^{4}.\end{aligned}\)}
	\StemEntry{\(\displaystyle \pi_{19}(S^{10})=\langle P(\iota_{21}), \nu^3_{10}, \mu_{10}, \eta_{10} \varepsilon _{11}\rangle\cong \mathbb Z_{(2)} \oplus (\mathbb Z/2)^{3}.\)}{\(\displaystyle \begin{aligned}[t]\pi_{19}(S^{10})&=\langle P(\iota_{21}), \nu^3_{10}, \mu_{10}, \eta_{10} \varepsilon _{11}\rangle\\&\cong \mathbb Z_{(2)} \oplus (\mathbb Z/2)^{3}.\end{aligned}\)}
	\StemEntry{\(\displaystyle \pi_{n+9}(S^{n})=\langle \nu^3_n, \mu_n, \eta_n \varepsilon _{n+1}\rangle\cong (\mathbb Z/2)^{3}\quad(n\ge 11).\)}{\(\displaystyle \begin{aligned}[t]\pi_{n+9}(S^{n})&=\langle \nu^3_n, \mu_n, \eta_n \varepsilon _{n+1}\rangle\\&\cong (\mathbb Z/2)^{3}\quad(n\ge 11).\end{aligned}\)}
	\StemEntry{\(\displaystyle \pi^S_{9}(S^0)=\langle \nu^3, \mu, \eta \varepsilon\rangle\cong \mathbb Z/2 \oplus \mathbb Z/2 \oplus \mathbb Z/2.\)}{\(\displaystyle \begin{aligned}[t]\pi^S_{9}(S^0)&=\langle \nu^3, \mu, \eta \varepsilon\rangle\\&\cong \mathbb Z/2 \oplus \mathbb Z/2 \oplus \mathbb Z/2.\end{aligned}\)}
	
	\subsection{Stems 10--11}\label{app:stem-10}
	\subsubsection*{Stem 10}
	\StemEntry{\(\displaystyle \pi_{12}(S^{2})=\langle \eta^2_2 \varepsilon _4, \eta_2 \mu_3\rangle\cong (\mathbb Z/2)^{2}.\)}{\(\displaystyle \begin{aligned}[t]\pi_{12}(S^{2})&=\langle \eta^2_2 \varepsilon _4, \eta_2 \mu_3\rangle\\&\cong (\mathbb Z/2)^{2}.\end{aligned}\)}
	\StemEntry{\(\displaystyle \pi_{13}(S^{3})=\langle \varepsilon ', \eta_3 \mu_4\rangle\cong \mathbb Z/4 \oplus \mathbb Z/2.\)}{\(\displaystyle \begin{aligned}[t]\pi_{13}(S^{3})&=\langle \varepsilon ', \eta_3 \mu_4\rangle\\&\cong \mathbb Z/4 \oplus \mathbb Z/2.\end{aligned}\)}
	\StemEntry{\(\displaystyle \pi_{14}(S^{4})=\langle \nu_4 \sigma', \Sigma \varepsilon ', \eta_4 \mu_5\rangle\cong \mathbb Z/8 \oplus \mathbb Z/4 \oplus \mathbb Z/2.\)}{\(\displaystyle \begin{aligned}[t]\pi_{14}(S^{4})&=\langle \nu_4 \sigma', \Sigma \varepsilon ', \eta_4 \mu_5\rangle\\&\cong \mathbb Z/8 \oplus \mathbb Z/4 \oplus \mathbb Z/2.\end{aligned}\)}
	\StemEntry{\(\displaystyle \pi_{n+10}(S^{n})=\langle \nu_n \sigma_{n+3}, \eta_n \mu_{n+1}\rangle\cong \mathbb Z/8 \oplus \mathbb Z/2\quad(n = 5, 6 \text{ or } 7).\)}{\(\displaystyle \begin{aligned}[t]\pi_{n+10}(S^{n})&=\langle \nu_n \sigma_{n+3}, \eta_n \mu_{n+1}\rangle\\&\cong \mathbb Z/8 \oplus \mathbb Z/2\quad(n = 5, 6 \text{ or } 7).\end{aligned}\)}
	\StemEntry{\(\displaystyle \pi_{18}(S^{8})=\langle \sigma_8 \nu_{15}, \nu_8 \sigma_{11}, \eta_8 \mu_9\rangle\cong \mathbb Z/8 \oplus \mathbb Z/8 \oplus \mathbb Z/2.\)}{\(\displaystyle \begin{aligned}[t]\pi_{18}(S^{8})&=\langle \sigma_8 \nu_{15}, \nu_8 \sigma_{11}, \eta_8 \mu_9\rangle\\&\cong \mathbb Z/8 \oplus \mathbb Z/8 \oplus \mathbb Z/2.\end{aligned}\)}
	\StemEntry{\(\displaystyle \pi_{19}(S^{9})=\langle \sigma_9 \nu_{16}, \eta_9 \mu_{10}\rangle\cong \mathbb Z/8 \oplus \mathbb Z/2.\)}{\(\displaystyle \begin{aligned}[t]\pi_{19}(S^{9})&=\langle \sigma_9 \nu_{16}, \eta_9 \mu_{10}\rangle\\&\cong \mathbb Z/8 \oplus \mathbb Z/2.\end{aligned}\)}
	\StemEntry{\(\displaystyle \pi_{20}(S^{10})=\langle \sigma_{10} \nu_{17}, \eta_{10} \mu_{11}\rangle\cong \mathbb Z/4 \oplus \mathbb Z/2.\)}{\(\displaystyle \begin{aligned}[t]\pi_{20}(S^{10})&=\langle \sigma_{10} \nu_{17}, \eta_{10} \mu_{11}\rangle\\&\cong \mathbb Z/4 \oplus \mathbb Z/2.\end{aligned}\)}
	\StemEntry{\(\displaystyle \pi_{21}(S^{11})=\langle \sigma_{11} \nu_{18}, \eta_{11} \mu_{12}\rangle\cong (\mathbb Z/2)^{2}.\)}{\(\displaystyle \begin{aligned}[t]\pi_{21}(S^{11})&=\langle \sigma_{11} \nu_{18}, \eta_{11} \mu_{12}\rangle\\&\cong (\mathbb Z/2)^{2}.\end{aligned}\)}
	\StemEntry{\(\displaystyle \pi_{n+10}(S^{n})=\langle \eta_n \mu_{n+1}\rangle\cong \mathbb Z/2\quad(n\ge 12).\)}{\(\displaystyle \begin{aligned}[t]\pi_{n+10}(S^{n})&=\langle \eta_n \mu_{n+1}\rangle\\&\cong \mathbb Z/2\quad(n\ge 12).\end{aligned}\)}
	\StemEntry{\(\displaystyle \pi^S_{10}(S^0)=\langle \eta \mu\rangle\cong \mathbb Z/2.\)}{\(\displaystyle \begin{aligned}[t]\pi^S_{10}(S^0)&=\langle \eta \mu\rangle\\&\cong \mathbb Z/2.\end{aligned}\)}
	
	\subsubsection*{Stem 11}
	\StemEntry{\(\displaystyle \pi_{13}(S^{2})=\langle \eta_2 \varepsilon ', \eta^2_2 \mu_4\rangle\cong \mathbb Z/4 \oplus \mathbb Z/2.\)}{\(\displaystyle \begin{aligned}[t]\pi_{13}(S^{2})&=\langle \eta_2 \varepsilon ', \eta^2_2 \mu_4\rangle\\&\cong \mathbb Z/4 \oplus \mathbb Z/2.\end{aligned}\)}
	\StemEntry{\(\displaystyle \pi_{14}(S^{3})=\langle \mu', \varepsilon _3 \nu_{11}, \nu' \varepsilon _6\rangle\cong \mathbb Z/4 \oplus (\mathbb Z/2)^{2}.\)}{\(\displaystyle \begin{aligned}[t]\pi_{14}(S^{3})&=\langle \mu', \varepsilon _3 \nu_{11}, \nu' \varepsilon _6\rangle\\&\cong \mathbb Z/4 \oplus (\mathbb Z/2)^{2}.\end{aligned}\)}
	\StemEntry{\(\displaystyle \pi_{15}(S^{4})=\langle \Sigma \mu', \nu_4 \sigma' \eta_{14}, \nu_4 \bar\nu_7, \nu_4 \varepsilon _7, \varepsilon _4 \nu_{12}, (\Sigma \nu') \varepsilon _7\rangle\cong \mathbb Z/4 \oplus (\mathbb Z/2)^{5}.\)}{\(\displaystyle \begin{aligned}[t]\pi_{15}(S^{4})&=\langle \Sigma \mu', \nu_4 \sigma' \eta_{14}, \nu_4 \bar\nu_7, \nu_4 \varepsilon _7, \varepsilon _4 \nu_{12}, (\Sigma \nu') \varepsilon _7\rangle\\&\cong \mathbb Z/4 \oplus (\mathbb Z/2)^{5}.\end{aligned}\)}
	\StemEntry{\(\displaystyle \pi_{16}(S^{5})=\langle \zeta_5, \nu_5 \bar\nu_8, \nu_5 \varepsilon _8\rangle\cong \mathbb Z/8 \oplus (\mathbb Z/2)^{2}.\)}{\(\displaystyle \begin{aligned}[t]\pi_{16}(S^{5})&=\langle \zeta_5, \nu_5 \bar\nu_8, \nu_5 \varepsilon _8\rangle\\&\cong \mathbb Z/8 \oplus (\mathbb Z/2)^{2}.\end{aligned}\)}
	\StemEntry{\(\displaystyle \pi_{17}(S^{6})=\langle \zeta_6, \bar\nu_6 \nu_{14}\rangle\cong \mathbb Z/8 \oplus \mathbb Z/4.\)}{\(\displaystyle \begin{aligned}[t]\pi_{17}(S^{6})&=\langle \zeta_6, \bar\nu_6 \nu_{14}\rangle\\&\cong \mathbb Z/8 \oplus \mathbb Z/4.\end{aligned}\)}
	\StemEntry{\(\displaystyle \pi_{n+11}(S^{n})=\langle \zeta_n, \bar\nu_n \nu_{n+8}\rangle\cong \mathbb Z/8 \oplus \mathbb Z/2\quad(n = 7, 8 \text{ or } 9).\)}{\(\displaystyle \begin{aligned}[t]\pi_{n+11}(S^{n})&=\langle \zeta_n, \bar\nu_n \nu_{n+8}\rangle\\&\cong \mathbb Z/8 \oplus \mathbb Z/2\quad(n = 7, 8 \text{ or } 9).\end{aligned}\)}
	\StemEntry{\(\displaystyle \pi_{n+11}(S^{n})=\langle \zeta_n\rangle\cong \mathbb Z/8\quad(n = 10, 11 \text{ or } n\ge 13).\)}{\(\displaystyle \begin{aligned}[t]\pi_{n+11}(S^{n})&=\langle \zeta_n\rangle\\&\cong \mathbb Z/8\quad(n = 10, 11 \text{ or } n\ge 13).\end{aligned}\)}
	\StemEntry{\(\displaystyle \pi_{23}(S^{12})=\langle P(\iota_{25}), \zeta_{12}\rangle\cong \mathbb Z_{(2)} \oplus \mathbb Z/8.\)}{\(\displaystyle \begin{aligned}[t]\pi_{23}(S^{12})&=\langle P(\iota_{25}), \zeta_{12}\rangle\\&\cong \mathbb Z_{(2)} \oplus \mathbb Z/8.\end{aligned}\)}
	\StemEntry{\(\displaystyle \pi^S_{11}(S^0)=\langle \zeta\rangle\cong \mathbb Z/8.\)}{\(\displaystyle \begin{aligned}[t]\pi^S_{11}(S^0)&=\langle \zeta\rangle\\&\cong \mathbb Z/8.\end{aligned}\)}
	
	\subsection{Stems 12--13}\label{app:stem-12}
	\subsubsection*{Stem 12}
	\StemEntry{\(\displaystyle \pi_{14}(S^{2})=\langle \eta_2 \mu', \eta_2 \nu' \bar\nu_6, \eta_2 \nu' \varepsilon _6\rangle\cong \mathbb Z/4 \oplus (\mathbb Z/2)^{2}.\)}{\(\displaystyle \begin{aligned}[t]\pi_{14}(S^{2})&=\langle \eta_2 \mu', \eta_2 \nu' \bar\nu_6, \eta_2 \nu' \varepsilon _6\rangle\\&\cong \mathbb Z/4 \oplus (\mathbb Z/2)^{2}.\end{aligned}\)}
	\StemEntry{\(\displaystyle \pi_{15}(S^{3})=\langle \nu' \mu_6, \nu' \eta_6 \varepsilon _7\rangle\cong (\mathbb Z/2)^{2}.\)}{\(\displaystyle \begin{aligned}[t]\pi_{15}(S^{3})&=\langle \nu' \mu_6, \nu' \eta_6 \varepsilon _7\rangle\\&\cong (\mathbb Z/2)^{2}.\end{aligned}\)}
	\StemEntry{\(\displaystyle \pi_{16}(S^{4})=\langle \nu_4 \sigma' \eta^2_{14}, \nu^4_4, \nu_4 \mu_7, \nu_4 \eta_7 \varepsilon _8, (\Sigma \nu') \mu_7, (\Sigma \nu') \eta_7 \varepsilon _8\rangle\cong (\mathbb Z/2)^{6}.\)}{\(\displaystyle \begin{aligned}[t]\pi_{16}(S^{4})&=\langle \nu_4 \sigma' \eta^2_{14}, \nu^4_4, \nu_4 \mu_7, \nu_4 \eta_7 \varepsilon _8, (\Sigma \nu') \mu_7,\\&\phantom{=\langle{}}(\Sigma \nu') \eta_7 \varepsilon _8\rangle\cong (\mathbb Z/2)^{6}.\end{aligned}\)}
	\StemEntry{\(\displaystyle \pi_{17}(S^{5})=\langle \nu^4_5, \nu_5 \mu_8, \nu_5 \eta_8 \varepsilon _9\rangle\cong (\mathbb Z/2)^{3}.\)}{\(\displaystyle \begin{aligned}[t]\pi_{17}(S^{5})&=\langle \nu^4_5, \nu_5 \mu_8, \nu_5 \eta_8 \varepsilon _9\rangle\\&\cong (\mathbb Z/2)^{3}.\end{aligned}\)}
	\StemEntry{\(\displaystyle \pi_{18}(S^{6})=\langle P(\sigma_{13})\rangle\cong \mathbb Z/16.\)}{\(\displaystyle \begin{aligned}[t]\pi_{18}(S^{6})&=\langle P(\sigma_{13})\rangle\\&\cong \mathbb Z/16.\end{aligned}\)}
	\StemEntry{\(\displaystyle \pi_{n+12}(S^{n})=0\quad(n = 7, 8 \text{ or } 9).\)}{\(\displaystyle \begin{aligned}[t]\pi_{n+12}(S^{n})&=0\quad(n = 7, 8 \text{ or } 9).\end{aligned}\)}
	\StemEntry{\(\displaystyle \pi_{22}(S^{10})=\langle [\e_{10},\nu_{10}]\rangle\cong \mathbb Z/4.\)}{\(\displaystyle \begin{aligned}[t]\pi_{22}(S^{10})&=\langle [\e_{10},\nu_{10}]\rangle\\&\cong \mathbb Z/4.\end{aligned}\)}
	\StemEntry{\(\displaystyle \pi_{23}(S^{11})=\langle \theta'\rangle\cong \mathbb Z/2.\)}{\(\displaystyle \begin{aligned}[t]\pi_{23}(S^{11})&=\langle \theta'\rangle\\&\cong \mathbb Z/2.\end{aligned}\)}
	\StemEntry{\(\displaystyle \pi_{24}(S^{12})=\langle \theta, \Sigma \theta'\rangle\cong (\mathbb Z/2)^{2}.\)}{\(\displaystyle \begin{aligned}[t]\pi_{24}(S^{12})&=\langle \theta, \Sigma \theta'\rangle\\&\cong (\mathbb Z/2)^{2}.\end{aligned}\)}
	\StemEntry{\(\displaystyle \pi_{25}(S^{13})=\langle \Sigma \theta\rangle\cong \mathbb Z/2.\)}{\(\displaystyle \begin{aligned}[t]\pi_{25}(S^{13})&=\langle \Sigma \theta\rangle\\&\cong \mathbb Z/2.\end{aligned}\)}
	\StemEntry{\(\displaystyle \pi_{n+12}(S^{n})=0\quad(n\ge 14).\)}{\(\displaystyle \begin{aligned}[t]\pi_{n+12}(S^{n})&=0\quad(n\ge 14).\end{aligned}\)}
	\StemEntry{\(\displaystyle \pi^S_{12}(S^0)=0.\)}{\(\displaystyle \pi^S_{12}(S^0)=0.\)}
	\subsubsection*{Stem 13}
	\StemEntry{\(\displaystyle \pi_{15}(S^{2})=\langle \eta_2 \nu' \mu_6, \eta_2 \nu' \eta_6 \varepsilon _7\rangle\cong (\mathbb Z/2)^{2}.\)}{\(\displaystyle \begin{aligned}[t]\pi_{15}(S^{2})&=\langle \eta_2 \nu' \mu_6, \eta_2 \nu' \eta_6 \varepsilon _7\rangle\\&\cong (\mathbb Z/2)^{2}.\end{aligned}\)}
	\StemEntry{\(\displaystyle \pi_{16}(S^{3})=\langle \nu' \eta_6 \mu_7\rangle\cong \mathbb Z/2.\)}{\(\displaystyle \begin{aligned}[t]\pi_{16}(S^{3})&=\langle \nu' \eta_6 \mu_7\rangle\\&\cong \mathbb Z/2.\end{aligned}\)}
	\StemEntry{\(\displaystyle \pi_{17}(S^{4})=\langle \nu^2_4\sigma_{10}, \nu_4 \eta_7 \mu_8, (\Sigma \nu') \eta_7 \mu_8\rangle\cong \mathbb Z/8 \oplus (\mathbb Z/2)^{2}.\)}{\(\displaystyle \begin{aligned}[t]\pi_{17}(S^{4})&=\langle \nu^2_4\sigma_{10}, \nu_4 \eta_7 \mu_8, (\Sigma \nu') \eta_7 \mu_8\rangle\\&\cong \mathbb Z/8 \oplus (\mathbb Z/2)^{2}.\end{aligned}\)}
	\StemEntry{\(\displaystyle \pi_{18}(S^{5})=\langle \nu_5 \sigma_8 \nu_{15}, \nu_5 \eta_8 \mu_9\rangle\cong (\mathbb Z/2)^{2}.\)}{\(\displaystyle \begin{aligned}[t]\pi_{18}(S^{5})&=\langle \nu_5 \sigma_8 \nu_{15}, \nu_5 \eta_8 \mu_9\rangle\\&\cong (\mathbb Z/2)^{2}.\end{aligned}\)}
	\StemEntry{\(\displaystyle \pi_{n+13}(S^{n})=\langle \nu_n \sigma_{n+3} \nu_{n+10}\rangle\cong \mathbb Z/2\quad(n = 6 \text{ or } 7).\)}{\(\displaystyle \begin{aligned}[t]\pi_{n+13}(S^{n})&=\langle \nu_n \sigma_{n+3} \nu_{n+10}\rangle\\&\cong \mathbb Z/2\quad(n = 6 \text{ or } 7).\end{aligned}\)}
	\StemEntry{\(\displaystyle \pi_{21}(S^{8})=\langle \sigma_8 \nu^2_{15}, \nu_8 \sigma_{11} \nu_{18}\rangle\cong (\mathbb Z/2)^{2}.\)}{\(\displaystyle \begin{aligned}[t]\pi_{21}(S^{8})&=\langle \sigma_8 \nu^2_{15}, \nu_8 \sigma_{11} \nu_{18}\rangle\\&\cong (\mathbb Z/2)^{2}.\end{aligned}\)}
	\StemEntry{\(\displaystyle \pi_{n+13}(S^{n})=\langle \sigma_n \nu^2_{n+7}\rangle\cong \mathbb Z/2\quad(n = 9 \text{ or } 10).\)}{\(\displaystyle \begin{aligned}[t]\pi_{n+13}(S^{n})&=\langle \sigma_n \nu^2_{n+7}\rangle\\&\cong \mathbb Z/2\quad(n = 9 \text{ or } 10).\end{aligned}\)}
	\StemEntry{\(\displaystyle \pi_{24}(S^{11})=\langle \theta' \eta_{23}, \sigma_{11} \nu^2_{18}\rangle\cong (\mathbb Z/2)^{2}.\)}{\(\displaystyle \begin{aligned}[t]\pi_{24}(S^{11})&=\langle \theta' \eta_{23}, \sigma_{11} \nu^2_{18}\rangle\\&\cong (\mathbb Z/2)^{2}.\end{aligned}\)}
	\StemEntry{\(\displaystyle \pi_{25}(S^{12})=\langle \theta \eta_{24}, (\Sigma \theta') \eta_{24}\rangle\cong (\mathbb Z/2)^{2}.\)}{\(\displaystyle \begin{aligned}[t]\pi_{25}(S^{12})&=\langle \theta \eta_{24}, (\Sigma \theta') \eta_{24}\rangle\\&\cong (\mathbb Z/2)^{2}.\end{aligned}\)}
	\StemEntry{\(\displaystyle \pi_{26}(S^{13})=\langle (\Sigma \theta) \eta_{25}\rangle\cong \mathbb Z/2.\)}{\(\displaystyle \begin{aligned}[t]\pi_{26}(S^{13})&=\langle (\Sigma \theta) \eta_{25}\rangle\\&\cong \mathbb Z/2.\end{aligned}\)}
	\StemEntry{\(\displaystyle \pi_{27}(S^{14})=\langle P(\iota_{29})\rangle\cong \mathbb Z_{(2)}.\)}{\(\displaystyle \begin{aligned}[t]\pi_{27}(S^{14})&=\langle P(\iota_{29})\rangle\\&\cong \mathbb Z_{(2)}.\end{aligned}\)}
	\StemEntry{\(\displaystyle \pi_{n+13}(S^{n})=0\quad(n\ge 15).\)}{\(\displaystyle \begin{aligned}[t]\pi_{n+13}(S^{n})&=0\quad(n\ge 15).\end{aligned}\)}
	\StemEntry{\(\displaystyle \pi^S_{13}(S^0)=0.\)}{\(\displaystyle \pi^S_{13}(S^0)=0.\)}
	\subsection{Stem 14}\label{app:stem-14}
	\StemEntry{\(\displaystyle \pi_{16}(S^{2})=\langle \eta_2 \nu' \eta_6 \mu_7\rangle\cong \mathbb Z/2.\)}{\(\displaystyle \begin{aligned}[t]\pi_{16}(S^{2})&=\langle \eta_2 \nu' \eta_6 \mu_7\rangle\\&\cong \mathbb Z/2.\end{aligned}\)}
	\StemEntry{\(\displaystyle \pi_{17}(S^{3})=\langle \varepsilon _3 \nu^2_{11}\rangle\cong \mathbb Z/2.\)}{\(\displaystyle \begin{aligned}[t]\pi_{17}(S^{3})&=\langle \varepsilon _3 \nu^2_{11}\rangle\\&\cong \mathbb Z/2.\end{aligned}\)}
	\StemEntry{\(\displaystyle \pi_{18}(S^{4})=\langle \varepsilon _4 \nu^2_{12}, \nu_4 \zeta_7, \nu_4 \bar\nu_7 \nu_{15}\rangle\cong \mathbb Z/2 \oplus \mathbb Z/8 \oplus \mathbb Z/2.\)}{\(\displaystyle \begin{aligned}[t]\pi_{18}(S^{4})&=\langle \varepsilon _4 \nu^2_{12}, \nu_4 \zeta_7, \nu_4 \bar\nu_7 \nu_{15}\rangle\\&\cong \mathbb Z/2 \oplus \mathbb Z/8 \oplus \mathbb Z/2.\end{aligned}\)}
	\StemEntry{\(\displaystyle \pi_{19}(S^{5})=\langle \nu_5 \zeta_8, \nu_5 \bar\nu_8 \nu_{16}\rangle\cong (\mathbb Z/2)^{2}.\)}{\(\displaystyle \begin{aligned}[t]\pi_{19}(S^{5})&=\langle \nu_5 \zeta_8, \nu_5 \bar\nu_8 \nu_{16}\rangle\\&\cong (\mathbb Z/2)^{2}.\end{aligned}\)}
	\StemEntry{\(\displaystyle \pi_{20}(S^{6})=\langle \sigma'' \sigma_{13}, \bar\nu_6 \nu^2_{14}\rangle\cong \mathbb Z/4 \oplus \mathbb Z/2.\)}{\(\displaystyle \begin{aligned}[t]\pi_{20}(S^{6})&=\langle \sigma'' \sigma_{13}, \bar\nu_6 \nu^2_{14}\rangle\\&\cong \mathbb Z/4 \oplus \mathbb Z/2.\end{aligned}\)}
	\StemEntry{\(\displaystyle \pi_{21}(S^{7})=\langle \sigma' \sigma_{14}, \kappa_7\rangle\cong \mathbb Z/8 \oplus \mathbb Z/4.\)}{\(\displaystyle \begin{aligned}[t]\pi_{21}(S^{7})&=\langle \sigma' \sigma_{14}, \kappa_7\rangle\\&\cong \mathbb Z/8 \oplus \mathbb Z/4.\end{aligned}\)}
	\StemEntry{\(\displaystyle \pi_{22}(S^{8})=\langle \sigma^2_8, (\Sigma \sigma') \sigma_{15}, \kappa_8\rangle\cong \mathbb Z/16 \oplus \mathbb Z/8 \oplus \mathbb Z/4.\)}{\(\displaystyle \begin{aligned}[t]\pi_{22}(S^{8})&=\langle \sigma^2_8, (\Sigma \sigma') \sigma_{15}, \kappa_8\rangle\\&\cong \mathbb Z/16 \oplus \mathbb Z/8 \oplus \mathbb Z/4.\end{aligned}\)}
	\StemEntry{\(\displaystyle \pi_{23}(S^{9})=\langle \sigma^2_9, \kappa_9\rangle\cong \mathbb Z/16 \oplus \mathbb Z/4.\)}{\(\displaystyle \begin{aligned}[t]\pi_{23}(S^{9})&=\langle \sigma^2_9, \kappa_9\rangle\\&\cong \mathbb Z/16 \oplus \mathbb Z/4.\end{aligned}\)}
	\StemEntry{\(\displaystyle \pi_{28}(S^{14})=\langle \sigma^2_{14}, \kappa_{14}\rangle\cong \mathbb Z/8 \oplus \mathbb Z/2.\)}{\(\displaystyle \begin{aligned}[t]\pi_{28}(S^{14})&=\langle \sigma^2_{14}, \kappa_{14}\rangle\\&\cong \mathbb Z/8 \oplus \mathbb Z/2.\end{aligned}\)}
	\StemEntry{\(\displaystyle \pi_{29}(S^{15})=\langle \sigma^2_{15}, \kappa_{15}\rangle\cong \mathbb Z/4 \oplus \mathbb Z/2.\)}{\(\displaystyle \begin{aligned}[t]\pi_{29}(S^{15})&=\langle \sigma^2_{15}, \kappa_{15}\rangle\\&\cong \mathbb Z/4 \oplus \mathbb Z/2.\end{aligned}\)}
	\StemEntry{\(\displaystyle \pi_{n+14}(S^{n})=\langle \sigma^2_n, \kappa_n\rangle\cong (\mathbb Z/2)^{2}\quad(n\ge 16).\)}{\(\displaystyle \begin{aligned}[t]\pi_{n+14}(S^{n})&=\langle \sigma^2_n, \kappa_n\rangle\\&\cong (\mathbb Z/2)^{2}\quad(n\ge 16).\end{aligned}\)}
	
	\StemEntry{\(\displaystyle \pi_{n+14}(S^{n})=\langle \sigma^2_n, \kappa_n\rangle\cong \mathbb Z/16 \oplus \mathbb Z/2\quad(n = 10, 11, 13).\)}{\(\displaystyle \begin{aligned}[t]\pi_{n+14}(S^{n})&=\langle \sigma^2_n, \kappa_n\rangle\\&\cong \mathbb Z/16 \oplus \mathbb Z/2\quad(n = 10, 11, 13).\end{aligned}\)}
	\StemEntry{\(\displaystyle \pi_{26}(S^{12})=\langle \sigma^2_{12}, \kappa_{12}, P(\nu_{25})\rangle\cong \mathbb Z/16 \oplus \mathbb Z/2 \oplus \mathbb Z/4.\)}{\(\displaystyle \begin{aligned}[t]\pi_{26}(S^{12})&=\langle \sigma^2_{12}, \kappa_{12}, P(\nu_{25})\rangle\\&\cong \mathbb Z/16 \oplus \mathbb Z/2 \oplus \mathbb Z/4.\end{aligned}\)}
	\StemEntry{\(\displaystyle \pi^S_{14}(S^0)=\langle \sigma^2, \kappa\rangle\cong \mathbb Z/2 \oplus \mathbb Z/2.\)}{\(\displaystyle \begin{aligned}[t]\pi^S_{14}(S^0)&=\langle \sigma^2, \kappa\rangle\\&\cong \mathbb Z/2 \oplus \mathbb Z/2.\end{aligned}\)}
	\subsection{Stem 15}\label{app:stem-15}
	\StemEntry{\(\displaystyle \pi_{17}(S^{2})=\langle \eta_2 \varepsilon _3 \nu^2_{11}\rangle\cong \mathbb Z/2.\)}{\(\displaystyle \begin{aligned}[t]\pi_{17}(S^{2})&=\langle \eta_2 \varepsilon _3 \nu^2_{11}\rangle\\&\cong \mathbb Z/2.\end{aligned}\)}
	\StemEntry{\(\displaystyle \pi_{n+15}(S^{n})=\langle \bar\varepsilon _n\rangle\cong \mathbb Z/2\quad(n = 3 \text{ or } 4).\)}{\(\displaystyle \begin{aligned}[t]\pi_{n+15}(S^{n})&=\langle \bar\varepsilon _n\rangle\\&\cong \mathbb Z/2\quad(n = 3 \text{ or } 4).\end{aligned}\)}
	\StemEntry{\(\displaystyle \pi_{20}(S^{5})=\langle \rho^{\mathrm{IV}}, \bar\varepsilon _5\rangle\cong (\mathbb Z/2)^{2}.\)}{\(\displaystyle \begin{aligned}[t]\pi_{20}(S^{5})&=\langle \rho^{\mathrm{IV}}, \bar\varepsilon _5\rangle\\&\cong (\mathbb Z/2)^{2}.\end{aligned}\)}
	\StemEntry{\(\displaystyle \pi_{21}(S^{6})=\langle \rho''', \bar\varepsilon _6\rangle\cong \mathbb Z/4 \oplus \mathbb Z/2.\)}{\(\displaystyle \begin{aligned}[t]\pi_{21}(S^{6})&=\langle \rho''', \bar\varepsilon _6\rangle\\&\cong \mathbb Z/4 \oplus \mathbb Z/2.\end{aligned}\)}
	\StemEntry{\(\displaystyle \pi_{22}(S^{7})=\langle \rho'', \sigma' \bar\nu_{14}, \sigma' \varepsilon _{14}, \bar\varepsilon _7\rangle\cong \mathbb Z/8 \oplus (\mathbb Z/2)^{3}.\)}{\(\displaystyle \begin{aligned}[t]\pi_{22}(S^{7})&=\langle \rho'', \sigma' \bar\nu_{14}, \sigma' \varepsilon _{14}, \bar\varepsilon _7\rangle\\&\cong \mathbb Z/8 \oplus (\mathbb Z/2)^{3}.\end{aligned}\)}
	\StemEntry{\(\displaystyle \pi_{23}(S^{8})=\langle \Sigma \rho'', \sigma_8 \bar\nu_{15}, \sigma_8 \varepsilon _{15}, (\Sigma \sigma') \bar\nu_{15}, (\Sigma \sigma') \varepsilon _{15}, \bar\varepsilon _8\rangle\cong \mathbb Z/8 \oplus (\mathbb Z/2)^{5}.\)}{\(\displaystyle \begin{aligned}[t]\pi_{23}(S^{8})&=\langle \Sigma \rho'', \sigma_8 \bar\nu_{15}, \sigma_8 \varepsilon _{15}, (\Sigma \sigma') \bar\nu_{15}, (\Sigma \sigma') \varepsilon _{15},\\&\phantom{=\langle{}}\bar\varepsilon _8\rangle\cong \mathbb Z/8 \oplus (\mathbb Z/2)^{5}.\end{aligned}\)}
	\StemEntry{\(\displaystyle \pi_{24}(S^{9})=\langle \rho', \sigma_9 \bar\nu_{16}, \sigma_9 \varepsilon _{16}, \bar\varepsilon _9\rangle\cong \mathbb Z/16 \oplus (\mathbb Z/2)^{3}.\)}{\(\displaystyle \begin{aligned}[t]\pi_{24}(S^{9})&=\langle \rho', \sigma_9 \bar\nu_{16}, \sigma_9 \varepsilon _{16}, \bar\varepsilon _9\rangle\\&\cong \mathbb Z/16 \oplus (\mathbb Z/2)^{3}.\end{aligned}\)}
	\StemEntry{\(\displaystyle \pi_{25}(S^{10})=\langle \Sigma \rho', \sigma_{10} \bar\nu_{17}, \bar\varepsilon _{10}\rangle\cong \mathbb Z/16 \oplus (\mathbb Z/2)^{2}.\)}{\(\displaystyle \begin{aligned}[t]\pi_{25}(S^{10})&=\langle \Sigma \rho', \sigma_{10} \bar\nu_{17}, \bar\varepsilon _{10}\rangle\\&\cong \mathbb Z/16 \oplus (\mathbb Z/2)^{2}.\end{aligned}\)}
	\StemEntry{\(\displaystyle \pi_{26}(S^{11})=\langle \Sigma ^2 \rho', \bar\varepsilon _{11}\rangle\cong \mathbb Z/16 \oplus \mathbb Z/2.\)}{\(\displaystyle \begin{aligned}[t]\pi_{26}(S^{11})&=\langle \Sigma ^2 \rho', \bar\varepsilon _{11}\rangle\\&\cong \mathbb Z/16 \oplus \mathbb Z/2.\end{aligned}\)}
	\StemEntry{\(\displaystyle \pi_{27}(S^{12})=\langle \Sigma ^3 \rho', \bar\varepsilon _{12}\rangle\cong \mathbb Z/16 \oplus \mathbb Z/2.\)}{\(\displaystyle \begin{aligned}[t]\pi_{27}(S^{12})&=\langle \Sigma ^3 \rho', \bar\varepsilon _{12}\rangle\\&\cong \mathbb Z/16 \oplus \mathbb Z/2.\end{aligned}\)}
	\StemEntry{\(\displaystyle \pi_{n+15}(S^{n})=\langle \rho_n, \bar\varepsilon _n\rangle\cong \mathbb Z/32 \oplus \mathbb Z/2\quad(n \ge 13 \text{ but } n\ne16).\)}{\(\displaystyle \begin{aligned}[t]\pi_{n+15}(S^{n})&=\langle \rho_n, \bar\varepsilon _n\rangle\\&\cong \mathbb Z/32 \oplus \mathbb Z/2\quad(n \ge 13 \text{ but } n\ne16).\end{aligned}\)}
	\StemEntry{\(\displaystyle \pi_{31}(S^{16})=\langle \rho_{16}, \bar\varepsilon _{16}, P(\iota_{33})\rangle\cong \mathbb Z/32 \oplus \mathbb Z/2 \oplus \mathbb Z_{(2)}.\)}{\(\displaystyle \begin{aligned}[t]\pi_{31}(S^{16})&=\langle \rho_{16}, \bar\varepsilon _{16}, P(\iota_{33})\rangle\\&\cong \mathbb Z/32 \oplus \mathbb Z/2 \oplus \mathbb Z_{(2)}.\end{aligned}\)}
	\StemEntry{\(\displaystyle \pi^S_{15}(S^0)=\langle \rho, \bar\varepsilon\rangle\cong \mathbb Z/32 \oplus \mathbb Z/2;\quad (\bar{\lt}=\ca\cn).\)}{\(\displaystyle \begin{aligned}[t]\pi^S_{15}(S^0)&=\langle \rho, \bar\varepsilon\rangle\\&\cong \mathbb Z/32 \oplus \mathbb Z/2;\quad (\bar{\lt}=\ca\cn).\end{aligned}\)}
	
	\subsection{Stem 16}\label{app:stem-16}
	\StemEntry{\(\displaystyle \pi_{18}(S^{2})=\langle \eta_2 \bar\varepsilon _3\rangle\cong \mathbb Z/2.\)}{\(\displaystyle \begin{aligned}[t]\pi_{18}(S^{2})&=\langle \eta_2 \bar\varepsilon _3\rangle\\&\cong \mathbb Z/2.\end{aligned}\)}
	\StemEntry{\(\displaystyle \pi_{n+16}(S^{n})=\langle \mu_n \sigma_{n+9}, \eta_n \bar\varepsilon _{n+1}\rangle\cong (\mathbb Z/2)^{2}\quad(n=3,5).\)}{\(\displaystyle \begin{aligned}[t]\pi_{n+16}(S^{n})&=\langle \mu_n \sigma_{n+9}, \eta_n \bar\varepsilon _{n+1}\rangle\\&\cong (\mathbb Z/2)^{2}\quad(n=3,5).\end{aligned}\)}
	\StemEntry{\(\displaystyle \pi_{20}(S^{4})=\langle \nu^2_4 \sigma_{10} \nu_{17}, \mu_4 \sigma_{13}, \eta_4 \bar\varepsilon _5\rangle\cong (\mathbb Z/2)^{3}.\)}{\(\displaystyle \begin{aligned}[t]\pi_{20}(S^{4})&=\langle \nu^2_4 \sigma_{10} \nu_{17}, \mu_4 \sigma_{13},\\&\phantom{=\langle{}}\eta_4 \bar\varepsilon _5\rangle\cong (\mathbb Z/2)^{3}.\end{aligned}\)}
	\StemEntry{\(\displaystyle \pi_{22}(S^{6})=\langle \zeta', \mu_6 \sigma_{15}, \eta_6 \bar\varepsilon _7\rangle\cong \mathbb Z/8 \oplus (\mathbb Z/2)^{2}.\)}{\(\displaystyle \begin{aligned}[t]\pi_{22}(S^{6})&=\langle \zeta', \mu_6 \sigma_{15}, \eta_6 \bar\varepsilon _7\rangle\\&\cong \mathbb Z/8 \oplus (\mathbb Z/2)^{2}.\end{aligned}\)}
	\StemEntry{\(\displaystyle \pi_{23}(S^{7})=\langle \sigma' \mu_{14}, \Sigma \zeta', \mu_7 \sigma_{16}, \eta_7 \bar\varepsilon _8\rangle\cong (\mathbb Z/2)^{4}.\)}{\(\displaystyle \begin{aligned}[t]\pi_{23}(S^{7})&=\langle \sigma' \mu_{14}, \Sigma \zeta', \mu_7 \sigma_{16}, \eta_7 \bar\varepsilon _8\rangle\\&\cong (\mathbb Z/2)^{4}.\end{aligned}\)}
	\StemEntry{\(\displaystyle \pi_{24}(S^{8})=\langle \sigma_8 \nu^3_{15}, \sigma_8 \mu_{15}, \sigma_8 \eta_{15} \varepsilon _{16}, (\Sigma \sigma') \mu_{15}, \Sigma ^2 \zeta', \mu_8 \sigma_{17}, \eta_8 \bar\varepsilon _9\rangle\cong (\mathbb Z/2)^{7}.\)}{\(\displaystyle \begin{aligned}[t]\pi_{24}(S^{8})&=\langle \sigma_8 \nu^3_{15}, \sigma_8 \mu_{15}, \sigma_8 \eta_{15} \varepsilon _{16}, (\Sigma \sigma') \mu_{15}, \Sigma ^2 \zeta', \mu_8 \sigma_{17},\\&\phantom{=\langle{}}\eta_8 \bar\varepsilon _9\rangle\cong (\mathbb Z/2)^{7}.\end{aligned}\)}
	\StemEntry{\(\displaystyle \pi_{25}(S^{9})=\langle \sigma_9 \nu^3_{16}, \sigma_9 \mu_{16}, \sigma_9 \eta_{16} \varepsilon _{17}, \mu_9 \sigma_{18}\rangle\cong (\mathbb Z/2)^{4}.\)}{\(\displaystyle \begin{aligned}[t]\pi_{25}(S^{9})&=\langle \sigma_9 \nu^3_{16}, \sigma_9 \mu_{16}, \sigma_9 \eta_{16} \varepsilon _{17},\\&\phantom{=\langle{}}\mu_9 \sigma_{18}\rangle\cong (\mathbb Z/2)^{4}.\end{aligned}\)}
	\StemEntry{\(\displaystyle \pi_{26}(S^{10})=\langle P(\sigma_{21}), \sigma_{10} \mu_{17}\rangle\cong \mathbb Z/16 \oplus \mathbb Z/2.\)}{\(\displaystyle \begin{aligned}[t]\pi_{26}(S^{10})&=\langle P(\sigma_{21}), \sigma_{10} \mu_{17}\rangle\\&\cong \mathbb Z/16 \oplus \mathbb Z/2.\end{aligned}\)}
	\StemEntry{\(\displaystyle \pi_{n+16}(S^{n})=\langle \sigma_n \mu_{n+7}\rangle\cong \mathbb Z/2\quad(n = 11, 12, 13).\)}{\(\displaystyle \begin{aligned}[t]\pi_{n+16}(S^{n})&=\langle \sigma_n \mu_{n+7}\rangle\\&\cong \mathbb Z/2\quad(n = 11, 12, 13).\end{aligned}\)}
	\StemEntry{\(\displaystyle \pi_{30}(S^{14})=\langle \omega_{14}, \sigma_{14} \mu_{21}\rangle\cong \mathbb Z/8 \oplus \mathbb Z/2.\)}{\(\displaystyle \begin{aligned}[t]\pi_{30}(S^{14})&=\langle \omega_{14}, \sigma_{14} \mu_{21}\rangle\\&\cong \mathbb Z/8 \oplus \mathbb Z/2.\end{aligned}\)}
	\StemEntry{\(\displaystyle \pi_{31}(S^{15})=\langle \eta^{*\prime}, \omega_{15}, \sigma_{15} \mu_{22}\rangle\cong (\mathbb Z/2)^{3}.\)}{\(\displaystyle \begin{aligned}[t]\pi_{31}(S^{15})&=\langle \eta^{*\prime}, \omega_{15}, \sigma_{15} \mu_{22}\rangle\\&\cong (\mathbb Z/2)^{3}.\end{aligned}\)}
	\StemEntry{\(\displaystyle \pi_{32}(S^{16})=\langle \eta^{*}_{16}, \Sigma \eta^{*\prime}, \omega_{16}, \sigma_{16} \mu_{23}\rangle\cong (\mathbb Z/2)^{4}.\)}{\(\displaystyle \begin{aligned}[t]\pi_{32}(S^{16})&=\langle \eta^{*}_{16}, \Sigma \eta^{*\prime}, \omega_{16}, \sigma_{16} \mu_{23}\rangle\\&\cong (\mathbb Z/2)^{4}.\end{aligned}\)}
	\StemEntry{\(\displaystyle \pi_{33}(S^{17})=\langle \eta^{*}_{17}, \omega_{17}, \sigma_{17} \mu_{24}\rangle\cong (\mathbb Z/2)^{3}.\)}{\(\displaystyle \begin{aligned}[t]\pi_{33}(S^{17})&=\langle \eta^{*}_{17}, \omega_{17}, \sigma_{17} \mu_{24}\rangle\\&\cong (\mathbb Z/2)^{3}.\end{aligned}\)}
	\StemEntry{\(\displaystyle \pi_{n+16}(S^{n})=\langle \omega_n, \sigma_n \mu_{n+7}\rangle\cong (\mathbb Z/2)^{2}\quad(n\ge 18).\)}{\(\displaystyle \begin{aligned}[t]\pi_{n+16}(S^{n})&=\langle \omega_n, \sigma_n \mu_{n+7}\rangle\\&\cong (\mathbb Z/2)^{2}\quad(n\ge 18).\end{aligned}\)}
	\StemEntry{\(\displaystyle \pi^S_{16}(S^0)=\langle \omega, \sigma \mu\rangle\cong \mathbb Z/2 \oplus \mathbb Z/2.\)}{\(\displaystyle \begin{aligned}[t]\pi^S_{16}(S^0)&=\langle \omega, \sigma \mu\rangle\\&\cong \mathbb Z/2 \oplus \mathbb Z/2.\end{aligned}\)}
	
	\subsection{Stem 17}\label{app:stem-17}
	\StemEntry{\(\displaystyle \pi_{19}(S^{2})=\langle \eta_2 \mu_3 \sigma_{12}, \eta^2_2 \bar\varepsilon _4\rangle\cong (\mathbb Z/2)^{2}.\)}{\(\displaystyle \begin{aligned}[t]\pi_{19}(S^{2})&=\langle \eta_2 \mu_3 \sigma_{12}, \eta^2_2 \bar\varepsilon _4\rangle\\&\cong (\mathbb Z/2)^{2}.\end{aligned}\)}
	\StemEntry{\(\displaystyle \pi_{20}(S^{3})=\langle \bar\varepsilon ', \bar\mu_3, \eta_3 \mu_4 \sigma_{13}\rangle\cong \mathbb Z/4 \oplus (\mathbb Z/2)^{2}.\)}{\(\displaystyle \begin{aligned}[t]\pi_{20}(S^{3})&=\langle \bar\varepsilon ', \bar\mu_3, \eta_3 \mu_4 \sigma_{13}\rangle\\&\cong \mathbb Z/4 \oplus (\mathbb Z/2)^{2}.\end{aligned}\)}
	\StemEntry{\(\displaystyle \pi_{21}(S^{4})=\langle \nu_4 \sigma' \sigma_{14}, \nu_4 \kappa_7, \Sigma \bar\varepsilon ', \bar\mu_4, \eta_4 \mu_5 \sigma_{14}\rangle\cong \mathbb Z/8 \oplus \mathbb Z/4 \oplus \mathbb Z/4 \oplus (\mathbb Z/2)^{2}.\)}{\(\displaystyle \begin{aligned}[t]\pi_{21}(S^{4})&=\langle \nu_4 \sigma' \sigma_{14}, \nu_4 \kappa_7, \Sigma \bar\varepsilon ', \bar\mu_4, \eta_4 \mu_5 \sigma_{14}\rangle\\&\cong \mathbb Z/8 \oplus \mathbb Z/4 \oplus \mathbb Z/4 \oplus (\mathbb Z/2)^{2}.\end{aligned}\)}
	\StemEntry{\(\displaystyle \pi_{22}(S^{5})=\langle \nu_5 \kappa_8, \bar\mu_5, \eta_5 \mu_6 \sigma_{15}\rangle\cong \mathbb Z/4 \oplus (\mathbb Z/2)^{2}.\)}{\(\displaystyle \begin{aligned}[t]\pi_{22}(S^{5})&=\langle \nu_5 \kappa_8, \bar\mu_5, \eta_5 \mu_6 \sigma_{15}\rangle\\&\cong \mathbb Z/4 \oplus (\mathbb Z/2)^{2}.\end{aligned}\)}
	\StemEntry{\(\displaystyle \pi_{23}(S^{6})=\langle P(\Sigma \theta), \nu_6 \kappa_9, \bar\mu_6, \eta_6 \mu_7 \sigma_{16}\rangle\cong (\mathbb Z/2)^{4}.\)}{\(\displaystyle \begin{aligned}[t]\pi_{23}(S^{6})&=\langle P(\Sigma \theta), \nu_6 \kappa_9, \bar\mu_6, \eta_6 \mu_7 \sigma_{16}\rangle\\&\cong (\mathbb Z/2)^{4}.\end{aligned}\)}
	\StemEntry{\(\displaystyle \pi_{24}(S^{7})=\langle \sigma' \eta_{14} \mu_{15}, \nu_7 \kappa_{10}, \bar\mu_7, \eta_7 \mu_8 \sigma_{17}\rangle\cong (\mathbb Z/2)^{4}.\)}{\(\displaystyle \begin{aligned}[t]\pi_{24}(S^{7})&=\langle \sigma' \eta_{14} \mu_{15}, \nu_7 \kappa_{10}, \bar\mu_7, \eta_7 \mu_8 \sigma_{17}\rangle\\&\cong (\mathbb Z/2)^{4}.\end{aligned}\)}
	\StemEntry{\(\displaystyle \pi_{25}(S^{8})=\langle \sigma_8 \eta_{15} \mu_{16}, (\Sigma \sigma') \eta_{15} \mu_{16}, \nu_8 \kappa_{11}, \bar\mu_8, \eta_8 \mu_9 \sigma_{18}\rangle\cong (\mathbb Z/2)^{5}.\)}{\(\displaystyle \begin{aligned}[t]\pi_{25}(S^{8})&=\langle \sigma_8 \eta_{15} \mu_{16}, (\Sigma \sigma') \eta_{15} \mu_{16}, \nu_8 \kappa_{11}, \bar\mu_8,\\&\phantom{=\langle{}}\eta_8 \mu_9 \sigma_{18}\rangle\cong (\mathbb Z/2)^{5}.\end{aligned}\)}
	\StemEntry{\(\displaystyle \pi_{26}(S^{9})=\langle \sigma_9 \eta_{16} \mu_{17}, \nu_9 \kappa_{12}, \bar\mu_9, \eta_9 \mu_{10} \sigma_{19}\rangle\cong (\mathbb Z/2)^{4}.\)}{\(\displaystyle \begin{aligned}[t]\pi_{26}(S^{9})&=\langle \sigma_9 \eta_{16} \mu_{17}, \nu_9 \kappa_{12}, \bar\mu_9, \eta_9 \mu_{10} \sigma_{19}\rangle\\&\cong (\mathbb Z/2)^{4}.\end{aligned}\)}
	\StemEntry{\(\displaystyle \pi_{n+17}(S^{n})=\langle \sigma_n \eta_{n+7} \mu_{n+8}, \nu_n \kappa_{n+3}, \bar\mu_n\rangle\cong (\mathbb Z/2)^{3}\quad(n = 10, 11).\)}{\(\displaystyle \begin{aligned}[t]\pi_{n+17}(S^{n})&=\langle \sigma_n \eta_{n+7} \mu_{n+8}, \nu_n \kappa_{n+3}, \bar\mu_n\rangle\\&\cong (\mathbb Z/2)^{3}\quad(n = 10, 11).\end{aligned}\)}
	\StemEntry{\(\displaystyle \pi_{n+17}(S^{n})=\langle \varepsilon ^{*}_n, \sigma_n \eta_{n+7} \mu_{n+8}, \nu_n \kappa_{n+3}, \bar\mu_n\rangle\cong (\mathbb Z/2)^{4}\quad(n =12, 13, 14).\)}{\(\displaystyle \begin{aligned}[t]\pi_{n+17}(S^{n})&=\langle \varepsilon ^{*}_n, \sigma_n \eta_{n+7} \mu_{n+8}, \nu_n \kappa_{n+3}, \bar\mu_n\rangle\\&\cong (\mathbb Z/2)^{4}\quad(n =12, 13, 14).\end{aligned}\)}
	\StemEntry{\(\displaystyle \pi_{32}(S^{15})=\langle \eta^{*\prime} \eta_{31}, \varepsilon ^{*}_{15}, \sigma_{15} \eta_{22} \mu_{23}, \nu_{15} \kappa_{18}, \bar\mu_{15}\rangle\cong (\mathbb Z/2)^{5}.\)}{\(\displaystyle \begin{aligned}[t]\pi_{32}(S^{15})&=\langle \eta^{*\prime} \eta_{31}, \varepsilon ^{*}_{15}, \sigma_{15} \eta_{22} \mu_{23}, \nu_{15} \kappa_{18},\\&\phantom{=\langle{}}\bar\mu_{15}\rangle\cong (\mathbb Z/2)^{5}.\end{aligned}\)}
	\StemEntry{\(\displaystyle \pi_{33}(S^{16})=\langle \eta^{*}_{16} \eta_{32}, (\Sigma \eta^{*\prime}) \eta_{32}, \varepsilon ^{*}_{16}, \sigma_{16} \eta_{23} \mu_{24}, \nu_{16} \kappa_{19}, \bar\mu_{16}\rangle\cong (\mathbb Z/2)^{6}.\)}{\(\displaystyle \begin{aligned}[t]\pi_{33}(S^{16})&=\langle \eta^{*}_{16} \eta_{32}, (\Sigma \eta^{*\prime}) \eta_{32}, \varepsilon ^{*}_{16}, \sigma_{16} \eta_{23} \mu_{24}, \nu_{16} \kappa_{19},\\&\phantom{=\langle{}}\bar\mu_{16}\rangle\cong (\mathbb Z/2)^{6}.\end{aligned}\)}
	\StemEntry{\(\displaystyle \pi_{34}(S^{17})=\langle \eta^{*}_{17} \eta_{33}, \varepsilon ^{*}_{17}, \sigma_{17} \eta_{24} \mu_{25}, \nu_{17} \kappa_{20}, \bar\mu_{17}\rangle\cong (\mathbb Z/2)^{5}.\)}{\(\displaystyle \begin{aligned}[t]\pi_{34}(S^{17})&=\langle \eta^{*}_{17} \eta_{33}, \varepsilon ^{*}_{17}, \sigma_{17} \eta_{24} \mu_{25}, \nu_{17} \kappa_{20},\\&\phantom{=\langle{}}\bar\mu_{17}\rangle\cong (\mathbb Z/2)^{5}.\end{aligned}\)}
	\StemEntry{\(\displaystyle \pi_{35}(S^{18})=\langle P(\iota_{37}), \varepsilon ^{*}_{18}, \sigma_{18} \eta_{25} \mu_{26}, \nu_{18} \kappa_{21}, \bar\mu_{18}\rangle\cong \mathbb Z_{(2)} \oplus (\mathbb Z/2)^{4}.\)}{\(\displaystyle \begin{aligned}[t]\pi_{35}(S^{18})&=\langle P(\iota_{37}), \varepsilon ^{*}_{18}, \sigma_{18} \eta_{25} \mu_{26}, \nu_{18} \kappa_{21},\\&\phantom{=\langle{}}\bar\mu_{18}\rangle\cong \mathbb Z_{(2)} \oplus (\mathbb Z/2)^{4}.\end{aligned}\)}
	\StemEntry{\(\displaystyle \pi_{n+17}(S^{n})=\langle \varepsilon ^{*}_n, \sigma_n \eta_{n+7} \mu_{n+8}, \nu_n \kappa_{n+3}, \bar\mu_n\rangle\cong (\mathbb Z/2)^{4}\quad(n\ge 19).\)}{\(\displaystyle \begin{aligned}[t]\pi_{n+17}(S^{n})&=\langle \varepsilon ^{*}_n, \sigma_n \eta_{n+7} \mu_{n+8}, \nu_n \kappa_{n+3}, \bar\mu_n\rangle\\&\cong (\mathbb Z/2)^{4}\quad(n\ge 19).\end{aligned}\)}
	\StemEntry{\(\displaystyle \pi^S_{17}(S^0)=\langle \varepsilon ^{*}, \sigma \eta \mu, \nu \kappa, \bar\mu\rangle\cong \mathbb Z/2 \oplus \mathbb Z/2 \oplus \mathbb Z/2 \oplus \mathbb Z/2;\quad ( \varepsilon ^{*}=\ca\w).\)}{\(\displaystyle \begin{aligned}[t]\pi^S_{17}(S^0)&=\langle \varepsilon ^{*}, \sigma \eta \mu, \nu \kappa, \bar\mu\rangle\\&\cong \mathbb Z/2 \oplus \mathbb Z/2 \oplus \mathbb Z/2 \oplus \mathbb Z/2;\quad ( \varepsilon ^{*}=\ca\w).\end{aligned}\)}
	
	\subsection{Stem 18}\label{app:stem-18}
	\StemEntry{\(\displaystyle \pi_{20}(S^{2})=\langle \eta_2 \bar\varepsilon ', \eta_2 \bar\mu_3, \eta^2_2 \mu_4 \sigma_{13}\rangle\cong \mathbb Z/4 \oplus (\mathbb Z/2)^{2}.\)}{\(\displaystyle \begin{aligned}[t]\pi_{20}(S^{2})&=\langle \eta_2 \bar\varepsilon ', \eta_2 \bar\mu_3, \eta^2_2 \mu_4 \sigma_{13}\rangle\\&\cong \mathbb Z/4 \oplus (\mathbb Z/2)^{2}.\end{aligned}\)}
	\StemEntry{\(\displaystyle \pi_{21}(S^{3})=\langle \mu' \sigma_{14}, \nu' \bar\varepsilon _6, \eta_3 \bar\mu_4\rangle\cong \mathbb Z/4 \oplus (\mathbb Z/2)^{2}.\)}{\(\displaystyle \begin{aligned}[t]\pi_{21}(S^{3})&=\langle \mu' \sigma_{14}, \nu' \bar\varepsilon _6, \eta_3 \bar\mu_4\rangle\\&\cong \mathbb Z/4 \oplus (\mathbb Z/2)^{2}.\end{aligned}\)}
	\StemEntry{\(\displaystyle \pi_{22}(S^{4})=\langle \nu_4 \rho'', \nu_4 \sigma' \bar\nu_{14}, \nu_4 \sigma' \varepsilon _{14}, \nu_4 \bar\varepsilon _7, (\Sigma \mu' )\sigma_{15}, (\Sigma \nu') \bar\varepsilon _7, \eta_4 \bar\mu_5\rangle\cong \mathbb Z/8 \oplus \mathbb Z/2 \oplus \mathbb Z/2 \oplus \mathbb Z/2 \oplus \mathbb Z/4 \oplus (\mathbb Z/2)^{2}.\)}{\(\displaystyle \begin{aligned}[t]\pi_{22}(S^{4})&=\langle \nu_4 \rho'', \nu_4 \sigma' \bar\nu_{14}, \nu_4 \sigma' \varepsilon _{14}, \nu_4 \bar\varepsilon _7, (\Sigma \mu' )\sigma_{15}, (\Sigma \nu') \bar\varepsilon _7,\\&\phantom{=\langle{}}\eta_4 \bar\mu_5\rangle\cong \mathbb Z/8 \oplus \mathbb Z/2 \oplus \mathbb Z/2 \oplus \mathbb Z/2 \oplus \mathbb Z/4 \oplus (\mathbb Z/2)^{2}.\end{aligned}\)}
	\StemEntry{\(\displaystyle \pi_{23}(S^{5})=\langle \zeta_5 \sigma_{16}, \nu_5 \bar\varepsilon _8, \eta_5 \bar\mu_6\rangle\cong \mathbb Z/8 \oplus (\mathbb Z/2)^{2}.\)}{\(\displaystyle \begin{aligned}[t]\pi_{23}(S^{5})&=\langle \zeta_5 \sigma_{16}, \nu_5 \bar\varepsilon _8, \eta_5 \bar\mu_6\rangle\\&\cong \mathbb Z/8 \oplus (\mathbb Z/2)^{2}.\end{aligned}\)}
	\StemEntry{\(\displaystyle \pi_{24}(S^{6})=\langle P(\Sigma \theta) \eta_{23}, \zeta_6 \sigma_{17}, \eta_6 \bar\mu_7\rangle\cong \mathbb Z/2 \oplus \mathbb Z/8 \oplus \mathbb Z/2.\)}{\(\displaystyle \begin{aligned}[t]\pi_{24}(S^{6})&=\langle P(\Sigma \theta) \eta_{23}, \zeta_6 \sigma_{17},\\&\phantom{=\langle{}}\eta_6 \bar\mu_7\rangle\cong \mathbb Z/2 \oplus \mathbb Z/8 \oplus \mathbb Z/2.\end{aligned}\)}
	\StemEntry{\(\displaystyle \pi_{25}(S^{7})=\langle \zeta_7 \sigma_{18}, \eta_7 \bar\mu_8\rangle\cong \mathbb Z/8 \oplus \mathbb Z/2.\)}{\(\displaystyle \begin{aligned}[t]\pi_{25}(S^{7})&=\langle \zeta_7 \sigma_{18}, \eta_7 \bar\mu_8\rangle\\&\cong \mathbb Z/8 \oplus \mathbb Z/2.\end{aligned}\)}
	\StemEntry{\(\displaystyle \pi_{26}(S^{8})=\langle \sigma_8 \zeta_{15}, \zeta_8 \sigma_{19}, \eta_8 \bar\mu_9\rangle\cong \mathbb Z/8 \oplus \mathbb Z/8 \oplus \mathbb Z/2.\)}{\(\displaystyle \begin{aligned}[t]\pi_{26}(S^{8})&=\langle \sigma_8 \zeta_{15}, \zeta_8 \sigma_{19}, \eta_8 \bar\mu_9\rangle\\&\cong \mathbb Z/8 \oplus \mathbb Z/8 \oplus \mathbb Z/2.\end{aligned}\)}

	\StemEntry{\(\displaystyle \pi_{27}(S^{9})=\langle \sigma_9 \zeta_{16}, \eta_9 \bar\mu_{10}\rangle\cong \mathbb Z/8 \oplus \mathbb Z/2.\)}{\(\displaystyle \begin{aligned}[t]\pi_{27}(S^{9})&=\langle \sigma_9 \zeta_{16}, \eta_9 \bar\mu_{10}\rangle\\&\cong \mathbb Z/8 \oplus \mathbb Z/2.\end{aligned}\)}\vspace{-0.8em}\;\\	In $\xi''\pm\lambda'',$   the sign is chosen so that the element has order $2$
	\cite[Theorem, pp.~231--232]{Oda}.\vspace{-0.4em}
	\StemEntry{\(\displaystyle \pi_{28}(S^{10})=\langle \xi'', \xi''\pm\lambda'', \eta_{10}\bar\mu_{11}\rangle\cong \mathbb Z/8 \oplus (\mathbb Z/2)^{2}.\)}{\(\displaystyle \begin{aligned}[t]\pi_{28}(S^{10})&=\langle \xi'', \xi''\pm\lambda'', \eta_{10}\bar\mu_{11}\rangle\\&\cong \mathbb Z/8 \oplus (\mathbb Z/2)^{2}.\end{aligned}\)}
	\StemEntry{\(\displaystyle \pi_{29}(S^{11})=\langle \xi', \xi'+\lambda', \eta_{11}\bar\mu_{12}\rangle\cong \mathbb Z/8 \oplus \mathbb Z/4 \oplus \mathbb Z/2.\)}{\(\displaystyle \begin{aligned}[t]\pi_{29}(S^{11})&=\langle \xi', \xi'+\lambda', \eta_{11}\bar\mu_{12}\rangle\\&\cong \mathbb Z/8 \oplus \mathbb Z/4 \oplus \mathbb Z/2.\end{aligned}\)}
	\StemEntry{\(\displaystyle \pi_{30}(S^{12})=\langle \xi_{12}, \Sigma\xi'+4\xi_{12}, \Sigma\xi'+\Sigma\lambda', \eta_{12}\bar\mu_{13}\rangle\cong \mathbb Z/32 \oplus \mathbb Z/4 \oplus \mathbb Z/4 \oplus \mathbb Z/2.\)}{\(\displaystyle \begin{aligned}[t]\pi_{30}(S^{12})&=\langle \xi_{12}, \Sigma\xi'+4\xi_{12}, \Sigma\xi'+\Sigma\lambda', \eta_{12}\bar\mu_{13}\rangle\\&\cong \mathbb Z/32 \oplus \mathbb Z/4 \oplus \mathbb Z/4 \oplus \mathbb Z/2.\end{aligned}\)}
	
	\StemEntry{\(\displaystyle \pi_{n+18}(S^{n})=\langle \xi_n, \Sigma ^{n-13}\lambda, \eta_n \bar\mu_{n+1}\rangle\cong \mathbb Z/8 \oplus \mathbb Z/8 \oplus \mathbb Z/2\quad(n = 13, 14, 15).\)}{\(\displaystyle \begin{aligned}[t]\pi_{n+18}(S^{n})&=\langle \xi_n, \Sigma ^{n-13}\lambda, \eta_n \bar\mu_{n+1}\rangle\\&\cong \mathbb Z/8 \oplus \mathbb Z/8 \oplus \mathbb Z/2\quad(n = 13, 14, 15).\end{aligned}\)}
	\StemEntry{\(\displaystyle \pi_{34}(S^{16})=\langle \nu^{*}_{16}, \xi_{16}, \Sigma ^3\lambda, \eta_{16} \bar\mu_{17}\rangle\cong \mathbb Z/8 \oplus \mathbb Z/8 \oplus \mathbb Z/8 \oplus \mathbb Z/2.\)}{\(\displaystyle \begin{aligned}[t]\pi_{34}(S^{16})&=\langle \nu^{*}_{16}, \xi_{16}, \Sigma ^3\lambda, \eta_{16} \bar\mu_{17}\rangle\\&\cong \mathbb Z/8 \oplus \mathbb Z/8 \oplus \mathbb Z/8 \oplus \mathbb Z/2.\end{aligned}\)}
	\StemEntry{\(\displaystyle \pi_{35}(S^{17})=\langle \nu^{*}_{17}, \xi_{17}, \eta_{17} \bar\mu_{18}\rangle\cong \mathbb Z/8 \oplus \mathbb Z/8 \oplus \mathbb Z/2.\)}{\(\displaystyle \begin{aligned}[t]\pi_{35}(S^{17})&=\langle \nu^{*}_{17}, \xi_{17}, \eta_{17} \bar\mu_{18}\rangle\\&\cong \mathbb Z/8 \oplus \mathbb Z/8 \oplus \mathbb Z/2.\end{aligned}\)}
	\StemEntry{\(\displaystyle \pi_{36}(S^{18})=\langle \nu^{*}_{18}, \nu^{*}_{18} + \xi_{18}, \eta_{18} \bar\mu_{19}\rangle\cong \mathbb Z/8 \oplus \mathbb Z/4 \oplus \mathbb Z/2.\)}{\(\displaystyle \begin{aligned}[t]\pi_{36}(S^{18})&=\langle \nu^{*}_{18}, \nu^{*}_{18} + \xi_{18}, \eta_{18} \bar\mu_{19}\rangle\\&\cong \mathbb Z/8 \oplus \mathbb Z/4 \oplus \mathbb Z/2.\end{aligned}\)}
	\StemEntry{\(\displaystyle \pi_{37}(S^{19})=\langle \nu^{*}_{19}, \nu^{*}_{19} + \xi_{19}, \eta_{19} \bar\mu_{20}\rangle\cong \mathbb Z/8 \oplus (\mathbb Z/2)^{2}.\)}{\(\displaystyle \begin{aligned}[t]\pi_{37}(S^{19})&=\langle \nu^{*}_{19}, \nu^{*}_{19} + \xi_{19}, \eta_{19} \bar\mu_{20}\rangle\\&\cong \mathbb Z/8 \oplus (\mathbb Z/2)^{2}.\end{aligned}\)}
	\StemEntry{\(\displaystyle \pi_{n+18}(S^{n})=\langle \nu^{*}_n, \eta_n \bar\mu_{n+1}\rangle\cong \mathbb Z/8 \oplus \mathbb Z/2\quad(n\ge 20).\)}{\(\displaystyle \begin{aligned}[t]\pi_{n+18}(S^{n})&=\langle \nu^{*}_n, \eta_n \bar\mu_{n+1}\rangle\\&\cong \mathbb Z/8 \oplus \mathbb Z/2\quad(n\ge 20).\end{aligned}\)}
	
	\StemEntry{\(\displaystyle \pi^S_{18}(S^0)=\langle \nu^{*}, \eta \bar\mu\rangle\cong \mathbb Z/8 \oplus \mathbb Z/2.\)}{\(\displaystyle \begin{aligned}[t]\pi^S_{18}(S^0)&=\langle \nu^{*}, \eta \bar\mu\rangle\\&\cong \mathbb Z/8 \oplus \mathbb Z/2.\end{aligned}\)}
	\subsection{Stem 19}\label{app:stem-19}
	\StemEntry{\(\displaystyle \pi_{21}(S^{2})=\langle \eta_2 \mu' \sigma_{14}, \eta_2 \nu' \bar\varepsilon _6, \eta^2_2 \bar\mu_4\rangle\cong \mathbb Z/4 \oplus (\mathbb Z/2)^{2}.\)}{\(\displaystyle \begin{aligned}[t]\pi_{21}(S^{2})&=\langle \eta_2 \mu' \sigma_{14}, \eta_2 \nu' \bar\varepsilon _6, \eta^2_2 \bar\mu_4\rangle\\&\cong \mathbb Z/4 \oplus (\mathbb Z/2)^{2}.\end{aligned}\)}
	\StemEntry{\(\displaystyle \pi_{22}(S^{3})=\langle \bar\mu', \nu' \mu_6 \sigma_{15}\rangle\cong \mathbb Z/4 \oplus \mathbb Z/2.\)}{\(\displaystyle \begin{aligned}[t]\pi_{22}(S^{3})&=\langle \bar\mu', \nu' \mu_6 \sigma_{15}\rangle\\&\cong \mathbb Z/4 \oplus \mathbb Z/2.\end{aligned}\)}
	\StemEntry{\(\displaystyle \pi_{23}(S^{4})=\langle \nu_4 \sigma' \mu_{14}, \nu_4 \Sigma \zeta', \nu_4 \mu_7 \sigma_{16}, \nu_4 \eta_7 \bar\varepsilon _8, \Sigma \bar\mu', (\Sigma \nu') \mu_7 \sigma_{16}\rangle\cong (\mathbb Z/2)^{4} \oplus \mathbb Z/4 \oplus \mathbb Z/2.\)}{\(\displaystyle \begin{aligned}[t]\pi_{23}(S^{4})&=\langle \nu_4 \sigma' \mu_{14}, \nu_4 \Sigma \zeta', \nu_4 \mu_7 \sigma_{16}, \nu_4 \eta_7 \bar\varepsilon _8, \Sigma \bar\mu',\\&\phantom{=\langle{}}(\Sigma \nu') \mu_7 \sigma_{16}\rangle\cong (\mathbb Z/2)^{4} \oplus \mathbb Z/4 \oplus \mathbb Z/2.\end{aligned}\)}
	\StemEntry{\(\displaystyle \pi_{24}(S^{5})=\langle \bar\zeta_5, \nu_5 \mu_8 \sigma_{17}\rangle\cong \mathbb Z/8 \oplus \mathbb Z/2.\)}{\(\displaystyle \begin{aligned}[t]\pi_{24}(S^{5})&=\langle \bar\zeta_5, \nu_5 \mu_8 \sigma_{17}\rangle\\&\cong \mathbb Z/8 \oplus \mathbb Z/2.\end{aligned}\)}
	\StemEntry{\(\displaystyle \pi_{25}(S^{6})=\langle \bar\zeta_6, \bar\sigma_6\rangle\cong \mathbb Z/8 \oplus \mathbb Z/32.\)}{\(\displaystyle \begin{aligned}[t]\pi_{25}(S^{6})&=\langle \bar\zeta_6, \bar\sigma_6\rangle\\&\cong \mathbb Z/8 \oplus \mathbb Z/32.\end{aligned}\)}
	\StemEntry{\(\displaystyle \pi_{n+19}(S^{n})=\langle \bar\zeta_n, \bar\sigma_n\rangle\cong \mathbb Z/8 \oplus \mathbb Z/2\quad(7\leq n\leq10).\)}{\(\displaystyle \begin{aligned}[t]\pi_{n+19}(S^{n})&=\langle \bar\zeta_n, \bar\sigma_n\rangle\\&\cong \mathbb Z/8 \oplus \mathbb Z/2\quad(7\leq n\leq10).\end{aligned}\)}

	\StemEntry{\(\displaystyle \pi_{30}(S^{11})=\langle \lambda' \eta_{29}, \xi' \eta_{29}, \bar\sigma_{11}, \bar\zeta_{11}\rangle\cong (\mathbb Z/2)^{3} \oplus \mathbb Z/8.\)}{\(\displaystyle \begin{aligned}[t]\pi_{30}(S^{11})&=\langle \lambda' \eta_{29}, \xi' \eta_{29}, \bar\sigma_{11}, \bar\zeta_{11}\rangle\\&\cong (\mathbb Z/2)^{3} \oplus \mathbb Z/8.\end{aligned}\)}
	\StemEntry{\(\displaystyle \pi_{31}(S^{12})=\langle \omega', \xi_{12} \eta_{30}, (\Sigma \lambda') \eta_{30}, (\Sigma \xi') \eta_{30}, \bar\sigma_{12}, \bar\zeta_{12}\rangle\cong (\mathbb Z/2)^{5} \oplus \mathbb Z/8.\)}{\(\displaystyle \begin{aligned}[t]\pi_{31}(S^{12})&=\langle \omega', \xi_{12} \eta_{30}, (\Sigma \lambda') \eta_{30}, (\Sigma \xi') \eta_{30}, \bar\sigma_{12},\\&\phantom{=\langle{}}\bar\zeta_{12}\rangle\cong (\mathbb Z/2)^{5} \oplus \mathbb Z/8.\end{aligned}\)}
	\StemEntry{\(\displaystyle \pi_{32}(S^{13})=\langle \Sigma \omega', \xi_{13} \eta_{31}, \bar\sigma_{13}, \bar\zeta_{13}\rangle\cong (\mathbb Z/2)^{3} \oplus \mathbb Z/8.\)}{\(\displaystyle \begin{aligned}[t]\pi_{32}(S^{13})&=\langle \Sigma \omega', \xi_{13} \eta_{31}, \bar\sigma_{13}, \bar\zeta_{13}\rangle\\&\cong (\mathbb Z/2)^{3} \oplus \mathbb Z/8.\end{aligned}\)}
	\StemEntry{\(\displaystyle \pi_{33}(S^{14})=\langle \omega_{14} \nu_{30}, \bar\sigma_{14}, \bar\zeta_{14}\rangle\cong \mathbb Z/4 \oplus \mathbb Z/2 \oplus \mathbb Z/8.\)}{\(\displaystyle \begin{aligned}[t]\pi_{33}(S^{14})&=\langle \omega_{14} \nu_{30}, \bar\sigma_{14}, \bar\zeta_{14}\rangle\\&\cong \mathbb Z/4 \oplus \mathbb Z/2 \oplus \mathbb Z/8.\end{aligned}\)}
	\StemEntry{\(\displaystyle \pi_{n+19}(S^{n})=\langle \omega_n \nu_{n+16}, \bar\sigma_n, \bar\zeta_n\rangle\cong (\mathbb Z/2)^{2} \oplus \mathbb Z/8\quad(n = 15, 16, 17).\)}{\(\displaystyle \begin{aligned}[t]\pi_{n+19}(S^{n})&=\langle \omega_n \nu_{n+16}, \bar\sigma_n, \bar\zeta_n\rangle\\&\cong (\mathbb Z/2)^{2} \oplus \mathbb Z/8\quad(n = 15, 16, 17).\end{aligned}\)}
	\StemEntry{\(\displaystyle \pi_{39}(S^{20})=\langle P(\iota_{41}), \bar\sigma_{20}, \bar\zeta_{20}\rangle\cong \mathbb Z_{(2)} \oplus \mathbb Z/2 \oplus \mathbb Z/8.\)}{\(\displaystyle \begin{aligned}[t]\pi_{39}(S^{20})&=\langle P(\iota_{41}), \bar\sigma_{20}, \bar\zeta_{20}\rangle\\&\cong \mathbb Z_{(2)} \oplus \mathbb Z/2 \oplus \mathbb Z/8.\end{aligned}\)}
	\StemEntry{\(\displaystyle \pi_{n+19}(S^{n})=\langle \bar\sigma_n, \bar\zeta_n\rangle\cong \mathbb Z/2 \oplus \mathbb Z/8\quad(n\ge 18 \text{ but } n \neq 20).\)}{\(\displaystyle \begin{aligned}[t]\pi_{n+19}(S^{n})&=\langle \bar\sigma_n, \bar\zeta_n\rangle\\&\cong \mathbb Z/2 \oplus \mathbb Z/8\quad(n\ge 18 \text{ but } n \neq 20).\end{aligned}\)}
	\StemEntry{\(\displaystyle \pi^S_{19}(S^0)=\langle \bar\sigma, \bar\zeta\rangle\cong \mathbb Z/2 \oplus \mathbb Z/8.\)}{\(\displaystyle \begin{aligned}[t]\pi^S_{19}(S^0)&=\langle \bar\sigma, \bar\zeta\rangle\\&\cong \mathbb Z/2 \oplus \mathbb Z/8.\end{aligned}\)}

	\section{The   Mimura-Toda Stem}
	\subsection{Stem 20}\label{app:stem-20}
	\StemEntry{\(\displaystyle \pi_{22}(S^{2})=\langle \eta_2\bar{\mu}', \eta_2\nu'\mu_6\sigma _{15}\rangle\cong \mathbb Z/4 \oplus \mathbb Z/2.\)}{\(\displaystyle \begin{aligned}[t]\pi_{22}(S^{2})&=\langle \eta_2\bar{\mu}', \eta_2\nu'\mu_6\sigma _{15}\rangle\\&\cong \mathbb Z/4 \oplus \mathbb Z/2.\end{aligned}\)}
	\StemEntry{\(\displaystyle \pi_{23}(S^{3})=\langle \nu'\bar{\mu}_6, \nu'\eta_6\mu_7\sigma _{16}\rangle\cong (\mathbb Z/2)^{2}.\)}{\(\displaystyle \begin{aligned}[t]\pi_{23}(S^{3})&=\langle \nu'\bar{\mu}_6, \nu'\eta_6\mu_7\sigma _{16}\rangle\\&\cong (\mathbb Z/2)^{2}.\end{aligned}\)}
	\StemEntry{\(\displaystyle \pi_{24}(S^{4})=\langle (\Sigma \nu')\bar{\mu}_7, (\Sigma \nu')\eta_7\mu_8\sigma _{17}, \nu_4\sigma '\eta_{14}\mu_{15}, \nu^2_4\kappa_{10}, \nu_4\bar{\mu}_7, \nu_4\eta_7\mu_8\sigma _{17}\rangle\cong (\mathbb Z/2)^{6}.\)}{\(\displaystyle \begin{aligned}[t]\pi_{24}(S^{4})&=\langle (\Sigma \nu')\bar{\mu}_7, (\Sigma \nu')\eta_7\mu_8\sigma _{17}, \nu_4\sigma '\eta_{14}\mu_{15}, \nu^2_4\kappa_{10}, \nu_4\bar{\mu}_7,\\&\phantom{=\langle{}}\nu_4\eta_7\mu_8\sigma _{17}\rangle\cong (\mathbb Z/2)^{6}.\end{aligned}\)}
	\StemEntry{\(\displaystyle \pi_{25}(S^{5})=\langle \nu^2_5\kappa_{11}, \nu_5\bar{\mu}_8, \nu_5\eta_8\mu_9\sigma _{18}\rangle\cong (\mathbb Z/2)^{3}.\)}{\(\displaystyle \begin{aligned}[t]\pi_{25}(S^{5})&=\langle \nu^2_5\kappa_{11}, \nu_5\bar{\mu}_8, \nu_5\eta_8\mu_9\sigma _{18}\rangle\\&\cong (\mathbb Z/2)^{3}.\end{aligned}\)}
	\StemEntry{\(\displaystyle \pi_{26}(S^{6})=\langle P(\rho_{13}), \bar{\kappa}'\rangle\cong \mathbb Z/32\oplus \mathbb Z/4.\)}{\(\displaystyle \begin{aligned}[t]\pi_{26}(S^{6})&=\langle P(\rho_{13}), \bar{\kappa}'\rangle\\&\cong \mathbb Z/32\oplus \mathbb Z/4.\end{aligned}\)}
	\StemEntry{\(\displaystyle \pi_{20+n}(S^{n})=\langle \bar{\kappa}_n\rangle\cong \mathbb Z/8\quad(n=7,8,9).\)}{\(\displaystyle \begin{aligned}[t]\pi_{20+n}(S^{n})&=\langle \bar{\kappa}_n\rangle\\&\cong \mathbb Z/8\quad(n=7,8,9).\end{aligned}\)}
	\StemEntry{\(\displaystyle \pi_{30}(S^{10})=\langle \bar{\kappa}_{10}, \beta'\rangle\cong \mathbb Z/8 \oplus \mathbb Z/8.\)}{\(\displaystyle \begin{aligned}[t]\pi_{30}(S^{10})&=\langle \bar{\kappa}_{10}, \beta'\rangle\\&\cong \mathbb Z/8 \oplus \mathbb Z/8.\end{aligned}\)}
	\StemEntry{\(\displaystyle \pi_{31}(S^{11})=\langle \bar{\kappa}_{11}, \theta'\varepsilon _{23}, \beta''\rangle\cong \mathbb Z/8\oplus (\mathbb Z/2)^{2}.\)}{\(\displaystyle \begin{aligned}[t]\pi_{31}(S^{11})&=\langle \bar{\kappa}_{11}, \theta'\varepsilon _{23}, \beta''\rangle\\&\cong \mathbb Z/8\oplus (\mathbb Z/2)^{2}.\end{aligned}\)}
	\StemEntry{\(\displaystyle \pi_{32}(S^{12})=\langle \bar{\kappa}_{12}, (\Sigma \theta')\varepsilon _{24}, P(\mu_{25}), \theta\varepsilon _{24}, \theta\bar{\nu}_{24}, \beta'''\rangle\cong \mathbb Z/8 \oplus (\mathbb Z/2)^{5}.\)}{\(\displaystyle \begin{aligned}[t]\pi_{32}(S^{12})&=\langle \bar{\kappa}_{12}, (\Sigma \theta')\varepsilon _{24}, P(\mu_{25}), \theta\varepsilon _{24}, \theta\bar{\nu}_{24}, \beta'''\rangle\\&\cong \mathbb Z/8 \oplus (\mathbb Z/2)^{5}.\end{aligned}\)}
	\StemEntry{\(\displaystyle \pi_{33}(S^{13})=\langle \bar{\kappa}_{13}, (\Sigma \theta)\varepsilon _{25}, (\Sigma \theta)\bar{\nu}_{25}, \Sigma \beta'''\rangle\cong \mathbb Z/8 \oplus (\mathbb Z/2)^{3}.\)}{\(\displaystyle \begin{aligned}[t]\pi_{33}(S^{13})&=\langle \bar{\kappa}_{13}, (\Sigma \theta)\varepsilon _{25}, (\Sigma \theta)\bar{\nu}_{25}, \Sigma \beta'''\rangle\\&\cong \mathbb Z/8 \oplus (\mathbb Z/2)^{3}.\end{aligned}\)}
	\StemEntry{\(\displaystyle \pi_{34}(S^{14})=\langle P(\sigma _{29}), \bar{\kappa}_{14}\rangle\cong \mathbb Z/16\oplus \mathbb Z/8.\)}{\(\displaystyle \begin{aligned}[t]\pi_{34}(S^{14})&=\langle P(\sigma _{29}), \bar{\kappa}_{14}\rangle\\&\cong \mathbb Z/16\oplus \mathbb Z/8.\end{aligned}\)}
	\StemEntry{\(\displaystyle \pi_{20+n}(S^{n})=\langle \bar{\kappa}_{n}\rangle\cong \mathbb Z/8\quad(n=15,16,17).\)}{\(\displaystyle \begin{aligned}[t]\pi_{20+n}(S^{n})&=\langle \bar{\kappa}_{n}\rangle\\&\cong \mathbb Z/8\quad(n=15,16,17).\end{aligned}\)}
	\StemEntry{\(\displaystyle \pi_{38}(S^{18})=\langle \bar{\kappa}_{18}, P(\nu_{37})\rangle\cong \mathbb Z/8\oplus \mathbb Z/4.\)}{\(\displaystyle \begin{aligned}[t]\pi_{38}(S^{18})&=\langle \bar{\kappa}_{18}, P(\nu_{37})\rangle\\&\cong \mathbb Z/8\oplus \mathbb Z/4.\end{aligned}\)}
	\StemEntry{\(\displaystyle \pi_{39}(S^{19})=\langle \bar{\kappa}_{19}, \bar{\beta}\rangle\cong \mathbb Z/8 \oplus \mathbb Z/2.\)}{\(\displaystyle \begin{aligned}[t]\pi_{39}(S^{19})&=\langle \bar{\kappa}_{19}, \bar{\beta}\rangle\\&\cong \mathbb Z/8 \oplus \mathbb Z/2.\end{aligned}\)}
	\StemEntry{\(\displaystyle \pi_{40}(S^{20})=\langle \bar{\kappa}_{20}, P(\eta_{41}), \bar{\bar{\beta}}\rangle\cong \mathbb Z/8 \oplus (\mathbb Z/2)^{2}.\)}{\(\displaystyle \begin{aligned}[t]\pi_{40}(S^{20})&=\langle \bar{\kappa}_{20}, P(\eta_{41}), \bar{\bar{\beta}}\rangle\\&\cong \mathbb Z/8 \oplus (\mathbb Z/2)^{2}.\end{aligned}\)}
	\StemEntry{\(\displaystyle \pi_{41}(S^{21})=\langle \bar{\kappa}_{21}, P(\iota_{43})\rangle\cong \mathbb Z/8 \oplus \mathbb Z/2.\)}{\(\displaystyle \begin{aligned}[t]\pi_{41}(S^{21})&=\langle \bar{\kappa}_{21}, P(\iota_{43})\rangle\\&\cong \mathbb Z/8 \oplus \mathbb Z/2.\end{aligned}\)}
	\StemEntry{\(\displaystyle \pi_{42}(S^{22})=\langle \bar{\kappa}_{22}\rangle\cong \mathbb Z/8.\)}{\(\displaystyle \begin{aligned}[t]\pi_{42}(S^{22})&=\langle \bar{\kappa}_{22}\rangle\\&\cong \mathbb Z/8.\end{aligned}\)}
	\StemEntry{\(\displaystyle \pi^S_{20}(S^0)=\langle \bar{\kappa}\rangle\cong \mathbb Z/8.\)}{\(\displaystyle \begin{aligned}[t]\pi^S_{20}(S^0)&=\langle \bar{\kappa}\rangle\\&\cong \mathbb Z/8.\end{aligned}\)}
	
	\section{The   Mimura Stems}
	\subsection{Stem 21}\label{app:stem-21}
	\StemEntry{\(\displaystyle \pi_{23}(S^{2})=\langle \eta_2\nu'\bar{\mu}_6, \eta_2\nu'\eta_6\mu_7\sigma _{16}\rangle\cong (\mathbb Z/2)^{2}.\)}{\(\displaystyle \begin{aligned}[t]\pi_{23}(S^{2})&=\langle \eta_2\nu'\bar{\mu}_6, \eta_2\nu'\eta_6\mu_7\sigma _{16}\rangle\\&\cong (\mathbb Z/2)^{2}.\end{aligned}\)}
	\StemEntry{\(\displaystyle \pi_{24}(S^{3})=\langle \nu'\eta_6\bar{\mu}_7\rangle\cong \mathbb Z/2.\)}{\(\displaystyle \begin{aligned}[t]\pi_{24}(S^{3})&=\langle \nu'\eta_6\bar{\mu}_7\rangle\\&\cong \mathbb Z/2.\end{aligned}\)}
	\StemEntry{\(\displaystyle \pi_{25}(S^{4})=\langle (\Sigma \nu')\eta_7\bar{\mu}_8, \nu_4\eta_7\bar{\mu}_8, \nu_4\zeta_7\sigma _{18}\rangle\cong (\mathbb Z/2)^{2} \oplus \mathbb Z/8.\)}{\(\displaystyle \begin{aligned}[t]\pi_{25}(S^{4})&=\langle (\Sigma \nu')\eta_7\bar{\mu}_8, \nu_4\eta_7\bar{\mu}_8,\\&\phantom{=\langle{}}\nu_4\zeta_7\sigma _{18}\rangle\cong (\mathbb Z/2)^{2} \oplus \mathbb Z/8.\end{aligned}\)}
	\StemEntry{\(\displaystyle \pi_{26}(S^{5})=\langle \alpha, \nu_5\eta_8\bar{\mu}_9\rangle\cong (\mathbb Z/2)^{2}.\)}{\(\displaystyle \begin{aligned}[t]\pi_{26}(S^{5})&=\langle \alpha, \nu_5\eta_8\bar{\mu}_9\rangle\\&\cong (\mathbb Z/2)^{2}.\end{aligned}\)}
	\StemEntry{\(\displaystyle \pi_{27}(S^{6})=\langle \eta_6\bar{\kappa}_7\rangle\cong \mathbb Z/2.\)}{\(\displaystyle \begin{aligned}[t]\pi_{27}(S^{6})&=\langle \eta_6\bar{\kappa}_7\rangle\\&\cong \mathbb Z/2.\end{aligned}\)}
	\StemEntry{\(\displaystyle \pi_{28}(S^{7})=\langle \eta_7 \bar{\kappa}_8, \sigma '\kappa_{14}\rangle\cong (\mathbb Z/2)^{2}.\)}{\(\displaystyle \begin{aligned}[t]\pi_{28}(S^{7})&=\langle \eta_7 \bar{\kappa}_8, \sigma '\kappa_{14}\rangle\\&\cong (\mathbb Z/2)^{2}.\end{aligned}\)}
	\StemEntry{\(\displaystyle \pi_{29}(S^{8})=\langle \eta_8\bar{\kappa}_9, (\Sigma \sigma ')\kappa_{15}, \sigma _8\kappa_{15}, \sigma ^3_8\rangle\cong (\mathbb Z/2)^{3} \oplus \mathbb Z/4.\)}{\(\displaystyle \begin{aligned}[t]\pi_{29}(S^{8})&=\langle \eta_8\bar{\kappa}_9, (\Sigma \sigma ')\kappa_{15}, \sigma _8\kappa_{15},\\&\phantom{=\langle{}}\sigma ^3_8\rangle\cong (\mathbb Z/2)^{3} \oplus \mathbb Z/4.\end{aligned}\)}
	\StemEntry{\(\displaystyle \pi_{21+n}(S^{n})=\langle \eta_n\bar{\kappa}_{n+1}, \sigma _n\kappa_{n+7}, \sigma ^3_n\rangle\cong (\mathbb Z/2)^{3}\quad(n=9, 10).\)}{\(\displaystyle \begin{aligned}[t]\pi_{21+n}(S^{n})&=\langle \eta_n\bar{\kappa}_{n+1}, \sigma _n\kappa_{n+7}, \sigma ^3_n\rangle\\&\cong (\mathbb Z/2)^{3}\quad(n=9, 10).\end{aligned}\)}
	\StemEntry{\(\displaystyle \pi_{32}(S^{11})=\langle \eta_{11}\bar{\kappa}_{12}, \sigma _{11}\kappa_{18}, \sigma ^3_{11}, \theta'\mu_{23}\rangle\cong (\mathbb Z/2)^{4}.\)}{\(\displaystyle \begin{aligned}[t]\pi_{32}(S^{11})&=\langle \eta_{11}\bar{\kappa}_{12}, \sigma _{11}\kappa_{18}, \sigma ^3_{11}, \theta'\mu_{23}\rangle\\&\cong (\mathbb Z/2)^{4}.\end{aligned}\)}
	\StemEntry{\(\displaystyle \pi_{33}(S^{12})=\langle \eta_{12}\bar{\kappa}_{13}, \sigma _{12}\kappa_{19}, \sigma ^3_{12}, (\Sigma \theta')\mu_{24}, \theta \mu_{24}\rangle\cong (\mathbb Z/2)^{5}.\)}{\(\displaystyle \begin{aligned}[t]\pi_{33}(S^{12})&=\langle \eta_{12}\bar{\kappa}_{13}, \sigma _{12}\kappa_{19}, \sigma ^3_{12}, (\Sigma \theta')\mu_{24},\\&\phantom{=\langle{}}\theta \mu_{24}\rangle\cong (\mathbb Z/2)^{5}.\end{aligned}\)}
	\StemEntry{\(\displaystyle \pi_{34}(S^{13})=\langle \eta_{13}\bar{\kappa}_{14}, \sigma ^3_{13}, (\Sigma \theta)\mu_{25}, \lambda\nu_{31}\rangle\cong (\mathbb Z/2)^{3} \oplus \mathbb Z/4.\)}{\(\displaystyle \begin{aligned}[t]\pi_{34}(S^{13})&=\langle \eta_{13}\bar{\kappa}_{14}, \sigma ^3_{13}, (\Sigma \theta)\mu_{25}, \lambda\nu_{31}\rangle\\&\cong (\mathbb Z/2)^{3} \oplus \mathbb Z/4.\end{aligned}\)}
	\StemEntry{\(\displaystyle \pi_{35}(S^{14})=\langle \eta_{14}\bar{\kappa}_{15}, \sigma ^3_{14}, (\Sigma \lambda)\nu_{32}\rangle\cong (\mathbb Z/2)^{2} \oplus \mathbb Z/4.\)}{\(\displaystyle \begin{aligned}[t]\pi_{35}(S^{14})&=\langle \eta_{14}\bar{\kappa}_{15}, \sigma ^3_{14}, (\Sigma \lambda)\nu_{32}\rangle\\&\cong (\mathbb Z/2)^{2} \oplus \mathbb Z/4.\end{aligned}\)}
	\StemEntry{\(\displaystyle \pi_{36}(S^{15})=\langle \eta_{15}\bar{\kappa}_{16}, \sigma ^3_{15}, (\Sigma ^2\lambda)\nu_{33}\rangle\cong (\mathbb Z/2)^{3}.\)}{\(\displaystyle \begin{aligned}[t]\pi_{36}(S^{15})&=\langle \eta_{15}\bar{\kappa}_{16}, \sigma ^3_{15}, (\Sigma ^2\lambda)\nu_{33}\rangle\\&\cong (\mathbb Z/2)^{3}.\end{aligned}\)}
	\StemEntry{\(\displaystyle \pi_{37}(S^{16})=\langle \eta_{16}\bar{\kappa}_{17}, \sigma ^3_{16}, (\Sigma ^3\lambda)\nu_{34}, \nu^\ast_{16}\nu_{34}\rangle\cong (\mathbb Z/2)^{4}.\)}{\(\displaystyle \begin{aligned}[t]\pi_{37}(S^{16})&=\langle \eta_{16}\bar{\kappa}_{17}, \sigma ^3_{16}, (\Sigma ^3\lambda)\nu_{34}, \nu^\ast_{16}\nu_{34}\rangle\\&\cong (\mathbb Z/2)^{4}.\end{aligned}\)}
	\StemEntry{\(\displaystyle \pi_{38}(S^{17})=\langle \eta_{17}\bar{\kappa}_{18}, \sigma ^3_{17}, \nu^\ast_{17}\nu_{35}\rangle\cong (\mathbb Z/2)^{3}.\)}{\(\displaystyle \begin{aligned}[t]\pi_{38}(S^{17})&=\langle \eta_{17}\bar{\kappa}_{18}, \sigma ^3_{17}, \nu^\ast_{17}\nu_{35}\rangle\\&\cong (\mathbb Z/2)^{3}.\end{aligned}\)}
	\StemEntry{\(\displaystyle \pi_{39}(S^{18})=\langle \eta_{18}\bar{\kappa}_{19}, \sigma ^3_{18}, \nu^\ast_{18}\nu_{36}\rangle\cong (\mathbb Z/2)^{3}.\)}{\(\displaystyle \begin{aligned}[t]\pi_{39}(S^{18})&=\langle \eta_{18}\bar{\kappa}_{19}, \sigma ^3_{18}, \nu^\ast_{18}\nu_{36}\rangle\\&\cong (\mathbb Z/2)^{3}.\end{aligned}\)}
	\StemEntry{\(\displaystyle \pi_{40}(S^{19})=\langle \eta_{19}\bar{\kappa}_{20}, \sigma ^3_{19}, \nu^\ast_{19}\nu_{37}, \bar{\beta}\eta_{39}\rangle\cong (\mathbb Z/2)^{4}.\)}{\(\displaystyle \begin{aligned}[t]\pi_{40}(S^{19})&=\langle \eta_{19}\bar{\kappa}_{20}, \sigma ^3_{19}, \nu^\ast_{19}\nu_{37}, \bar{\beta}\eta_{39}\rangle\\&\cong (\mathbb Z/2)^{4}.\end{aligned}\)}
	\StemEntry{\(\displaystyle \pi_{41}(S^{20})=\langle \eta_{20}\bar{\kappa}_{21}, \sigma ^3_{20}, (\Sigma \bar{\beta})\eta_{40}, \bar{\bar{\beta}}\eta_{40}\rangle\cong (\mathbb Z/2)^{4}.\)}{\(\displaystyle \begin{aligned}[t]\pi_{41}(S^{20})&=\langle \eta_{20}\bar{\kappa}_{21}, \sigma ^3_{20}, (\Sigma \bar{\beta})\eta_{40}, \bar{\bar{\beta}}\eta_{40}\rangle\\&\cong (\mathbb Z/2)^{4}.\end{aligned}\)}
	\StemEntry{\(\displaystyle \pi_{42}(S^{21})=\langle \eta_{21}\bar{\kappa}_{22}, \sigma ^3_{21}, (\Sigma \bar{\bar{\beta}})\eta_{41}\rangle\cong (\mathbb Z/2)^{3}.\)}{\(\displaystyle \begin{aligned}[t]\pi_{42}(S^{21})&=\langle \eta_{21}\bar{\kappa}_{22}, \sigma ^3_{21}, (\Sigma \bar{\bar{\beta}})\eta_{41}\rangle\\&\cong (\mathbb Z/2)^{3}.\end{aligned}\)}
	\StemEntry{\(\displaystyle \pi_{43}(S^{22})=\langle P(\iota_{45}), \eta_{22}\bar{\kappa}_{23}, \sigma ^3_{22}\rangle\cong \mathbb Z_{(2)} \oplus (\mathbb Z/2)^{2}.\)}{\(\displaystyle \begin{aligned}[t]\pi_{43}(S^{22})&=\langle P(\iota_{45}), \eta_{22}\bar{\kappa}_{23}, \sigma ^3_{22}\rangle\\&\cong \mathbb Z_{(2)} \oplus (\mathbb Z/2)^{2}.\end{aligned}\)}
	\StemEntry{\(\displaystyle \pi_{44}(S^{23})=\langle \eta_{23}\bar{\kappa}_{24}, \sigma ^3_{23}\rangle\cong (\mathbb Z/2)^{2}.\)}{\(\displaystyle \begin{aligned}[t]\pi_{44}(S^{23})&=\langle \eta_{23}\bar{\kappa}_{24}, \sigma ^3_{23}\rangle\\&\cong (\mathbb Z/2)^{2}.\end{aligned}\)}
	\StemEntry{\(\displaystyle \pi^S_{21}(S^0)=\langle \eta\bar{\kappa}, \sigma ^3\rangle\cong \mathbb Z/2 \oplus \mathbb Z/2.\)}{\(\displaystyle \begin{aligned}[t]\pi^S_{21}(S^0)&=\langle \eta\bar{\kappa}, \sigma ^3\rangle\\&\cong \mathbb Z/2 \oplus \mathbb Z/2.\end{aligned}\)}
	\subsection{Stem 22}\label{app:stem-22}
	\StemEntry{\(\displaystyle \pi_{24}(S^{2})=\langle \eta_2\nu'\eta_6\bar{\mu}_7\rangle\cong \mathbb Z/2.\)}{\(\displaystyle \begin{aligned}[t]\pi_{24}(S^{2})&=\langle \eta_2\nu'\eta_6\bar{\mu}_7\rangle\\&\cong \mathbb Z/2.\end{aligned}\)}
	\StemEntry{\(\displaystyle \pi_{25}(S^{3})=\langle \varepsilon _3\kappa_{11}\rangle\cong \mathbb Z/2.\)}{\(\displaystyle \begin{aligned}[t]\pi_{25}(S^{3})&=\langle \varepsilon _3\kappa_{11}\rangle\\&\cong \mathbb Z/2.\end{aligned}\)}
	
	\StemEntry
	{\(\displaystyle
		\pi_{26}(S^{4})
		=\langle \nu_{4}\bar{\zeta}_{7}, \varepsilon_{4}\kappa_{12}, \nu_{4}\bar{\sigma}_{7}\rangle
		\cong \mathbb Z/8 \oplus (\mathbb Z/2)^{2}.
		\)}
	{\(\displaystyle
		\begin{aligned}[t]
			\pi_{26}(S^{4})
			&=\langle \nu_{4}\bar{\zeta}_{7}, \varepsilon_{4}\kappa_{12}, \nu_{4}\bar{\sigma}_{7}\rangle\\
			&\cong \mathbb Z/8 \oplus (\mathbb Z/2)^{2}.
		\end{aligned}
		\)}
	
	\StemEntry
	{\(\displaystyle
		\pi_{27}(S^{5})
		=\langle \nu_{5}\bar{\zeta}_{8}, \varepsilon_{5}\kappa_{13}, \nu_{5}\bar{\sigma}_{8}\rangle
		\cong (\mathbb Z/2)^{3}.
		\)}
	{\(\displaystyle
		\begin{aligned}[t]
			\pi_{27}(S^{5})
			&=\langle \nu_{5}\bar{\zeta}_{8}, \varepsilon_{5}\kappa_{13}, \nu_{5}\bar{\sigma}_{8}\rangle\\
			&\cong (\mathbb Z/2)^{3}.
		\end{aligned}
		\)}

	\StemEntry{\(\displaystyle \pi_{28}(S^{6})=\langle \rho'''\sigma _{21}, \bar{\nu}_6\kappa_{14}, \varepsilon _6\kappa_{14}, \nu_6\bar{\sigma }_9\rangle\cong \mathbb Z/4 \oplus (\mathbb Z/2)^{3}.\)}{\(\displaystyle \begin{aligned}[t]\pi_{28}(S^{6})&=\langle \rho'''\sigma _{21}, \bar{\nu}_6\kappa_{14}, \varepsilon _6\kappa_{14}, \nu_6\bar{\sigma }_9\rangle\\&\cong \mathbb Z/4 \oplus (\mathbb Z/2)^{3}.\end{aligned}\)}
	\StemEntry{\(\displaystyle \pi_{29}(S^{7})=\langle \sigma '\rho_{14}, \bar{\nu}_7\kappa_{15}, \varepsilon _7\kappa_{15}, \nu_7\bar{\sigma }_{10}\rangle\cong \mathbb Z/8 \oplus (\mathbb Z/2)^{3}.\)}{\(\displaystyle \begin{aligned}[t]\pi_{29}(S^{7})&=\langle \sigma '\rho_{14}, \bar{\nu}_7\kappa_{15}, \varepsilon _7\kappa_{15}, \nu_7\bar{\sigma }_{10}\rangle\\&\cong \mathbb Z/8 \oplus (\mathbb Z/2)^{3}.\end{aligned}\)}
	\StemEntry{\(\displaystyle \pi_{30}(S^{8})=\langle (\Sigma \sigma ')\rho_{15}, \bar{\nu}_8\kappa_{16}, \varepsilon _8\kappa_{16}, \nu_8\bar{\sigma }_{11}, \sigma _8\rho_{15}, \sigma _8\bar{\varepsilon }_{15}\rangle\cong \mathbb Z/8 \oplus \mathbb Z/2 \oplus \mathbb Z/2 \oplus \mathbb Z/2 \oplus \mathbb Z/32 \oplus \mathbb Z/2.\)}{\(\displaystyle \begin{aligned}[t]\pi_{30}(S^{8})&=\langle (\Sigma \sigma ')\rho_{15}, \bar{\nu}_8\kappa_{16}, \varepsilon _8\kappa_{16}, \nu_8\bar{\sigma }_{11}, \sigma _8\rho_{15},\\&\phantom{=\langle{}}\sigma _8\bar{\varepsilon }_{15}\rangle\cong \mathbb Z/8 \oplus \mathbb Z/2 \oplus \mathbb Z/2 \oplus \mathbb Z/2 \oplus \mathbb Z/32 \oplus \mathbb Z/2.\end{aligned}\)}
	\StemEntry{\(\displaystyle \pi_{31}(S^{9})=\langle \sigma_9\rho_{16}, \varepsilon _9\kappa_{17}, \nu_9\bar{\sigma }_{12}, \sigma_9\bar{\varepsilon}_{16}\rangle\cong \mathbb Z/16 \oplus (\mathbb Z/2)^{3}.\)}{\(\displaystyle \begin{aligned}[t]\pi_{31}(S^{9})&=\langle \sigma_9\rho_{16}, \varepsilon _9\kappa_{17}, \nu_9\bar{\sigma }_{12}, \sigma_9\bar{\varepsilon}_{16}\rangle\\&\cong \mathbb Z/16 \oplus (\mathbb Z/2)^{3}.\end{aligned}\)}
	\StemEntry{\(\displaystyle \pi_{22+n}(S^{n})=\langle \sigma _{n}\rho_{n+7}, \varepsilon _{n}\kappa_{n+8}, \nu_n\bar{\sigma}_{n+3}\rangle\cong \mathbb Z/16 \oplus (\mathbb Z/2)^{2}\quad(n=10, 11).\)}{\(\displaystyle \begin{aligned}[t]\pi_{22+n}(S^{n})&=\langle \sigma _{n}\rho_{n+7}, \varepsilon _{n}\kappa_{n+8}, \nu_n\bar{\sigma}_{n+3}\rangle\\&\cong \mathbb Z/16 \oplus (\mathbb Z/2)^{2}\quad(n=10, 11).\end{aligned}\)}\vspace{-0.8em}\;\\
	The sign is chosen so that $\sigma_{12}\rho_{19}\pm2{\sigma^{\ast\prime\prime\prime}}$ has order $4.$ 
	\StemEntry{\(\displaystyle \pi_{34}(S^{12})=\langle {\sigma ^{\ast\prime\prime\prime}}, \sigma _{12}\rho_{19}\pm2{\sigma ^{\ast\prime\prime\prime}}, \varepsilon _{12}\kappa_{20}, \nu_{12}\bar{\sigma }_{15}\rangle\cong \mathbb Z/32 \oplus \mathbb Z/4 \oplus (\mathbb Z/2)^{2}.\)}{\(\displaystyle \begin{aligned}[t]\pi_{34}(S^{12})&=\langle {\sigma ^{\ast\prime\prime\prime}}, \sigma _{12}\rho_{19}\pm2{\sigma ^{\ast\prime\prime\prime}}, \varepsilon _{12}\kappa_{20},\\&\phantom{=\langle{}}\nu_{12}\bar{\sigma }_{15}\rangle\cong \mathbb Z/32 \oplus \mathbb Z/4 \oplus (\mathbb Z/2)^{2}.\end{aligned}\)}
	\StemEntry{\(\displaystyle \pi_{35}(S^{13})=\langle \rho_{13}\sigma _{28}, \varepsilon _{13}\kappa_{21}, \nu_{13}\bar{\sigma }_{16}\rangle\cong \mathbb Z/16 \oplus (\mathbb Z/2)^{2}.\)}{\(\displaystyle \begin{aligned}[t]\pi_{35}(S^{13})&=\langle \rho_{13}\sigma _{28}, \varepsilon _{13}\kappa_{21}, \nu_{13}\bar{\sigma }_{16}\rangle\\&\cong \mathbb Z/16 \oplus (\mathbb Z/2)^{2}.\end{aligned}\)}
	\StemEntry{\(\displaystyle \pi_{36}(S^{14})=\langle {\sigma ^{\ast\prime\prime}}, \omega_{14}\nu^2_{30}, \varepsilon _{14}\kappa_{22}, \nu_{14}\bar{\sigma }_{17}\rangle\cong \mathbb Z/16 \oplus (\mathbb Z/2)^{3}.\)}{\(\displaystyle \begin{aligned}[t]\pi_{36}(S^{14})&=\langle {\sigma ^{\ast\prime\prime}}, \omega_{14}\nu^2_{30}, \varepsilon _{14}\kappa_{22}, \nu_{14}\bar{\sigma }_{17}\rangle\\&\cong \mathbb Z/16 \oplus (\mathbb Z/2)^{3}.\end{aligned}\)}
	\StemEntry{\(\displaystyle \pi_{37}(S^{15})=\langle {\sigma ^{\ast\prime}}, \omega_{15}\nu^2_{31}, \varepsilon _{15}\kappa_{23}, \nu_{15}\bar{\sigma }_{18}\rangle\cong \mathbb Z/16 \oplus (\mathbb Z/2)^{3}.\)}{\(\displaystyle \begin{aligned}[t]\pi_{37}(S^{15})&=\langle {\sigma ^{\ast\prime}}, \omega_{15}\nu^2_{31}, \varepsilon _{15}\kappa_{23}, \nu_{15}\bar{\sigma }_{18}\rangle\\&\cong \mathbb Z/16 \oplus (\mathbb Z/2)^{3}.\end{aligned}\)}
	\StemEntry{\(\displaystyle \pi_{38}(S^{16})=\langle \sigma ^\ast_{16}, \Sigma {\sigma ^{\ast\prime}}, \omega_{16}\nu^2_{32}, \varepsilon _{16}\kappa_{24}, \nu_{16}\bar{\sigma }_{19}\rangle\cong \mathbb Z/16 \oplus \mathbb Z/16 \oplus (\mathbb Z/2)^{3}.\)}{\(\displaystyle \begin{aligned}[t]\pi_{38}(S^{16})&=\langle \sigma ^\ast_{16}, \Sigma {\sigma ^{\ast\prime}}, \omega_{16}\nu^2_{32}, \varepsilon _{16}\kappa_{24}, \nu_{16}\bar{\sigma }_{19}\rangle\\&\cong \mathbb Z/16 \oplus \mathbb Z/16 \oplus (\mathbb Z/2)^{3}.\end{aligned}\)}
	\StemEntry{\(\displaystyle \pi_{39}(S^{17})=\langle \sigma ^\ast_{17}, \omega_{17}\nu^2_{33}, \varepsilon _{17}\kappa_{25}, \nu_{17}\bar{\sigma }_{20}\rangle\cong \mathbb Z/16 \oplus (\mathbb Z/2)^{3}.\)}{\(\displaystyle \begin{aligned}[t]\pi_{39}(S^{17})&=\langle \sigma ^\ast_{17}, \omega_{17}\nu^2_{33}, \varepsilon _{17}\kappa_{25}, \nu_{17}\bar{\sigma }_{20}\rangle\\&\cong \mathbb Z/16 \oplus (\mathbb Z/2)^{3}.\end{aligned}\)}
	\StemEntry{\(\displaystyle \pi_{n+22}(S^{n})=\langle \sigma ^\ast_n, \varepsilon _n\kappa_{n+8}, \nu_n\bar{\sigma }_{n+3}\rangle\cong \mathbb Z/16 \oplus (\mathbb Z/2)^{2}\quad(n=18, 19).\)}{\(\displaystyle \begin{aligned}[t]\pi_{n+22}(S^{n})&=\langle \sigma ^\ast_n, \varepsilon _n\kappa_{n+8}, \nu_n\bar{\sigma }_{n+3}\rangle\\&\cong \mathbb Z/16 \oplus (\mathbb Z/2)^{2}\quad(n=18, 19).\end{aligned}\)}
	\StemEntry{\(\displaystyle \pi_{42}(S^{20})=\langle \sigma ^\ast_{20}, \varepsilon _{20}\kappa_{28}, \nu_{20}\bar{\sigma }_{23}, P(\nu_{41})+2\sigma ^\ast_{20}\rangle\cong \mathbb Z/16 \oplus \mathbb Z/2 \oplus \mathbb Z/2 \oplus \mathbb Z/4.\)}{\(\displaystyle \begin{aligned}[t]\pi_{42}(S^{20})&=\langle \sigma ^\ast_{20}, \varepsilon _{20}\kappa_{28}, \nu_{20}\bar{\sigma }_{23}, P(\nu_{41})+2\sigma ^\ast_{20}\rangle\\&\cong \mathbb Z/16 \oplus \mathbb Z/2 \oplus \mathbb Z/2 \oplus \mathbb Z/4.\end{aligned}\)}
	\StemEntry{\(\displaystyle \pi_{43}(S^{21})=\langle \sigma ^\ast_{21}, \varepsilon _{21}\kappa_{29}, \nu_{21}\bar{\sigma }_{24}\rangle\cong \mathbb Z/8\oplus (\mathbb Z/2)^{2}.\)}{\(\displaystyle \begin{aligned}[t]\pi_{43}(S^{21})&=\langle \sigma ^\ast_{21}, \varepsilon _{21}\kappa_{29}, \nu_{21}\bar{\sigma }_{24}\rangle\\&\cong \mathbb Z/8\oplus (\mathbb Z/2)^{2}.\end{aligned}\)}
	\StemEntry{\(\displaystyle \pi_{44}(S^{22})=\langle \sigma ^\ast_{22}, \varepsilon _{22}\kappa_{30}, \nu_{22}\bar{\sigma }_{25}\rangle\cong \mathbb Z/4\oplus (\mathbb Z/2)^{2}.\)}{\(\displaystyle \begin{aligned}[t]\pi_{44}(S^{22})&=\langle \sigma ^\ast_{22}, \varepsilon _{22}\kappa_{30}, \nu_{22}\bar{\sigma }_{25}\rangle\\&\cong \mathbb Z/4\oplus (\mathbb Z/2)^{2}.\end{aligned}\)}
	\StemEntry{\(\displaystyle \pi_{45}(S^{23})=\langle \sigma ^\ast_{23}, \varepsilon _{23}\kappa_{31}, \nu_{23}\bar{\sigma }_{26}\rangle\cong (\mathbb Z/2)^{3}.\)}{\(\displaystyle \begin{aligned}[t]\pi_{45}(S^{23})&=\langle \sigma ^\ast_{23}, \varepsilon _{23}\kappa_{31}, \nu_{23}\bar{\sigma }_{26}\rangle\\&\cong (\mathbb Z/2)^{3}.\end{aligned}\)}
	\StemEntry{\(\displaystyle \pi_{46}(S^{24})=\langle \varepsilon _{24}\kappa_{32}, \nu_{24}\bar{\sigma }_{27}\rangle\cong (\mathbb Z/2)^{2}.\)}{\(\displaystyle \begin{aligned}[t]\pi_{46}(S^{24})&=\langle \varepsilon _{24}\kappa_{32}, \nu_{24}\bar{\sigma }_{27}\rangle\\&\cong (\mathbb Z/2)^{2}.\end{aligned}\)}
	\StemEntry{\(\displaystyle \pi^S_{22}(S^0)=\langle \varepsilon \kappa, \nu\bar{\sigma }\rangle\cong \mathbb Z/2 \oplus \mathbb Z/2.\)}{\(\displaystyle \begin{aligned}[t]\pi^S_{22}(S^0)&=\langle \varepsilon \kappa, \nu\bar{\sigma }\rangle\\&\cong \mathbb Z/2 \oplus \mathbb Z/2.\end{aligned}\)}
	
	\section{The   Mimura-Mori-Oda Stems}
	\subsection{Stem 23}\label{app:stem-23}
	For the ordered generators in dimension $6,$   set
	$b_6=P(\xi_{13})+yP(\lambda),$   with $y\in\{0,1\}$ chosen so that
	$4b_6=0.$ 
	\StemEntry{\(\displaystyle \pi_{25}(S^{2})=\langle \eta_2\varepsilon _3\kappa_{11}\rangle\cong \mathbb Z/2.\)}{\(\displaystyle \begin{aligned}[t]\pi_{25}(S^{2})&=\langle \eta_2\varepsilon _3\kappa_{11}\rangle\\&\cong \mathbb Z/2.\end{aligned}\)}
	\StemEntry{\(\displaystyle \pi_{26}(S^{3})=\langle \bar{\alpha}\rangle\cong \mathbb Z/4.\)}{\(\displaystyle \begin{aligned}[t]\pi_{26}(S^{3})&=\langle \bar{\alpha}\rangle\\&\cong \mathbb Z/4.\end{aligned}\)}
	\StemEntry{\(\displaystyle \pi_{27}(S^{4})=\langle \Sigma \bar{\alpha}, \nu_4\bar{\kappa}_7\rangle\cong \mathbb Z/4 \oplus \mathbb Z/8.\)}{\(\displaystyle \begin{aligned}[t]\pi_{27}(S^{4})&=\langle \Sigma \bar{\alpha}, \nu_4\bar{\kappa}_7\rangle\\&\cong \mathbb Z/4 \oplus \mathbb Z/8.\end{aligned}\)}
	\StemEntry{\(\displaystyle \pi_{28}(S^{5})=\langle \nu_5\bar{\kappa}_8, \bar{\rho}''', \phi_5\rangle\cong \mathbb Z/8 \oplus (\mathbb Z/2)^{2}.\)}{\(\displaystyle \begin{aligned}[t]\pi_{28}(S^{5})&=\langle \nu_5\bar{\kappa}_8, \bar{\rho}''', \phi_5\rangle\\&\cong \mathbb Z/8 \oplus (\mathbb Z/2)^{2}.\end{aligned}\)}
	\StemEntry{\(\displaystyle \pi_{29}(S^{6})=\langle \nu_6\bar{\kappa}_9, P(\lambda), \bar{\rho}'', b_6, \phi_6\rangle\cong \mathbb Z/8 \oplus \mathbb Z/8 \oplus \mathbb Z/4 \oplus \mathbb Z/4 \oplus \mathbb Z/2.\)}{\(\displaystyle \begin{aligned}[t]\pi_{29}(S^{6})&=\langle \nu_6\bar{\kappa}_9, P(\lambda), \bar{\rho}'', b_6, \phi_6\rangle\\&\cong \mathbb Z/8 \oplus \mathbb Z/8 \oplus \mathbb Z/4 \oplus \mathbb Z/4 \oplus \mathbb Z/2.\end{aligned}\)}
	\StemEntry{\(\displaystyle \pi_{30}(S^{7})=\langle \nu_7\bar{\kappa}_{10}, \bar{\rho}', \phi_7, \bar{\kappa}_7\nu_{27}-\nu_7\bar{\kappa}_{10}, \sigma '\sigma _{14}\mu_{21}, \sigma '\omega_{14}\rangle\cong \mathbb Z/8 \oplus \mathbb Z/8 \oplus (\mathbb Z/2)^{4}.\)}{\(\displaystyle \begin{aligned}[t]\pi_{30}(S^{7})&=\langle \nu_7\bar{\kappa}_{10}, \bar{\rho}', \phi_7, \bar{\kappa}_7\nu_{27}-\nu_7\bar{\kappa}_{10}, \sigma '\sigma _{14}\mu_{21},\\&\phantom{=\langle{}}\sigma '\omega_{14}\rangle\cong \mathbb Z/8 \oplus \mathbb Z/8 \oplus (\mathbb Z/2)^{4}.\end{aligned}\)}
	\StemEntry{\(\displaystyle \pi_{31}(S^{8})=\langle \nu_8\bar{\kappa}_{11}, \Sigma \bar{\rho}', \phi_8, \bar{\kappa}_8\nu_{28}-\nu_8\bar{\kappa}_{11}, (\Sigma \sigma ')\sigma _{15}\mu_{22}, (\Sigma \sigma ')\omega_{15}, \sigma ^2_8\mu_{22}, \sigma _8\omega_{15}, \sigma _8{\eta^{\ast\prime}}\rangle\cong \mathbb Z/8 \oplus \mathbb Z/8 \oplus (\mathbb Z/2)^{7}.\)}{\(\displaystyle \begin{aligned}[t]\pi_{31}(S^{8})&=\langle \nu_8\bar{\kappa}_{11}, \Sigma \bar{\rho}', \phi_8, \bar{\kappa}_8\nu_{28}-\nu_8\bar{\kappa}_{11}, (\Sigma \sigma ')\sigma _{15}\mu_{22}, (\Sigma \sigma ')\omega_{15}, \sigma ^2_8\mu_{22},\\&\phantom{=\langle{}}\sigma _8\omega_{15}, \sigma _8{\eta^{\ast\prime}}\rangle\cong \mathbb Z/8 \oplus \mathbb Z/8 \oplus (\mathbb Z/2)^{7}.\end{aligned}\)}
	\StemEntry{\(\displaystyle \pi_{32}(S^{9})=\langle \bar{\rho}_9, \nu_9\ct_{12}, \phi_9, \bar{\kappa}_9\nu_{29}-\nu_9\bar{\kappa}_{12}, \sigma ^2_9\mu_{23}, \sigma _9\omega_{16}\rangle\cong \mathbb Z/16 \oplus \mathbb Z/8 \oplus (\mathbb Z/2)^{4}.\)}{\(\displaystyle \begin{aligned}[t]\pi_{32}(S^{9})&=\langle \bar{\rho}_9, \nu_9\ct_{12}, \phi_9, \bar{\kappa}_9\nu_{29}-\nu_9\bar{\kappa}_{12}, \sigma ^2_9\mu_{23},\\&\phantom{=\langle{}}\sigma _9\omega_{16}\rangle\cong \mathbb Z/16 \oplus \mathbb Z/8 \oplus (\mathbb Z/2)^{4}.\end{aligned}\)}
	\StemEntry{\(\displaystyle \pi_{33}(S^{10})=\langle \bar{\rho}_{10}, \nu_{10}\bar{\kappa}_{13}, \phi_{10}, \bar{\kappa}_{10}\nu_{30}-\nu_{10}\bar{\kappa}_{13}, \psi_{10}\rangle\cong \mathbb Z/16 \oplus \mathbb Z/8 \oplus \mathbb Z/2 \oplus \mathbb Z/2 \oplus \mathbb Z/4.\)}{\(\displaystyle \begin{aligned}[t]\pi_{33}(S^{10})&=\langle \bar{\rho}_{10}, \nu_{10}\bar{\kappa}_{13}, \phi_{10}, \bar{\kappa}_{10}\nu_{30}-\nu_{10}\bar{\kappa}_{13},\\&\phantom{=\langle{}}\psi_{10}\rangle\cong \mathbb Z/16 \oplus \mathbb Z/8 \oplus \mathbb Z/2 \oplus \mathbb Z/2 \oplus \mathbb Z/4.\end{aligned}\)}
	\StemEntry{\(\displaystyle \pi_{n+23}(S^{n})=\langle \bar{\rho}_n, \nu_n\bar{\kappa}_{n+3}, \phi_n, \psi_n\rangle\cong \mathbb Z/16 \oplus \mathbb Z/8 \oplus (\mathbb Z/2)^{2}\quad(n=11, 12, 13, 14).\)}{\(\displaystyle \begin{aligned}[t]\pi_{n+23}(S^{n})&=\langle \bar{\rho}_n, \nu_n\bar{\kappa}_{n+3}, \phi_n, \psi_n\rangle\\&\cong \mathbb Z/16 \oplus \mathbb Z/8 \oplus (\mathbb Z/2)^{2}\quad(n=11, 12, 13, 14).\end{aligned}\)}
	\StemEntry{\(\displaystyle \pi_{38}(S^{15})=\langle \bar{\rho}_{15}, \nu_{15}\bar{\kappa}_{18}, \phi_{15}, \psi_{15}, \bar{\varepsilon }^{\ast\prime}, \bar{\nu}^{\ast\prime}\rangle\cong \mathbb Z/16 \oplus \mathbb Z/8 \oplus (\mathbb Z/2)^{4}.\)}{\(\displaystyle \begin{aligned}[t]\pi_{38}(S^{15})&=\langle \bar{\rho}_{15}, \nu_{15}\bar{\kappa}_{18}, \phi_{15}, \psi_{15}, \bar{\varepsilon }^{\ast\prime}, \bar{\nu}^{\ast\prime}\rangle\\&\cong \mathbb Z/16 \oplus \mathbb Z/8 \oplus (\mathbb Z/2)^{4}.\end{aligned}\)}
	\StemEntry{\(\displaystyle \pi_{39}(S^{16})=\langle \bar{\rho}_{16}, \nu_{16}\bar{\kappa}_{19}, \phi_{16}, \psi_{16}, {\Sigma \bar{\varepsilon }^{\ast\prime}}, {\Sigma \bar{\nu}^{\ast\prime}}, \bar{\varepsilon }^\ast_{16}, \bar{\nu}^\ast_{16}\rangle\cong \mathbb Z/16 \oplus \mathbb Z/8 \oplus (\mathbb Z/2)^{6}.\)}{\(\displaystyle \begin{aligned}[t]\pi_{39}(S^{16})&=\langle \bar{\rho}_{16}, \nu_{16}\bar{\kappa}_{19}, \phi_{16}, \psi_{16}, {\Sigma \bar{\varepsilon }^{\ast\prime}}, {\Sigma \bar{\nu}^{\ast\prime}}, \bar{\varepsilon }^\ast_{16}, \bar{\nu}^\ast_{16}\rangle\\&\cong \mathbb Z/16 \oplus \mathbb Z/8 \oplus (\mathbb Z/2)^{6}.\end{aligned}\)}
	\StemEntry{\(\displaystyle \pi_{40}(S^{17})=\langle \bar{\rho}_{17}, \nu_{17}\bar{\kappa}_{20}, \phi_{17}, \psi_{17}, \bar{\varepsilon }^\ast_{17}, \bar{\nu}^\ast_{17}\rangle\cong \mathbb Z/16 \oplus \mathbb Z/8 \oplus (\mathbb Z/2)^{4}.\)}{\(\displaystyle \begin{aligned}[t]\pi_{40}(S^{17})&=\langle \bar{\rho}_{17}, \nu_{17}\bar{\kappa}_{20}, \phi_{17}, \psi_{17}, \bar{\varepsilon }^\ast_{17}, \bar{\nu}^\ast_{17}\rangle\\&\cong \mathbb Z/16 \oplus \mathbb Z/8 \oplus (\mathbb Z/2)^{4}.\end{aligned}\)}
	\StemEntry{\(\displaystyle \pi_{41}(S^{18})=\langle \bar{\rho}_{18}, \nu_{18}\bar{\kappa}_{21}, \phi_{18}, \psi_{18}, \bar{\varepsilon }^\ast_{18}\rangle\cong \mathbb Z/16 \oplus \mathbb Z/8 \oplus (\mathbb Z/2)^{3}.\)}{\(\displaystyle \begin{aligned}[t]\pi_{41}(S^{18})&=\langle \bar{\rho}_{18}, \nu_{18}\bar{\kappa}_{21}, \phi_{18}, \psi_{18}, \bar{\varepsilon }^\ast_{18}\rangle\\&\cong \mathbb Z/16 \oplus \mathbb Z/8 \oplus (\mathbb Z/2)^{3}.\end{aligned}\)}
	\StemEntry{\(\displaystyle \pi_{n+23}(S^{n})=\langle \bar{\rho}_{n}, \nu_{n}\bar{\kappa}_{n+3}, \phi_{n}, \psi_n\rangle\cong \mathbb Z/16 \oplus \mathbb Z/8 \oplus (\mathbb Z/2)^{2}\quad(n=19, 20, 21).\)}{\(\displaystyle \begin{aligned}[t]\pi_{n+23}(S^{n})&=\langle \bar{\rho}_{n}, \nu_{n}\bar{\kappa}_{n+3}, \phi_{n}, \psi_n\rangle\\&\cong \mathbb Z/16 \oplus \mathbb Z/8 \oplus (\mathbb Z/2)^{2}\quad(n=19, 20, 21).\end{aligned}\)}
	\StemEntry{\(\displaystyle \pi_{n+23}(S^{n})=\langle \bar{\rho}_{n}, \nu_{n}\bar{\kappa}_{n+3}, \phi_{n}\rangle\cong \mathbb Z/16 \oplus \mathbb Z/8 \oplus \mathbb Z/2\quad(n=22, 23).\)}{\(\displaystyle \begin{aligned}[t]\pi_{n+23}(S^{n})&=\langle \bar{\rho}_{n}, \nu_{n}\bar{\kappa}_{n+3}, \phi_{n}\rangle\\&\cong \mathbb Z/16 \oplus \mathbb Z/8 \oplus \mathbb Z/2\quad(n=22, 23).\end{aligned}\)}
	\StemEntry{\(\displaystyle \pi_{47}(S^{24})=\langle \bar{\rho}_{24}, \nu_{24}\bar{\kappa}_{27}, \phi_{24}, P(\iota_{49})\rangle\cong \mathbb Z/16 \oplus \mathbb Z/8 \oplus \mathbb Z/2 \oplus \mathbb Z_{(2)}.\)}{\(\displaystyle \begin{aligned}[t]\pi_{47}(S^{24})&=\langle \bar{\rho}_{24}, \nu_{24}\bar{\kappa}_{27}, \phi_{24}, P(\iota_{49})\rangle\\&\cong \mathbb Z/16 \oplus \mathbb Z/8 \oplus \mathbb Z/2 \oplus \mathbb Z_{(2)}.\end{aligned}\)}
	\StemEntry{\(\displaystyle \pi_{n+23}(S^{n})=\langle \bar{\rho}_n, \nu_n\bar{\kappa}_{n+3}, \phi_n\rangle\cong \mathbb Z/16 \oplus \mathbb Z/8 \oplus \mathbb Z/2\quad(n\ge25).\)}{\(\displaystyle \begin{aligned}[t]\pi_{n+23}(S^{n})&=\langle \bar{\rho}_n, \nu_n\bar{\kappa}_{n+3}, \phi_n\rangle\\&\cong \mathbb Z/16 \oplus \mathbb Z/8 \oplus \mathbb Z/2\quad(n\ge25).\end{aligned}\)}
	\StemEntry{\(\displaystyle \pi^S_{23}(S^0)=\langle \bar{\rho}, \nu\bar{\kappa}, \phi\rangle\cong \mathbb Z/16 \oplus \mathbb Z/8 \oplus \mathbb Z/2.\)}{\(\displaystyle \begin{aligned}[t]\pi^S_{23}(S^0)&=\langle \bar{\rho}, \nu\bar{\kappa}, \phi\rangle\\&\cong \mathbb Z/16 \oplus \mathbb Z/8 \oplus \mathbb Z/2.\end{aligned}\)}
	\subsection{Stem 24}\label{app:stem-24}
	\StemEntry{\(\displaystyle \pi_{26}(S^{2})=\langle \eta_2\bar{\alpha}\rangle\cong \mathbb Z/4.\)}{\(\displaystyle \begin{aligned}[t]\pi_{26}(S^{2})&=\langle \eta_2\bar{\alpha}\rangle\\&\cong \mathbb Z/4.\end{aligned}\)}
	\StemEntry{\(\displaystyle \pi_{27}(S^{3})=\langle \delta_3, \bar{\mu}_3\sigma _{20}, \varepsilon '\kappa_{13}\rangle\cong (\mathbb Z/2)^{3}.\)}{\(\displaystyle \begin{aligned}[t]\pi_{27}(S^{3})&=\langle \delta_3, \bar{\mu}_3\sigma _{20}, \varepsilon '\kappa_{13}\rangle\\&\cong (\mathbb Z/2)^{3}.\end{aligned}\)}
	\StemEntry{\(\displaystyle \pi_{28}(S^{4})=\langle \delta_4, \bar{\mu}_4\sigma _{21}, (\Sigma \varepsilon ')\kappa_{14}, \nu_4\eta_7\bar{\kappa}_8, \nu_4\sigma '\kappa_{14}\rangle\cong (\mathbb Z/2)^{5}.\)}{\(\displaystyle \begin{aligned}[t]\pi_{28}(S^{4})&=\langle \delta_4, \bar{\mu}_4\sigma _{21}, (\Sigma \varepsilon ')\kappa_{14}, \nu_4\eta_7\bar{\kappa}_8, \nu_4\sigma '\kappa_{14}\rangle\\&\cong (\mathbb Z/2)^{5}.\end{aligned}\)}
	\StemEntry{\(\displaystyle \pi_{29}(S^{5})=\langle \delta_5, \bar{\mu}_5\sigma _{22}, \nu_5\eta_8\bar{\kappa}_9, \nu_5\sigma _8\kappa_{15}\rangle\cong (\mathbb Z/2)^{4}.\)}{\(\displaystyle \begin{aligned}[t]\pi_{29}(S^{5})&=\langle \delta_5, \bar{\mu}_5\sigma _{22}, \nu_5\eta_8\bar{\kappa}_9, \nu_5\sigma _8\kappa_{15}\rangle\\&\cong (\mathbb Z/2)^{4}.\end{aligned}\)}
	\StemEntry{\(\displaystyle \pi_{30}(S^{6})=\langle \delta_6, \bar{\mu}_6\sigma _{23}, \nu_6\sigma _9\kappa_{16}, \bar{\zeta}_6', \bar{\sigma }_6', P(\lambda\eta_{31}), P(\xi_{13}\eta_{31})\rangle\cong (\mathbb Z/2)^{3} \oplus \mathbb Z/8 \oplus (\mathbb Z/2)^{3}.\)}{\(\displaystyle \begin{aligned}[t]\pi_{30}(S^{6})&=\langle \delta_6, \bar{\mu}_6\sigma _{23}, \nu_6\sigma _9\kappa_{16}, \bar{\zeta}_6', \bar{\sigma }_6', P(\lambda\eta_{31}), P(\xi_{13}\eta_{31})\rangle\\&\cong (\mathbb Z/2)^{3} \oplus \mathbb Z/8 \oplus (\mathbb Z/2)^{3}.\end{aligned}\)}
	\StemEntry{\(\displaystyle \pi_{31}(S^{7})=\langle \delta_7, \bar{\mu}_7\sigma _{24}, \nu_7\sigma _{10}\kappa_{17}, \bar{\zeta}'_7, \bar{\sigma }_7', \sigma '\bar{\mu}_{14}, \sigma '\omega_{14}\eta_{30}\rangle\cong (\mathbb Z/2)^{7}.\)}{\(\displaystyle \begin{aligned}[t]\pi_{31}(S^{7})&=\langle \delta_7, \bar{\mu}_7\sigma _{24}, \nu_7\sigma _{10}\kappa_{17}, \bar{\zeta}'_7, \bar{\sigma }_7', \sigma '\bar{\mu}_{14}, \sigma '\omega_{14}\eta_{30}\rangle\\&\cong (\mathbb Z/2)^{7}.\end{aligned}\)}
	\StemEntry{\(\displaystyle \pi_{32}(S^{8})=\langle \delta_8, \bar{\mu}_8\sigma _{25}, \nu_8\sigma _{11}\kappa_{18}, \bar{\zeta}'_8, \bar{\sigma }_8', (\Sigma \sigma ')\bar{\mu}_{15}, (\Sigma \sigma ')\omega_{15}\eta_{31}, \sigma _8\bar{\mu}_{15}, \sigma ^2_8\eta_{22}\mu_{23}, \sigma _8\nu_{15}\kappa_{18}, \sigma _8\varepsilon ^\ast_{15}, \sigma _8\eta^{\ast\prime}\eta_{31}\rangle\cong (\mathbb Z/2)^{12}.\)}{\(\displaystyle \begin{aligned}[t]\pi_{32}(S^{8})&=\langle \delta_8, \bar{\mu}_8\sigma _{25}, \nu_8\sigma _{11}\kappa_{18}, \bar{\zeta}'_8, \bar{\sigma }_8', (\Sigma \sigma ')\bar{\mu}_{15}, (\Sigma \sigma ')\omega_{15}\eta_{31}, \sigma _8\bar{\mu}_{15}, \sigma ^2_8\eta_{22}\mu_{23}, \sigma _8\nu_{15}\kappa_{18},\\&\phantom{=\langle{}}\sigma _8\varepsilon ^\ast_{15}, \sigma _8\eta^{\ast\prime}\eta_{31}\rangle\cong (\mathbb Z/2)^{12}.\end{aligned}\)}
	\StemEntry{\(\displaystyle \pi_{33}(S^{9})=\langle \delta_9, \bar{\mu}_9\sigma _{26}, \bar{\sigma }_9', \sigma _9\bar{\mu}_{16}, \sigma ^2_9\eta_{23} \mu_{24}, \sigma _9\nu_{16}\kappa_{19}, \sigma _9\varepsilon ^\ast_{16}\rangle\cong (\mathbb Z/2)^{7}.\)}{\(\displaystyle \begin{aligned}[t]\pi_{33}(S^{9})&=\langle \delta_9, \bar{\mu}_9\sigma _{26}, \bar{\sigma }_9', \sigma _9\bar{\mu}_{16}, \sigma ^2_9\eta_{23} \mu_{24}, \sigma _9\nu_{16}\kappa_{19},\\&\phantom{=\langle{}}\sigma _9\varepsilon ^\ast_{16}\rangle\cong (\mathbb Z/2)^{7}.\end{aligned}\)}
	\StemEntry{\(\displaystyle \pi_{34}(S^{10})=\langle \delta_{10}, \bar{\mu}_{10}\sigma _{27}, \bar{\sigma}'_{10}, \tilde{\varepsilon }_{10}, P(\rho_{21})\rangle\cong (\mathbb Z/2)^{3} \oplus \mathbb Z/4 \oplus \mathbb Z/32.\)}{\(\displaystyle \begin{aligned}[t]\pi_{34}(S^{10})&=\langle \delta_{10}, \bar{\mu}_{10}\sigma _{27}, \bar{\sigma}'_{10}, \tilde{\varepsilon }_{10}, P(\rho_{21})\rangle\\&\cong (\mathbb Z/2)^{3} \oplus \mathbb Z/4 \oplus \mathbb Z/32.\end{aligned}\)}
	\StemEntry{\(\displaystyle \pi_{35}(S^{11})=\langle \delta_{11}, \bar{\mu}_{11}\sigma _{28}, \bar{\sigma }_{11}', \tilde{\varepsilon }_{11}\rangle\cong (\mathbb Z/2)^{3} \oplus \mathbb Z/4.\)}{\(\displaystyle \begin{aligned}[t]\pi_{35}(S^{11})&=\langle \delta_{11}, \bar{\mu}_{11}\sigma _{28}, \bar{\sigma }_{11}', \tilde{\varepsilon }_{11}\rangle\\&\cong (\mathbb Z/2)^{3} \oplus \mathbb Z/4.\end{aligned}\)}
	\StemEntry{\(\displaystyle \pi_{n+24}(S^{n})=\langle \delta_{n}, \bar{\mu}_{n}\sigma _{n+17}, \bar{\sigma }_{n}', \tilde{\varepsilon }_{n}\rangle\cong (\mathbb Z/2)^{4}\quad(n=12, 13).\)}{\(\displaystyle \begin{aligned}[t]\pi_{n+24}(S^{n})&=\langle \delta_{n}, \bar{\mu}_{n}\sigma _{n+17}, \bar{\sigma }_{n}', \tilde{\varepsilon }_{n}\rangle\\&\cong (\mathbb Z/2)^{4}\quad(n=12, 13).\end{aligned}\)}
	\StemEntry{\(\displaystyle \pi_{38}(S^{14})=\langle \delta_{14}, \bar{\mu}_{14}\sigma _{31}, \bar{\sigma }_{14}', \tilde{\varepsilon }_{14}, \zeta^\ast\rangle\cong (\mathbb Z/2)^{4} \oplus \mathbb Z/8.\)}{\(\displaystyle \begin{aligned}[t]\pi_{38}(S^{14})&=\langle \delta_{14}, \bar{\mu}_{14}\sigma _{31}, \bar{\sigma }_{14}', \tilde{\varepsilon }_{14}, \zeta^\ast\rangle\\&\cong (\mathbb Z/2)^{4} \oplus \mathbb Z/8.\end{aligned}\)}
	\StemEntry{\(\displaystyle \pi_{39}(S^{15})=\langle \delta_{15}, \bar{\mu}_{15}\sigma _{32}, \bar{\sigma }_{15}', \tilde{\varepsilon }_{15}, \Sigma \zeta^\ast, \mu^{\ast\prime}\rangle\cong (\mathbb Z/2)^{6}.\)}{\(\displaystyle \begin{aligned}[t]\pi_{39}(S^{15})&=\langle \delta_{15}, \bar{\mu}_{15}\sigma _{32}, \bar{\sigma }_{15}', \tilde{\varepsilon }_{15}, \Sigma \zeta^\ast, \mu^{\ast\prime}\rangle\\&\cong (\mathbb Z/2)^{6}.\end{aligned}\)}
	\StemEntry{\(\displaystyle \pi_{40}(S^{16})=\langle \delta_{16}, \bar{\mu}_{16}\sigma _{33}, \bar{\sigma }_{16}', \tilde{\varepsilon }_{16}, \Sigma ^2\zeta^\ast, \Sigma \mu^{\ast\prime}, \mu^\ast_{16}, \eta_{16}^\ast \varepsilon _{32}, \eta_{16}^\ast\bar{\nu}_{32}\rangle\cong (\mathbb Z/2)^{9}.\)}{\(\displaystyle \begin{aligned}[t]\pi_{40}(S^{16})&=\langle \delta_{16}, \bar{\mu}_{16}\sigma _{33}, \bar{\sigma }_{16}', \tilde{\varepsilon }_{16}, \Sigma ^2\zeta^\ast, \Sigma \mu^{\ast\prime}, \mu^\ast_{16}, \eta_{16}^\ast \varepsilon _{32}, \eta_{16}^\ast\bar{\nu}_{32}\rangle\\&\cong (\mathbb Z/2)^{9}.\end{aligned}\)}
	\StemEntry{\(\displaystyle \pi_{41}(S^{17})=\langle \delta_{17}, \bar{\mu}_{17}\sigma _{34}, \bar{\sigma }_{17}', \mu^\ast_{17}, \eta_{17}^\ast \varepsilon _{33}, \eta_{17}^\ast\bar{\nu}_{33}\rangle\cong (\mathbb Z/2)^{6}.\)}{\(\displaystyle \begin{aligned}[t]\pi_{41}(S^{17})&=\langle \delta_{17}, \bar{\mu}_{17}\sigma _{34}, \bar{\sigma }_{17}', \mu^\ast_{17}, \eta_{17}^\ast \varepsilon _{33}, \eta_{17}^\ast\bar{\nu}_{33}\rangle\\&\cong (\mathbb Z/2)^{6}.\end{aligned}\)}
	\StemEntry{\(\displaystyle \pi_{42}(S^{18})=\langle \delta_{18}, \bar{\mu}_{18}\sigma _{35}, \bar{\sigma }_{18}', P(\sigma _{37})\rangle\cong (\mathbb Z/2)^{3} \oplus \mathbb Z/16.\)}{\(\displaystyle \begin{aligned}[t]\pi_{42}(S^{18})&=\langle \delta_{18}, \bar{\mu}_{18}\sigma _{35}, \bar{\sigma }_{18}', P(\sigma _{37})\rangle\\&\cong (\mathbb Z/2)^{3} \oplus \mathbb Z/16.\end{aligned}\)}
	\StemEntry{\(\displaystyle \pi_{43}(S^{19})=\langle \delta_{19}, \bar{\mu}_{19}\sigma _{36}, \bar{\sigma }_{19}'\rangle\cong (\mathbb Z/2)^{3}.\)}{\(\displaystyle \begin{aligned}[t]\pi_{43}(S^{19})&=\langle \delta_{19}, \bar{\mu}_{19}\sigma _{36}, \bar{\sigma }_{19}'\rangle\\&\cong (\mathbb Z/2)^{3}.\end{aligned}\)}
	\StemEntry{\(\displaystyle \pi_{n+24}(S^{n})=\langle \delta_{n}, \bar{\mu}_{n}\sigma _{n+17}\rangle\cong (\mathbb Z/2)^{2}\quad(n=20, 21).\)}{\(\displaystyle \begin{aligned}[t]\pi_{n+24}(S^{n})&=\langle \delta_{n}, \bar{\mu}_{n}\sigma _{n+17}\rangle\\&\cong (\mathbb Z/2)^{2}\quad(n=20, 21).\end{aligned}\)}
	\StemEntry{\(\displaystyle \pi_{46}(S^{22})=\langle \delta_{22}, \bar{\mu}_{22}\sigma _{39}, P(\nu_{45})\rangle\cong (\mathbb Z/2)^{2} \oplus \mathbb Z/4.\)}{\(\displaystyle \begin{aligned}[t]\pi_{46}(S^{22})&=\langle \delta_{22}, \bar{\mu}_{22}\sigma _{39}, P(\nu_{45})\rangle\\&\cong (\mathbb Z/2)^{2} \oplus \mathbb Z/4.\end{aligned}\)}
	\StemEntry{\(\displaystyle \pi_{47}(S^{23})=\langle \delta_{23}, \bar{\mu}_{23}\sigma _{40}, \tilde{\eta}'\rangle\cong (\mathbb Z/2)^{3}.\)}{\(\displaystyle \begin{aligned}[t]\pi_{47}(S^{23})&=\langle \delta_{23}, \bar{\mu}_{23}\sigma _{40}, \tilde{\eta}'\rangle\\&\cong (\mathbb Z/2)^{3}.\end{aligned}\)}
	\StemEntry{\(\displaystyle \pi_{48}(S^{24})=\langle \delta_{24}, \bar{\mu}_{24}\sigma _{41}, \Sigma \tilde{\eta}', \tilde{\eta}\rangle\cong (\mathbb Z/2)^{4}.\)}{\(\displaystyle \begin{aligned}[t]\pi_{48}(S^{24})&=\langle \delta_{24}, \bar{\mu}_{24}\sigma _{41}, \Sigma \tilde{\eta}', \tilde{\eta}\rangle\\&\cong (\mathbb Z/2)^{4}.\end{aligned}\)}
	\StemEntry{\(\displaystyle \pi_{49}(S^{25})=\langle \delta_{25}, \bar{\mu}_{25}\sigma _{42}, \Sigma \tilde{\eta}\rangle\cong (\mathbb Z/2)^{3}.\)}{\(\displaystyle \begin{aligned}[t]\pi_{49}(S^{25})&=\langle \delta_{25}, \bar{\mu}_{25}\sigma _{42}, \Sigma \tilde{\eta}\rangle\\&\cong (\mathbb Z/2)^{3}.\end{aligned}\)}
	\StemEntry{\(\displaystyle \pi_{n+24}(S^{n})=\langle \delta_{n}, \bar{\mu}_{n}\sigma _{n+17}\rangle\cong (\mathbb Z/2)^{2}\quad(n \ge 26).\)}{\(\displaystyle \begin{aligned}[t]\pi_{n+24}(S^{n})&=\langle \delta_{n}, \bar{\mu}_{n}\sigma _{n+17}\rangle\\&\cong (\mathbb Z/2)^{2}\quad(n \ge 26).\end{aligned}\)}
	\StemEntry{\(\displaystyle \pi^S_{24}(S^0)=\langle \delta, \bar{\mu}\sigma\rangle\cong \mathbb Z/2 \oplus \mathbb Z/2.\)}{\(\displaystyle \begin{aligned}[t]\pi^S_{24}(S^0)&=\langle \delta, \bar{\mu}\sigma\rangle\\&\cong \mathbb Z/2 \oplus \mathbb Z/2.\end{aligned}\)}
	
	\section{The   Oda Stems}
	\subsection*{Oda's notation for generators}
	Oda slightly modified the notation for certain generators. For ease of comparison,
	we record all the relations in \cite[pp.~52--53, (2.1)--(2.5)]{Oda1979},
	with $\Sigma$ denoting suspension. In (2.1)--(2.4) below, the symbols on the left
	are Oda's, while those on the right come from the earlier literature.
	The terms following $\bmod$ specify the subgroup modulo which each congruence
	is taken. The equalities in (2.5) give exact identifications.
	
	\begingroup
	\setlength{\LTleft}{\fill}
	\setlength{\LTright}{\fill}
	\setlength{\LTpre}{4pt}
	\setlength{\LTpost}{6pt}
	\renewcommand{\arraystretch}{1.15}
	
	\begin{longtable}{@{}r@{\hspace{0.4em}}l@{}}
		
		\multicolumn{2}{@{}c@{}}{\textit{Oda's relations (2.1)}}\\*[3pt]
		
		\(\displaystyle A_1\)
		& \(\displaystyle
		\equiv\theta\m \Sigma\theta';
		\qquad C_1=\theta';
		\)\\
		
		\(\displaystyle A_2\)
		& \(\displaystyle
		\equiv\bar{\bar{\beta}}\m \Sigma\bar\beta;
		\qquad C_2=\bar\beta;
		\)\\
		
		\(\displaystyle B_1\)
		& \(\displaystyle
		\equiv\eta_{16}^{*}+\omega_{16}
		\m \Sigma\eta^{*\prime},\ \sigma_{16}\mu_{23};
		\)\\
		
		\(\displaystyle B_2\)
		& \(\displaystyle
		\equiv\tilde\eta\m \Sigma\tilde\eta';
		\qquad D_2=\tilde\eta';
		\)\\
		
		\(\displaystyle D_1\)
		& \(\displaystyle
		\equiv\eta^{*\prime}\m \sigma_{15}\mu_{22};
		\)\\
		
		\(\displaystyle G_0\)
		& \(\displaystyle
		\equiv\zeta'\m  2\zeta',\ \eta_6\bar\varepsilon_7;
		\)\\
		
		\(\displaystyle G_1\)
		& \(\displaystyle
		\equiv\zeta^{*}\m  2\zeta^{*},\
		\omega_{14}\varepsilon_{30};
		\)\\
		
		\(\displaystyle F_1\)
		& \(\displaystyle
		\equiv\beta'\m  2\beta'.
		\)\\[8pt]
		
		\multicolumn{2}{@{}c@{}}{\textit{Oda's relations (2.2)}}\\*[3pt]
		
		\(\displaystyle A_1^{\mathrm{I}}\)
		& \(\displaystyle
		\equiv\omega'\m \xi_{12}\eta_{30},\
		(\Sigma\lambda')\eta_{30},\
		(\Sigma\xi')\eta_{30};
		\)\\
		
		\(\displaystyle C_1^{\mathrm{I}}\)
		& \(\displaystyle
		\equiv\lambda'\eta_{29}\m \xi'\eta_{29};
		\)\\
		
		\(\displaystyle B_1^{\mathrm{I}}\)
		& \(\displaystyle
		\equiv\bar\varepsilon_{16}^{*}
		\m \bar\varepsilon_{16}^{*}+\bar\nu_{16}^{*},\
		\Sigma\bar\varepsilon^{*\prime},\
		\Sigma\bar\nu^{*\prime},\
		4\nu_{16}\bar\kappa_{19};
		\)\\
		
		\(\displaystyle D_1^{\mathrm{I}}\)
		& \(\displaystyle
		\equiv\bar\varepsilon^{*\prime}
		\m \bar\varepsilon^{*\prime}+\bar\nu^{*\prime},\
		2\nu_{15}\bar\kappa_{18}.
		\)\\[8pt]
		
		\multicolumn{2}{@{}c@{}}{\textit{Oda's relations (2.3)}}\\*[3pt]
		
		\(\displaystyle A_1^{(1)}\)
		& \(\displaystyle
		\equiv\beta'''\m \theta\bar\nu_{24},\
		\theta\varepsilon_{24},\
		\Sigma\beta'',\
		(\Sigma\theta')\varepsilon_{24};
		\)\\
		
		\(\displaystyle C_1^{(1)}\)
		& \(\displaystyle
		\equiv\beta''\m \theta'\varepsilon_{23};
		\)\\
		
		\(\displaystyle B_1^{(1)}\)
		& \(\displaystyle
		\equiv\mu_{16}^{*}
		\m \eta_{16}^{*}\varepsilon_{32},\
		\eta_{16}^{*}\bar\nu_{32},\
		\Sigma\pi_{39}(S^{15});
		\)\\
		
		\(\displaystyle D_1^{(1)}\)
		& \(\displaystyle
		\equiv\mu^{*\prime}\m \Sigma\pi_{38}(S^{14});
		\)\\
		
		\(\displaystyle G_0^{(1)}\)
		& \(\displaystyle
		\equiv\bar\zeta'_6\m  2\bar\zeta'_6.
		\)\\[8pt]
		
		\multicolumn{2}{@{}c@{}}{\textit{Oda's relations (2.4)}}\\*[3pt]
		
		\(\displaystyle A_1\sigma_{24}\)
		& \(\displaystyle
		\equiv\xi_{12}\eta_{30}
		\m (\Sigma\xi')\eta_{30},\
		(\Sigma\lambda')\eta_{30};
		\)\\
		
		\(\displaystyle C_1\sigma_{23}\)
		& \(\displaystyle
		=\xi'\eta_{29};
		\)\\
		
		\(\displaystyle B_1\sigma_{32}\)
		& \(\displaystyle
		\equiv\bar\varepsilon_{16}^{*}+\bar\nu_{16}^{*}
		\m \Sigma\pi_{38}(S^{15});
		\)\\
		
		\(\displaystyle D_1\sigma_{31}\)
		& \(\displaystyle
		=\bar\varepsilon^{*\prime}+\bar\nu^{*\prime}.
		\)\\[8pt]
		
		\multicolumn{2}{@{}c@{}}{\textit{Oda's identifications (2.5)}}\\*[3pt]
		
		\(\displaystyle \alpha_0^{\mathrm{IV}}\)
		& \(\displaystyle
		=\sigma_9,\qquad
		\alpha_0'''=\sigma',\qquad
		\alpha_0''=\sigma'',\qquad
		\alpha_0'=\sigma''';
		\)\\*
		
		\(\displaystyle \alpha_1^{\mathrm{IV}}\)
		& \(\displaystyle
		=\rho',\qquad
		\alpha_1'''=\rho'',\qquad
		\alpha_1''=\rho''',\qquad
		\alpha_1'=\rho^{\mathrm{IV}};
		\)\\*
		
		\(\displaystyle \alpha_2^{\mathrm{IV}}\)
		& \(\displaystyle
		=\bar\rho_9,\qquad
		\alpha_2'''=\bar\rho',\qquad
		\alpha_2''=\bar\rho'',\qquad
		\alpha_2'=\bar\rho'''.
		\)
		
	\end{longtable}
	\endgroup

	\subsection{Stem 25}\label{app:stem-25}
	For the ordered bases in stems $25$ and $28,$   choose $t\in\mathbb Z$ such that
	$2\phi'''=\nu_7^2\bar\sigma_{13}+2t\Sigma^2\phi'',$  
	and write $u_7=\phi'''-t\Sigma^2\phi''.$ 
	\StemEntry{\(\displaystyle \pi_{27}(S^{2})=\langle \eta_2 \varepsilon ' \kappa_{13}, \eta_2 \delta_3, \eta_2 \bar\mu_3 \sigma_{20}\rangle\cong (\mathbb Z/2)^{3}.\)}{\(\displaystyle \begin{aligned}[t]\pi_{27}(S^{2})&=\langle \eta_2 \varepsilon ' \kappa_{13}, \eta_2 \delta_3, \eta_2 \bar\mu_3 \sigma_{20}\rangle\\&\cong (\mathbb Z/2)^{3}.\end{aligned}\)}
	\StemEntry{\(\displaystyle \pi_{28}(S^{3})=\langle \phi', \nu' \varepsilon _6 \kappa_{14}, \mu_{3,3}, \eta_3 \bar\mu_4 \sigma_{21}\rangle\cong \mathbb Z/4 \oplus (\mathbb Z/2)^{3}.\)}{\(\displaystyle \begin{aligned}[t]\pi_{28}(S^{3})&=\langle \phi', \nu' \varepsilon _6 \kappa_{14}, \mu_{3,3}, \eta_3 \bar\mu_4 \sigma_{21}\rangle\\&\cong \mathbb Z/4 \oplus (\mathbb Z/2)^{3}.\end{aligned}\)}
	\StemEntry{\(\displaystyle \pi_{29}(S^{4})=\langle \nu_4 \sigma' \rho_{14}, \nu_4 \bar\nu_7 \kappa_{15}, \nu_4 \varepsilon _7 \kappa_{15}, \nu^2_4 \bar\sigma_{10}, \Sigma \phi', \Sigma \nu' \varepsilon _7 \kappa_{15}, \mu_{3,4}, \eta_4 \bar\mu_5 \sigma_{22}\rangle\cong \mathbb Z/8 \oplus \mathbb Z/2 \oplus \mathbb Z/2 \oplus \mathbb Z/2 \oplus \mathbb Z/4 \oplus (\mathbb Z/2)^{3}.\)}{\(\displaystyle \begin{aligned}[t]\pi_{29}(S^{4})&=\langle \nu_4 \sigma' \rho_{14}, \nu_4 \bar\nu_7 \kappa_{15}, \nu_4 \varepsilon _7 \kappa_{15}, \nu^2_4 \bar\sigma_{10}, \Sigma \phi', \Sigma \nu' \varepsilon _7 \kappa_{15}, \mu_{3,4},\\&\phantom{=\langle{}}\eta_4 \bar\mu_5 \sigma_{22}\rangle\cong \mathbb Z/8 \oplus \mathbb Z/2 \oplus \mathbb Z/2 \oplus \mathbb Z/2 \oplus \mathbb Z/4 \oplus (\mathbb Z/2)^{3}.\end{aligned}\)}
	\StemEntry{\(\displaystyle \pi_{30}(S^{5})=\langle \phi'', \nu_5 \varepsilon _8 \kappa_{16}, \nu^2_5 \bar\sigma_{11}, \mu_{3,5}, \eta_5 \bar\mu_6 \sigma_{23}\rangle\cong \mathbb Z/8 \oplus (\mathbb Z/2)^{4}.\)}{\(\displaystyle \begin{aligned}[t]\pi_{30}(S^{5})&=\langle \phi'', \nu_5 \varepsilon _8 \kappa_{16}, \nu^2_5 \bar\sigma_{11}, \mu_{3,5}, \eta_5 \bar\mu_6 \sigma_{23}\rangle\\&\cong \mathbb Z/8 \oplus (\mathbb Z/2)^{4}.\end{aligned}\)}
	\StemEntry{\(\displaystyle \pi_{31}(S^{6})=\langle P(\bar\kappa_{13}), P((\Sigma A_1)\varepsilon_{25}), P(\Sigma A_1^{(1)}), \Sigma \phi'', \nu^2_6\bar\sigma_{12}, \mu_{3,6}, \eta_6 \bar\mu_7 \sigma_{24}\rangle\cong \mathbb Z/4 \oplus \mathbb Z/2 \oplus \mathbb Z/2 \oplus \mathbb Z/8 \oplus (\mathbb Z/2)^{3}.\)}{\(\displaystyle \begin{aligned}[t]\pi_{31}(S^{6})&=\langle P(\bar\kappa_{13}), P((\Sigma A_1)\varepsilon_{25}), P(\Sigma A_1^{(1)}), \Sigma \phi'', \nu^2_6\bar\sigma_{12}, \mu_{3,6},\\&\phantom{=\langle{}}\eta_6 \bar\mu_7 \sigma_{24}\rangle\cong \mathbb Z/4 \oplus \mathbb Z/2 \oplus \mathbb Z/2 \oplus \mathbb Z/8 \oplus (\mathbb Z/2)^{3}.\end{aligned}\)}
	\StemEntry{\(\displaystyle \pi_{32}(S^{7})=\langle u_7, \Sigma ^2\phi'', \sigma' \eta_{14} \bar\mu_{15}, \mu_{3,7}, \eta_7 \bar\mu_8 \sigma_{25}\rangle\cong \mathbb Z/4 \oplus \mathbb Z/8 \oplus (\mathbb Z/2)^{3}.\)}{\(\displaystyle \begin{aligned}[t]\pi_{32}(S^{7})&=\langle u_7, \Sigma ^2\phi'', \sigma' \eta_{14} \bar\mu_{15}, \mu_{3,7}, \eta_7 \bar\mu_8 \sigma_{25}\rangle\\&\cong \mathbb Z/4 \oplus \mathbb Z/8 \oplus (\mathbb Z/2)^{3}.\end{aligned}\)}
	\StemEntry{\(\displaystyle \pi_{33}(S^{8})=\langle \sigma_8 \xi_{15}, \sigma_8 \Sigma ^2\lambda, \sigma_8 \eta_{15} \bar\mu_{16}, \Sigma u_7, \Sigma ^3\phi'', (\Sigma \sigma') \eta_{15} \bar\mu_{16}, \mu_{3,8}, \eta_8 \bar\mu_9 \sigma_{26}\rangle\cong \mathbb Z/8 \oplus \mathbb Z/8 \oplus \mathbb Z/2 \oplus \mathbb Z/4 \oplus \mathbb Z/8 \oplus (\mathbb Z/2)^{3}.\)}{\(\displaystyle \begin{aligned}[t]\pi_{33}(S^{8})&=\langle \sigma_8 \xi_{15}, \sigma_8 \Sigma ^2\lambda, \sigma_8 \eta_{15} \bar\mu_{16}, \Sigma u_7, \Sigma ^3\phi'', (\Sigma \sigma') \eta_{15} \bar\mu_{16}, \mu_{3,8},\\&\phantom{=\langle{}}\eta_8 \bar\mu_9 \sigma_{26}\rangle\cong \mathbb Z/8 \oplus \mathbb Z/8 \oplus \mathbb Z/2 \oplus \mathbb Z/4 \oplus \mathbb Z/8 \oplus (\mathbb Z/2)^{3}.\end{aligned}\)}
	\StemEntry{\(\displaystyle \pi_{34}(S^{9})=\langle \sigma_9 \xi_{16}, \sigma_9 \eta_{16} \bar\mu_{17}, \sigma_9 \nu^*_{16}, \mu_{3,9}, \eta_9 \bar\mu_{10} \sigma_{27}\rangle\cong \mathbb Z/8 \oplus \mathbb Z/2 \oplus \mathbb Z/4 \oplus (\mathbb Z/2)^{2}.\)}{\(\displaystyle \begin{aligned}[t]\pi_{34}(S^{9})&=\langle \sigma_9 \xi_{16}, \sigma_9 \eta_{16} \bar\mu_{17}, \sigma_9 \nu^*_{16}, \mu_{3,9},\\&\phantom{=\langle{}}\eta_9 \bar\mu_{10} \sigma_{27}\rangle\cong \mathbb Z/8 \oplus \mathbb Z/2 \oplus \mathbb Z/4 \oplus (\mathbb Z/2)^{2}.\end{aligned}\)}
	\StemEntry{\(\displaystyle \pi_{35}(S^{10})=\langle \sigma_{10} \xi_{17}, \sigma_{10} \nu^*_{17}, \mu_{3,10}, \eta_{10} \bar\mu_{11} \sigma_{28}\rangle\cong \mathbb Z/4 \oplus \mathbb Z/4 \oplus (\mathbb Z/2)^{2}.\)}{\(\displaystyle \begin{aligned}[t]\pi_{35}(S^{10})&=\langle \sigma_{10} \xi_{17}, \sigma_{10} \nu^*_{17}, \mu_{3,10}, \eta_{10} \bar\mu_{11} \sigma_{28}\rangle\\&\cong \mathbb Z/4 \oplus \mathbb Z/4 \oplus (\mathbb Z/2)^{2}.\end{aligned}\)}
	\StemEntry{\(\displaystyle \pi_{36}(S^{11})=\langle K_1, \sigma_{11} \xi_{18}, \sigma_{11} \xi_{18} + \sigma_{11} \nu^*_{18}, \mu_{3,11}, \eta_{11} \bar\mu_{12} \sigma_{29}\rangle\cong \mathbb Z/2 \oplus \mathbb Z/4 \oplus (\mathbb Z/2)^{3}.\)}{\(\displaystyle \begin{aligned}[t]\pi_{36}(S^{11})&=\langle K_1, \sigma_{11} \xi_{18}, \sigma_{11} \xi_{18} + \sigma_{11} \nu^*_{18}, \mu_{3,11},\\&\phantom{=\langle{}}\eta_{11} \bar\mu_{12} \sigma_{29}\rangle\cong \mathbb Z/2 \oplus \mathbb Z/4 \oplus (\mathbb Z/2)^{3}.\end{aligned}\)}
	\StemEntry{\(\displaystyle \pi_{37}(S^{12})=\langle \xi_{12} \sigma_{30}, \Sigma K_1, \sigma_{12} \xi_{19} + \sigma_{12} \nu^*_{19}, \mu_{3,12}, \eta_{12} \bar\mu_{13} \sigma_{30}\rangle\cong \mathbb Z/8 \oplus (\mathbb Z/2)^{4}.\)}{\(\displaystyle \begin{aligned}[t]\pi_{37}(S^{12})&=\langle \xi_{12} \sigma_{30}, \Sigma K_1, \sigma_{12} \xi_{19} + \sigma_{12} \nu^*_{19}, \mu_{3,12},\\&\phantom{=\langle{}}\eta_{12} \bar\mu_{13} \sigma_{30}\rangle\cong \mathbb Z/8 \oplus (\mathbb Z/2)^{4}.\end{aligned}\)}
	\StemEntry{\(\displaystyle \pi_{n+25}(S^{n})=\langle \xi_n \sigma_{n+18}, \mu_{3,n}, \eta_n \bar\mu_{n+1} \sigma_{n+18}\rangle\cong \mathbb Z/8 \oplus (\mathbb Z/2)^{2}\quad(n = 13 , 14).\)}{\(\displaystyle \begin{aligned}[t]\pi_{n+25}(S^{n})&=\langle \xi_n \sigma_{n+18}, \mu_{3,n}, \eta_n \bar\mu_{n+1} \sigma_{n+18}\rangle\\&\cong \mathbb Z/8 \oplus (\mathbb Z/2)^{2}\quad(n = 13 , 14).\end{aligned}\)}
	\StemEntry{\(\displaystyle \pi_{40}(S^{15})=\langle D_1 \mu_{31}, \xi_{15} \sigma_{33}, \mu_{3,15}, \eta_{15} \bar\mu_{16} \sigma_{33}\rangle\cong \mathbb Z/2 \oplus \mathbb Z/8 \oplus (\mathbb Z/2)^{2}.\)}{\(\displaystyle \begin{aligned}[t]\pi_{40}(S^{15})&=\langle D_1 \mu_{31}, \xi_{15} \sigma_{33}, \mu_{3,15}, \eta_{15} \bar\mu_{16} \sigma_{33}\rangle\\&\cong \mathbb Z/2 \oplus \mathbb Z/8 \oplus (\mathbb Z/2)^{2}.\end{aligned}\)}
	\StemEntry{\(\displaystyle \pi_{41}(S^{16})=\langle B_1 \mu_{32}, (\Sigma D_1) \mu_{32}, \xi_{16} \sigma_{34}, \mu_{3,16}, \eta_{16} \bar\mu_{17} \sigma_{34}\rangle\cong (\mathbb Z/2)^{2} \oplus \mathbb Z/8 \oplus (\mathbb Z/2)^{2}.\)}{\(\displaystyle \begin{aligned}[t]\pi_{41}(S^{16})&=\langle B_1 \mu_{32}, (\Sigma D_1) \mu_{32}, \xi_{16} \sigma_{34}, \mu_{3,16},\\&\phantom{=\langle{}}\eta_{16} \bar\mu_{17} \sigma_{34}\rangle\cong (\mathbb Z/2)^{2} \oplus \mathbb Z/8 \oplus (\mathbb Z/2)^{2}.\end{aligned}\)}
	\StemEntry{\(\displaystyle \pi_{42}(S^{17})=\langle (\Sigma B_1) \mu_{33}, \xi_{17} \sigma_{35}, \mu_{3,17}, \eta_{17} \bar\mu_{18} \sigma_{35}\rangle\cong \mathbb Z/2 \oplus \mathbb Z/8 \oplus (\mathbb Z/2)^{2}.\)}{\(\displaystyle \begin{aligned}[t]\pi_{42}(S^{17})&=\langle (\Sigma B_1) \mu_{33}, \xi_{17} \sigma_{35}, \mu_{3,17},\\&\phantom{=\langle{}}\eta_{17} \bar\mu_{18} \sigma_{35}\rangle\cong \mathbb Z/2 \oplus \mathbb Z/8 \oplus (\mathbb Z/2)^{2}.\end{aligned}\)}
	\StemEntry{\(\displaystyle \pi_{43}(S^{18})=\langle \xi_{18} \sigma_{36}, \mu_{3,18}, \eta_{18} \bar\mu_{19} \sigma_{36}\rangle\cong \mathbb Z/4 \oplus (\mathbb Z/2)^{2}.\)}{\(\displaystyle \begin{aligned}[t]\pi_{43}(S^{18})&=\langle \xi_{18} \sigma_{36}, \mu_{3,18}, \eta_{18} \bar\mu_{19} \sigma_{36}\rangle\\&\cong \mathbb Z/4 \oplus (\mathbb Z/2)^{2}.\end{aligned}\)}
	\StemEntry{\(\displaystyle \pi_{44}(S^{19})=\langle \xi_{19} \sigma_{37}, \mu_{3,19}, \eta_{19} \bar\mu_{20} \sigma_{37}\rangle\cong (\mathbb Z/2)^{3}.\)}{\(\displaystyle \begin{aligned}[t]\pi_{44}(S^{19})&=\langle \xi_{19} \sigma_{37}, \mu_{3,19}, \eta_{19} \bar\mu_{20} \sigma_{37}\rangle\\&\cong (\mathbb Z/2)^{3}.\end{aligned}\)}
	\StemEntry{\(\displaystyle \pi_{n+25}(S^{n})=\langle \mu_{3,n}, \eta_n \bar\mu_{n+1} \sigma_{n+18}\rangle\cong (\mathbb Z/2)^{2}\quad(n = 20 , 21, 22).\)}{\(\displaystyle \begin{aligned}[t]\pi_{n+25}(S^{n})&=\langle \mu_{3,n}, \eta_n \bar\mu_{n+1} \sigma_{n+18}\rangle\\&\cong (\mathbb Z/2)^{2}\quad(n = 20 , 21, 22).\end{aligned}\)}
	\StemEntry{\(\displaystyle \pi_{48}(S^{23})=\langle D_2 \eta_{47}, \mu_{3,23}, \eta_{23} \bar\mu_{24} \sigma_{41}\rangle\cong (\mathbb Z/2)^{3}.\)}{\(\displaystyle \begin{aligned}[t]\pi_{48}(S^{23})&=\langle D_2 \eta_{47}, \mu_{3,23}, \eta_{23} \bar\mu_{24} \sigma_{41}\rangle\\&\cong (\mathbb Z/2)^{3}.\end{aligned}\)}
	\StemEntry{\(\displaystyle \pi_{49}(S^{24})=\langle B_2 \eta_{48}, (\Sigma D_2) \eta_{48}, \mu_{3,24}, \eta_{24} \bar\mu_{25} \sigma_{42}\rangle\cong (\mathbb Z/2)^{4}.\)}{\(\displaystyle \begin{aligned}[t]\pi_{49}(S^{24})&=\langle B_2 \eta_{48}, (\Sigma D_2) \eta_{48}, \mu_{3,24},\\&\phantom{=\langle{}}\eta_{24} \bar\mu_{25} \sigma_{42}\rangle\cong (\mathbb Z/2)^{4}.\end{aligned}\)}
	\StemEntry{\(\displaystyle \pi_{50}(S^{25})=\langle (\Sigma B_2) \eta_{49}, \mu_{3,25}, \eta_{25} \bar\mu_{26} \sigma_{43}\rangle\cong (\mathbb Z/2)^{3}.\)}{\(\displaystyle \begin{aligned}[t]\pi_{50}(S^{25})&=\langle (\Sigma B_2) \eta_{49}, \mu_{3,25}, \eta_{25} \bar\mu_{26} \sigma_{43}\rangle\\&\cong (\mathbb Z/2)^{3}.\end{aligned}\)}
	\StemEntry{\(\displaystyle \pi_{51}(S^{26})=\langle P(\iota_{53}), \mu_{3,26}, \eta_{26} \bar\mu_{27} \sigma_{44}\rangle\cong \mathbb Z_{(2)} \oplus (\mathbb Z/2)^{2}.\)}{\(\displaystyle \begin{aligned}[t]\pi_{51}(S^{26})&=\langle P(\iota_{53}), \mu_{3,26}, \eta_{26} \bar\mu_{27} \sigma_{44}\rangle\\&\cong \mathbb Z_{(2)} \oplus (\mathbb Z/2)^{2}.\end{aligned}\)}
	\StemEntry{\(\displaystyle \pi_{n+25}(S^{n})=\langle \mu_{3,n}, \eta_{n} \bar\mu_{n+1} \sigma_{n+18}\rangle\cong (\mathbb Z/2)^{2}\quad(n\ge 27).\)}{\(\displaystyle \begin{aligned}[t]\pi_{n+25}(S^{n})&=\langle \mu_{3,n}, \eta_{n} \bar\mu_{n+1} \sigma_{n+18}\rangle\\&\cong (\mathbb Z/2)^{2}\quad(n\ge 27).\end{aligned}\)}
	\StemEntry{\(\displaystyle \pi^S_{25}(S^0)=\langle \mu_{3}, \eta \bar\mu \sigma\rangle\cong \mathbb Z/2 \oplus \mathbb Z/2.\)}{\(\displaystyle \begin{aligned}[t]\pi^S_{25}(S^0)&=\langle \mu_{3}, \eta \bar\mu \sigma\rangle\\&\cong \mathbb Z/2 \oplus \mathbb Z/2.\end{aligned}\)}
	\subsection{Stem 26}\label{app:stem-26}
	\StemEntry{\(\displaystyle \pi_{28}(S^{2})=\langle \eta_2 \phi', \eta_2 \nu' \varepsilon _6 \kappa_{14}, \eta^2_2 \bar\mu_4 \sigma_{21}, \eta_2 \mu_{3,3}\rangle\cong \mathbb Z/4 \oplus (\mathbb Z/2)^{3}.\)}{\(\displaystyle \begin{aligned}[t]\pi_{28}(S^{2})&=\langle \eta_2 \phi', \eta_2 \nu' \varepsilon _6 \kappa_{14}, \eta^2_2 \bar\mu_4 \sigma_{21},\\&\phantom{=\langle{}}\eta_2 \mu_{3,3}\rangle\cong \mathbb Z/4 \oplus (\mathbb Z/2)^{3}.\end{aligned}\)}
	\StemEntry{\(\displaystyle \pi_{29}(S^{3})=\langle \bar\mu' \sigma_{22}, \nu' \phi_6, \bar\alpha \nu_{26}, \eta_3 \mu_{3,4}\rangle\cong \mathbb Z/4 \oplus (\mathbb Z/2)^{3}.\)}{\(\displaystyle \begin{aligned}[t]\pi_{29}(S^{3})&=\langle \bar\mu' \sigma_{22}, \nu' \phi_6, \bar\alpha \nu_{26}, \eta_3 \mu_{3,4}\rangle\\&\cong \mathbb Z/4 \oplus (\mathbb Z/2)^{3}.\end{aligned}\)}
	\StemEntry{\(\displaystyle \pi_{30}(S^{4})=\langle \nu^2_4 \bar\kappa_{10}, \nu_4 \alpha'''_2, \nu_4 \sigma' \omega_{14}, \nu_4 \sigma' \sigma_{14} \mu_{21}, \nu_4 (\bar\kappa_7 \nu_{27} - \nu_7 \bar\kappa_{10}), \nu_4 \phi_7, (\Sigma \bar\mu') \sigma_{23}, (\Sigma \nu') \phi_7, (\Sigma \bar\alpha )\nu_{27}, \eta_4 \mu_{3,5}\rangle\cong \mathbb Z/8 \oplus \mathbb Z/8 \oplus \mathbb Z/2 \oplus \mathbb Z/2 \oplus \mathbb Z/2 \oplus \mathbb Z/2 \oplus \mathbb Z/4 \oplus (\mathbb Z/2)^{3}.\)}{\(\displaystyle \begin{aligned}[t]\pi_{30}(S^{4})&=\langle \nu^2_4 \bar\kappa_{10}, \nu_4 \alpha'''_2, \nu_4 \sigma' \omega_{14}, \nu_4 \sigma' \sigma_{14} \mu_{21}, \nu_4 (\bar\kappa_7 \nu_{27} - \nu_7 \bar\kappa_{10}), \nu_4 \phi_7, (\Sigma \bar\mu') \sigma_{23}, (\Sigma \nu') \phi_7,\\&\phantom{=\langle{}}(\Sigma \bar\alpha )\nu_{27}, \eta_4 \mu_{3,5}\rangle\cong \mathbb Z/8 \oplus \mathbb Z/8 \oplus \mathbb Z/2 \oplus \mathbb Z/2 \oplus \mathbb Z/2 \oplus \mathbb Z/2 \oplus \mathbb Z/4 \oplus (\mathbb Z/2)^{3}.\end{aligned}\)}
	\StemEntry{\(\displaystyle \pi_{31}(S^{5})=\langle \bar\zeta_5 \sigma_{24}, \phi_5 \nu_{28}, \nu^2_5 \bar\kappa_{11}, \nu_5 \bar\kappa_8 \nu_{28}, \nu_5 \phi_8, \eta_5 \mu_{3,6}\rangle\cong \mathbb Z/8 \oplus (\mathbb Z/2)^{5}.\)}{\(\displaystyle \begin{aligned}[t]\pi_{31}(S^{5})&=\langle \bar\zeta_5 \sigma_{24}, \phi_5 \nu_{28}, \nu^2_5 \bar\kappa_{11}, \nu_5 \bar\kappa_8 \nu_{28}, \nu_5 \phi_8,\\&\phantom{=\langle{}}\eta_5 \mu_{3,6}\rangle\cong \mathbb Z/8 \oplus (\mathbb Z/2)^{5}.\end{aligned}\)}
	\StemEntry{\(\displaystyle \pi_{32}(S^{6})=\langle P((\Sigma A_1) \mu_{25}), P(\lambda \nu_{31}), \bar\sigma_6 \sigma_{25}, \bar\zeta_6\sigma_{25}, \phi_6 \nu_{29}, \nu^2_6\bar\kappa_{12}, \eta_6 \mu_{3,7}\rangle\cong \mathbb Z/2 \oplus \mathbb Z/4 \oplus \mathbb Z/4 \oplus \mathbb Z/8 \oplus (\mathbb Z/2)^{3}.\)}{\(\displaystyle \begin{aligned}[t]\pi_{32}(S^{6})&=\langle P((\Sigma A_1) \mu_{25}), P(\lambda \nu_{31}), \bar\sigma_6 \sigma_{25}, \bar\zeta_6\sigma_{25}, \phi_6 \nu_{29}, \nu^2_6\bar\kappa_{12},\\&\phantom{=\langle{}}\eta_6 \mu_{3,7}\rangle\cong \mathbb Z/2 \oplus \mathbb Z/4 \oplus \mathbb Z/4 \oplus \mathbb Z/8 \oplus (\mathbb Z/2)^{3}.\end{aligned}\)}
	\StemEntry{\(\displaystyle \pi_{33}(S^{7})=\langle \sigma' \omega_{14} \nu_{30}, \bar\kappa_7 \nu^2_{27}, \bar\sigma_7 \sigma_{26}, \bar\zeta_7 \sigma_{26}, \phi_7 \nu_{30}, \nu^2_7 \bar\kappa_{13}, \eta_7 \mu_{3,8}\rangle\cong (\mathbb Z/2)^{3} \oplus \mathbb Z/8 \oplus (\mathbb Z/2)^{3}.\)}{\(\displaystyle \begin{aligned}[t]\pi_{33}(S^{7})&=\langle \sigma' \omega_{14} \nu_{30}, \bar\kappa_7 \nu^2_{27}, \bar\sigma_7 \sigma_{26}, \bar\zeta_7 \sigma_{26}, \phi_7 \nu_{30}, \nu^2_7 \bar\kappa_{13},\\&\phantom{=\langle{}}\eta_7 \mu_{3,8}\rangle\cong (\mathbb Z/2)^{3} \oplus \mathbb Z/8 \oplus (\mathbb Z/2)^{3}.\end{aligned}\)}
	\StemEntry{\(\displaystyle \pi_{34}(S^{8})=\langle \sigma_8 \omega_{15} \nu_{31}, \sigma_8 \bar\sigma_{15}, \sigma_8 \bar\zeta_{15}, (\Sigma \sigma') \omega_{15} \nu_{31}, \bar\kappa_8 \nu^2_{28}, \bar\sigma_8 \sigma_{27}, \bar\zeta_8 \sigma_{27}, \phi_8 \nu_{31}, \nu^2_8 \bar\kappa_{14}, \eta_8 \mu_{3,9}\rangle\cong (\mathbb Z/2)^{2} \oplus \mathbb Z/8 \oplus \mathbb Z/2 \oplus \mathbb Z/2 \oplus \mathbb Z/2 \oplus \mathbb Z/8 \oplus (\mathbb Z/2)^{3}.\)}{\(\displaystyle \begin{aligned}[t]\pi_{34}(S^{8})&=\langle \sigma_8 \omega_{15} \nu_{31}, \sigma_8 \bar\sigma_{15}, \sigma_8 \bar\zeta_{15}, (\Sigma \sigma') \omega_{15} \nu_{31}, \bar\kappa_8 \nu^2_{28}, \bar\sigma_8 \sigma_{27}, \bar\zeta_8 \sigma_{27}, \phi_8 \nu_{31},\\&\phantom{=\langle{}}\nu^2_8 \bar\kappa_{14}, \eta_8 \mu_{3,9}\rangle\cong (\mathbb Z/2)^{2} \oplus \mathbb Z/8 \oplus \mathbb Z/2 \oplus \mathbb Z/2 \oplus \mathbb Z/2 \oplus \mathbb Z/8 \oplus (\mathbb Z/2)^{3}.\end{aligned}\)}
	\StemEntry{\(\displaystyle \pi_{35}(S^{9})=\langle \sigma_9 \omega_{16} \nu_{32}, \sigma_9 \bar\sigma_{16}, \sigma_9\bar\zeta_{16}, \bar\kappa_9 \nu^2_{29}, \bar\sigma_9 \sigma_{28}, \nu^2_9 \bar\kappa_{15}, \eta_9 \mu_{3,10}\rangle\cong (\mathbb Z/2)^{2} \oplus \mathbb Z/8 \oplus (\mathbb Z/2)^{4}.\)}{\(\displaystyle \begin{aligned}[t]\pi_{35}(S^{9})&=\langle \sigma_9 \omega_{16} \nu_{32}, \sigma_9 \bar\sigma_{16}, \sigma_9\bar\zeta_{16}, \bar\kappa_9 \nu^2_{29}, \bar\sigma_9 \sigma_{28}, \nu^2_9 \bar\kappa_{15},\\&\phantom{=\langle{}}\eta_9 \mu_{3,10}\rangle\cong (\mathbb Z/2)^{2} \oplus \mathbb Z/8 \oplus (\mathbb Z/2)^{4}.\end{aligned}\)}
	\StemEntry{\(\displaystyle \pi_{36}(S^{10})=\langle \tau'', R_1, \sigma_{10} \bar\sigma_{17}, \nu^2_{10} \bar\kappa_{16}, \eta_{10} \mu_{3,11}\rangle\cong \mathbb Z/8 \oplus \mathbb Z/4 \oplus (\mathbb Z/2)^{3}.\)}{\(\displaystyle \begin{aligned}[t]\pi_{36}(S^{10})&=\langle \tau'', R_1, \sigma_{10} \bar\sigma_{17}, \nu^2_{10} \bar\kappa_{16}, \eta_{10} \mu_{3,11}\rangle\\&\cong \mathbb Z/8 \oplus \mathbb Z/4 \oplus (\mathbb Z/2)^{3}.\end{aligned}\)}
	\StemEntry{\(\displaystyle \pi_{37}(S^{11})=\langle \tau''', C_1 \kappa_{23}, \sigma_{11} \bar\sigma_{18}, \nu^2_{11} \bar\kappa_{17}, \eta_{11} \mu_{3,12}\rangle\cong \mathbb Z/8 \oplus (\mathbb Z/2)^{4}.\)}{\(\displaystyle \begin{aligned}[t]\pi_{37}(S^{11})&=\langle \tau''', C_1 \kappa_{23}, \sigma_{11} \bar\sigma_{18}, \nu^2_{11} \bar\kappa_{17}, \eta_{11} \mu_{3,12}\rangle\\&\cong \mathbb Z/8 \oplus (\mathbb Z/2)^{4}.\end{aligned}\)}
	\StemEntry{\(\displaystyle \pi_{38}(S^{12})=\langle \tau^{\mathrm{IV}}, \Sigma \tau'''-8\tau^{\mathrm{IV}}, A_1 \kappa_{24}, \Sigma C_1 \kappa_{24}, \sigma_{12} \bar\sigma_{19}, \nu^2_{12} \bar\kappa_{18}, \eta_{12} \mu_{3,13}\rangle\cong \mathbb Z/64 \oplus \mathbb Z/4 \oplus (\mathbb Z/2)^{5}.\)}{\(\displaystyle \begin{aligned}[t]\pi_{38}(S^{12})&=\langle \tau^{\mathrm{IV}}, \Sigma \tau'''-8\tau^{\mathrm{IV}}, A_1 \kappa_{24}, \Sigma C_1 \kappa_{24}, \sigma_{12} \bar\sigma_{19}, \nu^2_{12} \bar\kappa_{18},\\&\phantom{=\langle{}}\eta_{12} \mu_{3,13}\rangle\cong \mathbb Z/64 \oplus \mathbb Z/4 \oplus (\mathbb Z/2)^{5}.\end{aligned}\)}
	\StemEntry{\(\displaystyle \pi_{39}(S^{13})=\langle \Sigma \tau^{\mathrm{IV}}, (\Sigma A_1) \kappa_{25}, \sigma_{13} \bar\sigma_{20}, \nu^2_{13} \bar\kappa_{19}, \eta_{13} \mu_{3,14}\rangle\cong \mathbb Z/8 \oplus (\mathbb Z/2)^{4}.\)}{\(\displaystyle \begin{aligned}[t]\pi_{39}(S^{13})&=\langle \Sigma \tau^{\mathrm{IV}}, (\Sigma A_1) \kappa_{25}, \sigma_{13} \bar\sigma_{20}, \nu^2_{13} \bar\kappa_{19},\\&\phantom{=\langle{}}\eta_{13} \mu_{3,14}\rangle\cong \mathbb Z/8 \oplus (\mathbb Z/2)^{4}.\end{aligned}\)}
	\StemEntry{\(\displaystyle \pi_{n+26}(S^{n})=\langle \Sigma ^{n-12} \tau^{\mathrm{IV}}, \nu^2_n \bar\kappa_{n+6}, \eta_n \mu_{3,n+1}\rangle\cong \mathbb Z/8 \oplus (\mathbb Z/2)^{2}\quad(n = 14 , 15).\)}{\(\displaystyle \begin{aligned}[t]\pi_{n+26}(S^{n})&=\langle \Sigma ^{n-12} \tau^{\mathrm{IV}}, \nu^2_n \bar\kappa_{n+6}, \eta_n \mu_{3,n+1}\rangle\\&\cong \mathbb Z/8 \oplus (\mathbb Z/2)^{2}\quad(n = 14 , 15).\end{aligned}\)}
	\StemEntry{\(\displaystyle \pi_{42}(S^{16})=\langle M^{(1)}_1, \Sigma ^4 \tau^{\mathrm{IV}}, \nu^2_{16} \bar\kappa_{22}, \eta_{16} \mu_{3,17}\rangle\cong \mathbb Z/8 \oplus \mathbb Z/8 \oplus (\mathbb Z/2)^{2}.\)}{\(\displaystyle \begin{aligned}[t]\pi_{42}(S^{16})&=\langle M^{(1)}_1, \Sigma ^4 \tau^{\mathrm{IV}}, \nu^2_{16} \bar\kappa_{22}, \eta_{16} \mu_{3,17}\rangle\\&\cong \mathbb Z/8 \oplus \mathbb Z/8 \oplus (\mathbb Z/2)^{2}.\end{aligned}\)}
	\StemEntry{\(\displaystyle \pi_{43}(S^{17})=\langle \Sigma M^{(1)}_1, \nu^2_{17} \bar\kappa_{23}, \eta_{17} \mu_{3,18}\rangle\cong \mathbb Z/8 \oplus (\mathbb Z/2)^{2}.\)}{\(\displaystyle \begin{aligned}[t]\pi_{43}(S^{17})&=\langle \Sigma M^{(1)}_1, \nu^2_{17} \bar\kappa_{23}, \eta_{17} \mu_{3,18}\rangle\\&\cong \mathbb Z/8 \oplus (\mathbb Z/2)^{2}.\end{aligned}\)}
	\StemEntry{\(\displaystyle \pi_{44}(S^{18})=\langle M'''_2, \nu^2_{18} \bar\kappa_{24}, \eta_{18} \mu_{3,19}\rangle\cong \mathbb Z/8 \oplus (\mathbb Z/2)^{2}.\)}{\(\displaystyle \begin{aligned}[t]\pi_{44}(S^{18})&=\langle M'''_2, \nu^2_{18} \bar\kappa_{24}, \eta_{18} \mu_{3,19}\rangle\\&\cong \mathbb Z/8 \oplus (\mathbb Z/2)^{2}.\end{aligned}\)}
	\StemEntry{\(\displaystyle \pi_{45}(S^{19})=\langle M''_2, \nu^2_{19} \bar\kappa_{25}, \eta_{19} \mu_{3,20}\rangle\cong \mathbb Z/8 \oplus (\mathbb Z/2)^{2}.\)}{\(\displaystyle \begin{aligned}[t]\pi_{45}(S^{19})&=\langle M''_2, \nu^2_{19} \bar\kappa_{25}, \eta_{19} \mu_{3,20}\rangle\\&\cong \mathbb Z/8 \oplus (\mathbb Z/2)^{2}.\end{aligned}\)}
	\StemEntry{\(\displaystyle \pi_{46}(S^{20})=\langle P(\sigma_{41}), \Sigma M''_2-2P(\sigma_{41}), \nu^2_{20} \bar\kappa_{26}, \eta_{20} \mu_{3,21}\rangle\cong \mathbb Z/16 \oplus \mathbb Z/4 \oplus (\mathbb Z/2)^{2}.\)}{\(\displaystyle \begin{aligned}[t]\pi_{46}(S^{20})&=\langle P(\sigma_{41}), \Sigma M''_2-2P(\sigma_{41}), \nu^2_{20} \bar\kappa_{26},\\&\phantom{=\langle{}}\eta_{20} \mu_{3,21}\rangle\cong \mathbb Z/16 \oplus \mathbb Z/4 \oplus (\mathbb Z/2)^{2}.\end{aligned}\)}
	\StemEntry{\(\displaystyle \pi_{n+26}(S^{n})=\langle \Sigma ^{n-21}M'_2, \nu^2_n \bar\kappa_{n+6}, \eta_n \mu_{3,n+1}\rangle\cong \mathbb Z/8 \oplus (\mathbb Z/2)^{2}\quad(n = 21, 22 , 23).\)}{\(\displaystyle \begin{aligned}[t]\pi_{n+26}(S^{n})&=\langle \Sigma ^{n-21}M'_2, \nu^2_n \bar\kappa_{n+6}, \eta_n \mu_{3,n+1}\rangle\\&\cong \mathbb Z/8 \oplus (\mathbb Z/2)^{2}\quad(n = 21, 22 , 23).\end{aligned}\)}
	\StemEntry{\(\displaystyle \pi_{50}(S^{24})=\langle M_2, \Sigma ^3M'_2, \nu^2_{24} \bar\kappa_{30}, \eta_{24} \mu_{3,25}\rangle\cong \mathbb Z/8 \oplus \mathbb Z/8 \oplus (\mathbb Z/2)^{2}.\)}{\(\displaystyle \begin{aligned}[t]\pi_{50}(S^{24})&=\langle M_2, \Sigma ^3M'_2, \nu^2_{24} \bar\kappa_{30}, \eta_{24} \mu_{3,25}\rangle\\&\cong \mathbb Z/8 \oplus \mathbb Z/8 \oplus (\mathbb Z/2)^{2}.\end{aligned}\)}
	\StemEntry{\(\displaystyle \pi_{51}(S^{25})=\langle \Sigma M_2, \nu^2_{25} \bar\kappa_{31}, \eta_{25} \mu_{3,26}\rangle\cong \mathbb Z/8 \oplus (\mathbb Z/2)^{2}.\)}{\(\displaystyle \begin{aligned}[t]\pi_{51}(S^{25})&=\langle \Sigma M_2, \nu^2_{25} \bar\kappa_{31}, \eta_{25} \mu_{3,26}\rangle\\&\cong \mathbb Z/8 \oplus (\mathbb Z/2)^{2}.\end{aligned}\)}
	\StemEntry{\(\displaystyle \pi_{52}(S^{26})=\langle \Sigma ^2M_2, \nu^2_{26} \bar\kappa_{32}, \eta_{26} \mu_{3,27}\rangle\cong \mathbb Z/4 \oplus (\mathbb Z/2)^{2}.\)}{\(\displaystyle \begin{aligned}[t]\pi_{52}(S^{26})&=\langle \Sigma ^2M_2, \nu^2_{26} \bar\kappa_{32}, \eta_{26} \mu_{3,27}\rangle\\&\cong \mathbb Z/4 \oplus (\mathbb Z/2)^{2}.\end{aligned}\)}
	\StemEntry{\(\displaystyle \pi_{53}(S^{27})=\langle \Sigma ^3M_2, \nu^2_{27} \bar\kappa_{33}, \eta_{27} \mu_{3,28}\rangle\cong (\mathbb Z/2)^{3}.\)}{\(\displaystyle \begin{aligned}[t]\pi_{53}(S^{27})&=\langle \Sigma ^3M_2, \nu^2_{27} \bar\kappa_{33}, \eta_{27} \mu_{3,28}\rangle\\&\cong (\mathbb Z/2)^{3}.\end{aligned}\)}
	\StemEntry{\(\displaystyle \pi_{n+26}(S^{n})=\langle \nu^2_n \bar\kappa_{n+6}, \eta_n \mu_{3,n+1}\rangle\cong (\mathbb Z/2)^{2}\quad(n\ge 28).\)}{\(\displaystyle \begin{aligned}[t]\pi_{n+26}(S^{n})&=\langle \nu^2_n \bar\kappa_{n+6}, \eta_n \mu_{3,n+1}\rangle\\&\cong (\mathbb Z/2)^{2}\quad(n\ge 28).\end{aligned}\)}
	\StemEntry{\(\displaystyle \pi^S_{26}(S^0)=\langle \nu^2 \bar\kappa, \eta \mu_{3}\rangle\cong \mathbb Z/2 \oplus \mathbb Z/2.\)}{\(\displaystyle \begin{aligned}[t]\pi^S_{26}(S^0)&=\langle \nu^2 \bar\kappa, \eta \mu_{3}\rangle\\&\cong \mathbb Z/2 \oplus \mathbb Z/2.\end{aligned}\)}
	\subsection{Stem 27}\label{app:stem-27}
	\StemEntry{\(\displaystyle \pi_{29}(S^{2})=\langle \eta_2 \bar\mu' \sigma_{22}, \eta_2 \nu' \phi_6, \eta_2 \bar\alpha \nu_{26}, \eta^2_2\mu_{3,4}\rangle\cong \mathbb Z/4 \oplus (\mathbb Z/2)^{3}.\)}{\(\displaystyle \begin{aligned}[t]\pi_{29}(S^{2})&=\langle \eta_2 \bar\mu' \sigma_{22}, \eta_2 \nu' \phi_6, \eta_2 \bar\alpha \nu_{26},\\&\phantom{=\langle{}}\eta^2_2\mu_{3,4}\rangle\cong \mathbb Z/4 \oplus (\mathbb Z/2)^{3}.\end{aligned}\)}
	\StemEntry{\(\displaystyle \pi_{30}(S^{3})=\langle \nu' \delta_6, \nu' \bar\mu_6 \sigma_{23}, \varepsilon _3 \bar\sigma_{11}, \mu'_3\rangle\cong (\mathbb Z/2)^{3} \oplus \mathbb Z/4.\)}{\(\displaystyle \begin{aligned}[t]\pi_{30}(S^{3})&=\langle \nu' \delta_6, \nu' \bar\mu_6 \sigma_{23}, \varepsilon _3 \bar\sigma_{11}, \mu'_3\rangle\\&\cong (\mathbb Z/2)^{3} \oplus \mathbb Z/4.\end{aligned}\)}
	\StemEntry{\(\displaystyle \pi_{31}(S^{4})=\langle \nu_4 \sigma' \omega_{14} \eta_{30}, \nu_4 \sigma' \bar\mu_{14}, \nu_4 \bar\sigma'_7, \nu_4 \bar\zeta'_7, \nu^2_4 \sigma_{10} \kappa_{17}, \nu_4 \bar\mu_7\sigma_{24}, \nu_4 \delta_7, (\Sigma \nu') \delta_7, (\Sigma \nu') \bar\mu_7 \sigma_{24}, \varepsilon _4 \bar\sigma_{12}, \Sigma \mu'_3\rangle\cong (\mathbb Z/2)^{10} \oplus \mathbb Z/4.\)}{\(\displaystyle \begin{aligned}[t]\pi_{31}(S^{4})&=\langle \nu_4 \sigma' \omega_{14} \eta_{30}, \nu_4 \sigma' \bar\mu_{14}, \nu_4 \bar\sigma'_7, \nu_4 \bar\zeta'_7, \nu^2_4 \sigma_{10} \kappa_{17}, \nu_4 \bar\mu_7\sigma_{24}, \nu_4 \delta_7, (\Sigma \nu') \delta_7, (\Sigma \nu') \bar\mu_7 \sigma_{24},\\&\phantom{=\langle{}}\varepsilon _4 \bar\sigma_{12}, \Sigma \mu'_3\rangle\cong (\mathbb Z/2)^{10} \oplus \mathbb Z/4.\end{aligned}\)}
	\StemEntry{\(\displaystyle \pi_{32}(S^{5})=\langle \sigma''' \bar\kappa_{12}, \nu_5 \bar\sigma'_8, \nu_5 \bar\mu_8 \sigma_{25}, \nu_5 \delta_8, \zeta_{3,5}\rangle\cong (\mathbb Z/2)^{4} \oplus \mathbb Z/8.\)}{\(\displaystyle \begin{aligned}[t]\pi_{32}(S^{5})&=\langle \sigma''' \bar\kappa_{12}, \nu_5 \bar\sigma'_8, \nu_5 \bar\mu_8 \sigma_{25}, \nu_5 \delta_8, \zeta_{3,5}\rangle\\&\cong (\mathbb Z/2)^{4} \oplus \mathbb Z/8.\end{aligned}\)}
	\StemEntry{\(\displaystyle \pi_{33}(S^{6})=\langle P(\rho_{13} \sigma_{28}), \sigma'' \bar\kappa_{13}, \bar\nu_6 \bar\sigma_{14}, \zeta_{3,6}\rangle\cong \mathbb Z/16 \oplus \mathbb Z/4 \oplus \mathbb Z/2 \oplus \mathbb Z/8.\)}{\(\displaystyle \begin{aligned}[t]\pi_{33}(S^{6})&=\langle P(\rho_{13} \sigma_{28}), \sigma'' \bar\kappa_{13}, \bar\nu_6 \bar\sigma_{14},\\&\phantom{=\langle{}}\zeta_{3,6}\rangle\cong \mathbb Z/16 \oplus \mathbb Z/4 \oplus \mathbb Z/2 \oplus \mathbb Z/8.\end{aligned}\)}
	\StemEntry{\(\displaystyle \pi_{34}(S^{7})=\langle \sigma' \bar\kappa_{14}, \bar\nu_7 \bar\sigma_{15}, \zeta_{3,7}\rangle\cong \mathbb Z/8 \oplus \mathbb Z/2 \oplus \mathbb Z/8.\)}{\(\displaystyle \begin{aligned}[t]\pi_{34}(S^{7})&=\langle \sigma' \bar\kappa_{14}, \bar\nu_7 \bar\sigma_{15}, \zeta_{3,7}\rangle\\&\cong \mathbb Z/8 \oplus \mathbb Z/2 \oplus \mathbb Z/8.\end{aligned}\)}
	\StemEntry{\(\displaystyle \pi_{35}(S^{8})=\langle \sigma_8 \bar\kappa_{15}, (\Sigma \sigma') \bar\kappa_{15}, \bar\nu_8 \bar\sigma_{16}, \zeta_{3,8}\rangle\cong \mathbb Z/8 \oplus \mathbb Z/8 \oplus \mathbb Z/2 \oplus \mathbb Z/8.\)}{\(\displaystyle \begin{aligned}[t]\pi_{35}(S^{8})&=\langle \sigma_8 \bar\kappa_{15}, (\Sigma \sigma') \bar\kappa_{15}, \bar\nu_8 \bar\sigma_{16},\\&\phantom{=\langle{}}\zeta_{3,8}\rangle\cong \mathbb Z/8 \oplus \mathbb Z/8 \oplus \mathbb Z/2 \oplus \mathbb Z/8.\end{aligned}\)}
	\StemEntry{\(\displaystyle \pi_{36}(S^{9})=\langle \sigma_9 \bar\kappa_{16}, \bar\nu_9 \bar\sigma_{17}, \zeta_{3,9}\rangle\cong \mathbb Z/8 \oplus \mathbb Z/2 \oplus \mathbb Z/8.\)}{\(\displaystyle \begin{aligned}[t]\pi_{36}(S^{9})&=\langle \sigma_9 \bar\kappa_{16}, \bar\nu_9 \bar\sigma_{17}, \zeta_{3,9}\rangle\\&\cong \mathbb Z/8 \oplus \mathbb Z/2 \oplus \mathbb Z/8.\end{aligned}\)}
	\StemEntry{\(\displaystyle \pi_{37}(S^{10})=\langle P(\nu^*_{21}), \sigma_{10} \bar\kappa_{17}, \zeta_{3,10}\rangle\cong \mathbb Z/4 \oplus \mathbb Z/8 \oplus \mathbb Z/8.\)}{\(\displaystyle \begin{aligned}[t]\pi_{37}(S^{10})&=\langle P(\nu^*_{21}), \sigma_{10} \bar\kappa_{17}, \zeta_{3,10}\rangle\\&\cong \mathbb Z/4 \oplus \mathbb Z/8 \oplus \mathbb Z/8.\end{aligned}\)}
	\StemEntry{\(\displaystyle \pi_{38}(S^{11})=\langle C^{\mathrm{II}}_1, C^{(1)}_1 \sigma_{31}, \sigma_{11} \bar\kappa_{18}, \zeta_{3,11}\rangle\cong (\mathbb Z/2)^{2} \oplus \mathbb Z/8 \oplus \mathbb Z/8.\)}{\(\displaystyle \begin{aligned}[t]\pi_{38}(S^{11})&=\langle C^{\mathrm{II}}_1, C^{(1)}_1 \sigma_{31}, \sigma_{11} \bar\kappa_{18}, \zeta_{3,11}\rangle\\&\cong (\mathbb Z/2)^{2} \oplus \mathbb Z/8 \oplus \mathbb Z/8.\end{aligned}\)}
	\StemEntry{\(\displaystyle \pi_{39}(S^{12})=\langle A^{\mathrm{II}}_1, A^{(1)}_1 \sigma_{32}, \Sigma C^{\mathrm{II}}_1, \Sigma C^{(1)}_1 \sigma_{32}, \sigma_{12} \bar\kappa_{19}, \zeta_{3,12}\rangle\cong (\mathbb Z/2)^{4} \oplus \mathbb Z/4 \oplus \mathbb Z/8.\)}{\(\displaystyle \begin{aligned}[t]\pi_{39}(S^{12})&=\langle A^{\mathrm{II}}_1, A^{(1)}_1 \sigma_{32}, \Sigma C^{\mathrm{II}}_1, \Sigma C^{(1)}_1 \sigma_{32}, \sigma_{12} \bar\kappa_{19},\\&\phantom{=\langle{}}\zeta_{3,12}\rangle\cong (\mathbb Z/2)^{4} \oplus \mathbb Z/4 \oplus \mathbb Z/8.\end{aligned}\)}
	\StemEntry{\(\displaystyle \pi_{40}(S^{13})=\langle \Sigma A^{\mathrm{II}}_1, \Sigma A^{(1)}_1 \sigma_{33}, \sigma_{13} \bar\kappa_{20}, \zeta_{3,13}\rangle\cong (\mathbb Z/2)^{2} \oplus \mathbb Z/4 \oplus \mathbb Z/8.\)}{\(\displaystyle \begin{aligned}[t]\pi_{40}(S^{13})&=\langle \Sigma A^{\mathrm{II}}_1, \Sigma A^{(1)}_1 \sigma_{33}, \sigma_{13} \bar\kappa_{20},\\&\phantom{=\langle{}}\zeta_{3,13}\rangle\cong (\mathbb Z/2)^{2} \oplus \mathbb Z/4 \oplus \mathbb Z/8.\end{aligned}\)}
	\StemEntry{\(\displaystyle \pi_{41}(S^{14})=\langle \Sigma ^2A^{\mathrm{II}}_1, \sigma_{14} \bar\kappa_{21}, \zeta_{3,14}\rangle\cong (\mathbb Z/2)^{2} \oplus \mathbb Z/8.\)}{\(\displaystyle \begin{aligned}[t]\pi_{41}(S^{14})&=\langle \Sigma ^2A^{\mathrm{II}}_1, \sigma_{14} \bar\kappa_{21}, \zeta_{3,14}\rangle\\&\cong (\mathbb Z/2)^{2} \oplus \mathbb Z/8.\end{aligned}\)}
	\StemEntry{\(\displaystyle \pi_{n+27}(S^{n})=\langle \zeta_{3,n}\rangle\cong \mathbb Z/8\quad(n = 15, 16, 17, 18).\)}{\(\displaystyle \begin{aligned}[t]\pi_{n+27}(S^{n})&=\langle \zeta_{3,n}\rangle\\&\cong \mathbb Z/8\quad(n = 15, 16, 17, 18).\end{aligned}\)}
	\StemEntry{\(\displaystyle \pi_{46}(S^{19})=\langle C^{\mathrm{I}}_2, C_2 \sigma_{39}, \zeta_{3,19}\rangle\cong (\mathbb Z/2)^{2} \oplus \mathbb Z/8.\)}{\(\displaystyle \begin{aligned}[t]\pi_{46}(S^{19})&=\langle C^{\mathrm{I}}_2, C_2 \sigma_{39}, \zeta_{3,19}\rangle\\&\cong (\mathbb Z/2)^{2} \oplus \mathbb Z/8.\end{aligned}\)}
	\StemEntry{\(\displaystyle \pi_{47}(S^{20})=\langle A^{\mathrm{I}}_2, A_2 \sigma_{40}, \Sigma C^{\mathrm{I}}_2, (\Sigma C_2) \sigma_{40}, \zeta_{3,20}\rangle\cong (\mathbb Z/2)^{4} \oplus \mathbb Z/8.\)}{\(\displaystyle \begin{aligned}[t]\pi_{47}(S^{20})&=\langle A^{\mathrm{I}}_2, A_2 \sigma_{40}, \Sigma C^{\mathrm{I}}_2, (\Sigma C_2) \sigma_{40}, \zeta_{3,20}\rangle\\&\cong (\mathbb Z/2)^{4} \oplus \mathbb Z/8.\end{aligned}\)}
	\StemEntry{\(\displaystyle \pi_{48}(S^{21})=\langle \Sigma A^{\mathrm{I}}_2, \Sigma A_2 \sigma_{41}, \zeta_{3,21}\rangle\cong (\mathbb Z/2)^{2} \oplus \mathbb Z/8.\)}{\(\displaystyle \begin{aligned}[t]\pi_{48}(S^{21})&=\langle \Sigma A^{\mathrm{I}}_2, \Sigma A_2 \sigma_{41}, \zeta_{3,21}\rangle\\&\cong (\mathbb Z/2)^{2} \oplus \mathbb Z/8.\end{aligned}\)}
	\StemEntry{\(\displaystyle \pi_{49}(S^{22})=\langle V_2, \zeta_{3,22}\rangle\cong \mathbb Z/4 \oplus \mathbb Z/8.\)}{\(\displaystyle \begin{aligned}[t]\pi_{49}(S^{22})&=\langle V_2, \zeta_{3,22}\rangle\\&\cong \mathbb Z/4 \oplus \mathbb Z/8.\end{aligned}\)}
	\StemEntry{\(\displaystyle \pi_{n+27}(S^{n})=\langle \Sigma ^{n-22} V_2, \zeta_{3,n}\rangle\cong \mathbb Z/2 \oplus \mathbb Z/8\quad(n = 23, 24, 25).\)}{\(\displaystyle \begin{aligned}[t]\pi_{n+27}(S^{n})&=\langle \Sigma ^{n-22} V_2, \zeta_{3,n}\rangle\\&\cong \mathbb Z/2 \oplus \mathbb Z/8\quad(n = 23, 24, 25).\end{aligned}\)}
	\StemEntry{\(\displaystyle \pi_{55}(S^{28})=\langle P(\iota_{57}), \zeta_{3,28}\rangle\cong \mathbb Z_{(2)} \oplus \mathbb Z/8.\)}{\(\displaystyle \begin{aligned}[t]\pi_{55}(S^{28})&=\langle P(\iota_{57}), \zeta_{3,28}\rangle\\&\cong \mathbb Z_{(2)} \oplus \mathbb Z/8.\end{aligned}\)}
	\StemEntry{\(\displaystyle \pi_{n+27}(S^{n})=\langle \zeta_{3,n}\rangle\cong \mathbb Z/8\quad(n=26,27\text{ or }n\ge29).\)}{\(\displaystyle \begin{aligned}[t]\pi_{n+27}(S^{n})&=\langle \zeta_{3,n}\rangle\\&\cong \mathbb Z/8\quad(n=26,27\text{ or }n\ge29).\end{aligned}\)}
	\StemEntry{\(\displaystyle \pi^S_{27}(S^0)=\langle \zeta_{3}\rangle\cong \mathbb Z/8.\)}{\(\displaystyle \begin{aligned}[t]\pi^S_{27}(S^0)&=\langle \zeta_{3}\rangle\\&\cong \mathbb Z/8.\end{aligned}\)}
	\subsection{Stem 28}\label{app:stem-28}
	\StemEntry{\(\displaystyle \pi_{30}(S^{2})=\langle \eta_2 \nu' \delta_6, \eta_2 \nu' \bar\mu_6 \sigma_{23}, \eta_2 \varepsilon _3 \bar\sigma_{11}, \eta_2 \mu'_3\rangle\cong (\mathbb Z/2)^{3} \oplus \mathbb Z/4.\)}{\(\displaystyle \begin{aligned}[t]\pi_{30}(S^{2})&=\langle \eta_2 \nu' \delta_6, \eta_2 \nu' \bar\mu_6 \sigma_{23}, \eta_2 \varepsilon _3 \bar\sigma_{11},\\&\phantom{=\langle{}}\eta_2 \mu'_3\rangle\cong (\mathbb Z/2)^{3} \oplus \mathbb Z/4.\end{aligned}\)}
	\StemEntry{\(\displaystyle \pi_{31}(S^{3})=\langle \nu' \eta_6 \bar\mu_7 \sigma_{24}, \nu' \mu_{3,6}, \phi' \nu_{28}, \varepsilon _3 \bar\kappa_{11}\rangle\cong (\mathbb Z/2)^{4}.\)}{\(\displaystyle \begin{aligned}[t]\pi_{31}(S^{3})&=\langle \nu' \eta_6 \bar\mu_7 \sigma_{24}, \nu' \mu_{3,6}, \phi' \nu_{28},\\&\phantom{=\langle{}}\varepsilon _3 \bar\kappa_{11}\rangle\cong (\mathbb Z/2)^{4}.\end{aligned}\)}
	\StemEntry{\(\displaystyle \pi_{32}(S^{4})=\langle \nu_4(\Sigma^2\phi''), \nu_4u_7, \nu_4 \sigma' \eta_{14} \bar\mu_{15}, \nu_4 \mu_{3,7}, \nu_4 \eta_7 \bar\mu_8 \sigma_{25}, (\Sigma \nu') \eta_7 \bar\mu_8 \sigma_{25}, (\Sigma \nu') \mu_{3,7}, (\Sigma \phi') \nu_{29}, \varepsilon _4 \bar\kappa_{12}\rangle\cong \mathbb Z/8 \oplus \mathbb Z/4 \oplus (\mathbb Z/2)^{7}.\)}{\(\displaystyle \begin{aligned}[t]\pi_{32}(S^{4})&=\langle \nu_4(\Sigma^2\phi''), \nu_4u_7, \nu_4 \sigma' \eta_{14} \bar\mu_{15}, \nu_4 \mu_{3,7}, \nu_4 \eta_7 \bar\mu_8 \sigma_{25}, (\Sigma \nu') \eta_7 \bar\mu_8 \sigma_{25}, (\Sigma \nu') \mu_{3,7},\\&\phantom{=\langle{}}(\Sigma \phi') \nu_{29}, \varepsilon _4 \bar\kappa_{12}\rangle\cong \mathbb Z/8 \oplus \mathbb Z/4 \oplus (\mathbb Z/2)^{7}.\end{aligned}\)}
	\StemEntry{\(\displaystyle \pi_{33}(S^{5})=\langle \nu_5 (\Sigma \phi'''), \nu_5 \mu_{3,8}, \nu_5 \eta_8 \bar\mu_9 \sigma_{26}, \varepsilon _5 \bar\kappa_{13}\rangle\cong \mathbb Z/4 \oplus (\mathbb Z/2)^{3}.\)}{\(\displaystyle \begin{aligned}[t]\pi_{33}(S^{5})&=\langle \nu_5 (\Sigma \phi'''), \nu_5 \mu_{3,8}, \nu_5 \eta_8 \bar\mu_9 \sigma_{26},\\&\phantom{=\langle{}}\varepsilon _5 \bar\kappa_{13}\rangle\cong \mathbb Z/4 \oplus (\mathbb Z/2)^{3}.\end{aligned}\)}
	\StemEntry{\(\displaystyle \pi_{34}(S^{6})=\langle P(\bar\rho_{13}), \bar\nu_6 \bar\kappa_{14}, \nu_6( \Sigma ^2\phi'''), \varepsilon _6 \bar\kappa_{14}\rangle\cong \mathbb Z/16 \oplus \mathbb Z/8 \oplus (\mathbb Z/2)^{2}.\)}{\(\displaystyle \begin{aligned}[t]\pi_{34}(S^{6})&=\langle P(\bar\rho_{13}), \bar\nu_6 \bar\kappa_{14}, \nu_6( \Sigma ^2\phi'''),\\&\phantom{=\langle{}}\varepsilon _6 \bar\kappa_{14}\rangle\cong \mathbb Z/16 \oplus \mathbb Z/8 \oplus (\mathbb Z/2)^{2}.\end{aligned}\)}
	\StemEntry{\(\displaystyle \pi_{35}(S^{7})=\langle \phi''' \nu_{32}, \sigma' \eta_{14} \bar\kappa_{15}, \bar\nu_7 \bar\kappa_{15}, \varepsilon _7 \bar\kappa_{15}\rangle\cong (\mathbb Z/2)^{4}.\)}{\(\displaystyle \begin{aligned}[t]\pi_{35}(S^{7})&=\langle \phi''' \nu_{32}, \sigma' \eta_{14} \bar\kappa_{15}, \bar\nu_7 \bar\kappa_{15},\\&\phantom{=\langle{}}\varepsilon _7 \bar\kappa_{15}\rangle\cong (\mathbb Z/2)^{4}.\end{aligned}\)}
	\StemEntry{\(\displaystyle \pi_{36}(S^{8})=\langle \sigma_8 (\Sigma ^2\lambda) \nu_{33}, \sigma^4_8, \sigma_8 \eta_{15} \bar\kappa_{16}, (\Sigma \phi''') \nu_{33}, (\Sigma \sigma' )\eta_{15} \bar\kappa_{16}, \bar\nu_8 \bar\kappa_{16}, \varepsilon _8 \bar\kappa_{16}\rangle\cong (\mathbb Z/2)^{7}.\)}{\(\displaystyle \begin{aligned}[t]\pi_{36}(S^{8})&=\langle \sigma_8 (\Sigma ^2\lambda) \nu_{33}, \sigma^4_8, \sigma_8 \eta_{15} \bar\kappa_{16}, (\Sigma \phi''') \nu_{33}, (\Sigma \sigma' )\eta_{15} \bar\kappa_{16},\\&\phantom{=\langle{}}\bar\nu_8 \bar\kappa_{16}, \varepsilon _8 \bar\kappa_{16}\rangle\cong (\mathbb Z/2)^{7}.\end{aligned}\)}
	\StemEntry{\(\displaystyle \pi_{37}(S^{9})=\langle \sigma^4_9, \sigma_9 \eta_{16} \bar\kappa_{17}, \sigma_9 \nu^*_{16} \nu_{34}, \bar\nu_9 \bar\kappa_{17}, \varepsilon _9 \bar\kappa_{17}\rangle\cong (\mathbb Z/2)^{5}.\)}{\(\displaystyle \begin{aligned}[t]\pi_{37}(S^{9})&=\langle \sigma^4_9, \sigma_9 \eta_{16} \bar\kappa_{17}, \sigma_9 \nu^*_{16} \nu_{34}, \bar\nu_9 \bar\kappa_{17},\\&\phantom{=\langle{}}\varepsilon _9 \bar\kappa_{17}\rangle\cong (\mathbb Z/2)^{5}.\end{aligned}\)}
	\StemEntry{\(\displaystyle \pi_{38}(S^{10})=\langle F^{(1)}_1, \sigma^4_{10}, \sigma_{10} \nu^*_{17} \nu_{35}, \bar\nu_{10} \bar\kappa_{18}, \varepsilon _{10} \bar\kappa_{18}\rangle\cong \mathbb Z/8 \oplus (\mathbb Z/2)^{4}.\)}{\(\displaystyle \begin{aligned}[t]\pi_{38}(S^{10})&=\langle F^{(1)}_1, \sigma^4_{10}, \sigma_{10} \nu^*_{17} \nu_{35}, \bar\nu_{10} \bar\kappa_{18},\\&\phantom{=\langle{}}\varepsilon _{10} \bar\kappa_{18}\rangle\cong \mathbb Z/8 \oplus (\mathbb Z/2)^{4}.\end{aligned}\)}
	\StemEntry{\(\displaystyle \pi_{39}(S^{11})=\langle C_1 \omega_{23}, C^{(2)}_1, \Sigma F^{(1)}_1, \sigma^4_{11}, \bar\nu_{11} \bar\kappa_{19}, \varepsilon _{11} \bar\kappa_{19}\rangle\cong (\mathbb Z/2)^{6}.\)}{\(\displaystyle \begin{aligned}[t]\pi_{39}(S^{11})&=\langle C_1 \omega_{23}, C^{(2)}_1, \Sigma F^{(1)}_1, \sigma^4_{11}, \bar\nu_{11} \bar\kappa_{19}, \varepsilon _{11} \bar\kappa_{19}\rangle\\&\cong (\mathbb Z/2)^{6}.\end{aligned}\)}
	\StemEntry{\(\displaystyle \pi_{40}(S^{12})=\langle A^{(2)}_1, A_1 \omega_{24}, A_1 \sigma_{24} \mu_{31}, \Sigma C^{(2)}_1, (\Sigma C_1) \omega_{24}, \Sigma ^2F^{(1)}_1, \bar\nu_{12} \bar\kappa_{20}, \varepsilon _{12} \bar\kappa_{20}\rangle\cong (\mathbb Z/2)^{8}.\)}{\(\displaystyle \begin{aligned}[t]\pi_{40}(S^{12})&=\langle A^{(2)}_1, A_1 \omega_{24}, A_1 \sigma_{24} \mu_{31}, \Sigma C^{(2)}_1, (\Sigma C_1) \omega_{24}, \Sigma ^2F^{(1)}_1, \bar\nu_{12} \bar\kappa_{20},\\&\phantom{=\langle{}}\varepsilon _{12} \bar\kappa_{20}\rangle\cong (\mathbb Z/2)^{8}.\end{aligned}\)}
	\StemEntry{\(\displaystyle \pi_{41}(S^{13})=\langle \Sigma A^{(2)}_1, (\Sigma A_1 )\omega_{25}, (\Sigma A_1) \sigma_{25} \mu_{32}, \varepsilon _{13} \bar\kappa_{21}\rangle\cong (\mathbb Z/2)^{4}.\)}{\(\displaystyle \begin{aligned}[t]\pi_{41}(S^{13})&=\langle \Sigma A^{(2)}_1, (\Sigma A_1 )\omega_{25}, (\Sigma A_1) \sigma_{25} \mu_{32},\\&\phantom{=\langle{}}\varepsilon _{13} \bar\kappa_{21}\rangle\cong (\mathbb Z/2)^{4}.\end{aligned}\)}
	\StemEntry{\(\displaystyle \pi_{42}(S^{14})=\langle P(\rho_{29}), \varepsilon _{14} \bar\kappa_{22}\rangle\cong \mathbb Z/32 \oplus \mathbb Z/2.\)}{\(\displaystyle \begin{aligned}[t]\pi_{42}(S^{14})&=\langle P(\rho_{29}), \varepsilon _{14} \bar\kappa_{22}\rangle\\&\cong \mathbb Z/32 \oplus \mathbb Z/2.\end{aligned}\)}
	\StemEntry{\(\displaystyle \pi_{n+28}(S^{n})=\langle \varepsilon _n \bar\kappa_{n+8}\rangle\cong \mathbb Z/2\quad(n = 15, 16, 17).\)}{\(\displaystyle \begin{aligned}[t]\pi_{n+28}(S^{n})&=\langle \varepsilon _n \bar\kappa_{n+8}\rangle\\&\cong \mathbb Z/2\quad(n = 15, 16, 17).\end{aligned}\)}
	\StemEntry{\(\displaystyle \pi_{46}(S^{18})=\langle F_2, \varepsilon _{18} \bar\kappa_{26}\rangle\cong \mathbb Z/8 \oplus \mathbb Z/2.\)}{\(\displaystyle \begin{aligned}[t]\pi_{46}(S^{18})&=\langle F_2, \varepsilon _{18} \bar\kappa_{26}\rangle\\&\cong \mathbb Z/8 \oplus \mathbb Z/2.\end{aligned}\)}
	\StemEntry{\(\displaystyle \pi_{47}(S^{19})=\langle C^{(1)}_2, \Sigma F_2, \varepsilon _{19} \bar\kappa_{27}\rangle\cong (\mathbb Z/2)^{3}.\)}{\(\displaystyle \begin{aligned}[t]\pi_{47}(S^{19})&=\langle C^{(1)}_2, \Sigma F_2, \varepsilon _{19} \bar\kappa_{27}\rangle\\&\cong (\mathbb Z/2)^{3}.\end{aligned}\)}
	\StemEntry{\(\displaystyle \pi_{48}(S^{20})=\langle A^{(1)}_2, A_2 \varepsilon _{40}, A_2 \eta_{40} \sigma_{41}, \Sigma C^{(1)}_2, \Sigma ^2F_2, \varepsilon _{20} \bar\kappa_{28}\rangle\cong (\mathbb Z/2)^{6}.\)}{\(\displaystyle \begin{aligned}[t]\pi_{48}(S^{20})&=\langle A^{(1)}_2, A_2 \varepsilon _{40}, A_2 \eta_{40} \sigma_{41}, \Sigma C^{(1)}_2, \Sigma ^2F_2,\\&\phantom{=\langle{}}\varepsilon _{20} \bar\kappa_{28}\rangle\cong (\mathbb Z/2)^{6}.\end{aligned}\)}
	\StemEntry{\(\displaystyle \pi_{49}(S^{21})=\langle \Sigma A^{(1)}_2, (\Sigma A_2) \varepsilon _{41}, (\Sigma A_2 )\eta_{41} \sigma_{42}, \varepsilon _{21} \bar\kappa_{29}\rangle\cong (\mathbb Z/2)^{4}.\)}{\(\displaystyle \begin{aligned}[t]\pi_{49}(S^{21})&=\langle \Sigma A^{(1)}_2, (\Sigma A_2) \varepsilon _{41}, (\Sigma A_2 )\eta_{41} \sigma_{42},\\&\phantom{=\langle{}}\varepsilon _{21} \bar\kappa_{29}\rangle\cong (\mathbb Z/2)^{4}.\end{aligned}\)}
	\StemEntry{\(\displaystyle \pi_{50}(S^{22})=\langle P(\sigma_{45}), \varepsilon _{22} \bar\kappa_{30}\rangle\cong \mathbb Z/16 \oplus \mathbb Z/2.\)}{\(\displaystyle \begin{aligned}[t]\pi_{50}(S^{22})&=\langle P(\sigma_{45}), \varepsilon _{22} \bar\kappa_{30}\rangle\\&\cong \mathbb Z/16 \oplus \mathbb Z/2.\end{aligned}\)}
	\StemEntry{\(\displaystyle \pi_{n+28}(S^{n})=\langle \varepsilon_n\bar\kappa_{n+8}\rangle\cong \mathbb Z/2\quad(n = 23, 24, 25).\)}{\(\displaystyle \begin{aligned}[t]\pi_{n+28}(S^{n})&=\langle \varepsilon_n\bar\kappa_{n+8}\rangle\\&\cong \mathbb Z/2\quad(n = 23, 24, 25).\end{aligned}\)}
	\StemEntry{\(\displaystyle \pi_{54}(S^{26})=\langle P(\nu_{53}), \varepsilon _{26} \bar\kappa_{34}\rangle\cong \mathbb Z/4 \oplus \mathbb Z/2.\)}{\(\displaystyle \begin{aligned}[t]\pi_{54}(S^{26})&=\langle P(\nu_{53}), \varepsilon _{26} \bar\kappa_{34}\rangle\\&\cong \mathbb Z/4 \oplus \mathbb Z/2.\end{aligned}\)}
	\StemEntry{\(\displaystyle \pi_{55}(S^{27})=\langle C_3, \varepsilon _{27} \bar\kappa_{35}\rangle\cong (\mathbb Z/2)^{2}.\)}{\(\displaystyle \begin{aligned}[t]\pi_{55}(S^{27})&=\langle C_3, \varepsilon _{27} \bar\kappa_{35}\rangle\\&\cong (\mathbb Z/2)^{2}.\end{aligned}\)}
	\StemEntry{\(\displaystyle \pi_{56}(S^{28})=\langle A_3, \Sigma C_3, \varepsilon _{28} \bar\kappa_{36}\rangle\cong (\mathbb Z/2)^{3}.\)}{\(\displaystyle \begin{aligned}[t]\pi_{56}(S^{28})&=\langle A_3, \Sigma C_3, \varepsilon _{28} \bar\kappa_{36}\rangle\\&\cong (\mathbb Z/2)^{3}.\end{aligned}\)}
	\StemEntry{\(\displaystyle \pi_{57}(S^{29})=\langle \Sigma A_3, \varepsilon _{29} \bar\kappa_{37}\rangle\cong (\mathbb Z/2)^{2}.\)}{\(\displaystyle \begin{aligned}[t]\pi_{57}(S^{29})&=\langle \Sigma A_3, \varepsilon _{29} \bar\kappa_{37}\rangle\\&\cong (\mathbb Z/2)^{2}.\end{aligned}\)}
	\StemEntry{\(\displaystyle \pi_{n+28}(S^{n})=\langle \varepsilon _n \bar\kappa_{n+8}\rangle\cong \mathbb Z/2\quad(n\ge 30).\)}{\(\displaystyle \begin{aligned}[t]\pi_{n+28}(S^{n})&=\langle \varepsilon _n \bar\kappa_{n+8}\rangle\\&\cong \mathbb Z/2\quad(n\ge 30).\end{aligned}\)}
	\StemEntry{\(\displaystyle \pi^S_{28}(S^0)=\langle \varepsilon  \bar\kappa\rangle\cong \mathbb Z/2.\)}{\(\displaystyle \begin{aligned}[t]\pi^S_{28}(S^0)&=\langle \varepsilon  \bar\kappa\rangle\\&\cong \mathbb Z/2.\end{aligned}\)}
	\subsection{Stem 29}\label{app:stem-29}
	\StemEntry{\(\displaystyle \pi_{31}(S^{2})=\langle \eta_2 \nu' \eta_6 \bar\mu_7 \sigma_{24}, \eta_2 \nu' \mu_{3,6}, \eta_2 \phi' \nu_{28}, \eta_2 \varepsilon _3 \bar\kappa_{11}\rangle\cong (\mathbb Z/2)^{4}.\)}{\(\displaystyle \begin{aligned}[t]\pi_{31}(S^{2})&=\langle \eta_2 \nu' \eta_6 \bar\mu_7 \sigma_{24}, \eta_2 \nu' \mu_{3,6}, \eta_2 \phi' \nu_{28},\\&\phantom{=\langle{}}\eta_2 \varepsilon _3 \bar\kappa_{11}\rangle\cong (\mathbb Z/2)^{4}.\end{aligned}\)}
	\StemEntry{\(\displaystyle \pi_{32}(S^{3})=\langle \bar\alpha \nu^2_{26}, \nu' \eta_6 \mu_{3,7}, \eta_3 \varepsilon _4 \bar\kappa_{12}\rangle\cong (\mathbb Z/2)^{3}.\)}{\(\displaystyle \begin{aligned}[t]\pi_{32}(S^{3})&=\langle \bar\alpha \nu^2_{26}, \nu' \eta_6 \mu_{3,7}, \eta_3 \varepsilon _4 \bar\kappa_{12}\rangle\\&\cong (\mathbb Z/2)^{3}.\end{aligned}\)}
	\StemEntry{\(\displaystyle \pi_{33}(S^{4})=\langle \nu_4 \sigma' \omega_{14} \nu_{30}, \nu_4 \bar\kappa_7 \nu^2_{27}, \nu_4 \bar\sigma_7 \sigma_{26}, \nu_4 \bar\zeta_7 \sigma_{26}, \nu_4 \phi_7 \nu_{30}, \nu^3_4 \bar\kappa_{13}, \nu_4 \eta_7 \mu_{3,8}, (\Sigma \bar\alpha) \nu^2_{27}, (\Sigma \nu' )\eta_7 \mu_{3,8}, \eta_4 \varepsilon _5 \bar\kappa_{13}\rangle\cong (\mathbb Z/2)^{3} \oplus \mathbb Z/8 \oplus (\mathbb Z/2)^{6}.\)}{\(\displaystyle \begin{aligned}[t]\pi_{33}(S^{4})&=\langle \nu_4 \sigma' \omega_{14} \nu_{30}, \nu_4 \bar\kappa_7 \nu^2_{27}, \nu_4 \bar\sigma_7 \sigma_{26}, \nu_4 \bar\zeta_7 \sigma_{26}, \nu_4 \phi_7 \nu_{30}, \nu^3_4 \bar\kappa_{13}, \nu_4 \eta_7 \mu_{3,8}, (\Sigma \bar\alpha) \nu^2_{27},\\&\phantom{=\langle{}}(\Sigma \nu' )\eta_7 \mu_{3,8}, \eta_4 \varepsilon _5 \bar\kappa_{13}\rangle\cong (\mathbb Z/2)^{3} \oplus \mathbb Z/8 \oplus (\mathbb Z/2)^{6}.\end{aligned}\)}
	\StemEntry{\(\displaystyle \pi_{34}(S^{5})=\langle \phi_5 \nu^2_{28}, \nu_5 \bar\kappa_8 \nu^2_{28}, \nu_5 \bar\sigma_8 \sigma_{27}, \nu^3_5 \bar\kappa_{14}, \nu_5 \eta_8 \mu_{3,9}, \eta_5 \varepsilon _6 \bar\kappa_{14}\rangle\cong (\mathbb Z/2)^{6}.\)}{\(\displaystyle \begin{aligned}[t]\pi_{34}(S^{5})&=\langle \phi_5 \nu^2_{28}, \nu_5 \bar\kappa_8 \nu^2_{28}, \nu_5 \bar\sigma_8 \sigma_{27}, \nu^3_5 \bar\kappa_{14}, \nu_5 \eta_8 \mu_{3,9},\\&\phantom{=\langle{}}\eta_5 \varepsilon _6 \bar\kappa_{14}\rangle\cong (\mathbb Z/2)^{6}.\end{aligned}\)}
	\vspace{-0.8em}\;\\
	In the following entry, $\bar\varepsilon'\in\pi_{35}(S^{6})$ denotes Oda's element \cite[p.~130]{Oda1979}; it is distinct from the element $\bar\varepsilon'\in\pi_{20}(S^{3})$ in Subsection~\ref{app:stem-17}.
	
	\StemEntry{\(\displaystyle \pi_{35}(S^{6})=\langle \delta', \bar\varepsilon ', \phi_6 \nu^2_{29}, \eta_6 \varepsilon _7 \bar\kappa_{15}\rangle\cong \mathbb Z/4 \oplus \mathbb Z/4 \oplus (\mathbb Z/2)^{2}.\)}{\(\displaystyle \begin{aligned}[t]\pi_{35}(S^{6})&=\langle \delta', \bar\varepsilon ', \phi_6 \nu^2_{29}, \eta_6 \varepsilon _7 \bar\kappa_{15}\rangle\\&\cong \mathbb Z/4 \oplus \mathbb Z/4 \oplus (\mathbb Z/2)^{2}.\end{aligned}\)}
	\StemEntry{\(\displaystyle \pi_{36}(S^{7})=\langle \delta'', \sigma' \varepsilon _{14} \kappa_{22}, \sigma' \omega_{14} \nu^2_{30}, \phi_7 \nu^2_{30}, \eta_7 \varepsilon _8 \bar\kappa_{16}\rangle\cong \mathbb Z/8 \oplus (\mathbb Z/2)^{4}.\)}{\(\displaystyle \begin{aligned}[t]\pi_{36}(S^{7})&=\langle \delta'', \sigma' \varepsilon _{14} \kappa_{22}, \sigma' \omega_{14} \nu^2_{30}, \phi_7 \nu^2_{30},\\&\phantom{=\langle{}}\eta_7 \varepsilon _8 \bar\kappa_{16}\rangle\cong \mathbb Z/8 \oplus (\mathbb Z/2)^{4}.\end{aligned}\)}
	\StemEntry{\(\displaystyle \pi_{37}(S^{8})=\langle \sigma_8 \sigma^{*\prime}, \sigma_8 \omega_{15} \nu^2_{31}, \sigma_8 \varepsilon _{15} \kappa_{23}, \sigma_8 \nu_{15} \bar\sigma_{18}, \Sigma \delta'', (\Sigma \sigma') \varepsilon _{15} \kappa_{23}, (\Sigma \sigma' )\omega_{15} \nu^2_{31}, \phi_8 \nu^2_{31}, \eta_8 \varepsilon _9 \bar\kappa_{17}\rangle\cong \mathbb Z/16 \oplus \mathbb Z/2 \oplus \mathbb Z/2 \oplus \mathbb Z/2 \oplus \mathbb Z/8 \oplus (\mathbb Z/2)^{4}.\)}{\(\displaystyle \begin{aligned}[t]\pi_{37}(S^{8})&=\langle \sigma_8 \sigma^{*\prime}, \sigma_8 \omega_{15} \nu^2_{31}, \sigma_8 \varepsilon _{15} \kappa_{23}, \sigma_8 \nu_{15} \bar\sigma_{18}, \Sigma \delta'', (\Sigma \sigma') \varepsilon _{15} \kappa_{23}, (\Sigma \sigma' )\omega_{15} \nu^2_{31},\\&\phantom{=\langle{}}\phi_8 \nu^2_{31}, \eta_8 \varepsilon _9 \bar\kappa_{17}\rangle\cong \mathbb Z/16 \oplus \mathbb Z/2 \oplus \mathbb Z/2 \oplus \mathbb Z/2 \oplus \mathbb Z/8 \oplus (\mathbb Z/2)^{4}.\end{aligned}\)}
	\StemEntry{\(\displaystyle \pi_{38}(S^{9})=\langle \sigma_9 \sigma^*_{16}, \sigma_9 \omega_{16} \nu^2_{32}, \sigma_9 \varepsilon _{16} \kappa_{24}, \sigma_9 \nu_{16} \bar\sigma_{19}, \eta_9 \varepsilon _{10} \bar\kappa_{18}\rangle\cong \mathbb Z/16 \oplus (\mathbb Z/2)^{4}.\)}{\(\displaystyle \begin{aligned}[t]\pi_{38}(S^{9})&=\langle \sigma_9 \sigma^*_{16}, \sigma_9 \omega_{16} \nu^2_{32}, \sigma_9 \varepsilon _{16} \kappa_{24}, \sigma_9 \nu_{16} \bar\sigma_{19},\\&\phantom{=\langle{}}\eta_9 \varepsilon _{10} \bar\kappa_{18}\rangle\cong \mathbb Z/16 \oplus (\mathbb Z/2)^{4}.\end{aligned}\)}
	\StemEntry{\(\displaystyle \pi_{39}(S^{10})=\langle P(\bar\kappa_{21}), P(\Sigma A_2), \sigma_{10} \sigma^*_{17}, \sigma_{10} \nu_{17} \bar\sigma_{20}\rangle\cong \mathbb Z/8 \oplus \mathbb Z/2 \oplus \mathbb Z/16 \oplus \mathbb Z/2.\)}{\(\displaystyle \begin{aligned}[t]\pi_{39}(S^{10})&=\langle P(\bar\kappa_{21}), P(\Sigma A_2), \sigma_{10} \sigma^*_{17},\\&\phantom{=\langle{}}\sigma_{10} \nu_{17} \bar\sigma_{20}\rangle\cong \mathbb Z/8 \oplus \mathbb Z/2 \oplus \mathbb Z/16 \oplus \mathbb Z/2.\end{aligned}\)}
	\StemEntry{\(\displaystyle \pi_{40}(S^{11})=\langle C_1 \bar\mu_{23}, \sigma_{11} \sigma^*_{18}, \sigma_{11} \nu_{18} \bar\sigma_{21}\rangle\cong \mathbb Z/2 \oplus \mathbb Z/16 \oplus \mathbb Z/2.\)}{\(\displaystyle \begin{aligned}[t]\pi_{40}(S^{11})&=\langle C_1 \bar\mu_{23}, \sigma_{11} \sigma^*_{18}, \sigma_{11} \nu_{18} \bar\sigma_{21}\rangle\\&\cong \mathbb Z/2 \oplus \mathbb Z/16 \oplus \mathbb Z/2.\end{aligned}\)}
	\StemEntry{\(\displaystyle \pi_{41}(S^{12})=\langle P(\nu^*_{25}) + 2\sigma_{12} \sigma^*_{19}, A_1 \bar\mu_{24}, (\Sigma C_1) \bar\mu_{24}, \sigma_{12} \sigma^*_{19}\rangle\cong \mathbb Z/4 \oplus \mathbb Z/2 \oplus \mathbb Z/2 \oplus \mathbb Z/16.\)}{\(\displaystyle \begin{aligned}[t]\pi_{41}(S^{12})&=\langle P(\nu^*_{25}) + 2\sigma_{12} \sigma^*_{19}, A_1 \bar\mu_{24}, (\Sigma C_1) \bar\mu_{24},\\&\phantom{=\langle{}}\sigma_{12} \sigma^*_{19}\rangle\cong \mathbb Z/4 \oplus \mathbb Z/2 \oplus \mathbb Z/2 \oplus \mathbb Z/16.\end{aligned}\)}
	\StemEntry{\(\displaystyle \pi_{42}(S^{13})=\langle (\Sigma A_1 )\bar\mu_{25}, \sigma_{13} \sigma^*_{20}\rangle\cong \mathbb Z/2 \oplus \mathbb Z/8.\)}{\(\displaystyle \begin{aligned}[t]\pi_{42}(S^{13})&=\langle (\Sigma A_1 )\bar\mu_{25}, \sigma_{13} \sigma^*_{20}\rangle\\&\cong \mathbb Z/2 \oplus \mathbb Z/8.\end{aligned}\)}
	\StemEntry{\(\displaystyle \pi_{43}(S^{14})=\langle \sigma_{14} \sigma^*_{21}\rangle\cong \mathbb Z/4.\)}{\(\displaystyle \begin{aligned}[t]\pi_{43}(S^{14})&=\langle \sigma_{14} \sigma^*_{21}\rangle\\&\cong \mathbb Z/4.\end{aligned}\)}
	\StemEntry{\(\displaystyle \pi_{44}(S^{15})=\langle L_1, \sigma_{15} \sigma^*_{22}\rangle\cong (\mathbb Z/2)^{2}.\)}{\(\displaystyle \begin{aligned}[t]\pi_{44}(S^{15})&=\langle L_1, \sigma_{15} \sigma^*_{22}\rangle\\&\cong (\mathbb Z/2)^{2}.\end{aligned}\)}
	\StemEntry{\(\displaystyle \pi_{45}(S^{16})=\langle \sigma^*_{16} \sigma_{38}, P_1, \Sigma L_1, \sigma_{16} \sigma^*_{23}\rangle\cong (\mathbb Z/2)^{4}.\)}{\(\displaystyle \begin{aligned}[t]\pi_{45}(S^{16})&=\langle \sigma^*_{16} \sigma_{38}, P_1, \Sigma L_1, \sigma_{16} \sigma^*_{23}\rangle\\&\cong (\mathbb Z/2)^{4}.\end{aligned}\)}
	\StemEntry{\(\displaystyle \pi_{n+29}(S^{n})=\langle \sigma^*_n \sigma_{n+22}, \Sigma ^{n-16}P_1\rangle\cong (\mathbb Z/2)^{2}\quad(n = 17 \text{ or } 18).\)}{\(\displaystyle \begin{aligned}[t]\pi_{n+29}(S^{n})&=\langle \sigma^*_n \sigma_{n+22}, \Sigma ^{n-16}P_1\rangle\\&\cong (\mathbb Z/2)^{2}\quad(n = 17 \text{ or } 18).\end{aligned}\)}
	\StemEntry{\(\displaystyle \pi_{48}(S^{19})=\langle C_2 \mu_{39}, \sigma^*_{19} \sigma_{41}, \Sigma ^3P_1\rangle\cong (\mathbb Z/2)^{3}.\)}{\(\displaystyle \begin{aligned}[t]\pi_{48}(S^{19})&=\langle C_2 \mu_{39}, \sigma^*_{19} \sigma_{41}, \Sigma ^3P_1\rangle\\&\cong (\mathbb Z/2)^{3}.\end{aligned}\)}
	\StemEntry{\(\displaystyle \pi_{49}(S^{20})=\langle A_2 \mu_{40}, (\Sigma C_2) \mu_{40}, \sigma^*_{20} \sigma_{42}, \Sigma ^4P_1\rangle\cong (\mathbb Z/2)^{4}.\)}{\(\displaystyle \begin{aligned}[t]\pi_{49}(S^{20})&=\langle A_2 \mu_{40}, (\Sigma C_2) \mu_{40}, \sigma^*_{20} \sigma_{42},\\&\phantom{=\langle{}}\Sigma ^4P_1\rangle\cong (\mathbb Z/2)^{4}.\end{aligned}\)}
	\StemEntry{\(\displaystyle \pi_{50}(S^{21})=\langle M'_2 \nu_{47}, (\Sigma A_2 )\mu_{41}, \sigma^*_{21} \sigma_{43}\rangle\cong \mathbb Z/4 \oplus (\mathbb Z/2)^{2}.\)}{\(\displaystyle \begin{aligned}[t]\pi_{50}(S^{21})&=\langle M'_2 \nu_{47}, (\Sigma A_2 )\mu_{41},\\&\phantom{=\langle{}}\sigma^*_{21} \sigma_{43}\rangle\cong \mathbb Z/4 \oplus (\mathbb Z/2)^{2}.\end{aligned}\)}
	\StemEntry{\(\displaystyle \pi_{51}(S^{22})=\langle (\Sigma M'_2) \nu_{48}, \sigma^*_{22} \sigma_{44}\rangle\cong \mathbb Z/4 \oplus \mathbb Z/2.\)}{\(\displaystyle \begin{aligned}[t]\pi_{51}(S^{22})&=\langle (\Sigma M'_2) \nu_{48}, \sigma^*_{22} \sigma_{44}\rangle\\&\cong \mathbb Z/4 \oplus \mathbb Z/2.\end{aligned}\)}
	\StemEntry{\(\displaystyle \pi_{52}(S^{23})=\langle (\Sigma ^2M'_2) \nu_{49}, \sigma^*_{23} \sigma_{45}\rangle\cong (\mathbb Z/2)^{2}.\)}{\(\displaystyle \begin{aligned}[t]\pi_{52}(S^{23})&=\langle (\Sigma ^2M'_2) \nu_{49},\\&\phantom{=\langle{}}\sigma^*_{23} \sigma_{45}\rangle\cong (\mathbb Z/2)^{2}.\end{aligned}\)}
	\StemEntry{\(\displaystyle \pi_{53}(S^{24})=\langle M_2 \nu_{50}, (\Sigma ^3M'_2) \nu_{50}\rangle\cong (\mathbb Z/2)^{2}.\)}{\(\displaystyle \begin{aligned}[t]\pi_{53}(S^{24})&=\langle M_2 \nu_{50}, (\Sigma ^3M'_2) \nu_{50}\rangle\\&\cong (\mathbb Z/2)^{2}.\end{aligned}\)}
	\StemEntry{\(\displaystyle \pi_{n+29}(S^{n})=\langle (\Sigma ^{n-24}M_2) \nu_{n+26}\rangle\cong \mathbb Z/2\quad(n = 25 \text{ or } 26).\)}{\(\displaystyle \begin{aligned}[t]\pi_{n+29}(S^{n})&=\langle (\Sigma ^{n-24}M_2) \nu_{n+26}\rangle\\&\cong \mathbb Z/2\quad(n = 25 \text{ or } 26).\end{aligned}\)}
	\StemEntry{\(\displaystyle \pi_{56}(S^{27})=\langle C_3 \eta_{55}, (\Sigma ^3M_2) \nu_{53}\rangle\cong (\mathbb Z/2)^{2}.\)}{\(\displaystyle \begin{aligned}[t]\pi_{56}(S^{27})&=\langle C_3 \eta_{55}, (\Sigma ^3M_2) \nu_{53}\rangle\\&\cong (\mathbb Z/2)^{2}.\end{aligned}\)}
	\StemEntry{\(\displaystyle \pi_{57}(S^{28})=\langle A_3 \eta_{56}, (\Sigma C_3) \eta_{56}\rangle\cong (\mathbb Z/2)^{2}.\)}{\(\displaystyle \begin{aligned}[t]\pi_{57}(S^{28})&=\langle A_3 \eta_{56}, (\Sigma C_3) \eta_{56}\rangle\\&\cong (\mathbb Z/2)^{2}.\end{aligned}\)}
	\StemEntry{\(\displaystyle \pi_{58}(S^{29})=\langle (\Sigma A_3 )\eta_{57}\rangle\cong \mathbb Z/2.\)}{\(\displaystyle \begin{aligned}[t]\pi_{58}(S^{29})&=\langle (\Sigma A_3 )\eta_{57}\rangle\\&\cong \mathbb Z/2.\end{aligned}\)}
	\StemEntry{\(\displaystyle \pi_{59}(S^{30})=\langle P(\iota_{61})\rangle\cong \mathbb Z_{(2)}.\)}{\(\displaystyle \begin{aligned}[t]\pi_{59}(S^{30})&=\langle P(\iota_{61})\rangle\\&\cong \mathbb Z_{(2)}.\end{aligned}\)}
	\StemEntry{\(\displaystyle \pi_{n+29}(S^{n})=0\quad(n\ge 31).\)}{\(\displaystyle \begin{aligned}[t]\pi_{n+29}(S^{n})&=0\quad(n\ge 31).\end{aligned}\)}
	\StemEntry{\(\displaystyle \pi^S_{29}(S^0)=0.\)}{\(\displaystyle \pi^S_{29}(S^0)=0.\)}
	\subsection{Stem 30}\label{app:stem-30}
	\StemEntry{\(\displaystyle \pi_{32}(S^{2})=\langle \eta_2 \bar\alpha \nu^2_{26}, \eta_2 \nu' \eta_6 \mu_{3,7}, \eta^2_2 \varepsilon _4 \bar\kappa_{12}\rangle\cong (\mathbb Z/2)^{3}.\)}{\(\displaystyle \begin{aligned}[t]\pi_{32}(S^{2})&=\langle \eta_2 \bar\alpha \nu^2_{26}, \eta_2 \nu' \eta_6 \mu_{3,7},\\&\phantom{=\langle{}}\eta^2_2 \varepsilon _4 \bar\kappa_{12}\rangle\cong (\mathbb Z/2)^{3}.\end{aligned}\)}
	\StemEntry{\(\displaystyle \pi_{33}(S^{3})=\langle \varepsilon ' \bar\kappa_{13}, \varepsilon _3 \nu_{11} \bar\sigma_{14}\rangle\cong \mathbb Z/4 \oplus \mathbb Z/2.\)}{\(\displaystyle \begin{aligned}[t]\pi_{33}(S^{3})&=\langle \varepsilon ' \bar\kappa_{13}, \varepsilon _3 \nu_{11} \bar\sigma_{14}\rangle\\&\cong \mathbb Z/4 \oplus \mathbb Z/2.\end{aligned}\)}
	\StemEntry{\(\displaystyle \pi_{34}(S^{4})=\langle \nu_4 \sigma' \bar\kappa_{14}, \nu_4 \zeta_{3,7}, \nu_4 \bar\nu_7 \bar\sigma_{15}, (\Sigma \varepsilon ') \bar\kappa_{14}, \varepsilon _4 \nu_{12} \bar\sigma_{15}\rangle\cong \mathbb Z/8 \oplus \mathbb Z/8 \oplus \mathbb Z/2 \oplus \mathbb Z/4 \oplus \mathbb Z/2.\)}{\(\displaystyle \begin{aligned}[t]\pi_{34}(S^{4})&=\langle \nu_4 \sigma' \bar\kappa_{14}, \nu_4 \zeta_{3,7}, \nu_4 \bar\nu_7 \bar\sigma_{15}, (\Sigma \varepsilon ') \bar\kappa_{14},\\&\phantom{=\langle{}}\varepsilon _4 \nu_{12} \bar\sigma_{15}\rangle\cong \mathbb Z/8 \oplus \mathbb Z/8 \oplus \mathbb Z/2 \oplus \mathbb Z/4 \oplus \mathbb Z/2.\end{aligned}\)}
	\StemEntry{\(\displaystyle \pi_{35}(S^{5})=\langle \nu_5 \sigma_8 \bar\kappa_{15}, \phi_5 \sigma_{28}, \nu_5 \zeta_{3,8}, \nu_5 \bar\nu_8 \bar\sigma_{16}\rangle\cong \mathbb Z/8 \oplus (\mathbb Z/2)^{3}.\)}{\(\displaystyle \begin{aligned}[t]\pi_{35}(S^{5})&=\langle \nu_5 \sigma_8 \bar\kappa_{15}, \phi_5 \sigma_{28}, \nu_5 \zeta_{3,8},\\&\phantom{=\langle{}}\nu_5 \bar\nu_8 \bar\sigma_{16}\rangle\cong \mathbb Z/8 \oplus (\mathbb Z/2)^{3}.\end{aligned}\)}
	\StemEntry{\(\displaystyle \pi_{36}(S^{6})=\langle P(\xi_{13} \sigma_{31}), \sigma'' \bar\rho_{13}, \nu_6 \sigma_9 \bar\kappa_{16}, \phi_6 \sigma_{29}\rangle\cong \mathbb Z/4 \oplus \mathbb Z/4 \oplus \mathbb Z/8 \oplus \mathbb Z/2.\)}{\(\displaystyle \begin{aligned}[t]\pi_{36}(S^{6})&=\langle P(\xi_{13} \sigma_{31}), \sigma'' \bar\rho_{13}, \nu_6 \sigma_9 \bar\kappa_{16},\\&\phantom{=\langle{}}\phi_6 \sigma_{29}\rangle\cong \mathbb Z/4 \oplus \mathbb Z/4 \oplus \mathbb Z/8 \oplus \mathbb Z/2.\end{aligned}\)}
	\StemEntry{\(\displaystyle \pi_{37}(S^{7})=\langle \sigma' \bar\rho_{14}, \sigma' \phi_{14}, \sigma' \psi_{14}, \nu_7 \sigma_{10} \bar\kappa_{17}, \phi_7 \sigma_{30}\rangle\cong \mathbb Z/8 \oplus \mathbb Z/2 \oplus \mathbb Z/2 \oplus \mathbb Z/8 \oplus \mathbb Z/2.\)}{\(\displaystyle \begin{aligned}[t]\pi_{37}(S^{7})&=\langle \sigma' \bar\rho_{14}, \sigma' \phi_{14}, \sigma' \psi_{14}, \nu_7 \sigma_{10} \bar\kappa_{17},\\&\phantom{=\langle{}}\phi_7 \sigma_{30}\rangle\cong \mathbb Z/8 \oplus \mathbb Z/2 \oplus \mathbb Z/2 \oplus \mathbb Z/8 \oplus \mathbb Z/2.\end{aligned}\)}
	\StemEntry{\(\displaystyle \pi_{38}(S^{8})=\langle \sigma_8 \bar\rho_{15}, \sigma_8 \nu_{15} \bar\kappa_{18}, \sigma_8 \phi_{15}, \sigma_8 \psi_{15}, \sigma_8 D^{\mathrm{I}}_1, \sigma_8 D_1 \sigma_{31}, (\Sigma \sigma') \bar\rho_{15}, (\Sigma \sigma') \phi_{15}, (\Sigma \sigma') \psi_{15}, \nu_8 \sigma_{11} \bar\kappa_{18}, \phi_8 \sigma_{31}\rangle\cong \mathbb Z/16 \oplus \mathbb Z/8 \oplus \mathbb Z/2 \oplus \mathbb Z/2 \oplus \mathbb Z/2 \oplus \mathbb Z/2 \oplus \mathbb Z/8 \oplus \mathbb Z/2 \oplus \mathbb Z/2 \oplus \mathbb Z/8 \oplus \mathbb Z/2.\)}{\(\displaystyle \begin{aligned}[t]\pi_{38}(S^{8})&=\langle \sigma_8 \bar\rho_{15}, \sigma_8 \nu_{15} \bar\kappa_{18}, \sigma_8 \phi_{15}, \sigma_8 \psi_{15}, \sigma_8 D^{\mathrm{I}}_1, \sigma_8 D_1 \sigma_{31}, (\Sigma \sigma') \bar\rho_{15}, (\Sigma \sigma') \phi_{15}, (\Sigma \sigma') \psi_{15},\\&\phantom{=\langle{}}\nu_8 \sigma_{11} \bar\kappa_{18}, \phi_8 \sigma_{31}\rangle\cong \mathbb Z/16 \oplus \mathbb Z/8 \oplus \mathbb Z/2 \oplus \mathbb Z/2 \oplus \mathbb Z/2 \oplus \mathbb Z/2 \oplus \mathbb Z/8 \oplus \mathbb Z/2 \oplus \mathbb Z/2 \oplus \mathbb Z/8 \oplus \mathbb Z/2.\end{aligned}\)}
	\StemEntry{\(\displaystyle \pi_{39}(S^{9})=\langle \sigma_9 \bar\rho_{16}, \sigma_9 \nu_{16} \bar\kappa_{19}, \sigma_9 \phi_{16}, \sigma_9 \psi_{16}, \phi_9 \sigma_{32}\rangle\cong \mathbb Z/16 \oplus \mathbb Z/8 \oplus (\mathbb Z/2)^{3}.\)}{\(\displaystyle \begin{aligned}[t]\pi_{39}(S^{9})&=\langle \sigma_9 \bar\rho_{16}, \sigma_9 \nu_{16} \bar\kappa_{19}, \sigma_9 \phi_{16}, \sigma_9 \psi_{16},\\&\phantom{=\langle{}}\phi_9 \sigma_{32}\rangle\cong \mathbb Z/16 \oplus \mathbb Z/8 \oplus (\mathbb Z/2)^{3}.\end{aligned}\)}
	\StemEntry{\(\displaystyle \pi_{40}(S^{10})=\langle P((\Sigma A_2 )\eta_{41}), \psi_{10} \sigma_{33}, \sigma_{10} \bar\rho_{17}, \sigma_{10} \nu_{17} \bar\kappa_{20}, \sigma_{10} \psi_{17}\rangle\cong \mathbb Z/2 \oplus \mathbb Z/4 \oplus \mathbb Z/16 \oplus \mathbb Z/4 \oplus \mathbb Z/2.\)}{\(\displaystyle \begin{aligned}[t]\pi_{40}(S^{10})&=\langle P((\Sigma A_2 )\eta_{41}), \psi_{10} \sigma_{33}, \sigma_{10} \bar\rho_{17}, \sigma_{10} \nu_{17} \bar\kappa_{20},\\&\phantom{=\langle{}}\sigma_{10} \psi_{17}\rangle\cong \mathbb Z/2 \oplus \mathbb Z/4 \oplus \mathbb Z/16 \oplus \mathbb Z/4 \oplus \mathbb Z/2.\end{aligned}\)}
	\StemEntry{\(\displaystyle \pi_{41}(S^{11})=\langle \psi_{11} \sigma_{34}, \sigma_{11} \bar\rho_{18}, \sigma_{11} \nu_{18} \bar\kappa_{21}, \sigma_{11} \psi_{18}\rangle\cong \mathbb Z/2 \oplus \mathbb Z/16 \oplus (\mathbb Z/2)^{2}.\)}{\(\displaystyle \begin{aligned}[t]\pi_{41}(S^{11})&=\langle \psi_{11} \sigma_{34}, \sigma_{11} \bar\rho_{18}, \sigma_{11} \nu_{18} \bar\kappa_{21},\\&\phantom{=\langle{}}\sigma_{11} \psi_{18}\rangle\cong \mathbb Z/2 \oplus \mathbb Z/16 \oplus (\mathbb Z/2)^{2}.\end{aligned}\)}\vspace{-0.8em}\;\\
	The sign is chosen so that $\sigma_{12}\bar\rho_{19}\pm2\theta^{\mathrm{VII}}$ has order $4.$ 
	\StemEntry{\(\displaystyle \pi_{42}(S^{12})=\langle \theta^{\mathrm{VII}}, \sigma_{12} \bar\rho_{19} \pm 2\theta^{\mathrm{VII}}, \psi_{12} \sigma_{35}, \sigma_{12} \psi_{19}\rangle\cong \mathbb Z/32 \oplus \mathbb Z/4 \oplus (\mathbb Z/2)^{2}.\)}{\(\displaystyle \begin{aligned}[t]\pi_{42}(S^{12})&=\langle \theta^{\mathrm{VII}}, \sigma_{12} \bar\rho_{19} \pm 2\theta^{\mathrm{VII}}, \psi_{12} \sigma_{35},\\&\phantom{=\langle{}}\sigma_{12} \psi_{19}\rangle\cong \mathbb Z/32 \oplus \mathbb Z/4 \oplus (\mathbb Z/2)^{2}.\end{aligned}\)}
	\StemEntry{\(\displaystyle \pi_{43}(S^{13})=\langle \rho^2_{13}, \psi_{13} \sigma_{36}\rangle\cong \mathbb Z/32 \oplus \mathbb Z/2.\)}{\(\displaystyle \begin{aligned}[t]\pi_{43}(S^{13})&=\langle \rho^2_{13}, \psi_{13} \sigma_{36}\rangle\\&\cong \mathbb Z/32 \oplus \mathbb Z/2.\end{aligned}\)}
	\StemEntry{\(\displaystyle \pi_{44}(S^{14})=\langle \theta^{\mathrm{VI}}, \omega_{14} \kappa_{30}, \psi_{14} \sigma_{37}\rangle\cong \mathbb Z/64 \oplus (\mathbb Z/2)^{2}.\)}{\(\displaystyle \begin{aligned}[t]\pi_{44}(S^{14})&=\langle \theta^{\mathrm{VI}}, \omega_{14} \kappa_{30}, \psi_{14} \sigma_{37}\rangle\\&\cong \mathbb Z/64 \oplus (\mathbb Z/2)^{2}.\end{aligned}\)}
	\StemEntry{\(\displaystyle \pi_{45}(S^{15})=\langle \theta^{\mathrm{V}}, \omega_{15} \kappa_{31}, \psi_{15} \sigma_{38}\rangle\cong \mathbb Z/64 \oplus (\mathbb Z/2)^{2}.\)}{\(\displaystyle \begin{aligned}[t]\pi_{45}(S^{15})&=\langle \theta^{\mathrm{V}}, \omega_{15} \kappa_{31}, \psi_{15} \sigma_{38}\rangle\\&\cong \mathbb Z/64 \oplus (\mathbb Z/2)^{2}.\end{aligned}\)}\vspace{-0.8em}\;\\
	The sign is chosen so that $(\Sigma\theta^{\mathrm{V}})\pm2\theta^{\mathrm{IV}}$ has order $16.$ 
	\StemEntry{\(\displaystyle \pi_{46}(S^{16})=\langle \theta^{\mathrm{IV}}, (\Sigma \theta^{\mathrm{V}}) \pm 2\theta^{\mathrm{IV}}, B_1 \kappa_{32}, \omega_{16} \kappa_{32}, \psi_{16} \sigma_{39}\rangle\cong \mathbb Z/128 \oplus \mathbb Z/16 \oplus (\mathbb Z/2)^{3}.\)}{\(\displaystyle \begin{aligned}[t]\pi_{46}(S^{16})&=\langle \theta^{\mathrm{IV}}, (\Sigma \theta^{\mathrm{V}}) \pm 2\theta^{\mathrm{IV}}, B_1 \kappa_{32}, \omega_{16} \kappa_{32},\\&\phantom{=\langle{}}\psi_{16} \sigma_{39}\rangle\cong \mathbb Z/128 \oplus \mathbb Z/16 \oplus (\mathbb Z/2)^{3}.\end{aligned}\)}
	\StemEntry{\(\displaystyle \pi_{47}(S^{17})=\langle \Sigma \theta^{\mathrm{IV}}, (\Sigma B_1) \kappa_{33}, \psi_{17} \sigma_{40}\rangle\cong \mathbb Z/64 \oplus (\mathbb Z/2)^{2}.\)}{\(\displaystyle \begin{aligned}[t]\pi_{47}(S^{17})&=\langle \Sigma \theta^{\mathrm{IV}}, (\Sigma B_1) \kappa_{33}, \psi_{17} \sigma_{40}\rangle\\&\cong \mathbb Z/64 \oplus (\mathbb Z/2)^{2}.\end{aligned}\)}
	\StemEntry{\(\displaystyle \pi_{n+30}(S^{n})=\langle \Sigma ^{n-16} \theta^{\mathrm{IV}}\rangle\cong \mathbb Z/64\quad(n = 18 \text{ or } 19).\)}{\(\displaystyle \begin{aligned}[t]\pi_{n+30}(S^{n})&=\langle \Sigma ^{n-16} \theta^{\mathrm{IV}}\rangle\\&\cong \mathbb Z/64\quad(n = 18 \text{ or } 19).\end{aligned}\)}\vspace{-0.8em}\;\\
	Choose an odd integer $x$ such that $4(\Sigma^4\theta^{\mathrm{IV}}+2x\theta''')=0.$ 
	\StemEntry{\(\displaystyle \pi_{50}(S^{20})=\langle \theta''', \Sigma ^4\theta^{\mathrm{IV}} + 2x\theta'''\rangle\cong \mathbb Z/128 \oplus \mathbb Z/4.\)}{\(\displaystyle \begin{aligned}[t]\pi_{50}(S^{20})&=\langle \theta''', \Sigma ^4\theta^{\mathrm{IV}} + 2x\theta'''\rangle\\&\cong \mathbb Z/128 \oplus \mathbb Z/4.\end{aligned}\)}
	\StemEntry{\(\displaystyle \pi_{51}(S^{21})=\langle \Sigma \theta'''\rangle\cong \mathbb Z/64.\)}{\(\displaystyle \begin{aligned}[t]\pi_{51}(S^{21})&=\langle \Sigma \theta'''\rangle\\&\cong \mathbb Z/64.\end{aligned}\)}
	\StemEntry{\(\displaystyle \pi_{52}(S^{22})=\langle \theta'', V_2 \nu_{49}\rangle\cong \mathbb Z/64 \oplus \mathbb Z/2.\)}{\(\displaystyle \begin{aligned}[t]\pi_{52}(S^{22})&=\langle \theta'', V_2 \nu_{49}\rangle\\&\cong \mathbb Z/64 \oplus \mathbb Z/2.\end{aligned}\)}
	\StemEntry{\(\displaystyle \pi_{53}(S^{23})=\langle \theta'_{23}, (\Sigma V_2) \nu_{50}\rangle\cong \mathbb Z/64 \oplus \mathbb Z/2.\)}{\(\displaystyle \begin{aligned}[t]\pi_{53}(S^{23})&=\langle \theta'_{23}, (\Sigma V_2) \nu_{50}\rangle\\&\cong \mathbb Z/64 \oplus \mathbb Z/2.\end{aligned}\)}
	\StemEntry{\(\displaystyle \pi_{54}(S^{24})=\langle P(\sigma_{49}) - 4\theta'_{24}, \theta'_{24}, (\Sigma ^2V_2) \nu_{51}\rangle\cong \mathbb Z/8 \oplus \mathbb Z/64 \oplus \mathbb Z/2.\)}{\(\displaystyle \begin{aligned}[t]\pi_{54}(S^{24})&=\langle P(\sigma_{49}) - 4\theta'_{24}, \theta'_{24},\\&\phantom{=\langle{}}(\Sigma ^2V_2) \nu_{51}\rangle\cong \mathbb Z/8 \oplus \mathbb Z/64 \oplus \mathbb Z/2.\end{aligned}\)}
	\StemEntry{\(\displaystyle \pi_{55}(S^{25})=\langle \theta'_{25}, (\Sigma ^3V_2 )\nu_{52}\rangle\cong \mathbb Z/32 \oplus \mathbb Z/2.\)}{\(\displaystyle \begin{aligned}[t]\pi_{55}(S^{25})&=\langle \theta'_{25}, (\Sigma ^3V_2 )\nu_{52}\rangle\\&\cong \mathbb Z/32 \oplus \mathbb Z/2.\end{aligned}\)}
	\StemEntry{\(\displaystyle \pi_{n+30}(S^{n})=\langle \theta'_n\rangle\cong \mathbb Z/32\quad(n = 26 \text{ or } 27).\)}{\(\displaystyle \begin{aligned}[t]\pi_{n+30}(S^{n})&=\langle \theta'_n\rangle\\&\cong \mathbb Z/32\quad(n = 26 \text{ or } 27).\end{aligned}\)}
	\StemEntry{\(\displaystyle \pi_{58}(S^{28})=\langle P(\nu_{57}) -4\theta'_{28}, \theta'_{28}\rangle\cong \mathbb Z/4 \oplus \mathbb Z/32.\)}{\(\displaystyle \begin{aligned}[t]\pi_{58}(S^{28})&=\langle P(\nu_{57}) -4\theta'_{28},\\&\phantom{=\langle{}}\theta'_{28}\rangle\cong \mathbb Z/4 \oplus \mathbb Z/32.\end{aligned}\)}
	\StemEntry{\(\displaystyle \pi_{59}(S^{29})=\langle \theta'_{29}\rangle\cong \mathbb Z/16.\)}{\(\displaystyle \begin{aligned}[t]\pi_{59}(S^{29})&=\langle \theta'_{29}\rangle\\&\cong \mathbb Z/16.\end{aligned}\)}
	\StemEntry{\(\displaystyle \pi_{60}(S^{30})=\langle \theta'_{30}\rangle\cong \mathbb Z/8.\)}{\(\displaystyle \begin{aligned}[t]\pi_{60}(S^{30})&=\langle \theta'_{30}\rangle\\&\cong \mathbb Z/8.\end{aligned}\)}
	\StemEntry{\(\displaystyle \pi_{61}(S^{31})=\langle \theta'_{31}\rangle\cong \mathbb Z/4.\)}{\(\displaystyle \begin{aligned}[t]\pi_{61}(S^{31})&=\langle \theta'_{31}\rangle\\&\cong \mathbb Z/4.\end{aligned}\)}
	\StemEntry{\(\displaystyle \pi_{n+30}(S^{n})=\langle \theta'_n\rangle\cong \mathbb Z/2\quad(n\ge32).\)}{\(\displaystyle \begin{aligned}[t]\pi_{n+30}(S^{n})&=\langle \theta'_n\rangle\\&\cong \mathbb Z/2\quad(n\ge32).\end{aligned}\)}
	\StemEntry{\(\displaystyle \pi^S_{30}(S^0)=\langle \theta'\rangle\cong \mathbb Z/2.\)}{\(\displaystyle \begin{aligned}[t]\pi^S_{30}(S^0)&=\langle \theta'\rangle\\&\cong \mathbb Z/2.\end{aligned}\)}
	\subsection{Stem 31}\label{app:stem-31}
	\StemEntry{\(\displaystyle \pi_{33}(S^{2})=\langle \eta_2 \varepsilon ' \bar\kappa_{13}, \eta_2 \varepsilon _3 \nu_{11} \bar\sigma_{14}\rangle\cong \mathbb Z/4 \oplus \mathbb Z/2.\)}{\(\displaystyle \begin{aligned}[t]\pi_{33}(S^{2})&=\langle \eta_2 \varepsilon ' \bar\kappa_{13}, \eta_2 \varepsilon _3 \nu_{11} \bar\sigma_{14}\rangle\\&\cong \mathbb Z/4 \oplus \mathbb Z/2.\end{aligned}\)}
	\StemEntry{\(\displaystyle \pi_{34}(S^{3})=\langle \delta_3 \sigma_{27}, \varepsilon _3 \nu_{11} \bar\kappa_{14}, \nu' \varepsilon _6 \bar\kappa_{14}, \phi' \nu^2_{28}\rangle\cong (\mathbb Z/2)^{4}.\)}{\(\displaystyle \begin{aligned}[t]\pi_{34}(S^{3})&=\langle \delta_3 \sigma_{27}, \varepsilon _3 \nu_{11} \bar\kappa_{14}, \nu' \varepsilon _6 \bar\kappa_{14},\\&\phantom{=\langle{}}\phi' \nu^2_{28}\rangle\cong (\mathbb Z/2)^{4}.\end{aligned}\)}
	\StemEntry{\(\displaystyle \pi_{35}(S^{4})=\langle \nu_4 \phi''' \nu_{32}, \nu_4 \sigma' \eta_{14} \bar\kappa_{15}, \nu_4 \bar\nu_7 \bar\kappa_{15}, \nu_4 \varepsilon _7 \bar\kappa_{15}, \delta_4 \sigma_{28}, \varepsilon _4\nu_{12} \bar\kappa_{15}, (\Sigma \nu') \varepsilon _7 \bar\kappa_{15}, (\Sigma \phi') \nu^2_{29}\rangle\cong (\mathbb Z/2)^{8}.\)}{\(\displaystyle \begin{aligned}[t]\pi_{35}(S^{4})&=\langle \nu_4 \phi''' \nu_{32}, \nu_4 \sigma' \eta_{14} \bar\kappa_{15}, \nu_4 \bar\nu_7 \bar\kappa_{15}, \nu_4 \varepsilon _7 \bar\kappa_{15}, \delta_4 \sigma_{28}, \varepsilon _4\nu_{12} \bar\kappa_{15},\\&\phantom{=\langle{}}(\Sigma \nu') \varepsilon _7 \bar\kappa_{15}, (\Sigma \phi') \nu^2_{29}\rangle\cong (\mathbb Z/2)^{8}.\end{aligned}\)}
	\StemEntry{\(\displaystyle \pi_{36}(S^{5})=\langle \nu_5 (\Sigma \phi''') \nu_{33}, \nu_5 \bar\nu_8 \bar\kappa_{16}, \nu_5 \varepsilon _8 \bar\kappa_{16}, \delta_5 \sigma_{29}, \alpha'_3\rangle\cong (\mathbb Z/2)^{5}.\)}{\(\displaystyle \begin{aligned}[t]\pi_{36}(S^{5})&=\langle \nu_5 (\Sigma \phi''') \nu_{33}, \nu_5 \bar\nu_8 \bar\kappa_{16}, \nu_5 \varepsilon _8 \bar\kappa_{16}, \delta_5 \sigma_{29},\\&\phantom{=\langle{}}\alpha'_3\rangle\cong (\mathbb Z/2)^{5}.\end{aligned}\)}
	\StemEntry{\(\displaystyle \pi_{37}(S^{6})=\langle P(\Sigma \tau^{\mathrm{IV}}), P((\Sigma A_1) \kappa_{25}), \bar\nu_6 \nu_{14} \bar\kappa_{17}, \nu_6 (\Sigma ^2\phi''') \nu_{34}, \delta_6 \sigma_{30}, \alpha''_3\rangle\cong \mathbb Z/8 \oplus \mathbb Z/2 \oplus \mathbb Z/4 \oplus \mathbb Z/2 \oplus \mathbb Z/2 \oplus \mathbb Z/4.\)}{\(\displaystyle \begin{aligned}[t]\pi_{37}(S^{6})&=\langle P(\Sigma \tau^{\mathrm{IV}}), P((\Sigma A_1) \kappa_{25}), \bar\nu_6 \nu_{14} \bar\kappa_{17}, \nu_6 (\Sigma ^2\phi''') \nu_{34},\\&\phantom{=\langle{}}\delta_6 \sigma_{30}, \alpha''_3\rangle\cong \mathbb Z/8 \oplus \mathbb Z/2 \oplus \mathbb Z/4 \oplus \mathbb Z/2 \oplus \mathbb Z/2 \oplus \mathbb Z/4.\end{aligned}\)}
	\StemEntry{\(\displaystyle \pi_{38}(S^{7})=\langle \sigma' \bar\sigma'_{14}, \sigma' \bar\mu_{14} \sigma_{31}, \bar\nu_7 \nu_{15} \bar\kappa_{18}, \delta_7 \sigma_{31}, \alpha'''_3\rangle\cong (\mathbb Z/2)^{4} \oplus \mathbb Z/8.\)}{\(\displaystyle \begin{aligned}[t]\pi_{38}(S^{7})&=\langle \sigma' \bar\sigma'_{14}, \sigma' \bar\mu_{14} \sigma_{31}, \bar\nu_7 \nu_{15} \bar\kappa_{18}, \delta_7 \sigma_{31},\\&\phantom{=\langle{}}\alpha'''_3\rangle\cong (\mathbb Z/2)^{4} \oplus \mathbb Z/8.\end{aligned}\)}
	\StemEntry{\(\displaystyle \pi_{39}(S^{8})=\langle \sigma_8 \delta_{15}, \sigma_8 \bar\mu_{15} \sigma_{32}, \sigma_8 \bar\sigma'_{15}, \sigma_8\tilde{\lt}_{15}, \sigma_8 (\Sigma G_1), \sigma_8 D^{(1)}_1, (\Sigma \sigma') \bar\sigma'_{15}, (\Sigma \sigma') \bar\mu_{15} \sigma_{32}, \bar\nu_8 \nu_{16} \bar\kappa_{19}, \delta_8 \sigma_{32}, \Sigma \alpha'''_3\rangle\cong (\mathbb Z/2)^{10} \oplus \mathbb Z/8.\)}{\(\displaystyle \begin{aligned}[t]\pi_{39}(S^{8})&=\langle \sigma_8 \delta_{15}, \sigma_8 \bar\mu_{15} \sigma_{32}, \sigma_8 \bar\sigma'_{15}, \sigma_8\tilde{\lt}_{15}, \sigma_8 (\Sigma G_1), \sigma_8 D^{(1)}_1, (\Sigma \sigma') \bar\sigma'_{15}, (\Sigma \sigma') \bar\mu_{15} \sigma_{32}, \bar\nu_8 \nu_{16} \bar\kappa_{19},\\&\phantom{=\langle{}}\delta_8 \sigma_{32}, \Sigma \alpha'''_3\rangle\cong (\mathbb Z/2)^{10} \oplus \mathbb Z/8.\end{aligned}\)}
	\StemEntry{\(\displaystyle \pi_{40}(S^{9})=\langle \sigma_9\delta_{16}, \sigma_9\bar\mu_{16}\sigma_{33}, \sigma_9\bar\sigma'_{16}, \delta_9\sigma_{33}, \alpha^{\mathrm{IV}}_3\rangle\cong (\mathbb Z/2)^{4} \oplus \mathbb Z/16.\)}{\(\displaystyle \begin{aligned}[t]\pi_{40}(S^{9})&=\langle \sigma_9\delta_{16}, \sigma_9\bar\mu_{16}\sigma_{33}, \sigma_9\bar\sigma'_{16}, \delta_9\sigma_{33}, \alpha^{\mathrm{IV}}_3\rangle\\&\cong (\mathbb Z/2)^{4} \oplus \mathbb Z/16.\end{aligned}\)}
	\StemEntry{\(\displaystyle \pi_{41}(S^{10})=\langle P(\sigma^{*}_{21}), \Sigma\alpha^{\mathrm{IV}}_3, \kappa^{*}_{10}\rangle\cong \mathbb Z/8 \oplus \mathbb Z/16 \oplus \mathbb Z/4.\)}{\(\displaystyle \begin{aligned}[t]\pi_{41}(S^{10})&=\langle P(\sigma^{*}_{21}), \Sigma\alpha^{\mathrm{IV}}_3, \kappa^{*}_{10}\rangle\\&\cong \mathbb Z/8 \oplus \mathbb Z/16 \oplus \mathbb Z/4.\end{aligned}\)}
	\StemEntry{\(\displaystyle \pi_{42}(S^{11})=\langle \kappa^{*}_{11}, \Sigma ^2 \alpha^{\mathrm{IV}}_3\rangle\cong \mathbb Z/2 \oplus \mathbb Z/16.\)}{\(\displaystyle \begin{aligned}[t]\pi_{42}(S^{11})&=\langle \kappa^{*}_{11}, \Sigma ^2 \alpha^{\mathrm{IV}}_3\rangle\\&\cong \mathbb Z/2 \oplus \mathbb Z/16.\end{aligned}\)}
	\StemEntry{\(\displaystyle \pi_{43}(S^{12})=\langle P(\bar\kappa_{25}), \kappa^{*}_{12}, \Sigma ^3 \alpha^{\mathrm{IV}}_3\rangle\cong \mathbb Z/4 \oplus \mathbb Z/2 \oplus \mathbb Z/16.\)}{\(\displaystyle \begin{aligned}[t]\pi_{43}(S^{12})&=\langle P(\bar\kappa_{25}), \kappa^{*}_{12}, \Sigma ^3 \alpha^{\mathrm{IV}}_3\rangle\\&\cong \mathbb Z/4 \oplus \mathbb Z/2 \oplus \mathbb Z/16.\end{aligned}\)}
	\StemEntry{\(\displaystyle \pi_{44}(S^{13})=\langle \kappa^{*}_{13}, \alpha^{\mathrm{V}}_3\rangle\cong \mathbb Z/2 \oplus \mathbb Z/32.\)}{\(\displaystyle \begin{aligned}[t]\pi_{44}(S^{13})&=\langle \kappa^{*}_{13}, \alpha^{\mathrm{V}}_3\rangle\\&\cong \mathbb Z/2 \oplus \mathbb Z/32.\end{aligned}\)}
	\StemEntry{\(\displaystyle \pi_{45}(S^{14})=\langle \omega^{*}_{14}, \kappa^{*}_{14}, \rho_{3,14}\rangle\cong \mathbb Z/8 \oplus \mathbb Z/2 \oplus \mathbb Z/64.\)}{\(\displaystyle \begin{aligned}[t]\pi_{45}(S^{14})&=\langle \omega^{*}_{14}, \kappa^{*}_{14}, \rho_{3,14}\rangle\\&\cong \mathbb Z/8 \oplus \mathbb Z/2 \oplus \mathbb Z/64.\end{aligned}\)}
	\StemEntry{\(\displaystyle \pi_{46}(S^{15})=\langle D^{\mathrm{II}}_1, D^{(1)}_1 \sigma_{39}, \kappa^{*\prime}, \omega^{*}_{15}, \kappa^{*}_{15}, \rho_{3,15}\rangle\cong (\mathbb Z/2)^{5}\oplus \mathbb Z/64.\)}{\(\displaystyle \begin{aligned}[t]\pi_{46}(S^{15})&=\langle D^{\mathrm{II}}_1, D^{(1)}_1 \sigma_{39}, \kappa^{*\prime}, \omega^{*}_{15}, \kappa^{*}_{15}, \rho_{3,15}\rangle\\&\cong (\mathbb Z/2)^{5}\oplus \mathbb Z/64.\end{aligned}\)}
	\StemEntry{\(\displaystyle \pi_{47}(S^{16})=\langle B^{\mathrm{II}}_1, B^{(1)}_1 \sigma_{40}, \Sigma D^{\mathrm{II}}_1, \Sigma D^{(1)}_1 \sigma_{40}, \Sigma \kappa^{*\prime}, \omega^{*}_{16}, \kappa^{*}_{16}, \rho_{3,16}\rangle\cong (\mathbb Z/2)^{7} \oplus \mathbb Z/64.\)}{\(\displaystyle \begin{aligned}[t]\pi_{47}(S^{16})&=\langle B^{\mathrm{II}}_1, B^{(1)}_1 \sigma_{40}, \Sigma D^{\mathrm{II}}_1, \Sigma D^{(1)}_1 \sigma_{40}, \Sigma \kappa^{*\prime}, \omega^{*}_{16}, \kappa^{*}_{16},\\&\phantom{=\langle{}}\rho_{3,16}\rangle\cong (\mathbb Z/2)^{7} \oplus \mathbb Z/64.\end{aligned}\)}
	\StemEntry{\(\displaystyle \pi_{48}(S^{17})=\langle \Sigma B^{\mathrm{II}}_1, (\Sigma B^{(1)}_1) \sigma_{41}, \Sigma ^2 \kappa^{*\prime}, \omega^{*}_{17}, \kappa^{*}_{17}, \rho_{3,17}\rangle\cong (\mathbb Z/2)^{5} \oplus \mathbb Z/64.\)}{\(\displaystyle \begin{aligned}[t]\pi_{48}(S^{17})&=\langle \Sigma B^{\mathrm{II}}_1, (\Sigma B^{(1)}_1) \sigma_{41}, \Sigma ^2 \kappa^{*\prime}, \omega^{*}_{17}, \kappa^{*}_{17},\\&\phantom{=\langle{}}\rho_{3,17}\rangle\cong (\mathbb Z/2)^{5} \oplus \mathbb Z/64.\end{aligned}\)}
	\StemEntry{\(\displaystyle \pi_{49}(S^{18})=\langle \Sigma ^2 B^{\mathrm{II}}_1, \Sigma ^3 \kappa^{*\prime}, \omega^{*}_{18}, \kappa^{*}_{18}, \rho_{3,18}\rangle\cong (\mathbb Z/2)^{4} \oplus \mathbb Z/64.\)}{\(\displaystyle \begin{aligned}[t]\pi_{49}(S^{18})&=\langle \Sigma ^2 B^{\mathrm{II}}_1, \Sigma ^3 \kappa^{*\prime}, \omega^{*}_{18}, \kappa^{*}_{18}, \rho_{3,18}\rangle\\&\cong (\mathbb Z/2)^{4} \oplus \mathbb Z/64.\end{aligned}\)}
	\StemEntry{\(\displaystyle \pi_{n+31}(S^{n})=\langle \omega^{*}_n, \kappa^{*}_n, \rho_{3,n}\rangle\cong (\mathbb Z/2)^{2} \oplus \mathbb Z/64\quad(n = 19, 20, 21 \text{ or } 22).\)}{\(\displaystyle \begin{aligned}[t]\pi_{n+31}(S^{n})&=\langle \omega^{*}_n, \kappa^{*}_n, \rho_{3,n}\rangle\\&\cong (\mathbb Z/2)^{2} \oplus \mathbb Z/64\quad(n = 19, 20, 21 \text{ or } 22).\end{aligned}\)}
	\StemEntry{\(\displaystyle \pi_{54}(S^{23})=\langle D^{\mathrm{I}}_2, D_2 \sigma_{47}, \omega^{*}_{23}, \kappa^{*}_{23}, \rho_{3,23}\rangle\cong (\mathbb Z/2)^{4} \oplus \mathbb Z/64.\)}{\(\displaystyle \begin{aligned}[t]\pi_{54}(S^{23})&=\langle D^{\mathrm{I}}_2, D_2 \sigma_{47}, \omega^{*}_{23}, \kappa^{*}_{23}, \rho_{3,23}\rangle\\&\cong (\mathbb Z/2)^{4} \oplus \mathbb Z/64.\end{aligned}\)}
	\StemEntry{\(\displaystyle \pi_{55}(S^{24})=\langle B^{\mathrm{I}}_2, B_2 \sigma_{48}, \Sigma D^{\mathrm{I}}_2, (\Sigma D_2) \sigma_{48}, \omega^{*}_{24}, \kappa^{*}_{24}, \rho_{3,24}\rangle\cong (\mathbb Z/2)^{6} \oplus \mathbb Z/64.\)}{\(\displaystyle \begin{aligned}[t]\pi_{55}(S^{24})&=\langle B^{\mathrm{I}}_2, B_2 \sigma_{48}, \Sigma D^{\mathrm{I}}_2, (\Sigma D_2) \sigma_{48}, \omega^{*}_{24}, \kappa^{*}_{24}, \rho_{3,24}\rangle\\&\cong (\mathbb Z/2)^{6} \oplus \mathbb Z/64.\end{aligned}\)}
	\StemEntry{\(\displaystyle \pi_{56}(S^{25})=\langle \Sigma B^{\mathrm{I}}_2, (\Sigma B_2 )\sigma_{49}, \omega^{*}_{25}, \kappa^{*}_{25}, \rho_{3,25}\rangle\cong (\mathbb Z/2)^{4} \oplus \mathbb Z/64.\)}{\(\displaystyle \begin{aligned}[t]\pi_{56}(S^{25})&=\langle \Sigma B^{\mathrm{I}}_2, (\Sigma B_2 )\sigma_{49}, \omega^{*}_{25}, \kappa^{*}_{25}, \rho_{3,25}\rangle\\&\cong (\mathbb Z/2)^{4} \oplus \mathbb Z/64.\end{aligned}\)}
	\StemEntry{\(\displaystyle \pi_{57}(S^{26})=\langle \Sigma ^2B^{\mathrm{I}}_2, \omega^{*}_{26}, \kappa^{*}_{26}, \rho_{3,26}\rangle\cong (\mathbb Z/2)^{3} \oplus \mathbb Z/64.\)}{\(\displaystyle \begin{aligned}[t]\pi_{57}(S^{26})&=\langle \Sigma ^2B^{\mathrm{I}}_2, \omega^{*}_{26}, \kappa^{*}_{26}, \rho_{3,26}\rangle\\&\cong (\mathbb Z/2)^{3} \oplus \mathbb Z/64.\end{aligned}\)}
	\StemEntry{\(\displaystyle \pi_{n+31}(S^{n})=\langle \omega^{*}_n, \kappa^{*}_n, \rho_{3,n}\rangle\cong (\mathbb Z/2)^{2} \oplus \mathbb Z/64\quad(n\geq27\text{ but }n\ne32).\)}{\(\displaystyle \begin{aligned}[t]\pi_{n+31}(S^{n})&=\langle \omega^{*}_n, \kappa^{*}_n, \rho_{3,n}\rangle\\&\cong (\mathbb Z/2)^{2} \oplus \mathbb Z/64\quad(n\geq27\text{ but }n\ne32).\end{aligned}\)}
	\StemEntry{\(\displaystyle \pi_{63}(S^{32})=\langle P(\iota_{65}), \omega^{*}_{32}, \kappa^{*}_{32}, \rho_{3,32}\rangle\cong \mathbb Z_{(2)} \oplus \mathbb Z/2 \oplus \mathbb Z/2 \oplus \mathbb Z/64.\)}{\(\displaystyle \begin{aligned}[t]\pi_{63}(S^{32})&=\langle P(\iota_{65}), \omega^{*}_{32}, \kappa^{*}_{32}, \rho_{3,32}\rangle\\&\cong \mathbb Z_{(2)} \oplus \mathbb Z/2 \oplus \mathbb Z/2 \oplus \mathbb Z/64.\end{aligned}\)}
	\StemEntry{\(\displaystyle \pi^S_{31}(S^0)=\langle \omega^{*}, \kappa^{*}, \rho_{3}\rangle\cong \mathbb Z/2 \oplus \mathbb Z/2 \oplus \mathbb Z/64.\)}{\(\displaystyle \begin{aligned}[t]\pi^S_{31}(S^0)&=\langle \omega^{*}, \kappa^{*}, \rho_{3}\rangle\\&\cong \mathbb Z/2 \oplus \mathbb Z/2 \oplus \mathbb Z/64.\end{aligned}\)}
	
	\section{The   Miyauchi-Mukai Stem}
	\subsection{Stem 32}\label{app:stem-32}
	\StemEntry{\(\displaystyle \pi_{34}(S^{2})=\langle \eta_2\delta_3\sigma_{27}, \eta_2\varepsilon_3\nu_{11}\bar{\kappa}_{14}, \eta_2\nu'\varepsilon_6\bar{\kappa}_{14}, \eta_2\phi'\nu^2_{28}\rangle\cong (\mathbb Z/2)^{4}.\)}{\(\displaystyle \begin{aligned}[t]\pi_{34}(S^{2})&=\langle \eta_2\delta_3\sigma_{27}, \eta_2\varepsilon_3\nu_{11}\bar{\kappa}_{14}, \eta_2\nu'\varepsilon_6\bar{\kappa}_{14},\\&\phantom{=\langle{}}\eta_2\phi'\nu^2_{28}\rangle\cong (\mathbb Z/2)^{4}.\end{aligned}\)}
	\StemEntry{\(\displaystyle \pi_{35}(S^{3})=\langle \phi'\sigma_{28}, \nu'\eta_6\varepsilon_7\bar{\kappa}_{15}, \mu_{3,3}\sigma_{28}\rangle\cong \mathbb Z/4\oplus(\mathbb Z/2)^{2}.\)}{\(\displaystyle \begin{aligned}[t]\pi_{35}(S^{3})&=\langle \phi'\sigma_{28}, \nu'\eta_6\varepsilon_7\bar{\kappa}_{15},\\&\phantom{=\langle{}}\mu_{3,3}\sigma_{28}\rangle\cong \mathbb Z/4\oplus(\mathbb Z/2)^{2}.\end{aligned}\)}
	\StemEntry{\(\displaystyle \pi_{36}(S^{4})=\langle \nu_4\delta'', \nu_4\sigma'\varepsilon_{14}\kappa_{22}, \nu_4\sigma'\omega_{14}\nu^2_{30}, \nu_4\phi_7\nu^2_{30}, \nu_4\eta_7\varepsilon_8\bar{\kappa}_{16}, (\Sigma \phi')\sigma_{29}, (\Sigma \nu')\eta_7\varepsilon_8\bar{\kappa}_{16}, \mu_{3,4}\sigma_{29}\rangle\cong \mathbb Z/8\oplus\mathbb Z/2 \oplus \mathbb Z/2 \oplus \mathbb Z/2 \oplus \mathbb Z/2\oplus\mathbb Z/4\oplus(\mathbb Z/2)^{2}.\)}{\(\displaystyle \begin{aligned}[t]\pi_{36}(S^{4})&=\langle \nu_4\delta'', \nu_4\sigma'\varepsilon_{14}\kappa_{22}, \nu_4\sigma'\omega_{14}\nu^2_{30}, \nu_4\phi_7\nu^2_{30}, \nu_4\eta_7\varepsilon_8\bar{\kappa}_{16}, (\Sigma \phi')\sigma_{29},\\&\phantom{=\langle{}}(\Sigma \nu')\eta_7\varepsilon_8\bar{\kappa}_{16}, \mu_{3,4}\sigma_{29}\rangle\cong \mathbb Z/8\oplus\mathbb Z/2 \oplus \mathbb Z/2 \oplus \mathbb Z/2 \oplus \mathbb Z/2\oplus\mathbb Z/4\oplus(\mathbb Z/2)^{2}.\end{aligned}\)}
	\StemEntry{\(\displaystyle \pi_{37}(S^{5})=\langle \phi''\sigma_{30}, \nu_5\eta_8\varepsilon_9\bar{\kappa}_{17}, \mu_{3,5}\sigma_{30}\rangle\cong \mathbb Z/8\oplus(\mathbb Z/2)^{2}.\)}{\(\displaystyle \begin{aligned}[t]\pi_{37}(S^{5})&=\langle \phi''\sigma_{30}, \nu_5\eta_8\varepsilon_9\bar{\kappa}_{17},\\&\phantom{=\langle{}}\mu_{3,5}\sigma_{30}\rangle\cong \mathbb Z/8\oplus(\mathbb Z/2)^{2}.\end{aligned}\)}
	\StemEntry{\(\displaystyle \pi_{38}(S^{6})=\langle G^{(2)}_0, P(\sigma_{13}\bar{\kappa}_{20}), P(\Sigma A^{\mathrm{II}}_1), P((\Sigma A^{(1)}_1)\sigma_{33}), (\Sigma \phi'')\sigma_{31}, \mu_{3,6}\sigma_{31}\rangle\cong \mathbb Z/8\oplus\mathbb Z/4\oplus\mathbb Z/2 \oplus \mathbb Z/2\oplus\mathbb Z/8\oplus\mathbb Z/2.\)}{\(\displaystyle \begin{aligned}[t]\pi_{38}(S^{6})&=\langle G^{(2)}_0, P(\sigma_{13}\bar{\kappa}_{20}), P(\Sigma A^{\mathrm{II}}_1), P((\Sigma A^{(1)}_1)\sigma_{33}),\\&\phantom{=\langle{}}(\Sigma \phi'')\sigma_{31}, \mu_{3,6}\sigma_{31}\rangle\cong \mathbb Z/8\oplus\mathbb Z/4\oplus\mathbb Z/2 \oplus \mathbb Z/2\oplus\mathbb Z/8\oplus\mathbb Z/2.\end{aligned}\)}
	\StemEntry{\(\displaystyle \pi_{39}(S^{7})=\langle \sigma'\mu_{3,14}, \sigma'\eta_{14}\sigma_{15}\bar{\mu}_{22}, (\Sigma ^2\phi'')\sigma_{32}, \mu_{3,7}\sigma_{32}\rangle\cong (\mathbb Z/2)^{2}\oplus\mathbb Z/8\oplus\mathbb Z/2.\)}{\(\displaystyle \begin{aligned}[t]\pi_{39}(S^{7})&=\langle \sigma'\mu_{3,14}, \sigma'\eta_{14}\sigma_{15}\bar{\mu}_{22}, (\Sigma ^2\phi'')\sigma_{32},\\&\phantom{=\langle{}}\mu_{3,7}\sigma_{32}\rangle\cong (\mathbb Z/2)^{2}\oplus\mathbb Z/8\oplus\mathbb Z/2.\end{aligned}\)}
	\StemEntry{\(\displaystyle \pi_{40}(S^{8})=\langle \sigma_8\xi_{15}\sigma_{33}, \sigma_8\eta^{\ast\prime}\mu_{31}, \sigma_8\mu_{3,15}, \sigma_8\eta_{15}\bar{\mu}_{16}\sigma_{33}, (\Sigma \sigma')\mu_{3,15}, (\Sigma \sigma')\eta_{15}\sigma_{16}\bar{\mu}_{23}, (\Sigma ^3\phi'')\sigma_{33}, \mu_{3,8}\sigma_{33}\rangle\cong \mathbb Z/8\oplus\mathbb Z/2 \oplus \mathbb Z/2 \oplus \mathbb Z/2 \oplus \mathbb Z/2 \oplus \mathbb Z/2\oplus\mathbb Z/8\oplus\mathbb Z/2.\)}{\(\displaystyle \begin{aligned}[t]\pi_{40}(S^{8})&=\langle \sigma_8\xi_{15}\sigma_{33}, \sigma_8\eta^{\ast\prime}\mu_{31}, \sigma_8\mu_{3,15}, \sigma_8\eta_{15}\bar{\mu}_{16}\sigma_{33}, (\Sigma \sigma')\mu_{3,15}, (\Sigma \sigma')\eta_{15}\sigma_{16}\bar{\mu}_{23},\\&\phantom{=\langle{}}(\Sigma ^3\phi'')\sigma_{33}, \mu_{3,8}\sigma_{33}\rangle\cong \mathbb Z/8\oplus\mathbb Z/2 \oplus \mathbb Z/2 \oplus \mathbb Z/2 \oplus \mathbb Z/2 \oplus \mathbb Z/2\oplus\mathbb Z/8\oplus\mathbb Z/2.\end{aligned}\)}
	\StemEntry{\(\displaystyle \pi_{41}(S^{9})=\langle \sigma_9\xi_{16}\sigma_{34}, \sigma_9\mu_{3,16}, \sigma^2_9\eta_{23}\bar{\mu}_{24}, \mu_{3,9}\sigma_{34}\rangle\cong \mathbb Z/8\oplus(\mathbb Z/2)^{3}.\)}{\(\displaystyle \begin{aligned}[t]\pi_{41}(S^{9})&=\langle \sigma_9\xi_{16}\sigma_{34}, \sigma_9\mu_{3,16}, \sigma^2_9\eta_{23}\bar{\mu}_{24},\\&\phantom{=\langle{}}\mu_{3,9}\sigma_{34}\rangle\cong \mathbb Z/8\oplus(\mathbb Z/2)^{3}.\end{aligned}\)}
	\StemEntry{\(\displaystyle \pi_{42}(S^{10})=\langle [\iota_{10},\bar{\rho}_{10}], \overline{\nu\bar{\kappa}}_{10}, \sigma_{10}\xi_{17}\sigma_{35}, \sigma_{10}\mu_{3,17}\rangle\cong \mathbb Z/16\oplus\mathbb Z/8\oplus\mathbb Z/4\oplus\mathbb Z/2.\)}{\(\displaystyle \begin{aligned}[t]\pi_{42}(S^{10})&=\langle [\iota_{10},\bar{\rho}_{10}], \overline{\nu\bar{\kappa}}_{10}, \sigma_{10}\xi_{17}\sigma_{35},\\&\phantom{=\langle{}}\sigma_{10}\mu_{3,17}\rangle\cong \mathbb Z/16\oplus\mathbb Z/8\oplus\mathbb Z/4\oplus\mathbb Z/2.\end{aligned}\)}
	\StemEntry{\(\displaystyle \pi_{43}(S^{11})=\langle \overline{\nu\bar{\kappa}}_{11}, \sigma_{11}\mu_{3,18}, \sigma_{11}\xi_{18}\sigma_{36}, \lambda'\kappa_{29}\rangle\cong (\mathbb Z/2)^{4}.\)}{\(\displaystyle \begin{aligned}[t]\pi_{43}(S^{11})&=\langle \overline{\nu\bar{\kappa}}_{11}, \sigma_{11}\mu_{3,18}, \sigma_{11}\xi_{18}\sigma_{36},\\&\phantom{=\langle{}}\lambda'\kappa_{29}\rangle\cong (\mathbb Z/2)^{4}.\end{aligned}\)}
	\StemEntry{\(\displaystyle \pi_{44}(S^{12})=\langle \theta\bar{\kappa}_{24}, \xi_{12}\sigma^2_{30}, \overline{\nu\bar{\kappa}}_{12}, \sigma_{12}\mu_{3,19}, (\Sigma \lambda')\kappa_{30}\rangle\cong (\mathbb Z/2)^{5}.\)}{\(\displaystyle \begin{aligned}[t]\pi_{44}(S^{12})&=\langle \theta\bar{\kappa}_{24}, \xi_{12}\sigma^2_{30}, \overline{\nu\bar{\kappa}}_{12}, \sigma_{12}\mu_{3,19},\\&\phantom{=\langle{}}(\Sigma \lambda')\kappa_{30}\rangle\cong (\mathbb Z/2)^{5}.\end{aligned}\)}
	\StemEntry{\(\displaystyle \pi_{45}(S^{13})=\langle \lambda\kappa_{31}, (\Sigma \theta)\bar{\kappa}_{25}, \xi_{13}\sigma^2_{31}, \overline{\nu\bar{\kappa}}_{13}, \sigma_{13}\mu_{3,20}\rangle\cong (\mathbb Z/2)^{5}.\)}{\(\displaystyle \begin{aligned}[t]\pi_{45}(S^{13})&=\langle \lambda\kappa_{31}, (\Sigma \theta)\bar{\kappa}_{25}, \xi_{13}\sigma^2_{31}, \overline{\nu\bar{\kappa}}_{13},\\&\phantom{=\langle{}}\sigma_{13}\mu_{3,20}\rangle\cong (\mathbb Z/2)^{5}.\end{aligned}\)}
	\StemEntry{\(\displaystyle \pi_{46}(S^{14})=\langle \bar{\zeta}^*, \bar{\bar{\sigma}}_{14}, \overline{\nu\bar{\kappa}}_{14}, (\Sigma \lambda)\kappa_{32}, \xi_{14}\sigma^2_{32}, \sigma_{14}\mu_{3,21}\rangle\cong \mathbb Z/8\oplus(\mathbb Z/2)^{5}.\)}{\(\displaystyle \begin{aligned}[t]\pi_{46}(S^{14})&=\langle \bar{\zeta}^*, \bar{\bar{\sigma}}_{14}, \overline{\nu\bar{\kappa}}_{14}, (\Sigma \lambda)\kappa_{32}, \xi_{14}\sigma^2_{32}, \sigma_{14}\mu_{3,21}\rangle\\&\cong \mathbb Z/8\oplus(\mathbb Z/2)^{5}.\end{aligned}\)}
	\StemEntry{\(\displaystyle \pi_{47}(S^{15})=\langle \overline{\nu\bar{\kappa}}_{15}, \bar{\bar{\sigma}}_{15}, \eta^{\ast\prime}\eta^*_{31}, \mu^{\ast\prime\prime}, \Sigma \bar{\zeta}^*, (\Sigma ^2\lambda)\kappa_{33}, \xi_{15}\sigma^2_{33}, \sigma_{15}\mu_{3,22}\rangle\cong (\mathbb Z/2)^{8}.\)}{\(\displaystyle \begin{aligned}[t]\pi_{47}(S^{15})&=\langle \overline{\nu\bar{\kappa}}_{15}, \bar{\bar{\sigma}}_{15}, \eta^{\ast\prime}\eta^*_{31}, \mu^{\ast\prime\prime}, \Sigma \bar{\zeta}^*, (\Sigma ^2\lambda)\kappa_{33}, \xi_{15}\sigma^2_{33},\\&\phantom{=\langle{}}\sigma_{15}\mu_{3,22}\rangle\cong (\mathbb Z/2)^{8}.\end{aligned}\)}
	\StemEntry{\(\displaystyle \pi_{48}(S^{16})=\langle \overline{\nu\bar{\kappa}}_{16}, \bar{\bar{\sigma}}_{16}, \eta^{*2}_{16}, \mu^{**}_{16}, \nu^*_{16}\kappa_{34}, \eta^*_{16}\eta_{32}\rho_{33}, (\Sigma \eta^{\ast\prime})\eta^*_{32}, \Sigma \mu^{\ast\prime\prime}, \Sigma ^2\bar{\zeta}^*, (\Sigma ^3\lambda)\kappa_{34}, \xi_{16}\sigma^2_{34}, \sigma_{16}\mu_{3,23}\rangle\cong (\mathbb Z/2)^{12}.\)}{\(\displaystyle \begin{aligned}[t]\pi_{48}(S^{16})&=\langle \overline{\nu\bar{\kappa}}_{16}, \bar{\bar{\sigma}}_{16}, \eta^{*2}_{16}, \mu^{**}_{16}, \nu^*_{16}\kappa_{34}, \eta^*_{16}\eta_{32}\rho_{33}, (\Sigma \eta^{\ast\prime})\eta^*_{32}, \Sigma \mu^{\ast\prime\prime}, \Sigma ^2\bar{\zeta}^*, (\Sigma ^3\lambda)\kappa_{34}, \xi_{16}\sigma^2_{34},\\&\phantom{=\langle{}}\sigma_{16}\mu_{3,23}\rangle\cong (\mathbb Z/2)^{12}.\end{aligned}\)}
	\StemEntry{\(\displaystyle \pi_{49}(S^{17})=\langle \overline{\nu\bar{\kappa}}_{17}, \bar{\bar{\sigma}}_{17}, \eta^{*2}_{17}, \mu^{**}_{17}, \nu^*_{17}\kappa_{35}, \eta^*_{17}\eta_{33}\rho_{34}, \xi_{17}\sigma^2_{35}, \sigma_{17}\mu_{3,24}\rangle\cong (\mathbb Z/2)^{8}.\)}{\(\displaystyle \begin{aligned}[t]\pi_{49}(S^{17})&=\langle \overline{\nu\bar{\kappa}}_{17}, \bar{\bar{\sigma}}_{17}, \eta^{*2}_{17}, \mu^{**}_{17}, \nu^*_{17}\kappa_{35}, \eta^*_{17}\eta_{33}\rho_{34}, \xi_{17}\sigma^2_{35},\\&\phantom{=\langle{}}\sigma_{17}\mu_{3,24}\rangle\cong (\mathbb Z/2)^{8}.\end{aligned}\)}
	\StemEntry{\(\displaystyle \pi_{50}(S^{18})=\langle [\iota_{18},\rho_{18}], \bar{\bar{\varepsilon}}_{18}, \overline{\nu\bar{\kappa}}_{18}, \bar{\bar{\sigma}}_{18}, \xi_{18}\sigma^2_{36}, \sigma_{18}\mu_{3,25}\rangle\cong \mathbb Z/32\oplus\mathbb Z/4\oplus(\mathbb Z/2)^{4}.\)}{\(\displaystyle \begin{aligned}[t]\pi_{50}(S^{18})&=\langle [\iota_{18},\rho_{18}], \bar{\bar{\varepsilon}}_{18}, \overline{\nu\bar{\kappa}}_{18}, \bar{\bar{\sigma}}_{18}, \xi_{18}\sigma^2_{36}, \sigma_{18}\mu_{3,25}\rangle\\&\cong \mathbb Z/32\oplus\mathbb Z/4\oplus(\mathbb Z/2)^{4}.\end{aligned}\)}
	\StemEntry{\(\displaystyle \pi_{51}(S^{19})=\langle \bar{\bar{\varepsilon}}_{19}, \overline{\nu\bar{\kappa}}_{19}, \bar{\bar{\sigma}}_{19}, \xi_{19}\sigma^2_{37}, \sigma_{19}\mu_{3,26}\rangle\cong \mathbb Z/4\oplus(\mathbb Z/2)^{4}.\)}{\(\displaystyle \begin{aligned}[t]\pi_{51}(S^{19})&=\langle \bar{\bar{\varepsilon}}_{19}, \overline{\nu\bar{\kappa}}_{19}, \bar{\bar{\sigma}}_{19}, \xi_{19}\sigma^2_{37}, \sigma_{19}\mu_{3,26}\rangle\\&\cong \mathbb Z/4\oplus(\mathbb Z/2)^{4}.\end{aligned}\)}
	\StemEntry{\(\displaystyle \pi_{n+32}(S^{n})=\langle \bar{\bar{\varepsilon}}_n, \overline{\nu\bar{\kappa}}_n, \bar{\bar{\sigma}}_n, \sigma_n\mu_{3,n+7}\rangle\cong (\mathbb Z/2)^{4}\quad(n=20, 21).\)}{\(\displaystyle \begin{aligned}[t]\pi_{n+32}(S^{n})&=\langle \bar{\bar{\varepsilon}}_n, \overline{\nu\bar{\kappa}}_n, \bar{\bar{\sigma}}_n, \sigma_n\mu_{3,n+7}\rangle\\&\cong (\mathbb Z/2)^{4}\quad(n=20, 21).\end{aligned}\)}
	\StemEntry{\(\displaystyle \pi_{54}(S^{22})=\langle \bar{\zeta}^{\ast\prime}, \overline{\nu\bar{\kappa}}_{22}, \bar{\bar{\varepsilon}}_{22}, \bar{\bar{\sigma}}_{22}, \sigma_{22}\mu_{3,29}\rangle\cong \mathbb Z/8\oplus(\mathbb Z/2)^{4}.\)}{\(\displaystyle \begin{aligned}[t]\pi_{54}(S^{22})&=\langle \bar{\zeta}^{\ast\prime}, \overline{\nu\bar{\kappa}}_{22}, \bar{\bar{\varepsilon}}_{22}, \bar{\bar{\sigma}}_{22}, \sigma_{22}\mu_{3,29}\rangle\\&\cong \mathbb Z/8\oplus(\mathbb Z/2)^{4}.\end{aligned}\)}
	\StemEntry{\(\displaystyle \pi_{55}(S^{23})=\langle \tilde{\eta}^{\ast\prime}, \theta'_{23}\eta^2_{53}, \bar{\bar{\varepsilon}}_{23}, \overline{\nu\bar{\kappa}}_{23}, \bar{\bar{\sigma}}_{23}, \sigma_{23}\mu_{3,30}\rangle\cong (\mathbb Z/2)^{6}.\)}{\(\displaystyle \begin{aligned}[t]\pi_{55}(S^{23})&=\langle \tilde{\eta}^{\ast\prime}, \theta'_{23}\eta^2_{53}, \bar{\bar{\varepsilon}}_{23}, \overline{\nu\bar{\kappa}}_{23}, \bar{\bar{\sigma}}_{23}, \sigma_{23}\mu_{3,30}\rangle\\&\cong (\mathbb Z/2)^{6}.\end{aligned}\)}
	\StemEntry{\(\displaystyle \pi_{56}(S^{24})=\langle \theta'_{24}\eta^2_{54}, \Sigma \tilde{\eta}^{\ast\prime}, \tilde{\eta}\eta_{48}\sigma_{49}, \tilde{\eta}\bar{\nu}_{48}, \tilde{\eta}^*, \bar{\bar{\varepsilon}}_{24}, \overline{\nu\bar{\kappa}}_{24}, \bar{\bar{\sigma}}_{24}, \sigma_{24}\mu_{3,31}\rangle\cong (\mathbb Z/2)^{9}.\)}{\(\displaystyle \begin{aligned}[t]\pi_{56}(S^{24})&=\langle \theta'_{24}\eta^2_{54}, \Sigma \tilde{\eta}^{\ast\prime}, \tilde{\eta}\eta_{48}\sigma_{49}, \tilde{\eta}\bar{\nu}_{48}, \tilde{\eta}^*, \bar{\bar{\varepsilon}}_{24}, \overline{\nu\bar{\kappa}}_{24}, \bar{\bar{\sigma}}_{24}, \sigma_{24}\mu_{3,31}\rangle\\&\cong (\mathbb Z/2)^{9}.\end{aligned}\)}
	\StemEntry{\(\displaystyle \pi_{57}(S^{25})=\langle (\Sigma \tilde{\eta})\eta_{49}\sigma_{50}, (\Sigma \tilde{\eta})\bar{\nu}_{49}, \Sigma \tilde{\eta}^*, \overline{\nu\bar{\kappa}}_{25}, \bar{\bar{\sigma}}_{25}, \sigma_{25}\mu_{3,32}\rangle\cong (\mathbb Z/2)^{6}.\)}{\(\displaystyle \begin{aligned}[t]\pi_{57}(S^{25})&=\langle (\Sigma \tilde{\eta})\eta_{49}\sigma_{50}, (\Sigma \tilde{\eta})\bar{\nu}_{49}, \Sigma \tilde{\eta}^*, \overline{\nu\bar{\kappa}}_{25}, \bar{\bar{\sigma}}_{25},\\&\phantom{=\langle{}}\sigma_{25}\mu_{3,32}\rangle\cong (\mathbb Z/2)^{6}.\end{aligned}\)}
	\StemEntry{\(\displaystyle \pi_{58}(S^{26})=\langle [\iota_{26},\sigma_{26}], \overline{\nu\bar{\kappa}}_{26}, \bar{\bar{\sigma}}_{26}, \sigma_{26}\mu_{3,33}\rangle\cong \mathbb Z/16\oplus(\mathbb Z/2)^{3}.\)}{\(\displaystyle \begin{aligned}[t]\pi_{58}(S^{26})&=\langle [\iota_{26},\sigma_{26}], \overline{\nu\bar{\kappa}}_{26}, \bar{\bar{\sigma}}_{26}, \sigma_{26}\mu_{3,33}\rangle\\&\cong \mathbb Z/16\oplus(\mathbb Z/2)^{3}.\end{aligned}\)}
	\StemEntry{\(\displaystyle \pi_{n+32}(S^{n})=\langle \overline{\nu\bar{\kappa}}_n, \bar{\bar{\sigma}}_n, \sigma_n\mu_{3,n+7}\rangle\cong (\mathbb Z/2)^{3}\quad(n=27, 28, 29).\)}{\(\displaystyle \begin{aligned}[t]\pi_{n+32}(S^{n})&=\langle \overline{\nu\bar{\kappa}}_n, \bar{\bar{\sigma}}_n, \sigma_n\mu_{3,n+7}\rangle\\&\cong (\mathbb Z/2)^{3}\quad(n=27, 28, 29).\end{aligned}\)}
	\StemEntry{\(\displaystyle \pi_{62}(S^{30})=\langle \eta_{4,30}, \overline{\nu\bar{\kappa}}_{30}, \bar{\bar{\sigma}}_{30}, \sigma_{30}\mu_{3,37}\rangle\cong \mathbb Z/8\oplus(\mathbb Z/2)^{3}.\)}{\(\displaystyle \begin{aligned}[t]\pi_{62}(S^{30})&=\langle \eta_{4,30}, \overline{\nu\bar{\kappa}}_{30}, \bar{\bar{\sigma}}_{30}, \sigma_{30}\mu_{3,37}\rangle\\&\cong \mathbb Z/8\oplus(\mathbb Z/2)^{3}.\end{aligned}\)}
	\StemEntry{\(\displaystyle \pi_{63}(S^{31})=\langle \eta_{4,31}, \theta^{\prime\ast}, \overline{\nu\bar{\kappa}}_{31}, \bar{\bar{\sigma}}_{31}, \sigma_{31}\mu_{3,38}\rangle\cong (\mathbb Z/2)^{5}.\)}{\(\displaystyle \begin{aligned}[t]\pi_{63}(S^{31})&=\langle \eta_{4,31}, \theta^{\prime\ast}, \overline{\nu\bar{\kappa}}_{31}, \bar{\bar{\sigma}}_{31}, \sigma_{31}\mu_{3,38}\rangle\\&\cong (\mathbb Z/2)^{5}.\end{aligned}\)}
	\StemEntry{\(\displaystyle \pi_{64}(S^{32})=\langle \eta_{4,32}, \theta''_{32}, [\iota_{32},\eta_{32}], \overline{\nu\bar{\kappa}}_{32}, \bar{\bar{\sigma}}_{32}, \sigma_{32}\mu_{3,39}\rangle\cong (\mathbb Z/2)^{6}.\)}{\(\displaystyle \begin{aligned}[t]\pi_{64}(S^{32})&=\langle \eta_{4,32}, \theta''_{32}, [\iota_{32},\eta_{32}], \overline{\nu\bar{\kappa}}_{32}, \bar{\bar{\sigma}}_{32}, \sigma_{32}\mu_{3,39}\rangle\\&\cong (\mathbb Z/2)^{6}.\end{aligned}\)}
	\StemEntry{\(\displaystyle \pi_{65}(S^{33})=\langle \eta_{4,33}, \theta''_{33}, \overline{\nu\bar{\kappa}}_{33}, \bar{\bar{\sigma}}_{33}, \sigma_{33}\mu_{3,40}\rangle\cong (\mathbb Z/2)^{5}.\)}{\(\displaystyle \begin{aligned}[t]\pi_{65}(S^{33})&=\langle \eta_{4,33}, \theta''_{33}, \overline{\nu\bar{\kappa}}_{33}, \bar{\bar{\sigma}}_{33}, \sigma_{33}\mu_{3,40}\rangle\\&\cong (\mathbb Z/2)^{5}.\end{aligned}\)}
	\StemEntry{\(\displaystyle \pi_{n+32}(S^{n})=\langle \eta_{4,n}, \overline{\nu\bar{\kappa}}_n, \bar{\bar{\sigma}}_n, \sigma_n\mu_{3,n+7}\rangle\cong (\mathbb Z/2)^{4}\quad(n\ge 34).\)}{\(\displaystyle \begin{aligned}[t]\pi_{n+32}(S^{n})&=\langle \eta_{4,n}, \overline{\nu\bar{\kappa}}_n, \bar{\bar{\sigma}}_n, \sigma_n\mu_{3,n+7}\rangle\\&\cong (\mathbb Z/2)^{4}\quad(n\ge 34).\end{aligned}\)}
	\StemEntry{\(\displaystyle \pi^S_{32}(S^0)=\langle \eta_{4}, \overline{\nu\bar{\kappa}}, \bar{\bar{\sigma}}, \sigma\mu_{3}\rangle\cong \mathbb Z/2 \oplus \mathbb Z/2 \oplus \mathbb Z/2 \oplus \mathbb Z/2.\)}{\(\displaystyle \begin{aligned}[t]\pi^S_{32}(S^0)&=\langle \eta_{4}, \overline{\nu\bar{\kappa}}, \bar{\bar{\sigma}}, \sigma\mu_{3}\rangle\\&\cong \mathbb Z/2 \oplus \mathbb Z/2 \oplus \mathbb Z/2 \oplus \mathbb Z/2.\end{aligned}\)}
	
	\section{Known Unstable Groups in Stem 33}\label{app:stem-33}
	% n=2,...,5: Oda (1979), Proposition 9.20, pp.144--145.
	% n=6,7,8: Yang--Wu (2025), Corollary 1.2, p.54.
	% n=9: Miyauchi--Mukai, 32STEM, final formula.
	\StemEntry{\(\displaystyle \pi_{35}(S^{2})=\langle \eta_2\nu'\eta_6\varepsilon_7\bar\kappa_{15}, \eta_2\phi'\sigma_{28}, \eta_2\mu_{3,3}\sigma_{28}\rangle\cong \mathbb Z/2 \oplus \mathbb Z/4 \oplus \mathbb Z/2.\)}{\(\displaystyle \begin{aligned}[t]\pi_{35}(S^{2})&=\langle \eta_2\nu'\eta_6\varepsilon_7\bar\kappa_{15}, \eta_2\phi'\sigma_{28},\\&\phantom{=\langle{}}\eta_2\mu_{3,3}\sigma_{28}\rangle\cong \mathbb Z/2 \oplus \mathbb Z/4 \oplus \mathbb Z/2.\end{aligned}\)}
	\StemEntry{\(\displaystyle \pi_{36}(S^{3})=\langle \mu_{4,3}, \nu'\phi_6\sigma_{29}, \eta_3\mu_{3,4}\sigma_{29}\rangle\cong (\mathbb Z/2)^{3}.\)}{\(\displaystyle \begin{aligned}[t]\pi_{36}(S^{3})&=\langle \mu_{4,3}, \nu'\phi_6\sigma_{29}, \eta_3\mu_{3,4}\sigma_{29}\rangle\\&\cong (\mathbb Z/2)^{3}.\end{aligned}\)}
	\StemEntry{\(\displaystyle \pi_{37}(S^{4})=\langle \nu_4\sigma'\bar\rho_{14}, \nu_4\sigma'\phi_{14}, \nu_4\sigma'\psi_{14}, \nu_4^2\sigma_{10}\bar\kappa_{17}, \nu_4\phi_7\sigma_{30}, \mu_{4,4}, (\Sigma\nu')\phi_7\sigma_{30}, \eta_4\mu_{3,5}\sigma_{30}\rangle\cong \mathbb Z/8 \oplus \mathbb Z/2 \oplus \mathbb Z/2 \oplus \mathbb Z/8 \oplus (\mathbb Z/2)^{4}.\)}{\(\displaystyle \begin{aligned}[t]\pi_{37}(S^{4})&=\langle \nu_4\sigma'\bar\rho_{14}, \nu_4\sigma'\phi_{14}, \nu_4\sigma'\psi_{14}, \nu_4^2\sigma_{10}\bar\kappa_{17}, \nu_4\phi_7\sigma_{30}, \mu_{4,4}, (\Sigma\nu')\phi_7\sigma_{30},\\&\phantom{=\langle{}}\eta_4\mu_{3,5}\sigma_{30}\rangle\cong \mathbb Z/8 \oplus \mathbb Z/2 \oplus \mathbb Z/2 \oplus \mathbb Z/8 \oplus (\mathbb Z/2)^{4}.\end{aligned}\)}
	\StemEntry{\(\displaystyle \pi_{38}(S^{5})=\langle \nu_5\sigma_8\nu_{15}\bar\kappa_{18}, \nu_5\phi_8\sigma_{31}, \mu_{4,5}, \eta_5\mu_{3,6}\sigma_{31}\rangle\cong (\mathbb Z/2)^{4}.\)}{\(\displaystyle \begin{aligned}[t]\pi_{38}(S^{5})&=\langle \nu_5\sigma_8\nu_{15}\bar\kappa_{18}, \nu_5\phi_8\sigma_{31}, \mu_{4,5}, \eta_5\mu_{3,6}\sigma_{31}\rangle\\&\cong (\mathbb Z/2)^{4}.\end{aligned}\)}
	\StemEntry{\(\displaystyle \pi_{39}(S^{6})=\langle \kappa'_6, P((\Sigma A_1)\omega_{25}), P((\Sigma A_1)\sigma_{25}\mu_{32}), \bar\sigma_6\sigma_{25}^2, P(\Sigma A_1^{(2)}), \mu_{4,6}, \eta_6\mu_{3,7}\sigma_{32}\rangle\cong \mathbb Z/4 \oplus (\mathbb Z/2)^{6}.\)}{\(\displaystyle \begin{aligned}[t]\pi_{39}(S^{6})&=\langle \kappa'_6, P((\Sigma A_1)\omega_{25}), P((\Sigma A_1)\sigma_{25}\mu_{32}), \bar\sigma_6\sigma_{25}^2, P(\Sigma A_1^{(2)}), \mu_{4,6},\\&\phantom{=\langle{}}\eta_6\mu_{3,7}\sigma_{32}\rangle\cong \mathbb Z/4 \oplus (\mathbb Z/2)^{6}.\end{aligned}\)}
	\StemEntry{\(\displaystyle \pi_{40}(S^{7})=\langle \kappa'_7, \sigma'\eta_{14}\mu_{3,15}, \bar\sigma_7\sigma_{26}^2, \mu_{4,7}, \eta_7\mu_{3,8}\sigma_{33}\rangle\cong \mathbb Z/4 \oplus (\mathbb Z/2)^{4}.\)}{\(\displaystyle \begin{aligned}[t]\pi_{40}(S^{7})&=\langle \kappa'_7, \sigma'\eta_{14}\mu_{3,15}, \bar\sigma_7\sigma_{26}^2, \mu_{4,7}, \eta_7\mu_{3,8}\sigma_{33}\rangle\\&\cong \mathbb Z/4 \oplus (\mathbb Z/2)^{4}.\end{aligned}\)}
	\StemEntry{\(\displaystyle \pi_{41}(S^{8})=\langle \sigma_8(\Sigma^3\tau^{\mathrm{IV}}), \kappa'_8, \sigma_8\eta_{15}\mu_{3,16}, \sigma_8\nu_{15}^2\bar\kappa_{21}, (\Sigma\sigma')\eta_{15}\mu_{3,16}, \bar\sigma_8\sigma_{27}^2, \mu_{4,8}, \eta_8\mu_{3,9}\sigma_{34}\rangle\cong \mathbb Z/8 \oplus \mathbb Z/4 \oplus (\mathbb Z/2)^{6}.\)}{\(\displaystyle \begin{aligned}[t]\pi_{41}(S^{8})&=\langle \sigma_8(\Sigma^3\tau^{\mathrm{IV}}), \kappa'_8, \sigma_8\eta_{15}\mu_{3,16}, \sigma_8\nu_{15}^2\bar\kappa_{21}, (\Sigma\sigma')\eta_{15}\mu_{3,16}, \bar\sigma_8\sigma_{27}^2, \mu_{4,8},\\&\phantom{=\langle{}}\eta_8\mu_{3,9}\sigma_{34}\rangle\cong \mathbb Z/8 \oplus \mathbb Z/4 \oplus (\mathbb Z/2)^{6}.\end{aligned}\)}
	\StemEntry{\(\displaystyle \pi_{42}(S^{9})=\langle \sigma_9\eta_{16}\mu_{3,17}, \sigma_9\nu_{16}^2\bar\kappa_{22}, \bar\sigma_9\sigma_{28}^2, \mu_{4,9}, \eta_9\mu_{3,10}\sigma_{35}, \kappa'_9\rangle\cong (\mathbb Z/2)^{6}.\)}{\(\displaystyle \begin{aligned}[t]\pi_{42}(S^{9})&=\langle \sigma_9\eta_{16}\mu_{3,17}, \sigma_9\nu_{16}^2\bar\kappa_{22}, \bar\sigma_9\sigma_{28}^2, \mu_{4,9}, \eta_9\mu_{3,10}\sigma_{35},\\&\phantom{=\langle{}}\kappa'_9\rangle\cong (\mathbb Z/2)^{6}.\end{aligned}\)}
	\endgroup

	\backmatter


\begin{thebibliography}{99}
		\addcontentsline{toc}{chapter}{Bibliography}
		\urlstyle{same}
		\small\setlength{\itemsep}{3pt}\setlength{\parskip}{0pt}
		\bibitem[Wh41]{Jh1}J. H. C. Whitehead, On adding relations to homotopy groups,
		Ann. of Math. (2) \textbf{42} (1941), no.~2, 409--428.
		\bibitem[Wh42]{Wh42}G.~W.~Whitehead, Homotopy properties of the real orthogonal groups, Ann. of Math.
		\textbf{43} (1942), 132--146.
		\bibitem[BM49]{BlMa1}
		A. Blakers and W. S. Massey, The homotopy groups of a triad, Proc. Natl. Acad. Sci. USA
		\textbf{35}(6) (1949), 322--328.
		\bibitem[Wh50]{Wh}G.~W.~Whitehead, A generalization of the Hopf invariant, Ann. of Math. \textbf{51}
		(1950), 192--237.
		\bibitem[Hi51]{Hil51}P. J. Hilton, Suspension theorems and the generalized Hopf invariant, Proc. London
		Math. Soc. (3) \textbf{1} (1951), 462--493.
		\bibitem[St51]{St}
		N. E. Steenrod, The topology of fibre bundles, Princeton University Press, 1951.
		
		
		\bibitem[Ma51]{Massey1951}
		W.~S.~Massey, A new cohomology invariant of topological spaces,
		Bull. Amer. Math. Soc. \textbf{57} (1951), 74.
		
		\bibitem[To52]{Toda1952}
		H.~Toda, Generalized Whitehead products and homotopy groups
		of spheres,
		J. Inst. Polytech. Osaka City Univ. Ser. A
		\textbf{3} (1952), no.~1--2, 43--82.
		
		
		
		
		
		
		
		\bibitem[BS53]{BS}
		A.~Borel and J.-P.~Serre, Groupes de {L}ie et puissances r\'{e}duites de {S}teenrod,
		Amer. J. Math. \textbf{75} (1953), 409--448.
		\bibitem[HW53]{HW53}
		P.~J.~Hilton and J.~H.~C.~Whitehead, Note on the Whitehead product, Ann. of Math.
		\textbf{58} (1953), 429--442.
		
		\bibitem[Su53]{Su1953}
		M.~Sugawara:
		On the homotopy groups of rotation groups,
		Math. J. Okayama Univ. \textbf{3} (1953), 11--21.
		
		
		\bibitem[Wh53]{JHCWh}
		J.~H.~C.~Whitehead, On certain theorems of G.~W.~Whitehead, Ann. of Math. \textbf{58}
		(1953), 418--428.
		\bibitem[Ja54]{Ja1}
		I. M. James, On the homotopy groups of certain pairs and triads, Quart. J. Math.
		\textbf{5}(1) (1954), 260--270.
		
		\bibitem[Su54]{Su1954}
		M.~Sugawara:
		Some remarks on homotopy groups of rotation groups,
		Math. J. Okayama Univ. \textbf{3} (1954), 129--133.
		
		
		
		\bibitem[To55]{Toda1955}
		H.~Toda, Sur les groupes d'homotopie des sph\`eres,
		C. R. Acad. Sci. Paris \textbf{240} (1955), 42--44.
		
		
		\bibitem[Ja56]{J}
		I. M.~James, On the suspension triad, Ann. of Math. \textbf{63}(2) (1956), 191--247.
		\bibitem[BB58]{BB}
		W. D. Barcus and M. G. Barratt, On the homotopy classification of the extensions of a
		fixed map, Trans. Amer. Math. Soc. \textbf{88} (1958), 57--74.
		\bibitem[Bo59]{Bo}
		R. Bott, The stable homotopy of the classical groups, Ann. of Math. \textbf{70} (1959),
		313--337.
		
		\bibitem[To59]{To1}
		H. Toda, A topological proof of theorems of Bott and Borel-Hirzebruch for homotopy
		groups of unitary groups, Mem. Coll. Sci. Univ. Kyoto Ser. A Math. \textbf{32}(1)
		(1959), 103--119.
		\bibitem[Ke60]{Ke}
		M. A. Kervaire, Some nonstable homotopy groups of Lie groups, Illinois J. Math.
		\textbf{4} (1960), 161--169.
		\bibitem[Ba61]{Ba}
		M. G. Barratt, Note on a formula due to Toda, J. London Math. Soc. \textbf{36} (1961),
		95--96.
		\bibitem[Ha61]{Ha}
		B.~Harris, On the homotopy groups of the classical groups, Ann. of Math. \textbf{74}
		(1961), 407--413.
	
		\bibitem[Sp62]{Spanier1962}
		E. H. Spanier,
		Secondary operations on mappings and cohomology,
		\textit{Ann. of Math.} (2) \textbf{75} (1962), 260--282.
		\bibitem[Toda]{T}
		H. Toda, Composition methods in homotopy groups of spheres, Ann. of Math. Studies,
		vol.~49, Princeton University Press, 1962.
		
		\bibitem[Ba63]{Ba2}
		M.~G.~Barratt, Homotopy operations and homotopy groups, AMS Summer topology institute
		(Seattle) (1963).
		\bibitem[20STEM]{MT1}
		M. Mimura and H. Toda, The $(n+20)$-th homotopy groups of $n$-spheres, J. Math. Kyoto
		Univ. \textbf{3} (1963), 37--58.
		\bibitem[Og63]{Oguchi1963}
		K. \^{O}guchi, A generalization of secondary composition and its applications. J. Fac. Sci.
		Univ. Tokyo \textbf{10} (1963), 29--79. Unofficial scanned copy: \href{https://drive.google.com/file/d/1GVuDR3Kq8pR0IgLa9YGWRn8CSVsD3Q70/view?usp=drive_link}{download}.
		\bibitem[To63]{Toda1963}
		H.~Toda, Order of the identity class of a suspension space,
		Ann. of Math. \textbf{78} (1963), 300--325.
		
		
		\bibitem[Mi64]{Mimura1964}
		M. Mimura, On the generalized Hopf homomorphism and the higher composition,
		Part I, J. Math. Kyoto Univ. \textbf{4} (1964), 171--190.
		
		
		
		
		
		
		
		
		
		\bibitem[MT64a]{MT2}
		M. Mimura and H. Toda, Homotopy groups of $\mathit{SU}(3), \mathit{SU}(4)$ and
		$\mathit{Sp}(2),$ J. Math. Kyoto Univ. \textbf{3} (1964), 217--250. 
		
		
		\bibitem[MT64b]{MT3}
		M.~Mimura and H.~Toda, Homotopy groups of symplectic groups, J. Math. Kyoto Univ.
		\textbf{3} (1964), 251--273. 
		
		\bibitem[Og64]{Og}
		K. \^{O}guchi, Generators of 2-primary components of homotopy groups of
		spheres, unitary groups and symplectic groups, J. Fac. Sci. Univ. of Tokyo \textbf{11}
		(1964), 65--111. Unofficial scanned copy: \href{https://drive.google.com/file/d/1OhXyy3QZYfIOP3x1YLYWbhbnpYkq6-Xo/view?usp=drivesdk}{download}.
		\bibitem[HM65]{HoM}
		C. S. Hoo and M. E. Mahowald, Some homotopy groups of Stiefel manifolds, Bull. Amer.
		Math. Soc. \textbf{71} (1965), 661--667.
		
		
		
		\bibitem[2122STEM]{Mi1}
		M. Mimura, On the generalized Hopf homomorphism and the higher composition,
		Part II. $\pi_{n+i}(S^n)$ for $i=21$ and $22,$  
		J. Math. Kyoto Univ. \textbf{4} (1965), 301--326.
		
		
		
		\bibitem[MM67]{imj}
		M.~E.~Mahowald and M.~C.~Tangora,
		Some differentials in the Adams spectral sequence,
		\textit{Topology} \textbf{6} (1967), 349--369.
		
		\bibitem[Mimura1967]{Mi2}
		M. Mimura, The homotopy groups of Lie groups of low rank, J. Math. Kyoto Univ.
		\textbf{6} (1967), 131--176.
		
		\bibitem[Co68]{Cohen1968}
		J.~Cohen, The decomposition of stable homotopy,
		Ann. of Math. \textbf{87} (1968), 305--320.
		
		\bibitem[R1516]{Ka}
		H. Kachi, On the homotopy groups of rotation groups $\mathit{R}_n,$ J. Fac. Sci.,
		Shinshu Univ. \textbf{3} (1968), 13--33. Unofficial scanned copy: \href{https://drive.google.com/file/d/117aRjyJegc7EZTkMNt20FGyRKL7KMz_4/view?usp=drive_link}{download}.
		
		
		\bibitem[Gr73]{Gray}
		B. Gray,
		On the homotopy groups of mapping cones,
		Proc. London Math. Soc. (3) 26 (1973), 497--520.
		\bibitem[2324STEM]{MMO}
		M. Mimura, M. Mori and N. Oda, Determination of $2$-components of the $23$ and
		$24$-stems in homotopy groups of spheres, Mem. Fac. Sci. Kyushu Univ. A \textbf{29}
		(1975), 1--42.
		\bibitem[Ja76]{Ja76}
		I. M. James, The topology of Stiefel manifolds, London Math. Soc. Lecture Note Ser.,
		vol.~24, Cambridge University Press, 1976.
		\bibitem[Fu78]{F}
		Y. Furukawa, The homotopy groups $\Pi_i(\mathit{SU}(3)),$ $\Pi_i(\mathit{SU}(4))$ and
		$\Pi_i(\mathit{Sp}(2))$ for $i=24, 25, 26$ and $27,$ The Bulletin of Aichi Univ. of
		Education \textbf{27}(Natural Science) (1978), 71--79.  Unofficial scanned copy: \href{https://drive.google.com/file/d/15qYCNSfyIQGgB8MAOsAv3eDZ5JfOBzCY/view?usp=drive_link}{download}.
		\bibitem[Oda1978]{Oda1978}
		N. Oda, Some homotopy groups of SU(3), SU(4) and Sp(2), Fukuoka Univ. Sci. Rep.
		\textbf{8}(2) (1978), 77--89. Unofficial scanned copy: \href{https://drive.google.com/file/d/1BalH9D8SQJMhv48c3NuvKKyIi2Ih3Qi6/view?usp=drive_link}{download}.
		\bibitem[Wh78]{GW}
		G. W. Whitehead, Elements of homotopy theory, Graduate Texts in Mathematics, Springer,
		Berlin, 1978.
		
		\bibitem[CMN79]{CMN1979}
		F.~R. Cohen, J.~C. Moore, and J.~A. Neisendorfer,
		Torsion in homotopy groups,
		Ann. of Math. 109 (1979), 121--168.
		
		
		
		
		
		
		
		
		
		
		
		
		
		
		
		\bibitem[Oda1979]{Oda1979}
		N. Oda, Unstable homotopy groups of spheres, Bull. Inst. Adv. Res. Fukuoka Univ.
		\textbf{44} (1979), 49--152.  Unofficial scanned copy: \href{https://drive.google.com/file/d/1iNXikrnfeoHmT-HtJbHPb-f1PKqN3f5l/view?usp=drivesdk}{download}.
		\bibitem[MO81]{MimuraOdaPeriodic}
		M.~Mimura and N.~Oda, Periodic families in homotopy groups, J. Math. Kyoto Univ.
		\textbf{21} (1981), 171--187.
		\bibitem[Os82]{Os}
		H.~\={O}shima, Some James numbers of Stiefel manifolds, Math. Proc. Cambridge Philos.
		Soc. \textbf{92}(1) (1982), 139--161.
		\bibitem[Od85]{Oda2}
		N. Oda, Some Toda brackets and the relations of H.~$\overline{\rm O}$shima in the
		homotopy groups of spheres, Fukuoka Univ. Sci. Rep. \textbf{15} (1985), 5--11. Unofficial scanned copy: \href{https://drive.google.com/file/d/1yepKBVPn3M7ebPOk9X_N9jNr8k1lTjvE/view?usp=drive_link}{download}.
		\bibitem[Od86]{Oda}
		N. Oda, On the orders of the generators in the 18-stem of the homotopy groups of
		spheres, Adv. Stud. Pure Math. \textbf{9} (1986), 231--236. Unofficial scanned copy: \href{https://drive.google.com/file/d/1es3vhqloPowi2uP1SD10g1Q8IBS4eRPe/view?usp=drive_link}{download}.
		
		
		\bibitem[Lu92]{Lun}
		A. T. Lundell, 
		Concise tables of James numbers and some homotopy of classical Lie groups
		and associated homogeneous spaces,
		in \textit{Algebraic Topology: Homotopy and Group Cohomology},
		Lecture Notes in Math., vol.~1509,
		Springer, 1992, pp.~250--272.
		\url{https://doi.org/10.1007/BFb0087515}.
		
		
		\bibitem[Ja95]{Hand}
		I. M. James (ed.), Handbook of algebraic topology, Elsevier, 1995.
		\bibitem[Ad96]{A}
		J. F. Adams, Lectures on exceptional Lie groups, Univ. Chicago Press, Chicago, IL,
		1996.
		
		\bibitem[R1718]{KM1}
		H. Kachi and J. Mukai, Some homotopy groups of the rotation group
		$\mathit{R}_n,$
		Hiroshima J. Math. \textbf{29} (1999), 327--345. Unofficial scanned copy: \href{https://drive.google.com/file/d/1DWfFFh1oHV4w_tJmPnPWAdcR4VO0B3EU/view?usp=drive_link}{download}.
		\bibitem[R1920]{KM2}
		H. Kachi and J. Mukai, The $19$ and $20$-th homotopy groups of the rotation groups
		$\mathit{R}_n,$ Math. J. Okayama Univ. \textbf{42} (2000), 89--113. Unofficial scanned copy: \href{https://drive.google.com/file/d/16P5o6oMCVZmX8yOGT2Kkfi_MbVNEpmOG/view?usp=drive_link}{download}.
		\bibitem[MO01]{MaO}
		H.~J.~Marcum and N.~Oda, Some classical and matrix Toda brackets in the $13$- and
		$15$-stems, Kyushu J. Math. \textbf{55} (2001), 405--428.
		
		
		
		\bibitem[Wu03]{Wu03}
		J. Wu, \emph{Homotopy theory of the suspensions of the projective plane},
		Mem. Amer. Math. Soc. \textbf{162} (2003), no.~769.
		
		
		
		\bibitem[R2122]{HKM}
		Y. Hirato, H. Kachi and J. Mukai, $21$-st and $22$-nd homotopy groups of the $n$-th
		rotation group, J. Fac. Sci., Shinshu Univ. \textbf{41} (2006), 1--28. Unofficial scanned copy: \href{https://drive.google.com/file/d/1omhJqo71I8-yfDHzL66kis5VXSHpuxAk/view?usp=drive_link}{download}.
		\bibitem[GM08]{GM}
		M. Golasi\'nski and J. Mukai, Gottlieb groups of spheres, Topology \textbf{47} (2008),
		399--430.
		\bibitem[Mu08]{Mukaiwhiteheadji} J. Mukai, Determination of the order of the $P$-image by Toda brackets,
		Geom. Topol. Monogr. \textbf{13} (2008), 355--383.
		\url{https://doi.org/10.2140/gtm.2008.13.355}
		
		
		\bibitem[Ro09]{Rotman} J. J. Rotman, An introduction to homological algebra, Springer New York  (2009).
		\bibitem[MM12]{MiM}
		T.~Miyauchi and J.~Mukai, Relations in the $24$-th homotopy groups of spheres, Preprint
		(2012), arXiv:1207.6170. 
		\bibitem[CW13]{CW2013}
		W. Chen and J. Wu, Decomposition of loop spaces and periodic problem on $\pi_{*}$,
		Algebr. Geom. Topol. 13 (2013), 3245--3260.
		
		
		\bibitem[IMM15]{IMM}
		T. Inoue, T. Miyauchi and J. Mukai,
		The 2-components of the 31-stem homotopy groups of the 9 and 10-spheres,
		J. Fac. Sci. Shinshu Univ. \textbf{46} (2015), 1--19. Unofficial scanned copy: \href{https://drive.google.com/file/d/1f2UNaa-WaDPsX2RS2tr4c9oe4wdXeYrz/view?usp=drive_link}{download}.
		\bibitem[OO15]{O2015}
		H.~{\=O}shima and K.~{\=O}shima,
		Quasi tertiary compositions and a Toda bracket in homotopy groups of
		$\mathrm{SU}(3),$  
		Math. J. Okayama Univ. \textbf{57} (2015), 13--78.
		\bibitem[Ha17]{HatcherVB}
		A.~Hatcher, \emph{Vector Bundles and K-Theory}, version 2.2, November 2017. \url{https://pi.math.cornell.edu/~hatcher/VBKT/VB.pdf}.
		\bibitem[32STEM]{MiM1}
		T. Miyauchi and J. Mukai, Determination of the $2$-primary components of the $32$-stem
		homotopy groups of $S^n,$   Bol. Soc. Mat. Mex. (3) \textbf{23} (2017), 319--387. 
		\bibitem[MMO18]{MMOh}
		T.~Miyauchi, J.~Mukai and M.~Ohara, On Mimura's extension problem, Topology Appl.
		\textbf{243} (2018), 135--145.
		\bibitem[R2324]{r2324}
		Y.~Hirato, J.~Lee, T.~Miyauchi, J.~Mukai and M.~Ohara, The $23$-rd and $24$-th homotopy
		groups of the $n$-th rotation group, arXiv:2107.00154 (2021)\uprightsemicolon{} incorporated in
		this book. 
		\url{https://arxiv.org/abs/2107.00154}.
		
		\bibitem[R252627]{r2567}
		Y.~Hirato, J.~Lee, T.~Miyauchi, J.~Mukai and M.~Ohara, The $25$-th, $26$-th and $27$-th
		homotopy groups of the $n$-th rotation group, (preprint in 2022)\uprightsemicolon{} incorporated in
		this book. 
		\bibitem[OO22]{O2022}
		H.~\={O}shima and K.~\={O}shima,
		A system of unstable higher Toda brackets,
		\textit{J. Math. Soc. Japan}
		\textbf{74} (2022), no.~4, 1079--1105.
		\url{https://doi.org/10.2969/jmsj/86558655}.
		
		\bibitem[R2829]{r289}J.~Yang,   T.~Miyauchi  and J.~Mukai, The $28$-th and  $29$-th
		homotopy groups of the $n$-th rotation group, (preprint in 2023)\uprightsemicolon{} incorporated in
		this book.  
		\bibitem[YMM24]{YMM}
		J.~Yang, T.~Miyauchi and J.~Mukai, An unstable approach to the May--Lawrence matrix
		Toda bracket and the 2nd James--Hopf invariant, Topology Appl. \textbf{355} (2024),
		109026.
		\url{https://doi.org/10.1016/j.topol.2024.109026}.
		\bibitem[YMW25]{hp2}
		J. Yang, J. Mukai, and J. Wu,
		On the homotopy groups of the suspended quaternionic projective plane and applications,
		Algebr. Geom. Topol. 25 (2025), no.~5, 2981--3033.
		\bibitem[YW25]{YangWuHHA}
		J.~Yang and J.~Wu, Oda's extension problem in the $33$-stem homotopy group of the
		$6$-sphere, Homology, Homotopy and Applications \textbf{27}(2) (2025), 53--60.
		\url{https://doi.org/10.4310/HHA.2025.v27.n2.a3}. 
		
		
	\end{thebibliography}
\end{document}